\documentclass[10pt,a4paper]{report}
\usepackage[utf8]{inputenc}
\usepackage[T1]{fontenc}
\usepackage{geometry}
\usepackage{amsthm}
\usepackage[francais,english]{babel}
\usepackage{amssymb}
\usepackage{lmodern}
\usepackage{amsmath,amsfonts}
\usepackage{mathrsfs} 
\usepackage{dsfont}
\usepackage{stmaryrd}
 
\usepackage{graphicx}
\usepackage{fullpage}
\usepackage[nottoc,notlot,notlof]{tocbibind}
\usepackage[colorlinks=true,urlcolor=black,linkcolor=black,citecolor=black]{hyperref}
\numberwithin{equation}{section}
\usepackage{mathtools,array}
\usepackage{verbatimbox}
\usepackage{color}
\usepackage{bm}
\usepackage{tikz-cd}
\usepackage{pst-node}
\usetikzlibrary{matrix}
\usepackage{mathtools}
\usepackage{verbatimbox}
\usepackage{diagbox}
\usepackage{array}
\newcolumntype{C}{>{$\displaystyle} c <{$}}
\usepackage{makecell}
\usepackage{float}
\usepackage{tabu}
\usepackage{cancel}
\makeatletter
\def\env@dmatrix{\hskip -\arraycolsep
	\let\@ifnextchar\new@ifnextchar
	\def\arraystretch{2}%
	\array{*{\c@MaxMatrixCols}{>{\displaystyle}c}}%
}

\newenvironment{dmatrix}{\left(\env@dmatrix}{\endmatrix\right)}

\makeatother

\begin{document}

	\renewcommand{\thefootnote}{\fnsymbol{footnote}}
	
	\title{Computer-Assisted Proof of the Main Theorem of\\ 'The Classification of Branched Willmore Spheres in the $3$-Sphere and the $4$-Sphere'}
	\author{Alexis Michelat\footnote{Department of Mathematics, ETH Zentrum, CH-8093 Zürich, Switzerland.}\; and Tristan \selectlanguage{french}Rivière$^{*}$\selectlanguage{english}\setcounter{footnote}{0}}
	\date{\today}
	
	\maketitle
	
	\vspace{-0.5em}

	\begin{abstract}
		We provide a computer-assisted proof of the holomorphy of the quartic and the octic meromorphic differentials arising in the main Theorem 4.11 of our paper 'The Classification of Branched Willmore spheres in the $3$-Sphere and the $4$-Sphere' (arXiv:1706.01405), using the free mathematical software Sage. 
		\vspace{1cm}
		\begin{center}
		\end{center}
	\end{abstract}
	
	\tableofcontents
	
	\theoremstyle{plain}
	\newtheorem*{theorem*}{Theorem}
	\newtheorem{theorem}{Theorem}[section]
	\newenvironment{theorembis}[1]
	{\renewcommand{\thetheorem}{\ref{#1}$'$}%
		\addtocounter{theorem}{-1}%
		\begin{theorem}}
		{\end{theorem}}
	\renewcommand*{\thetheorem}{\Alph{theorem}}
	\newtheorem{lemme}[theorem]{Lemma}
	\newtheorem{propdef}[theorem]{Définition-Proposition}
	\newtheorem{prop}[theorem]{Proposition}
	\newtheorem{cor}[theorem]{Corollary}
	\theoremstyle{definition}
	\newtheorem*{definition}{Definition}
	\newtheorem{defi}[theorem]{Definition}
	\newtheorem{rem}[theorem]{Remark}
	\newtheorem{rems}[theorem]{Remarks}
	\newtheorem{exemple}[theorem]{Example}
	\renewcommand\hat[1]{%
		\savestack{\tmpbox}{\stretchto{%
				\scaleto{%
					\scalerel*[\widthof{\ensuremath{#1}}]{\kern-.6pt\bigwedge\kern-.6pt}%
					{\rule[-\textheight/2]{1ex}{\textheight}}
				}{\textheight}%
			}{0.5ex}}%
		\stackon[1pt]{#1}{\tmpbox}
	}
	\parskip 1ex
	\newcommand{\totimes}{\ensuremath{\,\dot{\otimes}\,}}
	\newcommand{\vc}[3]{\overset{#2}{\underset{#3}{#1}}}
	\newcommand{\conv}[1]{\ensuremath{\underset{#1}{\longrightarrow}}}
	\newcommand{\A}{\ensuremath{\vec{A}}}
	\newcommand{\B}{\ensuremath{\vec{B}}}
	\newcommand{\C}{\ensuremath{\vec{C}}}
	\newcommand{\D}{\ensuremath{\vec{D}}}
	\newcommand{\E}{\ensuremath{\vec{E}}}
	\newcommand{\Q}{\ensuremath{\vec{Q}}}
	\newcommand{\z}{\ensuremath{\bar{z}}}
	\newcommand{\R}{\ensuremath{\mathbb{R}}}
	\newcommand{\N}{\ensuremath{\mathbb{N}}}
	\newcommand{\Z}{\ensuremath{\mathbb{Z}}}
	\newcommand{\p}[1]{\ensuremath{\partial_{#1}}}
	\newcommand{\Res}{\ensuremath{\mathrm{Res}}}
	\newcommand{\lp}[2]{\ensuremath{\mathrm{L}^{#1}(#2)}}
	\renewcommand{\wp}[3]{\ensuremath{\left\Vert #1\right\Vert_{\mathrm{W}^{#2}(#3)}}}
	\newcommand{\np}[3]{\ensuremath{\left\Vert #1\right\Vert_{\mathrm{L}^{#2}(#3)}}}
	\newcommand{\nc}[3]{\ensuremath{\left\Vert #1\right\Vert_{C^{#2}(#3)}}}
	\newcommand{\h}{\ensuremath{\vec{h}}}
	\renewcommand{\Re}{\ensuremath{\mathrm{Re}\,}}
	\renewcommand{\Im}{\ensuremath{\mathrm{Im}\,}}
	\newcommand{\diam}{\ensuremath{\mathrm{diam}\,}}
	\newcommand{\leb}{\ensuremath{\mathscr{L}}}
	\newcommand{\supp}{\ensuremath{\mathrm{supp}\,}}
	\renewcommand{\phi}{\ensuremath{\vec{\Phi}}}
	\renewcommand{\H}{\ensuremath{\vec{H}}}
	\renewcommand{\L}{\ensuremath{\mathrm{L}}}
	\renewcommand{\epsilon}{\ensuremath{\varepsilon}}
	\renewcommand{\bar}{\ensuremath{\overline}}
	\newcommand{\s}[2]{\ensuremath{\langle #1,#2\rangle}}
	\newcommand{\pwedge}[2]{\ensuremath{\,#1\wedge#2\,}}
	\newcommand{\bs}[2]{\ensuremath{\left\langle #1,#2\right\rangle}}
	\newcommand{\scal}[2]{\ensuremath{\langle #1,#2\rangle}}
	\newcommand{\sg}[2]{\ensuremath{\left\langle #1,#2\right\rangle_{\mkern-3mu g}}}
	\newcommand{\n}{\ensuremath{\vec{n}}}
	\newcommand{\ens}[1]{\ensuremath{\left\{ #1\right\}}}
	\newcommand{\lie}[2]{\ensuremath{\left[#1,#2\right]}}
	\newcommand{\g}{\ensuremath{g}}
	\newcommand{\e}{\ensuremath{\vec{e}}}
	\newcommand{\ig}{\ensuremath{|\vec{\mathbb{I}}_{\phi}|}}
	\newcommand{\ik}{\ensuremath{\left|\mathbb{I}_{\phi_k}\right|}}
	\newcommand{\w}{\ensuremath{\vec{w}}}
	\newcommand{\vg}{\ensuremath{\mathrm{vol}_g}}
	\newcommand{\im}{\ensuremath{\mathrm{W}^{2,2}_{\iota}(\Sigma,N^n)}}
	\newcommand{\imm}{\ensuremath{\mathrm{W}^{2,2}_{\iota}(\Sigma,\R^3)}}
	\newcommand{\timm}[1]{\ensuremath{\mathrm{W}^{2,2}_{#1}(\Sigma,T\R^3)}}
	\newcommand{\tim}[1]{\ensuremath{\mathrm{W}^{2,2}_{#1}(\Sigma,TN^n)}}
	\renewcommand{\d}[1]{\ensuremath{\partial_{x_{#1}}}}
	\newcommand{\dg}{\ensuremath{\mathrm{div}_{g}}}
	\renewcommand{\Res}{\ensuremath{\mathrm{Res}}}
	\newcommand{\un}[2]{\ensuremath{\bigcup\limits_{#1}^{#2}}}
	\newcommand{\res}{\mathbin{\vrule height 1.6ex depth 0pt width
			0.13ex\vrule height 0.13ex depth 0pt width 1.3ex}}
	\newcommand{\ala}[5]{\ensuremath{e^{-6\lambda}\left(e^{2\lambda_{#1}}\alpha_{#2}^{#3}-\mu\alpha_{#2}^{#1}\right)\left\langle \nabla_{\vec{e}_{#4}}\vec{w},\vec{\mathbb{I}}_{#5}\right\rangle}}
	\setlength\boxtopsep{1pt}
	\setlength\boxbottomsep{1pt}
	\newcommand\norm[1]{%
		\setbox1\hbox{$#1$}%
		\setbox2\hbox{\addvbuffer{\usebox1}}%
		\stretchrel{\lvert}{\usebox2}\stretchrel*{\lvert}{\usebox2}%
	}
	\allowdisplaybreaks
	\newcommand*\mcup{\mathbin{\mathpalette\mcapinn\relax}}
	\newcommand*\mcapinn[2]{\vcenter{\hbox{$\mathsurround=0pt
				\ifx\displaystyle#1\textstyle\else#1\fi\bigcup$}}}
	\def\Xint#1{\mathchoice
		{\XXint\displaystyle\textstyle{#1}}%
		{\XXint\textstyle\scriptstyle{#1}}%
		{\XXint\scriptstyle\scriptscriptstyle{#1}}%
		{\XXint\scriptscriptstyle\scriptscriptstyle{#1}}%
		\!\int}
	\def\XXint#1#2#3{{\setbox0=\hbox{$#1{#2#3}{\int}$ }
			\vcenter{\hbox{$#2#3$ }}\kern-.58\wd0}}
	\def\ddashint{\Xint=}
	\newcommand{\dashint}[1]{\ensuremath{{\Xint-}_{\mkern-10mu #1}}}
	\newcommand\ccancel[1]{\renewcommand\CancelColor{\color{red}}\cancel{#1}}
	\newcommand\colorcancel[2]{\renewcommand\CancelColor{\color{#2}}\cancel{#1}}
	
	\newpage
	
	\chapter{Introduction and organisation of the paper}
	
	\section{Description of the results}
	
	For a broader picture on the motivations, we refer the reader to \cite{classification}.
	
	By a classical theorem of Robert Bryant (see \cite{bryant}), there exists a holomorphic quartic differential $\mathscr{Q}_{\phi}$ associated to any Willmore immersion $\phi:\Sigma^2\rightarrow S^3$ from a closed Riemann immersion $\Sigma^2$ with the following property : if $\mathscr{Q}_{\phi}=0$, then there exists a stereographic projection $\pi:S^3\setminus\ens{p}\rightarrow\R^3$ from some $p\in S^3$ such that the mean curvature of the composition
	\begin{align*}
		\pi\circ \phi:\Sigma^2\setminus\phi^{-1}(\ens{p})\rightarrow\R^3
	\end{align*}
	vanishes identically. In particular, if $\Sigma^2$ is a topological sphere, we thereby deduce by Riemann-Roch theorem that $\mathscr{Q}_{\phi}$ must always vanish identically. At this point, this is not hard to show the following theorem.
	\begin{theorem*}[Bryant, \cite{bryant}]
		Let $\phi:S^2\rightarrow S^3$ a conformal Willmore immersion. Then $\phi$ is the inverse stereographic projection of a complete minimal immersion $\vec{\Psi}:S^2\setminus\ens{p_1,\cdots,p_m}\rightarrow \R^3$ with embedded planar ends.
	\end{theorem*} 

    Let us now consider branched Willmore immersion $\phi:\Sigma^2\rightarrow S^3$. It was showed by Tobias Lamm and Huy The Nguyen (\cite{lamm}) that the quartic differential $\mathscr{Q}_{\phi}$ is only meromorphic, and may have poles up to order equal to $2$ at each branch point. Let $\ens{p_1,\cdots,p_m}\subset \Sigma^2$ these branch points. That $\mathscr{Q}_{\phi}$ has poles of order at most $2$ at $p_1,\cdots,p_m$ is equivalent to say that it is a \emph{holomorphic} section of the holomorphic line bundle
    \begin{align*}
    	\mathscr{L}=K_{\Sigma^2}^4\otimes\mathscr{O}(2\,p_1+\cdots+2\,p_m),
    \end{align*}
    where $K_{\Sigma^2}=T^{\ast}\Sigma^2$ is the canonical bundle of the compact connected Riemann surface $\Sigma^2$. If $\Sigma^2$ has genus $g$, then the degree (or the first Chern class) of $\mathscr{L}$ is given by
    \begin{align*}
    	\mathrm{deg}(\mathscr{L})&=4\,\mathrm{deg}(\mathscr{L})+2\,\mathrm{deg}(\mathscr{O}(p_1+\cdots+p_m))\\
    	&=4(2g-2)+2m.
    \end{align*}
    In particular, if $g\geq 1$ or $m\geq 4$, then $\mathrm{deg}(\mathscr{L})\geq 0$ and by Riemann-Roch theorem, the space $H^0(\Sigma^2,\mathscr{L})$ of holomorphic sections of $\mathscr{L}$ has positive dimension.
    
    However, if $\Sigma^2$ has genus $0$, and $m\leq 3$, we deduce that
    \begin{align*}
    	\mathrm{deg}(\mathscr{L})=-8+2m\leq -2<0,
    \end{align*}
    so all holomorphic sections of $\mathscr{L}$ must vanish identically. In particular, we have $\mathscr{Q}_{\phi}=0$, and this easily implies the following theorem.
    
    \begin{theorem*}[\label{lamm}Lamm, Nguyen, \cite{lamm}]
    	Let $\phi: S^2\rightarrow S^3$ a branched Willmore immersion with less than $3$ branch points. Then $\phi$ is the inverse stereographic projection of a  complete minimal immersion $\vec{\Psi}:S^2\setminus\ens{p_1,\cdots,p_m}\rightarrow \R^3$ with finite total curvature.
    \end{theorem*}

    By performing a precise Taylor expansion of $\mathscr{Q}_{\phi}$ and obtaining informations on the first (\cite{classification}) and second residue (\cite{blow-up}) of \emph{variational} branched Willmore immersions, we were able to show that $\mathscr{Q}_{\phi}$ has poles of order at most $1$. Furthermore, for $\theta_0\leq 3$, we showed in \cite{classification} that poles are completely removable. 
    This shows that 
    \begin{align*}
    	\mathscr{Q}_{\phi}\in H^0(S^2,K_{S^2}^4\otimes \mathscr{O}(p_1+\cdots+p_m))
    \end{align*}
    and we have
    \begin{align*}
    	\mathrm{deg}\left(K_{S^2}^4\otimes \mathscr{O}(p_1+\cdots+p_m)\right)=-8+m<0
    \end{align*}
    for $m\leq 7$.    
    Therefore, we obtain in particular the following  improvement of \eqref{lamm}.
    \begin{theorem}[\cite{classification}]
    	Let $\phi: S^2\rightarrow S^3$ a \emph{variational} branched Willmore immersion with less than $7$ branch points. Then $\phi$ is the inverse stereographic projection of a  a complete minimal immersion $\vec{\Psi}:S^2\setminus\ens{p_1,\cdots,p_n}\rightarrow \R^3$ with finite total curvature and \emph{zero flux}.
    \end{theorem}
     
    However, as we shall see this is almost impossible to check by hand that $\mathscr{Q}_{\phi}$ has no poles at branch points of multiplicity $\theta_0\geq 4$. The goal of this paper is to show that there are indeed no poles. This permits to obtain the following result, which is a combination of \cite{classification} and of the forthcoming computations.
    
    \begin{theorem}
    	Let $n\geq 3$ and $\phi:D^2\rightarrow S^n$ a branched Willmore sphere with a unique branch point at the origin of multiplicity $\theta_0\geq 1$. Assume that $\phi$ is a \emph{true} branched Willmore disk for $1\leq \theta_0\leq 3$, and provided $\theta_0\geq 2$, the second residue $r(0)\in\ens{0,\cdots,\theta_0-1}$ satisfies $r(0)\leq \theta_0-2$.
    	If the quartic differential $\mathscr{Q}_{\phi}$ defined by 
    	\begin{align*}
    	\mathscr{Q}_{\phi}&=g^{-1}\otimes\left(\partial^N\bar{\partial}^N\h_0\totimes\h_0-\partial^N\h_0\totimes\bar{\partial}^N\h_0\right)+\frac{1}{4}\left(1+|\H|^2\right)\,\h_0\totimes\h_0
    	\end{align*} 
    	is \emph{meromorphic}. Then $\mathscr{Q}_{\phi}$ is holomorphic, and the form of degree $8$
    	\small
    	\begin{align*}
    		&\mathscr{O}_{\phi}=g^{-2}\otimes\bigg\{\frac{1}{4}(\partial^N\bar{\partial}^N\h_0\,\dot{\otimes}\,\partial^N\bar{\partial}^N\h_0)\otimes (\h_0\,\dot{\otimes}\,\h_0)+\frac{1}{4}(\partial^N\h_0\,\dot{\otimes}\, \partial^N\h_0)\otimes (\bar{\partial}^N\h_0\,\dot{\otimes}\, \bar{\partial}^N\h_0)\nonumber\\
    		&-\frac{1}{2}(\partial^N\bar{\partial}^N\h_0\totimes\partial^N\h_0)\otimes(\bar{\partial}^N\h_0\totimes\h_0)-\frac{1}{2}(\partial^N\bar{\partial}^N\h_0\totimes\bar{\partial}^N\h_0)\otimes(\partial^N\h_0\totimes\h_0)+\frac{1}{2}(\partial^N\bar{\partial}^N\h_0\totimes\h_0)\otimes(\partial^N\h_0\totimes\bar{\partial}^N\h_0)\bigg\}\nonumber\\
    		&+\frac{1}{4}(1+|\H|^2)\,g^{-1}\otimes\left\{\frac{1}{2}(\partial^N\bar{\partial}^N\h_0\totimes \h_0)\otimes (\h_0\totimes\h_0)-(\partial^N\h_0\totimes\h_0)\otimes (\bar{\partial}^N\h_0\totimes\h_0)+\frac{1}{2}(\partial^N\h_0\totimes\bar{\partial}^N\h_0)\otimes (\h_0\totimes\h_0)\right\}\nonumber\\
    		&+\frac{1}{64}\left(1+|\H|^2\right)^2\,\left(\h_0\totimes\h_0\right)^2.
    	\end{align*}
    	\normalsize
    	is bounded, \textit{i.e.}
    	\begin{align}
    		\mathscr{O}_{\phi}\in L^{\infty}(D^2).
    	\end{align}
    	In particular, $\mathscr{O}_{\phi}$ is holomorphic once it is meromorphic.
    \end{theorem}
    
    Therefore, we deduce the following special case.
    
    \begin{theorem}[\cite{classification}]
    	\emph{Variational} branched Willmore sphere in $S^3$ are inverse stereographic projection of complete minimal surfaces in $\R^3$ with finite total curvature and \emph{zero flux}.
    \end{theorem}

Indeed, by \cite{classification}, \emph{variational} branched Willmore spheres are \emph{true} and the condition on the second residue is satisfied by \cite{blow-up}, \cite{blow-up2}.

Finally, as the analysis is codimension free, the computations involving only the meromorphy of the quartic form will imply also that the poles of order \textit{a priori} equal to $4$ of Montiel's octic form $\mathscr{O}_{\phi}$ are completely removable.

\begin{theorem}
	Let $\phi:S^2\rightarrow S^4$ be \emph{variational} branched Willmore surface. Then $\phi$ is either the inverse stereographic projection of a complete branched minimal surface in $\mathbb{R}^4$ with finite total curvature and zero flux or the image by the Penrose twistor fibration of a (singular) algebraic curve $C\subset \mathbb{C}\mathbb{P}^3$.
	
	Furthermore, the two possibilities coincide if and only if the algebraic curve $C\subset \mathbb{C}\mathbb{P}^3$ lies in some hypersurface $H\simeq  \mathbb{C}\mathbb{P}^2\subset\mathbb{C}\mathbb{P}^3$. 
\end{theorem}

    \section{Description of the argument}
    
    The notations are consistent with \cite{classification} and will not be reintroduced again here.
    
    First, recall that the quartic form has the following structure for immersions $\phi:\Sigma^2\rightarrow\R^3$
    	\begin{align*}
    	\mathscr{Q}_{\phi}&=g^{-1}\otimes\left(\partial^N\bar{\partial}^N\h_0\totimes\h_0-\partial^N\h_0\totimes\bar{\partial}^N\h_0\right)+\frac{1}{4}|\H|^2\,\h_0\totimes\h_0\\
    	&=g^{-1}\otimes\left(\partial\bar{\partial}\h_0\totimes\h_0-\partial\h_0\totimes\bar{\partial}\h_0\right)+\left(\frac{1}{4}|\H|^2+|\h_0|^2_{WP}\right)\h_0\totimes\h_0+\s{\H}{\h_0}^2.
    	\end{align*} 
    	In particular, the expression makes sense for immersions $\phi:\Sigma^2\rightarrow \R^n$ for any $n\geq 3$, and if $n\geq 4$, then the quartic differential $\mathscr{Q}_{\phi}$ need not be holomorphic. Therefore, if we compute its Taylor expansion for arbitrary codimension immersions, there shall be some non-meromorphic components appearing, which must therefore vanish once we suppose that $\mathscr{Q}_{\phi}$ is meromorphic. Furthermore, one of the main achievements of this paper is to show that the relations furnished by these cancellations are non-trivial, and show the holomorphy of $\mathscr{Q}_{\phi}$, once we combine them with the conservation laws.
    	
    	Indeed, we show that for any branched Willmore disk $\phi:D^2\rightarrow \R^n$ with a unique branch point of multiplicity $\theta_0\geq 3$ at $0$, there exists $\vec{A}_1,\vec{C}_1\in \mathbb{C}^n$ such that
    	\begin{align*}
    		\mathscr{Q}_{\phi}&=(\theta_0-1)(\theta_0-2)\s{\vec{A}_1}{\vec{C}_1}\frac{dz^4}{z}-6(\theta_0-2)\left(|\vec{A}_1|^2\s{\vec{A}_1}{\vec{C}_1}-\s{\bar{\vec{A}_1}}{\vec{C}_1}\s{\vec{A}_1}{\vec{A}_1}\right)\z\, dz^4\\
    		&+\frac{3(\theta_0-2)}{2\theta_0}\left(|\vec{C}_1|^2\s{\vec{A}_1}{\vec{A}_1}-\s{\vec{A}_1}{\bar{\vec{C}_1}}\s{\vec{A}_1}{\vec{C}_1}\right)z^{\theta_0}\,\z^{2-\theta_0}\,dz^4+O(|z|^3).
    	\end{align*}
     If we suppose that $\mathscr{Q}_{\phi}$ is meromorphic, then we obtain (see \eqref{system})
    \begin{align}\label{system2}
    \left\{\begin{alignedat}{1}
    |\vec{A}_1|^2\s{\vec{A}_1}{\vec{C}_1}&=\s{\bar{\vec{A}_1}}{\vec{C}_1}\s{\vec{A}_1}{\vec{A}_1}\\
    |\vec{C}_1|^2\s{\vec{A}_1}{\vec{A}_1}&=\s{\vec{A}_1}{\bar{\vec{C}_1}}\s{\vec{A}_1}{\vec{C}_1},
    \end{alignedat}\right.
    \end{align}
    Remarking that is a linear system in $(\s{\vec{A}_1}{\vec{C}_1},\s{\vec{A}_1}{\vec{A}_1})$, we can recast \eqref{system2} as
    \begin{align}\label{salvationmatrix2}
    \begin{dmatrix}
    |\vec{A}_1|^2 & -\s{\bar{\vec{A}_1}}{\vec{C}_1}\\
    -\s{\vec{A}_1}{\bar{\vec{C}_1}} & |\vec{C}_1|^2
    \end{dmatrix}
    \begin{dmatrix}
    \s{\vec{A}_1}{\vec{C}_1}\\
    \s{\vec{A}_1}{\vec{A}_1}
    \end{dmatrix}=0.
    \end{align}
    Thanks of Cauchy-Schwarz inequality, we obtain
    \begin{align}\label{determinant2}
    \det 	\begin{dmatrix}
    |\vec{A}_1|^2 & -\s{\bar{\vec{A}_1}}{\vec{C}_1}\\
    -\s{\vec{A}_1}{\bar{\vec{C}_1}} & |\vec{C}_1|^2
    \end{dmatrix}=|\vec{A}_1|^2|\vec{C}_1|^2-|\s{\vec{A}_1}{\bar{\vec{C}_1}}|^2\geq 0.
    \end{align}
    Therefore, if the determinant is positive, we obtain
    \begin{align*}
    \s{\vec{A}_1}{\vec{C}_1}=0,
    \end{align*}
    and the holomorphy of the quartic form, and if the determinant vanishes,
    \begin{align}\label{ccancel}
    \vec{A}_1\;\,\text{and}\;\, \vec{C}_1\;\,\text{are proportional}.
    \end{align}
    
    Furthermore, in the case $\theta_0=4$, we obtain
    \begin{align*}
    \mathscr{Q}_{\phi}&=6\s{\vec{A}_1}{\vec{C}_1}\frac{dz^4}{z}-12\left(|\vec{A}_1|^2\s{\vec{A}_1}{\vec{C}_1}-\s{\bar{\vec{A}_1}}{\vec{C}_1}\s{\vec{A}_1}{\vec{A}_1}\right)\z\, dz^4\\
    &+\frac{3}{4}\left(|\vec{C}_1|^2\s{\vec{A}_1}{\vec{A}_1}-\s{\vec{A}_1}{\bar{\vec{C}_1}}\s{\vec{A}_1}{\vec{C}_1}\right)z^{\theta_0}\,\z^{2-\theta_0}\,dz^4-\frac{3}{8}\s{\vec{A}_1}{\vec{C}_1}\bar{\s{\vec{C}_1}{\vec{C}_1}}\frac{\z^4}{z}\log|z|+O(|z|^4).
    \end{align*}
    so we obtain the additional 
    \begin{align}\label{additional}
    	\s{\vec{A}_1}{\vec{C}_1}\bar{\s{\vec{C}_1}{\vec{C}_1}}=0,
    \end{align}
    which implies by the preceding argument, that
    \begin{align*}
    	\s{\vec{A}_1}{\vec{C}_1}=0,
    \end{align*}
    along with the holomorphy of the quartic form. So we can suppose that $\theta_0\geq 5$ in the following.
    
    Now, the final argument is to combine the invariance under inversions of the Willmore energy, which translates into the conservation law
    \begin{align*}
    	d\,\Im\left(\mathscr{I}_{\phi}\left(\partial\H+|\H|^2\partial\phi+2\,g^{-1}\otimes\s{\H}{\h_0}\otimes\bar{\partial}\phi\right)-g^{-1}\otimes\left(\bar{\partial}|\phi|^2\otimes\h_0-2\,\s{\phi}{\h_0}\otimes\bar{\partial}\phi\right)\right)=0
    \end{align*}
    where for all $\vec{X}\in\mathbb{C}^n$, we have
    \begin{align*}
    	\mathscr{I}_{\phi}=|\phi|^2\vec{X}-2\s{\phi}{\vec{X}}\phi
    \end{align*}
    and gives as coefficient in $z^{\theta_0+3}dx_1\wedge dx_2$, by supposing that 
    \begin{align*}
    	|\vec{A}_1|^2\s{\vec{A}_1}{\vec{C}_1}=\s{\bar{\vec{A}_1}}{\vec{C}_1}\s{\vec{A}_1}{\vec{A}_1}
    \end{align*}
    the following quantity (see \eqref{fautycroire})
    \begin{align}\label{alphaomega}
    	\frac{4(\theta_0-4)}{\theta_0^2(\theta_0-3)}|\vec{A}_1|^2\s{\vec{A}_1}{\vec{C}_1}=0.
    \end{align}
    Finally, as \eqref{alphaomega} trivially shows that
    \begin{align*}
    	\s{\vec{A}_1}{\vec{C}_1}=0
    \end{align*}
    we also obtain the holomorphy of the quartic form $\mathscr{Q}_{\phi}$.
    
    Finally, let us mention that there is a result on the removability of poles of a $8$-differential for branched Willmore spheres with values in $S^4$. The removability of these poles does not necessitate additional analysis as one can see in chapter \ref{chapteroctic}.

    \section{Heuristic argument}
    
    We can see easily that for a true Willmore sphere with branch points of multiplicity $\theta_0\leq 3$ the quartic form $\mathscr{Q}_{\phi}$ is holomorphic, as the coefficient $c\in \mathbb{C}$ in
    \begin{align*}
    	\mathscr{Q}_{\phi}=\frac{c\;dz^4}{z}+O(1)
    \end{align*}
    is a function of the first residue $\vec{\gamma}_0$, which vanishes whenever the first residue does. For branch points of multiplicity $\theta_0\geq 4$, we can even remove the hypothesis as $\s{\vec{A}_1}{\vec{C}_1}$ is independent of $\vec{\gamma}_0$. In this case, the coefficient has no geometric meaning, and this fact is general, as for any holomorphic function $h:D^2\rightarrow\mathbb{C}$ such that $h(0)\neq 0$
    \begin{align*}
    	\frac{h(z)\,dz^4}{z},
    \end{align*}
    there exists a local coordinate $w$ around $0$ such that
    \begin{align}
    	\frac{h(z)\,dz^4}{z}=\frac{dw^4}{w}.
    \end{align}
    We refer to  \cite{bryantoverflow} for more details on this.
    
	\section{On notations}
	
	In order to obtain an extremely readable format for the code, a function
	\begin{align*}
	\vec{F}(z)=\lambda_1\vec{A}z^{a}\z^b+\lambda_2\vec{B}z^{c}\z^d\log^p|z|,\quad \lambda_1,\lambda_2\in\mathbb{C},\;\,\vec{A},\vec{B}\in\mathbb{C}^n,\;\, a,b,c,d\in\Z,\;\, p\in\N.
	\end{align*} 
	will be written as the following matrix
	\begin{align*}
	\vec{F}=\begin{pmatrix}
	\lambda_1 & \vec{A} & a & b & 0\\
	\lambda_2 & \vec{B} & c & d & p
	\end{pmatrix}
	\end{align*}
	and in the case $p=0$, we will simply remove the last column of zeroes and write instead
	\begin{align*}
	\vec{F}=\begin{pmatrix}
	\lambda_1 & \vec{A} & a & b\\
	\lambda_2 & \vec{B} & c & d 
	\end{pmatrix}
	\end{align*}
	All program shall deal with matrices of this type, and as logarithm terms only appear at branch points of multiplicity $\theta_0\leq 4$, this will be very convenient to remove this last column to speed up computations. The error in $O(|z|^{\ast})$ will be most of the time omitted for simplicity of notations. 
	
	We used this notations as in  the Sage code, all variables are complex or real numbers, so to keep track of the order scalar product arise, this notations is quite convenient.
	
	Finally, we will indicate with a {\color{red} red cross} scalar products which vanish, while the scalar products which cancel will be indicated by a {\color{blue} blue cross}. For example, we have
	\begin{align*}
		\s{\vec{A}_0}{\vec{A}_0}=0,\quad \s{\vec{A}_1}{\vec{A}_1}+2\s{\vec{A}_0}{\vec{A}_2}=0,
	\end{align*}
	so we will write in expressions involving several scalar products
	\begin{align*}
		\ccancel{\s{\vec{A}_0}{\vec{A}_0}},\quad \colorcancel{\s{\vec{A}_1}{\vec{A}_1}}{blue}+\colorcancel{2\s{\vec{A}_0}{\vec{A}_2}}{blue}
	\end{align*}
	to indicate to the reader which cancellation we have used.
	
	\textbf{Finally, remark that there are hyperlinks to any equation in this document.}
	
	\section{The code}
	
	The Python code from Sage is available on the personal website of the first author at the following address:
	
	\href{https://people.math.ethz.ch/~alexism/publications}{https://people.math.ethz.ch/~alexism/publications}
	 
	\textbf{The link is available under the name \emph{Sage source code}, which can be found after the link to this present paper.}
	 
	The Sage version we used is the version 7.3 (\cite{sage}), and at the end of the project, ETH's Sage version moved to version 7.6, which did not create compatibility issues in the last computations, which were performed with this more recent version. 
	
	\begin{rem}
		The code was runned on ETH's supercomputers, and some programms needed around 3 hours to finish (on the model described below), especially for the last order development of the quartic form and the last order development of the conservation laws. The last order of the Weingarten tensor also took some time (less than one hour though), so we do not know how much time it would take (if it can terminate) on a personal computer. Nevertheless, the first two steps of the code (in the two files \texttt{Code\_Classification\_1.txt} and \texttt{Code\_Classification\_2.txt}) are not very resource demanding, and one shall be able to run it on a personal computer. As an indication, the strongest machine we had access to had the following characteristics:
	\begin{align*}
	  \left\{\begin{alignedat}{1}
	  	&\text{Type : Superserver\; $1027$TR-TF}\\
	  	&\text{CPU :	$48\,  \times$ Intel(R) Xeon(R) CPU E$5-2697$ v$2$ @ $2.70$GHz}\\
	  	&\text{RAM :	$256$GB}
	  	\end{alignedat}\right.
	\end{align*}
	For more informations on the hardware aspect, we refer the reader to the following link:
	
	\href{https://blogs.ethz.ch/isgdmath/central-clients/}{https://blogs.ethz.ch/isgdmath/central-clients/}
	\end{rem}
	
	\chapter{The direct verification of the proof}\label{I}

	\section{Step $1$ verification : upper regularity and poles of order $2$}
	
	First, recall that by the very first development of the four tensors is
	\begin{align}\label{1general}
	\left\{\begin{alignedat}{1}
	\p{z}\phi&=\A_0 z^{\theta_0-1}+O(|z|^{\theta_{0}})\\
	g&=|z|^{2\theta_{0}-2}+O(|z|^{2\theta_{0}-1})\\
	\H&=\Re\left(\frac{\vec{C}_1}{z^{\theta_0-2}}\right)+O(|z|^{3-\theta_0})\\
	\h_0&=O(|z|^{\theta_0-1}).
	\end{alignedat}\right.
	\end{align}
	One can check that each development made in \cite{classification} has an error in $O(|z|^{\alpha+1-\epsilon})$, where $\alpha$ is the maximal order of products of $z$ and $\z$.  From now on, we will not write these errors to simplify notations.
	With Sage, we define
	\begin{verbatim}
	dzphi = matrix([1,a0,n-1,0])
	g = matrix([1,n-1,n-1])
	ginv = matrix([1,1-n,1-n])
	H = matrix([[1/2,c1,1-n,0],[1/2,c1b,0,1-n]]) 
	\end{verbatim}
	with obvious notations (here \texttt{n} corresponds to the multiplicity $\theta_0$), and using the \texttt{latex()} function with the function
	\begin{verbatim}
	latex(dzphi),latex(g),latex(H)
	\end{verbatim}
	yields
	\begin{align*}
	\left\{\begin{alignedat}{1}
	&\left(\begin{array}{rrrr}
	1 & A_{0} & \theta_{0} - 1 & 0
	\end{array}\right),\\ &\left(\begin{array}{rrr}
	1 & \theta_{0} - 1 & \theta_{0} - 1
	\end{array}\right),\\& \left(\begin{array}{rrrr}
	\frac{1}{2} & C_{1} & -\theta_{0} + 2 & 0 \\
	\frac{1}{2} & \overline{C_{1}} & 0 & -\theta_{0} + 2
	\end{array}\right)
	\end{alignedat}
	\right.
	\end{align*}
	which we interpret as
	\begin{align*}
	\left\{\begin{alignedat}{1}
	\p{z}\phi&=\left(\begin{array}{rrrr}
	1 & \A_{0} & \theta_{0} - 1 & 0
	\end{array}\right)\\
	g&=\left(\begin{array}{rrr}
	1 & \theta_{0} - 1 & \theta_{0} - 1
	\end{array}\right)\\
	\H&=\left(\begin{array}{rrrr}
	\frac{1}{2} & \vec{C}_{1} & -\theta_{0} + 2 & 0 \\
	\frac{1}{2} & \overline{\vec{C}_{1}} & 0 & -\theta_{0} + 2
	\end{array}\right)
	\end{alignedat}\right.
	\end{align*}
	In these notations, the first column is a numerical constant in $\mathbb{C}$, the second the vector coefficient in $\mathbb{C}^n$ (if the matrix has $4$ columns), and the two last columns are the power in $z$ and $\z$ respectively.  
	
	\begin{rem}\label{logarithm}
		From step 5, as logarithm will appear in the developments, the scalars will have $4$ columns and the vectors $5$. This shall not create confusions as we will consistently come back and force between this matrix notation and the more traditional one from \cite{classification}.
	\end{rem}
	
	We will systematically omit all $dz, dz^2,dz^3,dz^4,dz^8,|dz|^2$ factors in the following expressions. In particular, scalars have $3$ columns, and vectors $4$.
	
	We will always write first the result given by Sage, then the transcription in correct mathematical notations, so that one can check consistently each step. We will also indicate to which equation one needs to compare with.

	Assuming that the second residue satisfies $r(0)\leq \theta_0-2$, we start from the expansion obtained in \cite{classification}
	\begin{align}\label{3rddev}
	\left\{\begin{alignedat}{1}
	\p{z}\phi&=\vec{A}_0z^{\theta_0-1}+\vec{A}_1z^{\theta_0}+\vec{A}_2z^{\theta_0+1}+\frac{1}{4\theta_0}\vec{C}_1z\,\z^{\theta_0}+\frac{1}{8}\bar{\vec{C}_1}z^{\theta_0-1}\z^2+O(|z|^{\theta_0+2-\epsilon})\\
	\H&=\Re\left(\frac{\vec{C}_1}{z^{\theta_0-2}}\right)+O(|z|^{3-\theta_0-\epsilon}),
	\end{alignedat}\right.
	\end{align}
	obtained by direct integration of
	\begin{align*}
		\p{\z}\left(\p{z}\phi\right)&=\frac{1}{2}e^{2\lambda}\H=\frac{1}{2}|z|^{2\theta_0-2}\left(1+O(|z|)\right)\times\left( \Re\left(\frac{\vec{C}_1}{z^{\theta_0-2}}\right)+O(|z|^{3-\theta_0-\epsilon})\right)\\
		&=\frac{1}{2}\Re\left(\vec{C}_1z\z^{\theta_0-1}\right)+O(|z|^{\theta_0+1-\epsilon})\\
		&=\frac{1}{4}\vec{C}_1z\z^{\theta_0-1}+\frac{1}{4}\bar{\vec{C}_1}z^{\theta_0-1}\z+O(|z|^{\theta_0+1-\epsilon}),
	\end{align*}
	using the first order expansions of $e^{2\lambda}$ and $\H$ in \eqref{1general} (and using Proposition $6.7.$ of \cite{classification} to integrate it).
	
	As $\s{\p{z}\phi}{\p{z}\phi}=0$, we find directly
	\begin{align}\label{sagecancel2}
	\s{\A_0}{\A_0}=\s{\vec{A}_0}{\vec{A}_1}=0, \qquad \s{\vec{A}_0}{\vec{C}_1}=\s{\bar{\vec{A}_0}}{\vec{C}_1}=0.
	\end{align}
	Now define
	\begin{align*}
	\left\{\begin{alignedat}{1}
		\alpha_0&=2\s{\bar{\vec{A}_0}}{\vec{A}_1},\\
		\alpha_1&=2\s{\bar{\vec{A}_0}}{\vec{A}_2}
		\end{alignedat}\right.
	\end{align*}
	As $|\vec{A}_0|^2=\dfrac{1}{2}$, \eqref{3rddev} implies that
	\begin{align*}
		e^{2\lambda}=2|\p{z}\phi|^2=|z|^{2\theta_0-2}\left(1+2\,\Re\left(\alpha_0z+\alpha_1z^2\right)+2|\vec{A}_1|^2|z|^2+O(|z|^{3-\epsilon})\right)
	\end{align*}
	which is coded as
	\begin{align}\label{newdevg}
		e^{2\lambda}=\begin{dmatrix}
		1 & \theta_{0} - 1 & \theta_{0} - 1 \\
		2 \, {\left| A_{1} \right|}^{2} & \theta_{0} & \theta_{0} \\
		\alpha_{0} & \theta_{0} & \theta_{0} - 1 
		\end{dmatrix}
		\begin{dmatrix}
		\alpha_{1} & \theta_{0} + 1 & \theta_{0} - 1 \\
		\overline{\alpha_{0}} & \theta_{0} - 1 & \theta_{0} \\
		\overline{\alpha_{1}} & \theta_{0} - 1 & \theta_{0} + 1
		\end{dmatrix}
	\end{align}
	as expected.
	We now compute
	\begin{align}\label{newdevh0}
		\h_0=\begin{dmatrix}
		-\frac{\theta_{0} - 2}{2 \, \theta_{0}} & C_{1} & 0 & \theta_{0} \\
		4 & A_{2} & \theta_{0} & 0 \\
		-2 \, \alpha_{0} & A_{1} & \theta_{0} & 0 \\
		2 \, \alpha_{0}^{2} - 4 \, \alpha_{1} & A_{0} & \theta_{0} & 0 \\
		2 & A_{1} & \theta_{0} - 1 & 0 \\
		-2 \, \alpha_{0} & A_{0} & \theta_{0} - 1 & 0 \\
		-4 \, {\left| A_{1} \right|}^{2} + 2 \, \alpha_{0} \overline{\alpha_{0}} & A_{0} & \theta_{0} - 1 & 1
		\end{dmatrix}+O(|z|^{\theta_0+1-\epsilon})
	\end{align}
	to be compared with
	\begin{align*}
		\h_0&=2\left(\vec{A}_1-\alpha_0\vec{A}_0+\Big(2|\vec{A}_1|^2-|\alpha_0|^2\Big)\vec{A}_0\z\right)z^{\theta_0-1}+2\left(2\vec{A}_2-\alpha_0\vec{A}_1-(2\alpha_1-\alpha_0^2)\vec{A}_0\right)z^{\theta_0}\\
		&-\frac{(\theta_0-2)}{2\theta_0}\vec{C}_1\z^{\theta_0}+O(|z|^{\theta_0+1-\epsilon})
	\end{align*}
	so we see that both developments coincide to this point.
	We also check that
	\begin{align*}
		\mathscr{Q}_{\phi}=(\theta_0-1)(\theta_0-2)\s{\vec{A}_1}{\vec{C}_1}\frac{dz^4}{z}+O(1)
	\end{align*}
	as in the code, we have
	\begin{align*}
		\mathscr{Q}_{\phi}=\begin{dmatrix}
		-{\left({\left(\alpha_{0} \theta_{0}^{2} - 3 \, \alpha_{0} \theta_{0} + 2 \, \alpha_{0}\right)} A_{0} - {\left(\theta_{0}^{2} - 3 \, \theta_{0} + 2\right)} A_{1}\right)} C_{1} & -1 & 0
		\end{dmatrix}+O(1)
	\end{align*}
	and as $\s{\vec{A}_0}{\vec{C}_1}=0$, and $(\theta_0-1)(\theta_0-2)=\theta_0^2-3\theta_0+2$, we indeed recover the same formula.

	\section{Step $2$ verification : higher development for $\theta_0\geq 4$}
	
	\textbf{From now we suppose that $\theta_0\geq 4$.}
	
	Let us sum up the different developments we obtained so far. From now, we can be much faster and check directly the important steps, that is the development of $\vec{Q}$, $\H$, $\p{z}\phi$, $e^{2\lambda}$ and $\h_0$.
	
	As previously, $\Q$ is such that
	\begin{align*}
	\partial \Q=-|\H|^2\partial \phi-2\,g^{-1}\otimes \s{\H}{\h_0}\otimes \bar{\partial}\phi
	\end{align*}
	and $\Q$ has no anti-holomorphic components of high singularity.
	Now, as $\H=O(|z|^{2-\theta_0})$ and $\p{z}\phi=O(|z|^{\theta_0-1})$, we have
	\begin{align}\label{step20}
	|\H|^2\partial\phi=O(|z|^{3-\theta_0})
	\end{align}
	and (see \texttt{del1H3})
	\begin{align}\label{step21}
	g^{-1}\otimes \s{\H}{\h_0}\otimes\bar{\partial}\phi=\left(\begin{array}{rrrr}
	\A_{1}\cdot \C_{1} & \overline{\A_{0}} & -\theta_{0} + 2 & 0 \\
	\A_{1} \cdot\overline{\C_{1}} & \overline{\A_{0}} & 0 & -\theta_{0}
	+ 2
	\end{array}\right)+O(|z|^{3-\theta_0})
	\end{align}
	which checks. 
	
	While it is not used in this step, we will need the development of $|\H|^2\partial\phi$
	in the next one, so we check that
	\begin{align}\label{delH41}
	|\H|^2\partial\phi=\left(\begin{array}{rrrr}
	\frac{1}{4} \, C_{1}^{2} & A_{0} & -\theta_{0} + 3 & 0 \\
	\frac{1}{2} \, C_{1} \overline{C_{1}} & A_{0} & 1 &
	-\theta_{0} + 2 \\
	\frac{1}{4} \, \overline{C_{1}}^{2} & A_{0} & \theta_{0} - 1
	& -2 \, \theta_{0} + 4.
	\end{array}\right)
	\end{align}

	In particular, thanks of \eqref{step20} and \eqref{step21}, we see that $|\H|^2\partial\phi$ is an error, and we obtain 
	\begin{align*}
		\vec{Q}=\begin{dmatrix}
		 \frac{2 \, {\left({\left(\alpha_{0} {\left(\theta_{0} - 1\right)} - \alpha_{0} \theta_{0}\right)} A_{0} C_{1} + A_{1} C_{1}\right)}}{\theta_{0} - 3} & \overline{A_{0}} & -\theta_{0} + 3 & 0 \\
		-2 \, {\left(\alpha_{0} {\left(\theta_{0} - 1\right)} - \alpha_{0} \theta_{0}\right)} A_{0} \overline{C_{1}} - 2 \, A_{1} \overline{C_{1}} & \overline{A_{0}} & 1 & -\theta_{0} + 2
		\end{dmatrix}
	\end{align*}
	and as $\s{\vec{A}_0}{\vec{C}_1}=\s{\bar{\vec{A}_0}}{\vec{C}_1}=0$, this reduces to
	\begin{align*}
	\vec{Q}=\left(\begin{array}{rrrr}
	\frac{2 \, \A_{1}\cdot \C_{1}}{\theta_{0} - 3} & \overline{\A_{0}} &
	-\theta_{0} + 3 & 0 \\
	-2 \, \A_{1}\cdot \overline{\C_{1}} & \overline{\A_{0}} & 1 &
	-\theta_{0} + 2
	\end{array}\right)
	\end{align*}
	and as 
	\begin{align*}
		\partial\left(\H-2i\vec{L}+\vec{\gamma}_0\log|z|-\vec{Q}\right)=0
	\end{align*}
	there exists $\vec{C}_2\in\mathbb{C}^n$ such that
	\begin{align*}
		\H-2i\vec{L}+\vec{\gamma}_0\log|z|&=\frac{\bar{\vec{C}_1}}{\z^{\theta_0-2}}+\frac{\bar{\vec{C}_2}}{\z^{\theta_0-3}}+\vec{Q}\\
		&=\frac{\bar{\vec{C}_1}}{\z^{\theta_0-2}}+\frac{\bar{\vec{C}_2}}{\z^{\theta_0-3}}+\begin{dmatrix}
				 \frac{2 \, {\left({\left(\alpha_{0} {\left(\theta_{0} - 1\right)} - \alpha_{0} \theta_{0}\right)} A_{0} C_{1} + A_{1} C_{1}\right)}}{\theta_{0} - 3} & \overline{A_{0}} & -\theta_{0} + 3 & 0 \\
		-2 \, {\left(\alpha_{0} {\left(\theta_{0} - 1\right)} - \alpha_{0} \theta_{0}\right)} A_{0} \overline{C_{1}} - 2 \, A_{1} \overline{C_{1}} & \overline{A_{0}} & 1 & -\theta_{0} + 2
		\end{dmatrix}
	\end{align*}
	so taking the real part we get
	\begin{align*}
		\H=\Re\left(\frac{\vec{C}_1}{z^{\theta_0-2}}\right)+\begin{dmatrix}
		 \frac{{\left(\alpha_{0} {\left(\theta_{0} - 1\right)} - \alpha_{0} \theta_{0}\right)} A_{0} C_{1} + A_{1} C_{1}}{\theta_{0} - 3} & \overline{A_{0}} & -\theta_{0} + 3 & 0 \\
		\frac{{\left({\left(\theta_{0} - 1\right)} \overline{\alpha_{0}} - \theta_{0} \overline{\alpha_{0}}\right)} \overline{A_{0}} \overline{C_{1}} + \overline{A_{1}} \overline{C_{1}}}{\theta_{0} - 3} & A_{0} & 0 & -\theta_{0} + 3 \\
		-{\left(\alpha_{0} {\left(\theta_{0} - 1\right)} - \alpha_{0} \theta_{0}\right)} A_{0} \overline{C_{1}} - A_{1} \overline{C_{1}} & \overline{A_{0}} & 1 & -\theta_{0} + 2 \\
		-{\left({\left(\theta_{0} - 1\right)} \overline{\alpha_{0}} - \theta_{0} \overline{\alpha_{0}}\right)} C_{1} \overline{A_{0}} - C_{1} \overline{A_{1}} & A_{0} & -\theta_{0} + 2 & 1
		\end{dmatrix}
	\end{align*}
	so using again $\s{\vec{A}_0}{\vec{C}_1}=\s{\bar{\vec{A}_0}}{\vec{C}_1}=0$, we finally obtain
	\small
	\begin{align*}
	\H=\begin{dmatrix}
	\frac{1}{2} & \overline{\C_{1}} & 0 & \theta_{0} - 2 \\
	\frac{1}{2} & \C_{1} & \theta_{0} - 2 & 0 \\
	\frac{1}{2} & \overline{\D_{2}} & 0 & -\theta_{0} + 3 \\
	\frac{1}{2} & \D_{2} & -\theta_{0} + 3 & 0 
	\end{dmatrix}
	\quad
	\begin{dmatrix}
	\frac{\A_{1}\cdot \C_{1}}{\theta_{0} - 3} & \overline{\A_{0}} &
	-\theta_{0} + 3 & 0\\
	\frac{\overline{\A_{1}}\cdot \overline{\C_{1}}}{\theta_{0} - 3} & A_{0}
	& 0 & -\theta_{0} + 3 \\
	-\A_{1} \cdot\overline{\C_{1}} & \overline{\A_{0}} & 1 & -\theta_{0}
	+ 2 \\
	-\C_{1} \cdot\overline{\A_{1}} & \A_{0} & -\theta_{0} + 2 & 1
	\end{dmatrix}
	\end{align*}
	\normalsize
	so for
	\begin{align}\label{BC2}
	\left\{\begin{alignedat}{1}
	\vec{C}_2&=\vec{D}_2+\frac{2}{\theta_0-3}\s{\vec{A}_1}{\vec{C}_1}\bar{\vec{A}_0},\\ \vec{B}_1&=-2\s{\bar{\vec{A}_1}}{\vec{C}_1}\vec{A}_0
	\end{alignedat}\right.
	\end{align}
	we obtain
	\begin{align}\label{devH0}
	\H&=\Re\left(\frac{\vec{C}_1}{z^{\theta_0-2}}+\frac{\vec{C}_2}{z^{\theta_0-3}}+\vec{B}_1\frac{\z}{z^{\theta_0-2}}\right)+O(|z|^{4-\theta_0-\epsilon}).
	\end{align}
	
	\section{Upper development for $\theta_0\geq 5$}
	
	\textbf{From now on, all computations are exclusively computer generated.}
	
	Reinserting the new development of $\H$ allows one to obtain for some $\vec{D}_4\in\mathbb{C}^n$
	\small
	\begin{align}\label{newHcode}
		\H&=\Re\left(\frac{\vec{C}_1}{z^{\theta_0-2}}+\frac{\vec{C}_2}{z^{\theta_0-3}}+\frac{\vec{D}_3}{z^{\theta_0-4}}+\frac{\vec{D}_4}{z^{\theta_0-5}}\right)\\
		&+\begin{dmatrix}
		-\frac{{\left(\overline{A_{0}} \overline{\alpha_{0}} - \overline{A_{1}}\right)} \overline{C_{1}}}{\theta_{0} - 3} & A_{0} & 0 & -\theta_{0} + 3 \\
		\frac{2 \, {\left({\left(\overline{\alpha_{0}}^{2} - \overline{\alpha_{1}}\right)} \overline{A_{0}} - \overline{A_{1}} \overline{\alpha_{0}} + \overline{A_{2}}\right)} \overline{C_{1}} - {\left(\overline{A_{0}} \overline{\alpha_{0}} - \overline{A_{1}}\right)} \overline{C_{2}}}{\theta_{0} - 4} & A_{0} & 0 & -\theta_{0} + 4 \\
		\frac{\overline{C_{1}}^{2}}{8 \, {\left(\theta_{0} - 4\right)}} & \overline{A_{0}} & 0 & -\theta_{0} + 4 \\
		{\left(A_{0} \alpha_{0} - A_{1}\right)} \overline{C_{1}} & \overline{A_{0}} & 1 & -\theta_{0} + 2 \\
		-\frac{{\left(\overline{A_{0}} \overline{\alpha_{0}} - \overline{A_{1}}\right)} \overline{B_{1}} + {\left(2 \, {\left| A_{1} \right|}^{2} \overline{A_{0}} - 2 \, \alpha_{0} \overline{A_{0}} \overline{\alpha_{0}} + \alpha_{0} \overline{A_{1}}\right)} \overline{C_{1}}}{\theta_{0} - 3} & A_{0} & 1 & -\theta_{0} + 3 \\
		-\frac{{\left(\overline{A_{0}} \overline{\alpha_{0}} - \overline{A_{1}}\right)} \overline{C_{1}}}{\theta_{0} - 3} & A_{1} & 1 & -\theta_{0} + 3 \\
		{\left(A_{0} \alpha_{0} - A_{1}\right)} \overline{C_{1}} & \overline{A_{1}} & 1 & -\theta_{0} + 3 \\
		{\left(2 \, A_{0} {\left| A_{1} \right|}^{2} - 2 \, A_{0} \alpha_{0} \overline{\alpha_{0}} + A_{1} \overline{\alpha_{0}}\right)} \overline{C_{1}} + {\left(A_{0} \alpha_{0} - A_{1}\right)} \overline{C_{2}} & \overline{A_{0}} & 1 & -\theta_{0} + 3 \\
		-\frac{C_{1} {\left(\theta_{0} + 2\right)} \overline{C_{1}}}{8 \, \theta_{0}} & A_{0} & 2 & -\theta_{0} + 2 \\
		\frac{1}{2} \, {\left(A_{0} \alpha_{0} - A_{1}\right)} \overline{B_{1}} - {\left({\left(\alpha_{0}^{2} - \alpha_{1}\right)} A_{0} - A_{1} \alpha_{0} + A_{2}\right)} \overline{C_{1}} & \overline{A_{0}} & 2 & -\theta_{0} + 2 \\
		-\frac{C_{1}^{2}}{4 \, \theta_{0}} & \overline{A_{0}} & -2 \, \theta_{0} + 4 & \theta_{0} \\
		{\left(\overline{A_{0}} \overline{\alpha_{0}} - \overline{A_{1}}\right)} C_{1} & A_{0} & -\theta_{0} + 2 & 1 \\
		-\frac{C_{1} {\left(\theta_{0} + 2\right)} \overline{C_{1}}}{8 \, \theta_{0}} & \overline{A_{0}} & -\theta_{0} + 2 & 2 \\
		\frac{1}{2} \, {\left(\overline{A_{0}} \overline{\alpha_{0}} - \overline{A_{1}}\right)} B_{1} - {\left({\left(\overline{\alpha_{0}}^{2} - \overline{\alpha_{1}}\right)} \overline{A_{0}} - \overline{A_{1}} \overline{\alpha_{0}} + \overline{A_{2}}\right)} C_{1} & A_{0} & -\theta_{0} + 2 & 2 \\
		-\frac{{\left(A_{0} \alpha_{0} - A_{1}\right)} C_{1}}{\theta_{0} - 3} & \overline{A_{0}} & -\theta_{0} + 3 & 0 \\
		-\frac{{\left(A_{0} \alpha_{0} - A_{1}\right)} B_{1} + {\left(2 \, A_{0} {\left| A_{1} \right|}^{2} - 2 \, A_{0} \alpha_{0} \overline{\alpha_{0}} + A_{1} \overline{\alpha_{0}}\right)} C_{1}}{\theta_{0} - 3} & \overline{A_{0}} & -\theta_{0} + 3 & 1 \\
		-\frac{{\left(A_{0} \alpha_{0} - A_{1}\right)} C_{1}}{\theta_{0} - 3} & \overline{A_{1}} & -\theta_{0} + 3 & 1 \\
		{\left(\overline{A_{0}} \overline{\alpha_{0}} - \overline{A_{1}}\right)} C_{1} & A_{1} & -\theta_{0} + 3 & 1 \\
		{\left(2 \, {\left| A_{1} \right|}^{2} \overline{A_{0}} - 2 \, \alpha_{0} \overline{A_{0}} \overline{\alpha_{0}} + \alpha_{0} \overline{A_{1}}\right)} C_{1} + {\left(\overline{A_{0}} \overline{\alpha_{0}} - \overline{A_{1}}\right)} C_{2} & A_{0} & -\theta_{0} + 3 & 1 \\
		\frac{2 \, {\left({\left(\alpha_{0}^{2} - \alpha_{1}\right)} A_{0} - A_{1} \alpha_{0} + A_{2}\right)} C_{1} - {\left(A_{0} \alpha_{0} - A_{1}\right)} C_{2}}{\theta_{0} - 4} & \overline{A_{0}} & -\theta_{0} + 4 & 0 \\
		\frac{C_{1}^{2}}{8 \, {\left(\theta_{0} - 4\right)}} & A_{0} & -\theta_{0} + 4 & 0 \\
		-\frac{\overline{C_{1}}^{2}}{4 \, \theta_{0}} & A_{0} & \theta_{0} & -2 \, \theta_{0} + 4
		\end{dmatrix}
	\end{align}
	\normalsize
	so for some $\vec{C}_3,\vec{B}_2,\vec{B}_3,\vec{E}_1\in \mathbb{C}^n$ we have
	\begin{align*}
		\H=\Re\left(\frac{\vec{C}_1}{z^{\theta_0-2}}+\frac{\vec{C}_2}{z^{\theta_0-3}}+\frac{\vec{C}_3}{z^{\theta_0-4}}+\vec{B}_1\frac{\z}{z^{\theta_0-2}}+\vec{B}_2\frac{\z^2}{z^{\theta_0-2}}+\vec{B}_3\frac{\z}{z^{\theta_0-3}}\right)+O(|z|^{5-\theta_0}).
	\end{align*} 
	or on the code
	\begin{align}\label{newH}
		\H=\begin{dmatrix}
		\frac{1}{2} & C_{1} & -\theta_{0} + 2 & 0 \\
		\frac{1}{2} & C_{2} & -\theta_{0} + 3 & 0 \\
		\frac{1}{2} & C_{3} & -\theta_{0} + 4 & 0 \\
		\frac{1}{2} & B_{1} & -\theta_{0} + 2 & 1 \\
		\frac{1}{2} & B_{2} & -\theta_{0} + 2 & 2 \\
		\frac{1}{2} & B_{3} & -\theta_{0} + 3 & 1 \\
		\frac{1}{2} & E_{1} & -2 \, \theta_{0} + 4 & \theta_{0} 
		\end{dmatrix}
		\begin{dmatrix}
		\frac{1}{2} & \overline{C_{1}} & 0 & -\theta_{0} + 2 \\
		\frac{1}{2} & \overline{C_{2}} & 0 & -\theta_{0} + 3 \\
		\frac{1}{2} & \overline{C_{3}} & 0 & -\theta_{0} + 4 \\
		\frac{1}{2} & \overline{B_{1}} & 1 & -\theta_{0} + 2 \\
		\frac{1}{2} & \overline{B_{2}} & 2 & -\theta_{0} + 2 \\
		\frac{1}{2} & \overline{B_{3}} & 1 & -\theta_{0} + 3 \\
		\frac{1}{2} & \overline{E_{1}} & \theta_{0} & -2 \, \theta_{0} + 4
		\end{dmatrix}
	\end{align}
	We will only need the precise expressions of $\vec{B}_3$ and $\vec{E}_1$. First, recall that \eqref{BC2} we have
	\begin{align*}
		\vec{B}_1=-2\s{\bar{\vec{A}_1}}{\vec{C}_1}\vec{A}_0
	\end{align*}
	while by \eqref{sagecancel2}, we have
	\begin{align*}
		\s{\vec{A}_0}{\vec{A}_0}=\s{\vec{A}_0}{\vec{A}_1}=\s{\vec{A}_0}{\vec{C}_1}=\s{\vec{A}_0}{\bar{\vec{C}_1}}=0.
	\end{align*}
	Therefore, comparing \eqref{newH} and \eqref{newHcode}, we obtain the expression
	\begin{align*}
		\frac{1}{2}\vec{B}_3&=\begin{dmatrix}
		-\frac{\ccancel{{\left(A_{0} \alpha_{0} - A_{1}\right)} B_{1}} + {\left(2 \, \ccancel{A_{0} {\left| A_{1} \right|}^{2}} - \ccancel{2 \, A_{0} \alpha_{0} \overline{\alpha_{0}}} + A_{1} \overline{\alpha_{0}}\right)} C_{1}}{\theta_{0} - 3} & \overline{A_{0}} & -\theta_{0} + 3 & 1 \\
		-\frac{{\left(\ccancel{A_{0} \alpha_{0}} - A_{1}\right)} C_{1}}{\theta_{0} - 3} & \overline{A_{1}} & -\theta_{0} + 3 & 1 \\
		{\left(\ccancel{\overline{A_{0}} \overline{\alpha_{0}}} - \overline{A_{1}}\right)} C_{1} & A_{1} & -\theta_{0} + 3 & 1 \\
		{\left(\ccancel{2 \, {\left| A_{1} \right|}^{2} \overline{A_{0}}} - \ccancel{\ccancel{2 \, \alpha_{0} \overline{A_{0}} \overline{\alpha_{0}}}} + \alpha_{0} \overline{A_{1}}\right)} C_{1} + {\left(\ccancel{\overline{A_{0}} \overline{\alpha_{0}}} - \overline{A_{1}}\right)} C_{2} & A_{0} & -\theta_{0} + 3 & 1 \\
		\end{dmatrix}\\
		&=\frac{\bar{\alpha_0}}{\theta_0-3}\s{\vec{A}_1}{\vec{C}_1}\bar{\vec{A}_0}+\frac{1}{\theta_0-3}\s{\vec{A}_1}{\vec{C}_1}\bar{\vec{A}_1}-\s{\bar{\vec{A}_1}}{\vec{C}_1}\vec{A}_1+\left(\alpha_0\s{\bar{\vec{A}_1}}{\vec{C}_1}-\s{\bar{\vec{A}_1}}{\vec{C}_2}\right)\vec{A}_0
	\end{align*}
	while we trivially have
	\begin{align*}
		\frac{1}{2}\vec{E}_1&=\begin{dmatrix}
		-\frac{C_{1}^{2}}{4 \, \theta_{0}} & \overline{A_{0}} & -2 \, \theta_{0} + 4 & \theta_{0} 
		\end{dmatrix}\\
		&=-\frac{1}{4\theta_0}\s{\vec{C}_1}{\vec{C}_1}\bar{\vec{A}_0}
	\end{align*}
	Now, we have as $\s{\bar{\vec{A}_0}}{\vec{B}_1}=-\s{\bar{\vec{A}_1}}{\vec{C}_1}$
	\begin{align*}
		\frac{1}{2}\vec{B}_2&=\begin{dmatrix}
		-\frac{C_{1} {\left(\theta_{0} + 2\right)} \overline{C_{1}}}{8 \, \theta_{0}} & \overline{A_{0}} & -\theta_{0} + 2 & 2 \\
		\frac{1}{2} \, {\left(\overline{A_{0}} \overline{\alpha_{0}} - \overline{A_{1}}\right)} B_{1} - {\left({\left(\overline{\alpha_{0}}^{2} - \overline{\alpha_{1}}\right)} \overline{A_{0}} - \overline{A_{1}} \overline{\alpha_{0}} + \overline{A_{2}}\right)} C_{1} & A_{0} & -\theta_{0} + 2 & 2 
		\end{dmatrix}\\
		&=-\frac{(\theta_0+2)}{8\theta_0}|\vec{C}_1|^2\bar{\vec{A}_0}+\left(-\frac{1}{2}\bar{\alpha_0}\s{\bar{\vec{A}_1}}{\vec{C}_1}+\bar{\alpha_0}\s{\bar{\vec{A}_1}}{\vec{C}_1}-\s{\bar{\vec{A}_2}}{\vec{C}_1}\right)\vec{A}_0\\
		&=-\frac{(\theta_0+2)}{8\theta_0}|\vec{C}_1|^2\bar{\vec{A}_0}+\left(\frac{1}{2}\bar{\alpha_0}\s{\bar{\vec{A}_1}}{\vec{C}_1}-\s{\bar{\vec{A}_2}}{\vec{C}_1}\right)\vec{A}_0
	\end{align*}
	so we obtain
	\begin{align}\label{b3c1}
		\left\{
		\begin{alignedat}{1}
		\vec{B}_2&=-\frac{(\theta_0+2)}{4\theta_0}|\vec{C}_1|^2\bar{\vec{A}_0}+\left(\bar{\alpha_0}\s{\bar{\vec{A}_1}}{\vec{C}_1}-2\s{\bar{\vec{A}_2}}{\vec{C}_1}\right)\vec{A}_0\\
		\vec{B}_3&=-2\s{\bar{\vec{A}_1}}{\vec{C}_1}\vec{A}_1+\frac{2}{\theta_0-3}\s{\vec{A}_1}{\vec{C}_1}\bar{\vec{A}_1}+\frac{2\bar{\alpha_0}}{\theta_0-3}\s{\vec{A}_1}{\vec{C}_1}\bar{\vec{A}_0}+2\left(\alpha_0\s{\bar{\vec{A}_1}}{\vec{C}_1}-\s{\bar{\vec{A}_1}}{\vec{C}_2}\right)\vec{A}_0\\
		\vec{E}_1&=-\frac{1}{2\theta_0}\s{\vec{C}_1}{\vec{C}_1}\bar{\vec{A}_0}.
		\end{alignedat}\right.
	\end{align}
	Integrating the equation
	\begin{align*}
		\p{\z}\left(\p{z}\phi\right)=\frac{e^{2\lambda}}{2}\H,
	\end{align*}
	one obtains for some $\vec{A}_3,\vec{A}_4\in\mathbb{C}^n$ the development
	\small
	\begin{align}
		\p{z}\phi=\begin{dmatrix}
		1 & A_{0} & \theta_{0} - 1 & 0 \\
		1 & A_{1} & \theta_{0} & 0 \\
		1 & A_{2} & \theta_{0} + 1 & 0 \\
		1 & A_{3} & \theta_{0} + 2 & 0 \\
		1 & A_{4} & \theta_{0} + 3 & 0 \\
		\frac{1}{4 \, \theta_{0}} & C_{1} & 1 & \theta_{0} \\
		\frac{1}{4 \, {\left(\theta_{0} + 1\right)}} & B_{1} & 1 & \theta_{0} + 1 \\
		\frac{\overline{\alpha_{0}}}{4 \, {\left(\theta_{0} + 1\right)}} & C_{1} & 1 & \theta_{0} + 1 \\
		\frac{1}{4 \, {\left(\theta_{0} + 2\right)}} & B_{2} & 1 & \theta_{0} + 2 \\
		\frac{\overline{\alpha_{1}}}{4 \, {\left(\theta_{0} + 2\right)}} & C_{1} & 1 & \theta_{0} + 2 \\
		\frac{\overline{\alpha_{0}}}{4 \, {\left(\theta_{0} + 2\right)}} & B_{1} & 1 & \theta_{0} + 2 \\
		\frac{1}{4 \, \theta_{0}} & C_{2} & 2 & \theta_{0} \\
		\frac{\alpha_{0}}{4 \, \theta_{0}} & C_{1} & 2 & \theta_{0} \\
		\frac{1}{4 \, {\left(\theta_{0} + 1\right)}} & B_{3} & 2 & \theta_{0} + 1 \\
		\frac{\overline{\alpha_{0}}}{4 \, {\left(\theta_{0} + 1\right)}} & C_{2} & 2 & \theta_{0} + 1 \\
		\frac{\alpha_{0}}{4 \, {\left(\theta_{0} + 1\right)}} & B_{1} & 2 & \theta_{0} + 1 \\
		\frac{{\left| A_{1} \right|}^{2}}{2 \, {\left(\theta_{0} + 1\right)}} & C_{1} & 2 & \theta_{0} + 1 
		\end{dmatrix}
		\begin{dmatrix}
		\frac{1}{4 \, \theta_{0}} & C_{3} & 3 & \theta_{0} \\
		\frac{\alpha_{1}}{4 \, \theta_{0}} & C_{1} & 3 & \theta_{0} \\
		\frac{\alpha_{0}}{4 \, \theta_{0}} & C_{2} & 3 & \theta_{0} \\
		\frac{1}{8 \, \theta_{0}} & E_{1} & -\theta_{0} + 3 & 2 \, \theta_{0} \\
		\frac{1}{8} & \overline{B_{1}} & \theta_{0} & 2 \\
		\frac{1}{8} \, \alpha_{0} & \overline{C_{1}} & \theta_{0} & 2 \\
		\frac{1}{12} & \overline{B_{3}} & \theta_{0} & 3 \\
		\frac{1}{12} \, \overline{\alpha_{0}} & \overline{B_{1}} & \theta_{0} & 3 \\
		\frac{1}{12} \, \alpha_{0} & \overline{C_{2}} & \theta_{0} & 3 \\
		\frac{1}{6} \, {\left| A_{1} \right|}^{2} & \overline{C_{1}} & \theta_{0} & 3 \\
		\frac{1}{8} & \overline{C_{1}} & \theta_{0} - 1 & 2 \\
		\frac{1}{12} & \overline{C_{2}} & \theta_{0} - 1 & 3 \\
		\frac{1}{12} \, \overline{\alpha_{0}} & \overline{C_{1}} & \theta_{0} - 1 & 3 \\
		\frac{1}{16} & \overline{C_{3}} & \theta_{0} - 1 & 4 \\
		\frac{1}{16} \, \overline{\alpha_{1}} & \overline{C_{1}} & \theta_{0} - 1 & 4 \\
		\frac{1}{16} \, \overline{\alpha_{0}} & \overline{C_{2}} & \theta_{0} - 1 & 4 \\
		\frac{1}{8} & \overline{B_{2}} & \theta_{0} + 1 & 2 \\
		\frac{1}{8} \, \alpha_{1} & \overline{C_{1}} & \theta_{0} + 1 & 2 \\
		\frac{1}{8} \, \alpha_{0} & \overline{B_{1}} & \theta_{0} + 1 & 2 \\
		-\frac{1}{4 \, {\left(\theta_{0} - 4\right)}} & \overline{E_{1}} & 2 \, \theta_{0} - 1 & -\theta_{0} + 4
		\end{dmatrix}+O(|z|^{\theta_0+4-\epsilon})
	\end{align}
	\normalsize
	
	Throwing all terms of degree more than equal to $\theta_0+2$ yields
	\begin{align*}
		\p{z}\phi=\begin{dmatrix}
		1 & A_{0} & \theta_{0} - 1 & 0 \\
		1 & A_{1} & \theta_{0} & 0 \\
		1 & A_{2} & \theta_{0} + 1 & 0 \\
		\frac{1}{4 \, \theta_{0}} & C_{1} & 1 & \theta_{0} \\
		\frac{1}{8} & \overline{C_{1}} & \theta_{0} - 1 & 2
		\end{dmatrix}
	\end{align*}
	so we recover \eqref{3rddev}. The next step is to look at relations given by the conformality of $\phi$. We have
	\begin{align*}
		0=\s{\p{z}\phi}{\p{z}\phi}=\begin{dmatrix}
		 A_{0}^{2} & 2 \, \theta_{0} - 2 & 0 \\
		2 \, A_{0} A_{1} & 2 \, \theta_{0} - 1 & 0 \\
		A_{1}^{2} + 2 \, A_{0} A_{2} & 2 \, \theta_{0} & 0 \\
		2 \, A_{1} A_{2} + 2 \, A_{0} A_{3} & 2 \, \theta_{0} + 1 & 0 \\
		A_{2}^{2} + 2 \, A_{1} A_{3} + 2 \, A_{0} A_{4} & 2 \, \theta_{0} + 2 & 0 \\
		\frac{A_{0} C_{1}}{2 \, \theta_{0}} & \theta_{0} & \theta_{0} \\
		\frac{A_{0} C_{1} \overline{\alpha_{0}} + A_{0} B_{1}}{2 \, {\left(\theta_{0} + 1\right)}} & \theta_{0} & \theta_{0} + 1 \\
		\frac{8 \, A_{0} B_{1} \theta_{0} \overline{\alpha_{0}} + 8 \, A_{0} C_{1} \theta_{0} \overline{\alpha_{1}} + 8 \, A_{0} B_{2} \theta_{0} + C_{1} {\left(\theta_{0} + 2\right)} \overline{C_{1}}}{16 \, {\left(\theta_{0}^{2} + 2 \, \theta_{0}\right)}} & \theta_{0} & \theta_{0} + 2 \\
		\frac{{\left(A_{0} \alpha_{0} + A_{1}\right)} C_{1} + A_{0} C_{2}}{2 \, \theta_{0}} & \theta_{0} + 1 & \theta_{0} \\
		\frac{A_{0} C_{2} \overline{\alpha_{0}} + {\left(A_{0} \alpha_{0} + A_{1}\right)} B_{1} + A_{0} B_{3} + {\left(2 \, A_{0} {\left| A_{1} \right|}^{2} + A_{1} \overline{\alpha_{0}}\right)} C_{1}}{2 \, {\left(\theta_{0} + 1\right)}} & \theta_{0} + 1 & \theta_{0} + 1 \\
		\frac{{\left(A_{1} \alpha_{0} + A_{0} \alpha_{1} + A_{2}\right)} C_{1} + {\left(A_{0} \alpha_{0} + A_{1}\right)} C_{2} + A_{0} C_{3}}{2 \, \theta_{0}} & \theta_{0} + 2 & \theta_{0} \\
		\frac{4 \, A_{0} E_{1} \theta_{0} + C_{1}^{2}}{16 \, \theta_{0}^{2}} & 2 & 2 \, \theta_{0} \\
		\frac{1}{4} \, A_{0} \overline{B_{1}} + \frac{1}{4} \, {\left(A_{0} \alpha_{0} + A_{1}\right)} \overline{C_{1}} & 2 \, \theta_{0} - 1 & 2 \\
		\frac{1}{6} \, A_{0} \overline{B_{1}} \overline{\alpha_{0}} + \frac{1}{6} \, A_{0} \overline{B_{3}} + \frac{1}{6} \, {\left(2 \, A_{0} {\left| A_{1} \right|}^{2} + A_{1} \overline{\alpha_{0}}\right)} \overline{C_{1}} + \frac{1}{6} \, {\left(A_{0} \alpha_{0} + A_{1}\right)} \overline{C_{2}} & 2 \, \theta_{0} - 1 & 3 \\
		\frac{1}{4} \, A_{0} \overline{C_{1}} & 2 \, \theta_{0} - 2 & 2 \\
		\frac{1}{6} \, A_{0} \overline{C_{1}} \overline{\alpha_{0}} + \frac{1}{6} \, A_{0} \overline{C_{2}} & 2 \, \theta_{0} - 2 & 3 \\
		\frac{1}{8} \, A_{0} \overline{C_{2}} \overline{\alpha_{0}} + \frac{1}{8} \, A_{0} \overline{C_{1}} \overline{\alpha_{1}} + \frac{1}{64} \, \overline{C_{1}}^{2} + \frac{1}{8} \, A_{0} \overline{C_{3}} & 2 \, \theta_{0} - 2 & 4 \\
		\frac{1}{4} \, {\left(A_{0} \alpha_{0} + A_{1}\right)} \overline{B_{1}} + \frac{1}{4} \, A_{0} \overline{B_{2}} + \frac{1}{4} \, {\left(A_{1} \alpha_{0} + A_{0} \alpha_{1} + A_{2}\right)} \overline{C_{1}} & 2 \, \theta_{0} & 2 \\
		-\frac{A_{0} \overline{E_{1}}}{2 \, {\left(\theta_{0} - 4\right)}} & 3 \, \theta_{0} - 2 & -\theta_{0} + 4
		\end{dmatrix}
	\end{align*}
		Most of these relations are trivial or uninteresting, and we will only need besides the previous ones (as $\s{\vec{A}_0}{\vec{C}_1}=\s{\vec{A}_0}{\bar{\vec{C}_1}}=0$)
	\begin{align}\label{a1c1}
	\s{\vec{A}_1}{\vec{A}_1}+2\s{\vec{A}_0}{\vec{A}_2}=0,\quad \s{\vec{A}_1}{\vec{C}_1}+\s{\vec{A}_0}{\vec{C}_2}=0,\quad \s{\vec{A}_0}{\bar{\vec{C}_2}}=0,
	\end{align}
	furnished by the lines $(2\theta_0, 0)$, $(\theta_0+1, \theta_0)$, and $(2\theta_0-2, 3)$. Therefore, we deduce that
	\begin{align*}
		0=\frac{1}{8} \, A_{0} \overline{C_{2}} \overline{\alpha_{0}} + \frac{1}{8} \, A_{0} \overline{C_{1}} \overline{\alpha_{1}} + \frac{1}{64} \, \overline{C_{1}}^{2} + \frac{1}{8} \, A_{0} \overline{C_{3}}=\frac{1}{64}\bar{\s{\vec{C}_1}{\vec{C}_1}}+\frac{1}{8}\s{\vec{A}_0}{\bar{\vec{C}_3}}
	\end{align*}
	so
	\begin{align*}
		\s{\vec{C}_1}{\vec{C}_1}+8\s{\bar{\vec{A}_0}}{\vec{C}_3}=0.
	\end{align*}
	
	We now stress out all relations we shall need in the sequel
	\begin{align}\label{a1c12}
	\left\{\begin{alignedat}{1}
		&\s{\vec{A}_0}{\vec{A}_0}=\s{\vec{A}_0}{\vec{A}_1}=\s{\vec{A}_0}{\vec{C}_1}=\s{\vec{A}_0}{\bar{\vec{C}_1}}=\s{\vec{A}_0}{\bar{\vec{C}_2}}=0\\
		&\s{\vec{A}_1}{\vec{A}_1}+2\s{\vec{A}_0}{\vec{A}_2}=0,\quad \s{\vec{A}_1}{\vec{C}_1}+\s{\vec{A}_0}{\vec{C}_2}=0,\quad \s{\vec{C}_1}{\vec{C}_1}+8\s{\bar{\vec{A}_0}}{\vec{C}_3}=0.
		\end{alignedat}\right.
	\end{align}
    As $e^{2\lambda}=2\s{\p{z}\phi}{\p{\z}\phi}$, we compute
    \small
    \begin{align*}
    	e^{2\lambda}=\begin{dmatrix}
    	 2 \, A_{0} \overline{A_{0}} & \theta_{0} - 1 & \theta_{0} - 1 \\
    	2 \, A_{0} \overline{A_{1}} & \theta_{0} - 1 & \theta_{0} \\
    	2 \, A_{0} \overline{A_{2}} + \frac{1}{4} \, \overline{A_{0}} \overline{C_{1}} & \theta_{0} - 1 & \theta_{0} + 1 \\
    	2 \, A_{0} \overline{A_{3}} + \frac{1}{12} \, {\left(2 \, \overline{A_{0}} \overline{\alpha_{0}} + 3 \, \overline{A_{1}}\right)} \overline{C_{1}} + \frac{1}{6} \, \overline{A_{0}} \overline{C_{2}} & \theta_{0} - 1 & \theta_{0} + 2 \\
    	2 \, A_{0} \overline{A_{4}} + \frac{1}{24} \, {\left(4 \, \overline{A_{1}} \overline{\alpha_{0}} + 3 \, \overline{A_{0}} \overline{\alpha_{1}} + 6 \, \overline{A_{2}}\right)} \overline{C_{1}} + \frac{1}{24} \, {\left(3 \, \overline{A_{0}} \overline{\alpha_{0}} + 4 \, \overline{A_{1}}\right)} \overline{C_{2}} + \frac{1}{8} \, \overline{A_{0}} \overline{C_{3}} & \theta_{0} - 1 & \theta_{0} + 3 \\
    	\frac{\ccancel{A_{0} \overline{C_{1}}}}{2 \, \theta_{0}} & 2 \, \theta_{0} - 1 & 1 \\
    	\frac{A_{0} \theta_{0} \overline{B_{1}} + {\left(A_{0} \alpha_{0} \theta_{0} + A_{1} {\left(\theta_{0} + 1\right)}\right)} \overline{C_{1}}}{2 \, {\left(\theta_{0}^{2} + \theta_{0}\right)}} & 2 \, \theta_{0} & 1 \\
    	\mu_1 & 2 \, \theta_{0} + 1 & 1 \\
    	\frac{\ccancel{A_{0} \overline{C_{1}} \overline{\alpha_{0}}} + \ccancel{A_{0} \overline{C_{2}}}}{2 \, \theta_{0}} & 2 \, \theta_{0} - 1 & 2 \\
    	\frac{A_{0} \theta_{0} \overline{B_{1}} \overline{\alpha_{0}} + A_{0} \theta_{0} \overline{B_{3}} + {\left(2 \, A_{0} \theta_{0} {\left| A_{1} \right|}^{2} + {\left(\theta_{0} \overline{\alpha_{0}} + \overline{\alpha_{0}}\right)} A_{1}\right)} \overline{C_{1}} + {\left(A_{0} \alpha_{0} \theta_{0} + A_{1} {\left(\theta_{0} + 1\right)}\right)} \overline{C_{2}}}{2 \, {\left(\theta_{0}^{2} + \theta_{0}\right)}} & 2 \, \theta_{0} & 2 \\
    	\frac{8 \, {\left(\theta_{0} \overline{\alpha_{1}} - 4 \, \overline{\alpha_{1}}\right)} A_{0} \overline{C_{1}} + {\left(\theta_{0} - 4\right)} \overline{C_{1}}^{2} + 8 \, {\left(\theta_{0} \overline{\alpha_{0}} - 4 \, \overline{\alpha_{0}}\right)} A_{0} \overline{C_{2}} + 8 \, A_{0} {\left(\theta_{0} - 4\right)} \overline{C_{3}} - 8 \, \theta_{0} \overline{A_{0}} \overline{E_{1}}}{16 \, {\left(\theta_{0}^{2} - 4 \, \theta_{0}\right)}} & 2 \, \theta_{0} - 1 & 3 \\
    	\frac{\ccancel{A_{0} \overline{E_{1}}}}{4 \, \theta_{0}} & 3 \, \theta_{0} - 1 & -\theta_{0} + 3 \\
    	\frac{1}{4} \, A_{0} C_{1} \overline{\alpha_{0}} + \frac{1}{4} \, A_{0} B_{1} + 2 \, A_{2} \overline{A_{1}} & \theta_{0} + 1 & \theta_{0} \\
    	\frac{1}{6} \, A_{0} C_{2} \overline{\alpha_{0}} + \frac{1}{12} \, {\left(2 \, A_{0} \alpha_{0} + 3 \, A_{1}\right)} B_{1} + \frac{1}{6} \, A_{0} B_{3} + \frac{1}{12} \, {\left(4 \, A_{0} {\left| A_{1} \right|}^{2} + 3 \, A_{1} \overline{\alpha_{0}}\right)} C_{1} + 2 \, A_{3} \overline{A_{1}} & \theta_{0} + 2 & \theta_{0} \\
    	\frac{1}{4} \, A_{0} C_{1} + 2 \, A_{2} \overline{A_{0}} & \theta_{0} + 1 & \theta_{0} - 1 \\
    	\frac{1}{12} \, {\left(2 \, A_{0} \alpha_{0} + 3 \, A_{1}\right)} C_{1} + \frac{1}{6} \, A_{0} C_{2} + 2 \, A_{3} \overline{A_{0}} & \theta_{0} + 2 & \theta_{0} - 1 \\
    	\frac{1}{24} \, {\left(4 \, A_{1} \alpha_{0} + 3 \, A_{0} \alpha_{1} + 6 \, A_{2}\right)} C_{1} + \frac{1}{24} \, {\left(3 \, A_{0} \alpha_{0} + 4 \, A_{1}\right)} C_{2} + \frac{1}{8} \, A_{0} C_{3} + 2 \, A_{4} \overline{A_{0}} & \theta_{0} + 3 & \theta_{0} - 1 \\
    	\mu_2 & \theta_{0} + 1 & \theta_{0} + 1 \\
    	\frac{C_{1}^{2} {\left(\theta_{0} - 4\right)} - 8 \, A_{0} E_{1} \theta_{0} + 8 \, {\left(\alpha_{1} \theta_{0} - 4 \, \alpha_{1}\right)} C_{1} \overline{A_{0}} + 8 \, {\left(\alpha_{0} \theta_{0} - 4 \, \alpha_{0}\right)} C_{2} \overline{A_{0}} + 8 \, C_{3} {\left(\theta_{0} - 4\right)} \overline{A_{0}}}{16 \, {\left(\theta_{0}^{2} - 4 \, \theta_{0}\right)}} & 3 & 2 \, \theta_{0} - 1 \\
    	2 \, A_{1} \overline{A_{0}} & \theta_{0} & \theta_{0} - 1 \\
    	2 \, A_{1} \overline{A_{1}} & \theta_{0} & \theta_{0} \\
    	\frac{1}{4} \, \alpha_{0} \overline{A_{0}} \overline{C_{1}} + 2 \, A_{1} \overline{A_{2}} + \frac{1}{4} \, \overline{A_{0}} \overline{B_{1}} & \theta_{0} & \theta_{0} + 1 \\
    	\frac{1}{6} \, \alpha_{0} \overline{A_{0}} \overline{C_{2}} + 2 \, A_{1} \overline{A_{3}} + \frac{1}{12} \, {\left(2 \, \overline{A_{0}} \overline{\alpha_{0}} + 3 \, \overline{A_{1}}\right)} \overline{B_{1}} + \frac{1}{6} \, \overline{A_{0}} \overline{B_{3}} + \frac{1}{12} \, {\left(4 \, {\left| A_{1} \right|}^{2} \overline{A_{0}} + 3 \, \alpha_{0} \overline{A_{1}}\right)} \overline{C_{1}} & \theta_{0} & \theta_{0} + 2 \\
    	\frac{\ccancel{C_{1} \overline{A_{0}}}}{2 \, \theta_{0}} & 1 & 2 \, \theta_{0} - 1 \\
    	\frac{B_{1} \theta_{0} \overline{A_{0}} + {\left(\theta_{0} \overline{A_{0}} \overline{\alpha_{0}} + {\left(\theta_{0} + 1\right)} \overline{A_{1}}\right)} C_{1}}{2 \, {\left(\theta_{0}^{2} + \theta_{0}\right)}} & 1 & 2 \, \theta_{0} \\
    	\mu_3 & 1 & 2 \, \theta_{0} + 1 \\
    	\frac{\ccancel{C_{1} \alpha_{0} \overline{A_{0}}} + \ccancel{C_{2} \overline{A_{0}}}}{2 \, \theta_{0}} & 2 & 2 \, \theta_{0} - 1 \\
    	\frac{B_{1} \alpha_{0} \theta_{0} \overline{A_{0}} + B_{3} \theta_{0} \overline{A_{0}} + {\left(2 \, \theta_{0} {\left| A_{1} \right|}^{2} \overline{A_{0}} + {\left(\alpha_{0} \theta_{0} + \alpha_{0}\right)} \overline{A_{1}}\right)} C_{1} + {\left(\theta_{0} \overline{A_{0}} \overline{\alpha_{0}} + {\left(\theta_{0} + 1\right)} \overline{A_{1}}\right)} C_{2}}{2 \, {\left(\theta_{0}^{2} + \theta_{0}\right)}} & 2 & 2 \, \theta_{0} \\
    	\frac{\ccancel{E_{1} \overline{A_{0}}}}{4 \, \theta_{0}} & -\theta_{0} + 3 & 3 \, \theta_{0} - 1
    	\end{dmatrix}
     \end{align*}
    \normalsize
    where
    \begin{align*}
    	\mu_1&=\frac{1}{2 \, {\left(\theta_{0}^{3} + 3 \, \theta_{0}^{2} + 2 \, \theta_{0}\right)}}\bigg\{{\left(\theta_{0}^{2} + \theta_{0}\right)} A_{0} \overline{B_{2}} + {\left({\left(\alpha_{0} \theta_{0}^{2} + \alpha_{0} \theta_{0}\right)} A_{0} + {\left(\theta_{0}^{2} + 2 \, \theta_{0}\right)} A_{1}\right)} \overline{B_{1}}\\
    	& + {\left({\left(\alpha_{1} \theta_{0}^{2} + \alpha_{1} \theta_{0}\right)} A_{0} + {\left(\alpha_{0} \theta_{0}^{2} + 2 \, \alpha_{0} \theta_{0}\right)} A_{1} + {\left(\theta_{0}^{2} + 3 \, \theta_{0} + 2\right)} A_{2}\right)} \overline{C_{1}}\bigg\} \\
    	\mu_2&=\frac{1}{32 \, \theta_{0}^{2}}\bigg\{8 \, \alpha_{0} \theta_{0}^{2} \overline{A_{0}} \overline{B_{1}} + 8 \, A_{0} B_{1} \theta_{0}^{2} \overline{\alpha_{0}} + 8 \, A_{0} C_{1} \theta_{0}^{2} \overline{\alpha_{1}} + 8 \, A_{0} B_{2} \theta_{0}^{2} + 64 \, A_{2} \theta_{0}^{2} \overline{A_{2}} + 8 \, \theta_{0}^{2} \overline{A_{0}} \overline{B_{2}}\\
    	& + {\left(8 \, \alpha_{1} \theta_{0}^{2} \overline{A_{0}} + {\left(\theta_{0}^{2} + 4\right)} C_{1}\right)} \overline{C_{1}}\bigg\}\\
    	\mu_3&=\frac{1}{2 \, {\left(\theta_{0}^{3} + 3 \, \theta_{0}^{2} + 2 \, \theta_{0}\right)}}\bigg\{{\left(\theta_{0}^{2} + \theta_{0}\right)} B_{2} \overline{A_{0}} + {\left({\left(\theta_{0}^{2} \overline{\alpha_{0}} + \theta_{0} \overline{\alpha_{0}}\right)} \overline{A_{0}} + {\left(\theta_{0}^{2} + 2 \, \theta_{0}\right)} \overline{A_{1}}\right)} B_{1}\\
    	& + {\left({\left(\theta_{0}^{2} \overline{\alpha_{1}} + \theta_{0} \overline{\alpha_{1}}\right)} \overline{A_{0}} + {\left(\theta_{0}^{2} \overline{\alpha_{0}} + 2 \, \theta_{0} \overline{\alpha_{0}}\right)} \overline{A_{1}} + {\left(\theta_{0}^{2} + 3 \, \theta_{0} + 2\right)} \overline{A_{2}}\right)} C_{1}\bigg\}.
    \end{align*}
    As $\vec{E}_1\in \mathrm{Span}(\bar{\vec{A}_0})$ thanks of \eqref{b3c1}, we have
    \begin{align*}
    	\s{\vec{E}_1}{\bar{\vec{A}_0}}=0
    \end{align*}
    so there is no coefficient in $\Re\left(\ast \dfrac{\z^{3\theta_0-1}}{z^{\theta_0-3}}\right)$ in the Taylor expansion of $e^{2\lambda}$. Also, thanks of \eqref{a1c12}, we know that
    \begin{align*}
    	\s{\bar{\vec{A}_0}}{\vec{C}_1}=0
    \end{align*}
    the coefficient in $\Re(\ast z\z^{2\theta_0-1})$ also vanish.
    Also, by \eqref{a1c12}, we have
    \begin{align*}
    	\s{\bar{\vec{A}_0}}{\vec{C}_2}=0
    \end{align*}
    so the coefficient in $\Re(\ast z^2\z^{2\theta_0-2})$ also vanishes. 
    
    Finally, we have
    \small
    \begin{align*}
    	e^{2\lambda}=\begin{dmatrix}
    	 2 \, A_{0} \overline{A_{0}} & \theta_{0} - 1 & \theta_{0} - 1 \\
    	2 \, A_{0} \overline{A_{1}} & \theta_{0} - 1 & \theta_{0} \\
    	2 \, A_{0} \overline{A_{2}} + \frac{1}{4} \, \overline{A_{0}} \overline{C_{1}} & \theta_{0} - 1 & \theta_{0} + 1 \\
    	2 \, A_{0} \overline{A_{3}} + \frac{1}{12} \, {\left(2 \, \overline{A_{0}} \overline{\alpha_{0}} + 3 \, \overline{A_{1}}\right)} \overline{C_{1}} + \frac{1}{6} \, \overline{A_{0}} \overline{C_{2}} & \theta_{0} - 1 & \theta_{0} + 2 \\
    	2 \, A_{0} \overline{A_{4}} + \frac{1}{24} \, {\left(4 \, \overline{A_{1}} \overline{\alpha_{0}} + 3 \, \overline{A_{0}} \overline{\alpha_{1}} + 6 \, \overline{A_{2}}\right)} \overline{C_{1}} + \frac{1}{24} \, {\left(3 \, \overline{A_{0}} \overline{\alpha_{0}} + 4 \, \overline{A_{1}}\right)} \overline{C_{2}} + \frac{1}{8} \, \overline{A_{0}} \overline{C_{3}} & \theta_{0} - 1 & \theta_{0} + 3 \\
    	\frac{A_{0} \theta_{0} \overline{B_{1}} + {\left(A_{0} \alpha_{0} \theta_{0} + A_{1} {\left(\theta_{0} + 1\right)}\right)} \overline{C_{1}}}{2 \, {\left(\theta_{0}^{2} + \theta_{0}\right)}} & 2 \, \theta_{0} & 1 \\
    	\mu_1 & 2 \, \theta_{0} + 1 & 1 \\
    	\frac{A_{0} \theta_{0} \overline{B_{1}} \overline{\alpha_{0}} + A_{0} \theta_{0} \overline{B_{3}} + {\left(2 \, A_{0} \theta_{0} {\left| A_{1} \right|}^{2} + {\left(\theta_{0} \overline{\alpha_{0}} + \overline{\alpha_{0}}\right)} A_{1}\right)} \overline{C_{1}} + {\left(A_{0} \alpha_{0} \theta_{0} + A_{1} {\left(\theta_{0} + 1\right)}\right)} \overline{C_{2}}}{2 \, {\left(\theta_{0}^{2} + \theta_{0}\right)}} & 2 \, \theta_{0} & 2 \\
    	\frac{8 \, {\left(\theta_{0} \overline{\alpha_{1}} - 4 \, \overline{\alpha_{1}}\right)} A_{0} \overline{C_{1}} + {\left(\theta_{0} - 4\right)} \overline{C_{1}}^{2} + 8 \, {\left(\theta_{0} \overline{\alpha_{0}} - 4 \, \overline{\alpha_{0}}\right)} A_{0} \overline{C_{2}} + 8 \, A_{0} {\left(\theta_{0} - 4\right)} \overline{C_{3}} - 8 \, \theta_{0} \overline{A_{0}} \overline{E_{1}}}{16 \, {\left(\theta_{0}^{2} - 4 \, \theta_{0}\right)}} & 2 \, \theta_{0} - 1 & 3 \\
    	\frac{1}{4} \, A_{0} C_{1} \overline{\alpha_{0}} + \frac{1}{4} \, A_{0} B_{1} + 2 \, A_{2} \overline{A_{1}} & \theta_{0} + 1 & \theta_{0} \\
    	\frac{1}{6} \, A_{0} C_{2} \overline{\alpha_{0}} + \frac{1}{12} \, {\left(2 \, A_{0} \alpha_{0} + 3 \, A_{1}\right)} B_{1} + \frac{1}{6} \, A_{0} B_{3} + \frac{1}{12} \, {\left(4 \, A_{0} {\left| A_{1} \right|}^{2} + 3 \, A_{1} \overline{\alpha_{0}}\right)} C_{1} + 2 \, A_{3} \overline{A_{1}} & \theta_{0} + 2 & \theta_{0} \\
    	\frac{1}{4} \, A_{0} C_{1} + 2 \, A_{2} \overline{A_{0}} & \theta_{0} + 1 & \theta_{0} - 1 \\
    	\frac{1}{12} \, {\left(2 \, A_{0} \alpha_{0} + 3 \, A_{1}\right)} C_{1} + \frac{1}{6} \, A_{0} C_{2} + 2 \, A_{3} \overline{A_{0}} & \theta_{0} + 2 & \theta_{0} - 1 \\
    	\frac{1}{24} \, {\left(4 \, A_{1} \alpha_{0} + 3 \, A_{0} \alpha_{1} + 6 \, A_{2}\right)} C_{1} + \frac{1}{24} \, {\left(3 \, A_{0} \alpha_{0} + 4 \, A_{1}\right)} C_{2} + \frac{1}{8} \, A_{0} C_{3} + 2 \, A_{4} \overline{A_{0}} & \theta_{0} + 3 & \theta_{0} - 1 \\
    	\mu_2 & \theta_{0} + 1 & \theta_{0} + 1 \\
    	\frac{C_{1}^{2} {\left(\theta_{0} - 4\right)} - 8 \, A_{0} E_{1} \theta_{0} + 8 \, {\left(\alpha_{1} \theta_{0} - 4 \, \alpha_{1}\right)} C_{1} \overline{A_{0}} + 8 \, {\left(\alpha_{0} \theta_{0} - 4 \, \alpha_{0}\right)} C_{2} \overline{A_{0}} + 8 \, C_{3} {\left(\theta_{0} - 4\right)} \overline{A_{0}}}{16 \, {\left(\theta_{0}^{2} - 4 \, \theta_{0}\right)}} & 3 & 2 \, \theta_{0} - 1 \\
    	2 \, A_{1} \overline{A_{0}} & \theta_{0} & \theta_{0} - 1 \\
    	2 \, A_{1} \overline{A_{1}} & \theta_{0} & \theta_{0} \\
    	\frac{1}{4} \, \alpha_{0} \overline{A_{0}} \overline{C_{1}} + 2 \, A_{1} \overline{A_{2}} + \frac{1}{4} \, \overline{A_{0}} \overline{B_{1}} & \theta_{0} & \theta_{0} + 1 \\
    	\frac{1}{6} \, \alpha_{0} \overline{A_{0}} \overline{C_{2}} + 2 \, A_{1} \overline{A_{3}} + \frac{1}{12} \, {\left(2 \, \overline{A_{0}} \overline{\alpha_{0}} + 3 \, \overline{A_{1}}\right)} \overline{B_{1}} + \frac{1}{6} \, \overline{A_{0}} \overline{B_{3}} + \frac{1}{12} \, {\left(4 \, {\left| A_{1} \right|}^{2} \overline{A_{0}} + 3 \, \alpha_{0} \overline{A_{1}}\right)} \overline{C_{1}} & \theta_{0} & \theta_{0} + 2 \\
    	\frac{B_{1} \theta_{0} \overline{A_{0}} + {\left(\theta_{0} \overline{A_{0}} \overline{\alpha_{0}} + {\left(\theta_{0} + 1\right)} \overline{A_{1}}\right)} C_{1}}{2 \, {\left(\theta_{0}^{2} + \theta_{0}\right)}} & 1 & 2 \, \theta_{0} \\
    	\mu_3 & 1 & 2 \, \theta_{0} + 1 \\
    	\frac{B_{1} \alpha_{0} \theta_{0} \overline{A_{0}} + B_{3} \theta_{0} \overline{A_{0}} + {\left(2 \, \theta_{0} {\left| A_{1} \right|}^{2} \overline{A_{0}} + {\left(\alpha_{0} \theta_{0} + \alpha_{0}\right)} \overline{A_{1}}\right)} C_{1} + {\left(\theta_{0} \overline{A_{0}} \overline{\alpha_{0}} + {\left(\theta_{0} + 1\right)} \overline{A_{1}}\right)} C_{2}}{2 \, {\left(\theta_{0}^{2} + \theta_{0}\right)}} & 2 & 2 \, \theta_{0} \\
    	\end{dmatrix}
    \end{align*}
    \normalsize
    for some uninteresting $\mu_3\in\C$.
    First, recall that by \eqref{b3c1}
    \begin{align*}
    	\vec{E}_1=-\frac{1}{2\theta_0}\s{\vec{C}_1}{\vec{C}_1}\bar{\vec{A}_0}
    \end{align*}
    while by \eqref{a1c12}
    \begin{align*}
    \s{\vec{C}_1}{\vec{C}_1}+8\s{\bar{\vec{A}_0}}{\vec{C}_3}=0,
    \end{align*}
    so we deduce that
    \begin{align}\label{defalpha7}
    	&\frac{\colorcancel{C_{1}^{2} {\left(\theta_{0} - 4\right)}}{blue} - 8 \, A_{0} E_{1} \theta_{0} + 8 \, {\left(\alpha_{1} \theta_{0} - 4 \, \alpha_{1}\right)} \ccancel{C_{1} \overline{A_{0}}} + 8 \, {\left(\alpha_{0} \theta_{0} - 4 \, \alpha_{0}\right)} \ccancel{C_{2} \overline{A_{0}}} + \colorcancel{8 \, C_{3} {\left(\theta_{0} - 4\right)} \overline{A_{0}}}{blue}}{16 \, {\left(\theta_{0}^{2} - 4 \, \theta_{0}\right)}}=-\frac{\theta_0}{2\theta_0(\theta_0-4)}\s{\vec{A}_0}{\vec{E}_1}\nonumber\\
    	&=-\frac{1}{2(\theta_0-4)}\bs{\vec{A}_0}{-\frac{1}{2\theta_0}\s{\vec{C}_1}{\vec{C}_1}\bar{\vec{A}_0}}=\frac{1}{8\theta_0(\theta_0-4)}\s{\vec{C}_1}{\vec{C}_1}.
    \end{align}
    as $|\vec{A}_0|^2=\dfrac{1}{2}$.
    
    In particular, recalling that $|\vec{A}_0|^2=\dfrac{1}{2}$, we see that there exists some constants
    \begin{align*}
    	\beta\in\R,\quad \alpha_2,\alpha_3,\alpha_4,\alpha_5,\alpha_6,\alpha_7,\alpha_8,\alpha_9\in \mathbb{C}
    \end{align*}
	such that
	\begin{align*}
		e^{2\lambda}&=|z|^{2\theta_0-2}+2|\vec{A}_1|^2|z|^{2\theta_0}+\beta|z|^{2\theta_0+2}+2\,\Re\bigg(\alpha_0z^{\theta_0}\z^{\theta_0-1}+\alpha_1z^{\theta_0+1}\z^{\theta_0-1}+\alpha_2z\z^{2\theta_0}+\alpha_3z^{\theta_0+2}\z^{\theta_0-1}+\alpha_4z^{\theta_0+3}\z^{\theta_0-1}\\
		&+\alpha_5z^{\theta_0+1}\z^{\theta_0}+\alpha_6z^{\theta_0+2}\z^{\theta_0}+\alpha_7z^3\z^{2\theta_0-1}+\alpha_8z^2\z^{2\theta_0}+\alpha_9z\z^{2\theta_0+1}\bigg)+O(|z|^{2\theta_0+3-\epsilon})
	\end{align*}
	so that
    \begin{align*}
    	e^{2\lambda}=\begin{dmatrix}
    	1 & \theta_{0} - 1 & \theta_{0} - 1 \\
    	2 \, {\left| A_{1} \right|}^{2} & \theta_{0} & \theta_{0} \\
    	\beta & \theta_{0} + 1 & \theta_{0} + 1 \\
    	\alpha_{0} & \theta_{0} & \theta_{0} - 1 \\
    	\alpha_{1} & \theta_{0} + 1 & \theta_{0} - 1 \\
    	\alpha_{2} & 1 & 2 \, \theta_{0} \\
    	\alpha_{3} & \theta_{0} + 2 & \theta_{0} - 1 \\
    	\alpha_{4} & \theta_{0} + 3 & \theta_{0} - 1 \\
    	\alpha_{5} & \theta_{0} + 1 & \theta_{0} \\
    	\alpha_{6} & \theta_{0} + 2 & \theta_{0} \\
    	\alpha_{7} & 3 & 2 \, \theta_{0} - 1 \\
    	\alpha_{8} & 2 & 2 \, \theta_{0} \\
    	\alpha_{9} & 1 & 2 \, \theta_{0} + 1 
    	\end{dmatrix}
    	\begin{dmatrix}
    	\overline{\alpha_{0}} & \theta_{0} - 1 & \theta_{0} \\
    	\overline{\alpha_{1}} & \theta_{0} - 1 & \theta_{0} + 1 \\
    	\overline{\alpha_{2}} & 2 \, \theta_{0} & 1 \\
    	\overline{\alpha_{3}} & \theta_{0} - 1 & \theta_{0} + 2 \\
    	\overline{\alpha_{4}} & \theta_{0} - 1 & \theta_{0} + 3 \\
    	\overline{\alpha_{5}} & \theta_{0} & \theta_{0} + 1 \\
    	\overline{\alpha_{6}} & \theta_{0} & \theta_{0} + 2 \\
    	\overline{\alpha_{7}} & 2 \, \theta_{0} - 1 & 3 \\
    	\overline{\alpha_{8}} & 2 \, \theta_{0} & 2 \\
    	\overline{\alpha_{9}} & 2 \, \theta_{0} + 1 & 1
    	\end{dmatrix}
    \end{align*}
    In these coefficients, the only interesting ones are $|\vec{A}_1|^2$, $\alpha_0$ and $\alpha_2$. As we need the exact formula for $\alpha_2$, we first recall that
    \begin{align*}
        \vec{B}_1=-2\s{\bar{\vec{A}_1}}{\vec{C}_1}\vec{A}_0,\quad |\vec{A}_0|^2=\frac{1}{2},\quad \s{\bar{\vec{A}_0}}{\vec{C}_1}=0.
    \end{align*}
    We deduce that 
    \begin{align*}
    \s{\bar{\vec{A}_0}}{\vec{B}_1}=\s{\bar{\vec{A}_0}}{-2\s{\bar{\vec{A}_1}}{\vec{C}_1}\vec{A}_0}=-\s{\bar{\vec{A}_1}}{\vec{C}_1}
    \end{align*}
    and finally
    \begin{align}\label{alpha2}
    	\alpha_2&=\frac{B_{1} \theta_{0} \overline{A_{0}} + {\left(\ccancel{\theta_{0} \overline{A_{0}} \overline{\alpha_{0}}} + {\left(\theta_{0} + 1\right)} \overline{A_{1}}\right)} C_{1}}{2 \, {\left(\theta_{0}^{2} + \theta_{0}\right)}}=\frac{1}{2\theta_0(\theta_0+1)}\left(-\theta_0\s{\bar{\vec{A}_1}}{\vec{C}_1}+(\theta_0+1)\s{\bar{\vec{A}_1}}{\vec{C}_1}\right)\nonumber\\
    	&=\frac{1}{2\theta_0(\theta_0+1)}\s{\bar{\vec{A}_1}}{\vec{C}_1}
    \end{align}
    while by \eqref{defalpha7}, we have
    \begin{align}
    	\alpha_7=\frac{1}{8\theta_0(\theta_0-4)}\s{\vec{C}_1}{\vec{C}_1}.
    \end{align}
    Now, we have
    \begin{align*}
    	\alpha_3&=\frac{1}{12} \, {\left(\ccancel{2 \, A_{0} \alpha_{0}} + 3 \, A_{1}\right)} C_{1} + \frac{1}{6} \, A_{0} C_{2} + 2 \, A_{3} \overline{A_{0}}\\
    	&=\frac{1}{4}\s{\vec{A}_1}{\vec{C}_1}+\frac{1}{6}\s{\vec{A}_0}{\vec{C}_2}+2\s{\bar{\vec{A}_0}}{\vec{A}_3}=\frac{1}{12}\s{\vec{A}_1}{\vec{C}_1}+2\s{\bar{\vec{A}_0}}{\vec{A}_3}
    \end{align*}
    Finally, we have (using the result not proved yet that $\alpha_0=0$, but $\alpha_6$ does not enter into play in  this proof at the reader may check)
    \begin{align*}
    	\alpha_6&=\frac{1}{6} \, A_{0} C_{2} \ccancel{\overline{\alpha_{0}}} + \frac{1}{12} \, {\left(\ccancel{2 \, A_{0} \alpha_{0}} + \ccancel{3 \, A_{1}}\right)} B_{1} + \ccancel{\frac{1}{6} \, A_{0} B_{3}} + \frac{1}{12} \, {\left(\ccancel{4 \, A_{0} {\left| A_{1} \right|}^{2}} + 3 \, A_{1} \ccancel{\overline{\alpha_{0}}}\right)} C_{1} + 2 \, A_{3} \overline{A_{1}}\\
    	&=2\s{\bar{\vec{A}_1}}{\vec{A}_3}.
    \end{align*}
    Now, we also have as  $\s{\vec{A}_1}{\vec{C}_1}+\s{\vec{A}_0}{\vec{C}_2}=0$ by conformality
    \begin{align*}
    	0={\left(\colorcancel{A_{1} \alpha_{0}}{blue} + \ccancel{A_{0} \alpha_{1}} + A_{2}\right)} C_{1} + {\left(\colorcancel{A_{0} \alpha_{0}}{blue} + A_{1}\right)} C_{2} + A_{0} C_{3}=\s{\vec{A}_2}{\vec{C}_1}+\s{\vec{A}_1}{\vec{C}_2}+\s{\vec{A}_0}{\vec{C}_3}
    \end{align*}
    and
    the development
    \begin{align*}
    	\alpha_4&=\frac{1}{24} \, {\left(4 \, A_{1} \alpha_{0} + \ccancel{3 \, A_{0} \alpha_{1}} + 6 \, A_{2}\right)} C_{1} + \frac{1}{24} \, {\left(3 \, A_{0} \alpha_{0} + 4 \, A_{1}\right)} C_{2} + \frac{1}{8} \, A_{0} C_{3} + 2 \, A_{4} \overline{A_{0}}\\
    	&=\frac{1}{6}\alpha_0\s{\vec{A}_1}{\vec{C}_1}+\frac{1}{4}\s{\vec{A}_2}{\vec{C}_1}-\frac{1}{8}\alpha_0\s{\vec{A}_1}{\vec{C}_1}+\frac{1}{6}\s{\vec{A}_1}{\vec{C}_2}-\frac{1}{8}\left(\s{\vec{A}_2}{\vec{C}_1}+\s{\vec{A}_1}{\vec{C}_2}\right)+2\s{\bar{\vec{A}_0}}{\vec{A}_4}\\
    	&=\frac{1}{24}\alpha_0\s{\vec{A}_1}{\vec{C}_1}+\frac{1}{24}\left(3\s{\vec{A}_2}{\vec{C}_1}+\s{\vec{A}_1}{\vec{C}_2}\right)+2\s{\bar{\vec{A}_0}}{\vec{A}_4}.
    \end{align*}
    
    \begin{align}\label{alpha27}
    	\left\{
    	\begin{alignedat}{1}
    	\alpha_1&=2\s{\bar{\vec{A}_0}}{\vec{A}_2}\\
    	\alpha_2&=\frac{1}{2\theta_0(\theta_0+1)}\s{\bar{\vec{A}_1}}{\vec{C}_1}\\
    	\alpha_3&=\frac{1}{12}\s{\vec{A}_1}{\vec{C}_1}+2\s{\bar{\vec{A}_0}}{\vec{A}_3}\\
    	\alpha_4&=\frac{1}{24}\alpha_0\s{\vec{A}_1}{\vec{C}_1}+\frac{1}{24}\left(3\s{\vec{A}_2}{\vec{C}_1}+\s{\vec{A}_1}{\vec{C}_2}\right)+2\s{\bar{\vec{A}_0}}{\vec{A}_4}\\
    	\alpha_5&=2\s{\bar{\vec{A}_1}}{\vec{A}_2}\\
    	\alpha_6&=2\s{\bar{\vec{A}_1}}{\vec{A}_3}\\
    	\alpha_7&=\frac{1}{8\theta_0(\theta_0-4)}\s{\vec{C}_1}{\vec{C}_1}
    	\end{alignedat}\right.
    \end{align}
    We now obtain
    \begin{align*}
    	\h_0=\begin{dmatrix}
    	2 & A_{1} & \theta_{0} - 1 & 0 \\
    	4 & A_{2} & \theta_{0} & 0 \\
    	6 & A_{3} & \theta_{0} + 1 & 0 \\
    	8 & A_{4} & \theta_{0} + 2 & 0 \\
    	-\frac{\theta_{0} - 2}{2 \, \theta_{0}} & C_{1} & 0 & \theta_{0} \\
    	-\frac{\theta_{0} - 2}{2 \, {\left(\theta_{0} + 1\right)}} & B_{1} & 0 & \theta_{0} + 1 \\
    	-\frac{\theta_{0} \overline{\alpha_{0}} - 2 \, \overline{\alpha_{0}}}{2 \, {\left(\theta_{0} + 1\right)}} & C_{1} & 0 & \theta_{0} + 1 \\
    	-\frac{\theta_{0} - 2}{2 \, {\left(\theta_{0} + 2\right)}} & B_{2} & 0 & \theta_{0} + 2 \\
    	-\frac{\theta_{0} \overline{\alpha_{1}} - 2 \, \overline{\alpha_{1}}}{2 \, {\left(\theta_{0} + 2\right)}} & C_{1} & 0 & \theta_{0} + 2 \\
    	-\frac{\theta_{0} \overline{\alpha_{0}} - 2 \, \overline{\alpha_{0}}}{2 \, {\left(\theta_{0} + 2\right)}} & B_{1} & 0 & \theta_{0} + 2 \\
    	-\frac{\theta_{0} - 3}{2 \, \theta_{0}} & C_{2} & 1 & \theta_{0} \\
    	-\frac{\alpha_{0} \theta_{0} - 2 \, \alpha_{0}}{2 \, \theta_{0}} & C_{1} & 1 & \theta_{0} \\
    	-\frac{\theta_{0} - 3}{2 \, {\left(\theta_{0} + 1\right)}} & B_{3} & 1 & \theta_{0} + 1 \\
    	-\frac{\theta_{0} \overline{\alpha_{0}} - 3 \, \overline{\alpha_{0}}}{2 \, {\left(\theta_{0} + 1\right)}} & C_{2} & 1 & \theta_{0} + 1 \\
    	-\frac{\alpha_{0} \theta_{0} - 2 \, \alpha_{0}}{2 \, {\left(\theta_{0} + 1\right)}} & B_{1} & 1 & \theta_{0} + 1 \\
    	-\frac{2 \, {\left(\theta_{0}^{2} - 2 \, \theta_{0} + 1\right)} {\left| A_{1} \right|}^{2} - \alpha_{0} \overline{\alpha_{0}}}{2 \, {\left(\theta_{0}^{2} + \theta_{0}\right)}} & C_{1} & 1 & \theta_{0} + 1 \\
    	-\frac{\theta_{0} - 4}{2 \, \theta_{0}} & C_{3} & 2 & \theta_{0} \\
    	-\frac{\alpha_{1} \theta_{0} - 2 \, \alpha_{1}}{2 \, \theta_{0}} & C_{1} & 2 & \theta_{0} \\
    	-\frac{\alpha_{0} \theta_{0} - 3 \, \alpha_{0}}{2 \, \theta_{0}} & C_{2} & 2 & \theta_{0} \\
    	-\frac{\theta_{0} - 2}{2 \, \theta_{0}} & E_{1} & -\theta_{0} + 2 & 2 \, \theta_{0} \\
    	\frac{1}{4} & \overline{B_{1}} & \theta_{0} - 1 & 2 \\
    	\frac{1}{6} & \overline{B_{3}} & \theta_{0} - 1 & 3 \\
    	\frac{1}{6} \, \overline{\alpha_{0}} & \overline{B_{1}} & \theta_{0} - 1 & 3 \\
    	-\frac{1}{6} \, {\left| A_{1} \right|}^{2} + \frac{1}{12} \, \alpha_{0} \overline{\alpha_{0}} & \overline{C_{1}} & \theta_{0} - 1 & 3 \\
    	\frac{1}{2} & \overline{B_{2}} & \theta_{0} & 2 \\
    	\frac{1}{4} \, \alpha_{0} & \overline{B_{1}} & \theta_{0} & 2 
    	    	\end{dmatrix}
    	\end{align*}
    	\begin{align*}
    	\begin{dmatrix}
    	-\frac{\theta_{0}}{2 \, {\left(\theta_{0} - 4\right)}} & \overline{E_{1}} & 2 \, \theta_{0} - 2 & -\theta_{0} + 4 \\
    	-4 \, {\left| A_{1} \right|}^{2} + 2 \, \alpha_{0} \overline{\alpha_{0}} & A_{0} & \theta_{0} - 1 & 1 \\
    	-4 \, {\left| A_{1} \right|}^{2} + 2 \, \alpha_{0} \overline{\alpha_{0}} & A_{1} & \theta_{0} & 1 \\
    	-4 \, {\left| A_{1} \right|}^{2} + 2 \, \alpha_{0} \overline{\alpha_{0}} & A_{2} & \theta_{0} + 1 & 1 \\
    	\mu_4& A_{0} & \theta_{0} & 2 \\
    	-2 \, \alpha_{0} & A_{0} & \theta_{0} - 1 & 0 \\
    	-2 \, \alpha_{0} & A_{1} & \theta_{0} & 0 \\
    	-2 \, \alpha_{0} & A_{2} & \theta_{0} + 1 & 0 \\
    	-2 \, \alpha_{0} & A_{3} & \theta_{0} + 2 & 0 \\
    	2 \, \alpha_{0}^{2} - 4 \, \alpha_{1} & A_{0} & \theta_{0} & 0 \\
    	2 \, \alpha_{0}^{2} - 4 \, \alpha_{1} & A_{1} & \theta_{0} + 1 & 0 \\
    	2 \, \alpha_{0}^{2} - 4 \, \alpha_{1} & A_{2} & \theta_{0} + 2 & 0 \\
    	2 \, \alpha_{2} \theta_{0} - 4 \, \alpha_{2} & A_{0} & 0 & \theta_{0} + 1 \\
    	2 \, \alpha_{2} \theta_{0} - 4 \, \alpha_{2} & A_{1} & 1 & \theta_{0} + 1 \\
    	-4 \, \alpha_{0}^{3} \theta_{0} + 2 \, \alpha_{0}^{3} + 6 \, \alpha_{0} \alpha_{1} - 6 \, \alpha_{3} & A_{0} & \theta_{0} + 1 & 0 \\
    	-4 \, \alpha_{0}^{3} \theta_{0} + 2 \, \alpha_{0}^{3} + 6 \, \alpha_{0} \alpha_{1} - 6 \, \alpha_{3} & A_{1} & \theta_{0} + 2 & 0 \\
    	2 \, \alpha_{0}^{4} + 4 \, \alpha_{0}^{2} \alpha_{1} + 4 \, \alpha_{1}^{2} + 8 \, \alpha_{0} \alpha_{3} - 4 \, {\left(\alpha_{0}^{4} + 3 \, \alpha_{0}^{2} \alpha_{1}\right)} \theta_{0} - 8 \, \alpha_{4} & A_{0} & \theta_{0} + 2 & 0 \\
    	-12 \, \alpha_{0}^{2} \theta_{0} \overline{\alpha_{0}} + 8 \, \alpha_{0} {\left| A_{1} \right|}^{2} + 8 \, \alpha_{0}^{2} \overline{\alpha_{0}} + 4 \, \alpha_{1} \overline{\alpha_{0}} - 4 \, \alpha_{5} & A_{0} & \theta_{0} & 1 \\
    	-12 \, \alpha_{0}^{2} \theta_{0} \overline{\alpha_{0}} + 8 \, \alpha_{0} {\left| A_{1} \right|}^{2} + 8 \, \alpha_{0}^{2} \overline{\alpha_{0}} + 4 \, \alpha_{1} \overline{\alpha_{0}} - 4 \, \alpha_{5} & A_{1} & \theta_{0} + 1 & 1 \\
    	\mu_5 & A_{0} & \theta_{0} + 1 & 1 \\
    	2 \, \alpha_{7} \theta_{0} - 8 \, \alpha_{7} & A_{0} & 2 & \theta_{0} \\
    	6 \, \alpha_{0} \alpha_{2} - 2 \, {\left(\alpha_{0} \alpha_{2} - \alpha_{8}\right)} \theta_{0} - 6 \, \alpha_{8} & A_{0} & 1 & \theta_{0} + 1 \\
    	-2 \, {\left(\alpha_{2} \overline{\alpha_{0}} - \alpha_{9}\right)} \theta_{0} + 4 \, \alpha_{2} \overline{\alpha_{0}} - 4 \, \alpha_{9} & A_{0} & 0 & \theta_{0} + 2 \\
    	-2 \, \theta_{0} \overline{\alpha_{2}} - 2 \, \overline{\alpha_{2}} & A_{0} & 2 \, \theta_{0} - 1 & -\theta_{0} + 2 \\
    	-2 \, \theta_{0} \overline{\alpha_{2}} - 2 \, \overline{\alpha_{2}} & A_{1} & 2 \, \theta_{0} & -\theta_{0} + 2 \\
    	-4 \, \theta_{0} \overline{\alpha_{0}}^{3} + 4 \, \overline{\alpha_{0}}^{3} & A_{0} & \theta_{0} - 2 & 3 \\
    	-4 \, \theta_{0} \overline{\alpha_{0}}^{3} + 4 \, \overline{\alpha_{0}}^{3} & A_{1} & \theta_{0} - 1 & 3 \\
    	4 \, \overline{\alpha_{0}}^{4} + 12 \, \overline{\alpha_{0}}^{2} \overline{\alpha_{1}} - 4 \, {\left(\overline{\alpha_{0}}^{4} + 3 \, \overline{\alpha_{0}}^{2} \overline{\alpha_{1}}\right)} \theta_{0} & A_{0} & \theta_{0} - 2 & 4 
    	\end{dmatrix}
    	\end{align*}
    	\begin{align*}
    	\begin{dmatrix}
    	-12 \, \alpha_{0} \theta_{0} \overline{\alpha_{0}}^{2} + 4 \, {\left| A_{1} \right|}^{2} \overline{\alpha_{0}} + 10 \, \alpha_{0} \overline{\alpha_{0}}^{2} + 2 \, \alpha_{0} \overline{\alpha_{1}} - 2 \, \overline{\alpha_{5}} & A_{0} & \theta_{0} - 1 & 2 \\
    	-12 \, \alpha_{0} \theta_{0} \overline{\alpha_{0}}^{2} + 4 \, {\left| A_{1} \right|}^{2} \overline{\alpha_{0}} + 10 \, \alpha_{0} \overline{\alpha_{0}}^{2} + 2 \, \alpha_{0} \overline{\alpha_{1}} - 2 \, \overline{\alpha_{5}} & A_{1} & \theta_{0} & 2 \\
    \mu_6& A_{0} & \theta_{0} - 1 & 3 \\
    	-2 \, \theta_{0} \overline{\alpha_{7}} & A_{0} & 2 \, \theta_{0} - 2 & -\theta_{0} + 4 \\
    	2 \, {\left(\overline{\alpha_{0}} \overline{\alpha_{2}} - \overline{\alpha_{8}}\right)} \theta_{0} + 2 \, \overline{\alpha_{0}} \overline{\alpha_{2}} - 2 \, \overline{\alpha_{8}} & A_{0} & 2 \, \theta_{0} - 1 & -\theta_{0} + 3 \\
    	2 \, {\left(\alpha_{0} \overline{\alpha_{2}} - \overline{\alpha_{9}}\right)} \theta_{0} + 4 \, \alpha_{0} \overline{\alpha_{2}} - 4 \, \overline{\alpha_{9}} & A_{0} & 2 \, \theta_{0} & -\theta_{0} + 2
    	\end{dmatrix}
    \end{align*}
    \normalsize
    where
    \begin{align*}
    	\mu_4&=8 \, {\left| A_{1} \right|}^{4} + 18 \, \alpha_{0}^{2} \overline{\alpha_{0}}^{2} - 16 \, {\left(3 \, \alpha_{0} \theta_{0} \overline{\alpha_{0}} - 2 \, \alpha_{0} \overline{\alpha_{0}}\right)} {\left| A_{1} \right|}^{2} + 8 \, \alpha_{1} \overline{\alpha_{0}}^{2} - 12 \, {\left(2 \, \alpha_{0}^{2} \overline{\alpha_{0}}^{2} + \alpha_{1} \overline{\alpha_{0}}^{2} + \alpha_{0}^{2} \overline{\alpha_{1}}\right)} \theta_{0}\\
    	& + 4 \, \alpha_{5} \overline{\alpha_{0}} + 4 \, {\left(2 \, \alpha_{0}^{2} + \alpha_{1}\right)} \overline{\alpha_{1}} + 4 \, \alpha_{0} \overline{\alpha_{5}} - 4 \, \beta \\
    	\mu_5&=10 \, \alpha_{0}^{3} \overline{\alpha_{0}} - 12 \, {\left(2 \, \alpha_{0}^{2} \theta_{0} - \alpha_{0}^{2} - \alpha_{1}\right)} {\left| A_{1} \right|}^{2} + 12 \, \alpha_{0} \alpha_{1} \overline{\alpha_{0}} + 6 \, \alpha_{0} \alpha_{5} - 8 \, {\left(2 \, \alpha_{0}^{3} \overline{\alpha_{0}} + 3 \, \alpha_{0} \alpha_{1} \overline{\alpha_{0}}\right)} \theta_{0} + 6 \, \alpha_{3} \overline{\alpha_{0}} - 6 \, \alpha_{6}\\
    	\mu_6&=	14 \, \alpha_{0} \overline{\alpha_{0}}^{3} - 4 \, {\left(6 \, \theta_{0} \overline{\alpha_{0}}^{2} - 5 \, \overline{\alpha_{0}}^{2} - \overline{\alpha_{1}}\right)} {\left| A_{1} \right|}^{2} + 20 \, \alpha_{0} \overline{\alpha_{0}} \overline{\alpha_{1}} - 8 \, {\left(2 \, \alpha_{0} \overline{\alpha_{0}}^{3} + 3 \, \alpha_{0} \overline{\alpha_{0}} \overline{\alpha_{1}}\right)} \theta_{0} + 2 \, \alpha_{0} \overline{\alpha_{3}} + 2 \, \overline{\alpha_{0}} \overline{\alpha_{5}} - 2 \, \overline{\alpha_{6}} 
    \end{align*}
    Throwing all terms of order larger or equal that $\theta_0+1$, we obtain
    \begin{align*}
    	\h_0=\begin{dmatrix}
    	2 & A_{1} & \theta_{0} - 1 & 0 \\
    	4 & A_{2} & \theta_{0} & 0 \\
    	-\frac{\theta_{0} - 2}{2 \, \theta_{0}} & C_{1} & 0 & \theta_{0} \\
    	-4 \, {\left| A_{1} \right|}^{2} + 2 \, \alpha_{0} \overline{\alpha_{0}} & A_{0} & \theta_{0} - 1 & 1 \\
    	-2 \, \alpha_{0} & A_{0} & \theta_{0} - 1 & 0 \\
    	-2 \, \alpha_{0} & A_{1} & \theta_{0} & 0 \\
    	2 \, \alpha_{0}^{2} - 4 \, \alpha_{1} & A_{0} & \theta_{0} & 0
    	\end{dmatrix}
    \end{align*}
    which agrees with the previous development \eqref{newdevh0}. 
    \section{The conservation law associated to the invariance by inversions}
    
    Finally, we have thanks of the third conservation law
    \begin{align*}
    d\,\Im\left(|\phi|^2\vec{\alpha}-2\s{{\phi}}{\vec{\alpha}}\phi-g^{-1}\otimes\left(\h_0\otimes\bar{\partial}|\phi|^2-2\s{\h_0}{\phi}\otimes\bar{\partial}\phi\right)\right)=0
    \end{align*}
    if
    \begin{align}
    \vec{\alpha}=\partial\H+|\H|^2\partial\phi+2\,g^{-1}\otimes\s{\H}{\h_0}\otimes\bar{\partial}\phi.
    \end{align}
    As the form
    \begin{align*}
    \vec{\beta}=|\phi|^2\vec{\alpha}-2\s{{\phi}}{\vec{\alpha}}\phi-g^{-1}\otimes\left(\h_0\otimes\bar{\partial}|\phi|^2-2\s{\h_0}{\phi}\otimes\bar{\partial}\phi\right)
    \end{align*}
    is a $\mathbb{C}^n$-valued $1$ differential of type $(1,0)$ (in $\mathrm{Span}_{\mathbb{C}^n}(dz)$), we write for some $\vec{F}:D^2\rightarrow\mathbb{C}^n$,
    \begin{align}\label{deff}
    \vec{\beta}=\vec{F}(z)dz,
    \end{align}
    and we have $\partial\vec{\beta}=0$. Therefore, we have
    \begin{align*}
    d\vec{\beta}=\partial\vec{\beta}+\bar{\partial}\vec{\beta}=\bar{\partial}\vec{\beta}
    \end{align*}
    and by \eqref{deff}, we deduce by linearity of imaginary part that
    \begin{align*}
    d\,\Im\left(\vec{\beta}\right)=\Im\left(\bar{\partial}\vec{\beta}\right)=\Im\left(\p{\z}\vec{F}(z)\,d\z\wedge dz\right)=-2\,\Re\left(\p{\z}\vec{F}(z)\right)dx_1\wedge dx_2=0
    \end{align*}
    as 
    \begin{align*}
    d\z\wedge dz=\left(dx_1-idx_2\right)\wedge\left(dx_1+idx_2\right)=2i\,dx_1\wedge dx_2,\quad \text{and}\;\; \Im(i\,\cdot\,)=-\Re\left(\,\cdot\,\right)
    \end{align*}
    In particular, we have $d\,\Im\left(\vec{\beta}\right)=0$ if and only if
    \begin{align}\label{pzbar}
    \Re\left(\p{\z}\vec{F}(z)\right)=0.
    \end{align}
    
    We first compute
    \small
    \begin{align*}
    	\vec{\alpha}=\begin{dmatrix}
    	-\frac{1}{2} \, \theta_{0} + 1 & C_{1} & -\theta_{0} + 1 & 0 \\
    	-\frac{1}{2} \, \theta_{0} + \frac{3}{2} & C_{2} & -\theta_{0} + 2 & 0 \\
    	-\frac{1}{2} \, \theta_{0} + 2 & C_{3} & -\theta_{0} + 3 & 0 \\
    	-\frac{1}{2} \, \theta_{0} + 1 & B_{1} & -\theta_{0} + 1 & 1 \\
    	-\frac{1}{2} \, \theta_{0} + 1 & B_{2} & -\theta_{0} + 1 & 2 \\
    	-\frac{1}{2} \, \theta_{0} + \frac{3}{2} & B_{3} & -\theta_{0} + 2 & 1 \\
    	-\theta_{0} + 2 & E_{1} & -2 \, \theta_{0} + 3 & \theta_{0} \\
    	\frac{1}{2} & \overline{B_{1}} & 0 & -\theta_{0} + 2 \\
    	1 & \overline{B_{2}} & 1 & -\theta_{0} + 2 \\
    	\frac{1}{2} & \overline{B_{3}} & 0 & -\theta_{0} + 3 \\
    	\frac{1}{2} \, \theta_{0} & \overline{E_{1}} & \theta_{0} - 1 & -2 \, \theta_{0} + 4 \\
    	\frac{1}{4} \, C_{1}^{2} & A_{0} & -\theta_{0} + 3 & 0 \\
    	\frac{1}{2} \, C_{1} \overline{C_{1}} & A_{0} & 1 & -\theta_{0} + 2 \\
    	\frac{1}{4} \, \overline{C_{1}}^{2} & A_{0} & \theta_{0} - 1 & -2 \, \theta_{0} + 4 \\
    	-2 \, A_{0} C_{1} \alpha_{0} + 2 \, A_{1} C_{1} & \overline{A_{0}} & -\theta_{0} + 2 & 0 \\
    	-2 \, A_{0} C_{1} \alpha_{0} + 2 \, A_{1} C_{1} & \overline{A_{1}} & -\theta_{0} + 2 & 1 \\
    	2 \, A_{0} C_{1} \alpha_{0}^{2} + 2 \, {\left(\alpha_{0}^{2} - 2 \, \alpha_{1}\right)} A_{0} C_{1} - 4 \, A_{1} C_{1} \alpha_{0} - 2 \, A_{0} C_{2} \alpha_{0} + 4 \, A_{2} C_{1} + 2 \, A_{1} C_{2} & \overline{A_{0}} & -\theta_{0} + 3 & 0 \\
    	2 \, A_{0} C_{1} \alpha_{0} \overline{\alpha_{0}} - 2 \, {\left(2 \, {\left| A_{1} \right|}^{2} - \alpha_{0} \overline{\alpha_{0}}\right)} A_{0} C_{1} - 2 \, A_{0} B_{1} \alpha_{0} - 2 \, A_{1} C_{1} \overline{\alpha_{0}} + 2 \, A_{1} B_{1} & \overline{A_{0}} & -\theta_{0} + 2 & 1 \\
    	-\frac{C_{1}^{2} {\left(\theta_{0} - 2\right)}}{2 \, \theta_{0}} & \overline{A_{0}} & -2 \, \theta_{0} + 3 & \theta_{0} \\
    	-2 \, A_{0} \alpha_{0} \overline{C_{1}} + 2 \, A_{1} \overline{C_{1}} & \overline{A_{0}} & 0 & -\theta_{0} + 2 \\
    	-2 \, A_{0} \alpha_{0} \overline{C_{1}} + 2 \, A_{1} \overline{C_{1}} & \overline{A_{1}} & 0 & -\theta_{0} + 3 \\
    	2 \, A_{0} \alpha_{0}^{2} \overline{C_{1}} - 2 \, A_{0} \alpha_{0} \overline{B_{1}} + 2 \, {\left(\alpha_{0}^{2} - 2 \, \alpha_{1}\right)} A_{0} \overline{C_{1}} - 4 \, A_{1} \alpha_{0} \overline{C_{1}} + 2 \, A_{1} \overline{B_{1}} + 4 \, A_{2} \overline{C_{1}} & \overline{A_{0}} & 1 & -\theta_{0} + 2 \\
    	2 \, A_{0} \alpha_{0} \overline{C_{1}} \overline{\alpha_{0}} - 2 \, {\left(2 \, {\left| A_{1} \right|}^{2} - \alpha_{0} \overline{\alpha_{0}}\right)} A_{0} \overline{C_{1}} - 2 \, A_{0} \alpha_{0} \overline{C_{2}} - 2 \, A_{1} \overline{C_{1}} \overline{\alpha_{0}} + 2 \, A_{1} \overline{C_{2}} & \overline{A_{0}} & 0 & -\theta_{0} + 3 \\
    	-\frac{C_{1} {\left(\theta_{0} - 2\right)} \overline{C_{1}}}{2 \, \theta_{0}} & \overline{A_{0}} & -\theta_{0} + 1 & 2
    	\end{dmatrix}
    \end{align*}
    \normalsize
    As $\phi$ is Willmore we have with the command
    \begin{align*} \texttt{latex(matrix\_sort(matrix\_full\_simplify(real\_part(diffzb(alpha)))))
    }
\end{align*}
the identity
    \begin{align*}
    	0=\Re\left(\p{\z}\vec{\alpha}\right)=\begin{dmatrix}
    	-\frac{1}{2} \, \theta_{0} + 1 & \overline{B_{1}} & 0 & -\theta_{0} + 1 \\
    	{\left({\left(\alpha_{0} \theta_{0} - 2 \, \alpha_{0}\right)} A_{0} - A_{1} {\left(\theta_{0} - 2\right)}\right)} \overline{C_{1}} & \overline{A_{0}} & 0 & -\theta_{0} + 1 \\
    	-\frac{1}{2} \, \theta_{0} + \frac{3}{2} & \overline{B_{3}} & 0 & -\theta_{0} + 2 \\
    	-{\left(\overline{A_{0}} \overline{\alpha_{0}} - \overline{A_{1}}\right)} \overline{C_{1}} & A_{1} & 0 & -\theta_{0} + 2 \\
    	{\left({\left(\alpha_{0} \theta_{0} - 3 \, \alpha_{0}\right)} A_{0} - A_{1} {\left(\theta_{0} - 3\right)}\right)} \overline{C_{1}} & \overline{A_{1}} & 0 & -\theta_{0} + 2 \\
    	\nu_1 & \overline{A_{0}} & 0 & -\theta_{0} + 2 \\
    	-{\left(\overline{A_{0}} \overline{\alpha_{0}} - \overline{A_{1}}\right)} \overline{B_{1}} - {\left(2 \, {\left| A_{1} \right|}^{2} \overline{A_{0}} - 2 \, \alpha_{0} \overline{A_{0}} \overline{\alpha_{0}} + \alpha_{0} \overline{A_{1}}\right)} \overline{C_{1}} & A_{0} & 0 & -\theta_{0} + 2 \\
    	-\theta_{0} + 2 & \overline{B_{2}} & 1 & -\theta_{0} + 1 \\
    	-\frac{{\left(\theta_{0}^{2} - 4\right)} C_{1} \overline{C_{1}}}{4 \, \theta_{0}} & A_{0} & 1 & -\theta_{0} + 1 \\
    	\nu_2 & \overline{A_{0}} & 1 & -\theta_{0} + 1 \\
    	-\theta_{0}^{2} + 2 \, \theta_{0} & E_{1} & -2 \, \theta_{0} + 3 & \theta_{0} - 1 \\
    	-\frac{1}{2} \, C_{1}^{2} {\left(\theta_{0} - 2\right)} & \overline{A_{0}} & -2 \, \theta_{0} + 3 & \theta_{0} - 1 \\
    	-\frac{1}{2} \, \theta_{0} + 1 & B_{1} & -\theta_{0} + 1 & 0 \\
    	{\left({\left(\theta_{0} \overline{\alpha_{0}} - 2 \, \overline{\alpha_{0}}\right)} \overline{A_{0}} - {\left(\theta_{0} - 2\right)} \overline{A_{1}}\right)} C_{1} & A_{0} & -\theta_{0} + 1 & 0 \\
    	-\theta_{0} + 2 & B_{2} & -\theta_{0} + 1 & 1 \\
    	-\frac{{\left(\theta_{0}^{2} - 4\right)} C_{1} \overline{C_{1}}}{4 \, \theta_{0}} & \overline{A_{0}} & -\theta_{0} + 1 & 1 \\
    	\nu_3 & A_{0} & -\theta_{0} + 1 & 1 \\
    	-\frac{1}{2} \, \theta_{0} + \frac{3}{2} & B_{3} & -\theta_{0} + 2 & 0 \\
    	{\left({\left(\theta_{0} \overline{\alpha_{0}} - 3 \, \overline{\alpha_{0}}\right)} \overline{A_{0}} - {\left(\theta_{0} - 3\right)} \overline{A_{1}}\right)} C_{1} & A_{1} & -\theta_{0} + 2 & 0 \\
    	-{\left(A_{0} \alpha_{0} - A_{1}\right)} C_{1} & \overline{A_{1}} & -\theta_{0} + 2 & 0 \\
    	\nu_4 & A_{0} & -\theta_{0} + 2 & 0 \\
    	-{\left(A_{0} \alpha_{0} - A_{1}\right)} B_{1} - {\left(2 \, A_{0} {\left| A_{1} \right|}^{2} - 2 \, A_{0} \alpha_{0} \overline{\alpha_{0}} + A_{1} \overline{\alpha_{0}}\right)} C_{1} & \overline{A_{0}} & -\theta_{0} + 2 & 0 \\
    	-\theta_{0}^{2} + 2 \, \theta_{0} & \overline{E_{1}} & \theta_{0} - 1 & -2 \, \theta_{0} + 3 \\
    	-\frac{1}{2} \, {\left(\theta_{0} - 2\right)} \overline{C_{1}}^{2} & A_{0} & \theta_{0} - 1 & -2 \, \theta_{0} + 3
    	\end{dmatrix}
    \end{align*}
    where
    \begin{align*}
    	\nu_1&={\left(2 \, A_{0} {\left(\theta_{0} - 3\right)} {\left| A_{1} \right|}^{2} - 2 \, {\left(\alpha_{0} \theta_{0} \overline{\alpha_{0}} - 3 \, \alpha_{0} \overline{\alpha_{0}}\right)} A_{0} + {\left(\theta_{0} \overline{\alpha_{0}} - 3 \, \overline{\alpha_{0}}\right)} A_{1}\right)} \overline{C_{1}} + {\left({\left(\alpha_{0} \theta_{0} - 3 \, \alpha_{0}\right)} A_{0} - A_{1} {\left(\theta_{0} - 3\right)}\right)} \overline{C_{2}}\\
    	\nu_2&={\left({\left(\alpha_{0} \theta_{0} - 2 \, \alpha_{0}\right)} A_{0} - A_{1} {\left(\theta_{0} - 2\right)}\right)} \overline{B_{1}} + 2 \, {\left({\left(2 \, \alpha_{0}^{2} - {\left(\alpha_{0}^{2} - \alpha_{1}\right)} \theta_{0} - 2 \, \alpha_{1}\right)} A_{0} + {\left(\alpha_{0} \theta_{0} - 2 \, \alpha_{0}\right)} A_{1} - A_{2} {\left(\theta_{0} - 2\right)}\right)} \overline{C_{1}}\\
    	\nu_3&={\left({\left(\theta_{0} \overline{\alpha_{0}} - 2 \, \overline{\alpha_{0}}\right)} \overline{A_{0}} - {\left(\theta_{0} - 2\right)} \overline{A_{1}}\right)} B_{1} - 2 \, {\left({\left({\left(\overline{\alpha_{0}}^{2} - \overline{\alpha_{1}}\right)} \theta_{0} - 2 \, \overline{\alpha_{0}}^{2} + 2 \, \overline{\alpha_{1}}\right)} \overline{A_{0}} - {\left(\theta_{0} \overline{\alpha_{0}} - 2 \, \overline{\alpha_{0}}\right)} \overline{A_{1}} + {\left(\theta_{0} - 2\right)} \overline{A_{2}}\right)} C_{1}\\
    	\nu_4&={\left(2 \, {\left(\theta_{0} - 3\right)} {\left| A_{1} \right|}^{2} \overline{A_{0}} - 2 \, {\left(\alpha_{0} \theta_{0} \overline{\alpha_{0}} - 3 \, \alpha_{0} \overline{\alpha_{0}}\right)} \overline{A_{0}} + {\left(\alpha_{0} \theta_{0} - 3 \, \alpha_{0}\right)} \overline{A_{1}}\right)} C_{1} + {\left({\left(\theta_{0} \overline{\alpha_{0}} - 3 \, \overline{\alpha_{0}}\right)} \overline{A_{0}} - {\left(\theta_{0} - 3\right)} \overline{A_{1}}\right)} C_{2}.
    \end{align*}
    Taking one order less in the expansion, we get
    \begin{align*}
    	\vec{\alpha}=\begin{dmatrix}
    	-\frac{1}{2} \, \theta_{0} + 1 & C_{1} & -\theta_{0} + 1 & 0 \\
    	-\frac{1}{2} \, \theta_{0} + \frac{3}{2} & C_{2} & -\theta_{0} + 2 & 0 \\
    	-\frac{1}{2} \, \theta_{0} + 1 & B_{1} & -\theta_{0} + 1 & 1 \\
    	\frac{1}{2} & \overline{B_{1}} & 0 & -\theta_{0} + 2 \\
    	-2 \, A_{0} C_{1} \alpha_{0} + 2 \, A_{1} C_{1} & \overline{A_{0}} & -\theta_{0} + 2 & 0 \\
    	-2 \, A_{0} \alpha_{0} \overline{C_{1}} + 2 \, A_{1} \overline{C_{1}} & \overline{A_{0}} & 0 & -\theta_{0} + 2
    	\end{dmatrix}
    	=\begin{dmatrix}
    	-\frac{1}{2} \, \theta_{0} + 1 & C_{1} & -\theta_{0} + 1 & 0 \\
    	-\frac{1}{2} \, \theta_{0} + \frac{3}{2} & C_{2} & -\theta_{0} + 2 & 0 \\
    	-\frac{1}{2} \, \theta_{0} + 1 & B_{1} & -\theta_{0} + 1 & 1 \\
    	\frac{1}{2} & \overline{B_{1}} & 0 & -\theta_{0} + 2 \\
    	 2 \, A_{1} C_{1} & \overline{A_{0}} & -\theta_{0} + 2 & 0 \\
    	 2 \, A_{1} \overline{C_{1}} & \overline{A_{0}} & 0 & -\theta_{0} + 2
    	\end{dmatrix}
    \end{align*}
    as $\s{\vec{A}_0}{\vec{C}_1}=\s{\bar{\vec{A}_0}}{\vec{C}_1}=0$.

    We will prove thanks of the invariance by inversions that
    \begin{align*}
    	\alpha_0=2\s{\bar{\vec{A}_0}}{\vec{A}_1}=0
    \end{align*}
    What follows is valid for all $\theta_0\geq 4$, as we do not use the complete developments of tensors.
    We compute
    \begin{align}\label{relinversion}
    	0=\Re\left(\p{\z}\vec{F}(z)\right)=
    	\begin{dmatrix}
    	\frac{4 \, {\left(\theta_{0}^{2} \overline{\alpha_{0}}^{3} + 2 \, \theta_{0} \overline{\alpha_{0}}^{3} - 3 \, \overline{\alpha_{0}}^{3}\right)} \overline{A_{0}}^{2}}{\theta_{0}} & A_{0} & -1 & \theta_{0} + 2 \\
    	-\frac{4 \, {\left(\theta_{0}^{2} \overline{\alpha_{0}}^{3} + 2 \, \theta_{0} \overline{\alpha_{0}}^{3} - 3 \, \overline{\alpha_{0}}^{3}\right)} A_{0} \overline{A_{0}}}{\theta_{0}} & \overline{A_{0}} & -1 & \theta_{0} + 2 \\
    	2 \, \overline{A_{0}}^{2} \overline{\alpha_{0}} - 4 \, \overline{A_{0}} \overline{A_{1}} & A_{1} & 0 & \theta_{0} \\
    	-2 \, A_{0} \alpha_{0} \overline{A_{0}} + 2 \, A_{1} \overline{A_{0}} & \overline{A_{1}} & 0 & \theta_{0} \\
    	4 \, {\left| A_{1} \right|}^{2} \overline{A_{0}}^{2} - 4 \, \alpha_{0} \overline{A_{0}}^{2} \overline{\alpha_{0}} + 4 \, \alpha_{0} \overline{A_{0}} \overline{A_{1}} & A_{0} & 0 & \theta_{0} \\
    	-4 \, A_{0} {\left| A_{1} \right|}^{2} \overline{A_{0}} + 4 \, A_{0} \alpha_{0} \overline{A_{0}} \overline{\alpha_{0}} - 2 \, A_{1} \overline{A_{0}} \overline{\alpha_{0}} - 2 \, {\left(A_{0} \alpha_{0} - A_{1}\right)} \overline{A_{1}} & \overline{A_{0}} & 0 & \theta_{0} \\
    	-2 \, \overline{A_{0}}^{2} & A_{1} & 0 & \theta_{0} - 1 \\
    	2 \, \alpha_{0} \overline{A_{0}}^{2} & A_{0} & 0 & \theta_{0} - 1 \\
    	-2 \, A_{0} \alpha_{0} \overline{A_{0}} + 2 \, A_{1} \overline{A_{0}} & \overline{A_{0}} & 0 & \theta_{0} - 1 \\
    	\lambda_1 & \overline{A_{0}} & 0 & \theta_{0} + 1 \\
    	\lambda_2 & A_{0} & 0 & \theta_{0} + 1 \\
    	-\frac{{\left(\theta_{0}^{2} + 2 \, \theta_{0} - 4\right)} \overline{A_{0}}^{2}}{4 \, \theta_{0}^{2}} & \overline{B_{1}} & 0 & \theta_{0} + 1 \\
    	-2 \, {\left(\overline{\alpha_{0}}^{2} - \overline{\alpha_{1}}\right)} \overline{A_{0}}^{2} + 4 \, \overline{A_{0}} \overline{A_{1}} \overline{\alpha_{0}} - 2 \, \overline{A_{1}}^{2} - 4 \, \overline{A_{0}} \overline{A_{2}} & A_{1} & 0 & \theta_{0} + 1 \\
    	-2 \, A_{0} \alpha_{0} \overline{A_{0}} + 2 \, A_{1} \overline{A_{0}} & \overline{A_{2}} & 0 & \theta_{0} + 1 \\
    	-4 \, A_{0} {\left| A_{1} \right|}^{2} \overline{A_{0}} + 4 \, A_{0} \alpha_{0} \overline{A_{0}} \overline{\alpha_{0}} - 2 \, A_{1} \overline{A_{0}} \overline{\alpha_{0}} - 2 \, {\left(A_{0} \alpha_{0} - A_{1}\right)} \overline{A_{1}} & \overline{A_{1}} & 0 & \theta_{0} + 1 
    	    	\end{dmatrix}
    	\end{align}
    	\begin{align*}
    	\begin{dmatrix}
    	\lambda_3 & \overline{A_{0}} & 1 & \theta_{0} \\
    	\lambda_4 & A_{0} & 1 & \theta_{0} \\
    	-\frac{A_{0} \overline{A_{0}} \overline{\alpha_{0}} - A_{0} \overline{A_{1}}}{2 \, \theta_{0}} & C_{1} & 1 & \theta_{0} \\
    	\frac{4 \, {\left(\alpha_{0}^{2} - \alpha_{1}\right)} A_{0} \theta_{0} \overline{A_{0}} - 4 \, A_{1} \alpha_{0} \theta_{0} \overline{A_{0}} + 4 \, A_{2} \theta_{0} \overline{A_{0}} - A_{0} C_{1}}{\theta_{0}} & \overline{A_{1}} & 1 & \theta_{0} \\
    	-\frac{A_{0} {\left(\theta_{0} - 2\right)} \overline{A_{0}}}{\theta_{0}^{2}} & B_{1} & 1 & \theta_{0} \\
    	4 \, \overline{A_{0}}^{2} \overline{\alpha_{0}} - 8 \, \overline{A_{0}} \overline{A_{1}} & A_{2} & 1 & \theta_{0} \\
    	8 \, {\left| A_{1} \right|}^{2} \overline{A_{0}}^{2} - 8 \, \alpha_{0} \overline{A_{0}}^{2} \overline{\alpha_{0}} + 8 \, \alpha_{0} \overline{A_{0}} \overline{A_{1}} & A_{1} & 1 & \theta_{0} \\
    	-\frac{8 \, {\left(\alpha_{0}^{2} - \alpha_{1}\right)} \theta_{0} \overline{A_{0}}^{2} - C_{1} {\left(\theta_{0} - 2\right)} \overline{A_{0}}}{2 \, \theta_{0}} & A_{0} & 1 & \theta_{0} - 1 \\
    	-4 \, \overline{A_{0}}^{2} & A_{2} & 1 & \theta_{0} - 1 \\
    	4 \, \alpha_{0} \overline{A_{0}}^{2} & A_{1} & 1 & \theta_{0} - 1 \\
    	4 \, {\left(\alpha_{0}^{2} - \alpha_{1}\right)} A_{0} \overline{A_{0}} - 4 \, A_{1} \alpha_{0} \overline{A_{0}} + 4 \, A_{2} \overline{A_{0}} & \overline{A_{0}} & 1 & \theta_{0} - 1 \\
    	\lambda_5 & A_{0} & 2 & \theta_{0} - 1 \\
    	\lambda_6 & \overline{A_{0}} & 2 & \theta_{0} - 1 \\
    	-\frac{3 \, {\left(4 \, {\left(\alpha_{0}^{2} + {\left(\alpha_{0}^{2} - \alpha_{1}\right)} \theta_{0} - \alpha_{1}\right)} \overline{A_{0}}^{2} + C_{1} \overline{A_{0}}\right)}}{2 \, {\left(\theta_{0} + 1\right)}} & A_{1} & 2 & \theta_{0} - 1 \\
    	-6 \, \overline{A_{0}}^{2} & A_{3} & 2 & \theta_{0} - 1 \\
    	6 \, \alpha_{0} \overline{A_{0}}^{2} & A_{2} & 2 & \theta_{0} - 1 \\
    	-\frac{1}{4} \, A_{0} \alpha_{0} \overline{A_{0}} + \frac{1}{4} \, A_{1} \overline{A_{0}} & C_{1} & 2 & \theta_{0} - 1 \\
    	\frac{{\left(\theta_{0} - 2\right)} \overline{A_{0}}^{2}}{2 \, \theta_{0}} & C_{1} & -\theta_{0} + 1 & 2 \, \theta_{0} - 1 \\
    	-\frac{2 \, {\left(\alpha_{2} \theta_{0}^{3} - 2 \, \alpha_{2} \theta_{0}^{2}\right)} \overline{A_{0}}^{2} - {\left({\left(\theta_{0} \overline{\alpha_{0}} - 2 \, \overline{\alpha_{0}}\right)} \overline{A_{0}}^{3} - {\left(\theta_{0} - 2\right)} \overline{A_{0}}^{2} \overline{A_{1}}\right)} C_{1}}{\theta_{0}^{2}} & A_{0} & -\theta_{0} + 1 & 2 \, \theta_{0} \\
    	-\frac{{\left(\theta_{0} \overline{\alpha_{0}} - 2 \, \overline{\alpha_{0}}\right)} \overline{A_{0}}^{2} - 2 \, {\left(\theta_{0}^{2} - \theta_{0} - 2\right)} \overline{A_{0}} \overline{A_{1}}}{2 \, {\left(\theta_{0}^{2} + \theta_{0}\right)}} & C_{1} & -\theta_{0} + 1 & 2 \, \theta_{0} \\
    	\lambda_7 & \overline{A_{0}} & -\theta_{0} + 1 & 2 \, \theta_{0} \\
    	\frac{{\left(\theta_{0}^{3} - 3 \, \theta_{0}^{2} + \theta_{0} + 2\right)} \overline{A_{0}}^{2}}{2 \, {\left(\theta_{0}^{3} + \theta_{0}^{2}\right)}} & B_{1} & -\theta_{0} + 1 & 2 \, \theta_{0} \\
    	-\frac{C_{1} {\left(\theta_{0} - 2\right)} \overline{A_{0}}}{2 \, {\left(\theta_{0}^{2} + \theta_{0}\right)}} & \overline{A_{1}} & -\theta_{0} + 1 & 2 \, \theta_{0} \\
    	\frac{2 \, {\left(A_{0} \alpha_{0} \overline{A_{0}}^{2} - A_{1} \overline{A_{0}}^{2}\right)} C_{1}}{\theta_{0}} & \overline{A_{0}} & -\theta_{0} + 2 & 2 \, \theta_{0} - 1 \\
    	\frac{{\left(\theta_{0} - 3\right)} \overline{A_{0}}^{2}}{2 \, \theta_{0}} & C_{2} & -\theta_{0} + 2 & 2 \, \theta_{0} - 1 \\
    	-2 \, A_{0} \overline{A_{0}} \overline{\alpha_{0}} + 2 \, A_{0} \overline{A_{1}} & A_{1} & \theta_{0} & 0 \\
    	2 \, A_{0}^{2} \alpha_{0} - 4 \, A_{0} A_{1} & \overline{A_{1}} & \theta_{0} & 0 
    	\end{dmatrix}
    	\end{align*}
    	\begin{align*}
    	\begin{dmatrix}
    	-4 \, A_{0} {\left| A_{1} \right|}^{2} \overline{A_{0}} + 4 \, A_{0} \alpha_{0} \overline{A_{0}} \overline{\alpha_{0}} - 2 \, A_{1} \overline{A_{0}} \overline{\alpha_{0}} - 2 \, {\left(A_{0} \alpha_{0} - A_{1}\right)} \overline{A_{1}} & A_{0} & \theta_{0} & 0 \\
    	4 \, A_{0}^{2} {\left| A_{1} \right|}^{2} - 4 \, A_{0}^{2} \alpha_{0} \overline{\alpha_{0}} + 4 \, A_{0} A_{1} \overline{\alpha_{0}} & \overline{A_{0}} & \theta_{0} & 0 \\
    	\lambda_8 & A_{0} & \theta_{0} & 1 \\
    	\lambda_9 & \overline{A_{0}} & \theta_{0} & 1 \\
    	-\frac{A_{0} \alpha_{0} \overline{A_{0}} - A_{1} \overline{A_{0}}}{2 \, \theta_{0}} & \overline{C_{1}} & \theta_{0} & 1 \\
    	\frac{4 \, {\left(\overline{\alpha_{0}}^{2} - \overline{\alpha_{1}}\right)} A_{0} \theta_{0} \overline{A_{0}} - 4 \, A_{0} \theta_{0} \overline{A_{1}} \overline{\alpha_{0}} + 4 \, A_{0} \theta_{0} \overline{A_{2}} - \overline{A_{0}} \overline{C_{1}}}{\theta_{0}} & A_{1} & \theta_{0} & 1 \\
    	-\frac{A_{0} {\left(\theta_{0} - 2\right)} \overline{A_{0}}}{\theta_{0}^{2}} & \overline{B_{1}} & \theta_{0} & 1 \\
    	4 \, A_{0}^{2} \alpha_{0} - 8 \, A_{0} A_{1} & \overline{A_{2}} & \theta_{0} & 1 \\
    	8 \, A_{0}^{2} {\left| A_{1} \right|}^{2} - 8 \, A_{0}^{2} \alpha_{0} \overline{\alpha_{0}} + 8 \, A_{0} A_{1} \overline{\alpha_{0}} & \overline{A_{1}} & \theta_{0} & 1 \\
    	-2 \, A_{0}^{2} & \overline{A_{1}} & \theta_{0} - 1 & 0 \\
    	-2 \, A_{0} \overline{A_{0}} \overline{\alpha_{0}} + 2 \, A_{0} \overline{A_{1}} & A_{0} & \theta_{0} - 1 & 0 \\
    	2 \, A_{0}^{2} \overline{\alpha_{0}} & \overline{A_{0}} & \theta_{0} - 1 & 0 \\
    	-\frac{8 \, {\left(\overline{\alpha_{0}}^{2} - \overline{\alpha_{1}}\right)} A_{0}^{2} \theta_{0} - A_{0} {\left(\theta_{0} - 2\right)} \overline{C_{1}}}{2 \, \theta_{0}} & \overline{A_{0}} & \theta_{0} - 1 & 1 \\
    	-4 \, A_{0}^{2} & \overline{A_{2}} & \theta_{0} - 1 & 1 \\
    	4 \, A_{0}^{2} \overline{\alpha_{0}} & \overline{A_{1}} & \theta_{0} - 1 & 1 \\
    	4 \, {\left(\overline{\alpha_{0}}^{2} - \overline{\alpha_{1}}\right)} A_{0} \overline{A_{0}} - 4 \, A_{0} \overline{A_{1}} \overline{\alpha_{0}} + 4 \, A_{0} \overline{A_{2}} & A_{0} & \theta_{0} - 1 & 1 \\
    	\lambda_{10} & \overline{A_{0}} & \theta_{0} - 1 & 2 \\
    	\lambda_{11} & A_{0} & \theta_{0} - 1 & 2 \\
    	-\frac{3 \, {\left(4 \, {\left({\left(\overline{\alpha_{0}}^{2} - \overline{\alpha_{1}}\right)} \theta_{0} + \overline{\alpha_{0}}^{2} - \overline{\alpha_{1}}\right)} A_{0}^{2} + A_{0} \overline{C_{1}}\right)}}{2 \, {\left(\theta_{0} + 1\right)}} & \overline{A_{1}} & \theta_{0} - 1 & 2 \\
    	-6 \, A_{0}^{2} & \overline{A_{3}} & \theta_{0} - 1 & 2 \\
    	-\frac{1}{4} \, A_{0} \overline{A_{0}} \overline{\alpha_{0}} + \frac{1}{4} \, A_{0} \overline{A_{1}} & \overline{C_{1}} & \theta_{0} - 1 & 2 \\
    	6 \, A_{0}^{2} \overline{\alpha_{0}} & \overline{A_{2}} & \theta_{0} - 1 & 2 \\
    	\lambda_{12} & A_{0} & \theta_{0} + 1 & 0 \\
    	\lambda_{13} & \overline{A_{0}} & \theta_{0} + 1 & 0 \\
    	-\frac{{\left(\theta_{0}^{2} + 2 \, \theta_{0} - 4\right)} A_{0}^{2}}{4 \, \theta_{0}^{2}} & B_{1} & \theta_{0} + 1 & 0 
    	    	    	\end{dmatrix}
    	\end{align*}
    	\begin{align}
    	\begin{dmatrix}
    	-2 \, A_{0} \overline{A_{0}} \overline{\alpha_{0}} + 2 \, A_{0} \overline{A_{1}} & A_{2} & \theta_{0} + 1 & 0 \\
    	-2 \, {\left(\alpha_{0}^{2} - \alpha_{1}\right)} A_{0}^{2} + 4 \, A_{0} A_{1} \alpha_{0} - 2 \, A_{1}^{2} - 4 \, A_{0} A_{2} & \overline{A_{1}} & \theta_{0} + 1 & 0 \\
    	-4 \, A_{0} {\left| A_{1} \right|}^{2} \overline{A_{0}} + 4 \, A_{0} \alpha_{0} \overline{A_{0}} \overline{\alpha_{0}} - 2 \, A_{1} \overline{A_{0}} \overline{\alpha_{0}} - 2 \, {\left(A_{0} \alpha_{0} - A_{1}\right)} \overline{A_{1}} & A_{1} & \theta_{0} + 1 & 0 \\
    	-\frac{4 \, {\left(\alpha_{0}^{3} \theta_{0}^{2} + 2 \, \alpha_{0}^{3} \theta_{0} - 3 \, \alpha_{0}^{3}\right)} A_{0} \overline{A_{0}}}{\theta_{0}} & A_{0} & \theta_{0} + 2 & -1 \\
    	\frac{4 \, {\left(\alpha_{0}^{3} \theta_{0}^{2} + 2 \, \alpha_{0}^{3} \theta_{0} - 3 \, \alpha_{0}^{3}\right)} A_{0}^{2}}{\theta_{0}} & \overline{A_{0}} & \theta_{0} + 2 & -1 \\
    	\frac{A_{0}^{2} {\left(\theta_{0} - 2\right)}}{2 \, \theta_{0}} & \overline{C_{1}} & 2 \, \theta_{0} - 1 & -\theta_{0} + 1 \\
    	\frac{2 \, {\left(A_{0}^{2} \overline{A_{0}} \overline{\alpha_{0}} - A_{0}^{2} \overline{A_{1}}\right)} \overline{C_{1}}}{\theta_{0}} & A_{0} & 2 \, \theta_{0} - 1 & -\theta_{0} + 2 \\
    	\frac{A_{0}^{2} {\left(\theta_{0} - 3\right)}}{2 \, \theta_{0}} & \overline{C_{2}} & 2 \, \theta_{0} - 1 & -\theta_{0} + 2 \\
    	-\frac{2 \, {\left(\theta_{0}^{3} \overline{\alpha_{2}} - 2 \, \theta_{0}^{2} \overline{\alpha_{2}}\right)} A_{0}^{2} - {\left({\left(\alpha_{0} \theta_{0} - 2 \, \alpha_{0}\right)} A_{0}^{3} - A_{0}^{2} A_{1} {\left(\theta_{0} - 2\right)}\right)} \overline{C_{1}}}{\theta_{0}^{2}} & \overline{A_{0}} & 2 \, \theta_{0} & -\theta_{0} + 1 \\
    	-\frac{{\left(\alpha_{0} \theta_{0} - 2 \, \alpha_{0}\right)} A_{0}^{2} - 2 \, {\left(\theta_{0}^{2} - \theta_{0} - 2\right)} A_{0} A_{1}}{2 \, {\left(\theta_{0}^{2} + \theta_{0}\right)}} & \overline{C_{1}} & 2 \, \theta_{0} & -\theta_{0} + 1 \\
    	\lambda_{14} & A_{0} & 2 \, \theta_{0} & -\theta_{0} + 1 \\
    	-\frac{A_{0} {\left(\theta_{0} - 2\right)} \overline{C_{1}}}{2 \, {\left(\theta_{0}^{2} + \theta_{0}\right)}} & A_{1} & 2 \, \theta_{0} & -\theta_{0} + 1 \\
    	\frac{{\left(\theta_{0}^{3} - 3 \, \theta_{0}^{2} + \theta_{0} + 2\right)} A_{0}^{2}}{2 \, {\left(\theta_{0}^{3} + \theta_{0}^{2}\right)}} & \overline{B_{1}} & 2 \, \theta_{0} & -\theta_{0} + 1
    	\end{dmatrix}
    	\end{align}
    	where
    	\small
    	\begin{align*}
    		\lambda_{1}&=\frac{1}{4 \, \theta_{0}^{2}}\bigg\{8 \, {\left(\overline{\alpha_{0}}^{2} - \overline{\alpha_{1}}\right)} A_{1} \theta_{0}^{2} \overline{A_{0}} + 16 \, {\left(2 \, A_{0} \theta_{0}^{2} \overline{A_{0}} \overline{\alpha_{0}} - A_{0} \theta_{0}^{2} \overline{A_{1}}\right)} {\left| A_{1} \right|}^{2}\\
    		& - 8 \, {\left(6 \, \alpha_{0} \theta_{0}^{3} \overline{\alpha_{0}}^{2} - 6 \, \alpha_{0} \theta_{0} \overline{\alpha_{0}}^{2} + {\left(3 \, \alpha_{0} \overline{\alpha_{0}}^{2} - 2 \, \alpha_{0} \overline{\alpha_{1}} + \overline{\alpha_{5}}\right)} \theta_{0}^{2}\right)} A_{0} \overline{A_{0}} + {\left(\theta_{0}^{2} + 2 \, \theta_{0} - 8\right)} \overline{A_{0}} \overline{B_{1}}\\
    		& + 8 \, {\left(2 \, A_{0} \alpha_{0} \theta_{0}^{2} \overline{\alpha_{0}} - A_{1} \theta_{0}^{2} \overline{\alpha_{0}}\right)} \overline{A_{1}} - 8 \, {\left(A_{0} \alpha_{0} \theta_{0}^{2} - A_{1} \theta_{0}^{2}\right)} \overline{A_{2}} + 4 \, {\left({\left(\alpha_{0} \theta_{0} + 2 \, \alpha_{0}\right)} A_{0} \overline{A_{0}}^{2} - A_{1} {\left(\theta_{0} + 2\right)} \overline{A_{0}}^{2}\right)} \overline{C_{1}}\bigg\}\\
    		\lambda_{2}&=-\frac{2}{\theta_{0}} \, \bigg\{4 \, \alpha_{0} \theta_{0} \overline{A_{0}} \overline{A_{1}} \overline{\alpha_{0}} - \alpha_{0} \theta_{0} \overline{A_{1}}^{2} - 2 \, \alpha_{0} \theta_{0} \overline{A_{0}} \overline{A_{2}} + 4 \, {\left(\theta_{0} \overline{A_{0}}^{2} \overline{\alpha_{0}} - \theta_{0} \overline{A_{0}} \overline{A_{1}}\right)} {\left| A_{1} \right|}^{2}\\
    		& - {\left(6 \, \alpha_{0} \theta_{0}^{2} \overline{\alpha_{0}}^{2} - 6 \, \alpha_{0} \overline{\alpha_{0}}^{2} + {\left(3 \, \alpha_{0} \overline{\alpha_{0}}^{2} - 2 \, \alpha_{0} \overline{\alpha_{1}} + \overline{\alpha_{5}}\right)} \theta_{0}\right)} \overline{A_{0}}^{2}\bigg\}\\
    		\lambda_{3}&=\frac{1}{2 \, \theta_{0}^{2}}\bigg\{16 \, A_{1} \alpha_{0} \theta_{0}^{2} \overline{A_{0}} \overline{\alpha_{0}} + 4 \, {\left(\alpha_{2} \theta_{0}^{3} + \alpha_{2} \theta_{0}^{2}\right)} A_{0}^{2} + A_{0} B_{1} {\left(\theta_{0} - 4\right)} + 16 \, {\left(2 \, A_{0} \alpha_{0} \theta_{0}^{2} \overline{A_{0}} - A_{1} \theta_{0}^{2} \overline{A_{0}}\right)} {\left| A_{1} \right|}^{2}\\
    		& - 8 \, {\left(3 \, \alpha_{0}^{2} \theta_{0}^{3} \overline{\alpha_{0}} + 3 \, \alpha_{0}^{2} \theta_{0} \overline{\alpha_{0}} - {\left(3 \, \alpha_{0}^{2} \overline{\alpha_{0}} + 2 \, \alpha_{1} \overline{\alpha_{0}} - \alpha_{5}\right)} \theta_{0}^{2}\right)} A_{0} \overline{A_{0}} - 8 \, {\left(\theta_{0}^{2} \overline{A_{0}} \overline{\alpha_{0}} - \theta_{0}^{2} \overline{A_{1}}\right)} A_{2}\\
    		& + 2 \, {\left(4 \, A_{0}^{2} \overline{A_{0}} \overline{\alpha_{0}} - 4 \, A_{0}^{2} \overline{A_{1}} + A_{0} \theta_{0} \overline{\alpha_{0}}\right)} C_{1} + 8 \, {\left({\left(\alpha_{0}^{2} - \alpha_{1}\right)} A_{0} \theta_{0}^{2} - A_{1} \alpha_{0} \theta_{0}^{2}\right)} \overline{A_{1}}\bigg\}\\
    		\lambda_{4}&=-\frac{1}{2 \, \theta_{0}^{2}}\bigg\{32 \, \alpha_{0} \theta_{0}^{2} {\left| A_{1} \right|}^{2} \overline{A_{0}}^{2} + 16 \, {\left(\alpha_{0}^{2} - \alpha_{1}\right)} \theta_{0}^{2} \overline{A_{0}} \overline{A_{1}} + 4 \, {\left(\alpha_{2} \theta_{0}^{3} + \alpha_{2} \theta_{0}^{2}\right)} A_{0} \overline{A_{0}} - {\left(\theta_{0}^{2} - 4\right)} B_{1} \overline{A_{0}}\\
    		& - 8 \, {\left(3 \, \alpha_{0}^{2} \theta_{0}^{3} \overline{\alpha_{0}} + 3 \, \alpha_{0}^{2} \theta_{0} \overline{\alpha_{0}} - {\left(3 \, \alpha_{0}^{2} \overline{\alpha_{0}} + 2 \, \alpha_{1} \overline{\alpha_{0}} - \alpha_{5}\right)} \theta_{0}^{2}\right)} \overline{A_{0}}^{2} + {\left(\theta_{0} \overline{A_{0}} \overline{\alpha_{0}} - {\left(\theta_{0}^{2} - \theta_{0}\right)} \overline{A_{1}}\right)} C_{1}\bigg\}\\
    		\lambda_{5}&=\frac{C_{2} {\left(\theta_{0} - 3\right)} \overline{A_{0}} + 4 \, {\left(2 \, \alpha_{0}^{3} \theta_{0}^{2} + 6 \, \alpha_{0}^{3} - {\left(5 \, \alpha_{0}^{3} + 6 \, \alpha_{0} \alpha_{1} - 3 \, \alpha_{3}\right)} \theta_{0}\right)} \overline{A_{0}}^{2} + {\left(4 \, A_{0} \alpha_{0} \overline{A_{0}}^{2} + \alpha_{0} \theta_{0} \overline{A_{0}} - 4 \, A_{1} \overline{A_{0}}^{2}\right)} C_{1}}{2 \, \theta_{0}}\\
    		\lambda_{6}&=\frac{1}{4 \, \theta_{0}}\bigg\{24 \, {\left(\alpha_{0}^{2} - \alpha_{1}\right)} A_{1} \theta_{0} \overline{A_{0}} - 24 \, A_{2} \alpha_{0} \theta_{0} \overline{A_{0}} - 8 \, {\left(2 \, \alpha_{0}^{3} \theta_{0}^{2} + 6 \, \alpha_{0}^{3} - {\left(5 \, \alpha_{0}^{3} + 6 \, \alpha_{0} \alpha_{1} - 3 \, \alpha_{3}\right)} \theta_{0}\right)} A_{0} \overline{A_{0}} + 24 \, A_{3} \theta_{0} \overline{A_{0}}\\
    		& - {\left(A_{0} \alpha_{0} \theta_{0} - A_{1} \theta_{0}\right)} C_{1}\bigg\}\\
    		\lambda_{7}&=\frac{1}{2 \, {\left(\theta_{0}^{3} + \theta_{0}^{2}\right)}}\bigg\{4 \, {\left(\alpha_{2} \theta_{0}^{4} - \alpha_{2} \theta_{0}^{3} - 2 \, \alpha_{2} \theta_{0}^{2}\right)} A_{0} \overline{A_{0}} + {\left(3 \, \theta_{0}^{2} - 4 \, \theta_{0} - 4\right)} B_{1} \overline{A_{0}}\\
    		& - {\left(4 \, {\left(\theta_{0}^{2} \overline{\alpha_{0}} - \theta_{0} \overline{\alpha_{0}} - 2 \, \overline{\alpha_{0}}\right)} A_{0} \overline{A_{0}}^{2} - 4 \, {\left(\theta_{0}^{2} - \theta_{0} - 2\right)} A_{0} \overline{A_{0}} \overline{A_{1}} - {\left(\theta_{0}^{2} \overline{\alpha_{0}} - 2 \, \theta_{0} \overline{\alpha_{0}}\right)} \overline{A_{0}}\right)} C_{1}\bigg\}\\
    		\lambda_{8}&=\frac{1}{2 \, \theta_{0}^{2}}\bigg\{8 \, {\left(\overline{\alpha_{0}}^{2} - \overline{\alpha_{1}}\right)} A_{1} \theta_{0}^{2} \overline{A_{0}} + 16 \, {\left(2 \, A_{0} \theta_{0}^{2} \overline{A_{0}} \overline{\alpha_{0}} - A_{0} \theta_{0}^{2} \overline{A_{1}}\right)} {\left| A_{1} \right|}^{2}\\
    		& - 8 \, {\left(3 \, \alpha_{0} \theta_{0}^{3} \overline{\alpha_{0}}^{2} + 3 \, \alpha_{0} \theta_{0} \overline{\alpha_{0}}^{2} - {\left(3 \, \alpha_{0} \overline{\alpha_{0}}^{2} + 2 \, \alpha_{0} \overline{\alpha_{1}} - \overline{\alpha_{5}}\right)} \theta_{0}^{2}\right)} A_{0} \overline{A_{0}} + 4 \, {\left(\theta_{0}^{3} \overline{\alpha_{2}} + \theta_{0}^{2} \overline{\alpha_{2}}\right)} \overline{A_{0}}^{2} + {\left(\theta_{0} - 4\right)} \overline{A_{0}} \overline{B_{1}}\\
    		& + 8 \, {\left(2 \, A_{0} \alpha_{0} \theta_{0}^{2} \overline{\alpha_{0}} - A_{1} \theta_{0}^{2} \overline{\alpha_{0}}\right)} \overline{A_{1}} - 8 \, {\left(A_{0} \alpha_{0} \theta_{0}^{2} - A_{1} \theta_{0}^{2}\right)} \overline{A_{2}} + 2 \, {\left(4 \, A_{0} \alpha_{0} \overline{A_{0}}^{2} + \alpha_{0} \theta_{0} \overline{A_{0}} - 4 \, A_{1} \overline{A_{0}}^{2}\right)} \overline{C_{1}}\bigg\}\\
    		\lambda_{9}&=-\frac{1}{2 \, \theta_{0}^{2}}\bigg\{32 \, A_{0}^{2} \theta_{0}^{2} {\left| A_{1} \right|}^{2} \overline{\alpha_{0}} + 16 \, {\left(\overline{\alpha_{0}}^{2} - \overline{\alpha_{1}}\right)} A_{0} A_{1} \theta_{0}^{2} - 8 \, {\left(3 \, \alpha_{0} \theta_{0}^{3} \overline{\alpha_{0}}^{2} + 3 \, \alpha_{0} \theta_{0} \overline{\alpha_{0}}^{2} - {\left(3 \, \alpha_{0} \overline{\alpha_{0}}^{2} + 2 \, \alpha_{0} \overline{\alpha_{1}} - \overline{\alpha_{5}}\right)} \theta_{0}^{2}\right)} A_{0}^{2}\\
    		& + 4 \, {\left(\theta_{0}^{3} \overline{\alpha_{2}} + \theta_{0}^{2} \overline{\alpha_{2}}\right)} A_{0} \overline{A_{0}} - {\left(\theta_{0}^{2} - 4\right)} A_{0} \overline{B_{1}} + {\left(A_{0} \alpha_{0} \theta_{0} - {\left(\theta_{0}^{2} - \theta_{0}\right)} A_{1}\right)} \overline{C_{1}}\bigg\}\\
    		\lambda_{10}&=\frac{4 \, {\left(2 \, \theta_{0}^{2} \overline{\alpha_{0}}^{3} + 6 \, \overline{\alpha_{0}}^{3} - {\left(5 \, \overline{\alpha_{0}}^{3} + 6 \, \overline{\alpha_{0}} \overline{\alpha_{1}} - 3 \, \overline{\alpha_{3}}\right)} \theta_{0}\right)} A_{0}^{2} + A_{0} {\left(\theta_{0} - 3\right)} \overline{C_{2}} + {\left(4 \, A_{0}^{2} \overline{A_{0}} \overline{\alpha_{0}} - 4 \, A_{0}^{2} \overline{A_{1}} + A_{0} \theta_{0} \overline{\alpha_{0}}\right)} \overline{C_{1}}}{2 \, \theta_{0}}\\
    		\lambda_{11}&=\frac{1}{4 \, \theta_{0}}\bigg\{24 \, {\left(\overline{\alpha_{0}}^{2} - \overline{\alpha_{1}}\right)} A_{0} \theta_{0} \overline{A_{1}} - 24 \, A_{0} \theta_{0} \overline{A_{2}} \overline{\alpha_{0}} - 8 \, {\left(2 \, \theta_{0}^{2} \overline{\alpha_{0}}^{3} + 6 \, \overline{\alpha_{0}}^{3} - {\left(5 \, \overline{\alpha_{0}}^{3} + 6 \, \overline{\alpha_{0}} \overline{\alpha_{1}} - 3 \, \overline{\alpha_{3}}\right)} \theta_{0}\right)} A_{0} \overline{A_{0}}\\
    		& + 24 \, A_{0} \theta_{0} \overline{A_{3}} - {\left(\theta_{0} \overline{A_{0}} \overline{\alpha_{0}} - \theta_{0} \overline{A_{1}}\right)} \overline{C_{1}}\bigg\}\\
    		\lambda_{12}&=\frac{1}{4 \, \theta_{0}^{2}}\bigg\{16 \, A_{1} \alpha_{0} \theta_{0}^{2} \overline{A_{0}} \overline{\alpha_{0}} + {\left(\theta_{0}^{2} + 2 \, \theta_{0} - 8\right)} A_{0} B_{1} + 16 \, {\left(2 \, A_{0} \alpha_{0} \theta_{0}^{2} \overline{A_{0}} - A_{1} \theta_{0}^{2} \overline{A_{0}}\right)} {\left| A_{1} \right|}^{2}\\
    		& - 8 \, {\left(6 \, \alpha_{0}^{2} \theta_{0}^{3} \overline{\alpha_{0}} - 6 \, \alpha_{0}^{2} \theta_{0} \overline{\alpha_{0}} + {\left(3 \, \alpha_{0}^{2} \overline{\alpha_{0}} - 2 \, \alpha_{1} \overline{\alpha_{0}} + \alpha_{5}\right)} \theta_{0}^{2}\right)} A_{0} \overline{A_{0}} - 8 \, {\left(\theta_{0}^{2} \overline{A_{0}} \overline{\alpha_{0}} - \theta_{0}^{2} \overline{A_{1}}\right)} A_{2}\\
    		& + 4 \, {\left({\left(\theta_{0} \overline{\alpha_{0}} + 2 \, \overline{\alpha_{0}}\right)} A_{0}^{2} \overline{A_{0}} - A_{0}^{2} {\left(\theta_{0} + 2\right)} \overline{A_{1}}\right)} C_{1} + 8 \, {\left({\left(\alpha_{0}^{2} - \alpha_{1}\right)} A_{0} \theta_{0}^{2} - A_{1} \alpha_{0} \theta_{0}^{2}\right)} \overline{A_{1}}\bigg\}\\
    		\lambda_{13}&=-\frac{2}{\theta_{0}} \, \bigg\{4 \, A_{0} A_{1} \alpha_{0} \theta_{0} \overline{\alpha_{0}} - A_{1}^{2} \theta_{0} \overline{\alpha_{0}} - 2 \, A_{0} A_{2} \theta_{0} \overline{\alpha_{0}} - {\left(6 \, \alpha_{0}^{2} \theta_{0}^{2} \overline{\alpha_{0}} - 6 \, \alpha_{0}^{2} \overline{\alpha_{0}} + {\left(3 \, \alpha_{0}^{2} \overline{\alpha_{0}} - 2 \, \alpha_{1} \overline{\alpha_{0}} + \alpha_{5}\right)} \theta_{0}\right)} A_{0}^{2} \\
    		&+ 4 \, {\left(A_{0}^{2} \alpha_{0} \theta_{0} - A_{0} A_{1} \theta_{0}\right)} {\left| A_{1} \right|}^{2}\bigg\}\\	
    		\lambda_{14}&=\frac{1}{2 \, {\left(\theta_{0}^{3} + \theta_{0}^{2}\right)}}\bigg\{4 \, {\left(\theta_{0}^{4} \overline{\alpha_{2}} - \theta_{0}^{3} \overline{\alpha_{2}} - 2 \, \theta_{0}^{2} \overline{\alpha_{2}}\right)} A_{0} \overline{A_{0}} + {\left(3 \, \theta_{0}^{2} - 4 \, \theta_{0} - 4\right)} A_{0} \overline{B_{1}}\\
    		& - {\left(4 \, {\left(\alpha_{0} \theta_{0}^{2} - \alpha_{0} \theta_{0} - 2 \, \alpha_{0}\right)} A_{0}^{2} \overline{A_{0}} - 4 \, {\left(\theta_{0}^{2} - \theta_{0} - 2\right)} A_{0} A_{1} \overline{A_{0}} - {\left(\alpha_{0} \theta_{0}^{2} - 2 \, \alpha_{0} \theta_{0}\right)} A_{0}\right)} \overline{C_{1}}\bigg\}
    	\end{align*}
    	\normalsize
    	
    	The only interesting coefficient is 
    	\begin{align*}
    		\frac{\z^{\theta_0+2}}{z}
    	\end{align*}
    	which is as $\s{\vec{A}_0}{\vec{A}_0}=0$ and $|\vec{A}_0|^2=\dfrac{1}{2}$
    	\begin{align}\label{id1}
    		\Omega_0&=\begin{dmatrix}
    		\frac{4 \, {\left(\theta_{0}^{2} \overline{\alpha_{0}}^{3} + 2 \, \theta_{0} \overline{\alpha_{0}}^{3} - 3 \, \overline{\alpha_{0}}^{3}\right)} \ccancel{\overline{A_{0}}^{2}}}{\theta_{0}} & A_{0} & -1 & \theta_{0} + 2 \\
    		-\frac{4 \, {\left(\theta_{0}^{2} \overline{\alpha_{0}}^{3} + 2 \, \theta_{0} \overline{\alpha_{0}}^{3} - 3 \, \overline{\alpha_{0}}^{3}\right)} A_{0} \overline{A_{0}}}{\theta_{0}} & \overline{A_{0}} & -1 & \theta_{0} + 2
    		\end{dmatrix}\nonumber\\
    		&=-\frac{2}{\theta_0}\left(\theta_0^2+2\theta_0-3\right)\bar{\alpha_0}^3\bar{\vec{A}_0}\nonumber\\
    		&=-\frac{2}{\theta_0}(\theta_0+3)(\theta_0-1)\bar{\alpha_0}^3\bar{\vec{A}_0}\nonumber\\
    		&=0
    	\end{align}
    	as $\theta_0^2+2\theta_0-1=(\theta_0+3)(\theta_0-1)$. As $\theta_0\geq 4$, and $\vec{A}_0\neq 0$ by the very definition of a branch point of multiplicity $\theta_0\geq 1$, we obtain
    	\begin{align}\label{alpha0}
    		\alpha_0=2\s{\bar{\vec{A}_0}}{\vec{A}_1}=0.
    	\end{align}
    	However, the rest of these expressions all vanish identically by previous results.
    	Therefore, we have
    	\begin{align*}
    		\left\{
    		\begin{alignedat}{1}
    		\vec{B}_2&=-\frac{(\theta_0+2)}{4\theta_0}|\vec{C}_1|^2\bar{\vec{A}_0}+\left(\ccancel{\frac{1}{2}\bar{\alpha_0}\s{\bar{\vec{A}_1}}{\vec{C}_1}}-2\s{\bar{\vec{A}_2}}{\vec{C}_1}\right)\vec{A}_0\\
    		\vec{B}_3&=-2\s{\bar{\vec{A}_1}}{\vec{C}_1}\vec{A}_1+\frac{2}{\theta_0-3}\s{\vec{A}_1}{\vec{C}_1}\bar{\vec{A}_1}+\ccancel{\frac{2\bar{\alpha_0}}{\theta_0-3}\s{\vec{A}_1}{\vec{C}_1}\bar{\vec{A}_0}}+2\left(\ccancel{\alpha_0\s{\bar{\vec{A}_1}}{\vec{C}_1}}-\s{\bar{\vec{A}_1}}{\vec{C}_2}\right)\vec{A}_0\\
    		\vec{E}_1&=-\frac{1}{2\theta_0}\s{\vec{C}_1}{\vec{C}_1}\bar{\vec{A}_0}.
    		\end{alignedat}\right.
    	\end{align*}
    	so
    	\begin{align}\label{newb2e1}
    	\left\{
    	\begin{alignedat}{1}
    	\vec{B}_1&=-2\s{\bar{\vec{A}_1}}{\vec{C}_1}\vec{A}_0\\
    	\vec{B}_2&=-\frac{(\theta_0+2)}{4\theta_0}|\vec{C}_1|^2\bar{\vec{A}_0}-2\s{\bar{\vec{A}_2}}{\vec{C}_1}\vec{A}_0\\
    	\vec{B}_3&=-2\s{\bar{\vec{A}_1}}{\vec{C}_1}\vec{A}_1+\frac{2}{\theta_0-3}\s{\vec{A}_1}{\vec{C}_1}\bar{\vec{A}_1}-2\s{\bar{\vec{A}_1}}{\vec{C}_2}\vec{A}_0\\
    	\vec{E}_1&=-\frac{1}{2\theta_0}\s{\vec{C}_1}{\vec{C}_1}\bar{\vec{A}_0}.
    	\end{alignedat}\right.
    	\end{align}

    \section{Relations given by the meromorphy}
    
    We have    
    \small
    \begin{align*}
    \mathscr{Q}_{\phi}=\begin{dmatrix}
    16 \, A_{0} A_{2} {\left| A_{1} \right|}^{2} - 8 \, A_{0} A_{1}
    \alpha_{5} - 8 \, {\left(2 \, A_{0}^{2} \alpha_{1} - A_{1}^{2}\right)}
    {\left| A_{1} \right|}^{2} & \theta_{0} - 1 & -\theta_{0} + 1 \\
    48 \, A_{1} A_{2} {\left| A_{1} \right|}^{2} + 48 \, A_{0} A_{3} {\left|
    	A_{1} \right|}^{2} - 24 \, A_{0} A_{1} \alpha_{6} - 48 \, {\left(A_{0}
    	A_{1} \alpha_{1} + A_{0}^{2} \alpha_{3}\right)} {\left| A_{1}
    	\right|}^{2} & \theta_{0} & -\theta_{0} + 1 \\
    \omega_{0} & 0 & 1 \\
    \omega_1 &
    \theta_{0} - 1 & -\theta_{0} + 2 \\
    \omega_2 & 2 \, \theta_{0} - 1 & -2 \, \theta_{0} + 2 \\
    {\left(\theta_{0}^{2} - 5 \, \theta_{0} + 6\right)} A_{1} C_{2} - 2 \,
    {\left({\left(\alpha_{1} \theta_{0}^{2} - 2 \, \alpha_{1}
    		\theta_{0}\right)} A_{0} - {\left(\theta_{0}^{2} - 2 \,
    		\theta_{0}\right)} A_{2}\right)} C_{1} & 0 & 0 \\
\omega_3 & 1 & 0 \\
    2 \, {\left(2 \, \theta_{0}^{2} - 7 \, \theta_{0} + 6\right)} A_{1}
    E_{1} & -\theta_{0} + 1 & \theta_{0} \\
    \omega_4 & -1 & 2 \\
    4 \, {\left(\theta_{0}^{3} \overline{\alpha_{7}} - 5 \, \theta_{0}^{2}
    	\overline{\alpha_{7}} + 4 \, \theta_{0} \overline{\alpha_{7}}\right)}
    A_{0} A_{1} + {\left(\theta_{0}^{2} - \theta_{0}\right)} A_{1}
    \overline{E_{1}} & 2 \, \theta_{0} - 3 & -2 \, \theta_{0} + 4 \\
    \omega_5 & -1 & 1 \\
    -8 \, {\left(\theta_{0}^{3} \overline{\alpha_{2}} - \theta_{0}
    	\overline{\alpha_{2}}\right)} A_{0}^{2} {\left| A_{1} \right|}^{2} + 4
    \, {\left(\theta_{0}^{3} \overline{\alpha_{8}} - 2 \, \theta_{0}^{2}
    	\overline{\alpha_{8}} - 3 \, \theta_{0} \overline{\alpha_{8}}\right)}
    A_{0} A_{1} & 2 \, \theta_{0} - 2 & -2 \, \theta_{0} + 3 \\
    4 \, {\left(\theta_{0}^{3} \overline{\alpha_{2}} - \theta_{0}^{2}
    	\overline{\alpha_{2}} - 2 \, \theta_{0} \overline{\alpha_{2}}\right)}
    A_{0} A_{1} & 2 \, \theta_{0} - 2 & -2 \, \theta_{0} + 2 \\
    {\left(\theta_{0}^{2} - 3 \, \theta_{0} + 2\right)} A_{1} C_{1} & -1
    & 0
    \end{dmatrix}
    \end{align*}
    \normalsize
    where
    \begin{align*}
    &	{\omega}_0=	-\frac{1}{\theta_{0}}\bigg\{2 \, {\left(\theta_{0}^{3} - 6 \, \theta_{0}^{2} + 11 \,
    	\theta_{0} - 6\right)} A_{0} C_{2} {\left| A_{1} \right|}^{2} - 8 \,
    {\left(\alpha_{1} \alpha_{2} \theta_{0}^{4} - \alpha_{1} \alpha_{2}
    	\theta_{0}^{3} - 2 \, \alpha_{1} \alpha_{2} \theta_{0}^{2}\right)}
    A_{0}^{2}\\
    & + 4 \, {\left(\alpha_{8} \theta_{0}^{4} - 4 \, \alpha_{8}
    	\theta_{0}^{3} + \alpha_{8} \theta_{0}^{2} + 6 \, \alpha_{8}
    	\theta_{0}\right)} A_{0} A_{1} + 4 \, {\left(\alpha_{2} \theta_{0}^{4} -
    	\alpha_{2} \theta_{0}^{3} - 2 \, \alpha_{2} \theta_{0}^{2}\right)}
    A_{1}^{2} + 8 \, {\left(\alpha_{2} \theta_{0}^{4} - \alpha_{2}
    	\theta_{0}^{3} - 2 \, \alpha_{2} \theta_{0}^{2}\right)} A_{0} A_{2}\\
    & -
    {\left(\theta_{0}^{3} - 5 \, \theta_{0}^{2} + 6 \, \theta_{0}\right)}
    A_{1} B_{3} + 2 \, {\left({\left(\alpha_{1} \theta_{0}^{3} - 2 \,
    		\alpha_{1} \theta_{0}^{2}\right)} A_{0} - {\left(\theta_{0}^{3} - 2 \,
    		\theta_{0}^{2}\right)} A_{2}\right)} B_{1}\\
    & + 2 \,
    {\left({\left(\theta_{0}^{3} - 2 \, \theta_{0}^{2} + 3 \, \theta_{0} -
    		6\right)} A_{1} {\left| A_{1} \right|}^{2} + {\left(\alpha_{5}
    		\theta_{0}^{3} - 3 \, \alpha_{5} \theta_{0}^{2} + 2 \, \alpha_{5}
    		\theta_{0}\right)} A_{0}\right)} C_{1}\bigg\}\\
    	\omega_{1}&=-\frac{2}{\theta_{0}}\bigg\{ \,16 \, A_{0} A_{1} \theta_{0} {\left| A_{1}\right|}^{4} + 8 \, A_{0}^{2} \alpha_{1} \theta_{0}\overline{\alpha_{5}} - 8 \, {\left(\alpha_{1} \overline{\alpha_{1}} -\beta\right)} A_{0} A_{1} \theta_{0} - 4 \, A_{1}^{2} \theta_{0}\overline{\alpha_{5}} - 8 \, A_{0} A_{2} \theta_{0}
    			\overline{\alpha_{5}}\\
    			& + {\left(2 \, \theta_{0}^{4} \overline{\alpha_{2}}- 5 \, \theta_{0}^{3} \overline{\alpha_{2}} + 5 \, \theta_{0}\overline{\alpha_{2}} - 2 \, \overline{\alpha_{2}}\right)} A_{0} C_{1} -
    			A_{1} \theta_{0} \overline{B_{2}} - {\left(A_{0} \alpha_{1} \theta_{0} -A_{2} \theta_{0}\right)} \overline{B_{1}}\bigg\}\\
    	\omega_{2}&=-8 \, {\left(\alpha_{1} \theta_{0}^{3} \overline{\alpha_{2}} - 2 \,
    		\alpha_{1} \theta_{0}^{2} \overline{\alpha_{2}} - \alpha_{1} \theta_{0}
    		\overline{\alpha_{2}} + 2 \, \alpha_{1} \overline{\alpha_{2}}\right)}
    	A_{0}^{2} + 4 \, {\left(\theta_{0}^{3} \overline{\alpha_{9}} +
    		\theta_{0}^{2} \overline{\alpha_{9}} - 4 \, \theta_{0}
    		\overline{\alpha_{9}} - 4 \, \overline{\alpha_{9}}\right)} A_{0} A_{1}\\
    	& +
    	4 \, {\left(\theta_{0}^{3} \overline{\alpha_{2}} - 2 \, \theta_{0}^{2}
    		\overline{\alpha_{2}} - \theta_{0} \overline{\alpha_{2}} + 2 \,
    		\overline{\alpha_{2}}\right)} A_{1}^{2} + 8 \, {\left(\theta_{0}^{3}
    		\overline{\alpha_{2}} - 2 \, \theta_{0}^{2} \overline{\alpha_{2}} -
    		\theta_{0} \overline{\alpha_{2}} + 2 \, \overline{\alpha_{2}}\right)}
    	A_{0} A_{2}\\
    	\omega_{3}&=    -4 \, {\left(\alpha_{7} \theta_{0}^{3} - 7 \, \alpha_{7} \theta_{0}^{2}
    		+ 12 \, \alpha_{7} \theta_{0}\right)} A_{0} A_{1} +
    	{\left(\theta_{0}^{2} - 7 \, \theta_{0} + 12\right)} A_{1} C_{3}\\
    	& -
    	{\left(3 \, {\left(\alpha_{3} \theta_{0}^{2} - \alpha_{3} \theta_{0} - 2
    			\, \alpha_{3}\right)} A_{0} + 2 \, {\left(\alpha_{1} \theta_{0}^{2} - 4
    			\, \alpha_{1}\right)} A_{1} - 3 \, {\left(\theta_{0}^{2} - \theta_{0} -
    			2\right)} A_{3}\right)} C_{1}
    		\\
    		& - 2 \, {\left({\left(\alpha_{1}
    			\theta_{0}^{2} - 4 \, \alpha_{1} \theta_{0} + 3 \, \alpha_{1}\right)}
    		A_{0} - {\left(\theta_{0}^{2} - 4 \, \theta_{0} + 3\right)}
    		A_{2}\right)} C_{2}\\
    	\omega_{4}&=\frac{1}{8 \, {\left(\theta_{0}^{2} +
    			\theta_{0}\right)}}\bigg\{64 \, {\left(\alpha_{2} \theta_{0}^{5} - 2 \, \alpha_{2}
    			\theta_{0}^{4} - \alpha_{2} \theta_{0}^{3} + 2 \, \alpha_{2}
    			\theta_{0}^{2}\right)} A_{0}^{2} {\left| A_{1} \right|}^{2} - 16 \,
    		{\left(\theta_{0}^{4} - 3 \, \theta_{0}^{3} + 2 \,
    			\theta_{0}^{2}\right)} A_{0} B_{1} {\left| A_{1} \right|}^{2}\\
    		& - 32 \,
    		{\left(\alpha_{9} \theta_{0}^{5} - 5 \, \alpha_{9} \theta_{0}^{3} + 4 \,
    			\alpha_{9} \theta_{0}\right)} A_{0} A_{1} + 8 \, {\left(\theta_{0}^{4} -
    			2 \, \theta_{0}^{3} - \theta_{0}^{2} + 2 \, \theta_{0}\right)} A_{1}
    		B_{2} \\
    		&- {\left(8 \, {\left(\theta_{0}^{4} \overline{\alpha_{5}} - 4 \,
    				\theta_{0}^{3} \overline{\alpha_{5}} + 3 \, \theta_{0}^{2}
    				\overline{\alpha_{5}} + 4 \, \theta_{0} \overline{\alpha_{5}} - 4 \,
    				\overline{\alpha_{5}}\right)} A_{0} - {\left(\theta_{0}^{4} - 4 \,
    				\theta_{0}^{3} + 3 \, \theta_{0}^{2} + 4 \, \theta_{0} - 4\right)}
    			\overline{B_{1}}\right)} C_{1}\bigg\}\\
    	\omega_{5}&=-\frac{2}{\theta_{0}} \bigg\{ \, {\left(\theta_{0}^{3} - 4 \, \theta_{0}^{2} + 5 \,
    			\theta_{0} - 2\right)} A_{0} C_{1} {\left| A_{1} \right|}^{2} + 4 \,
    		{\left(\alpha_{2} \theta_{0}^{4} - 2 \, \alpha_{2} \theta_{0}^{3} -
    			\alpha_{2} \theta_{0}^{2} + 2 \, \alpha_{2} \theta_{0}\right)} A_{0}
    		A_{1} - {\left(\theta_{0}^{3} - 3 \, \theta_{0}^{2} + 2 \,
    			\theta_{0}\right)} A_{1} B_{1}\bigg\}
    \end{align*}
    The only non-trivial component is as $\s{\vec{A}_1}{\vec{A}_1}+2\s{\vec{A}_0}{\vec{A}_2}=0$
    \begin{align*}
    	\omega_0&=-\frac{1}{\theta_{0}}\bigg\{2 \, {\left(\theta_{0}^{3} - 6 \, \theta_{0}^{2} + 11 \,
    		\theta_{0} - 6\right)} A_{0} C_{2} {\left| A_{1} \right|}^{2} - \ccancel{8 \,
    	{\left(\alpha_{1} \alpha_{2} \theta_{0}^{4} - \alpha_{1} \alpha_{2}
    		\theta_{0}^{3} - 2 \, \alpha_{1} \alpha_{2} \theta_{0}^{2}\right)}
    	A_{0}^{2}}\\
    	& + \ccancel{4 \, {\left(\alpha_{8} \theta_{0}^{4} - 4 \, \alpha_{8}
    		\theta_{0}^{3} + \alpha_{8} \theta_{0}^{2} + 6 \, \alpha_{8}
    		\theta_{0}\right)} A_{0} A_{1}} + \colorcancel{4 \, {\left(\alpha_{2} \theta_{0}^{4} -
    		\alpha_{2} \theta_{0}^{3} - 2 \, \alpha_{2} \theta_{0}^{2}\right)}
    	A_{1}^{2}}{blue} + \colorcancel{8 \, {\left(\alpha_{2} \theta_{0}^{4} - \alpha_{2}
    		\theta_{0}^{3} - 2 \, \alpha_{2} \theta_{0}^{2}\right)} A_{0} A_{2}}{blue}\\
    	& -
    	{\left(\theta_{0}^{3} - 5 \, \theta_{0}^{2} + 6 \, \theta_{0}\right)}
    	A_{1} B_{3} + 2 \, {\left(\ccancel{{\left(\alpha_{1} \theta_{0}^{3} - 2 \,
    			\alpha_{1} \theta_{0}^{2}\right)} A_{0}} - {\left(\theta_{0}^{3} - 2 \,
    			\theta_{0}^{2}\right)} A_{2}\right)} B_{1}\\
    	& + 2 \,
    	{\left({\left(\theta_{0}^{3} - 2 \, \theta_{0}^{2} + 3 \, \theta_{0} -
    			6\right)} A_{1} {\left| A_{1} \right|}^{2} + \ccancel{{\left(\alpha_{5}
    			\theta_{0}^{3} - 3 \, \alpha_{5} \theta_{0}^{2} + 2 \, \alpha_{5}
    			\theta_{0}\right)} A_{0}}\right)} C_{1}\bigg\}\\
    		&=-\frac{1}{\theta_0}\bigg\{-2\left(\theta_0^3-6\theta_0^2+11\theta_0-6\right)A_1 C_1|A_1|^2-\left(\theta_0^3-5\theta_0^2+6\theta_0\right)A_1B_3-2(\theta_0^3-2\theta_0^2)A_2B_1\\
    		&-2\left(\theta_0^3-2\theta_0^2+3\theta_0-6\right)A_1C_1|A_1|^2\bigg\}.
    \end{align*}
    Recall that
    \begin{align*}
    \left\{\begin{alignedat}{1}
        \vec{B}_1&=-2\s{\bar{\vec{A}_1}}{\vec{C}_1}\vec{A}_0\\
    	\vec{B}_2&=-\frac{(\theta_0+2)}{4\theta_0}|\vec{C}_1|^2\bar{\vec{A}_0}-2\s{\bar{\vec{A}_2}}{\vec{C}_1}\vec{A}_0\\
    	\vec{B}_3&=-2\s{\bar{\vec{A}_1}}{\vec{C}_1}\vec{A}_1+\frac{2}{\theta_0-3}\s{\vec{A}_1}{\vec{C}_1}\bar{\vec{A}_1}-2\s{\bar{\vec{A}_1}}{\vec{C}_2}\vec{A}_0\\
    	\vec{E}_1&=-\frac{1}{2\theta_0}\s{\vec{C}_1}{\vec{C}_1}\bar{\vec{A}_0}.
    \end{alignedat}\right.
\end{align*}
We first have
\begin{align*}
	\s{\vec{A}_1}{\vec{B}_3}&=-2\s{\bar{\vec{A}_1}}{\vec{C}_1}\s{\vec{A}_1}{\vec{A}_1}+\frac{2}{\theta_0-3}|\vec{A}_1|^2\s{\vec{A}_1}{\vec{C}_1}\\
	\s{\vec{A}_2}{\vec{B}_1}&=-2\s{\bar{\vec{A}_1}}{\vec{C}_1}\s{\vec{A}_0}{\vec{A}_2}=\s{\bar{\vec{A}_1}}{\vec{C}_1}\s{\vec{A}_1}{\vec{A}_1}
\end{align*}
and finally
\begin{align*}
	\omega_0&=-\frac{1}{\theta_0}\bigg\{-2\left(\theta_0^3-6\theta_0^2+11\theta_0-6\right)A_1 C_1|A_1|^2-\left(\theta_0^3-5\theta_0^2+6\theta_0\right)A_1B_3-2(\theta_0^3-2\theta_0^2)A_2B_1\\
	&-2\left(\theta_0^3-2\theta_0^2+3\theta_0-6\right)A_1C_1|A_1|^2\bigg\}\\
	&=-\frac{1}{\theta_0}\bigg\{\left(-2\left(\theta_0^3-6\theta_0^2+11\theta_0-6\right)-\left(\theta_0^3-5\theta_0^2+6\theta_0\right)\frac{2}{\theta_0-3}+2\left(\theta_0^3-2\theta_0^2+3\theta_0-6\right)\right)|\vec{A}_1|^2\s{\vec{A}_1}{\vec{C}_1}\\
	&+2\left(-(\theta_0^3-5\theta_0^2+6\theta_0)+(\theta_0^3-2\theta_0)^2\right)\s{\bar{\vec{A}_1}}{\vec{C}_1}\s{\vec{A}_1}{\vec{A}_1}\bigg\}\\
	&=-\frac{1}{\theta_0}\left\{6\theta_0(\theta_0-2)|\vec{A}_1|^2\s{\vec{A}_1}{\vec{C}_1}-6\theta(\theta_0-2)\s{\bar{\vec{A}_1}}{\vec{C}_1}\right\}\\
	&=-6(\theta_0-2)\left(|\vec{A}_1|^2\s{\vec{A}_1}{\vec{C}_1}-\s{\bar{\vec{A}_1}}{\vec{C}_1}\s{\vec{A}_1}{\vec{A}_1}\right)\\
	&=0
\end{align*}
    so
    \begin{align}\label{omega1}
    	|\vec{A}_1|^2\s{\vec{A}_1}{\vec{C}_1}=\s{\bar{\vec{A}_1}}{\vec{C}_1}\s{\vec{A}_1}{\vec{A}_1}
    \end{align}
    In particular, if $\vec{A}_1\neq 0$, we have by Cauchy-Schwarz inequality $|\s{\vec{A}_1}{\vec{A}_1}|\leq |\vec{A}_1|^2$, so we obtain
    \begin{align}\label{cauchy}
    	\left|\s{\vec{A}_1}{\vec{C}_1}\right|=\left|\s{\bar{\vec{A}_1}}{\vec{C}_1}\frac{\s{\vec{A}_1}{\vec{A}_1}}{|\vec{A}_1|^2}\right|=|\s{\bar{\vec{A}_1}}{\vec{C}_1}|\frac{|\s{\vec{A}_1}{\vec{A}_1}|}{|\vec{A}_1|^2}\leq |\s{\bar{\vec{A}_1}}{\vec{C}_1}|.
    \end{align}
    As $\vec{A}_1=0$ would finish the proof, we can safely assume that \eqref{cauchy} holds.
    
    \newpage

    \section{Another conservation law}
    
    We check that the scaling invariance is trivial, as it is algebraically trivial. Here, if
    \begin{align*}
    	2i\,\partial\vec{L}=\vec{\alpha}=\partial\H+|\H|^2\partial\phi+2\,g^{-1}\otimes\s{\H}{\h_0}+\partial\left(\vec{\gamma_0}\log|z|\right),
    \end{align*}
    it corresponds to
    \begin{align*}
    	0=\s{\p{z}\vec{L}}{{\p{\z}\phi}}=\begin{dmatrix}
    	\left(\overline{A_{0}} \overline{B_{1}}\right) \left(\frac{1}{2} \, \right) & 0 & 1 \\
    	\left(\overline{A_{0}} \overline{A_{0}}\right) \left(\left(2 \, A_{1} \overline{C_{1}}\right)\right) & 0 & 1 \\
    	\left(\overline{A_{0}} \overline{B_{3}}\right) \left(\frac{1}{2} \, \right) & 0 & 2 \\
    	\left(\overline{A_{1}} \overline{B_{1}}\right) \left(\frac{1}{2} \, \right) & 0 & 2 \\
    	\left(\overline{A_{0}} \overline{A_{1}}\right) \left(\left(2 \, A_{1} \overline{C_{1}}\right)\right) & 0 & 2 \\
    	\left(\overline{A_{0}} \overline{A_{1}}\right) \left(\left(2 \, A_{1} \overline{C_{1}}\right)\right) & 0 & 2 \\
    	\left(\overline{A_{0}} \overline{A_{0}}\right) \left({\left(-4 \, A_{0} {\left| A_{1} \right|}^{2} \overline{C_{1}} + 2 \, A_{1} \overline{C_{2}}\right)}\right) & 0 & 2 \\
    	\left(\overline{A_{0}} \overline{B_{2}}\right) \left(\right) & 1 & 1 \\
    	\left(A_{0} \overline{A_{0}}\right) \left(\left(\frac{1}{2} \, C_{1} \overline{C_{1}}\right)\right) & 1 & 1 \\
    	\left(\overline{A_{0}} \overline{A_{0}}\right) \left({\left(-4 \, A_{0} \alpha_{1} \overline{C_{1}} + 2 \, A_{1} \overline{B_{1}} + 4 \, A_{2} \overline{C_{1}}\right)}\right) & 1 & 1 \\
    	\left({\left(-\frac{1}{2} \, \theta_{0} + 1\right)} \left(\frac{1}{4 \, \theta_{0}}\right)\right) \left(C_{1} \overline{C_{1}}\right) & 1 & 1 \\
    	\left(E_{1} \overline{A_{0}}\right) \left({\left(-\theta_{0} + 2\right)}\right) & -2 \, \theta_{0} + 3 & 2 \, \theta_{0} - 1 \\
    	\left(\overline{A_{0}} \overline{A_{0}}\right) \left(\left(-\frac{C_{1}^{2} {\left(\theta_{0} - 2\right)}}{2 \, \theta_{0}}\right)\right) & -2 \, \theta_{0} + 3 & 2 \, \theta_{0} - 1 \\
    	\left(C_{1} \overline{A_{1}}\right) \left({\left(-\frac{1}{2} \, \theta_{0} + 1\right)}\right) & -\theta_{0} + 1 & \theta_{0} \\
    	\left(B_{1} \overline{A_{0}}\right) \left({\left(-\frac{1}{2} \, \theta_{0} + 1\right)}\right) & -\theta_{0} + 1 & \theta_{0} \\
    	\left(C_{1} \overline{A_{0}}\right) \left({\left(-\frac{1}{2} \, \theta_{0} + 1\right)}\right) & -\theta_{0} + 1 & \theta_{0} - 1 \\
    	\left(C_{1} \overline{A_{2}}\right) \left({\left(-\frac{1}{2} \, \theta_{0} + 1\right)}\right) & -\theta_{0} + 1 & \theta_{0} + 1 \\
    	\left(B_{1} \overline{A_{1}}\right) \left({\left(-\frac{1}{2} \, \theta_{0} + 1\right)}\right) & -\theta_{0} + 1 & \theta_{0} + 1 \\
    	\left(B_{2} \overline{A_{0}}\right) \left({\left(-\frac{1}{2} \, \theta_{0} + 1\right)}\right) & -\theta_{0} + 1 & \theta_{0} + 1 \\
    	\left(\overline{A_{0}} \overline{A_{0}}\right) \left(\left(-\frac{C_{1} {\left(\theta_{0} - 2\right)} \overline{C_{1}}}{2 \, \theta_{0}}\right)\right) & -\theta_{0} + 1 & \theta_{0} + 1 
    	\end{dmatrix}
    	\end{align*}
    	\begin{align*}
    	\begin{dmatrix}
    	\left(C_{2} \overline{A_{1}}\right) \left({\left(-\frac{1}{2} \, \theta_{0} + \frac{3}{2}\right)}\right) & -\theta_{0} + 2 & \theta_{0} \\
    	\left(B_{3} \overline{A_{0}}\right) \left({\left(-\frac{1}{2} \, \theta_{0} + \frac{3}{2}\right)}\right) & -\theta_{0} + 2 & \theta_{0} \\
    	\left(\overline{A_{0}} \overline{A_{1}}\right) \left(\left(2 \, A_{1} C_{1}\right)\right) & -\theta_{0} + 2 & \theta_{0} \\
    	\left(\overline{A_{0}} \overline{A_{1}}\right) \left(\left(2 \, A_{1} C_{1}\right)\right) & -\theta_{0} + 2 & \theta_{0} \\
    	\left(\overline{A_{0}} \overline{A_{0}}\right) \left({\left(-4 \, A_{0} C_{1} {\left| A_{1} \right|}^{2} + 2 \, A_{1} B_{1}\right)}\right) & -\theta_{0} + 2 & \theta_{0} \\
    	\left(C_{2} \overline{A_{0}}\right) \left({\left(-\frac{1}{2} \, \theta_{0} + \frac{3}{2}\right)}\right) & -\theta_{0} + 2 & \theta_{0} - 1 \\
    	\left(\overline{A_{0}} \overline{A_{0}}\right) \left(\left(2 \, A_{1} C_{1}\right)\right) & -\theta_{0} + 2 & \theta_{0} - 1 \\
    	\left(C_{1} C_{1}\right) \left(\frac{1}{8} \, {\left(-\frac{1}{2} \, \theta_{0} + 1\right)}\right) & -\theta_{0} + 3 & \theta_{0} - 1 \\
    	\left(C_{3} \overline{A_{0}}\right) \left({\left(-\frac{1}{2} \, \theta_{0} + 2\right)}\right) & -\theta_{0} + 3 & \theta_{0} - 1 \\
    	\left(A_{0} \overline{A_{0}}\right) \left(\left(\frac{1}{4} \, C_{1}^{2}\right)\right) & -\theta_{0} + 3 & \theta_{0} - 1 \\
    	\left(\overline{A_{0}} \overline{A_{0}}\right) \left({\left(-4 \, A_{0} C_{1} \alpha_{1} + 4 \, A_{2} C_{1} + 2 \, A_{1} C_{2}\right)}\right) & -\theta_{0} + 3 & \theta_{0} - 1 \\
    	\left(\overline{A_{0}} \overline{E_{1}}\right) \left(\left(\frac{1}{2} \, \theta_{0}\right)\right) & \theta_{0} - 1 & -\theta_{0} + 3 \\
    	\left(A_{0} \overline{A_{0}}\right) \left(\left(\frac{1}{4} \, \overline{C_{1}}^{2}\right)\right) & \theta_{0} - 1 & -\theta_{0} + 3
    	\end{dmatrix}
    \end{align*}
    We fist look at
    \begin{align*}
    	\Gamma_1=\begin{dmatrix}
    	\left(\overline{A_{0}} \overline{B_{3}}\right) \left(\frac{1}{2} \, \right) & 0 & 2 \\
    	\left(\overline{A_{1}} \overline{B_{1}}\right) \left(\frac{1}{2} \, \right) & 0 & 2 \\
    	\left(\overline{A_{0}} \overline{A_{1}}\right) \left(\left(2 \, A_{1} \overline{C_{1}}\right)\right) & 0 & 2 \\
    	\left(\overline{A_{0}} \overline{A_{1}}\right) \left(\left(2 \, A_{1} \overline{C_{1}}\right)\right) & 0 & 2 \\
    	\left(\overline{A_{0}} \overline{A_{0}}\right) \left({\left(-4 \, A_{0} {\left| A_{1} \right|}^{2} \overline{C_{1}} + 2 \, A_{1} \overline{C_{2}}\right)}\right) & 0 & 2
    	\end{dmatrix}=0
    \end{align*}
    As $\s{\vec{A}_0}{\vec{A}_0}=0$ and $\s{\vec{A}_0}{\vec{A}_1}=\s{\bar{\vec{A}_0}}{\vec{A}_1}=0$, and $\vec{B}_1\in\mathrm{Span}(\vec{A}_0)$, we have
    \begin{align*}
    	\Gamma_1=\frac{1}{2}\bar{\s{\vec{A}_0}{\vec{B}_3}}=0
    \end{align*}
    but
    \begin{align*}
    	\vec{B}_3&=-2\s{\bar{\vec{A}_1}}{\vec{C}_1}\vec{A}_1+\frac{2}{\theta_0-3}\s{\vec{A}_1}{\vec{C}_1}\bar{\vec{A}_1}-2\s{\bar{\vec{A}_1}}{\vec{C}_2}\vec{A}_0
    \end{align*}
    so
    \begin{align*}
    	\s{\vec{A}_0}{\vec{B}_3}=0
    \end{align*}
    Now, we have 
    \begin{align*}
    	\s{\bar{\vec{A}_0}}{\vec{B}_3}=-\s{\bar{\vec{A}_1}}{\vec{C}_2}.
    \end{align*}
    so everything is in fact trivial.

	\section{Next order developments of tensors}
	
	Up to harmonic factors, we have
	\small
	\begin{align*}
		\H=\begin{dmatrix}
		\frac{\overline{A_{1}} \overline{C_{1}}}{\theta_{0} - 3} & A_{0} & 0 & -\theta_{0} + 3 \\
		-\frac{2 \, {\left(\overline{A_{0}} \overline{\alpha_{1}} - \overline{A_{2}}\right)} \overline{C_{1}} - \overline{A_{1}} \overline{C_{2}}}{\theta_{0} - 4} & A_{0} & 0 & -\theta_{0} + 4 \\
		\frac{\overline{C_{1}}^{2}}{8 \, {\left(\theta_{0} - 4\right)}} & \overline{A_{0}} & 0 & -\theta_{0} + 4 \\
		-\frac{3 \, {\left(\overline{A_{1}} \overline{\alpha_{1}} + \overline{A_{0}} \overline{\alpha_{3}} - \overline{A_{3}}\right)} \overline{C_{1}} + 2 \, {\left(\overline{A_{0}} \overline{\alpha_{1}} - \overline{A_{2}}\right)} \overline{C_{2}} - \overline{A_{1}} \overline{C_{3}}}{\theta_{0} - 5} & A_{0} & 0 & -\theta_{0} + 5 \\
		\frac{\overline{C_{1}} \overline{C_{2}}}{4 \, {\left(\theta_{0} - 5\right)}} & \overline{A_{0}} & 0 & -\theta_{0} + 5 \\
		\frac{\overline{C_{1}}^{2}}{8 \, {\left(\theta_{0} - 5\right)}} & \overline{A_{1}} & 0 & -\theta_{0} + 5 \\
		\frac{\overline{A_{1}} \overline{C_{1}}}{8 \, {\left(\theta_{0} - 5\right)}} & \overline{C_{1}} & 0 & -\theta_{0} + 5 \\
		-A_{1} \overline{C_{1}} & \overline{A_{0}} & 1 & -\theta_{0} + 2 \\
		-\frac{2 \, {\left| A_{1} \right|}^{2} \overline{A_{0}} \overline{C_{1}} - \overline{A_{1}} \overline{B_{1}}}{\theta_{0} - 3} & A_{0} & 1 & -\theta_{0} + 3 \\
		\frac{\overline{A_{1}} \overline{C_{1}}}{\theta_{0} - 3} & A_{1} & 1 & -\theta_{0} + 3 \\
		-A_{1} \overline{C_{1}} & \overline{A_{1}} & 1 & -\theta_{0} + 3 \\
		2 \, A_{0} {\left| A_{1} \right|}^{2} \overline{C_{1}} - A_{1} \overline{C_{2}} & \overline{A_{0}} & 1 & -\theta_{0} + 3 \\
		-\frac{2 \, {\left(\overline{A_{0}} \overline{\alpha_{1}} - \overline{A_{2}}\right)} \overline{C_{1}} - \overline{A_{1}} \overline{C_{2}}}{\theta_{0} - 4} & A_{1} & 1 & -\theta_{0} + 4 \\
		-\frac{2 \, {\left| A_{1} \right|}^{2} \overline{A_{0}} \overline{C_{2}} + 2 \, {\left(\overline{A_{0}} \overline{\alpha_{1}} - \overline{A_{2}}\right)} \overline{B_{1}} - \overline{A_{1}} \overline{B_{3}} + 2 \, {\left(2 \, {\left| A_{1} \right|}^{2} \overline{A_{1}} + \overline{A_{0}} \overline{\alpha_{5}}\right)} \overline{C_{1}}}{\theta_{0} - 4} & A_{0} & 1 & -\theta_{0} + 4 \\
		\frac{16 \, A_{0} {\left(\theta_{0} - 4\right)} {\left| A_{1} \right|}^{2} \overline{C_{2}} - 8 \, A_{1} {\left(\theta_{0} - 4\right)} \overline{C_{3}} + {\left(8 \, {\left(\theta_{0} \overline{\alpha_{5}} - 4 \, \overline{\alpha_{5}}\right)} A_{0} + 8 \, {\left(\theta_{0} \overline{\alpha_{1}} - 4 \, \overline{\alpha_{1}}\right)} A_{1} - {\left(\theta_{0} - 6\right)} \overline{B_{1}}\right)} \overline{C_{1}}}{8 \, {\left(\theta_{0} - 4\right)}} & \overline{A_{0}} & 1 & -\theta_{0} + 4 \\
		-A_{1} \overline{C_{1}} & \overline{A_{2}} & 1 & -\theta_{0} + 4 \\
		2 \, A_{0} {\left| A_{1} \right|}^{2} \overline{C_{1}} - A_{1} \overline{C_{2}} & \overline{A_{1}} & 1 & -\theta_{0} + 4 \\
		-\frac{C_{1} {\left(\theta_{0} + 2\right)} \overline{C_{1}}}{8 \, \theta_{0}} & A_{0} & 2 & -\theta_{0} + 2 \\
		-\frac{1}{2} \, A_{1} \overline{B_{1}} + {\left(A_{0} \alpha_{1} - A_{2}\right)} \overline{C_{1}} & \overline{A_{0}} & 2 & -\theta_{0} + 2 \\
		-\frac{2 \, {\left| A_{1} \right|}^{2} \overline{A_{0}} \overline{C_{1}} - \overline{A_{1}} \overline{B_{1}}}{\theta_{0} - 3} & A_{1} & 2 & -\theta_{0} + 3 \\
		-\frac{16 \, \theta_{0} {\left| A_{1} \right|}^{2} \overline{A_{0}} \overline{B_{1}} - 8 \, \theta_{0} \overline{A_{1}} \overline{B_{2}} + {\left(\theta_{0}^{2} - \theta_{0} - 6\right)} C_{1} \overline{C_{2}} + {\left(8 \, \alpha_{5} \theta_{0} \overline{A_{0}} + 8 \, \alpha_{1} \theta_{0} \overline{A_{1}} + {\left(\theta_{0}^{2} - 2 \, \theta_{0} - 4\right)} B_{1}\right)} \overline{C_{1}}}{8 \, {\left(\theta_{0}^{2} - 3 \, \theta_{0}\right)}} & A_{0} & 2 & -\theta_{0} + 3 \\
		\frac{\overline{A_{1}} \overline{C_{1}}}{\theta_{0} - 3} & A_{2} & 2 & -\theta_{0} + 3 \\
		-\frac{1}{2} \, A_{1} \overline{B_{1}} + {\left(A_{0} \alpha_{1} - A_{2}\right)} \overline{C_{1}} & \overline{A_{1}} & 2 & -\theta_{0} + 3 \\
		A_{0} {\left| A_{1} \right|}^{2} \overline{B_{1}} - \frac{1}{2} \, A_{1} \overline{B_{3}} + {\left(2 \, A_{1} {\left| A_{1} \right|}^{2} + A_{0} \alpha_{5}\right)} \overline{C_{1}} + {\left(A_{0} \alpha_{1} - A_{2}\right)} \overline{C_{2}} & \overline{A_{0}} & 2 & -\theta_{0} + 3 
		\end{dmatrix}
		\end{align*}
		\footnotesize
		\begin{align*}
		\begin{dmatrix}
		-\frac{{\left(\theta_{0}^{2} + 4 \, \theta_{0} + 3\right)} C_{2} \overline{C_{1}} - {\left(12 \, {\left(\theta_{0}^{2} \overline{\alpha_{2}} + \theta_{0} \overline{\alpha_{2}}\right)} \overline{A_{0}} - {\left(\theta_{0}^{2} + 4 \, \theta_{0}\right)} \overline{B_{1}}\right)} C_{1}}{12 \, {\left(\theta_{0}^{2} + \theta_{0}\right)}} & A_{0} & 3 & -\theta_{0} + 2 \\
		-\frac{A_{1} C_{1}}{12 \, \theta_{0}} & \overline{C_{1}} & 3 & -\theta_{0} + 2 \\
		-\frac{C_{1} {\left(\theta_{0} + 3\right)} \overline{C_{1}}}{12 \, \theta_{0}} & A_{1} & 3 & -\theta_{0} + 2 \\
		-\frac{1}{24} \, A_{1} \overline{C_{1}} & C_{1} & 3 & -\theta_{0} + 2 \\
		\frac{1}{3} \, {\left(\theta_{0} \overline{\alpha_{2}} + \overline{\alpha_{2}}\right)} A_{0} C_{1} + \frac{2}{3} \, {\left(A_{0} \alpha_{1} - A_{2}\right)} \overline{B_{1}} - \frac{1}{3} \, A_{1} \overline{B_{2}} + {\left(A_{1} \alpha_{1} + A_{0} \alpha_{3} - A_{3}\right)} \overline{C_{1}} & \overline{A_{0}} & 3 & -\theta_{0} + 2 \\
		-\frac{C_{1}^{2}}{4 \, \theta_{0}} & \overline{A_{0}} & -2 \, \theta_{0} + 4 & \theta_{0} \\
		\frac{{\left(\alpha_{2} \theta_{0} + \alpha_{2}\right)} C_{1} \overline{A_{0}} - E_{1} \overline{A_{1}}}{\theta_{0} + 1} & A_{0} & -2 \, \theta_{0} + 4 & \theta_{0} + 1 \\
		-\frac{C_{1} \overline{A_{1}}}{4 \, {\left(\theta_{0}^{2} + \theta_{0}\right)}} & C_{1} & -2 \, \theta_{0} + 4 & \theta_{0} + 1 \\
		\frac{{\left(4 \, {\left(\alpha_{2} \theta_{0}^{2} + \alpha_{2} \theta_{0}\right)} A_{0} - B_{1} {\left(4 \, \theta_{0} + 1\right)}\right)} C_{1}}{8 \, {\left(\theta_{0}^{2} + \theta_{0}\right)}} & \overline{A_{0}} & -2 \, \theta_{0} + 4 & \theta_{0} + 1 \\
		-\frac{C_{1}^{2} {\left(2 \, \theta_{0} + 1\right)}}{8 \, {\left(\theta_{0}^{2} + \theta_{0}\right)}} & \overline{A_{1}} & -2 \, \theta_{0} + 4 & \theta_{0} + 1 \\
		-\frac{C_{1} C_{2} {\left(2 \, \theta_{0} - 5\right)} - 2 \, A_{1} E_{1} \theta_{0}}{2 \, {\left(2 \, \theta_{0}^{2} - 5 \, \theta_{0}\right)}} & \overline{A_{0}} & -2 \, \theta_{0} + 5 & \theta_{0} \\
		-C_{1} \overline{A_{1}} & A_{0} & -\theta_{0} + 2 & 1 \\
		-\frac{C_{1} {\left(\theta_{0} + 2\right)} \overline{C_{1}}}{8 \, \theta_{0}} & \overline{A_{0}} & -\theta_{0} + 2 & 2 \\
		{\left(\overline{A_{0}} \overline{\alpha_{1}} - \overline{A_{2}}\right)} C_{1} - \frac{1}{2} \, B_{1} \overline{A_{1}} & A_{0} & -\theta_{0} + 2 & 2 \\
		-\frac{{\left(\theta_{0}^{2} + 4 \, \theta_{0} + 3\right)} C_{1} \overline{C_{2}} - {\left(12 \, {\left(\alpha_{2} \theta_{0}^{2} + \alpha_{2} \theta_{0}\right)} A_{0} - {\left(\theta_{0}^{2} + 4 \, \theta_{0}\right)} B_{1}\right)} \overline{C_{1}}}{12 \, {\left(\theta_{0}^{2} + \theta_{0}\right)}} & \overline{A_{0}} & -\theta_{0} + 2 & 3 \\
		-\frac{\overline{A_{1}} \overline{C_{1}}}{12 \, \theta_{0}} & C_{1} & -\theta_{0} + 2 & 3 \\
		-\frac{C_{1} {\left(\theta_{0} + 3\right)} \overline{C_{1}}}{12 \, \theta_{0}} & \overline{A_{1}} & -\theta_{0} + 2 & 3 \\
		-\frac{1}{24} \, C_{1} \overline{A_{1}} & \overline{C_{1}} & -\theta_{0} + 2 & 3 \\
		\frac{1}{3} \, {\left(\alpha_{2} \theta_{0} + \alpha_{2}\right)} \overline{A_{0}} \overline{C_{1}} + \frac{2}{3} \, {\left(\overline{A_{0}} \overline{\alpha_{1}} - \overline{A_{2}}\right)} B_{1} + {\left(\overline{A_{1}} \overline{\alpha_{1}} + \overline{A_{0}} \overline{\alpha_{3}} - \overline{A_{3}}\right)} C_{1} - \frac{1}{3} \, B_{2} \overline{A_{1}} & A_{0} & -\theta_{0} + 2 & 3 \\
		\frac{A_{1} C_{1}}{\theta_{0} - 3} & \overline{A_{0}} & -\theta_{0} + 3 & 0 \\
		-\frac{2 \, A_{0} C_{1} {\left| A_{1} \right|}^{2} - A_{1} B_{1}}{\theta_{0} - 3} & \overline{A_{0}} & -\theta_{0} + 3 & 1 \\
		\frac{A_{1} C_{1}}{\theta_{0} - 3} & \overline{A_{1}} & -\theta_{0} + 3 & 1 \\
		-C_{1} \overline{A_{1}} & A_{1} & -\theta_{0} + 3 & 1 \\
		2 \, C_{1} {\left| A_{1} \right|}^{2} \overline{A_{0}} - C_{2} \overline{A_{1}} & A_{0} & -\theta_{0} + 3 & 1 \\
		-\frac{2 \, A_{0} C_{1} {\left| A_{1} \right|}^{2} - A_{1} B_{1}}{\theta_{0} - 3} & \overline{A_{1}} & -\theta_{0} + 3 & 2 \\
		-\frac{16 \, A_{0} B_{1} \theta_{0} {\left| A_{1} \right|}^{2} - 8 \, A_{1} B_{2} \theta_{0} + {\left(\theta_{0}^{2} - \theta_{0} - 6\right)} C_{2} \overline{C_{1}} + {\left(8 \, A_{1} \theta_{0} \overline{\alpha_{1}} + 8 \, A_{0} \theta_{0} \overline{\alpha_{5}} + {\left(\theta_{0}^{2} - 2 \, \theta_{0} - 4\right)} \overline{B_{1}}\right)} C_{1}}{8 \, {\left(\theta_{0}^{2} - 3 \, \theta_{0}\right)}} & \overline{A_{0}} & -\theta_{0} + 3 & 2 \\
		\frac{A_{1} C_{1}}{\theta_{0} - 3} & \overline{A_{2}} & -\theta_{0} + 3 & 2 \\
		{\left(\overline{A_{0}} \overline{\alpha_{1}} - \overline{A_{2}}\right)} C_{1} - \frac{1}{2} \, B_{1} \overline{A_{1}} & A_{1} & -\theta_{0} + 3 & 2 \\
		B_{1} {\left| A_{1} \right|}^{2} \overline{A_{0}} + {\left(2 \, {\left| A_{1} \right|}^{2} \overline{A_{1}} + \overline{A_{0}} \overline{\alpha_{5}}\right)} C_{1} + {\left(\overline{A_{0}} \overline{\alpha_{1}} - \overline{A_{2}}\right)} C_{2} - \frac{1}{2} \, B_{3} \overline{A_{1}} & A_{0} & -\theta_{0} + 3 & 2 
		\end{dmatrix}
		\end{align*}
		\small
		\begin{align*}
		\begin{dmatrix}
		-\frac{2 \, {\left(A_{0} \alpha_{1} - A_{2}\right)} C_{1} - A_{1} C_{2}}{\theta_{0} - 4} & \overline{A_{0}} & -\theta_{0} + 4 & 0 \\
		\frac{C_{1}^{2}}{8 \, {\left(\theta_{0} - 4\right)}} & A_{0} & -\theta_{0} + 4 & 0 \\
		-\frac{2 \, {\left(A_{0} \alpha_{1} - A_{2}\right)} C_{1} - A_{1} C_{2}}{\theta_{0} - 4} & \overline{A_{1}} & -\theta_{0} + 4 & 1 \\
		-\frac{2 \, A_{0} C_{2} {\left| A_{1} \right|}^{2} + 2 \, {\left(A_{0} \alpha_{1} - A_{2}\right)} B_{1} - A_{1} B_{3} + 2 \, {\left(2 \, A_{1} {\left| A_{1} \right|}^{2} + A_{0} \alpha_{5}\right)} C_{1}}{\theta_{0} - 4} & \overline{A_{0}} & -\theta_{0} + 4 & 1 \\
		\frac{16 \, C_{2} {\left(\theta_{0} - 4\right)} {\left| A_{1} \right|}^{2} \overline{A_{0}} - 8 \, C_{3} {\left(\theta_{0} - 4\right)} \overline{A_{1}} - {\left(B_{1} {\left(\theta_{0} - 6\right)} - 8 \, {\left(\alpha_{5} \theta_{0} - 4 \, \alpha_{5}\right)} \overline{A_{0}} - 8 \, {\left(\alpha_{1} \theta_{0} - 4 \, \alpha_{1}\right)} \overline{A_{1}}\right)} C_{1}}{8 \, {\left(\theta_{0} - 4\right)}} & A_{0} & -\theta_{0} + 4 & 1 \\
		-C_{1} \overline{A_{1}} & A_{2} & -\theta_{0} + 4 & 1 \\
		2 \, C_{1} {\left| A_{1} \right|}^{2} \overline{A_{0}} - C_{2} \overline{A_{1}} & A_{1} & -\theta_{0} + 4 & 1 \\
		-\frac{3 \, {\left(A_{1} \alpha_{1} + A_{0} \alpha_{3} - A_{3}\right)} C_{1} + 2 \, {\left(A_{0} \alpha_{1} - A_{2}\right)} C_{2} - A_{1} C_{3}}{\theta_{0} - 5} & \overline{A_{0}} & -\theta_{0} + 5 & 0 \\
		\frac{C_{1} C_{2}}{4 \, {\left(\theta_{0} - 5\right)}} & A_{0} & -\theta_{0} + 5 & 0 \\
		\frac{C_{1}^{2}}{8 \, {\left(\theta_{0} - 5\right)}} & A_{1} & -\theta_{0} + 5 & 0 \\
		\frac{A_{1} C_{1}}{8 \, {\left(\theta_{0} - 5\right)}} & C_{1} & -\theta_{0} + 5 & 0 \\
		-\frac{\overline{C_{1}}^{2}}{4 \, \theta_{0}} & A_{0} & \theta_{0} & -2 \, \theta_{0} + 4 \\
		-\frac{{\left(2 \, \theta_{0} - 5\right)} \overline{C_{1}} \overline{C_{2}} - 2 \, \theta_{0} \overline{A_{1}} \overline{E_{1}}}{2 \, {\left(2 \, \theta_{0}^{2} - 5 \, \theta_{0}\right)}} & A_{0} & \theta_{0} & -2 \, \theta_{0} + 5 \\
		\frac{{\left(\theta_{0} \overline{\alpha_{2}} + \overline{\alpha_{2}}\right)} A_{0} \overline{C_{1}} - A_{1} \overline{E_{1}}}{\theta_{0} + 1} & \overline{A_{0}} & \theta_{0} + 1 & -2 \, \theta_{0} + 4 \\
		-\frac{A_{1} \overline{C_{1}}}{4 \, {\left(\theta_{0}^{2} + \theta_{0}\right)}} & \overline{C_{1}} & \theta_{0} + 1 & -2 \, \theta_{0} + 4 \\
		\frac{{\left(4 \, {\left(\theta_{0}^{2} \overline{\alpha_{2}} + \theta_{0} \overline{\alpha_{2}}\right)} \overline{A_{0}} - {\left(4 \, \theta_{0} + 1\right)} \overline{B_{1}}\right)} \overline{C_{1}}}{8 \, {\left(\theta_{0}^{2} + \theta_{0}\right)}} & A_{0} & \theta_{0} + 1 & -2 \, \theta_{0} + 4 \\
		-\frac{{\left(2 \, \theta_{0} + 1\right)} \overline{C_{1}}^{2}}{8 \, {\left(\theta_{0}^{2} + \theta_{0}\right)}} & A_{1} & \theta_{0} + 1 & -2 \, \theta_{0} + 4
		\end{dmatrix}
	\end{align*}
	\normalsize
	which we translate for some $\vec{C}_4,\vec{B}_4,\vec{B}_5,\vec{B}_6,\vec{E}_2,\vec{E}_3\in\mathbb{C}^n$ as
	\begin{align*}
		\H=\begin{dmatrix}
		\frac{1}{2} & C_{1} & -\theta_{0} + 2 & 0 \\
		\frac{1}{2} & C_{2} & -\theta_{0} + 3 & 0 \\
		\frac{1}{2} & C_{3} & -\theta_{0} + 4 & 0 \\
		\frac{1}{2} & C_{4} & -\theta_{0} + 5 & 0 \\
		\frac{1}{2} & B_{1} & -\theta_{0} + 2 & 1 \\
		\frac{1}{2} & B_{2} & -\theta_{0} + 2 & 2 \\
		\frac{1}{2} & B_{3} & -\theta_{0} + 3 & 1 \\
		\frac{1}{2} & B_{4} & -\theta_{0} + 2 & 3 \\
		\frac{1}{2} & B_{5} & -\theta_{0} + 3 & 2 \\
		\frac{1}{2} & B_{6} & -\theta_{0} + 4 & 1 \\
		\frac{1}{2} & E_{1} & -2 \, \theta_{0} + 4 & \theta_{0} \\
		\frac{1}{2} & E_{2} & -2 \, \theta_{0} + 4 & \theta_{0} + 1 \\
		\frac{1}{2} & E_{3} & -2 \, \theta_{0} + 5 & \theta_{0} 
		\end{dmatrix}
		\begin{dmatrix}
		\frac{1}{2} & \overline{C_{1}} & 0 & -\theta_{0} + 2 \\
		\frac{1}{2} & \overline{C_{2}} & 0 & -\theta_{0} + 3 \\
		\frac{1}{2} & \overline{C_{3}} & 0 & -\theta_{0} + 4 \\
		\frac{1}{2} & \overline{C_{4}} & 0 & -\theta_{0} + 5 \\
		\frac{1}{2} & \overline{B_{1}} & 1 & -\theta_{0} + 2 \\
		\frac{1}{2} & \overline{B_{2}} & 2 & -\theta_{0} + 2 \\
		\frac{1}{2} & \overline{B_{3}} & 1 & -\theta_{0} + 3 \\
		\frac{1}{2} & \overline{B_{4}} & 3 & -\theta_{0} + 2 \\
		\frac{1}{2} & \overline{B_{5}} & 2 & -\theta_{0} + 3 \\
		\frac{1}{2} & \overline{B_{6}} & 1 & -\theta_{0} + 4 \\
		\frac{1}{2} & \overline{E_{1}} & \theta_{0} & -2 \, \theta_{0} + 4 \\
		\frac{1}{2} & \overline{E_{2}} & \theta_{0} + 1 & -2 \, \theta_{0} + 4 \\
		\frac{1}{2} & \overline{E_{3}} & \theta_{0} & -2 \, \theta_{0} + 5
		\end{dmatrix}
	\end{align*}
	We have as
	\begin{align*}
		\alpha_2=\frac{1}{2\theta_0(\theta_0+1)}\s{\bar{\vec{A}_1}}{\vec{C}_1},\quad \alpha_7=\frac{1}{8\theta_0(\theta_0-4)}\s{\vec{C}_1}{\vec{C}_1}
	\end{align*}
	word
	\begin{align*}
		\frac{1}{2}\vec{E}_2&=\begin{dmatrix}
		\frac{\ccancel{{\left(\alpha_{2} \theta_{0} + \alpha_{2}\right)} C_{1} \overline{A_{0}}} - \ccancel{E_{1} \overline{A_{1}}}}{\theta_{0} + 1} & A_{0} & -2 \, \theta_{0} + 4 & \theta_{0} + 1 \\
		-\frac{C_{1} \overline{A_{1}}}{4 \, {\left(\theta_{0}^{2} + \theta_{0}\right)}} & C_{1} & -2 \, \theta_{0} + 4 & \theta_{0} + 1 \\
		\frac{{\left(\ccancel{4 \, {\left(\alpha_{2} \theta_{0}^{2} + \alpha_{2} \theta_{0}\right)}} A_{0} - \ccancel{B_{1} {\left(4 \, \theta_{0} + 1\right)}}\right)} C_{1}}{8 \, {\left(\theta_{0}^{2} + \theta_{0}\right)}} & \overline{A_{0}} & -2 \, \theta_{0} + 4 & \theta_{0} + 1 \\
		-\frac{C_{1}^{2} {\left(2 \, \theta_{0} + 1\right)}}{8 \, {\left(\theta_{0}^{2} + \theta_{0}\right)}} & \overline{A_{1}} & -2 \, \theta_{0} + 4 & \theta_{0} + 1 
		\end{dmatrix}\\
		&=-\frac{\alpha_2}{2}\vec{C}_1-\frac{(2\theta_0+1)(\theta_0-4)}{\theta_0+1}\alpha_7\bar{\vec{A}_1}
	\end{align*}
	so
	\begin{align}
		\vec{E}_2=-\alpha_2\vec{C}_1-\frac{2(2\theta_0+1)(\theta_0-4)}{\theta_0+1}\alpha_7\bar{\vec{A}_1}
	\end{align}
	Finally,
	\begin{align*}
		\frac{1}{2}\vec{E}_3&=
		\begin{dmatrix}-\frac{C_{1} C_{2} {\left(2 \, \theta_{0} - 5\right)} - \ccancel{2 \, A_{1} E_{1} \theta_{0}}}{2 \, {\left(2 \, \theta_{0}^{2} - 5 \, \theta_{0}\right)}} & \overline{A_{0}} & -2 \, \theta_{0} + 5 & \theta_{0} 
		\end{dmatrix}\\
		&=-\frac{1}{2\theta_0}\s{\vec{C}_1}{\vec{C}_2}\bar{\vec{A}_0}
	\end{align*}
	So
	\begin{align*}
		\left\{
		\begin{alignedat}{1}
		\vec{E}_2&=-\alpha_2\,\vec{C}_1-\frac{2(2\theta_0+1)(\theta_0-4)}{\theta_0+1}\alpha_7\bar{\vec{A}_1}\\
		\vec{E}_3&=-\frac{1}{\theta_0}\s{\vec{C}_1}{\vec{C}_2}\bar{\vec{A}_0}
		\end{alignedat}\right.
	\end{align*}
	Now, we have as $\vec{B}_1\in\mathrm{Span}(\vec{A}_0)$, and $\vec{B}_2\in\mathrm{Span}(\vec{A}_0,\bar{\vec{A}_0})$
	\begin{align*}
		\frac{1}{2}\vec{B}_4&=\begin{dmatrix}
		-\frac{{\left(\theta_{0}^{2} + 4 \, \theta_{0} + 3\right)} C_{1} \overline{C_{2}} - {\left(\ccancel{12 \, {\left(\alpha_{2} \theta_{0}^{2} + \alpha_{2} \theta_{0}\right)} A_{0}} - \ccancel{{\left(\theta_{0}^{2} + 4 \, \theta_{0}\right)} B_{1}}\right)} \overline{C_{1}}}{12 \, {\left(\theta_{0}^{2} + \theta_{0}\right)}} & \overline{A_{0}} & -\theta_{0} + 2 & 3 \\
		-\frac{\overline{A_{1}} \overline{C_{1}}}{12 \, \theta_{0}} & C_{1} & -\theta_{0} + 2 & 3 \\
		-\frac{C_{1} {\left(\theta_{0} + 3\right)} \overline{C_{1}}}{12 \, \theta_{0}} & \overline{A_{1}} & -\theta_{0} + 2 & 3 \\
		-\frac{1}{24} \, C_{1} \overline{A_{1}} & \overline{C_{1}} & -\theta_{0} + 2 & 3 \\
		\frac{1}{3} \, \ccancel{{\left(\alpha_{2} \theta_{0} + \alpha_{2}\right)} \overline{A_{0}} \overline{C_{1}}} + \frac{2}{3} \, {\left(\overline{A_{0}} \overline{\alpha_{1}} - \overline{A_{2}}\right)} B_{1} + {\left(\overline{A_{1}} \overline{\alpha_{1}} + \ccancel{\overline{A_{0}} \overline{\alpha_{3}}} - \overline{A_{3}}\right)} C_{1} - \frac{1}{3}\ccancel{ \, B_{2} \overline{A_{1}}} & A_{0} & -\theta_{0} + 2 & 3 
		\end{dmatrix}\\
		&=\begin{dmatrix}
		-\frac{{\left(\theta_{0}+3\right)} C_{1} \overline{C_{2}}}{12 \, {\theta_{0}}} & \overline{A_{0}} & -\theta_{0} + 2 & 3 \\
		-\frac{\bar{\zeta_2}}{12 \, \theta_{0}} & C_{1} & -\theta_{0} + 2 & 3 \\
		- \frac{{\left(\theta_{0} + 3\right)} |{C_{1}}|^2}{12 \, \theta_{0}} & \overline{A_{1}} & -\theta_{0} + 2 & 3 \\
		-\frac{\theta_0(\theta_0+1)}{12}\alpha_2 & \overline{C_{1}} & -\theta_{0} + 2 & 3 \\
		\frac{2}{3} \, {\left(\overline{A_{0}} \overline{\alpha_{1}} - \overline{A_{2}}\right)} B_{1} + {\left(\overline{A_{1}} \overline{\alpha_{1}} - \overline{A_{3}}\right)} C_{1}   & A_{0} & -\theta_{0} + 2 & 3 
		\end{dmatrix}
	\end{align*}
	as
	\begin{align*}
	\left\{
	\begin{alignedat}{1}
	    \alpha_2&=\frac{1}{2\theta_0(\theta_0+1)}\s{\bar{\vec{A}_1}}{\vec{C}_1}\\
		\zeta_2&=\s{\vec{A}_1}{\vec{C}_1}.
		\end{alignedat}\right.
	\end{align*}
	As $\vec{B}_1=-2\s{\bar{\vec{A}_1}}{\vec{C}_1}\vec{A}_0$, and
	\begin{align*}
		&\s{\bar{\alpha_1}\bar{\vec{A}_0}-\bar{\vec{A}_2}}{\vec{A}_0}=\frac{1}{2}\left(\bar{\alpha_1}-\s{\vec{A}_0}{\bar{\vec{A}_2}}\right)=0\\
		&\s{\bar{\alpha_1}\bar{\vec{A}_0}-\bar{\vec{A}_2}}{\vec{B}_1}=0
	\end{align*}
	so finally
	\begin{align*}
		\vec{B}_4=-\frac{(\theta_0+3)}{6\theta_0}\s{\vec{C}_1}{\bar{\vec{C}_2}}\bar{\vec{A}_0}-\frac{\bar{\zeta_2}}{6\theta_0}\vec{C}_1-\frac{(\theta_0+3)}{6\theta_0}|\vec{C}_1|^2\bar{\vec{A}_1}-\frac{\theta_0(\theta_0+1)}{6}\alpha_2\bar{\vec{C}_1}+2\left(\bar{\alpha_1}\s{\bar{\vec{A}_1}}{\vec{C}_1}-\s{\bar{\vec{A}_3}}{\vec{C}_1}\right)\vec{A}_0.
	\end{align*}
	Then, we have
	\small
	\begin{align*}
		&\frac{1}{2}\vec{B}_5=\\
		&\begin{dmatrix}
				-\frac{\ccancel{2 \, A_{0} C_{1} {\left| A_{1} \right|}^{2}} - \ccancel{A_{1} B_{1}}}{\theta_{0} - 3} & \overline{A_{1}} & -\theta_{0} + 3 & 2 \\
		-\frac{\ccancel{16 \, A_{0} B_{1} \theta_{0} {\left| A_{1} \right|}^{2}} - \ccancel{8 \, A_{1} B_{2} \theta_{0}} + {\left(\theta_{0}^{2} - \theta_{0} - 6\right)} C_{2} \overline{C_{1}} + {\left(8 \, A_{1} \theta_{0} \overline{\alpha_{1}} + \ccancel{8 \, A_{0} \theta_{0} \overline{\alpha_{5}}} + \ccancel{{\left(\theta_{0}^{2} - 2 \, \theta_{0} - 4\right)} \overline{B_{1}}}\right)} C_{1}}{8 \, {\left(\theta_{0}^{2} - 3 \, \theta_{0}\right)}} & \overline{A_{0}} & -\theta_{0} + 3 & 2 \\
		\frac{A_{1} C_{1}}{\theta_{0} - 3} & \overline{A_{2}} & -\theta_{0} + 3 & 2 \\
		{\left(\ccancel{\overline{A_{0}} \overline{\alpha_{1}}} - \overline{A_{2}}\right)} C_{1} - \ccancel{\frac{1}{2} \, B_{1} \overline{A_{1}}} & A_{1} & -\theta_{0} + 3 & 2 \\
		B_{1} {\left| A_{1} \right|}^{2} \overline{A_{0}} + {\left(2 \, {\left| A_{1} \right|}^{2} \overline{A_{1}} + \ccancel{\overline{A_{0}} \overline{\alpha_{5}}}\right)} C_{1} + {\left(\ccancel{\overline{A_{0}} \overline{\alpha_{1}}} - \overline{A_{2}}\right)} C_{2} - \frac{1}{2} \, B_{3} \overline{A_{1}} & A_{0} & -\theta_{0} + 3 & 2 
		\end{dmatrix}\\
		&=\begin{dmatrix}
		-\frac{{\left(\theta_{0}^{2} - \theta_{0} - 6\right)} C_{2} \overline{C_{1}} + 8 \,  \theta_{0} \overline{\alpha_{1}} \zeta_2}{8 \, {\left(\theta_{0}^{2} - 3 \, \theta_{0}\right)}} & \overline{A_{0}} & -\theta_{0} + 3 & 2 \\
		\frac{\zeta_2}{\theta_{0} - 3} & \overline{A_{2}} & -\theta_{0} + 3 & 2 \\
		- \overline{A_{2}} C_{1}& A_{1} & -\theta_{0} + 3 & 2 \\
		B_{1} {\left| A_{1} \right|}^{2} \overline{A_{0}} +2 \, {\left| A_{1} \right|}^{2} \overline{A_{1}} C_{1} - \overline{A_{2}} C_{2} - \frac{1}{2} \, B_{3} \overline{A_{1}} & A_{0} & -\theta_{0} + 3 & 2 
		\end{dmatrix}\\
	\end{align*}
	\normalsize
	As
	\begin{align*}
	\left\{\begin{alignedat}{1}
		\vec{B}_1&=-2\s{\bar{\vec{A}_1}}{\vec{C}_1}\vec{A}_0\\
		\vec{B}_2&=-\frac{(\theta_0+2)}{4\theta_0}|\vec{C}_1|^2\bar{\vec{A}_0}-2\s{\bar{\vec{A}_2}}{\vec{C}_1}\vec{A}_0\\
		\vec{B}_3&=-2\s{\bar{\vec{A}_1}}{\vec{C}_1}\vec{A}_1+\frac{2}{\theta_0-3}\s{\vec{A}_1}{\vec{C}_1}\bar{\vec{A}_1}-2\s{\bar{\vec{A}_1}}{\vec{C}_2}\vec{A}_0\\
		\vec{E}_1&=-\frac{1}{2\theta_0}\s{\vec{C}_1}{\vec{C}_1}\bar{\vec{A}_0}.
	\end{alignedat}\right.
\end{align*}
we have
\begin{align*}
	&\s{\bar{\vec{A}_0}}{\vec{B}_1}=-\s{\bar{\vec{A}_1}}{\vec{C}_1}\\
	&\s{\bar{\vec{A}_1}}{\vec{B}_3}=-2|\vec{A}_1|^2\s{\bar{\vec{A}_1}}{\vec{C}_1}+\frac{2}{\theta_0-3}\s{\vec{A}_1}{\vec{C}_1}\bar{\s{\vec{A}_1}{\vec{A}_1}}\\
	&|\vec{A}_1|^2\s{\bar{\vec{A}_0}}{\vec{B}_1}-\frac{1}{2}\s{\bar{\vec{A}_1}}{\vec{B}_3}=-\colorcancel{|\vec{A}_1|^2\s{\bar{\vec{A}_1}}{\vec{C}_1}}{blue}-\frac{1}{2}\left(-\colorcancel{2|\vec{A}_1|^2\s{\bar{\vec{A}_1}}{\vec{C}_1}}{blue}+\frac{2}{\theta_0-3}\s{\vec{A}_1}{\vec{C}_1}\bar{\s{\vec{A}_1}{\vec{A}_1}}\right)\\
	&=-\frac{1}{\theta_0-3}\s{\vec{A}_1}{\vec{C}_1}\bar{\s{\vec{A}_1}{\vec{A}_1}}\\
	&=-\frac{1}{\theta_0-3}\bar{\zeta_0}\zeta_2
\end{align*}
where
\begin{align*}
	\zeta_0=\s{\vec{A}_1}{\vec{A}_1}
\end{align*}
	so we finally have as
	\begin{align*}
		\alpha_2=\frac{1}{2\theta_0(\theta_0+1)}\s{\bar{\vec{A}_1}}{\vec{C}_1}
	\end{align*}
	the identity
	\begin{align*}
		\frac{1}{2}\vec{B}_5=\begin{dmatrix}
		-\frac{{\left(\theta_{0}^{2} - \theta_{0} - 6\right)} C_{2} \overline{C_{1}} + 8 \,  \theta_{0} \overline{\alpha_{1}} \zeta_2}{8 \, {\left(\theta_{0}^{2} - 3 \, \theta_{0}\right)}} & \overline{A_{0}} & -\theta_{0} + 3 & 2 \\
		\frac{\zeta_2}{\theta_{0} - 3} & \overline{A_{2}} & -\theta_{0} + 3 & 2 \\
		- \overline{A_{2}} C_{1}& A_{1} & -\theta_{0} + 3 & 2 \\
		4\theta_0(\theta_0+1) \, {\left| A_{1} \right|}^{2} \alpha_2 - \overline{A_{2}} C_{2} -\frac{1}{\theta_0-3}\bar{\zeta_0}\zeta_2 & A_{0} & -\theta_{0} + 3 & 2 
		\end{dmatrix}\\
	\end{align*}
	and finally
	\begin{align*}
		\vec{B}_5&=-\left(\frac{(\theta_0+2)}{4}\s{\bar{\vec{C}_1}}{\vec{C}_2}+\frac{2}{\theta_0-3}\bar{\alpha_1}\zeta_2\right)\bar{\vec{A}_0}+\frac{2\,\zeta_2}{\theta_0-3}\bar{\vec{A}_2}-2\s{\bar{\vec{A}_2}}{\vec{C}_1}\vec{A}_1\\
		&+\left(8\theta_0(\theta_0+1)|\vec{A}_1|^2\alpha_2-2\s{\bar{\vec{A}_2}}{\vec{C}_2}-\frac{2}{\theta_0-3}\bar{\zeta_0}\zeta_2\right)\vec{A}_0.
	\end{align*}
	Finally, we have
	\small
	\begin{align*}
		&\frac{1}{2}\vec{B}_6=\\
		&\begin{dmatrix}
		-\frac{2 \, {\left(\ccancel{A_{0} \alpha_{1}} - A_{2}\right)} C_{1} - A_{1} C_{2}}{\theta_{0} - 4} & \overline{A_{1}} & -\theta_{0} + 4 & 1 \\
		-\frac{2 \, A_{0} C_{2} {\left| A_{1} \right|}^{2} + 2 \, {\left(\ccancel{A_{0} \alpha_{1}} - A_{2}\right)} B_{1} - A_{1} B_{3} + 2 \, {\left(2 \, A_{1} {\left| A_{1} \right|}^{2} + \ccancel{A_{0} \alpha_{5}}\right)} C_{1}}{\theta_{0} - 4} & \overline{A_{0}} & -\theta_{0} + 4 & 1 \\
		\frac{\ccancel{16 \, C_{2} {\left(\theta_{0} - 4\right)} {\left| A_{1} \right|}^{2} \overline{A_{0}}} - 8 \, C_{3} {\left(\theta_{0} - 4\right)} \overline{A_{1}} - {\left(\ccancel{B_{1} {\left(\theta_{0} - 6\right)}} - \ccancel{8 \, {\left(\alpha_{5} \theta_{0} - 4 \, \alpha_{5}\right)} \overline{A_{0}}} - 8 \, {\left(\alpha_{1} \theta_{0} - 4 \, \alpha_{1}\right)} \overline{A_{1}}\right)} C_{1}}{8 \, {\left(\theta_{0} - 4\right)}} & A_{0} & -\theta_{0} + 4 & 1 \\
		-C_{1} \overline{A_{1}} & A_{2} & -\theta_{0} + 4 & 1 \\
		\ccancel{2 \, C_{1} {\left| A_{1} \right|}^{2} \overline{A_{0}}} - C_{2} \overline{A_{1}} & A_{1} & -\theta_{0} + 4 & 1 
		\end{dmatrix}\\
		\normalsize
		&=\begin{dmatrix}
		\frac{2 A_{2}C_{1}+ A_{1} C_{2}}{\theta_{0} - 4} & \overline{A_{1}} & -\theta_{0} + 4 & 1 \\
		-\frac{-2 A_{2} B_{1} - A_{1} B_{3} + 2 \, {\left| A_{1} \right|}^{2} A_{1} C_{1}}{\theta_{0} - 4} & \overline{A_{0}} & -\theta_{0} + 4 & 1 \\
		- \overline{A_{1}} C_3 +\alpha_{1}\overline{A_{1}}C_{1}& A_{0} & -\theta_{0} + 4 & 1 \\
		-2\theta_0(\theta_0+1)\alpha_2 & A_{2} & -\theta_{0} + 4 & 1 \\
		- C_{2} \overline{A_{1}} & A_{1} & -\theta_{0} + 4 & 1 
		\end{dmatrix}
	\end{align*}
	We have 
	\begin{align*}
		2\s{\vec{A}_2}{\vec{B}_1}+\s{\vec{A}_1}{\vec{B}_3}&=-4\s{\bar{\vec{A}_1}}{\vec{C}_1}\s{\vec{A}_0}{\vec{A}_2}-2\s{\bar{\vec{A}_1}}{\vec{C}_1}\s{\vec{A}_1}{\vec{A}_1}+\frac{2}{\theta_0-3}|\vec{A}_1|^2\s{\vec{A}_1}{\vec{C}_1}\\
		&=\frac{2}{\theta_0-3}|\vec{A}_1|^2\s{\vec{A}_1}{\vec{C}_1}
	\end{align*}
	so we finally have
	\begin{align}
		-\frac{-2 A_{2} B_{1} - A_{1} B_{3} + 2 \, {\left| A_{1} \right|}^{2} A_{1} C_{1}}{\theta_{0} - 4}&=-\frac{1}{\theta_0-4}\left(2-\frac{2}{\theta_0-3}\right)|\vec{A}_1|^2\s{\vec{A}_1}{\vec{C}_1}\\
		&=-\frac{2}{\theta_0-3}|\vec{A}_1|^2\s{\vec{A}_1}{\vec{C}_1}
	\end{align}
	and
	\begin{align*}
		\vec{B}_6=\frac{2\,\zeta_5}{\theta_0-4}\bar{\vec{A}_1}-\frac{4}{\theta_0-3}|\vec{A}_1|^2\zeta_2\bar{\vec{A}_0}+\left(-2\s{\bar{\vec{A}_1}}{\vec{C}_3}+4\theta_0(\theta_0+1)\alpha_1\alpha_2\right)\vec{A}_0-4\theta_0(\theta_0+1)\alpha_2\vec{A}_2-2\s{\bar{\vec{A}_1}}{C_2}\vec{A}_1.
	\end{align*}
	Finally,
	\begin{align*}
		\left\{
		\begin{alignedat}{1}
		\vec{B}_4&=-\frac{(\theta_0+3)}{6\theta_0}\s{\vec{C}_1}{\bar{\vec{C}_2}}\bar{\vec{A}_0}-\frac{\bar{\zeta_2}}{6\theta_0}\vec{C}_1-\frac{(\theta_0+3)}{6\theta_0}|\vec{C}_1|^2\bar{\vec{A}_1}-\frac{\theta_0(\theta_0+1)}{6}\alpha_2\bar{\vec{C}_1}+2\left(\bar{\alpha_1}\s{\bar{\vec{A}_1}}{\vec{C}_1}-\s{\bar{\vec{A}_3}}{\vec{C}_1}\right)\vec{A}_0.\\
		\vec{B}_5&=-\left(\frac{(\theta_0+2)}{4}\s{\bar{\vec{C}_1}}{\vec{C}_2}+\frac{2}{\theta_0-3}\bar{\alpha_1}\zeta_2\right)\bar{\vec{A}_0}+\frac{2\,\zeta_2}{\theta_0-3}\bar{\vec{A}_2}-2\s{\bar{\vec{A}_2}}{\vec{C}_1}\vec{A}_1\\
		&+\left(8\theta_0(\theta_0+1)|\vec{A}_1|^2\alpha_2-2\s{\bar{\vec{A}_2}}{\vec{C}_2}-\frac{2}{\theta_0-3}\bar{\zeta_0}\zeta_2\right)\vec{A}_0\\
		\vec{B}_6&=\frac{2\,\zeta_5}{\theta_0-4}\bar{\vec{A}_1}-\frac{4}{\theta_0-3}|\vec{A}_1|^2\zeta_2\bar{\vec{A}_0}+\left(-2\s{\bar{\vec{A}_1}}{\vec{C}_3}+4\theta_0(\theta_0+1)\alpha_1\alpha_2\right)\vec{A}_0-4\theta_0(\theta_0+1)\alpha_2\vec{A}_2-2\s{\bar{\vec{A}_1}}{C_2}\vec{A}_1.
		\end{alignedat}\right.
	\end{align*}
	\begin{align*}
		\p{z}\phi=\begin{dmatrix}
		1 & A_{0} & \theta_{0} - 1 & 0 \\
		1 & A_{1} & \theta_{0} & 0 \\
		1 & A_{2} & \theta_{0} + 1 & 0 \\
		1 & A_{3} & \theta_{0} + 2 & 0 \\
		1 & A_{4} & \theta_{0} + 3 & 0 \\
		1 & A_{5} & \theta_{0} + 4 & 0 \\
		\frac{1}{4 \, \theta_{0}} & C_{1} & 1 & \theta_{0} \\
		\frac{1}{4 \, {\left(\theta_{0} + 1\right)}} & B_{1} & 1 & \theta_{0} + 1 \\
		\frac{1}{4 \, {\left(\theta_{0} + 2\right)}} & B_{2} & 1 & \theta_{0} + 2 \\
		\frac{\overline{\alpha_{1}}}{4 \, {\left(\theta_{0} + 2\right)}} & C_{1} & 1 & \theta_{0} + 2 \\
		\frac{1}{4 \, {\left(\theta_{0} + 3\right)}} & B_{4} & 1 & \theta_{0} + 3 \\
		\frac{\overline{\alpha_{3}}}{4 \, {\left(\theta_{0} + 3\right)}} & C_{1} & 1 & \theta_{0} + 3 \\
		\frac{\overline{\alpha_{1}}}{4 \, {\left(\theta_{0} + 3\right)}} & B_{1} & 1 & \theta_{0} + 3 \\
		\frac{\alpha_{2}}{4 \, {\left(\theta_{0} + 3\right)}} & \overline{C_{1}} & 1 & \theta_{0} + 3 \\
		\frac{1}{4 \, \theta_{0}} & C_{2} & 2 & \theta_{0} \\
		\frac{1}{4 \, {\left(\theta_{0} + 1\right)}} & B_{3} & 2 & \theta_{0} + 1 \\
		\frac{{\left| A_{1} \right|}^{2}}{2 \, {\left(\theta_{0} + 1\right)}} & C_{1} & 2 & \theta_{0} + 1 \\
		\frac{1}{4 \, {\left(\theta_{0} + 2\right)}} & B_{5} & 2 & \theta_{0} + 2 \\
		\frac{\overline{\alpha_{5}}}{4 \, {\left(\theta_{0} + 2\right)}} & C_{1} & 2 & \theta_{0} + 2 \\
		\frac{\overline{\alpha_{1}}}{4 \, {\left(\theta_{0} + 2\right)}} & C_{2} & 2 & \theta_{0} + 2 \\
		\frac{{\left| A_{1} \right|}^{2}}{2 \, {\left(\theta_{0} + 2\right)}} & B_{1} & 2 & \theta_{0} + 2 \\
		\frac{1}{4 \, \theta_{0}} & C_{3} & 3 & \theta_{0} \\
		\frac{\alpha_{1}}{4 \, \theta_{0}} & C_{1} & 3 & \theta_{0} \\
		\frac{1}{4 \, {\left(\theta_{0} + 1\right)}} & B_{6} & 3 & \theta_{0} + 1 \\
		\frac{\alpha_{5}}{4 \, {\left(\theta_{0} + 1\right)}} & C_{1} & 3 & \theta_{0} + 1 \\
		\frac{\alpha_{1}}{4 \, {\left(\theta_{0} + 1\right)}} & B_{1} & 3 & \theta_{0} + 1 
						\end{dmatrix}
		\begin{dmatrix}
		\frac{{\left| A_{1} \right|}^{2}}{2 \, {\left(\theta_{0} + 1\right)}} & C_{2} & 3 & \theta_{0} + 1 \\
		\frac{1}{4 \, \theta_{0}} & C_{4} & 4 & \theta_{0} \\
		\frac{\alpha_{3}}{4 \, \theta_{0}} & C_{1} & 4 & \theta_{0} \\
		\frac{\alpha_{1}}{4 \, \theta_{0}} & C_{2} & 4 & \theta_{0} \\
		\frac{1}{4 \, {\left(2 \, \theta_{0} + 1\right)}} & E_{2} & -\theta_{0} + 3 & 2 \, \theta_{0} + 1 \\
		\frac{\alpha_{2}}{4 \, {\left(2 \, \theta_{0} + 1\right)}} & C_{1} & -\theta_{0} + 3 & 2 \, \theta_{0} + 1 \\
		\frac{1}{8 \, \theta_{0}} & E_{1} & -\theta_{0} + 3 & 2 \, \theta_{0} \\
		\frac{1}{8 \, \theta_{0}} & E_{3} & -\theta_{0} + 4 & 2 \, \theta_{0} \\
		\frac{1}{8} & \overline{B_{1}} & \theta_{0} & 2 \\
		\frac{1}{12} & \overline{B_{3}} & \theta_{0} & 3 \\
		\frac{1}{6} \, {\left| A_{1} \right|}^{2} & \overline{C_{1}} & \theta_{0} & 3 \\
		\frac{1}{16} & \overline{B_{6}} & \theta_{0} & 4 \\
		\frac{1}{16} \, \overline{\alpha_{5}} & \overline{C_{1}} & \theta_{0} & 4 \\
		\frac{1}{16} \, \overline{\alpha_{1}} & \overline{B_{1}} & \theta_{0} & 4 \\
		\frac{1}{8} \, {\left| A_{1} \right|}^{2} & \overline{C_{2}} & \theta_{0} & 4 \\
		\frac{1}{8} & \overline{C_{1}} & \theta_{0} - 1 & 2 \\
		\frac{1}{12} & \overline{C_{2}} & \theta_{0} - 1 & 3 \\
		\frac{1}{16} & \overline{C_{3}} & \theta_{0} - 1 & 4 \\
		\frac{1}{16} \, \overline{\alpha_{1}} & \overline{C_{1}} & \theta_{0} - 1 & 4 \\
		\frac{1}{20} & \overline{C_{4}} & \theta_{0} - 1 & 5 \\
		\frac{1}{20} \, \overline{\alpha_{3}} & \overline{C_{1}} & \theta_{0} - 1 & 5 \\
		\frac{1}{20} \, \overline{\alpha_{1}} & \overline{C_{2}} & \theta_{0} - 1 & 5 \\
		\frac{1}{8} & \overline{B_{2}} & \theta_{0} + 1 & 2 \\
		\frac{1}{8} \, \alpha_{1} & \overline{C_{1}} & \theta_{0} + 1 & 2 \\
		\frac{1}{12} & \overline{B_{5}} & \theta_{0} + 1 & 3 \\
		\frac{1}{12} \, \alpha_{5} & \overline{C_{1}} & \theta_{0} + 1 & 3 \\
		\frac{1}{12} \, \alpha_{1} & \overline{C_{2}} & \theta_{0} + 1 & 3 
						\end{dmatrix}
		\end{align*}
		\begin{align}\label{rastophi}
		\begin{dmatrix}
		\frac{1}{6} \, {\left| A_{1} \right|}^{2} & \overline{B_{1}} & \theta_{0} + 1 & 3 \\
		\frac{1}{8} & \overline{B_{4}} & \theta_{0} + 2 & 2 \\
		\frac{1}{8} \, \overline{\alpha_{2}} & C_{1} & \theta_{0} + 2 & 2 \\
		\frac{1}{8} \, \alpha_{3} & \overline{C_{1}} & \theta_{0} + 2 & 2 \\
		\frac{1}{8} \, \alpha_{1} & \overline{B_{1}} & \theta_{0} + 2 & 2 \\
		-\frac{1}{4 \, {\left(\theta_{0} - 4\right)}} & \overline{E_{1}} & 2 \, \theta_{0} - 1 & -\theta_{0} + 4 \\
		-\frac{1}{4 \, {\left(\theta_{0} - 5\right)}} & \overline{E_{3}} & 2 \, \theta_{0} - 1 & -\theta_{0} + 5 \\
		-\frac{1}{4 \, {\left(\theta_{0} - 4\right)}} & \overline{E_{2}} & 2 \, \theta_{0} & -\theta_{0} + 4 \\
		-\frac{\overline{\alpha_{2}}}{4 \, {\left(\theta_{0} - 4\right)}} & \overline{C_{1}} & 2 \, \theta_{0} & -\theta_{0} + 4
		\end{dmatrix}
	\end{align}
	\small
	\begin{align*}
		&0=\s{\p{z}\phi}{\p{z}\phi}=\\
&	\begin{dmatrix}
	 \frac{4 \, {\left(\theta_{0}^{2} + \theta_{0}\right)} A_{0} E_{2} + {\left(4 \, {\left(\alpha_{2} \theta_{0}^{2} + \alpha_{2} \theta_{0}\right)} A_{0} + B_{1} {\left(2 \, \theta_{0} + 1\right)}\right)} C_{1}}{8 \, {\left(2 \, \theta_{0}^{3} + 3 \, \theta_{0}^{2} + \theta_{0}\right)}} & 2 & 2 \, \theta_{0} + 1 \\
	\frac{4 \, A_{0} E_{1} \theta_{0} + C_{1}^{2}}{16 \, \theta_{0}^{2}} & 2 & 2 \, \theta_{0} \\
	\frac{2 \, A_{1} E_{1} \theta_{0} + 2 \, A_{0} E_{3} \theta_{0} + C_{1} C_{2}}{8 \, \theta_{0}^{2}} & 3 & 2 \, \theta_{0} \\
	\frac{A_{0} C_{1}}{2 \, \theta_{0}} & \theta_{0} & \theta_{0} \\
	\frac{A_{0} B_{1}}{2 \, {\left(\theta_{0} + 1\right)}} & \theta_{0} & \theta_{0} + 1 \\
	\frac{8 \, A_{0} C_{1} \theta_{0} \overline{\alpha_{1}} + 8 \, A_{0} B_{2} \theta_{0} + C_{1} {\left(\theta_{0} + 2\right)} \overline{C_{1}}}{16 \, {\left(\theta_{0}^{2} + 2 \, \theta_{0}\right)}} & \theta_{0} & \theta_{0} + 2 \\
	\lambda_1 & \theta_{0} & \theta_{0} + 3 \\
	\frac{A_{1} C_{1} + A_{0} C_{2}}{2 \, \theta_{0}} & \theta_{0} + 1 & \theta_{0} \\
	\frac{2 \, A_{0} C_{1} {\left| A_{1} \right|}^{2} + A_{1} B_{1} + A_{0} B_{3}}{2 \, {\left(\theta_{0} + 1\right)}} & \theta_{0} + 1 & \theta_{0} + 1 \\
	\lambda_2 & \theta_{0} + 1 & \theta_{0} + 2 \\
	\frac{{\left(A_{0} \alpha_{1} + A_{2}\right)} C_{1} + A_{1} C_{2} + A_{0} C_{3}}{2 \, \theta_{0}} & \theta_{0} + 2 & \theta_{0} \\
	\frac{2 \, A_{0} C_{2} {\left| A_{1} \right|}^{2} + {\left(A_{0} \alpha_{1} + A_{2}\right)} B_{1} + A_{1} B_{3} + A_{0} B_{6} + {\left(2 \, A_{1} {\left| A_{1} \right|}^{2} + A_{0} \alpha_{5}\right)} C_{1}}{2 \, {\left(\theta_{0} + 1\right)}} & \theta_{0} + 2 & \theta_{0} + 1 \\
	\frac{{\left(A_{1} \alpha_{1} + A_{0} \alpha_{3} + A_{3}\right)} C_{1} + {\left(A_{0} \alpha_{1} + A_{2}\right)} C_{2} + A_{1} C_{3} + A_{0} C_{4}}{2 \, \theta_{0}} & \theta_{0} + 3 & \theta_{0} \\
	A_{0}^{2} & 2 \, \theta_{0} - 2 & 0 \\
	\frac{1}{4} \, A_{0} \overline{C_{1}} & 2 \, \theta_{0} - 2 & 2 \\
	\frac{1}{6} \, A_{0} \overline{C_{2}} & 2 \, \theta_{0} - 2 & 3 \\
	\frac{1}{8} \, A_{0} \overline{C_{1}} \overline{\alpha_{1}} + \frac{1}{64} \, \overline{C_{1}}^{2} + \frac{1}{8} \, A_{0} \overline{C_{3}} & 2 \, \theta_{0} - 2 & 4 \\
	\frac{1}{10} \, A_{0} \overline{C_{1}} \overline{\alpha_{3}} + \frac{1}{240} \, {\left(24 \, A_{0} \overline{\alpha_{1}} + 5 \, \overline{C_{1}}\right)} \overline{C_{2}} + \frac{1}{10} \, A_{0} \overline{C_{4}} & 2 \, \theta_{0} - 2 & 5 \\
	2 \, A_{0} A_{1} & 2 \, \theta_{0} - 1 & 0 \\
	\frac{1}{4} \, A_{0} \overline{B_{1}} + \frac{1}{4} \, A_{1} \overline{C_{1}} & 2 \, \theta_{0} - 1 & 2 \\
	\frac{1}{3} \, A_{0} {\left| A_{1} \right|}^{2} \overline{C_{1}} + \frac{1}{6} \, A_{0} \overline{B_{3}} + \frac{1}{6} \, A_{1} \overline{C_{2}} & 2 \, \theta_{0} - 1 & 3 \\
	\frac{1}{4} \, A_{0} {\left| A_{1} \right|}^{2} \overline{C_{2}} + \frac{1}{8} \, A_{0} \overline{B_{1}} \overline{\alpha_{1}} + \frac{1}{8} \, A_{0} \overline{B_{6}} + \frac{1}{32} \, {\left(4 \, A_{1} \overline{\alpha_{1}} + 4 \, A_{0} \overline{\alpha_{5}} + \overline{B_{1}}\right)} \overline{C_{1}} + \frac{1}{8} \, A_{1} \overline{C_{3}} & 2 \, \theta_{0} - 1 & 4 \\
	2 \, A_{1} A_{2} + 2 \, A_{0} A_{3} & 2 \, \theta_{0} + 1 & 0 \\
	\frac{1}{4} \, A_{0} C_{1} \overline{\alpha_{2}} + \frac{1}{4} \, {\left(A_{0} \alpha_{1} + A_{2}\right)} \overline{B_{1}} + \frac{1}{4} \, A_{1} \overline{B_{2}} + \frac{1}{4} \, A_{0} \overline{B_{4}} + \frac{1}{4} \, {\left(A_{1} \alpha_{1} + A_{0} \alpha_{3} + A_{3}\right)} \overline{C_{1}} & 2 \, \theta_{0} + 1 & 2 \\
	A_{2}^{2} + 2 \, A_{1} A_{3} + 2 \, A_{0} A_{4} & 2 \, \theta_{0} + 2 & 0 
	\end{dmatrix}
	\end{align*}
	\begin{align*}
	\begin{dmatrix}
	2 \, A_{2} A_{3} + 2 \, A_{1} A_{4} + 2 \, A_{0} A_{5} & 2 \, \theta_{0} + 3 & 0 \\
	-\frac{A_{0} \overline{E_{1}}}{2 \, {\left(\theta_{0} - 4\right)}} & 3 \, \theta_{0} - 2 & -\theta_{0} + 4 \\
	-\frac{A_{0} \overline{E_{3}}}{2 \, {\left(\theta_{0} - 5\right)}} & 3 \, \theta_{0} - 2 & -\theta_{0} + 5 \\
	-\frac{A_{0} \overline{C_{1}} \overline{\alpha_{2}} + A_{1} \overline{E_{1}} + A_{0} \overline{E_{2}}}{2 \, {\left(\theta_{0} - 4\right)}} & 3 \, \theta_{0} - 1 & -\theta_{0} + 4 \\
	A_{1}^{2} + 2 \, A_{0} A_{2} & 2 \, \theta_{0} & 0 \\
	\frac{1}{4} \, A_{1} \overline{B_{1}} + \frac{1}{4} \, A_{0} \overline{B_{2}} + \frac{1}{4} \, {\left(A_{0} \alpha_{1} + A_{2}\right)} \overline{C_{1}} & 2 \, \theta_{0} & 2 \\
	\frac{1}{3} \, A_{0} {\left| A_{1} \right|}^{2} \overline{B_{1}} + \frac{1}{6} \, A_{1} \overline{B_{3}} + \frac{1}{6} \, A_{0} \overline{B_{5}} + \frac{1}{6} \, {\left(2 \, A_{1} {\left| A_{1} \right|}^{2} + A_{0} \alpha_{5}\right)} \overline{C_{1}} + \frac{1}{6} \, {\left(A_{0} \alpha_{1} + A_{2}\right)} \overline{C_{2}} & 2 \, \theta_{0} & 3
	\end{dmatrix}
	\end{align*}
	\normalsize
    where
    \begin{align*}
    	\lambda_1&=\frac{1}{48 \, {\left(\theta_{0}^{3} + 4 \, \theta_{0}^{2} + 3 \, \theta_{0}\right)}}\bigg\{24 \, {\left(\theta_{0}^{2} \overline{\alpha_{1}} + \theta_{0} \overline{\alpha_{1}}\right)} A_{0} B_{1} + 24 \, {\left(\theta_{0}^{2} + \theta_{0}\right)} A_{0} B_{4} + 24 \, {\left(\theta_{0}^{2} \overline{\alpha_{3}} + \theta_{0} \overline{\alpha_{3}}\right)} A_{0} C_{1}\\
    	& + 2 \, {\left(\theta_{0}^{2} + 4 \, \theta_{0} + 3\right)} C_{1} \overline{C_{2}} + 3 \, {\left(8 \, {\left(\alpha_{2} \theta_{0}^{2} + \alpha_{2} \theta_{0}\right)} A_{0} + {\left(\theta_{0}^{2} + 3 \, \theta_{0}\right)} B_{1}\right)} \overline{C_{1}}\bigg\}\\
    	\lambda_2&=\frac{1}{16 \, {\left(\theta_{0}^{2} + 2 \, \theta_{0}\right)}}\bigg\{16 \, A_{0} B_{1} \theta_{0} {\left| A_{1} \right|}^{2} + 8 \, A_{1} B_{2} \theta_{0} + 8 \, A_{0} B_{5} \theta_{0} + {\left(8 \, A_{1} \theta_{0} \overline{\alpha_{1}} + 8 \, A_{0} \theta_{0} \overline{\alpha_{5}} + {\left(\theta_{0} + 2\right)} \overline{B_{1}}\right)} C_{1}\\
    	& + {\left(8 \, A_{0} \theta_{0} \overline{\alpha_{1}} + {\left(\theta_{0} + 2\right)} \overline{C_{1}}\right)} C_{2}\bigg\}
    \end{align*}
    or
	\small
	\begin{align*}
	e^{2\lambda}=\begin{dmatrix}
	\frac{C_{1} \overline{A_{0}}}{2 \, \theta_{0}} & 1 & 2 \, \theta_{0} - 1 \\
	\frac{{\left(\theta_{0}^{2} + \theta_{0}\right)} B_{2} \overline{A_{0}} + {\left(\theta_{0}^{2} + 2 \, \theta_{0}\right)} B_{1} \overline{A_{1}} + {\left({\left(\theta_{0}^{2} \overline{\alpha_{1}} + \theta_{0} \overline{\alpha_{1}}\right)} \overline{A_{0}} + {\left(\theta_{0}^{2} + 3 \, \theta_{0} + 2\right)} \overline{A_{2}}\right)} C_{1}}{2 \, {\left(\theta_{0}^{3} + 3 \, \theta_{0}^{2} + 2 \, \theta_{0}\right)}} & 1 & 2 \, \theta_{0} + 1 \\
	\mu_1 & 1 & 2 \, \theta_{0} + 2 \\
	\frac{B_{1} \theta_{0} \overline{A_{0}} + C_{1} {\left(\theta_{0} + 1\right)} \overline{A_{1}}}{2 \, {\left(\theta_{0}^{2} + \theta_{0}\right)}} & 1 & 2 \, \theta_{0} \\
	\frac{C_{2} \overline{A_{0}}}{2 \, \theta_{0}} & 2 & 2 \, \theta_{0} - 1 \\
	\mu_2 & 2 & 2 \, \theta_{0} + 1 \\
	\frac{2 \, C_{1} \theta_{0} {\left| A_{1} \right|}^{2} \overline{A_{0}} + B_{3} \theta_{0} \overline{A_{0}} + C_{2} {\left(\theta_{0} + 1\right)} \overline{A_{1}}}{2 \, {\left(\theta_{0}^{2} + \theta_{0}\right)}} & 2 & 2 \, \theta_{0} \\
	\frac{C_{1}^{2} {\left(\theta_{0} - 4\right)} - 8 \, A_{0} E_{1} \theta_{0} + 8 \, {\left(\alpha_{1} \theta_{0} - 4 \, \alpha_{1}\right)} C_{1} \overline{A_{0}} + 8 \, C_{3} {\left(\theta_{0} - 4\right)} \overline{A_{0}}}{16 \, {\left(\theta_{0}^{2} - 4 \, \theta_{0}\right)}} & 3 & 2 \, \theta_{0} - 1 \\
	\mu_3 & 3 & 2 \, \theta_{0} \\
	\mu_4 & 4 & 2 \, \theta_{0} - 1 \\
	\frac{E_{1} \overline{A_{0}}}{4 \, \theta_{0}} & -\theta_{0} + 3 & 3 \, \theta_{0} - 1 \\
	\frac{2 \, C_{1} \alpha_{2} \theta_{0} \overline{A_{0}} + 2 \, E_{2} \theta_{0} \overline{A_{0}} + E_{1} {\left(2 \, \theta_{0} + 1\right)} \overline{A_{1}}}{4 \, {\left(2 \, \theta_{0}^{2} + \theta_{0}\right)}} & -\theta_{0} + 3 & 3 \, \theta_{0} \\
	\frac{E_{3} \overline{A_{0}}}{4 \, \theta_{0}} & -\theta_{0} + 4 & 3 \, \theta_{0} - 1 \\
	2 \, A_{1} \overline{A_{1}} & \theta_{0} & \theta_{0} \\
	2 \, A_{1} \overline{A_{0}} & \theta_{0} & \theta_{0} - 1 \\
	2 \, A_{1} \overline{A_{2}} + \frac{1}{4} \, \overline{A_{0}} \overline{B_{1}} & \theta_{0} & \theta_{0} + 1 \\
	\frac{1}{3} \, {\left| A_{1} \right|}^{2} \overline{A_{0}} \overline{C_{1}} + 2 \, A_{1} \overline{A_{3}} + \frac{1}{4} \, \overline{A_{1}} \overline{B_{1}} + \frac{1}{6} \, \overline{A_{0}} \overline{B_{3}} & \theta_{0} & \theta_{0} + 2 \\
	\frac{1}{4} \, {\left| A_{1} \right|}^{2} \overline{A_{0}} \overline{C_{2}} + 2 \, A_{1} \overline{A_{4}} + \frac{1}{8} \, {\left(\overline{A_{0}} \overline{\alpha_{1}} + 2 \, \overline{A_{2}}\right)} \overline{B_{1}} + \frac{1}{6} \, \overline{A_{1}} \overline{B_{3}} + \frac{1}{8} \, \overline{A_{0}} \overline{B_{6}} + \frac{1}{24} \, {\left(8 \, {\left| A_{1} \right|}^{2} \overline{A_{1}} + 3 \, \overline{A_{0}} \overline{\alpha_{5}}\right)} \overline{C_{1}} & \theta_{0} & \theta_{0} + 3 \\
	2 \, A_{0} \overline{A_{1}} & \theta_{0} - 1 & \theta_{0} \\
	2 \, A_{0} \overline{A_{0}} & \theta_{0} - 1 & \theta_{0} - 1 \\
	2 \, A_{0} \overline{A_{2}} + \frac{1}{4} \, \overline{A_{0}} \overline{C_{1}} & \theta_{0} - 1 & \theta_{0} + 1 \\
	2 \, A_{0} \overline{A_{3}} + \frac{1}{4} \, \overline{A_{1}} \overline{C_{1}} + \frac{1}{6} \, \overline{A_{0}} \overline{C_{2}} & \theta_{0} - 1 & \theta_{0} + 2 \\
	2 \, A_{0} \overline{A_{4}} + \frac{1}{8} \, {\left(\overline{A_{0}} \overline{\alpha_{1}} + 2 \, \overline{A_{2}}\right)} \overline{C_{1}} + \frac{1}{6} \, \overline{A_{1}} \overline{C_{2}} + \frac{1}{8} \, \overline{A_{0}} \overline{C_{3}} & \theta_{0} - 1 & \theta_{0} + 3 \\
	2 \, A_{0} \overline{A_{5}} + \frac{1}{40} \, {\left(5 \, \overline{A_{1}} \overline{\alpha_{1}} + 4 \, \overline{A_{0}} \overline{\alpha_{3}} + 10 \, \overline{A_{3}}\right)} \overline{C_{1}} + \frac{1}{30} \, {\left(3 \, \overline{A_{0}} \overline{\alpha_{1}} + 5 \, \overline{A_{2}}\right)} \overline{C_{2}} + \frac{1}{8} \, \overline{A_{1}} \overline{C_{3}} + \frac{1}{10} \, \overline{A_{0}} \overline{C_{4}} & \theta_{0} - 1 & \theta_{0} + 4 \\
	\frac{1}{4} \, A_{0} B_{1} + 2 \, A_{2} \overline{A_{1}} & \theta_{0} + 1 & \theta_{0} \\
	\frac{1}{4} \, A_{0} C_{1} + 2 \, A_{2} \overline{A_{0}} & \theta_{0} + 1 & \theta_{0} - 1 \\
	\frac{8 \, A_{0} C_{1} \theta_{0}^{2} \overline{\alpha_{1}} + 8 \, A_{0} B_{2} \theta_{0}^{2} + 64 \, A_{2} \theta_{0}^{2} \overline{A_{2}} + 8 \, \theta_{0}^{2} \overline{A_{0}} \overline{B_{2}} + {\left(8 \, \alpha_{1} \theta_{0}^{2} \overline{A_{0}} + {\left(\theta_{0}^{2} + 4\right)} C_{1}\right)} \overline{C_{1}}}{32 \, \theta_{0}^{2}} & \theta_{0} + 1 & \theta_{0} + 1 
		\end{dmatrix}
	\end{align*}
	\begin{align*}
	\begin{dmatrix}
	\mu_5 & \theta_{0} + 1 & \theta_{0} + 2 \\
	\frac{1}{3} \, A_{0} C_{1} {\left| A_{1} \right|}^{2} + \frac{1}{4} \, A_{1} B_{1} + \frac{1}{6} \, A_{0} B_{3} + 2 \, A_{3} \overline{A_{1}} & \theta_{0} + 2 & \theta_{0} \\
	\frac{1}{4} \, A_{1} C_{1} + \frac{1}{6} \, A_{0} C_{2} + 2 \, A_{3} \overline{A_{0}} & \theta_{0} + 2 & \theta_{0} - 1 \\
	\mu_6 & \theta_{0} + 2 & \theta_{0} + 1 \\
	\frac{1}{4} \, A_{0} C_{2} {\left| A_{1} \right|}^{2} + \frac{1}{8} \, {\left(A_{0} \alpha_{1} + 2 \, A_{2}\right)} B_{1} + \frac{1}{6} \, A_{1} B_{3} + \frac{1}{8} \, A_{0} B_{6} + \frac{1}{24} \, {\left(8 \, A_{1} {\left| A_{1} \right|}^{2} + 3 \, A_{0} \alpha_{5}\right)} C_{1} + 2 \, A_{4} \overline{A_{1}} & \theta_{0} + 3 & \theta_{0} \\
	\frac{1}{8} \, {\left(A_{0} \alpha_{1} + 2 \, A_{2}\right)} C_{1} + \frac{1}{6} \, A_{1} C_{2} + \frac{1}{8} \, A_{0} C_{3} + 2 \, A_{4} \overline{A_{0}} & \theta_{0} + 3 & \theta_{0} - 1 \\
	\frac{1}{40} \, {\left(5 \, A_{1} \alpha_{1} + 4 \, A_{0} \alpha_{3} + 10 \, A_{3}\right)} C_{1} + \frac{1}{30} \, {\left(3 \, A_{0} \alpha_{1} + 5 \, A_{2}\right)} C_{2} + \frac{1}{8} \, A_{1} C_{3} + \frac{1}{10} \, A_{0} C_{4} + 2 \, A_{5} \overline{A_{0}} & \theta_{0} + 4 & \theta_{0} - 1 \\
	\frac{A_{0} \overline{C_{1}}}{2 \, \theta_{0}} & 2 \, \theta_{0} - 1 & 1 \\
	\frac{A_{0} \overline{C_{2}}}{2 \, \theta_{0}} & 2 \, \theta_{0} - 1 & 2 \\
	\frac{8 \, {\left(\theta_{0} \overline{\alpha_{1}} - 4 \, \overline{\alpha_{1}}\right)} A_{0} \overline{C_{1}} + {\left(\theta_{0} - 4\right)} \overline{C_{1}}^{2} + 8 \, A_{0} {\left(\theta_{0} - 4\right)} \overline{C_{3}} - 8 \, \theta_{0} \overline{A_{0}} \overline{E_{1}}}{16 \, {\left(\theta_{0}^{2} - 4 \, \theta_{0}\right)}} & 2 \, \theta_{0} - 1 & 3 \\
	\mu_7 & 2 \, \theta_{0} - 1 & 4 \\
	\frac{{\left(\theta_{0}^{2} + 2 \, \theta_{0}\right)} A_{1} \overline{B_{1}} + {\left(\theta_{0}^{2} + \theta_{0}\right)} A_{0} \overline{B_{2}} + {\left({\left(\alpha_{1} \theta_{0}^{2} + \alpha_{1} \theta_{0}\right)} A_{0} + {\left(\theta_{0}^{2} + 3 \, \theta_{0} + 2\right)} A_{2}\right)} \overline{C_{1}}}{2 \, {\left(\theta_{0}^{3} + 3 \, \theta_{0}^{2} + 2 \, \theta_{0}\right)}} & 2 \, \theta_{0} + 1 & 1 \\
	\mu_8 & 2 \, \theta_{0} + 1 & 2 \\
	\mu_9 & 2 \, \theta_{0} + 2 & 1 \\
	\frac{A_{0} \overline{E_{1}}}{4 \, \theta_{0}} & 3 \, \theta_{0} - 1 & -\theta_{0} + 3 \\
	\frac{A_{0} \overline{E_{3}}}{4 \, \theta_{0}} & 3 \, \theta_{0} - 1 & -\theta_{0} + 4 \\
	\frac{A_{0} \theta_{0} \overline{B_{1}} + A_{1} {\left(\theta_{0} + 1\right)} \overline{C_{1}}}{2 \, {\left(\theta_{0}^{2} + \theta_{0}\right)}} & 2 \, \theta_{0} & 1 \\
	\frac{2 \, A_{0} \theta_{0} {\left| A_{1} \right|}^{2} \overline{C_{1}} + A_{0} \theta_{0} \overline{B_{3}} + A_{1} {\left(\theta_{0} + 1\right)} \overline{C_{2}}}{2 \, {\left(\theta_{0}^{2} + \theta_{0}\right)}} & 2 \, \theta_{0} & 2 \\
	\mu_{10} & 2 \, \theta_{0} & 3 \\
	\frac{2 \, A_{0} \theta_{0} \overline{C_{1}} \overline{\alpha_{2}} + A_{1} {\left(2 \, \theta_{0} + 1\right)} \overline{E_{1}} + 2 \, A_{0} \theta_{0} \overline{E_{2}}}{4 \, {\left(2 \, \theta_{0}^{2} + \theta_{0}\right)}} & 3 \, \theta_{0} & -\theta_{0} + 3
	\end{dmatrix}
\end{align*}
where
\begin{align*}
&\mu_{1}=\frac{1}{2 \, {\left(\theta_{0}^{4} + 6 \, \theta_{0}^{3} + 11 \, \theta_{0}^{2} + 6 \, \theta_{0}\right)}}\bigg\{{\left(\theta_{0}^{3} + 3 \, \theta_{0}^{2} + 2 \, \theta_{0}\right)} B_{4} \overline{A_{0}} + {\left(\theta_{0}^{3} + 4 \, \theta_{0}^{2} + 3 \, \theta_{0}\right)} B_{2} \overline{A_{1}} + {\left(\alpha_{2} \theta_{0}^{3} + 3 \, \alpha_{2} \theta_{0}^{2} + 2 \, \alpha_{2} \theta_{0}\right)} \overline{A_{0}} \overline{C_{1}}\\
& + {\left({\left(\theta_{0}^{3} \overline{\alpha_{1}} + 3 \, \theta_{0}^{2} \overline{\alpha_{1}} + 2 \, \theta_{0} \overline{\alpha_{1}}\right)} \overline{A_{0}} + {\left(\theta_{0}^{3} + 5 \, \theta_{0}^{2} + 6 \, \theta_{0}\right)} \overline{A_{2}}\right)} B_{1} + \bigg({\left(\theta_{0}^{3} \overline{\alpha_{3}} + 3 \, \theta_{0}^{2} \overline{\alpha_{3}} + 2 \, \theta_{0} \overline{\alpha_{3}}\right)} \overline{A_{0}}\\
& + {\left(\theta_{0}^{3} \overline{\alpha_{1}} + 4 \, \theta_{0}^{2} \overline{\alpha_{1}} + 3 \, \theta_{0} \overline{\alpha_{1}}\right)} \overline{A_{1}} + {\left(\theta_{0}^{3} + 6 \, \theta_{0}^{2} + 11 \, \theta_{0} + 6\right)} \overline{A_{3}}\bigg) C_{1}\bigg\}\\
&\mu_{2}=\frac{1}{2 \, {\left(\theta_{0}^{3} + 3 \, \theta_{0}^{2} + 2 \, \theta_{0}\right)}}\bigg\{2 \, {\left(\theta_{0}^{2} + \theta_{0}\right)} B_{1} {\left| A_{1} \right|}^{2} \overline{A_{0}} + {\left(\theta_{0}^{2} + \theta_{0}\right)} B_{5} \overline{A_{0}} + {\left(\theta_{0}^{2} + 2 \, \theta_{0}\right)} B_{3} \overline{A_{1}}\\
& + {\left(2 \, {\left(\theta_{0}^{2} + 2 \, \theta_{0}\right)} {\left| A_{1} \right|}^{2} \overline{A_{1}} + {\left(\theta_{0}^{2} \overline{\alpha_{5}} + \theta_{0} \overline{\alpha_{5}}\right)} \overline{A_{0}}\right)} C_{1} + {\left({\left(\theta_{0}^{2} \overline{\alpha_{1}} + \theta_{0} \overline{\alpha_{1}}\right)} \overline{A_{0}} + {\left(\theta_{0}^{2} + 3 \, \theta_{0} + 2\right)} \overline{A_{2}}\right)} C_{2}\bigg\}\\
&\mu_{3}=\frac{1}{16 \, {\left(\theta_{0}^{3} - 3 \, \theta_{0}^{2} - 4 \, \theta_{0}\right)}}\bigg\{16 \, {\left(\theta_{0}^{2} - 4 \, \theta_{0}\right)} C_{2} {\left| A_{1} \right|}^{2} \overline{A_{0}} - 8 \, {\left(\theta_{0}^{2} + \theta_{0}\right)} A_{0} E_{2} + 8 \, {\left(\alpha_{1} \theta_{0}^{2} - 4 \, \alpha_{1} \theta_{0}\right)} B_{1} \overline{A_{0}} + 8 \, {\left(\theta_{0}^{2} - 4 \, \theta_{0}\right)} B_{6} \overline{A_{0}} \\
&+ 8 \, {\left(\theta_{0}^{2} - 3 \, \theta_{0} - 4\right)} C_{3} \overline{A_{1}} - {\left(8 \, {\left(\alpha_{2} \theta_{0}^{2} + \alpha_{2} \theta_{0}\right)} A_{0} - {\left(2 \, \theta_{0}^{2} - 7 \, \theta_{0} - 4\right)} B_{1} - 8 \, {\left(\alpha_{5} \theta_{0}^{2} - 4 \, \alpha_{5} \theta_{0}\right)} \overline{A_{0}} - 8 \, {\left(\alpha_{1} \theta_{0}^{2} - 3 \, \alpha_{1} \theta_{0} - 4 \, \alpha_{1}\right)} \overline{A_{1}}\right)} C_{1}\bigg\}\\
&\mu_{4}=-\frac{1}{48 \, {\left(\theta_{0}^{3} - 9 \, \theta_{0}^{2} + 20 \, \theta_{0}\right)}}\bigg\{24 \, {\left(\theta_{0}^{2} - 5 \, \theta_{0}\right)} A_{1} E_{1} + 24 \, {\left(\theta_{0}^{2} - 4 \, \theta_{0}\right)} A_{0} E_{3} - 24 \, {\left(\alpha_{3} \theta_{0}^{2} - 9 \, \alpha_{3} \theta_{0} + 20 \, \alpha_{3}\right)} C_{1} \overline{A_{0}}\\
& - 24 \, {\left(\theta_{0}^{2} - 9 \, \theta_{0} + 20\right)} C_{4} \overline{A_{0}} - {\left(5 \, {\left(\theta_{0}^{2} - 9 \, \theta_{0} + 20\right)} C_{1} + 24 \, {\left(\alpha_{1} \theta_{0}^{2} - 9 \, \alpha_{1} \theta_{0} + 20 \, \alpha_{1}\right)} \overline{A_{0}}\right)} C_{2}\bigg\}\\
&\mu_{5}=\frac{1}{96 \, {\left(\theta_{0}^{3} + \theta_{0}^{2}\right)}}\bigg\{32 \, {\left(\theta_{0}^{3} + \theta_{0}^{2}\right)} {\left| A_{1} \right|}^{2} \overline{A_{0}} \overline{B_{1}} + 24 \, {\left(\theta_{0}^{3} \overline{\alpha_{1}} + \theta_{0}^{2} \overline{\alpha_{1}}\right)} A_{0} B_{1} + 24 \, {\left(\theta_{0}^{3} + \theta_{0}^{2}\right)} A_{0} B_{4} + 24 \, {\left(\theta_{0}^{3} \overline{\alpha_{3}} + \theta_{0}^{2} \overline{\alpha_{3}}\right)} A_{0} C_{1}\\
& + 192 \, {\left(\theta_{0}^{3} + \theta_{0}^{2}\right)} A_{2} \overline{A_{3}} + 24 \, {\left(\theta_{0}^{3} + \theta_{0}^{2}\right)} \overline{A_{1}} \overline{B_{2}} + 16 \, {\left(\theta_{0}^{3} + \theta_{0}^{2}\right)} \overline{A_{0}} \overline{B_{5}} + \bigg(24 \, {\left(\alpha_{2} \theta_{0}^{3} + \alpha_{2} \theta_{0}^{2}\right)} A_{0} + 3 \, {\left(\theta_{0}^{3} + \theta_{0}^{2} + 4 \, \theta_{0}\right)} B_{1} \\
&+ 16 \, {\left(\alpha_{5} \theta_{0}^{3} + \alpha_{5} \theta_{0}^{2}\right)} \overline{A_{0}} + 24 \, {\left(\alpha_{1} \theta_{0}^{3} + \alpha_{1} \theta_{0}^{2}\right)} \overline{A_{1}}\bigg) \overline{C_{1}} + 2 \, {\left({\left(\theta_{0}^{3} + \theta_{0}^{2} + 6 \, \theta_{0} + 6\right)} C_{1} + 8 \, {\left(\alpha_{1} \theta_{0}^{3} + \alpha_{1} \theta_{0}^{2}\right)} \overline{A_{0}}\right)} \overline{C_{2}}\bigg\}\\
\mu_{6}&=\frac{1}{96 \, {\left(\theta_{0}^{3} + \theta_{0}^{2}\right)}}\bigg\{32 \, {\left(\theta_{0}^{3} + \theta_{0}^{2}\right)} A_{0} B_{1} {\left| A_{1} \right|}^{2} + 24 \, {\left(\theta_{0}^{3} + \theta_{0}^{2}\right)} A_{1} B_{2} + 16 \, {\left(\theta_{0}^{3} + \theta_{0}^{2}\right)} A_{0} B_{5} + 192 \, {\left(\theta_{0}^{3} + \theta_{0}^{2}\right)} A_{3} \overline{A_{2}}\\
& + 24 \, {\left(\alpha_{1} \theta_{0}^{3} + \alpha_{1} \theta_{0}^{2}\right)} \overline{A_{0}} \overline{B_{1}} + 24 \, {\left(\theta_{0}^{3} + \theta_{0}^{2}\right)} \overline{A_{0}} \overline{B_{4}} + 24 \, {\left(\alpha_{3} \theta_{0}^{3} + \alpha_{3} \theta_{0}^{2}\right)} \overline{A_{0}} \overline{C_{1}} + \bigg(16 \, {\left(\theta_{0}^{3} \overline{\alpha_{5}} + \theta_{0}^{2} \overline{\alpha_{5}}\right)} A_{0}\\
& + 24 \, {\left(\theta_{0}^{3} \overline{\alpha_{1}} + \theta_{0}^{2} \overline{\alpha_{1}}\right)} A_{1} + 24 \, {\left(\theta_{0}^{3} \overline{\alpha_{2}} + \theta_{0}^{2} \overline{\alpha_{2}}\right)} \overline{A_{0}} + 3 \, {\left(\theta_{0}^{3} + \theta_{0}^{2} + 4 \, \theta_{0}\right)} \overline{B_{1}}\bigg) C_{1}\\
& + 2 \, {\left(8 \, {\left(\theta_{0}^{3} \overline{\alpha_{1}} + \theta_{0}^{2} \overline{\alpha_{1}}\right)} A_{0} + {\left(\theta_{0}^{3} + \theta_{0}^{2} + 6 \, \theta_{0} + 6\right)} \overline{C_{1}}\right)} C_{2}\bigg\}\\
&\mu_{7}=\frac{1}{48 \, {\left(\theta_{0}^{3} - 9 \, \theta_{0}^{2} + 20 \, \theta_{0}\right)}}\bigg\{24 \, {\left(\theta_{0}^{2} \overline{\alpha_{3}} - 9 \, \theta_{0} \overline{\alpha_{3}} + 20 \, \overline{\alpha_{3}}\right)} A_{0} \overline{C_{1}} + 24 \, {\left(\theta_{0}^{2} - 9 \, \theta_{0} + 20\right)} A_{0} \overline{C_{4}} - 24 \, {\left(\theta_{0}^{2} - 5 \, \theta_{0}\right)} \overline{A_{1}} \overline{E_{1}}\\
& - 24 \, {\left(\theta_{0}^{2} - 4 \, \theta_{0}\right)} \overline{A_{0}} \overline{E_{3}} + {\left(24 \, {\left(\theta_{0}^{2} \overline{\alpha_{1}} - 9 \, \theta_{0} \overline{\alpha_{1}} + 20 \, \overline{\alpha_{1}}\right)} A_{0} + 5 \, {\left(\theta_{0}^{2} - 9 \, \theta_{0} + 20\right)} \overline{C_{1}}\right)} \overline{C_{2}}\bigg\}\\
&\mu_{8}=\frac{1}{2 \, {\left(\theta_{0}^{3} + 3 \, \theta_{0}^{2} + 2 \, \theta_{0}\right)}}\bigg\{2 \, {\left(\theta_{0}^{2} + \theta_{0}\right)} A_{0} {\left| A_{1} \right|}^{2} \overline{B_{1}} + {\left(\theta_{0}^{2} + 2 \, \theta_{0}\right)} A_{1} \overline{B_{3}} + {\left(\theta_{0}^{2} + \theta_{0}\right)} A_{0} \overline{B_{5}}\\
& + {\left(2 \, {\left(\theta_{0}^{2} + 2 \, \theta_{0}\right)} A_{1} {\left| A_{1} \right|}^{2} + {\left(\alpha_{5} \theta_{0}^{2} + \alpha_{5} \theta_{0}\right)} A_{0}\right)} \overline{C_{1}} + {\left({\left(\alpha_{1} \theta_{0}^{2} + \alpha_{1} \theta_{0}\right)} A_{0} + {\left(\theta_{0}^{2} + 3 \, \theta_{0} + 2\right)} A_{2}\right)} \overline{C_{2}}\bigg\}\\
&\mu_{9}=\frac{1}{2 \, {\left(\theta_{0}^{4} + 6 \, \theta_{0}^{3} + 11 \, \theta_{0}^{2} + 6 \, \theta_{0}\right)}}\bigg\{{\left(\theta_{0}^{3} \overline{\alpha_{2}} + 3 \, \theta_{0}^{2} \overline{\alpha_{2}} + 2 \, \theta_{0} \overline{\alpha_{2}}\right)} A_{0} C_{1} + {\left(\theta_{0}^{3} + 4 \, \theta_{0}^{2} + 3 \, \theta_{0}\right)} A_{1} \overline{B_{2}} + {\left(\theta_{0}^{3} + 3 \, \theta_{0}^{2} + 2 \, \theta_{0}\right)} A_{0} \overline{B_{4}}\\
& + {\left({\left(\alpha_{1} \theta_{0}^{3} + 3 \, \alpha_{1} \theta_{0}^{2} + 2 \, \alpha_{1} \theta_{0}\right)} A_{0} + {\left(\theta_{0}^{3} + 5 \, \theta_{0}^{2} + 6 \, \theta_{0}\right)} A_{2}\right)} \overline{B_{1}} + \bigg({\left(\alpha_{3} \theta_{0}^{3} + 3 \, \alpha_{3} \theta_{0}^{2} + 2 \, \alpha_{3} \theta_{0}\right)} A_{0}\\
& + {\left(\alpha_{1} \theta_{0}^{3} + 4 \, \alpha_{1} \theta_{0}^{2} + 3 \, \alpha_{1} \theta_{0}\right)} A_{1} + {\left(\theta_{0}^{3} + 6 \, \theta_{0}^{2} + 11 \, \theta_{0} + 6\right)} A_{3}\bigg) \overline{C_{1}}\bigg\}\\
&\mu_{10}=\frac{1}{16 \, {\left(\theta_{0}^{3} - 3 \, \theta_{0}^{2} - 4 \, \theta_{0}\right)}}\bigg\{16 \, {\left(\theta_{0}^{2} - 4 \, \theta_{0}\right)} A_{0} {\left| A_{1} \right|}^{2} \overline{C_{2}} + 8 \, {\left(\theta_{0}^{2} \overline{\alpha_{1}} - 4 \, \theta_{0} \overline{\alpha_{1}}\right)} A_{0} \overline{B_{1}} + 8 \, {\left(\theta_{0}^{2} - 4 \, \theta_{0}\right)} A_{0} \overline{B_{6}}\\
& + 8 \, {\left(\theta_{0}^{2} - 3 \, \theta_{0} - 4\right)} A_{1} \overline{C_{3}} - 8 \, {\left(\theta_{0}^{2} + \theta_{0}\right)} \overline{A_{0}} \overline{E_{2}} + \bigg(8 \, {\left(\theta_{0}^{2} \overline{\alpha_{5}} - 4 \, \theta_{0} \overline{\alpha_{5}}\right)} A_{0} + 8 \, {\left(\theta_{0}^{2} \overline{\alpha_{1}} - 3 \, \theta_{0} \overline{\alpha_{1}} - 4 \, \overline{\alpha_{1}}\right)} A_{1}\\
& - 8 \, {\left(\theta_{0}^{2} \overline{\alpha_{2}} + \theta_{0} \overline{\alpha_{2}}\right)} \overline{A_{0}} + {\left(2 \, \theta_{0}^{2} - 7 \, \theta_{0} - 4\right)} \overline{B_{1}}\bigg) \overline{C_{1}}\bigg\}
\end{align*}
The precise expression of these $\mu$ is irrelevant.
The only strange new powers of order $2\theta_0+3$ are
\begin{align*}
	\begin{dmatrix}
		\frac{2 \, E_{2} \theta_{0} \overline{A_{0}} + E_{1} {\left(2 \, \theta_{0} + 1\right)} \overline{A_{1}}}{4 \, {\left(2 \, \theta_{0}^{2} + \theta_{0}\right)}} & -\theta_{0} + 3 & 3 \, \theta_{0} \\
	\frac{E_{3} \overline{A_{0}}}{4 \, \theta_{0}} & -\theta_{0} + 4 & 3 \, \theta_{0} - 1 \\
	\end{dmatrix}
\end{align*}
and their conjugate.
As $\vec{E}_1\in\mathrm{Span}(\bar{\vec{A}_0})$, $\vec{E}_2\in\mathrm{Span}(\bar{\vec{A}_1},\vec{C}_1)$ and $\vec{E}_3\in\mathrm{Span}(\bar{\vec{A}_0})$ , they vanish and there exists $
\alpha_{10}, \alpha_{11}, \alpha_{12}, \alpha_{13}, \alpha_{14}, \alpha_{15}, \alpha_{16}, \in\mathbb{C}$ such that
\begin{align}\label{endmetric}
	&e^{2\lambda}=|z|^{2\theta_0-2}+2|\vec{A}_1|^2|z|^{2\theta_0}+\beta|z|^{2\theta_0+2}+2\,\Re\bigg(\alpha_0z^{\theta_0}\z^{\theta_0-1}+\alpha_1z^{\theta_0+1}\z^{\theta_0-2}+\alpha_2z\z^{2\theta_0}+\alpha_3z^{\theta_0+2}\z^{\theta_0-1}+\alpha_4z^{\theta_0+3}\z^{\theta_0-1}\nonumber\\
	&+\alpha_5z^{\theta_0+1}\z^{\theta_0}+\alpha_6z^{\theta_0+2}\z^{\theta_0}+\alpha_7z^3\z^{2\theta_0-1}+\alpha_8z^2\z^{2\theta_0}+\alpha_9z\z^{2\theta_0+1}\nonumber\\
	&+\alpha_{10}z\z^{2\theta_0+2}+\alpha_{11}z^2\z^{2\theta_0+1}+\alpha_{12}z^{3}\z^{2\theta_0}+\alpha_{13}z^4\z^{2\theta_0-1}+\alpha_{14}z^{\theta_0-1}\z^{\theta_0+4}+\alpha_{15}z^{\theta_0}\z^{\theta_0+3}+\alpha_{16}z^{\theta_0+1}\z^{\theta_0+2}\bigg)+O(|z|^{2\theta_0+4-\epsilon})
\end{align}
which gives in Sage
	\begin{align*}
	e^{2\lambda}=\begin{dmatrix}
	1 & \theta_{0} - 1 & \theta_{0} - 1 \\
	2 \, {\left| A_{1} \right|}^{2} & \theta_{0} & \theta_{0} \\
	\beta & \theta_{0} + 1 & \theta_{0} + 1 \\
	\alpha_{1} & \theta_{0} + 1 & \theta_{0} - 1 \\
	\alpha_{2} & 1 & 2 \, \theta_{0} \\
	\alpha_{3} & \theta_{0} + 2 & \theta_{0} - 1 \\
	\alpha_{4} & \theta_{0} + 3 & \theta_{0} - 1 \\
	\alpha_{5} & \theta_{0} + 1 & \theta_{0} \\
	\alpha_{6} & \theta_{0} + 2 & \theta_{0} \\
	\alpha_{7} & 3 & 2 \, \theta_{0} - 1 \\
	\alpha_{8} & 2 & 2 \, \theta_{0} \\
	\alpha_{9} & 1 & 2 \, \theta_{0} + 1 \\
	\alpha_{10} & 1 & 2 \, \theta_{0} + 2 \\
	\alpha_{11} & 2 & 2 \, \theta_{0} + 1 \\
	\alpha_{12} & 3 & 2 \, \theta_{0} \\
	\alpha_{13} & 4 & 2 \, \theta_{0} - 1 \\
	\alpha_{14} & \theta_{0} - 1 & \theta_{0} + 4 \\
	\alpha_{15} & \theta_{0} & \theta_{0} + 3 \\
	\alpha_{16} & \theta_{0} + 1 & \theta_{0} + 2 
	\end{dmatrix}
	\begin{dmatrix}
	\overline{\alpha_{1}} & \theta_{0} - 1 & \theta_{0} + 1 \\
	\overline{\alpha_{2}} & 2 \, \theta_{0} & 1 \\
	\overline{\alpha_{3}} & \theta_{0} - 1 & \theta_{0} + 2 \\
	\overline{\alpha_{4}} & \theta_{0} - 1 & \theta_{0} + 3 \\
	\overline{\alpha_{5}} & \theta_{0} & \theta_{0} + 1 \\
	\overline{\alpha_{6}} & \theta_{0} & \theta_{0} + 2 \\
	\overline{\alpha_{7}} & 2 \, \theta_{0} - 1 & 3 \\
	\overline{\alpha_{8}} & 2 \, \theta_{0} & 2 \\
	\overline{\alpha_{9}} & 2 \, \theta_{0} + 1 & 1 \\
	\overline{\alpha_{10}} & 2 \, \theta_{0} + 2 & 1 \\
	\overline{\alpha_{11}} & 2 \, \theta_{0} + 1 & 2 \\
	\overline{\alpha_{12}} & 2 \, \theta_{0} & 3 \\
	\overline{\alpha_{13}} & 2 \, \theta_{0} - 1 & 4 \\
	\overline{\alpha_{14}} & \theta_{0} + 4 & \theta_{0} - 1 \\
	\overline{\alpha_{15}} & \theta_{0} + 3 & \theta_{0} \\
	\overline{\alpha_{16}} & \theta_{0} + 2 & \theta_{0} + 1
	\end{dmatrix}
	\end{align*}

\begin{align*}
	\h_0=\begin{dmatrix}
	-\frac{\theta_{0} - 2}{2 \, \theta_{0}} & C_{1} & 0 & \theta_{0} \\
	-\frac{\theta_{0} - 2}{2 \, {\left(\theta_{0} + 1\right)}} & B_{1} & 0 & \theta_{0} + 1 \\
	2 \, \alpha_{2} \theta_{0} - 4 \, \alpha_{2} & A_{0} & 0 & \theta_{0} + 1 \\
	-\frac{\theta_{0} - 2}{2 \, {\left(\theta_{0} + 2\right)}} & B_{2} & 0 & \theta_{0} + 2 \\
	-\frac{\theta_{0} \overline{\alpha_{1}} - 2 \, \overline{\alpha_{1}}}{2 \, {\left(\theta_{0} + 2\right)}} & C_{1} & 0 & \theta_{0} + 2 \\
	2 \, \alpha_{9} \theta_{0} - 4 \, \alpha_{9} & A_{0} & 0 & \theta_{0} + 2 \\
	-\frac{\theta_{0} - 2}{2 \, {\left(\theta_{0} + 3\right)}} & B_{4} & 0 & \theta_{0} + 3 \\
	-\frac{\theta_{0} \overline{\alpha_{3}} - 2 \, \overline{\alpha_{3}}}{2 \, {\left(\theta_{0} + 3\right)}} & C_{1} & 0 & \theta_{0} + 3 \\
	-\frac{\theta_{0} \overline{\alpha_{1}} - 2 \, \overline{\alpha_{1}}}{2 \, {\left(\theta_{0} + 3\right)}} & B_{1} & 0 & \theta_{0} + 3 \\
	\frac{\alpha_{2} \theta_{0}^{2} - \alpha_{2} \theta_{0} - 2 \, \alpha_{2}}{4 \, {\left(\theta_{0} + 3\right)}} & \overline{C_{1}} & 0 & \theta_{0} + 3 \\
	-2 \, {\left(\alpha_{2} \overline{\alpha_{1}} - \alpha_{10}\right)} \theta_{0} + 4 \, \alpha_{2} \overline{\alpha_{1}} - 4 \, \alpha_{10} & A_{0} & 0 & \theta_{0} + 3 \\
	-\frac{\theta_{0} - 3}{2 \, \theta_{0}} & C_{2} & 1 & \theta_{0} \\
	-\frac{\theta_{0} - 3}{2 \, {\left(\theta_{0} + 1\right)}} & B_{3} & 1 & \theta_{0} + 1 \\
	-\frac{{\left(\theta_{0}^{2} - 2 \, \theta_{0} + 1\right)} {\left| A_{1} \right|}^{2}}{\theta_{0}^{2} + \theta_{0}} & C_{1} & 1 & \theta_{0} + 1 \\
	2 \, \alpha_{8} \theta_{0} - 6 \, \alpha_{8} & A_{0} & 1 & \theta_{0} + 1 \\
	2 \, \alpha_{2} \theta_{0} - 4 \, \alpha_{2} & A_{1} & 1 & \theta_{0} + 1 \\
	-\frac{\theta_{0} - 3}{2 \, {\left(\theta_{0} + 2\right)}} & B_{5} & 1 & \theta_{0} + 2 \\
	-\frac{\theta_{0} \overline{\alpha_{1}} - 3 \, \overline{\alpha_{1}}}{2 \, {\left(\theta_{0} + 2\right)}} & C_{2} & 1 & \theta_{0} + 2 \\
	-\frac{\theta_{0}^{2} \overline{\alpha_{5}} - 2 \, \theta_{0} \overline{\alpha_{5}} + 2 \, \overline{\alpha_{5}}}{2 \, {\left(\theta_{0}^{2} + 2 \, \theta_{0}\right)}} & C_{1} & 1 & \theta_{0} + 2 \\
	-\frac{{\left(\theta_{0}^{2} - \theta_{0} - 1\right)} {\left| A_{1} \right|}^{2}}{\theta_{0}^{2} + 3 \, \theta_{0} + 2} & B_{1} & 1 & \theta_{0} + 2 \\
	2 \, \alpha_{9} \theta_{0} - 4 \, \alpha_{9} & A_{1} & 1 & \theta_{0} + 2 \\
	-4 \, {\left(\alpha_{2} \theta_{0} - 3 \, \alpha_{2}\right)} {\left| A_{1} \right|}^{2} + 2 \, \alpha_{11} \theta_{0} - 6 \, \alpha_{11} & A_{0} & 1 & \theta_{0} + 2 \\
	-\frac{\theta_{0} - 4}{2 \, \theta_{0}} & C_{3} & 2 & \theta_{0} \\
	-\frac{\alpha_{1} \theta_{0} - 2 \, \alpha_{1}}{2 \, \theta_{0}} & C_{1} & 2 & \theta_{0} \\
	2 \, \alpha_{7} \theta_{0} - 8 \, \alpha_{7} & A_{0} & 2 & \theta_{0} 
			\end{dmatrix}
	\end{align*}
	\begin{align*}
	\begin{dmatrix}
	-\frac{\theta_{0} - 4}{2 \, {\left(\theta_{0} + 1\right)}} & B_{6} & 2 & \theta_{0} + 1 \\
	-\frac{\alpha_{1} \theta_{0} - 2 \, \alpha_{1}}{2 \, {\left(\theta_{0} + 1\right)}} & B_{1} & 2 & \theta_{0} + 1 \\
	-\frac{\alpha_{5} \theta_{0}^{2} - 2 \, \alpha_{5} \theta_{0} + 2 \, \alpha_{5}}{2 \, {\left(\theta_{0}^{2} + \theta_{0}\right)}} & C_{1} & 2 & \theta_{0} + 1 \\
	-\frac{{\left(\theta_{0}^{2} - 3 \, \theta_{0} + 1\right)} {\left| A_{1} \right|}^{2}}{\theta_{0}^{2} + \theta_{0}} & C_{2} & 2 & \theta_{0} + 1 \\
	2 \, \alpha_{8} \theta_{0} - 6 \, \alpha_{8} & A_{1} & 2 & \theta_{0} + 1 \\
	2 \, \alpha_{2} \theta_{0} - 4 \, \alpha_{2} & A_{2} & 2 & \theta_{0} + 1 \\
	8 \, \alpha_{1} \alpha_{2} - 2 \, {\left(\alpha_{1} \alpha_{2} - \alpha_{12}\right)} \theta_{0} - 8 \, \alpha_{12} & A_{0} & 2 & \theta_{0} + 1 \\
	-\frac{\theta_{0} - 5}{2 \, \theta_{0}} & C_{4} & 3 & \theta_{0} \\
	-\frac{\alpha_{3} \theta_{0} - 2 \, \alpha_{3}}{2 \, \theta_{0}} & C_{1} & 3 & \theta_{0} \\
	-\frac{\alpha_{1} \theta_{0} - 3 \, \alpha_{1}}{2 \, \theta_{0}} & C_{2} & 3 & \theta_{0} \\
	2 \, \alpha_{7} \theta_{0} - 8 \, \alpha_{7} & A_{1} & 3 & \theta_{0} \\
	2 \, \alpha_{13} \theta_{0} - 10 \, \alpha_{13} & A_{0} & 3 & \theta_{0} \\
	-\frac{\theta_{0} - 2}{2 \, \theta_{0} + 1} & E_{2} & -\theta_{0} + 2 & 2 \, \theta_{0} + 1 \\
	\frac{\alpha_{2} \theta_{0} - 2 \, \alpha_{2}}{2 \, {\left(2 \, \theta_{0}^{2} + \theta_{0}\right)}} & C_{1} & -\theta_{0} + 2 & 2 \, \theta_{0} + 1 \\
	-\frac{\theta_{0} - 2}{2 \, \theta_{0}} & E_{1} & -\theta_{0} + 2 & 2 \, \theta_{0} \\
	-\frac{2 \, \theta_{0} - 5}{4 \, \theta_{0}} & E_{3} & -\theta_{0} + 3 & 2 \, \theta_{0} \\
	4 & A_{2} & \theta_{0} & 0 \\
	-4 \, \alpha_{1} & A_{0} & \theta_{0} & 0 \\
	-4 \, \alpha_{5} & A_{0} & \theta_{0} & 1 \\
	-4 \, {\left| A_{1} \right|}^{2} & A_{1} & \theta_{0} & 1 \\
	\frac{1}{2} & \overline{B_{2}} & \theta_{0} & 2 \\
	-2 \, \overline{\alpha_{5}} & A_{1} & \theta_{0} & 2 \\
	8 \, {\left| A_{1} \right|}^{4} + 4 \, \alpha_{1} \overline{\alpha_{1}} - 4 \, \beta & A_{0} & \theta_{0} & 2 \\
	\frac{1}{3} & \overline{B_{5}} & \theta_{0} & 3 \\
	-\frac{1}{6} \, \alpha_{5} & \overline{C_{1}} & \theta_{0} & 3 \\
	\frac{1}{6} \, {\left| A_{1} \right|}^{2} & \overline{B_{1}} & \theta_{0} & 3 \\
	8 \, {\left| A_{1} \right|}^{2} \overline{\alpha_{5}} + 4 \, \alpha_{5} \overline{\alpha_{1}} + 4 \, \alpha_{1} \overline{\alpha_{3}} - 4 \, \alpha_{16} & A_{0} & \theta_{0} & 3 \\
	4 \, {\left| A_{1} \right|}^{2} \overline{\alpha_{1}} - 2 \, \overline{\alpha_{6}} & A_{1} & \theta_{0} & 3 
			\end{dmatrix}
	\end{align*}
	\begin{align*}
	\begin{dmatrix}
	2 & A_{1} & \theta_{0} - 1 & 0 \\
	-4 \, {\left| A_{1} \right|}^{2} & A_{0} & \theta_{0} - 1 & 1 \\
	\frac{1}{4} & \overline{B_{1}} & \theta_{0} - 1 & 2 \\
	-2 \, \overline{\alpha_{5}} & A_{0} & \theta_{0} - 1 & 2 \\
	\frac{1}{6} & \overline{B_{3}} & \theta_{0} - 1 & 3 \\
	-\frac{1}{6} \, {\left| A_{1} \right|}^{2} & \overline{C_{1}} & \theta_{0} - 1 & 3 \\
	4 \, {\left| A_{1} \right|}^{2} \overline{\alpha_{1}} - 2 \, \overline{\alpha_{6}} & A_{0} & \theta_{0} - 1 & 3 \\
	\frac{1}{8} & \overline{B_{6}} & \theta_{0} - 1 & 4 \\
	-\frac{1}{8} \, \overline{\alpha_{5}} & \overline{C_{1}} & \theta_{0} - 1 & 4 \\
	\frac{1}{8} \, \overline{\alpha_{1}} & \overline{B_{1}} & \theta_{0} - 1 & 4 \\
	-\frac{1}{12} \, {\left| A_{1} \right|}^{2} & \overline{C_{2}} & \theta_{0} - 1 & 4 \\
	4 \, {\left| A_{1} \right|}^{2} \overline{\alpha_{3}} + 2 \, \overline{\alpha_{1}} \overline{\alpha_{5}} - 2 \, \alpha_{15} & A_{0} & \theta_{0} - 1 & 4 \\
	6 & A_{3} & \theta_{0} + 1 & 0 \\
	-6 \, \alpha_{3} & A_{0} & \theta_{0} + 1 & 0 \\
	-4 \, \alpha_{1} & A_{1} & \theta_{0} + 1 & 0 \\
	-4 \, \alpha_{5} & A_{1} & \theta_{0} + 1 & 1 \\
	-4 \, {\left| A_{1} \right|}^{2} & A_{2} & \theta_{0} + 1 & 1 \\
	12 \, \alpha_{1} {\left| A_{1} \right|}^{2} - 6 \, \alpha_{6} & A_{0} & \theta_{0} + 1 & 1 \\
	\frac{3}{4} & \overline{B_{4}} & \theta_{0} + 1 & 2 \\
	\frac{\theta_{0} \overline{\alpha_{2}} - 2 \, \overline{\alpha_{2}}}{4 \, \theta_{0}} & C_{1} & \theta_{0} + 1 & 2 \\
	-2 \, \overline{\alpha_{5}} & A_{2} & \theta_{0} + 1 & 2 \\
	\frac{1}{4} \, \alpha_{1} & \overline{B_{1}} & \theta_{0} + 1 & 2 \\
	12 \, \alpha_{5} {\left| A_{1} \right|}^{2} + 6 \, \alpha_{3} \overline{\alpha_{1}} + 6 \, \alpha_{1} \overline{\alpha_{5}} - 6 \, \overline{\alpha_{16}} & A_{0} & \theta_{0} + 1 & 2 \\
	8 \, {\left| A_{1} \right|}^{4} + 4 \, \alpha_{1} \overline{\alpha_{1}} - 4 \, \beta & A_{1} & \theta_{0} + 1 & 2 \\
	8 & A_{4} & \theta_{0} + 2 & 0 \\
	-6 \, \alpha_{3} & A_{1} & \theta_{0} + 2 & 0 \\
	-4 \, \alpha_{1} & A_{2} & \theta_{0} + 2 & 0 \\
	4 \, \alpha_{1}^{2} - 8 \, \alpha_{4} & A_{0} & \theta_{0} + 2 & 0 
			\end{dmatrix}
	\end{align*}
	\begin{align*}
	\begin{dmatrix}
	-4 \, \alpha_{5} & A_{2} & \theta_{0} + 2 & 1 \\
	-4 \, {\left| A_{1} \right|}^{2} & A_{3} & \theta_{0} + 2 & 1 \\
	16 \, \alpha_{3} {\left| A_{1} \right|}^{2} + 8 \, \alpha_{1} \alpha_{5} - 8 \, \overline{\alpha_{15}} & A_{0} & \theta_{0} + 2 & 1 \\
	12 \, \alpha_{1} {\left| A_{1} \right|}^{2} - 6 \, \alpha_{6} & A_{1} & \theta_{0} + 2 & 1 \\
	10 & A_{5} & \theta_{0} + 3 & 0 \\
	-6 \, \alpha_{3} & A_{2} & \theta_{0} + 3 & 0 \\
	-4 \, \alpha_{1} & A_{3} & \theta_{0} + 3 & 0 \\
	10 \, \alpha_{1} \alpha_{3} - 10 \, \overline{\alpha_{14}} & A_{0} & \theta_{0} + 3 & 0 \\
	4 \, \alpha_{1}^{2} - 8 \, \alpha_{4} & A_{1} & \theta_{0} + 3 & 0 \\
	-\frac{\theta_{0}}{2 \, {\left(\theta_{0} - 4\right)}} & \overline{E_{1}} & 2 \, \theta_{0} - 2 & -\theta_{0} + 4 \\
	-2 \, \theta_{0} \overline{\alpha_{7}} & A_{0} & 2 \, \theta_{0} - 2 & -\theta_{0} + 4 \\
	-\frac{\theta_{0}}{2 \, {\left(\theta_{0} - 5\right)}} & \overline{E_{3}} & 2 \, \theta_{0} - 2 & -\theta_{0} + 5 \\
	-2 \, \theta_{0} \overline{\alpha_{13}} & A_{0} & 2 \, \theta_{0} - 2 & -\theta_{0} + 5 \\
	-2 \, \theta_{0} \overline{\alpha_{2}} - 2 \, \overline{\alpha_{2}} & A_{0} & 2 \, \theta_{0} - 1 & -\theta_{0} + 2 \\
	-2 \, \theta_{0} \overline{\alpha_{8}} - 2 \, \overline{\alpha_{8}} & A_{0} & 2 \, \theta_{0} - 1 & -\theta_{0} + 3 \\
	-\frac{\theta_{0} + 1}{2 \, {\left(\theta_{0} - 4\right)}} & \overline{E_{2}} & 2 \, \theta_{0} - 1 & -\theta_{0} + 4 \\
	-\frac{\theta_{0}^{2} \overline{\alpha_{2}} - \theta_{0} \overline{\alpha_{2}} - 2 \, \overline{\alpha_{2}}}{4 \, {\left(\theta_{0} - 4\right)}} & \overline{C_{1}} & 2 \, \theta_{0} - 1 & -\theta_{0} + 4 \\
	-2 \, \theta_{0} \overline{\alpha_{7}} & A_{1} & 2 \, \theta_{0} - 1 & -\theta_{0} + 4 \\
	2 \, {\left(\overline{\alpha_{1}} \overline{\alpha_{2}} - \overline{\alpha_{12}}\right)} \theta_{0} + 2 \, \overline{\alpha_{1}} \overline{\alpha_{2}} - 2 \, \overline{\alpha_{12}} & A_{0} & 2 \, \theta_{0} - 1 & -\theta_{0} + 4 \\
	-2 \, \theta_{0} \overline{\alpha_{9}} - 4 \, \overline{\alpha_{9}} & A_{1} & 2 \, \theta_{0} + 1 & -\theta_{0} + 2 \\
	-2 \, \theta_{0} \overline{\alpha_{2}} - 2 \, \overline{\alpha_{2}} & A_{2} & 2 \, \theta_{0} + 1 & -\theta_{0} + 2 \\
	2 \, {\left(\alpha_{1} \overline{\alpha_{2}} - \overline{\alpha_{10}}\right)} \theta_{0} + 6 \, \alpha_{1} \overline{\alpha_{2}} - 6 \, \overline{\alpha_{10}} & A_{0} & 2 \, \theta_{0} + 1 & -\theta_{0} + 2 \\
	-2 \, \theta_{0} \overline{\alpha_{9}} - 4 \, \overline{\alpha_{9}} & A_{0} & 2 \, \theta_{0} & -\theta_{0} + 2 \\
	-2 \, \theta_{0} \overline{\alpha_{2}} - 2 \, \overline{\alpha_{2}} & A_{1} & 2 \, \theta_{0} & -\theta_{0} + 2 \\
	-2 \, \theta_{0} \overline{\alpha_{8}} - 2 \, \overline{\alpha_{8}} & A_{1} & 2 \, \theta_{0} & -\theta_{0} + 3 \\
	4 \, {\left(\theta_{0} \overline{\alpha_{2}} + 2 \, \overline{\alpha_{2}}\right)} {\left| A_{1} \right|}^{2} - 2 \, \theta_{0} \overline{\alpha_{11}} - 4 \, \overline{\alpha_{11}} & A_{0} & 2 \, \theta_{0} & -\theta_{0} + 3
	\end{dmatrix}
	\end{align*}

	\newpage
	
	\section{Order $2$ terms in the quartic form}
	\normalsize
	We have
	\begin{align*}
		&g_0^{-1}\otimes\vec{Q}(\h_0)|_{\theta_0,2-\theta_0}=
		\frac{2}{\theta_{0}} \, \bigg\{32 \, A_{1}^{2} \theta_{0} {\left| A_{1} \right|}^{4} + 48 \, A_{0} A_{1} \alpha_{5} \theta_{0} {\left| A_{1} \right|}^{2} - 24 \, A_{0}^{2} \alpha_{3} \theta_{0} \overline{\alpha_{5}} + 8 \, {\left(3 \, \alpha_{3} \overline{\alpha_{1}} + \alpha_{1} \overline{\alpha_{5}} - 3 \, \overline{\alpha_{16}}\right)} A_{0} A_{1} \theta_{0}\\
		& + 16 \, {\left(\alpha_{1} \overline{\alpha_{1}} - \beta\right)} A_{1}^{2} \theta_{0} - 8 \, A_{1} A_{2} \theta_{0} \overline{\alpha_{5}} + 24 \, A_{0} A_{3} \theta_{0} \overline{\alpha_{5}} - 2 \, {\left(\theta_{0}^{4} \overline{\alpha_{2}} - 4 \, \theta_{0}^{3} \overline{\alpha_{2}} + 2 \, \theta_{0}^{2} \overline{\alpha_{2}} + 4 \, \theta_{0} \overline{\alpha_{2}} - 3 \, \overline{\alpha_{2}}\right)} A_{0} C_{2}\\
		& + 3 \, A_{1} \theta_{0} \overline{B_{4}} - {\left(2 \, {\left(\theta_{0}^{4} \overline{\alpha_{9}} - \theta_{0}^{3} \overline{\alpha_{9}} - 4 \, \theta_{0}^{2} \overline{\alpha_{9}} + 4 \, \theta_{0} \overline{\alpha_{9}}\right)} A_{0} + {\left(2 \, \theta_{0}^{4} \overline{\alpha_{2}} - 4 \, \theta_{0}^{3} \overline{\alpha_{2}} - 2 \, \theta_{0}^{2} \overline{\alpha_{2}} + 3 \, \theta_{0} \overline{\alpha_{2}} + 2 \, \overline{\alpha_{2}}\right)} A_{1}\right)} C_{1}\\
		& + 3 \, {\left(A_{1} \alpha_{1} \theta_{0} + A_{0} \alpha_{3} \theta_{0} - A_{3} \theta_{0}\right)} \overline{B_{1}}\bigg\}
	\end{align*}
	Then
	\begin{align*}
	g_0^{-1}\otimes Q(\h_0)=
	\begin{dmatrix}
		 {\left(\theta_{0}^{2} - 5 \, \theta_{0} + 6\right)} A_{1} C_{2} - 2 \, {\left({\left(\alpha_{1} \theta_{0}^{2} - 2 \, \alpha_{1} \theta_{0}\right)} A_{0} - {\left(\theta_{0}^{2} - 2 \, \theta_{0}\right)} A_{2}\right)} C_{1} & 0 & 0 \\
		16 \, A_{0} A_{2} {\left| A_{1} \right|}^{2} - 8 \, A_{0} A_{1} \alpha_{5} - 8 \, {\left(2 \, A_{0}^{2} \alpha_{1} + A_{1}^{2}\right)} {\left| A_{1} \right|}^{2} & \theta_{0} - 1 & -\theta_{0} + 1 \\
		\omega_{-1} & -1 & 1 \\
		4 \, {\left(\theta_{0}^{3} \overline{\alpha_{2}} - \theta_{0}^{2} \overline{\alpha_{2}} - 2 \, \theta_{0} \overline{\alpha_{2}}\right)} A_{0} A_{1} & 2 \, \theta_{0} - 2 & -2 \, \theta_{0} + 2 \\
		{\left(\theta_{0}^{2} - 3 \, \theta_{0} + 2\right)} A_{1} C_{1} & -1 & 0
				\end{dmatrix}
		\end{align*}
		where
		\begin{align*}
			\omega_{-1}=-\frac{2 \, {\left(\theta_{0}^{3} - 4 \, \theta_{0}^{2} + 5 \, \theta_{0} - 2\right)} A_{0} C_{1} {\left| A_{1} \right|}^{2} + 4 \, {\left(\alpha_{2} \theta_{0}^{4} - 2 \, \alpha_{2} \theta_{0}^{3} - \alpha_{2} \theta_{0}^{2} + 2 \, \alpha_{2} \theta_{0}\right)} A_{0} A_{1} - {\left(\theta_{0}^{3} - 3 \, \theta_{0}^{2} + 2 \, \theta_{0}\right)} A_{1} B_{1}}{\theta_{0}}=0
		\end{align*}
		\begin{align*}
		e^{-2u}-1=\begin{dmatrix}
		-2 \, {\left| A_{1} \right|}^{2} & 1 & 1 \\
		-\alpha_{1} & 2 & 0 \\
		-\alpha_{2} & -\theta_{0} + 2 & \theta_{0} + 1 \\
		-\alpha_{3} & 3 & 0 \\
		-\alpha_{5} & 2 & 1 \\
		-\overline{\alpha_{1}} & 0 & 2 \\
		-\overline{\alpha_{2}} & \theta_{0} + 1 & -\theta_{0} + 2 \\
		-\overline{\alpha_{3}} & 0 & 3 \\
		-\overline{\alpha_{5}} & 1 & 2
				\end{dmatrix},\quad
		e^{-2\lambda}=\begin{dmatrix}
		1 & -\theta_{0} + 1 & -\theta_{0} + 1 \\
		-2 \, {\left| A_{1} \right|}^{2} & -\theta_{0} + 2 & -\theta_{0} + 2 \\
		-\alpha_{1} & -\theta_{0} + 3 & -\theta_{0} + 1 \\
		-\alpha_{2} & -2 \, \theta_{0} + 3 & 2 \\
		-\alpha_{3} & -\theta_{0} + 4 & -\theta_{0} + 1 \\
		-\alpha_{5} & -\theta_{0} + 3 & -\theta_{0} + 2 \\
		-\overline{\alpha_{1}} & -\theta_{0} + 1 & -\theta_{0} + 3 \\
		-\overline{\alpha_{2}} & 2 & -2 \, \theta_{0} + 3 \\
		-\overline{\alpha_{3}} & -\theta_{0} + 1 & -\theta_{0} + 4 \\
		-\overline{\alpha_{5}} & -\theta_{0} + 2 & -\theta_{0} + 3
		\end{dmatrix}
		\end{align*}
		\small
		\begin{align*}
		&(g^{-1}-g_0^{-1})\otimes Q(\h_0)=\\
&\begin{dmatrix}
		\omega_1 & 1 & 1 \\
		-32 \, A_{0} A_{2} {\left| A_{1} \right|}^{4} + 16 \, A_{0} A_{1} \alpha_{5} {\left| A_{1} \right|}^{2} + 16 \, {\left(2 \, A_{0}^{2} \alpha_{1} + A_{1}^{2}\right)} {\left| A_{1} \right|}^{4} - {\left(\theta_{0}^{2} \overline{\alpha_{2}} - 3 \, \theta_{0} \overline{\alpha_{2}} + 2 \, \overline{\alpha_{2}}\right)} A_{1} C_{1} & \theta_{0} & -\theta_{0} + 2 \\
		\omega_2 & 0 & 2 \\
		-8 \, {\left(\theta_{0}^{3} \overline{\alpha_{2}} - \theta_{0}^{2} \overline{\alpha_{2}} - 2 \, \theta_{0} \overline{\alpha_{2}}\right)} A_{0} A_{1} {\left| A_{1} \right|}^{2} & 2 \, \theta_{0} - 1 & -2 \, \theta_{0} + 3 \\
		-2 \, {\left(\theta_{0}^{2} - 3 \, \theta_{0} + 2\right)} A_{1} C_{1} {\left| A_{1} \right|}^{2} & 0 & 1 \\
		\omega_3 & 2 & 0 \\
		-16 \, A_{0} A_{2} \alpha_{1} {\left| A_{1} \right|}^{2} + 8 \, A_{0} A_{1} \alpha_{1} \alpha_{5} + 8 \, {\left(2 \, A_{0}^{2} \alpha_{1}^{2} + A_{1}^{2} \alpha_{1}\right)} {\left| A_{1} \right|}^{2} & \theta_{0} + 1 & -\theta_{0} + 1 \\
		-4 \, {\left(\alpha_{1} \theta_{0}^{3} \overline{\alpha_{2}} - \alpha_{1} \theta_{0}^{2} \overline{\alpha_{2}} - 2 \, \alpha_{1} \theta_{0} \overline{\alpha_{2}}\right)} A_{0} A_{1} & 2 \, \theta_{0} & -2 \, \theta_{0} + 2 \\
		-{\left(\alpha_{1} \theta_{0}^{2} - 3 \, \alpha_{1} \theta_{0} + 2 \, \alpha_{1}\right)} A_{1} C_{1} & 1 & 0 \\
		-{\left(\alpha_{2} \theta_{0}^{2} - 3 \, \alpha_{2} \theta_{0} + 2 \, \alpha_{2}\right)} A_{1} C_{1} & -\theta_{0} + 1 & \theta_{0} + 1 \\
		-16 \, A_{0} A_{2} {\left| A_{1} \right|}^{2} \overline{\alpha_{1}} + 8 \, A_{0} A_{1} \alpha_{5} \overline{\alpha_{1}} + 8 \, {\left(2 \, A_{0}^{2} \alpha_{1} \overline{\alpha_{1}} + A_{1}^{2} \overline{\alpha_{1}}\right)} {\left| A_{1} \right|}^{2} & \theta_{0} - 1 & -\theta_{0} + 3 \\
		\omega_3 & -1 & 3 \\
		-4 \, {\left(\theta_{0}^{3} \overline{\alpha_{1}} \overline{\alpha_{2}} - \theta_{0}^{2} \overline{\alpha_{1}} \overline{\alpha_{2}} - 2 \, \theta_{0} \overline{\alpha_{1}} \overline{\alpha_{2}}\right)} A_{0} A_{1} & 2 \, \theta_{0} - 2 & -2 \, \theta_{0} + 4 \\
		-{\left(\theta_{0}^{2} \overline{\alpha_{1}} - 3 \, \theta_{0} \overline{\alpha_{1}} + 2 \, \overline{\alpha_{1}}\right)} A_{1} C_{1} & -1 & 2
		\end{dmatrix}
	\end{align*}
	where
	\begin{align*}
		\omega_1&=-\frac{1}{\theta_{0}}\Big\{2 \, {\left(\theta_{0}^{3} - 5 \, \theta_{0}^{2} + 6 \, \theta_{0}\right)} A_{1} C_{2} {\left| A_{1} \right|}^{2} - 4 \, {\left(\alpha_{1} \alpha_{2} \theta_{0}^{4} - 2 \, \alpha_{1} \alpha_{2} \theta_{0}^{3} - \alpha_{1} \alpha_{2} \theta_{0}^{2} + 2 \, \alpha_{1} \alpha_{2} \theta_{0}\right)} A_{0} A_{1}\\
		& + {\left(\alpha_{1} \theta_{0}^{3} - 3 \, \alpha_{1} \theta_{0}^{2} + 2 \, \alpha_{1} \theta_{0}\right)} A_{1} B_{1}\\
		& - \Big\{2 \, {\left(3 \, \alpha_{1} \theta_{0}^{3} - 8 \, \alpha_{1} \theta_{0}^{2} + 5 \, \alpha_{1} \theta_{0} - 2 \, \alpha_{1}\right)} A_{0} {\left| A_{1} \right|}^{2} - 4 \, {\left(\theta_{0}^{3} - 2 \, \theta_{0}^{2}\right)} A_{2} {\left| A_{1} \right|}^{2} - {\left(\alpha_{5} \theta_{0}^{3} - 3 \, \alpha_{5} \theta_{0}^{2} + 2 \, \alpha_{5} \theta_{0}\right)} A_{1}\Big\} C_{1}\Big\}\\
		\omega_2&=\frac{1}{\theta_{0}}\bigg\{8 \, {\left(\alpha_{2} \theta_{0}^{4} - 2 \, \alpha_{2} \theta_{0}^{3} - \alpha_{2} \theta_{0}^{2} + 2 \, \alpha_{2} \theta_{0}\right)} A_{0} A_{1} {\left| A_{1} \right|}^{2} - 2 \, {\left(\theta_{0}^{3} - 3 \, \theta_{0}^{2} + 2 \, \theta_{0}\right)} A_{1} B_{1} {\left| A_{1} \right|}^{2}\\
		& - {\left(\theta_{0}^{3} \overline{\alpha_{1}} - 5 \, \theta_{0}^{2} \overline{\alpha_{1}} + 6 \, \theta_{0} \overline{\alpha_{1}}\right)} A_{1} C_{2} + \Big\{4 \, {\left(\theta_{0}^{3} - 4 \, \theta_{0}^{2} + 5 \, \theta_{0} - 2\right)} A_{0} {\left| A_{1} \right|}^{4} + 2 \, {\left(\alpha_{1} \theta_{0}^{3} \overline{\alpha_{1}} - 2 \, \alpha_{1} \theta_{0}^{2} \overline{\alpha_{1}}\right)} A_{0}\\
		& - {\left(\theta_{0}^{3} \overline{\alpha_{5}} - 3 \, \theta_{0}^{2} \overline{\alpha_{5}} + 2 \, \theta_{0} \overline{\alpha_{5}}\right)} A_{1} - 2 \, {\left(\theta_{0}^{3} \overline{\alpha_{1}} - 2 \, \theta_{0}^{2} \overline{\alpha_{1}}\right)} A_{2}\Big\} C_{1}\bigg\}\\
		\omega_3&=-{\left(\alpha_{1} \theta_{0}^{2} - 5 \, \alpha_{1} \theta_{0} + 6 \, \alpha_{1}\right)} A_{1} C_{2} + {\left(2 \, {\left(\alpha_{1}^{2} \theta_{0}^{2} - 2 \, \alpha_{1}^{2} \theta_{0}\right)} A_{0} - {\left(\alpha_{3} \theta_{0}^{2} - 3 \, \alpha_{3} \theta_{0} + 2 \, \alpha_{3}\right)} A_{1} - 2 \, {\left(\alpha_{1} \theta_{0}^{2} - 2 \, \alpha_{1} \theta_{0}\right)} A_{2}\right)} C_{1}\\
		\omega_4&=\frac{1}{\theta_{0}}\bigg\{4 \, {\left(\alpha_{2} \theta_{0}^{4} \overline{\alpha_{1}} - 2 \, \alpha_{2} \theta_{0}^{3} \overline{\alpha_{1}} - \alpha_{2} \theta_{0}^{2} \overline{\alpha_{1}} + 2 \, \alpha_{2} \theta_{0} \overline{\alpha_{1}}\right)} A_{0} A_{1} - {\left(\theta_{0}^{3} \overline{\alpha_{1}} - 3 \, \theta_{0}^{2} \overline{\alpha_{1}} + 2 \, \theta_{0} \overline{\alpha_{1}}\right)} A_{1} B_{1}\\
		& + {\left(2 \, {\left(\theta_{0}^{3} \overline{\alpha_{1}} - 4 \, \theta_{0}^{2} \overline{\alpha_{1}} + 5 \, \theta_{0} \overline{\alpha_{1}} - 2 \, \overline{\alpha_{1}}\right)} A_{0} {\left| A_{1} \right|}^{2} - {\left(\theta_{0}^{3} \overline{\alpha_{3}} - 3 \, \theta_{0}^{2} \overline{\alpha_{3}} + 2 \, \theta_{0} \overline{\alpha_{3}}\right)} A_{1}\right)} C_{1}\bigg\}
	\end{align*}

	then
	\begin{align}\label{hh2}
	\frac{5}{4}|\H|^2\h_0\totimes\h_0&=\begin{dmatrix}
		 \left(4 \, A_{1}^{2}\right) \left(\frac{5}{16} \, C_{1}^{2}\right) & 2 & 0 \\
		\left(4 \, A_{1}^{2}\right) \left(\frac{5}{16} \, C_{1} \overline{C_{1}}\right) & \theta_{0} & -\theta_{0} + 2 \\
		\left(4 \, A_{1}^{2}\right) \left(\frac{5}{16} \, C_{1} \overline{C_{1}}\right) & \theta_{0} & -\theta_{0} + 2 \\
		\left(4 \, A_{1}^{2}\right) \left(\frac{5}{16} \, \overline{C_{1}}^{2}\right) & 2 \, \theta_{0} - 2 & -2 \, \theta_{0} + 4
	\end{dmatrix}\\
	\s{\H}{\h_0}^2&=\begin{dmatrix}
		\left(A_{1} C_{1}\right) \left(A_{1} C_{1}\right) & 2 & 0 \\
		\left(A_{1} C_{1}\right) \left(A_{1} \overline{C_{1}}\right) & \theta_{0} & -\theta_{0} + 2 \\
		\left(A_{1} C_{1}\right) \left(A_{1} \overline{C_{1}}\right) & \theta_{0} & -\theta_{0} + 2 \\
		\left(A_{1} \overline{C_{1}}\right) \left(A_{1} \overline{C_{1}}\right) & 2 \, \theta_{0} - 2 & -2 \, \theta_{0} + 4
		\end{dmatrix}
	\end{align}
	and finally
	\begin{align*}
		-K_g\h_0\totimes\h_0=\begin{dmatrix}
		16 \, A_{1}^{2} {\left| A_{1} \right|}^{2} & \theta_{0} - 1 & -\theta_{0} + 1 \\
		-64 \, A_{0} A_{1} \alpha_{1} {\left| A_{1} \right|}^{2} + 64 \, A_{1} A_{2} {\left| A_{1} \right|}^{2} + 16 \, A_{1}^{2} \alpha_{5} & \theta_{0} & -\theta_{0} + 1 \\
		\kappa_1 & \theta_{0} + 1 & -\theta_{0} + 1 \\
		-8 \, A_{1}^{2} \alpha_{2} {\left(\theta_{0} + 1\right)} {\left(\theta_{0} - 2\right)} - \frac{8 \, A_{1} C_{1} {\left(\theta_{0} - 2\right)} {\left| A_{1} \right|}^{2}}{\theta_{0}} & 0 & 1 \\
		\kappa_2 & 0 & 2 \\
		\kappa_3 & 1 & 1 \\
		64 \, A_{0}^{2} {\left| A_{1} \right|}^{6} - 16 \, A_{1}^{2} {\left| A_{1} \right|}^{2} \overline{\alpha_{1}} - 96 \, A_{0} A_{1} {\left| A_{1} \right|}^{2} \overline{\alpha_{5}} + 4 \, A_{1} {\left| A_{1} \right|}^{2} \overline{B_{1}} - 24 \, {\left(2 \, {\left| A_{1} \right|}^{2} \overline{\alpha_{1}} - \overline{\alpha_{6}}\right)} A_{1}^{2} & \theta_{0} - 1 & -\theta_{0} + 3 \\
		-64 \, A_{0} A_{1} {\left| A_{1} \right|}^{4} + 16 \, A_{1}^{2} \overline{\alpha_{5}} & \theta_{0} - 1 & -\theta_{0} + 2 \\
		\kappa_4 & \theta_{0} & -\theta_{0} + 2 \\
		32 \, A_{0} A_{1} {\left(\theta_{0} + 1\right)} {\left(\theta_{0} - 2\right)} {\left| A_{1} \right|}^{2} \overline{\alpha_{2}} - 32 \, {\left(\theta_{0} \overline{\alpha_{2}} + \overline{\alpha_{2}}\right)} A_{0} A_{1} {\left| A_{1} \right|}^{2} - 8 \, A_{1}^{2} {\left(\theta_{0} + 1\right)} {\left(\theta_{0} - 3\right)} \overline{\alpha_{8}} & 2 \, \theta_{0} - 1 & -2 \, \theta_{0} + 3 \\
		\frac{4 \, A_{1} C_{1} \alpha_{2} {\left(\theta_{0} + 1\right)} {\left(\theta_{0} - 2\right)}^{2}}{\theta_{0}} + \frac{C_{1}^{2} {\left(\theta_{0} - 2\right)}^{2} {\left| A_{1} \right|}^{2}}{\theta_{0}^{2}} & -\theta_{0} + 1 & \theta_{0} + 1 \\
		-8 \, A_{1}^{2} \alpha_{7} {\left(\theta_{0} - 4\right)} \theta_{0} & 2 & 0 \\
		-8 \, A_{1}^{2} {\left(\theta_{0} + 1\right)} {\left(\theta_{0} - 2\right)} \overline{\alpha_{2}} & 2 \, \theta_{0} - 1 & -2 \, \theta_{0} + 2 \\
		32 \, A_{0} A_{1} \alpha_{1} {\left(\theta_{0} + 1\right)} {\left(\theta_{0} - 2\right)} \overline{\alpha_{2}} - 32 \, A_{1} A_{2} {\left(\theta_{0} + 1\right)} {\left(\theta_{0} - 2\right)} \overline{\alpha_{2}} - 8 \, A_{1}^{2} {\left(\theta_{0} + 2\right)} {\left(\theta_{0} - 2\right)} \overline{\alpha_{9}} & 2 \, \theta_{0} & -2 \, \theta_{0} + 2 \\
		-8 \, A_{1}^{2} {\left(\theta_{0} - 4\right)} \theta_{0} \overline{\alpha_{7}} & 2 \, \theta_{0} - 2 & -2 \, \theta_{0} + 4
		\end{dmatrix}
	\end{align*}
	\begin{align*}
		\kappa_1&=64 \, A_{0}^{2} \alpha_{1}^{2} {\left| A_{1} \right|}^{2} - 80 \, A_{1}^{2} \alpha_{1} {\left| A_{1} \right|}^{2} - 128 \, A_{0} A_{2} \alpha_{1} {\left| A_{1} \right|}^{2} - 96 \, A_{0} A_{1} \alpha_{3} {\left| A_{1} \right|}^{2} - 64 \, A_{0} A_{1} \alpha_{1} \alpha_{5} + 64 \, A_{2}^{2} {\left| A_{1} \right|}^{2}\\
		& + 96 \, A_{1} A_{3} {\left| A_{1} \right|}^{2} - 24 \, {\left(2 \, \alpha_{1} {\left| A_{1} \right|}^{2} - \alpha_{6}\right)} A_{1}^{2} + 64 \, A_{1} A_{2} \alpha_{5}\\
		\kappa_2&=32 \, A_{0} A_{1} \alpha_{2} {\left(\theta_{0} + 1\right)} {\left(\theta_{0} - 2\right)} {\left| A_{1} \right|}^{2} + \frac{16 \, A_{0} C_{1} {\left(\theta_{0} - 2\right)} {\left| A_{1} \right|}^{4}}{\theta_{0}} - 8 \, A_{1}^{2} \alpha_{9} {\left(\theta_{0} + 2\right)} {\left(\theta_{0} - 2\right)} + 32 \, {\left(\alpha_{2} \theta_{0} - 2 \, \alpha_{2}\right)} A_{0} A_{1} {\left| A_{1} \right|}^{2}\\
		& - \frac{8 \, A_{1} B_{1} {\left(\theta_{0} - 2\right)} {\left| A_{1} \right|}^{2}}{\theta_{0} + 1} - \frac{8 \, A_{1} C_{1} {\left(\theta_{0} - 2\right)} \overline{\alpha_{5}}}{\theta_{0}}\\
		\kappa_3&=32 \, A_{0} A_{1} \alpha_{1} \alpha_{2} {\left(\theta_{0} + 1\right)} {\left(\theta_{0} - 2\right)} - 32 \, A_{1} A_{2} \alpha_{2} {\left(\theta_{0} + 1\right)} {\left(\theta_{0} - 2\right)} - 8 \, A_{1}^{2} \alpha_{8} {\left(\theta_{0} + 1\right)} {\left(\theta_{0} - 3\right)} + \frac{16 \, A_{0} C_{1} \alpha_{1} {\left(\theta_{0} - 2\right)} {\left| A_{1} \right|}^{2}}{\theta_{0}}\\
		& - \frac{16 \, A_{2} C_{1} {\left(\theta_{0} - 2\right)} {\left| A_{1} \right|}^{2}}{\theta_{0}} - \frac{8 \, A_{1} C_{2} {\left(\theta_{0} - 3\right)} {\left| A_{1} \right|}^{2}}{\theta_{0}} - \frac{8 \, A_{1} C_{1} \alpha_{5} {\left(\theta_{0} - 2\right)}}{\theta_{0}}\\
		\kappa_4&=128 \, A_{0}^{2} \alpha_{1} {\left| A_{1} \right|}^{4} - 96 \, A_{1}^{2} {\left| A_{1} \right|}^{4} - 128 \, A_{0} A_{2} {\left| A_{1} \right|}^{4} - 128 \, A_{0} A_{1} \alpha_{5} {\left| A_{1} \right|}^{2} + \frac{4 \, A_{1} C_{1} {\left(\theta_{0} + 1\right)} {\left(\theta_{0} - 2\right)}^{2} \overline{\alpha_{2}}}{\theta_{0}} - 64 \, A_{0} A_{1} \alpha_{1} \overline{\alpha_{5}} \\
		&- 32 \, {\left(2 \, {\left| A_{1} \right|}^{4} + \alpha_{1} \overline{\alpha_{1}} - \beta\right)} A_{1}^{2} + 64 \, A_{1} A_{2} \overline{\alpha_{5}}
	\end{align*}
	\newpage
	Finally, the coefficient in $z^{\theta_0}\z^{2-\theta_0}dz^4$ in the Taylor expansion of the quartic form $\mathscr{Q}_{\phi}$ is
	\begin{align}
		\Omega_0&=\frac{2}{\theta_{0}} \, \bigg\{32 \, A_{1}^{2} \theta_{0} {\left| A_{1} \right|}^{4} + \ccancel{48 \, A_{0} A_{1} \alpha_{5} \theta_{0} {\left| A_{1} \right|}^{2}} - \ccancel{24 \, A_{0}^{2} \alpha_{3} \theta_{0} \overline{\alpha_{5}}} + 8 \, {\left(3 \, \alpha_{3} \overline{\alpha_{1}} + \alpha_{1} \overline{\alpha_{5}} - 3 \, \overline{\alpha_{16}}\right)} \ccancel{A_{0} A_{1}} \theta_{0}\nonumber\\
		& + 16 \, {\left(\alpha_{1} \overline{\alpha_{1}} - \beta\right)} A_{1}^{2} \theta_{0} - 8 \, A_{1} A_{2} \theta_{0} \overline{\alpha_{5}} + 24 \, A_{0} A_{3} \theta_{0} \overline{\alpha_{5}} - 2 \, {\left(\theta_{0}^{4} \overline{\alpha_{2}} - 4 \, \theta_{0}^{3} \overline{\alpha_{2}} + 2 \, \theta_{0}^{2} \overline{\alpha_{2}} + 4 \, \theta_{0} \overline{\alpha_{2}} - 3 \, \overline{\alpha_{2}}\right)} A_{0} C_{2}\nonumber\\
		& + 3 \, A_{1} \theta_{0} \overline{B_{4}} - {\left(\ccancel{2 \, {\left(\theta_{0}^{4} \overline{\alpha_{9}} - \theta_{0}^{3} \overline{\alpha_{9}} - 4 \, \theta_{0}^{2} \overline{\alpha_{9}} + 4 \, \theta_{0} \overline{\alpha_{9}}\right)} A_{0}} + {\left(2 \, \theta_{0}^{4} \overline{\alpha_{2}} - 4 \, \theta_{0}^{3} \overline{\alpha_{2}} - 2 \, \theta_{0}^{2} \overline{\alpha_{2}} + 3 \, \theta_{0} \overline{\alpha_{2}} + 2 \, \overline{\alpha_{2}}\right)} A_{1}\right)} C_{1}\nonumber\\
		& + 3 \, {\left(\ccancel{A_{1} \alpha_{1} \theta_{0}} + A_{0} \alpha_{3} \theta_{0} - A_{3} \theta_{0}\right)} \overline{B_{1}}\bigg\}\nonumber\\
		&-32 \, A_{0} A_{2} {\left| A_{1} \right|}^{4} + \ccancel{16 \, A_{0} A_{1} \alpha_{5} {\left| A_{1} \right|}^{2}} + 16 \, {\left(\ccancel{2 \, A_{0}^{2} \alpha_{1}} + A_{1}^{2}\right)} {\left| A_{1} \right|}^{4} - {\left(\theta_{0}^{2} \overline{\alpha_{2}} - 3 \, \theta_{0} \overline{\alpha_{2}} + 2 \, \overline{\alpha_{2}}\right)} A_{1} C_{1}\nonumber\\
		&+\frac{5}{2}|\vec{C}_1|^2\s{\vec{A}_1}{\vec{A}_1}+2\s{\vec{A}_1}{\vec{C}_1}\s{\bar{\vec{A}_1}}{\bar{\vec{C}_1}}\nonumber\\
		&+\ccancel{128 \, A_{0}^{2} \alpha_{1} {\left| A_{1} \right|}^{4}} - 96 \, A_{1}^{2} {\left| A_{1} \right|}^{4} - 128 \, A_{0} A_{2} {\left| A_{1} \right|}^{4} - \ccancel{128 \, A_{0} A_{1} \alpha_{5} {\left| A_{1} \right|}^{2}} + \frac{4 \, A_{1} C_{1} {\left(\theta_{0} + 1\right)} {\left(\theta_{0} - 2\right)}^{2} \overline{\alpha_{2}}}{\theta_{0}} - \ccancel{64 \, A_{0} A_{1} \alpha_{1} \overline{\alpha_{5}}} \nonumber\\
		&- 32 \, {\left(2 \, {\left| A_{1} \right|}^{4} + \alpha_{1} \overline{\alpha_{1}} - \beta\right)} A_{1}^{2} + 64 \, A_{1} A_{2} \overline{\alpha_{5}}\\
		&=\frac{2}{\theta_{0}} \, \bigg\{32 \, A_{1}^{2} \theta_{0} {\left| A_{1} \right|}^{4} 
		 + 16 \, {\left(|\alpha_{1}|^2- \beta\right)} A_{1}^{2} \theta_{0} - 8 \, A_{1} A_{2} \theta_{0} \overline{\alpha_{5}} + 24 \, A_{0} A_{3} \theta_{0} \overline{\alpha_{5}} - 2 \, {\left(\theta_{0}^{4} - 4 \, \theta_{0}^{3} + 2 \, \theta_{0}^{2}  + 4 \, \theta_{0}  - 3 \, \right)}\overline{\alpha_{2}} A_{0} C_{2}\nonumber\\
		& + 3 \, A_{1} \theta_{0} \overline{B_{4}} - {\left({\left(2 \, \theta_{0}^{4}  - 4 \, \theta_{0}^{3}  - 2 \, \theta_{0}^{2}  + 3 \, \theta_{0} + 2 \, \right)}\overline{\alpha_{2}} A_{1}\right)} C_{1}
		 + 3 \, {\left(A_{0} \alpha_{3} \theta_{0} - A_{3} \theta_{0}\right)} \overline{B_{1}}\bigg\}\nonumber\\
		&-32 \, A_{0} A_{2} {\left| A_{1} \right|}^{4} + 16 \, {\left( A_{1}^{2}\right)} {\left| A_{1} \right|}^{4} - {\left(\theta_{0}^{2}  - 3 \, \theta_{0}+ 2 \, \right)}\overline{\alpha_{2}} A_{1} C_{1}\nonumber\\
		&+\frac{5}{2}|\vec{C}_1|^2\s{\vec{A}_1}{\vec{A}_1}+2\s{\vec{A}_1}{\vec{C}_1}\s{\bar{\vec{A}_1}}{\bar{\vec{C}_1}}\nonumber\\
		&- 96 \, A_{1}^{2} {\left| A_{1} \right|}^{4} - 128 \, A_{0} A_{2} {\left| A_{1} \right|}^{4} + \frac{4 \, A_{1} C_{1} {\left(\theta_{0} + 1\right)} {\left(\theta_{0} - 2\right)}^{2} \overline{\alpha_{2}}}{\theta_{0}}  \nonumber\\
		&- 32 \, {\left(2 \, {\left| A_{1} \right|}^{4} + |\alpha_{1}|^2 - \beta\right)} A_{1}^{2} + 64 \, A_{1} A_{2} \overline{\alpha_{5}}.
	\end{align}
	Let us first look at the coefficients
	\begin{align*}
		|\vec{A}_1|^4\s{\vec{A}_1}{\vec{A}_1}\,\quad \text{and}\;\, |\vec{A}_1|^4\s{\vec{A}_0}{\vec{A}_2}.
	\end{align*}
	We see that these are
	\begin{align*}
		&\frac{2}{\theta_0}\left(32\theta_0|\vec{A}_1|^4\s{\vec{A}_1}{\vec{A}_1}\right)-32|\vec{A}_1|^4\s{\vec{A}_0}{\vec{A}_2}+16|\vec{A}_1|^4\s{\vec{A}_1}{\vec{A}_1}-96|\vec{A}_1|^4\s{\vec{A}_1}{\vec{A}_1}-128|\vec{A}_1|^2\s{\vec{A}_1}{\vec{A}_1}-64|\vec{A}_1|^4\s{\vec{A}_1}{\vec{A}_1}\\
		&=-80|\vec{A}_1|^4\left(\s{\vec{A}_1}{\vec{A}_1}+2\s{\vec{A}_0}{\vec{A}_2}\right)=0.
	\end{align*}
	Then, the coefficient involving $|\alpha_1|^2-\beta$ is
	\begin{align*}
		\frac{2}{\theta_0}\left(16\theta_0|\alpha_1|^2-\beta\right)\s{\vec{A}_1}{\vec{A}_1}-32\left(|\alpha_1|^2-\beta\right)\s{\vec{A}_1}{\vec{A}_1}=0.
	\end{align*}
	Finally, the coefficient involving $\bar{\alpha_5}$ is
	\begin{align*}
		\frac{2}{\theta_0}\left(-8\theta_0\bar{\alpha_5}\s{\vec{A}_1}{\vec{A}_2}+24\theta_0\bar{\alpha_5}\s{\vec{A}_0}{\vec{A}_3}\right)+64\bar{\alpha_5}\s{\vec{A}_1}{\vec{A}_2}=48\bar{\alpha_5}\left(\s{\vec{A}_1}{\vec{A}_2}+\s{\vec{A}_0}{\vec{A}_3}\right)=0.
	\end{align*}
	Therefore, we deduce that (using $\s{\vec{A}_1}{\vec{C}_1}+\s{\vec{A}_0}{\vec{C}_2}=0$)
	\normalsize
	\begin{align}
	\small
		\Omega_0&=\frac{2}{\theta_{0}} \, \bigg\{\ccancel{32 \, A_{1}^{2} \theta_{0} {\left| A_{1} \right|}^{4}} 
		+ \ccancel{16 \, {\left(|\alpha_{1}|^2- \beta\right)} A_{1}^{2} \theta_{0}} - \ccancel{8 \, A_{1} A_{2} \theta_{0} \overline{\alpha_{5}}} + \ccancel{24 \, A_{0} A_{3} \theta_{0} \overline{\alpha_{5}}} - 2 \, {\left(\theta_{0}^{4} - 4 \, \theta_{0}^{3} + 2 \, \theta_{0}^{2}  + 4 \, \theta_{0}  - 3 \, \right)}\overline{\alpha_{2}} A_{0} C_{2}\nonumber\\
		& + 3 \, A_{1} \theta_{0} \overline{B_{4}} - {\left({\left(2 \, \theta_{0}^{4}  - 4 \, \theta_{0}^{3}  - 2 \, \theta_{0}^{2}  + 3 \, \theta_{0} + 2 \, \right)}\overline{\alpha_{2}} A_{1}\right)} C_{1}
		+ 3 \, {\left(A_{0} \alpha_{3} \theta_{0} - A_{3} \theta_{0}\right)} \overline{B_{1}}\bigg\}\nonumber\\
		&-\ccancel{32 \, A_{0} A_{2} {\left| A_{1} \right|}^{4}} + \ccancel{16 \, {\left( A_{1}^{2}\right)} {\left| A_{1} \right|}^{4}} - {\left(\theta_{0}^{2}  - 3 \, \theta_{0}+ 2 \, \right)}\overline{\alpha_{2}} A_{1} C_{1}\nonumber\\
		&+\frac{5}{2}|\vec{C}_1|^2\s{\vec{A}_1}{\vec{A}_1}+2\s{\vec{A}_1}{\vec{C}_1}\s{\bar{\vec{A}_1}}{\bar{\vec{C}_1}}\nonumber\\
		&- \ccancel{96 \, A_{1}^{2} {\left| A_{1} \right|}^{4}} - \ccancel{128 \, A_{0} A_{2} {\left| A_{1} \right|}^{4}} + \frac{4 \, A_{1} C_{1} {\left(\theta_{0} + 1\right)} {\left(\theta_{0} - 2\right)}^{2} \overline{\alpha_{2}}}{\theta_{0}}  \nonumber\\
		&- \ccancel{32 \, {\left(2 \, {\left| A_{1} \right|}^{4} + |\alpha_{1}|^2 - \beta\right)} A_{1}^{2}} + \ccancel{64 \, A_{1} A_{2} \overline{\alpha_{5}}}\nonumber\\
		\normalsize
		&=\frac{2}{\theta_{0}} \, \bigg\{
		   2 \, {\left(\theta_{0}^{4} - 4 \, \theta_{0}^{3} + 2 \, \theta_{0}^{2}  + 4 \, \theta_{0}  - 3 \, \right)}\overline{\alpha_{2}} A_{1} C_{1} + 3 \, A_{1} \theta_{0} \overline{B_{4}} - {\left({\left(2 \, \theta_{0}^{4}  - 4 \, \theta_{0}^{3}  - 2 \, \theta_{0}^{2}  + 3 \, \theta_{0} + 2 \, \right)}\overline{\alpha_{2}} A_{1}\right)} C_{1}\nonumber\\
		   &
		+ 3 \, {\left(A_{0} \alpha_{3} \theta_{0} - A_{3} \theta_{0}\right)} \overline{B_{1}}\bigg\}
		 - {\left(\theta_{0}^{2}  - 3 \, \theta_{0}+ 2 \, \right)}\overline{\alpha_{2}} A_{1} C_{1}
		+\frac{5}{2}|\vec{C}_1|^2\s{\vec{A}_1}{\vec{A}_1}+2\s{\vec{A}_1}{\vec{C}_1}\s{\bar{\vec{A}_1}}{\bar{\vec{C}_1}}\nonumber\\
		 &+ \frac{4 \, A_{1} C_{1} {\left(\theta_{0} + 1\right)} {\left(\theta_{0} - 2\right)}^{2} \overline{\alpha_{2}}}{\theta_{0}}  
	\end{align}
	Now, recall that
	\begin{align*}
		\frac{1}{2}\vec{B}_4=\begin{dmatrix}
		-\frac{{\left(\theta_{0}^{2} + 4 \, \theta_{0} + 3\right)} C_{1} \overline{C_{2}} - {\left(12 \, {\left(\alpha_{2} \theta_{0}^{2} + \alpha_{2} \theta_{0}\right)} A_{0} - {\left(\theta_{0}^{2} + 4 \, \theta_{0}\right)} B_{1}\right)} \overline{C_{1}}}{12 \, {\left(\theta_{0}^{2} + \theta_{0}\right)}} & \overline{A_{0}} & -\theta_{0} + 2 & 3 \\
		-\frac{\overline{A_{1}} \overline{C_{1}}}{12 \, \theta_{0}} & C_{1} & -\theta_{0} + 2 & 3 \\
		-\frac{C_{1} {\left(\theta_{0} + 3\right)} \overline{C_{1}}}{12 \, \theta_{0}} & \overline{A_{1}} & -\theta_{0} + 2 & 3 \\
		-\frac{1}{24} \, C_{1} \overline{A_{1}} & \overline{C_{1}} & -\theta_{0} + 2 & 3 \\
		\frac{1}{3} \, {\left(\alpha_{2} \theta_{0} + \alpha_{2}\right)} \overline{A_{0}} \overline{C_{1}} + \frac{2}{3} \, {\left(\overline{A_{0}} \overline{\alpha_{1}} - \overline{A_{2}}\right)} B_{1} + {\left(\overline{A_{1}} \overline{\alpha_{1}} + \overline{A_{0}} \overline{\alpha_{3}} - \overline{A_{3}}\right)} C_{1} - \frac{1}{3} \, B_{2} \overline{A_{1}} & A_{0} & -\theta_{0} + 2 & 3 
		\end{dmatrix}
	\end{align*}
	and for some unimportant $\lambda_1,\lambda_2\in\C$, we have
	\begin{align}
		\vec{B}_4=\lambda_1\vec{A}_0+\lambda_2\bar{\vec{A}_0}-\frac{1}{6\theta_0}\bar{\s{\vec{A}_1}{\vec{C}_1}}\vec{C}_1-\frac{(\theta_0+3)}{6\theta_0}|\vec{C}_1|^2\bar{\vec{A}_1}-\frac{1}{12}\s{\bar{\vec{A}_1}}{\vec{C}_1}\bar{\vec{C}_1}.
	\end{align}
	As $\s{\vec{A}_0}{\vec{A}_1}=\s{\vec{A}_0}{\bar{\vec{A}_1}}$, we obtain
	\begin{align}
		\s{\bar{\vec{A}_1}}{\vec{B}_4}=-\frac{(\theta_0+2)}{12\theta_0}\bar{\s{\vec{A}_1}{\vec{C}_1}}\vec{C}_1-\frac{(\theta_0+3)}{6\theta_0}|\vec{C}_1|^2\bar{\s{\vec{A}_1}{\vec{A}_1}}.
	\end{align}
	Now, as $\alpha_3\in\C$ is defined by
	\begin{align}
		\alpha_3=\frac{1}{12}\s{\vec{A}_1}{\vec{C}_1}+2\s{\bar{\vec{A}_0}}{\vec{A}_3},
	\end{align}
	we have
	\begin{align*}
		\s{\alpha_3\vec{A}_0-\vec{A}_3}{\bar{\vec{A}_0}}=\frac{\alpha_3}{2}-\s{\bar{\vec{A}_0}}{\vec{A}_3}=\frac{1}{24}\s{\vec{A}_1}{\vec{C}_1},
	\end{align*}
	so
	\begin{align}
		\s{\alpha_3\vec{A}_0-\vec{A}_3}{\bar{\vec{B}_1}}=-2\s{\vec{A}_1}{\bar{\vec{C}_1}}\s{\alpha_3\vec{A}_0-\vec{A}_3}{\bar{\vec{A}_0}}=-\frac{1}{12}\s{\vec{A}_1}{\vec{C}_1}\s{\vec{A}_1}{\bar{\vec{C}_1}}.
	\end{align}
	Finally, as 
	\begin{align}
		\alpha_{2}=\frac{1}{2\theta_0(\theta_0+1)}\s{\vec{A}_1}{\bar{\vec{C}_1}},
	\end{align}
    we get
    \begin{align}
    	\Omega_0&=\bigg\{\frac{2}{\theta_0}\left(\frac{(\theta_0^4-4\theta_0^3+2\theta_0^2+4\theta_0-3)}{\theta_0(\theta_0+1)}-\frac{(\theta_0+2)}{4}-\frac{2\theta_0^4-4\theta_0^3-2\theta_0^2+3\theta_0+2}{2\theta_0(\theta_0+1)}-\frac{\theta_0}{4}\right)\nonumber\\
    	&-\frac{\theta_0^2-3\theta_0+2}{2\theta_0(\theta_0+1)}+2+\frac{2(\theta_0-2)^2}{\theta_0^2}\bigg\}\s{\vec{A}_1}{\bar{\vec{C}_1}}\s{\vec{A}_1}{\vec{C}_1}
    	+\bigg\{\frac{2}{\theta_0}\left(-\frac{(\theta_0+3)}{2}\right)+\frac{5}{2}\bigg\}|\vec{C}_1|^2\s{\vec{A}_1}{\vec{A}_1}\nonumber\\
    	&=-\frac{3(\theta_0-2)}{2\theta_0}\s{\vec{A}_1}{\vec{C}_1}\s{\vec{A}_1}{\bar{\vec{C}_1}}+\frac{3(\theta_0-2)}{2\theta_0}|\vec{C}_1|^2\s{\vec{A}_1}{\vec{A}_1}\nonumber\\
    	&=\frac{3(\theta_0-2)}{2\theta_0}\left(|\vec{C}_1|^2\s{\vec{A}_1}{\vec{A}_1}-\s{\vec{A}_1}{\bar{\vec{C}_1}}\s{\vec{A}_1}{\vec{C}_1}\right)\nonumber\\
    	&=0
    \end{align}
    so
    \begin{align}
    	|\vec{C}_1|^2\s{\vec{A}_1}{\vec{A}_1}=\s{\vec{A}_1}{\bar{\vec{C}_1}}\s{\vec{A}_1}{\vec{C}_1}
    \end{align}

	 If we suppose that $\mathscr{Q}_{\phi}$ is meromorphic, then we obtain
	\begin{align}\label{system}
	\left\{\begin{alignedat}{1}
	|\vec{A}_1|^2\s{\vec{A}_1}{\vec{C}_1}&=\s{\bar{\vec{A}_1}}{\vec{C}_1}\s{\vec{A}_1}{\vec{A}_1}\\
	|\vec{C}_1|^2\s{\vec{A}_1}{\vec{A}_1}&=\s{\vec{A}_1}{\bar{\vec{C}_1}}\s{\vec{A}_1}{\vec{C}_1},
	\end{alignedat}\right.
	\end{align}
	Remarking that is a linear system in $(\s{\vec{A}_1}{\vec{C}_1},\s{\vec{A}_1}{\vec{A}_1})$, we can recast \eqref{system} as
	\begin{align}\label{salvationmatrix}
	\begin{dmatrix}
	|\vec{A}_1|^2 & -\s{\bar{\vec{A}_1}}{\vec{C}_1}\\
	-\s{\vec{A}_1}{\bar{\vec{C}_1}} & |\vec{C}_1|^2
	\end{dmatrix}
	\begin{dmatrix}
	\s{\vec{A}_1}{\vec{C}_1}\\
	\s{\vec{A}_1}{\vec{A}_1}
	\end{dmatrix}=0.
	\end{align}
	Thanks of Cauchy-Schwarz inequality, we obtain
	\begin{align}\label{determinant}
	\det 	\begin{dmatrix}
	|\vec{A}_1|^2 & -\s{\bar{\vec{A}_1}}{\vec{C}_1}\\
	-\s{\vec{A}_1}{\bar{\vec{C}_1}} & |\vec{C}_1|^2
	\end{dmatrix}=|\vec{A}_1|^2|\vec{C}_1|^2-|\s{\vec{A}_1}{\bar{\vec{C}_1}}|^2\geq 0.
	\end{align}
	Therefore, if the determinant is positive, we obtain
	\begin{align*}
	\s{\vec{A}_1}{\vec{C}_1}=0,
	\end{align*}
	and the holomorphy of the quartic form, and if the determinant vanishes,
	\begin{align}\label{ccancel}
	\vec{A}_1\;\,\text{and}\;\, \vec{C}_1\;\,\text{are proportional}.
	\end{align}

		\chapter{Return to the invariance by inversions}
		
		To obtain the next order development, we need to develop
		\begin{align*}
			\vec{\beta}=\mathscr{I}_{\phi}(\vec{\alpha})-g^{-1}\otimes\left(\h_0\otimes\bar{\partial}|\phi|^2-2\s{\phi}{\h_0}\otimes\bar{\partial}\phi\right)
		\end{align*}
		where
		\begin{align*}
			\vec{\alpha}&=\partial \H+|\H|^2\partial\phi+2\,g^{-1}\otimes\s{\H}{\h_0}\otimes\bar{\partial}\phi\\
			\mathscr{I}_{\phi}(\vec{X})&=|\phi|^2\vec{X}-2\s{\phi}{\vec{X}}\phi,
		\end{align*}
		up to an error of order ${\theta+4}$.

		\section{Computation of the order of development of tensors}
		
		We first need to see at which order we need to develop tensors to obtain the next order development of $\vec{\beta}$. We should obtain an error in $O(|z|^{\theta_0+4})$.
		
		To obtain such error, we need to develop $\vec{\alpha}$ and $\phi$ up to order $3$, for the component 
		\begin{align*}
		\mathscr{I}_{\phi}(\vec{\alpha})=\mathscr{I}_{\phi}\left(\partial\H+|\H|^2\partial\phi+2\,g^{-1}\otimes\s{\H}{\h_0}\otimes\bar{\partial}\phi\right)
		\end{align*}
		where for all $\vec{X}\in \mathbb{C}^n$, we have
		\begin{align*}
		\mathscr{I}_{\phi}(\vec{X})=|\phi|^2\vec{X}-2\s{\phi}{\vec{X}}\phi.
		\end{align*}
		As we perform this order $3$ development, we get
		\begin{align*}
		\vec{\alpha}=-\frac{(\theta_0-2)}{2}\frac{\vec{C}_1}{z^{\theta_0-1}}+\cdots+O(|z|^{4-\theta_0})
		\end{align*}
		and 
		\begin{align*}
		|\phi(z)|^2=\frac{1}{\theta_0^2}|z|^{2\theta_0}\left(1+\cdots+O(|z|^3)\right)
		\end{align*}
		so
		\begin{align*}
		|\phi|^2\vec{\alpha}=-\frac{(\theta_0-2)}{2\theta_0^2}\vec{C}_1 z\z^{\theta_0}+\cdots+O(|z|^{\theta_0+4})
		\end{align*}
		and likewise
		\begin{align*}
		\s{\phi}{\vec{\alpha}}{\phi}=\cdots +O(|z|^{\theta_0+4})
		\end{align*}
		so the developments of $\phi$ and $\vec{\alpha}$ are sufficient for the part
		\begin{align*}
		\mathscr{I}_{\phi}(\vec{\alpha})=|\phi|^2\vec{\alpha}-2\s{\phi}{\vec{\alpha}}\phi.
		\end{align*}
		Furthermore, we have
		\begin{align*}
		\s{\vec{\alpha}}{\phi}=\bs{-\frac{(\theta_0-2)}{2}\frac{\vec{C}_1}{z^{\theta_0-1}}+\cdots+O(|z|^{4-\theta_0})}{\Re\left(\frac{2}{\theta_0}\vec{A}_0z^{\theta_0}\right)+O(|z|^{\theta_0+3})}=\cdots+O(|z|^{4}).
		\end{align*}
		\textbf{We indicate the order of the errors in the partial developments so that one can check when we throw higher order terms with the \texttt{throw} function that the order is correct.}

		However, we need to develop $\phi$, $g$ and $\h_0$ up to order $4$ for the other part
		\begin{align*}
		-g^{-1}\otimes\left(\bar{\partial}|\phi|^2\otimes \h_0-2\s{\phi}{\h_0}\otimes\bar{\partial}\phi\right).
		\end{align*}
		Indeed, with
		\begin{align*}
		\h_0&=2\,\vec{A}_1z^{\theta_0-1}dz^2+\cdots+O(|z|^{\theta_0+3})\\
		g^{-1}&=|z|^{2-2\theta_0}\left(1+\cdots+O(|z|^{4})\right)\\
		|\phi|^2&=\frac{1}{\theta_0^2}|z|^{2\theta_0}\left(1+\cdots+O(|z|^{4})\right)
		\end{align*}
		we obtain
		\begin{align*}
		\bar{\partial}|\phi|^2=\frac{1}{\theta_0}|z|^{2\theta_0-2}\left(z\,d\z+\cdots+O(|z|^{5})\right)
		\end{align*}
		so
		\begin{align*}
		g^{-1}\otimes\h_0=2\vec{A}_1\z^{1-\theta_0}\frac{dz}{d\z}+\cdots+O(|z|^{5-\theta_0})
		\end{align*}
		and
		\begin{align*}
		g^{-1}\otimes\bar{\partial}|\phi|^2\otimes\h_0=\frac{2}{\theta_0}\vec{A}_1z^{\theta_0}\,dz+\cdots+O(|z|^{\theta_0+4}).
		\end{align*}
		We also have
		\begin{align*}
		g^{-1}\otimes\s{\h_0}{\phi}=\bs{2\vec{A}_1\z^{1-\theta_0}\frac{dz}{d\z}+\cdots+O(|z|^{5-\theta_0})}{\Re\left(\frac{2}{\theta_0}\vec{A}_0z^{\theta_0}\right)+O(|z|^{\theta_0+4})}=\cdots+O(|z|^{5}).
		\end{align*}
		
		\section{Development of $\phi$ and $|\phi|^2$ up to order $4$}
		
		We first recall that
		\begin{align*}
			\phi(z)&=\begin{dmatrix}
			\frac{1}{\theta_{0}} & \overline{A_{0}} & 0 & \theta_{0} \\
			\frac{1}{\theta_{0} + 1} & \overline{A_{1}} & 0 & \theta_{0} + 1 \\
			\frac{1}{\theta_{0} + 2} & \overline{A_{2}} & 0 & \theta_{0} + 2 \\
			\frac{1}{\theta_{0} + 3} & \overline{A_{3}} & 0 & \theta_{0} + 3 \\
			\frac{1}{\theta_{0}} & A_{0} & \theta_{0} & 0 \\
			\frac{1}{\theta_{0} + 1} & A_{1} & \theta_{0} + 1 & 0 \\
			\frac{1}{\theta_{0} + 2} & A_{2} & \theta_{0} + 2 & 0 \\
			\frac{1}{\theta_{0} + 3} & A_{3} & \theta_{0} + 3 & 0 
			\end{dmatrix}
			\begin{dmatrix}
			\frac{1}{8 \, \theta_{0}} & C_{1} & 2 & \theta_{0} \\
			\frac{1}{8 \, {\left(\theta_{0} + 1\right)}} & B_{1} & 2 & \theta_{0} + 1 \\
			\frac{1}{12 \, \theta_{0}} & C_{2} & 3 & \theta_{0} \\
			\frac{1}{8 \, {\left(\theta_{0} + 1\right)}} & \overline{B_{1}} & \theta_{0} + 1 & 2 \\
			\frac{1}{8 \, \theta_{0}} & \overline{C_{1}} & \theta_{0} & 2 \\
			\frac{1}{12 \, \theta_{0}} & \overline{C_{2}} & \theta_{0} & 3
			\end{dmatrix}+O(|z|^{\theta_0+4})\\
			&=\Re\bigg(\frac{2}{\theta_0}\vec{A}_0z^{\theta_0}+\frac{2}{\theta_0+1}\vec{A}_1z^{\theta_0+1}+\frac{2}{\theta_0+2}\vec{A}_2z^{\theta_0+2}+\frac{2}{\theta_0+3}\vec{A}_3z^{\theta_0+3}+\frac{1}{4\theta_0}\vec{C}_1z^{2}\z^{\theta_0}+\frac{1}{6\theta_0}\vec{C}_2z^3\z^{\theta_0}\\
			&+\frac{1}{4(\theta_0+1)}\vec{B}_1z^2\z^{\theta_0+1}\bigg)+O(|z|^{\theta_0+4})
		\end{align*}
		Recall that
		\begin{align*}
			\alpha_2=\frac{1}{2\theta_0(\theta_0+1)}\s{\bar{\vec{A}_1}}{\vec{C}_1}
		\end{align*}
		so
		\begin{align}\label{newb1}
			\vec{B}_1=-2\s{\bar{\vec{A}_1}}{\vec{C}_1}\vec{A}_0=-4\theta_0(\theta_0+1)\alpha_2\vec{A}_0.
		\end{align}
		Now we compute
		\begin{align*}
			|\phi(z)|^2=\begin{dmatrix}
			\frac{\ccancel{2 \, \overline{A_{0}} \overline{A_{1}}}}{\theta_{0}^{2} + \theta_{0}} & 0 & 2 \, \theta_{0} + 1 \\
			\frac{{\left(\theta_{0}^{2} + 2 \, \theta_{0}\right)} \overline{A_{1}}^{2} + 2 \, {\left(\theta_{0}^{2} + 2 \, \theta_{0} + 1\right)} \overline{A_{0}} \overline{A_{2}}}{\theta_{0}^{4} + 4 \, \theta_{0}^{3} + 5 \, \theta_{0}^{2} + 2 \, \theta_{0}} & 0 & 2 \, \theta_{0} + 2 \\
			\frac{2 \, {\left({\left(\theta_{0}^{2} + 3 \, \theta_{0}\right)} \overline{A_{1}} \overline{A_{2}} + {\left(\theta_{0}^{2} + 3 \, \theta_{0} + 2\right)} \overline{A_{0}} \overline{A_{3}}\right)}}{\theta_{0}^{4} + 6 \, \theta_{0}^{3} + 11 \, \theta_{0}^{2} + 6 \, \theta_{0}} & 0 & 2 \, \theta_{0} + 3 \\
			\frac{\ccancel{\overline{A_{0}}^{2}}}{\theta_{0}^{2}} & 0 & 2 \, \theta_{0} \\
			\frac{B_{1} \overline{A_{0}} + C_{1} \overline{A_{1}}}{4 \, {\left(\theta_{0}^{2} + \theta_{0}\right)}} & 2 & 2 \, \theta_{0} + 1 \\
			\frac{\ccancel{C_{1} \overline{A_{0}}}}{4 \, \theta_{0}^{2}} & 2 & 2 \, \theta_{0} \\
			\frac{\ccancel{C_{2} \overline{A_{0}}}}{6 \, \theta_{0}^{2}} & 3 & 2 \, \theta_{0} \\
			\frac{2 \, A_{0} \overline{A_{0}}}{\theta_{0}^{2}} & \theta_{0} & \theta_{0} \\
			\frac{\ccancel{2 \, A_{0} \overline{A_{1}}}}{\theta_{0}^{2} + \theta_{0}} & \theta_{0} & \theta_{0} + 1 \\
			\frac{8 \, A_{0} \theta_{0} \overline{A_{2}} + \ccancel{{\left(\theta_{0} + 2\right)} \overline{A_{0}} \overline{C_{1}}}}{4 \, {\left(\theta_{0}^{3} + 2 \, \theta_{0}^{2}\right)}} & \theta_{0} & \theta_{0} + 2 \\
			\frac{24 \, {\left(\theta_{0}^{2} + \theta_{0}\right)} A_{0} \overline{A_{3}} + 3 \, {\left(\theta_{0}^{2} + 3 \, \theta_{0}\right)} \overline{A_{1}} \overline{C_{1}} + 2 \, {\left(\theta_{0}^{2} + 4 \, \theta_{0} + 3\right)} \overline{A_{0}} \overline{C_{2}}}{12 \, {\left(\theta_{0}^{4} + 4 \, \theta_{0}^{3} + 3 \, \theta_{0}^{2}\right)}} & \theta_{0} & \theta_{0} + 3 \\
			\frac{\ccancel{2 \, A_{1} \overline{A_{0}}}}{\theta_{0}^{2} + \theta_{0}} & \theta_{0} + 1 & \theta_{0} \\
			\frac{2 \, A_{1} \overline{A_{1}}}{\theta_{0}^{2} + 2 \, \theta_{0} + 1} & \theta_{0} + 1 & \theta_{0} + 1 \\
			\frac{8 \, A_{1} \theta_{0} \overline{A_{2}} + \ccancel{{\left(\theta_{0} + 2\right)} \overline{A_{0}} \overline{B_{1}}}}{4 \, {\left(\theta_{0}^{3} + 3 \, \theta_{0}^{2} + 2 \, \theta_{0}\right)}} & \theta_{0} + 1 & \theta_{0} + 2 \\
			\frac{\ccancel{A_{0} C_{1} {\left(\theta_{0} + 2\right)}} + 8 \, A_{2} \theta_{0} \overline{A_{0}}}{4 \, {\left(\theta_{0}^{3} + 2 \, \theta_{0}^{2}\right)}} & \theta_{0} + 2 & \theta_{0} \\
			\frac{\ccancel{A_{0} B_{1} {\left(\theta_{0} + 2\right)}} + 8 \, A_{2} \theta_{0} \overline{A_{1}}}{4 \, {\left(\theta_{0}^{3} + 3 \, \theta_{0}^{2} + 2 \, \theta_{0}\right)}} & \theta_{0} + 2 & \theta_{0} + 1 \\
			\frac{3 \, {\left(\theta_{0}^{2} + 3 \, \theta_{0}\right)} A_{1} C_{1} + 2 \, {\left(\theta_{0}^{2} + 4 \, \theta_{0} + 3\right)} A_{0} C_{2} + 24 \, {\left(\theta_{0}^{2} + \theta_{0}\right)} A_{3} \overline{A_{0}}}{12 \, {\left(\theta_{0}^{4} + 4 \, \theta_{0}^{3} + 3 \, \theta_{0}^{2}\right)}} & \theta_{0} + 3 & \theta_{0} \\
			\frac{\ccancel{2 \, A_{0} A_{1}}}{\theta_{0}^{2} + \theta_{0}} & 2 \, \theta_{0} + 1 & 0 \\
			\frac{A_{0} \overline{B_{1}} + A_{1} \overline{C_{1}}}{4 \, {\left(\theta_{0}^{2} + \theta_{0}\right)}} & 2 \, \theta_{0} + 1 & 2 \\
			\frac{{\left(\theta_{0}^{2} + 2 \, \theta_{0}\right)} A_{1}^{2} + 2 \, {\left(\theta_{0}^{2} + 2 \, \theta_{0} + 1\right)} A_{0} A_{2}}{\theta_{0}^{4} + 4 \, \theta_{0}^{3} + 5 \, \theta_{0}^{2} + 2 \, \theta_{0}} & 2 \, \theta_{0} + 2 & 0 \\
			\frac{2 \, {\left({\left(\theta_{0}^{2} + 3 \, \theta_{0}\right)} A_{1} A_{2} + {\left(\theta_{0}^{2} + 3 \, \theta_{0} + 2\right)} A_{0} A_{3}\right)}}{\theta_{0}^{4} + 6 \, \theta_{0}^{3} + 11 \, \theta_{0}^{2} + 6 \, \theta_{0}} & 2 \, \theta_{0} + 3 & 0 \\
			\frac{\ccancel{A_{0}^{2}}}{\theta_{0}^{2}} & 2 \, \theta_{0} & 0 \\
			\frac{\ccancel{A_{0} \overline{C_{1}}}}{4 \, \theta_{0}^{2}} & 2 \, \theta_{0} & 2 \\
			\frac{\ccancel{A_{0} \overline{C_{2}}}}{6 \, \theta_{0}^{2}} & 2 \, \theta_{0} & 3
			\end{dmatrix}
		\end{align*}
		\begin{align*}
		=\begin{dmatrix}
		\frac{{\left(\theta_{0}^{2} + 2 \, \theta_{0}\right)} \overline{A_{1}}^{2} + 2 \, {\left(\theta_{0}^{2} + 2 \, \theta_{0} + 1\right)} \overline{A_{0}} \overline{A_{2}}}{\theta_{0}^{4} + 4 \, \theta_{0}^{3} + 5 \, \theta_{0}^{2} + 2 \, \theta_{0}} & 0 & 2 \, \theta_{0} + 2 &\textbf{(1)}\\
		\frac{2 \, {\left({\left(\theta_{0}^{2} + 3 \, \theta_{0}\right)} \overline{A_{1}} \overline{A_{2}} + {\left(\theta_{0}^{2} + 3 \, \theta_{0} + 2\right)} \overline{A_{0}} \overline{A_{3}}\right)}}{\theta_{0}^{4} + 6 \, \theta_{0}^{3} + 11 \, \theta_{0}^{2} + 6 \, \theta_{0}} & 0 & 2 \, \theta_{0} + 3 &\textbf{(2)}\\
		\frac{B_{1} \overline{A_{0}} + C_{1} \overline{A_{1}}}{4 \, {\left(\theta_{0}^{2} + \theta_{0}\right)}} & 2 & 2 \, \theta_{0} + 1 &\textbf{(3)}\\
		\frac{2 \, A_{0} \overline{A_{0}}}{\theta_{0}^{2}} & \theta_{0} & \theta_{0} &\textbf{(4)}\\
		\frac{8 \, A_{0} \overline{A_{2}}}{4 \, {\left(\theta_{0}^{2} + 2 \, \theta_{0}^{} \, \right)}} & \theta_{0} & \theta_{0} + 2 &\textbf{(5)}\\
		\frac{24 \, {\left(\theta_{0}^{2} + \theta_{0}\right)} A_{0} \overline{A_{3}} + 3 \, {\left(\theta_{0}^{2} + 3 \, \theta_{0}\right)} \overline{A_{1}} \overline{C_{1}} + 2 \, {\left(\theta_{0}^{2} + 4 \, \theta_{0} + 3\right)} \overline{A_{0}} \overline{C_{2}}}{12 \, {\left(\theta_{0}^{4} + 4 \, \theta_{0}^{3} + 3 \, \theta_{0}^{2}\right)}} & \theta_{0} & \theta_{0} + 3 &\textbf{(6)}\\
		\frac{2 \, A_{1} \overline{A_{1}}}{\theta_{0}^{2} + 2 \, \theta_{0} + 1} & \theta_{0} + 1 & \theta_{0} + 1 &\textbf{(7)}\\
		\frac{8 \, A_{1} \theta_{0} \overline{A_{2}}}{4 \, {\left(\theta_{0}^{3} + 3 \, \theta_{0}^{2} + 2 \, \theta_{0}\right)}} & \theta_{0} + 1 & \theta_{0} + 2 &\textbf{(8)}\\
		\frac{8 \, A_{2} \theta_{0} \overline{A_{0}}}{4 \, {\left(\theta_{0}^{3} + 2 \, \theta_{0}^{2}\right)}} & \theta_{0} + 2 & \theta_{0} &\textbf{(9)}\\
		\frac{ 8 \, A_{2} \theta_{0} \overline{A_{1}}}{4 \, {\left(\theta_{0}^{3} + 3 \, \theta_{0}^{2} + 2 \, \theta_{0}\right)}} & \theta_{0} + 2 & \theta_{0} + 1 &\textbf{(10)}\\
		\frac{3 \, {\left(\theta_{0}^{2} + 3 \, \theta_{0}\right)} A_{1} C_{1} + 2 \, {\left(\theta_{0}^{2} + 4 \, \theta_{0} + 3\right)} A_{0} C_{2} + 24 \, {\left(\theta_{0}^{2} + \theta_{0}\right)} A_{3} \overline{A_{0}}}{12 \, {\left(\theta_{0}^{4} + 4 \, \theta_{0}^{3} + 3 \, \theta_{0}^{2}\right)}} & \theta_{0} + 3 & \theta_{0} &\textbf{(11)}\\
		\frac{A_{0} \overline{B_{1}} + A_{1} \overline{C_{1}}}{4 \, {\left(\theta_{0}^{2} + \theta_{0}\right)}} & 2 \, \theta_{0} + 1 & 2 &\textbf{(12)}\\
		\frac{{\left(\theta_{0}^{2} + 2 \, \theta_{0}\right)} A_{1}^{2} + 2 \, {\left(\theta_{0}^{2} + 2 \, \theta_{0} + 1\right)} A_{0} A_{2}}{\theta_{0}^{4} + 4 \, \theta_{0}^{3} + 5 \, \theta_{0}^{2} + 2 \, \theta_{0}} & 2 \, \theta_{0} + 2 & 0 &\textbf{(13)}\\
		\frac{2 \, {\left({\left(\theta_{0}^{2} + 3 \, \theta_{0}\right)} A_{1} A_{2} + {\left(\theta_{0}^{2} + 3 \, \theta_{0} + 2\right)} A_{0} A_{3}\right)}}{\theta_{0}^{4} + 6 \, \theta_{0}^{3} + 11 \, \theta_{0}^{2} + 6 \, \theta_{0}} & 2 \, \theta_{0} + 3 & 0 &\textbf{(14)}
		\end{dmatrix}
		\end{align*}
		\normalsize
		Furthermore, we can further simplify this expression with
		\begin{align*}
			\s{\vec{A}_1}{\vec{A}_1}+2\s{\vec{A}_0}{\vec{A}_2}=0,\quad \s{\vec{A}_1}{\vec{A}_2}+\s{\vec{A}_0}{\vec{A}_3}=0,\quad \s{\vec{A}_1}{\vec{C}_1}+\s{\vec{A}_0}{\vec{C}_2}=0.
		\end{align*} 
		We now define $\zeta_0,\zeta_1,\zeta_2\in\mathbb{C}$ such that
		\begin{align}\label{defzeta012}
		\left\{
			\begin{alignedat}{1}
			\zeta_0&=\s{\vec{A}_1}{\vec{A}_1}\\
			\zeta_1&=\s{\vec{A}_1}{\vec{A}_2}\\
			\zeta_2&=\s{\vec{A}_1}{\vec{C}_1}
			\end{alignedat}\right.
        \end{align}
		which gives
		\begin{align}\label{zeta01}
			\textbf{(13)}&=\frac{{\left(\theta_{0}^{2} + 2 \, \theta_{0}\right)} A_{1}^{2} + 2 \, {\left(\theta_{0}^{2} + 2 \, \theta_{0} + 1\right)} A_{0} A_{2}}{\theta_{0}^{4} + 4 \, \theta_{0}^{3} + 5 \, \theta_{0}^{2} + 2 \, \theta_{0}}=\frac{2\s{\vec{A}_0}{\vec{A}_2}}{\theta_{0}^{4} + 4 \, \theta_{0}^{3} + 5 \, \theta_{0}^{2} + 2 \, \theta_{0}}=-\frac{\s{\vec{A}_1}{\vec{A}_1}}{\theta_{0}^{4} + 4 \, \theta_{0}^{3} + 5 \, \theta_{0}^{2} + 2 \, \theta_{0}}\nonumber\\
			&=\frac{-\zeta_0}{\theta_{0}^{4} + 4 \, \theta_{0}^{3} + 5 \, \theta_{0}^{2} + 2 \, \theta_{0}}\nonumber\\
			\textbf{(14)}&=\frac{2 \, {\left({\left(\theta_{0}^{2} + 3 \, \theta_{0}\right)} A_{1} A_{2} + {\left(\theta_{0}^{2} + 3 \, \theta_{0} + 2\right)} A_{0} A_{3}\right)}}{\theta_{0}^{4} + 6 \, \theta_{0}^{3} + 11 \, \theta_{0}^{2} + 6 \, \theta_{0}}=\frac{4\s{\vec{A}_0}{\vec{A}_3}}{\theta_{0}^{4} + 6 \, \theta_{0}^{3} + 11 \, \theta_{0}^{2} + 6 \, \theta_{0}}=\frac{-4\s{\vec{A}_1}{\vec{A}_2}}{\theta_{0}^{4} + 6 \, \theta_{0}^{3} + 11 \, \theta_{0}^{2} + 6 \, \theta_{0}}\nonumber\\
			&=\frac{-4\,\zeta_1}{\theta_{0}^{4} + 6 \, \theta_{0}^{3} + 11 \, \theta_{0}^{2} + 6 \, \theta_{0}}
		\end{align}
		Now, recall that
	     \begin{align}\label{alpha29}
	 \left\{
	 \begin{alignedat}{1}
	 \alpha_1&=2\s{\bar{\vec{A}_0}}{\vec{A}_2}\\
	 \alpha_2&=\frac{1}{2\theta_0(\theta_0+1)}\s{\bar{\vec{A}_1}}{\vec{C}_1}\\
	 \alpha_3&=\frac{1}{12}\s{\vec{A}_1}{\vec{C}_1}+2\s{\bar{\vec{A}_0}}{\vec{A}_3}\\
	 \alpha_5&=2\s{\bar{\vec{A}_1}}{\vec{A}_2}\\
	 \alpha_6&=2\s{\bar{\vec{A}_1}}{\vec{A}_3}\\
	 \alpha_7&=\frac{1}{8\theta_0(\theta_0-4)}\s{\vec{C}_1}{\vec{C}_1}
	 \end{alignedat}\right.
	 \end{align}
	 so
	 \begin{align}\label{zeta2}
	 	\textbf{(11)}&=\frac{3 \, {\left(\theta_{0}^{2} + 3 \, \theta_{0}\right)} A_{1} C_{1} + 2 \, {\left(\theta_{0}^{2} + 4 \, \theta_{0} + 3\right)} A_{0} C_{2} + 24 \, {\left(\theta_{0}^{2} + \theta_{0}\right)} A_{3} \overline{A_{0}}}{12 \, {\left(\theta_{0}^{4} + 4 \, \theta_{0}^{3} + 3 \, \theta_{0}^{2}\right)}}\nonumber\\
	 	&=\frac{3 \, {\left(\theta_{0}^{2} + 3 \, \theta_{0}\right)} A_{1} C_{1} - 2 \, {\left(\theta_{0}^{2} + 4 \, \theta_{0} + 3\right)} A_{1} C_{1} + 24 \, {\left(\theta_{0}^{2} + \theta_{0}\right)} A_{3} \overline{A_{0}}}{12 \, {\left(\theta_{0}^{4} + 4 \, \theta_{0}^{3} + 3 \, \theta_{0}^{2}\right)}}\nonumber\\
	 	&=\frac{ \, \left(3\theta_{0}^{2} + 9 \, \theta_{0}-2\,\theta_0^2-8\theta_0-6\right) A_{1} C_{1}  + 12 \, {\left(\theta_{0}(\theta_{0}+1)\right)} (\alpha_3-\frac{1}{12}A_1\,C_1)}{12 \, {\left(\theta_{0}^{4} + 4 \, \theta_{0}^{3} + 3 \, \theta_{0}^{2}\right)}}\nonumber\\
	 	&=\frac{-6\,A_1C_1+12\theta_0(\theta_0+1)\alpha_3}{12 \, {\left(\theta_{0}^{4} + 4 \, \theta_{0}^{3} + 3 \, \theta_{0}^{2}\right)}}=\frac{-A_1C_1+2\theta_0(\theta_0+1)\alpha_3}{2{\left(\theta_{0}^{4} + 4 \, \theta_{0}^{3} + 3 \, \theta_{0}^{2}\right)}}=\frac{-\zeta_2+2\theta_0(\theta_0+1)\alpha_3}{2{\left(\theta_{0}^{4} + 4 \, \theta_{0}^{3} + 3 \, \theta_{0}^{2}\right)}}
	 \end{align}
	 Then
	 \begin{align*}
	 	\s{\vec{A}_0}{\bar{\vec{B}_1}}=\s{\vec{A}_0}{-2\s{\vec{A}_1}{\bar{\vec{C}_1}}\vec{A}_0}=-\s{\vec{A}_1}{\bar{\vec{C}_1}}
	 \end{align*}
	 so
	 \begin{align*}
	 	\textbf{(12)}&=
	 	\frac{A_{0} \overline{B_{1}} + A_{1} \overline{C_{1}}}{4 \, {\left(\theta_{0}^{2} + \theta_{0}\right)}}=0 
	 \end{align*}
	 Now we directly have as $\alpha_5=2\s{\bar{\vec{A}_1}}{\vec{A}_2}$ the identity
	 \begin{align*}
	    \textbf{(10)}=\frac{ 8 \, A_{2} \theta_{0} \overline{A_{1}}}{4 \, {\left(\theta_{0}^{3} + 3 \, \theta_{0}^{2} + 2 \, \theta_{0}\right)}}=\frac{\alpha_5}{\theta_0^2+3\theta_0+2}=\frac{\alpha_5}{(\theta_0+1)(\theta_0+2)}.
	 \end{align*}
	 Likewise,
	 \begin{align*}
	 	\textbf{(9)}=\frac{8 \, A_{2} \theta_{0} \overline{A_{0}}}{4 \, {\left(\theta_{0}^{3} + 2 \, \theta_{0}^{2}\right)}}=\frac{\alpha_1}{\theta_0(\theta_0+2)}.
	 \end{align*}
	 Therefore, we finally obtain
	 \begin{align*}
	 	&|\phi(z)|^2=\begin{dmatrix}
	 	\frac{-\bar{\zeta_0}}{\theta_{0}^{4} + 4 \, \theta_{0}^{3} + 5 \, \theta_{0}^{2} + 2 \, \theta_{0}} & 0 & 2 \, \theta_{0} + 2 &\textbf{(1)}\\
	 	\frac{-4\bar{\zeta_1}}{\theta_{0}^{4} + 6 \, \theta_{0}^{3} + 11 \, \theta_{0}^{2} + 6 \, \theta_{0}} & 0 & 2 \, \theta_{0} + 3 &\textbf{(2)}\\
	 	\frac{1}{\theta_{0}^{2}} & \theta_{0} & \theta_{0} &\textbf{(4)}\\
	 	\frac{\bar{\alpha_1}}{{\theta_{0}(\theta_0 + 2)}} & \theta_{0} & \theta_{0} + 2 &\textbf{(5)}\\
	 	\frac{-\bar{\zeta_2}+2\theta_0(\theta_0+1)\bar{\alpha_3}}{2{\left(\theta_{0}^{4} + 4 \, \theta_{0}^{3} + 3 \, \theta_{0}^{2}\right)}} & \theta_{0} & \theta_{0} + 3 &\textbf{(6)}\\
	 	\frac{2 \, |A_{1}|^2}{(\theta_{0}+1)^2} & \theta_{0} + 1 & \theta_{0} + 1 &\textbf{(7)}\\
	 	\frac{\bar{\alpha_5}}{(\theta_0+1)(\theta_0+2)} & \theta_{0} + 1 & \theta_{0} + 2 &\textbf{(8)}\\
	 	\frac{\alpha_1}{\theta_0(\theta_0+2)} & \theta_{0} + 2 & \theta_{0} &\textbf{(9)}\\
	 	\frac{\alpha_5}{(\theta_0+1)(\theta_0+2)} & \theta_{0} + 2 & \theta_{0} + 1 &\textbf{(10)}\\
	 	\frac{-\zeta_2+2\theta_0(\theta_0+1)\alpha_3}{2{\left(\theta_{0}^{4} + 4 \, \theta_{0}^{3} + 3 \, \theta_{0}^{2}\right)}} & \theta_{0} + 3 & \theta_{0} &\textbf{(11)}\\
	 	\frac{-\zeta_0}{\theta_{0}^{4} + 4 \, \theta_{0}^{3} + 5 \, \theta_{0}^{2} + 2 \, \theta_{0}} & 2 \, \theta_{0} + 2 & 0 &\textbf{(13)}\\
	 	\frac{-4\zeta_1}{\theta_{0}^{4} + 6 \, \theta_{0}^{3} + 11 \, \theta_{0}^{2} + 6 \, \theta_{0}} & 2 \, \theta_{0} + 3 & 0 &\textbf{(14)}
	 	\end{dmatrix}\\
	 	&=\frac{1}{\theta_0^2}|z|^{2\theta_0}+\frac{2|\vec{A}_1|^2}{(\theta_0+1)^2}|z|^{2\theta_0+2}+2\,\Re\bigg(\frac{\,\alpha_1}{\theta_0(\theta_0+2)}z^{\theta_0+2}\z^{\theta_0}+\frac{\alpha_5}{(\theta_0+1)(\theta_0+2)}z^{\theta_0+2}\z^{\theta_0+1}+\\
	 	&+\frac{-\zeta_2+2\theta_0(\theta_0+1)\alpha_3}{2(\theta_0^4+4\theta_0^3+3\theta_0^2)}z^{\theta_0+3}\z^{\theta_0}-\frac{\,\zeta_0}{\theta_0^4+4\theta_0^3+5\theta_0^2+2\theta_0}z^{2\theta_0+2}-\frac{4\zeta_1}{\theta_0^4+6\theta_0^3+11\theta_0^2+6\theta_0}z^{2\theta_0+3}\bigg)+O(|z|^{2\theta_0+4})\\
	 	&=|z|^{2\theta_0}\bigg\{\frac{1}{\theta_0^2}+\frac{2|\vec{A}_1|^2}{(\theta_0+1)^2}|z|^2+2\,\Re\bigg(\frac{\,\alpha_1}{\theta_0(\theta_0+2)}z^{2}+\frac{\alpha_5}{(\theta_0+1)(\theta_0+2)}z^{2}\z+
	 	+\frac{-\zeta_2+2\theta_0(\theta_0+1)\alpha_3}{2(\theta_0^4+4\theta_0^3+3\theta_0^2)}z^{3}\\
	 	&-\frac{\,\zeta_0}{\theta_0^4+4\theta_0^3+5\theta_0^2+2\theta_0}z^{\theta_0+2}\z^{-\theta_0}-\frac{4\zeta_1}{\theta_0^4+6\theta_0^3+11\theta_0^2+6\theta_0}z^{\theta_0+3}\z^{-\theta_0} \bigg)+O(|z|^4) \bigg\}
	 \end{align*}
	 to be compare with the right-hand size in the next equation
	 \begin{align*}
	 	|\phi(z)|^2=\begin{dmatrix}
	 	\frac{-\bar{\zeta_0}}{\theta_{0}^{4} + 4 \, \theta_{0}^{3} + 5 \, \theta_{0}^{2} + 2 \, \theta_{0}} & 0 & 2 \, \theta_{0} + 2 &\textbf{(1)}\\
	 	\frac{-4\bar{\zeta_1}}{\theta_{0}^{4} + 6 \, \theta_{0}^{3} + 11 \, \theta_{0}^{2} + 6 \, \theta_{0}} & 0 & 2 \, \theta_{0} + 3 &\textbf{(2)}\\
	 	\frac{1}{\theta_{0}^{2}} & \theta_{0} & \theta_{0} &\textbf{(4)}\\
	 	\frac{\bar{\alpha_1}}{{\theta_{0}(\theta_0 + 2)}} & \theta_{0} & \theta_{0} + 2 &\textbf{(5)}\\
	 	\frac{-\bar{\zeta_2}+2\theta_0(\theta_0+1)\bar{\alpha_3}}{2{\left(\theta_{0}^{4} + 4 \, \theta_{0}^{3} + 3 \, \theta_{0}^{2}\right)}} & \theta_{0} & \theta_{0} + 3 &\textbf{(6)}\\
	 	\frac{2 \, |A_{1}|^2}{(\theta_{0}+1)^2} & \theta_{0} + 1 & \theta_{0} + 1 &\textbf{(7)}\\
	 	\frac{\bar{\alpha_5}}{(\theta_0+1)(\theta_0+2)} & \theta_{0} + 1 & \theta_{0} + 2 &\textbf{(8)}\\
	 	\frac{\alpha_1}{\theta_0(\theta_0+2)} & \theta_{0} + 2 & \theta_{0} &\textbf{(9)}\\
	 	\frac{\alpha_5}{(\theta_0+1)(\theta_0+2)} & \theta_{0} + 2 & \theta_{0} + 1 &\textbf{(10)}\\
	 	\frac{-\zeta_2+2\theta_0(\theta_0+1)\alpha_3}{2{\left(\theta_{0}^{4} + 4 \, \theta_{0}^{3} + 3 \, \theta_{0}^{2}\right)}} & \theta_{0} + 3 & \theta_{0} &\textbf{(11)}\\
	 	\frac{-\zeta_0}{\theta_{0}^{4} + 4 \, \theta_{0}^{3} + 5 \, \theta_{0}^{2} + 2 \, \theta_{0}} & 2 \, \theta_{0} + 2 & 0 &\textbf{(13)}\\
	 	\frac{-4\zeta_1}{\theta_{0}^{4} + 6 \, \theta_{0}^{3} + 11 \, \theta_{0}^{2} + 6 \, \theta_{0}} & 2 \, \theta_{0} + 3 & 0 &\textbf{(14)}
	 	\end{dmatrix}
	 	=\begin{dmatrix}
	 	-\frac{\overline{\zeta_{0}}}{\theta_{0}^{4} + 4 \, \theta_{0}^{3} + 5 \, \theta_{0}^{2} + 2 \, \theta_{0}} & 0 & 2 \, \theta_{0} + 2 \\
	 	-\frac{4 \, \overline{\zeta_{1}}}{\theta_{0}^{4} + 6 \, \theta_{0}^{3} + 11 \, \theta_{0}^{2} + 6 \, \theta_{0}} & 0 & 2 \, \theta_{0} + 3 \\
	 	\frac{1}{\theta_{0}^{2}} & \theta_{0} & \theta_{0} \\
	 	\frac{\overline{\alpha_{1}}}{{\left(\theta_{0} + 2\right)} \theta_{0}} & \theta_{0} & \theta_{0} + 2 \\
	 	\frac{2 \, {\left(\theta_{0} + 1\right)} \theta_{0} \overline{\alpha_{3}} - \overline{\zeta_{2}}}{2 \, {\left(\theta_{0}^{4} + 4 \, \theta_{0}^{3} + 3 \, \theta_{0}^{2}\right)}} & \theta_{0} & \theta_{0} + 3 \\
	 	\frac{2 \, {\left| A_{1} \right|}^{2}}{{\left(\theta_{0} + 1\right)}^{2}} & \theta_{0} + 1 & \theta_{0} + 1 \\
	 	\frac{\overline{\alpha_{5}}}{{\left(\theta_{0} + 2\right)} {\left(\theta_{0} + 1\right)}} & \theta_{0} + 1 & \theta_{0} + 2 \\
	 	\frac{\alpha_{1}}{{\left(\theta_{0} + 2\right)} \theta_{0}} & \theta_{0} + 2 & \theta_{0} \\
	 	\frac{\alpha_{5}}{{\left(\theta_{0} + 2\right)} {\left(\theta_{0} + 1\right)}} & \theta_{0} + 2 & \theta_{0} + 1 \\
	 	\frac{2 \, \alpha_{3} {\left(\theta_{0} + 1\right)} \theta_{0} - \zeta_{2}}{2 \, {\left(\theta_{0}^{4} + 4 \, \theta_{0}^{3} + 3 \, \theta_{0}^{2}\right)}} & \theta_{0} + 3 & \theta_{0} \\
	 	-\frac{\zeta_{0}}{\theta_{0}^{4} + 4 \, \theta_{0}^{3} + 5 \, \theta_{0}^{2} + 2 \, \theta_{0}} & 2 \, \theta_{0} + 2 & 0 \\
	 	-\frac{4 \, \zeta_{1}}{\theta_{0}^{4} + 6 \, \theta_{0}^{3} + 11 \, \theta_{0}^{2} + 6 \, \theta_{0}} & 2 \, \theta_{0} + 3 & 0
	 	\end{dmatrix}
	 \end{align*}
	 Another test to check that our quantity is real.
	 \begin{align*}
	 	|\phi(z)|^2=\begin{dmatrix}
	 	\frac{-\bar{\zeta_0}}{\theta_{0}^{4} + 4 \, \theta_{0}^{3} + 5 \, \theta_{0}^{2} + 2 \, \theta_{0}} & 0 & 2 \, \theta_{0} + 2 &\textbf{(1)}\\
	 	\frac{-4\bar{\zeta_1}}{\theta_{0}^{4} + 6 \, \theta_{0}^{3} + 11 \, \theta_{0}^{2} + 6 \, \theta_{0}} & 0 & 2 \, \theta_{0} + 3 &\textbf{(2)}\\
	 	\frac{1}{\theta_{0}^{2}} & \theta_{0} & \theta_{0} &\textbf{(4)}\\
	 	\frac{\bar{\alpha_1}}{{\theta_{0}(\theta_0 + 2)}} & \theta_{0} & \theta_{0} + 2 &\textbf{(5)}\\
	 	\frac{-\bar{\zeta_2}+2\theta_0(\theta_0+1)\bar{\alpha_3}}{2{\left(\theta_{0}^{4} + 4 \, \theta_{0}^{3} + 3 \, \theta_{0}^{2}\right)}} & \theta_{0} & \theta_{0} + 3 &\textbf{(6)}\\
	 	\frac{2 \, |A_{1}|^2}{(\theta_{0}+1)^2} & \theta_{0} + 1 & \theta_{0} + 1 &\textbf{(7)}\\
	 	\frac{\bar{\alpha_5}}{(\theta_0+1)(\theta_0+2)} & \theta_{0} + 1 & \theta_{0} + 2 &\textbf{(8)}\\
	 	\frac{\alpha_1}{\theta_0(\theta_0+2)} & \theta_{0} + 2 & \theta_{0} &\textbf{(9)}\\
	 	\frac{\alpha_5}{(\theta_0+1)(\theta_0+2)} & \theta_{0} + 2 & \theta_{0} + 1 &\textbf{(10)}\\
	 	\frac{-\zeta_2+2\theta_0(\theta_0+1)\alpha_3}{2{\left(\theta_{0}^{4} + 4 \, \theta_{0}^{3} + 3 \, \theta_{0}^{2}\right)}} & \theta_{0} + 3 & \theta_{0} &\textbf{(11)}\\
	 	\frac{-\zeta_0}{\theta_{0}^{4} + 4 \, \theta_{0}^{3} + 5 \, \theta_{0}^{2} + 2 \, \theta_{0}} & 2 \, \theta_{0} + 2 & 0 &\textbf{(13)}\\
	 	\frac{-4\zeta_1}{\theta_{0}^{4} + 6 \, \theta_{0}^{3} + 11 \, \theta_{0}^{2} + 6 \, \theta_{0}} & 2 \, \theta_{0} + 3 & 0 &\textbf{(14)}
	 	\end{dmatrix}=\begin{dmatrix}
	 	-\frac{\overline{\zeta_{0}}}{\theta_{0}^{4} + 4 \, \theta_{0}^{3} + 5 \, \theta_{0}^{2} + 2 \, \theta_{0}} & 0 & 2 \, \theta_{0} + 2 \\
	 	-\frac{4 \, \overline{\zeta_{1}}}{\theta_{0}^{4} + 6 \, \theta_{0}^{3} + 11 \, \theta_{0}^{2} + 6 \, \theta_{0}} & 0 & 2 \, \theta_{0} + 3 \\
	 	\frac{1}{\theta_{0}^{2}} & \theta_{0} & \theta_{0} \\
	 	\frac{\overline{\alpha_{1}}}{\theta_{0}^{2} + 2 \, \theta_{0}} & \theta_{0} & \theta_{0} + 2 \\
	 	\frac{2 \, \theta_{0}^{2} \overline{\alpha_{3}} + 2 \, \theta_{0} \overline{\alpha_{3}} - \overline{\zeta_{2}}}{2 \, {\left(\theta_{0}^{4} + 4 \, \theta_{0}^{3} + 3 \, \theta_{0}^{2}\right)}} & \theta_{0} & \theta_{0} + 3 \\
	 	\frac{2 \, {\left| A_{1} \right|}^{2}}{\theta_{0}^{2} + 2 \, \theta_{0} + 1} & \theta_{0} + 1 & \theta_{0} + 1 \\
	 	\frac{\overline{\alpha_{5}}}{\theta_{0}^{2} + 3 \, \theta_{0} + 2} & \theta_{0} + 1 & \theta_{0} + 2 \\
	 	\frac{\alpha_{1}}{\theta_{0}^{2} + 2 \, \theta_{0}} & \theta_{0} + 2 & \theta_{0} \\
	 	\frac{\alpha_{5}}{\theta_{0}^{2} + 3 \, \theta_{0} + 2} & \theta_{0} + 2 & \theta_{0} + 1 \\
	 	\frac{2 \, \alpha_{3} \theta_{0}^{2} + 2 \, \alpha_{3} \theta_{0} - \zeta_{2}}{2 \, {\left(\theta_{0}^{4} + 4 \, \theta_{0}^{3} + 3 \, \theta_{0}^{2}\right)}} & \theta_{0} + 3 & \theta_{0} \\
	 	-\frac{\zeta_{0}}{\theta_{0}^{4} + 4 \, \theta_{0}^{3} + 5 \, \theta_{0}^{2} + 2 \, \theta_{0}} & 2 \, \theta_{0} + 2 & 0 \\
	 	-\frac{4 \, \zeta_{1}}{\theta_{0}^{4} + 6 \, \theta_{0}^{3} + 11 \, \theta_{0}^{2} + 6 \, \theta_{0}} & 2 \, \theta_{0} + 3 & 0
	 	\end{dmatrix}
	 \end{align*}
	 
	 \section{Development of $g^{-1}\otimes\h_0$ up to order $4$}
	 We directly compute
	 \begin{align*}
	 	g^{-1}\otimes\h_0=\begin{dmatrix}
	 	2 & A_{1} & 0 & -\theta_{0} + 1 \\
	 	-4 \, {\left| A_{1} \right|}^{2} & A_{0} & 0 & -\theta_{0} + 2 \\
	 	\frac{1}{4} & \overline{B_{1}} & 0 & -\theta_{0} + 3 \\
	 	-2 \, \overline{\alpha_{5}} & A_{0} & 0 & -\theta_{0} + 3 \\
	 	-2 \, \overline{\alpha_{1}} & A_{1} & 0 & -\theta_{0} + 3 \\
	 	\frac{1}{6} & \overline{B_{3}} & 0 & -\theta_{0} + 4 \\
	 	-2 \, \overline{\alpha_{3}} & A_{1} & 0 & -\theta_{0} + 4 \\
	 	-\frac{1}{6} \, {\left| A_{1} \right|}^{2} & \overline{C_{1}} & 0 & -\theta_{0} + 4 \\
	 	8 \, {\left| A_{1} \right|}^{2} \overline{\alpha_{1}} - 2 \, \overline{\alpha_{6}} & A_{0} & 0 & -\theta_{0} + 4 \\
	 	4 & A_{2} & 1 & -\theta_{0} + 1 \\
	 	-4 \, \alpha_{1} & A_{0} & 1 & -\theta_{0} + 1 \\
	 	-4 \, \alpha_{5} & A_{0} & 1 & -\theta_{0} + 2 \\
	 	-8 \, {\left| A_{1} \right|}^{2} & A_{1} & 1 & -\theta_{0} + 2 \\
	 	\frac{1}{2} & \overline{B_{2}} & 1 & -\theta_{0} + 3 \\
	 	-4 \, \overline{\alpha_{5}} & A_{1} & 1 & -\theta_{0} + 3 \\
	 	-4 \, \overline{\alpha_{1}} & A_{2} & 1 & -\theta_{0} + 3 \\
	 	16 \, {\left| A_{1} \right|}^{4} + 8 \, \alpha_{1} \overline{\alpha_{1}} - 4 \, \beta & A_{0} & 1 & -\theta_{0} + 3 \\
	 	6 & A_{3} & 2 & -\theta_{0} + 1 \\
	 	-6 \, \alpha_{3} & A_{0} & 2 & -\theta_{0} + 1 \\
	 	-6 \, \alpha_{1} & A_{1} & 2 & -\theta_{0} + 1 \\
	 	-6 \, \alpha_{5} & A_{1} & 2 & -\theta_{0} + 2 \\
	 	-12 \, {\left| A_{1} \right|}^{2} & A_{2} & 2 & -\theta_{0} + 2 \\
	 	24 \, \alpha_{1} {\left| A_{1} \right|}^{2} - 6 \, \alpha_{6} & A_{0} & 2 & -\theta_{0} + 2 
	 		 	\end{dmatrix}
	 	\begin{dmatrix}
	 	8 & A_{4} & 3 & -\theta_{0} + 1 \\
	 	-8 \, \alpha_{3} & A_{1} & 3 & -\theta_{0} + 1 \\
	 	-8 \, \alpha_{1} & A_{2} & 3 & -\theta_{0} + 1 \\
	 	8 \, \alpha_{1}^{2} - 8 \, \alpha_{4} & A_{0} & 3 & -\theta_{0} + 1 \\
	 	-\frac{\theta_{0} - 2}{2 \, \theta_{0}} & E_{1} & -2 \, \theta_{0} + 3 & \theta_{0} + 1 \\
	 	-\frac{\theta_{0} - 2}{2 \, \theta_{0}} & C_{1} & -\theta_{0} + 1 & 1 \\
	 	-\frac{\theta_{0} - 2}{2 \, {\left(\theta_{0} + 1\right)}} & B_{1} & -\theta_{0} + 1 & 2 \\
	 	2 \, \alpha_{2} \theta_{0} - 4 \, \alpha_{2} & A_{0} & -\theta_{0} + 1 & 2 \\
	 	-\frac{\theta_{0} - 2}{2 \, {\left(\theta_{0} + 2\right)}} & B_{2} & -\theta_{0} + 1 & 3 \\
	 	\frac{\theta_{0} \overline{\alpha_{1}} - 2 \, \overline{\alpha_{1}}}{\theta_{0}^{2} + 2 \, \theta_{0}} & C_{1} & -\theta_{0} + 1 & 3 \\
	 	2 \, \alpha_{9} \theta_{0} - 4 \, \alpha_{9} & A_{0} & -\theta_{0} + 1 & 3 \\
	 	-\frac{\theta_{0} - 3}{2 \, \theta_{0}} & C_{2} & -\theta_{0} + 2 & 1 \\
	 	-\frac{\theta_{0} - 3}{2 \, {\left(\theta_{0} + 1\right)}} & B_{3} & -\theta_{0} + 2 & 2 \\
	 	\frac{{\left(\theta_{0} - 3\right)} {\left| A_{1} \right|}^{2}}{\theta_{0}^{2} + \theta_{0}} & C_{1} & -\theta_{0} + 2 & 2 \\
	 	2 \, \alpha_{8} \theta_{0} - 6 \, \alpha_{8} & A_{0} & -\theta_{0} + 2 & 2 \\
	 	2 \, \alpha_{2} \theta_{0} - 6 \, \alpha_{2} & A_{1} & -\theta_{0} + 2 & 2 \\
	 	-\frac{\theta_{0} - 4}{2 \, \theta_{0}} & C_{3} & -\theta_{0} + 3 & 1 \\
	 	2 \, \alpha_{7} \theta_{0} - 8 \, \alpha_{7} & A_{0} & -\theta_{0} + 3 & 1 \\
	 	-2 \, \theta_{0} \overline{\alpha_{2}} - 2 \, \overline{\alpha_{2}} & A_{0} & \theta_{0} & -2 \, \theta_{0} + 3 \\
	 	-2 \, \theta_{0} \overline{\alpha_{8}} - 2 \, \overline{\alpha_{8}} & A_{0} & \theta_{0} & -2 \, \theta_{0} + 4 \\
	 	-\frac{\theta_{0}}{2 \, {\left(\theta_{0} - 4\right)}} & \overline{E_{1}} & \theta_{0} - 1 & -2 \, \theta_{0} + 5 \\
	 	-2 \, \theta_{0} \overline{\alpha_{7}} & A_{0} & \theta_{0} - 1 & -2 \, \theta_{0} + 5 \\
	 	-2 \, \theta_{0} \overline{\alpha_{9}} - 4 \, \overline{\alpha_{9}} & A_{0} & \theta_{0} + 1 & -2 \, \theta_{0} + 3 \\
	 	-2 \, \theta_{0} \overline{\alpha_{2}} - 4 \, \overline{\alpha_{2}} & A_{1} & \theta_{0} + 1 & -2 \, \theta_{0} + 3
	 	\end{dmatrix}
	 \end{align*}
	 and indeed we have a development of $g^{-1}\otimes \h_0$ up to an error in $O(|z|^{5-\theta_0})$.
		
		\section{Development of $\vec{\alpha}$ up to order $3$}
		Now we need to do some preformating, and we first compute
		\small
		\begin{align*}
			\vec{\alpha}=\begin{dmatrix}
			\frac{1}{2} & \overline{B_{1}} & 0 & -\theta_{0} + 2 \\
			2 \, A_{1} \overline{C_{1}} & \overline{A_{0}} & 0 & -\theta_{0} + 2 \\
			\frac{1}{2} & \overline{B_{3}} & 0 & -\theta_{0} + 3 \\
			2 \, A_{1} \overline{C_{1}} & \overline{A_{1}} & 0 & -\theta_{0} + 3 \\
			-\ccancel{4 \, A_{0} {\left| A_{1} \right|}^{2} \overline{C_{1}}} + 2 \, A_{1} \overline{C_{2}} & \overline{A_{0}} & 0 & -\theta_{0} + 3 \\
			1 & \overline{B_{2}} & 1 & -\theta_{0} + 2 \\
			\frac{1}{2} \, C_{1} \overline{C_{1}} & A_{0} & 1 & -\theta_{0} + 2 \\
			\ccancel{2 \, A_{1} \overline{B_{1}}} - 4 \, {\left(\ccancel{A_{0} \alpha_{1}} - A_{2}\right)} \overline{C_{1}} & \overline{A_{0}} & 1 & -\theta_{0} + 2 \\
			-\theta_{0} + 2 & E_{1} & -2 \, \theta_{0} + 3 & \theta_{0} \\
			-\frac{C_{1}^{2} {\left(\theta_{0} - 2\right)}}{2 \, \theta_{0}} & \overline{A_{0}} & -2 \, \theta_{0} + 3 & \theta_{0} \\
			-\frac{1}{2} \, \theta_{0} + 1 & C_{1} & -\theta_{0} + 1 & 0 \\
			-\frac{1}{2} \, \theta_{0} + 1 & B_{1} & -\theta_{0} + 1 & 1 \\
			-\frac{1}{2} \, \theta_{0} + 1 & B_{2} & -\theta_{0} + 1 & 2 \\
			-\frac{C_{1} {\left(\theta_{0} - 2\right)} \overline{C_{1}}}{2 \, \theta_{0}} & \overline{A_{0}} & -\theta_{0} + 1 & 2 \\
			-\frac{1}{2} \, \theta_{0} + \frac{3}{2} & C_{2} & -\theta_{0} + 2 & 0 \\
			2 \, A_{1} C_{1} & \overline{A_{0}} & -\theta_{0} + 2 & 0 \\
			-\frac{1}{2} \, \theta_{0} + \frac{3}{2} & B_{3} & -\theta_{0} + 2 & 1 \\
			2 \, A_{1} C_{1} & \overline{A_{1}} & -\theta_{0} + 2 & 1 \\
			-\ccancel{4 \, A_{0} C_{1} {\left| A_{1} \right|}^{2}} + \ccancel{2 \, A_{1} B_{1}} & \overline{A_{0}} & -\theta_{0} + 2 & 1 \\
			-\frac{1}{2} \, \theta_{0} + 2 & C_{3} & -\theta_{0} + 3 & 0 \\
			\frac{1}{4} \, C_{1}^{2} & A_{0} & -\theta_{0} + 3 & 0 \\
			-4 \, {\left(\ccancel{A_{0} \alpha_{1}} - A_{2}\right)} C_{1} + 2 \, A_{1} C_{2} & \overline{A_{0}} & -\theta_{0} + 3 & 0 \\
			\frac{1}{2} \, \theta_{0} & \overline{E_{1}} & \theta_{0} - 1 & -2 \, \theta_{0} + 4 \\
			\frac{1}{4} \, \overline{C_{1}}^{2} & A_{0} & \theta_{0} - 1 & -2 \, \theta_{0} + 4
			\end{dmatrix}
			=\begin{dmatrix}
			\frac{1}{2} & \overline{B_{1}} & 0 & -\theta_{0} + 2 & \textbf{(1)}\\
			2 \, A_{1} \overline{C_{1}} & \overline{A_{0}} & 0 & -\theta_{0} + 2 & \textbf{(2)}\\
			\frac{1}{2} & \overline{B_{3}} & 0 & -\theta_{0} + 3 & \textbf{(3)}\\
			2 \, A_{1} \overline{C_{1}} & \overline{A_{1}} & 0 & -\theta_{0} + 3 & \textbf{(4)}\\
			 2 \, A_{1} \overline{C_{2}} & \overline{A_{0}} & 0 & -\theta_{0} + 3 & \textbf{(5)}\\
			1 & \overline{B_{2}} & 1 & -\theta_{0} + 2 & \textbf{(6)}\\
			\frac{1}{2} \, C_{1} \overline{C_{1}} & A_{0} & 1 & -\theta_{0} + 2 & \textbf{(7)}\\
		4 \, A_{2} \overline{C_{1}} & \overline{A_{0}} & 1 & -\theta_{0} + 2 & \textbf{(8)}\\
			-\theta_{0} + 2 & E_{1} & -2 \, \theta_{0} + 3 & \theta_{0} & \textbf{(9)}\\
			-\frac{C_{1}^{2} {\left(\theta_{0} - 2\right)}}{2 \, \theta_{0}} & \overline{A_{0}} & -2 \, \theta_{0} + 3 & \theta_{0} & \textbf{(10)}\\
			-\frac{1}{2} \, \theta_{0} + 1 & C_{1} & -\theta_{0} + 1 & 0 & \textbf{(11)}\\
			-\frac{1}{2} \, \theta_{0} + 1 & B_{1} & -\theta_{0} + 1 & 1 & \textbf{(12)}\\
			-\frac{1}{2} \, \theta_{0} + 1 & B_{2} & -\theta_{0} + 1 & 2 & \textbf{(13)}\\
			-\frac{C_{1} {\left(\theta_{0} - 2\right)} \overline{C_{1}}}{2 \, \theta_{0}} & \overline{A_{0}} & -\theta_{0} + 1 & 2 & \textbf{(14)}\\
			-\frac{1}{2} \, \theta_{0} + \frac{3}{2} & C_{2} & -\theta_{0} + 2 & 0 & \textbf{(15)}\\
			2 \, A_{1} C_{1} & \overline{A_{0}} & -\theta_{0} + 2 & 0 & \textbf{(16)}\\
			-\frac{1}{2} \, \theta_{0} + \frac{3}{2} & B_{3} & -\theta_{0} + 2 & 1 & \textbf{(17)}\\
			2 \, A_{1} C_{1} & \overline{A_{1}} & -\theta_{0} + 2 & 1 & \textbf{(18)}\\
			-\frac{1}{2} \, \theta_{0} + 2 & C_{3} & -\theta_{0} + 3 & 0 & \textbf{(19)}\\
			\frac{1}{4} \, C_{1}^{2} & A_{0} & -\theta_{0} + 3 & 0 & \textbf{(20)}\\
			4 \, A_{2} C_{1} + 2 \, A_{1} C_{2} & \overline{A_{0}} & -\theta_{0} + 3 & 0 & \textbf{(21)}\\
			\frac{1}{2} \, \theta_{0} & \overline{E_{1}} & \theta_{0} - 1 & -2 \, \theta_{0} + 4 & \textbf{(22)}\\
			\frac{1}{4} \, \overline{C_{1}}^{2} & A_{0} & \theta_{0} - 1 & -2 \, \theta_{0} + 4& \textbf{(23)}
			\end{dmatrix}
		\end{align*}
		Most of the coefficients here cancel. Indeed, recall that
			\begin{align}\label{invref1}
			\left\{
			\begin{alignedat}{1}
			\vec{B}_1&=-2\s{\bar{\vec{A}_1}}{\vec{C}_1}\vec{A}_0\\
			\vec{B}_2&=-\frac{(\theta_0+2)}{4\theta_0}|\vec{C}_1|^2\bar{\vec{A}_0}-2\s{\bar{\vec{A}_2}}{\vec{C}_1}\vec{A}_0\\
			\vec{B}_3&=-2\s{\bar{\vec{A}_1}}{\vec{C}_1}\vec{A}_1+\frac{2}{\theta_0-3}\s{\vec{A}_1}{\vec{C}_1}\bar{\vec{A}_1}-2\s{\bar{\vec{A}_1}}{\vec{C}_2}\vec{A}_0\\
			\vec{E}_1&=-\frac{1}{2\theta_0}\s{\vec{C}_1}{\vec{C}_1}\bar{\vec{A}_0}.
			\end{alignedat}\right.
			\end{align}
		 Therefore, we have
		 \begin{align*}
		 	\begin{dmatrix}
		 	\frac{1}{2} & \overline{B_{1}} & 0 & -\theta_{0} + 2 & \textbf{(1)}\\
		 	2 \, A_{1} \overline{C_{1}} & \overline{A_{0}} & 0 & -\theta_{0} + 2& \textbf{(2)} 
		 	\end{dmatrix}
		 	&=\begin{dmatrix}
		 	- A_{1} \overline{C_{1}} & \overline{A_{0}} & 0 & -\theta_{0} + 2 & \textbf{(1)}\\
		 	2 \, A_{1} \overline{C_{1}} & \overline{A_{0}} & 0 & -\theta_{0} + 2 & \textbf{(2)}\\
		 	\end{dmatrix}
		 	=\begin{dmatrix}
		 	A_{1} \overline{C_{1}} & \overline{A_{0}} & 0 & -\theta_{0} + 2 & \textbf{(1)} 
		 	\end{dmatrix}\\
		 	&=\begin{dmatrix}
		 	2\theta_0(\theta_0+1)\bar{\alpha_2} & \bar{A_0} & 0 & -\theta_0+2 & \textbf{(1)}
		 	\end{dmatrix}
		 \end{align*}
		 and
		 \begin{align*}
		 	&\begin{dmatrix}
		 	\frac{1}{2} & \overline{B_{3}} & 0 & -\theta_{0} + 3 & \textbf{(3)}\\
		 	2 \, A_{1} \overline{C_{1}} & \overline{A_{1}} & 0 & -\theta_{0} + 3 & \textbf{(4)}\\
		 	2 \, A_{1} \overline{C_{2}} & \overline{A_{0}} & 0 & -\theta_{0} + 3 & \textbf{(5)}
		 	\end{dmatrix}
		 	=\begin{dmatrix}
		 	\frac{1}{\theta_0-3}\bar{\s{\vec{A}_1}{\vec{C}_1}} & A_1 & 0 & -\theta_0+3& \textbf{(3)}\\
		 	\, A_{1} \overline{C_{1}} & \overline{A_{1}} & 0 & -\theta_{0} + 3 & \textbf{(4)}\\
		 	 \, A_{1} \overline{C_{2}} & \overline{A_{0}} & 0 & -\theta_{0} + 3 & \textbf{(5)}
		 	\end{dmatrix}\\
		 	&		 	
		 	=\begin{dmatrix}
		 	\frac{\bar{\zeta_2}}{\theta_0-3} & A_1 & 0 & -\theta_0+3& \textbf{(3)}\\
		 	2\theta_0(\theta_0+1)\bar{\alpha_2} & \bar{A_1} & 0 & -\theta_0+3 & \textbf{(4)}\\
		 	\bar{\zeta_3}  & \bar{A_0} & 0 & -\theta_0+3& \textbf{(5)}
		 	\end{dmatrix}
		 \end{align*}
		 if 
		 \begin{align}\label{defzeta3}
		 	\zeta_3=\s{\bar{\vec{A}_1}}{\vec{C}_2}.
		 \end{align}
		 Then
		 \begin{align*}
		 	\begin{dmatrix}
		 	1 & \overline{B_{2}} & 1 & -\theta_{0} + 2 & \textbf{(6)}\\
		 	\frac{1}{2} \, C_{1} \overline{C_{1}} & A_{0} & 1 & -\theta_{0} + 2 & \textbf{(7)}\\
		 	4 \, A_{2} \overline{C_{1}} & \overline{A_{0}} & 1 & -\theta_{0} + 2& \textbf{(8)}
		 	\end{dmatrix}
		 	&=\left(-\frac{(\theta_0+2)}{4\theta_0}|\vec{C}_1|^2\vec{A}_0-2\s{\vec{A}_2}{\bar{\vec{C}_1}}\bar{\vec{A}_0}+\frac{1}{2}|\vec{C}_1|^2\vec{A}_0+4\s{\vec{A}_2}{\bar{\vec{C}_1}}\bar{\vec{A}_0}\right) z\z^{2-\theta_0}\\
		 	&=\left(\frac{(\theta_0-2)}{4\theta_0}|\vec{C}_1|^2\vec{A}_0+2\s{\vec{A}_2}{\bar{\vec{C}_1}}\bar{\vec{A}_0}\right)z\z^{2-\theta_0}\\
		 	&=\begin{dmatrix}
		 	\frac{(\theta_0-2)}{4\theta_0}|C_1|^2 & A_0 & 1 & -\theta_0+2& \textbf{(6)}\\
		 	2\,\bar{\zeta_4} & \bar{A_0} & 1 &-\theta_0+2& \textbf{(8)}
		 	\end{dmatrix}
		 \end{align*}
		 if
		 \begin{align*}
		 	\zeta_4=\s{\bar{\vec{A}_2}}{\vec{C}_1}.
		 \end{align*}
		 Now as $\vec{E}_1=-\dfrac{1}{2\theta_0}\s{\vec{C}_1}{\vec{C}_1}\bar{\vec{A}_0}$, we have
		 \begin{align*}
		 	\begin{dmatrix}
		 	-\theta_{0} + 2 & E_{1} & -2 \, \theta_{0} + 3 & \theta_{0} & \textbf{(9)}\\
		 	-\frac{C_{1}^{2} {\left(\theta_{0} - 2\right)}}{2 \, \theta_{0}} & \overline{A_{0}} & -2 \, \theta_{0} + 3 & \theta_{0} & \textbf{(10)}
		 	\end{dmatrix}=0,\quad
		 	\begin{dmatrix}
		 	\frac{1}{2} \, \theta_{0} & \overline{E_{1}} & \theta_{0} - 1 & -2 \, \theta_{0} + 4 & \textbf{(22)}\\
		 	\frac{1}{4} \, \overline{C_{1}}^{2} & A_{0} & \theta_{0} - 1 & -2 \, \theta_{0} + 4& \textbf{(23)}
		 	\end{dmatrix}
		 	=0
		 \end{align*}
		 Then, we have
		 \begin{align*}
		 	\begin{dmatrix}
		 	-\frac{1}{2} \, \theta_{0} + 1 & C_{1} & -\theta_{0} + 1 & 0 & \textbf{(11)}\\
		 	-\frac{1}{2} \, \theta_{0} + 1 & B_{1} & -\theta_{0} + 1 & 1 & \textbf{(12)}\\
		 	-\frac{1}{2} \, \theta_{0} + 1 & B_{2} & -\theta_{0} + 1 & 2 & \textbf{(13)}\\
		 	-\frac{C_{1} {\left(\theta_{0} - 2\right)} \overline{C_{1}}}{2 \, \theta_{0}} & \overline{A_{0}} & -\theta_{0} + 1 & 2 & \textbf{(14)}
		 	\end{dmatrix}
		 	=
		 	\begin{dmatrix}
		 	-\frac{(\theta_0-2)}{2} & C_{1} & -\theta_{0} + 1 & 0 & \textbf{(11)}\\
		 	-\frac{(\theta_0-2)}{2} & B_{1} & -\theta_{0} + 1 & 1 & \textbf{(12)}\\
		 	-\frac{(\theta_0-2)}{2} & B_{2} & -\theta_{0} + 1 & 2 & \textbf{(13)}\\
		 	-\frac{{\left(\theta_{0} - 2\right)} |{C_{1}}|^2}{2 \, \theta_{0}} & \overline{A_{0}} & -\theta_{0} + 1 & 2 & \textbf{(14)}\\
		 	\end{dmatrix}.
		 \end{align*}
		 Finally, if 
		 \begin{align}\label{defzeta5}
		 	\zeta_5=2\s{\vec{A}_2}{\vec{C}_1}+\s{\vec{A}_1}{\vec{C}_2}
		 \end{align}
		 recalling that
		 \begin{align*}
		 	\alpha_7=\frac{1}{8\theta_0(\theta_0-4)}\s{\vec{C}_1}{\vec{C}_1}
		 \end{align*}
		 we obtain
		 \begin{align*}
		 	\begin{dmatrix}
		 	-\frac{1}{2} \, \theta_{0} + \frac{3}{2} & C_{2} & -\theta_{0} + 2 & 0 & \textbf{(15)}\\
		 	2 \, A_{1} C_{1} & \overline{A_{0}} & -\theta_{0} + 2 & 0 & \textbf{(16)}\\
		 	-\frac{1}{2} \, \theta_{0} + \frac{3}{2} & B_{3} & -\theta_{0} + 2 & 1 & \textbf{(17)}\\
		 	2 \, A_{1} C_{1} & \overline{A_{1}} & -\theta_{0} + 2 & 1 & \textbf{(18)}\\
		 	-\frac{1}{2} \, \theta_{0} + 2 & C_{3} & -\theta_{0} + 3 & 0 & \textbf{(19)}\\
		 	\frac{1}{4} \, C_{1}^{2} & A_{0} & -\theta_{0} + 3 & 0 & \textbf{(20)}\\
		 	4 \, A_{2} C_{1} + 2 \, A_{1} C_{2} & \overline{A_{0}} & -\theta_{0} + 3 & 0 & \textbf{(21)}
		 	\end{dmatrix}
		 	=\begin{dmatrix}
		 	-\frac{\theta_{0}-3}{2} & C_{2} & -\theta_{0} + 2 & 0 & \textbf{(15)}\\
		 	2 \, \zeta_2 & \overline{A_{0}} & -\theta_{0} + 2 & 0 & \textbf{(16)}\\
		 	-\frac{\theta_{0}-3}{2} & B_{3} & -\theta_{0} + 2 & 1 & \textbf{(17)}\\
		 	2 \, \zeta_2 & \overline{A_{1}} & -\theta_{0} + 2 & 1 & \textbf{(18)}\\
		 	-\frac{\theta_0-4}{2} & C_{3} & -\theta_{0} + 3 & 0 & \textbf{(19)}\\
		 	2\theta_0(\theta_0-4)\alpha_7 & A_{0} & -\theta_{0} + 3 & 0 & \textbf{(20)}\\
		 	2\zeta_5 & \overline{A_{0}} & -\theta_{0} + 3 & 0 & \textbf{(21)}
		 	\end{dmatrix}
		 \end{align*}
		 We sum up her the new coefficients $\zeta$ that we have introduced
		 \begin{align}\label{defzeta05}
		 \left\{
		 \begin{alignedat}{1}
		 \zeta_0&=\s{\vec{A}_1}{\vec{A}_1}\\
		 \zeta_1&=\s{\vec{A}_1}{\vec{A}_2}\\
		 \zeta_2&=\s{\vec{A}_1}{\vec{C}_1}\\
		 \zeta_3&=\s{\bar{\vec{A}_1}}{\vec{C}_2}\\
		 \zeta_4&=\s{\bar{\vec{A}_2}}{\vec{C}_1}\\
		 \zeta_5&=2\s{\vec{A}_2}{\vec{C}_1}+\s{\vec{A}_1}{\vec{C}_2}
		 \end{alignedat}\right.
		 \end{align}
		 Finally, we obtain
		 \begin{align*}
		 	\vec{\alpha}=\begin{dmatrix}
		 	\frac{1}{2} & \overline{B_{1}} & 0 & -\theta_{0} + 2 &\textbf{(1)}\\
		 	2 \, A_{1} \overline{C_{1}} & \overline{A_{0}} & 0 & -\theta_{0} + 2 &\textbf{(2)}\\
		 	\frac{1}{2} & \overline{B_{3}} & 0 & -\theta_{0} + 3 &\textbf{(3)}\\
		 	2 \, A_{1} \overline{C_{1}} & \overline{A_{1}} & 0 & -\theta_{0} + 3 &\textbf{(4)}\\
		 	2 \, A_{1} \overline{C_{2}} & \overline{A_{0}} & 0 & -\theta_{0} + 3 &\textbf{(5)}\\
		 	1 & \overline{B_{2}} & 1 & -\theta_{0} + 2 &\textbf{(6)}\\
		 	\frac{1}{2} \, C_{1} \overline{C_{1}} & A_{0} & 1 & -\theta_{0} + 2 &\textbf{(7)}\\
		 	4 \, A_{2} \overline{C_{1}} & \overline{A_{0}} & 1 & -\theta_{0} + 2 &\textbf{(8)}\\
		 	\ccancel{-\theta_{0} + 2} & E_{1} & -2 \, \theta_{0} + 3 & \theta_{0} &\textbf{(9)}\\
		 	\ccancel{-\frac{C_{1}^{2} {\left(\theta_{0} - 2\right)}}{2 \, \theta_{0}}} & \overline{A_{0}} & -2 \, \theta_{0} + 3 & \theta_{0} &\textbf{(10)}\\
		 	-\frac{1}{2} \, \theta_{0} + 1 & C_{1} & -\theta_{0} + 1 & 0 &\textbf{(11)}\\
		 	-\frac{1}{2} \, \theta_{0} + 1 & B_{1} & -\theta_{0} + 1 & 1 &\textbf{(12)}\\
		 	-\frac{1}{2} \, \theta_{0} + 1 & B_{2} & -\theta_{0} + 1 & 2 &\textbf{(13)}\\
		 	-\frac{C_{1} {\left(\theta_{0} - 2\right)} \overline{C_{1}}}{2 \, \theta_{0}} & \overline{A_{0}} & -\theta_{0} + 1 & 2 &\textbf{(14)}\\
		 	-\frac{1}{2} \, \theta_{0} + \frac{3}{2} & C_{2} & -\theta_{0} + 2 & 0 &\textbf{(15)}\\
		 	2 \, A_{1} C_{1} & \overline{A_{0}} & -\theta_{0} + 2 & 0 &\textbf{(16)}\\
		 	-\frac{1}{2} \, \theta_{0} + \frac{3}{2} & B_{3} & -\theta_{0} + 2 & 1 &\textbf{(17)}\\
		 	2 \, A_{1} C_{1} & \overline{A_{1}} & -\theta_{0} + 2 & 1 &\textbf{(18)}\\
		 	-\frac{1}{2} \, \theta_{0} + 2 & C_{3} & -\theta_{0} + 3 & 0 &\textbf{(19)}\\
		 	\frac{1}{4} \, C_{1}^{2} & A_{0} & -\theta_{0} + 3 & 0 &\textbf{(20)}\\
		 	4 \, A_{2} C_{1} + 2 \, A_{1} C_{2} & \overline{A_{0}} & -\theta_{0} + 3 & 0 &\textbf{(21)}\\
		 	\ccancel{\frac{1}{2} \, \theta_{0}} & \overline{E_{1}} & \theta_{0} - 1 & -2 \, \theta_{0} + 4 &\textbf{(22)}\\
		 	\ccancel{\frac{1}{4} \, \overline{C_{1}}^{2}} & A_{0} & \theta_{0} - 1 & -2 \, \theta_{0} + 4 &\textbf{(23)}
		 	\end{dmatrix}
		 	&=\begin{dmatrix}
		 	2\theta_0(\theta_0+1)\bar{\alpha_2}  & \overline{A_{0}} & 0 & -\theta_{0} + 2 &\textbf{(1)}\\
		 	\frac{\bar{\zeta_2}}{\theta_0-3} & A_1 & 0 & -\theta_0+3&\textbf{(3)}\\
		 	2\theta_0(\theta_0+1)\bar{\alpha_2} & \bar{A_1} & 0 & -\theta_0+3 &\textbf{(4)}\\
		 	\bar{\zeta_3}  & \bar{A_0} & 0 & -\theta_0+3 &\textbf{(5)}\\
		 	\frac{(\theta_0-2)}{4\theta_0}|C_1|^2 & A_0 & 1 & -\theta_0+2 &\textbf{(6)}\\
		 	2\,\bar{\zeta_4} & \bar{A_0} & 1 &-\theta_0+2 &\textbf{(8)}\\
		 	-\frac{(\theta_0-2)}{2} & C_{1} & -\theta_{0} + 1 & 0 &\textbf{(11)}\\
		 	-\frac{(\theta_0-2)}{2} & B_{1} & -\theta_{0} + 1 & 1 &\textbf{(12)}\\
		 	-\frac{(\theta_0-2)}{2} & B_{2} & -\theta_{0} + 1 & 2 &\textbf{(13)}\\
		 	-\frac{{\left(\theta_{0} - 2\right)} |{C_{1}}|^2}{2 \, \theta_{0}} & \overline{A_{0}} & -\theta_{0} + 1 & 2 &\textbf{(14)}\\
		 	-\frac{(\theta_{0}-3)}{2} & C_{2} & -\theta_{0} + 2 & 0 &\textbf{(15)}\\
		 	2 \, \zeta_2 & \overline{A_{0}} & -\theta_{0} + 2 & 0 &\textbf{(16)}\\
		 	-\frac{(\theta_{0}-3)}{2} & B_{3} & -\theta_{0} + 2 & 1 &\textbf{(17)}\\
		 	2 \, \zeta_2 & \overline{A_{1}} & -\theta_{0} + 2 & 1 &\textbf{(18)}\\
		 	-\frac{(\theta_0-4)}{2} & C_{3} & -\theta_{0} + 3 & 0 &\textbf{(19)}\\
		 	2\theta_0(\theta_0-4)\alpha_7 & A_{0} & -\theta_{0} + 3 & 0 &\textbf{(20)}\\
		 	2\,\zeta_5 & \overline{A_{0}} & -\theta_{0} + 3 & 0 &\textbf{(21)}
		 	\end{dmatrix}
		 \end{align*}
		 which once translated again into code yields
		 \begin{align*}
		 	\vec{\alpha}=\begin{dmatrix}
		 	2\theta_0(\theta_0+1)\bar{\alpha_2}  & \overline{A_{0}} & 0 & -\theta_{0} + 2 &\textbf{(1)}\\
		 	\frac{\bar{\zeta_2}}{\theta_0-3} & A_1 & 0 & -\theta_0+3 &\textbf{(3)}\\
		 	2\theta_0(\theta_0+1)\bar{\alpha_2} & \bar{A_1} & 0 & -\theta_0+3 &\textbf{(4)}\\
		 	\bar{\zeta_3}  & \bar{A_0} & 0 & -\theta_0+3 &\textbf{(5)}\\
		 	\frac{(\theta_0-2)}{4\theta_0}|C_1|^2 & A_0 & 1 & -\theta_0+2 &\textbf{(6)}\\
		 	2\,\bar{\zeta_4} & \bar{A_0} & 1 &-\theta_0+2 &\textbf{(8)}\\
		 	-\frac{(\theta_0-2)}{2} & C_{1} & -\theta_{0} + 1 & 0 &\textbf{(11)}\\
		 	-\frac{(\theta_0-2)}{2} & B_{1} & -\theta_{0} + 1 & 1 &\textbf{(12)}\\
		 	-\frac{(\theta_0-2)}{2} & B_{2} & -\theta_{0} + 1 & 2 &\textbf{(13)}\\
		 	-\frac{{\left(\theta_{0} - 2\right)} |{C_{1}}|^2}{2 \, \theta_{0}} & \overline{A_{0}} & -\theta_{0} + 1 & 2 &\textbf{(14)}\\
		 	-\frac{(\theta_{0}-3)}{2} & C_{2} & -\theta_{0} + 2 & 0 &\textbf{(15)}\\
		 	2 \, \zeta_2 & \overline{A_{0}} & -\theta_{0} + 2 & 0 &\textbf{(16)}\\
		 	-\frac{(\theta_{0}-3)}{2} & B_{3} & -\theta_{0} + 2 & 1 &\textbf{(17)}\\
		 	2 \, \zeta_2 & \overline{A_{1}} & -\theta_{0} + 2 & 1 &\textbf{(18)}\\
		 	-\frac{(\theta_0-4)}{2} & C_{3} & -\theta_{0} + 3 & 0 &\textbf{(19)}\\
		 	2\theta_0(\theta_0-4)\alpha_7 & A_{0} & -\theta_{0} + 3 & 0 &\textbf{(20)}\\
		 	2\,\zeta_5 & \overline{A_{0}} & -\theta_{0} + 3 & 0 &\textbf{(21)} 
		 	\end{dmatrix}\begin{dmatrix}
		 	2 \, {\left(\theta_{0} + 1\right)} \theta_{0} \overline{\alpha_{2}} & \overline{A_{0}} & 0 & -\theta_{0} + 2 \\
		 	\frac{\overline{\zeta_{2}}}{\theta_{0} - 3} & A_{1} & 0 & -\theta_{0} + 3 \\
		 	2 \, {\left(\theta_{0} + 1\right)} \theta_{0} \overline{\alpha_{2}} & \overline{A_{1}} & 0 & -\theta_{0} + 3 \\
		 	\overline{\zeta_{3}} & \overline{A_{0}} & 0 & -\theta_{0} + 3 \\
		 	\frac{{\left(\theta_{0} - 2\right)} {\left| C_{1} \right|}^{2}}{4 \, \theta_{0}} & A_{0} & 1 & -\theta_{0} + 2 \\
		 	2 \, \overline{\zeta_{4}} & \overline{A_{0}} & 1 & -\theta_{0} + 2 \\
		 	-\frac{1}{2} \, \theta_{0} + 1 & C_{1} & -\theta_{0} + 1 & 0 \\
		 	-\frac{1}{2} \, \theta_{0} + 1 & B_{1} & -\theta_{0} + 1 & 1 \\
		 	-\frac{1}{2} \, \theta_{0} + 1 & B_{2} & -\theta_{0} + 1 & 2 \\
		 	-\frac{{\left(\theta_{0} - 2\right)} {\left| C_{1} \right|}^{2}}{2 \, \theta_{0}} & \overline{A_{0}} & -\theta_{0} + 1 & 2 \\
		 	-\frac{1}{2} \, \theta_{0} + \frac{3}{2} & C_{2} & -\theta_{0} + 2 & 0 \\
		 	2 \, \zeta_{2} & \overline{A_{0}} & -\theta_{0} + 2 & 0 \\
		 	-\frac{1}{2} \, \theta_{0} + \frac{3}{2} & B_{3} & -\theta_{0} + 2 & 1 \\
		 	2 \, \zeta_{2} & \overline{A_{1}} & -\theta_{0} + 2 & 1 \\
		 	-\frac{1}{2} \, \theta_{0} + 2 & C_{3} & -\theta_{0} + 3 & 0 \\
		 	2 \, \alpha_{7} {\left(\theta_{0} - 4\right)} \theta_{0} & A_{0} & -\theta_{0} + 3 & 0 \\
		 	2 \, \zeta_{5} & \overline{A_{0}} & -\theta_{0} + 3 & 0
		 	\end{dmatrix}
		 \end{align*}
		 Left \TeX\; right Sage.

		\small
		\begin{align*}
			\s{\vec{\alpha}}{\phi}=\begin{dmatrix}
			2 \, {\left(\theta_{0} \overline{\alpha_{2}} + \overline{\alpha_{2}}\right)} \ccancel{\overline{A_{0}}^{2}} & 0 & 2 &\textbf{(1)}\\
			\frac{{\left(\theta_{0} - 3\right)} \ccancel{\overline{A_{0}}^{2}} \overline{\zeta_{3}} + 2 \, {\left(2 \, \theta_{0}^{3} \overline{\alpha_{2}} - 5 \, \theta_{0}^{2} \overline{\alpha_{2}} - 3 \, \theta_{0} \overline{\alpha_{2}}\right)} \ccancel{\overline{A_{0}} \overline{A_{1}}} + \ccancel{A_{1} \overline{A_{0}}} \overline{\zeta_{2}}}{\theta_{0}^{2} - 3 \, \theta_{0}} & 0 & 3 &\textbf{(2)} \\
			-\frac{\ccancel{A_{0} C_{1}} {\left(\theta_{0} - 2\right)}}{2 \, \theta_{0}} & 1 & 0 &\textbf{(3)}\\
			-\frac{\ccancel{A_{0} B_{1}} {\left(\theta_{0} - 2\right)}}{2 \, \theta_{0}} & 1 & 1 &\textbf{(4)}\\
			-\frac{4 \, A_{0} {\left(\theta_{0} - 2\right)} {\left| C_{1} \right|}^{2} \overline{A_{0}} - 32 \, \theta_{0} \ccancel{\overline{A_{0}}^{2}} \overline{\zeta_{4}} + 8 \, {\left(\theta_{0}^{2} - 2 \, \theta_{0}\right)} A_{0} B_{2} + {\left(\theta_{0}^{2} - 2 \, \theta_{0}\right)} C_{1} \overline{C_{1}}}{16 \, \theta_{0}^{2}} & 1 & 2 &\textbf{(5)}\\
			\frac{4 \, A_{0} {\left(\theta_{0} + 1\right)} \zeta_{2} \overline{A_{0}} - {\left(\theta_{0}^{2} - 2 \, \theta_{0}\right)} A_{1} C_{1} - {\left(\theta_{0}^{2} - 2 \, \theta_{0} - 3\right)} A_{0} C_{2}}{2 \, {\left(\theta_{0}^{2} + \theta_{0}\right)}} & 2 & 0 &\textbf{(6)}\\
			\frac{\ccancel{4 \, A_{0} {\left(\theta_{0} + 1\right)} \zeta_{2} \overline{A_{1}}} - \ccancel{{\left(\theta_{0}^{2} - 2 \, \theta_{0}\right)} A_{1} B_{1}} - \ccancel{{\left(\theta_{0}^{2} - 2 \, \theta_{0} - 3\right)} A_{0} B_{3}}}{2 \, {\left(\theta_{0}^{2} + \theta_{0}\right)}} & 2 & 1 &\textbf{(7)}\\
			\lambda_1 & 3 & 0 &\textbf{(8)}\\
			-\frac{\ccancel{C_{1} {\left(\theta_{0} - 2\right)} \overline{A_{0}}}}{2 \, \theta_{0}} & -\theta_{0} + 1 & \theta_{0} &\textbf{(9)}\\
			-\frac{{\left(\theta_{0}^{2} - \theta_{0} - 2\right)} B_{1} \overline{A_{0}} + {\left(\theta_{0}^{2} - 2 \, \theta_{0}\right)} C_{1} \overline{A_{1}}}{2 \, {\left(\theta_{0}^{2} + \theta_{0}\right)}} & -\theta_{0} + 1 & \theta_{0} + 1 &\textbf{(10)}\\
			\lambda_2 & -\theta_{0} + 1 & \theta_{0} + 2 &\textbf{(11)}\\
			-\frac{\ccancel{C_{2} {\left(\theta_{0} - 3\right)} \overline{A_{0}}} - \ccancel{4 \, \zeta_{2} \overline{A_{0}}^{2}}}{2 \, \theta_{0}} & -\theta_{0} + 2 & \theta_{0} &\textbf{(12)}\\
			\frac{\ccancel{4 \, {\left(2 \, \theta_{0} + 1\right)} \zeta_{2} \overline{A_{0}} \overline{A_{1}}} - {\left(\theta_{0}^{2} - 2 \, \theta_{0} - 3\right)} B_{3} \overline{A_{0}} - {\left(\theta_{0}^{2} - 3 \, \theta_{0}\right)} C_{2} \overline{A_{1}}}{2 \, {\left(\theta_{0}^{2} + \theta_{0}\right)}} & -\theta_{0} + 2 & \theta_{0} + 1 &\textbf{(13)}\\
			-\frac{C_{1}^{2} {\left(\theta_{0} - 2\right)} - 32 \, {\left(\alpha_{7} \theta_{0}^{2} - 4 \, \alpha_{7} \theta_{0}\right)} A_{0} \overline{A_{0}} + 8 \, C_{3} {\left(\theta_{0} - 4\right)} \overline{A_{0}} - \ccancel{32 \, \zeta_{5} \overline{A_{0}}^{2}}}{16 \, \theta_{0}} & -\theta_{0} + 3 & \theta_{0} &\textbf{(14)}\\
			2 \, {\left(\theta_{0} \overline{\alpha_{2}} + \overline{\alpha_{2}}\right)} A_{0} \overline{A_{0}} & \theta_{0} & -\theta_{0} + 2 &\textbf{(15)}\\
			\frac{A_{0} {\left(\theta_{0} - 3\right)} \overline{A_{0}} \overline{\zeta_{3}} + \ccancel{2 \, {\left(\theta_{0}^{3} \overline{\alpha_{2}} - 2 \, \theta_{0}^{2} \overline{\alpha_{2}} - 3 \, \theta_{0} \overline{\alpha_{2}}\right)} A_{0} \overline{A_{1}}} + \ccancel{A_{0} A_{1} \overline{\zeta_{2}}}}{\theta_{0}^{2} - 3 \, \theta_{0}} & \theta_{0} & -\theta_{0} + 3 &\textbf{(16)}\\
			\frac{\ccancel{8 \, A_{1} \theta_{0}^{3} \overline{A_{0}} \overline{\alpha_{2}}} + \ccancel{A_{0}^{2} {\left(\theta_{0} - 2\right)} {\left| C_{1} \right|}^{2}} + 8 \, A_{0} \theta_{0} \overline{A_{0}} \overline{\zeta_{4}}}{4 \, \theta_{0}^{2}} & \theta_{0} + 1 & -\theta_{0} + 2 &\textbf{(17)}
			\end{dmatrix}
		\end{align*}
		where
		\begin{align*}
			&\lambda_1=\frac{1}{2 \, {\left(\theta_{0}^{3} + 3 \, \theta_{0}^{2} + 2 \, \theta_{0}\right)}}\bigg\{\ccancel{4 \, {\left(\theta_{0}^{2} + 2 \, \theta_{0}\right)} A_{1} \zeta_{2} \overline{A_{0}}} + 4 \, {\left(\theta_{0}^{2} + 3 \, \theta_{0} + 2\right)} A_{0} \zeta_{5} \overline{A_{0}} + \ccancel{4 \, {\left(\alpha_{7} \theta_{0}^{4} - \alpha_{7} \theta_{0}^{3} - 10 \, \alpha_{7} \theta_{0}^{2} - 8 \, \alpha_{7} \theta_{0}\right)} A_{0}^{2}}\\
			& - {\left(\theta_{0}^{3} - \theta_{0}^{2} - 2 \, \theta_{0}\right)} A_{2} C_{1} - {\left(\theta_{0}^{3} - \theta_{0}^{2} - 6 \, \theta_{0}\right)} A_{1} C_{2} - {\left(\theta_{0}^{3} - \theta_{0}^{2} - 10 \, \theta_{0} - 8\right)} A_{0} C_{3}\bigg\}\\
			&\lambda_2=-\frac{\ccancel{{\left(\theta_{0}^{3} + \theta_{0}^{2} - 4 \, \theta_{0} - 4\right)} {\left| C_{1} \right|}^{2} \overline{A_{0}}^{2}} + {\left(\theta_{0}^{4} + \theta_{0}^{3} - 4 \, \theta_{0}^{2} - 4 \, \theta_{0}\right)} B_{2} \overline{A_{0}} + \ccancel{{\left(\theta_{0}^{4} - 4 \, \theta_{0}^{2}\right)} B_{1} \overline{A_{1}}} + {\left(\theta_{0}^{4} - \theta_{0}^{3} - 2 \, \theta_{0}^{2}\right)} C_{1} \overline{A_{2}}}{2 \, {\left(\theta_{0}^{4} + 3 \, \theta_{0}^{3} + 2 \, \theta_{0}^{2}\right)}}.
		\end{align*}
		Therefore, we have
		\begin{align*}
			\textbf{(1)}=\textbf{(2)}=\textbf{(3)}=\textbf{(4)}=\textbf{(7)}=\textbf{(9)}=\textbf{(12)}.
		\end{align*}
		Now, recall that by \eqref{a1c12} and \eqref{alpha27}
		\begin{align*}
			\s{\vec{C}_1}{\vec{C}_1}+8\s{\bar{\vec{A}_0}}{\vec{C}_3}=0,\quad \alpha_7=\frac{1}{8\theta_0(\theta_0-4)}\s{\vec{C}_1}{\vec{C}_1}
		\end{align*}
		so as $|\vec{A}_0|^2=\dfrac{1}{2}$
		\begin{align*}
			&C_{1}^{2} {\left(\theta_{0} - 2\right)} - 32 \, {\left(\alpha_{7} \theta_{0}^{2} - 4 \, \alpha_{7} \theta_{0}\right)} A_{0} \overline{A_{0}} + 8 \, C_{3} {\left(\theta_{0} - 4\right)} \overline{A_{0}} - \ccancel{32 \, \zeta_{5} \overline{A_{0}}^{2}}\\
			&=C_{1}^{2} {\left(\theta_{0} - 2\right)} - 16 \, {\theta_0\left(\theta_{0} - 4 \,  \right)}\alpha_{7}   + 8 \, C_{3} {\left(\theta_{0} - 4\right)} \overline{A_{0}}\\
			&=(\theta_0-4)\s{\vec{C}_1}{\vec{C}_1}+8(\theta_0-4)\s{\bar{\vec{A}_0}}{\vec{C}_3}=0
		\end{align*}
		so 
		\begin{align}
			\textbf{(14)}=0.
		\end{align}
		Now, we compute
		\begin{align*}
			&4 \, A_{0} {\left(\theta_{0} - 2\right)} {\left| C_{1} \right|}^{2} \overline{A_{0}} - 32 \, \theta_{0} \ccancel{\overline{A_{0}}^{2}} \overline{\zeta_{4}} + 8 \, {\left(\theta_{0}^{2} - 2 \, \theta_{0}\right)} A_{0} B_{2} + {\left(\theta_{0}^{2} - 2 \, \theta_{0}\right)} C_{1} \overline{C_{1}}\\
			&=2(\theta_0-2)|\vec{C}_1|^2+8\theta_0(\theta_0-2)\s{\vec{A}_0}{\vec{B}_2}+\theta_0(\theta_0-2)|\vec{C}_1|^2\\
			&=(\theta_0-2)\bigg((\theta_0+2)|\vec{C}_1|^2+8\theta_0\bs{\vec{A}_0}{-\frac{(\theta_0+2)}{4\theta_0}|\vec{C}_1|^2\bar{\vec{A}_0}-2\s{\bar{\vec{A}_2}}{\vec{C}_1}\vec{A}_0}\bigg)\\
			&=0
		\end{align*}
		as $|\vec{A}_0|^2=\dfrac{1}{2}$. Therefore, we have
		\begin{align*}
			\textbf{(5)}=0.
		\end{align*}
		Now, we have as $\zeta_2=\s{\vec{A}_1}{\vec{C}_1}$ and $\s{\vec{A}_1}{\vec{C}_1}+\s{\vec{A}_0}{\vec{C}_2}=0$ the identity
		\begin{align*}
			&4 \, A_{0} {\left(\theta_{0} + 1\right)} \zeta_{2} \overline{A_{0}} - {\left(\theta_{0}^{2} - 2 \, \theta_{0}\right)} A_{1} C_{1} - {\left(\theta_{0}^{2} - 2 \, \theta_{0} - 3\right)} A_{0} C_{2}\\
			&=2(\theta_0+1)A_1C_1- {\left(\theta_{0}^{2} - 2 \, \theta_{0}\right)} A_{1} C_{1} + {\left(\theta_{0}^{2} - 2 \, \theta_{0} - 3\right)} A_{1} C_{1}\\
			&=(2\theta_0-1)\s{\vec{A}_1}{\vec{C}_1}\\
			&=(2\theta_0-1)\zeta_2.
		\end{align*}
		and
		\begin{align}\label{inv6}
			\textbf{(6)}=\frac{4 \, A_{0} {\left(\theta_{0} + 1\right)} \zeta_{2} \overline{A_{0}} - {\left(\theta_{0}^{2} - 2 \, \theta_{0}\right)} A_{1} C_{1} - {\left(\theta_{0}^{2} - 2 \, \theta_{0} - 3\right)} A_{0} C_{2}}{2 \, {\left(\theta_{0}^{2} + \theta_{0}\right)}}=\frac{(2\theta_0-1)\zeta_2}{2\theta_0(\theta_0+1)}
		\end{align}
		Now, recall
		\begin{align}\label{defzeta05bis}
		\left\{
		\begin{alignedat}{1}
		\zeta_0&=\s{\vec{A}_1}{\vec{A}_1}\\
		\zeta_1&=\s{\vec{A}_1}{\vec{A}_2}\\
		\zeta_2&=\s{\vec{A}_1}{\vec{C}_1}\\
		\zeta_3&=\s{\bar{\vec{A}_1}}{\vec{C}_2}\\
		\zeta_4&=\s{\bar{\vec{A}_2}}{\vec{C}_1}\\
		\zeta_5&=2\s{\vec{A}_2}{\vec{C}_1}+\s{\vec{A}_1}{\vec{C}_2}
		\end{alignedat}\right.
		\end{align}
		In particular, as
		\begin{align*}
			\s{\vec{A}_2}{\vec{C}_1}+\s{\vec{A}_1}{\vec{C}_2}+\s{\vec{A}_0}{\vec{C}_3}=0,
		\end{align*}
		we have
		\begin{align*}
			\lambda_1&=
			\frac{1}{2 \, {\left(\theta_{0}^{3} + 3 \, \theta_{0}^{2} + 2 \, \theta_{0}\right)}}\bigg\{\ccancel{4 \, {\left(\theta_{0}^{2} + 2 \, \theta_{0}\right)} A_{1} \zeta_{2} \overline{A_{0}}} + 4 \, {\left(\theta_{0}^{2} + 3 \, \theta_{0} + 2\right)} A_{0} \zeta_{5} \overline{A_{0}} + \ccancel{4 \, {\left(\alpha_{7} \theta_{0}^{4} - \alpha_{7} \theta_{0}^{3} - 10 \, \alpha_{7} \theta_{0}^{2} - 8 \, \alpha_{7} \theta_{0}\right)} A_{0}^{2}}\\
		    &- {\left(\theta_{0}^{3} - \theta_{0}^{2} - 2 \, \theta_{0}\right)} A_{2} C_{1} - {\left(\theta_{0}^{3} - \theta_{0}^{2} - 6 \, \theta_{0}\right)} A_{1} C_{2} - {\left(\theta_{0}^{3} - \theta_{0}^{2} - 10 \, \theta_{0} - 8\right)} A_{0} C_{3}\bigg\}\\
		    &=\frac{1}{2 \, {\left(\theta_{0}^{3} + 3 \, \theta_{0}^{2} + 2 \, \theta_{0}\right)}}\bigg\{2(\theta_0^2+3\theta_0+2)\zeta_5\\
		    &- {\left(\theta_{0}^{3} - \theta_{0}^{2} - 2 \, \theta_{0}\right)} A_{2} C_{1} - {\left(\theta_{0}^{3} - \theta_{0}^{2} - 6 \, \theta_{0}\right)} A_{1} C_{2} + {\left(\theta_{0}^{3} - \theta_{0}^{2} - 10 \, \theta_{0} - 8\right)} (A_{2} C_{1}+A_{1} C_{2}) \bigg\}\\
		    &=\frac{1}{2 \, {\left(\theta_{0}^{3} + 3 \, \theta_{0}^{2} + 2 \, \theta_{0}\right)}}\bigg\{2(\theta_0^2+3\theta_0+2)\zeta_5 -8(\theta_0+1)A_2C_1-4(\theta_0+2)A_1C_2\bigg\}\\
		    &=\frac{1}{2 \, {\left(\theta_{0}^{3} + 3 \, \theta_{0}^{2} + 2 \, \theta_{0}\right)}}\bigg\{2(\theta_0^2+3\theta_0+2)\zeta_5-4(\theta_0+1)\zeta_5-4A_1C_2 \bigg\}\\
		    &=\frac{1}{2 \, {\left(\theta_{0}^{3} + 3 \, \theta_{0}^{2} + 2 \, \theta_{0}\right)}}\bigg\{2\theta_0(\theta_0+1)\zeta_5-4\zeta_6 \bigg\}
		\end{align*}
		where
		\begin{align*}
			\zeta_6=\s{\vec{A}_1}{\vec{C}_2}.
		\end{align*}
		Therefore, 
		\begin{align*}
			\textbf{(8)}=\frac{1}{2 \, {\left(\theta_{0}^{3} + 3 \, \theta_{0}^{2} + 2 \, \theta_{0}\right)}}\bigg\{2\theta_0(\theta_0+1)\zeta_5-4\zeta_6 \bigg\}.
		\end{align*}
		Then, we compute as $\vec{B}_1=-2\s{\bar{\vec{A}_1}}{\vec{C}_1}\vec{A}_0$
		\begin{align*}
			\textbf{(10)}&=-\frac{{\left(\theta_{0}^{2} - \theta_{0} - 2\right)} B_{1} \overline{A_{0}} + {\left(\theta_{0}^{2} - 2 \, \theta_{0}\right)} C_{1} \overline{A_{1}}}{2 \, {\left(\theta_{0}^{2} + \theta_{0}\right)}}
			=\frac{\big({\left(\theta_{0}^{2} - \theta_{0} - 2\right)} \bar{A_1}C_1 - {\left(\theta_{0}^{2} - 2 \, \theta_{0}\right)} C_{1} \overline{A_{1}}\big)}{2 \, {\left(\theta_{0}^{2} + \theta_{0}\right)}}\\
			&=\frac{(\theta_0-2)\s{\bar{\vec{A}_1}}{\vec{C}_1}}{2\theta_0(\theta_0+1)}=(\theta_0-2)\alpha_2.
		\end{align*}
		Then, we have
		\begin{align*}
			\s{\bar{\vec{A}_0}}{\vec{B}_2}=\bs{\bar{\vec{A}_0}}{-\frac{(\theta_0+2)}{4\theta_0}|\vec{C}_1|^2\bar{\vec{A}_0}-2\s{\bar{\vec{A}_2}}{\vec{C}_1}\vec{A}_0}=-\s{\bar{\vec{A}_2}}{\vec{C}_1}
		\end{align*}
		so
		\begin{align*}
			\textbf{(11)}&=\lambda_2=-\frac{\ccancel{{\left(\theta_{0}^{3} + \theta_{0}^{2} - 4 \, \theta_{0} - 4\right)} {\left| C_{1} \right|}^{2} \overline{A_{0}}^{2}} + {\left(\theta_{0}^{4} + \theta_{0}^{3} - 4 \, \theta_{0}^{2} - 4 \, \theta_{0}\right)} B_{2} \overline{A_{0}} + \ccancel{{\left(\theta_{0}^{4} - 4 \, \theta_{0}^{2}\right)} B_{1} \overline{A_{1}}} + {\left(\theta_{0}^{4} - \theta_{0}^{3} - 2 \, \theta_{0}^{2}\right)} C_{1} \overline{A_{2}}}{2 \, {\left(\theta_{0}^{4} + 3 \, \theta_{0}^{3} + 2 \, \theta_{0}^{2}\right)}}\\
			&=\frac{ {\left(\theta_{0}^{4} + \theta_{0}^{3} - 4 \, \theta_{0}^{2} - 4 \, \theta_{0}\right)} \bar{A_2}\, C_1 - {\left(\theta_{0}^{4} - \theta_{0}^{3} - 2 \, \theta_{0}^{2}\right)} C_{1} \overline{A_{2}}}{2 \, {\left(\theta_{0}^{4} + 3 \, \theta_{0}^{3} + 2 \, \theta_{0}^{2}\right)}}\\
			&=\frac{ {\left(2\theta_{0}^{3} - 2 \, \theta_{0}^{2} - 4 \, \theta_{0}\right)} \bar{A_2}\, C_1 }{2 \, {\left(\theta_{0}^{4} + 3 \, \theta_{0}^{3} + 2 \, \theta_{0}^{2}\right)}}
			=\frac{ {\left(\theta_0+1)(\theta_0-2\right)} \bar{A_2}\, C_1 }{ \, {\left(\theta_{0}^{3} + 3 \, \theta_{0}^{2} + 2 \, \theta_{0}\right)}}=\frac{ {\left(\theta_0+1)(\theta_0-2\right)} \bar{A_2}\, C_1 }{\theta_0(\theta_{0}+1)(\theta_0+2)}=\frac{(\theta_0-2)\bar{A_2}C_1}{\theta_0(\theta_0+2)}\\
			&=\frac{(\theta_0-2)\zeta_4}{\theta_0(\theta_0+2)}.
		\end{align*}
		Then, we have
		\begin{align*}
			\s{\bar{\vec{A}_0}}{\vec{B}_3}&=\bs{\bar{\vec{A}_0}}{-2\s{\bar{\vec{A}_1}}{\vec{C}_1}\vec{A}_1+\frac{2}{\theta_0-3}\s{\vec{A}_1}{\vec{C}_1}\bar{\vec{A}_1}-2\s{\bar{\vec{A}_1}}{\vec{C}_2}\vec{A}_0}=-\s{\bar{\vec{A}_1}}{\vec{C}_2}
		\end{align*}
		so
		\begin{align*}
			\textbf{(13)}&=\frac{\ccancel{4 \, {\left(2 \, \theta_{0} + 1\right)} \zeta_{2} \overline{A_{0}} \overline{A_{1}}} - {\left(\theta_{0}^{2} - 2 \, \theta_{0} - 3\right)} B_{3} \overline{A_{0}} - {\left(\theta_{0}^{2} - 3 \, \theta_{0}\right)} C_{2} \overline{A_{1}}}{2 \, {\left(\theta_{0}^{2} + \theta_{0}\right)}}\\
			&=\frac{ {\left(\theta_{0}^{2} - 2 \, \theta_{0} - 3\right)} \bar{A_1}C_2 - {\left(\theta_{0}^{2} - 3 \, \theta_{0}\right)} C_{2} \overline{A_{1}}}{2 \, {\left(\theta_{0}^{2} + \theta_{0}\right)}}=\frac{(\theta_0-3)\bar{A_1}C_2}{2\theta_0(\theta_0+1)}\\
			&=\frac{(\theta_0-3)\zeta_3}{2\theta_0(\theta_0+1)}.
		\end{align*}
		Finally, we trivially have
		\begin{align}
			\textbf{(15)}=2 \, {\left(\theta_{0} \overline{\alpha_{2}} + \overline{\alpha_{2}}\right)} A_{0} \overline{A_{0}}=(\theta_0+1)\bar{\alpha_2}.
		\end{align}
		then
		\begin{align}
			\textbf{(16)}&=\frac{A_{0} {\left(\theta_{0} - 3\right)} \overline{A_{0}} \overline{\zeta_{3}} + \ccancel{2 \, {\left(\theta_{0}^{3} \overline{\alpha_{2}} - 2 \, \theta_{0}^{2} \overline{\alpha_{2}} - 3 \, \theta_{0} \overline{\alpha_{2}}\right)} A_{0} \overline{A_{1}}} + \ccancel{A_{0} A_{1} \overline{\zeta_{2}}}}{\theta_{0}^{2} - 3 \, \theta_{0}}
			=\frac{\bar{\zeta_3}}{2\theta_0}
		\end{align}
		while
		\begin{align}
			\textbf{(17)}&=\frac{\ccancel{8 \, A_{1} \theta_{0}^{3} \overline{A_{0}} \overline{\alpha_{2}}} + \ccancel{A_{0}^{2} {\left(\theta_{0} - 2\right)} {\left| C_{1} \right|}^{2}} + 8 \, A_{0} \theta_{0} \overline{A_{0}} \overline{\zeta_{4}}}{4 \, \theta_{0}^{2}}=\frac{\bar{\zeta_4}}{\theta_0}.
		\end{align}
		Finally
		\begin{align*}
			\s{\vec{\alpha}}{\phi}=\begin{dmatrix}
			\frac{(2\theta_0-1)\zeta_2}{2\theta_0(\theta_0+1)} & 2 & 0 &\textbf{(6)}\\
			\frac{1}{2 \, {\left(\theta_{0}^{3} + 3 \, \theta_{0}^{2} + 2 \, \theta_{0}\right)}}\bigg\{2\theta_0(\theta_0+1)\zeta_5-4\zeta_6 \bigg\} & 3 & 0 &\textbf{(8)}\\
		    (\theta_0-2)\alpha_2 & -\theta_{0} + 1 & \theta_{0} + 1 &\textbf{(10)}\\
			\frac{(\theta_0-2)\zeta_4}{\theta_0(\theta_0+2)} & -\theta_{0} + 1 & \theta_{0} + 2 &\textbf{(11)}\\
			\frac{(\theta_0-3)\zeta_3}{2\theta_0(\theta_0+1)} & -\theta_{0} + 2 & \theta_{0} + 1 &\textbf{(13)}\\
			(\theta_0+1)\bar{\alpha_2} & \theta_{0} & -\theta_{0} + 2 &\textbf{(15)}\\
			\frac{\bar{\zeta_3}}{2\theta_0} & \theta_{0} & -\theta_{0} + 3 &\textbf{(16)}\\
			\frac{\bar{\zeta_4}}{\theta_0} & \theta_{0} + 1 & -\theta_{0} + 2 &\textbf{(17)}
			\end{dmatrix}
			=
		\end{align*}
		In particular, we see easily that if
		\begin{align*}
			\s{\vec{\alpha}}{\phi}=f(z)dz,
		\end{align*}
		Indeed, see \texttt{test\_scaling} in the code. We can also check it directly.
		Indeed, we have
		\begin{align*}
			\p{\z}f(z)&=\begin{dmatrix}
			(\theta_0+1)(\theta_0-2)\alpha_2 & -\theta_{0} + 1 & \theta_{0}  &\textbf{(10)}\\
			\frac{(\theta_0-2)\zeta_4}{\theta_0} & -\theta_{0} + 1 & \theta_{0} +1 &\textbf{(11)}\\
			\frac{(\theta_0-3)\zeta_3}{2\theta_0} & -\theta_{0} + 2 & \theta_{0}  &\textbf{(13)}\\
			-(\theta_0-2)(\theta_0+1)\bar{\alpha_2} & \theta_{0} & -\theta_{0} + 1 &\textbf{(15)}\\
			\frac{-(\theta_0-3)\bar{\zeta_3}}{2\theta_0} & \theta_{0} & -\theta_{0} + 2 &\textbf{(16)}\\
			\frac{-(\theta_0-2)\bar{\zeta_4}}{\theta_0} & \theta_{0} + 1 & -\theta_{0} + 1 &\textbf{(17)}
			\end{dmatrix}\\
			&=2i\,\Im\left((\theta_0+1)(\theta_0-2)\alpha_2\frac{\z^{\theta_0}}{z}+\frac{(\theta_0-2)\zeta_4}{\theta_0}\frac{\z^{\theta_0+1}}{z^{\theta_0-1}}+\frac{(\theta_0-3)\zeta_3}{2\theta_0}\frac{\z^{\theta_0}}{z^{\theta_0-2}}\right)
		\end{align*}
		so $\p{\z}f(z)$ is purely imaginary, as expected.

		Then
		\begin{align*}
			\Re(\p{\z}f(z))=0,
		\end{align*}
		as expected. We also define
		\begin{align*}
			\zeta_7=\frac{1}{2 \, {\left(\theta_{0}^{3} + 3 \, \theta_{0}^{2} + 2 \, \theta_{0}\right)}}\bigg\{2\theta_0(\theta_0+1)\zeta_5-4\zeta_6 \bigg\}
		\end{align*}
		so that
		\begin{align*}
			\s{\phi}{\vec{\alpha}}=\begin{dmatrix}
			\frac{(2\theta_0-1)\zeta_2}{2\theta_0(\theta_0+1)} & 2 & 0 &\textbf{(6)}\\
			\zeta_7 & 3 & 0 &\textbf{(8)}\\
			(\theta_0-2)\alpha_2 & -\theta_{0} + 1 & \theta_{0} + 1 &\textbf{(10)}\\
			\frac{(\theta_0-2)\zeta_4}{\theta_0(\theta_0+2)} & -\theta_{0} + 1 & \theta_{0} + 2 &\textbf{(11)}\\
			\frac{(\theta_0-3)\zeta_3}{2\theta_0(\theta_0+1)} & -\theta_{0} + 2 & \theta_{0} + 1 &\textbf{(13)}\\
			(\theta_0+1)\bar{\alpha_2} & \theta_{0} & -\theta_{0} + 2 &\textbf{(15)}\\
			\frac{\bar{\zeta_3}}{2\theta_0} & \theta_{0} & -\theta_{0} + 3 &\textbf{(16)}\\
			\frac{\bar{\zeta_4}}{\theta_0} & \theta_{0} + 1 & -\theta_{0} + 2 &\textbf{(17)}
			\end{dmatrix}=\begin{dmatrix}
			\frac{{\left(2 \, \theta_{0} - 1\right)} \zeta_{2}}{2 \, {\left(\theta_{0} + 1\right)} \theta_{0}} & 2 & 0 \\
			\zeta_{7} & 3 & 0 \\
			\alpha_{2} {\left(\theta_{0} - 2\right)} & -\theta_{0} + 1 & \theta_{0} + 1 \\
			\frac{{\left(\theta_{0} - 2\right)} \zeta_{4}}{{\left(\theta_{0} + 2\right)} \theta_{0}} & -\theta_{0} + 1 & \theta_{0} + 2 \\
			\frac{{\left(\theta_{0} - 3\right)} \zeta_{3}}{2 \, {\left(\theta_{0} + 1\right)} \theta_{0}} & -\theta_{0} + 2 & \theta_{0} + 1 \\
			{\left(\theta_{0} + 1\right)} \overline{\alpha_{2}} & \theta_{0} & -\theta_{0} + 2 \\
			\frac{\overline{\zeta_{3}}}{2 \, \theta_{0}} & \theta_{0} & -\theta_{0} + 3 \\
			\frac{\overline{\zeta_{4}}}{\theta_{0}} & \theta_{0} + 1 & -\theta_{0} + 2
			\end{dmatrix}
		\end{align*}
		
		\section{Development of $\h_0$ up to order $4$}
		\small
		\begin{align*}
			&g^{-1}\otimes\s{\h_0}{\phi}=\\
			&\begin{dmatrix}
			\frac{\ccancel{2 \, A_{1} \overline{A_{0}}}}{\theta_{0}} & 0 & 1 \\
			-\frac{2 \, {\left(2 \, A_{0} {\left(\theta_{0} + 1\right)} {\left| A_{1} \right|}^{2} \overline{A_{0}} - A_{1} \theta_{0} \overline{A_{1}}\right)}}{\theta_{0}^{2} + \theta_{0}}=-\frac{2|A_1|^2}{\theta_0(\theta_0+1)} & 0 & 2 \\
			\mu_1 & 0 & 3 \\
			\mu_2 & 0 & 4 \\
			-\frac{\colorcancel{8 \, A_{0} \alpha_{1} \theta_{0} \overline{A_{0}}}{blue} + \ccancel{A_{0} C_{1} {\left(\theta_{0} - 2\right)}} - \colorcancel{8 \, A_{2} \theta_{0} \overline{A_{0}}}{blue}}{2 \, \theta_{0}^{2}} & 1 & 1 \\
			\mu_3 & 1 & 2 \\
			\mu_4 & 1 & 3 \\
			\mu_5 & 2 & 1 \\
			\mu_6 & 2 & 2 \\
			\mu_7 & 3 & 1 \\
			-\frac{\ccancel{E_{1} {\left(\theta_{0} - 2\right)} \overline{A_{0}}}}{2 \, \theta_{0}^{2}}& -2 \, \theta_{0} + 3 & 2 \, \theta_{0} + 1 \\
			-\frac{\ccancel{C_{1} {\left(\theta_{0} - 2\right)} \overline{A_{0}}}}{2 \, \theta_{0}^{2}} & -\theta_{0} + 1 & \theta_{0} + 1 \\
			\frac{4 \, {\left(\alpha_{2} \theta_{0}^{2} - \alpha_{2} \theta_{0} - 2 \, \alpha_{2}\right)} A_{0} \overline{A_{0}} -\colorcancel{ B_{1} {\left(\theta_{0} - 2\right)} \overline{A_{0}}}{blue} - \colorcancel{C_{1} {\left(\theta_{0} - 2\right)} \overline{A_{1}}}{blue}}{2 \, {\left(\theta_{0}^{2} + \theta_{0}\right)}} & -\theta_{0} + 1 & \theta_{0} + 2 \\
			\mu_8 & -\theta_{0} + 1 & \theta_{0} + 3 \\
			-\frac{\ccancel{C_{2} {\left(\theta_{0} - 3\right)} \overline{A_{0}}}}{2 \, \theta_{0}^{2}} & -\theta_{0} + 2 & \theta_{0} + 1 \\
			\mu_9 & -\theta_{0} + 2 & \theta_{0} + 2 \\
			-\frac{C_{1}^{2} {\left(\theta_{0} - 2\right)} + 8 \, A_{0} E_{1} {\left(\theta_{0} - 2\right)} - 32 \, {\left(\alpha_{7} \theta_{0}^{2} - 4 \, \alpha_{7} \theta_{0}\right)} A_{0} \overline{A_{0}} + 8 \, C_{3} {\left(\theta_{0} - 4\right)} \overline{A_{0}}}{16 \, \theta_{0}^{2}} & -\theta_{0} + 3 & \theta_{0} + 1 \\
			\frac{\ccancel{2 \, A_{0} A_{1}}}{\theta_{0}} & \theta_{0} & -\theta_{0} + 1 \\
			-\frac{\ccancel{4 \, A_{0}^{2} {\left| A_{1} \right|}^{2}}}{\theta_{0}} & \theta_{0} & -\theta_{0} + 2 \\
			-\frac{8 \, {\left(\theta_{0} \overline{\alpha_{2}} + \overline{\alpha_{2}}\right)} A_{0} \overline{A_{0}} + \ccancel{8 \, A_{0} A_{1} \overline{\alpha_{1}}} + \ccancel{8 \, A_{0}^{2} \overline{\alpha_{5}}} - \colorcancel{A_{0} \overline{B_{1}}}{blue} - \colorcancel{A_{1} \overline{C_{1}}}{blue}}{4 \, \theta_{0}} & \theta_{0} & -\theta_{0} + 3 \\
			\mu_{10}& \theta_{0} & -\theta_{0} + 4 \\
			-\frac{4 \, {\left(\theta_{0} \overline{\alpha_{7}} - 4 \, \overline{\alpha_{7}}\right)} A_{0} \overline{A_{0}} + \overline{A_{0}} \overline{E_{1}}}{2 \, {\left(\theta_{0} - 4\right)}} & \theta_{0} - 1 & -\theta_{0} + 5 \\
			-\frac{2 \, {\left(\ccancel{2 \, {\left(\alpha_{1} \theta_{0} + \alpha_{1}\right)} A_{0}^{2}} - 2 \, A_{0} A_{2} {\left(\theta_{0} + 1\right)} - A_{1}^{2} \theta_{0}\right)}}{\theta_{0}^{2} + \theta_{0}} & \theta_{0} + 1 & -\theta_{0} + 1 \\
			-\frac{4 \, {\left(\ccancel{A_{0} A_{1} {\left(3 \, \theta_{0} + 2\right)} {\left| A_{1} \right|}^{2}} + \ccancel{{\left(\alpha_{5} \theta_{0} + \alpha_{5}\right)} A_{0}^{2}}\right)}}{\theta_{0}^{2} + \theta_{0}} & \theta_{0} + 1 & -\theta_{0} + 2 \\
			\mu_{11} & \theta_{0} + 1 & -\theta_{0} + 3 \\
			\mu_{12} & \theta_{0} + 2 & -\theta_{0} + 1 
						\end{dmatrix}
			\end{align*}
			\begin{align*}
			\begin{dmatrix}
			\mu_{13} & \theta_{0} + 2 & -\theta_{0} + 2 \\
			\mu_{14} & \theta_{0} + 3 & -\theta_{0} + 1 \\
			-\frac{4 \, {\left(\theta_{0} \overline{\alpha_{7}} - 4 \, \overline{\alpha_{7}}\right)} \ccancel{A_{0}^{2}} + \ccancel{A_{0} \overline{E_{1}}}}{2 \, {\left(\theta_{0} - 4\right)}} & 2 \, \theta_{0} - 1 & -2 \, \theta_{0} + 5 \\
			-\frac{2 \, {\left({\left(\theta_{0} \overline{\alpha_{9}} + 2 \, \overline{\alpha_{9}}\right)} \ccancel{A_{0}^{2}} + 2 \, {\left(\theta_{0} \overline{\alpha_{2}} + \overline{\alpha_{2}}\right)} \ccancel{A_{0} A_{1}}\right)}}{\theta_{0}} & 2 \, \theta_{0} + 1 & -2 \, \theta_{0} + 3 \\
			-\frac{2 \, {\left(\theta_{0} \overline{\alpha_{2}} + \overline{\alpha_{2}}\right)} \ccancel{A_{0}^{2}}}{\theta_{0}} & 2 \, \theta_{0} & -2 \, \theta_{0} + 3 \\
			-\frac{2 \, {\left(\theta_{0} \overline{\alpha_{8}} + \overline{\alpha_{8}}\right)} \ccancel{A_{0}^{2}}}{\theta_{0}} & 2 \, \theta_{0} & -2 \, \theta_{0} + 4
			\end{dmatrix}
		\end{align*}
		\begin{align*}
			=\begin{dmatrix}
			-\frac{2|A_1|^2}{\theta_0(\theta_0+1)} & 0 & 2 & \textbf{(1)}\\
			\mu_1 & 0 & 3 & \textbf{(2)}\\
			\mu_2 & 0 & 4 & \textbf{(3)}\\
			\mu_3 & 1 & 2 & \textbf{(4)}\\
			\mu_4 & 1 & 3 & \textbf{(5)}\\
			\mu_5 & 2 & 1 & \textbf{(6)}\\
			\mu_6 & 2 & 2 & \textbf{(7)}\\
			\mu_7 & 3 & 1 & \textbf{(8)}\\
			\frac{(\theta_0-2)\alpha_2}{\theta_0} & -\theta_{0} + 1 & \theta_{0} + 2 & \textbf{(9)}\\
			\mu_8 & -\theta_{0} + 1 & \theta_{0} + 3 & \textbf{(10)}\\
			\mu_9 & -\theta_{0} + 2 & \theta_{0} + 2 & \textbf{(11)}\\
			-\frac{C_{1}^{2} {\left(\theta_{0} - 2\right)} + 8 \, A_{0} E_{1} {\left(\theta_{0} - 2\right)} - 32 \, {\left(\alpha_{7} \theta_{0}^{2} - 4 \, \alpha_{7} \theta_{0}\right)} A_{0} \overline{A_{0}} + 8 \, C_{3} {\left(\theta_{0} - 4\right)} \overline{A_{0}}}{16 \, \theta_{0}^{2}} & -\theta_{0} + 3 & \theta_{0} + 1 & \textbf{(12)}\\
			-\frac{\, (\theta_{0}+1) \overline{\alpha_{2}}} { \, \theta_{0}} & \theta_{0} & -\theta_{0} + 3 & \textbf{(13)}\\
			\mu_{10}& \theta_{0} & -\theta_{0} + 4 & \textbf{(14)}\\
			-\frac{4 \, {\left(\theta_{0} \overline{\alpha_{7}} - 4 \, \overline{\alpha_{7}}\right)} A_{0} \overline{A_{0}} + \overline{A_{0}} \overline{E_{1}}}{2 \, {\left(\theta_{0} - 4\right)}} & \theta_{0} - 1 & -\theta_{0} + 5 & \textbf{(15)}\\
			-\frac{2{A}_1^2}{\theta_{0}(\theta_{0}+1)} & \theta_{0} + 1 & -\theta_{0} + 1 & \textbf{(16)}\\
			\mu_{11} & \theta_{0} + 1 & -\theta_{0} + 3 & \textbf{(17)}\\
			\mu_{12} & \theta_{0} + 2 & -\theta_{0} + 1 & \textbf{(18)}\\
			\mu_{13} & \theta_{0} + 2 & -\theta_{0} + 2 & \textbf{(19)}\\
			\mu_{14} & \theta_{0} + 3 & -\theta_{0} + 1 & \textbf{(20)}
			\end{dmatrix}
		\end{align*}
		where
		\begin{align}\label{endend}
			&\mu_{1}=-\frac{1}{4 \, {\left(\theta_{0}^{3} + 3 \, \theta_{0}^{2} + 2 \, \theta_{0}\right)}}\bigg\{16 \, {\left(\theta_{0}^{2} + 2 \, \theta_{0}\right)} A_{0} {\left| A_{1} \right|}^{2} \overline{A_{1}} + 8 \, {\left(\theta_{0}^{2} \overline{\alpha_{5}} + 3 \, \theta_{0} \overline{\alpha_{5}} + 2 \, \overline{\alpha_{5}}\right)} A_{0} \overline{A_{0}} + 8 \, {\left(\theta_{0}^{2} \overline{\alpha_{1}} + 3 \, \theta_{0} \overline{\alpha_{1}} + 2 \, \overline{\alpha_{1}}\right)} A_{1} \overline{A_{0}}\nonumber\\
			& - 8 \, {\left(\theta_{0}^{2} + \theta_{0}\right)} A_{1} \overline{A_{2}} - {\left(\theta_{0}^{2} + 3 \, \theta_{0} + 2\right)} \overline{A_{0}} \overline{B_{1}}\bigg\}\nonumber\\
			&\mu_{2}=\frac{1}{12 \, {\left(\theta_{0}^{4} + 6 \, \theta_{0}^{3} + 11 \, \theta_{0}^{2} + 6 \, \theta_{0}\right)}}\bigg\{96 \, {\left(\theta_{0}^{3} \overline{\alpha_{1}} + 6 \, \theta_{0}^{2} \overline{\alpha_{1}} + 11 \, \theta_{0} \overline{\alpha_{1}} + 6 \, \overline{\alpha_{1}}\right)} A_{0} {\left| A_{1} \right|}^{2} \overline{A_{0}} - 48 \, {\left(\theta_{0}^{3} + 4 \, \theta_{0}^{2} + 3 \, \theta_{0}\right)} A_{0} {\left| A_{1} \right|}^{2} \overline{A_{2}}\nonumber\\
			& - 2 \, {\left(\theta_{0}^{3} + 6 \, \theta_{0}^{2} + 11 \, \theta_{0} + 6\right)} {\left| A_{1} \right|}^{2} \overline{A_{0}} \overline{C_{1}} - 24 \, {\left(\theta_{0}^{3} \overline{\alpha_{6}} + 6 \, \theta_{0}^{2} \overline{\alpha_{6}} + 11 \, \theta_{0} \overline{\alpha_{6}} + 6 \, \overline{\alpha_{6}}\right)} A_{0} \overline{A_{0}}\nonumber\\
			& - 24 \, {\left(\theta_{0}^{3} \overline{\alpha_{3}} + 6 \, \theta_{0}^{2} \overline{\alpha_{3}} + 11 \, \theta_{0} \overline{\alpha_{3}} + 6 \, \overline{\alpha_{3}}\right)} A_{1} \overline{A_{0}} + 24 \, {\left(\theta_{0}^{3} + 3 \, \theta_{0}^{2} + 2 \, \theta_{0}\right)} A_{1} \overline{A_{3}} + 3 \, {\left(\theta_{0}^{3} + 5 \, \theta_{0}^{2} + 6 \, \theta_{0}\right)} \overline{A_{1}} \overline{B_{1}}\nonumber\\
			& + 2 \, {\left(\theta_{0}^{3} + 6 \, \theta_{0}^{2} + 11 \, \theta_{0} + 6\right)} \overline{A_{0}} \overline{B_{3}} - 24 \, {\left({\left(\theta_{0}^{3} \overline{\alpha_{5}} + 5 \, \theta_{0}^{2} \overline{\alpha_{5}} + 6 \, \theta_{0} \overline{\alpha_{5}}\right)} A_{0} + {\left(\theta_{0}^{3} \overline{\alpha_{1}} + 5 \, \theta_{0}^{2} \overline{\alpha_{1}} + 6 \, \theta_{0} \overline{\alpha_{1}}\right)} A_{1}\right)} \overline{A_{1}}\bigg\}\nonumber\\
			&\mu_{3}=-\frac{16 \, A_{1} {\left(\theta_{0} + 1\right)} {\left| A_{1} \right|}^{2} \overline{A_{0}} + 8 \, A_{0} \alpha_{1} \theta_{0} \overline{A_{1}} - 4 \, {\left(\alpha_{2} \theta_{0}^{2} - \alpha_{2} \theta_{0} - 2 \, \alpha_{2}\right)} A_{0}^{2} + A_{0} B_{1} {\left(\theta_{0} - 2\right)} + 8 \, {\left(\alpha_{5} \theta_{0} + \alpha_{5}\right)} A_{0} \overline{A_{0}} - 8 \, A_{2} \theta_{0} \overline{A_{1}}}{2 \, {\left(\theta_{0}^{2} + \theta_{0}\right)}}\nonumber\\
			&\mu_{4}=\frac{1}{16 \, {\left(\theta_{0}^{4} + 3 \, \theta_{0}^{3} + 2 \, \theta_{0}^{2}\right)}}\bigg\{256 \, {\left(\theta_{0}^{3} + 3 \, \theta_{0}^{2} + 2 \, \theta_{0}\right)} A_{0} {\left| A_{1} \right|}^{4} \overline{A_{0}} - 128 \, {\left(\theta_{0}^{3} + 2 \, \theta_{0}^{2}\right)} A_{1} {\left| A_{1} \right|}^{2} \overline{A_{1}}\nonumber\\
			& + 32 \, {\left(\alpha_{9} \theta_{0}^{4} + \alpha_{9} \theta_{0}^{3} - 4 \, \alpha_{9} \theta_{0}^{2} - 4 \, \alpha_{9} \theta_{0}\right)} A_{0}^{2} - 8 \, {\left(\theta_{0}^{3} - \theta_{0}^{2} - 2 \, \theta_{0}\right)} A_{0} B_{2} + 16 \, {\left(\theta_{0}^{2} \overline{\alpha_{1}} - \theta_{0} \overline{\alpha_{1}} - 2 \, \overline{\alpha_{1}}\right)} A_{0} C_{1}\nonumber\\
			& + 64 \, {\left({\left(2 \, \alpha_{1} \overline{\alpha_{1}} - \beta\right)} \theta_{0}^{3} + 3 \, {\left(2 \, \alpha_{1} \overline{\alpha_{1}} - \beta\right)} \theta_{0}^{2} + 2 \, {\left(2 \, \alpha_{1} \overline{\alpha_{1}} - \beta\right)} \theta_{0}\right)} A_{0} \overline{A_{0}} - 64 \, {\left(\theta_{0}^{3} \overline{\alpha_{5}} + 3 \, \theta_{0}^{2} \overline{\alpha_{5}} + 2 \, \theta_{0} \overline{\alpha_{5}}\right)} A_{1} \overline{A_{0}}\nonumber\\
			& - 64 \, {\left(\theta_{0}^{3} \overline{\alpha_{1}} + 3 \, \theta_{0}^{2} \overline{\alpha_{1}} + 2 \, \theta_{0} \overline{\alpha_{1}}\right)} A_{2} \overline{A_{0}} - 64 \, {\left(\alpha_{5} \theta_{0}^{3} + 2 \, \alpha_{5} \theta_{0}^{2}\right)} A_{0} \overline{A_{1}} + 8 \, {\left(\theta_{0}^{3} + 3 \, \theta_{0}^{2} + 2 \, \theta_{0}\right)} \overline{A_{0}} \overline{B_{2}} - {\left(\theta_{0}^{3} + \theta_{0}^{2} - 4 \, \theta_{0} - 4\right)} C_{1} \overline{C_{1}}\nonumber\\
			& - 64 \, {\left({\left(\alpha_{1} \theta_{0}^{3} + \alpha_{1} \theta_{0}^{2}\right)} A_{0} - {\left(\theta_{0}^{3} + \theta_{0}^{2}\right)} A_{2}\right)} \overline{A_{2}}\bigg\}\nonumber\\
			&\mu_{5}=-\frac{{\left(\theta_{0}^{2} - 5 \, \theta_{0}\right)} A_{1} C_{1} + 2 \, {\left(\theta_{0}^{2} - 2 \, \theta_{0} - 3\right)} A_{0} C_{2} + 24 \, {\left(\alpha_{3} \theta_{0}^{2} + \alpha_{3} \theta_{0}\right)} A_{0} \overline{A_{0}} + 24 \, {\left(\alpha_{1} \theta_{0}^{2} + \alpha_{1} \theta_{0}\right)} A_{1} \overline{A_{0}} - 24 \, {\left(\theta_{0}^{2} + \theta_{0}\right)} A_{3} \overline{A_{0}}}{4 \, {\left(\theta_{0}^{3} + \theta_{0}^{2}\right)}}\nonumber\\
			&\mu_{6}=-\frac{1}{4 \, {\left(\theta_{0}^{4} + 2 \, \theta_{0}^{3} + \theta_{0}^{2}\right)}}\bigg\{2 \, {\left(\theta_{0}^{3} + 5 \, \theta_{0} + 6\right)} A_{0} C_{1} {\left| A_{1} \right|}^{2} - 96 \, {\left(\alpha_{1} \theta_{0}^{3} + 2 \, \alpha_{1} \theta_{0}^{2} + \alpha_{1} \theta_{0}\right)} A_{0} {\left| A_{1} \right|}^{2} \overline{A_{0}}\\
			& + 48 \, {\left(\theta_{0}^{3} + 2 \, \theta_{0}^{2} + \theta_{0}\right)} A_{2} {\left| A_{1} \right|}^{2} \overline{A_{0}} - 8 \, {\left(\alpha_{8} \theta_{0}^{4} - \alpha_{8} \theta_{0}^{3} - 5 \, \alpha_{8} \theta_{0}^{2} - 3 \, \alpha_{8} \theta_{0}\right)} A_{0}^{2} + {\left(\theta_{0}^{3} - 5 \, \theta_{0}^{2}\right)} A_{1} B_{1} + 2 \, {\left(\theta_{0}^{3} - 2 \, \theta_{0}^{2} - 3 \, \theta_{0}\right)} A_{0} B_{3}\nonumber\\
			& + 24 \, {\left(\alpha_{6} \theta_{0}^{3} + 2 \, \alpha_{6} \theta_{0}^{2} + \alpha_{6} \theta_{0}\right)} A_{0} \overline{A_{0}} - 24 \, {\left(\theta_{0}^{3} + \theta_{0}^{2}\right)} A_{3} \overline{A_{1}} - 8 \, \bigg({\left(2 \, \alpha_{2} \theta_{0}^{4} - 2 \, \alpha_{2} \theta_{0}^{3} - 7 \, \alpha_{2} \theta_{0}^{2} - 3 \, \alpha_{2} \theta_{0}\right)} A_{0}\nonumber\\
			& - 3 \, {\left(\alpha_{5} \theta_{0}^{3} + 2 \, \alpha_{5} \theta_{0}^{2} + \alpha_{5} \theta_{0}\right)} \overline{A_{0}}\bigg) A_{1} + 24 \, {\left({\left(\alpha_{3} \theta_{0}^{3} + \alpha_{3} \theta_{0}^{2}\right)} A_{0} + {\left(\alpha_{1} \theta_{0}^{3} + \alpha_{1} \theta_{0}^{2}\right)} A_{1}\right)} \overline{A_{1}}\bigg\}\nonumber\\
			&\mu_{7}=\frac{1}{6 \, {\left(\theta_{0}^{4} + 3 \, \theta_{0}^{3} + 2 \, \theta_{0}^{2}\right)}}\bigg\{12 \, {\left(\alpha_{7} \theta_{0}^{4} - \alpha_{7} \theta_{0}^{3} - 10 \, \alpha_{7} \theta_{0}^{2} - 8 \, \alpha_{7} \theta_{0}\right)} A_{0}^{2} - 2 \, {\left(\theta_{0}^{3} - 3 \, \theta_{0}^{2} - 10 \, \theta_{0}\right)} A_{1} C_{2}\nonumber\\
			& - 3 \, {\left(\theta_{0}^{3} - \theta_{0}^{2} - 10 \, \theta_{0} - 8\right)} A_{0} C_{3} + 48 \, {\left({\left(\alpha_{1}^{2} - \alpha_{4}\right)} \theta_{0}^{3} + 3 \, {\left(\alpha_{1}^{2} - \alpha_{4}\right)} \theta_{0}^{2} + 2 \, {\left(\alpha_{1}^{2} - \alpha_{4}\right)} \theta_{0}\right)} A_{0} \overline{A_{0}}\\
			& - 48 \, {\left(\alpha_{3} \theta_{0}^{3} + 3 \, \alpha_{3} \theta_{0}^{2} + 2 \, \alpha_{3} \theta_{0}\right)} A_{1} \overline{A_{0}} - 48 \, {\left(\alpha_{1} \theta_{0}^{3} + 3 \, \alpha_{1} \theta_{0}^{2} + 2 \, \alpha_{1} \theta_{0}\right)} A_{2} \overline{A_{0}} + 48 \, {\left(\theta_{0}^{3} + 3 \, \theta_{0}^{2} + 2 \, \theta_{0}\right)} A_{4} \overline{A_{0}}\nonumber \\
			&- 3 \, {\left({\left(\alpha_{1} \theta_{0}^{3} + 3 \, \alpha_{1} \theta_{0}^{2} + 2 \, \alpha_{1} \theta_{0}\right)} A_{0} - 4 \, {\left(\theta_{0}^{2} + \theta_{0}\right)} A_{2}\right)} C_{1}\bigg\}\nonumber\\
			&\mu_{8}=\frac{1}{2 \, {\left(\theta_{0}^{5} + 4 \, \theta_{0}^{4} + 5 \, \theta_{0}^{3} + 2 \, \theta_{0}^{2}\right)}}\bigg\{4 \, {\left(\alpha_{9} \theta_{0}^{5} + 2 \, \alpha_{9} \theta_{0}^{4} - 3 \, \alpha_{9} \theta_{0}^{3} - 8 \, \alpha_{9} \theta_{0}^{2} - 4 \, \alpha_{9} \theta_{0}\right)} A_{0} \overline{A_{0}} - {\left(\theta_{0}^{4} - 3 \, \theta_{0}^{2} - 2 \, \theta_{0}\right)} B_{2} \overline{A_{0}}\nonumber\\
			& + 4 \, {\left(\alpha_{2} \theta_{0}^{5} + \alpha_{2} \theta_{0}^{4} - 4 \, \alpha_{2} \theta_{0}^{3} - 4 \, \alpha_{2} \theta_{0}^{2}\right)} A_{0} \overline{A_{1}} - {\left(\theta_{0}^{4} - 4 \, \theta_{0}^{2}\right)} B_{1} \overline{A_{1}} + {\left(2 \, {\left(\theta_{0}^{3} \overline{\alpha_{1}} - 3 \, \theta_{0} \overline{\alpha_{1}} - 2 \, \overline{\alpha_{1}}\right)} \overline{A_{0}} - {\left(\theta_{0}^{4} - 3 \, \theta_{0}^{2} - 2 \, \theta_{0}\right)} \overline{A_{2}}\right)} C_{1}\bigg\}\nonumber\\
			&\mu_{9}=\frac{1}{2 \, {\left(\theta_{0}^{3} + \theta_{0}^{2}\right)}}\bigg\{2 \, C_{1} {\left(\theta_{0} - 3\right)} {\left| A_{1} \right|}^{2} \overline{A_{0}} + 4 \, {\left(\alpha_{8} \theta_{0}^{3} - 2 \, \alpha_{8} \theta_{0}^{2} - 3 \, \alpha_{8} \theta_{0}\right)} A_{0} \overline{A_{0}} + 4 \, {\left(\alpha_{2} \theta_{0}^{3} - 2 \, \alpha_{2} \theta_{0}^{2} - 3 \, \alpha_{2} \theta_{0}\right)} A_{1} \overline{A_{0}}\\
			& - {\left(\theta_{0}^{2} - 3 \, \theta_{0}\right)} B_{3} \overline{A_{0}} - {\left(\theta_{0}^{2} - 3 \, \theta_{0}\right)} C_{2} \overline{A_{1}}\bigg\}\nonumber\\
			&\mu_{10}=\frac{48 \, A_{0}^{2} {\left| A_{1} \right|}^{2} \overline{\alpha_{1}} - 4 \, A_{0} {\left| A_{1} \right|}^{2} \overline{C_{1}} - 12 \, A_{0} \theta_{0} \overline{A_{1}} \overline{\alpha_{2}} - 12 \, {\left(\theta_{0} \overline{\alpha_{8}} + \overline{\alpha_{8}}\right)} A_{0} \overline{A_{0}} - 12 \, A_{0} A_{1} \overline{\alpha_{3}} - 12 \, A_{0}^{2} \overline{\alpha_{6}} + A_{0} \overline{B_{3}} + A_{1} \overline{C_{2}}}{6 \, \theta_{0}}\nonumber\\
			&\mu_{11}=\frac{1}{2 \, {\left(\theta_{0}^{2} + \theta_{0}\right)}}\bigg\{32 \, A_{0}^{2} {\left(\theta_{0} + 1\right)} {\left| A_{1} \right|}^{4} - 4 \, A_{1}^{2} \theta_{0} \overline{\alpha_{1}} + 8 \, {\left({\left(2 \, \alpha_{1} \overline{\alpha_{1}} - \beta\right)} \theta_{0} + 2 \, \alpha_{1} \overline{\alpha_{1}} - \beta\right)} A_{0}^{2} - 8 \, {\left(\theta_{0} \overline{\alpha_{1}} + \overline{\alpha_{1}}\right)} A_{0} A_{2} \nonumber\\
			&- 4 \, {\left(\theta_{0}^{2} \overline{\alpha_{9}} + 3 \, \theta_{0} \overline{\alpha_{9}} + 2 \, \overline{\alpha_{9}}\right)} A_{0} \overline{A_{0}} + A_{1} \theta_{0} \overline{B_{1}} + A_{0} {\left(\theta_{0} + 1\right)} \overline{B_{2}} - 4 \, \bigg({\left(3 \, \theta_{0} \overline{\alpha_{5}} + 2 \, \overline{\alpha_{5}}\right)} A_{0} \nonumber\\
				&+ {\left(\theta_{0}^{2} \overline{\alpha_{2}} + 3 \, \theta_{0} \overline{\alpha_{2}} + 2 \, \overline{\alpha_{2}}\right)} \overline{A_{0}}\bigg) A_{1} - {\left({\left(\alpha_{1} \theta_{0} + \alpha_{1}\right)} A_{0} - A_{2} {\left(\theta_{0} + 1\right)}\right)} \overline{C_{1}}\bigg\}\nonumber\\
			&\mu_{12}=-\frac{2 \, {\left(3 \, {\left(\alpha_{3} \theta_{0}^{2} + 3 \, \alpha_{3} \theta_{0} + 2 \, \alpha_{3}\right)} A_{0}^{2} + {\left(5 \, \alpha_{1} \theta_{0}^{2} + 13 \, \alpha_{1} \theta_{0} + 6 \, \alpha_{1}\right)} A_{0} A_{1} - {\left(3 \, \theta_{0}^{2} + 5 \, \theta_{0}\right)} A_{1} A_{2} - 3 \, {\left(\theta_{0}^{2} + 3 \, \theta_{0} + 2\right)} A_{0} A_{3}\right)}}{\theta_{0}^{3} + 3 \, \theta_{0}^{2} + 2 \, \theta_{0}}\nonumber\\
			&\mu_{13}=-\frac{2}{\theta_{0}^{3} + 3 \, \theta_{0}^{2} + 2 \, \theta_{0}}\bigg\{4 \, {\left(2 \, \theta_{0}^{2} + 5 \, \theta_{0} + 3\right)} A_{0} A_{2} {\left| A_{1} \right|}^{2} + 3 \, {\left(\alpha_{6} \theta_{0}^{2} + 3 \, \alpha_{6} \theta_{0} + 2 \, \alpha_{6}\right)} A_{0}^{2} + {\left(5 \, \alpha_{5} \theta_{0}^{2} + 13 \, \alpha_{5} \theta_{0} + 6 \, \alpha_{5}\right)} A_{0} A_{1}\nonumber\\
			& - 4 \, {\left(3 \, {\left(\alpha_{1} \theta_{0}^{2} + 3 \, \alpha_{1} \theta_{0} + 2 \, \alpha_{1}\right)} A_{0}^{2} - {\left(\theta_{0}^{2} + 2 \, \theta_{0}\right)} A_{1}^{2}\right)} {\left| A_{1} \right|}^{2}\bigg\}\nonumber\\
			&\mu_{14}=\frac{2}{{\theta_{0}^{4} + 6 \, \theta_{0}^{3} + 11 \, \theta_{0}^{2} + 6 \, \theta_{0}}}\bigg\{ \, 4 \, {\left({\left(\alpha_{1}^{2} - \alpha_{4}\right)} \theta_{0}^{3} + 6 \, {\left(\alpha_{1}^{2} - \alpha_{4}\right)} \theta_{0}^{2} + 6 \, \alpha_{1}^{2} + 11 \, {\left(\alpha_{1}^{2} - \alpha_{4}\right)} \theta_{0} - 6 \, \alpha_{4}\right)} A_{0}^{2} \nonumber\\
			&- {\left(7 \, \alpha_{3} \theta_{0}^{3} + 39 \, \alpha_{3} \theta_{0}^{2} + 62 \, \alpha_{3} \theta_{0} + 24 \, \alpha_{3}\right)} A_{0} A_{1} - 3 \, {\left(\alpha_{1} \theta_{0}^{3} + 5 \, \alpha_{1} \theta_{0}^{2} + 6 \, \alpha_{1} \theta_{0}\right)} A_{1}^{2}\nonumber\\
			& - 2 \, {\left(3 \, \alpha_{1} \theta_{0}^{3} + 16 \, \alpha_{1} \theta_{0}^{2} + 25 \, \alpha_{1} \theta_{0} + 12 \, \alpha_{1}\right)} A_{0} A_{2} + 2 \, {\left(\theta_{0}^{3} + 4 \, \theta_{0}^{2} + 3 \, \theta_{0}\right)} A_{2}^{2} + 2 \, {\left(2 \, \theta_{0}^{3} + 9 \, \theta_{0}^{2} + 10 \, \theta_{0}\right)} A_{1} A_{3}\nonumber\\
			& + 4 \, {\left(\theta_{0}^{3} + 6 \, \theta_{0}^{2} + 11 \, \theta_{0} + 6\right)} A_{0} A_{4}\bigg\}
		\end{align}
		\normalsize
		\textbf{We will not need all these $\mu$ at all.}
		
		Now, recall that
		\begin{align*}
			\alpha_7=\frac{1}{8\theta_0(\theta_0-4)}\s{\vec{C}_1}{\vec{C}_1},\quad \vec{E}_1=-\frac{1}{2\theta_0}\s{\vec{C}_1}{\vec{C}_1}\bar{\vec{A}_0},\quad \s{\vec{C}_1}{\vec{C}_1}+8\s{\bar{\vec{A}_0}}{\vec{C}_3}=0
		\end{align*}
		so that
		\begin{align*}
			\s{\vec{A}_0}{\vec{E}_1}=-\frac{1}{4\theta_0}\s{\vec{C}_1}{\vec{C}_1}=-2(\theta_0-4)\alpha_7
		\end{align*}
		Therefore, we have
		\begin{align}\label{invh015}
			\textbf{(12)}&=-\frac{C_{1}^{2} {\left(\theta_{0} - 2\right)} + 8 \, A_{0} E_{1} {\left(\theta_{0} - 2\right)} - 32 \, {\left(\alpha_{7} \theta_{0}^{2} - 4 \, \alpha_{7} \theta_{0}\right)} A_{0} \overline{A_{0}} + 8 \, C_{3} {\left(\theta_{0} - 4\right)} \overline{A_{0}}}{16 \, \theta_{0}^{2}}\nonumber\\
			&=-\frac{1}{16\theta_0^2}\left((\theta_0-2)\s{\vec{C}_1}{\vec{C}_1}-\frac{2(\theta_0-2)}{\theta_0}\s{\vec{C}_1}{\vec{C}_1}-16\theta_0(\theta_0-4)\alpha_7-(\theta_0-4)\s{\vec{C}_1}{\vec{C}_1}\right)\nonumber\\
			&=-\frac{1}{16\theta_0^2}\left((\theta_0-2)\s{\vec{C}_1}{\vec{C}_1}-\frac{2(\theta_0-2)}{\theta_0}\s{\vec{C}_1}{\vec{C}_1}-{2}\s{\vec{C}_1}{\vec{C}_1}-(\theta_0-4)\s{\vec{C}_1}{\vec{C}_1}\right)\nonumber\\
			&=\frac{(\theta_0-2)}{8\theta_0^3}\s{\vec{C}_1}{\vec{C}_1}=\frac{(\theta_0-2)(\theta_0-4)}{\theta_0^2}\alpha_7.
		\end{align}
		as
		\begin{align*}
			(\theta_0-2)-\frac{2(\theta_0-2)}{\theta_0}-\frac{2}{\theta_0}-(\theta_0-4)=2-2+\frac{4}{\theta_0}-\frac{2}{\theta_0}=\frac{2}{\theta_0}
		\end{align*}
		Then
		\begin{align*}
			\textbf{(15)}=
			-\frac{4 \, {\left(\theta_{0} \overline{\alpha_{7}} - 4 \, \overline{\alpha_{7}}\right)} A_{0} \overline{A_{0}} + \overline{A_{0}} \overline{E_{1}}}{2 \, {\left(\theta_{0} - 4\right)}}=-\frac{1}{2(\theta_0-4)}\left(2(\theta_0-4)\bar{\alpha_7}-2(\theta_0-4)\alpha_7\right)=0
		\end{align*}
	    Finally, we have (the $\mu$ coefficients do not seem to vanish so easily, even if there are numerous cancellations)
	    \begin{align}
	    	g^{-1}\otimes\s{\h_0}{\phi}=\begin{dmatrix}
	    	-\frac{2|A_1|^2}{\theta_0(\theta_0+1)} & 0 & 2 & \textbf{(1)}\\
	    	\mu_1 & 0 & 3 & \textbf{(2)}\\
	    	\mu_2 & 0 & 4 & \textbf{(3)}\\
	    	\mu_3 & 1 & 2 & \textbf{(4)}\\
	    	\mu_4 & 1 & 3 & \textbf{(5)}\\
	    	\mu_5 & 2 & 1 & \textbf{(6)}\\
	    	\mu_6 & 2 & 2 & \textbf{(7)}\\
	    	\mu_7 & 3 & 1 & \textbf{(8)}\\
	    	\frac{(\theta_0-2)\alpha_2}{\theta_0} & -\theta_{0} + 1 & \theta_{0} + 2 & \textbf{(9)}\\
	    	\mu_8 & -\theta_{0} + 1 & \theta_{0} + 3 & \textbf{(10)}\\
	    	\mu_9 & -\theta_{0} + 2 & \theta_{0} + 2 & \textbf{(11)}\\
	    	\frac{(\theta_0-2)(\theta_0-4)\alpha_7}{\theta_0^2} & -\theta_{0} + 3 & \theta_{0} + 1 & \textbf{(12)}\\
	    	-\frac{\, (\theta_{0}+1) \overline{\alpha_{2}}} { \, \theta_{0}} & \theta_{0} & -\theta_{0} + 3 & \textbf{(13)}\\
	    	\mu_{10}& \theta_{0} & -\theta_{0} + 4 & \textbf{(14)}\\
	    	-\frac{2\zeta_0}{\theta_{0}(\theta_{0}+1)} & \theta_{0} + 1 & -\theta_{0} + 1 & \textbf{(16)}\\
	    	\mu_{11} & \theta_{0} + 1 & -\theta_{0} + 3 & \textbf{(17)}\\
	    	\mu_{12} & \theta_{0} + 2 & -\theta_{0} + 1 & \textbf{(18)}\\
	    	\mu_{13} & \theta_{0} + 2 & -\theta_{0} + 2 & \textbf{(19)}\\
	    	\mu_{14} & \theta_{0} + 3 & -\theta_{0} + 1 & \textbf{(20)}
	    	\end{dmatrix}
	    	=\begin{dmatrix}
	    	-\frac{2 \, {\left| A_{1} \right|}^{2}}{{\left(\theta_{0} + 1\right)} \theta_{0}} & 0 & 2 \\
	    	\mu_{1} & 0 & 3 \\
	    	\mu_{2} & 0 & 4 \\
	    	\mu_{3} & 1 & 2 \\
	    	\mu_{4} & 1 & 3 \\
	    	\mu_{5} & 2 & 1 \\
	    	\mu_{6} & 2 & 2 \\
	    	\mu_{7} & 3 & 1 \\
	    	\frac{\alpha_{2} {\left(\theta_{0} - 2\right)}}{\theta_{0}} & -\theta_{0} + 1 & \theta_{0} + 2 \\
	    	\mu_{8} & -\theta_{0} + 1 & \theta_{0} + 3 \\
	    	\mu_{9} & -\theta_{0} + 2 & \theta_{0} + 2 \\
	    	\frac{\alpha_{7} {\left(\theta_{0} - 2\right)} {\left(\theta_{0} - 4\right)}}{\theta_{0}^{2}} & -\theta_{0} + 3 & \theta_{0} + 1 \\
	    	-\frac{{\left(\theta_{0} + 1\right)} \overline{\alpha_{2}}}{\theta_{0}} & \theta_{0} & -\theta_{0} + 3 \\
	    	\mu_{10} & \theta_{0} & -\theta_{0} + 4 \\
	    	-\frac{2 \, \zeta_{0}}{{\left(\theta_{0} + 1\right)} \theta_{0}} & \theta_{0} + 1 & -\theta_{0} + 1 \\
	    	\mu_{11} & \theta_{0} + 1 & -\theta_{0} + 3 \\
	    	\mu_{12} & \theta_{0} + 2 & -\theta_{0} + 1 \\
	    	\mu_{13} & \theta_{0} + 2 & -\theta_{0} + 2 \\
	    	\mu_{14} & \theta_{0} + 3 & -\theta_{0} + 1	    	
	    	\end{dmatrix}
	    \end{align}
	    Code check the right. Observe that Sage enjoys mirror symmetry.

	    \section{Second order development of the cancellation law}
	    Finally,
	    if
	    \begin{align*}
	    \vec{\beta}=\mathscr{I}_{\phi}(\vec{\alpha})-g^{-1}\otimes\left(\bar{\partial}|\phi|^2\h_0-2\s{\phi}{\h_0}\otimes\bar{\partial}\phi\right)=\vec{F}(z)dz
	    \end{align*}
	    where as usual
	    \begin{align*}
	    	\vec{\alpha}=\partial\H+|\H|^2\partial\phi+2\,g^{-1}\otimes\s{\H}{\h_0}
	    \end{align*}
        we have
        \begin{align*}
        	\Re(\p{\z}\vec{F}(z))=0,
        \end{align*}
        and
        \footnotesize
        \begin{align*}
            &\Re(\p{\z}\vec{F}(z))=\\
            &\begin{dmatrix}
            \frac{\mu_{1} \theta_{0}^{2} + 2 \, \mu_{1} \theta_{0} + 2 \, \overline{\alpha_{5}}}{\theta_{0}} & \overline{A_{0}} & 0 & \theta_{0} + 1 \\
            \frac{\mu_{1} \theta_{0}^{3} + 5 \, \mu_{1} \theta_{0}^{2} + 2 \, {\left(3 \, \mu_{1} + \overline{\alpha_{5}}\right)} \theta_{0} + 6 \, \overline{\alpha_{5}}}{\theta_{0}^{2} + 2 \, \theta_{0}} & \overline{A_{1}} & 0 & \theta_{0} + 2 \\
            \frac{\theta_{0}^{3} \overline{\mu_{13}} - 4 \, {\left(\theta_{0} + 3\right)} {\left| A_{1} \right|}^{2} \overline{\zeta_{0}} + 3 \, \theta_{0}^{2} \overline{\mu_{13}} + 2 \, \theta_{0} \overline{\mu_{13}}}{\theta_{0}^{3} + 3 \, \theta_{0}^{2} + 2 \, \theta_{0}} & A_{0} & 0 & \theta_{0} + 2 \\
            \frac{\theta_{0}^{3} \overline{\mu_{12}} + 3 \, \theta_{0}^{2} \overline{\mu_{12}} + 2 \, \theta_{0} \overline{\mu_{12}} + 4 \, {\left(2 \, \theta_{0} + 3\right)} \overline{\zeta_{1}}}{\theta_{0}^{3} + 3 \, \theta_{0}^{2} + 2 \, \theta_{0}} & A_{1} & 0 & \theta_{0} + 2 \\
            \frac{\mu_{2} \theta_{0}^{4} + 6 \, \mu_{2} \theta_{0}^{3} + {\left(11 \, \mu_{2} + 3 \, \overline{\alpha_{6}}\right)} \theta_{0}^{2} - 2 \, {\left(3 \, \theta_{0}^{2} \overline{\alpha_{1}} + 13 \, \theta_{0} \overline{\alpha_{1}} + 12 \, \overline{\alpha_{1}}\right)} {\left| A_{1} \right|}^{2} + 3 \, {\left(2 \, \mu_{2} + 3 \, \overline{\alpha_{6}}\right)} \theta_{0} + 6 \, \overline{\alpha_{6}}}{\theta_{0}^{3} + 3 \, \theta_{0}^{2} + 2 \, \theta_{0}} & \overline{A_{0}} & 0 & \theta_{0} + 2 \\
            -\frac{\theta_{0} - 1}{2 \, \theta_{0}^{2}} & B_{1} & 1 & \theta_{0} \\
            \frac{\mu_{3} \theta_{0}^{2} + \mu_{3} \theta_{0} + 2 \, \alpha_{5}}{\theta_{0}} & \overline{A_{0}} & 1 & \theta_{0} \\
            -\frac{2 \, {\left(\alpha_{2} \theta_{0}^{2} - \alpha_{2}\right)}}{\theta_{0}} & A_{0} & 1 & \theta_{0} \\
            -\frac{\theta_{0} - 1}{\theta_{0}^{2}} & B_{2} & 1 & \theta_{0} + 1 \\
            \frac{\mu_{3} \theta_{0}^{3} + 3 \, \mu_{3} \theta_{0}^{2} + 2 \, {\left(\alpha_{5} + \mu_{3}\right)} \theta_{0} + 4 \, \alpha_{5}}{\theta_{0}^{2} + \theta_{0}} & \overline{A_{1}} & 1 & \theta_{0} + 1 \\
            -\frac{\alpha_{9} \theta_{0}^{4} + {\left(\alpha_{9} - 2 \, \overline{\mu_{11}}\right)} \theta_{0}^{3} - 2 \, {\left(2 \, \alpha_{9} + \overline{\mu_{11}}\right)} \theta_{0}^{2} + 4 \, \alpha_{1} \theta_{0} \overline{\zeta_{0}} - 4 \, \alpha_{9} \theta_{0} + {\left(\theta_{0}^{2} - \theta_{0} - 2\right)} \zeta_{4}}{\theta_{0}^{3} + \theta_{0}^{2}} & A_{0} & 1 & \theta_{0} + 1 \\
            \frac{4 \, \mu_{4} \theta_{0}^{5} + 12 \, \mu_{4} \theta_{0}^{4} - 32 \, {\left(\theta_{0}^{3} + 2 \, \theta_{0}^{2}\right)} {\left| A_{1} \right|}^{4} - 8 \, {\left(2 \, \alpha_{1} \overline{\alpha_{1}} - 2 \, \beta - \mu_{4}\right)} \theta_{0}^{3} - 16 \, {\left(\alpha_{1} \overline{\alpha_{1}} - \beta\right)} \theta_{0}^{2} - {\left(\theta_{0}^{3} - 3 \, \theta_{0} - 2\right)} {\left| C_{1} \right|}^{2}}{4 \, {\left(\theta_{0}^{4} + \theta_{0}^{3}\right)}} & \overline{A_{0}} & 1 & \theta_{0} + 1 \\
            -\frac{2 \, \theta_{0} - 3}{4 \, \theta_{0}^{2}} & B_{3} & 2 & \theta_{0} \\
            -\frac{2 \, \alpha_{8} \theta_{0}^{3} - 2 \, {\left(2 \, \alpha_{8} + 3 \, \overline{\mu_{10}}\right)} \theta_{0}^{2} - 6 \, \alpha_{8} \theta_{0} + {\left(\theta_{0} - 3\right)} \zeta_{3}}{2 \, \theta_{0}^{2}} & A_{0} & 2 & \theta_{0} \\
            \frac{\mu_{6} \theta_{0}^{2} - 6 \, \alpha_{1} {\left| A_{1} \right|}^{2} + \mu_{6} \theta_{0} + 3 \, \alpha_{6}}{\theta_{0}} & \overline{A_{0}} & 2 & \theta_{0} \\
            -\frac{2 \, \alpha_{2} \theta_{0}^{2} - \alpha_{2} \theta_{0} - 3 \, \alpha_{2}}{\theta_{0}} & A_{1} & 2 & \theta_{0} \\
            \frac{2 \, \mu_{5} \theta_{0}^{4} - 4 \, \mu_{5} \theta_{0}^{3} - 6 \, \mu_{5} \theta_{0}^{2} + {\left(5 \, \theta_{0} - 12\right)} \zeta_{2}}{2 \, {\left(\theta_{0}^{3} - 3 \, \theta_{0}^{2}\right)}} & \overline{A_{1}} & 2 & \theta_{0} \\
            \frac{2 \, \mu_{5} \theta_{0}^{3} + 2 \, \mu_{5} \theta_{0}^{2} + 3 \, \zeta_{2}}{2 \, {\left(\theta_{0}^{2} + \theta_{0}\right)}} & \overline{A_{0}} & 2 & \theta_{0} - 1 \\
            \frac{\mu_{7} \theta_{0}^{2} - \theta_{0} \zeta_{7} + \zeta_{5}}{\theta_{0}} & \overline{A_{0}} & 3 & \theta_{0} - 1 \\
            -\frac{{\left(\alpha_{9} - 2 \, \mu_{8}\right)} \theta_{0}^{3} + 2 \, {\left(\alpha_{9} - 3 \, \mu_{8}\right)} \theta_{0}^{2} - 4 \, {\left(\alpha_{9} + \mu_{8}\right)} \theta_{0} + {\left(\theta_{0} - 2\right)} \zeta_{4} - 8 \, \alpha_{9}}{\theta_{0}^{2} + 2 \, \theta_{0}} & \overline{A_{0}} & -\theta_{0} + 1 & 2 \, \theta_{0} + 1 \\
            -\frac{2 \, {\left(\alpha_{8} - 2 \, \mu_{9}\right)} \theta_{0}^{3} - 2 \, {\left(\alpha_{8} + 3 \, \mu_{9}\right)} \theta_{0}^{2} - 2 \, {\left(5 \, \alpha_{8} + \mu_{9}\right)} \theta_{0} + {\left(\theta_{0} - 3\right)} \zeta_{3} - 6 \, \alpha_{8}}{2 \, {\left(\theta_{0}^{2} + \theta_{0}\right)}} & \overline{A_{0}} & -\theta_{0} + 2 & 2 \, \theta_{0} \\
            \frac{\theta_{0} - 4}{4 \, \theta_{0}} & E_{1} & -\theta_{0} + 3 & 2 \, \theta_{0} - 1 \\
            \frac{\alpha_{7} \theta_{0}^{2} - 8 \, \alpha_{7} \theta_{0} + 16 \, \alpha_{7}}{\theta_{0}} & \overline{A_{0}} & -\theta_{0} + 3 & 2 \, \theta_{0} - 1 \\
            -\frac{\theta_{0} - 1}{2 \, \theta_{0}^{2}} & \overline{B_{1}} & \theta_{0} & 1 \\
            \frac{\theta_{0}^{2} \overline{\mu_{3}} + \theta_{0} \overline{\mu_{3}} + 2 \, \overline{\alpha_{5}}}{\theta_{0}} & A_{0} & \theta_{0} & 1 \\
            -\frac{2 \, {\left(\theta_{0}^{2} \overline{\alpha_{2}} - \overline{\alpha_{2}}\right)}}{\theta_{0}} & \overline{A_{0}} & \theta_{0} & 1 
            \end{dmatrix}
            \end{align*}
            \begin{align*}
            \begin{dmatrix}
            -\frac{2 \, \theta_{0} - 3}{4 \, \theta_{0}^{2}} & \overline{B_{3}} & \theta_{0} & 2 \\
            -\frac{2 \, \theta_{0}^{3} \overline{\alpha_{8}} - 2 \, {\left(3 \, \mu_{10} + 2 \, \overline{\alpha_{8}}\right)} \theta_{0}^{2} - 6 \, \theta_{0} \overline{\alpha_{8}} + {\left(\theta_{0} - 3\right)} \overline{\zeta_{3}}}{2 \, \theta_{0}^{2}} & \overline{A_{0}} & \theta_{0} & 2 \\
            -\frac{2 \, \theta_{0}^{2} \overline{\alpha_{2}} - \theta_{0} \overline{\alpha_{2}} - 3 \, \overline{\alpha_{2}}}{\theta_{0}} & \overline{A_{1}} & \theta_{0} & 2 \\
            -\frac{6 \, {\left| A_{1} \right|}^{2} \overline{\alpha_{1}} - \theta_{0}^{2} \overline{\mu_{6}} - \theta_{0} \overline{\mu_{6}} - 3 \, \overline{\alpha_{6}}}{\theta_{0}} & A_{0} & \theta_{0} & 2 \\
            \frac{2 \, \theta_{0}^{4} \overline{\mu_{5}} - 4 \, \theta_{0}^{3} \overline{\mu_{5}} - 6 \, \theta_{0}^{2} \overline{\mu_{5}} + {\left(5 \, \theta_{0} - 12\right)} \overline{\zeta_{2}}}{2 \, {\left(\theta_{0}^{3} - 3 \, \theta_{0}^{2}\right)}} & A_{1} & \theta_{0} & 2 \\
            \frac{2 \, \theta_{0}^{3} \overline{\mu_{5}} + 2 \, \theta_{0}^{2} \overline{\mu_{5}} + 3 \, \overline{\zeta_{2}}}{2 \, {\left(\theta_{0}^{2} + \theta_{0}\right)}} & A_{0} & \theta_{0} - 1 & 2 \\
            \frac{\theta_{0}^{2} \overline{\mu_{7}} - \theta_{0} \overline{\zeta_{7}} + \overline{\zeta_{5}}}{\theta_{0}} & A_{0} & \theta_{0} - 1 & 3 \\
            \frac{\theta_{0}^{2} \overline{\mu_{1}} + 2 \, \theta_{0} \overline{\mu_{1}} + 2 \, \alpha_{5}}{\theta_{0}} & A_{0} & \theta_{0} + 1 & 0 \\
            -\frac{\theta_{0} - 1}{\theta_{0}^{2}} & \overline{B_{2}} & \theta_{0} + 1 & 1 \\
            \frac{\theta_{0}^{3} \overline{\mu_{3}} + 3 \, \theta_{0}^{2} \overline{\mu_{3}} + 2 \, \theta_{0} {\left(\overline{\alpha_{5}} + \overline{\mu_{3}}\right)} + 4 \, \overline{\alpha_{5}}}{\theta_{0}^{2} + \theta_{0}} & A_{1} & \theta_{0} + 1 & 1 \\
            -\frac{\theta_{0}^{4} \overline{\alpha_{9}} - {\left(2 \, \mu_{11} - \overline{\alpha_{9}}\right)} \theta_{0}^{3} - 2 \, {\left(\mu_{11} + 2 \, \overline{\alpha_{9}}\right)} \theta_{0}^{2} + 4 \, \theta_{0} \zeta_{0} \overline{\alpha_{1}} - 4 \, \theta_{0} \overline{\alpha_{9}} + {\left(\theta_{0}^{2} - \theta_{0} - 2\right)} \overline{\zeta_{4}}}{\theta_{0}^{3} + \theta_{0}^{2}} & \overline{A_{0}} & \theta_{0} + 1 & 1 \\
            \frac{4 \, \theta_{0}^{5} \overline{\mu_{4}} - 32 \, {\left(\theta_{0}^{3} + 2 \, \theta_{0}^{2}\right)} {\left| A_{1} \right|}^{4} + 12 \, \theta_{0}^{4} \overline{\mu_{4}} - 8 \, {\left(2 \, \alpha_{1} \overline{\alpha_{1}} - 2 \, \beta - \overline{\mu_{4}}\right)} \theta_{0}^{3} - 16 \, {\left(\alpha_{1} \overline{\alpha_{1}} - \beta\right)} \theta_{0}^{2} - {\left(\theta_{0}^{3} - 3 \, \theta_{0} - 2\right)} {\left| C_{1} \right|}^{2}}{4 \, {\left(\theta_{0}^{4} + \theta_{0}^{3}\right)}} & A_{0} & \theta_{0} + 1 & 1 \\
            \frac{\theta_{0}^{3} \overline{\mu_{1}} + 5 \, \theta_{0}^{2} \overline{\mu_{1}} + 2 \, {\left(\alpha_{5} + 3 \, \overline{\mu_{1}}\right)} \theta_{0} + 6 \, \alpha_{5}}{\theta_{0}^{2} + 2 \, \theta_{0}} & A_{1} & \theta_{0} + 2 & 0 \\
            \frac{\mu_{13} \theta_{0}^{3} - 4 \, {\left(\theta_{0} + 3\right)} \zeta_{0} {\left| A_{1} \right|}^{2} + 3 \, \mu_{13} \theta_{0}^{2} + 2 \, \mu_{13} \theta_{0}}{\theta_{0}^{3} + 3 \, \theta_{0}^{2} + 2 \, \theta_{0}} & \overline{A_{0}} & \theta_{0} + 2 & 0 \\
            \frac{\mu_{12} \theta_{0}^{3} + 3 \, \mu_{12} \theta_{0}^{2} + 2 \, \mu_{12} \theta_{0} + 4 \, {\left(2 \, \theta_{0} + 3\right)} \zeta_{1}}{\theta_{0}^{3} + 3 \, \theta_{0}^{2} + 2 \, \theta_{0}} & \overline{A_{1}} & \theta_{0} + 2 & 0 \\
            \frac{\theta_{0}^{4} \overline{\mu_{2}} + 6 \, \theta_{0}^{3} \overline{\mu_{2}} + {\left(3 \, \alpha_{6} + 11 \, \overline{\mu_{2}}\right)} \theta_{0}^{2} - 2 \, {\left(3 \, \alpha_{1} \theta_{0}^{2} + 13 \, \alpha_{1} \theta_{0} + 12 \, \alpha_{1}\right)} {\left| A_{1} \right|}^{2} + 3 \, {\left(3 \, \alpha_{6} + 2 \, \overline{\mu_{2}}\right)} \theta_{0} + 6 \, \alpha_{6}}{\theta_{0}^{3} + 3 \, \theta_{0}^{2} + 2 \, \theta_{0}} & A_{0} & \theta_{0} + 2 & 0 \\
            \frac{\theta_{0} - 4}{4 \, \theta_{0}} & \overline{E_{1}} & 2 \, \theta_{0} - 1 & -\theta_{0} + 3 \\
            \frac{\theta_{0}^{2} \overline{\alpha_{7}} - 8 \, \theta_{0} \overline{\alpha_{7}} + 16 \, \overline{\alpha_{7}}}{\theta_{0}} & A_{0} & 2 \, \theta_{0} - 1 & -\theta_{0} + 3 \\
            -\frac{\theta_{0}^{3} {\left(\overline{\alpha_{9}} - 2 \, \overline{\mu_{8}}\right)} + 2 \, \theta_{0}^{2} {\left(\overline{\alpha_{9}} - 3 \, \overline{\mu_{8}}\right)} - 4 \, \theta_{0} {\left(\overline{\alpha_{9}} + \overline{\mu_{8}}\right)} + {\left(\theta_{0} - 2\right)} \overline{\zeta_{4}} - 8 \, \overline{\alpha_{9}}}{\theta_{0}^{2} + 2 \, \theta_{0}} & A_{0} & 2 \, \theta_{0} + 1 & -\theta_{0} + 1 \\
            -\frac{2 \, \theta_{0}^{3} {\left(\overline{\alpha_{8}} - 2 \, \overline{\mu_{9}}\right)} - 2 \, \theta_{0}^{2} {\left(\overline{\alpha_{8}} + 3 \, \overline{\mu_{9}}\right)} - 2 \, \theta_{0} {\left(5 \, \overline{\alpha_{8}} + \overline{\mu_{9}}\right)} + {\left(\theta_{0} - 3\right)} \overline{\zeta_{3}} - 6 \, \overline{\alpha_{8}}}{2 \, {\left(\theta_{0}^{2} + \theta_{0}\right)}} & A_{0} & 2 \, \theta_{0} & -\theta_{0} + 2
            \end{dmatrix}        	
        \end{align*}
        \normalsize
        One can check that all relations of this order are trivial.

    	\chapter{Last order in the inversion}
    	\section{Development of $\phi$}
    	\small
    	\begin{align*}
    		&|\phi(z)|^2=\\
    		&\begin{dmatrix}
    		\frac{2 \, \overline{A_{0}} \overline{A_{1}}}{\theta_{0}^{2} + \theta_{0}} & 0 & 2 \, \theta_{0} + 1 &\textbf{(1)}\\
    		\frac{{\left(\theta_{0}^{2} + 2 \, \theta_{0}\right)} \overline{A_{1}}^{2} + 2 \, {\left(\theta_{0}^{2} + 2 \, \theta_{0} + 1\right)} \overline{A_{0}} \overline{A_{2}}}{\theta_{0}^{4} + 4 \, \theta_{0}^{3} + 5 \, \theta_{0}^{2} + 2 \, \theta_{0}} & 0 & 2 \, \theta_{0} + 2 &\textbf{(2)}\\
    		\frac{2 \, {\left({\left(\theta_{0}^{2} + 3 \, \theta_{0}\right)} \overline{A_{1}} \overline{A_{2}} + {\left(\theta_{0}^{2} + 3 \, \theta_{0} + 2\right)} \overline{A_{0}} \overline{A_{3}}\right)}}{\theta_{0}^{4} + 6 \, \theta_{0}^{3} + 11 \, \theta_{0}^{2} + 6 \, \theta_{0}} & 0 & 2 \, \theta_{0} + 3 &\textbf{(3)}\\
    		\lambda_1 & 0 & 2 \, \theta_{0} + 4 &\textbf{(4)}\\
    		\frac{\overline{A_{0}}^{2}}{\theta_{0}^{2}} & 0 & 2 \, \theta_{0} &\textbf{(5)}\\
    		\frac{B_{1} \overline{A_{0}} + C_{1} \overline{A_{1}}}{4 \, {\left(\theta_{0}^{2} + \theta_{0}\right)}} & 2 & 2 \, \theta_{0} + 1 &\textbf{(6)}\\
    		\frac{{\left(\theta_{0}^{2} + 2 \, \theta_{0} + 1\right)} B_{2} \overline{A_{0}} + {\left(\theta_{0}^{2} + 2 \, \theta_{0}\right)} B_{1} \overline{A_{1}} + {\left({\left(\theta_{0}^{2} \overline{\alpha_{1}} + 2 \, \theta_{0} \overline{\alpha_{1}} + \overline{\alpha_{1}}\right)} \overline{A_{0}} + {\left(\theta_{0}^{2} + 2 \, \theta_{0} + 1\right)} \overline{A_{2}}\right)} C_{1}}{4 \, {\left(\theta_{0}^{4} + 4 \, \theta_{0}^{3} + 5 \, \theta_{0}^{2} + 2 \, \theta_{0}\right)}} & 2 & 2 \, \theta_{0} + 2 &\textbf{(7)}\\
    		\frac{C_{1} \overline{A_{0}}}{4 \, \theta_{0}^{2}} & 2 & 2 \, \theta_{0} &\textbf{(8)}\\
    		\frac{2 \, C_{1} {\left| A_{1} \right|}^{2} \overline{A_{0}} + B_{3} \overline{A_{0}} + C_{2} \overline{A_{1}}}{6 \, {\left(\theta_{0}^{2} + \theta_{0}\right)}} & 3 & 2 \, \theta_{0} + 1 &\textbf{(9)}\\
    		\frac{C_{2} \overline{A_{0}}}{6 \, \theta_{0}^{2}} & 3 & 2 \, \theta_{0} &\textbf{(10)}\\
    		\frac{C_{1}^{2} {\left(\theta_{0} - 4\right)} + 8 \, {\left(\alpha_{1} \theta_{0} - 4 \, \alpha_{1}\right)} C_{1} \overline{A_{0}} + 8 \, C_{3} {\left(\theta_{0} - 4\right)} \overline{A_{0}} - 16 \, A_{0} E_{1}}{64 \, {\left(\theta_{0}^{3} - 4 \, \theta_{0}^{2}\right)}} & 4 & 2 \, \theta_{0} &\textbf{(11)}\\
    		-\frac{E_{1} \overline{A_{0}}}{4 \, {\left(\theta_{0}^{3} - 4 \, \theta_{0}^{2}\right)}} & -\theta_{0} + 4 & 3 \, \theta_{0} &\textbf{(12)}\\
    		\frac{2 \, A_{0} \overline{A_{0}}}{\theta_{0}^{2}} & \theta_{0} & \theta_{0} &\textbf{(13)}\\
    		\frac{2 \, A_{0} \overline{A_{1}}}{\theta_{0}^{2} + \theta_{0}} & \theta_{0} & \theta_{0} + 1 &\textbf{(14)}\\
    		\frac{8 \, A_{0} \theta_{0} \overline{A_{2}} + {\left(\theta_{0} + 2\right)} \overline{A_{0}} \overline{C_{1}}}{4 \, {\left(\theta_{0}^{3} + 2 \, \theta_{0}^{2}\right)}} & \theta_{0} & \theta_{0} + 2 &\textbf{(15)}\\
    		\frac{24 \, {\left(\theta_{0}^{2} + \theta_{0}\right)} A_{0} \overline{A_{3}} + 3 \, {\left(\theta_{0}^{2} + 3 \, \theta_{0}\right)} \overline{A_{1}} \overline{C_{1}} + 2 \, {\left(\theta_{0}^{2} + 4 \, \theta_{0} + 3\right)} \overline{A_{0}} \overline{C_{2}}}{12 \, {\left(\theta_{0}^{4} + 4 \, \theta_{0}^{3} + 3 \, \theta_{0}^{2}\right)}} & \theta_{0} & \theta_{0} + 3 &\textbf{(16)}\\
    		\lambda_2 & \theta_{0} & \theta_{0} + 4 &\textbf{(17)}\\
    		\frac{2 \, A_{1} \overline{A_{0}}}{\theta_{0}^{2} + \theta_{0}} & \theta_{0} + 1 & \theta_{0} &\textbf{(18)}
    		    		\end{dmatrix}
    		\end{align*}
    		\begin{align*}
    		\begin{dmatrix}
    		\frac{2 \, A_{1} \overline{A_{1}}}{\theta_{0}^{2} + 2 \, \theta_{0} + 1} & \theta_{0} + 1 & \theta_{0} + 1 &\textbf{(19)}\\
    		\frac{8 \, A_{1} \theta_{0} \overline{A_{2}} + {\left(\theta_{0} + 2\right)} \overline{A_{0}} \overline{B_{1}}}{4 \, {\left(\theta_{0}^{3} + 3 \, \theta_{0}^{2} + 2 \, \theta_{0}\right)}} & \theta_{0} + 1 & \theta_{0} + 2 &\textbf{(20)}\\
    		\frac{4 \, {\left(\theta_{0}^{2} + 4 \, \theta_{0} + 3\right)} {\left| A_{1} \right|}^{2} \overline{A_{0}} \overline{C_{1}} + 24 \, {\left(\theta_{0}^{2} + \theta_{0}\right)} A_{1} \overline{A_{3}} + 3 \, {\left(\theta_{0}^{2} + 3 \, \theta_{0}\right)} \overline{A_{1}} \overline{B_{1}} + 2 \, {\left(\theta_{0}^{2} + 4 \, \theta_{0} + 3\right)} \overline{A_{0}} \overline{B_{3}}}{12 \, {\left(\theta_{0}^{4} + 5 \, \theta_{0}^{3} + 7 \, \theta_{0}^{2} + 3 \, \theta_{0}\right)}} & \theta_{0} + 1 & \theta_{0} + 3 &\textbf{(21)}\\
    		\frac{A_{0} C_{1} {\left(\theta_{0} + 2\right)} + 8 \, A_{2} \theta_{0} \overline{A_{0}}}{4 \, {\left(\theta_{0}^{3} + 2 \, \theta_{0}^{2}\right)}} & \theta_{0} + 2 & \theta_{0} &\textbf{(22)}\\
    		\frac{A_{0} B_{1} {\left(\theta_{0} + 2\right)} + 8 \, A_{2} \theta_{0} \overline{A_{1}}}{4 \, {\left(\theta_{0}^{3} + 3 \, \theta_{0}^{2} + 2 \, \theta_{0}\right)}} & \theta_{0} + 2 & \theta_{0} + 1 &\textbf{(23)}\\
    		\lambda_3 & \theta_{0} + 2 & \theta_{0} + 2 &\textbf{(24)}\\
    		\frac{3 \, {\left(\theta_{0}^{2} + 3 \, \theta_{0}\right)} A_{1} C_{1} + 2 \, {\left(\theta_{0}^{2} + 4 \, \theta_{0} + 3\right)} A_{0} C_{2} + 24 \, {\left(\theta_{0}^{2} + \theta_{0}\right)} A_{3} \overline{A_{0}}}{12 \, {\left(\theta_{0}^{4} + 4 \, \theta_{0}^{3} + 3 \, \theta_{0}^{2}\right)}} & \theta_{0} + 3 & \theta_{0} &\textbf{(25)}\\
    		\frac{4 \, {\left(\theta_{0}^{2} + 4 \, \theta_{0} + 3\right)} A_{0} C_{1} {\left| A_{1} \right|}^{2} + 3 \, {\left(\theta_{0}^{2} + 3 \, \theta_{0}\right)} A_{1} B_{1} + 2 \, {\left(\theta_{0}^{2} + 4 \, \theta_{0} + 3\right)} A_{0} B_{3} + 24 \, {\left(\theta_{0}^{2} + \theta_{0}\right)} A_{3} \overline{A_{1}}}{12 \, {\left(\theta_{0}^{4} + 5 \, \theta_{0}^{3} + 7 \, \theta_{0}^{2} + 3 \, \theta_{0}\right)}} & \theta_{0} + 3 & \theta_{0} + 1 &\textbf{(26)}\\
    		\lambda_4 & \theta_{0} + 4 & \theta_{0} &\textbf{(27)}\\
    		\frac{2 \, A_{0} A_{1}}{\theta_{0}^{2} + \theta_{0}} & 2 \, \theta_{0} + 1 & 0 &\textbf{(28)}\\
    		\frac{A_{0} \overline{B_{1}} + A_{1} \overline{C_{1}}}{4 \, {\left(\theta_{0}^{2} + \theta_{0}\right)}} & 2 \, \theta_{0} + 1 & 2 &\textbf{(29)}\\
    		\frac{2 \, A_{0} {\left| A_{1} \right|}^{2} \overline{C_{1}} + A_{0} \overline{B_{3}} + A_{1} \overline{C_{2}}}{6 \, {\left(\theta_{0}^{2} + \theta_{0}\right)}} & 2 \, \theta_{0} + 1 & 3 &\textbf{(30)}\\
    		\frac{{\left(\theta_{0}^{2} + 2 \, \theta_{0}\right)} A_{1}^{2} + 2 \, {\left(\theta_{0}^{2} + 2 \, \theta_{0} + 1\right)} A_{0} A_{2}}{\theta_{0}^{4} + 4 \, \theta_{0}^{3} + 5 \, \theta_{0}^{2} + 2 \, \theta_{0}} & 2 \, \theta_{0} + 2 & 0 &\textbf{(31)}\\
    		\frac{{\left(\theta_{0}^{2} + 2 \, \theta_{0}\right)} A_{1} \overline{B_{1}} + {\left(\theta_{0}^{2} + 2 \, \theta_{0} + 1\right)} A_{0} \overline{B_{2}} + {\left({\left(\alpha_{1} \theta_{0}^{2} + 2 \, \alpha_{1} \theta_{0} + \alpha_{1}\right)} A_{0} + {\left(\theta_{0}^{2} + 2 \, \theta_{0} + 1\right)} A_{2}\right)} \overline{C_{1}}}{4 \, {\left(\theta_{0}^{4} + 4 \, \theta_{0}^{3} + 5 \, \theta_{0}^{2} + 2 \, \theta_{0}\right)}} & 2 \, \theta_{0} + 2 & 2 &\textbf{(32)}\\
    		\frac{2 \, {\left({\left(\theta_{0}^{2} + 3 \, \theta_{0}\right)} A_{1} A_{2} + {\left(\theta_{0}^{2} + 3 \, \theta_{0} + 2\right)} A_{0} A_{3}\right)}}{\theta_{0}^{4} + 6 \, \theta_{0}^{3} + 11 \, \theta_{0}^{2} + 6 \, \theta_{0}} & 2 \, \theta_{0} + 3 & 0 &\textbf{(33)}\\
    		\lambda_5 & 2 \, \theta_{0} + 4 & 0 &\textbf{(34)}\\
    		\frac{A_{0}^{2}}{\theta_{0}^{2}} & 2 \, \theta_{0} & 0 &\textbf{(35)}\\
    		\frac{A_{0} \overline{C_{1}}}{4 \, \theta_{0}^{2}} & 2 \, \theta_{0} & 2 &\textbf{(36)}\\
    		\frac{A_{0} \overline{C_{2}}}{6 \, \theta_{0}^{2}} & 2 \, \theta_{0} & 3 &\textbf{(37)}\\
    		\frac{8 \, {\left(\theta_{0} \overline{\alpha_{1}} - 4 \, \overline{\alpha_{1}}\right)} A_{0} \overline{C_{1}} + {\left(\theta_{0} - 4\right)} \overline{C_{1}}^{2} + 8 \, A_{0} {\left(\theta_{0} - 4\right)} \overline{C_{3}} - 16 \, \overline{A_{0}} \overline{E_{1}}}{64 \, {\left(\theta_{0}^{3} - 4 \, \theta_{0}^{2}\right)}} & 2 \, \theta_{0} & 4 &\textbf{(38)}\\
    		-\frac{A_{0} \overline{E_{1}}}{4 \, {\left(\theta_{0}^{3} - 4 \, \theta_{0}^{2}\right)}} & 3 \, \theta_{0} & -\theta_{0} + 4&\textbf{(39)}
    		\end{dmatrix}
    	\end{align*}
    	where
    	\begin{align}\label{lambda}
    	    \lambda_1&=\frac{{\left(\theta_{0}^{4} + 8 \, \theta_{0}^{3} + 19 \, \theta_{0}^{2} + 12 \, \theta_{0}\right)} \overline{A_{2}}^{2} + 2 \, {\left(\theta_{0}^{4} + 8 \, \theta_{0}^{3} + 20 \, \theta_{0}^{2} + 16 \, \theta_{0}\right)} \overline{A_{1}} \overline{A_{3}} + 2 \, {\left(\theta_{0}^{4} + 8 \, \theta_{0}^{3} + 23 \, \theta_{0}^{2} + 28 \, \theta_{0} + 12\right)} \overline{A_{0}} \overline{A_{4}}}{\theta_{0}^{6} + 12 \, \theta_{0}^{5} + 55 \, \theta_{0}^{4} + 120 \, \theta_{0}^{3} + 124 \, \theta_{0}^{2} + 48 \, \theta_{0}}\\
    		\lambda_2&=\frac{1}{24 \, {\left(\theta_{0}^{5} + 7 \, \theta_{0}^{4} + 14 \, \theta_{0}^{3} + 8 \, \theta_{0}^{2}\right)}}\bigg\{48 \, {\left(\theta_{0}^{3} + 3 \, \theta_{0}^{2} + 2 \, \theta_{0}\right)} A_{0} \overline{A_{4}} + 4 \, {\left(\theta_{0}^{3} + 6 \, \theta_{0}^{2} + 8 \, \theta_{0}\right)} \overline{A_{1}} \overline{C_{2}} + 3 \, {\left(\theta_{0}^{3} + 7 \, \theta_{0}^{2} + 14 \, \theta_{0} + 8\right)} \overline{A_{0}} \overline{C_{3}}\\
    		& + 3 \, {\left({\left(\theta_{0}^{3} \overline{\alpha_{1}} + 7 \, \theta_{0}^{2} \overline{\alpha_{1}} + 14 \, \theta_{0} \overline{\alpha_{1}} + 8 \, \overline{\alpha_{1}}\right)} \overline{A_{0}} + 2 \, {\left(\theta_{0}^{3} + 5 \, \theta_{0}^{2} + 4 \, \theta_{0}\right)} \overline{A_{2}}\right)} \overline{C_{1}}\bigg\}\\
    		\lambda_3&=\frac{1}{32 \, {\left(\theta_{0}^{4} + 4 \, \theta_{0}^{3} + 4 \, \theta_{0}^{2}\right)}}\bigg\{64 \, A_{2} \theta_{0}^{2} \overline{A_{2}} + 8 \, {\left(\theta_{0}^{2} + 2 \, \theta_{0}\right)} A_{0} B_{2} + 8 \, {\left(\theta_{0}^{2} \overline{\alpha_{1}} + 2 \, \theta_{0} \overline{\alpha_{1}}\right)} A_{0} C_{1} + 8 \, {\left(\theta_{0}^{2} + 2 \, \theta_{0}\right)} \overline{A_{0}} \overline{B_{2}}\\
    		& + {\left({\left(\theta_{0}^{2} + 4 \, \theta_{0} + 4\right)} C_{1} + 8 \, {\left(\alpha_{1} \theta_{0}^{2} + 2 \, \alpha_{1} \theta_{0}\right)} \overline{A_{0}}\right)} \overline{C_{1}}\bigg\}\\
    		\lambda_4&=\frac{1}{24 \, {\left(\theta_{0}^{5} + 7 \, \theta_{0}^{4} + 14 \, \theta_{0}^{3} + 8 \, \theta_{0}^{2}\right)}}\bigg\{4 \, {\left(\theta_{0}^{3} + 6 \, \theta_{0}^{2} + 8 \, \theta_{0}\right)} A_{1} C_{2} + 3 \, {\left(\theta_{0}^{3} + 7 \, \theta_{0}^{2} + 14 \, \theta_{0} + 8\right)} A_{0} C_{3} + 48 \, {\left(\theta_{0}^{3} + 3 \, \theta_{0}^{2} + 2 \, \theta_{0}\right)} A_{4} \overline{A_{0}}\\
    		& + 3 \, {\left({\left(\alpha_{1} \theta_{0}^{3} + 7 \, \alpha_{1} \theta_{0}^{2} + 14 \, \alpha_{1} \theta_{0} + 8 \, \alpha_{1}\right)} A_{0} + 2 \, {\left(\theta_{0}^{3} + 5 \, \theta_{0}^{2} + 4 \, \theta_{0}\right)} A_{2}\right)} C_{1}\bigg\}\\
    		\lambda_5&=\frac{{\left(\theta_{0}^{4} + 8 \, \theta_{0}^{3} + 19 \, \theta_{0}^{2} + 12 \, \theta_{0}\right)} A_{2}^{2} + 2 \, {\left(\theta_{0}^{4} + 8 \, \theta_{0}^{3} + 20 \, \theta_{0}^{2} + 16 \, \theta_{0}\right)} A_{1} A_{3} + 2 \, {\left(\theta_{0}^{4} + 8 \, \theta_{0}^{3} + 23 \, \theta_{0}^{2} + 28 \, \theta_{0} + 12\right)} A_{0} A_{4}}{\theta_{0}^{6} + 12 \, \theta_{0}^{5} + 55 \, \theta_{0}^{4} + 120 \, \theta_{0}^{3} + 124 \, \theta_{0}^{2} + 48 \, \theta_{0}}
    	\end{align}
    	Recall by REF that
    		\begin{align*}
    	&|\phi(z)|^2=\begin{dmatrix}
    	\frac{-\bar{\zeta_0}}{\theta_{0}^{4} + 4 \, \theta_{0}^{3} + 5 \, \theta_{0}^{2} + 2 \, \theta_{0}} & 0 & 2 \, \theta_{0} + 2 &\textbf{(1)}\\
    	\frac{-4\bar{\zeta_1}}{\theta_{0}^{4} + 6 \, \theta_{0}^{3} + 11 \, \theta_{0}^{2} + 6 \, \theta_{0}} & 0 & 2 \, \theta_{0} + 3 &\textbf{(2)}\\
    	\frac{1}{\theta_{0}^{2}} & \theta_{0} & \theta_{0} &\textbf{(4)}\\
    	\frac{\bar{\alpha_1}}{{\theta_{0}(\theta_0 + 2)}} & \theta_{0} & \theta_{0} + 2 &\textbf{(5)}\\
    	\frac{-\bar{\zeta_2}+2\theta_0(\theta_0+1)\bar{\alpha_3}}{2{\left(\theta_{0}^{4} + 4 \, \theta_{0}^{3} + 3 \, \theta_{0}^{2}\right)}} & \theta_{0} & \theta_{0} + 3 &\textbf{(6)}\\
    	\frac{2 \, |A_{1}|^2}{(\theta_{0}+1)^2} & \theta_{0} + 1 & \theta_{0} + 1 &\textbf{(7)}\\
    	\frac{\bar{\alpha_5}}{(\theta_0+1)(\theta_0+2)} & \theta_{0} + 1 & \theta_{0} + 2 &\textbf{(8)}\\
    	\frac{\alpha_1}{\theta_0(\theta_0+2)} & \theta_{0} + 2 & \theta_{0} &\textbf{(9)}\\
    	\frac{\alpha_5}{(\theta_0+1)(\theta_0+2)} & \theta_{0} + 2 & \theta_{0} + 1 &\textbf{(10)}\\
    	\frac{-\zeta_2+2\theta_0(\theta_0+1)\alpha_3}{2{\left(\theta_{0}^{4} + 4 \, \theta_{0}^{3} + 3 \, \theta_{0}^{2}\right)}} & \theta_{0} + 3 & \theta_{0} &\textbf{(11)}\\
    	\frac{-\zeta_0}{\theta_{0}^{4} + 4 \, \theta_{0}^{3} + 5 \, \theta_{0}^{2} + 2 \, \theta_{0}} & 2 \, \theta_{0} + 2 & 0 &\textbf{(13)}\\
    	\frac{-4\zeta_1}{\theta_{0}^{4} + 6 \, \theta_{0}^{3} + 11 \, \theta_{0}^{2} + 6 \, \theta_{0}} & 2 \, \theta_{0} + 3 & 0 &\textbf{(14)}
    	\end{dmatrix}\\
    	&=\frac{1}{\theta_0^2}|z|^{2\theta_0}+\frac{2|\vec{A}_1|^2}{(\theta_0+1)^2}|z|^{2\theta_0+2}+2\,\Re\bigg(\frac{\,\alpha_1}{\theta_0(\theta_0+2)}z^{\theta_0+2}\z^{\theta_0}+\frac{\alpha_5}{(\theta_0+1)(\theta_0+2)}z^{\theta_0+2}\z^{\theta_0+1}+\\
    	&+\frac{-\zeta_2+2\theta_0(\theta_0+1)\alpha_3}{2(\theta_0^4+4\theta_0^3+3\theta_0^2)}z^{\theta_0+3}\z^{\theta_0}-\frac{\,\zeta_0}{\theta_0^4+4\theta_0^3+5\theta_0^2+2\theta_0}z^{2\theta_0+2}-\frac{4\zeta_1}{\theta_0^4+6\theta_0^3+11\theta_0^2+6\theta_0}z^{2\theta_0+3}\bigg)+O(|z|^{2\theta_0+4})\\
    	&=|z|^{2\theta_0}\bigg\{\frac{1}{\theta_0^2}+\frac{2|\vec{A}_1|^2}{(\theta_0+1)^2}|z|^2+2\,\Re\bigg(\frac{\,\alpha_1}{\theta_0(\theta_0+2)}z^{2}+\frac{\alpha_5}{(\theta_0+1)(\theta_0+2)}z^{2}\z+
    	+\frac{-\zeta_2+2\theta_0(\theta_0+1)\alpha_3}{2(\theta_0^4+4\theta_0^3+3\theta_0^2)}z^{3}\\
    	&-\frac{\,\zeta_0}{\theta_0^4+4\theta_0^3+5\theta_0^2+2\theta_0}z^{\theta_0+2}\z^{-\theta_0}-\frac{4\zeta_1}{\theta_0^4+6\theta_0^3+11\theta_0^2+6\theta_0}z^{\theta_0+3}\z^{-\theta_0} \bigg)+O(|z|^4) \bigg\}
    	\end{align*}
    	where
    	\begin{align*}
    	\left\{\begin{alignedat}{1}
    	\zeta_0&=\s{\vec{A}_1}{\vec{A}_1}\\
    	\zeta_1&=\s{\vec{A}_1}{\vec{A}_2}\\
    	\zeta_2&=\s{\vec{A}_1}{\vec{C}_1}
    	\end{alignedat}\right.
    	\end{align*}
    	so we need only consider powers of order equal to $2\theta_0+4$, which correspond to the lines
    	\begin{align*}
    		\textbf{(4)}, \textbf{(7)}, \textbf{(9)}, \textbf{(11)}, \textbf{(12)}, \textbf{(17)}, \textbf{(21)}, \textbf{(24)}, \textbf{(26)}, \textbf{(27)}, \textbf{(30)}, \textbf{(32)}, \textbf{(34)}, \textbf{(38)}, \textbf{(39)}. 
    	\end{align*}
    	and this is thanks of
    				\begin{align}\label{invref1}
    	\left\{
    	\begin{alignedat}{1}
    	\vec{B}_1&=-2\s{\bar{\vec{A}_1}}{\vec{C}_1}\vec{A}_0\\
    	\vec{B}_2&=-\frac{(\theta_0+2)}{4\theta_0}|\vec{C}_1|^2\bar{\vec{A}_0}-2\s{\bar{\vec{A}_2}}{\vec{C}_1}\vec{A}_0\\
    	\vec{B}_3&=-2\s{\bar{\vec{A}_1}}{\vec{C}_1}\vec{A}_1+\frac{2}{\theta_0-3}\s{\vec{A}_1}{\vec{C}_1}\bar{\vec{A}_1}-2\s{\bar{\vec{A}_1}}{\vec{C}_2}\vec{A}_0\\
    	\vec{E}_1&=-\frac{1}{2\theta_0}\s{\vec{C}_1}{\vec{C}_1}\bar{\vec{A}_0}.
    	\end{alignedat}\right.
    	\end{align}
    	some words
    	\begin{align*}
    	&\begin{dmatrix}
    		\lambda_1 & 0 & 2 \, \theta_{0} + 4 &\textbf{(4)}\\
    		\frac{\colorcancel{{\left(\theta_{0}^{2} + 2 \, \theta_{0} + 1\right)} B_{2} \overline{A_{0}}}{blue} + \ccancel{{\left(\theta_{0}^{2} + 2 \, \theta_{0}\right)} B_{1} \overline{A_{1}}} + {\left(\ccancel{{\left(\theta_{0}^{2} \overline{\alpha_{1}} + 2 \, \theta_{0} \overline{\alpha_{1}} + \overline{\alpha_{1}}\right)} \overline{A_{0}}} + \colorcancel{{\left(\theta_{0}^{2} + 2 \, \theta_{0} + 1\right)} \overline{A_{2}}}{blue}\right)} C_{1}}{4 \, {\left(\theta_{0}^{4} + 4 \, \theta_{0}^{3} + 5 \, \theta_{0}^{2} + 2 \, \theta_{0}\right)}} & 2 & 2 \, \theta_{0} + 2 &\textbf{(7)}\\
    		\frac{\ccancel{2 \, C_{1} {\left| A_{1} \right|}^{2} \overline{A_{0}}} + \colorcancel{B_{3} \overline{A_{0}}}{blue} + \colorcancel{C_{2} \overline{A_{1}}}{blue}}{6 \, {\left(\theta_{0}^{2} + \theta_{0}\right)}} & 3 & 2 \, \theta_{0} + 1 &\textbf{(9)}\\
    		\frac{\colorcancel{C_{1}^{2} {\left(\theta_{0} - 4\right)}}{blue} + \ccancel{8 \, {\left(\alpha_{1} \theta_{0} - 4 \, \alpha_{1}\right)} C_{1} \overline{A_{0}}} + \colorcancel{8 \, C_{3} {\left(\theta_{0} - 4\right)} \overline{A_{0}}}{blue} - 16 \, A_{0} E_{1}}{64 \, {\left(\theta_{0}^{3} - 4 \, \theta_{0}^{2}\right)}} & 4 & 2 \, \theta_{0} &\textbf{(11)}\\
    		-\frac{\ccancel{E_{1} \overline{A_{0}}}}{4 \, {\left(\theta_{0}^{3} - 4 \, \theta_{0}^{2}\right)}} & -\theta_{0} + 4 & 3 \, \theta_{0} &\textbf{(12)}\\
    		\lambda_2 & \theta_{0} & \theta_{0} + 4 &\textbf{(17)}\\
    		\frac{\ccancel{4 \, {\left(\theta_{0}^{2} + 4 \, \theta_{0} + 3\right)} {\left| A_{1} \right|}^{2} \overline{A_{0}} \overline{C_{1}}} + 24 \, {\left(\theta_{0}^{2} + \theta_{0}\right)} A_{1} \overline{A_{3}} + \ccancel{3 \, {\left(\theta_{0}^{2} + 3 \, \theta_{0}\right)} \overline{A_{1}} \overline{B_{1}}} + \ccancel{2 \, {\left(\theta_{0}^{2} + 4 \, \theta_{0} + 3\right)} \overline{A_{0}} \overline{B_{3}}}}{12 \, {\left(\theta_{0}^{4} + 5 \, \theta_{0}^{3} + 7 \, \theta_{0}^{2} + 3 \, \theta_{0}\right)}} & \theta_{0} + 1 & \theta_{0} + 3 &\textbf{(21)}\\
    		\lambda_3 & \theta_{0} + 2 & \theta_{0} + 2 &\textbf{(24)}\\
    		\frac{\ccancel{4 \, {\left(\theta_{0}^{2} + 4 \, \theta_{0} + 3\right)} A_{0} C_{1} {\left| A_{1} \right|}^{2}} + \ccancel{3 \, {\left(\theta_{0}^{2} + 3 \, \theta_{0}\right)} A_{1} B_{1}} + \ccancel{2 \, {\left(\theta_{0}^{2} + 4 \, \theta_{0} + 3\right)} A_{0} B_{3}} + 24 \, {\left(\theta_{0}^{2} + \theta_{0}\right)} A_{3} \overline{A_{1}}}{12 \, {\left(\theta_{0}^{4} + 5 \, \theta_{0}^{3} + 7 \, \theta_{0}^{2} + 3 \, \theta_{0}\right)}} & \theta_{0} + 3 & \theta_{0} + 1 &\textbf{(26)}\\
    		\lambda_4 & \theta_{0} + 4 & \theta_{0} &\textbf{(27)}\\
    		\frac{\ccancel{2 \, A_{0} {\left| A_{1} \right|}^{2} \overline{C_{1}}} + \colorcancel{A_{0} \overline{B_{3}}}{blue} + \colorcancel{A_{1} \overline{C_{2}}}{blue}}{6 \, {\left(\theta_{0}^{2} + \theta_{0}\right)}} & 2 \, \theta_{0} + 1 & 3 &\textbf{(30)}\\
    		\frac{\ccancel{{\left(\theta_{0}^{2} + 2 \, \theta_{0}\right)} A_{1} \overline{B_{1}}} + \colorcancel{{\left(\theta_{0}^{2} + 2 \, \theta_{0} + 1\right)} A_{0} \overline{B_{2}}}{blue} + {\left(\ccancel{{\left(\alpha_{1} \theta_{0}^{2} + 2 \, \alpha_{1} \theta_{0} + \alpha_{1}\right)} A_{0}} + \colorcancel{{\left(\theta_{0}^{2} + 2 \, \theta_{0} + 1\right)} A_{2}}{blue}\right)} \overline{C_{1}}}{4 \, {\left(\theta_{0}^{4} + 4 \, \theta_{0}^{3} + 5 \, \theta_{0}^{2} + 2 \, \theta_{0}\right)}} & 2 \, \theta_{0} + 2 & 2 &\textbf{(32)}\\
    		\lambda_5 & 2 \, \theta_{0} + 4 & 0 &\textbf{(34)}\\
    		\frac{\ccancel{8 \, {\left(\theta_{0} \overline{\alpha_{1}} - 4 \, \overline{\alpha_{1}}\right)} A_{0} \overline{C_{1}}} + \colorcancel{{\left(\theta_{0} - 4\right)} \overline{C_{1}}^{2}}{blue} + \colorcancel{8 \, A_{0} {\left(\theta_{0} - 4\right)} \overline{C_{3}}}{blue} - 16 \, \overline{A_{0}} \overline{E_{1}}}{64 \, {\left(\theta_{0}^{3} - 4 \, \theta_{0}^{2}\right)}} & 2 \, \theta_{0} & 4 &\textbf{(38)}\\
    		-\frac{\ccancel{A_{0} \overline{E_{1}}}}{4 \, {\left(\theta_{0}^{3} - 4 \, \theta_{0}^{2}\right)}} & 3 \, \theta_{0} & -\theta_{0} + 4&\textbf{(39)} 
    		\end{dmatrix}\\
    		&=\begin{dmatrix}
    		\lambda_1 & 0 & 2 \, \theta_{0} + 4 &\textbf{(4)}\\
    		\frac{- 16 \, A_{0} E_{1}}{64 \, {\left(\theta_{0}^{3} - 4 \, \theta_{0}^{2}\right)}} & 4 & 2 \, \theta_{0} &\textbf{(11)}\\
    		\lambda_2 & \theta_{0} & \theta_{0} + 4 &\textbf{(17)}\\
    		\frac{24 \, {\left(\theta_{0}^{2} + \theta_{0}\right)} A_{1} \overline{A_{3}}}{12 \, {\left(\theta_{0}^{4} + 5 \, \theta_{0}^{3} + 7 \, \theta_{0}^{2} + 3 \, \theta_{0}\right)}} & \theta_{0} + 1 & \theta_{0} + 3 &\textbf{(21)}\\
    		\lambda_3 & \theta_{0} + 2 & \theta_{0} + 2 &\textbf{(24)}\\
    		\frac{24 \, {\left(\theta_{0}^{2} + \theta_{0}\right)} A_{3} \overline{A_{1}}}{12 \, {\left(\theta_{0}^{4} + 5 \, \theta_{0}^{3} + 7 \, \theta_{0}^{2} + 3 \, \theta_{0}\right)}} & \theta_{0} + 3 & \theta_{0} + 1 &\textbf{(26)}\\
    		\lambda_4 & \theta_{0} + 4 & \theta_{0} &\textbf{(27)}\\
    		\lambda_5 & 2 \, \theta_{0} + 4 & 0 &\textbf{(34)}\\
    		\frac{- 16 \, \overline{A_{0}} \overline{E_{1}}}{64 \, {\left(\theta_{0}^{3} - 4 \, \theta_{0}^{2}\right)}} & 2 \, \theta_{0} & 4 &\textbf{(38)}
    		\end{dmatrix}
    	\end{align*}
    	\normalsize
    	$\lambda_1$ does not simplify much and as $\phi$ is real, we have $\lambda_5=\bar{\lambda_1}$. Then, we also have $\lambda_4=\bar{\lambda_2}$, and finally, 
    	\begin{align*}
    		\lambda_3&=\frac{1}{32 \, {\left(\theta_{0}^{4} + 4 \, \theta_{0}^{3} + 4 \, \theta_{0}^{2}\right)}}\bigg\{64 \, A_{2} \theta_{0}^{2} \overline{A_{2}} + 8 \, {\left(\theta_{0}^{2} + 2 \, \theta_{0}\right)} A_{0} B_{2} + \ccancel{8 \, {\left(\theta_{0}^{2} \overline{\alpha_{1}} + 2 \, \theta_{0} \overline{\alpha_{1}}\right)} A_{0} C_{1}} + 8 \, {\left(\theta_{0}^{2} + 2 \, \theta_{0}\right)} \overline{A_{0}} \overline{B_{2}}\\
    		& + {\left({\left(\theta_{0}^{2} + 4 \, \theta_{0} + 4\right)} C_{1} + \ccancel{8 \, {\left(\alpha_{1} \theta_{0}^{2} + 2 \, \alpha_{1} \theta_{0}\right)} \overline{A_{0}}}\right)} \overline{C_{1}}\bigg\}\\
    		&=\frac{1}{32\theta_0^2(\theta_0+1)^2}\bigg\{64\,\theta_0^2|\vec{A}_2|^2+(\theta_0+2)^2|\vec{C}_1|^2+16\theta_0(\theta_0+2)\Re\left(\s{\vec{A}_0}{\vec{B}_2}\right) \bigg\}
    	\end{align*}
    	Now, we have
    	\begin{align*}
    		\s{\vec{A}_0}{\vec{B}_2}&=\bs{\vec{A}_0}{-\frac{(\theta_0+2)}{4\theta_0}|\vec{C}_1|^2\bar{\vec{A}_0}-2\s{\bar{\vec{A}_2}}{\vec{C}_1}\vec{A}_0}=-\frac{(\theta_0+2)}{8\theta_0}|\vec{C}_1|^2\in\R
    	\end{align*}
    	so 
    	\begin{align*}
    		16\theta_0(\theta_0+2)\Re\left(\s{\vec{A}_0}{\vec{B}_2}\right)=16\theta_0(\theta_0+2)\s{\vec{A}_0}{\vec{B}_2}=16\theta_0(\theta_0+2)\left(-\frac{(\theta_0+2)}{8\theta_0}|\vec{C}_1|^2\right)=-2(\theta_0+2)^2|\vec{C}_1|^2
    	\end{align*}
    	and finally
    	\begin{align*}
    		\lambda_3=\frac{1}{32\theta_0^2(\theta_0+2)^2}\left(64\theta_0^2|\vec{A}_2|^2-(\theta_0+2)^2|\vec{C}_1|^2\right)
    	\end{align*}
    	Now, recall that
    	\begin{align*}
    	&\alpha_7=\frac{1}{8\theta_0(\theta_0-4)}\s{\vec{C}_1}{\vec{C}_1}\\
    	&\vec{E}_1=-\frac{1}{2\theta_0}\s{\vec{C}_1}{\vec{C}_1}\bar{\vec{A}_0}=-4(\theta_0-4)\alpha_7\vec{A}_0
    	\end{align*}
    	so as $|\vec{A}_0|^2=\dfrac{1}{2}$
    	\begin{align*}
    		\s{\vec{A}_0}{\vec{E}_1}=-2(\theta_0-4)\alpha_7
    	\end{align*}
    	and
    	\begin{align*}
    		\textbf{(11)}=\frac{- 16 \, A_{0} E_{1}}{64 \, {\left(\theta_{0}^{3} - 4 \, \theta_{0}^{2}\right)}}=-\frac{1}{4\theta_0^2(\theta_0-4)}\left(-2(\theta_0-4)\alpha_7\right)=\frac{\alpha_7}{2\theta_0^2}.
    	\end{align*}
    	To keep consistent notations, we define
    	\begin{align*}
    		\left\{\begin{alignedat}{1}
    		\zeta_8&=\lambda_4\\
    		\zeta_9&=\lambda_5\\
    		\zeta_{10}&=\lambda_3=\frac{1}{32\theta_0^2(\theta_0+2)^2}\left(64\,\theta_0^2|\vec{A}_2|^2-(\theta_0+2)^2\vec{C}_1|^2\right)\\
    		\zeta_{11}&=\s{\bar{\vec{A}_1}}{\vec{A}_3}
    		\end{alignedat}\right.
       	\end{align*}
    	where $\lambda_4,\lambda_5$ are given by \eqref{lambda},
    	so that
    	\begin{align*}
    		\begin{dmatrix}
    		\lambda_1 & 0 & 2 \, \theta_{0} + 4 &\textbf{(4)}\\
    		\frac{- 16 \, A_{0} E_{1}}{64 \, {\left(\theta_{0}^{3} - 4 \, \theta_{0}^{2}\right)}} & 4 & 2 \, \theta_{0} &\textbf{(11)}\\
    		\lambda_2 & \theta_{0} & \theta_{0} + 4 &\textbf{(17)}\\
    		\frac{24 \, {\left(\theta_{0}^{2} + \theta_{0}\right)} A_{1} \overline{A_{3}}}{12 \, {\left(\theta_{0}^{4} + 5 \, \theta_{0}^{3} + 7 \, \theta_{0}^{2} + 3 \, \theta_{0}\right)}} & \theta_{0} + 1 & \theta_{0} + 3 &\textbf{(21)}\\
    		\lambda_3 & \theta_{0} + 2 & \theta_{0} + 2 &\textbf{(24)}\\
    		\frac{24 \, {\left(\theta_{0}^{2} + \theta_{0}\right)} A_{3} \overline{A_{1}}}{12 \, {\left(\theta_{0}^{4} + 5 \, \theta_{0}^{3} + 7 \, \theta_{0}^{2} + 3 \, \theta_{0}\right)}} & \theta_{0} + 3 & \theta_{0} + 1 &\textbf{(26)}\\
    		\lambda_4 & \theta_{0} + 4 & \theta_{0} &\textbf{(27)}\\
    		\lambda_5 & 2 \, \theta_{0} + 4 & 0 &\textbf{(34)}\\
    		\frac{- 16 \, \overline{A_{0}} \overline{E_{1}}}{64 \, {\left(\theta_{0}^{3} - 4 \, \theta_{0}^{2}\right)}} & 2 \, \theta_{0} & 4 &\textbf{(38)}
    		\end{dmatrix}
    		=\begin{dmatrix}
    		\bar{\zeta_9} & 0 & 2 \, \theta_{0} + 4 &\textbf{(4)}\\
    		\frac{\alpha_7}{2\theta_0^2} & 4 & 2 \, \theta_{0} &\textbf{(11)}\\
    		\bar{\zeta_{8}} & \theta_{0} & \theta_{0} + 4 &\textbf{(17)}\\
    		\frac{2\,\bar{\zeta_{11}}}{(\theta_0+1)(\theta_0+3)} & \theta_{0} + 1 & \theta_{0} + 3 &\textbf{(21)}\\
    		\zeta_{10} & \theta_{0} + 2 & \theta_{0} + 2 &\textbf{(24)}\\
    		\frac{2\zeta_{11}}{(\theta_0+1)(\theta_0+3)} & \theta_{0} + 3 & \theta_{0} + 1 &\textbf{(26)}\\
    		\zeta_{8} & \theta_{0} + 4 & \theta_{0} &\textbf{(27)}\\
    		\zeta_{9} & 2 \, \theta_{0} + 4 & 0 &\textbf{(34)}\\
    		\frac{\bar{\alpha_7}}{2\theta_0^2} & 2 \, \theta_{0} & 4 &\textbf{(38)}
    		\end{dmatrix}
    	\end{align*}
    	as
    	\begin{align*}
    		\, {\left(\theta_{0}^{4} + 5 \, \theta_{0}^{3} + 7 \, \theta_{0}^{2} + 3 \, \theta_{0}\right)}=\theta_0(\theta_0+1)^2(\theta_0+3)
    	\end{align*}
    	so that
    	\begin{align*}
    		\frac{24 \, {\left(\theta_{0}^{2} + \theta_{0}\right)} A_{3} \overline{A_{1}}}{12 \, {\left(\theta_{0}^{4} + 5 \, \theta_{0}^{3} + 7 \, \theta_{0}^{2} + 3 \, \theta_{0}\right)}}=\frac{2\theta_0(\theta_0+1)\,\zeta_{11}}{\theta_0(\theta_0+1)^2(\theta_0+3)}=\frac{2\zeta_{11}}{(\theta_0+1)(\theta_0+3)}
    	\end{align*}
    	Now we compare the Sage transcription on the right to the expression we computed on the left
    	\begin{align*}
    		\begin{dmatrix}
    		\bar{\zeta_9} & 0 & 2 \, \theta_{0} + 4 &\textbf{(4)}\\
    		\frac{\alpha_7}{2\theta_0^2} & 4 & 2 \, \theta_{0} &\textbf{(11)}\\
    		\bar{\zeta_{8}} & \theta_{0} & \theta_{0} + 4 &\textbf{(17)}\\
    		\frac{2\,\bar{\zeta_{11}}}{(\theta_0+1)(\theta_0+3)} & \theta_{0} + 1 & \theta_{0} + 3 &\textbf{(21)}\\
    		\zeta_{10} & \theta_{0} + 2 & \theta_{0} + 2 &\textbf{(24)}\\
    		\frac{2\zeta_{11}}{(\theta_0+1)(\theta_0+3)} & \theta_{0} + 3 & \theta_{0} + 1 &\textbf{(26)}\\
    		\zeta_{8} & \theta_{0} + 4 & \theta_{0} &\textbf{(27)}\\
    		\zeta_{9} & 2 \, \theta_{0} + 4 & 0 &\textbf{(34)}\\
    		\frac{\bar{\alpha_7}}{2\theta_0^2} & 2 \, \theta_{0} & 4 &\textbf{(38)}
    		\end{dmatrix}
    		,\;\,\begin{dmatrix}
    		\overline{\zeta_{9}} & 0 & 2 \, \theta_{0} + 4 \\
    		\frac{\alpha_{7}}{2 \, \theta_{0}^{2}} & 4 & 2 \, \theta_{0} \\
    		\overline{\zeta_{8}} & \theta_{0} & \theta_{0} + 4 \\
    		\frac{2 \, \overline{\zeta_{11}}}{{\left(\theta_{0} + 3\right)} {\left(\theta_{0} + 1\right)}} & \theta_{0} + 1 & \theta_{0} + 3 \\
    		\zeta_{10} & \theta_{0} + 2 & \theta_{0} + 2 \\
    		\frac{2 \, \zeta_{11}}{{\left(\theta_{0} + 3\right)} {\left(\theta_{0} + 1\right)}} & \theta_{0} + 3 & \theta_{0} + 1 \\
    		\zeta_{8} & \theta_{0} + 4 & \theta_{0} \\
    		\zeta_{9} & 2 \, \theta_{0} + 4 & 0 \\
    		\frac{\overline{\alpha_{7}}}{2 \, \theta_{0}^{2}} & 2 \, \theta_{0} & 4
    		\end{dmatrix}
    	\end{align*}
    	and we see that both expressions coincide. Adding up with the previous computations, we obtain
    	\begin{align*}
    		|\phi(z)|^2=\begin{dmatrix}
    		-\frac{\overline{\zeta_{0}}}{\theta_{0}^{4} + 4 \, \theta_{0}^{3} + 5 \, \theta_{0}^{2} + 2 \, \theta_{0}} & 0 & 2 \, \theta_{0} + 2 \\
    		-\frac{4 \, \overline{\zeta_{1}}}{\theta_{0}^{4} + 6 \, \theta_{0}^{3} + 11 \, \theta_{0}^{2} + 6 \, \theta_{0}} & 0 & 2 \, \theta_{0} + 3 \\
    		\frac{1}{\theta_{0}^{2}} & \theta_{0} & \theta_{0} \\
    		\frac{\overline{\alpha_{1}}}{{\left(\theta_{0} + 2\right)} \theta_{0}} & \theta_{0} & \theta_{0} + 2 \\
    		\frac{2 \, {\left(\theta_{0} + 1\right)} \theta_{0} \overline{\alpha_{3}} - \overline{\zeta_{2}}}{2 \, {\left(\theta_{0}^{4} + 4 \, \theta_{0}^{3} + 3 \, \theta_{0}^{2}\right)}} & \theta_{0} & \theta_{0} + 3 \\
    		\frac{2 \, {\left| A_{1} \right|}^{2}}{{\left(\theta_{0} + 1\right)}^{2}} & \theta_{0} + 1 & \theta_{0} + 1 \\
    		\frac{\overline{\alpha_{5}}}{{\left(\theta_{0} + 2\right)} {\left(\theta_{0} + 1\right)}} & \theta_{0} + 1 & \theta_{0} + 2 \\
    		\frac{\alpha_{1}}{{\left(\theta_{0} + 2\right)} \theta_{0}} & \theta_{0} + 2 & \theta_{0} \\
    		\frac{\alpha_{5}}{{\left(\theta_{0} + 2\right)} {\left(\theta_{0} + 1\right)}} & \theta_{0} + 2 & \theta_{0} + 1 \\
    		\frac{2 \, \alpha_{3} {\left(\theta_{0} + 1\right)} \theta_{0} - \zeta_{2}}{2 \, {\left(\theta_{0}^{4} + 4 \, \theta_{0}^{3} + 3 \, \theta_{0}^{2}\right)}} & \theta_{0} + 3 & \theta_{0} \\
    		-\frac{\zeta_{0}}{\theta_{0}^{4} + 4 \, \theta_{0}^{3} + 5 \, \theta_{0}^{2} + 2 \, \theta_{0}} & 2 \, \theta_{0} + 2 & 0 \\
    		-\frac{4 \, \zeta_{1}}{\theta_{0}^{4} + 6 \, \theta_{0}^{3} + 11 \, \theta_{0}^{2} + 6 \, \theta_{0}} & 2 \, \theta_{0} + 3 & 0 
    		\end{dmatrix}
    		\begin{dmatrix}
    		\overline{\zeta_{9}} & 0 & 2 \, \theta_{0} + 4 \\
    		\frac{\alpha_{7}}{2 \, \theta_{0}^{2}} & 4 & 2 \, \theta_{0} \\
    		\overline{\zeta_{8}} & \theta_{0} & \theta_{0} + 4 \\
    		\frac{2 \, \overline{\zeta_{11}}}{{\left(\theta_{0} + 3\right)} {\left(\theta_{0} + 1\right)}} & \theta_{0} + 1 & \theta_{0} + 3 \\
    		\zeta_{10} & \theta_{0} + 2 & \theta_{0} + 2 \\
    		\frac{2 \, \zeta_{11}}{{\left(\theta_{0} + 3\right)} {\left(\theta_{0} + 1\right)}} & \theta_{0} + 3 & \theta_{0} + 1 \\
    		\zeta_{8} & \theta_{0} + 4 & \theta_{0} \\
    		\zeta_{9} & 2 \, \theta_{0} + 4 & 0 \\
    		\frac{\overline{\alpha_{7}}}{2 \, \theta_{0}^{2}} & 2 \, \theta_{0} & 4
    		\end{dmatrix}
    	\end{align*}
    	
    	\section{Development of $\vec{\alpha}$}
    	
    		\begin{align*}
    			&\vec{\alpha}=\\
    			&\begin{dmatrix}
    			\frac{1}{2} & \overline{B_{1}} & 0 & -\theta_{0} + 2 & \textbf{(1)}\\
    			2 \, A_{1} \overline{C_{1}} & \overline{A_{0}} & 0 & -\theta_{0} + 2  & \textbf{(2)}\\
    			\frac{1}{2} & \overline{B_{3}} & 0 & -\theta_{0} + 3  & \textbf{(3)}\\
    			2 \, A_{1} \overline{C_{1}} & \overline{A_{1}} & 0 & -\theta_{0} + 3  & \textbf{(4)}\\
    			-\ccancel{4 \, A_{0} {\left| A_{1} \right|}^{2} \overline{C_{1}}} + 2 \, A_{1} \overline{C_{2}} & \overline{A_{0}} & 0 & -\theta_{0} + 3  & \textbf{(5)}\\
    			\frac{1}{2} & \overline{B_{6}} & 0 & -\theta_{0} + 4  & \textbf{(6)}\\
    			2 \, A_{1} \overline{C_{1}} & \overline{A_{2}} & 0 & -\theta_{0} + 4  & \textbf{(7)}\\
    			-\ccancel{4 \, A_{0} {\left| A_{1} \right|}^{2} \overline{C_{2}}} - 2 \, A_{1} \overline{C_{1}} \overline{\alpha_{1}} - \ccancel{2 \, A_{0} \overline{C_{1}} \overline{\alpha_{5}}} + \frac{1}{4} \, \ccancel{\overline{B_{1}} \overline{C_{1}}} + 2 \, A_{1} \overline{C_{3}} & \overline{A_{0}} & 0 & -\theta_{0} + 4  & \textbf{(8)}\\
    			-\ccancel{4 \, A_{0} {\left| A_{1} \right|}^{2} \overline{C_{1}}} + 2 \, A_{1} \overline{C_{2}} & \overline{A_{1}} & 0 & -\theta_{0} + 4  & \textbf{(9)}\\
    			1 & \overline{B_{2}} & 1 & -\theta_{0} + 2  & \textbf{(10)}\\
    			\frac{1}{2} \, C_{1} \overline{C_{1}} & A_{0} & 1 & -\theta_{0} + 2  & \textbf{(11)}\\
    			-\ccancel{4 \, A_{0} \alpha_{1} \overline{C_{1}}} + \ccancel{2 \, A_{1} \overline{B_{1}}} + 4 \, A_{2} \overline{C_{1}} & \overline{A_{0}} & 1 & -\theta_{0} + 2  & \textbf{(12)}\\
    			1 & \overline{B_{5}} & 1 & -\theta_{0} + 3  & \textbf{(13)}\\
    			\frac{1}{2} \, \ccancel{B_{1} \overline{C_{1}}} + \frac{1}{2} \, C_{1} \overline{C_{2}} & A_{0} & 1 & -\theta_{0} + 3  & \textbf{(14)}\\
    			-\ccancel{4 \, A_{0} \alpha_{1} \overline{C_{1}}} + \ccancel{2 \, A_{1} \overline{B_{1}}} + 4 \, A_{2} \overline{C_{1}} & \overline{A_{1}} & 1 & -\theta_{0} + 3  & \textbf{(15)}\\
    			-4 \, A_{0} {\left| A_{1} \right|}^{2} \overline{B_{1}} - 8 \, A_{1} {\left| A_{1} \right|}^{2} \overline{C_{1}} - \ccancel{4 \, A_{0} \alpha_{5} \overline{C_{1}}} - \ccancel{4 \, A_{0} \alpha_{1} \overline{C_{2}}} + 2 \, A_{1} \overline{B_{3}} + 4 \, A_{2} \overline{C_{2}} & \overline{A_{0}} & 1 & -\theta_{0} + 3  & \textbf{(16)}\\
    			\frac{3}{2} & \overline{B_{4}} & 2 & -\theta_{0} + 2  & \textbf{(17)}\\
    			\frac{A_{1} C_{1}}{2 \, \theta_{0}} & \overline{C_{1}} & 2 & -\theta_{0} + 2  & \textbf{(18)}\\
    			\frac{1}{2} \, C_{1} \overline{C_{1}} & A_{1} & 2 & -\theta_{0} + 2  & \textbf{(19)}\\
    			\frac{1}{4} \, A_{1} \overline{C_{1}} & C_{1} & 2 & -\theta_{0} + 2  & \textbf{(20)}\\
    			\frac{1}{2} \, \ccancel{C_{1} \overline{B_{1}}} + \frac{1}{2} \, C_{2} \overline{C_{1}} & A_{0} & 2 & -\theta_{0} + 2  & \textbf{(21)}\\
    			-\ccancel{2 \, {\left(\theta_{0} \overline{\alpha_{2}} + \overline{\alpha_{2}}\right)} A_{0} C_{1}} - 4 \, A_{0} \alpha_{1} \overline{B_{1}} - 6 \, A_{1} \alpha_{1} \overline{C_{1}} - \ccancel{6 \, A_{0} \alpha_{3} \overline{C_{1}}} + 4 \, A_{2} \overline{B_{1}} + \ccancel{2 \, A_{1} \overline{B_{2}}} + 6 \, A_{3} \overline{C_{1}} & \overline{A_{0}} & 2 & -\theta_{0} + 2 & \textbf{(22)} 
    			\end{dmatrix}
    			\end{align*}
    			\begin{align*}
    			\begin{dmatrix}
    			-\theta_{0} + 2 & E_{1} & -2 \, \theta_{0} + 3 & \theta_{0}  & \textbf{(23)}\\
    			-\frac{C_{1}^{2} {\left(\theta_{0} - 2\right)}}{2 \, \theta_{0}} & \overline{A_{0}} & -2 \, \theta_{0} + 3 & \theta_{0}  & \textbf{(24)}\\
    			-\theta_{0} + 2 & E_{2} & -2 \, \theta_{0} + 3 & \theta_{0} + 1  & \textbf{(25)}\\
    			-\frac{C_{1}^{2} {\left(\theta_{0} - 2\right)}}{2 \, \theta_{0}} & \overline{A_{1}} & -2 \, \theta_{0} + 3 & \theta_{0} + 1  & \textbf{(26)}\\
    			\ccancel{2 \, {\left(\alpha_{2} \theta_{0} - 2 \, \alpha_{2}\right)} A_{0} C_{1}} - \frac{\ccancel{B_{1} C_{1} {\left(\theta_{0} - 2\right)}}}{2 \, {\left(\theta_{0} + 1\right)}} - \frac{\ccancel{B_{1} C_{1} {\left(\theta_{0} - 2\right)}}}{2 \, \theta_{0}} & \overline{A_{0}} & -2 \, \theta_{0} + 3 & \theta_{0} + 1  & \textbf{(27)}\\
    			-\theta_{0} + \frac{5}{2} & E_{3} & -2 \, \theta_{0} + 4 & \theta_{0}  & \textbf{(28)}\\
    			\ccancel{2 \, A_{1} E_{1}} - \frac{C_{1} C_{2} {\left(\theta_{0} - 2\right)}}{2 \, \theta_{0}} - \frac{C_{1} C_{2} {\left(\theta_{0} - 3\right)}}{2 \, \theta_{0}} & \overline{A_{0}} & -2 \, \theta_{0} + 4 & \theta_{0}  & \textbf{(29)}\\
    			-\frac{1}{2} \, \theta_{0} + 1 & C_{1} & -\theta_{0} + 1 & 0  & \textbf{(30)}\\
    			-\frac{1}{2} \, \theta_{0} + 1 & B_{1} & -\theta_{0} + 1 & 1  & \textbf{(31)}\\
    			-\frac{1}{2} \, \theta_{0} + 1 & B_{2} & -\theta_{0} + 1 & 2  & \textbf{(32)}\\
    			-\frac{C_{1} {\left(\theta_{0} - 2\right)} \overline{C_{1}}}{2 \, \theta_{0}} & \overline{A_{0}} & -\theta_{0} + 1 & 2  & \textbf{(33)}\\
    			-\frac{1}{2} \, \theta_{0} + 1 & B_{4} & -\theta_{0} + 1 & 3  & \textbf{(34)}\\
    			-\frac{C_{1} {\left(\theta_{0} - 2\right)} \overline{C_{1}}}{2 \, \theta_{0}} & \overline{A_{1}} & -\theta_{0} + 1 & 3  & \textbf{(35)}\\
    			\ccancel{2 \, {\left(\alpha_{2} \theta_{0} - 2 \, \alpha_{2}\right)} A_{0} \overline{C_{1}}} - \frac{\ccancel{B_{1} {\left(\theta_{0} - 2\right)} \overline{C_{1}}}}{2 \, {\left(\theta_{0} + 1\right)}} - \frac{C_{1} {\left(\theta_{0} - 2\right)} \overline{C_{2}}}{2 \, \theta_{0}} & \overline{A_{0}} & -\theta_{0} + 1 & 3  & \textbf{(36)}\\
    			-\frac{1}{2} \, \theta_{0} + \frac{3}{2} & C_{2} & -\theta_{0} + 2 & 0  & \textbf{(37)}\\
    			2 \, A_{1} C_{1} & \overline{A_{0}} & -\theta_{0} + 2 & 0  & \textbf{(38)}\\
    			-\frac{1}{2} \, \theta_{0} + \frac{3}{2} & B_{3} & -\theta_{0} + 2 & 1  & \textbf{(39)}\\
    			2 \, A_{1} C_{1} & \overline{A_{1}} & -\theta_{0} + 2 & 1  & \textbf{(40)}\\
    			-\ccancel{4 \, A_{0} C_{1} {\left| A_{1} \right|}^{2}} + \ccancel{2 \, A_{1} B_{1}} & \overline{A_{0}} & -\theta_{0} + 2 & 1  & \textbf{(41)}\\
    			-\frac{1}{2} \, \theta_{0} + \frac{3}{2} & B_{5} & -\theta_{0} + 2 & 2  & \textbf{(42)}\\
    			2 \, A_{1} C_{1} & \overline{A_{2}} & -\theta_{0} + 2 & 2  & \textbf{(43)}\\
    			-\ccancel{4 \, A_{0} C_{1} {\left| A_{1} \right|}^{2}} + \ccancel{2 \, A_{1} B_{1}} & \overline{A_{1}} & -\theta_{0} + 2 & 2  & \textbf{(44)}\\
    			\zeta_{12} & \overline{A_{0}} & -\theta_{0} + 2 & 2  & \textbf{(45)}\\
    			-\frac{1}{2} \, \theta_{0} + 2 & C_{3} & -\theta_{0} + 3 & 0  & \textbf{(46)}\\
    			\frac{1}{4} \, C_{1}^{2} & A_{0} & -\theta_{0} + 3 & 0  & \textbf{(47)}\\
    			-\ccancel{4 \, A_{0} C_{1} \alpha_{1}} + 4 \, A_{2} C_{1} + 2 \, A_{1} C_{2} & \overline{A_{0}} & -\theta_{0} + 3 & 0 & \textbf{(48)} 
    			\end{dmatrix}
    			\end{align*}
    			\begin{align*}
    			\begin{dmatrix}
    			-\frac{1}{2} \, \theta_{0} + 2 & B_{6} & -\theta_{0} + 3 & 1  & \textbf{(49)}\\
    			\frac{1}{2} \, \ccancel{B_{1} C_{1}} & A_{0} & -\theta_{0} + 3 & 1  & \textbf{(50)}\\
    			-\ccancel{4 \, A_{0} C_{1} \alpha_{1}} + 4 \, A_{2} C_{1} + 2 \, A_{1} C_{2} & \overline{A_{1}} & -\theta_{0} + 3 & 1  & \textbf{(51)}\\
    			-8 \, A_{1} C_{1} {\left| A_{1} \right|}^{2} - 4 \, A_{0} C_{2} {\left| A_{1} \right|}^{2} - \ccancel{4 \, A_{0} B_{1} \alpha_{1}} - \ccancel{4 \, A_{0} C_{1} \alpha_{5}} + 4 \, A_{2} B_{1} + \ccancel{2 \, A_{1} B_{3}} & \overline{A_{0}} & -\theta_{0} + 3 & 1  & \textbf{(52)}\\
    			-\frac{1}{2} \, \theta_{0} + \frac{5}{2} & C_{4} & -\theta_{0} + 4 & 0  & \textbf{(53)}\\
    			\frac{1}{2} \, C_{1} C_{2} & A_{0} & -\theta_{0} + 4 & 0  & \textbf{(54)}\\
    			\frac{1}{4} \, C_{1}^{2} & A_{1} & -\theta_{0} + 4 & 0  & \textbf{(55)}\\
    			\frac{1}{4} \, A_{1} C_{1} & C_{1} & -\theta_{0} + 4 & 0  & \textbf{(56)}\\
    			-6 \, A_{1} C_{1} \alpha_{1} - 4 \, A_{0} C_{2} \alpha_{1} - \ccancel{6 \, A_{0} C_{1} \alpha_{3}} + 6 \, A_{3} C_{1} + 4 \, A_{2} C_{2} + 2 \, A_{1} C_{3} & \overline{A_{0}} & -\theta_{0} + 4 & 0  & \textbf{(57)}\\
    			\frac{A_{1} \overline{C_{1}}}{2 \, \theta_{0}} & \overline{C_{1}} & \theta_{0} & -2 \, \theta_{0} + 4  & \textbf{(58)}\\
    			\frac{1}{2} \, \theta_{0} + \frac{1}{2} & \overline{E_{2}} & \theta_{0} & -2 \, \theta_{0} + 4  & \textbf{(59)}\\
    			\frac{1}{4} \, \overline{C_{1}}^{2} & A_{1} & \theta_{0} & -2 \, \theta_{0} + 4  & \textbf{(60)}\\
    			\frac{1}{2} \, \ccancel{\overline{B_{1}} \overline{C_{1}}} & A_{0} & \theta_{0} & -2 \, \theta_{0} + 4  & \textbf{(61)}\\
    			-\ccancel{2 \, {\left(\theta_{0} \overline{\alpha_{2}} + \overline{\alpha_{2}}\right)} A_{0} \overline{C_{1}}} + \ccancel{2 \, A_{1} \overline{E_{1}}} & \overline{A_{0}} & \theta_{0} & -2 \, \theta_{0} + 4  & \textbf{(62)}\\
    			\frac{1}{2} \, \theta_{0} & \overline{E_{1}} & \theta_{0} - 1 & -2 \, \theta_{0} + 4  & \textbf{(63)}\\
    			\frac{1}{4} \, \overline{C_{1}}^{2} & A_{0} & \theta_{0} - 1 & -2 \, \theta_{0} + 4  & \textbf{(64)}\\
    			\frac{1}{2} \, \theta_{0} & \overline{E_{3}} & \theta_{0} - 1 & -2 \, \theta_{0} + 5  & \textbf{(65)}\\
    			\frac{1}{2} \, \overline{C_{1}} \overline{C_{2}} & A_{0} & \theta_{0} - 1 & -2 \, \theta_{0} + 5 & \textbf{(66)}
    			\end{dmatrix}
    		\end{align*}
    		where
    		\begin{align*}
    			\zeta_{12}&=-\ccancel{4 \, A_{0} B_{1} {\left| A_{1} \right|}^{2}} - 2 \, A_{1} C_{1} \overline{\alpha_{1}} - \ccancel{2 \, A_{0} C_{1} \overline{\alpha_{5}}} + \ccancel{2 \, A_{1} B_{2}} + \frac{1}{4} \,\ccancel{ C_{1} \overline{B_{1}}} - \frac{\ccancel{C_{1} {\left(\theta_{0} - 2\right)} \overline{B_{1}}}}{2 \, \theta_{0}} - \frac{C_{2} {\left(\theta_{0} - 3\right)} \overline{C_{1}}}{2 \, \theta_{0}}\\
    			&=-\frac{(\theta_0-3)}{2\theta_0}\s{\bar{\vec{C}_1}}{\vec{C}_2}
    		\end{align*}
    		We have as
    		\begin{align*}
    		\alpha_2=\frac{1}{2\theta_0(\theta_0+1)}\s{\bar{\vec{A}_1}}{\vec{C}_1},\quad \alpha_7=\frac{1}{8\theta_0(\theta_0-4)}\s{\vec{C}_1}{\vec{C}_1}
    		\end{align*}
    		word
    		\begin{align*}
    		\frac{1}{2}\vec{E}_2&=\begin{dmatrix}
    		\frac{\ccancel{{\left(\alpha_{2} \theta_{0} + \alpha_{2}\right)} C_{1} \overline{A_{0}}} - \ccancel{E_{1} \overline{A_{1}}}}{\theta_{0} + 1} & A_{0} & -2 \, \theta_{0} + 4 & \theta_{0} + 1 \\
    		-\frac{C_{1} \overline{A_{1}}}{4 \, {\left(\theta_{0}^{2} + \theta_{0}\right)}} & C_{1} & -2 \, \theta_{0} + 4 & \theta_{0} + 1 \\
    		\frac{{\left(\ccancel{4 \, {\left(\alpha_{2} \theta_{0}^{2} + \alpha_{2} \theta_{0}\right)}} A_{0} - \ccancel{B_{1} {\left(4 \, \theta_{0} + 1\right)}}\right)} C_{1}}{8 \, {\left(\theta_{0}^{2} + \theta_{0}\right)}} & \overline{A_{0}} & -2 \, \theta_{0} + 4 & \theta_{0} + 1 \\
    		-\frac{C_{1}^{2} {\left(2 \, \theta_{0} + 1\right)}}{8 \, {\left(\theta_{0}^{2} + \theta_{0}\right)}} & \overline{A_{1}} & -2 \, \theta_{0} + 4 & \theta_{0} + 1 
    		\end{dmatrix}\\
    		&=-\frac{\alpha_2}{2}\vec{C}_1-\frac{(2\theta_0+1)(\theta_0-4)}{\theta_0+1}\alpha_7\bar{\vec{A}_1}
    		\end{align*}
    		so
    		\begin{align}
    		\vec{E}_2=-\alpha_2\vec{C}_1-\frac{2(2\theta_0+1)(\theta_0-4)}{\theta_0+1}\alpha_7\bar{\vec{A}_1}
    		\end{align}
    		Finally,
    		\begin{align*}
    		\frac{1}{2}\vec{E}_3&=
    		\begin{dmatrix}-\frac{C_{1} C_{2} {\left(2 \, \theta_{0} - 5\right)} - \ccancel{2 \, A_{1} E_{1} \theta_{0}}}{2 \, {\left(2 \, \theta_{0}^{2} - 5 \, \theta_{0}\right)}} & \overline{A_{0}} & -2 \, \theta_{0} + 5 & \theta_{0} 
    		\end{dmatrix}\\
    		&=-\frac{1}{2\theta_0}\s{\vec{C}_1}{\vec{C}_2}\bar{\vec{A}_0}
    		\end{align*}
    		So
    		\begin{align*}
    		\left\{
    		\begin{alignedat}{1}
    		\vec{E}_2&=-\alpha_2\,\vec{C}_1-\frac{2(2\theta_0+1)(\theta_0-4)}{\theta_0+1}\alpha_7\bar{\vec{A}_1}\\
    		\vec{E}_3&=-\frac{1}{\theta_0}\s{\vec{C}_1}{\vec{C}_2}\bar{\vec{A}_0}
    		\end{alignedat}\right.
    		\end{align*}
    		we have
    		Thanks of the previous development we need only compute the terms of degree $4-\theta_0$, and they correspond to the powers
    		\begin{align*}
    			\begin{dmatrix}
    			\frac{1}{2} \, \theta_{0} & \overline{E_{3}} & \theta_{0} - 1 & -2 \, \theta_{0} + 5  & \textbf{(65)}\\
    			\frac{1}{2} \, \overline{C_{1}} \overline{C_{2}} & A_{0} & \theta_{0} - 1 & -2 \, \theta_{0} + 5 & \textbf{(66)}
    			\end{dmatrix}\\
    			&=\frac{\theta_0}{2}\left(-\frac{1}{\theta_0}\bar{\s{\vec{C}_1}{\vec{C}_2}}\bar{\vec{A}_0}\right)+\frac{1}{2}\bar{\s{\vec{C}_1}{\vec{C}_2}}\vec{A}_0=0
    		\end{align*}
    		Then, we have
    		\begin{align*}
    			&\begin{dmatrix}
    			\frac{A_{1} \overline{C_{1}}}{2 \, \theta_{0}} & \overline{C_{1}} & \theta_{0} & -2 \, \theta_{0} + 4  & \textbf{(58)}\\
    			\frac{1}{2} \, \theta_{0} + \frac{1}{2} & \overline{E_{2}} & \theta_{0} & -2 \, \theta_{0} + 4  & \textbf{(59)}\\
    			\frac{1}{4} \, \overline{C_{1}}^{2} & A_{1} & \theta_{0} & -2 \, \theta_{0} + 4  & \textbf{(60)}\\
    			\frac{1}{2} \, \ccancel{\overline{B_{1}} \overline{C_{1}}} & A_{0} & \theta_{0} & -2 \, \theta_{0} + 4  & \textbf{(61)}\\
    			-\ccancel{2 \, {\left(\theta_{0} \overline{\alpha_{2}} + \overline{\alpha_{2}}\right)} A_{0} \overline{C_{1}}} + \ccancel{2 \, A_{1} \overline{E_{1}}} & \overline{A_{0}} & \theta_{0} & -2 \, \theta_{0} + 4  & \textbf{(62)}
    			\end{dmatrix}\\
    			&=(\theta_0+1)\bar{\alpha_2}\bar{\vec{C}_1}+\frac{(\theta_0+1)}{2}\bar{\vec{E}_2}+2\theta_0(\theta_0-4)\bar{\alpha_7}\vec{A}_1\\
    			&=(\theta_0+1)\bar{\alpha_2}\bar{\vec{C}_1}+\frac{(\theta_0+1)}{2}\left(-\bar{\alpha_2}\bar{\vec{C}_1}-\frac{2(2\theta_0+1)(\theta_0-4)}{\theta_0+1}\bar{\alpha_7}\vec{A}_1\right)+2\theta_0(\theta_0-4)\bar{\alpha_7}\vec{A}_1\\
    			&=\frac{(\theta_0+1)}{2}\bar{\alpha_2}\bar{\vec{C}_1}-(2\theta_0+1)(\theta_0-4)\bar{\alpha_7}\vec{A}_1+2\theta_0(\theta_0-4)\bar{\alpha_7}\vec{A}_1\\
    			&=\frac{(\theta_0+1)}{2}\bar{\alpha_2}\bar{\vec{C}_1}-(\theta_0-4)\bar{\alpha_7}\vec{A_1}
    		\end{align*}
    		Then, we compute
    		\begin{align*}
    		&	\begin{dmatrix}
    				-\theta_{0} + 2 & E_{2} & -2 \, \theta_{0} + 3 & \theta_{0} + 1  & \textbf{(25)}\\
    			-\frac{C_{1}^{2} {\left(\theta_{0} - 2\right)}}{2 \, \theta_{0}} & \overline{A_{1}} & -2 \, \theta_{0} + 3 & \theta_{0} + 1  & \textbf{(26)}\\
    			\ccancel{2 \, {\left(\alpha_{2} \theta_{0} - 2 \, \alpha_{2}\right)} A_{0} C_{1}} - \frac{\ccancel{B_{1} C_{1} {\left(\theta_{0} - 2\right)}}}{2 \, {\left(\theta_{0} + 1\right)}} - \frac{\ccancel{B_{1} C_{1} {\left(\theta_{0} - 2\right)}}}{2 \, \theta_{0}} & \overline{A_{0}} & -2 \, \theta_{0} + 3 & \theta_{0} + 1 & \textbf{(27)} 
    			\end{dmatrix}\\
    			&=-(\theta_0-2)\left(-\alpha_2\vec{C}_1-\frac{2(2\theta_0+1)(\theta_0-4)}{\theta_0+1}\alpha_7\bar{\vec{A}_1}\right)-4(\theta_0-2)(\theta_0-4)\alpha_7\bar{\vec{A}_1}\\
    			&=(\theta_0-2)\alpha_2\vec{C}_1-\frac{2(\theta_0-2)(\theta_0-4)}{\theta_0+1}\alpha_7\bar{\vec{A}_1}
    			\end{align*}
    			and finally
    			\begin{align*}
    			&\begin{dmatrix}
    			-\theta_{0} + \frac{5}{2} & E_{3} & -2 \, \theta_{0} + 4 & \theta_{0}  & \textbf{(28)}\\
    			\ccancel{2 \, A_{1} E_{1}} - \frac{C_{1} C_{2} {\left(\theta_{0} - 2\right)}}{2 \, \theta_{0}} - \frac{C_{1} C_{2} {\left(\theta_{0} - 3\right)}}{2 \, \theta_{0}} & \overline{A_{0}} & -2 \, \theta_{0} + 4 & \theta_{0} & \textbf{(29)} 
    			\end{dmatrix}\\
    			&=-\frac{(2\theta_0-5)}{2}\vec{E}_3-\frac{(2\theta_0-5)}{2\theta_0}\s{\vec{C}_1}{\vec{C}_2}\bar{\vec{A}_0}\\
    			&=-\frac{(2\theta_0-5)}{2}\left(-\frac{1}{\theta_0}\s{\vec{C}_1}{\vec{C}_2}\bar{\vec{A}_0}\right)-\frac{(2\theta_0-5)}{2\theta_0}\s{\vec{C}_1}{\vec{C}_2}\bar{\vec{A}_0}\\
    			&=0
    		\end{align*}
    		Finally, we see that the new powers of order $4-\theta_0$ are
    		\begin{align*}
    			\begin{dmatrix}
    			0&-\theta_0+4 & \textbf{(6)}-\textbf{(9)}\\
    			1&-\theta_0+3& \textbf{(13)}-\textbf{(16)}\\
    			2&-\theta_0+2 & \textbf{(17)}-\textbf{(22)}\\
    			-2\theta_0+3& \theta_0+1& \textbf{(25)}-\textbf{(27)}\\
    			-\theta_0+1& 3& \textbf{(34)}-\textbf{(36)}\\
    			-\theta_0+2& 2& \textbf{(42)}-\textbf{(45)}\\
    			-\theta_0+3& 1& \textbf{(49)}-\textbf{(52)}\\
    			-\theta_0+4& 0& \textbf{(54)}-\textbf{(57)}\\
    			\theta_0& -2\theta_0+4& \textbf{(58)}-\textbf{(62)}
    			\end{dmatrix}
    		\end{align*}
    		so for some $\vec{B}_7,\vec{B}_8,\vec{B}_9,\vec{B}_{10},\vec{B}_{11},\vec{B}_{12},\vec{B}_{13}\in\mathbb{C}^n$, the new powers are
    		\begin{align*}
    			\begin{dmatrix}
    			1 & B_7 & 0 & -\theta_0+4 &\textbf{(6)}-\textbf{(9)}\\
    			1 & {B}_8 & 1 & -\theta_0+3& \textbf{(13)}-\textbf{(16)}\\
    			1 & {B}_9 & 2 & -\theta_0+2&\textbf{(17)}-\textbf{(22)}\\
    			(\theta_0-2)\alpha_2 & C_1 & -2\theta_0+3 & \theta_0+1&\textbf{(25)}-\textbf{(27)}\\
    			-\frac{2(\theta_0-2)(\theta_0-4)}{\theta_0+1}\alpha_{7} & \bar{A_1} & -2\theta_0+3 & \theta_0+1&\textbf{(25)}-\textbf{(27)}\\
    			1 & B_{10} & -\theta_0+1 &3 &\textbf{(34)}-\textbf{(36)}\\
    			1 & B_{11} & -\theta_0+2 &2&\textbf{(42)}-\textbf{(45)}\\
    			1 & B_{12} & -\theta_0+3 & 1&\textbf{(49)}-\textbf{(52)}\\
    			1 & B_{13} & -\theta_0+4 & 0&\textbf{(54)}-\textbf{(57)}\\
    			\frac{(\theta_0+1)}{2}\bar{\alpha_2} & \bar{C_1} & \theta_0 & -2\theta_0+4&\textbf{(58)}-\textbf{(62)}\\
    			-(\theta_0-4)\bar{\alpha_7} & A_1 & \theta_0 & -2\theta_0+4 &\textbf{(58)}-\textbf{(62)}.
    			\end{dmatrix}
    		\end{align*}
    		to be compared with the Sage version
    		\small
    		\begin{align*}
    		\begin{dmatrix}
    		1 & B_7 & 0 & -\theta_0+4 &\textbf{(6)}-\textbf{(9)}\\
    		1 & {B}_8 & 1 & -\theta_0+3& \textbf{(13)}-\textbf{(16)}\\
    		1 & {B}_9 & 2 & -\theta_0+2&\textbf{(17)}-\textbf{(22)}\\
    		(\theta_0-2)\alpha_2 & C_1 & -2\theta_0+3 & \theta_0+1&\textbf{(25)}-\textbf{(27)}\\
    		-\frac{2(\theta_0-2)(\theta_0-4)}{\theta_0+1}\alpha_{7} & \bar{A_1} & -2\theta_0+3 & \theta_0+1&\textbf{(25)}-\textbf{(27)}\\
    		1 & B_{10} & -\theta_0+1 &3 &\textbf{(34)}-\textbf{(36)}\\
    		1 & B_{11} & -\theta_0+2 &2&\textbf{(42)}-\textbf{(45)}\\
    		1 & B_{12} & -\theta_0+3 & 1&\textbf{(49)}-\textbf{(52)}\\
    		1 & B_{13} & -\theta_0+4 & 0&\textbf{(54)}-\textbf{(57)}\\
    		\frac{(\theta_0+1)}{2}\bar{\alpha_2} & \bar{C_1} & \theta_0 & -2\theta_0+4&\textbf{(58)}-\textbf{(62)}\\
    		-(\theta_0-4)\bar{\alpha_7} & A_1 & \theta_0 & -2\theta_0+4 &\textbf{(58)}-\textbf{(62)}.
    		\end{dmatrix}
    			\begin{dmatrix}
    			1 & B_{7} & 0 & -\theta_{0} + 4 \\
    			1 & B_{8} & 1 & -\theta_{0} + 3 \\
    			1 & B_{9} & 2 & -\theta_{0} + 2 \\
    			\alpha_{2} {\left(\theta_{0} - 2\right)} & C_{1} & -2 \, \theta_{0} + 3 & \theta_{0} + 1 \\
    			-\frac{2 \, \alpha_{7} {\left(\theta_{0} - 2\right)} {\left(\theta_{0} - 4\right)}}{\theta_{0} + 1} & \overline{A_{1}} & -2 \, \theta_{0} + 3 & \theta_{0} + 1 \\
    			1 & B_{10} & -\theta_{0} + 1 & 3 \\
    			1 & B_{11} & -\theta_{0} + 2 & 2 \\
    			1 & B_{12} & -\theta_{0} + 3 & 1 \\
    			1 & B_{13} & -\theta_{0} + 4 & 0 \\
    			\frac{1}{2} \, {\left(\theta_{0} + 1\right)} \overline{\alpha_{2}} & \overline{C_{1}} & \theta_{0} & -2 \, \theta_{0} + 4 \\
    			-{\left(\theta_{0} - 4\right)} \overline{\alpha_{7}} & A_{1} & \theta_{0} & -2 \, \theta_{0} + 4
    			\end{dmatrix}
    		\end{align*}
    		Using the previous development of $\vec{\alpha}$ up to order $4-\theta_0$, we obtain
    		\begin{align*}
    			\vec{\alpha}=\begin{dmatrix}
    			2 \, {\left(\theta_{0} + 1\right)} \theta_{0} \overline{\alpha_{2}} & \overline{A_{0}} & 0 & -\theta_{0} + 2 \\
    			\frac{\overline{\zeta_{2}}}{\theta_{0} - 3} & A_{1} & 0 & -\theta_{0} + 3 \\
    			2 \, {\left(\theta_{0} + 1\right)} \theta_{0} \overline{\alpha_{2}} & \overline{A_{1}} & 0 & -\theta_{0} + 3 \\
    			\overline{\zeta_{3}} & \overline{A_{0}} & 0 & -\theta_{0} + 3 \\
    			\frac{{\left(\theta_{0} - 2\right)} {\left| C_{1} \right|}^{2}}{4 \, \theta_{0}} & A_{0} & 1 & -\theta_{0} + 2 \\
    			2 \, \overline{\zeta_{4}} & \overline{A_{0}} & 1 & -\theta_{0} + 2 \\
    			-\frac{1}{2} \, \theta_{0} + 1 & C_{1} & -\theta_{0} + 1 & 0 \\
    			-\frac{1}{2} \, \theta_{0} + 1 & B_{1} & -\theta_{0} + 1 & 1 \\
    			-\frac{1}{2} \, \theta_{0} + 1 & B_{2} & -\theta_{0} + 1 & 2 \\
    			-\frac{{\left(\theta_{0} - 2\right)} {\left| C_{1} \right|}^{2}}{2 \, \theta_{0}} & \overline{A_{0}} & -\theta_{0} + 1 & 2 \\
    			-\frac{1}{2} \, \theta_{0} + \frac{3}{2} & C_{2} & -\theta_{0} + 2 & 0 \\
    			2 \, \zeta_{2} & \overline{A_{0}} & -\theta_{0} + 2 & 0 \\
    			-\frac{1}{2} \, \theta_{0} + \frac{3}{2} & B_{3} & -\theta_{0} + 2 & 1 \\
    			2 \, \zeta_{2} & \overline{A_{1}} & -\theta_{0} + 2 & 1 \\
    			-\frac{1}{2} \, \theta_{0} + 2 & C_{3} & -\theta_{0} + 3 & 0 \\
    			2 \, \alpha_{7} {\left(\theta_{0} - 4\right)} \theta_{0} & A_{0} & -\theta_{0} + 3 & 0 \\
    			2 \, \zeta_{5} & \overline{A_{0}} & -\theta_{0} + 3 & 0 
    			\end{dmatrix}
    			\begin{dmatrix}
    			1 & B_{7} & 0 & -\theta_{0} + 4 \\
    			1 & B_{8} & 1 & -\theta_{0} + 3 \\
    			1 & B_{9} & 2 & -\theta_{0} + 2 \\
    			\alpha_{2} {\left(\theta_{0} - 2\right)} & C_{1} & -2 \, \theta_{0} + 3 & \theta_{0} + 1 \\
    			-\frac{2 \, \alpha_{7} {\left(\theta_{0} - 2\right)} {\left(\theta_{0} - 4\right)}}{\theta_{0} + 1} & \overline{A_{1}} & -2 \, \theta_{0} + 3 & \theta_{0} + 1 \\
    			1 & B_{10} & -\theta_{0} + 1 & 3 \\
    			1 & B_{11} & -\theta_{0} + 2 & 2 \\
    			1 & B_{12} & -\theta_{0} + 3 & 1 \\
    			1 & B_{13} & -\theta_{0} + 4 & 0 \\
    			\frac{1}{2} \, {\left(\theta_{0} + 1\right)} \overline{\alpha_{2}} & \overline{C_{1}} & \theta_{0} & -2 \, \theta_{0} + 4 \\
    			-{\left(\theta_{0} - 4\right)} \overline{\alpha_{7}} & A_{1} & \theta_{0} & -2 \, \theta_{0} + 4
    			\end{dmatrix}
    		\end{align*}
    		
    		\section{Development of $\s{\vec{\alpha}}{\phi}$}
    		
    		We have
    		\footnotesize
    		\begin{align*}
    			&\s{\vec{\alpha}}{\phi}=\\
    			&\begin{dmatrix}
    			2 \, {\left(\theta_{0} \overline{\alpha_{2}} + \overline{\alpha_{2}}\right)} \overline{A_{0}}^{2} & 0 & 2 & \textbf{(1)}\\
    			\frac{{\left(\theta_{0} - 3\right)} \overline{A_{0}}^{2} \overline{\zeta_{3}} + 2 \, {\left(2 \, \theta_{0}^{3} \overline{\alpha_{2}} - 5 \, \theta_{0}^{2} \overline{\alpha_{2}} - 3 \, \theta_{0} \overline{\alpha_{2}}\right)} \overline{A_{0}} \overline{A_{1}} + A_{1} \overline{A_{0}} \overline{\zeta_{2}}}{\theta_{0}^{2} - 3 \, \theta_{0}} & 0 & 3 & \textbf{(2)}\\
    			\zeta_{13} & 0 & 4 & \textbf{(3)}\\
    			-\frac{A_{0} C_{1} {\left(\theta_{0} - 2\right)}}{2 \, \theta_{0}} & 1 & 0 & \textbf{(4)}\\
    			-\frac{A_{0} B_{1} {\left(\theta_{0} - 2\right)}}{2 \, \theta_{0}} & 1 & 1 & \textbf{(5)}\\
    			-\frac{4 \, A_{0} {\left(\theta_{0} - 2\right)} {\left| C_{1} \right|}^{2} \overline{A_{0}} - 32 \, \theta_{0} \overline{A_{0}}^{2} \overline{\zeta_{4}} + 8 \, {\left(\theta_{0}^{2} - 2 \, \theta_{0}\right)} A_{0} B_{2} + {\left(\theta_{0}^{2} - 2 \, \theta_{0}\right)} C_{1} \overline{C_{1}}}{16 \, \theta_{0}^{2}} & 1 & 2 & \textbf{(6)}\\
    			\zeta_{14} & 1 & 3 & \textbf{(7)}\\
    			\frac{4 \, A_{0} {\left(\theta_{0} + 1\right)} \zeta_{2} \overline{A_{0}} - {\left(\theta_{0}^{2} - 2 \, \theta_{0}\right)} A_{1} C_{1} - {\left(\theta_{0}^{2} - 2 \, \theta_{0} - 3\right)} A_{0} C_{2}}{2 \, {\left(\theta_{0}^{2} + \theta_{0}\right)}} & 2 & 0 & \textbf{(8)}\\
    			\frac{4 \, A_{0} {\left(\theta_{0} + 1\right)} \zeta_{2} \overline{A_{1}} - {\left(\theta_{0}^{2} - 2 \, \theta_{0}\right)} A_{1} B_{1} - {\left(\theta_{0}^{2} - 2 \, \theta_{0} - 3\right)} A_{0} B_{3}}{2 \, {\left(\theta_{0}^{2} + \theta_{0}\right)}} & 2 & 1 & \textbf{(9)}\\
    			\zeta_{15} & 2 & 2 & \textbf{(10)}\\
    			\zeta_{7} & 3 & 0 & \textbf{(11)}\\
    			\zeta_{16} & 3 & 1 & \textbf{(12)}\\
    			\zeta_{17} & 4 & 0 & \textbf{(13)}\\
    			\frac{{\left(\alpha_{2} \theta_{0}^{2} - \alpha_{2} \theta_{0} - 2 \, \alpha_{2}\right)} C_{1} \overline{A_{0}} - 2 \, {\left(\alpha_{7} \theta_{0}^{2} - 6 \, \alpha_{7} \theta_{0} + 8 \, \alpha_{7}\right)} \overline{A_{0}} \overline{A_{1}}}{\theta_{0}^{2} + \theta_{0}} & -2 \, \theta_{0} + 3 & 2 \, \theta_{0} + 1 & \textbf{(14)}\\
    			-\frac{C_{1} {\left(\theta_{0} - 2\right)} \overline{A_{0}}}{2 \, \theta_{0}} & -\theta_{0} + 1 & \theta_{0} & \textbf{(15)}\\
    			-\frac{{\left(\theta_{0}^{2} - \theta_{0} - 2\right)} B_{1} \overline{A_{0}} + {\left(\theta_{0}^{2} - 2 \, \theta_{0}\right)} C_{1} \overline{A_{1}}}{2 \, {\left(\theta_{0}^{2} + \theta_{0}\right)}} & -\theta_{0} + 1 & \theta_{0} + 1 & \textbf{(16)}\\
    			\lambda_2 & -\theta_{0} + 1 & \theta_{0} + 2 & \textbf{(17)}\\
    			\zeta_{18} & -\theta_{0} + 1 & \theta_{0} + 3 & \textbf{(18)}\\
    			-\frac{C_{2} {\left(\theta_{0} - 3\right)} \overline{A_{0}} - 4 \, \zeta_{2} \overline{A_{0}}^{2}}{2 \, \theta_{0}} & -\theta_{0} + 2 & \theta_{0} & \textbf{(19)}\\
    			\frac{4 \, {\left(2 \, \theta_{0} + 1\right)} \zeta_{2} \overline{A_{0}} \overline{A_{1}} - {\left(\theta_{0}^{2} - 2 \, \theta_{0} - 3\right)} B_{3} \overline{A_{0}} - {\left(\theta_{0}^{2} - 3 \, \theta_{0}\right)} C_{2} \overline{A_{1}}}{2 \, {\left(\theta_{0}^{2} + \theta_{0}\right)}} & -\theta_{0} + 2 & \theta_{0} + 1 & \textbf{(20)}\\
    			\zeta_{19} & -\theta_{0} + 2 & \theta_{0} + 2 & \textbf{(21)}\\
    			-\frac{C_{1}^{2} {\left(\theta_{0} - 2\right)} - 32 \, {\left(\alpha_{7} \theta_{0}^{2} - 4 \, \alpha_{7} \theta_{0}\right)} A_{0} \overline{A_{0}} + 8 \, C_{3} {\left(\theta_{0} - 4\right)} \overline{A_{0}} - 32 \, \zeta_{5} \overline{A_{0}}^{2}}{16 \, \theta_{0}} & -\theta_{0} + 3 & \theta_{0} & \textbf{(22)}\\
    			\zeta_{20} & -\theta_{0} + 3 & \theta_{0} + 1 & \textbf{(23)}\\
    			-\frac{C_{1} C_{2} {\left(5 \, \theta_{0} - 13\right)} - 12 \, C_{1} \zeta_{2} \overline{A_{0}} - 48 \, B_{13} \overline{A_{0}}}{48 \, \theta_{0}} & -\theta_{0} + 4 & \theta_{0} & \textbf{(24)}
    			\end{dmatrix}
    			\end{align*}
    			\begin{align}\label{ordre4}
    			\begin{dmatrix}
    			2 \, {\left(\theta_{0} \overline{\alpha_{2}} + \overline{\alpha_{2}}\right)} A_{0} \overline{A_{0}} & \theta_{0} & -\theta_{0} + 2 & \textbf{(25)}\\
    			\frac{A_{0} {\left(\theta_{0} - 3\right)} \overline{A_{0}} \overline{\zeta_{3}} + 2 \, {\left(\theta_{0}^{3} \overline{\alpha_{2}} - 2 \, \theta_{0}^{2} \overline{\alpha_{2}} - 3 \, \theta_{0} \overline{\alpha_{2}}\right)} A_{0} \overline{A_{1}} + A_{0} A_{1} \overline{\zeta_{2}}}{\theta_{0}^{2} - 3 \, \theta_{0}} & \theta_{0} & -\theta_{0} + 3 & \textbf{(26)}\\
    			-\frac{4 \, {\left(\theta_{0} \overline{\alpha_{7}} - 4 \, \overline{\alpha_{7}}\right)} A_{1} \overline{A_{0}} - {\left(\theta_{0}^{2} \overline{\alpha_{2}} + 3 \, \theta_{0} \overline{\alpha_{2}} + 2 \, \overline{\alpha_{2}}\right)} \overline{A_{0}} \overline{C_{1}} - 4 \, A_{0} B_{7}}{4 \, \theta_{0}} & \theta_{0} & -\theta_{0} + 4 & \textbf{(27)}\\
    			\frac{8 \, A_{1} \theta_{0}^{3} \overline{A_{0}} \overline{\alpha_{2}} + A_{0}^{2} {\left(\theta_{0} - 2\right)} {\left| C_{1} \right|}^{2} + 8 \, A_{0} \theta_{0} \overline{A_{0}} \overline{\zeta_{4}}}{4 \, \theta_{0}^{2}} & \theta_{0} + 1 & -\theta_{0} + 2 & \textbf{(28)}\\
    			\frac{A_{1}^{2} \theta_{0} \overline{\zeta_{2}} + {\left(\theta_{0}^{2} - 3 \, \theta_{0}\right)} A_{1} \overline{A_{0}} \overline{\zeta_{3}} + {\left(\theta_{0}^{2} - 2 \, \theta_{0} - 3\right)} A_{0} B_{8} + 2 \, {\left(\theta_{0}^{4} \overline{\alpha_{2}} - 2 \, \theta_{0}^{3} \overline{\alpha_{2}} - 3 \, \theta_{0}^{2} \overline{\alpha_{2}}\right)} A_{1} \overline{A_{1}}}{\theta_{0}^{3} - 2 \, \theta_{0}^{2} - 3 \, \theta_{0}} & \theta_{0} + 1 & -\theta_{0} + 3& \textbf{(29)} \\
    			\frac{{\left(\theta_{0}^{2} - 4\right)} A_{0} A_{1} {\left| C_{1} \right|}^{2} + 8 \, {\left(\theta_{0}^{2} + 2 \, \theta_{0}\right)} A_{1} \overline{A_{0}} \overline{\zeta_{4}} + 4 \, {\left(\theta_{0}^{2} + 3 \, \theta_{0} + 2\right)} A_{0} B_{9} + 8 \, {\left(\theta_{0}^{4} \overline{\alpha_{2}} + 2 \, \theta_{0}^{3} \overline{\alpha_{2}} + \theta_{0}^{2} \overline{\alpha_{2}}\right)} A_{2} \overline{A_{0}}}{4 \, {\left(\theta_{0}^{3} + 3 \, \theta_{0}^{2} + 2 \, \theta_{0}\right)}} & \theta_{0} + 2 & -\theta_{0} + 2 & \textbf{(30)}\\
    			-\frac{2 \, {\left(\theta_{0} \overline{\alpha_{7}} - 4 \, \overline{\alpha_{7}}\right)} A_{0} A_{1} - {\left(\theta_{0} \overline{\alpha_{2}} + \overline{\alpha_{2}}\right)} A_{0} \overline{C_{1}}}{2 \, \theta_{0}} & 2 \, \theta_{0} & -2 \, \theta_{0} + 4& \textbf{(31)}
    			\end{dmatrix}
    		\end{align}
    		\small
    		where
    		\begin{align*}
    			&\zeta_{13}=\frac{1}{\theta_{0}^{4} - 7 \, \theta_{0}^{2} - 6 \, \theta_{0}}\bigg\{{\left(\theta_{0}^{3} - 7 \, \theta_{0} - 6\right)} B_{7} \overline{A_{0}} + 2 \, {\left(\theta_{0}^{5} \overline{\alpha_{2}} - 7 \, \theta_{0}^{3} \overline{\alpha_{2}} - 6 \, \theta_{0}^{2} \overline{\alpha_{2}}\right)} \overline{A_{1}}^{2} + 2 \, {\left(\theta_{0}^{5} \overline{\alpha_{2}} - \theta_{0}^{4} \overline{\alpha_{2}} - 5 \, \theta_{0}^{3} \overline{\alpha_{2}} - 3 \, \theta_{0}^{2} \overline{\alpha_{2}}\right)} \overline{A_{0}} \overline{A_{2}}\\
    			& + {\left({\left(\theta_{0}^{2} + 2 \, \theta_{0}\right)} A_{1} \overline{\zeta_{2}} + {\left(\theta_{0}^{3} - \theta_{0}^{2} - 6 \, \theta_{0}\right)} \overline{A_{0}} \overline{\zeta_{3}}\right)} \overline{A_{1}}\bigg\}\\
    			&\zeta_{14}=\frac{1}{48 \, {\left(\theta_{0}^{2} + \theta_{0}\right)}}\bigg\{12 \, A_{0} {\left(\theta_{0} - 2\right)} {\left| C_{1} \right|}^{2} \overline{A_{1}} + 96 \, \theta_{0} \overline{A_{0}} \overline{A_{1}} \overline{\zeta_{4}} + 48 \, A_{0} B_{10} {\left(\theta_{0} + 1\right)} + 48 \, B_{8} {\left(\theta_{0} + 1\right)} \overline{A_{0}} - 3 \, {\left(\theta_{0}^{2} - \theta_{0} - 2\right)} B_{1} \overline{C_{1}}\\
    			& - 2 \, {\left(\theta_{0}^{2} - \theta_{0} - 2\right)} C_{1} \overline{C_{2}}\bigg\}\\
    			&\zeta_{15}=-\frac{1}{16 \, {\left(\theta_{0}^{2} + \theta_{0}\right)}}\bigg\{8 \, A_{1} {\left(\theta_{0} - 2\right)} {\left| C_{1} \right|}^{2} \overline{A_{0}} - 4 \, {\left(\theta_{0} + 1\right)} \zeta_{2} \overline{A_{0}} \overline{C_{1}} + 8 \, {\left(\theta_{0}^{2} - 2 \, \theta_{0}\right)} A_{1} B_{2} - 16 \, A_{0} B_{11} {\left(\theta_{0} + 1\right)}\\
    			& - 16 \, B_{9} {\left(\theta_{0} + 1\right)} \overline{A_{0}} + {\left(\theta_{0}^{2} - 2 \, \theta_{0} - 3\right)} C_{2} \overline{C_{1}} - {\left(4 \, {\left(\theta_{0}^{3} \overline{\alpha_{2}} + 2 \, \theta_{0}^{2} \overline{\alpha_{2}} + \theta_{0} \overline{\alpha_{2}}\right)} \overline{A_{0}} - {\left(\theta_{0}^{2} - 2 \, \theta_{0}\right)} \overline{B_{1}}\right)} C_{1}\bigg\}\\
    			&\zeta_{16}=\frac{4 \, {\left(\theta_{0}^{2} + 2 \, \theta_{0}\right)} A_{1} \zeta_{2} \overline{A_{1}} - {\left(\theta_{0}^{3} - \theta_{0}^{2} - 2 \, \theta_{0}\right)} A_{2} B_{1} + 2 \, {\left(\theta_{0}^{2} + 3 \, \theta_{0} + 2\right)} A_{0} B_{12} - {\left(\theta_{0}^{3} - \theta_{0}^{2} - 6 \, \theta_{0}\right)} A_{1} B_{3}}{2 \, {\left(\theta_{0}^{3} + 3 \, \theta_{0}^{2} + 2 \, \theta_{0}\right)}}\\
    			& - {\left(\theta_{0}^{3} - \theta_{0}^{2} - 2 \, \theta_{0}\right)} A_{2} C_{1} - {\left(\theta_{0}^{3} - \theta_{0}^{2} - 6 \, \theta_{0}\right)} A_{1} C_{2} - {\left(\theta_{0}^{3} - \theta_{0}^{2} - 10 \, \theta_{0} - 8\right)} A_{0} C_{3}\bigg\}\\
    			&\zeta_{17}=\frac{1}{2 \, {\left(\theta_{0}^{4} + 6 \, \theta_{0}^{3} + 11 \, \theta_{0}^{2} + 6 \, \theta_{0}\right)}}\bigg\{4 \, {\left(\theta_{0}^{3} + 4 \, \theta_{0}^{2} + 3 \, \theta_{0}\right)} A_{2} \zeta_{2} \overline{A_{0}} + 2 \, {\left(\theta_{0}^{3} + 6 \, \theta_{0}^{2} + 11 \, \theta_{0} + 6\right)} A_{0} B_{13}\\
    			& - {\left(\theta_{0}^{4} + \theta_{0}^{3} - 4 \, \theta_{0}^{2} - 4 \, \theta_{0}\right)} A_{3} C_{1} - {\left(\theta_{0}^{4} + \theta_{0}^{3} - 9 \, \theta_{0}^{2} - 9 \, \theta_{0}\right)} A_{2} C_{2} - {\left(\theta_{0}^{4} + \theta_{0}^{3} - 14 \, \theta_{0}^{2} - 24 \, \theta_{0}\right)} A_{1} C_{3}\\
    			& + 4 \, {\left({\left(\theta_{0}^{3} + 5 \, \theta_{0}^{2} + 6 \, \theta_{0}\right)} \zeta_{5} \overline{A_{0}} + {\left(\alpha_{7} \theta_{0}^{5} + \alpha_{7} \theta_{0}^{4} - 14 \, \alpha_{7} \theta_{0}^{3} - 24 \, \alpha_{7} \theta_{0}^{2}\right)} A_{0}\right)} A_{1}\bigg\}\\
    			&\lambda_2=-\frac{{\left(\theta_{0}^{3} + \theta_{0}^{2} - 4 \, \theta_{0} - 4\right)} {\left| C_{1} \right|}^{2} \overline{A_{0}}^{2} + {\left(\theta_{0}^{4} + \theta_{0}^{3} - 4 \, \theta_{0}^{2} - 4 \, \theta_{0}\right)} B_{2} \overline{A_{0}} + {\left(\theta_{0}^{4} - 4 \, \theta_{0}^{2}\right)} B_{1} \overline{A_{1}} + {\left(\theta_{0}^{4} - \theta_{0}^{3} - 2 \, \theta_{0}^{2}\right)} C_{1} \overline{A_{2}}}{2 \, {\left(\theta_{0}^{4} + 3 \, \theta_{0}^{3} + 2 \, \theta_{0}^{2}\right)}}\\
    			&\zeta_{18}=-\frac{1}{2 \, {\left(\theta_{0}^{4} + 6 \, \theta_{0}^{3} + 11 \, \theta_{0}^{2} + 6 \, \theta_{0}\right)}}\bigg\{{\left(\theta_{0}^{3} + 3 \, \theta_{0}^{2} - 4 \, \theta_{0} - 12\right)} {\left| C_{1} \right|}^{2} \overline{A_{0}} \overline{A_{1}} - 2 \, {\left(\theta_{0}^{3} + 6 \, \theta_{0}^{2} + 11 \, \theta_{0} + 6\right)} B_{10} \overline{A_{0}}\\
    			& + {\left(\theta_{0}^{4} + 3 \, \theta_{0}^{3} - 4 \, \theta_{0}^{2} - 12 \, \theta_{0}\right)} B_{2} \overline{A_{1}} + {\left(\theta_{0}^{4} + 2 \, \theta_{0}^{3} - 5 \, \theta_{0}^{2} - 6 \, \theta_{0}\right)} B_{1} \overline{A_{2}} + {\left(\theta_{0}^{4} + \theta_{0}^{3} - 4 \, \theta_{0}^{2} - 4 \, \theta_{0}\right)} C_{1} \overline{A_{3}}\bigg\}\\
    			&\zeta_{19}=\frac{4 \, {\left(\theta_{0}^{2} + 2 \, \theta_{0}\right)} \zeta_{2} \overline{A_{1}}^{2} + 4 \, {\left(\theta_{0}^{2} + \theta_{0}\right)} \zeta_{2} \overline{A_{0}} \overline{A_{2}} + 2 \, {\left(\theta_{0}^{2} + 3 \, \theta_{0} + 2\right)} B_{11} \overline{A_{0}} - {\left(\theta_{0}^{3} - \theta_{0}^{2} - 6 \, \theta_{0}\right)} B_{3} \overline{A_{1}} - {\left(\theta_{0}^{3} - 2 \, \theta_{0}^{2} - 3 \, \theta_{0}\right)} C_{2} \overline{A_{2}}}{2 \, {\left(\theta_{0}^{3} + 3 \, \theta_{0}^{2} + 2 \, \theta_{0}\right)}}\\
    			&\zeta_{20}=\frac{1}{16 \, {\left(\theta_{0}^{2} + \theta_{0}\right)}}\bigg\{16 \, B_{12} {\left(\theta_{0} + 1\right)} \overline{A_{0}} - 8 \, {\left(\theta_{0}^{2} - 4 \, \theta_{0}\right)} C_{3} \overline{A_{1}} + {\left(16 \, {\left(\alpha_{2} \theta_{0}^{2} - \alpha_{2} \theta_{0} - 2 \, \alpha_{2}\right)} A_{0} - {\left(2 \, \theta_{0}^{2} - 3 \, \theta_{0} - 2\right)} B_{1}\right)} C_{1}\\
    			& + 32 \, {\left(\theta_{0} \zeta_{5} \overline{A_{0}} + {\left(\alpha_{7} \theta_{0}^{3} - 5 \, \alpha_{7} \theta_{0}^{2} + 6 \, \alpha_{7} \theta_{0} - 8 \, \alpha_{7}\right)} A_{0}\right)} \overline{A_{1}}\bigg\}
    		\end{align*}
    		\normalsize
    		\begin{align*}
    		 \zeta_{16}&=\frac{4 \, {\left(\theta_{0}^{2} + 2 \, \theta_{0}\right)} A_{1} \zeta_{2} \overline{A_{1}} - {\left(\theta_{0}^{3} - \theta_{0}^{2} - 2 \, \theta_{0}\right)} A_{2} B_{1} + 2 \, {\left(\theta_{0}^{2} + 3 \, \theta_{0} + 2\right)} A_{0} B_{12} - {\left(\theta_{0}^{3} - \theta_{0}^{2} - 6 \, \theta_{0}\right)} A_{1} B_{3}}{2 \, {\left(\theta_{0}^{3} + 3 \, \theta_{0}^{2} + 2 \, \theta_{0}\right)}}\\
    		 &=\frac{4 \, {\left(\theta_{0}^{2} + 2 \, \theta_{0}\right)} A_{1} \zeta_{2} \overline{A_{1}} - {\left(\theta_{0}^{3} - \theta_{0}^{2} - 2 \, \theta_{0}\right)} A_{2} B_{1} + 2 \, {\left(\theta_{0}^{2} + 3 \, \theta_{0} + 2\right)} A_{0} B_{12} - {\left(\theta_{0}^{3} - \theta_{0}^{2} - 6 \, \theta_{0}\right)} A_{1} B_{3}}{2 \, {\left(\theta_{0}^{3} + 3 \, \theta_{0}^{2} + 2 \, \theta_{0}\right)}} 
    		\end{align*}
    		Now, recall that
    		\begin{align*}
    			\s{\vec{\alpha}}{\phi}=\begin{dmatrix}
    			\frac{{\left(2 \, \theta_{0} - 1\right)} \zeta_{2}}{2 \, {\left(\theta_{0} + 1\right)} \theta_{0}} & 2 & 0 \\
    			\zeta_{7} & 3 & 0 \\
    			\alpha_{2} {\left(\theta_{0} - 2\right)} & -\theta_{0} + 1 & \theta_{0} + 1 \\
    			\frac{{\left(\theta_{0} - 2\right)} \zeta_{4}}{{\left(\theta_{0} + 2\right)} \theta_{0}} & -\theta_{0} + 1 & \theta_{0} + 2 \\
    			\frac{{\left(\theta_{0} - 3\right)} \zeta_{3}}{2 \, {\left(\theta_{0} + 1\right)} \theta_{0}} & -\theta_{0} + 2 & \theta_{0} + 1 \\
    			{\left(\theta_{0} + 1\right)} \overline{\alpha_{2}} & \theta_{0} & -\theta_{0} + 2 \\
    			\frac{\overline{\zeta_{3}}}{2 \, \theta_{0}} & \theta_{0} & -\theta_{0} + 3 \\
    			\frac{\overline{\zeta_{4}}}{\theta_{0}} & \theta_{0} + 1 & -\theta_{0} + 2
    			\end{dmatrix}
    		\end{align*}
    		so we need only consider the order $4$ in \eqref{ordre4}. Now, recall that
    		\begin{align*}
    			\left\{
    			\begin{alignedat}{1}
    			\vec{B}_1&=-2\s{\bar{\vec{A}_1}}{\vec{C}_1}\vec{A}_0\\
    			\vec{B}_2&=-\frac{(\theta_0+2)}{4\theta_0}|\vec{C}_1|^2\bar{\vec{A}_0}+\left(\bar{\alpha_0}\s{\bar{\vec{A}_1}}{\vec{C}_1}-2\s{\bar{\vec{A}_2}}{\vec{C}_1}\right)\vec{A}_0\\
    			\vec{B}_3&=-2\s{\bar{\vec{A}_1}}{\vec{C}_1}\vec{A}_1+\frac{2}{\theta_0-3}\s{\vec{A}_1}{\vec{C}_1}\bar{\vec{A}_1}+\frac{2\bar{\alpha_0}}{\theta_0-3}\s{\vec{A}_1}{\vec{C}_1}\bar{\vec{A}_0}+2\left(\alpha_0\s{\bar{\vec{A}_1}}{\vec{C}_1}-\s{\bar{\vec{A}_1}}{\vec{C}_2}\right)\vec{A}_0\\
    			\vec{B}_4&=-\frac{(\theta_0+3)}{6\theta_0}\s{\vec{C}_1}{\bar{\vec{C}_2}}\bar{\vec{A}_0}-\frac{\bar{\zeta_2}}{6\theta_0}\vec{C}_1-\frac{(\theta_0+3)}{6\theta_0}|\vec{C}_1|^2\bar{\vec{A}_1}-\frac{\theta_0(\theta_0+1)}{6}\alpha_2\bar{\vec{C}_1}+2\left(\bar{\alpha_1}\s{\bar{\vec{A}_1}}{\vec{C}_1}-\s{\bar{\vec{A}_3}}{\vec{C}_1}\right)\vec{A}_0.\\
    			\vec{B}_5&=-\left(\frac{(\theta_0+2)}{4}\s{\bar{\vec{C}_1}}{\vec{C}_2}+\frac{2}{\theta_0-3}\bar{\alpha_1}\zeta_2\right)\bar{\vec{A}_0}+\frac{2\,\zeta_2}{\theta_0-3}\bar{\vec{A}_2}-2\s{\bar{\vec{A}_2}}{\vec{C}_1}\vec{A}_1\\
    			&+\left(8\theta_0(\theta_0+1)|\vec{A}_1|^2\alpha_2-2\s{\bar{\vec{A}_2}}{\vec{C}_2}-\frac{2}{\theta_0-3}\bar{\zeta_0}\zeta_2\right)\vec{A}_0\\
    			\vec{B}_6&=\frac{2\,\zeta_5}{\theta_0-4}\bar{\vec{A}_1}-\frac{4}{\theta_0-3}|\vec{A}_1|^2\zeta_2\bar{\vec{A}_0}+\left(-2\s{\bar{\vec{A}_1}}{\vec{C}_3}+4\theta_0(\theta_0+1)\alpha_1\alpha_2\right)\vec{A}_0-4\theta_0(\theta_0+1)\alpha_2\vec{A}_2-2\s{\bar{\vec{A}_1}}{C_2}\vec{A}_1\\
    			\vec{E}_1&=-\frac{1}{2\theta_0}\s{\vec{C}_1}{\vec{C}_1}\bar{\vec{A}_0}\\
    			\vec{E}_2&=-\alpha_2\,\vec{C}_1-\frac{2(2\theta_0+1)(\theta_0-4)}{\theta_0+1}\alpha_7\bar{\vec{A}_1}\\
    			\vec{E}_3&=-\frac{1}{\theta_0}\s{\vec{C}_1}{\vec{C}_2}\bar{\vec{A}_0}
    			\end{alignedat}\right.
    		\end{align*}
    		The new powers of order $4$ are
    		\begin{align*}
    			\textbf{(3)}, \textbf{(7)}, \textbf{(10)}, \textbf{(12)}-\textbf{(14)}, \textbf{(18)}, \textbf{(21)}, \textbf{(23)}-\textbf{(24)}, \textbf{(27)}, \textbf{(29)}-\textbf{(31)}. 
    		\end{align*}
    		or
    		\footnotesize
    		\begin{align*}
    			&\begin{dmatrix}
    			\zeta_{13} & 0 & 4 & \textbf{(3)}\\
    			\zeta_{14} & 1 & 3 & \textbf{(7)}\\
    			\zeta_{15} & 2 & 2 & \textbf{(10)}\\
    			\zeta_{16} & 3 & 1 & \textbf{(12)}\\
    			\zeta_{17} & 4 & 0 & \textbf{(13)}\\
    			\frac{\ccancel{{\left(\alpha_{2} \theta_{0}^{2} - \alpha_{2} \theta_{0} - 2 \, \alpha_{2}\right)} C_{1} \overline{A_{0}}} - \ccancel{2 \, {\left(\alpha_{7} \theta_{0}^{2} - 6 \, \alpha_{7} \theta_{0} + 8 \, \alpha_{7}\right)} \overline{A_{0}} \overline{A_{1}}}}{\theta_{0}^{2} + \theta_{0}} & -2 \, \theta_{0} + 3 & 2 \, \theta_{0} + 1 & \textbf{(14)}\\
    			\zeta_{18} & -\theta_{0} + 1 & \theta_{0} + 3 & \textbf{(18)}\\
    			\zeta_{19} & -\theta_{0} + 2 & \theta_{0} + 2 & \textbf{(21)}\\
    			\zeta_{20} & -\theta_{0} + 3 & \theta_{0} + 1 & \textbf{(23)}\\
    			-\frac{C_{1} C_{2} {\left(5 \, \theta_{0} - 13\right)} - \ccancel{12 \, C_{1} \zeta_{2} \overline{A_{0}}} - 48 \, B_{13} \overline{A_{0}}}{48 \, \theta_{0}} & -\theta_{0} + 4 & \theta_{0} & \textbf{(24)}\\
    			-\frac{\ccancel{4 \, {\left(\theta_{0} \overline{\alpha_{7}} - 4 \, \overline{\alpha_{7}}\right)} A_{1} \overline{A_{0}}} - \ccancel{{\left(\theta_{0}^{2} \overline{\alpha_{2}} + 3 \, \theta_{0} \overline{\alpha_{2}} + 2 \, \overline{\alpha_{2}}\right)} \overline{A_{0}} \overline{C_{1}}} - 4 \, A_{0} B_{7}}{4 \, \theta_{0}} & \theta_{0} & -\theta_{0} + 4 & \textbf{(27)}\\
    			\frac{A_{1}^{2} \theta_{0} \overline{\zeta_{2}} + \ccancel{{\left(\theta_{0}^{2} - 3 \, \theta_{0}\right)} A_{1} \overline{A_{0}} \overline{\zeta_{3}}} + {\left(\theta_{0}^{2} - 2 \, \theta_{0} - 3\right)} A_{0} B_{8} + 2 \, {\left(\theta_{0}^{4} \overline{\alpha_{2}} - 2 \, \theta_{0}^{3} \overline{\alpha_{2}} - 3 \, \theta_{0}^{2} \overline{\alpha_{2}}\right)} A_{1} \overline{A_{1}}}{\theta_{0}^{3} - 2 \, \theta_{0}^{2} - 3 \, \theta_{0}} & \theta_{0} + 1 & -\theta_{0} + 3& \textbf{(29)} \\
    			\frac{\ccancel{{\left(\theta_{0}^{2} - 4\right)} A_{0} A_{1} {\left| C_{1} \right|}^{2}} + \ccancel{8 \, {\left(\theta_{0}^{2} + 2 \, \theta_{0}\right)} A_{1} \overline{A_{0}} \overline{\zeta_{4}} }+ 4 \, {\left(\theta_{0}^{2} + 3 \, \theta_{0} + 2\right)} A_{0} B_{9} + 8 \, {\left(\theta_{0}^{4} \overline{\alpha_{2}} + 2 \, \theta_{0}^{3} \overline{\alpha_{2}} + \theta_{0}^{2} \overline{\alpha_{2}}\right)} A_{2} \overline{A_{0}}}{4 \, {\left(\theta_{0}^{3} + 3 \, \theta_{0}^{2} + 2 \, \theta_{0}\right)}} & \theta_{0} + 2 & -\theta_{0} + 2 & \textbf{(30)}\\
    			-\frac{\ccancel{2 \, {\left(\theta_{0} \overline{\alpha_{7}} - 4 \, \overline{\alpha_{7}}\right)} A_{0} A_{1}} - \ccancel{{\left(\theta_{0} \overline{\alpha_{2}} + \overline{\alpha_{2}}\right)} A_{0} \overline{C_{1}}}}{2 \, \theta_{0}} & 2 \, \theta_{0} & -2 \, \theta_{0} + 4& \textbf{(31)}
    			\end{dmatrix}\\
    			&=
    			\begin{dmatrix}
    			\zeta_{13} & 0 & 4 & \textbf{(3)}\\
    			\zeta_{14} & 1 & 3 & \textbf{(7)}\\
    			\zeta_{15} & 2 & 2 & \textbf{(10)}\\
    			\zeta_{16} & 3 & 1 & \textbf{(12)}\\
    			\zeta_{17} & 4 & 0 & \textbf{(13)}\\
    			\zeta_{18} & -\theta_{0} + 1 & \theta_{0} + 3 & \textbf{(18)}\\
    			\zeta_{19} & -\theta_{0} + 2 & \theta_{0} + 2 & \textbf{(21)}\\
    			\zeta_{20} & -\theta_{0} + 3 & \theta_{0} + 1 & \textbf{(23)}\\
    			-\frac{C_{1} C_{2} {\left(5 \, \theta_{0} - 13\right)}- 48 \, B_{13} \overline{A_{0}}}{48 \, \theta_{0}} & -\theta_{0} + 4 & \theta_{0} & \textbf{(24)}\\
    			\frac{A_{0} B_{7}}{ \, \theta_{0}} & \theta_{0} & -\theta_{0} + 4 & \textbf{(27)}\\
    			\frac{A_{1}^{2} \theta_{0} \overline{\zeta_{2}}+ {\left(\theta_{0}^{2} - 2 \, \theta_{0} - 3\right)} A_{0} B_{8} + 2 \, {\left(\theta_{0}^{4} \overline{\alpha_{2}} - 2 \, \theta_{0}^{3} \overline{\alpha_{2}} - 3 \, \theta_{0}^{2} \overline{\alpha_{2}}\right)} A_{1} \overline{A_{1}}}{\theta_{0}^{3} - 2 \, \theta_{0}^{2} - 3 \, \theta_{0}} & \theta_{0} + 1 & -\theta_{0} + 3& \textbf{(29)} \\
    			\frac{ 4 \, {\left(\theta_{0}^{2} + 3 \, \theta_{0} + 2\right)} A_{0} B_{9} + 8 \, {\left(\theta_{0}^{4} \overline{\alpha_{2}} + 2 \, \theta_{0}^{3} \overline{\alpha_{2}} + \theta_{0}^{2} \overline{\alpha_{2}}\right)} A_{2} \overline{A_{0}}}{4 \, {\left(\theta_{0}^{3} + 3 \, \theta_{0}^{2} + 2 \, \theta_{0}\right)}} & \theta_{0} + 2 & -\theta_{0} + 2 & \textbf{(30)}
    			\end{dmatrix}
    		\end{align*}
    		\normalsize
    		so for 
    		\begin{align*}
    			\left\{\begin{alignedat}{1}
    			\zeta_{21}&=-\frac{1}{48\theta_0}\left((5\theta_0-13)\s{\vec{C}_1}{\vec{C}_2}+\s{\bar{\vec{A}_0}}{\vec{B}_{13}}\right)\\
    			\zeta_{22}&=\frac{1}{\theta_0}\s{\vec{A}_0}{\vec{B}_7}\\
    			\zeta_{23}&=\frac{\zeta_0\bar{\zeta}_2}{(\theta_0+1)(\theta_0-3)}+\s{\vec{A}_0}{\vec{B}_8}+2|\vec{A}_1|^2\bar{\alpha_2}\\
    			\zeta_{24}&=\frac{1}{\theta_0}\s{\vec{A}_0}{\vec{B}_9}+\frac{\theta_0(\theta_0+2)}{\theta_0+1}\alpha_1\bar{\alpha_2}
    			\end{alignedat}\right.
    		\end{align*}
    		we obtain as matrix of powers of order $4$
    		\begin{align*}
    			\begin{dmatrix}
    			\zeta_{13} & 0 & 4 & \textbf{(3)}\\
    			\zeta_{14} & 1 & 3 & \textbf{(7)}\\
    			\zeta_{15} & 2 & 2 & \textbf{(10)}\\
    			\zeta_{16} & 3 & 1 & \textbf{(12)}\\
    			\zeta_{17} & 4 & 0 & \textbf{(13)}\\
    			\zeta_{18} & -\theta_{0} + 1 & \theta_{0} + 3 & \textbf{(18)}\\
    			\zeta_{19} & -\theta_{0} + 2 & \theta_{0} + 2 & \textbf{(21)}\\
    			\zeta_{20} & -\theta_{0} + 3 & \theta_{0} + 1 & \textbf{(23)}\\
    			\zeta_{21} & -\theta_{0} + 4 & \theta_{0} & \textbf{(24)}\\
    			\zeta_{22}& \theta_{0} & -\theta_{0} + 4 & \textbf{(27)}\\
    			\zeta_{23} & \theta_{0} + 1 & -\theta_{0} + 3& \textbf{(29)} \\
    			\zeta_{24} & \theta_{0} + 2 & -\theta_{0} + 2 & \textbf{(30)}
    			\end{dmatrix}
    		\end{align*}
    		At this point, no risk of mistakes...
    		
    		\TeX\; on the left and Sage on the right
    		\begin{align*}
    		    \begin{dmatrix}
    		    \zeta_{13} & 0 & 4 & \textbf{(3)}\\
    		    \zeta_{14} & 1 & 3 & \textbf{(7)}\\
    		    \zeta_{15} & 2 & 2 & \textbf{(10)}\\
    		    \zeta_{16} & 3 & 1 & \textbf{(12)}\\
    		    \zeta_{17} & 4 & 0 & \textbf{(13)}\\
    		    \zeta_{18} & -\theta_{0} + 1 & \theta_{0} + 3 & \textbf{(18)}\\
    		    \zeta_{19} & -\theta_{0} + 2 & \theta_{0} + 2 & \textbf{(21)}\\
    		    \zeta_{20} & -\theta_{0} + 3 & \theta_{0} + 1 & \textbf{(23)}\\
    		    \zeta_{21} & -\theta_{0} + 4 & \theta_{0} & \textbf{(24)}\\
    		    \zeta_{22}& \theta_{0} & -\theta_{0} + 4 & \textbf{(27)}\\
    		    \zeta_{23} & \theta_{0} + 1 & -\theta_{0} + 3& \textbf{(29)} \\
    		    \zeta_{24} & \theta_{0} + 2 & -\theta_{0} + 2 & \textbf{(30)}
    		    \end{dmatrix}
    			\begin{dmatrix}
    			\zeta_{13} & 0 & 4 \\
    			\zeta_{14} & 1 & 3 \\
    			\zeta_{15} & 2 & 2 \\
    			\zeta_{16} & 3 & 1 \\
    			\zeta_{17} & 4 & 0 \\
    			\zeta_{18} & -\theta_{0} + 1 & \theta_{0} + 3 \\
    			\zeta_{19} & -\theta_{0} + 2 & \theta_{0} + 2 \\
    			\zeta_{20} & -\theta_{0} + 3 & \theta_{0} + 1 \\
    			\zeta_{21} & -\theta_{0} + 4 & \theta_{0} \\
    			\zeta_{22} & \theta_{0} & -\theta_{0} + 4 \\
    			\zeta_{23} & \theta_{0} + 1 & -\theta_{0} + 3 \\
    			\zeta_{24} & \theta_{0} + 2 & -\theta_{0} + 2
    			\end{dmatrix}
    		\end{align*}
    		
    		Finally, we have
    		\begin{align*}
    			\s{\alpha}{\phi}=\begin{dmatrix}
    			\frac{{\left(2 \, \theta_{0} - 1\right)} \zeta_{2}}{2 \, {\left(\theta_{0} + 1\right)} \theta_{0}} & 2 & 0 \\
    			\zeta_{7} & 3 & 0 \\
    			\alpha_{2} {\left(\theta_{0} - 2\right)} & -\theta_{0} + 1 & \theta_{0} + 1 \\
    			\frac{{\left(\theta_{0} - 2\right)} \zeta_{4}}{{\left(\theta_{0} + 2\right)} \theta_{0}} & -\theta_{0} + 1 & \theta_{0} + 2 \\
    			\frac{{\left(\theta_{0} - 3\right)} \zeta_{3}}{2 \, {\left(\theta_{0} + 1\right)} \theta_{0}} & -\theta_{0} + 2 & \theta_{0} + 1 \\
    			{\left(\theta_{0} + 1\right)} \overline{\alpha_{2}} & \theta_{0} & -\theta_{0} + 2 \\
    			\frac{\overline{\zeta_{3}}}{2 \, \theta_{0}} & \theta_{0} & -\theta_{0} + 3 \\
    			\frac{\overline{\zeta_{4}}}{\theta_{0}} & \theta_{0} + 1 & -\theta_{0} + 2 
    			\end{dmatrix}
    			\begin{dmatrix}
    			\zeta_{13} & 0 & 4 \\
    			\zeta_{14} & 1 & 3 \\
    			\zeta_{15} & 2 & 2 \\
    			\zeta_{16} & 3 & 1 \\
    			\zeta_{17} & 4 & 0 \\
    			\zeta_{18} & -\theta_{0} + 1 & \theta_{0} + 3 \\
    			\zeta_{19} & -\theta_{0} + 2 & \theta_{0} + 2 \\
    			\zeta_{20} & -\theta_{0} + 3 & \theta_{0} + 1 \\
    			\zeta_{21} & -\theta_{0} + 4 & \theta_{0} \\
    			\zeta_{22} & \theta_{0} & -\theta_{0} + 4 \\
    			\zeta_{23} & \theta_{0} + 1 & -\theta_{0} + 3 \\
    			\zeta_{24} & \theta_{0} + 2 & -\theta_{0} + 2
    			\end{dmatrix}
    		\end{align*}
    		
    	    \section{Next order development of $\h_0$}
    	    We first have
    	    \begin{align*}
    	    	g^{-1}\otimes\h_0=\begin{dmatrix}
    	    	2 & A_{1} & 0 & -\theta_{0} + 1 \\
    	    	-4 \, {\left| A_{1} \right|}^{2} & A_{0} & 0 & -\theta_{0} + 2 \\
    	    	\frac{1}{4} & \overline{B_{1}} & 0 & -\theta_{0} + 3 \\
    	    	-2 \, \overline{\alpha_{5}} & A_{0} & 0 & -\theta_{0} + 3 \\
    	    	-2 \, \overline{\alpha_{1}} & A_{1} & 0 & -\theta_{0} + 3 \\
    	    	\frac{1}{6} & \overline{B_{3}} & 0 & -\theta_{0} + 4 \\
    	    	-2 \, \overline{\alpha_{3}} & A_{1} & 0 & -\theta_{0} + 4 \\
    	    	-\frac{1}{6} \, {\left| A_{1} \right|}^{2} & \overline{C_{1}} & 0 & -\theta_{0} + 4 \\
    	    	8 \, {\left| A_{1} \right|}^{2} \overline{\alpha_{1}} - 2 \, \overline{\alpha_{6}} & A_{0} & 0 & -\theta_{0} + 4 \\
    	    	\frac{1}{8} & \overline{B_{6}} & 0 & -\theta_{0} + 5 \\
    	    	-\frac{1}{8} \, \overline{\alpha_{5}} & \overline{C_{1}} & 0 & -\theta_{0} + 5 \\
    	    	-\frac{1}{8} \, \overline{\alpha_{1}} & \overline{B_{1}} & 0 & -\theta_{0} + 5 \\
    	    	2 \, \overline{\alpha_{1}}^{2} - 2 \, \overline{\alpha_{4}} & A_{1} & 0 & -\theta_{0} + 5 \\
    	    	-\frac{1}{12} \, {\left| A_{1} \right|}^{2} & \overline{C_{2}} & 0 & -\theta_{0} + 5 \\
    	    	8 \, {\left| A_{1} \right|}^{2} \overline{\alpha_{3}} + 4 \, \overline{\alpha_{1}} \overline{\alpha_{5}} - 2 \, \alpha_{15} & A_{0} & 0 & -\theta_{0} + 5 \\
    	    	4 & A_{2} & 1 & -\theta_{0} + 1 \\
    	    	-4 \, \alpha_{1} & A_{0} & 1 & -\theta_{0} + 1 \\
    	    	-4 \, \alpha_{5} & A_{0} & 1 & -\theta_{0} + 2 \\
    	    	-8 \, {\left| A_{1} \right|}^{2} & A_{1} & 1 & -\theta_{0} + 2 \\
    	    	\frac{1}{2} & \overline{B_{2}} & 1 & -\theta_{0} + 3 \\
    	    	-4 \, \overline{\alpha_{5}} & A_{1} & 1 & -\theta_{0} + 3 \\
    	    	-4 \, \overline{\alpha_{1}} & A_{2} & 1 & -\theta_{0} + 3 \\
    	    	16 \, {\left| A_{1} \right|}^{4} + 8 \, \alpha_{1} \overline{\alpha_{1}} - 4 \, \beta & A_{0} & 1 & -\theta_{0} + 3 
    	    	\end{dmatrix}
    	    	\end{align*}
    	    	\begin{align*}
    	    	\begin{dmatrix}
    	    	\frac{1}{3} & \overline{B_{5}} & 1 & -\theta_{0} + 4 \\
    	    	-4 \, \overline{\alpha_{3}} & A_{2} & 1 & -\theta_{0} + 4 \\
    	    	-\frac{1}{6} \, \alpha_{5} & \overline{C_{1}} & 1 & -\theta_{0} + 4 \\
    	    	-\frac{1}{3} \, {\left| A_{1} \right|}^{2} & \overline{B_{1}} & 1 & -\theta_{0} + 4 \\
    	    	16 \, {\left| A_{1} \right|}^{2} \overline{\alpha_{5}} + 8 \, \alpha_{5} \overline{\alpha_{1}} + 8 \, \alpha_{1} \overline{\alpha_{3}} - 4 \, \alpha_{16} & A_{0} & 1 & -\theta_{0} + 4 \\
    	    	16 \, {\left| A_{1} \right|}^{2} \overline{\alpha_{1}} - 4 \, \overline{\alpha_{6}} & A_{1} & 1 & -\theta_{0} + 4 \\
    	    	6 & A_{3} & 2 & -\theta_{0} + 1 \\
    	    	-6 \, \alpha_{3} & A_{0} & 2 & -\theta_{0} + 1 \\
    	    	-6 \, \alpha_{1} & A_{1} & 2 & -\theta_{0} + 1 \\
    	    	-6 \, \alpha_{5} & A_{1} & 2 & -\theta_{0} + 2 \\
    	    	-12 \, {\left| A_{1} \right|}^{2} & A_{2} & 2 & -\theta_{0} + 2 \\
    	    	24 \, \alpha_{1} {\left| A_{1} \right|}^{2} - 6 \, \alpha_{6} & A_{0} & 2 & -\theta_{0} + 2 \\
    	    	\frac{3}{4} & \overline{B_{4}} & 2 & -\theta_{0} + 3 \\
    	    	\frac{3 \, {\left(\theta_{0} \overline{\alpha_{2}} - 2 \, \overline{\alpha_{2}}\right)}}{4 \, \theta_{0}} & C_{1} & 2 & -\theta_{0} + 3 \\
    	    	-6 \, \overline{\alpha_{5}} & A_{2} & 2 & -\theta_{0} + 3 \\
    	    	-6 \, \overline{\alpha_{1}} & A_{3} & 2 & -\theta_{0} + 3 \\
    	    	24 \, \alpha_{5} {\left| A_{1} \right|}^{2} + 12 \, \alpha_{3} \overline{\alpha_{1}} + 12 \, \alpha_{1} \overline{\alpha_{5}} - 6 \, \overline{\alpha_{16}} & A_{0} & 2 & -\theta_{0} + 3 \\
    	    	24 \, {\left| A_{1} \right|}^{4} + 12 \, \alpha_{1} \overline{\alpha_{1}} - 6 \, \beta & A_{1} & 2 & -\theta_{0} + 3 \\
    	    	8 & A_{4} & 3 & -\theta_{0} + 1 \\
    	    	-8 \, \alpha_{3} & A_{1} & 3 & -\theta_{0} + 1 \\
    	    	-8 \, \alpha_{1} & A_{2} & 3 & -\theta_{0} + 1 \\
    	    	8 \, \alpha_{1}^{2} - 8 \, \alpha_{4} & A_{0} & 3 & -\theta_{0} + 1 \\
    	    	-8 \, \alpha_{5} & A_{2} & 3 & -\theta_{0} + 2 \\
    	    	-16 \, {\left| A_{1} \right|}^{2} & A_{3} & 3 & -\theta_{0} + 2 \\
    	    	32 \, \alpha_{3} {\left| A_{1} \right|}^{2} + 16 \, \alpha_{1} \alpha_{5} - 8 \, \overline{\alpha_{15}} & A_{0} & 3 & -\theta_{0} + 2 \\
    	    	32 \, \alpha_{1} {\left| A_{1} \right|}^{2} - 8 \, \alpha_{6} & A_{1} & 3 & -\theta_{0} + 2 
    	    	\end{dmatrix}
    	    	\end{align*}
    	    	\begin{align*}
    	    	\begin{dmatrix}
    	    	10 & A_{5} & 4 & -\theta_{0} + 1 \\
    	    	-10 \, \alpha_{3} & A_{2} & 4 & -\theta_{0} + 1 \\
    	    	-10 \, \alpha_{1} & A_{3} & 4 & -\theta_{0} + 1 \\
    	    	20 \, \alpha_{1} \alpha_{3} - 10 \, \overline{\alpha_{14}} & A_{0} & 4 & -\theta_{0} + 1 \\
    	    	10 \, \alpha_{1}^{2} - 10 \, \alpha_{4} & A_{1} & 4 & -\theta_{0} + 1 \\
    	    	-\frac{\theta_{0} - 2}{2 \, \theta_{0}} & E_{1} & -2 \, \theta_{0} + 3 & \theta_{0} + 1 \\
    	    	-\frac{\theta_{0} - 2}{2 \, \theta_{0} + 1} & E_{2} & -2 \, \theta_{0} + 3 & \theta_{0} + 2 \\
    	    	\frac{\alpha_{2} \theta_{0}^{2} - \alpha_{2} \theta_{0} - 2 \, \alpha_{2}}{2 \, \theta_{0}^{2} + \theta_{0}} & C_{1} & -2 \, \theta_{0} + 3 & \theta_{0} + 2 \\
    	    	-\frac{2 \, \theta_{0} - 5}{4 \, \theta_{0}} & E_{3} & -2 \, \theta_{0} + 4 & \theta_{0} + 1 \\
    	    	-\frac{\theta_{0} - 2}{2 \, \theta_{0}} & C_{1} & -\theta_{0} + 1 & 1 \\
    	    	-\frac{\theta_{0} - 2}{2 \, {\left(\theta_{0} + 1\right)}} & B_{1} & -\theta_{0} + 1 & 2 \\
    	    	2 \, \alpha_{2} \theta_{0} - 4 \, \alpha_{2} & A_{0} & -\theta_{0} + 1 & 2 \\
    	    	-\frac{\theta_{0} - 2}{2 \, {\left(\theta_{0} + 2\right)}} & B_{2} & -\theta_{0} + 1 & 3 \\
    	    	\frac{\theta_{0} \overline{\alpha_{1}} - 2 \, \overline{\alpha_{1}}}{\theta_{0}^{2} + 2 \, \theta_{0}} & C_{1} & -\theta_{0} + 1 & 3 \\
    	    	2 \, \alpha_{9} \theta_{0} - 4 \, \alpha_{9} & A_{0} & -\theta_{0} + 1 & 3 \\
    	    	-\frac{\theta_{0} - 2}{2 \, {\left(\theta_{0} + 3\right)}} & B_{4} & -\theta_{0} + 1 & 4 \\
    	    	\frac{\alpha_{2} \theta_{0}^{2} - \alpha_{2} \theta_{0} - 2 \, \alpha_{2}}{4 \, {\left(\theta_{0} + 3\right)}} & \overline{C_{1}} & -\theta_{0} + 1 & 4 \\
    	    	\frac{3 \, {\left(\theta_{0} \overline{\alpha_{3}} - 2 \, \overline{\alpha_{3}}\right)}}{2 \, {\left(\theta_{0}^{2} + 3 \, \theta_{0}\right)}} & C_{1} & -\theta_{0} + 1 & 4 \\
    	    	\frac{\theta_{0} \overline{\alpha_{1}} - 2 \, \overline{\alpha_{1}}}{\theta_{0}^{2} + 4 \, \theta_{0} + 3} & B_{1} & -\theta_{0} + 1 & 4 \\
    	    	-2 \, {\left(2 \, \alpha_{2} \overline{\alpha_{1}} - \alpha_{10}\right)} \theta_{0} + 8 \, \alpha_{2} \overline{\alpha_{1}} - 4 \, \alpha_{10} & A_{0} & -\theta_{0} + 1 & 4 \\
    	    	-\frac{\theta_{0} - 3}{2 \, \theta_{0}} & C_{2} & -\theta_{0} + 2 & 1 \\
    	    	-\frac{\theta_{0} - 3}{2 \, {\left(\theta_{0} + 1\right)}} & B_{3} & -\theta_{0} + 2 & 2 \\
    	    	\frac{{\left(\theta_{0} - 3\right)} {\left| A_{1} \right|}^{2}}{\theta_{0}^{2} + \theta_{0}} & C_{1} & -\theta_{0} + 2 & 2 \\
    	    	2 \, \alpha_{8} \theta_{0} - 6 \, \alpha_{8} & A_{0} & -\theta_{0} + 2 & 2 \\
    	    	2 \, \alpha_{2} \theta_{0} - 6 \, \alpha_{2} & A_{1} & -\theta_{0} + 2 & 2 
    	    	\end{dmatrix}
    	    	\end{align*}
    	    	\begin{align*}
    	    	\begin{dmatrix}
    	    	-\frac{\theta_{0} - 3}{2 \, {\left(\theta_{0} + 2\right)}} & B_{5} & -\theta_{0} + 2 & 3 \\
    	    	\frac{\theta_{0} \overline{\alpha_{5}} - 3 \, \overline{\alpha_{5}}}{\theta_{0}^{2} + 2 \, \theta_{0}} & C_{1} & -\theta_{0} + 2 & 3 \\
    	    	\frac{\theta_{0} \overline{\alpha_{1}} - 3 \, \overline{\alpha_{1}}}{\theta_{0}^{2} + 2 \, \theta_{0}} & C_{2} & -\theta_{0} + 2 & 3 \\
    	    	\frac{{\left(\theta_{0} - 3\right)} {\left| A_{1} \right|}^{2}}{\theta_{0}^{2} + 3 \, \theta_{0} + 2} & B_{1} & -\theta_{0} + 2 & 3 \\
    	    	2 \, \alpha_{9} \theta_{0} - 6 \, \alpha_{9} & A_{1} & -\theta_{0} + 2 & 3 \\
    	    	-8 \, {\left(\alpha_{2} \theta_{0} - 3 \, \alpha_{2}\right)} {\left| A_{1} \right|}^{2} + 2 \, \alpha_{11} \theta_{0} - 6 \, \alpha_{11} & A_{0} & -\theta_{0} + 2 & 3 \\
    	    	-\frac{\theta_{0} - 4}{2 \, \theta_{0}} & C_{3} & -\theta_{0} + 3 & 1 \\
    	    	2 \, \alpha_{7} \theta_{0} - 8 \, \alpha_{7} & A_{0} & -\theta_{0} + 3 & 1 \\
    	    	-\frac{\theta_{0} - 4}{2 \, {\left(\theta_{0} + 1\right)}} & B_{6} & -\theta_{0} + 3 & 2 \\
    	    	\frac{\alpha_{5} \theta_{0} - 4 \, \alpha_{5}}{2 \, {\left(\theta_{0}^{2} + \theta_{0}\right)}} & C_{1} & -\theta_{0} + 3 & 2 \\
    	    	\frac{{\left(\theta_{0} - 4\right)} {\left| A_{1} \right|}^{2}}{\theta_{0}^{2} + \theta_{0}} & C_{2} & -\theta_{0} + 3 & 2 \\
    	    	2 \, \alpha_{8} \theta_{0} - 8 \, \alpha_{8} & A_{1} & -\theta_{0} + 3 & 2 \\
    	    	2 \, \alpha_{2} \theta_{0} - 8 \, \alpha_{2} & A_{2} & -\theta_{0} + 3 & 2 \\
    	    	16 \, \alpha_{1} \alpha_{2} - 2 \, {\left(2 \, \alpha_{1} \alpha_{2} - \alpha_{12}\right)} \theta_{0} - 8 \, \alpha_{12} & A_{0} & -\theta_{0} + 3 & 2 \\
    	    	-\frac{\theta_{0} - 5}{2 \, \theta_{0}} & C_{4} & -\theta_{0} + 4 & 1 \\
    	    	2 \, \alpha_{7} \theta_{0} - 10 \, \alpha_{7} & A_{1} & -\theta_{0} + 4 & 1 \\
    	    	2 \, \alpha_{13} \theta_{0} - 10 \, \alpha_{13} & A_{0} & -\theta_{0} + 4 & 1 \\
    	    	-2 \, \theta_{0} \overline{\alpha_{2}} - 2 \, \overline{\alpha_{2}} & A_{0} & \theta_{0} & -2 \, \theta_{0} + 3 \\
    	    	-2 \, \theta_{0} \overline{\alpha_{8}} - 2 \, \overline{\alpha_{8}} & A_{0} & \theta_{0} & -2 \, \theta_{0} + 4 \\
    	    	-\frac{\theta_{0} + 1}{2 \, {\left(\theta_{0} - 4\right)}} & \overline{E_{2}} & \theta_{0} & -2 \, \theta_{0} + 5 \\
    	    	-\frac{\theta_{0}^{2} \overline{\alpha_{2}} - \theta_{0} \overline{\alpha_{2}} - 2 \, \overline{\alpha_{2}}}{4 \, {\left(\theta_{0} - 4\right)}} & \overline{C_{1}} & \theta_{0} & -2 \, \theta_{0} + 5 \\
    	    	-2 \, \theta_{0} \overline{\alpha_{7}} - 2 \, \overline{\alpha_{7}} & A_{1} & \theta_{0} & -2 \, \theta_{0} + 5 \\
    	    	2 \, {\left(2 \, \overline{\alpha_{1}} \overline{\alpha_{2}} - \overline{\alpha_{12}}\right)} \theta_{0} + 4 \, \overline{\alpha_{1}} \overline{\alpha_{2}} - 2 \, \overline{\alpha_{12}} & A_{0} & \theta_{0} & -2 \, \theta_{0} + 5 
    	    	\end{dmatrix}
    	    	\end{align*}
    	    	\begin{align*}
    	    	\begin{dmatrix}
    	    	-\frac{\theta_{0}}{2 \, {\left(\theta_{0} - 4\right)}} & \overline{E_{1}} & \theta_{0} - 1 & -2 \, \theta_{0} + 5 \\
    	    	-2 \, \theta_{0} \overline{\alpha_{7}} & A_{0} & \theta_{0} - 1 & -2 \, \theta_{0} + 5 \\
    	    	-\frac{\theta_{0}}{2 \, {\left(\theta_{0} - 5\right)}} & \overline{E_{3}} & \theta_{0} - 1 & -2 \, \theta_{0} + 6 \\
    	    	-2 \, \theta_{0} \overline{\alpha_{13}} & A_{0} & \theta_{0} - 1 & -2 \, \theta_{0} + 6 \\
    	    	-2 \, \theta_{0} \overline{\alpha_{9}} - 4 \, \overline{\alpha_{9}} & A_{0} & \theta_{0} + 1 & -2 \, \theta_{0} + 3 \\
    	    	-2 \, \theta_{0} \overline{\alpha_{2}} - 4 \, \overline{\alpha_{2}} & A_{1} & \theta_{0} + 1 & -2 \, \theta_{0} + 3 \\
    	    	-2 \, \theta_{0} \overline{\alpha_{8}} - 4 \, \overline{\alpha_{8}} & A_{1} & \theta_{0} + 1 & -2 \, \theta_{0} + 4 \\
    	    	8 \, {\left(\theta_{0} \overline{\alpha_{2}} + 2 \, \overline{\alpha_{2}}\right)} {\left| A_{1} \right|}^{2} - 2 \, \theta_{0} \overline{\alpha_{11}} - 4 \, \overline{\alpha_{11}} & A_{0} & \theta_{0} + 1 & -2 \, \theta_{0} + 4 \\
    	    	-2 \, \theta_{0} \overline{\alpha_{9}} - 6 \, \overline{\alpha_{9}} & A_{1} & \theta_{0} + 2 & -2 \, \theta_{0} + 3 \\
    	    	-2 \, \theta_{0} \overline{\alpha_{2}} - 6 \, \overline{\alpha_{2}} & A_{2} & \theta_{0} + 2 & -2 \, \theta_{0} + 3 \\
    	    	2 \, {\left(2 \, \alpha_{1} \overline{\alpha_{2}} - \overline{\alpha_{10}}\right)} \theta_{0} + 12 \, \alpha_{1} \overline{\alpha_{2}} - 6 \, \overline{\alpha_{10}} & A_{0} & \theta_{0} + 2 & -2 \, \theta_{0} + 3
    	    	\end{dmatrix}
    	    \end{align*}
    	    
    		\small
    		\begin{align*}
    			&g^{-1}\otimes\s{\phi}{\h_0}=\\
    			&\begin{dmatrix}
    			\frac{2 \, A_{1} \overline{A_{0}}}{\theta_{0}} & 0 & 1 \\
    			-\frac{2 \, {\left(2 \, A_{0} {\left(\theta_{0} + 1\right)} {\left| A_{1} \right|}^{2} \overline{A_{0}} - A_{1} \theta_{0} \overline{A_{1}}\right)}}{\theta_{0}^{2} + \theta_{0}} & 0 & 2 \\
    			\nu_1 & 0 & 3 \\
    			\nu_2 & 0 & 4 \\
    			\nu_3 & 0 & 5 \\
    			-\frac{8 \, A_{0} \alpha_{1} \theta_{0} \overline{A_{0}} + A_{0} C_{1} {\left(\theta_{0} - 2\right)} - 8 \, A_{2} \theta_{0} \overline{A_{0}}}{2 \, \theta_{0}^{2}} & 1 & 1 \\
    			\nu_4 & 1 & 2 \\
    			\nu_5 & 1 & 3 \\
    			\nu_6 & 1 & 4 \\
    			\nu_7 & 2 & 1 \\
    			\nu_8 & 2 & 2 \\
    			\nu_9 & 2 & 3 \\
    			\nu_{10} & 3 & 1 \\
    			\nu_{11} & 3 & 2 \\
    			\nu_{12} & 4 & 1 \\
    			-\frac{E_{1} {\left(\theta_{0} - 2\right)} \overline{A_{0}}}{2 \, \theta_{0}^{2}} & -2 \, \theta_{0} + 3 & 2 \, \theta_{0} + 1 \\
    			\frac{2 \, {\left(\alpha_{2} \theta_{0}^{3} - 3 \, \alpha_{2} \theta_{0} - 2 \, \alpha_{2}\right)} C_{1} \overline{A_{0}} - 2 \, {\left(\theta_{0}^{3} - \theta_{0}^{2} - 2 \, \theta_{0}\right)} E_{2} \overline{A_{0}} - {\left(2 \, \theta_{0}^{3} - 3 \, \theta_{0}^{2} - 2 \, \theta_{0}\right)} E_{1} \overline{A_{1}}}{2 \, {\left(2 \, \theta_{0}^{4} + 3 \, \theta_{0}^{3} + \theta_{0}^{2}\right)}} & -2 \, \theta_{0} + 3 & 2 \, \theta_{0} + 2 \\
    			-\frac{E_{3} {\left(2 \, \theta_{0} - 5\right)} \overline{A_{0}}}{4 \, \theta_{0}^{2}} & -2 \, \theta_{0} + 4 & 2 \, \theta_{0} + 1 \\
    			-\frac{C_{1} {\left(\theta_{0} - 2\right)} \overline{A_{0}}}{2 \, \theta_{0}^{2}} & -\theta_{0} + 1 & \theta_{0} + 1 \\
    			\frac{4 \, {\left(\alpha_{2} \theta_{0}^{2} - \alpha_{2} \theta_{0} - 2 \, \alpha_{2}\right)} A_{0} \overline{A_{0}} - B_{1} {\left(\theta_{0} - 2\right)} \overline{A_{0}} - C_{1} {\left(\theta_{0} - 2\right)} \overline{A_{1}}}{2 \, {\left(\theta_{0}^{2} + \theta_{0}\right)}} & -\theta_{0} + 1 & \theta_{0} + 2 \\
    			\nu_{13} & -\theta_{0} + 1 & \theta_{0} + 3 \\
    			\nu_{14} & -\theta_{0} + 1 & \theta_{0} + 4 \\
    			-\frac{C_{2} {\left(\theta_{0} - 3\right)} \overline{A_{0}}}{2 \, \theta_{0}^{2}} & -\theta_{0} + 2 & \theta_{0} + 1 \\
    			\nu_{15} & -\theta_{0} + 2 & \theta_{0} + 2 \\
    			\nu_{16} & -\theta_{0} + 2 & \theta_{0} + 3 \\
    			-\frac{C_{1}^{2} {\left(\theta_{0} - 2\right)} + 8 \, A_{0} E_{1} {\left(\theta_{0} - 2\right)} - 32 \, {\left(\alpha_{7} \theta_{0}^{2} - 4 \, \alpha_{7} \theta_{0}\right)} A_{0} \overline{A_{0}} + 8 \, C_{3} {\left(\theta_{0} - 4\right)} \overline{A_{0}}}{16 \, \theta_{0}^{2}} & -\theta_{0} + 3 & \theta_{0} + 1 
    			    			\end{dmatrix}
    			\end{align*}
    			\begin{align*}
    			\begin{dmatrix}
    			\nu_{17} & -\theta_{0} + 3 & \theta_{0} + 2 \\
    			\nu_{18} & -\theta_{0} + 4 & \theta_{0} + 1 \\
    			\frac{2 \, A_{0} A_{1}}{\theta_{0}} & \theta_{0} & -\theta_{0} + 1 \\
    			-\frac{4 \, A_{0}^{2} {\left| A_{1} \right|}^{2}}{\theta_{0}} & \theta_{0} & -\theta_{0} + 2 \\
    			-\frac{8 \, {\left(\theta_{0} \overline{\alpha_{2}} + \overline{\alpha_{2}}\right)} A_{0} \overline{A_{0}} + 8 \, A_{0} A_{1} \overline{\alpha_{1}} + 8 \, A_{0}^{2} \overline{\alpha_{5}} - A_{0} \overline{B_{1}} - A_{1} \overline{C_{1}}}{4 \, \theta_{0}} & \theta_{0} & -\theta_{0} + 3 \\
    			\nu_{19} & \theta_{0} & -\theta_{0} + 4 \\
    			\nu_{20} & \theta_{0} & -\theta_{0} + 5 \\
    			-\frac{4 \, {\left(\theta_{0} \overline{\alpha_{7}} - 4 \, \overline{\alpha_{7}}\right)} A_{0} \overline{A_{0}} + \overline{A_{0}} \overline{E_{1}}}{2 \, {\left(\theta_{0} - 4\right)}} & \theta_{0} - 1 & -\theta_{0} + 5 \\
    			\nu_{21} & \theta_{0} - 1 & -\theta_{0} + 6 \\
    			-\frac{2 \, {\left(2 \, {\left(\alpha_{1} \theta_{0} + \alpha_{1}\right)} A_{0}^{2} - 2 \, A_{0} A_{2} {\left(\theta_{0} + 1\right)} - A_{1}^{2} \theta_{0}\right)}}{\theta_{0}^{2} + \theta_{0}} & \theta_{0} + 1 & -\theta_{0} + 1 \\
    			-\frac{4 \, {\left(A_{0} A_{1} {\left(3 \, \theta_{0} + 2\right)} {\left| A_{1} \right|}^{2} + {\left(\alpha_{5} \theta_{0} + \alpha_{5}\right)} A_{0}^{2}\right)}}{\theta_{0}^{2} + \theta_{0}} & \theta_{0} + 1 & -\theta_{0} + 2 \\
    			\nu_{22} & \theta_{0} + 1 & -\theta_{0} + 3 \\
    			\nu_{23} & \theta_{0} + 1 & -\theta_{0} + 4 \\
    			\nu_{24} & \theta_{0} + 2 & -\theta_{0} + 1 \\
    			\nu_{25} & \theta_{0} + 2 & -\theta_{0} + 2 \\
    			\nu_{26} & \theta_{0} + 2 & -\theta_{0} + 3 \\
    			\nu_{27} & \theta_{0} + 3 & -\theta_{0} + 1 \\
    			\nu_{28} & \theta_{0} + 3 & -\theta_{0} + 2 \\
    			\nu_{29} & \theta_{0} + 4 & -\theta_{0} + 1 \\
    			-\frac{4 \, {\left(\theta_{0} \overline{\alpha_{7}} - 4 \, \overline{\alpha_{7}}\right)} A_{0}^{2} + A_{0} \overline{E_{1}}}{2 \, {\left(\theta_{0} - 4\right)}} & 2 \, \theta_{0} - 1 & -2 \, \theta_{0} + 5 \\
    			-\frac{4 \, {\left(\theta_{0} \overline{\alpha_{13}} - 5 \, \overline{\alpha_{13}}\right)} A_{0}^{2} + A_{0} \overline{E_{3}}}{2 \, {\left(\theta_{0} - 5\right)}} & 2 \, \theta_{0} - 1 & -2 \, \theta_{0} + 6 \\
    			-\frac{2 \, {\left({\left(\theta_{0} \overline{\alpha_{9}} + 2 \, \overline{\alpha_{9}}\right)} A_{0}^{2} + 2 \, {\left(\theta_{0} \overline{\alpha_{2}} + \overline{\alpha_{2}}\right)} A_{0} A_{1}\right)}}{\theta_{0}} & 2 \, \theta_{0} + 1 & -2 \, \theta_{0} + 3 \\
    			\frac{2 \, {\left(4 \, {\left(\theta_{0} \overline{\alpha_{2}} + 2 \, \overline{\alpha_{2}}\right)} A_{0}^{2} {\left| A_{1} \right|}^{2} - {\left(\theta_{0} \overline{\alpha_{11}} + 2 \, \overline{\alpha_{11}}\right)} A_{0}^{2} - 2 \, {\left(\theta_{0} \overline{\alpha_{8}} + \overline{\alpha_{8}}\right)} A_{0} A_{1}\right)}}{\theta_{0}} & 2 \, \theta_{0} + 1 & -2 \, \theta_{0} + 4 \\
    			\nu_{30} & 2 \, \theta_{0} + 2 & -2 \, \theta_{0} + 3 \\
    			-\frac{2 \, {\left(\theta_{0} \overline{\alpha_{2}} + \overline{\alpha_{2}}\right)} A_{0}^{2}}{\theta_{0}} & 2 \, \theta_{0} & -2 \, \theta_{0} + 3 \\
    			-\frac{2 \, {\left(\theta_{0} \overline{\alpha_{8}} + \overline{\alpha_{8}}\right)} A_{0}^{2}}{\theta_{0}} & 2 \, \theta_{0} & -2 \, \theta_{0} + 4 \\
    			\nu_{31} & 2 \, \theta_{0} & -2 \, \theta_{0} + 5
    			\end{dmatrix}
    		\end{align*}
    		where
            \footnotesize

    	    \normalsize
    	    Now, recall that
    	    \begin{align*}
    	    	g^{-1}\otimes\s{\h_0}{\phi}=
    	    	\begin{dmatrix}
    	    	-\frac{2 \, {\left| A_{1} \right|}^{2}}{{\left(\theta_{0} + 1\right)} \theta_{0}} & 0 & 2 \\
    	    	\mu_{1} & 0 & 3 \\
    	    	\mu_{2} & 0 & 4 \\
    	    	\mu_{3} & 1 & 2 \\
    	    	\mu_{4} & 1 & 3 \\
    	    	\mu_{5} & 2 & 1 \\
    	    	\mu_{6} & 2 & 2 \\
    	    	\mu_{7} & 3 & 1 \\
    	    	\frac{\alpha_{2} {\left(\theta_{0} - 2\right)}}{\theta_{0}} & -\theta_{0} + 1 & \theta_{0} + 2 \\
    	    	\mu_{8} & -\theta_{0} + 1 & \theta_{0} + 3 
    	    	\end{dmatrix}
    	    	\begin{dmatrix}
    	    	\mu_{9} & -\theta_{0} + 2 & \theta_{0} + 2 \\
    	    	\frac{\alpha_{7} {\left(\theta_{0} - 2\right)} {\left(\theta_{0} - 4\right)}}{\theta_{0}^{2}} & -\theta_{0} + 3 & \theta_{0} + 1 \\
    	    	-\frac{{\left(\theta_{0} + 1\right)} \overline{\alpha_{2}}}{\theta_{0}} & \theta_{0} & -\theta_{0} + 3 \\
    	    	\mu_{10} & \theta_{0} & -\theta_{0} + 4 \\
    	    	-\frac{2 \, \zeta_{0}}{{\left(\theta_{0} + 1\right)} \theta_{0}} & \theta_{0} + 1 & -\theta_{0} + 1 \\
    	    	\mu_{11} & \theta_{0} + 1 & -\theta_{0} + 3 \\
    	    	\mu_{12} & \theta_{0} + 2 & -\theta_{0} + 1 \\
    	    	\mu_{13} & \theta_{0} + 2 & -\theta_{0} + 2 \\
    	    	\mu_{14} & \theta_{0} + 3 & -\theta_{0} + 1	
    	    	\end{dmatrix}
    	    \end{align*}
    	    so we need only look for powers of order $5$. They are
    	    \begin{align*}
    	    	\begin{dmatrix}
    	    	\nu_3 & 0 &5\\
    	    	\nu_6& 1 & 4\\
    	    	\nu_9&  2 & 3\\
    	    	\nu_{11} & 3 &2\\
    	    	\nu_{12} & 4 &1\\
    	    	\nu_{14}& -\theta_0+1 & \theta_0+4\\
    	    	\nu_{16}& -\theta_0+2 & \theta_0+3\\
    	    	\nu_{17}& -\theta_0+3 & \theta_0+2\\
    	    	\nu_{18}& -\theta_0+4 & \theta_0+1\\
    	    	\nu_{20}& \theta_0 & -\theta_0+5\\
    	    	\nu_{21} & \theta_0-1 & -\theta_0+6\\
    	    	\nu_{23}& \theta_0+1& -\theta_0+4\\
    	    	\nu_{26}& \theta_0+2& -\theta_0+3\\
    	    	\nu_{28}& \theta_0+3& -\theta_0+2\\
    	    	\nu_{29}& \theta_0+4&-\theta_0+1\\
    	    	\nu_{30}& 2\theta_0+2 &-2\theta_0+3\\
    	    	\nu_{31}& 2\theta_0& -2\theta_0+5 
    	    	\end{dmatrix}
    	    	=\begin{dmatrix}
    	    	\mu_{15} & 0 &5\\
    	    	\mu_{16}& 1 & 4\\
    	    	\mu_{17}&  2 & 3\\
    	    	\mu_{18} & 3 &2\\
    	    	\mu_{19} & 4 &1\\
    	    	\mu_{20}& -\theta_0+1 & \theta_0+4\\
    	    	\mu_{21}& -\theta_0+2 & \theta_0+3\\
    	    	\mu_{22}& -\theta_0+3 & \theta_0+2\\
    	    	\mu_{23}& -\theta_0+4 & \theta_0+1\\
    	    	\mu_{24}& \theta_0 & -\theta_0+5\\
    	    	\mu_{25} & \theta_0-1 & -\theta_0+6\\
    	    	\mu_{26}& \theta_0+1& -\theta_0+4\\
    	    	\mu_{27}& \theta_0+2& -\theta_0+3\\
    	    	\mu_{28}& \theta_0+3& -\theta_0+2\\
    	    	\mu_{29}& \theta_0+4&-\theta_0+1\\
    	    	\mu_{30}& 2\theta_0+2 &-2\theta_0+3\\
    	    	\mu_{31}& 2\theta_0& -2\theta_0+5 
    	    	\end{dmatrix}
    	    \end{align*}
    	    to preserve notations we write the new terms with $\ens{\mu_{j}}_{15\leq j \leq 31}$.
    	    In Sage, we have
    	    \begin{align}
    	    \begin{dmatrix}
    	    \mu_{15} & 0 &5\\
    	    \mu_{16}& 1 & 4\\
    	    \mu_{17}&  2 & 3\\
    	    \mu_{18} & 3 &2\\
    	    \mu_{19} & 4 &1\\
    	    \mu_{20}& -\theta_0+1 & \theta_0+4\\
    	    \mu_{21}& -\theta_0+2 & \theta_0+3\\
    	    \mu_{22}& -\theta_0+3 & \theta_0+2\\
    	    \mu_{23}& -\theta_0+4 & \theta_0+1\\
    	    \mu_{24}& \theta_0 & -\theta_0+5\\
    	    \mu_{25} & \theta_0-1 & -\theta_0+6\\
    	    \mu_{26}& \theta_0+1& -\theta_0+4\\
    	    \mu_{27}& \theta_0+2& -\theta_0+3\\
    	    \mu_{28}& \theta_0+3& -\theta_0+2\\
    	    \mu_{29}& \theta_0+4&-\theta_0+1\\
    	    \mu_{30}& 2\theta_0+2 &-2\theta_0+3\\
    	    \mu_{31}& 2\theta_0& -2\theta_0+5 
    	    \end{dmatrix}
    	    	\begin{dmatrix}
    	    	\mu_{15} & 0 & 5 \\
    	    	\mu_{16} & 1 & 4 \\
    	    	\mu_{17} & 2 & 3 \\
    	    	\mu_{18} & 3 & 2 \\
    	    	\mu_{19} & 4 & 1 \\
    	    	\mu_{20} & -\theta_{0} + 1 & \theta_{0} + 4 \\
    	    	\mu_{21} & -\theta_{0} + 2 & \theta_{0} + 3 \\
    	    	\mu_{22} & -\theta_{0} + 3 & \theta_{0} + 2 \\
    	    	\mu_{23} & -\theta_{0} + 4 & \theta_{0} + 1 \\
    	    	\mu_{24} & \theta_{0} & -\theta_{0} + 5 \\
    	    	\mu_{25} & \theta_{0} - 1 & -\theta_{0} + 6 \\
    	    	\mu_{26} & \theta_{0} + 1 & -\theta_{0} + 4 \\
    	    	\mu_{27} & \theta_{0} + 2 & -\theta_{0} + 3 \\
    	    	\mu_{28} & \theta_{0} + 3 & -\theta_{0} + 2 \\
    	    	\mu_{29} & \theta_{0} + 4 & -\theta_{0} + 1 \\
    	    	\mu_{30} & 2 \, \theta_{0} + 2 & -2 \, \theta_{0} + 3 \\
    	    	\mu_{31} & 2 \, \theta_{0} & -2 \, \theta_{0} + 5
    	    	\end{dmatrix}
    	    \end{align}
    	    where
    	    \begin{align*}
    	    	&\mu_{15}=\nu_3=-\frac{1}{24 \, {\left(\theta_{0}^{5} + 10 \, \theta_{0}^{4} + 35 \, \theta_{0}^{3} + 50 \, \theta_{0}^{2} + 24 \, \theta_{0}\right)}}\bigg\{96 \, {\left(\theta_{0}^{4} + 7 \, \theta_{0}^{3} + 14 \, \theta_{0}^{2} + 8 \, \theta_{0}\right)} A_{0} {\left| A_{1} \right|}^{2} \overline{A_{3}}\\
    	    	& + 2 \, {\left(\theta_{0}^{4} + 10 \, \theta_{0}^{3} + 35 \, \theta_{0}^{2} + 50 \, \theta_{0} + 24\right)} {\left| A_{1} \right|}^{2} \overline{A_{0}} \overline{C_{2}} - 192 \, \bigg({\left(\theta_{0}^{4} \overline{\alpha_{3}} + 10 \, \theta_{0}^{3} \overline{\alpha_{3}} + 35 \, \theta_{0}^{2} \overline{\alpha_{3}} + 50 \, \theta_{0} \overline{\alpha_{3}} + 24 \, \overline{\alpha_{3}}\right)} A_{0} \overline{A_{0}}\\
    	    	& + {\left(\theta_{0}^{4} \overline{\alpha_{1}} + 9 \, \theta_{0}^{3} \overline{\alpha_{1}} + 26 \, \theta_{0}^{2} \overline{\alpha_{1}} + 24 \, \theta_{0} \overline{\alpha_{1}}\right)} A_{0} \overline{A_{1}}\bigg) {\left| A_{1} \right|}^{2} - 48 \, \bigg({\left(2 \, \overline{\alpha_{1}} \overline{\alpha_{5}} - \alpha_{15}\right)} \theta_{0}^{4} + 10 \, {\left(2 \, \overline{\alpha_{1}} \overline{\alpha_{5}} - \alpha_{15}\right)} \theta_{0}^{3}\\
    	    	& + 35 \, {\left(2 \, \overline{\alpha_{1}} \overline{\alpha_{5}} - \alpha_{15}\right)} \theta_{0}^{2} + 50 \, {\left(2 \, \overline{\alpha_{1}} \overline{\alpha_{5}} - \alpha_{15}\right)} \theta_{0} + 48 \, \overline{\alpha_{1}} \overline{\alpha_{5}} - 24 \, \alpha_{15}\bigg) A_{0} \overline{A_{0}}\\
    	    	& - 48 \, {\left({\left(\overline{\alpha_{1}}^{2} - \overline{\alpha_{4}}\right)} \theta_{0}^{4} + 10 \, {\left(\overline{\alpha_{1}}^{2} - \overline{\alpha_{4}}\right)} \theta_{0}^{3} + 35 \, {\left(\overline{\alpha_{1}}^{2} - \overline{\alpha_{4}}\right)} \theta_{0}^{2} + 50 \, {\left(\overline{\alpha_{1}}^{2} - \overline{\alpha_{4}}\right)} \theta_{0} + 24 \, \overline{\alpha_{1}}^{2} - 24 \, \overline{\alpha_{4}}\right)} A_{1} \overline{A_{0}}\\
    	    	& - 48 \, {\left(\theta_{0}^{4} + 6 \, \theta_{0}^{3} + 11 \, \theta_{0}^{2} + 6 \, \theta_{0}\right)} A_{1} \overline{A_{4}} - 4 \, {\left(\theta_{0}^{4} + 9 \, \theta_{0}^{3} + 26 \, \theta_{0}^{2} + 24 \, \theta_{0}\right)} \overline{A_{1}} \overline{B_{3}}\\
    	    	& - 3 \, {\left(\theta_{0}^{4} + 10 \, \theta_{0}^{3} + 35 \, \theta_{0}^{2} + 50 \, \theta_{0} + 24\right)} \overline{A_{0}} \overline{B_{6}} + 48 \, \bigg({\left(\theta_{0}^{4} \overline{\alpha_{6}} + 9 \, \theta_{0}^{3} \overline{\alpha_{6}} + 26 \, \theta_{0}^{2} \overline{\alpha_{6}} + 24 \, \theta_{0} \overline{\alpha_{6}}\right)} A_{0}\\
    	    	& + {\left(\theta_{0}^{4} \overline{\alpha_{3}} + 9 \, \theta_{0}^{3} \overline{\alpha_{3}} + 26 \, \theta_{0}^{2} \overline{\alpha_{3}} + 24 \, \theta_{0} \overline{\alpha_{3}}\right)} A_{1}\bigg) \overline{A_{1}} + 48 \, \bigg({\left(\theta_{0}^{4} \overline{\alpha_{5}} + 8 \, \theta_{0}^{3} \overline{\alpha_{5}} + 19 \, \theta_{0}^{2} \overline{\alpha_{5}} + 12 \, \theta_{0} \overline{\alpha_{5}}\right)} A_{0}\\
    	    	& + {\left(\theta_{0}^{4} \overline{\alpha_{1}} + 8 \, \theta_{0}^{3} \overline{\alpha_{1}} + 19 \, \theta_{0}^{2} \overline{\alpha_{1}} + 12 \, \theta_{0} \overline{\alpha_{1}}\right)} A_{1}\bigg) \overline{A_{2}} + 3 \, \bigg({\left(\theta_{0}^{4} \overline{\alpha_{1}} + 10 \, \theta_{0}^{3} \overline{\alpha_{1}} + 35 \, \theta_{0}^{2} \overline{\alpha_{1}} + 50 \, \theta_{0} \overline{\alpha_{1}} + 24 \, \overline{\alpha_{1}}\right)} \overline{A_{0}}\\
    	    	& - 2 \, {\left(\theta_{0}^{4} + 8 \, \theta_{0}^{3} + 19 \, \theta_{0}^{2} + 12 \, \theta_{0}\right)} \overline{A_{2}}\bigg) \overline{B_{1}}\\
    	    	& + \bigg(4 \, {\left(\theta_{0}^{4} + 9 \, \theta_{0}^{3} + 26 \, \theta_{0}^{2} + 24 \, \theta_{0}\right)} {\left| A_{1} \right|}^{2} \overline{A_{1}} + 3 \, {\left(\theta_{0}^{4} \overline{\alpha_{5}} + 10 \, \theta_{0}^{3} \overline{\alpha_{5}} + 35 \, \theta_{0}^{2} \overline{\alpha_{5}} + 50 \, \theta_{0} \overline{\alpha_{5}} + 24 \, \overline{\alpha_{5}}\right)} \overline{A_{0}}\bigg) \overline{C_{1}}\bigg\}
    	    \end{align*}
    	    We will not need the precise expression of the other $\mu_{j}$ coefficients for $j\geq 16$, they are just written for completeness.

    	    \section{Final development of the tensors related to the invariance by inversions}
    	    \normalsize
    	    \begin{align*}
    	    	&|\phi|^2\vec{\alpha}-2\s{\vec{\alpha}}{\phi}\phi
    	    	=\begin{dmatrix}
    	    	-\frac{2 \, {\left(\theta_{0} \overline{\alpha_{2}} \overline{\zeta_{0}} + {\left(\theta_{0}^{2} + 3 \, \theta_{0} + 2\right)} \zeta_{13}\right)}}{\theta_{0}^{3} + 3 \, \theta_{0}^{2} + 2 \, \theta_{0}} & \overline{A_{0}} & 0 & \theta_{0} + 4 \\
    	    	-\frac{\theta_{0} - 2}{2 \, \theta_{0}^{2}} & C_{1} & 1 & \theta_{0} \\
    	    	-\frac{\theta_{0} - 2}{2 \, \theta_{0}^{2}} & B_{1} & 1 & \theta_{0} + 1 \\
    	    	-\frac{2 \, {\left(\alpha_{2} \theta_{0} - 2 \, \alpha_{2}\right)}}{\theta_{0}} & A_{0} & 1 & \theta_{0} + 1 \\
    	    	-\frac{\theta_{0} - 2}{2 \, \theta_{0}^{2}} & B_{2} & 1 & \theta_{0} + 2 \\
    	    	-\frac{{\left(\theta_{0} - 2\right)} {\left| C_{1} \right|}^{2}}{2 \, \theta_{0}^{3}} & \overline{A_{0}} & 1 & \theta_{0} + 2 \\
    	    	-\frac{\theta_{0} \overline{\alpha_{1}} - 2 \, \overline{\alpha_{1}}}{2 \, {\left(\theta_{0}^{2} + 2 \, \theta_{0}\right)}} & C_{1} & 1 & \theta_{0} + 2 \\
    	    	-\frac{2 \, {\left(\theta_{0} - 2\right)} \zeta_{4}}{\theta_{0}^{3} + 2 \, \theta_{0}^{2}} & A_{0} & 1 & \theta_{0} + 2 \\
    	    	\frac{1}{\theta_{0}^{2}} & B_{10} & 1 & \theta_{0} + 3 \\
    	    	-\frac{2 \, \zeta_{18}}{\theta_{0}} & A_{0} & 1 & \theta_{0} + 3 \\
    	    	-\frac{2 \, \zeta_{14}}{\theta_{0}} & \overline{A_{0}} & 1 & \theta_{0} + 3 \\
    	    	-\frac{\alpha_{2} \theta_{0} - 2 \, \alpha_{2}}{4 \, \theta_{0}} & \overline{C_{1}} & 1 & \theta_{0} + 3 \\
    	    	-\frac{\theta_{0} \overline{\alpha_{1}} - 2 \, \overline{\alpha_{1}}}{2 \, {\left(\theta_{0}^{2} + 2 \, \theta_{0}\right)}} & B_{1} & 1 & \theta_{0} + 3 \\
    	    	-\frac{2 \, \theta_{0}^{3} \overline{\alpha_{3}} - 2 \, \theta_{0}^{2} \overline{\alpha_{3}} - 4 \, \theta_{0} \overline{\alpha_{3}} - {\left(\theta_{0} - 2\right)} \overline{\zeta_{2}}}{4 \, {\left(\theta_{0}^{4} + 4 \, \theta_{0}^{3} + 3 \, \theta_{0}^{2}\right)}} & C_{1} & 1 & \theta_{0} + 3 \\
    	    	-\frac{\theta_{0} - 3}{2 \, \theta_{0}^{2}} & C_{2} & 2 & \theta_{0} \\
    	    	\frac{3 \, \zeta_{2}}{\theta_{0}^{3} + \theta_{0}^{2}} & \overline{A_{0}} & 2 & \theta_{0} \\
    	    	-\frac{\theta_{0} - 3}{2 \, \theta_{0}^{2}} & B_{3} & 2 & \theta_{0} + 1 \\
    	    	-\frac{2 \, {\left(\alpha_{2} \theta_{0} - 2 \, \alpha_{2}\right)}}{\theta_{0} + 1} & A_{1} & 2 & \theta_{0} + 1 \\
    	    	-\frac{{\left(\theta_{0} - 3\right)} \zeta_{3}}{\theta_{0}^{3} + \theta_{0}^{2}} & A_{0} & 2 & \theta_{0} + 1 \\
    	    	\frac{{\left(5 \, \theta_{0} + 2\right)} \zeta_{2}}{\theta_{0}^{4} + 2 \, \theta_{0}^{3} + \theta_{0}^{2}} & \overline{A_{1}} & 2 & \theta_{0} + 1 \\
    	    	-\frac{{\left(\theta_{0} - 2\right)} {\left| A_{1} \right|}^{2}}{\theta_{0}^{2} + 2 \, \theta_{0} + 1} & C_{1} & 2 & \theta_{0} + 1 
    	    	\end{dmatrix}
    	    	\end{align*}
    	    	\begin{align*}
    	    	\begin{dmatrix}
    	    	\frac{1}{\theta_{0}^{2}} & B_{11} & 2 & \theta_{0} + 2 \\
    	    	-\frac{2 \, \zeta_{19}}{\theta_{0}} & A_{0} & 2 & \theta_{0} + 2 \\
    	    	-\frac{\theta_{0} \overline{\alpha_{5}} - 2 \, \overline{\alpha_{5}}}{2 \, {\left(\theta_{0}^{2} + 3 \, \theta_{0} + 2\right)}} & C_{1} & 2 & \theta_{0} + 2 \\
    	    	-\frac{\theta_{0} \overline{\alpha_{1}} - 3 \, \overline{\alpha_{1}}}{2 \, {\left(\theta_{0}^{2} + 2 \, \theta_{0}\right)}} & C_{2} & 2 & \theta_{0} + 2 \\
    	    	-\frac{2 \, {\left({\left(\theta_{0} + 2\right)} \zeta_{15} - \zeta_{2} \overline{\alpha_{1}}\right)}}{\theta_{0}^{2} + 2 \, \theta_{0}} & \overline{A_{0}} & 2 & \theta_{0} + 2 \\
    	    	-\frac{2 \, {\left(\theta_{0} - 2\right)} \zeta_{4}}{\theta_{0}^{3} + 3 \, \theta_{0}^{2} + 2 \, \theta_{0}} & A_{1} & 2 & \theta_{0} + 2 \\
    	    	-\frac{{\left(2 \, \theta_{0} - 1\right)} \zeta_{2}}{\theta_{0}^{3} + 3 \, \theta_{0}^{2} + 2 \, \theta_{0}} & \overline{A_{2}} & 2 & \theta_{0} + 2 \\
    	    	-\frac{{\left(\theta_{0} - 2\right)} {\left| A_{1} \right|}^{2}}{\theta_{0}^{2} + 2 \, \theta_{0} + 1} & B_{1} & 2 & \theta_{0} + 2 \\
    	    	-\frac{\theta_{0} - 4}{2 \, \theta_{0}^{2}} & C_{3} & 3 & \theta_{0} \\
    	    	-\frac{2 \, {\left(\theta_{0} \zeta_{7} - \zeta_{5}\right)}}{\theta_{0}^{2}} & \overline{A_{0}} & 3 & \theta_{0} \\
    	    	\frac{2 \, {\left(\alpha_{7} \theta_{0} - 4 \, \alpha_{7}\right)}}{\theta_{0}} & A_{0} & 3 & \theta_{0} \\
    	    	-\frac{\alpha_{1} \theta_{0} - 2 \, \alpha_{1}}{2 \, {\left(\theta_{0}^{2} + 2 \, \theta_{0}\right)}} & C_{1} & 3 & \theta_{0} \\
    	    	\frac{1}{\theta_{0}^{2}} & B_{12} & 3 & \theta_{0} + 1 \\
    	    	-\frac{2 \, \zeta_{7}}{\theta_{0} + 1} & \overline{A_{1}} & 3 & \theta_{0} + 1 \\
    	    	-\frac{2 \, \zeta_{20}}{\theta_{0}} & A_{0} & 3 & \theta_{0} + 1 \\
    	    	-\frac{2 \, {\left(\alpha_{2} \theta_{0} - 2 \, \alpha_{2}\right)}}{\theta_{0} + 2} & A_{2} & 3 & \theta_{0} + 1 \\
    	    	-\frac{\alpha_{1} \theta_{0} - 2 \, \alpha_{1}}{2 \, {\left(\theta_{0}^{2} + 2 \, \theta_{0}\right)}} & B_{1} & 3 & \theta_{0} + 1 \\
    	    	-\frac{\alpha_{5} \theta_{0} - 2 \, \alpha_{5}}{2 \, {\left(\theta_{0}^{2} + 3 \, \theta_{0} + 2\right)}} & C_{1} & 3 & \theta_{0} + 1 \\
    	    	\frac{2 \, {\left(2 \, \theta_{0} \zeta_{2} {\left| A_{1} \right|}^{2} - {\left(\theta_{0}^{2} + 2 \, \theta_{0} + 1\right)} \zeta_{16}\right)}}{\theta_{0}^{3} + 2 \, \theta_{0}^{2} + \theta_{0}} & \overline{A_{0}} & 3 & \theta_{0} + 1 \\
    	    	-\frac{{\left(\theta_{0} - 3\right)} \zeta_{3}}{\theta_{0}^{3} + 2 \, \theta_{0}^{2} + \theta_{0}} & A_{1} & 3 & \theta_{0} + 1 \\
    	    	-\frac{{\left(\theta_{0} - 3\right)} {\left| A_{1} \right|}^{2}}{\theta_{0}^{2} + 2 \, \theta_{0} + 1} & C_{2} & 3 & \theta_{0} + 1 \\
    	    	\frac{1}{\theta_{0}^{2}} & B_{13} & 4 & \theta_{0} \\
    	    	-\frac{2 \, \zeta_{21}}{\theta_{0}} & A_{0} & 4 & \theta_{0} \\
    	    	-\frac{2 \, {\left({\left(\theta_{0} + 2\right)} \zeta_{17} - \alpha_{1} \zeta_{2}\right)}}{\theta_{0}^{2} + 2 \, \theta_{0}} & \overline{A_{0}} & 4 & \theta_{0} \\
    	    	-\frac{\alpha_{1} \theta_{0} - 3 \, \alpha_{1}}{2 \, {\left(\theta_{0}^{2} + 2 \, \theta_{0}\right)}} & C_{2} & 4 & \theta_{0} \\
    	    	-\frac{4 \, \alpha_{3} \theta_{0}^{2} - 8 \, \alpha_{3} \theta_{0} + {\left(2 \, \theta_{0} + 1\right)} \zeta_{2}}{8 \, {\left(\theta_{0}^{3} + 3 \, \theta_{0}^{2}\right)}} & C_{1} & 4 & \theta_{0} 
    	    	\end{dmatrix}
    	    	\end{align*}
    	    	\begin{align*}
    	    	\begin{dmatrix}
    	    	-\frac{2 \, {\left(\alpha_{2} \theta_{0} - 2 \, \alpha_{2}\right)}}{\theta_{0}} & \overline{A_{0}} & -\theta_{0} + 1 & 2 \, \theta_{0} + 1 \\
    	    	-\frac{2 \, {\left(\alpha_{2} \theta_{0} - 2 \, \alpha_{2}\right)}}{\theta_{0} + 1} & \overline{A_{1}} & -\theta_{0} + 1 & 2 \, \theta_{0} + 2 \\
    	    	\frac{{\left(\theta_{0} - 2\right)} \overline{\zeta_{0}}}{2 \, {\left(\theta_{0}^{4} + 4 \, \theta_{0}^{3} + 5 \, \theta_{0}^{2} + 2 \, \theta_{0}\right)}} & C_{1} & -\theta_{0} + 1 & 2 \, \theta_{0} + 2 \\
    	    	-\frac{2 \, {\left(\theta_{0} - 2\right)} \zeta_{4}}{\theta_{0}^{3} + 2 \, \theta_{0}^{2}} & \overline{A_{0}} & -\theta_{0} + 1 & 2 \, \theta_{0} + 2 \\
    	    	-\frac{2 \, \zeta_{18}}{\theta_{0}} & \overline{A_{0}} & -\theta_{0} + 1 & 2 \, \theta_{0} + 3 \\
    	    	-\frac{2 \, {\left(\alpha_{2} \theta_{0} - 2 \, \alpha_{2}\right)}}{\theta_{0} + 2} & \overline{A_{2}} & -\theta_{0} + 1 & 2 \, \theta_{0} + 3 \\
    	    	\frac{2 \, {\left(\theta_{0} - 2\right)} \overline{\zeta_{1}}}{\theta_{0}^{4} + 6 \, \theta_{0}^{3} + 11 \, \theta_{0}^{2} + 6 \, \theta_{0}} & C_{1} & -\theta_{0} + 1 & 2 \, \theta_{0} + 3 \\
    	    	\frac{{\left(\theta_{0} - 2\right)} \overline{\zeta_{0}}}{2 \, {\left(\theta_{0}^{4} + 4 \, \theta_{0}^{3} + 5 \, \theta_{0}^{2} + 2 \, \theta_{0}\right)}} & B_{1} & -\theta_{0} + 1 & 2 \, \theta_{0} + 3 \\
    	    	-\frac{2 \, {\left(\theta_{0} - 2\right)} \zeta_{4}}{\theta_{0}^{3} + 3 \, \theta_{0}^{2} + 2 \, \theta_{0}} & \overline{A_{1}} & -\theta_{0} + 1 & 2 \, \theta_{0} + 3 \\
    	    	-\frac{{\left(\theta_{0} - 3\right)} \zeta_{3}}{\theta_{0}^{3} + \theta_{0}^{2}} & \overline{A_{0}} & -\theta_{0} + 2 & 2 \, \theta_{0} + 1 \\
    	    	-\frac{2 \, {\left({\left(\theta_{0}^{3} + 4 \, \theta_{0}^{2} + 5 \, \theta_{0} + 2\right)} \zeta_{19} + \zeta_{2} \overline{\zeta_{0}}\right)}}{\theta_{0}^{4} + 4 \, \theta_{0}^{3} + 5 \, \theta_{0}^{2} + 2 \, \theta_{0}} & \overline{A_{0}} & -\theta_{0} + 2 & 2 \, \theta_{0} + 2 \\
    	    	\frac{{\left(\theta_{0} - 3\right)} \overline{\zeta_{0}}}{2 \, {\left(\theta_{0}^{4} + 4 \, \theta_{0}^{3} + 5 \, \theta_{0}^{2} + 2 \, \theta_{0}\right)}} & C_{2} & -\theta_{0} + 2 & 2 \, \theta_{0} + 2 \\
    	    	-\frac{{\left(\theta_{0} - 3\right)} \zeta_{3}}{\theta_{0}^{3} + 2 \, \theta_{0}^{2} + \theta_{0}} & \overline{A_{1}} & -\theta_{0} + 2 & 2 \, \theta_{0} + 2 \\
    	    	-\frac{\alpha_{2} \theta_{0}^{2} - 6 \, \alpha_{2} \theta_{0} + 8 \, \alpha_{2}}{4 \, \theta_{0}^{2}} & C_{1} & -\theta_{0} + 3 & 2 \, \theta_{0} + 1 \\
    	    	-\frac{2 \, \zeta_{20}}{\theta_{0}} & \overline{A_{0}} & -\theta_{0} + 3 & 2 \, \theta_{0} + 1 \\
    	    	-\frac{2 \, {\left(\alpha_{7} \theta_{0}^{2} - 6 \, \alpha_{7} \theta_{0} + 8 \, \alpha_{7}\right)}}{\theta_{0}^{3} + \theta_{0}^{2}} & \overline{A_{1}} & -\theta_{0} + 3 & 2 \, \theta_{0} + 1 \\
    	    	-\frac{2 \, \zeta_{21}}{\theta_{0}} & \overline{A_{0}} & -\theta_{0} + 4 & 2 \, \theta_{0} \\
    	    	\frac{\overline{\zeta_{2}}}{\theta_{0}^{3} - 3 \, \theta_{0}^{2}} & A_{1} & \theta_{0} & 3 \\
    	    	\frac{2 \, \overline{\alpha_{2}}}{\theta_{0}} & \overline{A_{1}} & \theta_{0} & 3 \\
    	    	\frac{1}{\theta_{0}^{2}} & B_{7} & \theta_{0} & 4 \\
    	    	-\frac{\overline{\zeta_{3}}}{\theta_{0}^{2} + \theta_{0}} & \overline{A_{1}} & \theta_{0} & 4 \\
    	    	-\frac{2 \, \zeta_{13}}{\theta_{0}} & A_{0} & \theta_{0} & 4 \\
    	    	-\frac{2 \, {\left(\theta_{0} \overline{\alpha_{2}} + \overline{\alpha_{2}}\right)}}{\theta_{0} + 2} & \overline{A_{2}} & \theta_{0} & 4 \\
    	    	\frac{2 \, {\left(\theta_{0}^{2} \overline{\alpha_{1}} \overline{\alpha_{2}} + \theta_{0} \overline{\alpha_{1}} \overline{\alpha_{2}} - {\left(\theta_{0} + 2\right)} \zeta_{22}\right)}}{\theta_{0}^{2} + 2 \, \theta_{0}} & \overline{A_{0}} & \theta_{0} & 4 
    	    	\end{dmatrix}
    	    	\end{align*}
    	    	\begin{align*}
    	    	\begin{dmatrix}
    	    	\frac{{\left(\theta_{0} - 2\right)} {\left| C_{1} \right|}^{2}}{4 \, \theta_{0}^{3}} & A_{0} & \theta_{0} + 1 & 2 \\
    	    	\frac{1}{\theta_{0}^{2}} & B_{8} & \theta_{0} + 1 & 3 \\
    	    	-\frac{2 \, \overline{\zeta_{4}}}{\theta_{0}^{2} + \theta_{0}} & \overline{A_{1}} & \theta_{0} + 1 & 3 \\
    	    	-\frac{2 \, \zeta_{14}}{\theta_{0}} & A_{0} & \theta_{0} + 1 & 3 \\
    	    	\frac{2 \, {\left(2 \, \theta_{0}^{2} {\left| A_{1} \right|}^{2} \overline{\alpha_{2}} - {\left(\theta_{0} + 1\right)} \zeta_{23}\right)}}{\theta_{0}^{2} + \theta_{0}} & \overline{A_{0}} & \theta_{0} + 1 & 3 \\
    	    	-\frac{{\left(2 \, \theta_{0} - 1\right)} \zeta_{2}}{\theta_{0}^{3} + \theta_{0}^{2}} & A_{0} & \theta_{0} + 2 & 0 \\
    	    	\frac{1}{\theta_{0}^{2}} & B_{9} & \theta_{0} + 2 & 2 \\
    	    	-\frac{2 \, \zeta_{15}}{\theta_{0}} & A_{0} & \theta_{0} + 2 & 2 \\
    	    	-\frac{\theta_{0} \overline{\alpha_{2}} + \overline{\alpha_{2}}}{4 \, \theta_{0}} & C_{1} & \theta_{0} + 2 & 2 \\
    	    	\frac{2 \, {\left(\alpha_{1} \theta_{0}^{2} \overline{\alpha_{2}} + \alpha_{1} \theta_{0} \overline{\alpha_{2}} - {\left(\theta_{0} + 2\right)} \zeta_{24}\right)}}{\theta_{0}^{2} + 2 \, \theta_{0}} & \overline{A_{0}} & \theta_{0} + 2 & 2 \\
    	    	-\frac{{\left(2 \, \theta_{0} - 1\right)} \zeta_{2}}{8 \, {\left(\theta_{0}^{3} + \theta_{0}^{2}\right)}} & \overline{C_{1}} & \theta_{0} + 2 & 2 \\
    	    	-\frac{2 \, \zeta_{7}}{\theta_{0}} & A_{0} & \theta_{0} + 3 & 0 \\
    	    	-\frac{{\left(2 \, \theta_{0} - 1\right)} \zeta_{2}}{\theta_{0}^{3} + 2 \, \theta_{0}^{2} + \theta_{0}} & A_{1} & \theta_{0} + 3 & 0 \\
    	    	\frac{{\left(\theta_{0} - 2\right)} \zeta_{0}}{2 \, {\left(\theta_{0}^{4} + 4 \, \theta_{0}^{3} + 5 \, \theta_{0}^{2} + 2 \, \theta_{0}\right)}} & C_{1} & \theta_{0} + 3 & 0 \\
    	    	-\frac{2 \, \zeta_{16}}{\theta_{0}} & A_{0} & \theta_{0} + 3 & 1 \\
    	    	\frac{{\left(\theta_{0} - 2\right)} \zeta_{0}}{2 \, {\left(\theta_{0}^{4} + 4 \, \theta_{0}^{3} + 5 \, \theta_{0}^{2} + 2 \, \theta_{0}\right)}} & B_{1} & \theta_{0} + 3 & 1 \\
    	    	-\frac{2 \, \zeta_{7}}{\theta_{0} + 1} & A_{1} & \theta_{0} + 4 & 0 \\
    	    	-\frac{2 \, \zeta_{17}}{\theta_{0}} & A_{0} & \theta_{0} + 4 & 0 \\
    	    	-\frac{2 \, \zeta_{0} \zeta_{2}}{\theta_{0}^{4} + 4 \, \theta_{0}^{3} + 5 \, \theta_{0}^{2} + 2 \, \theta_{0}} & \overline{A_{0}} & \theta_{0} + 4 & 0 \\
    	    	-\frac{{\left(2 \, \theta_{0} - 1\right)} \zeta_{2}}{\theta_{0}^{3} + 3 \, \theta_{0}^{2} + 2 \, \theta_{0}} & A_{2} & \theta_{0} + 4 & 0 \\
    	    	\frac{2 \, {\left(\theta_{0} - 2\right)} \zeta_{1}}{\theta_{0}^{4} + 6 \, \theta_{0}^{3} + 11 \, \theta_{0}^{2} + 6 \, \theta_{0}} & C_{1} & \theta_{0} + 4 & 0 \\
    	    	\frac{{\left(\theta_{0} - 3\right)} \zeta_{0}}{2 \, {\left(\theta_{0}^{4} + 4 \, \theta_{0}^{3} + 5 \, \theta_{0}^{2} + 2 \, \theta_{0}\right)}} & C_{2} & \theta_{0} + 4 & 0 
    	    	\end{dmatrix}
    	    	\end{align*}
    	    	\begin{align*}
    	    	\begin{dmatrix}
    	    	-\frac{2 \, \overline{\zeta_{4}}}{\theta_{0}^{2}} & A_{0} & 2 \, \theta_{0} + 1 & -\theta_{0} + 2 \\
    	    	-2 \, \overline{\alpha_{2}} & A_{1} & 2 \, \theta_{0} + 1 & -\theta_{0} + 2 \\
    	    	-\frac{\overline{\zeta_{3}}}{\theta_{0}^{2} + \theta_{0}} & A_{1} & 2 \, \theta_{0} + 1 & -\theta_{0} + 3 \\
    	    	-\frac{2 \, \zeta_{23}}{\theta_{0}} & A_{0} & 2 \, \theta_{0} + 1 & -\theta_{0} + 3 \\
    	    	-\frac{2 \, \overline{\zeta_{4}}}{\theta_{0}^{2} + \theta_{0}} & A_{1} & 2 \, \theta_{0} + 2 & -\theta_{0} + 2 \\
    	    	-\frac{2 \, \zeta_{24}}{\theta_{0}} & A_{0} & 2 \, \theta_{0} + 2 & -\theta_{0} + 2 \\
    	    	-\frac{2 \, {\left(\theta_{0} \overline{\alpha_{2}} + \overline{\alpha_{2}}\right)}}{\theta_{0} + 2} & A_{2} & 2 \, \theta_{0} + 2 & -\theta_{0} + 2 \\
    	    	-\frac{2 \, \zeta_{0} \overline{\alpha_{2}}}{\theta_{0}^{2} + 3 \, \theta_{0} + 2} & \overline{A_{0}} & 2 \, \theta_{0} + 2 & -\theta_{0} + 2 \\
    	    	-\frac{2 \, {\left(\theta_{0} \overline{\alpha_{2}} + \overline{\alpha_{2}}\right)}}{\theta_{0}} & A_{0} & 2 \, \theta_{0} & -\theta_{0} + 2 \\
    	    	-\frac{\overline{\zeta_{3}}}{\theta_{0}^{2}} & A_{0} & 2 \, \theta_{0} & -\theta_{0} + 3 \\
    	    	-\frac{\theta_{0} \overline{\alpha_{7}} - 4 \, \overline{\alpha_{7}}}{\theta_{0}^{2}} & A_{1} & 2 \, \theta_{0} & -\theta_{0} + 4 \\
    	    	-\frac{\theta_{0}^{2} \overline{\alpha_{2}} - \theta_{0} \overline{\alpha_{2}} - 2 \, \overline{\alpha_{2}}}{4 \, \theta_{0}^{2}} & \overline{C_{1}} & 2 \, \theta_{0} & -\theta_{0} + 4 \\
    	    	-\frac{2 \, \zeta_{22}}{\theta_{0}} & A_{0} & 2 \, \theta_{0} & -\theta_{0} + 4
    	    	\end{dmatrix}
    	    \end{align*}
    	    while
    	    \footnotesize
    	    \begin{align*}
    	    	&-g^{-1}\otimes\left(\bar{\partial}|\phi|^2\otimes\h_0-2\s{\phi}{\h_0}\otimes\bar{\partial}\phi\right)
    	    	=\begin{dmatrix}
    	    	-\frac{4 \, {\left| A_{1} \right|}^{2}}{\theta_{0}^{2} + \theta_{0}} & \overline{A_{0}} & 0 & \theta_{0} + 1 \\
    	    	\frac{4 \, \overline{\zeta_{0}}}{\theta_{0}^{3} + 3 \, \theta_{0}^{2} + 2 \, \theta_{0}} & A_{1} & 0 & \theta_{0} + 2 \\
    	    	-\frac{4 \, {\left| A_{1} \right|}^{2}}{\theta_{0}^{2} + \theta_{0}} & \overline{A_{1}} & 0 & \theta_{0} + 2 \\
    	    	2 \, \mu_{1} & \overline{A_{0}} & 0 & \theta_{0} + 2 \\
    	    	\frac{8 \, {\left(2 \, \theta_{0} + 3\right)} \overline{\zeta_{1}}}{\theta_{0}^{4} + 6 \, \theta_{0}^{3} + 11 \, \theta_{0}^{2} + 6 \, \theta_{0}} & A_{1} & 0 & \theta_{0} + 3 \\
    	    	-\frac{4 \, {\left| A_{1} \right|}^{2}}{\theta_{0}^{2} + \theta_{0}} & \overline{A_{2}} & 0 & \theta_{0} + 3 \\
    	    	2 \, \mu_{2} & \overline{A_{0}} & 0 & \theta_{0} + 3 \\
    	    	2 \, \mu_{1} & \overline{A_{1}} & 0 & \theta_{0} + 3 \\
    	    	-\frac{8 \, {\left| A_{1} \right|}^{2} \overline{\zeta_{0}}}{\theta_{0}^{3} + 3 \, \theta_{0}^{2} + 2 \, \theta_{0}} & A_{0} & 0 & \theta_{0} + 3 \\
    	    	\frac{\overline{\zeta_{0}}}{2 \, {\left(\theta_{0}^{3} + 3 \, \theta_{0}^{2} + 2 \, \theta_{0}\right)}} & \overline{B_{1}} & 0 & \theta_{0} + 4 \\
    	    	-\frac{4 \, {\left(\overline{\alpha_{1}} \overline{\zeta_{0}} + {\left(\theta_{0}^{4} + 5 \, \theta_{0}^{3} + 8 \, \theta_{0}^{2} + 4 \, \theta_{0}\right)} \overline{\zeta_{9}}\right)}}{\theta_{0}^{3} + 3 \, \theta_{0}^{2} + 2 \, \theta_{0}} & A_{1} & 0 & \theta_{0} + 4 \\
    	    	-\frac{4 \, {\left(4 \, {\left(2 \, \theta_{0} + 3\right)} {\left| A_{1} \right|}^{2} \overline{\zeta_{1}} + {\left(\theta_{0} \overline{\alpha_{5}} + 3 \, \overline{\alpha_{5}}\right)} \overline{\zeta_{0}}\right)}}{\theta_{0}^{4} + 6 \, \theta_{0}^{3} + 11 \, \theta_{0}^{2} + 6 \, \theta_{0}} & A_{0} & 0 & \theta_{0} + 4 \\
    	    	-\frac{4 \, {\left| A_{1} \right|}^{2}}{\theta_{0}^{2} + \theta_{0}} & \overline{A_{3}} & 0 & \theta_{0} + 4 \\
    	    	2 \, \mu_{2} & \overline{A_{1}} & 0 & \theta_{0} + 4 \\
    	    	2 \, \mu_{15} & \overline{A_{0}} & 0 & \theta_{0} + 4 \\
    	    	2 \, \mu_{1} & \overline{A_{2}} & 0 & \theta_{0} + 4 \\
    	    	\frac{\theta_{0} - 2}{2 \, \theta_{0}^{2}} & C_{1} & 1 & \theta_{0} \\
    	    	-\frac{2 \, {\left(\alpha_{2} \theta_{0} - 2 \, \alpha_{2}\right)}}{\theta_{0}} & A_{0} & 1 & \theta_{0} + 1 \\
    	    	\frac{\theta_{0} - 2}{2 \, {\left(\theta_{0}^{2} + \theta_{0}\right)}} & B_{1} & 1 & \theta_{0} + 1 \\
    	    	2 \, \mu_{3} & \overline{A_{0}} & 1 & \theta_{0} + 1 \\
    	    	\frac{8 \, \overline{\zeta_{0}}}{\theta_{0}^{3} + 3 \, \theta_{0}^{2} + 2 \, \theta_{0}} & A_{2} & 1 & \theta_{0} + 2 \\
    	    	\frac{\theta_{0} - 2}{2 \, {\left(\theta_{0}^{2} + 2 \, \theta_{0}\right)}} & B_{2} & 1 & \theta_{0} + 2 \\
    	    	\frac{\theta_{0} \overline{\alpha_{1}} - 2 \, \overline{\alpha_{1}}}{2 \, {\left(\theta_{0}^{2} + 2 \, \theta_{0}\right)}} & C_{1} & 1 & \theta_{0} + 2 \\
    	    	-\frac{2 \, {\left(\alpha_{9} \theta_{0}^{3} + \alpha_{9} \theta_{0}^{2} - 4 \, \alpha_{9} \theta_{0} + 4 \, \alpha_{1} \overline{\zeta_{0}} - 4 \, \alpha_{9}\right)}}{\theta_{0}^{3} + 3 \, \theta_{0}^{2} + 2 \, \theta_{0}} & A_{0} & 1 & \theta_{0} + 2 \\
    	    	2 \, \mu_{4} & \overline{A_{0}} & 1 & \theta_{0} + 2 \\
    	    	2 \, \mu_{3} & \overline{A_{1}} & 1 & \theta_{0} + 2 \\
    	    	\frac{\theta_{0} - 2}{2 \, {\left(\theta_{0}^{2} + 3 \, \theta_{0}\right)}} & B_{4} & 1 & \theta_{0} + 3 \\
    	    	\frac{\theta_{0} \overline{\alpha_{1}} - 2 \, \overline{\alpha_{1}}}{2 \, {\left(\theta_{0}^{2} + 3 \, \theta_{0}\right)}} & B_{1} & 1 & \theta_{0} + 3 \\
    	    	-\frac{\alpha_{2} \theta_{0}^{3} - 3 \, \alpha_{2} \theta_{0}^{2} - 4 \, \alpha_{2} \theta_{0} + 12 \, \alpha_{2}}{4 \, {\left(\theta_{0}^{3} + 3 \, \theta_{0}^{2}\right)}} & \overline{C_{1}} & 1 & \theta_{0} + 3 
    	    	\end{dmatrix}
    	    	\end{align*}
    	    	\begin{align*}
    	    	\begin{dmatrix}
    	    	\pi_1 & A_{0} & 1 & \theta_{0} + 3 \\
    	    	\frac{2 \, \theta_{0}^{4} \overline{\alpha_{3}} - 2 \, \theta_{0}^{3} \overline{\alpha_{3}} - 4 \, \theta_{0}^{2} \overline{\alpha_{3}} - {\left(\theta_{0}^{2} + \theta_{0} - 6\right)} \overline{\zeta_{2}}}{4 \, {\left(\theta_{0}^{5} + 4 \, \theta_{0}^{4} + 3 \, \theta_{0}^{3}\right)}} & C_{1} & 1 & \theta_{0} + 3 \\
    	    	\frac{16 \, {\left(2 \, \theta_{0} + 3\right)} \overline{\zeta_{1}}}{\theta_{0}^{4} + 6 \, \theta_{0}^{3} + 11 \, \theta_{0}^{2} + 6 \, \theta_{0}} & A_{2} & 1 & \theta_{0} + 3 \\
    	    	2 \, \mu_{4} & \overline{A_{1}} & 1 & \theta_{0} + 3 \\
    	    	2 \, \mu_{3} & \overline{A_{2}} & 1 & \theta_{0} + 3 \\
    	    	2 \, \mu_{16} & \overline{A_{0}} & 1 & \theta_{0} + 3 \\
    	    	-\frac{16 \, {\left| A_{1} \right|}^{2} \overline{\zeta_{0}}}{\theta_{0}^{3} + 3 \, \theta_{0}^{2} + 2 \, \theta_{0}} & A_{1} & 1 & \theta_{0} + 3 \\
    	    	\frac{\theta_{0} - 3}{2 \, \theta_{0}^{2}} & C_{2} & 2 & \theta_{0} \\
    	    	2 \, \mu_{5} & \overline{A_{0}} & 2 & \theta_{0} \\
    	    	-\frac{2 \, {\left(\alpha_{8} \theta_{0} - 3 \, \alpha_{8}\right)}}{\theta_{0}} & A_{0} & 2 & \theta_{0} + 1 \\
    	    	-\frac{2 \, {\left(\alpha_{2} \theta_{0} - 3 \, \alpha_{2}\right)}}{\theta_{0}} & A_{1} & 2 & \theta_{0} + 1 \\
    	    	\frac{\theta_{0} - 3}{2 \, {\left(\theta_{0}^{2} + \theta_{0}\right)}} & B_{3} & 2 & \theta_{0} + 1 \\
    	    	2 \, \mu_{6} & \overline{A_{0}} & 2 & \theta_{0} + 1 \\
    	    	2 \, \mu_{5} & \overline{A_{1}} & 2 & \theta_{0} + 1 \\
    	    	\frac{{\left(2 \, \theta_{0}^{2} - 7 \, \theta_{0} + 6\right)} {\left| A_{1} \right|}^{2}}{2 \, {\left(\theta_{0}^{3} + \theta_{0}^{2}\right)}} & C_{1} & 2 & \theta_{0} + 1 \\
    	    	\frac{12 \, \overline{\zeta_{0}}}{\theta_{0}^{3} + 3 \, \theta_{0}^{2} + 2 \, \theta_{0}} & A_{3} & 2 & \theta_{0} + 2 \\
    	    	\frac{\theta_{0} - 3}{2 \, {\left(\theta_{0}^{2} + 2 \, \theta_{0}\right)}} & B_{5} & 2 & \theta_{0} + 2 \\
    	    	\frac{\theta_{0} \overline{\alpha_{1}} - 3 \, \overline{\alpha_{1}}}{2 \, {\left(\theta_{0}^{2} + 2 \, \theta_{0}\right)}} & C_{2} & 2 & \theta_{0} + 2 \\
    	    	-\frac{2 \, {\left(\alpha_{9} \theta_{0}^{3} - 7 \, \alpha_{9} \theta_{0} + 6 \, \alpha_{1} \overline{\zeta_{0}} - 6 \, \alpha_{9}\right)}}{\theta_{0}^{3} + 3 \, \theta_{0}^{2} + 2 \, \theta_{0}} & A_{1} & 2 & \theta_{0} + 2 \\
    	    	-\frac{2 \, {\left(\alpha_{11} \theta_{0}^{3} - 2 \, {\left(\alpha_{2} \theta_{0}^{3} - 10 \, \alpha_{2} \theta_{0} - 12 \, \alpha_{2}\right)} {\left| A_{1} \right|}^{2} - 7 \, \alpha_{11} \theta_{0} + 6 \, \alpha_{3} \overline{\zeta_{0}} - 6 \, \alpha_{11}\right)}}{\theta_{0}^{3} + 3 \, \theta_{0}^{2} + 2 \, \theta_{0}} & A_{0} & 2 & \theta_{0} + 2 \\
    	    	\frac{\mu_{1} \theta_{0}^{4} + {\left(3 \, \mu_{1} + 2 \, \overline{\alpha_{5}}\right)} \theta_{0}^{3} + 2 \, {\left(\mu_{1} - 2 \, \overline{\alpha_{5}}\right)} \theta_{0}^{2} + 12 \, \overline{\alpha_{5}}}{4 \, {\left(\theta_{0}^{4} + 3 \, \theta_{0}^{3} + 2 \, \theta_{0}^{2}\right)}} & C_{1} & 2 & \theta_{0} + 2 \\
    	    	2 \, \mu_{6} & \overline{A_{1}} & 2 & \theta_{0} + 2 \\
    	    	2 \, \mu_{5} & \overline{A_{2}} & 2 & \theta_{0} + 2 \\
    	    	2 \, \mu_{17} & \overline{A_{0}} & 2 & \theta_{0} + 2 \\
    	    	\frac{{\left(2 \, \theta_{0}^{3} - 3 \, \theta_{0}^{2} - 7 \, \theta_{0} + 4\right)} {\left| A_{1} \right|}^{2}}{2 \, {\left(\theta_{0}^{4} + 4 \, \theta_{0}^{3} + 5 \, \theta_{0}^{2} + 2 \, \theta_{0}\right)}} & B_{1} & 2 & \theta_{0} + 2 \\
    	    	\frac{\theta_{0} - 4}{2 \, \theta_{0}^{2}} & C_{3} & 3 & \theta_{0} \\
    	    	-\frac{2 \, {\left(\alpha_{7} \theta_{0} - 4 \, \alpha_{7}\right)}}{\theta_{0}} & A_{0} & 3 & \theta_{0} \\
    	    	\frac{\alpha_{1} \theta_{0} - 2 \, \alpha_{1}}{2 \, {\left(\theta_{0}^{2} + 2 \, \theta_{0}\right)}} & C_{1} & 3 & \theta_{0} \\
    	    	2 \, \mu_{7} & \overline{A_{0}} & 3 & \theta_{0} 
    	    	\end{dmatrix}
    	    	\end{align*}
    	    	\begin{align*}
    	    	\begin{dmatrix}
    	    	-\frac{2 \, {\left(\alpha_{8} \theta_{0} - 4 \, \alpha_{8}\right)}}{\theta_{0}} & A_{1} & 3 & \theta_{0} + 1 \\
    	    	-\frac{2 \, {\left(\alpha_{2} \theta_{0} - 4 \, \alpha_{2}\right)}}{\theta_{0}} & A_{2} & 3 & \theta_{0} + 1 \\
    	    	\frac{\theta_{0} - 4}{2 \, {\left(\theta_{0}^{2} + \theta_{0}\right)}} & B_{6} & 3 & \theta_{0} + 1 \\
    	    	\frac{2 \, {\left({\left(\alpha_{1} \alpha_{2} - \alpha_{12}\right)} \theta_{0}^{2} - 16 \, \alpha_{1} \alpha_{2} - 2 \, {\left(\alpha_{1} \alpha_{2} - \alpha_{12}\right)} \theta_{0} + 8 \, \alpha_{12}\right)}}{\theta_{0}^{2} + 2 \, \theta_{0}} & A_{0} & 3 & \theta_{0} + 1 \\
    	    	\frac{\alpha_{1} \theta_{0} - 2 \, \alpha_{1}}{2 \, {\left(\theta_{0}^{2} + 3 \, \theta_{0} + 2\right)}} & B_{1} & 3 & \theta_{0} + 1 \\
    	    	\frac{\mu_{3} \theta_{0}^{4} + {\left(2 \, \alpha_{5} + 3 \, \mu_{3}\right)} \theta_{0}^{3} - 2 \, {\left(2 \, \alpha_{5} - \mu_{3}\right)} \theta_{0}^{2} + 16 \, \alpha_{5}}{4 \, {\left(\theta_{0}^{4} + 3 \, \theta_{0}^{3} + 2 \, \theta_{0}^{2}\right)}} & C_{1} & 3 & \theta_{0} + 1 \\
    	    	2 \, \mu_{7} & \overline{A_{1}} & 3 & \theta_{0} + 1 \\
    	    	2 \, \mu_{18} & \overline{A_{0}} & 3 & \theta_{0} + 1 \\
    	    	\frac{{\left(3 \, \theta_{0}^{2} - 13 \, \theta_{0} + 12\right)} {\left| A_{1} \right|}^{2}}{3 \, {\left(\theta_{0}^{3} + \theta_{0}^{2}\right)}} & C_{2} & 3 & \theta_{0} + 1 \\
    	    	\frac{\theta_{0} - 5}{2 \, \theta_{0}^{2}} & C_{4} & 4 & \theta_{0} \\
    	    	-\frac{2 \, {\left(\alpha_{7} \theta_{0} - 4 \, \alpha_{7}\right)}}{\theta_{0}} & A_{1} & 4 & \theta_{0} \\
    	    	-\frac{2 \, {\left(\alpha_{13} \theta_{0} - 5 \, \alpha_{13}\right)}}{\theta_{0}} & A_{0} & 4 & \theta_{0} \\
    	    	\frac{\alpha_{1} \theta_{0} - 3 \, \alpha_{1}}{2 \, {\left(\theta_{0}^{2} + 2 \, \theta_{0}\right)}} & C_{2} & 4 & \theta_{0} \\
    	    	\frac{\mu_{5} \theta_{0}^{4} + 2 \, {\left(\alpha_{3} + 2 \, \mu_{5}\right)} \theta_{0}^{3} - {\left(2 \, \alpha_{3} - 3 \, \mu_{5}\right)} \theta_{0}^{2} - 4 \, \alpha_{3} \theta_{0} - {\left(\theta_{0} - 2\right)} \zeta_{2}}{4 \, {\left(\theta_{0}^{4} + 4 \, \theta_{0}^{3} + 3 \, \theta_{0}^{2}\right)}} & C_{1} & 4 & \theta_{0} \\
    	    	2 \, \mu_{19} & \overline{A_{0}} & 4 & \theta_{0} \\
    	    	\frac{2 \, {\left(\alpha_{2} \theta_{0} - 2 \, \alpha_{2}\right)}}{\theta_{0}} & \overline{A_{0}} & -\theta_{0} + 1 & 2 \, \theta_{0} + 1 \\
    	    	\frac{2 \, {\left(\alpha_{2} \theta_{0} - 2 \, \alpha_{2}\right)}}{\theta_{0}} & \overline{A_{1}} & -\theta_{0} + 1 & 2 \, \theta_{0} + 2 \\
    	    	-\frac{{\left(\theta_{0} - 2\right)} \overline{\zeta_{0}}}{\theta_{0}^{4} + 3 \, \theta_{0}^{3} + 2 \, \theta_{0}^{2}} & C_{1} & -\theta_{0} + 1 & 2 \, \theta_{0} + 2 \\
    	    	2 \, \mu_{8} & \overline{A_{0}} & -\theta_{0} + 1 & 2 \, \theta_{0} + 2 \\
    	    	\frac{2 \, {\left(\alpha_{2} \theta_{0} - 2 \, \alpha_{2}\right)}}{\theta_{0}} & \overline{A_{2}} & -\theta_{0} + 1 & 2 \, \theta_{0} + 3 \\
    	    	-\frac{2 \, {\left(2 \, \theta_{0}^{2} - \theta_{0} - 6\right)} \overline{\zeta_{1}}}{\theta_{0}^{5} + 6 \, \theta_{0}^{4} + 11 \, \theta_{0}^{3} + 6 \, \theta_{0}^{2}} & C_{1} & -\theta_{0} + 1 & 2 \, \theta_{0} + 3 \\
    	    	-\frac{{\left(\theta_{0} - 2\right)} \overline{\zeta_{0}}}{\theta_{0}^{4} + 4 \, \theta_{0}^{3} + 5 \, \theta_{0}^{2} + 2 \, \theta_{0}} & B_{1} & -\theta_{0} + 1 & 2 \, \theta_{0} + 3 \\
    	    	\frac{4 \, {\left(\alpha_{2} \theta_{0} - 2 \, \alpha_{2}\right)} \overline{\zeta_{0}}}{\theta_{0}^{3} + 3 \, \theta_{0}^{2} + 2 \, \theta_{0}} & A_{0} & -\theta_{0} + 1 & 2 \, \theta_{0} + 3 \\
    	    	2 \, \mu_{8} & \overline{A_{1}} & -\theta_{0} + 1 & 2 \, \theta_{0} + 3 \\
    	    	2 \, \mu_{20} & \overline{A_{0}} & -\theta_{0} + 1 & 2 \, \theta_{0} + 3 \\
    	    	2 \, \mu_{9} & \overline{A_{0}} & -\theta_{0} + 2 & 2 \, \theta_{0} + 1 \\
    	    	-\frac{{\left(\theta_{0} - 3\right)} \overline{\zeta_{0}}}{\theta_{0}^{4} + 3 \, \theta_{0}^{3} + 2 \, \theta_{0}^{2}} & C_{2} & -\theta_{0} + 2 & 2 \, \theta_{0} + 2 \\
    	    	2 \, \mu_{9} & \overline{A_{1}} & -\theta_{0} + 2 & 2 \, \theta_{0} + 2 \\
    	    	2 \, \mu_{21} & \overline{A_{0}} & -\theta_{0} + 2 & 2 \, \theta_{0} + 2 
    	    	\end{dmatrix}
    	    	\end{align*}
    	    	\begin{align*}
    	    	\begin{dmatrix}
    	    	\frac{2 \, {\left(\alpha_{7} \theta_{0}^{2} - 6 \, \alpha_{7} \theta_{0} + 8 \, \alpha_{7}\right)}}{\theta_{0}^{2}} & \overline{A_{1}} & -\theta_{0} + 3 & 2 \, \theta_{0} + 1 \\
    	    	\frac{\theta_{0} - 2}{2 \, \theta_{0}^{2} + \theta_{0}} & E_{2} & -\theta_{0} + 3 & 2 \, \theta_{0} + 1 \\
    	    	\frac{2 \, \alpha_{2} \theta_{0}^{3} - 7 \, \alpha_{2} \theta_{0}^{2} + 2 \, \alpha_{2} \theta_{0} + 8 \, \alpha_{2}}{4 \, {\left(2 \, \theta_{0}^{3} + \theta_{0}^{2}\right)}} & C_{1} & -\theta_{0} + 3 & 2 \, \theta_{0} + 1 \\
    	    	2 \, \mu_{22} & \overline{A_{0}} & -\theta_{0} + 3 & 2 \, \theta_{0} + 1 \\
    	    	\frac{\theta_{0} - 2}{2 \, \theta_{0}^{2}} & E_{1} & -\theta_{0} + 3 & 2 \, \theta_{0} \\
    	    	\frac{2 \, {\left(\alpha_{7} \theta_{0}^{2} - 6 \, \alpha_{7} \theta_{0} + 8 \, \alpha_{7}\right)}}{\theta_{0}^{2}} & \overline{A_{0}} & -\theta_{0} + 3 & 2 \, \theta_{0} \\
    	    	\frac{2 \, \theta_{0} - 5}{4 \, \theta_{0}^{2}} & E_{3} & -\theta_{0} + 4 & 2 \, \theta_{0} \\
    	    	2 \, \mu_{23} & \overline{A_{0}} & -\theta_{0} + 4 & 2 \, \theta_{0} \\
    	    	-\frac{2}{\theta_{0}} & A_{1} & \theta_{0} & 0 \\
    	    	\frac{4 \, {\left| A_{1} \right|}^{2}}{\theta_{0}} & A_{0} & \theta_{0} & 1 \\
    	    	-\frac{1}{4 \, \theta_{0}} & \overline{B_{1}} & \theta_{0} & 2 \\
    	    	\frac{2 \, \overline{\alpha_{5}}}{\theta_{0}} & A_{0} & \theta_{0} & 2 \\
    	    	-\frac{2 \, {\left(\theta_{0} \overline{\alpha_{2}} + \overline{\alpha_{2}}\right)}}{\theta_{0}} & \overline{A_{0}} & \theta_{0} & 2 \\
    	    	-\frac{1}{6 \, \theta_{0}} & \overline{B_{3}} & \theta_{0} & 3 \\
    	    	\frac{\overline{\zeta_{2}}}{\theta_{0}^{3} + \theta_{0}^{2}} & A_{1} & \theta_{0} & 3 \\
    	    	-\frac{2 \, {\left(\theta_{0} \overline{\alpha_{2}} + \overline{\alpha_{2}}\right)}}{\theta_{0}} & \overline{A_{1}} & \theta_{0} & 3 \\
    	    	-\frac{2 \, {\left(2 \, {\left| A_{1} \right|}^{2} \overline{\alpha_{1}} - \overline{\alpha_{6}}\right)}}{\theta_{0}} & A_{0} & \theta_{0} & 3 \\
    	    	2 \, \mu_{10} & \overline{A_{0}} & \theta_{0} & 3 \\
    	    	\frac{{\left(\theta_{0}^{2} + \theta_{0} - 6\right)} {\left| A_{1} \right|}^{2}}{6 \, {\left(\theta_{0}^{3} + \theta_{0}^{2}\right)}} & \overline{C_{1}} & \theta_{0} & 3 \\
    	    	-\frac{1}{8 \, \theta_{0}} & \overline{B_{6}} & \theta_{0} & 4 \\
    	    	-\frac{\overline{\alpha_{1}}}{8 \, \theta_{0}} & \overline{B_{1}} & \theta_{0} & 4 \\
    	    	\frac{4 \, \mu_{1} + \overline{\alpha_{5}}}{8 \, \theta_{0}} & \overline{C_{1}} & \theta_{0} & 4 \\
    	    	-\frac{2 \, {\left({\left(\theta_{0}^{2} + 4 \, \theta_{0}\right)} \overline{\zeta_{8}} - \overline{\alpha_{4}}\right)}}{\theta_{0}} & A_{1} & \theta_{0} & 4 \\
    	    	-\frac{2 \, {\left(\theta_{0} \overline{\alpha_{2}} + \overline{\alpha_{2}}\right)}}{\theta_{0}} & \overline{A_{2}} & \theta_{0} & 4 \\
    	    	\pi_2 & A_{0} & \theta_{0} & 4 \\
    	    	2 \, \mu_{24} & \overline{A_{0}} & \theta_{0} & 4 \\
    	    	2 \, \mu_{10} & \overline{A_{1}} & \theta_{0} & 4 \\
    	    	\frac{{\left(\theta_{0}^{2} + \theta_{0} - 12\right)} {\left| A_{1} \right|}^{2}}{12 \, {\left(\theta_{0}^{3} + \theta_{0}^{2}\right)}} & \overline{C_{2}} & \theta_{0} & 4 \\
    	    	2 \, \mu_{25} & \overline{A_{0}} & \theta_{0} - 1 & 5 
    	    	\end{dmatrix}
    	    	\end{align*}
    	    	\begin{align*}
    	    	\begin{dmatrix}
    	    	-\frac{4}{\theta_{0}} & A_{2} & \theta_{0} + 1 & 0 \\
    	    	-\frac{4 \, \zeta_{0}}{\theta_{0}^{2} + \theta_{0}} & \overline{A_{0}} & \theta_{0} + 1 & 0 \\
    	    	\frac{4 \, \alpha_{1}}{\theta_{0}} & A_{0} & \theta_{0} + 1 & 0 \\
    	    	-\frac{4 \, \zeta_{0}}{\theta_{0}^{2} + \theta_{0}} & \overline{A_{1}} & \theta_{0} + 1 & 1 \\
    	    	\frac{4 \, \alpha_{5}}{\theta_{0}} & A_{0} & \theta_{0} + 1 & 1 \\
    	    	\frac{4 \, {\left(\theta_{0} + 2\right)} {\left| A_{1} \right|}^{2}}{\theta_{0}^{2} + \theta_{0}} & A_{1} & \theta_{0} + 1 & 1 \\
    	    	-\frac{1}{2 \, \theta_{0}} & \overline{B_{2}} & \theta_{0} + 1 & 2 \\
    	    	-\frac{4 \, \zeta_{0}}{\theta_{0}^{2} + \theta_{0}} & \overline{A_{2}} & \theta_{0} + 1 & 2 \\
    	    	\frac{2 \, {\left(\theta_{0} \overline{\alpha_{5}} + 2 \, \overline{\alpha_{5}}\right)}}{\theta_{0}^{2} + \theta_{0}} & A_{1} & \theta_{0} + 1 & 2 \\
    	    	-\frac{4 \, {\left(2 \, {\left(\theta_{0} + 2\right)} {\left| A_{1} \right|}^{4} + {\left(\alpha_{1} \overline{\alpha_{1}} - \beta\right)} \theta_{0} + \alpha_{1} \overline{\alpha_{1}} - \beta\right)}}{\theta_{0}^{2} + \theta_{0}} & A_{0} & \theta_{0} + 1 & 2 \\
    	    	2 \, \mu_{11} & \overline{A_{0}} & \theta_{0} + 1 & 2 \\
    	    	-\frac{1}{3 \, \theta_{0}} & \overline{B_{5}} & \theta_{0} + 1 & 3 \\
    	    	\frac{2 \, \overline{\zeta_{2}}}{\theta_{0}^{3} + \theta_{0}^{2}} & A_{2} & \theta_{0} + 1 & 3 \\
    	    	-\frac{4 \, \zeta_{0}}{\theta_{0}^{2} + \theta_{0}} & \overline{A_{3}} & \theta_{0} + 1 & 3 \\
    	    	\frac{\alpha_{5} + 3 \, \mu_{3}}{6 \, \theta_{0}} & \overline{C_{1}} & \theta_{0} + 1 & 3 \\
    	    	-\frac{4 \, {\left({\left(\theta_{0} \overline{\alpha_{1}} + 2 \, \overline{\alpha_{1}}\right)} {\left| A_{1} \right|}^{2} - \theta_{0} \overline{\alpha_{6}} + \theta_{0} \overline{\zeta_{11}} - \overline{\alpha_{6}}\right)}}{\theta_{0}^{2} + \theta_{0}} & A_{1} & \theta_{0} + 1 & 3 \\
    	    	-\frac{2 \, {\left(2 \, {\left(\alpha_{5} \overline{\alpha_{1}} + \alpha_{1} \overline{\alpha_{3}} - \alpha_{16}\right)} \theta_{0}^{2} + 4 \, {\left(\theta_{0}^{2} \overline{\alpha_{5}} + 2 \, \theta_{0} \overline{\alpha_{5}}\right)} {\left| A_{1} \right|}^{2} + 2 \, {\left(\alpha_{5} \overline{\alpha_{1}} + \alpha_{1} \overline{\alpha_{3}} - \alpha_{16}\right)} \theta_{0} + \alpha_{1} \overline{\zeta_{2}}\right)}}{\theta_{0}^{3} + \theta_{0}^{2}} & A_{0} & \theta_{0} + 1 & 3 \\
    	    	2 \, \mu_{26} & \overline{A_{0}} & \theta_{0} + 1 & 3 \\
    	    	2 \, \mu_{11} & \overline{A_{1}} & \theta_{0} + 1 & 3 \\
    	    	-\frac{{\left(\theta_{0}^{2} - \theta_{0} + 4\right)} {\left| A_{1} \right|}^{2}}{6 \, {\left(\theta_{0}^{3} + 2 \, \theta_{0}^{2} + \theta_{0}\right)}} & \overline{B_{1}} & \theta_{0} + 1 & 3 \\
    	    	-\frac{6}{\theta_{0}} & A_{3} & \theta_{0} + 2 & 0 \\
    	    	\frac{6 \, \alpha_{3}}{\theta_{0}} & A_{0} & \theta_{0} + 2 & 0 \\
    	    	\frac{4 \, {\left(\alpha_{1} \theta_{0} + 3 \, \alpha_{1}\right)}}{\theta_{0}^{2} + 2 \, \theta_{0}} & A_{1} & \theta_{0} + 2 & 0 \\
    	    	2 \, \mu_{12} & \overline{A_{0}} & \theta_{0} + 2 & 0 \\
    	    	\frac{4 \, {\left(\alpha_{5} \theta_{0} + 3 \, \alpha_{5}\right)}}{\theta_{0}^{2} + 2 \, \theta_{0}} & A_{1} & \theta_{0} + 2 & 1 \\
    	    	\frac{2 \, {\left(3 \, \alpha_{6} \theta_{0}^{2} - 2 \, {\left(3 \, \alpha_{1} \theta_{0}^{2} + 13 \, \alpha_{1} \theta_{0} + 12 \, \alpha_{1}\right)} {\left| A_{1} \right|}^{2} + 9 \, \alpha_{6} \theta_{0} + 6 \, \alpha_{6}\right)}}{\theta_{0}^{3} + 3 \, \theta_{0}^{2} + 2 \, \theta_{0}} & A_{0} & \theta_{0} + 2 & 1 \\
    	    	2 \, \mu_{13} & \overline{A_{0}} & \theta_{0} + 2 & 1 \\
    	    	2 \, \mu_{12} & \overline{A_{1}} & \theta_{0} + 2 & 1 \\
    	    	\frac{4 \, {\left(\theta_{0} + 3\right)} {\left| A_{1} \right|}^{2}}{\theta_{0}^{2} + \theta_{0}} & A_{2} & \theta_{0} + 2 & 1 
    	    	\end{dmatrix}
    	    	\end{align*}
    	    	\begin{align*}
    	    	\begin{dmatrix}
    	    	-\frac{3}{4 \, \theta_{0}} & \overline{B_{4}} & \theta_{0} + 2 & 2 \\
    	    	-\frac{\theta_{0}^{2} \overline{\alpha_{2}} + 4 \, \theta_{0} \overline{\alpha_{2}} - 6 \, \overline{\alpha_{2}}}{4 \, \theta_{0}^{2}} & C_{1} & \theta_{0} + 2 & 2 \\
    	    	\frac{\mu_{5}}{2 \, \theta_{0}} & \overline{C_{1}} & \theta_{0} + 2 & 2 \\
    	    	-\frac{\alpha_{1}}{4 \, {\left(\theta_{0} + 2\right)}} & \overline{B_{1}} & \theta_{0} + 2 & 2 \\
    	    	\frac{2 \, {\left(\theta_{0} \overline{\alpha_{5}} + 3 \, \overline{\alpha_{5}}\right)}}{\theta_{0}^{2} + \theta_{0}} & A_{2} & \theta_{0} + 2 & 2 \\
    	    	\pi_3 & A_{0} & \theta_{0} + 2 & 2 \\
    	    	\pi_4 & A_{1} & \theta_{0} + 2 & 2 \\
    	    	2 \, \mu_{27} & \overline{A_{0}} & \theta_{0} + 2 & 2 \\
    	    	2 \, \mu_{13} & \overline{A_{1}} & \theta_{0} + 2 & 2 \\
    	    	2 \, \mu_{12} & \overline{A_{2}} & \theta_{0} + 2 & 2 \\
    	    	-\frac{8}{\theta_{0}} & A_{4} & \theta_{0} + 3 & 0 \\
    	    	-\frac{\zeta_{0}}{2 \, {\left(\theta_{0}^{2} + \theta_{0}\right)}} & C_{1} & \theta_{0} + 3 & 0 \\
    	    	\frac{4 \, {\left(\alpha_{1} \theta_{0} + 4 \, \alpha_{1}\right)}}{\theta_{0}^{2} + 2 \, \theta_{0}} & A_{2} & \theta_{0} + 3 & 0 \\
    	    	-\frac{4 \, {\left(4 \, \alpha_{1}^{2} + {\left(\alpha_{1}^{2} - 2 \, \alpha_{4}\right)} \theta_{0} - 4 \, \alpha_{4}\right)}}{\theta_{0}^{2} + 2 \, \theta_{0}} & A_{0} & \theta_{0} + 3 & 0 \\
    	    	\frac{6 \, \alpha_{3} \theta_{0}^{2} + 30 \, \alpha_{3} \theta_{0} + 24 \, \alpha_{3} + \zeta_{2}}{\theta_{0}^{3} + 4 \, \theta_{0}^{2} + 3 \, \theta_{0}} & A_{1} & \theta_{0} + 3 & 0 \\
    	    	2 \, \mu_{14} & \overline{A_{0}} & \theta_{0} + 3 & 0 \\
    	    	-\frac{\zeta_{0}}{2 \, {\left(\theta_{0}^{2} + \theta_{0}\right)}} & B_{1} & \theta_{0} + 3 & 1 \\
    	    	\frac{4 \, {\left(\alpha_{5} \theta_{0} + 4 \, \alpha_{5}\right)}}{\theta_{0}^{2} + 2 \, \theta_{0}} & A_{2} & \theta_{0} + 3 & 1 \\
    	    	\pi_5 & A_{1} & \theta_{0} + 3 & 1 \\
    	    	\pi_6 & A_{0} & \theta_{0} + 3 & 1 \\
    	    	2 \, \mu_{28} & \overline{A_{0}} & \theta_{0} + 3 & 1 \\
    	    	2 \, \mu_{14} & \overline{A_{1}} & \theta_{0} + 3 & 1 \\
    	    	\frac{4 \, {\left(\theta_{0} + 4\right)} {\left| A_{1} \right|}^{2}}{\theta_{0}^{2} + \theta_{0}} & A_{3} & \theta_{0} + 3 & 1 \\
    	    	-\frac{10}{\theta_{0}} & A_{5} & \theta_{0} + 4 & 0 \\
    	    	-\frac{\zeta_{0}}{3 \, {\left(\theta_{0}^{2} + \theta_{0}\right)}} & C_{2} & \theta_{0} + 4 & 0 \\
    	    	\frac{4 \, {\left(\alpha_{1} \theta_{0} + 5 \, \alpha_{1}\right)}}{\theta_{0}^{2} + 2 \, \theta_{0}} & A_{3} & \theta_{0} + 4 & 0 \\
    	    	-\frac{2 \, {\left(10 \, \alpha_{1}^{2} + {\left(2 \, \alpha_{1}^{2} - 5 \, \alpha_{4}\right)} \theta_{0} + {\left(\theta_{0}^{3} + 2 \, \theta_{0}^{2}\right)} \zeta_{8} - 10 \, \alpha_{4}\right)}}{\theta_{0}^{2} + 2 \, \theta_{0}} & A_{1} & \theta_{0} + 4 & 0 \\
    	    	\frac{2 \, {\left(3 \, \alpha_{3} \theta_{0}^{2} + 18 \, \alpha_{3} \theta_{0} + 15 \, \alpha_{3} + \zeta_{2}\right)}}{\theta_{0}^{3} + 4 \, \theta_{0}^{2} + 3 \, \theta_{0}} & A_{2} & \theta_{0} + 4 & 0 \\
    	    	-\frac{2 \, {\left(5 \, {\left(\alpha_{1} \alpha_{3} - \overline{\alpha_{14}}\right)} \theta_{0}^{3} + 6 \, {\left(7 \, \alpha_{1} \alpha_{3} - 5 \, \overline{\alpha_{14}}\right)} \theta_{0}^{2} + 60 \, \alpha_{1} \alpha_{3} + {\left(97 \, \alpha_{1} \alpha_{3} - 55 \, \overline{\alpha_{14}}\right)} \theta_{0} + {\left(\alpha_{1} \theta_{0} + 2 \, \alpha_{1}\right)} \zeta_{2} - 30 \, \overline{\alpha_{14}}\right)}}{\theta_{0}^{4} + 6 \, \theta_{0}^{3} + 11 \, \theta_{0}^{2} + 6 \, \theta_{0}} & A_{0} & \theta_{0} + 4 & 0 \\
    	    	2 \, \mu_{29} & \overline{A_{0}} & \theta_{0} + 4 & 0 \\
    	    	\frac{1}{4} \, \mu_{12} & C_{1} & \theta_{0} + 4 & 0 
    	    	\end{dmatrix}
    	    	\end{align*}
    	    	\begin{align*}
    	    	\begin{dmatrix}
    	    	\frac{1}{2 \, {\left(\theta_{0} - 4\right)}} & \overline{E_{1}} & 2 \, \theta_{0} - 1 & -\theta_{0} + 4 \\
    	    	2 \, \overline{\alpha_{7}} & A_{0} & 2 \, \theta_{0} - 1 & -\theta_{0} + 4 \\
    	    	\frac{1}{2 \, {\left(\theta_{0} - 5\right)}} & \overline{E_{3}} & 2 \, \theta_{0} - 1 & -\theta_{0} + 5 \\
    	    	2 \, \overline{\alpha_{13}} & A_{0} & 2 \, \theta_{0} - 1 & -\theta_{0} + 5 \\
    	    	-\frac{\zeta_{0}}{\theta_{0}^{3} + \theta_{0}^{2}} & \overline{C_{1}} & 2 \, \theta_{0} + 1 & -\theta_{0} + 2 \\
    	    	\frac{2 \, {\left(\theta_{0} \overline{\alpha_{9}} + 2 \, \overline{\alpha_{9}}\right)}}{\theta_{0}} & A_{0} & 2 \, \theta_{0} + 1 & -\theta_{0} + 2 \\
    	    	\frac{2 \, {\left(\theta_{0} \overline{\alpha_{2}} + 2 \, \overline{\alpha_{2}}\right)}}{\theta_{0}} & A_{1} & 2 \, \theta_{0} + 1 & -\theta_{0} + 2 \\
    	    	-\frac{\zeta_{0}}{\theta_{0}^{3} + \theta_{0}^{2}} & \overline{C_{2}} & 2 \, \theta_{0} + 1 & -\theta_{0} + 3 \\
    	    	\frac{2 \, {\left(\theta_{0} \overline{\alpha_{8}} + 2 \, \overline{\alpha_{8}}\right)}}{\theta_{0}} & A_{1} & 2 \, \theta_{0} + 1 & -\theta_{0} + 3 \\
    	    	-\frac{2 \, {\left(2 \, {\left(\theta_{0} \overline{\alpha_{2}} + 4 \, \overline{\alpha_{2}}\right)} {\left| A_{1} \right|}^{2} - \theta_{0} \overline{\alpha_{11}} - 2 \, \overline{\alpha_{11}}\right)}}{\theta_{0}} & A_{0} & 2 \, \theta_{0} + 1 & -\theta_{0} + 3 \\
    	    	-\frac{\zeta_{0}}{\theta_{0}^{3} + 2 \, \theta_{0}^{2} + \theta_{0}} & \overline{B_{1}} & 2 \, \theta_{0} + 2 & -\theta_{0} + 2 \\
    	    	\frac{\mu_{12}}{2 \, \theta_{0}} & \overline{C_{1}} & 2 \, \theta_{0} + 2 & -\theta_{0} + 2 \\
    	    	\frac{2 \, {\left(\theta_{0} \overline{\alpha_{9}} + 3 \, \overline{\alpha_{9}}\right)}}{\theta_{0}} & A_{1} & 2 \, \theta_{0} + 2 & -\theta_{0} + 2 \\
    	    	\frac{2 \, {\left(\theta_{0} \overline{\alpha_{2}} + 3 \, \overline{\alpha_{2}}\right)}}{\theta_{0}} & A_{2} & 2 \, \theta_{0} + 2 & -\theta_{0} + 2 \\
    	    	-\frac{2 \, {\left({\left(\alpha_{1} \overline{\alpha_{2}} - \overline{\alpha_{10}}\right)} \theta_{0}^{2} + {\left(9 \, \alpha_{1} \overline{\alpha_{2}} - 5 \, \overline{\alpha_{10}}\right)} \theta_{0} + 12 \, \alpha_{1} \overline{\alpha_{2}} - 6 \, \overline{\alpha_{10}}\right)}}{\theta_{0}^{2} + 2 \, \theta_{0}} & A_{0} & 2 \, \theta_{0} + 2 & -\theta_{0} + 2 \\
    	    	2 \, \mu_{30} & \overline{A_{0}} & 2 \, \theta_{0} + 2 & -\theta_{0} + 2 \\
    	    	\frac{2 \, {\left(\theta_{0} \overline{\alpha_{2}} + \overline{\alpha_{2}}\right)}}{\theta_{0}} & A_{0} & 2 \, \theta_{0} & -\theta_{0} + 2 \\
    	    	\frac{2 \, {\left(\theta_{0} \overline{\alpha_{8}} + \overline{\alpha_{8}}\right)}}{\theta_{0}} & A_{0} & 2 \, \theta_{0} & -\theta_{0} + 3 \\
    	    	\frac{2 \, {\left(\theta_{0}^{2} \overline{\alpha_{7}} + \theta_{0} \overline{\alpha_{7}} - 2 \, \overline{\alpha_{7}}\right)}}{\theta_{0}^{2}} & A_{1} & 2 \, \theta_{0} & -\theta_{0} + 4 \\
    	    	-\frac{2 \, {\left({\left(\overline{\alpha_{1}} \overline{\alpha_{2}} - \overline{\alpha_{12}}\right)} \theta_{0} + \overline{\alpha_{1}} \overline{\alpha_{2}} - \overline{\alpha_{12}}\right)}}{\theta_{0}} & A_{0} & 2 \, \theta_{0} & -\theta_{0} + 4 \\
    	    	\frac{\theta_{0} + 1}{2 \, {\left(\theta_{0}^{2} - 4 \, \theta_{0}\right)}} & \overline{E_{2}} & 2 \, \theta_{0} & -\theta_{0} + 4 \\
    	    	\frac{\theta_{0}^{3} \overline{\alpha_{2}} - 3 \, \theta_{0}^{2} \overline{\alpha_{2}} + 4 \, \theta_{0} \overline{\alpha_{2}} + 8 \, \overline{\alpha_{2}}}{4 \, {\left(\theta_{0}^{3} - 4 \, \theta_{0}^{2}\right)}} & \overline{C_{1}} & 2 \, \theta_{0} & -\theta_{0} + 4 \\
    	    	2 \, \mu_{31} & \overline{A_{0}} & 2 \, \theta_{0} & -\theta_{0} + 4
    	    	\end{dmatrix}
    	    \end{align*}
    	    \small
    	    where
    	    \begin{align*}
    	    	\pi_{1}&=\frac{2}{\theta_{0}^{4} + 6 \, \theta_{0}^{3} + 11 \, \theta_{0}^{2} + 6 \, \theta_{0}}\bigg\{ {\left(\alpha_{2} \overline{\alpha_{1}} - \alpha_{10}\right)} \theta_{0}^{4} + 4 \, {\left(\alpha_{2} \overline{\alpha_{1}} - \alpha_{10}\right)} \theta_{0}^{3} - {\left(\alpha_{2} \overline{\alpha_{1}} - \alpha_{10}\right)} \theta_{0}^{2} - 16 \, {\left(\alpha_{2} \overline{\alpha_{1}} - \alpha_{10}\right)} \theta_{0} - 12 \, \alpha_{2} \overline{\alpha_{1}}\\
    	    	& - 4 \, {\left(\alpha_{5} \theta_{0} + 3 \, \alpha_{5}\right)} \overline{\zeta_{0}} - 8 \, {\left(2 \, \alpha_{1} \theta_{0} + 3 \, \alpha_{1}\right)} \overline{\zeta_{1}} + 12 \, \alpha_{10}\bigg\}\\
    	    	\pi_{2}&=-\frac{2}{\theta_{0}^{4} + 3 \, \theta_{0}^{3} + 2 \, \theta_{0}^{2}}\bigg\{{\left(\overline{\alpha_{1}} \overline{\alpha_{5}} - \alpha_{15}\right)} \theta_{0}^{3} + 3 \, {\left(\overline{\alpha_{1}} \overline{\alpha_{5}} - \alpha_{15}\right)} \theta_{0}^{2} + {\left(2 \, \theta_{0}^{3} \overline{\alpha_{3}} + 6 \, \theta_{0}^{2} \overline{\alpha_{3}} + 4 \, \theta_{0} \overline{\alpha_{3}} + {\left(\theta_{0} + 2\right)} \overline{\zeta_{2}}\right)} {\left| A_{1} \right|}^{2} \\
    	    	&+ 2 \, {\left(\overline{\alpha_{1}} \overline{\alpha_{5}} - \alpha_{15}\right)} \theta_{0} + 2 \, {\left(\theta_{0}^{2} \overline{\alpha_{2}} + \theta_{0} \overline{\alpha_{2}}\right)} \overline{\zeta_{0}}\bigg\}\\
    	    	\pi_{3}&=-\frac{2}{\theta_{0}^{3} + 3 \, \theta_{0}^{2} + 2 \, \theta_{0}}\bigg\{3 \, {\left(\alpha_{3} \overline{\alpha_{1}} + \alpha_{1} \overline{\alpha_{5}} - \overline{\alpha_{16}}\right)} \theta_{0}^{2} + 2 \, {\left(3 \, \alpha_{5} \theta_{0}^{2} + 13 \, \alpha_{5} \theta_{0} + 12 \, \alpha_{5}\right)} {\left| A_{1} \right|}^{2}\\
    	    	& + {\left(9 \, \alpha_{3} \overline{\alpha_{1}} + 13 \, \alpha_{1} \overline{\alpha_{5}} - 9 \, \overline{\alpha_{16}}\right)} \theta_{0} + 6 \, \alpha_{3} \overline{\alpha_{1}} + 12 \, \alpha_{1} \overline{\alpha_{5}} - 6 \, \overline{\alpha_{16}}\bigg\}\\
    	    	\pi_{4}&=-\frac{2 \, {\left(4 \, {\left(\theta_{0}^{2} + 5 \, \theta_{0} + 6\right)} {\left| A_{1} \right|}^{4} + {\left(2 \, \alpha_{1} \overline{\alpha_{1}} - 3 \, \beta\right)} \theta_{0}^{2} + {\left(8 \, \alpha_{1} \overline{\alpha_{1}} - 9 \, \beta\right)} \theta_{0} + {\left(\theta_{0}^{4} + 5 \, \theta_{0}^{3} + 8 \, \theta_{0}^{2} + 4 \, \theta_{0}\right)} \zeta_{10} + 6 \, \alpha_{1} \overline{\alpha_{1}} - 6 \, \beta\right)}}{\theta_{0}^{3} + 3 \, \theta_{0}^{2} + 2 \, \theta_{0}}\\
    	    	\pi_{5}&=\frac{4 \, {\left(2 \, \alpha_{6} \theta_{0}^{3} + 12 \, \alpha_{6} \theta_{0}^{2} - {\left(3 \, \alpha_{1} \theta_{0}^{3} + 25 \, \alpha_{1} \theta_{0}^{2} + 64 \, \alpha_{1} \theta_{0} + 48 \, \alpha_{1}\right)} {\left| A_{1} \right|}^{2} + 22 \, \alpha_{6} \theta_{0} - {\left(\theta_{0}^{3} + 3 \, \theta_{0}^{2} + 2 \, \theta_{0}\right)} \zeta_{11} + 12 \, \alpha_{6}\right)}}{\theta_{0}^{4} + 6 \, \theta_{0}^{3} + 11 \, \theta_{0}^{2} + 6 \, \theta_{0}}\\
    	    	\pi_{6}&=-\frac{2}{\theta_{0}^{4} + 6 \, \theta_{0}^{3} + 11 \, \theta_{0}^{2} + 6 \, \theta_{0}}\bigg\{4 \, {\left(\alpha_{1} \alpha_{5} - \overline{\alpha_{15}}\right)} \theta_{0}^{3}\\
    	    	& + 8 \, {\left(4 \, \alpha_{1} \alpha_{5} - 3 \, \overline{\alpha_{15}}\right)} \theta_{0}^{2} + {\left(8 \, \alpha_{3} \theta_{0}^{3} + 60 \, \alpha_{3} \theta_{0}^{2} + 136 \, \alpha_{3} \theta_{0} + {\left(\theta_{0} + 2\right)} \zeta_{2} + 96 \, \alpha_{3}\right)} {\left| A_{1} \right|}^{2} + 48 \, \alpha_{1} \alpha_{5}\\
    	    	& + 4 \, {\left(19 \, \alpha_{1} \alpha_{5} - 11 \, \overline{\alpha_{15}}\right)} \theta_{0} - 24 \, \overline{\alpha_{15}}\bigg\}
    	    \end{align*}
    	    Finally,
    	    \begin{align*}
    	    	&\Re(\p{\z}\vec{F}(z))=\\
    	    	&\begin{dmatrix}
    	    	\frac{\mu_{1} \theta_{0}^{2} + 2 \, \mu_{1} \theta_{0} + 2 \, \overline{\alpha_{5}}}{\theta_{0}} & \overline{A_{0}} & 0 & \theta_{0} + 1 \\
    	    	\frac{\mu_{1} \theta_{0}^{3} + 5 \, \mu_{1} \theta_{0}^{2} + 2 \, {\left(3 \, \mu_{1} + \overline{\alpha_{5}}\right)} \theta_{0} + 6 \, \overline{\alpha_{5}}}{\theta_{0}^{2} + 2 \, \theta_{0}} & \overline{A_{1}} & 0 & \theta_{0} + 2 \\
    	    	\frac{\theta_{0}^{3} \overline{\mu_{13}} - 4 \, {\left(\theta_{0} + 3\right)} {\left| A_{1} \right|}^{2} \overline{\zeta_{0}} + 3 \, \theta_{0}^{2} \overline{\mu_{13}} + 2 \, \theta_{0} \overline{\mu_{13}}}{\theta_{0}^{3} + 3 \, \theta_{0}^{2} + 2 \, \theta_{0}} & A_{0} & 0 & \theta_{0} + 2 \\
    	    	\frac{\theta_{0}^{3} \overline{\mu_{12}} + 3 \, \theta_{0}^{2} \overline{\mu_{12}} + 2 \, \theta_{0} \overline{\mu_{12}} + 4 \, {\left(2 \, \theta_{0} + 3\right)} \overline{\zeta_{1}}}{\theta_{0}^{3} + 3 \, \theta_{0}^{2} + 2 \, \theta_{0}} & A_{1} & 0 & \theta_{0} + 2 \\
    	    	\frac{\mu_{2} \theta_{0}^{4} + 6 \, \mu_{2} \theta_{0}^{3} + {\left(11 \, \mu_{2} + 3 \, \overline{\alpha_{6}}\right)} \theta_{0}^{2} - 2 \, {\left(3 \, \theta_{0}^{2} \overline{\alpha_{1}} + 13 \, \theta_{0} \overline{\alpha_{1}} + 12 \, \overline{\alpha_{1}}\right)} {\left| A_{1} \right|}^{2} + 3 \, {\left(2 \, \mu_{2} + 3 \, \overline{\alpha_{6}}\right)} \theta_{0} + 6 \, \overline{\alpha_{6}}}{\theta_{0}^{3} + 3 \, \theta_{0}^{2} + 2 \, \theta_{0}} & \overline{A_{0}} & 0 & \theta_{0} + 2 \\
    	    	\frac{3 \, \overline{\zeta_{0}}}{4 \, {\left(\theta_{0}^{3} + 4 \, \theta_{0}^{2} + 5 \, \theta_{0} + 2\right)}} & \overline{B_{1}} & 0 & \theta_{0} + 3 \\
    	    	\frac{\mu_{1} \theta_{0}^{3} + 6 \, \mu_{1} \theta_{0}^{2} + 2 \, {\left(4 \, \mu_{1} + \overline{\alpha_{5}}\right)} \theta_{0} + 8 \, \overline{\alpha_{5}}}{\theta_{0}^{2} + 2 \, \theta_{0}} & \overline{A_{2}} & 0 & \theta_{0} + 3 \\
    	    	\frac{\theta_{0}^{3} \overline{\mu_{14}} + 3 \, \theta_{0}^{2} \overline{\mu_{14}} + 2 \, \theta_{0} \overline{\mu_{14}} - 2 \, {\left(\theta_{0} \overline{\alpha_{1}} + 4 \, \overline{\alpha_{1}}\right)} \overline{\zeta_{0}} - 2 \, {\left(\theta_{0}^{5} + 9 \, \theta_{0}^{4} + 28 \, \theta_{0}^{3} + 36 \, \theta_{0}^{2} + 16 \, \theta_{0}\right)} \overline{\zeta_{9}}}{\theta_{0}^{3} + 3 \, \theta_{0}^{2} + 2 \, \theta_{0}} & A_{1} & 0 & \theta_{0} + 3 \\
    	    	\frac{\theta_{0}^{4} \overline{\mu_{28}} + 6 \, \theta_{0}^{3} \overline{\mu_{28}} - 8 \, {\left(2 \, \theta_{0}^{2} + 11 \, \theta_{0} + 12\right)} {\left| A_{1} \right|}^{2} \overline{\zeta_{1}} + 11 \, \theta_{0}^{2} \overline{\mu_{28}} + 6 \, \theta_{0} \overline{\mu_{28}} - 2 \, {\left(\theta_{0}^{2} \overline{\alpha_{5}} + 7 \, \theta_{0} \overline{\alpha_{5}} + 12 \, \overline{\alpha_{5}}\right)} \overline{\zeta_{0}}}{\theta_{0}^{4} + 6 \, \theta_{0}^{3} + 11 \, \theta_{0}^{2} + 6 \, \theta_{0}} & A_{0} & 0 & \theta_{0} + 3 \\
    	    	\lambda_1 & \overline{A_{1}} & 0 & \theta_{0} + 3 \\
    	    	\lambda_2 & \overline{A_{0}} & 0 & \theta_{0} + 3 \\
    	    	-\frac{\theta_{0} - 1}{2 \, \theta_{0}^{2}} & B_{1} & 1 & \theta_{0} \\
    	    	\frac{\mu_{3} \theta_{0}^{2} + \mu_{3} \theta_{0} + 2 \, \alpha_{5}}{\theta_{0}} & \overline{A_{0}} & 1 & \theta_{0} \\
    	    	-\frac{2 \, {\left(\alpha_{2} \theta_{0}^{2} - \alpha_{2}\right)}}{\theta_{0}} & A_{0} & 1 & \theta_{0} \\
    	    	-\frac{\theta_{0} - 1}{\theta_{0}^{2}} & B_{2} & 1 & \theta_{0} + 1 \\
    	    	\frac{\mu_{3} \theta_{0}^{3} + 3 \, \mu_{3} \theta_{0}^{2} + 2 \, {\left(\alpha_{5} + \mu_{3}\right)} \theta_{0} + 4 \, \alpha_{5}}{\theta_{0}^{2} + \theta_{0}} & \overline{A_{1}} & 1 & \theta_{0} + 1 \\
    	    	-\frac{\alpha_{9} \theta_{0}^{4} + {\left(\alpha_{9} - 2 \, \overline{\mu_{11}}\right)} \theta_{0}^{3} - 2 \, {\left(2 \, \alpha_{9} + \overline{\mu_{11}}\right)} \theta_{0}^{2} + 4 \, \alpha_{1} \theta_{0} \overline{\zeta_{0}} - 4 \, \alpha_{9} \theta_{0} + {\left(\theta_{0}^{2} - \theta_{0} - 2\right)} \zeta_{4}}{\theta_{0}^{3} + \theta_{0}^{2}} & A_{0} & 1 & \theta_{0} + 1 \\
    	    	\frac{4 \, \mu_{4} \theta_{0}^{5} + 12 \, \mu_{4} \theta_{0}^{4} - 32 \, {\left(\theta_{0}^{3} + 2 \, \theta_{0}^{2}\right)} {\left| A_{1} \right|}^{4} - 8 \, {\left(2 \, \alpha_{1} \overline{\alpha_{1}} - 2 \, \beta - \mu_{4}\right)} \theta_{0}^{3} - 16 \, {\left(\alpha_{1} \overline{\alpha_{1}} - \beta\right)} \theta_{0}^{2} - {\left(\theta_{0}^{3} - 3 \, \theta_{0} - 2\right)} {\left| C_{1} \right|}^{2}}{4 \, {\left(\theta_{0}^{4} + \theta_{0}^{3}\right)}} & \overline{A_{0}} & 1 & \theta_{0} + 1 \\
    	    	\frac{1}{\theta_{0}^{2}} & \overline{B_{9}} & 1 & \theta_{0} + 2 \\
    	    	\frac{\theta_{0} + 3}{2 \, \theta_{0}^{2}} & B_{10} & 1 & \theta_{0} + 2 \\
    	    	\frac{\theta_{0} - 5}{4 \, \theta_{0}} & B_{4} & 1 & \theta_{0} + 2 
    	    	    	    	\end{dmatrix}
    	    	\end{align*}
    	    	\begin{align*}
    	    	\begin{dmatrix}
    	    	-\frac{\theta_{0} \overline{\alpha_{1}} - \overline{\alpha_{1}}}{2 \, {\left(\theta_{0}^{2} + 2 \, \theta_{0}\right)}} & B_{1} & 1 & \theta_{0} + 2 \\
    	    	\frac{\mu_{3} \theta_{0}^{3} + 4 \, \mu_{3} \theta_{0}^{2} + {\left(2 \, \alpha_{5} + 3 \, \mu_{3}\right)} \theta_{0} + 6 \, \alpha_{5}}{\theta_{0}^{2} + \theta_{0}} & \overline{A_{2}} & 1 & \theta_{0} + 2 \\
    	    	\frac{2 \, {\left(\theta_{0}^{3} \overline{\mu_{13}} - 4 \, {\left(\theta_{0} + 3\right)} {\left| A_{1} \right|}^{2} \overline{\zeta_{0}} + 3 \, \theta_{0}^{2} \overline{\mu_{13}} + 2 \, \theta_{0} \overline{\mu_{13}}\right)}}{\theta_{0}^{3} + 3 \, \theta_{0}^{2} + 2 \, \theta_{0}} & A_{1} & 1 & \theta_{0} + 2 \\
    	    	\frac{2 \, {\left(\theta_{0}^{3} \overline{\mu_{12}} + 3 \, \theta_{0}^{2} \overline{\mu_{12}} + 2 \, \theta_{0} \overline{\mu_{12}} + 4 \, {\left(2 \, \theta_{0} + 3\right)} \overline{\zeta_{1}}\right)}}{\theta_{0}^{3} + 3 \, \theta_{0}^{2} + 2 \, \theta_{0}} & A_{2} & 1 & \theta_{0} + 2 \\
    	    	\lambda_3 & \overline{A_{1}} & 1 & \theta_{0} + 2 \\
    	    	\lambda_4 & \overline{A_{0}} & 1 & \theta_{0} + 2 \\
    	    	\lambda_5 & A_{0} & 1 & \theta_{0} + 2 \\
    	    	\frac{2 \, \theta_{0}^{3} \overline{\mu_{5}} + 2 \, \theta_{0}^{2} \overline{\mu_{5}} - {\left(\theta_{0}^{2} + \theta_{0} - 3\right)} \overline{\zeta_{2}}}{4 \, {\left(\theta_{0}^{4} + \theta_{0}^{3}\right)}} & C_{1} & 1 & \theta_{0} + 2 \\
    	    	-\frac{1}{4} \, \alpha_{2} \theta_{0} - \frac{1}{4} \, \alpha_{2} & \overline{C_{1}} & 1 & \theta_{0} + 2 \\
    	    	-\frac{2 \, \theta_{0} - 3}{4 \, \theta_{0}^{2}} & B_{3} & 2 & \theta_{0} \\
    	    	-\frac{2 \, \alpha_{8} \theta_{0}^{3} - 2 \, {\left(2 \, \alpha_{8} + 3 \, \overline{\mu_{10}}\right)} \theta_{0}^{2} - 6 \, \alpha_{8} \theta_{0} + {\left(\theta_{0} - 3\right)} \zeta_{3}}{2 \, \theta_{0}^{2}} & A_{0} & 2 & \theta_{0} \\
    	    	\frac{\mu_{6} \theta_{0}^{2} - 6 \, \alpha_{1} {\left| A_{1} \right|}^{2} + \mu_{6} \theta_{0} + 3 \, \alpha_{6}}{\theta_{0}} & \overline{A_{0}} & 2 & \theta_{0} \\
    	    	-\frac{2 \, \alpha_{2} \theta_{0}^{2} - \alpha_{2} \theta_{0} - 3 \, \alpha_{2}}{\theta_{0}} & A_{1} & 2 & \theta_{0} \\
    	    	\frac{2 \, \mu_{5} \theta_{0}^{4} - 4 \, \mu_{5} \theta_{0}^{3} - 6 \, \mu_{5} \theta_{0}^{2} + {\left(5 \, \theta_{0} - 12\right)} \zeta_{2}}{2 \, {\left(\theta_{0}^{3} - 3 \, \theta_{0}^{2}\right)}} & \overline{A_{1}} & 2 & \theta_{0} \\
    	    	\frac{2 \, \mu_{5} \theta_{0}^{3} + 2 \, \mu_{5} \theta_{0}^{2} + 3 \, \zeta_{2}}{2 \, {\left(\theta_{0}^{2} + \theta_{0}\right)}} & \overline{A_{0}} & 2 & \theta_{0} - 1 \\
    	    	\frac{3}{2 \, \theta_{0}^{2}} & \overline{B_{8}} & 2 & \theta_{0} + 1 \\
    	    	\frac{\theta_{0} + 2}{2 \, \theta_{0}^{2}} & B_{11} & 2 & \theta_{0} + 1 \\
    	    	\frac{\theta_{0} - 5}{4 \, \theta_{0}} & B_{5} & 2 & \theta_{0} + 1 \\
    	    	\frac{\mu_{6} \theta_{0}^{3} + 3 \, \mu_{6} \theta_{0}^{2} - 6 \, {\left(\alpha_{1} \theta_{0} + 2 \, \alpha_{1}\right)} {\left| A_{1} \right|}^{2} + 2 \, {\left(3 \, \alpha_{6} + \mu_{6}\right)} \theta_{0} - 6 \, \theta_{0} \zeta_{11} + 6 \, \alpha_{6}}{\theta_{0}^{2} + \theta_{0}} & \overline{A_{1}} & 2 & \theta_{0} + 1 \\
    	    	-\frac{\alpha_{9} \theta_{0}^{3} - 3 \, \theta_{0}^{2} \overline{\mu_{11}} - {\left(7 \, \alpha_{9} + 3 \, \overline{\mu_{11}}\right)} \theta_{0} + {\left(\theta_{0} + 1\right)} \zeta_{4} + 6 \, \alpha_{1} \overline{\zeta_{0}} - 6 \, \alpha_{9}}{\theta_{0}^{2} + \theta_{0}} & A_{1} & 2 & \theta_{0} + 1 \\
    	    	\lambda_6 & A_{0} & 2 & \theta_{0} + 1 
    	    	    	    	    	    	\end{dmatrix}
    	    	\end{align*}
    	    	\begin{align*}
    	    	\begin{dmatrix}
    	    	\frac{2 \, \mu_{5} \theta_{0}^{4} + 6 \, \mu_{5} \theta_{0}^{3} + 4 \, \mu_{5} \theta_{0}^{2} - {\left(2 \, \theta_{0}^{2} - \theta_{0} - 6\right)} \zeta_{2}}{2 \, {\left(\theta_{0}^{3} + \theta_{0}^{2}\right)}} & \overline{A_{2}} & 2 & \theta_{0} + 1 \\
    	    	\lambda_7 & \overline{A_{0}} & 2 & \theta_{0} + 1 \\
    	    	\frac{\mu_{1} \theta_{0}^{4} + 3 \, \mu_{1} \theta_{0}^{3} + 2 \, {\left(\mu_{1} + \overline{\alpha_{5}} + 3 \, \overline{\mu_{3}}\right)} \theta_{0}^{2} + 2 \, \theta_{0} {\left(\overline{\alpha_{5}} + 3 \, \overline{\mu_{3}}\right)} + 12 \, \overline{\alpha_{5}}}{8 \, {\left(\theta_{0}^{3} + \theta_{0}^{2}\right)}} & C_{1} & 2 & \theta_{0} + 1 \\
    	    	-\frac{{\left(2 \, \theta_{0} - 1\right)} {\left| A_{1} \right|}^{2}}{2 \, {\left(\theta_{0}^{2} + 2 \, \theta_{0} + 1\right)}} & B_{1} & 2 & \theta_{0} + 1 \\
    	    	\frac{2}{\theta_{0}^{2}} & \overline{B_{7}} & 3 & \theta_{0} \\
    	    	\frac{\theta_{0} + 1}{2 \, \theta_{0}^{2}} & B_{12} & 3 & \theta_{0} \\
    	    	\frac{\theta_{0} - 5}{4 \, \theta_{0}} & B_{6} & 3 & \theta_{0} \\
    	    	\frac{\mu_{7} \theta_{0}^{2} + \mu_{7} \theta_{0} - \theta_{0} \zeta_{7} - 4 \, {\left(\theta_{0}^{2} + 4 \, \theta_{0}\right)} \zeta_{8} + 4 \, \alpha_{4}}{\theta_{0}} & \overline{A_{1}} & 3 & \theta_{0} \\
    	    	-\frac{2 \, \alpha_{8} \theta_{0}^{2} - 2 \, {\left(3 \, \alpha_{8} + 4 \, \overline{\mu_{10}}\right)} \theta_{0} - 8 \, \alpha_{8} + \zeta_{3}}{2 \, \theta_{0}} & A_{1} & 3 & \theta_{0} \\
    	    	-\frac{\alpha_{1}}{2 \, {\left(\theta_{0} + 2\right)}} & B_{1} & 3 & \theta_{0} \\
    	    	\lambda_8 & A_{0} & 3 & \theta_{0} \\
    	    	\frac{\mu_{3} \theta_{0}^{4} + 3 \, \mu_{3} \theta_{0}^{3} + 2 \, {\left(\alpha_{5} + \mu_{3} + 4 \, \overline{\mu_{1}}\right)} \theta_{0}^{2} + 4 \, {\left(\alpha_{5} + 4 \, \overline{\mu_{1}}\right)} \theta_{0} + 16 \, \alpha_{5}}{8 \, {\left(\theta_{0}^{3} + 2 \, \theta_{0}^{2}\right)}} & C_{1} & 3 & \theta_{0} \\
    	    	\lambda_9 & \overline{A_{0}} & 3 & \theta_{0} \\
    	    	-2 \, \alpha_{2} \theta_{0} - 2 \, \alpha_{2} & A_{2} & 3 & \theta_{0} \\
    	    	\frac{\mu_{7} \theta_{0}^{2} - \theta_{0} \zeta_{7} + \zeta_{5}}{\theta_{0}} & \overline{A_{0}} & 3 & \theta_{0} - 1 \\
    	    	\frac{1}{2 \, \theta_{0}} & B_{13} & 4 & \theta_{0} - 1 \\
    	    	\frac{\theta_{0} - 5}{4 \, \theta_{0}} & C_{4} & 4 & \theta_{0} - 1 \\
    	    	\frac{\mu_{19} \theta_{0}^{2} + 2 \, \mu_{19} \theta_{0} - {\left(\theta_{0} + 2\right)} \zeta_{17} + \alpha_{1} \zeta_{2}}{\theta_{0} + 2} & \overline{A_{0}} & 4 & \theta_{0} - 1 \\
    	    	\frac{2 \, \mu_{5} \theta_{0}^{3} + 2 \, \mu_{5} \theta_{0}^{2} - {\left(2 \, \theta_{0} - 1\right)} \zeta_{2}}{16 \, {\left(\theta_{0}^{2} + \theta_{0}\right)}} & C_{1} & 4 & \theta_{0} - 1 \\
    	    	-\alpha_{7} \theta_{0} + 4 \, \alpha_{7} & A_{1} & 4 & \theta_{0} - 1 \\
    	    	-\alpha_{13} \theta_{0} + 5 \, \alpha_{13} - \zeta_{21} + 5 \, \overline{\mu_{25}} & A_{0} & 4 & \theta_{0} - 1 \\
    	    	-\frac{{\left(\alpha_{9} - 2 \, \mu_{8}\right)} \theta_{0}^{3} + 2 \, {\left(\alpha_{9} - 3 \, \mu_{8}\right)} \theta_{0}^{2} - 4 \, {\left(\alpha_{9} + \mu_{8}\right)} \theta_{0} + {\left(\theta_{0} - 2\right)} \zeta_{4} - 8 \, \alpha_{9}}{\theta_{0}^{2} + 2 \, \theta_{0}} & \overline{A_{0}} & -\theta_{0} + 1 & 2 \, \theta_{0} + 1 \\
    	    	-\frac{{\left(\alpha_{9} - 2 \, \mu_{8}\right)} \theta_{0}^{3} + {\left(3 \, \alpha_{9} - 7 \, \mu_{8}\right)} \theta_{0}^{2} - 2 \, {\left(2 \, \alpha_{9} + 3 \, \mu_{8}\right)} \theta_{0} + {\left(\theta_{0} - 2\right)} \zeta_{4} - 12 \, \alpha_{9}}{\theta_{0}^{2} + 2 \, \theta_{0}} & \overline{A_{1}} & -\theta_{0} + 1 & 2 \, \theta_{0} + 2 \\
    	    	\lambda_{10} & \overline{A_{0}} & -\theta_{0} + 1 & 2 \, \theta_{0} + 2 \\
    	    	-\frac{\theta_{0}^{4} \overline{\mu_{30}} + \theta_{0}^{3} \overline{\mu_{30}} - 4 \, \theta_{0}^{2} \overline{\mu_{30}} - 4 \, \theta_{0} \overline{\mu_{30}} - {\left(5 \, \alpha_{2} \theta_{0}^{2} - 4 \, \alpha_{2} \theta_{0} - 12 \, \alpha_{2}\right)} \overline{\zeta_{0}}}{\theta_{0}^{3} + 3 \, \theta_{0}^{2} + 2 \, \theta_{0}} & A_{0} & -\theta_{0} + 1 & 2 \, \theta_{0} + 2 
    	    	    	    	    	    	\end{dmatrix}
    	    	\end{align*}
    	    	\begin{align*}
    	    	\begin{dmatrix}
    	    	-\frac{\theta_{0}^{4} \overline{\mu_{12}} + \theta_{0}^{3} \overline{\mu_{12}} - 4 \, \theta_{0}^{2} \overline{\mu_{12}} - 4 \, \theta_{0} \overline{\mu_{12}} + 4 \, {\left(2 \, \theta_{0}^{2} - \theta_{0} - 6\right)} \overline{\zeta_{1}}}{4 \, {\left(\theta_{0}^{4} + 3 \, \theta_{0}^{3} + 2 \, \theta_{0}^{2}\right)}} & C_{1} & -\theta_{0} + 1 & 2 \, \theta_{0} + 2 \\
    	    	\frac{{\left(\theta_{0} - 2\right)} \overline{\zeta_{0}}}{4 \, {\left(\theta_{0}^{4} + 4 \, \theta_{0}^{3} + 5 \, \theta_{0}^{2} + 2 \, \theta_{0}\right)}} & B_{1} & -\theta_{0} + 1 & 2 \, \theta_{0} + 2 \\
    	    	-\frac{2 \, {\left(\alpha_{8} - 2 \, \mu_{9}\right)} \theta_{0}^{3} - 8 \, \mu_{9} \theta_{0}^{2} - 2 \, {\left(7 \, \alpha_{8} + 2 \, \mu_{9}\right)} \theta_{0} + {\left(\theta_{0} - 3\right)} \zeta_{3} - 12 \, \alpha_{8}}{2 \, {\left(\theta_{0}^{2} + \theta_{0}\right)}} & \overline{A_{1}} & -\theta_{0} + 2 & 2 \, \theta_{0} + 1 \\
    	    	\lambda_{11} & \overline{A_{0}} & -\theta_{0} + 2 & 2 \, \theta_{0} + 1 \\
    	    	-\frac{2 \, {\left(\alpha_{8} - 2 \, \mu_{9}\right)} \theta_{0}^{3} - 2 \, {\left(\alpha_{8} + 3 \, \mu_{9}\right)} \theta_{0}^{2} - 2 \, {\left(5 \, \alpha_{8} + \mu_{9}\right)} \theta_{0} + {\left(\theta_{0} - 3\right)} \zeta_{3} - 6 \, \alpha_{8}}{2 \, {\left(\theta_{0}^{2} + \theta_{0}\right)}} & \overline{A_{0}} & -\theta_{0} + 2 & 2 \, \theta_{0} \\
    	    	\frac{\theta_{0} - 4}{4 \, \theta_{0}} & E_{1} & -\theta_{0} + 3 & 2 \, \theta_{0} - 1 \\
    	    	\frac{\alpha_{7} \theta_{0}^{2} - 8 \, \alpha_{7} \theta_{0} + 16 \, \alpha_{7}}{\theta_{0}} & \overline{A_{0}} & -\theta_{0} + 3 & 2 \, \theta_{0} - 1 \\
    	    	\frac{\theta_{0} - 5}{4 \, \theta_{0}} & E_{2} & -\theta_{0} + 3 & 2 \, \theta_{0} \\
    	    	\frac{\alpha_{2} \theta_{0} - 5 \, \alpha_{2}}{4 \, \theta_{0}} & C_{1} & -\theta_{0} + 3 & 2 \, \theta_{0} \\
    	    	\frac{{\left(\alpha_{1} \alpha_{2} - \alpha_{12} + 2 \, \mu_{22}\right)} \theta_{0}^{2} - 4 \, \alpha_{1} \alpha_{2} - {\left(3 \, \alpha_{1} \alpha_{2} - 3 \, \alpha_{12} - \mu_{22}\right)} \theta_{0} - {\left(2 \, \theta_{0} + 1\right)} \zeta_{20} + {\left(\theta_{0} - 4\right)} \overline{\zeta_{22}} + 4 \, \alpha_{12}}{\theta_{0}} & \overline{A_{0}} & -\theta_{0} + 3 & 2 \, \theta_{0} \\
    	    	\frac{2 \, \alpha_{7} \theta_{0}^{3} - 17 \, \alpha_{7} \theta_{0}^{2} + 31 \, \alpha_{7} \theta_{0} + 20 \, \alpha_{7}}{2 \, {\left(\theta_{0}^{2} + \theta_{0}\right)}} & \overline{A_{1}} & -\theta_{0} + 3 & 2 \, \theta_{0} \\
    	    	-\theta_{0} \overline{\mu_{31}} + 4 \, \overline{\mu_{31}} & A_{0} & -\theta_{0} + 3 & 2 \, \theta_{0} \\
    	    	\frac{\theta_{0} - 5}{4 \, \theta_{0}} & E_{3} & -\theta_{0} + 4 & 2 \, \theta_{0} - 1 \\
    	    	-{\left(\alpha_{13} - 2 \, \mu_{23}\right)} \theta_{0} + 5 \, \alpha_{13} - 2 \, \zeta_{21} & \overline{A_{0}} & -\theta_{0} + 4 & 2 \, \theta_{0} - 1 \\
    	    	-\frac{\theta_{0} - 1}{2 \, \theta_{0}^{2}} & \overline{B_{1}} & \theta_{0} & 1 \\
    	    	\frac{\theta_{0}^{2} \overline{\mu_{3}} + \theta_{0} \overline{\mu_{3}} + 2 \, \overline{\alpha_{5}}}{\theta_{0}} & A_{0} & \theta_{0} & 1 \\
    	    	-\frac{2 \, {\left(\theta_{0}^{2} \overline{\alpha_{2}} - \overline{\alpha_{2}}\right)}}{\theta_{0}} & \overline{A_{0}} & \theta_{0} & 1 \\
    	    	-\frac{2 \, \theta_{0} - 3}{4 \, \theta_{0}^{2}} & \overline{B_{3}} & \theta_{0} & 2 \\
    	    	-\frac{2 \, \theta_{0}^{3} \overline{\alpha_{8}} - 2 \, {\left(3 \, \mu_{10} + 2 \, \overline{\alpha_{8}}\right)} \theta_{0}^{2} - 6 \, \theta_{0} \overline{\alpha_{8}} + {\left(\theta_{0} - 3\right)} \overline{\zeta_{3}}}{2 \, \theta_{0}^{2}} & \overline{A_{0}} & \theta_{0} & 2 \\
    	    	-\frac{2 \, \theta_{0}^{2} \overline{\alpha_{2}} - \theta_{0} \overline{\alpha_{2}} - 3 \, \overline{\alpha_{2}}}{\theta_{0}} & \overline{A_{1}} & \theta_{0} & 2 \\
    	    	-\frac{6 \, {\left| A_{1} \right|}^{2} \overline{\alpha_{1}} - \theta_{0}^{2} \overline{\mu_{6}} - \theta_{0} \overline{\mu_{6}} - 3 \, \overline{\alpha_{6}}}{\theta_{0}} & A_{0} & \theta_{0} & 2 \\
    	    	\frac{2 \, \theta_{0}^{4} \overline{\mu_{5}} - 4 \, \theta_{0}^{3} \overline{\mu_{5}} - 6 \, \theta_{0}^{2} \overline{\mu_{5}} + {\left(5 \, \theta_{0} - 12\right)} \overline{\zeta_{2}}}{2 \, {\left(\theta_{0}^{3} - 3 \, \theta_{0}^{2}\right)}} & A_{1} & \theta_{0} & 2 \\
    	    	\frac{2}{\theta_{0}^{2}} & B_{7} & \theta_{0} & 3 \\
    	    	\frac{\theta_{0} + 1}{2 \, \theta_{0}^{2}} & \overline{B_{12}} & \theta_{0} & 3 
    	    	    	    	    	    	\end{dmatrix}
    	    	\end{align*}
    	    	\begin{align*}
    	    	\begin{dmatrix}
    	    	-\frac{\overline{\alpha_{1}}}{2 \, {\left(\theta_{0} + 2\right)}} & \overline{B_{1}} & \theta_{0} & 3 \\
    	    	\frac{\theta_{0} - 5}{4 \, \theta_{0}} & \overline{B_{6}} & \theta_{0} & 3 \\
    	    	\frac{\theta_{0}^{2} \overline{\mu_{7}} + \theta_{0} \overline{\mu_{7}} - \theta_{0} \overline{\zeta_{7}} - 4 \, {\left(\theta_{0}^{2} + 4 \, \theta_{0}\right)} \overline{\zeta_{8}} + 4 \, \overline{\alpha_{4}}}{\theta_{0}} & A_{1} & \theta_{0} & 3 \\
    	    	-\frac{2 \, \theta_{0}^{2} \overline{\alpha_{8}} - 2 \, {\left(4 \, \mu_{10} + 3 \, \overline{\alpha_{8}}\right)} \theta_{0} - 8 \, \overline{\alpha_{8}} + \overline{\zeta_{3}}}{2 \, \theta_{0}} & \overline{A_{1}} & \theta_{0} & 3 \\
    	    	\lambda_{12} & \overline{A_{0}} & \theta_{0} & 3 \\
    	    	\frac{\theta_{0}^{4} \overline{\mu_{3}} + 3 \, \theta_{0}^{3} \overline{\mu_{3}} + 2 \, {\left(4 \, \mu_{1} + \overline{\alpha_{5}} + \overline{\mu_{3}}\right)} \theta_{0}^{2} + 4 \, {\left(4 \, \mu_{1} + \overline{\alpha_{5}}\right)} \theta_{0} + 16 \, \overline{\alpha_{5}}}{8 \, {\left(\theta_{0}^{3} + 2 \, \theta_{0}^{2}\right)}} & \overline{C_{1}} & \theta_{0} & 3 \\
    	    	\lambda_{13} & A_{0} & \theta_{0} & 3 \\
    	    	-2 \, \theta_{0} \overline{\alpha_{2}} - 2 \, \overline{\alpha_{2}} & \overline{A_{2}} & \theta_{0} & 3 \\
    	    	\frac{2 \, \theta_{0}^{3} \overline{\mu_{5}} + 2 \, \theta_{0}^{2} \overline{\mu_{5}} + 3 \, \overline{\zeta_{2}}}{2 \, {\left(\theta_{0}^{2} + \theta_{0}\right)}} & A_{0} & \theta_{0} - 1 & 2 \\
    	    	\frac{\theta_{0}^{2} \overline{\mu_{7}} - \theta_{0} \overline{\zeta_{7}} + \overline{\zeta_{5}}}{\theta_{0}} & A_{0} & \theta_{0} - 1 & 3 \\
    	    	\frac{1}{2 \, \theta_{0}} & \overline{B_{13}} & \theta_{0} - 1 & 4 \\
    	    	\frac{\theta_{0} - 5}{4 \, \theta_{0}} & \overline{C_{4}} & \theta_{0} - 1 & 4 \\
    	    	\frac{\theta_{0}^{2} \overline{\mu_{19}} + 2 \, \theta_{0} \overline{\mu_{19}} - {\left(\theta_{0} + 2\right)} \overline{\zeta_{17}} + \overline{\alpha_{1}} \overline{\zeta_{2}}}{\theta_{0} + 2} & A_{0} & \theta_{0} - 1 & 4 \\
    	    	\frac{2 \, \theta_{0}^{3} \overline{\mu_{5}} + 2 \, \theta_{0}^{2} \overline{\mu_{5}} - {\left(2 \, \theta_{0} - 1\right)} \overline{\zeta_{2}}}{16 \, {\left(\theta_{0}^{2} + \theta_{0}\right)}} & \overline{C_{1}} & \theta_{0} - 1 & 4 \\
    	    	-\theta_{0} \overline{\alpha_{7}} + 4 \, \overline{\alpha_{7}} & \overline{A_{1}} & \theta_{0} - 1 & 4 \\
    	    	-\theta_{0} \overline{\alpha_{13}} + 5 \, \mu_{25} + 5 \, \overline{\alpha_{13}} - \overline{\zeta_{21}} & \overline{A_{0}} & \theta_{0} - 1 & 4 \\
    	    	\frac{\theta_{0}^{2} \overline{\mu_{1}} + 2 \, \theta_{0} \overline{\mu_{1}} + 2 \, \alpha_{5}}{\theta_{0}} & A_{0} & \theta_{0} + 1 & 0 \\
    	    	-\frac{\theta_{0} - 1}{\theta_{0}^{2}} & \overline{B_{2}} & \theta_{0} + 1 & 1 \\
    	    	\frac{\theta_{0}^{3} \overline{\mu_{3}} + 3 \, \theta_{0}^{2} \overline{\mu_{3}} + 2 \, \theta_{0} {\left(\overline{\alpha_{5}} + \overline{\mu_{3}}\right)} + 4 \, \overline{\alpha_{5}}}{\theta_{0}^{2} + \theta_{0}} & A_{1} & \theta_{0} + 1 & 1 \\
    	    	-\frac{\theta_{0}^{4} \overline{\alpha_{9}} - {\left(2 \, \mu_{11} - \overline{\alpha_{9}}\right)} \theta_{0}^{3} - 2 \, {\left(\mu_{11} + 2 \, \overline{\alpha_{9}}\right)} \theta_{0}^{2} + 4 \, \theta_{0} \zeta_{0} \overline{\alpha_{1}} - 4 \, \theta_{0} \overline{\alpha_{9}} + {\left(\theta_{0}^{2} - \theta_{0} - 2\right)} \overline{\zeta_{4}}}{\theta_{0}^{3} + \theta_{0}^{2}} & \overline{A_{0}} & \theta_{0} + 1 & 1 \\
    	    	\frac{4 \, \theta_{0}^{5} \overline{\mu_{4}} - 32 \, {\left(\theta_{0}^{3} + 2 \, \theta_{0}^{2}\right)} {\left| A_{1} \right|}^{4} + 12 \, \theta_{0}^{4} \overline{\mu_{4}} - 8 \, {\left(2 \, \alpha_{1} \overline{\alpha_{1}} - 2 \, \beta - \overline{\mu_{4}}\right)} \theta_{0}^{3} - 16 \, {\left(\alpha_{1} \overline{\alpha_{1}} - \beta\right)} \theta_{0}^{2} - {\left(\theta_{0}^{3} - 3 \, \theta_{0} - 2\right)} {\left| C_{1} \right|}^{2}}{4 \, {\left(\theta_{0}^{4} + \theta_{0}^{3}\right)}} & A_{0} & \theta_{0} + 1 & 1 
    	    	    	    	    	    	\end{dmatrix}
    	    	\end{align*}
    	    	\begin{align*}
    	    	\begin{dmatrix}
    	    	\frac{3}{2 \, \theta_{0}^{2}} & B_{8} & \theta_{0} + 1 & 2 \\
    	    	\frac{\theta_{0} + 2}{2 \, \theta_{0}^{2}} & \overline{B_{11}} & \theta_{0} + 1 & 2 \\
    	    	\frac{\theta_{0} - 5}{4 \, \theta_{0}} & \overline{B_{5}} & \theta_{0} + 1 & 2 \\
    	    	\frac{\theta_{0}^{3} \overline{\mu_{6}} - 6 \, {\left(\theta_{0} \overline{\alpha_{1}} + 2 \, \overline{\alpha_{1}}\right)} {\left| A_{1} \right|}^{2} + 3 \, \theta_{0}^{2} \overline{\mu_{6}} + 2 \, \theta_{0} {\left(3 \, \overline{\alpha_{6}} + \overline{\mu_{6}}\right)} - 6 \, \theta_{0} \overline{\zeta_{11}} + 6 \, \overline{\alpha_{6}}}{\theta_{0}^{2} + \theta_{0}} & A_{1} & \theta_{0} + 1 & 2 \\
    	    	-\frac{\theta_{0}^{3} \overline{\alpha_{9}} - 3 \, \mu_{11} \theta_{0}^{2} - {\left(3 \, \mu_{11} + 7 \, \overline{\alpha_{9}}\right)} \theta_{0} + 6 \, \zeta_{0} \overline{\alpha_{1}} + {\left(\theta_{0} + 1\right)} \overline{\zeta_{4}} - 6 \, \overline{\alpha_{9}}}{\theta_{0}^{2} + \theta_{0}} & \overline{A_{1}} & \theta_{0} + 1 & 2 \\
    	    	\lambda_{14} & \overline{A_{0}} & \theta_{0} + 1 & 2 \\
    	    	\frac{2 \, \theta_{0}^{4} \overline{\mu_{5}} + 6 \, \theta_{0}^{3} \overline{\mu_{5}} + 4 \, \theta_{0}^{2} \overline{\mu_{5}} - {\left(2 \, \theta_{0}^{2} - \theta_{0} - 6\right)} \overline{\zeta_{2}}}{2 \, {\left(\theta_{0}^{3} + \theta_{0}^{2}\right)}} & A_{2} & \theta_{0} + 1 & 2 \\
    	    	\lambda_{15} & A_{0} & \theta_{0} + 1 & 2 \\
    	    	\frac{\theta_{0}^{4} \overline{\mu_{1}} + 3 \, \theta_{0}^{3} \overline{\mu_{1}} + 2 \, {\left(\alpha_{5} + 3 \, \mu_{3} + \overline{\mu_{1}}\right)} \theta_{0}^{2} + 2 \, {\left(\alpha_{5} + 3 \, \mu_{3}\right)} \theta_{0} + 12 \, \alpha_{5}}{8 \, {\left(\theta_{0}^{3} + \theta_{0}^{2}\right)}} & \overline{C_{1}} & \theta_{0} + 1 & 2 \\
    	    	-\frac{{\left(2 \, \theta_{0} - 1\right)} {\left| A_{1} \right|}^{2}}{2 \, {\left(\theta_{0}^{2} + 2 \, \theta_{0} + 1\right)}} & \overline{B_{1}} & \theta_{0} + 1 & 2 \\
    	    	\frac{\theta_{0}^{3} \overline{\mu_{1}} + 5 \, \theta_{0}^{2} \overline{\mu_{1}} + 2 \, {\left(\alpha_{5} + 3 \, \overline{\mu_{1}}\right)} \theta_{0} + 6 \, \alpha_{5}}{\theta_{0}^{2} + 2 \, \theta_{0}} & A_{1} & \theta_{0} + 2 & 0 \\
    	    	\frac{\mu_{13} \theta_{0}^{3} - 4 \, {\left(\theta_{0} + 3\right)} \zeta_{0} {\left| A_{1} \right|}^{2} + 3 \, \mu_{13} \theta_{0}^{2} + 2 \, \mu_{13} \theta_{0}}{\theta_{0}^{3} + 3 \, \theta_{0}^{2} + 2 \, \theta_{0}} & \overline{A_{0}} & \theta_{0} + 2 & 0 \\
    	    	\frac{\mu_{12} \theta_{0}^{3} + 3 \, \mu_{12} \theta_{0}^{2} + 2 \, \mu_{12} \theta_{0} + 4 \, {\left(2 \, \theta_{0} + 3\right)} \zeta_{1}}{\theta_{0}^{3} + 3 \, \theta_{0}^{2} + 2 \, \theta_{0}} & \overline{A_{1}} & \theta_{0} + 2 & 0 \\
    	    	\frac{\theta_{0}^{4} \overline{\mu_{2}} + 6 \, \theta_{0}^{3} \overline{\mu_{2}} + {\left(3 \, \alpha_{6} + 11 \, \overline{\mu_{2}}\right)} \theta_{0}^{2} - 2 \, {\left(3 \, \alpha_{1} \theta_{0}^{2} + 13 \, \alpha_{1} \theta_{0} + 12 \, \alpha_{1}\right)} {\left| A_{1} \right|}^{2} + 3 \, {\left(3 \, \alpha_{6} + 2 \, \overline{\mu_{2}}\right)} \theta_{0} + 6 \, \alpha_{6}}{\theta_{0}^{3} + 3 \, \theta_{0}^{2} + 2 \, \theta_{0}} & A_{0} & \theta_{0} + 2 & 0 \\
    	    	\frac{1}{\theta_{0}^{2}} & B_{9} & \theta_{0} + 2 & 1 \\
    	    	\frac{\theta_{0} + 3}{2 \, \theta_{0}^{2}} & \overline{B_{10}} & \theta_{0} + 2 & 1 \\
    	    	\frac{\theta_{0} - 5}{4 \, \theta_{0}} & \overline{B_{4}} & \theta_{0} + 2 & 1 \\
    	    	\frac{\theta_{0}^{3} \overline{\mu_{3}} + 4 \, \theta_{0}^{2} \overline{\mu_{3}} + \theta_{0} {\left(2 \, \overline{\alpha_{5}} + 3 \, \overline{\mu_{3}}\right)} + 6 \, \overline{\alpha_{5}}}{\theta_{0}^{2} + \theta_{0}} & A_{2} & \theta_{0} + 2 & 1 \\
    	    	-\frac{\alpha_{1} \theta_{0} - \alpha_{1}}{2 \, {\left(\theta_{0}^{2} + 2 \, \theta_{0}\right)}} & \overline{B_{1}} & \theta_{0} + 2 & 1 \\
    	    	\frac{2 \, {\left(\mu_{13} \theta_{0}^{3} - 4 \, {\left(\theta_{0} + 3\right)} \zeta_{0} {\left| A_{1} \right|}^{2} + 3 \, \mu_{13} \theta_{0}^{2} + 2 \, \mu_{13} \theta_{0}\right)}}{\theta_{0}^{3} + 3 \, \theta_{0}^{2} + 2 \, \theta_{0}} & \overline{A_{1}} & \theta_{0} + 2 & 1 \\
    	    	\frac{2 \, {\left(\mu_{12} \theta_{0}^{3} + 3 \, \mu_{12} \theta_{0}^{2} + 2 \, \mu_{12} \theta_{0} + 4 \, {\left(2 \, \theta_{0} + 3\right)} \zeta_{1}\right)}}{\theta_{0}^{3} + 3 \, \theta_{0}^{2} + 2 \, \theta_{0}} & \overline{A_{2}} & \theta_{0} + 2 & 1 \\
    	    	\lambda_{16} & A_{0} & \theta_{0} + 2 & 1 \\
    	    	\lambda_{17} & A_{1} & \theta_{0} + 2 & 1 \\
    	    	\lambda_{18} & \overline{A_{0}} & \theta_{0} + 2 & 1 \\
    	    	\frac{2 \, \mu_{5} \theta_{0}^{3} + 2 \, \mu_{5} \theta_{0}^{2} - {\left(\theta_{0}^{2} + \theta_{0} - 3\right)} \zeta_{2}}{4 \, {\left(\theta_{0}^{4} + \theta_{0}^{3}\right)}} & \overline{C_{1}} & \theta_{0} + 2 & 1 \\
    	    	-\frac{1}{4} \, \theta_{0} \overline{\alpha_{2}} - \frac{1}{4} \, \overline{\alpha_{2}} & C_{1} & \theta_{0} + 2 & 1 
    	    	    	    	    	    	\end{dmatrix}
    	    	\end{align*}
    	    	\begin{align*}
    	    	\begin{dmatrix}
    	    	\frac{3 \, \zeta_{0}}{4 \, {\left(\theta_{0}^{3} + 4 \, \theta_{0}^{2} + 5 \, \theta_{0} + 2\right)}} & B_{1} & \theta_{0} + 3 & 0 \\
    	    	\frac{\theta_{0}^{3} \overline{\mu_{1}} + 6 \, \theta_{0}^{2} \overline{\mu_{1}} + 2 \, {\left(\alpha_{5} + 4 \, \overline{\mu_{1}}\right)} \theta_{0} + 8 \, \alpha_{5}}{\theta_{0}^{2} + 2 \, \theta_{0}} & A_{2} & \theta_{0} + 3 & 0 \\
    	    	\frac{\mu_{14} \theta_{0}^{3} + 3 \, \mu_{14} \theta_{0}^{2} + 2 \, \mu_{14} \theta_{0} - 2 \, {\left(\alpha_{1} \theta_{0} + 4 \, \alpha_{1}\right)} \zeta_{0} - 2 \, {\left(\theta_{0}^{5} + 9 \, \theta_{0}^{4} + 28 \, \theta_{0}^{3} + 36 \, \theta_{0}^{2} + 16 \, \theta_{0}\right)} \zeta_{9}}{\theta_{0}^{3} + 3 \, \theta_{0}^{2} + 2 \, \theta_{0}} & \overline{A_{1}} & \theta_{0} + 3 & 0 \\
    	    	\frac{\mu_{28} \theta_{0}^{4} + 6 \, \mu_{28} \theta_{0}^{3} - 8 \, {\left(2 \, \theta_{0}^{2} + 11 \, \theta_{0} + 12\right)} \zeta_{1} {\left| A_{1} \right|}^{2} + 11 \, \mu_{28} \theta_{0}^{2} + 6 \, \mu_{28} \theta_{0} - 2 \, {\left(\alpha_{5} \theta_{0}^{2} + 7 \, \alpha_{5} \theta_{0} + 12 \, \alpha_{5}\right)} \zeta_{0}}{\theta_{0}^{4} + 6 \, \theta_{0}^{3} + 11 \, \theta_{0}^{2} + 6 \, \theta_{0}} & \overline{A_{0}} & \theta_{0} + 3 & 0 \\
    	    	\lambda_{19} & A_{1} & \theta_{0} + 3 & 0 \\
    	    	\lambda_{20} & A_{0} & \theta_{0} + 3 & 0 \\
    	    	\frac{\theta_{0} - 4}{4 \, \theta_{0}} & \overline{E_{1}} & 2 \, \theta_{0} - 1 & -\theta_{0} + 3 \\
    	    	\frac{\theta_{0}^{2} \overline{\alpha_{7}} - 8 \, \theta_{0} \overline{\alpha_{7}} + 16 \, \overline{\alpha_{7}}}{\theta_{0}} & A_{0} & 2 \, \theta_{0} - 1 & -\theta_{0} + 3 \\
    	    	\frac{\theta_{0} - 5}{4 \, \theta_{0}} & \overline{E_{3}} & 2 \, \theta_{0} - 1 & -\theta_{0} + 4 \\
    	    	-\theta_{0} {\left(\overline{\alpha_{13}} - 2 \, \overline{\mu_{23}}\right)} + 5 \, \overline{\alpha_{13}} - 2 \, \overline{\zeta_{21}} & A_{0} & 2 \, \theta_{0} - 1 & -\theta_{0} + 4 \\
    	    	-\frac{\theta_{0}^{3} {\left(\overline{\alpha_{9}} - 2 \, \overline{\mu_{8}}\right)} + 2 \, \theta_{0}^{2} {\left(\overline{\alpha_{9}} - 3 \, \overline{\mu_{8}}\right)} - 4 \, \theta_{0} {\left(\overline{\alpha_{9}} + \overline{\mu_{8}}\right)} + {\left(\theta_{0} - 2\right)} \overline{\zeta_{4}} - 8 \, \overline{\alpha_{9}}}{\theta_{0}^{2} + 2 \, \theta_{0}} & A_{0} & 2 \, \theta_{0} + 1 & -\theta_{0} + 1 \\
    	    	-\frac{2 \, \theta_{0}^{3} {\left(\overline{\alpha_{8}} - 2 \, \overline{\mu_{9}}\right)} - 8 \, \theta_{0}^{2} \overline{\mu_{9}} - 2 \, \theta_{0} {\left(7 \, \overline{\alpha_{8}} + 2 \, \overline{\mu_{9}}\right)} + {\left(\theta_{0} - 3\right)} \overline{\zeta_{3}} - 12 \, \overline{\alpha_{8}}}{2 \, {\left(\theta_{0}^{2} + \theta_{0}\right)}} & A_{1} & 2 \, \theta_{0} + 1 & -\theta_{0} + 2 \\
    	    	\lambda_{21} & A_{0} & 2 \, \theta_{0} + 1 & -\theta_{0} + 2 \\
    	    	-\frac{\theta_{0}^{3} {\left(\overline{\alpha_{9}} - 2 \, \overline{\mu_{8}}\right)} + \theta_{0}^{2} {\left(3 \, \overline{\alpha_{9}} - 7 \, \overline{\mu_{8}}\right)} - 2 \, \theta_{0} {\left(2 \, \overline{\alpha_{9}} + 3 \, \overline{\mu_{8}}\right)} + {\left(\theta_{0} - 2\right)} \overline{\zeta_{4}} - 12 \, \overline{\alpha_{9}}}{\theta_{0}^{2} + 2 \, \theta_{0}} & A_{1} & 2 \, \theta_{0} + 2 & -\theta_{0} + 1 \\
    	    	\lambda_{22} & A_{0} & 2 \, \theta_{0} + 2 & -\theta_{0} + 1 \\
    	    	-\frac{\mu_{30} \theta_{0}^{4} + \mu_{30} \theta_{0}^{3} - 4 \, \mu_{30} \theta_{0}^{2} - 4 \, \mu_{30} \theta_{0} - {\left(5 \, \theta_{0}^{2} \overline{\alpha_{2}} - 4 \, \theta_{0} \overline{\alpha_{2}} - 12 \, \overline{\alpha_{2}}\right)} \zeta_{0}}{\theta_{0}^{3} + 3 \, \theta_{0}^{2} + 2 \, \theta_{0}} & \overline{A_{0}} & 2 \, \theta_{0} + 2 & -\theta_{0} + 1 \\
    	    	-\frac{\mu_{12} \theta_{0}^{4} + \mu_{12} \theta_{0}^{3} - 4 \, \mu_{12} \theta_{0}^{2} - 4 \, \mu_{12} \theta_{0} + 4 \, {\left(2 \, \theta_{0}^{2} - \theta_{0} - 6\right)} \zeta_{1}}{4 \, {\left(\theta_{0}^{4} + 3 \, \theta_{0}^{3} + 2 \, \theta_{0}^{2}\right)}} & \overline{C_{1}} & 2 \, \theta_{0} + 2 & -\theta_{0} + 1 \\
    	    	\frac{{\left(\theta_{0} - 2\right)} \zeta_{0}}{4 \, {\left(\theta_{0}^{4} + 4 \, \theta_{0}^{3} + 5 \, \theta_{0}^{2} + 2 \, \theta_{0}\right)}} & \overline{B_{1}} & 2 \, \theta_{0} + 2 & -\theta_{0} + 1 \\
    	    	-\frac{2 \, \theta_{0}^{3} {\left(\overline{\alpha_{8}} - 2 \, \overline{\mu_{9}}\right)} - 2 \, \theta_{0}^{2} {\left(\overline{\alpha_{8}} + 3 \, \overline{\mu_{9}}\right)} - 2 \, \theta_{0} {\left(5 \, \overline{\alpha_{8}} + \overline{\mu_{9}}\right)} + {\left(\theta_{0} - 3\right)} \overline{\zeta_{3}} - 6 \, \overline{\alpha_{8}}}{2 \, {\left(\theta_{0}^{2} + \theta_{0}\right)}} & A_{0} & 2 \, \theta_{0} & -\theta_{0} + 2 \\
    	    	\frac{\theta_{0} - 5}{4 \, \theta_{0}} & \overline{E_{2}} & 2 \, \theta_{0} & -\theta_{0} + 3 \\
    	    	\frac{\theta_{0} \overline{\alpha_{2}} - 5 \, \overline{\alpha_{2}}}{4 \, \theta_{0}} & \overline{C_{1}} & 2 \, \theta_{0} & -\theta_{0} + 3 \\
    	    	\frac{{\left(\overline{\alpha_{1}} \overline{\alpha_{2}} - \overline{\alpha_{12}} + 2 \, \overline{\mu_{22}}\right)} \theta_{0}^{2} - {\left(3 \, \overline{\alpha_{1}} \overline{\alpha_{2}} - 3 \, \overline{\alpha_{12}} - \overline{\mu_{22}}\right)} \theta_{0} + {\left(\theta_{0} - 4\right)} \zeta_{22} - 4 \, \overline{\alpha_{1}} \overline{\alpha_{2}} - {\left(2 \, \theta_{0} + 1\right)} \overline{\zeta_{20}} + 4 \, \overline{\alpha_{12}}}{\theta_{0}} & A_{0} & 2 \, \theta_{0} & -\theta_{0} + 3 \\
    	    	\frac{2 \, \theta_{0}^{3} \overline{\alpha_{7}} - 17 \, \theta_{0}^{2} \overline{\alpha_{7}} + 31 \, \theta_{0} \overline{\alpha_{7}} + 20 \, \overline{\alpha_{7}}}{2 \, {\left(\theta_{0}^{2} + \theta_{0}\right)}} & A_{1} & 2 \, \theta_{0} & -\theta_{0} + 3 \\
    	    	-\mu_{31} \theta_{0} + 4 \, \mu_{31} & \overline{A_{0}} & 2 \, \theta_{0} & -\theta_{0} + 3
    	    	\end{dmatrix}
    	    \end{align*}
	where the $\lambda$ coefficients are given as follows.
	\begin{align*}
		&\lambda_{1}=\frac{1}{\theta_{0}^{4} + 6 \, \theta_{0}^{3} + 11 \, \theta_{0}^{2} + 6 \, \theta_{0}}\bigg\{\mu_{2} \theta_{0}^{5} + 10 \, \mu_{2} \theta_{0}^{4} + {\left(35 \, \mu_{2} + 4 \, \overline{\alpha_{6}}\right)} \theta_{0}^{3} + 2 \, {\left(25 \, \mu_{2} + 12 \, \overline{\alpha_{6}}\right)} \theta_{0}^{2}\\
		& - 2 \, {\left(3 \, \theta_{0}^{3} \overline{\alpha_{1}} + 25 \, \theta_{0}^{2} \overline{\alpha_{1}} + 64 \, \theta_{0} \overline{\alpha_{1}} + 48 \, \overline{\alpha_{1}}\right)} {\left| A_{1} \right|}^{2} + 4 \, {\left(6 \, \mu_{2} + 11 \, \overline{\alpha_{6}}\right)} \theta_{0} - 2 \, {\left(\theta_{0}^{3} + 3 \, \theta_{0}^{2} + 2 \, \theta_{0}\right)} \overline{\zeta_{11}} + 24 \, \overline{\alpha_{6}}\bigg\}\\
		\lambda_{2}&=\frac{1}{\theta_{0}^{4} + 6 \, \theta_{0}^{3} + 11 \, \theta_{0}^{2} + 6 \, \theta_{0}}\bigg\{\mu_{15} \theta_{0}^{5} + 10 \, \mu_{15} \theta_{0}^{4} - {\left(4 \, \overline{\alpha_{1}} \overline{\alpha_{5}} - 4 \, \alpha_{15} - 35 \, \mu_{15}\right)} \theta_{0}^{3} - 2 \, {\left(16 \, \overline{\alpha_{1}} \overline{\alpha_{5}} - 12 \, \alpha_{15} - 25 \, \mu_{15}\right)} \theta_{0}^{2}\\
		& - {\left(8 \, \theta_{0}^{3} \overline{\alpha_{3}} + 60 \, \theta_{0}^{2} \overline{\alpha_{3}} + 136 \, \theta_{0} \overline{\alpha_{3}} + {\left(\theta_{0} + 2\right)} \overline{\zeta_{2}} + 96 \, \overline{\alpha_{3}}\right)} {\left| A_{1} \right|}^{2} - 4 \, {\left(19 \, \overline{\alpha_{1}} \overline{\alpha_{5}} - 11 \, \alpha_{15} - 6 \, \mu_{15}\right)} \theta_{0}\\
		& - {\left(\theta_{0}^{4} + 10 \, \theta_{0}^{3} + 35 \, \theta_{0}^{2} + 50 \, \theta_{0} + 24\right)} \zeta_{13} - 48 \, \overline{\alpha_{1}} \overline{\alpha_{5}} - {\left(\theta_{0}^{3} \overline{\alpha_{2}} + 7 \, \theta_{0}^{2} \overline{\alpha_{2}} + 12 \, \theta_{0} \overline{\alpha_{2}}\right)} \overline{\zeta_{0}} - {\left(\theta_{0}^{3} + 6 \, \theta_{0}^{2} + 11 \, \theta_{0} + 6\right)} \overline{\zeta_{16}} + 24 \, \alpha_{15}\bigg\}\\
		\lambda_{3}&=\frac{1}{\theta_{0}^{3} + 3 \, \theta_{0}^{2} + 2 \, \theta_{0}}\bigg\{\mu_{4} \theta_{0}^{4} - 8 \, {\left(\theta_{0}^{2} + 5 \, \theta_{0} + 6\right)} {\left| A_{1} \right|}^{4} + 6 \, \mu_{4} \theta_{0}^{3} - {\left(4 \, \alpha_{1} \overline{\alpha_{1}} - 6 \, \beta - 11 \, \mu_{4}\right)} \theta_{0}^{2} - 2 \, {\left(8 \, \alpha_{1} \overline{\alpha_{1}} - 9 \, \beta - 3 \, \mu_{4}\right)} \theta_{0}\\
		& - 12 \, \alpha_{1} \overline{\alpha_{1}} - 2 \, {\left(\theta_{0}^{4} + 5 \, \theta_{0}^{3} + 8 \, \theta_{0}^{2} + 4 \, \theta_{0}\right)} \overline{\zeta_{10}} + 12 \, \beta\bigg\}\\
		\lambda_{4}&=\frac{1}{\theta_{0}^{3} + 3 \, \theta_{0}^{2} + 2 \, \theta_{0}}\bigg\{\mu_{16} \theta_{0}^{4} + 6 \, \mu_{16} \theta_{0}^{3} - {\left(6 \, \alpha_{5} \overline{\alpha_{1}} + 6 \, \alpha_{1} \overline{\alpha_{3}} - 6 \, \alpha_{16} - 11 \, \mu_{16}\right)} \theta_{0}^{2} - 4 \, {\left(3 \, \theta_{0}^{2} \overline{\alpha_{5}} + 13 \, \theta_{0} \overline{\alpha_{5}} + 12 \, \overline{\alpha_{5}}\right)} {\left| A_{1} \right|}^{2}\\
		& - 2 \, {\left(13 \, \alpha_{5} \overline{\alpha_{1}} + 9 \, \alpha_{1} \overline{\alpha_{3}} - 9 \, \alpha_{16} - 3 \, \mu_{16}\right)} \theta_{0} - {\left(\theta_{0}^{3} + 6 \, \theta_{0}^{2} + 11 \, \theta_{0} + 6\right)} \zeta_{14} - 24 \, \alpha_{5} \overline{\alpha_{1}} - 12 \, \alpha_{1} \overline{\alpha_{3}}\\
		& - 2 \, {\left(\theta_{0}^{2} + 3 \, \theta_{0} + 2\right)} \overline{\zeta_{15}} + 12 \, \alpha_{16}\bigg\}\\
		\lambda_{5}&=\frac{1}{\theta_{0}^{3} + 3 \, \theta_{0}^{2} + 2 \, \theta_{0}}\bigg\{{\left(\alpha_{2} \overline{\alpha_{1}} - \alpha_{10}\right)} \theta_{0}^{4} + 2 \, {\left(3 \, \alpha_{2} \overline{\alpha_{1}} - 2 \, \alpha_{10} + \overline{\mu_{27}}\right)} \theta_{0}^{3} + {\left(3 \, \alpha_{2} \overline{\alpha_{1}} + \alpha_{10} + 6 \, \overline{\mu_{27}}\right)} \theta_{0}^{2}\\
		& - 2 \, {\left(7 \, \alpha_{2} \overline{\alpha_{1}} - 8 \, \alpha_{10} - 2 \, \overline{\mu_{27}}\right)} \theta_{0} - {\left(\theta_{0}^{3} + 6 \, \theta_{0}^{2} + 11 \, \theta_{0} + 6\right)} \zeta_{18} - 12 \, \alpha_{2} \overline{\alpha_{1}} - 4 \, {\left(\alpha_{5} \theta_{0} + 3 \, \alpha_{5}\right)} \overline{\zeta_{0}} - 8 \, {\left(2 \, \alpha_{1} \theta_{0} + 3 \, \alpha_{1}\right)} \overline{\zeta_{1}}\\
		& - 2 \, {\left(\theta_{0}^{2} + 3 \, \theta_{0} + 2\right)} \overline{\zeta_{24}} + 12 \, \alpha_{10}\bigg\}\\
		\lambda_{6}&=-\frac{1}{\theta_{0}^{2} + \theta_{0}}\bigg\{\alpha_{11} \theta_{0}^{3} - 2 \, {\left(\alpha_{2} \theta_{0}^{3} + 3 \, \alpha_{2} \theta_{0}^{2} - 10 \, \alpha_{2} \theta_{0} - 12 \, \alpha_{2}\right)} {\left| A_{1} \right|}^{2} - 3 \, \theta_{0}^{2} \overline{\mu_{26}} - {\left(7 \, \alpha_{11} + 3 \, \overline{\mu_{26}}\right)} \theta_{0} + {\left(\theta_{0}^{2} + 3 \, \theta_{0} + 2\right)} \zeta_{19}\\
		& + 6 \, \alpha_{3} \overline{\zeta_{0}} + 3 \, {\left(\theta_{0} + 1\right)} \overline{\zeta_{23}} - 6 \, \alpha_{11}\bigg\}\\
		\lambda_{7}&=\frac{1}{\theta_{0}^{3} + \theta_{0}^{2}}\bigg\{\mu_{17} \theta_{0}^{4} + 3 \, \mu_{17} \theta_{0}^{3} - 2 \, {\left(3 \, \alpha_{3} \overline{\alpha_{1}} + 3 \, \alpha_{1} \overline{\alpha_{5}} - \mu_{17} - 3 \, \overline{\alpha_{16}}\right)} \theta_{0}^{2} - 12 \, {\left(\alpha_{5} \theta_{0}^{2} + 2 \, \alpha_{5} \theta_{0}\right)} {\left| A_{1} \right|}^{2}\\
		& - 6 \, {\left(\alpha_{3} \overline{\alpha_{1}} + \alpha_{1} \overline{\alpha_{5}} - \overline{\alpha_{16}}\right)} \theta_{0} - {\left(\theta_{0}^{3} + 3 \, \theta_{0}^{2} + 2 \, \theta_{0}\right)} \zeta_{15} + {\left(\theta_{0}^{2} \overline{\alpha_{1}} + \theta_{0} \overline{\alpha_{1}} - 3 \, \overline{\alpha_{1}}\right)} \zeta_{2} - 3 \, {\left(\theta_{0}^{2} + \theta_{0}\right)} \overline{\zeta_{14}}\bigg\}\\
		\lambda_{8}&=\frac{1}{\theta_{0}^{2} + 2 \, \theta_{0}}\bigg\{{\left(\alpha_{1} \alpha_{2} - \alpha_{12}\right)} \theta_{0}^{3} + {\left(3 \, \alpha_{1} \alpha_{2} + \alpha_{12} + 4 \, \overline{\mu_{24}}\right)} \theta_{0}^{2} - 16 \, \alpha_{1} \alpha_{2} - 2 \, {\left(7 \, \alpha_{1} \alpha_{2} - 5 \, \alpha_{12} - 4 \, \overline{\mu_{24}}\right)} \theta_{0}\\
		& - {\left(\theta_{0}^{2} + 3 \, \theta_{0} + 2\right)} \zeta_{20} - 4 \, {\left(\theta_{0} + 2\right)} \overline{\zeta_{22}} + 8 \, \alpha_{12}\bigg\}\\
		\lambda_{9}&=\frac{1}{\theta_{0}^{4} + 3 \, \theta_{0}^{3} + 2 \, \theta_{0}^{2}}\bigg\{\mu_{18} \theta_{0}^{5} + 4 \, \mu_{18} \theta_{0}^{4} - {\left(4 \, \alpha_{1} \alpha_{5} - 5 \, \mu_{18} - 4 \, \overline{\alpha_{15}}\right)} \theta_{0}^{3} - 2 \, {\left(6 \, \alpha_{1} \alpha_{5} - \mu_{18} - 6 \, \overline{\alpha_{15}}\right)} \theta_{0}^{2}\\
		& - 2 \, {\left(4 \, \alpha_{3} \theta_{0}^{3} + 12 \, \alpha_{3} \theta_{0}^{2} + 8 \, \alpha_{3} \theta_{0} - {\left(\theta_{0}^{3} + 2 \, \theta_{0}^{2} - 2 \, \theta_{0} - 4\right)} \zeta_{2}\right)} {\left| A_{1} \right|}^{2} - 8 \, {\left(\alpha_{1} \alpha_{5} - \overline{\alpha_{15}}\right)} \theta_{0} - 8 \, {\left(\alpha_{2} \theta_{0}^{2} + \alpha_{2} \theta_{0}\right)} \zeta_{0}\\
		& - {\left(\theta_{0}^{4} + 4 \, \theta_{0}^{3} + 5 \, \theta_{0}^{2} + 2 \, \theta_{0}\right)} \zeta_{16} - 4 \, {\left(\theta_{0}^{3} + 3 \, \theta_{0}^{2} + 2 \, \theta_{0}\right)} \overline{\zeta_{13}}\bigg\}\\
		\lambda_{10}&=\frac{1}{\theta_{0}^{2} + 2 \, \theta_{0}}\bigg\{{\left(\alpha_{2} \overline{\alpha_{1}} - \alpha_{10} + 2 \, \mu_{20}\right)} \theta_{0}^{3} + {\left(7 \, \alpha_{2} \overline{\alpha_{1}} - 3 \, \alpha_{10} + 7 \, \mu_{20}\right)} \theta_{0}^{2} - 2 \, {\left(3 \, \alpha_{2} \overline{\alpha_{1}} - 2 \, \alpha_{10} - 3 \, \mu_{20}\right)} \theta_{0}\\
		& - {\left(2 \, \theta_{0}^{2} + 7 \, \theta_{0} + 6\right)} \zeta_{18} - 24 \, \alpha_{2} \overline{\alpha_{1}} + {\left(\theta_{0}^{2} - 4\right)} \overline{\zeta_{24}} + 12 \, \alpha_{10}\bigg\}\\
		\lambda_{11}&=-\frac{1}{\theta_{0}^{3} + 3 \, \theta_{0}^{2} + 2 \, \theta_{0}}\bigg\{{\left(\alpha_{11} - 2 \, \mu_{21}\right)} \theta_{0}^{4} + 2 \, {\left(\alpha_{11} - 4 \, \mu_{21}\right)} \theta_{0}^{3} - {\left(7 \, \alpha_{11} + 10 \, \mu_{21}\right)} \theta_{0}^{2}\\
		& - 2 \, {\left(\alpha_{2} \theta_{0}^{4} + 4 \, \alpha_{2} \theta_{0}^{3} - 7 \, \alpha_{2} \theta_{0}^{2} - 34 \, \alpha_{2} \theta_{0} - 24 \, \alpha_{2}\right)} {\left| A_{1} \right|}^{2} - 4 \, {\left(5 \, \alpha_{11} + \mu_{21}\right)} \theta_{0} + 2 \, {\left(\theta_{0}^{3} + 4 \, \theta_{0}^{2} + 5 \, \theta_{0} + 2\right)} \zeta_{19} + 2 \, \zeta_{2} \overline{\zeta_{0}}\\
		& - {\left(\theta_{0}^{3} - 7 \, \theta_{0} - 6\right)} \overline{\zeta_{23}} - 12 \, \alpha_{11}\bigg\}\\
		\lambda_{12}&=\frac{1}{\theta_{0}^{2} + 2 \, \theta_{0}}\bigg\{{\left(\overline{\alpha_{1}} \overline{\alpha_{2}} - \overline{\alpha_{12}}\right)} \theta_{0}^{3} + {\left(3 \, \overline{\alpha_{1}} \overline{\alpha_{2}} + 4 \, \mu_{24} + \overline{\alpha_{12}}\right)} \theta_{0}^{2} - 2 \, {\left(7 \, \overline{\alpha_{1}} \overline{\alpha_{2}} - 4 \, \mu_{24} - 5 \, \overline{\alpha_{12}}\right)} \theta_{0} - 4 \, {\left(\theta_{0} + 2\right)} \zeta_{22}\\
		& - 16 \, \overline{\alpha_{1}} \overline{\alpha_{2}} - {\left(\theta_{0}^{2} + 3 \, \theta_{0} + 2\right)} \overline{\zeta_{20}} + 8 \, \overline{\alpha_{12}}\bigg\}\\
		\lambda_{13}&=\frac{1}{\theta_{0}^{4} + 3 \, \theta_{0}^{3} + 2 \, \theta_{0}^{2}}\bigg\{\theta_{0}^{5} \overline{\mu_{18}} + 4 \, \theta_{0}^{4} \overline{\mu_{18}} - {\left(4 \, \overline{\alpha_{1}} \overline{\alpha_{5}} - 4 \, \alpha_{15} - 5 \, \overline{\mu_{18}}\right)} \theta_{0}^{3} - 2 \, {\left(6 \, \overline{\alpha_{1}} \overline{\alpha_{5}} - 6 \, \alpha_{15} - \overline{\mu_{18}}\right)} \theta_{0}^{2}\\
		& - 2 \, {\left(4 \, \theta_{0}^{3} \overline{\alpha_{3}} + 12 \, \theta_{0}^{2} \overline{\alpha_{3}} + 8 \, \theta_{0} \overline{\alpha_{3}} - {\left(\theta_{0}^{3} + 2 \, \theta_{0}^{2} - 2 \, \theta_{0} - 4\right)} \overline{\zeta_{2}}\right)} {\left| A_{1} \right|}^{2} - 8 \, {\left(\overline{\alpha_{1}} \overline{\alpha_{5}} - \alpha_{15}\right)} \theta_{0} - 4 \, {\left(\theta_{0}^{3} + 3 \, \theta_{0}^{2} + 2 \, \theta_{0}\right)} \zeta_{13}\\
		& - 8 \, {\left(\theta_{0}^{2} \overline{\alpha_{2}} + \theta_{0} \overline{\alpha_{2}}\right)} \overline{\zeta_{0}} - {\left(\theta_{0}^{4} + 4 \, \theta_{0}^{3} + 5 \, \theta_{0}^{2} + 2 \, \theta_{0}\right)} \overline{\zeta_{16}}\bigg\}\\
		\lambda_{14}&=-\frac{1}{\theta_{0}^{2} + \theta_{0}}\bigg\{\theta_{0}^{3} \overline{\alpha_{11}} - 3 \, \mu_{26} \theta_{0}^{2} - 2 \, {\left(\theta_{0}^{3} \overline{\alpha_{2}} + 3 \, \theta_{0}^{2} \overline{\alpha_{2}} - 10 \, \theta_{0} \overline{\alpha_{2}} - 12 \, \overline{\alpha_{2}}\right)} {\left| A_{1} \right|}^{2} - {\left(3 \, \mu_{26} + 7 \, \overline{\alpha_{11}}\right)} \theta_{0} + 3 \, {\left(\theta_{0} + 1\right)} \zeta_{23}\\
		& + 6 \, \zeta_{0} \overline{\alpha_{3}} + {\left(\theta_{0}^{2} + 3 \, \theta_{0} + 2\right)} \overline{\zeta_{19}} - 6 \, \overline{\alpha_{11}}\bigg\}\\
		\lambda_{15}&=\frac{1}{\theta_{0}^{3} + \theta_{0}^{2}}\bigg\{\theta_{0}^{4} \overline{\mu_{17}} + 3 \, \theta_{0}^{3} \overline{\mu_{17}} - 2 \, {\left(3 \, \alpha_{5} \overline{\alpha_{1}} + 3 \, \alpha_{1} \overline{\alpha_{3}} - 3 \, \alpha_{16} - \overline{\mu_{17}}\right)} \theta_{0}^{2} - 12 \, {\left(\theta_{0}^{2} \overline{\alpha_{5}} + 2 \, \theta_{0} \overline{\alpha_{5}}\right)} {\left| A_{1} \right|}^{2}\\
		& - 6 \, {\left(\alpha_{5} \overline{\alpha_{1}} + \alpha_{1} \overline{\alpha_{3}} - \alpha_{16}\right)} \theta_{0} - 3 \, {\left(\theta_{0}^{2} + \theta_{0}\right)} \zeta_{14} - {\left(\theta_{0}^{3} + 3 \, \theta_{0}^{2} + 2 \, \theta_{0}\right)} \overline{\zeta_{15}} + {\left(\alpha_{1} \theta_{0}^{2} + \alpha_{1} \theta_{0} - 3 \, \alpha_{1}\right)} \overline{\zeta_{2}}\bigg\}\\
		\lambda_{16}&=\frac{1}{\theta_{0}^{3} + 3 \, \theta_{0}^{2} + 2 \, \theta_{0}}\bigg\{\theta_{0}^{4} \overline{\mu_{16}} + 6 \, \theta_{0}^{3} \overline{\mu_{16}} - {\left(6 \, \alpha_{3} \overline{\alpha_{1}} + 6 \, \alpha_{1} \overline{\alpha_{5}} - 6 \, \overline{\alpha_{16}} - 11 \, \overline{\mu_{16}}\right)} \theta_{0}^{2} - 4 \, {\left(3 \, \alpha_{5} \theta_{0}^{2} + 13 \, \alpha_{5} \theta_{0} + 12 \, \alpha_{5}\right)} {\left| A_{1} \right|}^{2}\\
		& - 2 \, {\left(9 \, \alpha_{3} \overline{\alpha_{1}} + 13 \, \alpha_{1} \overline{\alpha_{5}} - 9 \, \overline{\alpha_{16}} - 3 \, \overline{\mu_{16}}\right)} \theta_{0} - 2 \, {\left(\theta_{0}^{2} + 3 \, \theta_{0} + 2\right)} \zeta_{15} - 12 \, \alpha_{3} \overline{\alpha_{1}} - 24 \, \alpha_{1} \overline{\alpha_{5}}\\
		& - {\left(\theta_{0}^{3} + 6 \, \theta_{0}^{2} + 11 \, \theta_{0} + 6\right)} \overline{\zeta_{14}} + 12 \, \overline{\alpha_{16}}\bigg\}\\
		\lambda_{17}&=-\frac{1}{\theta_{0}^{3} + 3 \, \theta_{0}^{2} + 2 \, \theta_{0}}\bigg\{8 \, {\left(\theta_{0}^{2} + 5 \, \theta_{0} + 6\right)} {\left| A_{1} \right|}^{4} - \theta_{0}^{4} \overline{\mu_{4}} - 6 \, \theta_{0}^{3} \overline{\mu_{4}} + {\left(4 \, \alpha_{1} \overline{\alpha_{1}} - 6 \, \beta - 11 \, \overline{\mu_{4}}\right)} \theta_{0}^{2} + 2 \, {\left(8 \, \alpha_{1} \overline{\alpha_{1}} - 9 \, \beta - 3 \, \overline{\mu_{4}}\right)} \theta_{0}\\
		& + 2 \, {\left(\theta_{0}^{4} + 5 \, \theta_{0}^{3} + 8 \, \theta_{0}^{2} + 4 \, \theta_{0}\right)} \zeta_{10} + 12 \, \alpha_{1} \overline{\alpha_{1}} - 12 \, \beta\bigg\}\\
		\lambda_{18}&=\frac{1}{\theta_{0}^{3} + 3 \, \theta_{0}^{2} + 2 \, \theta_{0}}\bigg\{{\left(\alpha_{1} \overline{\alpha_{2}} - \overline{\alpha_{10}}\right)} \theta_{0}^{4} + 2 \, {\left(3 \, \alpha_{1} \overline{\alpha_{2}} + \mu_{27} - 2 \, \overline{\alpha_{10}}\right)} \theta_{0}^{3} + {\left(3 \, \alpha_{1} \overline{\alpha_{2}} + 6 \, \mu_{27} + \overline{\alpha_{10}}\right)} \theta_{0}^{2}\\
		& - 2 \, {\left(7 \, \alpha_{1} \overline{\alpha_{2}} - 2 \, \mu_{27} - 8 \, \overline{\alpha_{10}}\right)} \theta_{0} - 4 \, {\left(\theta_{0} \overline{\alpha_{5}} + 3 \, \overline{\alpha_{5}}\right)} \zeta_{0} - 8 \, {\left(2 \, \theta_{0} \overline{\alpha_{1}} + 3 \, \overline{\alpha_{1}}\right)} \zeta_{1} - 2 \, {\left(\theta_{0}^{2} + 3 \, \theta_{0} + 2\right)} \zeta_{24} - 12 \, \alpha_{1} \overline{\alpha_{2}}\\
		& - {\left(\theta_{0}^{3} + 6 \, \theta_{0}^{2} + 11 \, \theta_{0} + 6\right)} \overline{\zeta_{18}} + 12 \, \overline{\alpha_{10}}\bigg\}\\
		\lambda_{19}&=\frac{1}{\theta_{0}^{4} + 6 \, \theta_{0}^{3} + 11 \, \theta_{0}^{2} + 6 \, \theta_{0}}\bigg\{\theta_{0}^{5} \overline{\mu_{2}} + 10 \, \theta_{0}^{4} \overline{\mu_{2}} + {\left(4 \, \alpha_{6} + 35 \, \overline{\mu_{2}}\right)} \theta_{0}^{3} + 2 \, {\left(12 \, \alpha_{6} + 25 \, \overline{\mu_{2}}\right)} \theta_{0}^{2}\\
		& - 2 \, {\left(3 \, \alpha_{1} \theta_{0}^{3} + 25 \, \alpha_{1} \theta_{0}^{2} + 64 \, \alpha_{1} \theta_{0} + 48 \, \alpha_{1}\right)} {\left| A_{1} \right|}^{2} + 4 \, {\left(11 \, \alpha_{6} + 6 \, \overline{\mu_{2}}\right)} \theta_{0} - 2 \, {\left(\theta_{0}^{3} + 3 \, \theta_{0}^{2} + 2 \, \theta_{0}\right)} \zeta_{11} + 24 \, \alpha_{6}\bigg\}\\
		\lambda_{20}&=\frac{1}{\theta_{0}^{4} + 6 \, \theta_{0}^{3} + 11 \, \theta_{0}^{2} + 6 \, \theta_{0}}\bigg\{\theta_{0}^{5} \overline{\mu_{15}} + 10 \, \theta_{0}^{4} \overline{\mu_{15}} - {\left(4 \, \alpha_{1} \alpha_{5} - 4 \, \overline{\alpha_{15}} - 35 \, \overline{\mu_{15}}\right)} \theta_{0}^{3} - 2 \, {\left(16 \, \alpha_{1} \alpha_{5} - 12 \, \overline{\alpha_{15}} - 25 \, \overline{\mu_{15}}\right)} \theta_{0}^{2}\\
		& - {\left(8 \, \alpha_{3} \theta_{0}^{3} + 60 \, \alpha_{3} \theta_{0}^{2} + 136 \, \alpha_{3} \theta_{0} + {\left(\theta_{0} + 2\right)} \zeta_{2} + 96 \, \alpha_{3}\right)} {\left| A_{1} \right|}^{2} - 48 \, \alpha_{1} \alpha_{5} - 4 \, {\left(19 \, \alpha_{1} \alpha_{5} - 11 \, \overline{\alpha_{15}} - 6 \, \overline{\mu_{15}}\right)} \theta_{0} \\
		&- {\left(\alpha_{2} \theta_{0}^{3} + 7 \, \alpha_{2} \theta_{0}^{2} + 12 \, \alpha_{2} \theta_{0}\right)} \zeta_{0} - {\left(\theta_{0}^{3} + 6 \, \theta_{0}^{2} + 11 \, \theta_{0} + 6\right)} \zeta_{16} - {\left(\theta_{0}^{4} + 10 \, \theta_{0}^{3} + 35 \, \theta_{0}^{2} + 50 \, \theta_{0} + 24\right)} \overline{\zeta_{13}} + 24 \, \overline{\alpha_{15}}\bigg\}\\
		\lambda_{21}&=-\frac{1}{\theta_{0}^{3} + 3 \, \theta_{0}^{2} + 2 \, \theta_{0}}\bigg\{\theta_{0}^{4} {\left(\overline{\alpha_{11}} - 2 \, \overline{\mu_{21}}\right)} + 2 \, \theta_{0}^{3} {\left(\overline{\alpha_{11}} - 4 \, \overline{\mu_{21}}\right)} - \theta_{0}^{2} {\left(7 \, \overline{\alpha_{11}} + 10 \, \overline{\mu_{21}}\right)}\\
		& - 2 \, {\left(\theta_{0}^{4} \overline{\alpha_{2}} + 4 \, \theta_{0}^{3} \overline{\alpha_{2}} - 7 \, \theta_{0}^{2} \overline{\alpha_{2}} - 34 \, \theta_{0} \overline{\alpha_{2}} - 24 \, \overline{\alpha_{2}}\right)} {\left| A_{1} \right|}^{2} - {\left(\theta_{0}^{3} - 7 \, \theta_{0} - 6\right)} \zeta_{23} - 4 \, \theta_{0} {\left(5 \, \overline{\alpha_{11}} + \overline{\mu_{21}}\right)}\\
		& + 2 \, {\left(\theta_{0}^{3} + 4 \, \theta_{0}^{2} + 5 \, \theta_{0} + 2\right)} \overline{\zeta_{19}} + 2 \, \zeta_{0} \overline{\zeta_{2}} - 12 \, \overline{\alpha_{11}}\bigg\}\\
		\lambda_{22}&=\frac{1}{\theta_{0}^{2} + 2 \, \theta_{0}}\bigg\{{\left(\alpha_{1} \overline{\alpha_{2}} - \overline{\alpha_{10}} + 2 \, \overline{\mu_{20}}\right)} \theta_{0}^{3} + {\left(7 \, \alpha_{1} \overline{\alpha_{2}} - 3 \, \overline{\alpha_{10}} + 7 \, \overline{\mu_{20}}\right)} \theta_{0}^{2} - 2 \, {\left(3 \, \alpha_{1} \overline{\alpha_{2}} - 2 \, \overline{\alpha_{10}} - 3 \, \overline{\mu_{20}}\right)} \theta_{0} + {\left(\theta_{0}^{2} - 4\right)} \zeta_{24}\\
		& - 24 \, \alpha_{1} \overline{\alpha_{2}} - {\left(2 \, \theta_{0}^{2} + 7 \, \theta_{0} + 6\right)} \overline{\zeta_{18}} + 12 \, \overline{\alpha_{10}}\bigg\}
	\end{align*}
	\normalsize
	    \textbf{We shall only need the precise expression of $\lambda_2$.}
	    
	    Also, recall that
	    \small
	    \begin{align*}
	    &\zeta_{13}=\frac{1}{\theta_{0}^{4} - 7 \, \theta_{0}^{2} - 6 \, \theta_{0}}\bigg\{{\left(\theta_{0}^{3} - 7 \, \theta_{0} - 6\right)} B_{7} \overline{A_{0}} + 2 \, {\left(\theta_{0}^{5} \overline{\alpha_{2}} - 7 \, \theta_{0}^{3} \overline{\alpha_{2}} - 6 \, \theta_{0}^{2} \overline{\alpha_{2}}\right)} \overline{A_{1}}^{2} + 2 \, {\left(\theta_{0}^{5} \overline{\alpha_{2}} - \theta_{0}^{4} \overline{\alpha_{2}} - 5 \, \theta_{0}^{3} \overline{\alpha_{2}} - 3 \, \theta_{0}^{2} \overline{\alpha_{2}}\right)} \overline{A_{0}} \overline{A_{2}}\\
	    & + {\left({\left(\theta_{0}^{2} + 2 \, \theta_{0}\right)} A_{1} \overline{\zeta_{2}} + {\left(\theta_{0}^{3} - \theta_{0}^{2} - 6 \, \theta_{0}\right)} \overline{A_{0}} \overline{\zeta_{3}}\right)} \overline{A_{1}}\bigg\}\\
	    &\zeta_{14}=\frac{1}{48 \, {\left(\theta_{0}^{2} + \theta_{0}\right)}}\bigg\{12 \, A_{0} {\left(\theta_{0} - 2\right)} {\left| C_{1} \right|}^{2} \overline{A_{1}} + 96 \, \theta_{0} \overline{A_{0}} \overline{A_{1}} \overline{\zeta_{4}} + 48 \, A_{0} B_{10} {\left(\theta_{0} + 1\right)} + 48 \, B_{8} {\left(\theta_{0} + 1\right)} \overline{A_{0}}\\
	    & - 3 \, {\left(\theta_{0}^{2} - \theta_{0} - 2\right)} B_{1} \overline{C_{1}} - 2 \, {\left(\theta_{0}^{2} - \theta_{0} - 2\right)} C_{1} \overline{C_{2}}\bigg\}\\
	    &\zeta_{15}=-\frac{1}{16 \, {\left(\theta_{0}^{2} + \theta_{0}\right)}}\bigg\{8 \, A_{1} {\left(\theta_{0} - 2\right)} {\left| C_{1} \right|}^{2} \overline{A_{0}} - 4 \, {\left(\theta_{0} + 1\right)} \zeta_{2} \overline{A_{0}} \overline{C_{1}} + 8 \, {\left(\theta_{0}^{2} - 2 \, \theta_{0}\right)} A_{1} B_{2} - 16 \, A_{0} B_{11} {\left(\theta_{0} + 1\right)}\\
	    & - 16 \, B_{9} {\left(\theta_{0} + 1\right)} \overline{A_{0}} + {\left(\theta_{0}^{2} - 2 \, \theta_{0} - 3\right)} C_{2} \overline{C_{1}} - {\left(4 \, {\left(\theta_{0}^{3} \overline{\alpha_{2}} + 2 \, \theta_{0}^{2} \overline{\alpha_{2}} + \theta_{0} \overline{\alpha_{2}}\right)} \overline{A_{0}} - {\left(\theta_{0}^{2} - 2 \, \theta_{0}\right)} \overline{B_{1}}\right)} C_{1}\bigg\}\\
	    &\zeta_{16}=\frac{1}{2 \, {\left(\theta_{0}^{3} + 3 \, \theta_{0}^{2} + 2 \, \theta_{0}\right)}}\bigg\{4 \, {\left(\theta_{0}^{2} + 2 \, \theta_{0}\right)} A_{1} \zeta_{2} \overline{A_{0}} + 4 \, {\left(\theta_{0}^{2} + 3 \, \theta_{0} + 2\right)} A_{0} \zeta_{5} \overline{A_{0}} + 4 \, {\left(\alpha_{7} \theta_{0}^{4} - \alpha_{7} \theta_{0}^{3} - 10 \, \alpha_{7} \theta_{0}^{2} - 8 \, \alpha_{7} \theta_{0}\right)} A_{0}^{2}\\
	    & - {\left(\theta_{0}^{3} - \theta_{0}^{2} - 2 \, \theta_{0}\right)} A_{2} C_{1} - {\left(\theta_{0}^{3} - \theta_{0}^{2} - 6 \, \theta_{0}\right)} A_{1} C_{2} - {\left(\theta_{0}^{3} - \theta_{0}^{2} - 10 \, \theta_{0} - 8\right)} A_{0} C_{3}\bigg\}\\	
	    &\zeta_{17}=\frac{1}{2 \, {\left(\theta_{0}^{4} + 6 \, \theta_{0}^{3} + 11 \, \theta_{0}^{2} + 6 \, \theta_{0}\right)}}\bigg\{4 \, {\left(\theta_{0}^{3} + 4 \, \theta_{0}^{2} + 3 \, \theta_{0}\right)} A_{2} \zeta_{2} \overline{A_{0}} + 2 \, {\left(\theta_{0}^{3} + 6 \, \theta_{0}^{2} + 11 \, \theta_{0} + 6\right)} A_{0} B_{13}\\
	    & - {\left(\theta_{0}^{4} + \theta_{0}^{3} - 4 \, \theta_{0}^{2} - 4 \, \theta_{0}\right)} A_{3} C_{1} - {\left(\theta_{0}^{4} + \theta_{0}^{3} - 9 \, \theta_{0}^{2} - 9 \, \theta_{0}\right)} A_{2} C_{2} - {\left(\theta_{0}^{4} + \theta_{0}^{3} - 14 \, \theta_{0}^{2} - 24 \, \theta_{0}\right)} A_{1} C_{3}\\
	    & + 4 \, {\left({\left(\theta_{0}^{3} + 5 \, \theta_{0}^{2} + 6 \, \theta_{0}\right)} \zeta_{5} \overline{A_{0}} + {\left(\alpha_{7} \theta_{0}^{5} + \alpha_{7} \theta_{0}^{4} - 14 \, \alpha_{7} \theta_{0}^{3} - 24 \, \alpha_{7} \theta_{0}^{2}\right)} A_{0}\right)} A_{1}\bigg\}\\
	    &\zeta_{18}=-\frac{{\left(\theta_{0}^{3} + \theta_{0}^{2} - 4 \, \theta_{0} - 4\right)} {\left| C_{1} \right|}^{2} \overline{A_{0}}^{2} + {\left(\theta_{0}^{4} + \theta_{0}^{3} - 4 \, \theta_{0}^{2} - 4 \, \theta_{0}\right)} B_{2} \overline{A_{0}} + {\left(\theta_{0}^{4} - 4 \, \theta_{0}^{2}\right)} B_{1} \overline{A_{1}} + {\left(\theta_{0}^{4} - \theta_{0}^{3} - 2 \, \theta_{0}^{2}\right)} C_{1} \overline{A_{2}}}{2 \, {\left(\theta_{0}^{4} + 3 \, \theta_{0}^{3} + 2 \, \theta_{0}^{2}\right)}}\\
	    &\zeta_{19}=-\frac{1}{2 \, {\left(\theta_{0}^{4} + 6 \, \theta_{0}^{3} + 11 \, \theta_{0}^{2} + 6 \, \theta_{0}\right)}}\bigg\{{\left(\theta_{0}^{3} + 3 \, \theta_{0}^{2} - 4 \, \theta_{0} - 12\right)} {\left| C_{1} \right|}^{2} \overline{A_{0}} \overline{A_{1}} - 2 \, {\left(\theta_{0}^{3} + 6 \, \theta_{0}^{2} + 11 \, \theta_{0} + 6\right)} B_{10} \overline{A_{0}}\\
	    & + {\left(\theta_{0}^{4} + 3 \, \theta_{0}^{3} - 4 \, \theta_{0}^{2} - 12 \, \theta_{0}\right)} B_{2} \overline{A_{1}} + {\left(\theta_{0}^{4} + 2 \, \theta_{0}^{3} - 5 \, \theta_{0}^{2} - 6 \, \theta_{0}\right)} B_{1} \overline{A_{2}} + {\left(\theta_{0}^{4} + \theta_{0}^{3} - 4 \, \theta_{0}^{2} - 4 \, \theta_{0}\right)} C_{1} \overline{A_{3}}\bigg\}\\
	    &\zeta_{20}=\frac{4 \, {\left(\theta_{0}^{2} + 2 \, \theta_{0}\right)} \zeta_{2} \overline{A_{1}}^{2} + 4 \, {\left(\theta_{0}^{2} + \theta_{0}\right)} \zeta_{2} \overline{A_{0}} \overline{A_{2}} + 2 \, {\left(\theta_{0}^{2} + 3 \, \theta_{0} + 2\right)} B_{11} \overline{A_{0}} - {\left(\theta_{0}^{3} - \theta_{0}^{2} - 6 \, \theta_{0}\right)} B_{3} \overline{A_{1}} - {\left(\theta_{0}^{3} - 2 \, \theta_{0}^{2} - 3 \, \theta_{0}\right)} C_{2} \overline{A_{2}}}{2 \, {\left(\theta_{0}^{3} + 3 \, \theta_{0}^{2} + 2 \, \theta_{0}\right)}}\\
	    &\zeta_{21}=\frac{1}{16 \, {\left(\theta_{0}^{2} + \theta_{0}\right)}}\bigg\{16 \, B_{12} {\left(\theta_{0} + 1\right)} \overline{A_{0}} - 8 \, {\left(\theta_{0}^{2} - 4 \, \theta_{0}\right)} C_{3} \overline{A_{1}} + {\left(16 \, {\left(\alpha_{2} \theta_{0}^{2} - \alpha_{2} \theta_{0} - 2 \, \alpha_{2}\right)} A_{0} - {\left(2 \, \theta_{0}^{2} - 3 \, \theta_{0} - 2\right)} B_{1}\right)} C_{1}\\
	    & + 32 \, {\left(\theta_{0} \zeta_{5} \overline{A_{0}} + {\left(\alpha_{7} \theta_{0}^{3} - 5 \, \alpha_{7} \theta_{0}^{2} + 6 \, \alpha_{7} \theta_{0} - 8 \, \alpha_{7}\right)} A_{0}\right)} \overline{A_{1}}\bigg\}
	    \end{align*}
	    \normalsize
	    \textbf{Here, we will only need $\zeta_{13}$ and $\zeta_{16}$.}
	    	
	    \section{The coefficient in $(0, \theta_0+3)$ in the Taylor expansion of $\Re(\p{\z}\vec{F}(z))=0$}
	    	
	    The coefficient in $\z^{\theta_0+3}$ in the Taylor	expansion of $\Re(\p{\z}\vec{F}(z))=0$ is
	    \small
	    \begin{align}\label{deftheta0+3}
	    \begin{dmatrix}
	    \frac{3 \, \overline{\zeta_{0}}}{4 \, {\left(\theta_{0}^{3} + 4 \, \theta_{0}^{2} + 5 \, \theta_{0} + 2\right)}} & \overline{B_{1}} & 0 & \theta_{0} + 3 &\textbf{(1)}\\
	    \frac{\mu_{1} \theta_{0}^{3} + 6 \, \mu_{1} \theta_{0}^{2} + 2 \, {\left(4 \, \mu_{1} + \overline{\alpha_{5}}\right)} \theta_{0} + 8 \, \overline{\alpha_{5}}}{\theta_{0}^{2} + 2 \, \theta_{0}} & \overline{A_{2}} & 0 & \theta_{0} + 3  &\textbf{(2)}\\
	    \frac{\theta_{0}^{3} \overline{\mu_{14}} + 3 \, \theta_{0}^{2} \overline{\mu_{14}} + 2 \, \theta_{0} \overline{\mu_{14}} - 2 \, {\left(\theta_{0} \overline{\alpha_{1}} + 4 \, \overline{\alpha_{1}}\right)} \overline{\zeta_{0}} - 2 \, {\left(\theta_{0}^{5} + 9 \, \theta_{0}^{4} + 28 \, \theta_{0}^{3} + 36 \, \theta_{0}^{2} + 16 \, \theta_{0}\right)} \overline{\zeta_{9}}}{\theta_{0}^{3} + 3 \, \theta_{0}^{2} + 2 \, \theta_{0}} & A_{1} & 0 & \theta_{0} + 3  &\textbf{(3)}\\
	    \frac{\theta_{0}^{4} \overline{\mu_{28}} + 6 \, \theta_{0}^{3} \overline{\mu_{28}} - 8 \, {\left(2 \, \theta_{0}^{2} + 11 \, \theta_{0} + 12\right)} {\left| A_{1} \right|}^{2} \overline{\zeta_{1}} + 11 \, \theta_{0}^{2} \overline{\mu_{28}} + 6 \, \theta_{0} \overline{\mu_{28}} - 2 \, {\left(\theta_{0}^{2} \overline{\alpha_{5}} + 7 \, \theta_{0} \overline{\alpha_{5}} + 12 \, \overline{\alpha_{5}}\right)} \overline{\zeta_{0}}}{\theta_{0}^{4} + 6 \, \theta_{0}^{3} + 11 \, \theta_{0}^{2} + 6 \, \theta_{0}} & A_{0} & 0 & \theta_{0} + 3  &\textbf{(4)}\\
	    \lambda_1 & \overline{A_{1}} & 0 & \theta_{0} + 3  &\textbf{(5)}\\
	    \lambda_2 & \overline{A_{0}} & 0 & \theta_{0} + 3  &\textbf{(6)}\\
	    \end{dmatrix}
	    \end{align}
	    where
	    \begin{align*}
	    \lambda_1&=\frac{1}{\theta_{0}^{4} + 6 \, \theta_{0}^{3} + 11 \, \theta_{0}^{2} + 6 \, \theta_{0}}\bigg\{\mu_{2} \theta_{0}^{5} + 10 \, \mu_{2} \theta_{0}^{4} + {\left(35 \, \mu_{2} + 4 \, \overline{\alpha_{6}}\right)} \theta_{0}^{3} + 2 \, {\left(25 \, \mu_{2} + 12 \, \overline{\alpha_{6}}\right)} \theta_{0}^{2}\\
	    & - 2 \, {\left(3 \, \theta_{0}^{3} \overline{\alpha_{1}} + 25 \, \theta_{0}^{2} \overline{\alpha_{1}} + 64 \, \theta_{0} \overline{\alpha_{1}} + 48 \, \overline{\alpha_{1}}\right)} {\left| A_{1} \right|}^{2} + 4 \, {\left(6 \, \mu_{2} + 11 \, \overline{\alpha_{6}}\right)} \theta_{0} - 2 \, {\left(\theta_{0}^{3} + 3 \, \theta_{0}^{2} + 2 \, \theta_{0}\right)} \overline{\zeta_{11}} + 24 \, \overline{\alpha_{6}}\bigg\}\\
	    \lambda_2&=\frac{1}{\theta_{0}^{4} + 6 \, \theta_{0}^{3} + 11 \, \theta_{0}^{2} + 6 \, \theta_{0}}\bigg\{\mu_{15} \theta_{0}^{5} + 10 \, \mu_{15} \theta_{0}^{4} - {\left(4 \, \overline{\alpha_{1}} \overline{\alpha_{5}} - 4 \, \alpha_{15} - 35 \, \mu_{15}\right)} \theta_{0}^{3} - 2 \, {\left(16 \, \overline{\alpha_{1}} \overline{\alpha_{5}} - 12 \, \alpha_{15} - 25 \, \mu_{15}\right)} \theta_{0}^{2}\\
	    & - {\left(8 \, \theta_{0}^{3} \overline{\alpha_{3}} + 60 \, \theta_{0}^{2} \overline{\alpha_{3}} + 136 \, \theta_{0} \overline{\alpha_{3}} + {\left(\theta_{0} + 2\right)} \overline{\zeta_{2}} + 96 \, \overline{\alpha_{3}}\right)} {\left| A_{1} \right|}^{2} - 4 \, {\left(19 \, \overline{\alpha_{1}} \overline{\alpha_{5}} - 11 \, \alpha_{15} - 6 \, \mu_{15}\right)} \theta_{0}\\
	    & - {\left(\theta_{0}^{4} + 10 \, \theta_{0}^{3} + 35 \, \theta_{0}^{2} + 50 \, \theta_{0} + 24\right)} \zeta_{13} - 48 \, \overline{\alpha_{1}} \overline{\alpha_{5}} - {\left(\theta_{0}^{3} \overline{\alpha_{2}} + 7 \, \theta_{0}^{2} \overline{\alpha_{2}} + 12 \, \theta_{0} \overline{\alpha_{2}}\right)} \overline{\zeta_{0}} - {\left(\theta_{0}^{3} + 6 \, \theta_{0}^{2} + 11 \, \theta_{0} + 6\right)} \overline{\zeta_{16}} + 24 \, \alpha_{15}\bigg\}
	    \end{align*}    
	    \normalsize
	    Now, recall that by \eqref{endend} and as 
	    \begin{align*}
	    	|\vec{A}_0|^2=\dfrac{1}{2},\quad \alpha_5=2\s{\bar{\vec{A}_1}}{\vec{A}_2},
	    \end{align*}
	    we obtain
	    \small
	    \begin{align*}
	    \mu_{1}&=-\frac{1}{4 \, {\left(\theta_{0}^{3} + 3 \, \theta_{0}^{2} + 2 \, \theta_{0}\right)}}\bigg\{\ccancel{16 \, {\left(\theta_{0}^{2} + 2 \, \theta_{0}\right)} A_{0} {\left| A_{1} \right|}^{2} \overline{A_{1}}} + 8 \, {\left(\theta_{0}^{2} \overline{\alpha_{5}} + 3 \, \theta_{0} \overline{\alpha_{5}} + 2 \, \overline{\alpha_{5}}\right)} A_{0} \overline{A_{0}} + \ccancel{8 \, {\left(\theta_{0}^{2} \overline{\alpha_{1}} + 3 \, \theta_{0} \overline{\alpha_{1}} + 2 \, \overline{\alpha_{1}}\right)} A_{1} \overline{A_{0}}}\nonumber\\
	    & - 8 \, {\left(\theta_{0}^{2} + \theta_{0}\right)} A_{1} \overline{A_{2}} - \ccancel{{\left(\theta_{0}^{2} + 3 \, \theta_{0} + 2\right)} \overline{A_{0}} \overline{B_{1}}}\bigg\}\\
	    &=-\frac{1}{4 \, {\left(\theta_{0}^{3} + 3 \, \theta_{0}^{2} + 2 \, \theta_{0}\right)}}\bigg\{4(\theta_0^2+3\theta_0+2)\bar{\alpha_5}-4(\theta_0^2+\theta_0)\bar{\alpha_5}\bigg\}\\
	    &=-\frac{1}{4 \, {\left(\theta_{0}^{3} + 3 \, \theta_{0}^{2} + 2 \, \theta_{0}\right)}}\bigg\{4(2\theta_0+2)\bar{\alpha_5}\}\\
	    &=-\frac{2(\theta_0+1)}{ \, {\left(\theta_{0}^{3} + 3 \, \theta_{0}^{2} + 2 \, \theta_{0}\right)}}\bar{\alpha_5}\\
	    &=-\frac{2(\theta_0+1)}{\theta_0(\theta_0+1)(\theta_0+2)}\bar{\alpha_5}\\
	    &=-\frac{2}{\theta_0(\theta_0+2)}\bar{\alpha_5}
	    \end{align*}
	    \normalsize	
	    Then, we have
	    \begin{align*}
	    	\theta_0^3+6\theta_0^4+8\theta_0=\theta_0(\theta_0+2)(\theta_0+4),
	    \end{align*}
	    so
	    \begin{align*}
	    	\frac{\mu_{1} \theta_{0}^{3} + 6 \, \mu_{1} \theta_{0}^{2} + 2 \, {\left(4 \, \mu_{1} + \overline{\alpha_{5}}\right)} \theta_{0} + 8 \, \overline{\alpha_{5}}}{\theta_{0}^{2} + 2 \, \theta_{0}}&=\frac{1}{\theta_0(\theta_0+2)}\bigg\{\left(\theta_0^3+6\theta_0^2+8\right)\mu_1+2(\theta_0+4)\bar{\alpha_5}\bigg\}\\
	    	&=\frac{1}{\theta_0(\theta_0+2)}\bigg\{\theta_0(\theta_0+2)(\theta_0+4)\left(-\frac{2}{\theta_0(\theta_0+2)}\bar{\alpha_5}\right)+2(\theta_0+4)\bar{\alpha_5} \bigg\}\\
	    	&=0
	    \end{align*}
	    and in \eqref{deftheta0+3}, we have
	    \begin{align}\label{2=0}
	    	\textbf{(2)}=0
	    \end{align}	   
	    Therefore, we obtain by \eqref{deftheta0+3}
	    \small
	    \begin{align}\label{endendcancel}
	    	&\frac{3 \, \overline{\zeta_{0}}}{4 \, {\left(\theta_{0}^{3} + 4 \, \theta_{0}^{2} + 5 \, \theta_{0} + 2\right)}} \overline{B_{1}}+\frac{\theta_{0}^{3} \overline{\mu_{14}} + 3 \, \theta_{0}^{2} \overline{\mu_{14}} + 2 \, \theta_{0} \overline{\mu_{14}} - 2 \, {\left(\theta_{0} \overline{\alpha_{1}} + 4 \, \overline{\alpha_{1}}\right)} \overline{\zeta_{0}} - 2 \, {\left(\theta_{0}^{5} + 9 \, \theta_{0}^{4} + 28 \, \theta_{0}^{3} + 36 \, \theta_{0}^{2} + 16 \, \theta_{0}\right)} \overline{\zeta_{9}}}{\theta_{0}^{3} + 3 \, \theta_{0}^{2} + 2 \, \theta_{0}} A_{1}\nonumber\\
	    	&+\frac{\theta_{0}^{4} \overline{\mu_{28}} + 6 \, \theta_{0}^{3} \overline{\mu_{28}} - 8 \, {\left(2 \, \theta_{0}^{2} + 11 \, \theta_{0} + 12\right)} {\left| A_{1} \right|}^{2} \overline{\zeta_{1}} + 11 \, \theta_{0}^{2} \overline{\mu_{28}} + 6 \, \theta_{0} \overline{\mu_{28}} - 2 \, {\left(\theta_{0}^{2} \overline{\alpha_{5}} + 7 \, \theta_{0} \overline{\alpha_{5}} + 12 \, \overline{\alpha_{5}}\right)} \overline{\zeta_{0}}}{\theta_{0}^{4} + 6 \, \theta_{0}^{3} + 11 \, \theta_{0}^{2} + 6 \, \theta_{0}}  A_{0}\nonumber\\
	    	&+\lambda_1\bar{A_1}+\lambda_2\bar{A_0}=0
	    \end{align}
	    \normalsize
	    Furthermore, as
	    \begin{align*}
	    	\s{\vec{A}_0}{\vec{A}_0}=\s{\vec{A}_0}{\vec{A}_1}=\s{\vec{A}_0}{\bar{\vec{A}_1}}=0,
	    \end{align*}
	    we deduce that $\bar{\vec{A}_0}$ is linearly independent with $\vec{A}_0,\vec{A}_1$ and $\bar{\vec{A}_1}$, and as we also have
	    \begin{align*}
	    	\vec{B}_1=-2\s{\bar{\vec{A}_1}}{\vec{C}_1}\vec{A}_0\in \mathrm{Span}(\vec{A}_0)
	    \end{align*}
	    the linear relation in $\bar{\vec{A}_0}$ in \eqref{endendcancel} must be trivial, so that
	    \begin{align}\label{endendendcancel}
	    	\frac{3 \, \overline{\zeta_{0}}}{4 \, {\left(\theta_{0}^{3} + 4 \, \theta_{0}^{2} + 5 \, \theta_{0} + 2\right)}} \overline{B_{1}}+\lambda_2\bar{A_0}=0.
	    \end{align}
	    The rest of the proof will amount at computing this coefficient in front of $\bar{\vec{A}_0}$ in \eqref{endendendcancel}, which must cancel as $\vec{A}_0\neq 0$ by the very definition of a branch point of multiplicity $\theta_0\in \N$.
	    
	    We recall that
	    \begin{align*}
	    	\alpha_{15}&=\frac{1}{4} \, {\left| A_{1} \right|}^{2} \overline{A_{0}} \overline{C_{2}} + 2 \, A_{1} \overline{A_{4}} + \frac{1}{8} \, {\left(\overline{A_{0}} \overline{\alpha_{1}} + 2 \, \overline{A_{2}}\right)} \overline{B_{1}} + \frac{1}{6} \, \overline{A_{1}} \overline{B_{3}} + \frac{1}{8} \, \overline{A_{0}} \overline{B_{6}} + \frac{1}{24} \, {\left(8 \, {\left| A_{1} \right|}^{2} \overline{A_{1}} + 3 \, \overline{A_{0}} \overline{\alpha_{5}}\right)} \overline{C_{1}}
	    \end{align*}
	    while
	    \small
	    		\begin{align*}
	    \left\{
	    \begin{alignedat}{1}
	    \vec{B}_1&=-2\s{\bar{\vec{A}_1}}{\vec{C}_1}\vec{A}_0\\
	    \vec{B}_2&=-\frac{(\theta_0+2)}{4\theta_0}|\vec{C}_1|^2\bar{\vec{A}_0}+\left(\bar{\alpha_0}\s{\bar{\vec{A}_1}}{\vec{C}_1}-2\s{\bar{\vec{A}_2}}{\vec{C}_1}\right)\vec{A}_0\\
	    \vec{B}_3&=-2\s{\bar{\vec{A}_1}}{\vec{C}_1}\vec{A}_1+\frac{2}{\theta_0-3}\s{\vec{A}_1}{\vec{C}_1}\bar{\vec{A}_1}+\frac{2\bar{\alpha_0}}{\theta_0-3}\s{\vec{A}_1}{\vec{C}_1}\bar{\vec{A}_0}+2\left(\alpha_0\s{\bar{\vec{A}_1}}{\vec{C}_1}-\s{\bar{\vec{A}_1}}{\vec{C}_2}\right)\vec{A}_0\\
	    \vec{B}_4&=-\frac{(\theta_0+3)}{6\theta_0}\s{\vec{C}_1}{\bar{\vec{C}_2}}\bar{\vec{A}_0}-\frac{\bar{\zeta_2}}{6\theta_0}\vec{C}_1-\frac{(\theta_0+3)}{6\theta_0}|\vec{C}_1|^2\bar{\vec{A}_1}-\frac{\theta_0(\theta_0+1)}{6}\alpha_2\bar{\vec{C}_1}+2\left(\bar{\alpha_1}\s{\bar{\vec{A}_1}}{\vec{C}_1}-\s{\bar{\vec{A}_3}}{\vec{C}_1}\right)\vec{A}_0.\\
	    \vec{B}_5&=-\left(\frac{(\theta_0+2)}{4}\s{\bar{\vec{C}_1}}{\vec{C}_2}+\frac{2}{\theta_0-3}\bar{\alpha_1}\zeta_2\right)\bar{\vec{A}_0}+\frac{2\,\zeta_2}{\theta_0-3}\bar{\vec{A}_2}-2\s{\bar{\vec{A}_2}}{\vec{C}_1}\vec{A}_1\\
	    &+\left(8\theta_0(\theta_0+1)|\vec{A}_1|^2\alpha_2-2\s{\bar{\vec{A}_2}}{\vec{C}_2}-\frac{2}{\theta_0-3}\bar{\zeta_0}\zeta_2\right)\vec{A}_0\\
	    \vec{B}_6&=\frac{2\,\zeta_5}{\theta_0-4}\bar{\vec{A}_1}-\frac{4}{\theta_0-3}|\vec{A}_1|^2\zeta_2\bar{\vec{A}_0}+\left(-2\s{\bar{\vec{A}_1}}{\vec{C}_3}+4\theta_0(\theta_0+1)\alpha_1\alpha_2\right)\vec{A}_0-4\theta_0(\theta_0+1)\alpha_2\vec{A}_2-2\s{\bar{\vec{A}_1}}{C_2}\vec{A}_1\\
	    \vec{E}_1&=-\frac{1}{2\theta_0}\s{\vec{C}_1}{\vec{C}_1}\bar{\vec{A}_0}\\
	    \vec{E}_2&=-\alpha_2\,\vec{C}_1-\frac{2(2\theta_0+1)(\theta_0-4)}{\theta_0+1}\alpha_7\bar{\vec{A}_1}\\
	    \vec{E}_3&=-\frac{1}{\theta_0}\s{\vec{C}_1}{\vec{C}_2}\bar{\vec{A}_0}
	    \end{alignedat}\right.
	    \end{align*}
	    \normalsize
	    From now on, we assume that
	    \begin{align}\label{holomorphe}
	    	|\vec{A}_1|^2\s{\vec{A}_1}{\vec{C}_1}=\s{\bar{\vec{A}_1}}{\vec{C}_1}\s{\vec{A}_1}{\vec{A}_1}
	    \end{align}
	    thanks of the meromorphy of the quartic form (see \eqref{omega1}).
	    As
	    \begin{align}\label{h0l1}
	    \left\{
	    \begin{alignedat}{1}
	    	&\alpha_2=\frac{1}{2\theta_0(\theta_0+1)}\s{\bar{\vec{A}_1}}{\vec{C}_1},\\
	    	&\zeta_0=\s{\vec{A}_1}{\vec{A}_1}\\
	    	&\zeta_2=\s{\vec{A}_1}{\vec{C}_1}
	    	\end{alignedat}\right.
	    \end{align}
	    this also implies that
	    \begin{align}\label{hol2}
	    	|\vec{A}_1|^2\s{\vec{A}_1}{\vec{C}_1}=|\vec{A}_1|^2\zeta_2=2\theta_0(\theta_0+1)\zeta_0\alpha_2=2\theta_0(\theta_0+1)\s{\vec{A}_1}{\vec{A}_1}\alpha_2=\s{\vec{A}_1}{\vec{A}_1}\s{\bar{\vec{A}_1}}{\vec{C}_1}
	    \end{align}
	    In particular, we have as $\s{\vec{A}_1}{\vec{A}_1}+2\s{\vec{A}_0}{\vec{A}_2}=0$
	    \begin{align}\label{step0}
	    	\s{\vec{A}_2}{\vec{B}_1}&=\s{\vec{A}_2}{-2\s{\bar{\vec{A}_1}}{\vec{C}_1}\vec{A}_0}=-2\s{\bar{\vec{A}_1}}{\vec{C}_1}\s{\vec{A}_0}{\vec{A}_2}=\s{\bar{\vec{A}_1}}{\vec{C}_1}\s{\vec{A}_1}{\vec{A}_1}
	    	=|\vec{A}_1|^2\zeta_2.
	    \end{align}
	    Then, we have
	    \begin{align}\label{step1}
	    	\s{\vec{A}_1}{\vec{B}_3}&=-2\s{\bar{\vec{A}_1}}{\vec{C}_1}\s{\vec{A}_1}{\vec{A}_1}+\frac{2}{\theta_0-3}|\vec{A}_1|^2\s{\vec{A}_1}{\vec{C}_1}\nonumber\\
	    	&=\left(-2+\frac{2}{\theta_0-3}\right)|\vec{A}_1|^2\s{\vec{A}_1}{\vec{C}_1}\nonumber\\
	    	&=-2\frac{(\theta_0-4)}{\theta_0-3}|\vec{A}_1|^2\s{\vec{A}_1}{\vec{C}_1}\nonumber\\
	    	&=-\frac{2(\theta_0-4)}{\theta_0-3}|\vec{A}_1|^2\zeta_2	    \end{align}
	    	as
	    	\begin{align*}
	    		\zeta_2=\s{\vec{A}_1}{\vec{C}_1}.
	    	\end{align*}
	    	Now, we compute as $|\vec{A}_0|^2=\dfrac{1}{2}$ and $\s{\vec{A}_1}{\vec{A}_1}+2\s{\vec{A}_0}{\vec{A}_2}=0$ 
	    	\begin{align}\label{step2}
	    		\s{\vec{A}_0}{\vec{B}_6}&=\bs{\vec{A}_0}{\frac{2\,\zeta_5}{\theta_0-4}\bar{\vec{A}_1}-\frac{4}{\theta_0-3}|\vec{A}_1|^2\zeta_2\bar{\vec{A}_0}+\left(-2\s{\bar{\vec{A}_1}}{\vec{C}_3}+4\theta_0(\theta_0+1)\alpha_1\alpha_2\right)\vec{A}_0-4\theta_0(\theta_0+1)\alpha_2\vec{A}_2-2\s{\bar{\vec{A}_1}}{C_2}\vec{A}_1}\nonumber\\
	    		&=-\frac{2}{\theta_0-3}|\vec{A}_1|^2\zeta_2-4\theta_0(\theta_0+1)\alpha_2\s{\vec{A}_0}{\vec{A}_2}\nonumber\\
	    		&=-\frac{2}{\theta_0-3}|\vec{A}_1|^2\zeta_2+2\theta_0(\theta_0+1)\alpha_2\s{\vec{A}_1}{\vec{A}_1}\nonumber\\
	    		&=\left(-\frac{2}{\theta_0-3}+1\right)|\vec{A}_1|^2\zeta_2\nonumber\\
	    		&=\frac{(\theta_0-5)}{\theta_0-3}|\vec{A}_1|^2\zeta_2.
	    	\end{align}
	    	Finally, we have as $\s{\vec{A}_1}{\vec{C}_1}+\s{\vec{A}_0}{\vec{C}_2}=0$ and thanks of- \eqref{step0}, \eqref{step1} and \eqref{step2}
	    	\begin{align}\label{stepalpha15}
	    		\alpha_{15}&=\frac{1}{4} \, {\left| A_{1} \right|}^{2} \overline{A_{0}} \overline{C_{2}} + 2 \, A_{1} \overline{A_{4}} + \frac{1}{8} \, {\left(\ccancel{\overline{A_{0}} \overline{\alpha_{1}}} + 2 \, \overline{A_{2}}\right)} \overline{B_{1}} + \frac{1}{6} \, \overline{A_{1}} \overline{B_{3}} + \frac{1}{8} \, \overline{A_{0}} \overline{B_{6}} + \frac{1}{24} \, {\left(8 \, {\left| A_{1} \right|}^{2} \overline{A_{1}} + \ccancel{3 \, \overline{A_{0}} \overline{\alpha_{5}}}\right)} \overline{C_{1}}\nonumber\\
	    		&=-\frac{1}{4}|\vec{A}_1|^2\bar{\zeta_2}+2\s{\vec{A}_1}{\bar{\vec{A}_4}}+\frac{1}{4}\bar{\s{\vec{A}_2}{\vec{B}_1}}+\frac{1}{6}\bar{\s{\vec{A}_1}{\vec{B}_3}}+\frac{1}{8}\bar{\s{\vec{A}_0}{\vec{B}_6}}+\frac{1}{3}|\vec{A}_1|^2\bar{\zeta_2}\nonumber\\
	    		&=2\s{\vec{A}_1}{\bar{\vec{A}_4}}+\left(-\frac{1}{4}+\frac{1}{4}+\frac{1}{6}\times \left(-\frac{2(\theta_0-4)}{\theta_0-3}\right)+\frac{1}{8}\times\left(\frac{(\theta_0-5)}{\theta_0-3}\right)+\frac{1}{3}\right)|\vec{A}_1|^2\bar{\zeta_2}\nonumber\\
	    		&=2\s{\vec{A}_1}{\bar{\vec{A}_4}}+\frac{3\theta_0-7}{24(\theta_0-3)}|\vec{A}_1|^2\bar{\zeta_2}.
	    	\end{align}
	    	Then, we have
	    	\begin{align*}
	    		\zeta_{13}&=\frac{1}{\theta_{0}^{4} - 7 \, \theta_{0}^{2} - 6 \, \theta_{0}}\bigg\{{\left(\theta_{0}^{3} - 7 \, \theta_{0} - 6\right)} B_{7} \overline{A_{0}} + 2 \, {\left(\theta_{0}^{5} \overline{\alpha_{2}} - 7 \, \theta_{0}^{3} \overline{\alpha_{2}} - 6 \, \theta_{0}^{2} \overline{\alpha_{2}}\right)} \overline{A_{1}}^{2}\\
	    		& + 2 \, {\left(\theta_{0}^{5} \overline{\alpha_{2}} - \theta_{0}^{4} \overline{\alpha_{2}} - 5 \, \theta_{0}^{3} \overline{\alpha_{2}} - 3 \, \theta_{0}^{2} \overline{\alpha_{2}}\right)} \overline{A_{0}} \overline{A_{2}}
	    		 + {\left({\left(\theta_{0}^{2} + 2 \, \theta_{0}\right)} A_{1} \overline{\zeta_{2}} + {\left(\theta_{0}^{3} - \theta_{0}^{2} - 6 \, \theta_{0}\right)} \overline{A_{0}} \overline{\zeta_{3}}\right)} \overline{A_{1}}\bigg\}\\
	    	\end{align*}
	    	and
	    	\begin{align}
	    		\zeta_{16}&=\frac{4 \, {\left(\theta_{0}^{2} + 2 \, \theta_{0}\right)} A_{1} \zeta_{2} \overline{A_{1}} - {\left(\theta_{0}^{3} - \theta_{0}^{2} - 2 \, \theta_{0}\right)} A_{2} B_{1} + 2 \, {\left(\theta_{0}^{2} + 3 \, \theta_{0} + 2\right)} A_{0} B_{12} - {\left(\theta_{0}^{3} - \theta_{0}^{2} - 6 \, \theta_{0}\right)} A_{1} B_{3}}{2 \, {\left(\theta_{0}^{3} + 3 \, \theta_{0}^{2} + 2 \, \theta_{0}\right)}} 
	    	\end{align}
	    		\begin{align}\label{defzeta05ter}
	    		\left\{
	    		\begin{alignedat}{1}
	    		\zeta_0&=\s{\vec{A}_1}{\vec{A}_1}\\
	    		\zeta_1&=\s{\vec{A}_1}{\vec{A}_2}\\
	    		\zeta_2&=\s{\vec{A}_1}{\vec{C}_1}\\
	    		\zeta_3&=\s{\bar{\vec{A}_1}}{\vec{C}_2}\\
	    		\zeta_4&=\s{\bar{\vec{A}_2}}{\vec{C}_1}\\
	    		\zeta_5&=2\s{\vec{A}_2}{\vec{C}_1}+\s{\vec{A}_1}{\vec{C}_2}
	    		\end{alignedat}\right.
	    		\end{align}
	    		Then, we have
	    		\begin{align*}
	    			\vec{B}_7=\begin{dmatrix}
	    			\frac{1}{2} & \overline{B_{6}} & 0 & -\theta_{0} + 4  & \textbf{(6)}\\
	    			2 \, A_{1} \overline{C_{1}} & \overline{A_{2}} & 0 & -\theta_{0} + 4  & \textbf{(7)}\\
	    			-\ccancel{4 \, A_{0} {\left| A_{1} \right|}^{2} \overline{C_{2}}} - 2 \, A_{1} \overline{C_{1}} \overline{\alpha_{1}} - \ccancel{2 \, A_{0} \overline{C_{1}} \overline{\alpha_{5}}} + \frac{1}{4} \, \ccancel{\overline{B_{1}} \overline{C_{1}}} + 2 \, A_{1} \overline{C_{3}} & \overline{A_{0}} & 0 & -\theta_{0} + 4  & \textbf{(8)}\\
	    			-\ccancel{4 \, A_{0} {\left| A_{1} \right|}^{2} \overline{C_{1}}} + 2 \, A_{1} \overline{C_{2}} & \overline{A_{1}} & 0 & -\theta_{0} + 4  & \textbf{(9)}
	    			\end{dmatrix}
	    		\end{align*}
	    		so
	    		\begin{align*}
	    			\s{\bar{\vec{A}_0}}{\vec{B}_7}&=\frac{1}{2}\bar{\s{\vec{A}_0}{\vec{B}_6}}+2\s{\vec{A}_1}{\bar{\vec{C}_1}}\bar{\s{\vec{A}_0}{\vec{A}_2}}\\
	    			&=\frac{1}{2}\bar{\s{\vec{A}_0}{\vec{B}_6}}-\s{\vec{A}_1}{\bar{\vec{C}_1}}\bar{\s{\vec{A}_1}{\vec{A}_1}}\\
	    			&=\frac{1}{2}\bar{\s{\vec{A}_0}{\vec{B}_6}}-|\vec{A}_1|^2\bar{\zeta_2}.
	    		\end{align*}
	    		Then, we obtain by \eqref{step2}
	    		\begin{align*}
	    			\s{\vec{A}_0}{\vec{B}_6}=\frac{(\theta_0-5)}{\theta_0-3}|\vec{A}_1|^2\zeta_2
	    		\end{align*}
	    		so
	    		\begin{align}\label{step3}
	    			\s{\bar{\vec{A}_0}}{\vec{B}_7}&=\frac{(\theta_0-5)}{2(\theta_0-3)}|\vec{A}_1|^2\bar{\zeta_2}-|\vec{A}_1|^2\bar{\zeta_2}\nonumber\\
	    			&=-\frac{(\theta_0-1)}{2(\theta_0-3)}|\vec{A}_1|^2\bar{\zeta_2}
	    		\end{align}
	    		Finally, we obtain as $2\s{\vec{A}_0}{\vec{A}_2}=-\s{\vec{A}_1}{\vec{A}_1}$
	    		\begin{align}\label{step4}
	    			\zeta_{13}&=\frac{1}{\theta_{0}^{4} - 7 \, \theta_{0}^{2} - 6 \, \theta_{0}}\bigg\{{\left(\theta_{0}^{3} - 7 \, \theta_{0} - 6\right)} B_{7} \overline{A_{0}} + 2 \, {\left(\theta_{0}^{5} \overline{\alpha_{2}} - 7 \, \theta_{0}^{3} \overline{\alpha_{2}} - 6 \, \theta_{0}^{2} \overline{\alpha_{2}}\right)} \overline{A_{1}}^{2}\nonumber\\
	    			& + 2 \, {\left(\theta_{0}^{5} \overline{\alpha_{2}} - \theta_{0}^{4} \overline{\alpha_{2}} - 5 \, \theta_{0}^{3} \overline{\alpha_{2}} - 3 \, \theta_{0}^{2} \overline{\alpha_{2}}\right)} \overline{A_{0}} \overline{A_{2}}
	    			+ {\left({\left(\theta_{0}^{2} + 2 \, \theta_{0}\right)} A_{1} \overline{\zeta_{2}} + \ccancel{{\left(\theta_{0}^{3} - \theta_{0}^{2} - 6 \, \theta_{0}\right)} \overline{A_{0}} \overline{\zeta_{3}}}\right)} \overline{A_{1}}\bigg\}\nonumber\\
	    			&=\frac{1}{\theta_0(\theta_0+1)(\theta_0+2)(\theta_0-3)}\bigg\{(\theta_0+1)(\theta_0+2)(\theta_0-3)\left(-\frac{(\theta_0-1)}{2(\theta_0-3)}|\vec{A}_1|^2\bar{\zeta_2}\right)\nonumber\\
	    			&+\Big(2(\theta_0^5-7\theta_0^3-6\theta_0^2)-(\theta_0^5-\theta_0^4-5\theta_0^3-3\theta_0^2)\Big)\bar{\zeta_0}\bar{\alpha_2}+\theta_0(\theta_0+2)|\vec{A}_1|^2\bar{\zeta_2} \bigg\}\nonumber\\
	    			&=\frac{1}{\theta_0(\theta_0+1)(\theta_0+2)(\theta_0-3)}\bigg\{-\frac{1}{2}(\theta_0-1)(\theta_0+1)(\theta_0+2)|\vec{A}_1|^2\bar{\zeta_2}+\theta_0^2(\theta_0+1)(\theta_0+3)(\theta_0-3)\bar{\zeta_0}\bar{\alpha_2}\nonumber\\
	    			& +\theta_0(\theta_0+2)|\vec{A}_1|^2\bar{\zeta_2}\bigg\}\nonumber\\
	    			&=\frac{1}{\theta_0(\theta_0+1)(\theta_0+2)(\theta_0-3)}\bigg\{-\frac{1}{2}(\theta_0-1)(\theta_0+1)(\theta_0+2)+\frac{1}{2}\theta_0(\theta_0+3)(\theta_0-3)+\theta_0(\theta_0+2)
	    			\bigg\}|\vec{A}_1|^2\bar{\zeta_2}\nonumber\\
	    			&=\frac{-(2\theta_0-1)}{\theta_0(\theta_0+1)(\theta_0+2)(\theta_0-3)}|\vec{A}_1|^2\bar{\zeta_2}
	    		\end{align}
	    		as
	    		\begin{align*}
	    			&2(\theta_0^5-7\theta_0^3-6\theta_0^2)-(\theta_0^5-\theta_0^4-5\theta_0^3-3\theta_0^2)=\theta_0^2(\theta_0+1)(\theta_0+3)(\theta_0-3)\\
	    			&\theta_0(\theta_0+1)\bar{\zeta_0}\bar{\alpha_2}=\frac{1}{2}\bar{\s{\vec{A}_1}{\vec{A}_1}}\s{\vec{A}_1}{\bar{\vec{C}_1}}=\frac{1}{2}|\vec{A}_1|^2\bar{\zeta_2}\\
	    			&-\frac{1}{2}(\theta_0-1)(\theta_0+1)(\theta_0+2)+\frac{1}{2}\theta_0(\theta_0+3)(\theta_0-3)+\theta_0(\theta_0+2)=-2\theta_0+1
	    		\end{align*}
	    		Now we recall that
	    		\begin{align*}
	    			\vec{B}_{12}=\begin{dmatrix}
	    			-\frac{1}{2} \, \theta_{0} + 2 & B_{6} & -\theta_{0} + 3 & 1  & \textbf{(49)}\\
	    			\frac{1}{2} \, \ccancel{B_{1} C_{1}} & A_{0} & -\theta_{0} + 3 & 1  & \textbf{(50)}\\
	    			-\ccancel{4 \, A_{0} C_{1} \alpha_{1}} + 4 \, A_{2} C_{1} + 2 \, A_{1} C_{2} & \overline{A_{1}} & -\theta_{0} + 3 & 1  & \textbf{(51)}\\
	    			-8 \, A_{1} C_{1} {\left| A_{1} \right|}^{2} - 4 \, A_{0} C_{2} {\left| A_{1} \right|}^{2} - \ccancel{4 \, A_{0} B_{1} \alpha_{1}} - \ccancel{4 \, A_{0} C_{1} \alpha_{5}} + 4 \, A_{2} B_{1} + \ccancel{2 \, A_{1} B_{3}} & \overline{A_{0}} & -\theta_{0} + 3 & 1  & \textbf{(52)}
	    			\end{dmatrix}
	    		\end{align*}
	    		Furthermore, remark that $\s{\vec{A}_2}{-2\vec{A}_0}=\s{\vec{A}_1}{\vec{A}_1}$ so by \eqref{holomorphe}
	    		\begin{align*}
	    			&-8 \, A_{1} C_{1} {\left| A_{1} \right|}^{2} - 4 \, A_{0} C_{2} {\left| A_{1} \right|}^{2} - \ccancel{4 \, A_{0} B_{1} \alpha_{1}} - \ccancel{4 \, A_{0} C_{1} \alpha_{5}} + 4 \, A_{2} B_{1} + \ccancel{2 \, A_{1} B_{3}}\\
	    			&=-4|\vec{A}_1|^2\s{\vec{A}_1}{\vec{C}_1}+4\s{\vec{A}_2}{-2\s{\bar{\vec{A}_1}}{\vec{C}_1}\vec{A}_0}\\
	    			&=-4|\vec{A}_1|^2\s{\vec{A}_1}{\vec{C}_1}+4\s{\bar{\vec{A}_1}}{\vec{C}_1}\s{\vec{A}_1}{\vec{A}_1}\\
	    			&=0
	    		\end{align*}
	    		so
	    		\begin{align}\label{52}
	    			(\textbf{52})=0
	    		\end{align}
	    		and by \eqref{step2}, we obtain
	    		\begin{align}\label{step5}
	    			\s{\vec{A}_0}{\vec{B}_{12}}=-\frac{(\theta_0-4)}{2}\s{\vec{A}_0}{\vec{B}_6}=-\frac{(\theta_0-4)}{2}\left(\frac{(\theta_0-5)}{\theta_0-3}|\vec{A}_1|^2\zeta_2\right)=-\frac{(\theta_0-4)(\theta_0-5)}{2(\theta_0-3)}|\vec{A}_1|^2\zeta_2
	    		\end{align}
	    		We deduce by \eqref{step0}, \eqref{step1}, \eqref{step2} and \eqref{step5} that
	    		\begin{align}\label{lambda16}
	    			\zeta_{16}&=\frac{4 \, {\left(\theta_{0}^{2} + 2 \, \theta_{0}\right)} A_{1} \zeta_{2} \overline{A_{1}} - {\left(\theta_{0}^{3} - \theta_{0}^{2} - 2 \, \theta_{0}\right)} A_{2} B_{1} + 2 \, {\left(\theta_{0}^{2} + 3 \, \theta_{0} + 2\right)} A_{0} B_{12} - {\left(\theta_{0}^{3} - \theta_{0}^{2} - 6 \, \theta_{0}\right)} A_{1} B_{3}}{2 \, {\left(\theta_{0}^{3} + 3 \, \theta_{0}^{2} + 2 \, \theta_{0}\right)}}\nonumber\\
	    			&=\frac{1}{2\theta_0(\theta_0+1)(\theta_0+2)}\bigg\{4\theta_0(\theta_0+2)-\theta_0(\theta_0+1)(\theta_0-2)\times 1+2(\theta_0+1)(\theta_0+2)\times\left(-\frac{(\theta_0-4)(\theta_0-5)}{2(\theta_0-3)}\right)\nonumber\\
	    			& -\theta_0(\theta_0+2)(\theta_0-3)\times\left(-\frac{2(\theta_0-4)}{\theta_0-3}\right)\bigg\}|\vec{A}_1|^2\zeta_2\nonumber\\
	    			&=\frac{2(\theta_0^3-\theta_0^2-6\theta_0-10)}{\theta_0(\theta_0+1)(\theta_0+2)(\theta_0-3)}|\vec{A}_1|^2\zeta_2
	    		\end{align}
	    		as 
	    		\begin{align*}
	    			&4\theta_0(\theta_0+2)-\theta_0(\theta_0+1)(\theta_0-2)\times 1+2(\theta_0+1)(\theta_0+2)\times\left(-\frac{(\theta_0-4)(\theta_0-5)}{2(\theta_0-3)}\right)\\
	    			& -\theta_0(\theta_0+2)(\theta_0-3)\times\left(-\frac{2(\theta_0-4)}{\theta_0-3}\right)=4\theta_0(\theta_0+2)-\frac{40}{\theta_0-3}=\frac{4\left(\theta_0^3-\theta_0^2-6\theta_0-10\right)}{\theta_0-3}
	    		\end{align*}
	    		
	    		\small 
	    		\begin{align}\label{partmu15}
	    			\mu_{15}&=
	    			-\frac{1}{24 \, {\left(\theta_{0}^{5} + 10 \, \theta_{0}^{4} + 35 \, \theta_{0}^{3} + 50 \, \theta_{0}^{2} + 24 \, \theta_{0}\right)}}  \bigg\{96 \, {\left(\theta_{0}^{4} + 7 \, \theta_{0}^{3} + 14 \, \theta_{0}^{2} + 8 \, \theta_{0}\right)} A_{0} {\left| A_{1} \right|}^{2} \overline{A_{3}}\nonumber\\
	    			& + 2 \, {\left(\theta_{0}^{4} + 10 \, \theta_{0}^{3} + 35 \, \theta_{0}^{2} + 50 \, \theta_{0} + 24\right)} {\left| A_{1} \right|}^{2} \overline{A_{0}} \overline{C_{2}} - 192 \, \bigg\{{\left(\theta_{0}^{4} \overline{\alpha_{3}} + 10 \, \theta_{0}^{3} \overline{\alpha_{3}} + 35 \, \theta_{0}^{2} \overline{\alpha_{3}} + 50 \, \theta_{0} \overline{\alpha_{3}} + 24 \, \overline{\alpha_{3}}\right)} A_{0} \overline{A_{0}}\nonumber\\
	    			& + \ccancel{{\left(\theta_{0}^{4} \overline{\alpha_{1}} + 9 \, \theta_{0}^{3} \overline{\alpha_{1}} + 26 \, \theta_{0}^{2} \overline{\alpha_{1}} + 24 \, \theta_{0} \overline{\alpha_{1}}\right)} A_{0} \overline{A_{1}}}\bigg\} {\left| A_{1} \right|}^{2} - 48 \, \bigg\{{\left(2 \, \overline{\alpha_{1}} \overline{\alpha_{5}} - \alpha_{15}\right)} \theta_{0}^{4} + 10 \, {\left(2 \, \overline{\alpha_{1}} \overline{\alpha_{5}} - \alpha_{15}\right)} \theta_{0}^{3} + 35 \, {\left(2 \, \overline{\alpha_{1}} \overline{\alpha_{5}} - \alpha_{15}\right)} \theta_{0}^{2} \nonumber\\
	    			&+ 50 \, {\left(2 \, \overline{\alpha_{1}} \overline{\alpha_{5}} - \alpha_{15}\right)} \theta_{0} + 48 \, \overline{\alpha_{1}} \overline{\alpha_{5}} - 24 \, \alpha_{15}\bigg\} A_{0} \overline{A_{0}} - 48 \, \bigg\{{\left(\overline{\alpha_{1}}^{2} - \overline{\alpha_{4}}\right)} \theta_{0}^{4} + 10 \, {\left(\overline{\alpha_{1}}^{2} - \overline{\alpha_{4}}\right)} \theta_{0}^{3} + 35 \, {\left(\overline{\alpha_{1}}^{2} - \overline{\alpha_{4}}\right)} \theta_{0}^{2}\nonumber\\
	    			& + 50 \, {\left(\overline{\alpha_{1}}^{2} - \overline{\alpha_{4}}\right)} \theta_{0} + 24 \, \overline{\alpha_{1}}^{2} - 24 \, \overline{\alpha_{4}}\bigg\} \ccancel{A_{1} \overline{A_{0}}} - 48 \, {\left(\theta_{0}^{4} + 6 \, \theta_{0}^{3} + 11 \, \theta_{0}^{2} + 6 \, \theta_{0}\right)} A_{1} \overline{A_{4}} - 4 \, {\left(\theta_{0}^{4} + 9 \, \theta_{0}^{3} + 26 \, \theta_{0}^{2} + 24 \, \theta_{0}\right)} \overline{A_{1}} \overline{B_{3}}\\
	    			& - 3 \, {\left(\theta_{0}^{4} + 10 \, \theta_{0}^{3} + 35 \, \theta_{0}^{2} + 50 \, \theta_{0} + 24\right)} \overline{A_{0}} \overline{B_{6}}\nonumber\\
	    			& + 48 \, {\left(\ccancel{{\left(\theta_{0}^{4} \overline{\alpha_{6}} + 9 \, \theta_{0}^{3} \overline{\alpha_{6}} + 26 \, \theta_{0}^{2} \overline{\alpha_{6}} + 24 \, \theta_{0} \overline{\alpha_{6}}\right)} A_{0}} + {\left(\theta_{0}^{4} \overline{\alpha_{3}} + 9 \, \theta_{0}^{3} \overline{\alpha_{3}} + 26 \, \theta_{0}^{2} \overline{\alpha_{3}} + 24 \, \theta_{0} \overline{\alpha_{3}}\right)} A_{1}\right)} \overline{A_{1}}\nonumber\\
	    			& + 48 \, {\left({\left(\theta_{0}^{4} \overline{\alpha_{5}} + 8 \, \theta_{0}^{3} \overline{\alpha_{5}} + 19 \, \theta_{0}^{2} \overline{\alpha_{5}} + 12 \, \theta_{0} \overline{\alpha_{5}}\right)} A_{0} + {\left(\theta_{0}^{4} \overline{\alpha_{1}} + 8 \, \theta_{0}^{3} \overline{\alpha_{1}} + 19 \, \theta_{0}^{2} \overline{\alpha_{1}} + 12 \, \theta_{0} \overline{\alpha_{1}}\right)} A_{1}\right)} \overline{A_{2}}\nonumber\\
	    			& + 3 \, {\left(\ccancel{{\left(\theta_{0}^{4} \overline{\alpha_{1}} + 10 \, \theta_{0}^{3} \overline{\alpha_{1}} + 35 \, \theta_{0}^{2} \overline{\alpha_{1}} + 50 \, \theta_{0} \overline{\alpha_{1}} + 24 \, \overline{\alpha_{1}}\right)} \overline{A_{0}}} - 2 \, {\left(\theta_{0}^{4} + 8 \, \theta_{0}^{3} + 19 \, \theta_{0}^{2} + 12 \, \theta_{0}\right)} \overline{A_{2}}\right)} \overline{B_{1}}\nonumber\\
	    			& + {\left(4 \, {\left(\theta_{0}^{4} + 9 \, \theta_{0}^{3} + 26 \, \theta_{0}^{2} + 24 \, \theta_{0}\right)} {\left| A_{1} \right|}^{2} \overline{A_{1}} + \ccancel{3 \, {\left(\theta_{0}^{4} \overline{\alpha_{5}} + 10 \, \theta_{0}^{3} \overline{\alpha_{5}} + 35 \, \theta_{0}^{2} \overline{\alpha_{5}} + 50 \, \theta_{0} \overline{\alpha_{5}} + 24 \, \overline{\alpha_{5}}\right)} \overline{A_{0}}}\right)} \overline{C_{1}}\bigg\}\nonumber\\
	    			&=-\frac{1}{24 \, {\left(\theta_{0}^{5} + 10 \, \theta_{0}^{4} + 35 \, \theta_{0}^{3} + 50 \, \theta_{0}^{2} + 24 \, \theta_{0}\right)}}  \bigg\{96 \, {\left(\theta_{0}^{4} + 7 \, \theta_{0}^{3} + 14 \, \theta_{0}^{2} + 8 \, \theta_{0}\right)} A_{0} {\left| A_{1} \right|}^{2} \overline{A_{3}}\nonumber\\
	    			& + 2 \, {\left(\theta_{0}^{4} + 10 \, \theta_{0}^{3} + 35 \, \theta_{0}^{2} + 50 \, \theta_{0} + 24\right)} {\left| A_{1} \right|}^{2} \overline{A_{0}} \overline{C_{2}} - 192 \, \bigg\{{\left(\theta_{0}^{4} \overline{\alpha_{3}} + 10 \, \theta_{0}^{3} \overline{\alpha_{3}} + 35 \, \theta_{0}^{2} \overline{\alpha_{3}} + 50 \, \theta_{0} \overline{\alpha_{3}} + 24 \, \overline{\alpha_{3}}\right)} A_{0} \overline{A_{0}}\bigg\} {\left| A_{1} \right|}^{2}\nonumber\\
	    			& - 48 \, \bigg\{{\left(2 \, \overline{\alpha_{1}} \overline{\alpha_{5}} - \alpha_{15}\right)} \theta_{0}^{4} + 10 \, {\left(2 \, \overline{\alpha_{1}} \overline{\alpha_{5}} - \alpha_{15}\right)} \theta_{0}^{3} + 35 \, {\left(2 \, \overline{\alpha_{1}} \overline{\alpha_{5}} - \alpha_{15}\right)} \theta_{0}^{2}
	    			+ 50 \, {\left(2 \, \overline{\alpha_{1}} \overline{\alpha_{5}} - \alpha_{15}\right)} \theta_{0} + 48 \, \overline{\alpha_{1}} \overline{\alpha_{5}} - 24 \, \alpha_{15}\bigg\} A_{0} \overline{A_{0}}\nonumber\\
	    			& - 48 \, {\left(\theta_{0}^{4} + 6 \, \theta_{0}^{3} + 11 \, \theta_{0}^{2} + 6 \, \theta_{0}\right)} A_{1} \overline{A_{4}} - 4 \, {\left(\theta_{0}^{4} + 9 \, \theta_{0}^{3} + 26 \, \theta_{0}^{2} + 24 \, \theta_{0}\right)} \overline{A_{1}} \overline{B_{3}}\nonumber\\
	    			& - 3 \, {\left(\theta_{0}^{4} + 10 \, \theta_{0}^{3} + 35 \, \theta_{0}^{2} + 50 \, \theta_{0} + 24\right)} \overline{A_{0}} \overline{B_{6}}\nonumber\\
	    			& + 48  {\left(\theta_{0}^{4} \overline{\alpha_{3}} + 9 \, \theta_{0}^{3} \overline{\alpha_{3}} + 26 \, \theta_{0}^{2} \overline{\alpha_{3}} + 24 \, \theta_{0} \overline{\alpha_{3}}\right)} A_{1} \overline{A_{1}}\nonumber\\
	    			& + 48 \, {\left({\left(\theta_{0}^{4} \overline{\alpha_{5}} + 8 \, \theta_{0}^{3} \overline{\alpha_{5}} + 19 \, \theta_{0}^{2} \overline{\alpha_{5}} + 12 \, \theta_{0} \overline{\alpha_{5}}\right)} A_{0} + {\left(\theta_{0}^{4} \overline{\alpha_{1}} + 8 \, \theta_{0}^{3} \overline{\alpha_{1}} + 19 \, \theta_{0}^{2} \overline{\alpha_{1}} + 12 \, \theta_{0} \overline{\alpha_{1}}\right)} A_{1}\right)} \overline{A_{2}}\nonumber\\
	    			& -6 \, {\left(\theta_{0}^{4} + 8 \, \theta_{0}^{3} + 19 \, \theta_{0}^{2} + 12 \, \theta_{0}\right)} \overline{A_{2}} \overline{B_{1}}\nonumber\\
	    			& + 4 \, {\left(\theta_{0}^{4} + 9 \, \theta_{0}^{3} + 26 \, \theta_{0}^{2} + 24 \, \theta_{0}\right)} {\left| A_{1} \right|}^{2} \overline{A_{1}} \overline{C_{1}}\bigg\}     
	    		\end{align}
	    		\normalsize
	    		Now, recall that
	    		\begin{align*}
	    		\left\{
	    		\begin{alignedat}{1}
	    	        \alpha_1&=2\s{\bar{\vec{A}_0}}{\vec{A}_2}\\
	    			\alpha_3&=2\s{\bar{\vec{A}_0}}{\vec{A}_3}+\frac{1}{12}\s{\vec{A}_1}{\vec{C}_1}\\
	    			\alpha_5&=2\s{\bar{\vec{A}_1}}{\vec{A}_2}\\
	    			\alpha_{15}&=2\s{\vec{A}_1}{\bar{\vec{A}_4}}+\frac{3\theta_0-7}{24(\theta_0-3)}|\vec{A}_1|^2\bar{\zeta_2}.
	    			\end{alignedat}\right.
	    		\end{align*}
	    		so we have
	    		\begin{align*}
	    			&\s{\bar{\vec{A}_0}}{\vec{A}_3}=\frac{1}{2}\alpha_3-\frac{1}{24}\zeta_2\\
	    			&\s{\vec{A}_1}{\bar{\vec{A}_4}}=\frac{1}{2}\alpha_{15}-\frac{3\theta_0-7}{48(\theta_0-3)}|\vec{A}_1|^2\bar{\zeta_2}.
	    		\end{align*}
	    		We will compute separately each term in \eqref{partmu15}. We first have
	    		\begin{align}\label{mu151}
	    			\mathrm{(I)}&=96 \, {\left(\theta_{0}^{4} + 7 \, \theta_{0}^{3} + 14 \, \theta_{0}^{2} + 8 \, \theta_{0}\right)} A_{0} {\left| A_{1} \right|}^{2} \overline{A_{3}}=96|\vec{A}_1|^2{\left(\theta_{0}^{4} + 7 \, \theta_{0}^{3} + 14 \, \theta_{0}^{2} + 8 \, \theta_{0}\right)}\left(\frac{1}{2}\bar{\alpha_3}-\frac{1}{24}\bar{\zeta_2}\right)\nonumber\\
	    			&=48{\left(\theta_{0}^{4} + 7 \, \theta_{0}^{3} + 14 \, \theta_{0}^{2} + 8 \, \theta_{0}\right)}|\vec{A}_1|^2\bar{\alpha_3}-4{\left(\theta_{0}^{4} + 7 \, \theta_{0}^{3} + 14 \, \theta_{0}^{2} + 8 \, \theta_{0}\right)}|\vec{A}_1|^2\bar{\zeta_2}.
	    		\end{align}
	    		Now, as $\s{\vec{A}_1}{\vec{C}_1}+\s{\vec{A}_0}{\vec{C}_2}=0$, we deduce that
	    		\begin{align}\label{mu152}
	    			\mathrm{(II)}&=2 \, {\left(\theta_{0}^{4} + 10 \, \theta_{0}^{3} + 35 \, \theta_{0}^{2} + 50 \, \theta_{0} + 24\right)} {\left| A_{1} \right|}^{2} \overline{A_{0}} \overline{C_{2}}&=-2 \, {\left(\theta_{0}^{4} + 10 \, \theta_{0}^{3} + 35 \, \theta_{0}^{2} + 50 \, \theta_{0} + 24\right)} {\left| A_{1} \right|}^{2} \overline{A_{1}} \overline{C_{1}}\nonumber\\
	    			&=-2 \, {\left(\theta_{0}^{4} + 10 \, \theta_{0}^{3} + 35 \, \theta_{0}^{2} + 50 \, \theta_{0} + 24\right)} {\left| A_{1} \right|}^{2}\bar{\zeta_2}.
	    		\end{align}
	    		Then, we have as $|\vec{A}_0|^2=\dfrac{1}{2}$ the identity
	    		\begin{align}\label{mu153}
	    			\mathrm{(III)}&=- 192 \, \bigg\{{\left(\theta_{0}^{4} \overline{\alpha_{3}} + 10 \, \theta_{0}^{3} \overline{\alpha_{3}} + 35 \, \theta_{0}^{2} \overline{\alpha_{3}} + 50 \, \theta_{0} \overline{\alpha_{3}} + 24 \, \overline{\alpha_{3}}\right)} A_{0} \overline{A_{0}}\bigg\} {\left| A_{1} \right|}^{2}\nonumber\\
	    			&=-96 {\left(\theta_{0}^{4}  + 10 \, \theta_{0}^{3}  + 35 \, \theta_{0}^{2}  + 50 \, \theta_{0}  + 24 \, \right)}|\vec{A}_1|^2\bar{\alpha_3}
	    		\end{align}
	    		Similarly, we have
	    		\begin{align}\label{mu154}
	    			\mathrm{(IV)}&=- 48 \, \bigg\{{\left(2 \, \overline{\alpha_{1}} \overline{\alpha_{5}} - \alpha_{15}\right)} \theta_{0}^{4} + 10 \, {\left(2 \, \overline{\alpha_{1}} \overline{\alpha_{5}} - \alpha_{15}\right)} \theta_{0}^{3} + 35 \, {\left(2 \, \overline{\alpha_{1}} \overline{\alpha_{5}} - \alpha_{15}\right)} \theta_{0}^{2}
	    			+ 50 \, {\left(2 \, \overline{\alpha_{1}} \overline{\alpha_{5}} - \alpha_{15}\right)} \theta_{0}\nonumber\\
	    			& + 48 \, \overline{\alpha_{1}} \overline{\alpha_{5}} - 24 \, \alpha_{15}\bigg\} A_{0} \overline{A_{0}}\nonumber\\
	    			&=-24\left(\theta_0^4+10\,\theta_0^3+35\,\theta_0^2+50\,\theta_0+24\right)(2\bar{\alpha_1\alpha_5}-\alpha_{15})
	    		\end{align}
	    		Then, we have
	    		\begin{align}\label{mu155}
	    			\mathrm{(V)}&=- 48 \, {\left(\theta_{0}^{4} + 6 \, \theta_{0}^{3} + 11 \, \theta_{0}^{2} + 6 \, \theta_{0}\right)} A_{1} \overline{A_{4}}=- 48 \, {\left(\theta_{0}^{4} + 6 \, \theta_{0}^{3} + 11 \, \theta_{0}^{2} + 6 \, \theta_{0}\right)} \left(\frac{1}{2}\alpha_{15}-\frac{3\theta_0-7}{48(\theta_0-3)}|\vec{A}_1|^2\bar{\zeta_2}\right)\nonumber\\
	    			&=- 24 \, {\left(\theta_{0}^{4} + 6 \, \theta_{0}^{3} + 11 \, \theta_{0}^{2} + 6 \, \theta_{0}\right)} \alpha_{15}+\frac{ {\left(\theta_{0}^{4} + 6 \, \theta_{0}^{3} + 11 \, \theta_{0}^{2} + 6 \, \theta_{0}\right)} (3\theta_0-7)}{\theta_0-3}|\vec{A}_1|^2\bar{\zeta_2}
	    		\end{align}
	    		Now, recalling by \eqref{step1} that
	    		\begin{align*}
	    			\s{\vec{A}_1}{\vec{B}_3}=-\frac{2(\theta_0-4)}{\theta_0-3}|\vec{A}_1|^2\bar{\zeta_2},
	    		\end{align*}
	    		we obtain
	    		\begin{align}\label{mu156}
	                \mathrm{(VI)}&=- 4 \, {\left(\theta_{0}^{4} + 9 \, \theta_{0}^{3} + 26 \, \theta_{0}^{2} + 24 \, \theta_{0}\right)} \overline{A_{1}} \overline{B_{3}}=\frac{8{\left(\theta_{0}^{4} + 9 \, \theta_{0}^{3} + 26 \, \theta_{0}^{2} + 24 \, \theta_{0}\right)}(\theta_0-4)}{\theta_0-3}|\vec{A}_1|^2\bar{\zeta_2}.
	    		\end{align}
	    		Now, we have as by \eqref{step2}
	    		\begin{align*}
	    			\s{\vec{A}_0}{\vec{B}_6}=\frac{(\theta_0-5)}{\theta_0-3}|\vec{A}_1|^2\zeta_2,
	    		\end{align*}
	    		the identity
	    		\begin{align}\label{mu157}
	    			 \mathrm{(VII)}&=- 3 \, {\left(\theta_{0}^{4} + 10 \, \theta_{0}^{3} + 35 \, \theta_{0}^{2} + 50 \, \theta_{0} + 24\right)} \overline{A_{0}} \overline{B_{6}}=\frac{- 3 \, {\left(\theta_{0}^{4} + 10 \, \theta_{0}^{3} + 35 \, \theta_{0}^{2} + 50 \, \theta_{0} + 24\right)}(\theta_0-5)}{\theta_0-3}|\vec{A}_1|^2\bar{\zeta_2}
	    		\end{align}
	    		Then, we trivially have
	    		\begin{align}\label{mu158}
	    			 \mathrm{(VIII)}&=48  {\left(\theta_{0}^{4} \overline{\alpha_{3}} + 9 \, \theta_{0}^{3} \overline{\alpha_{3}} + 26 \, \theta_{0}^{2} \overline{\alpha_{3}} + 24 \, \theta_{0} \overline{\alpha_{3}}\right)} A_{1} \overline{A_{1}}=48{\left(\theta_{0}^{4} + 9 \, \theta_{0}^{3} + 26 \, \theta_{0}^{2} + 24 \, \theta_{0} \right)}|\vec{A}_1|^2\overline{\alpha_{3}}.
	    		\end{align}
	    		The next coefficient is as
	    		\begin{align*}
	    			\bar{\alpha_1}\s{\vec{A}_1}{\bar{\vec{A}_2}}=\bar{\alpha_5}\s{\vec{A}_0}{\bar{\vec{A}_2}}=\frac{1}{2}\bar{\alpha_1\alpha_5}
	    		\end{align*}
	    		word
	    		\begin{align}\label{mu159}
	    			\mathrm{(IX)}&=48 \, {\left({\left(\theta_{0}^{4} \overline{\alpha_{5}} + 8 \, \theta_{0}^{3} \overline{\alpha_{5}} + 19 \, \theta_{0}^{2} \overline{\alpha_{5}} + 12 \, \theta_{0} \overline{\alpha_{5}}\right)} A_{0} + {\left(\theta_{0}^{4} \overline{\alpha_{1}} + 8 \, \theta_{0}^{3} \overline{\alpha_{1}} + 19 \, \theta_{0}^{2} \overline{\alpha_{1}} + 12 \, \theta_{0} \overline{\alpha_{1}}\right)} A_{1}\right)} \overline{A_{2}}\nonumber\\
	    			&=48{\left(\theta_{0}^{4} + 8 \, \theta_{0}^{3} + 19 \, \theta_{0}^{2} + 12 \, \theta_{0} \right)}\bar{\alpha_1\alpha_5}
	    		\end{align}
	    		Then, we have by \eqref{step0}
	    		\begin{align}\label{mu1510}
	    			\mathrm{(X)}&=-6 \, {\left(\theta_{0}^{4} + 8 \, \theta_{0}^{3} + 19 \, \theta_{0}^{2} + 12 \, \theta_{0}\right)} \overline{A_{2}} \overline{B_{1}}=-6 \, {\left(\theta_{0}^{4} + 8 \, \theta_{0}^{3} + 19 \, \theta_{0}^{2} + 12 \, \theta_{0}\right)} |\vec{A}_1|^2\bar{\zeta_2}
	    		\end{align}
	    		and finally
	    		\begin{align}\label{mu1511}
	    			\mathrm{(XI)}&=4 \, {\left(\theta_{0}^{4} + 9 \, \theta_{0}^{3} + 26 \, \theta_{0}^{2} + 24 \, \theta_{0}\right)} {\left| A_{1} \right|}^{2} \overline{A_{1}} \overline{C_{1}}=4 \, {\left(\theta_{0}^{4} + 9 \, \theta_{0}^{3} + 26 \, \theta_{0}^{2} + 24 \, \theta_{0}\right)} {\left| A_{1} \right|}^{2} \bar{\zeta_2}
	    		\end{align}
	    		And by the very definition of $\mu_{15}$, we have
	    		\begin{align*}
	    			\mu_{15}&=-\frac{1}{24 \, {\left(\theta_{0}^{5} + 10 \, \theta_{0}^{4} + 35 \, \theta_{0}^{3} + 50 \, \theta_{0}^{2} + 24 \, \theta_{0}\right)}}\bigg\{\mathrm{(I)}+\mathrm{(II)}+\mathrm{(III)}+\mathrm{(IV)}+\mathrm{(V)}+\mathrm{(VI)}+\mathrm{(VII)}+\mathrm{(VIII)}\\
	    			&+\mathrm{(IX)}+\mathrm{(X)}+\mathrm{(XI)} \bigg\}
	    		\end{align*}
	    		Now, remark that
	    		\begin{align*}
	    			\frac{3\bar{\zeta_0}}{4(\theta_0^3+4\theta_0^2+5\theta_0+2)}\bar{\vec{B}_1}=-\frac{3\bar{\s{\vec{A}_1}{\vec{A}_1}}\s{\vec{A}_1}{\bar{\vec{C}_1}}}{2(\theta_0^3+4\theta_0^2+5\theta_0+2)}\bar{\vec{A}_0}=-\frac{3|\vec{A}_1|^2\bar{\zeta_2}}{2(\theta_0^3+4\theta_0^2+5\theta_0+2)}\bar{\vec{A}_0}
	    		\end{align*}
	    		As the coefficient in $\bar{\vec{A}_0}\z^{\theta_0+3}$ in the Taylor development of
	    		\begin{align*}
	    			\Re\left(\p{\z}\vec{F}(z)\right)=0
	    		\end{align*}
	    		is
	    		\begin{align*}
	    			\frac{3\bar{\zeta_0}}{2(\theta_0^3+4\theta_0^2+5\theta_0+2)}\bar{\vec{B}_1}+\lambda_2\bar{\vec{A}_0}		=\left(-\frac{3|\vec{A}_1|^2\bar{\zeta_2}}{2(\theta_0^3+4\theta_0^2+5\theta_0+2)}+\lambda_2\right)\bar{\vec{A}_0}
	    		\end{align*}
	    		while the coefficient in $\bar{\vec{A}_2}$ vanishes, and the vectors are
	    		\begin{align*}
	    			\vec{A}_0,\vec{A}_1,\bar{\vec{A}_1}
	    		\end{align*}
	    		while
	    		\begin{align*}
	    			\s{\vec{A}_0}{\vec{A}_0}=\s{\vec{A}_0}{\vec{A}_1}=\s{\vec{A}_0}{\bar{\vec{A}_1}}=0
	    		\end{align*}
	    		we deduce that
	    		\begin{align}\label{notquitelasthopebutstillalittlebit}
	    			\Omega=-\frac{3|\vec{A}_1|^2\bar{\zeta_2}}{2(\theta_0^3+4\theta_0^2+5\theta_0+2)}+\lambda_2=0
	    		\end{align}
	    		Now, we only need to compute
	    		\begin{align*}
	    			&\lambda_2=\frac{1}{\theta_{0}^{4} + 6 \, \theta_{0}^{3} + 11 \, \theta_{0}^{2} + 6 \, \theta_{0}}\bigg\{\mu_{15} \theta_{0}^{5} + 10 \, \mu_{15} \theta_{0}^{4} - {\left(4 \, \overline{\alpha_{1}} \overline{\alpha_{5}} - 4 \, \alpha_{15} - 35 \, \mu_{15}\right)} \theta_{0}^{3} - 2 \, {\left(16 \, \overline{\alpha_{1}} \overline{\alpha_{5}} - 12 \, \alpha_{15} - 25 \, \mu_{15}\right)} \theta_{0}^{2}\\
	    			& - {\left(8 \, \theta_{0}^{3} \overline{\alpha_{3}} + 60 \, \theta_{0}^{2} \overline{\alpha_{3}} + 136 \, \theta_{0} \overline{\alpha_{3}} + {\left(\theta_{0} + 2\right)} \overline{\zeta_{2}} + 96 \, \overline{\alpha_{3}}\right)} {\left| A_{1} \right|}^{2} - 4 \, {\left(19 \, \overline{\alpha_{1}} \overline{\alpha_{5}} - 11 \, \alpha_{15} - 6 \, \mu_{15}\right)} \theta_{0}\\
	    			& - {\left(\theta_{0}^{4} + 10 \, \theta_{0}^{3} + 35 \, \theta_{0}^{2} + 50 \, \theta_{0} + 24\right)} \zeta_{13} - 48 \, \overline{\alpha_{1}} \overline{\alpha_{5}} - {\left(\theta_{0}^{3} \overline{\alpha_{2}} + 7 \, \theta_{0}^{2} \overline{\alpha_{2}} + 12 \, \theta_{0} \overline{\alpha_{2}}\right)} \overline{\zeta_{0}} - {\left(\theta_{0}^{3} + 6 \, \theta_{0}^{2} + 11 \, \theta_{0} + 6\right)} \overline{\zeta_{16}} + 24 \, \alpha_{15}\bigg\}\\
	    			&=\frac{1}{\theta_{0}^{4} + 6 \, \theta_{0}^{3} + 11 \, \theta_{0}^{2} + 6 \, \theta_{0}}\bigg\{\mu_{15} \theta_{0}^{5} + 10 \, \mu_{15} \theta_{0}^{4} - {\left(4 \, \overline{\alpha_{1}} \overline{\alpha_{5}} - 4 \, \alpha_{15} - 35 \, \mu_{15}\right)} \theta_{0}^{3} - 2 \, {\left(16 \, \overline{\alpha_{1}} \overline{\alpha_{5}} - 12 \, \alpha_{15} - 25 \, \mu_{15}\right)} \theta_{0}^{2}\\
	    			& - {\left(8 \, \theta_{0}^{3} \overline{\alpha_{3}} + 60 \, \theta_{0}^{2} \overline{\alpha_{3}} + 136 \, \theta_{0} \overline{\alpha_{3}} + {\left(\theta_{0} + 2\right)} \overline{\zeta_{2}} + 96 \, \overline{\alpha_{3}}\right)} {\left| A_{1} \right|}^{2} - 4 \, {\left(19 \, \overline{\alpha_{1}} \overline{\alpha_{5}} - 11 \, \alpha_{15} - 6 \, \mu_{15}\right)} \theta_{0}\\
	    			& - {\left(\theta_{0}^{4} + 10 \, \theta_{0}^{3} + 35 \, \theta_{0}^{2} + 50 \, \theta_{0} + 24\right)} \zeta_{13} - 48 \, \overline{\alpha_{1}} \overline{\alpha_{5}} - \frac{{\left(\theta_{0}^{2}  + 7 \, \theta_{0}  + 12 \,  \right)}}{2(\theta_0+1)} |A_1|^2\bar{\zeta_2} - {\left(\theta_{0}^{3} + 6 \, \theta_{0}^{2} + 11 \, \theta_{0} + 6\right)} \overline{\zeta_{16}} + 24 \, \alpha_{15}\bigg\}
	    		\end{align*}
	    		Here the only change between the two lines is to write
	    		\begin{align*}
	    			{\left(\theta_{0}^{3} \overline{\alpha_{2}} + 7 \, \theta_{0}^{2} \overline{\alpha_{2}} + 12 \, \theta_{0} \overline{\alpha_{2}}\right)} \overline{\zeta_{0}}=\frac{{\left(\theta_{0}^{2}  + 7 \, \theta_{0} \overline{\alpha_{2}} + 12 \,  \overline{\alpha_{2}}\right)}}{2(\theta_0+1)} |A_1|^2\bar{\zeta_2}
	    		\end{align*}
	    		as
	    		\begin{align*}
	    			2\theta_0(\theta_0+1)\alpha_2\zeta_0=\s{\bar{\vec{A}_1}}{\vec{C}_1}\s{\vec{A}_1}{\vec{A}_1}=|\vec{A}_1|^2\s{\vec{A}_1}{\vec{C}_1}=|\vec{A}_1|^2\zeta_2
	    		\end{align*}
	    		We emphasize the new development of $\lambda_2$ here
	    		\begin{align}\label{deflambda2}
	    			&\lambda_2=\frac{1}{\theta_{0}^{4} + 6 \, \theta_{0}^{3} + 11 \, \theta_{0}^{2} + 6 \, \theta_{0}}\bigg\{\mu_{15} \theta_{0}^{5} + 10 \, \mu_{15} \theta_{0}^{4} - {\left(4 \, \overline{\alpha_{1}} \overline{\alpha_{5}} - 4 \, \alpha_{15} - 35 \, \mu_{15}\right)} \theta_{0}^{3} - 2 \, {\left(16 \, \overline{\alpha_{1}} \overline{\alpha_{5}} - 12 \, \alpha_{15} - 25 \, \mu_{15}\right)} \theta_{0}^{2}\nonumber\\
	    			& - {\left(8 \, \theta_{0}^{3} \overline{\alpha_{3}} + 60 \, \theta_{0}^{2} \overline{\alpha_{3}} + 136 \, \theta_{0} \overline{\alpha_{3}} + {\left(\theta_{0} + 2\right)} \overline{\zeta_{2}} + 96 \, \overline{\alpha_{3}}\right)} {\left| A_{1} \right|}^{2} - 4 \, {\left(19 \, \overline{\alpha_{1}} \overline{\alpha_{5}} - 11 \, \alpha_{15} - 6 \, \mu_{15}\right)} \theta_{0}\nonumber\\
	    			& - {\left(\theta_{0}^{4} + 10 \, \theta_{0}^{3} + 35 \, \theta_{0}^{2} + 50 \, \theta_{0} + 24\right)} \zeta_{13} - 48 \, \overline{\alpha_{1}} \overline{\alpha_{5}} - \frac{{\left(\theta_{0}^2  + 7 \, \theta_{0}  + 12 \,  \right)}}{2(\theta_0+1)} |A_1|^2\bar{\zeta_2} - {\left(\theta_{0}^{3} + 6 \, \theta_{0}^{2} + 11 \, \theta_{0} + 6\right)} \overline{\zeta_{16}} + 24 \, \alpha_{15}\bigg\}
	    		\end{align}
	    		Now, we will code each term and compare with the expression with obtained, to finally obtain the expression of the coefficient which we shall not name $\Omega$, as it did not bring us any luck so far.
	    		\textbf{Sage version}
	    		\begin{align*}
	    			&\lambda_2=\frac{1}{2 \, {\left(\theta_{0}^{4} + 6 \, \theta_{0}^{3} + 11 \, \theta_{0}^{2} + 6 \, \theta_{0}\right)}}\bigg\{2 \, \mu_{15} \theta_{0}^{5} + 20 \, \mu_{15} \theta_{0}^{4} - 2 \, {\left(4 \, \overline{\alpha_{1}} \overline{\alpha_{5}} - 4 \, \alpha_{15} - 35 \, \mu_{15}\right)} \theta_{0}^{3}\\
	    			& - 4 \, {\left(16 \, \overline{\alpha_{1}} \overline{\alpha_{5}} - 12 \, \alpha_{15} - 25 \, \mu_{15}\right)} \theta_{0}^{2} - 2 \, {\left(8 \, \theta_{0}^{3} \overline{\alpha_{3}} + 60 \, \theta_{0}^{2} \overline{\alpha_{3}} + 136 \, \theta_{0} \overline{\alpha_{3}} + {\left(\theta_{0} + 2\right)} \overline{\zeta_{2}} + 96 \, \overline{\alpha_{3}}\right)} {\left| A_{1} \right|}^{2}\\
	    			& - \frac{{\left(\theta_{0}^{2} + 7 \, \theta_{0} + 12\right)} {\left| A_{1} \right|}^{2} \overline{\zeta_{2}}}{\theta_{0} + 1} - 8 \, {\left(19 \, \overline{\alpha_{1}} \overline{\alpha_{5}} - 11 \, \alpha_{15} - 6 \, \mu_{15}\right)} \theta_{0} - 2 \, {\left(\theta_{0}^{4} + 10 \, \theta_{0}^{3} + 35 \, \theta_{0}^{2} + 50 \, \theta_{0} + 24\right)} \zeta_{13}\\
	    			& - 96 \, \overline{\alpha_{1}} \overline{\alpha_{5}} - 2 \, {\left(\theta_{0}^{3} + 6 \, \theta_{0}^{2} + 11 \, \theta_{0} + 6\right)} \overline{\zeta_{16}} + 48 \, \alpha_{15}\bigg\}
	    		\end{align*}
	    		Sage decided to factor the $\dfrac{1}{2}$ coming from the last line and changes the order of some factors, but one can nevertheless check that both expressions coincide.
	    		We also check that (the left -hand side is our expression and the right-hand side is the Sage expression - the attentive reader will check that Sage always puts spaces between numbers and symbolic characters)
	    		\begin{align*}
	    			-\frac{3|\vec{A}_1|^2\bar{\zeta_2}}{2(\theta_0^3+4\theta_0^2+5\theta_0+2)}=\Omega-\lambda_2=-\frac{3 \, {\left| A_{1} \right|}^{2} \overline{\zeta_{2}}}{2 \, {\left(\theta_{0}^{3} + 4 \, \theta_{0}^{2} + 5 \, \theta_{0} + 2\right)}}.
	    		\end{align*}
	    		Then, we have
	    		\begin{align*}
	    			\mu_{15}&=-\frac{1}{24 \, {\left(\theta_{0}^{5} + 10 \, \theta_{0}^{4} + 35 \, \theta_{0}^{3} + 50 \, \theta_{0}^{2} + 24 \, \theta_{0}\right)}}\bigg\{\mathrm{(I)}+\mathrm{(II)}+\mathrm{(III)}+\mathrm{(IV)}+\mathrm{(V)}+\mathrm{(VI)}+\mathrm{(VII)}+\mathrm{(VIII)}\\
	    			&+\mathrm{(IX)}+\mathrm{(X)}+\mathrm{(XI)} \bigg\}\\
	    			&=-\frac{1}{24 \, {\left(\theta_{0}^{5} + 10 \, \theta_{0}^{4} + 35 \, \theta_{0}^{3} + 50 \, \theta_{0}^{2} + 24 \, \theta_{0}\right)}}\bigg\{\mathit{roman}_{1} + \mathit{roman}_{10} + \mathit{roman}_{11} + \mathit{roman}_{2} + \mathit{roman}_{3} + \mathit{roman}_{4}\\
	    			& + \mathit{roman}_{5} + \mathit{roman}_{6} + \mathit{roman}_{7} + \mathit{roman}_{8} + \mathit{roman}_{9}\bigg\}
	    		\end{align*}
	    		where for all $1\leq j \leq n$, we have $ 			\mathit{roman}_j$ is the corresponding Roman numeral to the Arabic numeral $j$.
	    		Now, we have
	    		\begin{align*}
	    			\zeta_{13}&=\frac{-(2\theta_0-1)}{\theta_0(\theta_0+1)(\theta_0+2)(\theta_0-3)}|\vec{A}_1|^2\bar{\zeta_2}\\
	    			&=-\frac{{\left(2 \, \theta_{0} - 1\right)} {\left| A_{1} \right|}^{2} \overline{\zeta_{2}}}{{\left(\theta_{0} + 2\right)} {\left(\theta_{0} + 1\right)} {\left(\theta_{0} - 3\right)} \theta_{0}}
	    		\end{align*}
	    		and
	    		\begin{align*}
	    			\bar{\zeta_{16}}&=\frac{2(\theta_0^3-\theta_0^2-6\theta_0-10)}{\theta_0(\theta_0+1)(\theta_0+2)(\theta_0-3)}|\vec{A}_1|^2\bar{\zeta_2}\\
	    			&=\frac{2 \, {\left(\theta_{0}^{3} - \theta_{0}^{2} - 6 \, \theta_{0} - 10\right)} {\left| A_{1} \right|}^{2} \overline{\zeta_{2}}}{{\left(\theta_{0} + 2\right)} {\left(\theta_{0} + 1\right)} {\left(\theta_{0} - 3\right)} \theta_{0}}
	    		\end{align*}
	    		Finally, we have (the first line is the \TeX\;, and the second one the Sage version)
	    		\begin{align*}
	    			\mathrm{(I)}&=48{\left(\theta_{0}^{4} + 7 \, \theta_{0}^{3} + 14 \, \theta_{0}^{2} + 8 \, \theta_{0}\right)}|\vec{A}_1|^2\bar{\alpha_3}-4{\left(\theta_{0}^{4} + 7 \, \theta_{0}^{3} + 14 \, \theta_{0}^{2} + 8 \, \theta_{0}\right)}|\vec{A}_1|^2\bar{\zeta_2}\\
	    			&=48 \, {\left(\theta_{0}^{4} + 7 \, \theta_{0}^{3} + 14 \, \theta_{0}^{2} + 8 \, \theta_{0}\right)} {\left| A_{1} \right|}^{2} \overline{\alpha_{3}} - 4 \, {\left(\theta_{0}^{4} + 7 \, \theta_{0}^{3} + 14 \, \theta_{0}^{2} + 8 \, \theta_{0}\right)} {\left| A_{1} \right|}^{2} \overline{\zeta_{2}}\\
	    			\mathrm{(II)}&=-2 \, {\left(\theta_{0}^{4} + 10 \, \theta_{0}^{3} + 35 \, \theta_{0}^{2} + 50 \, \theta_{0} + 24\right)} {\left| A_{1} \right|}^{2}\bar{\zeta_2}\\
	    			&=-2 \, {\left(\theta_{0}^{4} + 10 \, \theta_{0}^{3} + 35 \, \theta_{0}^{2} + 50 \, \theta_{0} + 24\right)} {\left| A_{1} \right|}^{2}\bar{\zeta_2}\\
	    			\mathrm{(III)}&=-96 {\left(\theta_{0}^{4}  + 10 \, \theta_{0}^{3}  + 35 \, \theta_{0}^{2}  + 50 \, \theta_{0}  + 24 \, \right)}|\vec{A}_1|^2\bar{\alpha_3}\\
	    			&=-96 \, {\left(\theta_{0}^{4} + 10 \, \theta_{0}^{3} + 35 \, \theta_{0}^{2} + 50 \, \theta_{0} + 24\right)} {\left| A_{1} \right|}^{2} \overline{\alpha_{3}}\\
	    			\mathrm{(IV)}&=-24\left(\theta_0^4+10\,\theta_0^3+35\,\theta_0^2+50\,\theta_0+24\right)(2\bar{\alpha_1\alpha_5}-\alpha_{15})\\
	    			&=-24 \, {\left(\theta_{0}^{4} + 10 \, \theta_{0}^{3} + 35 \, \theta_{0}^{2} + 50 \, \theta_{0} + 24\right)} {\left(2 \, \overline{\alpha_{1}} \overline{\alpha_{5}} - \alpha_{15}\right)}\\
	    			\mathrm{(V)}&=- 24 \, {\left(\theta_{0}^{4} + 6 \, \theta_{0}^{3} + 11 \, \theta_{0}^{2} + 6 \, \theta_{0}\right)} \alpha_{15}+\frac{ {\left(\theta_{0}^{4} + 6 \, \theta_{0}^{3} + 11 \, \theta_{0}^{2} + 6 \, \theta_{0}\right)} (3\theta_0-7)}{\theta_0-3}|\vec{A}_1|^2\bar{\zeta_2}\\
	    			&=\frac{{\left(\theta_{0}^{4} + 6 \, \theta_{0}^{3} + 11 \, \theta_{0}^{2} + 6 \, \theta_{0}\right)} {\left(3 \, \theta_{0} - 7\right)} {\left| A_{1} \right|}^{2} \overline{\zeta_{2}}}{\theta_{0} - 3} - 24 \, {\left(\theta_{0}^{4} + 6 \, \theta_{0}^{3} + 11 \, \theta_{0}^{2} + 6 \, \theta_{0}\right)} \alpha_{15}\\
	    			\mathrm{(VI)}&=\frac{8{\left(\theta_{0}^{4} + 9 \, \theta_{0}^{3} + 26 \, \theta_{0}^{2} + 24 \, \theta_{0}\right)}(\theta_0-4)}{\theta_0-3}|\vec{A}_1|^2\bar{\zeta_2}\\
	    			&=\frac{8 \, {\left(\theta_{0}^{4} + 9 \, \theta_{0}^{3} + 26 \, \theta_{0}^{2} + 24 \, \theta_{0}\right)} {\left(\theta_{0} - 4\right)} {\left| A_{1} \right|}^{2} \overline{\zeta_{2}}}{\theta_{0} - 3}\\
	    			\mathrm{(VII)}&=\frac{- 3 \, {\left(\theta_{0}^{4} + 10 \, \theta_{0}^{3} + 35 \, \theta_{0}^{2} + 50 \, \theta_{0} + 24\right)}(\theta_0-5)}{\theta_0-3}|\vec{A}_1|^2\bar{\zeta_2}\\
	    			&=-\frac{3 \, {\left(\theta_{0}^{4} + 10 \, \theta_{0}^{3} + 35 \, \theta_{0}^{2} + 50 \, \theta_{0} + 24\right)} {\left(\theta_{0} - 5\right)} {\left| A_{1} \right|}^{2} \overline{\zeta_{2}}}{\theta_{0} - 3}
	    			\\
	    			\mathrm{(VIII)}&=48{\left(\theta_{0}^{4} + 9 \, \theta_{0}^{3} + 26 \, \theta_{0}^{2} + 24 \, \theta_{0} \right)}|\vec{A}_1|^2\overline{\alpha_{3}}\\
	    			&=48 \, {\left(\theta_{0}^{4} + 9 \, \theta_{0}^{3} + 26 \, \theta_{0}^{2} + 24 \, \theta_{0}\right)} {\left| A_{1} \right|}^{2} \overline{\alpha_{3}}\\
	    			\mathrm{(IX)}&=48{\left(\theta_{0}^{4} + 8 \, \theta_{0}^{3} + 19 \, \theta_{0}^{2} + 12 \, \theta_{0} \right)}\bar{\alpha_1\alpha_5}\\
	    			&=48 \, {\left(\theta_{0}^{4} + 8 \, \theta_{0}^{3} + 19 \, \theta_{0}^{2} + 12 \, \theta_{0}\right)} \overline{\alpha_{1}} \overline{\alpha_{5}}
	    			\\
	    			\mathrm{(X)}&=-6 \, {\left(\theta_{0}^{4} + 8 \, \theta_{0}^{3} + 19 \, \theta_{0}^{2} + 12 \, \theta_{0}\right)} |\vec{A}_1|^2\bar{\zeta_2}\\
	    			&=-6 \, {\left(\theta_{0}^{4} + 8 \, \theta_{0}^{3} + 19 \, \theta_{0}^{2} + 12 \, \theta_{0}\right)} {\left| A_{1} \right|}^{2} \overline{\zeta_{2}}
	    			\\
	    			\mathrm{(XI)}&=4 \, {\left(\theta_{0}^{4} + 9 \, \theta_{0}^{3} + 26 \, \theta_{0}^{2} + 24 \, \theta_{0}\right)} {\left| A_{1} \right|}^{2} \bar{\zeta_2}\\
	    			&=4 \, {\left(\theta_{0}^{4} + 9 \, \theta_{0}^{3} + 26 \, \theta_{0}^{2} + 24 \, \theta_{0}\right)} {\left| A_{1} \right|}^{2} \overline{\zeta_{2}}
	    		\end{align*}
	    		And finally
	    		\begin{align}\label{siseulement}
	    			\Omega=-\frac{2 \, {\left(\theta_{0} - 4\right)} {\left| A_{1} \right|}^{2} \overline{\zeta_{2}}}{\theta_{0}^{3} - 3 \, \theta_{0}^{2}}=0.
	    		\end{align}
	    		where $\zeta_2=\s{\vec{A}_1}{\vec{C}_1}$.
	    		This works for $\theta_0\geq 5$... Therefore, we obtain the relation
	    		\begin{align}\label{fautycroire}
	    			\bar{\Omega}=-\frac{2(\theta_0-4)}{\theta_0^2(\theta_0-3)}|\vec{A}_1|^2\s{\vec{A}_1}{\vec{C}_1}=0.
	    		\end{align}
	    		To get this relation, we have used the meromorphy of the quartic form to obtain
	    		\begin{align*}
	    			|\vec{A}_1|^2\s{\vec{A}_1}{\vec{C}_1}=\s{\bar{\vec{A}_1}}{\vec{C}_1}\s{\vec{A}_1}{\vec{A}_1}.
	    		\end{align*}
	    		\section{Conclusion}
	    		
	    		Thanks of \eqref{fautycroire}, we obtain
	    		\begin{align}\label{shinjiru}
	    			\s{\vec{A}_1}{\vec{C}_1}=0
	    		\end{align}
	    		and the holomorphy of the quartic form, as
	    		\begin{align*}
	    			\mathscr{Q}_{\phi}=(\theta_0-1)(\theta_0-2)\s{\vec{A}_1}{\vec{C}_1}\frac{dz^4}{z}+O(1).
	    		\end{align*}
	    		Furthermore, notice that as
	    		\begin{align*}
	    			&|\vec{A}_1|^2\s{\vec{A}_1}{\vec{C}_1}=\s{\bar{\vec{A}_1}}{\vec{C}_1}\s{\vec{A}_1}{\vec{C}_1}\\
	    			&|\vec{C}_1|^2\s{\vec{A}_1}{\vec{A}_1}=\s{\vec{A}_1}{\vec{C}_1}\s{\vec{A}_1}{\bar{\vec{C}_1}}=0,
	    		\end{align*}
	    		we have $\vec{C}_1=0$ or $\s{\vec{A}_1}{\vec{A}_1}=0$, and $\s{\vec{A}_1}{\vec{A}_1}=0$ or $\s{\bar{\vec{A}_1}}{\vec{C}_1}=0$.

	    		\chapter{Removability of the poles of the octic form}\label{chapteroctic}
	    		\normalsize
	    	As the holomorphy of the octic form $\mathscr{O}_{\phi}$ does not follow trivially from the asymptotic expansion of $\h_0$, we will need a normal derivative free expression, which is the content of the following proposition.

	    	\begin{prop}
	    		Let $\phi:\Sigma^2\rightarrow S^4$ a smooth immersion. Then we have
	    		\begin{align}\label{publiable}
	    		&\mathscr{O}_{\phi}=-\frac{1}{4}\,g^{-3}\otimes\partial\left(\h_0\totimes\bar{\h_0}\right)\otimes \left(\bar{\partial}\h_0\totimes\h_0\right)\otimes\left(\h_0\totimes\h_0\right)\nonumber\\
	    		&+g^{-2}\otimes\bigg\{\frac{1}{4}\left(\partial\bar{\partial}\h_0\,\dot{\otimes}\,\partial\bar{\partial}\h_0\right)\otimes \left(\h_0\,\dot{\otimes}\,\h_0\right)+\frac{1}{4}\left(\partial\h_0\,\dot{\otimes}\, \partial\h_0\right)\otimes \left(\bar{\partial}\h_0\,\dot{\otimes}\, \bar{\partial}\h_0\right)\nonumber\\
	    		&-\frac{1}{2}\left(\partial\bar{\partial}\h_0\totimes\partial\h_0\right)\otimes\left(\bar{\partial}\h_0\totimes\h_0\right)-\frac{1}{2}\left(\partial\bar{\partial}\h_0\totimes\bar{\partial}\h_0\right)\otimes\left(\partial\h_0\totimes\h_0\right)+\frac{1}{2}\left(\partial\bar{\partial}\h_0\totimes\h_0\right)\otimes\left(\partial\h_0\totimes\bar{\partial}\h_0\right)\nonumber\\
	    		&-\frac{1}{4}\Big\{2\s{\H}{\bar{\partial}\h_0}\otimes\left(\bar{\partial}\h_0\otimes\h_0\right)+\s{\H}{\partial\h_0}\otimes\partial\left(\h_0\totimes\bar{\h_0}\right)+\s{\H}{\h_0}\otimes\left(\bar{\partial}\h_0\totimes\bar{\partial}\h_0\right)\Big\}\otimes\left(\h_0\totimes\h_0\right)\nonumber\\
	    		&+\frac{1}{4}\,\bs{\H}{\bar{\partial}\left(\left(\h_0\totimes\h_0\right)\h_0\right)}\otimes\left(\bar{\partial}\h_0\totimes\h_0\right)
	    		+\frac{1}{4}\,\bs{\H}{\partial\left(\left(\h_0\totimes\bar{\h_0}\right)\h_0\right)}\otimes \left({\partial}\h_0\totimes\h_0\right)\bigg\}\nonumber\\
	    		&+\frac{1}{4}\,g^{-1}\,\otimes\bigg\{\bigg(\s{\H}{\partial\bar{\partial}\h_0}\otimes\s{\H}{\h_0}-\s{\H}{\partial\h_0}\otimes\s{\H}{\bar{\partial}\h_0}\bigg)\otimes\left(\h_0\totimes\h_0\right)\nonumber\\
	    		&-\left(\partial\bar{\partial}\h_0\totimes\h_0-\partial\h_0\totimes\bar{\partial}\h_0\right)\otimes\left(|\h_0|^{2}_{WP}\,\h_0\totimes\h_0+\s{\H}{\h_0}^2\right)\bigg\}\nonumber\\
	    		&+\frac{1}{4}|\h_0|^2_{WP}\,g^{-1}\otimes\left(\partial\bar{\partial}\h_0\totimes\h_0\right)\otimes\left(\h_{0}\totimes\h_0\right)
	    		-\frac{1}{4}\,|\h_0|^{2}_{WP}\,\bigg\{\left(\partial\h_0\totimes\partial\h_0\right)\otimes\s{\H}{\h_0}-\left(\h_0\totimes\h_0\right)\otimes\s{\H}{\h_0}^2\bigg\}\nonumber\\
	    		&+\frac{1}{16}\,\left\{\left(|\h_0|_{WP}^2\,\h_0\totimes\h_0\right)^2+|\H|^2\left(\h_0\totimes\h_0\right)\otimes\s{\H}{\h_0}+2|\h_0|_{WP}^2\left(\h_0\totimes\h_0\right)\otimes \s{\H}{\h_0}^2\right\}\nonumber\\
	    		&+\frac{1}{4}\left(1+|\H|^2\right)\,g^{-1}\otimes\left\{\frac{1}{2}\left(
	    		\partial\bar{\partial}\h_0\totimes\h_0\right)\otimes\left(\h_0\totimes\h_0\right)-\left(\partial\h_0\totimes\h_0\right)\otimes\left(\bar{\partial}\h_0\totimes\h_0\right) +\frac{1}{2}\left(\partial\h_0\totimes\bar{\partial}\h_0\right)\otimes\left(\h_0\totimes\h_0\right)\right\}\nonumber\\
	    		&+\frac{1}{64}\left(1+|\H|^2\right)^2\left(\h_0\totimes\h_0\right)\otimes\left(\h_0\totimes\h_0\right)-\frac{1}{8} \left(|\h_0|^{2}_{WP}\,\h_0\totimes\h_0+\s{\H}{\h_0}^2\right)^2
	    		\end{align}
	    	\end{prop}
	    	\begin{proof}
	    		We first remark that for any two tensors $\vec{\alpha},\vec{\beta}$, and any two differential operators $\D_1,\D_2\in \ens{\partial,\bar{\partial}}$, we have
	    		\begin{align}\label{denormalisation}
	    		\D_1^{\perp}\vec{\alpha}\totimes\D_2^{N}{\vec{\beta}}&=\left(\D_1-\D_1^\top\right)\vec{\alpha}\totimes \left(\D_2-\D_2^\top\right)\vec{\beta}\nonumber\\
	    		&=\D_1\vec{\alpha}\totimes\D_2\vec{\beta}-\D_1\vec{\alpha}\totimes\D_2^{\top}\vec{\beta}-\D_1^{\top}\vec{\alpha}\totimes\D_2\vec{\beta}+\D_1^{\top}\vec{\alpha}\totimes\D_2^{\top}\vec{\beta}\nonumber\\
	    		&=\D_1\vec{\alpha}\totimes\D_2\vec{\beta}-\D_1^{\top}\vec{\alpha}\totimes\D_2^{\top}\vec{\beta}
	    		\end{align}
	    		as $\D_1\vec{\alpha}\totimes\D_2^{\top}\vec{\beta}=\D_1^{\top}\vec{\alpha}\totimes\D_2^{\top}\vec{\beta}$. 
	    		Now, recall that
	    		\begin{align}\label{codazzi0}
	    		\left\{\begin{alignedat}{1}
	    		\partial^\top\h_0&=-\s{\H}{\h_0}\otimes \partial\phi-g^{-1}\otimes (\h_0\totimes\h_0)\otimes\bar{\partial}\phi\\
	    		\bar{\partial}^\top\h_0&=-|\h_0|^2_{WP}\;g\otimes \partial\phi-\s{\H}{\h_0}\otimes \bar{\partial}\phi
	    		\end{alignedat}\right.
	    		\end{align}
	    		As $\phi$ is conformal, $g=\partial\totimes\bar{\partial}\phi$ and $|\h_0|^2_{WP}=g^{-2}\otimes\left(\h_0\totimes\bar{\h_0}\right)$, we obtain
	    		\begin{align}\label{codazzi1}
	    		\left\{\begin{alignedat}{1}
	    		\partial^\top\h_0\totimes\bar{\partial}^\top\h_0&=\frac{1}{2}\,g\otimes\left(|\h_0|^2_{WP}\,\h_0\totimes\h_0+\s{\H}{\h_0}^2\right)\\
	    		\partial^\top\h_0\totimes\partial^\top\h_0&=\s{\H}{\h_0}\otimes\left(\h_0\totimes\h_0\right)\\
	    		\bar{\partial}^\top\h_0\totimes\bar{\partial}^\top\h_0&=|\h_0|^2_{WP}\,g^2\otimes\s{\H}{\h_0}=\s{\H}{\h_0}\otimes \left(\h_0\totimes\bar{\h_0}\right).
	    		\end{alignedat}\right.
	    		\end{align}
	    		Then we compute thanks to \eqref{codazzi0}
	    		\begin{align}\label{codazzi2}
	    		&\partial^N\bar{\partial}^N\h_0=\partial^N\left(\bar{\partial}\h_0+|\h_0|^{2}_{WP}\,g\otimes\partial\phi+\s{\H}{\h_0}\otimes\bar{\partial}\phi\right)\nonumber\\
	    		&=\partial^N\bar{\partial}\h_0+|\h_0|^2_{WP}\,g\otimes \partial^N\partial\phi+\s{\H}{\h_0}\otimes \partial^N\bar{\partial}\phi\nonumber\\
	    		&=\partial^N\bar{\partial}\h_0+\frac{1}{2}\,g\otimes\left(|\h_0|^{2}_{WP}\,\h_0+\s{\H}{\h_0}\H\right)
	    		\end{align}
	    		In particular, as $\h_0$ is normal, we have
	    		\begin{align}\label{codazzi3}
	    		\partial^N\bar{\partial}^N\h_0\totimes\h_0=\partial\bar{\partial}\h_0\totimes\h_0+\frac{1}{2}\,g\otimes\left(|\h_0|^{2}_{WP}\,\h_0\totimes\h_0+\s{\H}{\h_0}^2\right).
	    		\end{align}
	    		Then, by \eqref{codazzi1} and \eqref{denormalisation}, we have
	    		\begin{align}\label{codazzi4}
	    		&\partial^N\h_0\totimes\bar{\partial}^N\h_0=\partial\h_0\totimes\bar{\partial}\h_0-\partial^\top\h_0\totimes\bar{\partial}^\top\h_0=\partial\h_0\totimes\bar{\partial}\h_0-\frac{1}{2}\,g\otimes\left(|\h_0|^{2}_{WP}\,\h_0\totimes\h_0+\s{\H}{\h_0}^2\right).
	    		\end{align}
	    		Therefore, we finally obtain by \eqref{codazzi3} and \eqref{codazzi4}
	    		\begin{align*}
	    		&\left(
	    		\partial^N\bar{\partial}^N\h_0\totimes\h_0\right)\otimes\left(\h_0\totimes\h_0\right)+\left(\partial^N\h_0\totimes\bar{\partial}^N\h_0\right)\otimes\left(\h_0\totimes\h_0\right)\\
	    		&=\left(
	    		\partial\bar{\partial}\h_0\totimes\h_0\right)\otimes\left(\h_0\totimes\h_0\right)+ \left(\partial\h_0\totimes\bar{\partial}\h_0\right)\otimes\left(\h_0\totimes\h_0\right).
	    		\end{align*}
	    		which shows by normality of $\h_0$ that
	    		\begin{align}\label{removenormals}
	    		&g^{-1}\otimes\left\{\frac{1}{2}\left(
	    		\partial^N\bar{\partial}^N\h_0\totimes\h_0\right)\otimes\left(\h_0\totimes\h_0\right)-\left(\partial^N\h_0\totimes\h_0\right)\otimes\left(\bar{\partial}^N\h_0\totimes\h_0\right) +\frac{1}{2}\left(\partial^N\h_0\totimes\bar{\partial}^N\h_0\right)\otimes\left(\h_0\totimes\h_0\right)\right\}\nonumber\\
	    		&=\,g^{-1}\otimes\left\{\frac{1}{2}\left(
	    		\partial\bar{\partial}\h_0\totimes\h_0\right)\otimes\left(\h_0\totimes\h_0\right)-\left(\partial\h_0\totimes\h_0\right)\otimes\left(\bar{\partial}\h_0\totimes\h_0\right) +\frac{1}{2}\left(\partial\h_0\totimes\bar{\partial}\h_0\right)\otimes\left(\h_0\totimes\h_0\right)\right\}
	    		\end{align}
	    		so we can simply remove the normal derivatives.
	    		Then, we have almost by definition
	    		\begin{align}\label{almostdef}
	    		\partial^\top\bar{\partial}\h_0=2\,g^{-1}\otimes\left(\partial\bar{\partial}\h_0\totimes\bar{\partial}\phi\right)\partial\phi+2\,g^{-1}\otimes\left(\partial\bar{\partial}\h_0\totimes\partial\phi\right)\bar{\partial}\phi.
	    		\end{align}
	    		Recall now that
	    		\begin{align*}
	    		\h_0=2\,\vec{\mathbb{I}}(\p{z}\phi,\p{z}\phi)dz^2=2\left(\p{z}^2\phi-2(\p{z}\lambda)\p{z}\phi\right)dz^2.
	    		\end{align*}
	    		We have as $\h_0$ is normal the identity
	    		\begin{align*}
	    		\s{\bar{\partial}\h_0}{\partial\phi}=-\s{\h_0}{\partial\bar{\partial}\phi}=-\frac{1}{2}\,g\otimes\s{\H}{\h_0}
	    		\end{align*}
	    		so by Codazzi identity, we have
	    		\begin{align}\label{octic1}
	    		\s{\partial\bar{\partial}\h_0}{\partial\phi}&=-\frac{1}{2}\partial\left(g\otimes \s{\H}{\h_0}\right)-\s{\bar{\partial}\h_0}{\p{z}^2\phi}\otimes dz^2\nonumber\\
	    		&=-\frac{1}{2}\left(\left(\partial g\right)\otimes \s{\H}{\h_0}+g\otimes \s{\partial\H}{\h_0}+g\otimes\s{\H}{\partial\h_0}\right)-\s{\bar{\partial}\h_0}{\p{z}^2\phi}\otimes dz^2\nonumber\\
	    		&=-\frac{1}{2}\left(\left(\partial g\right)\otimes \s{\H}{\h_0}+ \s{\bar{\partial}\h_0}{\h_0}+g\otimes\s{\H}{\partial\h_0}\right)-\s{\bar{\partial}\h_0}{\p{z}^2\phi}\otimes dz^2
	    		\end{align}
	    		but
	    		\begin{align}\label{octic2}
	    		\s{\bar{\partial}\h_0}{\p{z}^2\phi}\otimes dz^2&=\frac{1}{2}\s{\bar{\partial}\h_0}{\h_0}+2(\p{z}\lambda)\s{\bar{\partial}\h_0}{\p{z}\phi}\otimes dz^2=\frac{1}{2}\,\s{\bar{\partial}\h_0}{\h_0}-(\p{z}\lambda)e^{2\lambda}\s{\H}{\h_0}\otimes dz^2\otimes d\z\nonumber\\
	    		&=\frac{1}{2}\,\s{\bar{\partial}\h_0}{\h_0}-\frac{1}{2}\,\p{z}(e^{2\lambda})dz^2\otimes d\z\otimes \s{\H}{\h_0}\nonumber\\
	    		&=\frac{1}{2}\,\s{\bar{\partial}\h_0}{\h_0}-\frac{1}{2}(\partial g)\otimes\s{\H}{\h_0}.
	    		\end{align}
	    		Therefore, by \eqref{octic1} and \eqref{octic2}, we have
	    		\begin{align}\label{tangentoctic1}
	    		\partial\bar{\partial}\h_0\totimes\partial\phi&=-\frac{1}{2}\left(\ccancel{\left(\partial g\right)\otimes \s{\H}{\h_0}}+ \s{\bar{\partial}\h_0}{\h_0}+g\otimes\s{\H}{\partial\h_0}\right)-\left(\frac{1}{2}\,\s{\bar{\partial}\h_0}{\h_0}-\frac{1}{2}\ccancel{(\partial g)\otimes\s{\H}{\h_0}}\right)\nonumber\\
	    		&=-\bar{\partial}\h_0\totimes\h_0-\frac{1}{2}\,g\otimes\s{\H}{\partial\h_0}.
	    		\end{align}
	    		Now, we have again by normality of $\h_0$ the identity
	    		\begin{align}\label{octic3}
	    		\s{\bar{\partial}\h_0}{\bar{\partial}\phi}=-\s{\h_0}{\bar{\partial}^N\bar{\partial}\phi}=-\frac{1}{2}\,\h_0\totimes\bar{\h_0}
	    		\end{align}
	    		so (as for any tensor $\partial\bar{\vec{\alpha}}=\bar{\bar{\partial}{\vec{\alpha}}}$)
	    		\begin{align}\label{tangentoctic2} 
	    		\partial\bar{\partial}\h_0\totimes\bar{\partial}\phi&=-\frac{1}{2}\left(\partial\h_0\totimes\bar{\h_0}+\h_0\totimes\partial\bar{\h_0}\right)-\s{\bar{\partial}\h_0}{\partial\bar{\partial}\phi}]\nonumber\\
	    		&=-\frac{1}{2}\left(\partial\h_0\totimes\bar{\h_0}+\h_0\totimes\bar{\bar{\partial}\h_0}\right)-\frac{1}{2}\,g\otimes\s{\H}{\bar{\partial}\h_0}.
	    		\end{align}
	    		Putting together \eqref{almostdef}, \eqref{tangentoctic1} and \eqref{tangentoctic2}, we obtain
	    		\begin{align}\label{tangentend}
	    		\partial^\top\bar{\partial}\h_0=-\left(g^{-1}\otimes\left(\partial\h_0\totimes\bar{\h_0}+\h_0\totimes\bar{\bar{\partial}\h_0}\right)+\s{\H}{\bar{\partial}\h_0}\right)\partial\phi-\left(2\,g^{-1}\otimes\left(\bar{\partial}\h_0\totimes\h_0\right)+\s{\H}{\partial\h_0}\right)\bar{\partial}\phi
	    		\end{align}
	    		This relation implies immediately that
	    		\begin{align}\label{partn2}
	    		&\partial^\top\bar{\partial}\h_0\totimes\partial^\top\bar{\partial}\h_0=g\otimes\left(g^{-1}\otimes\left(\partial\h_0\totimes\bar{\h_0}+\h_0\totimes\bar{\bar{\partial}\h_0}\right)+\s{\H}{\bar{\partial}\h_0}\right)\otimes
	    		\left(2\,g^{-1}\otimes\left(\bar{\partial}\h_0\totimes\h_0\right)+\s{\H}{\partial\h_0}\right).
	    		\end{align}
	    		By \eqref{denormalisation} and \eqref{partn2}, we have
	    		\begin{align}\label{clong}
	    		&\partial^N\bar{\partial}^N\h_0\totimes\partial^N\bar{\partial}^N\h_0=\left(\partial^N\bar{\partial}\h_0+\frac{1}{2}\,g\otimes\left(|\h_0|^{2}_{WP}\,\h_0+\s{\H}{\h_0}\H\right)\right)\totimes\left(\partial^N\bar{\partial}\h_0+\frac{1}{2}\,g\otimes\left(|\h_0|^{2}_{WP}\,\h_0+\s{\H}{\h_0}\H\right)\right)\nonumber\\
	    		&=\partial^N\bar{\partial}\h_0\totimes\partial^N\bar{\partial}\h_0+g\otimes\left(|\h_0|^2_{WP}\,\partial\bar{\partial}\h_0\totimes\h_0+\s{\H}{\partial\bar{\partial}\h_0}\totimes\s{\H}{\h_0}\right)\nonumber\\
	    		&+\frac{1}{4}\,g^2\otimes\left(|\h_0|_{WP}^4\,\h_0\totimes\h_0+|\H|^2\s{\H}{\h_0}+2|\h_0|_{WP}^2\s{\H}{\h_0}^2\right)\nonumber\\
	    		&=\partial\bar{\partial}\h_0\otimes\partial\bar{\partial}\h_0-\partial^\top\bar{\partial}\h_0\totimes\partial^\top\bar{\partial}\h_0+g\otimes\left(|\h_0|^2_{WP}\,\partial\bar{\partial}\h_0\totimes\h_0+\s{\H}{\partial\bar{\partial}\h_0}\totimes\s{\H}{\h_0}\right)\nonumber\\
	    		&+\frac{1}{4}\,g^2\otimes\left(|\h_0|_{WP}^4\,\h_0\totimes\h_0+|\H|^2\s{\H}{\h_0}+2|\h_0|_{WP}^2\s{\H}{\h_0}^2\right)\nonumber\\
	    		&=\partial\bar{\partial}\h_0\totimes\partial\bar{\partial}\h_0-g\otimes\left(g^{-1}\otimes\left(\partial\h_0\totimes\bar{\h_0}+\h_0\totimes{\partial}\bar{\h_0}\right)+\s{\H}{\bar{\partial}\h_0}\right)\otimes
	    		\left(2\,g^{-1}\otimes\left(\bar{\partial}\h_0\totimes\h_0\right)+\s{\H}{\partial\h_0}\right)\nonumber\\
	    		&+g\otimes\left(|\h_0|^2_{WP}\,\partial\bar{\partial}\h_0\totimes\h_0+\s{\H}{\partial\bar{\partial}\h_0}\totimes\s{\H}{\h_0}\right)
	    		+\frac{1}{4}\,g^2\otimes\left(|\h_0|_{WP}^4\,\h_0\totimes\h_0+|\H|^2\s{\H}{\h_0}+2|\h_0|_{WP}^2\s{\H}{\h_0}^2\right)
	    		\end{align}
	    		Now, we have
	    		\begin{align}\label{d2bar1}
	    		&g\otimes\left(g^{-1}\otimes\left(\partial\h_0\totimes\bar{\h_0}+\h_0\totimes\partial\bar{\h_0}\right)+\s{\H}{\bar{\partial}\h_0}\right)\otimes
	    		\left(2\,g^{-1}\otimes\left(\bar{\partial}\h_0\totimes\h_0\right)+\s{\H}{\partial\h_0}\right)\nonumber\\
	    		&=2\,g^{-1}\otimes\left(\partial\h_0\totimes\bar{\h_0}\right)\otimes\left(\bar{\partial}\h_0\totimes\h_0\right)+2\,g^{-1}\otimes\left(\partial\bar{\h_0}\totimes\h_0\right)\otimes \left(\bar{\partial}\h_0\totimes\h_0\right)+g\otimes\s{\H}{\partial\h_0}\totimes\s{\H}{\bar{\partial}\h_0}\nonumber\\
	    		&+2\s{\H}{\bar{\partial}\h_0}\otimes\left(\bar{\partial}\h_0\totimes\h_0\right)+\s{\H}{\partial\h_0}\otimes\left(\partial\h_0\totimes\bar{\h_0}\right)+\s{\H}{\partial\h_0}\otimes\left(\partial\bar{\h_0}\totimes\h_0\right)
	    		\end{align}
	    		As
	    		\begin{align}\label{d2bar2}
	    		g\otimes \left(|\h_0|^{2}_{WP}\partial\bar{\partial}\h_0\totimes\h_0\right)=g\otimes\left(g^{-2}\otimes\left(\h_0\totimes\bar{\h_0}\right)\totimes\left(\partial\bar{\partial}\h_0\totimes\h_0\right)\right)=g^{-1}\otimes\left(\partial\bar{\partial}\h_0\totimes\h_0\right)\otimes\left(\h_0\otimes\bar{\h_0}\right),
	    		\end{align}
	    		we finally obtain by \eqref{clong}, \eqref{d2bar1} and \eqref{d2bar2}
	    		\begin{align}\label{deldelbar2}
	    		&\partial^N\bar{\partial}^N\h_0\totimes\partial^N\bar{\partial}^N\h_0=\partial\bar{\partial}\h_0\totimes\partial\bar{\partial}\h_0\nonumber\\
	    		&+g^{-1}\otimes \bigg(\left(\partial\bar{\partial}\h_0\totimes\h_0\right)\otimes\left(\h_0\otimes\bar{\h_0}\right)-2\left(\partial\h_0\totimes\bar{\h_0}\right)\otimes\left(\bar{\partial}\h_0\totimes\h_0\right)-2\left(\partial\bar{\h_0}\totimes\h_0\right)\otimes\left(\bar{\partial}\h_0\totimes\h_0\right)\bigg)\nonumber\\
	    		&-\bigg(2\s{\H}{\bar{\partial}\h_0}\otimes\left(\bar{\partial}\h_0\totimes\h_0\right)+\s{\H}{\partial\h_0}\otimes\left(\partial\h_0\totimes\bar{\h_0}\right)+\s{\H}{\partial\h_0}\otimes\left(\partial\bar{\h_0}\totimes\h_0\right)\bigg)\nonumber\\
	    		&+g\otimes\bigg(\s{\H}{\partial\bar{\partial}\h_0}\otimes\s{\H}{\h_0}-\s{\H}{\partial\h_0}\otimes\s{\H}{\bar{\partial}\h_0}\bigg)\nonumber\\
	    		&+\frac{1}{4}\,g^2\otimes\left(|\h_0|_{WP}^4\,\h_0\totimes\h_0+|\H|^2\s{\H}{\h_0}+2|\h_0|_{WP}^2\s{\H}{\h_0}^2\right).
	    		\end{align}
	    		Remains only to compute
	    		\begin{align*}
	    		\partial^N\bar{\partial}^N\h_0\totimes\partial^N\h_0,\quad\text{and}\;\, \partial^N\bar{\partial}^N\h_0\totimes\bar{\partial}^N\h_0
	    		\end{align*}
	    		Thanks to \eqref{denormalisation}, we have
	    		\begin{align}\label{altcod1}
	    		&\partial^N\bar{\partial}^N\h_0\totimes\partial^N\h_0=\left(\partial^N\bar{\partial}\h_0+\frac{1}{2}\,g\otimes \left(|\h_0|^2_{WP}\,\h_0+\s{\H}{\h_0}\H\right)\right)\otimes\partial^N\h_0\nonumber\\
	    		&=\partial^N\bar{\partial}\h_0\totimes\partial^N\h_0+\frac{1}{2}\,g\otimes\left(|\h_0|^2_{WP}\,\partial\h_0\totimes\h_0+\s{\H}{\partial\h_0}\otimes\s{\H}{\h_0}\right)\nonumber\\
	    		&=\partial\bar{\partial}\h_0\totimes\partial\h_0-\partial^\top\bar{\partial}\h_0\totimes\partial^\top\h_0+\frac{1}{2}\,g\otimes\left(|\h_0|^2_{WP}\,\partial\h_0\totimes\h_0+\s{\H}{\partial\h_0}\otimes\s{\H}{\h_0}\right).
	    		\end{align}
	    		Now, by  \eqref{codazzi0} and \eqref{tangentend}, as $g=2\,\partial\phi\totimes\bar{\partial}\phi$, we have
	    		\begin{align}\label{altcod2}
	    		&\partial^\top\bar{\partial}\h_0\h_0\totimes\partial^\top\h_0=
	    		\bigg\{-\left(g^{-1}\otimes\left(\partial\h_0\totimes\bar{\h_0}+\h_0\totimes\bar{\bar{\partial}\h_0}\right)+\s{\H}{\bar{\partial}\h_0}\right)\partial\phi\nonumber\\
	    		&-\left(2\,g^{-1}\otimes\left(\bar{\partial}\h_0\totimes\h_0\right)+\s{\H}{\partial\h_0}\right)\bar{\partial}\phi\bigg\}\totimes
	    		\bigg\{-\s{\H}{\h_0}\otimes \partial\phi-g^{-1}\otimes (\h_0\totimes\h_0)\otimes\bar{\partial}\phi\bigg\}\nonumber\\
	    		&=\frac{1}{2}\,\left(g^{-1}\otimes\left(\partial\h_0\totimes\bar{\h_0}+\h_0\totimes\bar{\bar{\partial}\h_0}\right)+\s{\H}{\bar{\partial}\h_0}\right)\otimes\left(\h_0\totimes\h_0\right)\nonumber\\
	    		&+\s{\H}{\h_0}\otimes\left(\bar{\partial}\h_0\totimes\h_0\right)+\frac{1}{2}\,g\otimes\s{\H}{\partial\h_0}\otimes\s{\H}{\h_0}.
	    		\end{align}
	    		Putting together \eqref{altcod1} and \eqref{altcod2}, we have
	    		\begin{align}\label{del111}
	    		&\partial^N\bar{\partial}^N\h_0\totimes\partial^N\h_0=\partial\bar{\partial}\h_0\totimes\partial\h_0-\bigg(\frac{1}{2}\,\left(g^{-1}\otimes\left(\partial\h_0\totimes\bar{\h_0}+\h_0\totimes\bar{\bar{\partial}\h_0}\right)+\s{\H}{\bar{\partial}\h_0}\right)\otimes\left(\h_0\totimes\h_0\right)\nonumber\\
	    		&+\s{\H}{\h_0}\totimes\left(\bar{\partial}\h_0\totimes\h_0\right)+\ccancel{\frac{1}{2}\,g\otimes\s{\H}{\partial\h_0}\otimes\s{\H}{\h_0}}\bigg)+\frac{1}{2}\,g\otimes\left(|\h_0|^2_{WP}\,\partial\h_0\totimes\h_0+\ccancel{\s{\H}{\partial\h_0}\otimes\s{\H}{\h_0}}\right)\nonumber\\
	    		&=\partial\bar{\partial}\h_0\totimes\partial\h_0+\frac{1}{2}\,g^{-1}\otimes\left\{\left(\partial\h_0\totimes\h_0\right)\otimes\left(\h_0\totimes\bar{\h_0}\right)-\left(\partial\h_0\totimes\bar{\h_0}\right)\otimes\left(\h_0\totimes\h_0\right)-\left(\partial\bar{\h_0}\totimes\h_0\right)\otimes\left(\h_0\totimes\h_0\right)\right\}\nonumber\\
	    		&-\frac{1}{2}\s{\H}{\bar{\partial}\h_0}\otimes\left({\h_0\totimes\h_0}\right)-\s{\H}{\h_0}\otimes\left(\bar{\partial}\h_0\totimes\h_0\right)\nonumber\\
	    		&=\partial\bar{\partial}\h_0\totimes\partial\h_0+\frac{1}{2}\,g^{-1}\otimes\left\{\left(\partial\h_0\totimes\h_0\right)\otimes\left(\h_0\totimes\bar{\h_0}\right)-\left(\partial\h_0\totimes\bar{\h_0}\right)\otimes\left(\h_0\totimes\h_0\right)-\left(\partial\bar{\h_0}\totimes\h_0\right)\otimes\left(\h_0\totimes\h_0\right)\right\}\nonumber\\
	    		&-\frac{1}{2}\bs{\H}{\bar{\partial}\left(\left(\h_0\totimes\h_0\right)\h_0\right)}.
	    		\end{align}
	    		Finally, we have by \eqref{altcod1} (indeed, we just need to virtually replace $\partial\h_0$ by $\bar{\partial}\h_0$)
	    		\begin{align}\label{endcod1}
	    		\partial^N\bar{\partial}^N\h_0\totimes\bar{\partial}^N\h_0=\partial\bar{\partial}\h_0\totimes\bar{\partial}\h_0-\partial^\top\bar{\partial}\h_0\totimes\bar{\partial}^\top\h_0+\frac{1}{2}\,g\otimes\left(|\h_0|^2_{WP}\,\bar{\partial}\h_0\totimes\h_0+\s{\H}{\bar{\partial}\h_0}\otimes\s{\H}{\h_0}\right).
	    		\end{align}
	    		Now, by \eqref{codazzi0} and \eqref{tangentend}, we have
	    		\begin{align}\label{endcod2}
	    		&\partial^\top\bar{\partial}\h_0\totimes\bar{\partial}^\top\h_0=
	    		\bigg\{-\left(g^{-1}\otimes\left(\partial\h_0\totimes\bar{\h_0}+\h_0\totimes\bar{\bar{\partial}\h_0}\right)+\s{\H}{\bar{\partial}\h_0}\right)\partial\phi\nonumber\\
	    		&-\left(2\,g^{-1}\otimes\left(\bar{\partial}\h_0\totimes\h_0\right)+\s{\H}{\partial\h_0}\right)\bar{\partial}\phi\bigg\}\totimes
	    		\bigg\{-|\h_0|^2_{WP}\;g\otimes \partial\phi-\s{\H}{\h_0}\otimes \bar{\partial}\phi\bigg\}\nonumber\\
	    		&=\frac{1}{2}\,g^{-1}\otimes\s{\H}{\h_0}\otimes\left(\partial\h_0\totimes\bar{\h_0}+\partial\bar{\h_0}\totimes\h_0\right)+\frac{1}{2}\,g\otimes\s{\H}{\bar{\partial}\h_0}\otimes\s{\H}{\h_0}\nonumber\\
	    		&+g^{-1}\,\otimes\left(\bar{\partial}\h_0\totimes\h_0\right)\otimes\left(\h_0\totimes\bar{\h_0}\right)+\frac{1}{2}\s{\H}{\partial\h_0}\otimes\left(\h_0\totimes\bar{\h_0}\right)
	    		\end{align}
	    		as $|\h_0|^2_{WP}\,g^2=\h_0\totimes\bar{\h_0}$. Finally, by \eqref{endcod1} and \eqref{endcod2}, we have
	    		\begin{align}\label{endcod3}
	    		&\partial^N\bar{\partial}^N\h_0\totimes\bar{\partial}^N\h_0=\partial\bar{\partial}\h_0\totimes\bar{\partial}\h_0-\bigg\{\frac{1}{2}\,g^{-1}\otimes\s{\H}{\h_0}\otimes\left(\partial\h_0\totimes\bar{\h_0}+\partial\bar{\h_0}\totimes\h_0\right)+\ccancel{\frac{1}{2}\,g\otimes\s{\H}{\bar{\partial}\h_0}\otimes\s{\H}{\h_0}}\nonumber\nonumber\\
	    		&+g^{-1}\,\otimes\left(\bar{\partial}\h_0\totimes\h_0\right)\otimes\left(\h_0\totimes\bar{\h_0}\right)+\frac{1}{2}\s{\H}{\partial\h_0}\otimes\left(\h_0\totimes\bar{\h_0}\right)\bigg\}+\frac{1}{2}\,g\otimes\left(|\h_0|^2_{WP}\,\bar{\partial}\h_0\totimes\h_0+\ccancel{\s{\H}{\bar{\partial}\h_0}\otimes\s{\H}{\h_0}}\right)\nonumber\\
	    		&=\partial\bar{\partial}\h_0\totimes\bar{\partial}\h_0-\frac{1}{2}\,g^{-1}\otimes\left(\bar{\partial}\h_0\totimes\h_0\right)\otimes\left(\h_0\totimes\bar{\h_0}\right)-\frac{1}{2}\bs{\H}{\partial\left(\left(\h_0\totimes\bar{\h_0}\right)\h_0\right)}
	    		\end{align}
	    		as
	    		\begin{align*}
	    		&g\otimes\left(|\h_0|^2_{WP}\,\bar{\partial}\h_0\totimes\h_0\right)=g^{-1}\,\otimes\left(\bar{\partial}\h_0\totimes\h_0\right)\otimes\left(\h_0\totimes\bar{\h_0}\right)\\
	    		&\frac{1}{2}\,g^{-1}\otimes\s{\H}{\h_0}\otimes\left(\partial\h_0\totimes\bar{\h_0}+\partial\bar{\h_0}\totimes\h_0\right)+\frac{1}{2}\s{\H}{\partial\h_0}\otimes\left(\h_0\totimes\bar{\h_0}\right)=\frac{1}{2}\bs{\H}{\partial\left(\left(\h_0\totimes\bar{\h_0}\right)\h_0\right)}
	    		\end{align*}
	    		We see that in $\mathscr{O}_{\phi}$, the part that we need to express without normal derivatives is
	    		\begin{align}\label{opart}
	    		\widetilde{\mathscr{O}}_{\phi}&=g^{-2}\otimes\bigg\{\frac{1}{4}(\partial^N\bar{\partial}^N\h_0\,\dot{\otimes}\,\partial^N\bar{\partial}^N\h_0)\otimes (\h_0\,\dot{\otimes}\,\h_0)+\frac{1}{4}(\partial^N\h_0\,\dot{\otimes}\, \partial^N\h_0)\otimes (\bar{\partial}^N\h_0\,\dot{\otimes}\, \bar{\partial}^N\h_0)\nonumber\\
	    		&-\frac{1}{2}(\partial^N\bar{\partial}^N\h_0\totimes\partial^N\h_0)\otimes(\bar{\partial}^N\h_0\totimes\h_0)-\frac{1}{2}(\partial^N\bar{\partial}^N\h_0\totimes\bar{\partial}^N\h_0)\otimes(\partial^N\h_0\totimes\h_0)\nonumber
	    		\\
	    		&+\frac{1}{2}(\partial^N\bar{\partial}^N\h_0\totimes\h_0)\otimes(\partial^N\h_0\totimes\bar{\partial}^N\h_0)\bigg\}\nonumber\\
	    		&=g^{-2}\otimes\bigg\{\frac{1}{4}\mathrm{(I)}+\frac{1}{4}\mathrm{(II)}-\frac{1}{2}\mathrm{(III)}-\frac{1}{2}\mathrm{(IV)}+\frac{1}{2}\mathrm{(V)}\bigg\}.
	    		\end{align}
	    		We first have by \eqref{deldelbar2}
	    		\begin{align}\label{oI}
	    		&\mathrm{(I)}=\left(\partial^N\bar{\partial}^N\h_0\totimes\partial^N\bar{\partial}^N\h_0\right)\otimes\left(\h_0\totimes\h_0\right)=\left(\partial\bar{\partial}\h_0\totimes\partial\bar{\partial}\h_0\right)\otimes\left(\h_0\totimes\h_0\right)
	    		+g^{-1}\otimes \bigg(\left(\partial\bar{\partial}\h_0\totimes\h_0\right)\otimes\left(\h_0\otimes\bar{\h_0}\right)\nonumber\\
	    		&-2\left(\partial\h_0\totimes\bar{\h_0}\right)\otimes\left(\bar{\partial}\h_0\totimes\h_0\right)
	    		-2\left(\partial\bar{\h_0}\totimes\h_0\right)\otimes\left(\bar{\partial}\h_0\totimes\h_0\right)\bigg)\otimes\left(\h_0\totimes\h_0\right)\nonumber
	    		\\
	    		&-\bigg(2\s{\H}{\bar{\partial}\h_0}\otimes\left(\bar{\partial}\h_0\totimes\h_0\right)+\s{\H}{\partial\h_0}\otimes\left(\partial\h_0\totimes\bar{\h_0}\right)+\s{\H}{\partial\h_0}\otimes\left(\partial\bar{\h_0}\totimes\h_0\right)\bigg)\otimes\left(\h_0\totimes\h_0\right)\nonumber\\
	    		&+g\otimes\bigg(\s{\H}{\partial\bar{\partial}\h_0}\otimes\s{\H}{\h_0}-\s{\H}{\partial\h_0}\otimes\s{\H}{\bar{\partial}\h_0}\bigg)\otimes\left(\h_0\totimes\h_0\right)\nonumber\\
	    		&+\frac{1}{4}\,g^2\otimes\left(|\h_0|_{WP}^4\,\h_0\totimes\h_0+|\H|^2\s{\H}{\h_0}+2|\h_0|_{WP}^2\s{\H}{\h_0}^2\right)\otimes\left(\h_0\totimes\h_0\right).
	    		\end{align}
	    		Therefore, we have
	    		\begin{align}\label{ooI}
	    		&\frac{1}{4}\,g^{-2}\otimes \mathrm{(I)}=\frac{1}{4}\,g^{-2}\otimes\left(\partial\bar{\partial}\h_0\otimes\partial\bar{\partial}\h_0\right)\otimes\left(\h_0\totimes\h_0\right)+\frac{1}{4}|\h_0|^{2}_{WP}\,g^{-1}\otimes\left(\partial\bar{\partial}\h_0\otimes\h_0\right)\otimes\left(\h_0\totimes\h_0\right)\nonumber\\
	    		&-\frac{1}{2}\,g^{-3}\otimes \partial\left(\h_0\totimes\bar{\h_0}\right)\otimes\left(\bar{\partial}\h_0\totimes\h_0\right)\otimes\left(\h_0\totimes\h_0\right)\nonumber\\
	    		&-\frac{1}{4}\,g^{-2}\otimes\bigg\{2\s{\H}{\bar{\partial}\h_0}\otimes\left(\bar{\partial}\h_0\otimes\h_0\right)+\s{\H}{\partial\h_0}\otimes\partial\left(\h_0\totimes\bar{\h_0}\right)\bigg\}\otimes\left(\h_0\totimes\h_0\right)\nonumber\\
	    		&+\frac{1}{4}\,g^{-1}\otimes\bigg(\s{\H}{\partial\bar{\partial}\h_0}\otimes\s{\H}{\h_0}-\s{\H}{\partial\h_0}\otimes\s{\H}{\bar{\partial}\h_0}\bigg)\otimes\left(\h_0\totimes\h_0\right)\nonumber\\
	    		&+\frac{1}{16}\,\left\{\left(|\h_0|_{WP}^2\,\h_0\totimes\h_0\right)^2+|\H|^2\left(\h_0\totimes\h_0\right)\otimes\s{\H}{\h_0}+2|\h_0|_{WP}^2\left(\h_0\totimes\h_0\right)\otimes \s{\H}{\h_0}^2\right\}
	    		\end{align}
	    		Then, by \eqref{codazzi1}, we have
	    		\begin{align}\label{newcodazzi1}
	    		\left\{\begin{alignedat}{1}
	    		\partial^N\h_0\totimes\bar{\partial}^N\h_0&=\partial\h_0\totimes\bar{\partial}\h_0-\frac{1}{2}\,g\otimes\left(|\h_0|^2_{WP}\,\h_0\totimes\h_0+\s{\H}{\h_0}^2\right)\\
	    		\partial^N\h_0\totimes\partial^N\h_0&=\partial\h_0\totimes\partial\h_0-\s{\H}{\h_0}\otimes\left(\h_0\totimes\h_0\right)\\
	    		\bar{\partial}^N\h_0\totimes\bar{\partial}^N\h_0&=\bar{\partial}\h_0\totimes\bar{\partial}\h_0-\s{\H}{\h_0}\otimes \left(\h_0\totimes\bar{\h_0}\right).
	    		\end{alignedat}\right.
	    		\end{align}
	    		Therefore, we have
	    		\begin{align}\label{oII}
	    		&\mathrm{(II)}=\left(\partial^N\h_0\totimes\partial^N\h_0\right)\otimes\left(\bar{\partial}^N\h_0\totimes\bar{\partial}^N\h_0\right)=\left(\partial\h_0\totimes{\partial}\h_0\right)\otimes\left(\bar{\partial}\h_0\totimes\bar{\partial}\h_0\right)\nonumber\\
	    		&-\s{\H}{\h_0}\otimes\left\{\left(\partial\h_0\totimes\partial\h_0\right)\otimes\left(\h_0\totimes\bar{\h_0}\right)+\left(\bar{\partial}\h_0\totimes\bar{\partial}\h_0\right)\otimes\left(\h_0\totimes\h_0\right)\right\}
	    		+\s{\H}{\h_0}^2\otimes\left(\h_0\totimes\h_0\right)\otimes\left(\h_0\totimes\bar{\h_0}\right)
	    		\end{align}
	    		and
	    		\begin{align}\label{ooII}
	    		\frac{1}{4}\,g^{-2}\otimes\mathrm{(II)}&=\frac{1}{4}\,g^{-2}\otimes\left(\partial\h_0\totimes\partial\h_0\right)\otimes\left(\bar{\partial}\h_0\totimes\bar{\partial}\h_0\right)-\frac{1}{4}|\h_0|^{2}_{WP}\s{\H}{\h_0}\otimes\left(\partial\h_0\totimes\partial\h_0\right)\nonumber\\
	    		&-\frac{1}{4}\,g^{-2}\otimes \s{\H}{\h_0}\otimes\left(\bar{\partial}\h_0\totimes\bar{\partial}\h_0\right)\otimes\left(\h_0\totimes\h_0\right)
	    		+\frac{1}{4}|\h_0|^2_{WP}\left(\h_0\totimes\h_0\right)\otimes\s{\H}{\h_0}^2
	    		\end{align}
	    		Then, we have by \eqref{del111}
	    		\begin{align}\label{oIII}
	    		&\mathrm{(III)}=(\partial^N\bar{\partial}^N\h_0\totimes\partial^N\h_0)\otimes(\bar{\partial}^N\h_0\totimes\h_0)=\left(\partial\bar{\partial}\h_0\totimes\partial\h_0\right)\otimes\left(\bar{\partial}\h_0\totimes\h_0\right)\nonumber\\
	    		&+\frac{1}{2}\,g^{-1}\otimes\bigg\{\left(\partial\h_0\totimes\h_0\right)\otimes\left(\h_0\totimes\bar{\h_0}\right)-\left(\partial\h_0\totimes\bar{\h_0}\right)\otimes\left(\h_0\totimes\h_0\right)\nonumber\\
	    		&-\left(\partial\bar{\h_0}\totimes\h_0\right)\otimes\left(\h_0\totimes\h_0\right)\bigg\}\otimes\left(\bar{\partial}\h_0\totimes\h_0\right)
	    		-\frac{1}{2}\bs{\H}{\bar{\partial}\left(\left(\h_0\totimes\h_0\right)\h_0\right)}\otimes\left(\bar{\partial}\h_0\totimes\h_0\right)
	    		\end{align}
	    		so
	    		\begin{align}\label{ooIII}
	    		&-\frac{1}{2}\,g^{-2}\otimes\mathrm{(III)}=-\frac{1}{2}\,g^{-2}\otimes\left(\partial\bar{\partial}\h_0\otimes\partial\h_0\right)\totimes\left(\bar{\partial}\h_0\totimes\h_0\right)-\frac{1}{4}|\h_0|^2_{WP}\,g^{-1}\otimes \left(\partial\h_0\totimes\h_0\right)\otimes\left(\bar{\partial}\h_0\totimes\h_0\right)\nonumber\\
	    		&+\frac{1}{4}\,g^{-3}\otimes\partial\left(\h_0\totimes\bar{\h_0}\right)\otimes\left(\bar{\partial}\h_0\totimes\h_0\right)\otimes\left(\h_0\totimes\h_0\right)+\frac{1}{4}\,g^{-2}\otimes\bs{\H}{\bar{\partial}\left(\left(\h_0\totimes\h_0\right)\h_0\right)}\otimes\left(\bar{\partial}\h_0\totimes\h_0\right).
	    		\end{align}
	    		Now, we easily have by \eqref{endcod3}
	    		\begin{align}\label{oIV}
	    		&\mathrm{(IV)}=\left(\partial^N\bar{\partial}^N\h_0\totimes\bar{\partial}^N\h_0\right)\totimes\left({\partial}^N\h_0\totimes\h_0\right)=\left(\partial^N\bar{\partial}^N\h_0\totimes\bar{\partial}^N\h_0\right)\totimes\left({\partial}\h_0\totimes\h_0\right)\nonumber\\
	    		&=\left(\partial\bar{\partial}\h_0\totimes\bar{\partial}\h_0\right)\otimes \left({\partial}\h_0\totimes\h_0\right)-\frac{1}{2}\,g^{-1}\otimes\left({\partial}\h_0\totimes\h_0\right)\otimes\left(\bar{\partial}\h_0\totimes\h_0\right)\otimes\left(\h_0\totimes\bar{\h_0}\right)\nonumber\\
	    		&-\frac{1}{2}\bs{\H}{\partial\left(\left(\h_0\totimes\bar{\h_0}\right)\h_0\right)}\otimes \left({\partial}\h_0\totimes\h_0\right)
	    		\end{align}
	    		so
	    		\begin{align}\label{ooIV}
	    		-\frac{1}{2}\,g^{-2}\otimes \mathrm{(IV)}&=-\frac{1}{2}\,g^{-2}\otimes\left(\partial\bar{\partial}\h_0\totimes\bar{\partial}\h_0\right)\otimes \left({\partial}\h_0\totimes\h_0\right)+\frac{1}{4}|\h_0|^{2}_{WP}\,g^{-1}\otimes\left({\partial}\h_0\totimes\h_0\right)\otimes\left(\bar{\partial}\h_0\totimes\h_0\right)\nonumber\\
	    		&+\frac{1}{4}\,g^{-2}\otimes\bs{\H}{\partial\left(\left(\h_0\totimes\bar{\h_0}\right)\h_0\right)}\otimes \left({\partial}\h_0\totimes\h_0\right)
	    		\end{align}
	    		As by \eqref{newcodazzi1}, we have
	    		\begin{align*}
	    		\partial^N\h_0\totimes\bar{\partial}^N\h_0=\partial\h_0\totimes\bar{\partial}\h_0-\frac{1}{2}\,g\otimes\left(|\h_0|^2_{WP}\,\h_0\totimes\h_0+\s{\H}{\h_0}^2\right)
	    		\end{align*}
	    		and by \eqref{codazzi3}
	    		\begin{align*}
	    		\partial^N\bar{\partial}^N\h_0\totimes\h_0=\partial\bar{\partial}\h_0\totimes\h_0+\frac{1}{2}\,g\otimes\left(|\h_0|^{2}_{WP}\,\h_0\totimes\h_0+\s{\H}{\h_0}^2\right),
	    		\end{align*}
	    		we readily obtain that
	    		\begin{align}\label{oV}
	    		&\mathrm{(V)}=\left(\partial^N\bar{\partial}^N\h_0\totimes\h_0\right)\otimes\left(\partial^N\h_0\totimes\bar{\partial}^N\h_0\right)=\left(\partial\bar{\partial}\h_0\totimes\h_0\right)\otimes\left(\partial\h_0\otimes\bar{\partial}\h_0\right)\nonumber\\
	    		&-\frac{1}{2}\,g\otimes\left(\partial\bar{\partial}\h_0\totimes\h_0-\partial\h_0\totimes\bar{\partial}\h_0\right)\otimes\left(|\h_0|^{2}_{WP}\,\h_0\totimes\h_0+\s{\H}{\h_0}^2\right)-\frac{1}{4}\,g^2\otimes \left(|\h_0|^{2}_{WP}\,\h_0\totimes\h_0+\s{\H}{\h_0}^2\right)^2.
	    		\end{align}
	    		and
	    		\begin{align}\label{ooV}
	    		&\frac{1}{2}\,g^{-2}\otimes\mathrm{(V)}=\frac{1}{2}\,g^{-2}\otimes\left(\partial\bar{\partial}\h_0\totimes\h_0\right)\otimes\left(\partial\h_0\totimes\bar{\partial}\h_0\right)\nonumber\\
	    		&-\frac{1}{4}\,g^{-1}\otimes\left(\partial\bar{\partial}\h_0\totimes\h_0-\partial\h_0\totimes\bar{\partial}\h_0\right)\otimes\left(|\h_0|^{2}_{WP}\,\h_0\totimes\h_0+\s{\H}{\h_0}^2\right)-\frac{1}{8} \left(|\h_0|^{2}_{WP}\,\h_0\totimes\h_0+\s{\H}{\h_0}^2\right)^2
	    		\end{align}
	    		and we see that we recover the main term of the quartic form $\mathscr{Q}_{\phi}$.
	    		Putting together \eqref{ooI}, \eqref{ooII}, \eqref{ooIII}, \eqref{ooIV}, \eqref{ooV}, we obtain
	    		\begin{align}\label{superoctic}
	    		&\widetilde{\mathscr{O}}_{\phi}=g^{-2}\otimes\bigg\{\frac{1}{4}(\partial\bar{\partial}\h_0\,\dot{\otimes}\,\partial\bar{\partial}\h_0)\otimes (\h_0\,\dot{\otimes}\,\h_0)+\frac{1}{4}(\partial\h_0\,\dot{\otimes}\, \partial\h_0)\otimes (\bar{\partial}\h_0\,\dot{\otimes}\, \bar{\partial}\h_0)\nonumber\\
	    		&-\frac{1}{2}(\partial\bar{\partial}\h_0\totimes\partial\h_0)\otimes(\bar{\partial}\h_0\totimes\h_0)-\frac{1}{2}(\partial\bar{\partial}\h_0\totimes\bar{\partial}\h_0)\otimes(\partial\h_0\totimes\h_0)+\frac{1}{2}(\partial\bar{\partial}\h_0\totimes\h_0)\otimes(\partial\h_0\totimes\bar{\partial}\h_0)\bigg\}\nonumber\\
	    		&+\frac{1}{4}|\h_0|^{2}_{WP}\,g^{-1}\otimes\left(\partial\bar{\partial}\h_0\otimes\h_0\right)\otimes\left(\h_0\totimes\h_0\right)-\frac{1}{2}\,g^{-3}\otimes \partial\left(\h_0\totimes\bar{\h_0}\right)\otimes\left(\bar{\partial}\h_0\totimes\h_0\right)\otimes\left(\h_0\totimes\h_0\right)\nonumber\\
	    		&-\frac{1}{4}\,g^{-2}\otimes\bigg\{2\s{\H}{\bar{\partial}\h_0}\otimes\left(\bar{\partial}\h_0\otimes\h_0\right)+\s{\H}{\partial\h_0}\otimes\partial\left(\h_0\totimes\bar{\h_0}\right)\bigg\}\otimes\left(\h_0\totimes\h_0\right)\nonumber\\
	    		&+\frac{1}{4}\,g^{-1}\otimes\bigg(\s{\H}{\partial\bar{\partial}\h_0}\otimes\s{\H}{\h_0}-\s{\H}{\partial\h_0}\otimes\s{\H}{\bar{\partial}\h_0}\bigg)\otimes\left(\h_0\totimes\h_0\right)\nonumber\\
	    		&+\frac{1}{16}\,\left\{\left(|\h_0|_{WP}^2\,\h_0\totimes\h_0\right)^2+|\H|^2\left(\h_0\totimes\h_0\right)\otimes\s{\H}{\h_0}+2|\h_0|_{WP}^2\left(\h_0\totimes\h_0\right)\otimes \s{\H}{\h_0}^2\right\}\nonumber\\
	    		&-\frac{1}{4}|\h_0|^{2}_{WP}\s{\H}{\h_0}\otimes\left(\partial\h_0\totimes\partial\h_0\right)-\frac{1}{4}\,g^{-2}\otimes \s{\H}{\h_0}\otimes\left(\bar{\partial}\h_0\totimes\bar{\partial}\h_0\right)\otimes\left(\h_0\totimes\h_0\right)\nonumber\\
	    		&+\frac{1}{4}|\h_0|^2_{WP}\left(\h_0\totimes\h_0\right)\otimes\s{\H}{\h_0}^2
	    		-\ccancel{\frac{1}{4}|\h_0|^2_{WP}\,g^{-1}\otimes \left(\partial\h_0\totimes\h_0\right)\otimes\left(\bar{\partial}\h_0\totimes\h_0\right)}\nonumber\\
	    		&+\frac{1}{4}\,g^{-3}\otimes\partial\left(\h_0\totimes\bar{\h_0}\right)\otimes\left(\bar{\partial}\h_0\totimes\h_0\right)\otimes\left(\h_0\totimes\h_0\right)+\frac{1}{4}\,g^{-2}\otimes\bs{\H}{\bar{\partial}\left(\left(\h_0\totimes\h_0\right)\h_0\right)}\otimes\left(\bar{\partial}\h_0\totimes\h_0\right)\nonumber\\
	    		&+\ccancel{\frac{1}{4}|\h_0|^{2}_{WP}\,g^{-1}\otimes\left({\partial}\h_0\totimes\h_0\right)\otimes\left(\bar{\partial}\h_0\totimes\h_0\right)}
	    		+\frac{1}{4}\,g^{-2}\otimes\bs{\H}{\partial\left(\left(\h_0\totimes\bar{\h_0}\right)\h_0\right)}\otimes \left({\partial}\h_0\totimes\h_0\right)\nonumber\\
	    		&-\frac{1}{4}\,g^{-1}\otimes\left(\partial\bar{\partial}\h_0\totimes\h_0-\partial\h_0\totimes\bar{\partial}\h_0\right)\otimes\left(|\h_0|^{2}_{WP}\,\h_0\totimes\h_0+\s{\H}{\h_0}^2\right)-\frac{1}{8} \left(|\h_0|^{2}_{WP}\,\h_0\totimes\h_0+\s{\H}{\h_0}^2\right)^2
	    		\end{align}
	    		Finally, we obtain by \eqref{48formula}, \eqref{removenormals}, and \eqref{superoctic}
	    		\begin{align*}
	    		&\mathscr{O}_{\phi}=\widetilde{\mathscr{O}}_{\phi}
	    		+\frac{1}{4}\left(1+|\H|^2\right)\,g^{-1}\otimes\bigg\{\frac{1}{2}\left(
	    		\partial\bar{\partial}\h_0\totimes\h_0\right)\otimes\left(\h_0\totimes\h_0\right)-\left(\partial\h_0\totimes\h_0\right)\otimes\left(\bar{\partial}\h_0\totimes\h_0\right)\\ &+\frac{1}{2}\left(\partial\h_0\totimes\bar{\partial}\h_0\right)\otimes\left(\h_0\totimes\h_0\right)\bigg\}    +\frac{1}{64}\left(1+|\H|^2\right)^2\left(\h_0\totimes\h_0\right)^2
	    		\end{align*}
	    		which concludes the proof of the proposition.  
	    	\end{proof}
	    		We recall that
	    		\begin{align*}
	    		\h_0=\begin{dmatrix}
	    		-\frac{\theta_{0} - 2}{2 \, \theta_{0}} & C_{1} & 0 & \theta_{0} \\
	    		-\frac{\theta_{0} - 2}{2 \, {\left(\theta_{0} + 1\right)}} & B_{1} & 0 & \theta_{0} + 1 \\
	    		2 \, \alpha_{2} \theta_{0} - 4 \, \alpha_{2} & A_{0} & 0 & \theta_{0} + 1 \\
	    		-\frac{\theta_{0} - 3}{2 \, \theta_{0}} & C_{2} & 1 & \theta_{0} \\
	    		4 & A_{2} & \theta_{0} & 0 \\
	    		-4 \, \alpha_{1} & A_{0} & \theta_{0} & 0 \\
	    		-4 \, \alpha_{5} & A_{0} & \theta_{0} & 1 
	    		\end{dmatrix}
	    		\begin{dmatrix}
	    		-4 \, {\left| A_{1} \right|}^{2} & A_{1} & \theta_{0} & 1 \\
	    		2 & A_{1} & \theta_{0} - 1 & 0 \\
	    		-4 \, {\left| A_{1} \right|}^{2} & A_{0} & \theta_{0} - 1 & 1 \\
	    		\frac{1}{4} & \overline{B_{1}} & \theta_{0} - 1 & 2 \\
	    		-2 \, \overline{\alpha_{5}} & A_{0} & \theta_{0} - 1 & 2 \\
	    		6 & A_{3} & \theta_{0} + 1 & 0 \\
	    		-6 \, \alpha_{3} & A_{0} & \theta_{0} + 1 & 0 \\
	    		-4 \, \alpha_{1} & A_{1} & \theta_{0} + 1 & 0 \\
	    		-2 \, \theta_{0} \overline{\alpha_{2}} - 2 \, \overline{\alpha_{2}} & A_{0} & 2 \, \theta_{0} - 1 & -\theta_{0} + 2
	    		\end{dmatrix}
	    		\end{align*}
	    		Now, we compute
	    		\footnotesize
	    		\begin{align*}
	    		&\h_0\totimes\h_0=\begin{dmatrix}
	    		\frac{{\left(\theta_{0}^{2} - 4 \, \theta_{0} + 4\right)} C_{1}^{2}}{4 \, \theta_{0}^{2}} & 0 & 2 \, \theta_{0} \\
	    		-\frac{2 \, {\left(A_{1} C_{2} {\left(\theta_{0} - 3\right)} - 2 \, {\left({\left(\alpha_{1} \theta_{0} - 2 \, \alpha_{1}\right)} A_{0} - A_{2} {\left(\theta_{0} - 2\right)}\right)} C_{1}\right)}}{\theta_{0}} & \theta_{0} & \theta_{0} \\
	    		-\frac{2 \, A_{1} C_{1} {\left(\theta_{0} - 2\right)}}{\theta_{0}} & \theta_{0} - 1 & \theta_{0} \\
	    		\frac{2 \, {\left(2 \, {\left(\theta_{0}^{2} - \theta_{0} - 2\right)} A_{0} C_{1} {\left| A_{1} \right|}^{2} + 4 \, {\left(\alpha_{2} \theta_{0}^{3} - \alpha_{2} \theta_{0}^{2} - 2 \, \alpha_{2} \theta_{0}\right)} A_{0} A_{1} - {\left(\theta_{0}^{2} - 2 \, \theta_{0}\right)} A_{1} B_{1}\right)}}{\theta_{0}^{2} + \theta_{0}} & \theta_{0} - 1 & \theta_{0} + 1 \\
	    		4 \, A_{1}^{2} & 2 \, \theta_{0} - 2 & 0 \\
	    		-16 \, A_{0} A_{1} {\left| A_{1} \right|}^{2} & 2 \, \theta_{0} - 2 & 1 \\
	    		16 \, A_{0}^{2} {\left| A_{1} \right|}^{4} - 8 \, A_{0} A_{1} \overline{\alpha_{5}} + A_{1} \overline{B_{1}} & 2 \, \theta_{0} - 2 & 2 \\
	    		-16 \, A_{0} A_{1} \alpha_{1} + 16 \, A_{1} A_{2} & 2 \, \theta_{0} - 1 & 0 \\
	    		-32 \, A_{0} A_{2} {\left| A_{1} \right|}^{2} - 16 \, A_{0} A_{1} \alpha_{5} + 16 \, {\left(2 \, A_{0}^{2} \alpha_{1} - A_{1}^{2}\right)} {\left| A_{1} \right|}^{2} & 2 \, \theta_{0} - 1 & 1 \\
	    		-8 \, {\left(\theta_{0} \overline{\alpha_{2}} + \overline{\alpha_{2}}\right)} A_{0} A_{1} & 3 \, \theta_{0} - 2 & -\theta_{0} + 2 \\
	    		16 \, A_{0}^{2} \alpha_{1}^{2} - 16 \, A_{1}^{2} \alpha_{1} - 32 \, A_{0} A_{2} \alpha_{1} - 24 \, A_{0} A_{1} \alpha_{3} + 16 \, A_{2}^{2} + 24 \, A_{1} A_{3} & 2 \, \theta_{0} & 0
	    		\end{dmatrix}\\
	    		&\partial\h_0\totimes\h_0=\\
	    		&\begin{dmatrix}
	    		-\frac{{\left(\theta_{0}^{2} - 3 \, \theta_{0} + 2\right)} A_{1} C_{1}}{\theta_{0}} & \theta_{0} - 2 & \theta_{0} \\
	    		\frac{2 \, {\left(\theta_{0}^{3} - 2 \, \theta_{0}^{2} - \theta_{0} + 2\right)} A_{0} C_{1} {\left| A_{1} \right|}^{2} + 4 \, {\left(\alpha_{2} \theta_{0}^{4} - 2 \, \alpha_{2} \theta_{0}^{3} - \alpha_{2} \theta_{0}^{2} + 2 \, \alpha_{2} \theta_{0}\right)} A_{0} A_{1} - {\left(\theta_{0}^{3} - 3 \, \theta_{0}^{2} + 2 \, \theta_{0}\right)} A_{1} B_{1}}{\theta_{0}^{2} + \theta_{0}} & \theta_{0} - 2 & \theta_{0} + 1 \\
	    		-A_{1} C_{2} {\left(\theta_{0} - 3\right)} + 2 \, {\left({\left(\alpha_{1} \theta_{0} - 2 \, \alpha_{1}\right)} A_{0} - A_{2} {\left(\theta_{0} - 2\right)}\right)} C_{1} & \theta_{0} - 1 & \theta_{0} \\
	    		4 \, A_{1}^{2} {\left(\theta_{0} - 1\right)} & 2 \, \theta_{0} - 3 & 0 \\
	    		-16 \, A_{0} A_{1} {\left(\theta_{0} - 1\right)} {\left| A_{1} \right|}^{2} & 2 \, \theta_{0} - 3 & 1 \\
	    		16 \, A_{0}^{2} {\left(\theta_{0} - 1\right)} {\left| A_{1} \right|}^{4} - 8 \, {\left(\theta_{0} \overline{\alpha_{5}} - \overline{\alpha_{5}}\right)} A_{0} A_{1} + A_{1} {\left(\theta_{0} - 1\right)} \overline{B_{1}} & 2 \, \theta_{0} - 3 & 2 \\
	    		-8 \, {\left(2 \, \alpha_{1} \theta_{0} - \alpha_{1}\right)} A_{0} A_{1} + 8 \, A_{1} A_{2} {\left(2 \, \theta_{0} - 1\right)} & 2 \, \theta_{0} - 2 & 0 \\
	    		-16 \, A_{0} A_{2} {\left(2 \, \theta_{0} - 1\right)} {\left| A_{1} \right|}^{2} - 8 \, {\left(2 \, \alpha_{5} \theta_{0} - \alpha_{5}\right)} A_{0} A_{1} + 8 \, {\left(2 \, {\left(2 \, \alpha_{1} \theta_{0} - \alpha_{1}\right)} A_{0}^{2} - A_{1}^{2} {\left(2 \, \theta_{0} - 1\right)}\right)} {\left| A_{1} \right|}^{2} & 2 \, \theta_{0} - 2 & 1 \\
	    		16 \, A_{0}^{2} \alpha_{1}^{2} \theta_{0} - 16 \, A_{1}^{2} \alpha_{1} \theta_{0} - 32 \, A_{0} A_{2} \alpha_{1} \theta_{0} - 24 \, A_{0} A_{1} \alpha_{3} \theta_{0} + 16 \, A_{2}^{2} \theta_{0} + 24 \, A_{1} A_{3} \theta_{0} & 2 \, \theta_{0} - 1 & 0 \\
	    		-4 \, {\left(3 \, \theta_{0}^{2} \overline{\alpha_{2}} + \theta_{0} \overline{\alpha_{2}} - 2 \, \overline{\alpha_{2}}\right)} A_{0} A_{1} & 3 \, \theta_{0} - 3 & -\theta_{0} + 2
	    		\end{dmatrix}\\
	    		&\bar{\partial}\h_0\totimes \h_0=\\
	    		&\begin{dmatrix}
	    		\frac{{\left(\theta_{0}^{2} - 4 \, \theta_{0} + 4\right)} C_{1}^{2}}{4 \, \theta_{0}} & 0 & 2 \, \theta_{0} - 1 \\
	    		-A_{1} C_{2} {\left(\theta_{0} - 3\right)} + 2 \, {\left({\left(\alpha_{1} \theta_{0} - 2 \, \alpha_{1}\right)} A_{0} - A_{2} {\left(\theta_{0} - 2\right)}\right)} C_{1} & \theta_{0} & \theta_{0} - 1 \\
	    		\frac{2 \, {\left(\theta_{0}^{2} - \theta_{0} - 2\right)} A_{0} C_{1} {\left| A_{1} \right|}^{2} + 4 \, {\left(\alpha_{2} \theta_{0}^{3} - \alpha_{2} \theta_{0}^{2} - 2 \, \alpha_{2} \theta_{0}\right)} A_{0} A_{1} - {\left(\theta_{0}^{2} - 2 \, \theta_{0}\right)} A_{1} B_{1}}{\theta_{0}} & \theta_{0} - 1 & \theta_{0} \\
	    		-A_{1} C_{1} {\left(\theta_{0} - 2\right)} & \theta_{0} - 1 & \theta_{0} - 1 \\
	    		-8 \, A_{0} A_{1} {\left| A_{1} \right|}^{2} & 2 \, \theta_{0} - 2 & 0 \\
	    		16 \, A_{0}^{2} {\left| A_{1} \right|}^{4} - 8 \, A_{0} A_{1} \overline{\alpha_{5}} + A_{1} \overline{B_{1}} & 2 \, \theta_{0} - 2 & 1 \\
	    		-16 \, A_{0} A_{2} {\left| A_{1} \right|}^{2} - 8 \, A_{0} A_{1} \alpha_{5} + 8 \, {\left(2 \, A_{0}^{2} \alpha_{1} - A_{1}^{2}\right)} {\left| A_{1} \right|}^{2} & 2 \, \theta_{0} - 1 & 0 \\
	    		4 \, {\left(\theta_{0}^{2} \overline{\alpha_{2}} - \theta_{0} \overline{\alpha_{2}} - 2 \, \overline{\alpha_{2}}\right)} A_{0} A_{1} & 3 \, \theta_{0} - 2 & -\theta_{0} + 1
	    		\end{dmatrix}\\
	    		&\partial\h_0\totimes\partial\h_0=\\
	    		&\begin{dmatrix}
	    		-\frac{2 \, {\left(\theta_{0}^{2} - 4 \, \theta_{0} + 3\right)} A_{1} C_{2}}{\theta_{0}} & \theta_{0} - 2 & \theta_{0} \\
	    		4 \, {\left(\theta_{0}^{2} - 2 \, \theta_{0} + 1\right)} A_{1}^{2} & 2 \, \theta_{0} - 4 & 0 \\
	    		-16 \, {\left(\theta_{0}^{2} - 2 \, \theta_{0} + 1\right)} A_{0} A_{1} {\left| A_{1} \right|}^{2} & 2 \, \theta_{0} - 4 & 1 \\
	    		16 \, {\left(\theta_{0}^{2} - 2 \, \theta_{0} + 1\right)} A_{0}^{2} {\left| A_{1} \right|}^{4} - 8 \, {\left(\theta_{0}^{2} \overline{\alpha_{5}} - 2 \, \theta_{0} \overline{\alpha_{5}} + \overline{\alpha_{5}}\right)} A_{0} A_{1} + {\left(\theta_{0}^{2} - 2 \, \theta_{0} + 1\right)} A_{1} \overline{B_{1}} & 2 \, \theta_{0} - 4 & 2 \\
	    		-16 \, {\left(\alpha_{1} \theta_{0}^{2} - \alpha_{1} \theta_{0}\right)} A_{0} A_{1} + 16 \, {\left(\theta_{0}^{2} - \theta_{0}\right)} A_{1} A_{2} & 2 \, \theta_{0} - 3 & 0 \\
	    		-32 \, {\left(\theta_{0}^{2} - \theta_{0}\right)} A_{0} A_{2} {\left| A_{1} \right|}^{2} - 16 \, {\left(\alpha_{5} \theta_{0}^{2} - \alpha_{5} \theta_{0}\right)} A_{0} A_{1} + 16 \, {\left(2 \, {\left(\alpha_{1} \theta_{0}^{2} - \alpha_{1} \theta_{0}\right)} A_{0}^{2} - {\left(\theta_{0}^{2} - \theta_{0}\right)} A_{1}^{2}\right)} {\left| A_{1} \right|}^{2} & 2 \, \theta_{0} - 3 & 1 \\
	    		16 \, A_{0}^{2} \alpha_{1}^{2} \theta_{0}^{2} - 32 \, A_{0} A_{2} \alpha_{1} \theta_{0}^{2} + 16 \, A_{2}^{2} \theta_{0}^{2} - 24 \, {\left(\alpha_{3} \theta_{0}^{2} - \alpha_{3}\right)} A_{0} A_{1} - 16 \, {\left(\alpha_{1} \theta_{0}^{2} - \alpha_{1}\right)} A_{1}^{2} + 24 \, {\left(\theta_{0}^{2} - 1\right)} A_{1} A_{3} & 2 \, \theta_{0} - 2 & 0 \\
	    		-8 \, {\left(2 \, \theta_{0}^{3} \overline{\alpha_{2}} - \theta_{0}^{2} \overline{\alpha_{2}} - 2 \, \theta_{0} \overline{\alpha_{2}} + \overline{\alpha_{2}}\right)} A_{0} A_{1} & 3 \, \theta_{0} - 4 & -\theta_{0} + 2
	    		\end{dmatrix}\\
	    		&\bar{\partial}\h_0\totimes\bar{\partial}\h_0=\\
	    		&\begin{dmatrix}
	    		\frac{1}{4} \, {\left(\theta_{0}^{2} - 4 \, \theta_{0} + 4\right)} C_{1}^{2} & 0 & 2 \, \theta_{0} - 2 \\
	    		-\frac{1}{2} \, {\left(4 \, {\left(\alpha_{2} \theta_{0}^{3} - 3 \, \alpha_{2} \theta_{0}^{2} + 4 \, \alpha_{2}\right)} A_{0} - {\left(\theta_{0}^{2} - 4 \, \theta_{0} + 4\right)} B_{1}\right)} C_{1} & 0 & 2 \, \theta_{0} - 1 \\
	    		\frac{1}{2} \, {\left(\theta_{0}^{2} - 5 \, \theta_{0} + 6\right)} C_{1} C_{2} & 1 & 2 \, \theta_{0} - 2 \\
	    		4 \, A_{0} C_{2} {\left(\theta_{0} - 3\right)} {\left| A_{1} \right|}^{2} + 4 \, {\left(A_{1} {\left(\theta_{0} - 2\right)} {\left| A_{1} \right|}^{2} + {\left(\alpha_{5} \theta_{0} - 2 \, \alpha_{5}\right)} A_{0}\right)} C_{1} & \theta_{0} & \theta_{0} - 1 \\
	    		-16 \, {\left(\alpha_{2} \theta_{0}^{2} - \alpha_{2} \theta_{0} - 2 \, \alpha_{2}\right)} A_{0}^{2} {\left| A_{1} \right|}^{2} + 4 \, A_{0} B_{1} {\left(\theta_{0} - 2\right)} {\left| A_{1} \right|}^{2} + \frac{1}{2} \, {\left(8 \, {\left(\theta_{0} \overline{\alpha_{5}} - 2 \, \overline{\alpha_{5}}\right)} A_{0} - {\left(\theta_{0} - 2\right)} \overline{B_{1}}\right)} C_{1} & \theta_{0} - 1 & \theta_{0} \\
	    		4 \, A_{0} C_{1} {\left(\theta_{0} - 2\right)} {\left| A_{1} \right|}^{2} & \theta_{0} - 1 & \theta_{0} - 1 \\
	    		16 \, A_{0}^{2} {\left| A_{1} \right|}^{4} & 2 \, \theta_{0} - 2 & 0 \\
	    		32 \, A_{0}^{2} {\left| A_{1} \right|}^{2} \overline{\alpha_{5}} - 4 \, A_{0} {\left| A_{1} \right|}^{2} \overline{B_{1}} & 2 \, \theta_{0} - 2 & 1 \\
	    		32 \, A_{0} A_{1} {\left| A_{1} \right|}^{4} + 32 \, A_{0}^{2} \alpha_{5} {\left| A_{1} \right|}^{2} - 2 \, {\left(\theta_{0}^{3} \overline{\alpha_{2}} - 3 \, \theta_{0}^{2} \overline{\alpha_{2}} + 4 \, \overline{\alpha_{2}}\right)} A_{0} C_{1} & 2 \, \theta_{0} - 1 & 0 \\
	    		-16 \, {\left(\theta_{0}^{2} \overline{\alpha_{2}} - \theta_{0} \overline{\alpha_{2}} - 2 \, \overline{\alpha_{2}}\right)} A_{0}^{2} {\left| A_{1} \right|}^{2} & 3 \, \theta_{0} - 2 & -\theta_{0} + 1
	    		\end{dmatrix}\\
	    		&\partial\h_0\totimes\bar{\partial}\h_0=\\
	    		&\begin{dmatrix}
	    		2 \, {\left(\theta_{0}^{2} - 3 \, \theta_{0} + 2\right)} A_{0} C_{1} {\left| A_{1} \right|}^{2} + 4 \, {\left(\alpha_{2} \theta_{0}^{3} - 2 \, \alpha_{2} \theta_{0}^{2} - \alpha_{2} \theta_{0} + 2 \, \alpha_{2}\right)} A_{0} A_{1} - {\left(\theta_{0}^{2} - 3 \, \theta_{0} + 2\right)} A_{1} B_{1} & \theta_{0} - 2 & \theta_{0} \\
	    		-{\left(\theta_{0}^{2} - 3 \, \theta_{0} + 2\right)} A_{1} C_{1} & \theta_{0} - 2 & \theta_{0} - 1 \\
	    		-{\left(\theta_{0}^{2} - 4 \, \theta_{0} + 3\right)} A_{1} C_{2} + 2 \, {\left({\left(\alpha_{1} \theta_{0}^{2} - 2 \, \alpha_{1} \theta_{0}\right)} A_{0} - {\left(\theta_{0}^{2} - 2 \, \theta_{0}\right)} A_{2}\right)} C_{1} & \theta_{0} - 1 & \theta_{0} - 1 \\
	    		-8 \, A_{0} A_{1} {\left(\theta_{0} - 1\right)} {\left| A_{1} \right|}^{2} & 2 \, \theta_{0} - 3 & 0 \\
	    		16 \, A_{0}^{2} {\left(\theta_{0} - 1\right)} {\left| A_{1} \right|}^{4} - 8 \, {\left(\theta_{0} \overline{\alpha_{5}} - \overline{\alpha_{5}}\right)} A_{0} A_{1} + A_{1} {\left(\theta_{0} - 1\right)} \overline{B_{1}} & 2 \, \theta_{0} - 3 & 1 \\
	    		-16 \, A_{0} A_{2} \theta_{0} {\left| A_{1} \right|}^{2} - 8 \, {\left(\alpha_{5} \theta_{0} - \alpha_{5}\right)} A_{0} A_{1} + 8 \, {\left(2 \, A_{0}^{2} \alpha_{1} \theta_{0} - A_{1}^{2} {\left(\theta_{0} - 1\right)}\right)} {\left| A_{1} \right|}^{2} & 2 \, \theta_{0} - 2 & 0 \\
	    		4 \, {\left(\theta_{0}^{3} \overline{\alpha_{2}} - 2 \, \theta_{0}^{2} \overline{\alpha_{2}} - \theta_{0} \overline{\alpha_{2}} + 2 \, \overline{\alpha_{2}}\right)} A_{0} A_{1} & 3 \, \theta_{0} - 3 & -\theta_{0} + 1
	    		\end{dmatrix}\\
	    		&\partial\bar{\partial}\h_0\totimes\h_0=\\
	    		&\begin{dmatrix}
	    		\frac{2 \, {\left(\theta_{0}^{2} - 3 \, \theta_{0} + 2\right)} A_{0} C_{1} {\left| A_{1} \right|}^{2}}{\theta_{0}} & \theta_{0} - 2 & \theta_{0} \\
	    		-A_{1} C_{2} {\left(\theta_{0} - 3\right)} & \theta_{0} - 1 & \theta_{0} - 1 \\
	    		-8 \, A_{0} A_{1} {\left(\theta_{0} - 1\right)} {\left| A_{1} \right|}^{2} & 2 \, \theta_{0} - 3 & 0 \\
	    		16 \, A_{0}^{2} {\left(\theta_{0} - 1\right)} {\left| A_{1} \right|}^{4} - 8 \, {\left(\theta_{0} \overline{\alpha_{5}} - \overline{\alpha_{5}}\right)} A_{0} A_{1} + A_{1} {\left(\theta_{0} - 1\right)} \overline{B_{1}} & 2 \, \theta_{0} - 3 & 1 \\
	    		-16 \, A_{0} A_{2} {\left(\theta_{0} - 1\right)} {\left| A_{1} \right|}^{2} - 8 \, A_{0} A_{1} \alpha_{5} \theta_{0} + 8 \, {\left(2 \, {\left(\alpha_{1} \theta_{0} - \alpha_{1}\right)} A_{0}^{2} - A_{1}^{2} \theta_{0}\right)} {\left| A_{1} \right|}^{2} & 2 \, \theta_{0} - 2 & 0 \\
	    		4 \, {\left(2 \, \theta_{0}^{3} \overline{\alpha_{2}} - 3 \, \theta_{0}^{2} \overline{\alpha_{2}} - 3 \, \theta_{0} \overline{\alpha_{2}} + 2 \, \overline{\alpha_{2}}\right)} A_{0} A_{1} & 3 \, \theta_{0} - 3 & -\theta_{0} + 1
	    		\end{dmatrix}\\
	    		&\partial\bar{\partial}\h_0\totimes\partial\h_0=\\
	    		&\begin{dmatrix}
	    		-{\left(\theta_{0}^{2} - 4 \, \theta_{0} + 3\right)} A_{1} C_{2} & \theta_{0} - 2 & \theta_{0} - 1 \\
	    		-8 \, {\left(\theta_{0}^{2} - 2 \, \theta_{0} + 1\right)} A_{0} A_{1} {\left| A_{1} \right|}^{2} & 2 \, \theta_{0} - 4 & 0 \\
	    		16 \, {\left(\theta_{0}^{2} - 2 \, \theta_{0} + 1\right)} A_{0}^{2} {\left| A_{1} \right|}^{4} - 8 \, {\left(\theta_{0}^{2} \overline{\alpha_{5}} - 2 \, \theta_{0} \overline{\alpha_{5}} + \overline{\alpha_{5}}\right)} A_{0} A_{1} + {\left(\theta_{0}^{2} - 2 \, \theta_{0} + 1\right)} A_{1} \overline{B_{1}} & 2 \, \theta_{0} - 4 & 1 \\
	    		-16 \, {\left(\theta_{0}^{2} - \theta_{0}\right)} A_{0} A_{2} {\left| A_{1} \right|}^{2} - 8 \, {\left(\alpha_{5} \theta_{0}^{2} - \alpha_{5} \theta_{0}\right)} A_{0} A_{1} + 8 \, {\left(2 \, {\left(\alpha_{1} \theta_{0}^{2} - \alpha_{1} \theta_{0}\right)} A_{0}^{2} - {\left(\theta_{0}^{2} - \theta_{0}\right)} A_{1}^{2}\right)} {\left| A_{1} \right|}^{2} & 2 \, \theta_{0} - 3 & 0 \\
	    		4 \, {\left(2 \, \theta_{0}^{4} \overline{\alpha_{2}} - 5 \, \theta_{0}^{3} \overline{\alpha_{2}} + 5 \, \theta_{0} \overline{\alpha_{2}} - 2 \, \overline{\alpha_{2}}\right)} A_{0} A_{1} & 3 \, \theta_{0} - 4 & -\theta_{0} + 1
	    		\end{dmatrix}\\
	    		&\partial\bar{\partial}\h_0\totimes\bar{\partial}\h_0=\\
	    		&\begin{dmatrix}
	    		\frac{1}{4} \, {\left(\theta_{0}^{2} - 5 \, \theta_{0} + 6\right)} C_{1} C_{2} & 0 & 2 \, \theta_{0} - 2 \\
	    		\lambda_1 & \theta_{0} - 2 & \theta_{0} \\
	    		2 \, {\left(\theta_{0}^{2} - 3 \, \theta_{0} + 2\right)} A_{0} C_{1} {\left| A_{1} \right|}^{2} & \theta_{0} - 2 & \theta_{0} - 1 \\
	    		2 \, {\left(\theta_{0}^{2} - 3 \, \theta_{0}\right)} A_{0} C_{2} {\left| A_{1} \right|}^{2} + 2 \, {\left({\left(\theta_{0}^{2} - 2 \, \theta_{0}\right)} A_{1} {\left| A_{1} \right|}^{2} + {\left(\alpha_{5} \theta_{0}^{2} - 2 \, \alpha_{5} \theta_{0}\right)} A_{0}\right)} C_{1} & \theta_{0} - 1 & \theta_{0} - 1 \\
	    		16 \, A_{0}^{2} {\left(\theta_{0} - 1\right)} {\left| A_{1} \right|}^{4} & 2 \, \theta_{0} - 3 & 0 \\
	    		32 \, {\left(\theta_{0} \overline{\alpha_{5}} - \overline{\alpha_{5}}\right)} A_{0}^{2} {\left| A_{1} \right|}^{2} - 4 \, A_{0} {\left(\theta_{0} - 1\right)} {\left| A_{1} \right|}^{2} \overline{B_{1}} & 2 \, \theta_{0} - 3 & 1 \\
	    		16 \, A_{0} A_{1} {\left(2 \, \theta_{0} - 1\right)} {\left| A_{1} \right|}^{4} + 16 \, {\left(2 \, \alpha_{5} \theta_{0} - \alpha_{5}\right)} A_{0}^{2} {\left| A_{1} \right|}^{2} - {\left(2 \, \theta_{0}^{4} \overline{\alpha_{2}} - 7 \, \theta_{0}^{3} \overline{\alpha_{2}} + 3 \, \theta_{0}^{2} \overline{\alpha_{2}} + 8 \, \theta_{0} \overline{\alpha_{2}} - 4 \, \overline{\alpha_{2}}\right)} A_{0} C_{1} & 2 \, \theta_{0} - 2 & 0 \\
	    		-8 \, {\left(3 \, \theta_{0}^{3} \overline{\alpha_{2}} - 5 \, \theta_{0}^{2} \overline{\alpha_{2}} - 4 \, \theta_{0} \overline{\alpha_{2}} + 4 \, \overline{\alpha_{2}}\right)} A_{0}^{2} {\left| A_{1} \right|}^{2} & 3 \, \theta_{0} - 3 & -\theta_{0} + 1
	    		\end{dmatrix}\\
	    		&\lambda_1=-8 \, {\left(\alpha_{2} \theta_{0}^{3} - 2 \, \alpha_{2} \theta_{0}^{2} - \alpha_{2} \theta_{0} + 2 \, \alpha_{2}\right)} A_{0}^{2} {\left| A_{1} \right|}^{2} + 2 \, {\left(\theta_{0}^{2} - 3 \, \theta_{0} + 2\right)} A_{0} B_{1} {\left| A_{1} \right|}^{2}\\
	    		& + \frac{1}{4} \, {\left(8 \, {\left(\theta_{0}^{2} \overline{\alpha_{5}} - 3 \, \theta_{0} \overline{\alpha_{5}} + 2 \, \overline{\alpha_{5}}\right)} A_{0} - {\left(\theta_{0}^{2} - 3 \, \theta_{0} + 2\right)} \overline{B_{1}}\right)} C_{1}\\
	    		&\partial\bar{\partial}\h_0\totimes\partial\bar{\partial}\h_0=\begin{dmatrix}
	    		4 \, {\left(\theta_{0}^{2} - 4 \, \theta_{0} + 3\right)} A_{0} C_{2} {\left| A_{1} \right|}^{2} & \theta_{0} - 2 & \theta_{0} - 1 \\
	    		16 \, {\left(\theta_{0}^{2} - 2 \, \theta_{0} + 1\right)} A_{0}^{2} {\left| A_{1} \right|}^{4} & 2 \, \theta_{0} - 4 & 0 \\
	    		32 \, {\left(\theta_{0}^{2} \overline{\alpha_{5}} - 2 \, \theta_{0} \overline{\alpha_{5}} + \overline{\alpha_{5}}\right)} A_{0}^{2} {\left| A_{1} \right|}^{2} - 4 \, {\left(\theta_{0}^{2} - 2 \, \theta_{0} + 1\right)} A_{0} {\left| A_{1} \right|}^{2} \overline{B_{1}} & 2 \, \theta_{0} - 4 & 1 \\
	    		32 \, {\left(\theta_{0}^{2} - \theta_{0}\right)} A_{0} A_{1} {\left| A_{1} \right|}^{4} + 32 \, {\left(\alpha_{5} \theta_{0}^{2} - \alpha_{5} \theta_{0}\right)} A_{0}^{2} {\left| A_{1} \right|}^{2} & 2 \, \theta_{0} - 3 & 0 \\
	    		-16 \, {\left(2 \, \theta_{0}^{4} \overline{\alpha_{2}} - 5 \, \theta_{0}^{3} \overline{\alpha_{2}} + 5 \, \theta_{0} \overline{\alpha_{2}} - 2 \, \overline{\alpha_{2}}\right)} A_{0}^{2} {\left| A_{1} \right|}^{2} & 3 \, \theta_{0} - 4 & -\theta_{0} + 1
	    		\end{dmatrix}
	    		\end{align*}
	    		
	    		\normalsize
	    		
	    		We have for some constants $\alpha_1\in\mathbb{C}$ and
	    		$\vec{A}_0,\vec{A}_1,\vec{A}_2,\vec{C}_1\in \mathbb{C}^n$ the development
	    		\begin{align}\label{3form12}
	    		\left\{
	    		\begin{alignedat}{1}
	    		g&=|z|^{2\theta_0-2}|dz|^2\\
	    		\H&=\Re\left(\frac{\vec{C}_1}{z^{\theta_0-2}}\right)+O(|z|^{3-\theta_0})\\
	    		\h_0&=2\left(\vec{A}_1-2|\vec{A}_1|^2\vec{A}_0\z\right)z^{\theta_0-1}+4\left(\vec{A}_2-\alpha_1\vec{A}_0\right)z^{\theta_0}dz^2-\frac{(\theta_0-2)}{2\theta_0}\vec{C}_1\z^{\theta_0}dz^2+O(|z|^{\theta_0+1})
	    		\end{alignedat}\right.
	    		\end{align}
	    		In particular, we have
	    		\begin{align}\label{decayh0}
	    		\left\{
	    		\begin{alignedat}{1}
	    		&\partial\h_0=\left(2(\theta_0-1)\left(\vec{A}_1-2|\vec{A}_1|^2\vec{A}_0\z\right)z^{\theta_0-2}+4\left(\vec{A}_2-\alpha_1\vec{A}_0\right)z^{\theta_0-1}\right)dz^3+O(|z|^{\theta_0})=O(|z|^{\theta_0-2})\\
	    		&\bar{\partial}\h_0=\left(-4|\vec{A}_1|^2\vec{A}_0z^{\theta_0-1}-\frac{(\theta_0-2)}{2}\vec{C}_1\z^{\theta_0-1}\right)dz^2\totimes d\z+O(|z|^{\theta_0})=O(|z|^{\theta_0-1})\\
	    		&\partial\bar{\partial}\h_0=-4(\theta_0-1)|\vec{A}_1|^2\vec{A}_0z^{\theta_0-2}+O(|z|^{\theta_0-1})
	    		\end{alignedat}\right.
	    		\end{align}
	    		Now, recall that
	    		\begin{align*}
	    		\left\{
	    		\begin{alignedat}{1}
	    		&\s{\vec{A}_0}{\vec{A}_0}=\s{\vec{A}_0}{\vec{A}_1}=\s{\vec{A}_0}{\bar{\vec{A}_1}}=\s{\vec{A}_0}{\vec{C}_1}=\s{\vec{A}_0}{\bar{\vec{C}_1}}=0,\quad \s{\vec{A}_1}{\vec{A}_1}+2\s{\vec{A}_0}{\vec{A}_2}=0,\\
	    		& \s{\vec{A}_1}{\vec{C}_1}=\s{\vec{A}_0}{\vec{C}_2}=0,\quad |\vec{C}_1|^2\s{\vec{A}_1}{\vec{A}_1}=0
	    		\end{alignedat}\right.
	    		\end{align*}
	    		we obtain
	    		\small
	    		\begin{align}\label{h021}
	    		&\h_0\totimes\h_0=\begin{dmatrix}
	    		\frac{{\left(\theta_{0}^{2} - 4 \, \theta_{0} + 4\right)} C_{1}^{2}}{4 \, \theta_{0}^{2}} & 0 & 2 \, \theta_{0} \\
	    		-\frac{2 \, {\left(A_{1} C_{2} {\left(\theta_{0} - 3\right)} - 2 \, {\left(\ccancel{{\left(\alpha_{1} \theta_{0} - 2 \, \alpha_{1}\right)} A_{0}} - A_{2} {\left(\theta_{0} - 2\right)}\right)} C_{1}\right)}}{\theta_{0}} & \theta_{0} & \theta_{0} \\
	    		-\frac{2 \, \ccancel{A_{1} C_{1}} {\left(\theta_{0} - 2\right)}}{\theta_{0}} & \theta_{0} - 1 & \theta_{0} \\
	    		\frac{2 \, {\left(\ccancel{2 \, {\left(\theta_{0}^{2} - \theta_{0} - 2\right)} A_{0} C_{1} {\left| A_{1} \right|}^{2}} + \ccancel{4 \, {\left(\alpha_{2} \theta_{0}^{3} - \alpha_{2} \theta_{0}^{2} - 2 \, \alpha_{2} \theta_{0}\right)} A_{0} A_{1}} - \ccancel{{\left(\theta_{0}^{2} - 2 \, \theta_{0}\right)} A_{1} B_{1}}\right)}}{\theta_{0}^{2} + \theta_{0}} & \theta_{0} - 1 & \theta_{0} + 1 \\
	    		4 \, A_{1}^{2} & 2 \, \theta_{0} - 2 & 0 \\
	    		-\ccancel{16 \, A_{0} A_{1} {\left| A_{1} \right|}^{2}} & 2 \, \theta_{0} - 2 & 1 \\
	    		\ccancel{16 \, A_{0}^{2} {\left| A_{1} \right|}^{4}} - \ccancel{8 \, A_{0} A_{1} \overline{\alpha_{5}}} + \ccancel{A_{1} \overline{B_{1}}} & 2 \, \theta_{0} - 2 & 2 \\
	    		-\ccancel{16 \, A_{0} A_{1} \alpha_{1}} + 16 \, A_{1} A_{2} & 2 \, \theta_{0} - 1 & 0 \\
	    		-\colorcancel{32 \, A_{0} A_{2} {\left| A_{1} \right|}^{2}}{blue} - \ccancel{16 \, A_{0} A_{1} \alpha_{5}} + 16 \, {\left(\ccancel{2 \, A_{0}^{2} \alpha_{1}} - \colorcancel{A_{1}^{2}}{blue}\right)} {\left| A_{1} \right|}^{2} & 2 \, \theta_{0} - 1 & 1 \\
	    		-\ccancel{8 \, {\left(\theta_{0} \overline{\alpha_{2}} + \overline{\alpha_{2}}\right)} A_{0} A_{1}} & 3 \, \theta_{0} - 2 & -\theta_{0} + 2 \\
	    		\ccancel{16 \, A_{0}^{2} \alpha_{1}^{2}} - \colorcancel{16 \, A_{1}^{2} \alpha_{1}}{blue} - \colorcancel{32 \, A_{0} A_{2} \alpha_{1}}{blue} - \ccancel{24 \, A_{0} A_{1} \alpha_{3}} + 16 \, A_{2}^{2} + 24 \, A_{1} A_{3} & 2 \, \theta_{0} & 0
	    		\end{dmatrix}\\
	    		&=\begin{dmatrix}
	    		\frac{{\left(\theta_{0}^{2} - 4 \, \theta_{0} + 4\right)} C_{1}^{2}}{4 \, \theta_{0}^{2}} & 0 & 2 \, \theta_{0}& \textbf{(1)}\\
	    		-\frac{2 \, {\left(A_{1} C_{2} {\left(\theta_{0} - 3\right)} + 2 \,{\left(\theta_{0} - 2\right)} A_{2} C_{1}\right)}}{\theta_{0}} & \theta_{0} & \theta_{0}  & \textbf{(2)} \\
	    		4 \, A_{1}^{2} & 2 \, \theta_{0} - 2 & 0  & \textbf{(3)} \\
	    		16 \, A_{1} A_{2} & 2 \, \theta_{0} - 1 & 0  & \textbf{(4)} \\
	    		16 \, A_{2}^{2} + 24 \, A_{1} A_{3} & 2 \, \theta_{0} & 0 & \textbf{(5)} 
	    		\end{dmatrix}
	    		\end{align}
	    		Then, we have
	    		\begin{align}\label{h022}
	    		&\partial\h_0\totimes\h_0=\\
	    		&\begin{dmatrix}
	    		-\frac{\ccancel{{\left(\theta_{0}^{2} - 3 \, \theta_{0} + 2\right)} A_{1} C_{1}}{\theta_{0}}} & \theta_{0} - 2 & \theta_{0} \\
	    		\frac{\ccancel{2 \, {\left(\theta_{0}^{3} - 2 \, \theta_{0}^{2} - \theta_{0} + 2\right)} A_{0} C_{1} {\left| A_{1} \right|}^{2}} + \ccancel{4 \, {\left(\alpha_{2} \theta_{0}^{4} - 2 \, \alpha_{2} \theta_{0}^{3} - \alpha_{2} \theta_{0}^{2} + 2 \, \alpha_{2} \theta_{0}\right)} A_{0} A_{1}} - \ccancel{{\left(\theta_{0}^{3} - 3 \, \theta_{0}^{2} + 2 \, \theta_{0}\right)} A_{1} B_{1}}}{\theta_{0}^{2} + \theta_{0}} & \theta_{0} - 2 & \theta_{0} + 1 \\
	    		-A_{1} C_{2} {\left(\theta_{0} - 3\right)} + 2 \, {\left(\ccancel{{\left(\alpha_{1} \theta_{0} - 2 \, \alpha_{1}\right)} A_{0}} - A_{2} {\left(\theta_{0} - 2\right)}\right)} C_{1} & \theta_{0} - 1 & \theta_{0} \\
	    		4 \, A_{1}^{2} {\left(\theta_{0} - 1\right)} & 2 \, \theta_{0} - 3 & 0 \\
	    		-\ccancel{16 \, A_{0} A_{1} {\left(\theta_{0} - 1\right)} {\left| A_{1} \right|}^{2}} & 2 \, \theta_{0} - 3 & 1 \\
	    		\ccancel{16 \, A_{0}^{2} {\left(\theta_{0} - 1\right)} {\left| A_{1} \right|}^{4}} - \ccancel{8 \, {\left(\theta_{0} \overline{\alpha_{5}} - \overline{\alpha_{5}}\right)} A_{0} A_{1}} + \ccancel{A_{1} {\left(\theta_{0} - 1\right)} \overline{B_{1}}} & 2 \, \theta_{0} - 3 & 2 \\
	    		-\ccancel{8 \, {\left(2 \, \alpha_{1} \theta_{0} - \alpha_{1}\right)} A_{0} A_{1}} + 8 \, A_{1} A_{2} {\left(2 \, \theta_{0} - 1\right)} & 2 \, \theta_{0} - 2 & 0 \\
	    		-\colorcancel{16 \, A_{0} A_{2} {\left(2 \, \theta_{0} - 1\right)} {\left| A_{1} \right|}^{2}}{blue} - \ccancel{8 \, {\left(2 \, \alpha_{5} \theta_{0} - \alpha_{5}\right)} A_{0} A_{1}} + 8 \, {\left(\ccancel{2 \, {\left(2 \, \alpha_{1} \theta_{0} - \alpha_{1}\right)} A_{0}^{2}} - \colorcancel{A_{1}^{2} {\left(2 \, \theta_{0} - 1\right)}}{blue}\right)} {\left| A_{1} \right|}^{2} & 2 \, \theta_{0} - 2 & 1 \\
	    		\ccancel{16 \, A_{0}^{2} \alpha_{1}^{2} \theta_{0}} - \colorcancel{16 \, A_{1}^{2} \alpha_{1} \theta_{0}}{blue} - \colorcancel{32 \, A_{0} A_{2} \alpha_{1} \theta_{0}}{blue} - \ccancel{24 \, A_{0} A_{1} \alpha_{3} \theta_{0}} + 16 \, A_{2}^{2} \theta_{0} + 24 \, A_{1} A_{3} \theta_{0} & 2 \, \theta_{0} - 1 & 0 \\
	    		-\ccancel{4 \, {\left(3 \, \theta_{0}^{2} \overline{\alpha_{2}} + \theta_{0} \overline{\alpha_{2}} - 2 \, \overline{\alpha_{2}}\right)} A_{0} A_{1}} & 3 \, \theta_{0} - 3 & -\theta_{0} + 2
	    		\end{dmatrix}\nonumber\\
	    		&=\begin{dmatrix}
	    		-A_{1} C_{2} {\left(\theta_{0} - 3\right)} - 2 \,{\left(\theta_{0} - 2\right)} A_{2} C_{1} & \theta_{0} - 1 & \theta_{0} &\textbf{(1)}\\
	    		4 \, A_{1}^{2} {\left(\theta_{0} - 1\right)} & 2 \, \theta_{0} - 3 & 0 &\textbf{(2)}\\
	    		8 \, A_{1} A_{2} {\left(2 \, \theta_{0} - 1\right)} & 2 \, \theta_{0} - 2 & 0 &\textbf{(3)}\\
	    		16 \, A_{2}^{2} \theta_{0} + 24 \, A_{1} A_{3} \theta_{0} & 2 \, \theta_{0} - 1 & 0 &\textbf{(4)}\\
	    		\end{dmatrix}
	    		\end{align}
	    		Now, 
	    		\begin{align}\label{h023}
	    		&\bar{\partial}\h_0\totimes \h_0=\\
	    		&\begin{dmatrix}
	    		\frac{{\left(\theta_{0}^{2} - 4 \, \theta_{0} + 4\right)} C_{1}^{2}}{4 \, \theta_{0}} & 0 & 2 \, \theta_{0} - 1 \\
	    		-A_{1} C_{2} {\left(\theta_{0} - 3\right)} + 2 \, {\left(\ccancel{{\left(\alpha_{1} \theta_{0} - 2 \, \alpha_{1}\right)} A_{0}} - A_{2} {\left(\theta_{0} - 2\right)}\right)} C_{1} & \theta_{0} & \theta_{0} - 1 \\
	    		\frac{\ccancel{2 \, {\left(\theta_{0}^{2} - \theta_{0} - 2\right)} A_{0} C_{1} {\left| A_{1} \right|}^{2}} + \ccancel{4 \, {\left(\alpha_{2} \theta_{0}^{3} - \alpha_{2} \theta_{0}^{2} - 2 \, \alpha_{2} \theta_{0}\right)} A_{0} A_{1}} - \ccancel{{\left(\theta_{0}^{2} - 2 \, \theta_{0}\right)} A_{1} B_{1}}}{\theta_{0}} & \theta_{0} - 1 & \theta_{0} \\
	    		-\ccancel{A_{1} C_{1} {\left(\theta_{0} - 2\right)}} & \theta_{0} - 1 & \theta_{0} - 1 \\
	    		-\ccancel{8 \, A_{0} A_{1} {\left| A_{1} \right|}^{2}} & 2 \, \theta_{0} - 2 & 0 \\
	    		\ccancel{16 \, A_{0}^{2} {\left| A_{1} \right|}^{4}} - \ccancel{8 \, A_{0} A_{1} \overline{\alpha_{5}}} + \ccancel{A_{1} \overline{B_{1}}} & 2 \, \theta_{0} - 2 & 1 \\
	    		-\colorcancel{16 \, A_{0} A_{2} {\left| A_{1} \right|}^{2}}{blue} -\ccancel{ 8 \, A_{0} A_{1} \alpha_{5}} + 8 \, {\left(\ccancel{2 \, A_{0}^{2} \alpha_{1}} - \colorcancel{A_{1}^{2}}{blue}\right)} {\left| A_{1} \right|}^{2} & 2 \, \theta_{0} - 1 & 0 \\
	    		\ccancel{4 \, {\left(\theta_{0}^{2} \overline{\alpha_{2}} - \theta_{0} \overline{\alpha_{2}} - 2 \, \overline{\alpha_{2}}\right)} A_{0} A_{1}} & 3 \, \theta_{0} - 2 & -\theta_{0} + 1
	    		\end{dmatrix}\nonumber\\
	    		&=\begin{dmatrix}
	    		\frac{{(\theta_{0}-2)^2} C_{1}^{2}}{4 \, \theta_{0}} & 0 & 2 \, \theta_{0} - 1 &\textbf{(1)}\\
	    		-A_{1} C_{2} {\left(\theta_{0} - 3\right)} -2\, {\left(\theta_{0} - 2\right)} A_{2} C_{1} & \theta_{0} & \theta_{0} - 1 &\textbf{(2)} 
	    		\end{dmatrix}
	    		\end{align}
	    		Then
	    		\begin{align}\label{h024}
	    		&\partial\h_0\totimes\partial\h_0=\nonumber\\
	    		&\begin{dmatrix}
	    		-\frac{2 \, {\left(\theta_{0}^{2} - 4 \, \theta_{0} + 3\right)} A_{1} C_{2}}{\theta_{0}} & \theta_{0} - 2 & \theta_{0} \\
	    		4 \, {\left(\theta_{0}^{2} - 2 \, \theta_{0} + 1\right)} A_{1}^{2} & 2 \, \theta_{0} - 4 & 0 \\
	    		-\ccancel{16 \, {\left(\theta_{0}^{2} - 2 \, \theta_{0} + 1\right)} A_{0} A_{1} {\left| A_{1} \right|}^{2}} & 2 \, \theta_{0} - 4 & 1 \\
	    		\ccancel{16 \, {\left(\theta_{0}^{2} - 2 \, \theta_{0} + 1\right)} A_{0}^{2} {\left| A_{1} \right|}^{4}} - \ccancel{8 \, {\left(\theta_{0}^{2} \overline{\alpha_{5}} - 2 \, \theta_{0} \overline{\alpha_{5}} + \overline{\alpha_{5}}\right)} A_{0} A_{1}} + \ccancel{{\left(\theta_{0}^{2} - 2 \, \theta_{0} + 1\right)} A_{1} \overline{B_{1}}} & 2 \, \theta_{0} - 4 & 2 \\
	    		-\ccancel{16 \, {\left(\alpha_{1} \theta_{0}^{2} - \alpha_{1} \theta_{0}\right)} A_{0} A_{1}} + \ccancel{16 \, {\left(\theta_{0}^{2} - \theta_{0}\right)} A_{1} A_{2}} & 2 \, \theta_{0} - 3 & 0 \\
	    		-\colorcancel{32 \, {\left(\theta_{0}^{2} - \theta_{0}\right)} A_{0} A_{2} {\left| A_{1} \right|}^{2}}{blue} - \ccancel{16 \, {\left(\alpha_{5} \theta_{0}^{2} - \alpha_{5} \theta_{0}\right)} A_{0} A_{1}} + 16 \, {\left(\ccancel{2 \, {\left(\alpha_{1} \theta_{0}^{2} - \alpha_{1} \theta_{0}\right)} A_{0}^{2}} - \colorcancel{{\left(\theta_{0}^{2} - \theta_{0}\right)} A_{1}^{2}}{blue}\right)} {\left| A_{1} \right|}^{2} & 2 \, \theta_{0} - 3 & 1 \\
	    		\ccancel{16 \, A_{0}^{2} \alpha_{1}^{2} \theta_{0}^{2}} - \colorcancel{32 \, A_{0} A_{2} \alpha_{1} \theta_{0}^{2}}{blue} + 16 \, A_{2}^{2} \theta_{0}^{2} - \ccancel{24 \, {\left(\alpha_{3} \theta_{0}^{2} - \alpha_{3}\right)} A_{0} A_{1}} - \colorcancel{16 \, {\left(\alpha_{1} \theta_{0}^{2} - \alpha_{1}\right)} A_{1}^{2}}{blue} + 24 \, {\left(\theta_{0}^{2} - 1\right)} A_{1} A_{3} & 2 \, \theta_{0} - 2 & 0 \\
	    		-\ccancel{8 \, {\left(2 \, \theta_{0}^{3} \overline{\alpha_{2}} - \theta_{0}^{2} \overline{\alpha_{2}} - 2 \, \theta_{0} \overline{\alpha_{2}} + \overline{\alpha_{2}}\right)} A_{0} A_{1}} & 3 \, \theta_{0} - 4 & -\theta_{0} + 2
	    		\end{dmatrix}\nonumber\\
	    		&=\begin{dmatrix}
	    		-\frac{2 \, {\left(\theta_{0}^{2} - 4 \, \theta_{0} + 3\right)} A_{1} C_{2}}{\theta_{0}} & \theta_{0} - 2 & \theta_{0} &\textbf{(1)}\\
	    		4 \, {\left(\theta_{0}^{2} - 2 \, \theta_{0} + 1\right)} A_{1}^{2} & 2 \, \theta_{0} - 4 & 0 &\textbf{(2)}\\
	    		16 \, A_{2}^{2} \theta_{0}^{2} + 24 \, {\left(\theta_{0}^{2} - 1\right)} A_{1} A_{3} & 2 \, \theta_{0} - 2 & 0 &\textbf{(3)}
	    		\end{dmatrix}
	    		\end{align}
	    		Now
	    		\begin{align}\label{h025}
	    		&\bar{\partial}\h_0\totimes\bar{\partial}\h_0=\nonumber\\
	    		&\begin{dmatrix}
	    		\frac{1}{4} \, {\left(\theta_{0}^{2} - 4 \, \theta_{0} + 4\right)} C_{1}^{2} & 0 & 2 \, \theta_{0} - 2 \\
	    		-\frac{1}{2} \, {\left(\ccancel{4 \, {\left(\alpha_{2} \theta_{0}^{3} - 3 \, \alpha_{2} \theta_{0}^{2} + 4 \, \alpha_{2}\right)} A_{0}} - \ccancel{{\left(\theta_{0}^{2} - 4 \, \theta_{0} + 4\right)} B_{1}}\right)} C_{1} & 0 & 2 \, \theta_{0} - 1 \\
	    		\frac{1}{2} \, {\left(\theta_{0}^{2} - 5 \, \theta_{0} + 6\right)} C_{1} C_{2} & 1 & 2 \, \theta_{0} - 2 \\
	    		\ccancel{4 \, A_{0} C_{2} {\left(\theta_{0} - 3\right)} {\left| A_{1} \right|}^{2}} + 4 \, {\left(\ccancel{A_{1} {\left(\theta_{0} - 2\right)} {\left| A_{1} \right|}^{2}} + \ccancel{{\left(\alpha_{5} \theta_{0} - 2 \, \alpha_{5}\right)} A_{0}}\right)} C_{1} & \theta_{0} & \theta_{0} - 1 \\
	    		-\ccancel{16 \, {\left(\alpha_{2} \theta_{0}^{2} - \alpha_{2} \theta_{0} - 2 \, \alpha_{2}\right)} A_{0}^{2} {\left| A_{1} \right|}^{2}} + \ccancel{4 \, A_{0} B_{1} {\left(\theta_{0} - 2\right)} {\left| A_{1} \right|}^{2}} + \frac{1}{2} \, {\left(\ccancel{8 \, {\left(\theta_{0} \overline{\alpha_{5}} - 2 \, \overline{\alpha_{5}}\right)} A_{0}} - \ccancel{{\left(\theta_{0} - 2\right)} \overline{B_{1}}}\right)} C_{1} & \theta_{0} - 1 & \theta_{0} \\
	    		\ccancel{4 \, A_{0} C_{1} {\left(\theta_{0} - 2\right)} {\left| A_{1} \right|}^{2}} & \theta_{0} - 1 & \theta_{0} - 1 \\
	    		\ccancel{16 \, A_{0}^{2} {\left| A_{1} \right|}^{4}} & 2 \, \theta_{0} - 2 & 0 \\
	    		\ccancel{32 \, A_{0}^{2} {\left| A_{1} \right|}^{2} \overline{\alpha_{5}}} - 4 \, A_{0} {\left| A_{1} \right|}^{2} \overline{B_{1}} & 2 \, \theta_{0} - 2 & 1 \\
	    		\ccancel{32 \, A_{0} A_{1} {\left| A_{1} \right|}^{4}} + \ccancel{32 \, A_{0}^{2} \alpha_{5} {\left| A_{1} \right|}^{2}} - \ccancel{2 \, {\left(\theta_{0}^{3} \overline{\alpha_{2}} - 3 \, \theta_{0}^{2} \overline{\alpha_{2}} + 4 \, \overline{\alpha_{2}}\right)} A_{0} C_{1} }& 2 \, \theta_{0} - 1 & 0 \\
	    		-\ccancel{16 \, {\left(\theta_{0}^{2} \overline{\alpha_{2}} - \theta_{0} \overline{\alpha_{2}} - 2 \, \overline{\alpha_{2}}\right)} A_{0}^{2} {\left| A_{1} \right|}^{2}} & 3 \, \theta_{0} - 2 & -\theta_{0} + 1
	    		\end{dmatrix}\nonumber\\
	    		&=\begin{dmatrix}
	    		\frac{1}{4} \, {\left(\theta_{0}^{2} - 4 \, \theta_{0} + 4\right)} C_{1}^{2} & 0 & 2 \, \theta_{0} - 2 &\textbf{(1)}\\
	    		\frac{1}{2} \, {\left(\theta_{0}^{2} - 5 \, \theta_{0} + 6\right)} C_{1} C_{2} & 1 & 2 \, \theta_{0} - 2 &\textbf{(2)}\\
	    		8\theta_0(\theta_0+1)|\vec{A}_1|^2\bar{\alpha_2}& 2 \, \theta_{0} - 2 & 1 &\textbf{(3)}
	    		\end{dmatrix}
	    		\end{align}
	    		as
	    		\begin{align}\label{octic1}
	    		&\alpha_2=\frac{1}{2\theta_0(\theta_0+1)}\s{\bar{\vec{A}_1}}{\vec{C}_1}\nonumber\\
	    		&\s{\vec{A}_0}{\bar{\vec{B}_1}}=\s{\vec{A}_0}{-2\s{\vec{A}_1}{\bar{\vec{C}_1}}\bar{\vec{A}_0}}=-\s{\vec{A}_1}{\bar{\vec{C}_1}}=-2\theta_0(\theta_0+1)\bar{\alpha_2}.
	    		\end{align}
	    		Now, we have
	    		\begin{align}\label{h026}
	    		&\partial\h_0\totimes\bar{\partial}\h_0=\nonumber\\
	    		&\begin{dmatrix}
	    		\ccancel{2 \, {\left(\theta_{0}^{2} - 3 \, \theta_{0} + 2\right)} A_{0} C_{1} {\left| A_{1} \right|}^{2}} + \ccancel{4 \, {\left(\alpha_{2} \theta_{0}^{3} - 2 \, \alpha_{2} \theta_{0}^{2} - \alpha_{2} \theta_{0} + 2 \, \alpha_{2}\right)} A_{0} A_{1}} -\ccancel{ {\left(\theta_{0}^{2} - 3 \, \theta_{0} + 2\right)} A_{1} B_{1}} & \theta_{0} - 2 & \theta_{0} \\
	    		-\ccancel{{\left(\theta_{0}^{2} - 3 \, \theta_{0} + 2\right)} A_{1} C_{1}} & \theta_{0} - 2 & \theta_{0} - 1 \\
	    		-{\left(\theta_{0}^{2} - 4 \, \theta_{0} + 3\right)} A_{1} C_{2} + 2 \, {\left(\ccancel{{\left(\alpha_{1} \theta_{0}^{2} - 2 \, \alpha_{1} \theta_{0}\right)} A_{0}} - {\left(\theta_{0}^{2} - 2 \, \theta_{0}\right)} A_{2}\right)} C_{1} & \theta_{0} - 1 & \theta_{0} - 1 \\
	    		-\ccancel{8 \, A_{0} A_{1} {\left(\theta_{0} - 1\right)} {\left| A_{1} \right|}^{2}} & 2 \, \theta_{0} - 3 & 0 \\
	    		\ccancel{16 \, A_{0}^{2} {\left(\theta_{0} - 1\right)} {\left| A_{1} \right|}^{4}} - \ccancel{8 \, {\left(\theta_{0} \overline{\alpha_{5}} - \overline{\alpha_{5}}\right)} A_{0} A_{1}} + \ccancel{A_{1} {\left(\theta_{0} - 1\right)} \overline{B_{1}}} & 2 \, \theta_{0} - 3 & 1 \\
	    		-\colorcancel{16 \, A_{0} A_{2} \theta_{0} {\left| A_{1} \right|}^{2}}{blue} - \ccancel{8 \, {\left(\alpha_{5} \theta_{0} - \alpha_{5}\right)} A_{0} A_{1}} + 8 \, {\left(\ccancel{2 \, A_{0}^{2} \alpha_{1} \theta_{0}} - \colorcancel{A_{1}^{2} {\left(\theta_{0} - 1\right)}}{blue}\right)} {\left| A_{1} \right|}^{2} & 2 \, \theta_{0} - 2 & 0 \\
	    		\ccancel{4 \, {\left(\theta_{0}^{3} \overline{\alpha_{2}} - 2 \, \theta_{0}^{2} \overline{\alpha_{2}} - \theta_{0} \overline{\alpha_{2}} + 2 \, \overline{\alpha_{2}}\right)} A_{0} A_{1}} & 3 \, \theta_{0} - 3 & -\theta_{0} + 1
	    		\end{dmatrix}\nonumber\\
	    		&=\begin{dmatrix}
	    		-{\left(\theta_{0}^{2} - 4 \, \theta_{0} + 3\right)} A_{1} C_{2} - 2 \, {\left(\theta_{0}^{2} - 2 \, \theta_{0}\right)} A_{2} C_{1} & \theta_{0} - 1 & \theta_{0} - 1 
	    		\end{dmatrix}
	    		\end{align}
	    		Then
	    		\begin{align}\label{h027}
	    		\partial\bar{\partial}\h_0\totimes\h_0=\nonumber\\
	    		&\begin{dmatrix}
	    		\frac{\ccancel{2 \, {\left(\theta_{0}^{2} - 3 \, \theta_{0} + 2\right)} A_{0} C_{1} {\left| A_{1} \right|}^{2}}}{\theta_{0}} & \theta_{0} - 2 & \theta_{0} \\
	    		-A_{1} C_{2} {\left(\theta_{0} - 3\right)} & \theta_{0} - 1 & \theta_{0} - 1 \\
	    		-\ccancel{8 \, A_{0} A_{1} {\left(\theta_{0} - 1\right)} {\left| A_{1} \right|}^{2}} & 2 \, \theta_{0} - 3 & 0 \\
	    		\ccancel{16 \, A_{0}^{2} {\left(\theta_{0} - 1\right)} {\left| A_{1} \right|}^{4}} - \ccancel{8 \, {\left(\theta_{0} \overline{\alpha_{5}} - \overline{\alpha_{5}}\right)} A_{0} A_{1}} + \ccancel{A_{1} {\left(\theta_{0} - 1\right)} \overline{B_{1}}} & 2 \, \theta_{0} - 3 & 1 \\
	    		-\colorcancel{16 \, A_{0} A_{2} {\left(\theta_{0} - 1\right)} {\left| A_{1} \right|}^{2}}{blue} - \ccancel{8 \, A_{0} A_{1} \alpha_{5} \theta_{0}} + 8 \, {\left(\ccancel{2 \, {\left(\alpha_{1} \theta_{0} - \alpha_{1}\right)} A_{0}^{2}} - \colorcancel{A_{1}^{2} \theta_{0}}{blue}\right)} {\left| A_{1} \right|}^{2} & 2 \, \theta_{0} - 2 & 0 \\
	    		\ccancel{4 \, {\left(2 \, \theta_{0}^{3} \overline{\alpha_{2}} - 3 \, \theta_{0}^{2} \overline{\alpha_{2}} - 3 \, \theta_{0} \overline{\alpha_{2}} + 2 \, \overline{\alpha_{2}}\right)} A_{0} A_{1}} & 3 \, \theta_{0} - 3 & -\theta_{0} + 1
	    		\end{dmatrix}\nonumber\\
	    		&=\begin{dmatrix}
	    		-A_{1} C_{2} {\left(\theta_{0} - 3\right)} & \theta_{0} - 1 & \theta_{0} - 1
	    		\end{dmatrix}
	    		\end{align}
	    		while
	    		\begin{align}\label{h028}
	    		&\partial\bar{\partial}\h_0\totimes\partial\h_0=\nonumber\\
	    		&\begin{dmatrix}
	    		-{\left(\theta_{0}^{2} - 4 \, \theta_{0} + 3\right)} A_{1} C_{2} & \theta_{0} - 2 & \theta_{0} - 1 \\
	    		-\ccancel{8 \, {\left(\theta_{0}^{2} - 2 \, \theta_{0} + 1\right)} A_{0} A_{1} {\left| A_{1} \right|}^{2}} & 2 \, \theta_{0} - 4 & 0 \\
	    		\ccancel{16 \, {\left(\theta_{0}^{2} - 2 \, \theta_{0} + 1\right)} A_{0}^{2} {\left| A_{1} \right|}^{4}} - \ccancel{8 \, {\left(\theta_{0}^{2} \overline{\alpha_{5}} - 2 \, \theta_{0} \overline{\alpha_{5}} + \overline{\alpha_{5}}\right)} A_{0} A_{1}} + \ccancel{{\left(\theta_{0}^{2} - 2 \, \theta_{0} + 1\right)} A_{1} \overline{B_{1}}} & 2 \, \theta_{0} - 4 & 1 \\
	    		-\colorcancel{16 \, {\left(\theta_{0}^{2} - \theta_{0}\right)} A_{0} A_{2} {\left| A_{1} \right|}^{2}}{blue} - \ccancel{8 \, {\left(\alpha_{5} \theta_{0}^{2} - \alpha_{5} \theta_{0}\right)} A_{0} A_{1}} + 8 \, {\left(\ccancel{2 \, {\left(\alpha_{1} \theta_{0}^{2} - \alpha_{1} \theta_{0}\right)} A_{0}^{2}} -\colorcancel{{\left(\theta_{0}^{2} - \theta_{0}\right)} A_{1}^{2}}{blue}\right)} {\left| A_{1} \right|}^{2} & 2 \, \theta_{0} - 3 & 0 \\
	    		\ccancel{4 \, {\left(2 \, \theta_{0}^{4} \overline{\alpha_{2}} - 5 \, \theta_{0}^{3} \overline{\alpha_{2}} + 5 \, \theta_{0} \overline{\alpha_{2}} - 2 \, \overline{\alpha_{2}}\right)} A_{0} A_{1}} & 3 \, \theta_{0} - 4 & -\theta_{0} + 1
	    		\end{dmatrix}\nonumber\\
	    		&=\begin{dmatrix}
	    		-{\left(\theta_{0}^{2} - 4 \, \theta_{0} + 3\right)} A_{1} C_{2} & \theta_{0} - 2 & \theta_{0} - 1
	    		\end{dmatrix}
	    		\end{align}
	    		Now
	    		\begin{align}\label{h029}
	    		&\partial\bar{\partial}\h_0\totimes\bar{\partial}\h_0=\nonumber\\
	    		&\begin{dmatrix}
	    		\frac{1}{4} \, {\left(\theta_{0}^{2} - 5 \, \theta_{0} + 6\right)} C_{1} C_{2} & 0 & 2 \, \theta_{0} - 2 \\
	    		\ccancel{\lambda_1 }& \theta_{0} - 2 & \theta_{0} \\
	    		\ccancel{2 \, {\left(\theta_{0}^{2} - 3 \, \theta_{0} + 2\right)} A_{0} C_{1} {\left| A_{1} \right|}^{2}} & \theta_{0} - 2 & \theta_{0} - 1 \\
	    		\colorcancel{2 \, {\left(\theta_{0}^{2} - 3 \, \theta_{0}\right)} A_{0} C_{2} {\left| A_{1} \right|}^{2}}{blue} + 2 \, {\left(\colorcancel{{\left(\theta_{0}^{2} - 2 \, \theta_{0}\right)} A_{1} {\left| A_{1} \right|}^{2}}{blue} + \ccancel{{\left(\alpha_{5} \theta_{0}^{2} - 2 \, \alpha_{5} \theta_{0}\right)} A_{0}}\right)} C_{1} & \theta_{0} - 1 & \theta_{0} - 1 \\
	    		\ccancel{16 \, A_{0}^{2} {\left(\theta_{0} - 1\right)} {\left| A_{1} \right|}^{4}} & 2 \, \theta_{0} - 3 & 0 \\
	    		\ccancel{32 \, {\left(\theta_{0} \overline{\alpha_{5}} - \overline{\alpha_{5}}\right)} A_{0}^{2} {\left| A_{1} \right|}^{2}} - 4 \, A_{0} {\left(\theta_{0} - 1\right)} {\left| A_{1} \right|}^{2} \overline{B_{1}} & 2 \, \theta_{0} - 3 & 1 \\
	    		\ccancel{16 \, A_{0} A_{1} {\left(2 \, \theta_{0} - 1\right)} {\left| A_{1} \right|}^{4}} + \ccancel{16 \, {\left(2 \, \alpha_{5} \theta_{0} - \alpha_{5}\right)} A_{0}^{2} {\left| A_{1} \right|}^{2}}- \ccancel{{\left(2 \, \theta_{0}^{4} \overline{\alpha_{2}} - 7 \, \theta_{0}^{3} \overline{\alpha_{2}} + 3 \, \theta_{0}^{2} \overline{\alpha_{2}} + 8 \, \theta_{0} \overline{\alpha_{2}} - 4 \, \overline{\alpha_{2}}\right)} A_{0} C_{1}} & 2 \, \theta_{0} - 2 & 0 \\
	    		-\ccancel{8 \, {\left(3 \, \theta_{0}^{3} \overline{\alpha_{2}} - 5 \, \theta_{0}^{2} \overline{\alpha_{2}} - 4 \, \theta_{0} \overline{\alpha_{2}} + 4 \, \overline{\alpha_{2}}\right)} A_{0}^{2} {\left| A_{1} \right|}^{2}} & 3 \, \theta_{0} - 3 & -\theta_{0} + 1
	    		\end{dmatrix}\nonumber\\
	    		&=\begin{dmatrix}
	    		\frac{1}{4} \, {\left(\theta_{0}^{2} - 5 \, \theta_{0} + 6\right)} C_{1} C_{2} & 0 & 2 \, \theta_{0} - 2 \\
	    		8\theta_0(\theta_0+1)|\vec{A}_1|^2\bar{\alpha_2} & 2 \, \theta_{0} - 3 & 1
	    		\end{dmatrix}
	    		\end{align}
	    		as
	    		\begin{align*}
	    		&\lambda_1=-\ccancel{8 \, {\left(\alpha_{2} \theta_{0}^{3} - 2 \, \alpha_{2} \theta_{0}^{2} - \alpha_{2} \theta_{0} + 2 \, \alpha_{2}\right)} A_{0}^{2} {\left| A_{1} \right|}^{2}} + \ccancel{2 \, {\left(\theta_{0}^{2} - 3 \, \theta_{0} + 2\right)} A_{0} B_{1} {\left| A_{1} \right|}^{2}}\\
	    		& + \frac{1}{4} \, {\left(\ccancel{8 \, {\left(\theta_{0}^{2} \overline{\alpha_{5}} - 3 \, \theta_{0} \overline{\alpha_{5}} + 2 \, \overline{\alpha_{5}}\right)} A_{0}} - \ccancel{{\left(\theta_{0}^{2} - 3 \, \theta_{0} + 2\right)} \overline{B_{1}}}\right)} C_{1}=0
	    		\end{align*}
	    		Finally,
	    		\begin{align}\label{h0210}
	    		&\partial\bar{\partial}\h_0\totimes\partial\bar{\partial}\h_0=\begin{dmatrix}
	    		\ccancel{4 \, {\left(\theta_{0}^{2} - 4 \, \theta_{0} + 3\right)} A_{0} C_{2} {\left| A_{1} \right|}^{2}} & \theta_{0} - 2 & \theta_{0} - 1 \\
	    		\ccancel{16 \, {\left(\theta_{0}^{2} - 2 \, \theta_{0} + 1\right)} A_{0}^{2} {\left| A_{1} \right|}^{4}} & 2 \, \theta_{0} - 4 & 0 \\
	    		\ccancel{32 \, {\left(\theta_{0}^{2} \overline{\alpha_{5}} - 2 \, \theta_{0} \overline{\alpha_{5}} + \overline{\alpha_{5}}\right)} A_{0}^{2} {\left| A_{1} \right|}^{2}} - 4 \, {\left(\theta_{0}^{2} - 2 \, \theta_{0} + 1\right)} A_{0} {\left| A_{1} \right|}^{2} \overline{B_{1}} & 2 \, \theta_{0} - 4 & 1 \\
	    		\ccancel{32 \, {\left(\theta_{0}^{2} - \theta_{0}\right)} A_{0} A_{1} {\left| A_{1} \right|}^{4}} + \ccancel{32 \, {\left(\alpha_{5} \theta_{0}^{2} - \alpha_{5} \theta_{0}\right)} A_{0}^{2} {\left| A_{1} \right|}^{2}} & 2 \, \theta_{0} - 3 & 0 \\
	    		-\ccancel{16 \, {\left(2 \, \theta_{0}^{4} \overline{\alpha_{2}} - 5 \, \theta_{0}^{3} \overline{\alpha_{2}} + 5 \, \theta_{0} \overline{\alpha_{2}} - 2 \, \overline{\alpha_{2}}\right)} A_{0}^{2} {\left| A_{1} \right|}^{2}} & 3 \, \theta_{0} - 4 & -\theta_{0} + 1
	    		\end{dmatrix}\nonumber\\
	    		&=\begin{dmatrix}
	    		8\theta_0(\theta_0-1)^2(\theta_0+1)|\vec{A}_1|^2\bar{\alpha}_2 & 2 \, \theta_{0} - 4 & 1 
	    		\end{dmatrix}
	    		\end{align}
	    		thanks of \eqref{octic1}.
	    		\normalsize
	    		Finally, thanks of \eqref{h021}, \eqref{h022}, \eqref{h023}, \eqref{h024}, \eqref{h025}, \eqref{h026}, \eqref{h027}, \eqref{h028}, \eqref{h029}, \eqref{h0210}, we have (omitting the $dz$, $d\z$ to simplify notations)
	    		\begin{align}\label{h010}
	    		\left\{
	    		\begin{alignedat}{1}
	    		\h_0\totimes\h_0&=4\s{\vec{A}_1}{\vec{A}_1}z^{2\theta_0-2}+16\s{\vec{A}_1}{\vec{A}_2}z^{2\theta_0-1}+\left(16\s{\vec{A}_2}{\vec{A}_2}+24\s{\vec{A}_1}{\vec{A}_3}\right)z^{2\theta_0}\\
	    		&-\frac{2}{\theta_0}\left((\theta_0-3)\s{\vec{A}_1}{\vec{C}_2}+2(\theta_0-2)\s{\vec{A}_2}{\vec{C}_1}\right)|z|^{2\theta_0}+\frac{(\theta_0-2)^2}{4\theta_0^2}\s{\vec{C}_1}{\vec{C}_1}\z^{2\theta_0}+O(|z|^{2\theta_0+1})\\
	    		\partial\h_0\totimes\h_0&=4(\theta_0-1)\s{\vec{A}_1}{\vec{A}_1}z^{2\theta_0-3}+8(2\theta_0-1)\s{\vec{A}_1}{\vec{A}_2}z^{2\theta_0-2}+\theta_0\left(16\s{\vec{A}_2}{\vec{A}_2}+24\s{\vec{A}_1}{\vec{A}_2}\right)z^{2\theta_0-1}\\
	    		&-\left((\theta_0-3)\s{\vec{A}_1}{\vec{C}_2}+2(\theta_0-2)\s{\vec{A}_2}{\vec{C}_1}\right)z^{\theta_0-1}\z^{\theta_0}+O(|z|^{2\theta_0})\\
	    		\bar{\partial}\h_0\totimes\h_0&=-\left((\theta_0-3)\s{\vec{A}_1}{\vec{C}_2}+2(\theta_0-2)\s{\vec{A}_2}{\vec{C}_1}\right)z^{\theta_0}\z^{\theta_0-1}+\frac{(\theta_0-2)^2}{4\theta_0}\s{\vec{C}_1}{\vec{C}_1}\z^{2\theta_0-1}+O(|z|^{2\theta_0})\\
	    		\partial\h_0\totimes\partial\h_0&=4(\theta_0-1)^2\s{\vec{A}_1}{\vec{A}_1}z^{2\theta_0-4}+\left(16\theta_0^2\s{\vec{A}_2}{\vec{A}_2}+24(\theta_0-1)(\theta_0+1)\s{\vec{A}_1}{\vec{A}_3}\right)z^{2\theta_0-2}\\
	    		&-\frac{2(\theta_0-1)(\theta_0-3)}{\theta_0}\s{\vec{A}_1}{\vec{C}_2}z^{\theta_0-2}\z^{\theta_0}+O(|z|^{2\theta_0-1})\\
	    		\bar{\partial}\h_0\totimes\bar{\partial}\h_0&=\frac{1}{4}(\theta_0-2)^2\s{\vec{C}_1}{\vec{C}_1}\z^{2\theta_0-2}+\frac{1}{2}(\theta_0-2)(\theta_0-3)\s{\vec{C}_1}{\vec{C}_2}z\z^{2\theta_0-2}+8\theta_0(\theta_0+1)|\vec{A}_1|^2\bar{\alpha_2}z^{2\theta_0-2}\z\\
	    		&+O(|z|^{2\theta_0})\\
	    		\partial\h_0\totimes\bar{\partial}\h_0&=-\left((\theta_0-1)(\theta_0-3)\s{\vec{A}_1}{\vec{C}_2}+2\theta_0(\theta_0-2)\s{\vec{A}_2}{\vec{C}_1}\right)|z|^{2\theta_0-2}+O(|z|^{2\theta_0-1})\\
	    		\partial\bar{\partial}\h_0\otimes\h_0&=-(\theta_0-3)\s{\vec{A}_1}{\vec{C}_2}|z|^{2\theta_0-2}+O(|z|^{2\theta_0-1})\\
	    		\partial\bar{\partial}\h_0\totimes\partial\h_0&=-(\theta_0-1)(\theta_0-3)\s{\vec{A}_1}{\vec{C}_2}z^{\theta_0-2}\z^{\theta_0-1}+O(|z|^{2\theta_0-2})\\
	    		\partial\bar{\partial}\h_0\totimes\bar{\partial}\h_0&=\frac{1}{4}(\theta_0-2)(\theta_0-3)\s{\vec{C}_1}{\vec{C}_2}\z^{2\theta_0-2}+8\theta_0(\theta_0+1)|\vec{A}_1|^2\bar{\alpha_2}z^{2\theta_0-3}\z+O(|z|^{2\theta_0-1})\\
	    		\partial\bar{\partial}\h_0\totimes\partial\bar{\partial}\h_0&=8\theta_0(\theta_0-1)^2(\theta_0+1)|\vec{A}_1|^2\bar{\alpha_2}z^{2\theta_0-4}\z+O(|z|^{2\theta_0-2}).
	    		\end{alignedat}\right.
	    		\end{align} 
	    		We can easily check that the only potential singular term in $\mathscr{O}_{\phi}$ is
	    		\begin{align*}
	    		O(\h_0)&=g^{-2}\,\otimes\bigg\{\frac{1}{4}\left(\partial\bar{\partial}\h_0\,\dot{\otimes}\,\partial\bar{\partial}\h_0\right)\otimes \left(\h_0\,\dot{\otimes}\,\h_0\right)+\frac{1}{4}\left(\partial\h_0\,\dot{\otimes}\, \partial\h_0\right)\otimes \left(\bar{\partial}\h_0\,\dot{\otimes}\, \bar{\partial}\h_0\right)\nonumber\\
	    		&-\frac{1}{2}\left(\partial\bar{\partial}\h_0\totimes\partial\h_0\right)\otimes\left(\bar{\partial}\h_0\totimes\h_0\right)-\frac{1}{2}\left(\partial\bar{\partial}\h_0\totimes\bar{\partial}\h_0\right)\otimes\left(\partial\h_0\totimes\h_0\right)+\frac{1}{2}\left(\partial\bar{\partial}\h_0\totimes\h_0\right)\otimes\left(\partial\h_0\totimes\bar{\partial}\h_0\right) \bigg\}
	    		\end{align*}
	    		We have  by \eqref{h010}
	    		\begin{align}\label{endoctic1}
	    		&\left(\partial\bar{\partial}\h_0\totimes\partial\bar{\partial}\h_0\right)\otimes\left(\h_0\totimes\h_0\right)\nonumber\\
	    		&=\left(8\theta_0(\theta_0-1)^2(\theta_0+1)|\vec{A}_1|^2\bar{\alpha_2}z^{2\theta_0-4}\z+O(|z|^{2\theta_0-2})\right)\cdot \left(4\s{\vec{A}_1}{\vec{A}_1}z^{2\theta_0-2}+O(|z|^{2\theta_0-1})\right)\nonumber\\
	    		&=32\theta_0(\theta_0-1)^2(\theta_0+1)|\vec{A}_1|^2\bar{\alpha_2}\s{\vec{A}_1}{\vec{A}_1}z^{4\theta_0-6}\z+O(|z|^{4\theta_0-4})\nonumber\\
	    		&=16(\theta_0-1)^2|\vec{A}_1|^2\s{\vec{A}_1}{\bar{\vec{C}_1}}\s{\vec{A}_1}{\vec{A}_1}z^{4\theta_0-6}\z+O(|z|^{4\theta_0-4})
	    		\end{align}
	    		but
	    		\begin{align}\label{newcancel}
	    		|\vec{C}_1|^2\s{\vec{A}_1}{\vec{A}_1}=0
	    		\end{align}
	    		so either $\vec{C}_1=0$, or $\s{\vec{A}_1}{\vec{A}_1}=0$, so that in all cases
	    		\begin{align}\label{addrel}
	    		\s{\vec{A}_1}{\bar{\vec{C}_1}}\s{\vec{A}_1}{\vec{A}_1}=0,
	    		\end{align}
	    		and by \eqref{endoctic1} and \eqref{addrel} 
	    		\begin{align}\label{remove1}
	    		\left(\partial\bar{\partial}\h_0\totimes\partial\bar{\partial}\h_0\right)\otimes\left(\h_0\totimes\h_0\right)=O(|z|^{4\theta_0-4}).
	    		\end{align}
	    		Now, we have
	    		\begin{align}\label{remove2}
	    		&\left(\partial\h_0\totimes\partial\h_0\right)\otimes\left(\bar{\partial}\h_0\totimes\bar{\partial}\h_0\right)
	    		=\left(4(\theta_0-1)^2\s{\vec{A}_1}{\vec{A}_1}z^{2\theta_0-4} +O(|z|^{2\theta_0-2})\right)
	    		\times \nonumber\\
	    		&\bigg(\frac{1}{4}(\theta_0-2)^2\s{\vec{C}_1}{\vec{C}_1}\z^{2\theta_0-2}+\frac{1}{2}(\theta_0-2)(\theta_0-3)\s{\vec{C}_1}{\vec{C}_2}z\z^{2\theta_0-2}
	    		+8\theta_0(\theta_0+1)|\vec{A}_1|^2\bar{\alpha_2}z^{2\theta_0-2}\z+O(|z|^{2\theta_0})\bigg)\nonumber\\
	    		&=(\theta_0-1)^2(\theta_0-2)^2\s{\vec{A}_1}{\vec{A}_1}\s{\vec{C}_1}{\vec{C}_1}z^{2\theta_0-4}\z^{2\theta_0-2}+2(\theta_0-1)^2(\theta_0-2)(\theta_0-3)\s{\vec{A}_1}{\vec{A}_1}\s{\vec{C}_1}{\vec{C}_2}z^{2\theta_0-5}\z^{2\theta_0-2}\nonumber\\
	    		&+16(\theta_0-1)^2|\vec{A}_1|^2\s{\vec{A}_1}{\bar{\vec{C}_1}}\s{\vec{A}_1}{\vec{A}_1}z^{4\theta_0-6}\z\nonumber
	    		+O(|z|^{4\theta_0-4})\nonumber\\
	    		&=O(|z|^{4\theta_0-4})
	    		\end{align}
	    		thanks of \eqref{addrel} and as by \eqref{newcancel}, we have
	    		\begin{align*}
	    		\s{\vec{A}_1}{\vec{A}_1}\vec{C}_1=0,
	    		\end{align*}
	    		so that
	    		\begin{align}\label{newcancel2}
	    		\s{\vec{A}_1}{\vec{A}_1}\s{\vec{C}_1}{\vec{C}_1}=\s{\vec{A}_1}{\vec{A}_1}\s{\vec{C}_1}{\vec{C}_2}=0.
	    		\end{align}

	    		Then, we compute
	    		\begin{align}\label{remove3}
	    		&\left(\partial\bar{\partial}\h_0\totimes\partial\h_0\right)\otimes\left(\bar{\partial}\h_0\totimes\h_0\right)=\left(-(\theta_0-1)(\theta_0-3)\s{\vec{A}_1}{\vec{C}_2}z^{\theta_0-2}\z^{\theta_0-1}+O(|z|^{2\theta_0-2})\right)\times\nonumber\\
	    		& \left(-\left((\theta_0-3)\s{\vec{A}_1}{\vec{C}_2}+2(\theta_0-2)\s{\vec{A}_2}{\vec{C}_1}\right)z^{\theta_0}\z^{\theta_0-1}+\frac{(\theta_0-2)^2}{4\theta_0}\s{\vec{C}_1}{\vec{C}_1}\z^{2\theta_0-1}+O(|z|^{2\theta_0})\right)\nonumber\\
	    		&=O(|z|^{2\theta_0-3}\cdot |z|^{2\theta_0-1})=O(|z|^{4\theta_0-4}).
	    		\end{align}	    		
	    		The next term is thanks of \eqref{newcancel} and \eqref{newcancel2}
	    		\begin{align}\label{remove4}
	    		&\left(\partial\bar{\partial}\h_0\totimes\bar{\partial}\h_0\right)\otimes\left(\partial\h_0\totimes\h_0\right)=\left(\frac{1}{4}(\theta_0-2)(\theta_0-3)\s{\vec{C}_1}{\vec{C}_2}\z^{2\theta_0-2}+8\theta_0(\theta_0+1)|\vec{A}_1|^2\bar{\alpha_2}z^{2\theta_0-3}\z+O(|z|^{2\theta_0-1})\right)\times\nonumber\\
	    		&\left(4(\theta_0-1)\s{\vec{A}_1}{\vec{A}_1}z^{2\theta_0-3}+O(|z|^{2\theta_0-2}\right)\nonumber\\
	    		&=(\theta_0-1)(\theta_0-2)(\theta_0-3)\ccancel{\s{\vec{A}_1}{\vec{A}_1}\s{\vec{C}_1}{\vec{C}_2}}z^{2\theta_0-3}\z^{2\theta_0-2}+16(\theta_0-1)|\vec{A}_1|^2\ccancel{\s{\vec{A}_1}{\bar{\vec{C}_1}}\s{\vec{A}_1}{\vec{A}_1}}z^{4\theta_0-6}\z+O(|z|^{4\theta_0-4})\nonumber\\
	    		&=O(|z|^{4\theta_0-4})
	    		\end{align}
	    		Finally, we have
	    		\begin{align}\label{remove5}
	    		&\left(\partial\bar{\partial}\h_0\totimes\h_0\right)\otimes\left(\partial\h_0\totimes\bar{\partial}\h_0\right)=\left(-(\theta_0-3)\s{\vec{A}_1}{\vec{C}_2}|z|^{2\theta_0-2}+O(|z|^{2\theta_0-1}\right)\times\nonumber\\
	    		&\left(-\left((\theta_0-1)(\theta_0-3)\s{\vec{A}_1}{\vec{C}_2}+2\theta_0(\theta_0-2)\s{\vec{A}_2}{\vec{C}_1}|z|^{2\theta_0-2}+O(|z|^{2\theta_0-1})\right)\right)\nonumber\\
	    		&=O(|z|^{4\theta_0-4})
	    		\end{align}
	    		Finally, as $g=|z|^{2\theta_0-2}(1+O(|z|))$, we have
	    		\begin{align*}
	    		g^{-2}=\frac{|z|^{4-4\theta_0}}{|dz|^4}{(1+O(|z|^2))}
	    		\end{align*}
	    		and thanks of \eqref{remove1}, \eqref{remove2}, \eqref{remove3}, \eqref{remove4}, \eqref{remove5}, we have
	    		\begin{align}\label{octicrem1}
	    		&O(\h_0)=g^{-2}\,\otimes\bigg\{\frac{1}{4}\left(\partial\bar{\partial}\h_0\,\dot{\otimes}\,\partial\bar{\partial}\h_0\right)\otimes \left(\h_0\,\dot{\otimes}\,\h_0\right)+\frac{1}{4}\left(\partial\h_0\,\dot{\otimes}\, \partial\h_0\right)\otimes \left(\bar{\partial}\h_0\,\dot{\otimes}\, \bar{\partial}\h_0\right)\nonumber\\
	    		&-\frac{1}{2}\left(\partial\bar{\partial}\h_0\totimes\partial\h_0\right)\otimes\left(\bar{\partial}\h_0\totimes\h_0\right)-\frac{1}{2}\left(\partial\bar{\partial}\h_0\totimes\bar{\partial}\h_0\right)\otimes\left(\partial\h_0\totimes\h_0\right)+\frac{1}{2}\left(\partial\bar{\partial}\h_0\totimes\h_0\right)\otimes\left(\partial\h_0\totimes\bar{\partial}\h_0\right) \bigg\}\in L^{\infty}(D^2)
	    		\end{align} 
	    		Now, if we suppose that $\mathscr{O}_{\phi}$ is holomorphic, then
	    		   	
	    		Now, recall that
	    		\begin{align}\label{octic2}
	    		&\frac{1}{4}|\H|^2\,g^{-1}\otimes\bigg\{\frac{1}{2}\left(\partial^N\bar{\partial}^N\h_0\totimes\h_0\right)\otimes\left(\h_0\totimes\h_0\right)-\left(\partial^N\h_0\totimes\h_0\right)\otimes\left(\bar{\partial}^N\h_0\totimes\h_0 \right) +\frac{1}{2}\left(\partial^N\h_0\totimes\bar{\partial}^N\h_0\right)\otimes\left(\h_0\totimes\h_0\right)\bigg\}\nonumber\\
	    		&=\frac{1}{4}|\H|^2\,g^{-1}\otimes\bigg\{\frac{1}{2}\left(\partial\bar{\partial}\h_0\totimes\h_0\right)\otimes\left(\h_0\totimes\h_0\right)-\left(\partial\h_0\totimes\h_0\right)\otimes\left(\bar{\partial}\h_0\totimes\h_0 \right) +\frac{1}{2}\left(\partial\h_0\totimes\bar{\partial}\h_0\right)\otimes\left(\h_0\totimes\h_0\right)\bigg\}
	    		\end{align}
	    		We compute thanks of \eqref{h010}
	    		\begin{align*}
	    		&\left(\partial\bar{\partial}\h_0\totimes\h_0\right)\otimes\left(\h_0\totimes\h_0\right)=\left(-(\theta_0-3)\s{\vec{A}_1}{\vec{C}_2}|z|^{2\theta_0-2}+O(|z|^{2\theta_0-1})\right)\times \left(4\s{\vec{A}_1}{\vec{A}_1}z^{2\theta_0-2}+O(|z|^{2\theta_0-1})\right)\\
	    		&=O(|z|^{4\theta_0-4}).\\
	    		&\left(\partial\h_0\totimes\bar{\partial}\h_0\right)\otimes\left(\bar{\partial}\h_0\totimes\h_0\right)=\left(4(\theta_0-1)\s{\vec{A}_1}{\vec{A}_1}z^{2\theta_0-3}+O(|z|^{2\theta_0-2})\right)\\
	    		&\times\left(-\left((\theta_0-3)\s{\vec{A}_1}{\vec{C}_2}+2(\theta_0-2)\s{\vec{A}_2}{\vec{C}_1}\right)z^{\theta_0}\z^{\theta_0-1}+\frac{(\theta_0-2)^2}{4\theta_0}\s{\vec{C}_1}{\vec{C}_1}\z^{2\theta_0-1}+O(|z|^{2\theta_0})\right)=O(|z|^{4\theta_0-4})\\
	    		&\left(\partial\h_0\totimes\bar{\partial}\h_0\right)\otimes\left(\h_0\totimes\h_0\right)=O(|z|^{4\theta_0-4})
	    		\end{align*}
	    		and as
	    		\begin{align*}
	    		\H=\Re\left(\frac{\vec{C}_1}{z^{\theta_0-2}}\right)+O(|z|^{3-\theta_0})=O(|z|^{2-\theta_0})
	    		\end{align*}
	    		we have $|\H|^2=O(|z|^{4-2\theta_0})$, so that
	    		\begin{align*}
	    		\frac{1}{4}|\H|^2\,g^{-1}\otimes\bigg\{\frac{1}{2}\left(\partial\bar{\partial}\h_0\totimes\h_0\right)\otimes\left(\h_0\totimes\h_0\right)-\left(\partial\h_0\totimes\h_0\right)\otimes\left(\bar{\partial}\h_0\totimes\h_0 \right) +\frac{1}{2}\left(\partial\h_0\totimes\bar{\partial}\h_0\right)\otimes\left(\h_0\totimes\h_0\right)\bigg\}=O(|z|^{2})
	    		\end{align*}
	    		and \textit{a fortiori}
	    		\begin{align}\label{octicrem2}
	    		\frac{1}{4}|\H|^2\,g^{-1}\otimes\bigg\{\frac{1}{2}\left(\partial\bar{\partial}\h_0\totimes\h_0\right)\otimes\left(\h_0\totimes\h_0\right)-\left(\partial\h_0\totimes\h_0\right)\otimes\left(\bar{\partial}\h_0\totimes\h_0 \right) +\frac{1}{2}\left(\partial\h_0\totimes\bar{\partial}\h_0\right)\otimes\left(\h_0\totimes\h_0\right)\bigg\}\in L^{\infty}(D^2).
	    		\end{align}
	    		
	    		One can check the other bounds similarly, 
	    		which proves the holomorphy of $\mathscr{O}_{\phi}$ once it is meromorphic.
	    		
	    		\chapter{The special cases of low multiplicity $\theta_0=2,3,4$}
	    		
	    		\section{The case where $\theta_0=4$}
	    		
	    		In this case, we already have by the holomorphy of the quartic form the relation
	    		\begin{align}\label{4strange11}
	    			|\vec{A}_1|^2\s{\vec{A}_1}{\vec{C}_1}=\s{\bar{\vec{A}_1}}{\vec{C}_1}\s{\vec{A}_1}{\vec{A}_1}
	    		\end{align}
	    		and we can show that the coefficient in $z^{\theta_0}\z^{2-\theta_0}dz^4=z^4\z^{-2}dz^4$ in the quartic form is
	    		\begin{align*}
	    		\Omega_0=\frac{3}{4}\left(|\vec{C}_1|^2\s{\vec{A}_1}{\vec{A}_1}-\s{\vec{A}_1}{\bar{\vec{C}_1}}\s{\vec{A}_1}{\vec{C}_1}\right)=0
	    		\end{align*}
	    		so we obtain
	    		\begin{align}\label{4strange22}
	    		|\vec{C}_1|^2\s{\vec{A}_1}{\vec{A}_1}=\s{\vec{A}_1}{\bar{\vec{C}_1}}\s{\vec{A}_1}{\vec{C}_1}.
	    		\end{align}
	    		Coupling this equation with \eqref{4strange11}, we have recovered the system \begin{align}\label{system3}
	    		\left\{\begin{alignedat}{1}
	    		|\vec{A}_1|^2\s{\vec{A}_1}{\vec{C}_1}&=\s{\bar{\vec{A}_1}}{\vec{C}_1}\s{\vec{A}_1}{\vec{A}_1}\\
	    		|\vec{C}_1|^2\s{\vec{A}_1}{\vec{A}_1}&=\s{\vec{A}_1}{\bar{\vec{C}_1}}\s{\vec{A}_1}{\vec{C}_1},
	    		\end{alignedat}\right.
	    		\end{align}
	    		Remarking that is a linear system in $(\s{\vec{A}_1}{\vec{C}_1},\s{\vec{A}_1}{\vec{A}_1})$, we can recast \eqref{system3} as
	    		\begin{align}\label{salvationmatrix3}
	    		\begin{dmatrix}
	    		|\vec{A}_1|^2 & -\s{\bar{\vec{A}_1}}{\vec{C}_1}\\
	    		-\s{\vec{A}_1}{\bar{\vec{C}_1}} & |\vec{C}_1|^2
	    		\end{dmatrix}
	    		\begin{dmatrix}
	    		\s{\vec{A}_1}{\vec{C}_1}\\
	    		\s{\vec{A}_1}{\vec{A}_1}
	    		\end{dmatrix}=0.
	    		\end{align}
	    		Thanks of Cauchy-Schwarz inequality, we obtain
	    		\begin{align}\label{determinant3}
	    		\det 	\begin{dmatrix}
	    		|\vec{A}_1|^2 & -\s{\bar{\vec{A}_1}}{\vec{C}_1}\\
	    		-\s{\vec{A}_1}{\bar{\vec{C}_1}} & |\vec{C}_1|^2
	    		\end{dmatrix}=|\vec{A}_1|^2|\vec{C}_1|^2-|\s{\vec{A}_1}{\bar{\vec{C}_1}}|^2\geq  0.
	    		\end{align}
	    		and the determinant vanishes if and only if (by the equality case of the Cauchy-Schwarz identity) there exists $\lambda\in\mathbb{C}$ such that
	    		\begin{align}\label{alt1}
	    		\vec{A}_1=\lambda\,\vec{C}_1\quad \text{or}\;\, \vec{C}_1=\lambda \vec{A}_1.
	    		\end{align}
	    		If the determinant in \eqref{determinant3} is strictly positive, then we are done. Therefore, suppose now that \eqref{alt} holds. If $\lambda=0$, we are also done, so we can suppose that $\lambda\neq 0$ and that the relation
	    		\begin{align}\label{4endrel1}
	    		\vec{C}_1=\lambda\vec{A}_1
	    		\end{align}
	    		holds.
	    		Then, the conclusion will follow immediately by remarking that the only non-trivial coefficient in $(-1\;\, 4\;\, 1)$ or
	    		\begin{align*}
	    		\frac{\z^4}{z}\log|z|
	    		\end{align*}
	    		in the Taylor expansion of the quartic form is 
	    		\begin{align}
	    		\Omega_1=-\frac{3}{8}\s{\vec{A}_1}{\vec{C}_1}\bar{\s{\vec{C}_1}{\vec{C}_1}}=\frac{3\bar{\lambda}}{8}|\s{\vec{A}_1}{\vec{C}_1}|^2=0.
	    		\end{align}
	    		As $\lambda\neq 0$, we immediately obtain
	    		\begin{align}\label{4alphaomega1}
	    		\s{\vec{A}_1}{\vec{C}_1}=0.
	    		\end{align}
	    		Therefore, in all cases, \eqref{4alphaomega1} is always satisfied, and we are done in the case $\theta_0=4$.
	    		
	    		Now, let us present the details of this argument.

	    		As we will have to integrate functions of the type
	    		\begin{align*}
	    		z^\alpha\z^\beta\log^p|z|
	    		\end{align*}
	    		for $\alpha,b\in\Z$, with $\alpha\neq 1$ and $p\in \N$, let $a_0^p,\cdots,a_p^p\in \R$ and
	    		\begin{align*}
	    		f_p(z)=\frac{1}{\alpha+1}z^{\alpha+1}\z^\beta\left(\sum_{j=0}^p a_j^p\log^{p-j}|z|\right)
	    		\end{align*}
	    		such that
	    		\begin{align*}
	    		\p{z}f_p(z)=z^{\alpha}\z^{\beta}\log^p|z|.
	    		\end{align*}
	    		As
	    		\begin{align*}
	    		\p{z}\left(\frac{1}{\alpha+1}z^{\alpha+1}\z^{\beta}\log|z|\right)=z^{\alpha}\z^{\beta}\log|z|+\frac{1}{2(\alpha+1)}z^{\alpha}\z^{\beta},
	    		\end{align*}
	    		we deduce that
	    		\begin{align*}
	    		f_1(z)=\frac{1}{\alpha+1}\z^{\alpha+1}\z^{\beta}\left(\log|z|-\frac{1}{2(\alpha+1)}\right),
	    		\end{align*}
	    		so
	    		\begin{align*}
	    		a_{0}^1=1,\quad a_1^1=-\frac{1}{2(\alpha+1)}.
	    		\end{align*}
	    		In the general case, we have
	    		\begin{align*}
	    		\p{z}f_p(z)&=z^{\alpha}\z^{\beta}\left(\sum_{j=0}^{p}a_j^p\log^{p-j}(z)\right)+\frac{1}{2(\alpha+1)}z^{\alpha}\z^{\beta}\left(\sum_{j=0}^{p-1}(p-j)a_j^{p}\log^{p-j-1}|z|\right)\\
	    		&=z^{\alpha}\z^{\beta}a_0^p\log^p|z|+z^{\alpha}\z^{\beta}\sum_{j=1}^{p}\left(a_j^p+\frac{(p-j+1)}{2(\alpha+1)}a_{j-1}^p\right)\log^{p-j}|z|\\
	    		&=z^{\alpha}\z^{\beta}\log^p|z|
	    		\end{align*}
	    		if and only if
	    		\begin{align}
	    		\left\{\begin{alignedat}{1}
	    		a_0^p&=1\\
	    		a_{j}^p&=-\frac{(p-j+1)}{2(\alpha+1)}a_{j-1}^p,\quad j=1,\cdots,p.
	    		\end{alignedat}\right.
	    		\end{align}
	    		In particular, we obtain
	    		\begin{align}\label{intlog}
	    		a_j^p&=(-1)^j\frac{p!}{(p-j)!}\frac{1}{2^j(\alpha+1)^j},\quad j=0,\cdots, p.
	    		\end{align}
	    		As for $\alpha=-1$, we trivially have
	    		\begin{align*}
	    		f_{-1}(z)=\frac{2}{p+1}\z^{\beta}\log^{p+1}|z|,
	    		\end{align*}
	    		one can check that the function \texttt{intz} precisely gives these formulae. Indeed, fixing some arbitrary any $a_1,a_2,a_4,b_1,b_2,b_3,b_4\in \Z$ for example
	    		\begin{align*}
	    		h(z)=\lambda_1z^{a_1}\z^{b_1}\log^4|z|+\lambda_2\z^{b_2}\frac{\log^2|z|}{z}+\lambda_3 z^{a_3}\z^{b_3}\log|z|+\lambda_4 z^{a_4}z^{b_4}
	    		\end{align*}
	    		then a primitive $f$ of $h$ is given thanks of \eqref{intlog} by
	    		\begin{align}\label{trivialog}
	    		f(z)&=\frac{\lambda_1}{a_1+1}z^{a+1}\z^{b_1}\left(\log^4|z|-\frac{2}{a_1+1}\log^3|z|+\frac{3}{(a_1+1)^2}\log^2|z|-\frac{3}{(a_1+1)^3}\log|z|+\frac{3}{2(a_1+1)^4}\right)\nonumber\\
	    		&+\frac{2}{3}\lambda_2\z^{b_2}\log^3|z|+\frac{\lambda_3}{a_3+1}z^{a_3+1}\z^{b_3}\left(\log|z|-\frac{1}{2(a_3+1)}\right)+\frac{\lambda_4}{a_4+1}z^{a_4+1}\z^{b_4}.
	    		\end{align}
	    		and using the code
	    		\begin{verbatim}
	    		var('lambda1,lambda2,lambda3,lambda4',domain='complex'),var('a,b,p', domain='real')
	    		
	    		def intz(v):
	    		length=0
	    		for i in range(v.nrows()):
	    		if bool((v[i,v.ncols()-3]+1).is_zero()): # if the coefficient is z^{-1}\z^{b}\log^p|z|,\\
	    		the primitive has only one components
	    		length=length+1
	    		else:
	    		length=length+v[i,v.ncols()-1]+1 # if the coefficient is z^{a}\z^{b}\log^p|z|\\
	    		with a\neq -1, then the primitive has p+1 components   
	    		m=matrix(SR,length,v.ncols())
	    		n=0
	    		if v.ncols()==4:
	    		for i in range(v.nrows()):
	    		if (v[i,1]+1).is_zero():        #integration of 1/z
	    		m[i+n,0]=2*v[i,0]/(v[i,3]+1)
	    		m[i+n,1]=0   
	    		m[i+n,2]=v[i,2]
	    		m[i+n,3]=v[i,3]+1
	    		else:
	    		if v[i,3].is_zero():
	    		m[i+n,0]=v[i,0]/(v[i,1]+1)
	    		m[i+n,1]=v[i,1]+1
	    		m[i+n,2]=v[i,2]
	    		else: 
	    		a0=1/(v[i,1]+1)
	    		m[i+n,0]=a0*v[i,0]
	    		m[i+n,1]=v[i,1]+1
	    		m[i+n,2]=v[i,2]
	    		m[i+n,3]=v[i,3]
	    		n=n+1
	    		for j in range(v[i,3]):
	    		a0=-(v[i,3]-j)/(2*(v[i,1]+1))*a0
	    		m[i+n,0]=a0*v[i,0]
	    		m[i+n,1]=v[i,1]+1
	    		m[i+n,2]=v[i,2]
	    		m[i+n,3]=v[i,3]-j-1
	    		n=n+1
	    		n=n-1
	    		return m
	    		else:
	    		for i in range(v.nrows()):
	    		if (v[i,2]+1).is_zero():        #integration of 1/z
	    		m[i+n,0]=2*v[i,0]/(v[i,4]+1)
	    		m[i+n,1]=v[i,1]
	    		m[i+n,2]=0
	    		m[i+n,3]=v[i,3]
	    		m[i+n,4]=v[i,4]+1
	    		else:    
	    		if v[i,4].is_zero():
	    		m[i+n,0]=v[i,0]/(v[i,2]+1)
	    		m[i+n,1]=v[i,1]
	    		m[i+n,2]=v[i,2]+1
	    		m[i+n,3]=v[i,3]
	    		else:
	    		a0=1/(v[i,2]+1)
	    		m[i+n,0]=a0*v[i,0]
	    		m[i+n,1]=v[i,1]
	    		m[i+n,2]=v[i,2]+1
	    		m[i+n,3]=v[i,3]
	    		m[i+n,4]=v[i,4]
	    		n=n+1
	    		for j in range(v[i,4]):
	    		a0=-(v[i,4]-j)/(2*(v[i,2]+1))*a0
	    		m[i+n,0]=a0*v[i,0]
	    		m[i+n,1]=v[i,1]
	    		m[i+n,2]=v[i,2]+1
	    		m[i+n,3]=v[i,3]
	    		m[i+n,4]=v[i,4]-j-1
	    		n=n+1
	    		n=n-1                       
	    		return m
	    		
	    		m1=matrix([[lambda1,a,b,4],[lambda2,-1,b,2],[lambda3,a,b,1],[lambda4,a,b,0]])    
	    		
	    		latex(intz(m1))
	    		\end{verbatim}
	    		yields
	    		\small
	    		\begin{align*}
	    		f(z)=\begin{dmatrix}
	    		\frac{\lambda_{1}}{a_{1} + 1} & a_{1} + 1 & b_{1} & 4 \\
	    		-\frac{2 \, \lambda_{1}}{{\left(a_{1} + 1\right)}^{2}} & a_{1} + 1 & b_{1} & 3 \\
	    		\frac{3 \, \lambda_{1}}{{\left(a_{1} + 1\right)}^{3}} & a_{1} + 1 & b_{1} & 2 \\
	    		-\frac{3 \, \lambda_{1}}{{\left(a_{1} + 1\right)}^{4}} & a_{1} + 1 & b_{1} & 1 \\
	    		\frac{3 \, \lambda_{1}}{2 \, {\left(a_{1} + 1\right)}^{5}} & a_{1} + 1 & b_{1} & 0 
	    		\end{dmatrix}
	    		\begin{dmatrix}
	    		\frac{2}{3} \, \lambda_{2} & 0 & b_{2} & 3 \\
	    		\frac{\lambda_{3}}{a_{3} + 1} & a_{3} + 1 & b_{3} & 1 \\
	    		-\frac{\lambda_{3}}{2 \, {\left(a_{3} + 1\right)}^{2}} & a_{3} + 1 & b_{3} & 0 \\
	    		\frac{\lambda_{4}}{a_{4} + 1} & a_{4} + 1 & b_{4} & 0
	    		\end{dmatrix}
	    		\end{align*}
	    		\normalsize
	    		which coincides with \eqref{trivialog}.
	    		
	    		From now on, we assume that $\theta_0=4$. We remark that the developments \eqref{3rddev}, \eqref{newdevg}, \eqref{newdevh0} and \eqref{devH0} and the relation $\alpha_0=0$ (see \eqref{alpha0}) hold for all $\theta_0\geq 4$, so that
	    		\small
	    		\begin{align*}
	    		\left\{\begin{alignedat}{1}
	    		\p{z}\phi&=\begin{dmatrix}
	    		1 & A_{0} & \theta_{0} - 1 & 0 \\
	    		1 & A_{1} & \theta_{0} & 0 \\
	    		1 & A_{2} & \theta_{0} + 1 & 0 \\
	    		\frac{1}{4 \, \theta_{0}} & C_{1} & 1 & \theta_{0} \\
	    		\frac{1}{8} & \overline{C_{1}} & \theta_{0} - 1 & 2
	    		\end{dmatrix}\\
	    		\h_0&=\begin{dmatrix}
	    		2 & A_{1} & \theta_{0} - 1 & 0 \\
	    		-4 \, {\left| A_{1} \right|}^{2} & A_{0} & \theta_{0} - 1 &
	    		1 \\
	    		4 & A_{2} & \theta_{0} & 0 \\
	    		-\frac{\theta_{0} - 2}{2 \, \theta_{0}} & C_{1} & \theta_0-1 & 0\\
	    		- 4 \, \alpha_{1} & A_{0} & \theta_{0} & 0
	    		\theta_{0}
	    		\end{dmatrix}
	    		\end{alignedat}\right.
	    		\quad
	    		\left\{\begin{alignedat}{1}
	    		g&=\begin{dmatrix}
	    		1 & \theta_{0} - 1 & \theta_{0} - 1 \\
	    		2 \, {\left| A_{1} \right|}^{2} & \theta_{0} & \theta_{0}\\
	    		\alpha_1 & \theta_0 & \theta_0-1\\
	    		\bar{\alpha_1} & \theta_0-1 & \theta_0
	    		\end{dmatrix}\\
	    		\H&=\begin{dmatrix}
	    		\frac{1}{2} & C_{1} & -\theta_{0} + 2 & 0 \\
	    		\frac{1}{2} & \overline{C_{1}} & 0 & -\theta_{0} + 2 \\
	    		\frac{1}{2} & C_{2} & -\theta_{0} + 3 & 0 \\
	    		\frac{1}{2} & \overline{C_{2}} & 0 & -\theta_{0} + 3 \\
	    		\frac{1}{2} & B_{1} & -\theta_{0} + 2 & 1 \\
	    		\frac{1}{2} & \overline{B_{1}} & 1 & -\theta_{0} + 2
	    		\end{dmatrix}
	    		\end{alignedat}\right.
	    		\end{align*}
	    		\normalsize
	    		where by \eqref{BC2},
	    		\begin{align}\label{BC21}
	    		\left\{\begin{alignedat}{1}
	    		\vec{C}_2&=\vec{D}_2+\frac{2}{\theta_0-3}\s{\vec{A}_1}{\vec{C}_1}\bar{\vec{A}_0},\\ \vec{B}_1&=-2\s{\bar{\vec{A}_1}}{\vec{C}_1}\vec{A}_0
	    		\end{alignedat}\right.
	    		\end{align}
	    		Taking $\theta_0=4$, we obtain
	    		\begin{align}\label{coeffB1}
	    		\left\{\begin{alignedat}{1}
	    		\vec{C}_2&=\vec{D}_2+2\s{\vec{A}_1}{\vec{C}_1}\bar{\vec{A}_0},\\ \vec{B}_1&=-2\s{\bar{\vec{A}_1}}{\vec{C}_1}\vec{A}_0
	    		\end{alignedat}\right.
	    		\end{align}
	    		and 
	    		\begin{align*}
	    		\left\{\begin{alignedat}{2}
	    		&\p{z}\phi=\begin{dmatrix}
	    		1 & A_{0} & 3 & 0 & 0 \\
	    		1 & A_{1} & 4 & 0 & 0 \\
	    		1 & A_{2} & 5 & 0 & 0 \\
	    		\frac{1}{16} & C_{1} & 1 & 4 & 0 \\
	    		\frac{1}{8} & \overline{C_{1}} & 3 & 2 & 0
	    		\end{dmatrix}\quad
	    		&&g=\begin{dmatrix}
	    		1 & 3 & 3 & 0 \\
	    		2 \, {\left| A_{1} \right|}^{2} & 4 & 4 & 0 \\
	    		\alpha_{1} & 5 & 3 & 0 \\
	    		\overline{\alpha_{1}} & 3 & 5 & 0
	    		\end{dmatrix}\\
	    		&H=\begin{dmatrix}
	    		\frac{1}{2} & C_{1} & -2 & 0 & 0 \\
	    		\frac{1}{2} & \overline{C_{1}} & 0 & -2 & 0 \\
	    		\frac{1}{2} & C_{2} & -1 & 0 & 0 \\
	    		\frac{1}{2} & \overline{C_{2}} & 0 & -1 & 0 \\
	    		\frac{1}{2} & B_{1} & -2 & 1 & 0 \\
	    		\frac{1}{2} & \overline{B_{1}} & 1 & -2 & 0
	    		\end{dmatrix}
	    		&&\h_0=\begin{dmatrix}
	    		2 & A_{1} & 3 & 0 & 0 \\
	    		4 & A_{2} & 4 & 0 & 0 \\
	    		-\frac{1}{4} & C_{1} & 0 & 4 & 0 \\
	    		-4 \, {\left| A_{1} \right|}^{2} & A_{0} & 3 & 1 & 0 \\
	    		-4 \, \alpha_{1} & A_{0} & 4 & 0 & 0
	    		\end{dmatrix}
	    		\end{alignedat}\right.
	    		\end{align*}
	    		In some sense, the result is much more mysterious when we fix the parameter $\theta_0$, as the coefficients allowing the different cancellations which were polynomials in $\theta_0$ become random constants. This explains why the computer is so convenient in this special case.
	    		
	    		\textbf{In the following expressions, the fifth column will indicate the logarithmic power, as indicated in remark \ref{logarithm}.}
	    		
	    		If $\vec{Q}\in C^{\infty}(D^2\setminus\ens{0},\mathbb{C}^n)$ is as usual the anti-meromorphic free function such that
	    		\begin{align*}
	    		\partial\vec{Q}=-|\H|^2\partial\phi-2\,g^{-1}\otimes\s{\H}{\h_0}\otimes \bar{\partial}\phi,
	    		\end{align*}
	    		then we have
	    		\small
	    		\begin{align}\label{4Q1}
	    		\vec{Q}=\begin{dmatrix}
	    			-\frac{1}{16} \, C_{1}^{2} & \overline{A_{0}} & -4 & 4 & 0 \\
	    			-\frac{1}{8} \, C_{1} \overline{C_{1}} & \overline{A_{0}} & -2 & 2 & 0 \\
	    			2 \, A_{1} C_{1} & \overline{A_{0}} & -1 & 0 & 0 \\
	    			2 \, A_{1} C_{1} & \overline{A_{1}} & -1 & 1 & 0 \\
	    			-4 \, A_{0} C_{1} {\left| A_{1} \right|}^{2} + 2 \, A_{1} B_{1} & \overline{A_{0}} & -1 & 1 & 0 \\
	    			-\frac{1}{2} \, C_{1}^{2} & A_{0} & 0 & 0 & 1 \\
	    			8 \, A_{0} C_{1} \alpha_{1} - 8 \, A_{2} C_{1} - 4 \, A_{1} C_{2} & \overline{A_{0}} & 0 & 0 & 1 
	    			\end{dmatrix}
	    			\begin{dmatrix}
	    			-2 \, A_{1} \overline{C_{1}} & \overline{A_{0}} & 1 & -2 & 0 \\
	    			-2 \, A_{1} \overline{C_{1}} & \overline{A_{1}} & 1 & -1 & 0 \\
	    			4 \, A_{0} {\left| A_{1} \right|}^{2} \overline{C_{1}} - 2 \, A_{1} C_{2} & \overline{A_{0}} & 1 & -1 & 0 \\
	    			-\frac{1}{4} \, C_{1} \overline{C_{1}} & A_{0} & 2 & -2 & 0 \\
	    			2 \, A_{0} \alpha_{1} \overline{C_{1}} - A_{1} \overline{B_{1}} - 2 \, A_{2} \overline{C_{1}} & \overline{A_{0}} & 2 & -2 & 0 \\
	    			-\frac{1}{16} \, \overline{C_{1}}^{2} & A_{0} & 4 & -4 & 0
	    			\end{dmatrix}
	    		\end{align}
	    		\normalsize
	    		and as
	    		\begin{align*}
	    		\partial\left(\H-2i\vec{L}+\vec{\gamma}_0\log|z|\right)=-|\H|^2\partial\phi-2\,g^{-1}\otimes\s{\H}{\h_0}\otimes \bar{\partial}\phi,
	    		\end{align*}
	    		there exists $\vec{D}_2,\vec{D}_3\in \mathbb{C}^n$ such that
	    		\begin{align*}
	    		\vec{H}-2i\vec{L}+\vec{\gamma}_0\log|z|=\frac{\bar{\vec{C}_1}}{\z^2}+\frac{\bar{\vec{D}_2}}{\z}+\bar{\vec{D}_3}+\vec{Q}+O(|z|^{1-\epsilon})
	    		\end{align*}
	    		and
	    		\begin{align*}
	    		\H+\vec{\gamma}_0\log|z|=\Re\left(\frac{\vec{C}_1}{z^2}+\frac{\vec{D}_2}{z}+\vec{D}_3+\vec{Q}\right)+O(|z|^{1-\epsilon})
	    		\end{align*}
	    		so we obtain
	    		\small
	    		\begin{align}\label{4H1}
	    		&\H+\vec{\gamma}_0\log|z|=\Re\left(\frac{\vec{C}_1}{z^2}+\frac{\vec{D}_2}{z}+\vec{D}_3\right)+\\
	    		&\begin{dmatrix}
	    		-\frac{1}{16} \, C_{1}^{2} & \overline{A_{0}} & -4 & 4 & 0 \\
	    		-C_{1} \overline{A_{1}} & A_{0} & -2 & 1 & 0 \\
	    		-\frac{3}{16} \, C_{1} \overline{C_{1}} & \overline{A_{0}} & -2 & 2 & 0 \\
	    		\ccancel{C_{1} \overline{A_{0}} \overline{\alpha_{1}}} - \frac{1}{2} \ccancel{\, B_{1} \overline{A_{1}}} - C_{1} \overline{A_{2}} & A_{0} & -2 & 2 & 0 \\
	    		A_{1} C_{1} & \overline{A_{0}} & -1 & 0 & 0 \\
	    		A_{1} C_{1} & \overline{A_{1}} & -1 & 1 & 0 \\
	    		-\ccancel{2 \, A_{0} C_{1} {\left| A_{1} \right|}^{2}} + \ccancel{A_{1} B_{1}} & \overline{A_{0}} & -1 & 1 & 0 \\
	    		-C_{1} \overline{A_{1}} & A_{1} & -1 & 1 & 0 \\
	    		\ccancel{2 \, C_{1} {\left| A_{1} \right|}^{2} \overline{A_{0}}} - C_{2} \overline{A_{1}} & A_{0} & -1 & 1 & 0 \\
	    		\overline{A_{1}} \overline{C_{1}} & A_{0} & 0 & -1 & 0 \\
	    		\ccancel{4 \, \overline{A_{0}} \overline{C_{1}} \overline{\alpha_{1}}} - \frac{1}{4} \, C_{1}^{2} - 2 \, C_{2} \overline{A_{1}} - 4 \, \overline{A_{2}} \overline{C_{1}} & A_{0} & 0 & 0 & 1 \\
	    		\ccancel{4 \, A_{0} C_{1} \alpha_{1}} - 4 \, A_{2} C_{1} - 2 \, A_{1} C_{2} - \frac{1}{4} \, \overline{C_{1}}^{2} & \overline{A_{0}} & 0 & 0 & 1 
	    	    \end{dmatrix}
	    		\begin{dmatrix}
	    		-A_{1} \overline{C_{1}} & \overline{A_{0}} & 1 & -2 & 0 \\
	    		\overline{A_{1}} \overline{C_{1}} & A_{1} & 1 & -1 & 0 \\
	    		-\ccancel{2 \, {\left| A_{1} \right|}^{2} \overline{A_{0}} \overline{C_{1}}} + \ccancel{\overline{A_{1}} \overline{B_{1}}} & A_{0} & 1 & -1 & 0 \\
	    		-A_{1} \overline{C_{1}} & \overline{A_{1}} & 1 & -1 & 0 \\
	    		\ccancel{2 \, A_{0} {\left| A_{1} \right|}^{2} \overline{C_{1}}} - A_{1} C_{2} & \overline{A_{0}} & 1 & -1 & 0 \\
	    		-\frac{3}{16} \, C_{1} \overline{C_{1}} & A_{0} & 2 & -2 & 0 \\
	    		\ccancel{A_{0} \alpha_{1} \overline{C_{1}}} - \frac{1}{2} \ccancel{\, A_{1} \overline{B_{1}}} - A_{2} \overline{C_{1}} & \overline{A_{0}} & 2 & -2 & 0 \\
	    		-\frac{1}{16} \, \overline{C_{1}}^{2} & A_{0} & 4 & -4 & 0
	    		\end{dmatrix}
	    		\end{align}
	    		\normalsize
	    		There are some cancellations, here, as
	    		\begin{align*}
	    		\s{\vec{A}_0}{\vec{A}_1}=\s{\vec{A}_0}{\bar{\vec{A}_1}}=0,\quad \s{\vec{A}_0}{\vec{C}_1}=\s{\vec{A}_0}{\bar{\vec{C}_1}}=0
	    		\end{align*}
	    		and thanks of \eqref{coeffB1}, 
	    		\begin{align*}
	    		\s{\vec{A}_1}{\vec{B}_1}=\s{\bar{\vec{A}_1}}{\vec{B}_1}=0
	    		\end{align*}
	    		so with (we just need to multiply the corresponding coefficients in \eqref{4H1} by $2$)
	    		\begin{align}\label{4consts}
	    		\left\{\begin{alignedat}{1}
	    		\vec{C}_2&=\vec{D}_2+2\s{\vec{A}_1}{\vec{C}_1}\bar{\vec{A}_0}\\
	    		\vec{C}_3&=\Re(\vec{D}_3)\\
	    		\vec{B}_1&=-2\s{\bar{\vec{A}_1}}{\vec{C}_1}\vec{A}_0\\
	    		\vec{B}_2&=-\frac{3}{8}|\vec{C}_1|^2\bar{\vec{A}_0}-2\s{\bar{\vec{A}_2}}{\vec{C}_1}\vec{A}_0\\
	    		\vec{B}_3&=-2\s{\bar{\vec{A}_1}}{\vec{C}_1}\vec{A}_1+2\s{\vec{A}_1}{\vec{C}_1}\bar{\vec{A}_1}-2\s{\bar{\vec{A}_1}}{\vec{C}_2}\vec{A}_0\\
	    		\vec{E}_1&=-\frac{1}{8}\s{\vec{C}_1}{\vec{C}_1}\bar{\vec{A}_0}\\
	    		\vec{\gamma}_1&=-\vec{\gamma}_0-\Re\left\{\left(\frac{1}{2}\s{\vec{C}_1}{\vec{C}_1}+8\bar{\s{\vec{A}_2}{\vec{C}_1}}+4\bar{\s{\vec{A}_1}{\vec{C}_2}}\right)\vec{A}_0\right\}\in\R^n
	    		\end{alignedat}\right.
	    		\end{align}
	    		and we obtain 
	    		\begin{align}\label{4eqH1}
	    		\H&=
	    		\Re\left(\frac{\vec{C}_1}{z^{2}}+\frac{\vec{C}_2}{z^{}}+\vec{C}_3+\vec{B}_1\frac{\z}{z^{2}}+\vec{B}_2\frac{\z^2}{z^2}+\vec{B}_3\frac{\z}{z}+\vec{E}_1\frac{\z^{4}}{z^{4}}\right)+\vec{\gamma}_1\log|z|+O(|z|^{1-\epsilon}).
	    		\end{align}
	    		The notations were chosen to be consistent with the higher development of the previous chapters. Translated into code, this is
	    		\small
	    		\begin{align}\label{4H1code}
	    		&\H=\begin{dmatrix}
	    		\frac{1}{2} & C_{1} & -2 & 0 & 0 \\
	    		\frac{1}{2} & C_{2} & -1 & 0 & 0 \\
	    		1 & C_{3} & 0 & 0 & 0 \\
	    		\frac{1}{2} & B_{1} & -2 & 1 & 0 \\
	    		\frac{1}{2} & B_{2} & -2 & 2 & 0 \\
	    		\frac{1}{2} & B_{3} & -1 & 1 & 0 \\
	    		\frac{1}{2} & E_{1} & -4 & 4 & 0 
	    		\end{dmatrix}
	    		\begin{dmatrix}
	    		\frac{1}{2} & \overline{C_{1}} & 0 & -2 & 0 \\
	    		\frac{1}{2} & \overline{C_{2}} & 0 & -1 & 0 \\
	    		\frac{1}{2} & \overline{B_{1}} & 1 & -2 & 0 \\
	    		\frac{1}{2} & \overline{B_{2}} & 2 & -2 & 0 \\
	    		\frac{1}{2} & \overline{B_{3}} & 1 & -1 & 0 \\
	    		\frac{1}{2} & \overline{E_{1}} & 4 & -4 & 0 \\
	    		1 & \gamma_{1} & 0 & 0 & 1
	    		\end{dmatrix}
	    		\end{align}
	    		\normalsize 			    
	    		Now, as
	    		\begin{align*}
	    		\p{\z}\left(\p{z}\phi\right)=\frac{e^{2\lambda}}{2}\H,
	    		\end{align*}
	    		we deduce that
	    		\small
	    		\begin{align}\label{4dzphi1}
	    		&\p{z}\phi=\begin{dmatrix}
	    		1 & A_{0} & 3 & 0 & 0 \\
	    		1 & A_{1} & 4 & 0 & 0 \\
	    		1 & A_{2} & 5 & 0 & 0 \\
	    		1 & A_{3} & 6 & 0 & 0 \\
	    		1 & A_{4} & 7 & 0 & 0 \\
	    		\frac{1}{16} & C_{1} & 1 & 4 & 0 \\
	    		\frac{1}{16} & C_{2} & 2 & 4 & 0 \\
	    		\frac{1}{8} & C_{3} & 3 & 4 & 0 \\
	    		\frac{1}{20} & B_{1} & 1 & 5 & 0 
	    		\end{dmatrix}
	    		\begin{dmatrix}
	    		\frac{1}{24} & B_{2} & 1 & 6 & 0 \\
	    		\frac{1}{20} & B_{3} & 2 & 5 & 0 \\
	    		\frac{1}{32} & E_{1} & -1 & 8 & 0 \\
	    		\frac{1}{8} & \overline{C_{1}} & 3 & 2 & 0 \\
	    		\frac{1}{12} & \overline{C_{2}} & 3 & 3 & 0 \\
	    		\frac{1}{8} & \overline{B_{1}} & 4 & 2 & 0 \\
	    		\frac{1}{8} & \overline{B_{2}} & 5 & 2 & 0 \\
	    		\frac{1}{12} & \overline{B_{3}} & 4 & 3 & 0 \\
	    		\frac{1}{2} & \overline{E_{1}} & 7 & 0 & 1 
	    			    		\end{dmatrix}
	    		\begin{dmatrix}
	    		\frac{1}{8} & \gamma_{1} & 3 & 4 & 1 \\
	    		-\frac{1}{64} & \gamma_{1} & 3 & 4 & 0 \\
	    		\frac{1}{10} \, {\left| A_{1} \right|}^{2} & C_{1} & 2 & 5 & 0 \\
	    		\frac{1}{6} \, {\left| A_{1} \right|}^{2} & \overline{C_{1}} & 4 & 3 & 0 \\
	    		\frac{1}{16} \, \alpha_{1} & C_{1} & 3 & 4 & 0 \\
	    		\frac{1}{8} \, \alpha_{1} & \overline{C_{1}} & 5 & 2 & 0 \\
	    		\frac{1}{24} \, \overline{\alpha_{1}} & C_{1} & 1 & 6 & 0 \\
	    		\frac{1}{16} \, \overline{\alpha_{1}} & \overline{C_{1}} & 3 & 4 & 0
	    		\end{dmatrix}
	    		\end{align}
	    		\normalsize
	    		and by conformality of $\phi$, we obtain
	    		\small
	    		\begin{align}\label{4conf1}
	    		0=\s{\p{z}\phi}{\p{z}\phi}=\begin{dmatrix}
	    		 A_{0}^{2} & 6 & 0 & 0 \\
	    		2 \, A_{0} A_{1} & 7 & 0 & 0 \\
	    		A_{1}^{2} + 2 \, A_{0} A_{2} & 8 & 0 & 0 \\
	    		2 \, A_{1} A_{2} + 2 \, A_{0} A_{3} & 9 & 0 & 0 \\
	    		A_{2}^{2} + 2 \, A_{1} A_{3} + 2 \, A_{0} A_{4} & 10 & 0 & 0 \\
	    		\frac{1}{8} \, A_{0} C_{1} & 4 & 4 & 0 \\
	    		\frac{1}{8} \, A_{1} C_{1} + \frac{1}{8} \, A_{0} C_{2} & 5 & 4 & 0 \\
	    		\frac{1}{8} \, A_{0} C_{1} \alpha_{1} + \frac{1}{8} \, A_{0} \overline{C_{1}} \overline{\alpha_{1}} + \frac{1}{8} \, A_{2} C_{1} + \frac{1}{8} \, A_{1} C_{2} + \frac{1}{4} \, A_{0} C_{3} - \frac{1}{32} \, A_{0} \gamma_{1} + \frac{1}{64} \, \overline{C_{1}}^{2} & 6 & 4 & 0 \\
	    		\frac{1}{10} \, A_{0} B_{1} & 4 & 5 & 0 \\
	    		\frac{1}{12} \, A_{0} C_{1} \overline{\alpha_{1}} + \frac{1}{12} \, A_{0} B_{2} + \frac{1}{64} \, C_{1} \overline{C_{1}} & 4 & 6 & 0 \\
	    		\frac{1}{5} \, A_{0} C_{1} {\left| A_{1} \right|}^{2} + \frac{1}{10} \, A_{1} B_{1} + \frac{1}{10} \, A_{0} B_{3} & 5 & 5 & 0 \\
	    		\frac{1}{256} \, C_{1}^{2} + \frac{1}{16} \, A_{0} E_{1} & 2 & 8 & 0 \\
	    		\frac{1}{4} \, A_{0} \overline{C_{1}} & 6 & 2 & 0 \\
	    		\frac{1}{6} \, A_{0} \overline{C_{2}} & 6 & 3 & 0 \\
	    		\frac{1}{4} \, A_{0} \overline{B_{1}} + \frac{1}{4} \, A_{1} \overline{C_{1}} & 7 & 2 & 0 \\
	    		\frac{1}{4} \, A_{0} \alpha_{1} \overline{C_{1}} + \frac{1}{4} \, A_{1} \overline{B_{1}} + \frac{1}{4} \, A_{0} \overline{B_{2}} + \frac{1}{4} \, A_{2} \overline{C_{1}} & 8 & 2 & 0 \\
	    		\frac{1}{3} \, A_{0} {\left| A_{1} \right|}^{2} \overline{C_{1}} + \frac{1}{6} \, A_{0} \overline{B_{3}} + \frac{1}{6} \, A_{1} \overline{C_{2}} & 7 & 3 & 0 \\
	    		A_{0} \overline{E_{1}} & 10 & 0 & 1 \\
	    		\frac{1}{4} \, A_{0} \gamma_{1} & 6 & 4 & 1
	    		\end{dmatrix}
	    		\end{align} 
	    		\normalsize
	    		We already know that the only interesting relations are
	    		\begin{align}\label{4cancel1}
	    		\left\{\begin{alignedat}{1}
	    		&\s{\vec{A}_0}{\vec{A}_0}=\s{\vec{A}_0}{\vec{A}_1}=\s{\vec{A}_0}{\bar{\vec{A}_1}}=\s{\vec{A}_0}{\vec{\gamma}_1}=0,\\
	    		&\s{\vec{A}_1}{\vec{A}_1}+2\s{\vec{A}_0}{\vec{A}_2}=0,\quad \s{\vec{A}_1}{\vec{A}_2}+\s{\vec{A}_0}{\vec{A}_3}=0\\
	    		&\s{\vec{A}_0}{\vec{C}_1}=\s{\vec{A}_0}{\bar{\vec{C}_1}}=\s{\vec{A}_0}{\bar{\vec{C}_2}}=0,\quad \s{\vec{A}_1}{\vec{C}_1}+\s{\vec{A}_0}{\vec{C}_2}=0.
	    		\end{alignedat}\right.
	    		\end{align}
	    		where we added the previous $\s{\vec{A}_0}{\bar{\vec{A}_1}}=\s{\vec{A}_0}{\bar{\vec{A}_2}}=0$ coming from the upper regularity. The rest of the cancellations in \eqref{4conf1} is either trivial or useless.

	    		and 
	    		\small
	    		\begin{align}\label{4e2lambdatemp}
	    		&e^{2\lambda}=\\
	    		&\begin{dmatrix}
	    		2 \, A_{0} \overline{A_{0}} & 3 & 3 & 0 \\
	    		2 \, A_{0} \overline{A_{1}} & 3 & 4 & 0 \\
	    		2 \, A_{0} \overline{A_{2}} + \frac{1}{4} \, \overline{A_{0}} \overline{C_{1}} & 3 & 5 & 0 \\
	    		2 \, A_{0} \overline{A_{3}} + \frac{1}{4} \, \overline{A_{1}} \overline{C_{1}} + \frac{1}{6} \, \overline{A_{0}} \overline{C_{2}} & 3 & 6 & 0 \\
	    		\frac{1}{8} \, C_{1} \alpha_{1} \overline{A_{0}} + \frac{1}{8} \, \overline{A_{0}} \overline{C_{1}} \overline{\alpha_{1}} + \frac{1}{64} \, C_{1}^{2} + \frac{1}{4} \, C_{3} \overline{A_{0}} - \frac{1}{32} \, \gamma_{1} \overline{A_{0}} + 2 \, A_{0} \overline{A_{4}} + \frac{1}{4} \, \overline{A_{2}} \overline{C_{1}} + \frac{1}{6} \, \overline{A_{1}} \overline{C_{2}} & 3 & 7 & 0 \\
	    		\frac{1}{8} \, A_{0} \overline{C_{1}} & 7 & 1 & 0 \\
	    		\frac{1}{8} \, A_{0} \overline{C_{2}} & 7 & 2 & 0 \\
	    		\frac{1}{8} \, A_{0} C_{1} \alpha_{1} + \frac{1}{8} \, A_{0} \overline{C_{1}} \overline{\alpha_{1}} + \frac{1}{4} \, A_{2} C_{1} + \frac{1}{6} \, A_{1} C_{2} + \frac{1}{4} \, A_{0} C_{3} - \frac{1}{32} \, A_{0} \gamma_{1} + 2 \, A_{4} \overline{A_{0}} + \frac{1}{64} \, \overline{C_{1}}^{2} & 7 & 3 & 0 \\
	    		\frac{1}{10} \, A_{0} \overline{B_{1}} + \frac{1}{8} \, A_{1} \overline{C_{1}} & 8 & 1 & 0 \\
	    		\frac{1}{12} \, A_{0} \alpha_{1} \overline{C_{1}} + \frac{1}{10} \, A_{1} \overline{B_{1}} + \frac{1}{12} \, A_{0} \overline{B_{2}} + \frac{1}{8} \, A_{2} \overline{C_{1}} & 9 & 1 & 0 \\
	    		\frac{1}{5} \, A_{0} {\left| A_{1} \right|}^{2} \overline{C_{1}} + \frac{1}{10} \, A_{0} \overline{B_{3}} + \frac{1}{8} \, A_{1} \overline{C_{2}} & 8 & 2 & 0 \\
	    		\frac{1}{16} \, A_{0} \overline{E_{1}} & 11 & -1 & 0 \\
	    		\frac{1}{4} \, A_{0} C_{1} + 2 \, A_{2} \overline{A_{0}} & 5 & 3 & 0 \\
	    		\frac{1}{4} \, A_{1} C_{1} + \frac{1}{6} \, A_{0} C_{2} + 2 \, A_{3} \overline{A_{0}} & 6 & 3 & 0 \\
	    		\frac{1}{4} \, A_{0} B_{1} + 2 \, A_{2} \overline{A_{1}} & 5 & 4 & 0 \\
	    		\frac{1}{4} \, \alpha_{1} \overline{A_{0}} \overline{C_{1}} + \frac{1}{4} \, A_{0} C_{1} \overline{\alpha_{1}} + \frac{1}{4} \, A_{0} B_{2} + 2 \, A_{2} \overline{A_{2}} + \frac{1}{4} \, \overline{A_{0}} \overline{B_{2}} + \frac{5}{128} \, C_{1} \overline{C_{1}} & 5 & 5 & 0 \\
	    		\frac{1}{3} \, A_{0} C_{1} {\left| A_{1} \right|}^{2} + \frac{1}{4} \, A_{1} B_{1} + \frac{1}{6} \, A_{0} B_{3} + 2 \, A_{3} \overline{A_{1}} & 6 & 4 & 0 \\
	    		A_{0} E_{1} + \frac{1}{4} \, \gamma_{1} \overline{A_{0}} & 3 & 7 & 1 \\
	    		\frac{1}{4} \, A_{0} \gamma_{1} + \overline{A_{0}} \overline{E_{1}} & 7 & 3 & 1 \\
	    		2 \, A_{1} \overline{A_{0}} & 4 & 3 & 0 \\
	    		2 \, A_{1} \overline{A_{1}} & 4 & 4 & 0 \\
	    		2 \, A_{1} \overline{A_{2}} + \frac{1}{4} \, \overline{A_{0}} \overline{B_{1}} & 4 & 5 & 0 \\
	    		\frac{1}{3} \, {\left| A_{1} \right|}^{2} \overline{A_{0}} \overline{C_{1}} + 2 \, A_{1} \overline{A_{3}} + \frac{1}{4} \, \overline{A_{1}} \overline{B_{1}} + \frac{1}{6} \, \overline{A_{0}} \overline{B_{3}} & 4 & 6 & 0 \\
	    		\frac{1}{8} \, C_{1} \overline{A_{0}} & 1 & 7 & 0 \\
	    		\frac{1}{10} \, B_{1} \overline{A_{0}} + \frac{1}{8} \, C_{1} \overline{A_{1}} & 1 & 8 & 0 \\
	    		\frac{1}{12} \, C_{1} \overline{A_{0}} \overline{\alpha_{1}} + \frac{1}{12} \, B_{2} \overline{A_{0}} + \frac{1}{10} \, B_{1} \overline{A_{1}} + \frac{1}{8} \, C_{1} \overline{A_{2}} & 1 & 9 & 0 \\
	    		\frac{1}{8} \, C_{2} \overline{A_{0}} & 2 & 7 & 0 \\
	    		\frac{1}{5} \, C_{1} {\left| A_{1} \right|}^{2} \overline{A_{0}} + \frac{1}{10} \, B_{3} \overline{A_{0}} + \frac{1}{8} \, C_{2} \overline{A_{1}} & 2 & 8 & 0 \\
	    		\frac{1}{16} \, E_{1} \overline{A_{0}} & -1 & 11 & 0
	    		\end{dmatrix}
	    		\begin{dmatrix}
	    		\textbf{(1)}\\
	    		\textbf{(2)}\\
	    		\textbf{(3)}\\
	    		\textbf{(4)}\\
	    		\textbf{(5)}\\
	    		\textbf{(6)}\\
	    		\textbf{(7)}\\
	    		\textbf{(8)}\\
	    		\textbf{(9)}\\
	    		\textbf{(10)}\\
	    		\textbf{(11)}\\
	    		\textbf{(12)}\\
	    		\textbf{(13)}\\
	    		\textbf{(14)}\\
	    		\textbf{(15)}\\
	    		\textbf{(16)}\\
	    		\textbf{(17)}\\
	    		\textbf{(18)}\\
	    		\textbf{(19)}\\
	    		\textbf{(20)}\\
	    		\textbf{(21)}\\
	    		\textbf{(22)}\\
	    		\textbf{(23)}\\
	    		\textbf{(24)}\\
	    		\textbf{(25)}\\
	    		\textbf{(26)}\\
	    		\textbf{(27)}\\
	    		\textbf{(28)}\\
	    		\textbf{(29)}
	    		\end{dmatrix}
	    		\end{align}
	    		\normalsize
	    		We see easily that
	    		\begin{align*}
	    		(2),\;\, (6),\;\, (7),\;\, (12),\;\, (20),\;\, (24),\;\, (27),\;\, (29)
	    		\end{align*}
	    		all vanish. Therefore, for some 
	    		\begin{align*}
	    		\alpha_2,\alpha_3,\alpha_4,\alpha_5,\alpha_6,\alpha_7,\alpha_8,\alpha_9\in\mathbb{C},\quad \zeta_0\in \mathbb{C}, \quad \beta\in \R,
	    		\end{align*}
	    		we have $\left(\text{recall that}\;\, |\vec{A}_0|^2=\dfrac{1}{2}\right)$, with \TeX\;\, on the left, and Sage on the right
	    		\begin{align}\label{4e2lambda1}
	    		&e^{2\lambda}=\begin{dmatrix}
	    		1 & 3 & 3 & 0\\
	    		2|\vec{A}_1|^2 & 4 & 4 & 0 \\
	    		\beta & 5 & 5 & 0\\
	    		\alpha_1 & 5 & 3 & 0\\
	    		\alpha_2 & 1 & 8 & 0\\
	    		\alpha_3 & 6 & 3 & 0\\
	    		\alpha_4 & 5 & 4 & 0\\
	    		\alpha_5 & 6 & 4 & 0\\
	    		\alpha_6 & 8 & 2 & 0\\
	    		\alpha_7 & 9 & 1 & 0\\
	    		\zeta_0 & 7 & 3 & 1
	    		\end{dmatrix}
	    		\begin{dmatrix}
	    		\bar{\alpha_1} & 3 & 5 & 0\\
	    		\bar{\alpha_2} & 8 & 1 & 0\\
	    		\bar{\alpha_3} & 3 & 6 & 0\\
	    		\bar{\alpha_4} & 4 & 5 & 0\\
	    		\bar{\alpha_5} & 4 & 6 & 0\\
	    		\bar{\alpha_6} & 2 & 8 & 0\\
	    		\bar{\alpha_7} & 1 & 9 & 0\\
	    		\bar{\zeta_0} & 3 & 7 & 1\\
	    		\end{dmatrix}
	    		\begin{dmatrix}
	    		1 & 3 & 3 & 0 \\
	    		2 \, {\left| A_{1} \right|}^{2} & 4 & 4 & 0 \\
	    		\beta & 5 & 5 & 0 \\
	    		\alpha_{1} & 5 & 3 & 0 \\
	    		\alpha_{2} & 1 & 8 & 0 \\
	    		\alpha_{3} & 6 & 3 & 0 \\
	    		\alpha_{4} & 5 & 4 & 0 \\
	    		\alpha_{5} & 6 & 4 & 0 \\
	    		\alpha_{6} & 8 & 2 & 0 \\
	    		\alpha_{7} & 9 & 1 & 0 \\
	    		\zeta_{0} & 7 & 3 & 1 \\
	    		\end{dmatrix}
	    		\begin{dmatrix}
	    		\overline{\alpha_{1}} & 3 & 5 & 0 \\
	    		\overline{\alpha_{2}} & 8 & 1 & 0 \\
	    		\overline{\alpha_{3}} & 3 & 6 & 0 \\
	    		\overline{\alpha_{4}} & 4 & 5 & 0 \\
	    		\overline{\alpha_{5}} & 4 & 6 & 0 \\
	    		\overline{\alpha_{6}} & 2 & 8 & 0 \\
	    		\overline{\alpha_{7}} & 1 & 9 & 0 \\
	    		\overline{\zeta_{0}} & 3 & 7 & 1
	    		\end{dmatrix}
	    		\end{align}
	    		At this point, this is easy to see that the first non-trivial logarithm arising in the quartic form will be given by the expression in the beginning of the chapter, so we are done.
	    	    First, $\alpha_3$ is the coefficient in $(6\;\, 3\;\, 0)$ in \eqref{4e2lambdatemp}, which is the coefficient corresponding to the line $(14)$, that is
	    		\begin{align}\label{4alpha0}
	    		\alpha_3=\frac{1}{4} \, A_{1} C_{1} + \frac{1}{6} \, A_{0} C_{2} + 2 \, A_{3} \overline{A_{0}}=\frac{1}{12}\s{\vec{A}_1}{\vec{C}_1}+2\s{\bar{\vec{A}_0}}{\vec{A}_3},
	    		\end{align}
	    		where we have used $\s{\vec{A}_1}{\vec{C}_1}+\s{\vec{A}_0}{\vec{C}_2}=0$ in \eqref{4cancel1}. Then, $\alpha_4$ corresponds to the coefficient in $(5,\;\, 4,\;\, 0)$, \textit{i.e.} to the line $(15)$, giving
	    		\begin{align}\label{4alpha1}
	    		\alpha_4=\frac{1}{4} \, A_{0} B_{1} + 2 \, A_{2} \overline{A_{1}}=2\s{\bar{\vec{A}_1}}{\vec{A}_2}
	    		\end{align}
	    		as $\vec{B}_1\in\mathrm{Span}(\vec{A}_0)$. Finally, $\alpha_2$ corresponds to $(1,\;\, 8,\;\, 0),\;\, (25)$
	    		\begin{align}\label{4alpha2}
	    		\alpha_2=\frac{1}{10} \, B_{1} \overline{A_{0}} + \frac{1}{8} \, C_{1} \overline{A_{1}}=-\frac{1}{10}\s{\bar{\vec{A}_1}}{\vec{C}_1}+\frac{1}{8}\s{\bar{\vec{A}_1}}{\vec{C}_1}=\frac{1}{40}\s{\bar{\vec{A}_1}}{\vec{C}_1}
	    		\end{align}
	    		as
	    		\begin{align*}
	    		\s{\bar{\vec{A}_0}}{\vec{B}_1}=\bs{\bar{\vec{A}_0}}{-2\s{\bar{\vec{A}_1}}{\vec{C}_1}\vec{A}_0}=-\s{\bar{\vec{A}_1}}{\vec{C}_1}.
	    		\end{align*}
	    		Gathering \eqref{4alpha0}, \eqref{4alpha1} and \eqref{4alpha2}, we obtain
	    		\begin{align}\label{4alpha012}
	    		\left\{\begin{alignedat}{1}
	    		\alpha_2&=\frac{1}{40}\s{\bar{\vec{A}_1}}{\vec{C}_1}\\
	    		\alpha_3&=\frac{1}{12}\s{\vec{A}_1}{\vec{C}_1}+2\s{\bar{\vec{A}_0}}{\vec{A}_3}\\
	    		\alpha_4&=2\s{\bar{\vec{A}_1}}{\vec{A}_2}
	    		\end{alignedat}\right.
	    		\end{align}
	    		We now check that our transcription in Sage was made without mistakes. We have
	    		so
	    		\small
	    		\begin{align}
	    		\h_0&=\begin{dmatrix}
	    		 2 & A_{1} & 3 & 0 & 0 \\
	    		4 & A_{2} & 4 & 0 & 0 \\
	    		6 & A_{3} & 5 & 0 & 0 \\
	    		8 & A_{4} & 6 & 0 & 0 \\
	    		-\frac{1}{4} & C_{1} & 0 & 4 & 0 \\
	    		-\frac{1}{8} & C_{2} & 1 & 4 & 0 \\
	    		-\frac{1}{5} & B_{1} & 0 & 5 & 0 \\
	    		-\frac{1}{6} & B_{2} & 0 & 6 & 0 \\
	    		-\frac{1}{10} & B_{3} & 1 & 5 & 0 \\
	    		-\frac{1}{4} & E_{1} & -2 & 8 & 0 \\
	    		\frac{1}{4} & \overline{B_{1}} & 3 & 2 & 0 \\
	    		\frac{1}{2} & \overline{B_{2}} & 4 & 2 & 0 \\
	    		\frac{1}{6} & \overline{B_{3}} & 3 & 3 & 0 \\
	    		4 & \overline{E_{1}} & 6 & 0 & 1 \\
	    		\frac{1}{2} & \overline{E_{1}} & 6 & 0 & 0 
	    		\end{dmatrix}
	    		\begin{dmatrix}
	    		\frac{1}{8} & \gamma_{1} & 2 & 4 & 0 \\
	    		-\frac{9}{20} \, {\left| A_{1} \right|}^{2} & C_{1} & 1 & 5 & 0 \\
	    		-\frac{1}{6} \, {\left| A_{1} \right|}^{2} & \overline{C_{1}} & 3 & 3 & 0 \\
	    		-\frac{1}{4} \, \alpha_{1} & C_{1} & 2 & 4 & 0 \\
	    		-\frac{1}{6} \, \overline{\alpha_{1}} & C_{1} & 0 & 6 & 0 \\
	    		-4 \, {\left| A_{1} \right|}^{2} & A_{0} & 3 & 1 & 0 \\
	    		-4 \, {\left| A_{1} \right|}^{2} & A_{1} & 4 & 1 & 0 \\
	    		-4 \, {\left| A_{1} \right|}^{2} & A_{2} & 5 & 1 & 0 \\
	    		8 \, {\left| A_{1} \right|}^{4} + 4 \, \alpha_{1}  \overline{\alpha_{1}} - 4 \, \beta & A_{0} & 4 & 2 & 0 \\
	    		-4 \, \alpha_{1} & A_{0} & 4 & 0 & 0 \\
	    		-4 \, \alpha_{1} & A_{1} & 5 & 0 & 0 \\
	    		-4 \, \alpha_{1} & A_{2} & 6 & 0 & 0 \\
	    		4 \, \alpha_{2} & A_{0} & 0 & 5 & 0 \\
	    		4 \, \alpha_{2} & A_{1} & 1 & 5 & 0 
	    		\end{dmatrix}
	    		\begin{dmatrix}
	    		-6 \, \alpha_{3} & A_{0} & 5 & 0 & 0 \\
	    		-6 \, \alpha_{3} & A_{1} & 6 & 0 & 0 \\
	    		-4 \, \alpha_{4} & A_{0} & 4 & 1 & 0 \\
	    		-4 \, \alpha_{4} & A_{1} & 5 & 1 & 0 \\
	    		12 \, \alpha_{1} {\left| A_{1} \right|}^{2} - 6 \, \alpha_{5} & A_{0} & 5 & 1 & 0 \\
	    		-10 \, \alpha_{6} & A_{0} & 7 & -1 & 0 \\
	    		-12 \, \alpha_{7} & A_{0} & 8 & -2 & 0 \\
	    		-8 \, \zeta_{0} & A_{0} & 6 & 0 & 1 \\
	    		4 \, \alpha_{1}^{2} - \zeta_{0} & A_{0} & 6 & 0 & 0 \\
	    		-10 \, \overline{\alpha_{2}} & A_{0} & 7 & -2 & 0 \\
	    		-10 \, \overline{\alpha_{2}} & A_{1} & 8 & -2 & 0 \\
	    		-2 \, \overline{\alpha_{4}} & A_{0} & 3 & 2 & 0 \\
	    		-2 \, \overline{\alpha_{4}} & A_{1} & 4 & 2 & 0 \\
	    		4 \, {\left| A_{1} \right|}^{2} \overline{\alpha_{1}} - 2 \, \overline{\alpha_{5}} & A_{0} & 3 & 3 & 0 \\
	    		2 \, \overline{\alpha_{6}} & A_{0} & 1 & 5 & 0 \\
	    		4 \, \overline{\alpha_{7}} & A_{0} & 0 & 6 & 0 \\
	    		-\overline{\zeta_{0}} & A_{0} & 2 & 4 & 0
	    		\end{dmatrix}
	    		\end{align}
	    		\normalsize
	    		Let $Q:\mathcal{QD}(D^2,\mathbb{C}^n)\rightarrow\mathbb{C}$ be the function defined on the $\mathbb{C}$-vector space of $\mathbb{C}^n$-valued quadratic differential $\mathcal{QD}(D^2,\mathbb{C}^n)$, such that for all $\vec{\alpha}\in \mathcal{QD}(D^2,\mathbb{C}^n)$, we have
	    		\begin{align*}
	    		Q(\alpha)=\partial\bar{\partial}\vec{\alpha}\totimes\vec{\alpha}-\partial\vec{\alpha}\totimes\bar{\partial}\vec{\alpha}.
	    		\end{align*}
	    		Then we have
	    		\begin{align}
	    		\mathscr{Q}_{\phi}=g^{-1}\otimes Q(\h_0)+\left(\frac{1}{4}|\vec{H}|^2+|\h_0|^2_{WP}\right)\h_0\totimes\h_0+\s{\H}{\h_0}^2
	    		\end{align}
	    		As for all $\theta_0\geq 4$, we have
	    		\begin{align*}
	    		\H=O(|z|^{2-\theta_0}),\quad \h_0=O(|z|^{\theta_0-1})
	    		\end{align*}
	    		we have
	    		\begin{align*}
	    		\frac{1}{4}|\H|^2\h_0\totimes \h_0+\s{\H}{\h_0}^2=O(|z|^2)
	    		\end{align*}
	    		while
	    		\begin{align*}
	    		|\h_0|^2_{WP}\h_0\totimes\h_0=\left(|\H|^2-K_g\right)\h_0\totimes\h_0=-K_g\h_0\totimes\h_0+O(|z|^2)
	    		\end{align*}
	    		Then we compute
	    		\small
	    		\begin{align}
	    		&g^{-1}\otimes Q(\h_0)=\nonumber\\
	    		&\begin{dmatrix}
	    		 -16 \, A_{0} C_{1} \alpha_{1} + 16 \, A_{2} C_{1} + 2 \, A_{1} C_{2} & 0 & 0 & 0 \\
	    		-30 \, A_{0} C_{1} \alpha_{3} + 8 \, A_{0} A_{1} \overline{\zeta_{0}} + 30 \, A_{3} C_{1} + 6 \, A_{2} C_{2} - 6 \, {\left(4 \, A_{1} C_{1} + A_{0} C_{2}\right)} \alpha_{1} - A_{1} \gamma_{1} & 1 & 0 & 0 \\
	    		20 \, A_{1} E_{1} & -3 & 4 & 0 \\
	    		\omega_1 & 0 & 1 & 0 \\
	    		\lambda_2 & 3 & -2 & 0 \\
	    		192 \, A_{0}^{2} \alpha_{2} {\left| A_{1} \right|}^{2} - \frac{48}{5} \, A_{0} B_{1} {\left| A_{1} \right|}^{2} - 3 \, A_{0} C_{1} \overline{\alpha_{4}} - 144 \, A_{0} A_{1} \overline{\alpha_{7}} + 6 \, A_{1} B_{2} + \frac{3}{8} \, C_{1} \overline{B_{1}} & -1 & 2 & 0 \\
	    		-24 \, A_{0} A_{1} \zeta_{0} + 12 \, A_{1} \overline{E_{1}} & 5 & -4 & 0 \\
	    		-16 \, A_{0}^{2} \alpha_{1} {\left| A_{1} \right|}^{2} - 8 \, A_{1}^{2} {\left| A_{1} \right|}^{2} + 16 \, A_{0} A_{2} {\left| A_{1} \right|}^{2} - 8 \, A_{0} A_{1} \alpha_{4} & 3 & -3 & 0 \\
	    		16 \, A_{0} A_{1} \alpha_{1} {\left| A_{1} \right|}^{2} - 48 \, A_{0}^{2} \alpha_{3} {\left| A_{1} \right|}^{2} - 16 \, A_{1} A_{2} {\left| A_{1} \right|}^{2} + 48 \, A_{0} A_{3} {\left| A_{1} \right|}^{2} - 16 \, A_{1}^{2} \alpha_{4} - 24 \, A_{0} A_{1} \alpha_{5} & 4 & -3 & 0 \\
	    		-9 \, A_{0} C_{1} {\left| A_{1} \right|}^{2} - 120 \, A_{0} A_{1} \alpha_{2} + 6 \, A_{1} B_{1} & -1 & 1 & 0 \\
	    		-480 \, A_{0}^{2} {\left| A_{1} \right|}^{2} \overline{\alpha_{2}} + 80 \, A_{0} A_{1} \alpha_{6} & 6 & -5 & 0 \\
	    		240 \, A_{0} A_{1} \alpha_{7} - 40 \, {\left(6 \, A_{0}^{2} \alpha_{1} - 5 \, A_{1}^{2} - 6 \, A_{0} A_{2}\right)} \overline{\alpha_{2}} & 7 & -6 & 0 \\
	    		160 \, A_{0} A_{1} \overline{\alpha_{2}} & 6 & -6 & 0 \\
	    		6 \, A_{1} C_{1} & -1 & 0 & 0
	    		\end{dmatrix}
	    		\end{align}
	    		\normalsize
	    		where
	    		\begin{align}
	    		&\lambda_1=-15 \, A_{1} C_{1} {\left| A_{1} \right|}^{2} - 3 \, A_{0} C_{2} {\left| A_{1} \right|}^{2} - 16 \, A_{0} B_{1} \alpha_{1} - 12 \, A_{0} C_{1} \alpha_{4} - 40 \, A_{0} A_{1} \overline{\alpha_{6}} + 16 \, A_{2} B_{1} + 2 \, A_{1} B_{3}\nonumber\\
	    		& + 80 \, {\left(4 \, A_{0}^{2} \alpha_{1} - A_{1}^{2} - 4 \, A_{0} A_{2}\right)} \alpha_{2}\nonumber\\
	    		&\lambda_2=32 \, A_{0} A_{1} {\left| A_{1} \right|}^{4} + 16 \, A_{0} A_{1} \alpha_{1} \overline{\alpha_{1}} - 16 \, A_{0} A_{1} \beta + 2 \, A_{0} \alpha_{1} \overline{B_{1}} - 105 \, A_{0} C_{1} \overline{\alpha_{2}} - 2 \, A_{2} \overline{B_{1}} + 2 \, A_{1} \overline{B_{2}}\nonumber\\
	    		& - 8 \, {\left(2 \, A_{0}^{2} \alpha_{1} + A_{1}^{2} - 2 \, A_{0} A_{2}\right)} \overline{\alpha_{4}}
	    		\end{align}
	    		Then we have
	    		\small
	    		\begin{align}
	    		-K_g\h_0\totimes\h_0=\begin{dmatrix}
	    		 16 \, A_{1}^{2} {\left| A_{1} \right|}^{2} & 3 & -3 & 0 \\
	    		-64 \, A_{0} A_{1} \alpha_{1} {\left| A_{1} \right|}^{2} + 64 \, A_{1} A_{2} {\left| A_{1} \right|}^{2} + 16 \, A_{1}^{2} \alpha_{4} & 4 & -3 & 0 \\
	    		-4 \, A_{1} C_{1} {\left| A_{1} \right|}^{2} - 80 \, A_{1}^{2} \alpha_{2} & 0 & 1 & 0 \\
	    		-64 \, A_{0} A_{1} {\left| A_{1} \right|}^{4} + 16 \, A_{1}^{2} \overline{\alpha_{4}} & 3 & -2 & 0 \\
	    		-80 \, A_{1}^{2} \overline{\alpha_{2}} & 7 & -6 & 0
	    		\end{dmatrix}
	    		\end{align}
	    		\normalsize
	    		and we finally obtain as $|\H|^2\h_0\totimes\h_0=O(|z|^2)$ and $\s{\H}{\h_0}^2=O(|z|^2)$
	    		\small
	    		\begin{align}
	    		&\mathscr{Q}_{\phi}=\nonumber\\
	    		&\begin{dmatrix}
	    		 -16 \, A_{0} C_{1} \alpha_{1} + 16 \, A_{2} C_{1} + 2 \, A_{1} C_{2} & 0 & 0 & 0 \\
	    		-30 \, A_{0} C_{1} \alpha_{3} + 8 \, A_{0} A_{1} \overline{\zeta_{0}} + 30 \, A_{3} C_{1} + 6 \, A_{2} C_{2} - 6 \, {\left(4 \, A_{1} C_{1} + A_{0} C_{2}\right)} \alpha_{1} - A_{1} \gamma_{1} & 1 & 0 & 0 \\
	    		20 \, A_{1} E_{1} & -3 & 4 & 0 \\
	    		\mu_1 & 0 & 1 & 0 \\
	    		\mu_2 & 3 & -2 & 0 \\
	    		192 \, A_{0}^{2} \alpha_{2} {\left| A_{1} \right|}^{2} - \frac{48}{5} \, A_{0} B_{1} {\left| A_{1} \right|}^{2} - 3 \, A_{0} C_{1} \overline{\alpha_{4}} - 144 \, A_{0} A_{1} \overline{\alpha_{7}} + 6 \, A_{1} B_{2} + \frac{3}{8} \, C_{1} \overline{B_{1}} & -1 & 2 & 0 \\
	    		-24 \, A_{0} A_{1} \zeta_{0} + 12 \, A_{1} \overline{E_{1}} & 5 & -4 & 0 \\
	    		-16 \, A_{0}^{2} \alpha_{1} {\left| A_{1} \right|}^{2} + 8 \, A_{1}^{2} {\left| A_{1} \right|}^{2} + 16 \, A_{0} A_{2} {\left| A_{1} \right|}^{2} - 8 \, A_{0} A_{1} \alpha_{4} & 3 & -3 & 0 \\
	    		-48 \, A_{0} A_{1} \alpha_{1} {\left| A_{1} \right|}^{2} - 48 \, A_{0}^{2} \alpha_{3} {\left| A_{1} \right|}^{2} + 48 \, A_{1} A_{2} {\left| A_{1} \right|}^{2} + 48 \, A_{0} A_{3} {\left| A_{1} \right|}^{2} - 24 \, A_{0} A_{1} \alpha_{5} & 4 & -3 & 0 \\
	    		-9 \, A_{0} C_{1} {\left| A_{1} \right|}^{2} - 120 \, A_{0} A_{1} \alpha_{2} + 6 \, A_{1} B_{1} & -1 & 1 & 0 \\
	    		-480 \, A_{0}^{2} {\left| A_{1} \right|}^{2} \overline{\alpha_{2}} + 80 \, A_{0} A_{1} \alpha_{6} & 6 & -5 & 0 \\
	    		240 \, A_{0} A_{1} \alpha_{7} - 80 \, A_{1}^{2} \overline{\alpha_{2}} - 40 \, {\left(6 \, A_{0}^{2} \alpha_{1} - 5 \, A_{1}^{2} - 6 \, A_{0} A_{2}\right)} \overline{\alpha_{2}} & 7 & -6 & 0 \\
	    		160 \, A_{0} A_{1} \overline{\alpha_{2}} & 6 & -6 & 0 \\
	    		6 \, A_{1} C_{1} & -1 & 0 & 0
	    		\end{dmatrix}
	    		\end{align}
	    		\normalsize
	    		where
	    		\begin{align*}
	    		\mu_1&=-19 \, A_{1} C_{1} {\left| A_{1} \right|}^{2} - 3 \, A_{0} C_{2} {\left| A_{1} \right|}^{2} - 16 \, A_{0} B_{1} \alpha_{1} - 80 \, A_{1}^{2} \alpha_{2} - 12 \, A_{0} C_{1} \alpha_{4} - 40 \, A_{0} A_{1} \overline{\alpha_{6}} + 16 \, A_{2} B_{1} + 2 \, A_{1} B_{3}\\
	    		& + 80 \, {\left(4 \, A_{0}^{2} \alpha_{1} - A_{1}^{2} - 4 \, A_{0} A_{2}\right)} \alpha_{2}\\
	    		\mu_2&=-32 \, A_{0} A_{1} {\left| A_{1} \right|}^{4} + 16 \, A_{0} A_{1} \alpha_{1} \overline{\alpha_{1}} - 16 \, A_{0} A_{1} \beta + 2 \, A_{0} \alpha_{1} \overline{B_{1}} - 105 \, A_{0} C_{1} \overline{\alpha_{2}} + 16 \, A_{1}^{2} \overline{\alpha_{4}} - 2 \, A_{2} \overline{B_{1}} + 2 \, A_{1} \overline{B_{2}}\\
	    		& - 8 \, {\left(2 \, A_{0}^{2} \alpha_{1} + A_{1}^{2} - 2 \, A_{0} A_{2}\right)} \overline{\alpha_{4}}
	    		\end{align*}
	    		It is easy so see that all relations coming from the meromorphy of $\mathscr{Q}_{\phi}$ are trivial, besides
	    		\begin{align*}
	    		\mu_1=0.
	    		\end{align*}
	    		As 
	    		\begin{align*}
	    		\s{\vec{A}_0}{\vec{A}_1}=\s{\vec{A}_0}{\vec{C}_1}=\s{\vec{A}_1}{\vec{A}_1}+2\s{\vec{A}_0}{\vec{A}_2}=0,
	    		\end{align*}
	    		we have
	    		\begin{align*}
	    		\mu_1&=-19 \, A_{1} C_{1} {\left| A_{1} \right|}^{2} - 3 \, A_{0} C_{2} {\left| A_{1} \right|}^{2} - \ccancel{16 \, A_{0} B_{1} \alpha_{1}} - 80 \, A_{1}^{2} \alpha_{2} - \ccancel{12 \, A_{0} C_{1} \alpha_{4}} - \ccancel{40 \, A_{0} A_{1} \overline{\alpha_{6}}} + 16 \, A_{2} B_{1} + 2 \, A_{1} B_{3}\\
	    		& + 80 \, {\left(\ccancel{4 \, A_{0}^{2} \alpha_{1}} - A_{1}^{2} - 4 \, A_{0} A_{2}\right)} \alpha_{2}\\
	    		&=-19|\vec{A}_1|^2\s{\vec{A}_1}{\vec{C}_1}-3|\vec{A}_1|^2\s{\vec{A}_0}{\vec{C}_2}+16\s{\vec{A}_2}{\vec{B}_1}+2\s{\vec{A}_1}{\vec{B}_3}-15|\vec{A}_1|^2\s{\vec{A}_1}{\vec{C}_1}-160\left(\colorcancel{[]A_1^2}{blue}+\colorcancel{2\,A_0A_2}{blue}\right)\\
	    		&=-19|\vec{A}_1|^2\s{\vec{A}_1}{\vec{C}_1}-3|\vec{A}_1|^2\s{\vec{A}_0}{\vec{C}_2}+16\s{\vec{A}_2}{\vec{B}_1}+2\s{\vec{A}_1}{\vec{B}_3}-15|\vec{A}_1|^2\s{\vec{A}_1}{\vec{C}_1}\\
	    		&=-16|\vec{A}_1|^2\s{\vec{A}_1}{\vec{C}_1}+16\s{\vec{A}_2}{\vec{B}_1}+2\s{\vec{A}_1}{\vec{B}_3}
	    		\end{align*}
	    		as
	    		\begin{align*}
	    		\s{\vec{A}_1}{\vec{C}_1}+\s{\vec{A}_0}{\vec{C}_2}=0.
	    		\end{align*}
	    		Then, recall that
	    		\begin{align}\label{4consts1}
	    		\left\{\begin{alignedat}{1}
	    		\vec{C}_2&=\vec{D}_2+2\s{\vec{A}_1}{\vec{C}_1}\bar{\vec{A}_0}\\
	    		\vec{C}_3&=\vec{D}_3\\
	    		\vec{B}_1&=-2\s{\bar{\vec{A}_1}}{\vec{C}_1}\vec{A}_0\\
	    		\vec{B}_2&=-\frac{3}{8}|\vec{C}_1|^2\bar{\vec{A}_0}-2\s{\bar{\vec{A}_2}}{\vec{C}_1}\vec{A}_0\\
	    		\vec{B}_3&=-2\s{\bar{\vec{A}_1}}{\vec{C}_1}\vec{A}_1+2\s{\vec{A}_1}{\vec{C}_1}\bar{\vec{A}_1}-2\s{\bar{\vec{A}_1}}{\vec{C}_2}\vec{A}_0\\
	    		\vec{E}_1&=-\frac{1}{8}\s{\vec{C}_1}{\vec{C}_1}\bar{\vec{A}_0}\\
	    		\vec{\gamma}_1&=-\vec{\gamma}_0-\Re\left\{\left(\frac{1}{2}\s{\vec{C}_1}{\vec{C}_1}+8\bar{\s{\vec{A}_2}{\vec{C}_1}}+4\bar{\s{\vec{A}_1}{\vec{C}_2}}\right)\vec{A}_0\right\}\in\R^n
	    		\end{alignedat}\right.
	    		\end{align}
	    		In particular, we have
	    		\begin{align}\label{a12b13}
	    		\left\{\begin{alignedat}{1}
	    		\s{\vec{A}_2}{\vec{B}_1}&=-2\s{\bar{\vec{A}_1}}{\vec{C}_1}\s{\vec{A}_0}{\vec{A}_2}=\s{\bar{\vec{A}_1}}{\vec{C}_1}\s{\vec{A}_1}{\vec{A}_1}\\
	    		\s{\vec{A}_1}{\vec{B}_3}&=-2\s{\bar{\vec{A}_1}}{\vec{C}_1}\s{\vec{A}_1}{\vec{A}_1}+2|\vec{A}_1|^2\s{\vec{A}_1}{\vec{C}_1}
	    		\end{alignedat}\right.
	    		\end{align}
	    		and
	    		\begin{align*}
	    		\mu_1&=-16|\vec{A}_1|^2\s{\vec{A}_1}{\vec{C}_1}+16\s{\bar{\vec{A}_1}}{\vec{C}_1}\s{\vec{A}_1}{\vec{A}_1}+2\left(-2\s{\bar{\vec{A}_1}}{\vec{C}_1}\s{\vec{A}_1}{\vec{A}_1}+2|\vec{A}_1|^2\s{\vec{A}_1}{\vec{C}_1}\right)\\
	    		&=-12|\vec{A}_1|^2\s{\vec{A}_1}{\vec{C}_1}+12\s{\bar{\vec{A}_1}}{\vec{C}_1}\s{\vec{A}_1}{\vec{A}_1}=0
	    		\end{align*}
	    		so we recover as expected
	    		\begin{align}\label{4strange1}
	    		|\vec{A}_1|^2\s{\vec{A}_1}{\vec{C}_1}=\s{\bar{\vec{A}_1}}{\vec{C}_1}\s{\vec{A}_1}{\vec{A}_1}
	    		\end{align}
	    		In particular, looking at \eqref{a12b13}, we see that
	    		\begin{align}\label{4usefulb1b3}
	    		\s{\vec{A}_2}{\vec{B}_1}&=|\vec{A}_1|^2\s{\vec{A}_1}{\vec{C}_1},\quad 
	    		\s{\vec{A}_1}{\vec{B}_3}=0
	    		\end{align}

	    		Then we get for some $\vec{D}_4\in \C^n$
	    		\small
	    		\begin{align*}
	    			&\H+\vec{\gamma}_0\log|z|=\Re\left(\frac{\vec{C}_1}{z^2}+\frac{\vec{D}_2}{z}+\vec{D}_3+\vec{D}_4z\right)+\\
	    			&\begin{dmatrix}
	    			-\frac{1}{16} \, C_{1}^{2} & \overline{A_{0}} & -4 & 4 & 0 \\
	    			-\frac{9}{160} \, C_{1}^{2} & \overline{A_{1}} & -4 & 5 & 0 \\
	    			\frac{1}{2} \, A_{0} C_{1} \alpha_{2} - \frac{17}{160} \, B_{1} C_{1} & \overline{A_{0}} & -4 & 5 & 0 \\
	    			-\frac{1}{80} \, C_{1} \overline{A_{1}} & C_{1} & -4 & 5 & 0 \\
	    			C_{1} \alpha_{2} \overline{A_{0}} - \frac{1}{5} \, E_{1} \overline{A_{1}} & A_{0} & -4 & 5 & 0 \\
	    			-\frac{1}{8} \, C_{1} C_{2} + \frac{1}{3} \, A_{1} E_{1} & \overline{A_{0}} & -3 & 4 & 0 \\
	    			-C_{1} \overline{A_{1}} & A_{0} & -2 & 1 & 0 \\
	    			-\frac{3}{16} \, C_{1} \overline{C_{1}} & \overline{A_{0}} & -2 & 2 & 0 \\
	    			C_{1} \overline{A_{0}} \overline{\alpha_{1}} - \frac{1}{2} \, B_{1} \overline{A_{1}} - C_{1} \overline{A_{2}} & A_{0} & -2 & 2 & 0 \\
	    			-\frac{7}{48} \, C_{1} \overline{C_{1}} & \overline{A_{1}} & -2 & 3 & 0 \\
	    			A_{0} \alpha_{2} \overline{C_{1}} - \frac{2}{15} \, B_{1} \overline{C_{1}} - \frac{7}{48} \, C_{1} \overline{C_{2}} & \overline{A_{0}} & -2 & 3 & 0 \\
	    			-\frac{1}{48} \, \overline{A_{1}} \overline{C_{1}} & C_{1} & -2 & 3 & 0 \\
	    			\frac{5}{3} \, \alpha_{2} \overline{A_{0}} \overline{C_{1}} + C_{1} \overline{A_{0}} \overline{\alpha_{3}} - \frac{1}{3} \, B_{2} \overline{A_{1}} - \frac{2}{3} \, B_{1} \overline{A_{2}} - C_{1} \overline{A_{3}} + \frac{1}{3} \, {\left(2 \, B_{1} \overline{A_{0}} + 3 \, C_{1} \overline{A_{1}}\right)} \overline{\alpha_{1}} & A_{0} & -2 & 3 & 0 \\
	    			-\frac{1}{24} \, C_{1} \overline{A_{1}} & \overline{C_{1}} & -2 & 3 & 0 \\
	    			A_{1} C_{1} & \overline{A_{0}} & -1 & 0 & 0 \\
	    			A_{1} C_{1} & \overline{A_{1}} & -1 & 1 & 0 \\
	    			-2 \, A_{0} C_{1} {\left| A_{1} \right|}^{2} + A_{1} B_{1} & \overline{A_{0}} & -1 & 1 & 0 \\
	    			-C_{1} \overline{A_{1}} & A_{1} & -1 & 1 & 0 \\
	    			2 \, C_{1} {\left| A_{1} \right|}^{2} \overline{A_{0}} - C_{2} \overline{A_{1}} & A_{0} & -1 & 1 & 0 \\
	    			-2 \, A_{0} B_{1} {\left| A_{1} \right|}^{2} - A_{1} C_{1} \overline{\alpha_{1}} - A_{0} C_{1} \overline{\alpha_{4}} + A_{1} B_{2} - \frac{1}{8} \, C_{1} \overline{B_{1}} - \frac{3}{16} \, C_{2} \overline{C_{1}} & \overline{A_{0}} & -1 & 2 & 0 \\
	    			A_{1} C_{1} & \overline{A_{2}} & -1 & 2 & 0 \\
	    			-2 \, A_{0} C_{1} {\left| A_{1} \right|}^{2} + A_{1} B_{1} & \overline{A_{1}} & -1 & 2 & 0 \\
	    			B_{1} {\left| A_{1} \right|}^{2} \overline{A_{0}} + 2 \, C_{1} {\left| A_{1} \right|}^{2} \overline{A_{1}} + C_{2} \overline{A_{0}} \overline{\alpha_{1}} + C_{1} \overline{A_{0}} \overline{\alpha_{4}} - \frac{1}{2} \, B_{3} \overline{A_{1}} - C_{2} \overline{A_{2}} & A_{0} & -1 & 2 & 0 \\
	    			C_{1} \overline{A_{0}} \overline{\alpha_{1}} - \frac{1}{2} \, B_{1} \overline{A_{1}} - C_{1} \overline{A_{2}} & A_{1} & -1 & 2 & 0 \\
	    			\overline{A_{1}} \overline{C_{1}} & A_{0} & 0 & -1 & 0 \\
	    			4 \, \overline{A_{0}} \overline{C_{1}} \overline{\alpha_{1}} - \frac{1}{4} \, C_{1}^{2} - 4 \, \overline{A_{2}} \overline{C_{1}} - 2 \, \overline{A_{1}} \overline{C_{2}} & A_{0} & 0 & 0 & 1 \\
	    			4 \, A_{0} C_{1} \alpha_{1} - 4 \, A_{2} C_{1} - 2 \, A_{1} C_{2} - \frac{1}{4} \, \overline{C_{1}}^{2} & \overline{A_{0}} & 0 & 0 & 1 \\
	    			-\frac{1}{8} \, \overline{C_{1}}^{2} & \overline{A_{1}} & 0 & 1 & 0 
	    			\end{dmatrix}
	    			\begin{dmatrix}
	    			\textbf{(1)}\\
	    			\textbf{(2)}\\
	    			\textbf{(3)}\\
	    			\textbf{(4)}\\
	    			\textbf{(5)}\\
	    			\textbf{(6)}\\
	    			\textbf{(7)}\\
	    			\textbf{(8)}\\
	    			\textbf{(9)}\\
	    			\textbf{(10)}\\
	    			\textbf{(11)}\\
	    			\textbf{(12)}\\
	    			\textbf{(13)}\\
	    			\textbf{(14)}\\
	    			\textbf{(15)}\\
	    			\textbf{(16)}\\
	    			\textbf{(17)}\\
	    			\textbf{(18)}\\
	    			\textbf{(19)}\\
	    			\textbf{(20)}\\
	    			\textbf{(21)}\\
	    			\textbf{(22)}\\
	    			\textbf{(23)}\\
	    			\textbf{(24)}\\
	    			\textbf{(25)}\\
	    			\textbf{(26)}\\
	    			\textbf{(27)}\\
	    			\textbf{(28)}
	    			\end{dmatrix}
	    			\end{align*}
	    			\begin{align*}
	    			\begin{dmatrix}
	    			-\frac{1}{4} \, \overline{C_{1}} \overline{C_{2}} & \overline{A_{0}} & 0 & 1 & 0 \\
	    			-\frac{1}{8} \, \overline{A_{1}} \overline{C_{1}} & \overline{C_{1}} & 0 & 1 & 0 \\
	    			\lambda_1 & A_{0} & 0 & 1 & 0 \\
	    			-C_{1} \overline{A_{1}} & A_{2} & 0 & 1 & 0 \\
	    			2 \, C_{1} {\left| A_{1} \right|}^{2} \overline{A_{0}} - C_{2} \overline{A_{1}} & A_{1} & 0 & 1 & 0 \\
	    			-\frac{1}{2} \, B_{1} C_{1} - 2 \, \gamma_{1} \overline{A_{1}} & A_{0} & 0 & 1 & 1 \\
	    			8 \, A_{1} C_{1} {\left| A_{1} \right|}^{2} + 4 \, A_{0} C_{2} {\left| A_{1} \right|}^{2} + 4 \, A_{0} B_{1} \alpha_{1} + 4 \, A_{0} C_{1} \alpha_{4} - 4 \, A_{2} B_{1} - 2 \, A_{1} B_{3} & \overline{A_{0}} & 0 & 1 & 1 \\
	    			4 \, A_{0} C_{1} \alpha_{1} - 4 \, A_{2} C_{1} - 2 \, A_{1} C_{2} & \overline{A_{1}} & 0 & 1 & 1 \\
	    			-A_{1} \overline{C_{1}} & \overline{A_{0}} & 1 & -2 & 0 \\
	    			\overline{A_{1}} \overline{C_{1}} & A_{1} & 1 & -1 & 0 \\
	    			-2 \, {\left| A_{1} \right|}^{2} \overline{A_{0}} \overline{C_{1}} + \overline{A_{1}} \overline{B_{1}} & A_{0} & 1 & -1 & 0 \\
	    			-A_{1} \overline{C_{1}} & \overline{A_{1}} & 1 & -1 & 0 \\
	    			2 \, A_{0} {\left| A_{1} \right|}^{2} \overline{C_{1}} - A_{1} \overline{C_{2}} & \overline{A_{0}} & 1 & -1 & 0 \\
	    			-\frac{1}{8} \, C_{1}^{2} & A_{1} & 1 & 0 & 0 \\
	    			-\frac{1}{4} \, C_{1} C_{2} & A_{0} & 1 & 0 & 0 \\
	    			-\frac{1}{8} \, A_{1} C_{1} & C_{1} & 1 & 0 & 0 \\
	    			\lambda_2 & \overline{A_{0}} & 1 & 0 & 0 \\
	    			-A_{1} \overline{C_{1}} & \overline{A_{2}} & 1 & 0 & 0 \\
	    			2 \, A_{0} {\left| A_{1} \right|}^{2} \overline{C_{1}} - A_{1} \overline{C_{2}} & \overline{A_{1}} & 1 & 0 & 0 \\
	    			-2 \, A_{1} \gamma_{1} - \frac{1}{2} \, \overline{B_{1}} \overline{C_{1}} & \overline{A_{0}} & 1 & 0 & 1 \\
	    			8 \, {\left| A_{1} \right|}^{2} \overline{A_{1}} \overline{C_{1}} + 4 \, {\left| A_{1} \right|}^{2} \overline{A_{0}} \overline{C_{2}} + 4 \, \overline{A_{0}} \overline{B_{1}} \overline{\alpha_{1}} + 4 \, \overline{A_{0}} \overline{C_{1}} \overline{\alpha_{4}} - 4 \, \overline{A_{2}} \overline{B_{1}} - 2 \, \overline{A_{1}} \overline{B_{3}} & A_{0} & 1 & 0 & 1 \\
	    			4 \, \overline{A_{0}} \overline{C_{1}} \overline{\alpha_{1}} - 4 \, \overline{A_{2}} \overline{C_{1}} - 2 \, \overline{A_{1}} \overline{C_{2}} & A_{1} & 1 & 0 & 1 \\
	    			-\frac{3}{16} \, C_{1} \overline{C_{1}} & A_{0} & 2 & -2 & 0 \\
	    			A_{0} \alpha_{1} \overline{C_{1}} - \frac{1}{2} \, A_{1} \overline{B_{1}} - A_{2} \overline{C_{1}} & \overline{A_{0}} & 2 & -2 & 0 \\
	    			-2 \, {\left| A_{1} \right|}^{2} \overline{A_{0}} \overline{B_{1}} - \alpha_{4} \overline{A_{0}} \overline{C_{1}} - \alpha_{1} \overline{A_{1}} \overline{C_{1}} + \overline{A_{1}} \overline{B_{2}} - \frac{1}{8} \, B_{1} \overline{C_{1}} - \frac{3}{16} \, C_{1} \overline{C_{2}} & A_{0} & 2 & -1 & 0 \\
	    			\overline{A_{1}} \overline{C_{1}} & A_{2} & 2 & -1 & 0 \\
	    			-2 \, {\left| A_{1} \right|}^{2} \overline{A_{0}} \overline{C_{1}} + \overline{A_{1}} \overline{B_{1}} & A_{1} & 2 & -1 & 0 \\
	    			A_{0} {\left| A_{1} \right|}^{2} \overline{B_{1}} + 2 \, A_{1} {\left| A_{1} \right|}^{2} \overline{C_{1}} + A_{0} \alpha_{4} \overline{C_{1}} + A_{0} \alpha_{1} \overline{C_{2}} - \frac{1}{2} \, A_{1} \overline{B_{3}} - A_{2} \overline{C_{2}} & \overline{A_{0}} & 2 & -1 & 0 \\
	    			A_{0} \alpha_{1} \overline{C_{1}} - \frac{1}{2} \, A_{1} \overline{B_{1}} - A_{2} \overline{C_{1}} & \overline{A_{1}} & 2 & -1 & 0 
	    			\end{dmatrix}
	    			\begin{dmatrix}
	    			\textbf{(29)}\\
	    			\textbf{(30)}\\
	    			\textbf{(31)}\\
	    			\textbf{(32)}\\
	    			\textbf{(33)}\\
	    			\textbf{(34)}\\
	    			\textbf{(35)}\\
	    			\textbf{(36)}\\
	    			\textbf{(37)}\\
	    			\textbf{(38)}\\
	    			\textbf{(39)}\\
	    			\textbf{(40)}\\
	    			\textbf{(41)}\\
	    			\textbf{(42)}\\
	    			\textbf{(43)}\\
	    			\textbf{(44)}\\
	    			\textbf{(45)}\\
	    			\textbf{(46)}\\
	    			\textbf{(47)}\\
	    			\textbf{(48)}\\
	    			\textbf{(49)}\\
	    			\textbf{(50)}\\
	    			\textbf{(51)}\\
	    			\textbf{(52)}\\
	    			\textbf{(53)}\\
	    			\textbf{(54)}\\
	    			\textbf{(55)}\\
	    			\textbf{(56)}\\
	    			\textbf{(57)}
	    			\end{dmatrix}
	    			\end{align*}
	    			\begin{align*}
	    			\begin{dmatrix}
	    			-\frac{7}{48} \, C_{1} \overline{C_{1}} & A_{1} & 3 & -2 & 0 \\
	    			C_{1} \overline{A_{0}} \overline{\alpha_{2}} - \frac{2}{15} \, C_{1} \overline{B_{1}} - \frac{7}{48} \, C_{2} \overline{C_{1}} & A_{0} & 3 & -2 & 0 \\
	    			-\frac{1}{48} \, A_{1} C_{1} & \overline{C_{1}} & 3 & -2 & 0 \\
	    			A_{0} \alpha_{3} \overline{C_{1}} + \frac{5}{3} \, A_{0} C_{1} \overline{\alpha_{2}} + \frac{1}{3} \, {\left(2 \, A_{0} \overline{B_{1}} + 3 \, A_{1} \overline{C_{1}}\right)} \alpha_{1} - \frac{2}{3} \, A_{2} \overline{B_{1}} - \frac{1}{3} \, A_{1} \overline{B_{2}} - A_{3} \overline{C_{1}} & \overline{A_{0}} & 3 & -2 & 0 \\
	    			-\frac{1}{24} \, A_{1} \overline{C_{1}} & C_{1} & 3 & -2 & 0 \\
	    			-\frac{1}{16} \, \overline{C_{1}}^{2} & A_{0} & 4 & -4 & 0 \\
	    			-\frac{1}{8} \, \overline{C_{1}} \overline{C_{2}} + \frac{1}{3} \, \overline{A_{1}} \overline{E_{1}} & A_{0} & 4 & -3 & 0 \\
	    			-\frac{9}{160} \, \overline{C_{1}}^{2} & A_{1} & 5 & -4 & 0 \\
	    			\frac{1}{2} \, \overline{A_{0}} \overline{C_{1}} \overline{\alpha_{2}} - \frac{17}{160} \, \overline{B_{1}} \overline{C_{1}} & A_{0} & 5 & -4 & 0 \\
	    			-\frac{1}{80} \, A_{1} \overline{C_{1}} & \overline{C_{1}} & 5 & -4 & 0 \\
	    			A_{0} \overline{C_{1}} \overline{\alpha_{2}} - \frac{1}{5} \, A_{1} \overline{E_{1}} & \overline{A_{0}} & 5 & -4 & 0
	    			\end{dmatrix}
	    			\begin{dmatrix}
	    			\textbf{(58)}\\
	    			\textbf{(59)}\\
	    			\textbf{(60)}\\
	    			\textbf{(61)}\\
	    			\textbf{(62)}\\
	    			\textbf{(63)}\\
	    			\textbf{(64)}\\
	    			\textbf{(65)}\\
	    			\textbf{(66)}\\
	    			\textbf{(67)}\\
	    			\textbf{(68)}
	    			\end{dmatrix}
	    		\end{align*}
	    		where
	    		\begin{align*}
	    			\lambda_1&=2 \, C_{2} {\left| A_{1} \right|}^{2} \overline{A_{0}} + C_{1} \alpha_{4} \overline{A_{0}} + C_{1} \alpha_{1} \overline{A_{1}} + 3 \, \overline{A_{0}} \overline{C_{1}} \overline{\alpha_{3}} - \frac{1}{8} \, B_{1} C_{1} - 2 \, C_{3} \overline{A_{1}} + \gamma_{1} \overline{A_{1}} - 3 \, \overline{A_{3}} \overline{C_{1}} - 2 \, \overline{A_{2}} \overline{C_{2}}\\
	    			& + {\left(3 \, \overline{A_{1}} \overline{C_{1}} + 2 \, \overline{A_{0}} \overline{C_{2}}\right)} \overline{\alpha_{1}}\\
	    			\lambda_2&=2 \, A_{0} {\left| A_{1} \right|}^{2} \overline{C_{2}} + 3 \, A_{0} C_{1} \alpha_{3} + A_{1} \overline{C_{1}} \overline{\alpha_{1}} + A_{0} \overline{C_{1}} \overline{\alpha_{4}} - 3 \, A_{3} C_{1} - 2 \, A_{2} C_{2} - 2 \, A_{1} C_{3}\\
	    			& + {\left(3 \, A_{1} C_{1} + 2 \, A_{0} C_{2}\right)} \alpha_{1} + A_{1} \gamma_{1} - \frac{1}{8} \, \overline{B_{1}} \overline{C_{1}}
	    		\end{align*}
	    		\normalsize
	    		The new powers arising are
	    		\begin{align*}
	    		\begin{dmatrix}
	    		-4 & 5 & 0\\
	    		-3 & 4 & 0\\
	    		-2 & 3 & 0\
	    		\end{dmatrix}
	    		\begin{dmatrix}
	    		-1 & 2 & 0\\
	    		1 & 0 & 0\\
	    		1 & 0 & 1
	    		\end{dmatrix}
	    		\end{align*}
	    		so there exists $\vec{C}_4,\vec{B}_4,\vec{B}_5,\vec{E}_2,\vec{E}_3,\vec{\gamma}_2\in\mathbb{C}^n$ such that
	    		\begin{align}\label{4H2}
	    		\H&=\Re\left(\frac{\vec{C}_1}{z^2}+\frac{\vec{C}_2}{z}+\vec{
	    			C}_4z+\left(\vec{B}_1\z+\vec{B}_2\z^2+\vec{B}_4\z^3\right)\frac{1}{z^2}+\left(\vec{B}_3\z+\vec{B}_5\z^2\right)\frac{1}{z}+\left(\vec{E}_1\z^4+\vec{E}_2\z^{5}\right)\frac{1}{z^4}+\vec{E}_3\frac{\z^4}{z^3}\right)\nonumber\\   &+\vec{C}_3+\vec{\gamma}_1\log|z|+\Re\Big(\vec{\gamma}_2\, z\Big)\log|z|+O(|z|^{2-\epsilon}).
	    		\end{align}
	    		First, we check that the mean curvature given in the code coincides with \eqref{4H2} (recall that $\vec{C}_3\in \R^n$).
	    		\begin{align}\label{4devH2}
	    		\H=\begin{dmatrix}
	    		\frac{1}{2} & C_{1} & -2 & 0 & 0 \\
	    		\frac{1}{2} & C_{2} & -1 & 0 & 0 \\
	    		1 & C_{3} & 0 & 0 & 0 \\
	    		\frac{1}{2} & C_{4} & 1 & 0 & 0 \\
	    		\frac{1}{2} & B_{1} & -2 & 1 & 0 \\
	    		\frac{1}{2} & B_{2} & -2 & 2 & 0 \\
	    		\frac{1}{2} & B_{4} & -2 & 3 & 0 
	    		\end{dmatrix}
	    		\begin{dmatrix}
	    		\frac{1}{2} & B_{3} & -1 & 1 & 0 \\
	    		\frac{1}{2} & B_{5} & -1 & 2 & 0 \\
	    		\frac{1}{2} & E_{1} & -4 & 4 & 0 \\
	    		\frac{1}{2} & E_{2} & -4 & 5 & 0 \\
	    		\frac{1}{2} & E_{3} & -3 & 4 & 0 \\
	    		\frac{1}{2} & \gamma_{2} & 1 & 0 & 1 
	    		\end{dmatrix}
	    		\begin{dmatrix}
	    		\frac{1}{2} & \overline{C_{1}} & 0 & -2 & 0 \\
	    		\frac{1}{2} & \overline{C_{2}} & 0 & -1 & 0 \\
	    		\frac{1}{2} & \overline{C_{4}} & 0 & 1 & 0 \\
	    		\frac{1}{2} & \overline{B_{1}} & 1 & -2 & 0 \\
	    		\frac{1}{2} & \overline{B_{2}} & 2 & -2 & 0 \\
	    		\frac{1}{2} & \overline{B_{4}} & 3 & -2 & 0 
	    		\end{dmatrix}
	    		\begin{dmatrix}
	    		\frac{1}{2} & \overline{B_{3}} & 1 & -1 & 0 \\
	    		\frac{1}{2} & \overline{B_{5}} & 2 & -1 & 0 \\
	    		\frac{1}{2} & \overline{E_{1}} & 4 & -4 & 0 \\
	    		\frac{1}{2} & \overline{E_{2}} & 5 & -4 & 0 \\
	    		\frac{1}{2} & \overline{E_{3}} & 4 & -3 & 0 \\
	    		\frac{1}{2} & \overline{\gamma_{2}} & 0 & 1 & 1 \\
	    		\end{dmatrix}
	    		\end{align}
	    		
	    		We see that $\vec{\gamma}_2$ corresponds to \emph{twice} the lines
	    		\begin{align*}
	    		(48),\;\, (49)\;\, (50)
	    		\end{align*}
	    		so
	    		\begin{align*}
	    		\vec{\gamma}_2=-2\left(4\bar{\s{\vec{A}_2}{\vec{C}_1}}+2\bar{\s{\vec{A}_1}{\vec{C}_2}}\right)\vec{A}_1+4\s{\vec{A}_1}{\vec{\gamma}_0}\bar{\vec{A}_0}
	    		\end{align*}
	    		as by \eqref{a12b13} and \eqref{4strange1}, $\s{\vec{A}_1}{\vec{B}_3}=0$ and by \eqref{4cancel1}
	    		\begin{align*}
	    		&8 \, {\left| A_{1} \right|}^{2} \overline{A_{1}} \overline{C_{1}} + 4 \, {\left| A_{1} \right|}^{2} \overline{A_{0}} \overline{C_{2}} + 4 \, \overline{A_{0}} \overline{B_{1}} \overline{\alpha_{1}} + 4 \, \overline{A_{0}} \overline{C_{1}} \overline{\alpha_{4}} - 4 \, \overline{A_{2}} \overline{B_{1}} - 2 \, \overline{A_{1}} \overline{B_{3}}\\
	    		&=4|\vec{A}_1|^2\bar{\s{\vec{A}_1}{\vec{C}_1}}-4\s{\bar{\vec{A}_1}}{\vec{C}_1}\s{\vec{A}_1}{\vec{A}_1}=0.
	    		\end{align*}
	    		Now we have for some $\vec{D}_5\in\mathbb{C}^n$
	    		\footnotesize
	    		\begin{align*}
	    			&\H+\vec{\gamma}_0\log|z|=\Re\left(\frac{\vec{C}_1}{z^2}+\frac{\vec{D}_2}{z}+\vec{D}_3+\vec{D}_4z+\vec{D}_5z^2\right)+\\
	    			&\begin{dmatrix}
	    			-\frac{7}{96} \, C_{1} E_{1} & \overline{A_{0}} & -6 & 8 & 0 &\textbf{(1)}\\
	    			-\frac{1}{16} \, C_{1}^{2} & \overline{A_{0}} & -4 & 4 & 0 &\textbf{(2)}\\
	    			-\frac{9}{160} \, C_{1}^{2} & \overline{A_{1}} & -4 & 5 & 0 &\textbf{(3)}\\
	    			\frac{1}{2} \, A_{0} C_{1} \alpha_{2} - \frac{17}{160} \, B_{1} C_{1} & \overline{A_{0}} & -4 & 5 & 0 &\textbf{(4)}\\
	    			-\frac{1}{80} \, C_{1} \overline{A_{1}} & C_{1} & -4 & 5 & 0 &\textbf{(5)}\\
	    			C_{1} \alpha_{2} \overline{A_{0}} - \frac{1}{5} \, E_{1} \overline{A_{1}} & A_{0} & -4 & 5 & 0 &\textbf{(6)}\\
	    			\frac{1}{2} \, A_{0} B_{1} \alpha_{2} + \frac{1}{96} \, C_{1}^{2} \overline{\alpha_{1}} + \frac{1}{2} \, A_{0} C_{1} \overline{\alpha_{7}} - \frac{11}{240} \, B_{1}^{2} - \frac{3}{32} \, B_{2} C_{1} - \frac{7}{96} \, E_{1} \overline{C_{1}} & \overline{A_{0}} & -4 & 6 & 0 &\textbf{(7)}\\
	    			-\frac{5}{96} \, C_{1}^{2} & \overline{A_{2}} & -4 & 6 & 0 &\textbf{(8)}\\
	    			\frac{1}{2} \, A_{0} C_{1} \alpha_{2} - \frac{47}{480} \, B_{1} C_{1} & \overline{A_{1}} & -4 & 6 & 0 &\textbf{(9)}\\
	    			-\frac{1}{120} \, C_{1} \overline{A_{1}} & B_{1} & -4 & 6 & 0 &\textbf{(10)}\\
	    			\frac{1}{3} \, E_{1} \overline{A_{0}} \overline{\alpha_{1}} + C_{1} \overline{A_{0}} \overline{\alpha_{7}} + \frac{1}{6} \, {\left(5 \, B_{1} \overline{A_{0}} + 6 \, C_{1} \overline{A_{1}}\right)} \alpha_{2} - \frac{1}{6} \, E_{2} \overline{A_{1}} - \frac{1}{3} \, E_{1} \overline{A_{2}} & A_{0} & -4 & 6 & 0 &\textbf{(11)}\\
	    			\frac{1}{48} \, C_{1} \overline{A_{0}} \overline{\alpha_{1}} - \frac{1}{96} \, B_{1} \overline{A_{1}} - \frac{1}{48} \, C_{1} \overline{A_{2}} & C_{1} & -4 & 6 & 0 &\textbf{(12)}\\
	    			-\frac{1}{8} \, C_{1} C_{2} + \frac{1}{3} \, A_{1} E_{1} & \overline{A_{0}} & -3 & 4 & 0 &\textbf{(13)}\\
	    			-\frac{9}{80} \, C_{1} C_{2} + \frac{1}{3} \, A_{1} E_{1} & \overline{A_{1}} & -3 & 5 & 0 &\textbf{(14)}\\
	    			\frac{1}{120} \, C_{1}^{2} {\left| A_{1} \right|}^{2} - \frac{2}{3} \, A_{0} E_{1} {\left| A_{1} \right|}^{2} + \frac{1}{3} \, A_{0} C_{1} \overline{\alpha_{6}} - \frac{13}{120} \, B_{3} C_{1} - \frac{5}{48} \, B_{1} C_{2} + \frac{1}{3} \, A_{1} E_{2} + \frac{1}{3} \, {\left(A_{1} C_{1} + 2 \, A_{0} C_{2}\right)} \alpha_{2} & \overline{A_{0}} & -3 & 5 & 0 &\textbf{(15)}\\
	    			-\frac{1}{80} \, C_{1} \overline{A_{1}} & C_{2} & -3 & 5 & 0 &\textbf{(16)}\\
	    			\frac{1}{40} \, C_{1} {\left| A_{1} \right|}^{2} \overline{A_{0}} - \frac{1}{80} \, C_{2} \overline{A_{1}} & C_{1} & -3 & 5 & 0 &\textbf{(17)}\\
	    			\frac{2}{5} \, E_{1} {\left| A_{1} \right|}^{2} \overline{A_{0}} + C_{2} \alpha_{2} \overline{A_{0}} + C_{1} \overline{A_{0}} \overline{\alpha_{6}} - \frac{1}{5} \, E_{3} \overline{A_{1}} & A_{0} & -3 & 5 & 0 &\textbf{(18)}\\
	    			C_{1} \alpha_{2} \overline{A_{0}} - \frac{1}{5} \, E_{1} \overline{A_{1}} & A_{1} & -3 & 5 & 0 &\textbf{(19)}\\
	    			-C_{1} \overline{A_{1}} & A_{0} & -2 & 1 & 0 &\textbf{(20)}\\
	    			-\frac{3}{16} \, C_{1} \overline{C_{1}} & \overline{A_{0}} & -2 & 2 & 0 &\textbf{(21)}\\
	    			C_{1} \overline{A_{0}} \overline{\alpha_{1}} - \frac{1}{2} \, B_{1} \overline{A_{1}} - C_{1} \overline{A_{2}} & A_{0} & -2 & 2 & 0 &\textbf{(22)}\\
	    			-\frac{7}{48} \, C_{1} \overline{C_{1}} & \overline{A_{1}} & -2 & 3 & 0 &\textbf{(23)}\\
	    			A_{0} \alpha_{2} \overline{C_{1}} - \frac{2}{15} \, B_{1} \overline{C_{1}} - \frac{7}{48} \, C_{1} \overline{C_{2}} & \overline{A_{0}} & -2 & 3 & 0 &\textbf{(24)}\\
	    			-\frac{1}{48} \, \overline{A_{1}} \overline{C_{1}} & C_{1} & -2 & 3 & 0 &\textbf{(25)}\\
	    			\frac{5}{3} \, \alpha_{2} \overline{A_{0}} \overline{C_{1}} + C_{1} \overline{A_{0}} \overline{\alpha_{3}} - \frac{1}{3} \, B_{2} \overline{A_{1}} - \frac{2}{3} \, B_{1} \overline{A_{2}} - C_{1} \overline{A_{3}} + \frac{1}{3} \, {\left(2 \, B_{1} \overline{A_{0}} + 3 \, C_{1} \overline{A_{1}}\right)} \overline{\alpha_{1}} & A_{0} & -2 & 3 & 0 &\textbf{(26)}\\
	    			-\frac{1}{24} \, C_{1} \overline{A_{1}} & \overline{C_{1}} & -2 & 3 & 0 &\textbf{(27)}\\
	    			\frac{1}{32} \, \overline{A_{0}} \overline{C_{1}} \overline{\alpha_{1}} - \frac{1}{128} \, C_{1}^{2} - \frac{1}{32} \, \overline{A_{2}} \overline{C_{1}} - \frac{1}{64} \, \overline{A_{1}} \overline{C_{2}} & C_{1} & -2 & 4 & 0 &\textbf{(28)}\\
	    			\lambda_1 & A_{0} & -2 & 4 & 0 &\textbf{(29)}\\
	    			-\frac{1}{8} \, C_{1} \overline{C_{1}} & \overline{A_{2}} & -2 & 4 & 0 &\textbf{(30)}\\
	    			A_{0} \alpha_{2} \overline{C_{1}} - \frac{9}{80} \, B_{1} \overline{C_{1}} - \frac{1}{8} \, C_{1} \overline{C_{2}} & \overline{A_{1}} & -2 & 4 & 0 &\textbf{(31)}
	    			\end{dmatrix}
	    			\end{align*}
	    			\begin{align*}
	    			\begin{dmatrix}
	    			\lambda_2 & \overline{A_{0}} & -2 & 4 & 0 &\textbf{(32)}\\
	    			-\frac{1}{80} \, \overline{A_{1}} \overline{C_{1}} & B_{1} & -2 & 4 & 0 &\textbf{(33)}\\
	    			-\frac{1}{48} \, C_{1} \overline{A_{1}} & \overline{C_{2}} & -2 & 4 & 0 &\textbf{(34)}\\
	    			\frac{1}{16} \, C_{1} \overline{A_{0}} \overline{\alpha_{1}} - \frac{1}{32} \, B_{1} \overline{A_{1}} - \frac{1}{16} \, C_{1} \overline{A_{2}} & \overline{C_{1}} & -2 & 4 & 0 &\textbf{(35)}\\
	    			-\frac{1}{4} \, C_{1} \gamma_{1} & \overline{A_{0}} & -2 & 4 & 1 &\textbf{(36)}\\
	    			C_{1} \overline{A_{0}} \overline{\zeta_{0}} - \frac{1}{2} \, C_{1} E_{1} & A_{0} & -2 & 4 & 1 &\textbf{(37)}\\
	    			A_{1} C_{1} & \overline{A_{0}} & -1 & 0 & 0 &\textbf{(38)}\\
	    			A_{1} C_{1} & \overline{A_{1}} & -1 & 1 & 0 &\textbf{(39)}\\
	    			-2 \, A_{0} C_{1} {\left| A_{1} \right|}^{2} + A_{1} B_{1} & \overline{A_{0}} & -1 & 1 & 0 &\textbf{(40)}\\
	    			-C_{1} \overline{A_{1}} & A_{1} & -1 & 1 & 0 &\textbf{(41)}\\
	    			2 \, C_{1} {\left| A_{1} \right|}^{2} \overline{A_{0}} - C_{2} \overline{A_{1}} & A_{0} & -1 & 1 & 0 &\textbf{(42)}\\
	    			-2 \, A_{0} B_{1} {\left| A_{1} \right|}^{2} - A_{1} C_{1} \overline{\alpha_{1}} - A_{0} C_{1} \overline{\alpha_{4}} + A_{1} B_{2} - \frac{1}{8} \, C_{1} \overline{B_{1}} - \frac{3}{16} \, C_{2} \overline{C_{1}} & \overline{A_{0}} & -1 & 2 & 0 &\textbf{(43)}\\
	    			A_{1} C_{1} & \overline{A_{2}} & -1 & 2 & 0 &\textbf{(44)}\\
	    			-2 \, A_{0} C_{1} {\left| A_{1} \right|}^{2} + A_{1} B_{1} & \overline{A_{1}} & -1 & 2 & 0 &\textbf{(45)}\\
	    			B_{1} {\left| A_{1} \right|}^{2} \overline{A_{0}} + 2 \, C_{1} {\left| A_{1} \right|}^{2} \overline{A_{1}} + C_{2} \overline{A_{0}} \overline{\alpha_{1}} + C_{1} \overline{A_{0}} \overline{\alpha_{4}} - \frac{1}{2} \, B_{3} \overline{A_{1}} - C_{2} \overline{A_{2}} & A_{0} & -1 & 2 & 0 &\textbf{(46)}\\
	    			C_{1} \overline{A_{0}} \overline{\alpha_{1}} - \frac{1}{2} \, B_{1} \overline{A_{1}} - C_{1} \overline{A_{2}} & A_{1} & -1 & 2 & 0 &\textbf{(47)}\\
	    			-2 \, A_{0} B_{1} {\left| A_{1} \right|}^{2} - A_{1} C_{1} \overline{\alpha_{1}} - A_{0} C_{1} \overline{\alpha_{4}} + A_{1} B_{2} - \frac{1}{12} \, C_{1} \overline{B_{1}} - \frac{7}{48} \, C_{2} \overline{C_{1}} & \overline{A_{1}} & -1 & 3 & 0 &\textbf{(48)}\\
	    			\lambda_3 & \overline{A_{0}} & -1 & 3 & 0 &\textbf{(49)}\\
	    			A_{1} C_{1} & \overline{A_{3}} & -1 & 3 & 0 &\textbf{(50)}\\
	    			-\frac{1}{48} \, \overline{A_{1}} \overline{C_{1}} & C_{2} & -1 & 3 & 0 &\textbf{(51)}\\
	    			-2 \, A_{0} C_{1} {\left| A_{1} \right|}^{2} + A_{1} B_{1} & \overline{A_{2}} & -1 & 3 & 0 &\textbf{(52)}\\
	    			\frac{1}{24} \, {\left| A_{1} \right|}^{2} \overline{A_{0}} \overline{C_{1}} - \frac{1}{48} \, \overline{A_{1}} \overline{B_{1}} & C_{1} & -1 & 3 & 0 &\textbf{(53)}\\
	    			\lambda_4 & A_{0} & -1 & 3 & 0 &\textbf{(54)}\\
	    			\frac{5}{3} \, \alpha_{2} \overline{A_{0}} \overline{C_{1}} + C_{1} \overline{A_{0}} \overline{\alpha_{3}} - \frac{1}{3} \, B_{2} \overline{A_{1}} - \frac{2}{3} \, B_{1} \overline{A_{2}} - C_{1} \overline{A_{3}} + \frac{1}{3} \, {\left(2 \, B_{1} \overline{A_{0}} + 3 \, C_{1} \overline{A_{1}}\right)} \overline{\alpha_{1}} & A_{1} & -1 & 3 & 0 &\textbf{(55)}\\
	    			-\frac{1}{24} \, C_{1} \overline{A_{1}} & \overline{B_{1}} & -1 & 3 & 0 &\textbf{(56)}\\
	    			\frac{1}{12} \, C_{1} {\left| A_{1} \right|}^{2} \overline{A_{0}} - \frac{1}{24} \, C_{2} \overline{A_{1}} & \overline{C_{1}} & -1 & 3 & 0 &\textbf{(57)}\\
	    			\overline{A_{1}} \overline{C_{1}} & A_{0} & 0 & -1 & 0 &\textbf{(58)}\\
	    			4 \, \overline{A_{0}} \overline{C_{1}} \overline{\alpha_{1}} - \frac{1}{4} \, C_{1}^{2} - 4 \, \overline{A_{2}} \overline{C_{1}} - 2 \, \overline{A_{1}} \overline{C_{2}} & A_{0} & 0 & 0 & 1 &\textbf{(59)}\\
	    			4 \, A_{0} C_{1} \alpha_{1} - 4 \, A_{2} C_{1} - 2 \, A_{1} C_{2} - \frac{1}{4} \, \overline{C_{1}}^{2} & \overline{A_{0}} & 0 & 0 & 1 &\textbf{(60)}
	    			\end{dmatrix}
	    			\end{align*}
	    			\begin{align*}
	    			\begin{dmatrix}
	    			-\frac{1}{8} \, \overline{C_{1}}^{2} & \overline{A_{1}} & 0 & 1 & 0 &\textbf{(61)}\\
	    			-\frac{1}{4} \, \overline{C_{1}} \overline{C_{2}} & \overline{A_{0}} & 0 & 1 & 0 &\textbf{(62)}\\
	    			-\frac{1}{8} \, \overline{A_{1}} \overline{C_{1}} & \overline{C_{1}} & 0 & 1 & 0 &\textbf{(63)}\\
	    			\lambda_5 & A_{0} & 0 & 1 & 0 &\textbf{(64)}\\
	    			-C_{1} \overline{A_{1}} & A_{2} & 0 & 1 & 0 &\textbf{(65)}\\
	    			2 \, C_{1} {\left| A_{1} \right|}^{2} \overline{A_{0}} - C_{2} \overline{A_{1}} & A_{1} & 0 & 1 & 0 &\textbf{(66)}\\
	    			-\frac{1}{2} \, B_{1} C_{1} - 2 \, \gamma_{1} \overline{A_{1}} & A_{0} & 0 & 1 & 1 &\textbf{(67)}\\
	    			8 \, A_{1} C_{1} {\left| A_{1} \right|}^{2} + 4 \, A_{0} C_{2} {\left| A_{1} \right|}^{2} + 4 \, A_{0} B_{1} \alpha_{1} + 4 \, A_{0} C_{1} \alpha_{4} - 4 \, A_{2} B_{1} - 2 \, A_{1} B_{3} & \overline{A_{0}} & 0 & 1 & 1 &\textbf{(68)}\\
	    			4 \, A_{0} C_{1} \alpha_{1} - 4 \, A_{2} C_{1} - 2 \, A_{1} C_{2} & \overline{A_{1}} & 0 & 1 & 1 &\textbf{(69)}\\
	    			-\frac{1}{16} \, \overline{C_{1}}^{2} & \overline{A_{2}} & 0 & 2 & 0 &\textbf{(70)}\\
	    			-\frac{1}{8} \, \overline{C_{1}} \overline{C_{2}} & \overline{A_{1}} & 0 & 2 & 0 &\textbf{(71)}\\
	    			-\frac{1}{8} \, C_{2} \overline{B_{1}} - \frac{1}{8} \, C_{1} \overline{B_{2}} - \frac{1}{4} \, C_{3} \overline{C_{1}} + \frac{1}{16} \, \gamma_{1} \overline{C_{1}} - \frac{1}{16} \, \overline{C_{2}}^{2} & \overline{A_{0}} & 0 & 2 & 0 &\textbf{(72)}\\
	    			-\frac{3}{256} \, C_{1} \overline{C_{1}} & C_{1} & 0 & 2 & 0 &\textbf{(73)}\\
	    			\frac{1}{8} \, \overline{A_{0}} \overline{C_{1}} \overline{\alpha_{1}} - \frac{1}{256} \, C_{1}^{2} - \frac{1}{8} \, \overline{A_{2}} \overline{C_{1}} - \frac{1}{16} \, \overline{A_{1}} \overline{C_{2}} & \overline{C_{1}} & 0 & 2 & 0 &\textbf{(74)}\\
	    			-\frac{1}{24} \, \overline{A_{1}} \overline{C_{1}} & \overline{C_{2}} & 0 & 2 & 0 &\textbf{(75)}\\
	    			\lambda_6 & A_{0} & 0 & 2 & 0 &\textbf{(76)}\\
	    			B_{1} {\left| A_{1} \right|}^{2} \overline{A_{0}} + 2 \, C_{1} {\left| A_{1} \right|}^{2} \overline{A_{1}} + C_{2} \overline{A_{0}} \overline{\alpha_{1}} + C_{1} \overline{A_{0}} \overline{\alpha_{4}} - \frac{1}{2} \, B_{3} \overline{A_{1}} - C_{2} \overline{A_{2}} & A_{1} & 0 & 2 & 0 &\textbf{(77)}\\
	    			C_{1} \overline{A_{0}} \overline{\alpha_{1}} - \frac{1}{2} \, B_{1} \overline{A_{1}} - C_{1} \overline{A_{2}} & A_{2} & 0 & 2 & 0 &\textbf{(78)}\\
	    			-\frac{1}{64} \, C_{1}^{2} & \overline{C_{1}} & 0 & 2 & 1 &\textbf{(79)}\\
	    			2 \, \overline{A_{0}} \overline{C_{1}} \overline{\zeta_{0}} - \frac{1}{4} \, B_{1}^{2} - \frac{1}{2} \, B_{2} C_{1} + 2 \, {\left(\overline{A_{0}} \overline{\alpha_{1}} - \overline{A_{2}}\right)} \gamma_{1} - \frac{3}{2} \, E_{1} \overline{C_{1}} - \frac{1}{2} \, \overline{A_{1}} \overline{\gamma_{2}} & A_{0} & 0 & 2 & 1 &\textbf{(80)}\\
	    			\lambda_7 & \overline{A_{0}} & 0 & 2 & 1 &\textbf{(81)}\\
	    			8 \, A_{1} C_{1} {\left| A_{1} \right|}^{2} + 4 \, A_{0} C_{2} {\left| A_{1} \right|}^{2} + 4 \, A_{0} B_{1} \alpha_{1} + 4 \, A_{0} C_{1} \alpha_{4} - 4 \, A_{2} B_{1} - 2 \, A_{1} B_{3} & \overline{A_{1}} & 0 & 2 & 1 &\textbf{(82)}\\
	    			4 \, A_{0} C_{1} \alpha_{1} - 4 \, A_{2} C_{1} - 2 \, A_{1} C_{2} & \overline{A_{2}} & 0 & 2 & 1 &\textbf{(83)}\\
	    			-A_{1} \overline{C_{1}} & \overline{A_{0}} & 1 & -2 & 0 &\textbf{(84)}\\
	    			\overline{A_{1}} \overline{C_{1}} & A_{1} & 1 & -1 & 0 &\textbf{(85)}\\
	    			-2 \, {\left| A_{1} \right|}^{2} \overline{A_{0}} \overline{C_{1}} + \overline{A_{1}} \overline{B_{1}} & A_{0} & 1 & -1 & 0 &\textbf{(86)}\\
	    			-A_{1} \overline{C_{1}} & \overline{A_{1}} & 1 & -1 & 0 &\textbf{(87)}\\
	    			2 \, A_{0} {\left| A_{1} \right|}^{2} \overline{C_{1}} - A_{1} \overline{C_{2}} & \overline{A_{0}} & 1 & -1 & 0 &\textbf{(88)}\\
	    			-\frac{1}{8} \, C_{1}^{2} & A_{1} & 1 & 0 & 0 &\textbf{(89)}\\
	    			-\frac{1}{4} \, C_{1} C_{2} & A_{0} & 1 & 0 & 0 &\textbf{(90)}
	    			\end{dmatrix}
	    			\end{align*}
	    			\begin{align*}
	    			\begin{dmatrix}
	    			-\frac{1}{8} \, A_{1} C_{1} & C_{1} & 1 & 0 & 0 &\textbf{(91)}\\
	    			\lambda_8 & \overline{A_{0}} & 1 & 0 & 0 &\textbf{(92)}\\
	    			-A_{1} \overline{C_{1}} & \overline{A_{2}} & 1 & 0 & 0 &\textbf{(93)}\\
	    			2 \, A_{0} {\left| A_{1} \right|}^{2} \overline{C_{1}} - A_{1} \overline{C_{2}} & \overline{A_{1}} & 1 & 0 & 0 &\textbf{(94)}\\
	    			-2 \, A_{1} \gamma_{1} - \frac{1}{2} \, \overline{B_{1}} \overline{C_{1}} & \overline{A_{0}} & 1 & 0 & 1 &\textbf{(95)}\\
	    			8 \, {\left| A_{1} \right|}^{2} \overline{A_{1}} \overline{C_{1}} + 4 \, {\left| A_{1} \right|}^{2} \overline{A_{0}} \overline{C_{2}} + 4 \, \overline{A_{0}} \overline{B_{1}} \overline{\alpha_{1}} + 4 \, \overline{A_{0}} \overline{C_{1}} \overline{\alpha_{4}} - 4 \, \overline{A_{2}} \overline{B_{1}} - 2 \, \overline{A_{1}} \overline{B_{3}} & A_{0} & 1 & 0 & 1 &\textbf{(96)}\\
	    			4 \, \overline{A_{0}} \overline{C_{1}} \overline{\alpha_{1}} - 4 \, \overline{A_{2}} \overline{C_{1}} - 2 \, \overline{A_{1}} \overline{C_{2}} & A_{1} & 1 & 0 & 1 &\textbf{(97)}\\
	    			\lambda_9 & A_{1} & 1 & 1 & 0 &\textbf{(98)}\\
	    			\lambda_{10} & \overline{A_{1}} & 1 & 1 & 0 &\textbf{(99)}\\
	    			\lambda_{11} & A_{0} & 1 & 1 & 0 &\textbf{(100)}\\
	    			\lambda_{12} & \overline{A_{0}} & 1 & 1 & 0 &\textbf{(101)}\\
	    			-\frac{1}{8} \, A_{1} C_{1} & B_{1} & 1 & 1 & 0 &\textbf{(102)}\\
	    			-\frac{1}{8} \, \overline{A_{1}} \overline{C_{1}} & \overline{B_{1}} & 1 & 1 & 0 &\textbf{(103)}\\
	    			\frac{1}{4} \, A_{0} C_{1} {\left| A_{1} \right|}^{2} - \frac{1}{8} \, A_{1} B_{1} & C_{1} & 1 & 1 & 0 &\textbf{(104)}\\
	    			\frac{1}{4} \, {\left| A_{1} \right|}^{2} \overline{A_{0}} \overline{C_{1}} - \frac{1}{8} \, \overline{A_{1}} \overline{B_{1}} & \overline{C_{1}} & 1 & 1 & 0 &\textbf{(105)}\\
	    			-A_{1} \overline{C_{1}} & \overline{A_{3}} & 1 & 1 & 0 &\textbf{(106)}\\
	    			-C_{1} \overline{A_{1}} & A_{3} & 1 & 1 & 0 &\textbf{(107)}\\
	    			2 \, A_{0} {\left| A_{1} \right|}^{2} \overline{C_{1}} - A_{1} \overline{C_{2}} & \overline{A_{2}} & 1 & 1 & 0 &\textbf{(108)}\\
	    			2 \, C_{1} {\left| A_{1} \right|}^{2} \overline{A_{0}} - C_{2} \overline{A_{1}} & A_{2} & 1 & 1 & 0 &\textbf{(109)}\\
	    			4 \, A_{0} \gamma_{1} {\left| A_{1} \right|}^{2} - A_{1} \overline{\gamma_{2}} & \overline{A_{0}} & 1 & 1 & 1 &\textbf{(110)}\\
	    			4 \, \gamma_{1} {\left| A_{1} \right|}^{2} \overline{A_{0}} - \gamma_{2} \overline{A_{1}} & A_{0} & 1 & 1 & 1 &\textbf{(111)}\\
	    			-2 \, A_{1} \gamma_{1} & \overline{A_{1}} & 1 & 1 & 1 &\textbf{(112)}\\
	    			-2 \, \gamma_{1} \overline{A_{1}} & A_{1} & 1 & 1 & 1 &\textbf{(113)}\\
	    			-\frac{3}{16} \, C_{1} \overline{C_{1}} & A_{0} & 2 & -2 & 0 &\textbf{(114)}\\
	    			A_{0} \alpha_{1} \overline{C_{1}} - \frac{1}{2} \, A_{1} \overline{B_{1}} - A_{2} \overline{C_{1}} & \overline{A_{0}} & 2 & -2 & 0 &\textbf{(115)}\\
	    			-2 \, {\left| A_{1} \right|}^{2} \overline{A_{0}} \overline{B_{1}} - \alpha_{4} \overline{A_{0}} \overline{C_{1}} - \alpha_{1} \overline{A_{1}} \overline{C_{1}} + \overline{A_{1}} \overline{B_{2}} - \frac{1}{8} \, B_{1} \overline{C_{1}} - \frac{3}{16} \, C_{1} \overline{C_{2}} & A_{0} & 2 & -1 & 0 &\textbf{(116)}\\
	    			\overline{A_{1}} \overline{C_{1}} & A_{2} & 2 & -1 & 0 &\textbf{(117)}\\
	    			-2 \, {\left| A_{1} \right|}^{2} \overline{A_{0}} \overline{C_{1}} + \overline{A_{1}} \overline{B_{1}} & A_{1} & 2 & -1 & 0 &\textbf{(118)}\\
	    			A_{0} {\left| A_{1} \right|}^{2} \overline{B_{1}} + 2 \, A_{1} {\left| A_{1} \right|}^{2} \overline{C_{1}} + A_{0} \alpha_{4} \overline{C_{1}} + A_{0} \alpha_{1} \overline{C_{2}} - \frac{1}{2} \, A_{1} \overline{B_{3}} - A_{2} \overline{C_{2}} & \overline{A_{0}} & 2 & -1 & 0 &\textbf{(119)}\\
	    			A_{0} \alpha_{1} \overline{C_{1}} - \frac{1}{2} \, A_{1} \overline{B_{1}} - A_{2} \overline{C_{1}} & \overline{A_{1}} & 2 & -1 & 0&\textbf{(120)} 
	    			\end{dmatrix}
	    			\end{align*}
	    			\begin{align*}
	    			\begin{dmatrix}
	    			-\frac{1}{16} \, C_{1}^{2} & A_{2} & 2 & 0 & 0 &\textbf{(121)}\\
	    			-\frac{1}{8} \, C_{1} C_{2} & A_{1} & 2 & 0 & 0 &\textbf{(122)}\\
	    			-\frac{1}{16} \, C_{2}^{2} - \frac{1}{4} \, C_{1} C_{3} + \frac{1}{16} \, C_{1} \gamma_{1} - \frac{1}{8} \, B_{2} \overline{C_{1}} - \frac{1}{8} \, B_{1} \overline{C_{2}} & A_{0} & 2 & 0 & 0 &\textbf{(123)}\\
	    			-\frac{3}{256} \, C_{1} \overline{C_{1}} & \overline{C_{1}} & 2 & 0 & 0 &\textbf{(124)}\\
	    			\frac{1}{8} \, A_{0} C_{1} \alpha_{1} - \frac{1}{8} \, A_{2} C_{1} - \frac{1}{16} \, A_{1} C_{2} - \frac{1}{256} \, \overline{C_{1}}^{2} & C_{1} & 2 & 0 & 0 &\textbf{(125)}\\
	    			-\frac{1}{24} \, A_{1} C_{1} & C_{2} & 2 & 0 & 0 &\textbf{(126)}\\
	    			\lambda_{13} & \overline{A_{0}} & 2 & 0 & 0 &\textbf{(127)}\\
	    			A_{0} {\left| A_{1} \right|}^{2} \overline{B_{1}} + 2 \, A_{1} {\left| A_{1} \right|}^{2} \overline{C_{1}} + A_{0} \alpha_{4} \overline{C_{1}} + A_{0} \alpha_{1} \overline{C_{2}} - \frac{1}{2} \, A_{1} \overline{B_{3}} - A_{2} \overline{C_{2}} & \overline{A_{1}} & 2 & 0 & 0 &\textbf{(128)}\\
	    			A_{0} \alpha_{1} \overline{C_{1}} - \frac{1}{2} \, A_{1} \overline{B_{1}} - A_{2} \overline{C_{1}} & \overline{A_{2}} & 2 & 0 & 0 &\textbf{(129)}\\
	    			-\frac{1}{64} \, \overline{C_{1}}^{2} & C_{1} & 2 & 0 & 1 &\textbf{(130)}\\
	    			2 \, A_{0} C_{1} \zeta_{0} + 2 \, {\left(A_{0} \alpha_{1} - A_{2}\right)} \gamma_{1} - \frac{1}{2} \, A_{1} \gamma_{2} - \frac{1}{4} \, \overline{B_{1}}^{2} - \frac{1}{2} \, \overline{B_{2}} \overline{C_{1}} - \frac{3}{2} \, C_{1} \overline{E_{1}} & \overline{A_{0}} & 2 & 0 & 1 &\textbf{(131)}\\
	    			\lambda_{14} & A_{0} & 2 & 0 & 1 &\textbf{(132)}\\
	    			8 \, {\left| A_{1} \right|}^{2} \overline{A_{1}} \overline{C_{1}} + 4 \, {\left| A_{1} \right|}^{2} \overline{A_{0}} \overline{C_{2}} + 4 \, \overline{A_{0}} \overline{B_{1}} \overline{\alpha_{1}} + 4 \, \overline{A_{0}} \overline{C_{1}} \overline{\alpha_{4}} - 4 \, \overline{A_{2}} \overline{B_{1}} - 2 \, \overline{A_{1}} \overline{B_{3}} & A_{1} & 2 & 0 & 1 &\textbf{(133)}\\
	    			4 \, \overline{A_{0}} \overline{C_{1}} \overline{\alpha_{1}} - 4 \, \overline{A_{2}} \overline{C_{1}} - 2 \, \overline{A_{1}} \overline{C_{2}} & A_{2} & 2 & 0 & 1 &\textbf{(134)}\\
	    			-\frac{7}{48} \, C_{1} \overline{C_{1}} & A_{1} & 3 & -2 & 0 &\textbf{(135)}\\
	    			C_{1} \overline{A_{0}} \overline{\alpha_{2}} - \frac{2}{15} \, C_{1} \overline{B_{1}} - \frac{7}{48} \, C_{2} \overline{C_{1}} & A_{0} & 3 & -2 & 0 &\textbf{(136)}\\
	    			-\frac{1}{48} \, A_{1} C_{1} & \overline{C_{1}} & 3 & -2 & 0 &\textbf{(137)}\\
	    			A_{0} \alpha_{3} \overline{C_{1}} + \frac{5}{3} \, A_{0} C_{1} \overline{\alpha_{2}} + \frac{1}{3} \, {\left(2 \, A_{0} \overline{B_{1}} + 3 \, A_{1} \overline{C_{1}}\right)} \alpha_{1} - \frac{2}{3} \, A_{2} \overline{B_{1}} - \frac{1}{3} \, A_{1} \overline{B_{2}} - A_{3} \overline{C_{1}} & \overline{A_{0}} & 3 & -2 & 0 &\textbf{(138)}\\
	    			-\frac{1}{24} \, A_{1} \overline{C_{1}} & C_{1} & 3 & -2 & 0 &\textbf{(139)}\\
	    			-2 \, {\left| A_{1} \right|}^{2} \overline{A_{0}} \overline{B_{1}} - \alpha_{4} \overline{A_{0}} \overline{C_{1}} - \alpha_{1} \overline{A_{1}} \overline{C_{1}} + \overline{A_{1}} \overline{B_{2}} - \frac{1}{12} \, B_{1} \overline{C_{1}} - \frac{7}{48} \, C_{1} \overline{C_{2}} & A_{1} & 3 & -1 & 0 &\textbf{(140)}\\
	    			\lambda_{15} & A_{0} & 3 & -1 & 0 &\textbf{(141)}\\
	    			\overline{A_{1}} \overline{C_{1}} & A_{3} & 3 & -1 & 0 &\textbf{(142)}\\
	    			-\frac{1}{48} \, A_{1} C_{1} & \overline{C_{2}} & 3 & -1 & 0 &\textbf{(143)}\\
	    			-2 \, {\left| A_{1} \right|}^{2} \overline{A_{0}} \overline{C_{1}} + \overline{A_{1}} \overline{B_{1}} & A_{2} & 3 & -1 & 0 &\textbf{(144)}\\
	    			\frac{1}{24} \, A_{0} C_{1} {\left| A_{1} \right|}^{2} - \frac{1}{48} \, A_{1} B_{1} & \overline{C_{1}} & 3 & -1 & 0 &\textbf{(145)}\\
	    			\lambda_{16} & \overline{A_{0}} & 3 & -1 & 0 &\textbf{(146)}\\
	    			A_{0} \alpha_{3} \overline{C_{1}} + \frac{5}{3} \, A_{0} C_{1} \overline{\alpha_{2}} + \frac{1}{3} \, {\left(2 \, A_{0} \overline{B_{1}} + 3 \, A_{1} \overline{C_{1}}\right)} \alpha_{1} - \frac{2}{3} \, A_{2} \overline{B_{1}} - \frac{1}{3} \, A_{1} \overline{B_{2}} - A_{3} \overline{C_{1}} & \overline{A_{1}} & 3 & -1 & 0 &\textbf{(147)}\\
	    			-\frac{1}{24} \, A_{1} \overline{C_{1}} & B_{1} & 3 & -1 & 0 &\textbf{(148)}\\
	    			\frac{1}{12} \, A_{0} {\left| A_{1} \right|}^{2} \overline{C_{1}} - \frac{1}{24} \, A_{1} \overline{C_{2}} & C_{1} & 3 & -1 & 0 &\textbf{(149)}\\
	    			-\frac{1}{16} \, \overline{C_{1}}^{2} & A_{0} & 4 & -4 & 0 &\textbf{(150)}
	    			\end{dmatrix}
	    			\end{align*}
	    			\begin{align}\label{4h4}
	    			\begin{dmatrix}
	    			-\frac{1}{8} \, \overline{C_{1}} \overline{C_{2}} + \frac{1}{3} \, \overline{A_{1}} \overline{E_{1}} & A_{0} & 4 & -3 & 0 &\textbf{(151)}\\
	    			\frac{1}{32} \, A_{0} C_{1} \alpha_{1} - \frac{1}{32} \, A_{2} C_{1} - \frac{1}{64} \, A_{1} C_{2} - \frac{1}{128} \, \overline{C_{1}}^{2} & \overline{C_{1}} & 4 & -2 & 0 &\textbf{(152)}\\
	    			\lambda_{17} & \overline{A_{0}} & 4 & -2 & 0 &\textbf{(153)}\\
	    			-\frac{1}{8} \, C_{1} \overline{C_{1}} & A_{2} & 4 & -2 & 0 &\textbf{(154)}\\
	    			C_{1} \overline{A_{0}} \overline{\alpha_{2}} - \frac{9}{80} \, C_{1} \overline{B_{1}} - \frac{1}{8} \, C_{2} \overline{C_{1}} & A_{1} & 4 & -2 & 0 &\textbf{(155)}\\
	    			\lambda_{18} & A_{0} & 4 & -2 & 0 &\textbf{(156)}\\
	    			-\frac{1}{80} \, A_{1} C_{1} & \overline{B_{1}} & 4 & -2 & 0 &\textbf{(157)}\\
	    			-\frac{1}{48} \, A_{1} \overline{C_{1}} & C_{2} & 4 & -2 & 0 &\textbf{(158)}\\
	    			\frac{1}{16} \, A_{0} \alpha_{1} \overline{C_{1}} - \frac{1}{32} \, A_{1} \overline{B_{1}} - \frac{1}{16} \, A_{2} \overline{C_{1}} & C_{1} & 4 & -2 & 0 &\textbf{(159)}\\
	    			-\frac{1}{4} \, \gamma_{1} \overline{C_{1}} & A_{0} & 4 & -2 & 1 &\textbf{(160)}\\
	    			A_{0} \zeta_{0} \overline{C_{1}} - \frac{1}{2} \, \overline{C_{1}} \overline{E_{1}} & \overline{A_{0}} & 4 & -2 & 1 &\textbf{(161)}\\
	    			-\frac{9}{160} \, \overline{C_{1}}^{2} & A_{1} & 5 & -4 & 0 &\textbf{(162)}\\
	    			\frac{1}{2} \, \overline{A_{0}} \overline{C_{1}} \overline{\alpha_{2}} - \frac{17}{160} \, \overline{B_{1}} \overline{C_{1}} & A_{0} & 5 & -4 & 0&\textbf{(163)} \\
	    			-\frac{1}{80} \, A_{1} \overline{C_{1}} & \overline{C_{1}} & 5 & -4 & 0 &\textbf{(164)}\\
	    			A_{0} \overline{C_{1}} \overline{\alpha_{2}} - \frac{1}{5} \, A_{1} \overline{E_{1}} & \overline{A_{0}} & 5 & -4 & 0 &\textbf{(165)}\\
	    			-\frac{9}{80} \, \overline{C_{1}} \overline{C_{2}} + \frac{1}{3} \, \overline{A_{1}} \overline{E_{1}} & A_{1} & 5 & -3 & 0 &\textbf{(166)}\\
	    			\frac{1}{120} \, {\left| A_{1} \right|}^{2} \overline{C_{1}}^{2} - \frac{2}{3} \, {\left| A_{1} \right|}^{2} \overline{A_{0}} \overline{E_{1}} + \frac{1}{3} \, \alpha_{6} \overline{A_{0}} \overline{C_{1}} - \frac{13}{120} \, \overline{B_{3}} \overline{C_{1}} - \frac{5}{48} \, \overline{B_{1}} \overline{C_{2}} + \frac{1}{3} \, \overline{A_{1}} \overline{E_{2}} + \frac{1}{3} \, {\left(\overline{A_{1}} \overline{C_{1}} + 2 \, \overline{A_{0}} \overline{C_{2}}\right)} \overline{\alpha_{2}} & A_{0} & 5 & -3 & 0 &\textbf{(167)}\\
	    			-\frac{1}{80} \, A_{1} \overline{C_{1}} & \overline{C_{2}} & 5 & -3 & 0 &\textbf{(168)}\\
	    			\frac{1}{40} \, A_{0} {\left| A_{1} \right|}^{2} \overline{C_{1}} - \frac{1}{80} \, A_{1} \overline{C_{2}} & \overline{C_{1}} & 5 & -3 & 0&\textbf{(169)} \\
	    			\frac{2}{5} \, A_{0} {\left| A_{1} \right|}^{2} \overline{E_{1}} + A_{0} \alpha_{6} \overline{C_{1}} + A_{0} \overline{C_{2}} \overline{\alpha_{2}} - \frac{1}{5} \, A_{1} \overline{E_{3}} & \overline{A_{0}} & 5 & -3 & 0 &\textbf{(170)}\\
	    			A_{0} \overline{C_{1}} \overline{\alpha_{2}} - \frac{1}{5} \, A_{1} \overline{E_{1}} & \overline{A_{1}} & 5 & -3 & 0 &\textbf{(171)}\\
	    			\frac{1}{2} \, \alpha_{7} \overline{A_{0}} \overline{C_{1}} + \frac{1}{96} \, \alpha_{1} \overline{C_{1}}^{2} + \frac{1}{2} \, \overline{A_{0}} \overline{B_{1}} \overline{\alpha_{2}} - \frac{11}{240} \, \overline{B_{1}}^{2} - \frac{3}{32} \, \overline{B_{2}} \overline{C_{1}} - \frac{7}{96} \, C_{1} \overline{E_{1}} & A_{0} & 6 & -4 & 0 &\textbf{(172)}\\
	    			-\frac{5}{96} \, \overline{C_{1}}^{2} & A_{2} & 6 & -4 & 0 &\textbf{(173)}\\
	    			\frac{1}{2} \, \overline{A_{0}} \overline{C_{1}} \overline{\alpha_{2}} - \frac{47}{480} \, \overline{B_{1}} \overline{C_{1}} & A_{1} & 6 & -4 & 0 &\textbf{(174)}\\
	    			-\frac{1}{120} \, A_{1} \overline{C_{1}} & \overline{B_{1}} & 6 & -4 & 0 &\textbf{(175)}\\
	    			A_{0} \alpha_{7} \overline{C_{1}} + \frac{1}{3} \, A_{0} \alpha_{1} \overline{E_{1}} - \frac{1}{3} \, A_{2} \overline{E_{1}} - \frac{1}{6} \, A_{1} \overline{E_{2}} + \frac{1}{6} \, {\left(5 \, A_{0} \overline{B_{1}} + 6 \, A_{1} \overline{C_{1}}\right)} \overline{\alpha_{2}} & \overline{A_{0}} & 6 & -4 & 0 &\textbf{(176)}\\
	    			\frac{1}{48} \, A_{0} \alpha_{1} \overline{C_{1}} - \frac{1}{96} \, A_{1} \overline{B_{1}} - \frac{1}{48} \, A_{2} \overline{C_{1}} & \overline{C_{1}} & 6 & -4 & 0 &\textbf{(177)}\\
	    			-\frac{7}{96} \, \overline{C_{1}} \overline{E_{1}} & A_{0} & 8 & -6 & 0&\textbf{(178)}
	    			\end{dmatrix}
	    		\end{align}
	    		\small
	    		where
	    		\begin{align*}
	    			\lambda_{1}&=-C_{1} \overline{A_{0}} \overline{\alpha_{1}}^{2} + \frac{3}{2} \, \overline{A_{0}} \overline{C_{1}} \overline{\alpha_{7}} + \frac{1}{8} \, C_{1} E_{1} + \frac{1}{4} \, {\left(6 \, \overline{A_{1}} \overline{C_{1}} + 5 \, \overline{A_{0}} \overline{C_{2}}\right)} \alpha_{2} - \frac{1}{4} \, B_{4} \overline{A_{1}} - \frac{1}{2} \, B_{2} \overline{A_{2}} - \frac{3}{4} \, B_{1} \overline{A_{3}} - C_{1} \overline{A_{4}}\\
	    			& + \frac{1}{4} \, {\left(2 \, B_{2} \overline{A_{0}} + 3 \, B_{1} \overline{A_{1}} + 4 \, C_{1} \overline{A_{2}}\right)} \overline{\alpha_{1}} + \frac{1}{4} \, {\left(3 \, B_{1} \overline{A_{0}} + 4 \, C_{1} \overline{A_{1}}\right)} \overline{\alpha_{3}}\\
	    			\lambda_{2}&=-A_{0} E_{1} \alpha_{1} + A_{0} \alpha_{2} \overline{C_{2}} + \frac{1}{48} \, C_{1} \overline{C_{1}} \overline{\alpha_{1}} + A_{0} \overline{C_{1}} \overline{\alpha_{7}} - \frac{1}{4} \, A_{0} C_{1} \overline{\zeta_{0}} - \frac{1}{16} \, C_{2}^{2} - \frac{1}{4} \, C_{1} C_{3} + A_{2} E_{1} + \frac{1}{2} \, A_{1} E_{3} + \frac{1}{64} \, C_{1} \gamma_{1}\\
	    			&  - \frac{5}{48} \, B_{2} \overline{C_{1}} - \frac{9}{80} \, B_{1} \overline{C_{2}}\\
	    			\lambda_{3}&=-2 \, A_{0} B_{2} {\left| A_{1} \right|}^{2} - A_{1} C_{1} \overline{\alpha_{3}} - A_{0} B_{1} \overline{\alpha_{4}} - A_{0} C_{1} \overline{\alpha_{5}} + A_{0} \overline{C_{1}} \overline{\alpha_{6}} + A_{1} B_{4} + {\left(2 \, A_{0} \overline{B_{1}} + A_{1} \overline{C_{1}}\right)} \alpha_{2} - \frac{7}{120} \, B_{1} \overline{B_{1}}\\
	    			&  - \frac{1}{8} \, C_{1} \overline{B_{3}} - \frac{1}{120} \, {\left(7 \, C_{1} {\left| A_{1} \right|}^{2} + 16 \, B_{3}\right)} \overline{C_{1}} - \frac{7}{48} \, C_{2} \overline{C_{2}} + {\left(4 \, A_{0} C_{1} {\left| A_{1} \right|}^{2} - A_{1} B_{1}\right)} \overline{\alpha_{1}}\\
	    			\lambda_{4}&=\frac{2}{3} \, B_{2} {\left| A_{1} \right|}^{2} \overline{A_{0}} + \frac{4}{3} \, B_{1} {\left| A_{1} \right|}^{2} \overline{A_{1}} + 2 \, C_{1} {\left| A_{1} \right|}^{2} \overline{A_{2}} + \frac{5}{3} \, \alpha_{2} \overline{A_{0}} \overline{B_{1}} + C_{2} \overline{A_{0}} \overline{\alpha_{3}} + C_{1} \overline{A_{0}} \overline{\alpha_{5}} + \frac{5}{3} \, \overline{A_{0}} \overline{C_{1}} \overline{\alpha_{6}} - \frac{1}{3} \, B_{5} \overline{A_{1}}\\
	    			&  - \frac{2}{3} \, B_{3} \overline{A_{2}} - C_{2} \overline{A_{3}} - \frac{1}{3} \, {\left(12 \, C_{1} {\left| A_{1} \right|}^{2} \overline{A_{0}} - 2 \, B_{3} \overline{A_{0}} - 3 \, C_{2} \overline{A_{1}}\right)} \overline{\alpha_{1}} + \frac{1}{3} \, {\left(2 \, B_{1} \overline{A_{0}} + 3 \, C_{1} \overline{A_{1}}\right)} \overline{\alpha_{4}}\\
	    			\lambda_{5}&=2 \, C_{2} {\left| A_{1} \right|}^{2} \overline{A_{0}} + C_{1} \alpha_{4} \overline{A_{0}} + C_{1} \alpha_{1} \overline{A_{1}} + 3 \, \overline{A_{0}} \overline{C_{1}} \overline{\alpha_{3}} - \frac{1}{8} \, B_{1} C_{1} - 2 \, C_{3} \overline{A_{1}} + \gamma_{1} \overline{A_{1}} - 3 \, \overline{A_{3}} \overline{C_{1}} - 2 \, \overline{A_{2}} \overline{C_{2}}\\
	    			&  + {\left(3 \, \overline{A_{1}} \overline{C_{1}} + 2 \, \overline{A_{0}} \overline{C_{2}}\right)} \overline{\alpha_{1}}\\
	    			\lambda_{6}&=B_{3} {\left| A_{1} \right|}^{2} \overline{A_{0}} + 2 \, C_{2} {\left| A_{1} \right|}^{2} \overline{A_{1}} - 2 \, \overline{A_{0}} \overline{C_{1}} \overline{\alpha_{1}}^{2} + C_{1} \beta \overline{A_{0}} + C_{2} \overline{A_{0}} \overline{\alpha_{4}} - \frac{1}{4} \, \overline{A_{0}} \overline{C_{1}} \overline{\zeta_{0}} - \frac{1}{16} \, B_{1}^{2} - \frac{1}{8} \, {\left(32 \, {\left| A_{1} \right|}^{4} \overline{A_{0}} + B_{2}\right)} C_{1}\\
	    			&  + \frac{1}{2} \, {\left(B_{1} \overline{A_{1}} + 2 \, C_{1} \overline{A_{2}}\right)} \alpha_{1} + \frac{1}{2} \, {\left(B_{1} \overline{A_{0}} + 2 \, C_{1} \overline{A_{1}}\right)} \alpha_{4} - \frac{1}{2} \, {\left(\overline{A_{0}} \overline{\alpha_{1}} - \overline{A_{2}}\right)} \gamma_{1} - 2 \, C_{3} \overline{A_{2}} + \frac{3}{16} \, E_{1} \overline{C_{1}} - 2 \, \overline{A_{4}} \overline{C_{1}} - \frac{3}{2} \, \overline{A_{3}} \overline{C_{2}}\\
	    			&  - \frac{1}{2} \, \overline{A_{1}} \overline{C_{4}} - \frac{1}{2} \, {\left(4 \, C_{1} \alpha_{1} \overline{A_{0}} - 4 \, C_{3} \overline{A_{0}} - 4 \, \overline{A_{2}} \overline{C_{1}} - 3 \, \overline{A_{1}} \overline{C_{2}}\right)} \overline{\alpha_{1}} + \frac{1}{2} \, {\left(4 \, \overline{A_{1}} \overline{C_{1}} + 3 \, \overline{A_{0}} \overline{C_{2}}\right)} \overline{\alpha_{3}} + \frac{1}{8} \, \overline{A_{1}} \overline{\gamma_{2}}\\
	    			\lambda_{7}&=8 \, A_{1} B_{1} {\left| A_{1} \right|}^{2} + 4 \, A_{0} B_{3} {\left| A_{1} \right|}^{2} + 4 \, A_{0} B_{2} \alpha_{1} + 4 \, A_{0} B_{1} \alpha_{4} + 4 \, A_{0} C_{1} \beta + A_{0} \overline{C_{1}} \overline{\zeta_{0}} - 4 \, A_{2} B_{2} - 2 \, A_{1} B_{5}\\
	    			&  - \frac{1}{4} \, {\left(64 \, A_{0} {\left| A_{1} \right|}^{4} + \overline{B_{2}}\right)} C_{1} - \frac{1}{8} \, C_{2} \overline{B_{1}} - \frac{3}{8} \, \gamma_{1} \overline{C_{1}} - 2 \, {\left(4 \, A_{0} C_{1} \alpha_{1} - 2 \, A_{2} C_{1} - A_{1} C_{2}\right)} \overline{\alpha_{1}} + 2 \, {\left(2 \, A_{1} C_{1} + A_{0} C_{2}\right)} \overline{\alpha_{4}}\\
	    			\lambda_{8}&=2 \, A_{0} {\left| A_{1} \right|}^{2} \overline{C_{2}} + 3 \, A_{0} C_{1} \alpha_{3} + A_{1} \overline{C_{1}} \overline{\alpha_{1}} + A_{0} \overline{C_{1}} \overline{\alpha_{4}} - 3 \, A_{3} C_{1} - 2 \, A_{2} C_{2} - 2 \, A_{1} C_{3} + {\left(3 \, A_{1} C_{1} + 2 \, A_{0} C_{2}\right)} \alpha_{1} + A_{1} \gamma_{1}\\
	    			&  - \frac{1}{8} \, \overline{B_{1}} \overline{C_{1}}\\
	    			\lambda_{9}&=2 \, C_{2} {\left| A_{1} \right|}^{2} \overline{A_{0}} + C_{1} \alpha_{4} \overline{A_{0}} + C_{1} \alpha_{1} \overline{A_{1}} + 3 \, \overline{A_{0}} \overline{C_{1}} \overline{\alpha_{3}} - \frac{3}{8} \, B_{1} C_{1} - 2 \, C_{3} \overline{A_{1}} + \gamma_{1} \overline{A_{1}} - 3 \, \overline{A_{3}} \overline{C_{1}} - 2 \, \overline{A_{2}} \overline{C_{2}}\\
	    			&  + {\left(3 \, \overline{A_{1}} \overline{C_{1}} + 2 \, \overline{A_{0}} \overline{C_{2}}\right)} \overline{\alpha_{1}}\\
	    			\lambda_{10}&=2 \, A_{0} {\left| A_{1} \right|}^{2} \overline{C_{2}} + 3 \, A_{0} C_{1} \alpha_{3} + A_{1} \overline{C_{1}} \overline{\alpha_{1}} + A_{0} \overline{C_{1}} \overline{\alpha_{4}} - 3 \, A_{3} C_{1} - 2 \, A_{2} C_{2} - 2 \, A_{1} C_{3} + {\left(3 \, A_{1} C_{1} + 2 \, A_{0} C_{2}\right)} \alpha_{1} + A_{1} \gamma_{1}\\
	    			&  - \frac{3}{8} \, \overline{B_{1}} \overline{C_{1}}\\
	    			\lambda_{11}&=\frac{1}{12} \, C_{1}^{2} {\left| A_{1} \right|}^{2} + 4 \, C_{3} {\left| A_{1} \right|}^{2} \overline{A_{0}} - 2 \, \gamma_{1} {\left| A_{1} \right|}^{2} \overline{A_{0}} + 6 \, {\left| A_{1} \right|}^{2} \overline{A_{2}} \overline{C_{1}} + 4 \, {\left| A_{1} \right|}^{2} \overline{A_{1}} \overline{C_{2}} + C_{2} \alpha_{4} \overline{A_{0}} + C_{1} \alpha_{5} \overline{A_{0}} + C_{1} \alpha_{3} \overline{A_{1}} \\
	    			& + 3 \, \overline{A_{0}} \overline{B_{1}} \overline{\alpha_{3}} + 3 \, \overline{A_{0}} \overline{C_{1}} \overline{\alpha_{5}} - \frac{1}{3} \, B_{3} C_{1} - \frac{3}{8} \, B_{1} C_{2} - {\left(4 \, C_{1} {\left| A_{1} \right|}^{2} \overline{A_{0}} - C_{2} \overline{A_{1}}\right)} \alpha_{1} - C_{4} \overline{A_{1}} + \frac{1}{2} \, \gamma_{2} \overline{A_{1}} - 3 \, \overline{A_{3}} \overline{B_{1}} - 2 \, \overline{A_{2}} \overline{B_{3}}\\
	    			&  - {\left(12 \, {\left| A_{1} \right|}^{2} \overline{A_{0}} \overline{C_{1}} - 3 \, \overline{A_{1}} \overline{B_{1}} - 2 \, \overline{A_{0}} \overline{B_{3}}\right)} \overline{\alpha_{1}} + {\left(3 \, \overline{A_{1}} \overline{C_{1}} + 2 \, \overline{A_{0}} \overline{C_{2}}\right)} \overline{\alpha_{4}}\\
	    			\lambda_{12}&=6 \, A_{2} C_{1} {\left| A_{1} \right|}^{2} + 4 \, A_{1} C_{2} {\left| A_{1} \right|}^{2} + 4 \, A_{0} C_{3} {\left| A_{1} \right|}^{2} - 2 \, A_{0} \gamma_{1} {\left| A_{1} \right|}^{2} + \frac{1}{12} \, {\left| A_{1} \right|}^{2} \overline{C_{1}}^{2} + 3 \, A_{0} B_{1} \alpha_{3} + 3 \, A_{0} C_{1} \alpha_{5} + A_{1} \overline{C_{1}} \overline{\alpha_{3}}\\
	    			&  + A_{0} \overline{C_{2}} \overline{\alpha_{4}} + A_{0} \overline{C_{1}} \overline{\alpha_{5}} - 3 \, A_{3} B_{1} - 2 \, A_{2} B_{3} - {\left(12 \, A_{0} C_{1} {\left| A_{1} \right|}^{2} - 3 \, A_{1} B_{1} - 2 \, A_{0} B_{3}\right)} \alpha_{1} + {\left(3 \, A_{1} C_{1} + 2 \, A_{0} C_{2}\right)} \alpha_{4}\\
	    			&  - \frac{1}{3} \, \overline{B_{3}} \overline{C_{1}} - \frac{3}{8} \, \overline{B_{1}} \overline{C_{2}} - A_{1} \overline{C_{4}} - {\left(4 \, A_{0} {\left| A_{1} \right|}^{2} \overline{C_{1}} - A_{1} \overline{C_{2}}\right)} \overline{\alpha_{1}} + \frac{1}{2} \, A_{1} \overline{\gamma_{2}}\\
	    			\lambda_{13}&=-2 \, A_{0} C_{1} \alpha_{1}^{2} + A_{0} {\left| A_{1} \right|}^{2} \overline{B_{3}} + 2 \, A_{1} {\left| A_{1} \right|}^{2} \overline{C_{2}} - \frac{1}{4} \, A_{0} C_{1} \zeta_{0} + A_{0} \beta \overline{C_{1}} + A_{0} \alpha_{4} \overline{C_{2}} - 2 \, A_{4} C_{1} - \frac{3}{2} \, A_{3} C_{2} - 2 \, A_{2} C_{3}\\
	    			&                                - \frac{1}{2} \, A_{1} C_{4} + \frac{1}{2} \, {\left(4 \, A_{2} C_{1} + 3 \, A_{1} C_{2} + 4 \, A_{0} C_{3}\right)} \alpha_{1} + \frac{1}{2} \, {\left(4 \, A_{1} C_{1} + 3 \, A_{0} C_{2}\right)} \alpha_{3} - \frac{1}{2} \, {\left(A_{0} \alpha_{1} - A_{2}\right)} \gamma_{1} + \frac{1}{8} \, A_{1} \gamma_{2} - \frac{1}{16} \, \overline{B_{1}}^{2}\\
	    			&  - \frac{1}{8} \, {\left(32 \, A_{0} {\left| A_{1} \right|}^{4} + \overline{B_{2}}\right)} \overline{C_{1}} + \frac{3}{16} \, C_{1} \overline{E_{1}} - \frac{1}{2} \, {\left(4 \, A_{0} \alpha_{1} \overline{C_{1}} - A_{1} \overline{B_{1}} - 2 \, A_{2} \overline{C_{1}}\right)} \overline{\alpha_{1}} + \frac{1}{2} \, {\left(A_{0} \overline{B_{1}} + 2 \, A_{1} \overline{C_{1}}\right)} \overline{\alpha_{4}}\\	
	    			\lambda_{14}&=8 \, {\left| A_{1} \right|}^{2} \overline{A_{1}} \overline{B_{1}} + 4 \, {\left| A_{1} \right|}^{2} \overline{A_{0}} \overline{B_{3}} + C_{1} \zeta_{0} \overline{A_{0}} + 4 \, \beta \overline{A_{0}} \overline{C_{1}} + 4 \, \overline{A_{0}} \overline{B_{1}} \overline{\alpha_{4}} + 2 \, {\left(2 \, \overline{A_{2}} \overline{C_{1}} + \overline{A_{1}} \overline{C_{2}}\right)} \alpha_{1} + 2 \, {\left(2 \, \overline{A_{1}} \overline{C_{1}} + \overline{A_{0}} \overline{C_{2}}\right)} \alpha_{4}\\
	    			&  - \frac{3}{8} \, C_{1} \gamma_{1} - 4 \, \overline{A_{2}} \overline{B_{2}} - 2 \, \overline{A_{1}} \overline{B_{5}} - \frac{1}{4} \, {\left(64 \, {\left| A_{1} \right|}^{4} \overline{A_{0}} + B_{2}\right)} \overline{C_{1}} - \frac{1}{8} \, B_{1} \overline{C_{2}} - 4 \, {\left(2 \, \alpha_{1} \overline{A_{0}} \overline{C_{1}} - \overline{A_{0}} \overline{B_{2}}\right)} \overline{\alpha_{1}}\\
	    			\lambda_{15}&=-2 \, {\left| A_{1} \right|}^{2} \overline{A_{0}} \overline{B_{2}} + C_{1} \alpha_{6} \overline{A_{0}} - \alpha_{4} \overline{A_{0}} \overline{B_{1}} - \alpha_{5} \overline{A_{0}} \overline{C_{1}} - \alpha_{3} \overline{A_{1}} \overline{C_{1}} + {\left(4 \, {\left| A_{1} \right|}^{2} \overline{A_{0}} \overline{C_{1}} - \overline{A_{1}} \overline{B_{1}}\right)} \alpha_{1} - \frac{7}{120} \, B_{1} \overline{B_{1}} - \frac{2}{15} \, C_{1} \overline{B_{3}}\\
	    			&  + \overline{A_{1}} \overline{B_{4}} - \frac{1}{120} \, {\left(7 \, C_{1} {\left| A_{1} \right|}^{2} + 15 \, B_{3}\right)} \overline{C_{1}} - \frac{7}{48} \, C_{2} \overline{C_{2}} + {\left(2 \, B_{1} \overline{A_{0}} + C_{1} \overline{A_{1}}\right)} \overline{\alpha_{2}}\\
	    			\lambda_{16}&=\frac{4}{3} \, A_{1} {\left| A_{1} \right|}^{2} \overline{B_{1}} + \frac{2}{3} \, A_{0} {\left| A_{1} \right|}^{2} \overline{B_{2}} + 2 \, A_{2} {\left| A_{1} \right|}^{2} \overline{C_{1}} + \frac{5}{3} \, A_{0} C_{1} \alpha_{6} + A_{0} \alpha_{5} \overline{C_{1}} + A_{0} \alpha_{3} \overline{C_{2}} + \frac{5}{3} \, A_{0} B_{1} \overline{\alpha_{2}}\\
	    			&  - \frac{1}{3} \, {\left(12 \, A_{0} {\left| A_{1} \right|}^{2} \overline{C_{1}} - 2 \, A_{0} \overline{B_{3}} - 3 \, A_{1} \overline{C_{2}}\right)} \alpha_{1} + \frac{1}{3} \, {\left(2 \, A_{0} \overline{B_{1}} + 3 \, A_{1} \overline{C_{1}}\right)} \alpha_{4} - \frac{2}{3} \, A_{2} \overline{B_{3}} - \frac{1}{3} \, A_{1} \overline{B_{5}} - A_{3} \overline{C_{2}}\\
	    			\lambda_{17}&=-A_{0} \alpha_{1}^{2} \overline{C_{1}} + \frac{3}{2} \, A_{0} C_{1} \alpha_{7} + \frac{1}{4} \, {\left(3 \, A_{1} \overline{B_{1}} + 2 \, A_{0} \overline{B_{2}} + 4 \, A_{2} \overline{C_{1}}\right)} \alpha_{1} + \frac{1}{4} \, {\left(3 \, A_{0} \overline{B_{1}} + 4 \, A_{1} \overline{C_{1}}\right)} \alpha_{3} - \frac{3}{4} \, A_{3} \overline{B_{1}}\\
	    			&  - \frac{1}{2} \, A_{2} \overline{B_{2}} - \frac{1}{4} \, A_{1} \overline{B_{4}} - A_{4} \overline{C_{1}} + \frac{1}{8} \, \overline{C_{1}} \overline{E_{1}} + \frac{1}{4} \, {\left(6 \, A_{1} C_{1} + 5 \, A_{0} C_{2}\right)} \overline{\alpha_{2}}\\
	    			\lambda_{18}&=C_{1} \alpha_{7} \overline{A_{0}} + \frac{1}{48} \, C_{1} \alpha_{1} \overline{C_{1}} - \frac{1}{4} \, \zeta_{0} \overline{A_{0}} \overline{C_{1}} - \overline{A_{0}} \overline{E_{1}} \overline{\alpha_{1}} + C_{2} \overline{A_{0}} \overline{\alpha_{2}} - \frac{9}{80} \, C_{2} \overline{B_{1}} - \frac{5}{48} \, C_{1} \overline{B_{2}} - \frac{1}{4} \, C_{3} \overline{C_{1}} + \frac{1}{64} \, \gamma_{1} \overline{C_{1}} - \frac{1}{16} \, \overline{C_{2}}^{2}\\
	    			&  + \overline{A_{2}} \overline{E_{1}} + \frac{1}{2} \, \overline{A_{1}} \overline{E_{3}}
	    		\end{align*}
	    		\normalsize
	    		We see that the new powers are as $\s{\vec{C}_1}{\vec{E}_1}=0$ (recall that $\vec{E}_1\in\mathrm{Span}(\bar{\vec{A}_0})$)
	    		\begin{align*}
	    			\begin{dmatrix}
	    			-4 & 6 & 0\\
	    			-3 & 5 & 0\\
	    			-2 & 4 & 0\\
	    			-1 & 3 & 0
	    			\end{dmatrix}
	    			\begin{dmatrix}
	    			2 & 0 & 0\\
	    			2 & 0 & 1\\
	    			1 & 1 & 0\\
	    			1 & 1 & 1
	    			\end{dmatrix}
	    		\end{align*}
	    		and their conjugates, so there exists $\vec{C}_5,\vec{C}_6,\vec{B}_6,\vec{B}_7,\vec{E}_4,\vec{E}_5,\vec{\gamma}_2\in \mathbb{C}^n$ and $\vec{\gamma}_4\in \R^n$ such that
	    		\begin{align}\label{4H3}
	    			\H&=\Re\bigg(\frac{\vec{C}_1}{z^2}+\frac{\vec{C}_2}{z}+\vec{
	    				C}_4z+\vec{C}_5z^2+\left(\vec{B}_1\z+\vec{B}_2\z^2+\vec{B}_4\z^3+\vec{B}_6\z^4\right)\frac{1}{z^2}+\left(\vec{B}_3\z+\vec{B}_5\z^2+\vec{B}_7\z^3\right)\frac{1}{z}\nonumber\\
    				&+\left(\vec{E}_1\z^4+\vec{E}_2\z^{5}+\vec{E}_4\z^6\right)\frac{1}{z^4}+\left(\vec{E}_3\z^4+\vec{E}_5\z^5\right)\frac{1}{z^3}\bigg)\nonumber\\   &+\vec{C}_3+\vec{C}_6|z|^2+\vec{\gamma}_1\log|z|+\Re\Big(\vec{\gamma}_2\, z+\vec{\gamma}_3z^2\Big)\log|z|+\vec{\gamma}_4|z|^2\log|z|+O(|z|^{3-\epsilon}).
	    		\end{align}
	    		or
	    		\begin{align*}
	    			\H=\begin{dmatrix}
	    			\frac{1}{2} & C_{1} & -2 & 0 & 0 \\
	    			\frac{1}{2} & C_{2} & -1 & 0 & 0 \\
	    			1 & C_{3} & 0 & 0 & 0 \\
	    			\frac{1}{2} & C_{4} & 1 & 0 & 0 \\
	    			\frac{1}{2} & C_{5} & 2 & 0 & 0 \\
	    			1 & C_{6} & 1 & 1 & 0 \\
	    			\frac{1}{2} & B_{1} & -2 & 1 & 0 \\
	    			\frac{1}{2} & B_{2} & -2 & 2 & 0 \\
	    			\frac{1}{2} & B_{4} & -2 & 3 & 0 \\
	    			\frac{1}{2} & B_{6} & -2 & 4 & 0 \\
	    			\frac{1}{2} & B_{3} & -1 & 1 & 0 
	    			\end{dmatrix}
	    			\begin{dmatrix}
	    			\frac{1}{2} & B_{5} & -1 & 2 & 0 \\
	    			\frac{1}{2} & B_{7} & -1 & 3 & 0 \\
	    			\frac{1}{2} & E_{1} & -4 & 4 & 0 \\
	    			\frac{1}{2} & E_{2} & -4 & 5 & 0 \\
	    			\frac{1}{2} & E_{4} & -4 & 6 & 0 \\
	    			\frac{1}{2} & E_{3} & -3 & 4 & 0 \\
	    			\frac{1}{2} & E_{5} & -3 & 5 & 0 \\
	    			1 & \gamma_{1} & 0 & 0 & 1 \\
	    			\frac{1}{2} & \gamma_{2} & 1 & 0 & 1 \\
	    			\frac{1}{2} & \gamma_{3} & 2 & 0 & 1 \\
	    			1 & \gamma_{4} & 1 & 1 & 1 
	    			\end{dmatrix}
	    			\begin{dmatrix}
	    			\frac{1}{2} & \overline{C_{1}} & 0 & -2 & 0 \\
	    			\frac{1}{2} & \overline{C_{2}} & 0 & -1 & 0 \\
	    			\frac{1}{2} & \overline{C_{4}} & 0 & 1 & 0 \\
	    			\frac{1}{2} & \overline{C_{5}} & 0 & 2 & 0 \\
	    			\frac{1}{2} & \overline{B_{1}} & 1 & -2 & 0 \\
	    			\frac{1}{2} & \overline{B_{2}} & 2 & -2 & 0 \\
	    			\frac{1}{2} & \overline{B_{4}} & 3 & -2 & 0 \\
	    			\frac{1}{2} & \overline{B_{6}} & 4 & -2 & 0 \\
	    			\frac{1}{2} & \overline{B_{3}} & 1 & -1 & 0 
	    			\end{dmatrix}
	    			\begin{dmatrix}
	    			\frac{1}{2} & \overline{B_{5}} & 2 & -1 & 0 \\
	    			\frac{1}{2} & \overline{B_{7}} & 3 & -1 & 0 \\
	    			\frac{1}{2} & \overline{E_{1}} & 4 & -4 & 0 \\
	    			\frac{1}{2} & \overline{E_{2}} & 5 & -4 & 0 \\
	    			\frac{1}{2} & \overline{E_{4}} & 6 & -4 & 0 \\
	    			\frac{1}{2} & \overline{E_{3}} & 4 & -3 & 0 \\
	    			\frac{1}{2} & \overline{E_{5}} & 5 & -3 & 0 \\
	    			\frac{1}{2} & \overline{\gamma_{2}} & 0 & 1 & 1 \\
	    			\frac{1}{2} & \overline{\gamma_{3}} & 0 & 2 & 1
	    			\end{dmatrix}
	    		\end{align*}
	    		We can easily check that both developments coincide.
	    		Then, we have
	    		\small
	    		\begin{align*}
	    			\p{z}\phi=\begin{dmatrix}
	    			1 & A_{0} & 3 & 0 & 0 \\
	    			1 & A_{1} & 4 & 0 & 0 \\
	    			1 & A_{2} & 5 & 0 & 0 \\
	    			1 & A_{3} & 6 & 0 & 0 \\
	    			1 & A_{4} & 7 & 0 & 0 \\
	    			1 & A_{5} & 8 & 0 & 0 \\
	    			1 & A_{6} & 9 & 0 & 0 \\
	    			\frac{1}{16} & C_{1} & 1 & 4 & 0 \\
	    			\frac{1}{16} & C_{2} & 2 & 4 & 0 \\
	    			\frac{1}{8} & C_{3} & 3 & 4 & 0 \\
	    			\frac{1}{16} & C_{4} & 4 & 4 & 0 \\
	    			\frac{1}{16} & C_{5} & 5 & 4 & 0 \\
	    			\frac{1}{10} & C_{6} & 4 & 5 & 0 \\
	    			\frac{1}{20} & B_{1} & 1 & 5 & 0 \\
	    			\frac{1}{24} & B_{2} & 1 & 6 & 0 \\
	    			\frac{1}{28} & B_{4} & 1 & 7 & 0 \\
	    			\frac{1}{32} & B_{6} & 1 & 8 & 0 \\
	    			\frac{1}{20} & B_{3} & 2 & 5 & 0 \\
	    			\frac{1}{24} & B_{5} & 2 & 6 & 0 \\
	    			\frac{1}{28} & B_{7} & 2 & 7 & 0 \\
	    			\frac{1}{32} & E_{1} & -1 & 8 & 0 \\
	    			\frac{1}{36} & E_{2} & -1 & 9 & 0 \\
	    			\frac{1}{40} & E_{4} & -1 & 10 & 0 \\
	    			\frac{1}{32} & E_{3} & 0 & 8 & 0 \\
	    			\frac{1}{36} & E_{5} & 0 & 9 & 0 \\
	    			\frac{1}{8} & \gamma_{1} & 3 & 4 & 1 \\
	    			-\frac{1}{64} & \gamma_{1} & 3 & 4 & 0 \\
	    			\frac{1}{16} & \gamma_{2} & 4 & 4 & 1 
	    				    						\end{dmatrix}
	    			\begin{dmatrix}
	    			-\frac{1}{128} & \gamma_{2} & 4 & 4 & 0 \\
	    			\frac{1}{16} & \gamma_{3} & 5 & 4 & 1 \\
	    			-\frac{1}{128} & \gamma_{3} & 5 & 4 & 0 \\
	    			\frac{1}{10} & \gamma_{4} & 4 & 5 & 1 \\
	    			-\frac{1}{100} & \gamma_{4} & 4 & 5 & 0 \\
	    			\frac{1}{8} & \overline{C_{1}} & 3 & 2 & 0 \\
	    			\frac{1}{12} & \overline{C_{2}} & 3 & 3 & 0 \\
	    			\frac{1}{20} & \overline{C_{4}} & 3 & 5 & 0 \\
	    			\frac{1}{24} & \overline{C_{5}} & 3 & 6 & 0 \\
	    			\frac{1}{8} & \overline{B_{1}} & 4 & 2 & 0 \\
	    			\frac{1}{8} & \overline{B_{2}} & 5 & 2 & 0 \\
	    			\frac{1}{8} & \overline{B_{4}} & 6 & 2 & 0 \\
	    			\frac{1}{8} & \overline{B_{6}} & 7 & 2 & 0 \\
	    			\frac{1}{12} & \overline{B_{3}} & 4 & 3 & 0 \\
	    			\frac{1}{12} & \overline{B_{5}} & 5 & 3 & 0 \\
	    			\frac{1}{12} & \overline{B_{7}} & 6 & 3 & 0 \\
	    			\frac{1}{2} & \overline{E_{1}} & 7 & 0 & 1 \\
	    			\frac{1}{2} & \overline{E_{2}} & 8 & 0 & 1 \\
	    			\frac{1}{2} & \overline{E_{4}} & 9 & 0 & 1 \\
	    			\frac{1}{4} & \overline{E_{3}} & 7 & 1 & 0 \\
	    			\frac{1}{4} & \overline{E_{5}} & 8 & 1 & 0 \\
	    			\frac{1}{20} & \overline{\gamma_{2}} & 3 & 5 & 1 \\
	    			-\frac{1}{200} & \overline{\gamma_{2}} & 3 & 5 & 0 \\
	    			\frac{1}{24} & \overline{\gamma_{3}} & 3 & 6 & 1 \\
	    			-\frac{1}{288} & \overline{\gamma_{3}} & 3 & 6 & 0 \\
	    			\frac{1}{10} \, {\left| A_{1} \right|}^{2} & C_{1} & 2 & 5 & 0 \\
	    			\frac{1}{10} \, {\left| A_{1} \right|}^{2} & C_{2} & 3 & 5 & 0 \\
	    			\frac{1}{5} \, {\left| A_{1} \right|}^{2} & C_{3} & 4 & 5 & 0 
	    				    			\end{dmatrix}
	    			\begin{dmatrix}
	    			\frac{1}{12} \, {\left| A_{1} \right|}^{2} & B_{1} & 2 & 6 & 0 \\
	    			\frac{1}{14} \, {\left| A_{1} \right|}^{2} & B_{2} & 2 & 7 & 0 \\
	    			\frac{1}{12} \, {\left| A_{1} \right|}^{2} & B_{3} & 3 & 6 & 0 \\
	    			\frac{1}{18} \, {\left| A_{1} \right|}^{2} & E_{1} & 0 & 9 & 0 \\
	    			\frac{1}{5} \, {\left| A_{1} \right|}^{2} & \gamma_{1} & 4 & 5 & 1 \\
	    			-\frac{1}{50} \, {\left| A_{1} \right|}^{2} & \gamma_{1} & 4 & 5 & 0 \\
	    			\frac{1}{6} \, {\left| A_{1} \right|}^{2} & \overline{C_{1}} & 4 & 3 & 0 \\
	    			\frac{1}{8} \, {\left| A_{1} \right|}^{2} & \overline{C_{2}} & 4 & 4 & 0 \\
	    			\frac{1}{6} \, {\left| A_{1} \right|}^{2} & \overline{B_{1}} & 5 & 3 & 0 \\
	    			\frac{1}{6} \, {\left| A_{1} \right|}^{2} & \overline{B_{2}} & 6 & 3 & 0 \\
	    			\frac{1}{8} \, {\left| A_{1} \right|}^{2} & \overline{B_{3}} & 5 & 4 & 0 \\
	    			\frac{1}{2} \, {\left| A_{1} \right|}^{2} & \overline{E_{1}} & 8 & 1 & 0 \\
	    			\frac{1}{24} \, \beta & C_{1} & 3 & 6 & 0 \\
	    			\frac{1}{16} \, \beta & \overline{C_{1}} & 5 & 4 & 0 \\
	    			\frac{1}{16} \, \alpha_{1} & C_{1} & 3 & 4 & 0 \\
	    			\frac{1}{16} \, \alpha_{1} & C_{2} & 4 & 4 & 0 \\
	    			\frac{1}{8} \, \alpha_{1} & C_{3} & 5 & 4 & 0 \\
	    			\frac{1}{20} \, \alpha_{1} & B_{1} & 3 & 5 & 0 \\
	    			\frac{1}{24} \, \alpha_{1} & B_{2} & 3 & 6 & 0 \\
	    			\frac{1}{20} \, \alpha_{1} & B_{3} & 4 & 5 & 0 \\
	    			\frac{1}{32} \, \alpha_{1} & E_{1} & 1 & 8 & 0 \\
	    			\frac{1}{8} \, \alpha_{1} & \gamma_{1} & 5 & 4 & 1 \\
	    			-\frac{1}{64} \, \alpha_{1} & \gamma_{1} & 5 & 4 & 0 \\
	    			\frac{1}{8} \, \alpha_{1} & \overline{C_{1}} & 5 & 2 & 0 \\
	    			\frac{1}{12} \, \alpha_{1} & \overline{C_{2}} & 5 & 3 & 0 \\
	    			\frac{1}{8} \, \alpha_{1} & \overline{B_{1}} & 6 & 2 & 0 \\
	    			\frac{1}{8} \, \alpha_{1} & \overline{B_{2}} & 7 & 2 & 0 \\
	    			\frac{1}{12} \, \alpha_{1} & \overline{B_{3}} & 6 & 3 & 0
	    			\end{dmatrix}
	    			\end{align*}
	    			\begin{align}\label{4devphifinal}
	    			\begin{dmatrix}
	    			\frac{1}{2} \, \alpha_{1} & \overline{E_{1}} & 9 & 0 & 1 \\
	    			\frac{1}{36} \, \alpha_{2} & C_{1} & -1 & 9 & 0 \\
	    			\frac{1}{36} \, \alpha_{2} & C_{2} & 0 & 9 & 0 \\
	    			\frac{1}{40} \, \alpha_{2} & B_{1} & -1 & 10 & 0 \\
	    			\frac{1}{28} \, \alpha_{2} & \overline{C_{1}} & 1 & 7 & 0 \\
	    			\frac{1}{32} \, \alpha_{2} & \overline{C_{2}} & 1 & 8 & 0 \\
	    			\frac{1}{28} \, \alpha_{2} & \overline{B_{1}} & 2 & 7 & 0 \\
	    			\frac{1}{16} \, \alpha_{3} & C_{1} & 4 & 4 & 0 \\
	    			\frac{1}{16} \, \alpha_{3} & C_{2} & 5 & 4 & 0 \\
	    			\frac{1}{20} \, \alpha_{3} & B_{1} & 4 & 5 & 0 \\
	    			\frac{1}{8} \, \alpha_{3} & \overline{C_{1}} & 6 & 2 & 0 \\
	    			\frac{1}{12} \, \alpha_{3} & \overline{C_{2}} & 6 & 3 & 0 \\
	    			\frac{1}{8} \, \alpha_{3} & \overline{B_{1}} & 7 & 2 & 0 \\
	    			\frac{1}{20} \, \alpha_{4} & C_{1} & 3 & 5 & 0 \\
	    			\frac{1}{20} \, \alpha_{4} & C_{2} & 4 & 5 & 0 \\
	    			\frac{1}{24} \, \alpha_{4} & B_{1} & 3 & 6 & 0 \\
	    			\frac{1}{12} \, \alpha_{4} & \overline{C_{1}} & 5 & 3 & 0 \\
	    			\frac{1}{16} \, \alpha_{4} & \overline{C_{2}} & 5 & 4 & 0 \\
	    			\frac{1}{12} \, \alpha_{4} & \overline{B_{1}} & 6 & 3 & 0 \\
	    			\frac{1}{20} \, \alpha_{5} & C_{1} & 4 & 5 & 0 \\
	    			\frac{1}{12} \, \alpha_{5} & \overline{C_{1}} & 6 & 3 & 0 \\
	    			\frac{1}{12} \, \alpha_{6} & C_{1} & 6 & 3 & 0 \\
	    			\frac{1}{4} \, \alpha_{6} & \overline{C_{1}} & 8 & 1 & 0 \\
	    			\frac{1}{8} \, \alpha_{7} & C_{1} & 7 & 2 & 0 
	    				    			\end{dmatrix}
	    			\begin{dmatrix}
	    			\frac{1}{2} \, \alpha_{7} & \overline{C_{1}} & 9 & 0 & 1 \\
	    			\frac{1}{16} \, \zeta_{0} & C_{1} & 5 & 4 & 1 \\
	    			-\frac{1}{128} \, \zeta_{0} & C_{1} & 5 & 4 & 0 \\
	    			\frac{1}{8} \, \zeta_{0} & \overline{C_{1}} & 7 & 2 & 1 \\
	    			-\frac{1}{32} \, \zeta_{0} & \overline{C_{1}} & 7 & 2 & 0 \\
	    			\frac{1}{24} \, \overline{\alpha_{1}} & C_{1} & 1 & 6 & 0 \\
	    			\frac{1}{24} \, \overline{\alpha_{1}} & C_{2} & 2 & 6 & 0 \\
	    			\frac{1}{12} \, \overline{\alpha_{1}} & C_{3} & 3 & 6 & 0 \\
	    			\frac{1}{28} \, \overline{\alpha_{1}} & B_{1} & 1 & 7 & 0 \\
	    			\frac{1}{32} \, \overline{\alpha_{1}} & B_{2} & 1 & 8 & 0 \\
	    			\frac{1}{28} \, \overline{\alpha_{1}} & B_{3} & 2 & 7 & 0 \\
	    			\frac{1}{40} \, \overline{\alpha_{1}} & E_{1} & -1 & 10 & 0 \\
	    			\frac{1}{12} \, \overline{\alpha_{1}} & \gamma_{1} & 3 & 6 & 1 \\
	    			-\frac{1}{144} \, \overline{\alpha_{1}} & \gamma_{1} & 3 & 6 & 0 \\
	    			\frac{1}{16} \, \overline{\alpha_{1}} & \overline{C_{1}} & 3 & 4 & 0 \\
	    			\frac{1}{20} \, \overline{\alpha_{1}} & \overline{C_{2}} & 3 & 5 & 0 \\
	    			\frac{1}{16} \, \overline{\alpha_{1}} & \overline{B_{1}} & 4 & 4 & 0 \\
	    			\frac{1}{16} \, \overline{\alpha_{1}} & \overline{B_{2}} & 5 & 4 & 0 \\
	    			\frac{1}{20} \, \overline{\alpha_{1}} & \overline{B_{3}} & 4 & 5 & 0 \\
	    			\frac{1}{8} \, \overline{\alpha_{1}} & \overline{E_{1}} & 7 & 2 & 0 \\
	    			\frac{1}{8} \, \overline{\alpha_{2}} & C_{1} & 6 & 2 & 0 \\
	    			\frac{1}{8} \, \overline{\alpha_{2}} & C_{2} & 7 & 2 & 0 \\
	    			\frac{1}{12} \, \overline{\alpha_{2}} & B_{1} & 6 & 3 & 0 \\
	    			\frac{1}{2} \, \overline{\alpha_{2}} & \overline{C_{1}} & 8 & 0 & 1 
	    				    						\end{dmatrix}
	    			\begin{dmatrix}
	    			\frac{1}{4} \, \overline{\alpha_{2}} & \overline{C_{2}} & 8 & 1 & 0 \\
	    			\frac{1}{2} \, \overline{\alpha_{2}} & \overline{B_{1}} & 9 & 0 & 1 \\
	    			\frac{1}{28} \, \overline{\alpha_{3}} & C_{1} & 1 & 7 & 0 \\
	    			\frac{1}{28} \, \overline{\alpha_{3}} & C_{2} & 2 & 7 & 0 \\
	    			\frac{1}{32} \, \overline{\alpha_{3}} & B_{1} & 1 & 8 & 0 \\
	    			\frac{1}{20} \, \overline{\alpha_{3}} & \overline{C_{1}} & 3 & 5 & 0 \\
	    			\frac{1}{24} \, \overline{\alpha_{3}} & \overline{C_{2}} & 3 & 6 & 0 \\
	    			\frac{1}{20} \, \overline{\alpha_{3}} & \overline{B_{1}} & 4 & 5 & 0 \\
	    			\frac{1}{24} \, \overline{\alpha_{4}} & C_{1} & 2 & 6 & 0 \\
	    			\frac{1}{24} \, \overline{\alpha_{4}} & C_{2} & 3 & 6 & 0 \\
	    			\frac{1}{28} \, \overline{\alpha_{4}} & B_{1} & 2 & 7 & 0 \\
	    			\frac{1}{16} \, \overline{\alpha_{4}} & \overline{C_{1}} & 4 & 4 & 0 \\
	    			\frac{1}{20} \, \overline{\alpha_{4}} & \overline{C_{2}} & 4 & 5 & 0 \\
	    			\frac{1}{16} \, \overline{\alpha_{4}} & \overline{B_{1}} & 5 & 4 & 0 \\
	    			\frac{1}{28} \, \overline{\alpha_{5}} & C_{1} & 2 & 7 & 0 \\
	    			\frac{1}{20} \, \overline{\alpha_{5}} & \overline{C_{1}} & 4 & 5 & 0 \\
	    			\frac{1}{36} \, \overline{\alpha_{6}} & C_{1} & 0 & 9 & 0 \\
	    			\frac{1}{28} \, \overline{\alpha_{6}} & \overline{C_{1}} & 2 & 7 & 0 \\
	    			\frac{1}{40} \, \overline{\alpha_{7}} & C_{1} & -1 & 10 & 0 \\
	    			\frac{1}{32} \, \overline{\alpha_{7}} & \overline{C_{1}} & 1 & 8 & 0 \\
	    			\frac{1}{32} \, \overline{\zeta_{0}} & C_{1} & 1 & 8 & 1 \\
	    			-\frac{1}{512} \, \overline{\zeta_{0}} & C_{1} & 1 & 8 & 0 \\
	    			\frac{1}{24} \, \overline{\zeta_{0}} & \overline{C_{1}} & 3 & 6 & 1 \\
	    			-\frac{1}{288} \, \overline{\zeta_{0}} & \overline{C_{1}} & 3 & 6 & 0
	    			\end{dmatrix}
	    		\end{align}
	    		
	    		By conformality of $\phi$, we have
	    		\small
	    		\begin{align*}
	    		&0=\s{\p{z}\phi}{\p{z}\phi}=\\
	    		&\begin{dmatrix}
	    		A_{0}^{2} & 6 & 0 & 0 \\
	    		2 \, A_{0} A_{1} & 7 & 0 & 0 \\
	    		A_{1}^{2} + 2 \, A_{0} A_{2} & 8 & 0 & 0 \\
	    		2 \, A_{1} A_{2} + 2 \, A_{0} A_{3} & 9 & 0 & 0 \\
	    		A_{2}^{2} + 2 \, A_{1} A_{3} + 2 \, A_{0} A_{4} & 10 & 0 & 0 \\
	    		2 \, A_{2} A_{3} + 2 \, A_{1} A_{4} + 2 \, A_{0} A_{5} & 11 & 0 & 0 \\
	    		A_{3}^{2} + 2 \, A_{2} A_{4} + 2 \, A_{1} A_{5} + 2 \, A_{0} A_{6} & 12 & 0 & 0 \\
	    		\frac{1}{8} \, A_{0} C_{1} & 4 & 4 & 0 \\
	    		\frac{1}{8} \, A_{1} C_{1} + \frac{1}{8} \, A_{0} C_{2} & 5 & 4 & 0 \\
	    		\frac{1}{8} \, A_{0} C_{1} \alpha_{1} + \frac{1}{8} \, A_{0} \overline{C_{1}} \overline{\alpha_{1}} + \frac{1}{8} \, A_{2} C_{1} + \frac{1}{8} \, A_{1} C_{2} + \frac{1}{4} \, A_{0} C_{3} - \frac{1}{32} \, A_{0} \gamma_{1} + \frac{1}{64} \, \overline{C_{1}}^{2} & 6 & 4 & 0 \\
	    		\mu_1 & 7 & 4 & 0 \\
	    		\mu_2 & 8 & 4 & 0 \\
	    		\mu_3 & 7 & 5 & 0 \\
	    		\frac{1}{10} \, A_{0} B_{1} & 4 & 5 & 0 \\
	    		\frac{1}{12} \, A_{0} C_{1} \overline{\alpha_{1}} + \frac{1}{12} \, A_{0} B_{2} + \frac{1}{64} \, C_{1} \overline{C_{1}} & 4 & 6 & 0 \\
	    		\frac{1}{14} \, A_{0} \alpha_{2} \overline{C_{1}} + \frac{1}{14} \, A_{0} B_{1} \overline{\alpha_{1}} + \frac{1}{14} \, A_{0} C_{1} \overline{\alpha_{3}} + \frac{1}{14} \, A_{0} B_{4} + \frac{1}{80} \, B_{1} \overline{C_{1}} + \frac{1}{96} \, C_{1} \overline{C_{2}} & 4 & 7 & 0 \\
	    		\mu_4 & 4 & 8 & 0 \\
	    		\frac{1}{5} \, A_{0} C_{1} {\left| A_{1} \right|}^{2} + \frac{1}{10} \, A_{1} B_{1} + \frac{1}{10} \, A_{0} B_{3} & 5 & 5 & 0 \\
	    		\frac{1}{6} \, A_{0} B_{1} {\left| A_{1} \right|}^{2} + \frac{1}{12} \, A_{0} C_{1} \overline{\alpha_{4}} + \frac{1}{12} \, A_{1} B_{2} + \frac{1}{12} \, A_{0} B_{5} + \frac{1}{64} \, C_{1} \overline{B_{1}} + \frac{1}{64} \, C_{2} \overline{C_{1}} + \frac{1}{12} \, {\left(A_{1} C_{1} + A_{0} C_{2}\right)} \overline{\alpha_{1}} & 5 & 6 & 0 \\
	    		\mu_5 & 5 & 7 & 0 \\
	    		\frac{1}{256} \, C_{1}^{2} + \frac{1}{16} \, A_{0} E_{1} & 2 & 8 & 0 \\
	    		\frac{1}{18} \, A_{0} C_{1} \alpha_{2} + \frac{1}{160} \, B_{1} C_{1} + \frac{1}{18} \, A_{0} E_{2} & 2 & 9 & 0 \\
	    		\frac{1}{20} \, A_{0} B_{1} \alpha_{2} + \frac{1}{20} \, A_{0} C_{1} \overline{\alpha_{7}} + \frac{1}{400} \, B_{1}^{2} + \frac{1}{192} \, B_{2} C_{1} + \frac{1}{20} \, A_{0} E_{4} + \frac{1}{128} \, E_{1} \overline{C_{1}} + \frac{1}{960} \, {\left(5 \, C_{1}^{2} + 48 \, A_{0} E_{1}\right)} \overline{\alpha_{1}} & 2 & 10 & 0 \\
	    		\frac{1}{128} \, C_{1} C_{2} + \frac{1}{16} \, A_{1} E_{1} + \frac{1}{16} \, A_{0} E_{3} & 3 & 8 & 0 \\
	    		\mu_6 & 3 & 9 & 0 \\
	    		\frac{1}{4} \, A_{0} \gamma_{1} & 6 & 4 & 1 \\
	    		\frac{1}{4} \, A_{1} \gamma_{1} + \frac{1}{8} \, A_{0} \gamma_{2} & 7 & 4 & 1 \\
	    		\frac{1}{8} \, A_{0} C_{1} \zeta_{0} + \frac{1}{4} \, {\left(A_{0} \alpha_{1} + A_{2}\right)} \gamma_{1} + \frac{1}{8} \, A_{1} \gamma_{2} + \frac{1}{8} \, A_{0} \gamma_{3} + \frac{1}{16} \, C_{1} \overline{E_{1}} & 8 & 4 & 1 
	    		\end{dmatrix}
	    		\begin{dmatrix}
	    		\textbf{(1)}\\
	    		\textbf{(2)}\\
	    		\textbf{(3)}\\
	    		\textbf{(4)}\\
	    		\textbf{(5)}\\
	    		\textbf{(6)}\\
	    		\textbf{(7)}\\
	    		\textbf{(8)}\\
	    		\textbf{(9)}\\
	    		\textbf{(10)}\\
	    		\textbf{(11)}\\
	    		\textbf{(12)}\\
	    		\textbf{(13)}\\
	    		\textbf{(14)}\\
	    		\textbf{(15)}\\
	    		\textbf{(16)}\\
	    		\textbf{(17)}\\
	    		\textbf{(18)}\\
	    		\textbf{(19)}\\
	    		\textbf{(20)}\\
	    		\textbf{(21)}\\
	    		\textbf{(22)}\\
	    		\textbf{(23)}\\
	    		\textbf{(24)}\\
	    		\textbf{(25)}\\
	    		\textbf{(26)}\\
	    		\textbf{(27)}\\
	    		\textbf{(28)}\\
	    		\end{dmatrix}
	    		\end{align*}
	    		\begin{align}\label{4conf}
	    		\begin{dmatrix}
	    		\frac{2}{5} \, A_{0} \gamma_{1} {\left| A_{1} \right|}^{2} + \frac{1}{5} \, A_{0} \gamma_{4} + \frac{1}{10} \, A_{1} \overline{\gamma_{2}} & 7 & 5 & 1 \\
	    		\frac{1}{4} \, A_{0} \overline{C_{1}} & 6 & 2 & 0 \\
	    		\frac{1}{6} \, A_{0} \overline{C_{2}} & 6 & 3 & 0 \\
	    		\mu_7 & 6 & 5 & 0 \\
	    		\mu_8 & 6 & 6 & 0 \\
	    		\frac{1}{4} \, A_{0} \overline{B_{1}} + \frac{1}{4} \, A_{1} \overline{C_{1}} & 7 & 2 & 0 \\
	    		\frac{1}{4} \, A_{0} \alpha_{1} \overline{C_{1}} + \frac{1}{4} \, A_{1} \overline{B_{1}} + \frac{1}{4} \, A_{0} \overline{B_{2}} + \frac{1}{4} \, A_{2} \overline{C_{1}} & 8 & 2 & 0 \\
	    		\frac{1}{4} \, A_{0} \alpha_{3} \overline{C_{1}} + \frac{1}{4} \, A_{0} C_{1} \overline{\alpha_{2}} + \frac{1}{4} \, {\left(A_{0} \overline{B_{1}} + A_{1} \overline{C_{1}}\right)} \alpha_{1} + \frac{1}{4} \, A_{2} \overline{B_{1}} + \frac{1}{4} \, A_{1} \overline{B_{2}} + \frac{1}{4} \, A_{0} \overline{B_{4}} + \frac{1}{4} \, A_{3} \overline{C_{1}} & 9 & 2 & 0 \\
	    		\mu_9 & 10 & 2 & 0 \\
	    		\frac{1}{3} \, A_{0} {\left| A_{1} \right|}^{2} \overline{C_{1}} + \frac{1}{6} \, A_{0} \overline{B_{3}} + \frac{1}{6} \, A_{1} \overline{C_{2}} & 7 & 3 & 0 \\
	    		\frac{1}{3} \, A_{0} {\left| A_{1} \right|}^{2} \overline{B_{1}} + \frac{1}{3} \, A_{1} {\left| A_{1} \right|}^{2} \overline{C_{1}} + \frac{1}{6} \, A_{0} \alpha_{4} \overline{C_{1}} + \frac{1}{6} \, A_{0} \alpha_{1} \overline{C_{2}} + \frac{1}{6} \, A_{1} \overline{B_{3}} + \frac{1}{6} \, A_{0} \overline{B_{5}} + \frac{1}{6} \, A_{2} \overline{C_{2}} & 8 & 3 & 0 \\
	    		\mu_{10} & 9 & 3 & 0 \\
	    		A_{0} \overline{E_{1}} & 10 & 0 & 1 \\
	    		A_{0} \overline{C_{1}} \overline{\alpha_{2}} + A_{1} \overline{E_{1}} + A_{0} \overline{E_{2}} & 11 & 0 & 1 \\
	    		A_{0} \alpha_{7} \overline{C_{1}} + A_{0} \alpha_{1} \overline{E_{1}} + A_{2} \overline{E_{1}} + A_{1} \overline{E_{2}} + A_{0} \overline{E_{4}} + {\left(A_{0} \overline{B_{1}} + A_{1} \overline{C_{1}}\right)} \overline{\alpha_{2}} & 12 & 0 & 1 \\
	    		\frac{1}{2} \, A_{0} \overline{E_{3}} & 10 & 1 & 0 \\
	    		A_{0} {\left| A_{1} \right|}^{2} \overline{E_{1}} + \frac{1}{2} \, A_{0} \alpha_{6} \overline{C_{1}} + \frac{1}{2} \, A_{0} \overline{C_{2}} \overline{\alpha_{2}} + \frac{1}{2} \, A_{1} \overline{E_{3}} + \frac{1}{2} \, A_{0} \overline{E_{5}} & 11 & 1 & 0 \\
	    		\frac{1}{10} \, A_{0} \overline{\gamma_{2}} & 6 & 5 & 1 \\
	    		\frac{1}{12} \, A_{0} \overline{C_{1}} \overline{\zeta_{0}} + \frac{1}{96} \, {\left(16 \, A_{0} \overline{\alpha_{1}} + 3 \, \overline{C_{1}}\right)} \gamma_{1} + \frac{1}{12} \, A_{0} \overline{\gamma_{3}} & 6 & 6 & 1 \\
	    		\frac{1}{4} \, A_{0} \zeta_{0} \overline{C_{1}} + \frac{1}{8} \, \overline{C_{1}} \overline{E_{1}} & 10 & 2 & 1 \\
	    		\frac{1}{16} \, A_{0} C_{1} \overline{\zeta_{0}} + \frac{1}{64} \, C_{1} \gamma_{1} & 4 & 8 & 1 \\
	    		\frac{1}{256} \, C_{1} E_{1} & 0 & 12 & 0
	    		\end{dmatrix}
	    		\begin{dmatrix}
	    		\textbf{(29)}\\
	    		\textbf{(30)}\\
	    		\textbf{(31)}\\
	    		\textbf{(32)}\\
	    		\textbf{(33)}\\
	    		\textbf{(34)}\\
	    		\textbf{(35)}\\
	    		\textbf{(36)}\\
	    		\textbf{(37)}\\
	    		\textbf{(38)}\\
	    		\textbf{(39)}\\
	    		\textbf{(40)}\\
	    		\textbf{(41)}\\
	    		\textbf{(42)}\\
	    		\textbf{(43)}\\
	    		\textbf{(44)}\\
	    		\textbf{(45)}\\
	    		\textbf{(46)}\\
	    		\textbf{(47)}\\
	    		\textbf{(48)}\\
	    		\textbf{(49)}\\
	    		\textbf{(50)}
	    		\end{dmatrix}
	    		\end{align}
	    		where
	    		\begin{align*}
	    			\mu_{1}&=\frac{1}{4} \, A_{0} {\left| A_{1} \right|}^{2} \overline{C_{2}} + \frac{1}{8} \, A_{0} C_{1} \alpha_{3} + \frac{1}{8} \, A_{0} \overline{C_{1}} \overline{\alpha_{4}} + \frac{1}{8} \, A_{3} C_{1} + \frac{1}{8} \, A_{2} C_{2} + \frac{1}{4} \, A_{1} C_{3} + \frac{1}{8} \, A_{0} C_{4} + \frac{1}{8} \, {\left(A_{1} C_{1} + A_{0} C_{2}\right)} \alpha_{1} - \frac{1}{32} \, A_{1} \gamma_{1}\\
	    			& - \frac{1}{64} \, A_{0} \gamma_{2} + \frac{1}{32} \, \overline{B_{1}} \overline{C_{1}} + \frac{1}{8} \, {\left(A_{0} \overline{B_{1}} + A_{1} \overline{C_{1}}\right)} \overline{\alpha_{1}}\\
	    			\mu_{2}&=\frac{1}{4} \, A_{0} {\left| A_{1} \right|}^{2} \overline{B_{3}} + \frac{1}{4} \, A_{1} {\left| A_{1} \right|}^{2} \overline{C_{2}} - \frac{1}{64} \, A_{0} C_{1} \zeta_{0} + \frac{1}{8} \, A_{0} \beta \overline{C_{1}} + \frac{1}{8} \, A_{0} \alpha_{4} \overline{C_{2}} + \frac{1}{8} \, A_{4} C_{1} + \frac{1}{8} \, A_{3} C_{2} + \frac{1}{4} \, A_{2} C_{3} + \frac{1}{8} \, A_{1} C_{4}\\
	    			& + \frac{1}{8} \, A_{0} C_{5} + \frac{1}{32} \, {\left(4 \, A_{2} C_{1} + 4 \, A_{1} C_{2} + 8 \, A_{0} C_{3} + \overline{C_{1}}^{2}\right)} \alpha_{1} + \frac{1}{8} \, {\left(A_{1} C_{1} + A_{0} C_{2}\right)} \alpha_{3} - \frac{1}{32} \, {\left(A_{0} \alpha_{1} + A_{2}\right)} \gamma_{1}\\
	    			& - \frac{1}{64} \, A_{1} \gamma_{2} - \frac{1}{64} \, A_{0} \gamma_{3} + \frac{1}{64} \, \overline{B_{1}}^{2} + \frac{1}{32} \, \overline{B_{2}} \overline{C_{1}} + \frac{1}{8} \, {\left(A_{1} \overline{B_{1}} + A_{0} \overline{B_{2}} + A_{2} \overline{C_{1}}\right)} \overline{\alpha_{1}} + \frac{1}{8} \, {\left(A_{0} \overline{B_{1}} + A_{1} \overline{C_{1}}\right)} \overline{\alpha_{4}}\\
	    			\mu_{3}&=\frac{1}{5} \, A_{2} C_{1} {\left| A_{1} \right|}^{2} + \frac{1}{5} \, A_{1} C_{2} {\left| A_{1} \right|}^{2} + \frac{2}{5} \, A_{0} C_{3} {\left| A_{1} \right|}^{2} - \frac{1}{25} \, A_{0} \gamma_{1} {\left| A_{1} \right|}^{2} + \frac{1}{24} \, {\left| A_{1} \right|}^{2} \overline{C_{1}}^{2} + \frac{1}{10} \, A_{0} B_{1} \alpha_{3} + \frac{1}{10} \, A_{0} C_{1} \alpha_{5}\\
	    			& + \frac{1}{10} \, A_{0} \overline{C_{2}} \overline{\alpha_{4}} + \frac{1}{10} \, A_{0} \overline{C_{1}} \overline{\alpha_{5}} + \frac{1}{10} \, A_{3} B_{1} + \frac{1}{10} \, A_{2} B_{3} + \frac{1}{5} \, A_{0} C_{6} + \frac{1}{10} \, {\left(A_{1} B_{1} + A_{0} B_{3}\right)} \alpha_{1} + \frac{1}{10} \, {\left(A_{1} C_{1} + A_{0} C_{2}\right)} \alpha_{4}\\
	    			& - \frac{1}{50} \, A_{0} \gamma_{4} + \frac{1}{48} \, \overline{B_{3}} \overline{C_{1}} + \frac{1}{48} \, \overline{B_{1}} \overline{C_{2}} + \frac{1}{10} \, A_{1} \overline{C_{4}} + \frac{1}{10} \, {\left(A_{0} \overline{B_{3}} + A_{1} \overline{C_{2}}\right)} \overline{\alpha_{1}} + \frac{1}{10} \, {\left(A_{0} \overline{B_{1}} + A_{1} \overline{C_{1}}\right)} \overline{\alpha_{3}} - \frac{1}{100} \, A_{1} \overline{\gamma_{2}}\\
	    			\mu_{4}&=\frac{1}{16} \, A_{0} \alpha_{2} \overline{C_{2}} + \frac{1}{16} \, A_{0} B_{1} \overline{\alpha_{3}} + \frac{1}{16} \, A_{0} \overline{C_{1}} \overline{\alpha_{7}} - \frac{1}{256} \, A_{0} C_{1} \overline{\zeta_{0}} + \frac{1}{16} \, A_{0} B_{6} + \frac{1}{256} \, C_{2}^{2} + \frac{1}{64} \, C_{1} C_{3} + \frac{1}{16} \, A_{2} E_{1} + \frac{1}{16} \, A_{1} E_{3}\\
	    			& + \frac{1}{128} \, {\left(C_{1}^{2} + 8 \, A_{0} E_{1}\right)} \alpha_{1} - \frac{1}{512} \, C_{1} \gamma_{1} + \frac{1}{96} \, B_{2} \overline{C_{1}} + \frac{1}{120} \, B_{1} \overline{C_{2}} + \frac{1}{384} \, {\left(24 \, A_{0} B_{2} + 7 \, C_{1} \overline{C_{1}}\right)} \overline{\alpha_{1}}\\
	    			\mu_{5}&=\frac{1}{7} \, A_{0} B_{2} {\left| A_{1} \right|}^{2} + \frac{1}{14} \, A_{0} B_{1} \overline{\alpha_{4}} + \frac{1}{14} \, A_{0} C_{1} \overline{\alpha_{5}} + \frac{1}{14} \, A_{0} \overline{C_{1}} \overline{\alpha_{6}} + \frac{1}{14} \, A_{1} B_{4} + \frac{1}{14} \, A_{0} B_{7} + \frac{1}{14} \, {\left(A_{0} \overline{B_{1}} + A_{1} \overline{C_{1}}\right)} \alpha_{2} + \frac{1}{80} \, B_{1} \overline{B_{1}}\\
	    			& + \frac{1}{96} \, C_{1} \overline{B_{3}} + \frac{1}{240} \, {\left(11 \, C_{1} {\left| A_{1} \right|}^{2} + 3 \, B_{3}\right)} \overline{C_{1}} + \frac{1}{96} \, C_{2} \overline{C_{2}} + \frac{1}{14} \, {\left(A_{1} B_{1} + A_{0} B_{3}\right)} \overline{\alpha_{1}} + \frac{1}{14} \, {\left(A_{1} C_{1} + A_{0} C_{2}\right)} \overline{\alpha_{3}}\\
	    			\mu_{6}&=\frac{1}{80} \, C_{1}^{2} {\left| A_{1} \right|}^{2} + \frac{1}{9} \, A_{0} E_{1} {\left| A_{1} \right|}^{2} + \frac{1}{18} \, A_{0} C_{1} \overline{\alpha_{6}} + \frac{1}{160} \, B_{3} C_{1} + \frac{1}{160} \, B_{1} C_{2} + \frac{1}{18} \, A_{1} E_{2} + \frac{1}{18} \, A_{0} E_{5} + \frac{1}{18} \, {\left(A_{1} C_{1} + A_{0} C_{2}\right)} \alpha_{2}\\
	    			\mu_{7}&=\frac{1}{5} \, A_{1} C_{1} {\left| A_{1} \right|}^{2} + \frac{1}{5} \, A_{0} C_{2} {\left| A_{1} \right|}^{2} + \frac{1}{10} \, A_{0} B_{1} \alpha_{1} + \frac{1}{10} \, A_{0} C_{1} \alpha_{4} + \frac{1}{10} \, A_{0} \overline{C_{2}} \overline{\alpha_{1}} + \frac{1}{10} \, A_{0} \overline{C_{1}} \overline{\alpha_{3}} + \frac{1}{10} \, A_{2} B_{1} + \frac{1}{10} \, A_{1} B_{3}\\
	    			& + \frac{1}{48} \, \overline{C_{1}} \overline{C_{2}} + \frac{1}{10} \, A_{0} \overline{C_{4}} - \frac{1}{100} \, A_{0} \overline{\gamma_{2}}\\
	    			\mu_{8}&=\frac{1}{6} \, A_{1} B_{1} {\left| A_{1} \right|}^{2} + \frac{1}{6} \, A_{0} B_{3} {\left| A_{1} \right|}^{2} + \frac{1}{12} \, A_{0} B_{1} \alpha_{4} + \frac{1}{12} \, A_{0} C_{1} \beta + \frac{1}{12} \, A_{0} \overline{C_{2}} \overline{\alpha_{3}} - \frac{1}{144} \, A_{0} \overline{C_{1}} \overline{\zeta_{0}} + \frac{1}{12} \, A_{2} B_{2} + \frac{1}{12} \, A_{1} B_{5}\\
	    			& + \frac{1}{96} \, {\left(8 \, A_{0} B_{2} + 3 \, C_{1} \overline{C_{1}}\right)} \alpha_{1} - \frac{1}{2304} \, {\left(32 \, A_{0} \overline{\alpha_{1}} + 9 \, \overline{C_{1}}\right)} \gamma_{1} + \frac{1}{64} \, C_{2} \overline{B_{1}} + \frac{1}{64} \, C_{1} \overline{B_{2}} + \frac{1}{32} \, C_{3} \overline{C_{1}} + \frac{1}{144} \, \overline{C_{2}}^{2} + \frac{1}{12} \, A_{0} \overline{C_{5}} \\
	    			&+ \frac{1}{192} \, {\left(16 \, A_{2} C_{1} + 16 \, A_{1} C_{2} + 32 \, A_{0} C_{3} + 3 \, \overline{C_{1}}^{2}\right)} \overline{\alpha_{1}} + \frac{1}{12} \, {\left(A_{1} C_{1} + A_{0} C_{2}\right)} \overline{\alpha_{4}} - \frac{1}{144} \, A_{0} \overline{\gamma_{3}}\\
	    			\mu_{9}&=\frac{1}{4} \, A_{0} C_{1} \alpha_{7} - \frac{1}{16} \, A_{0} \zeta_{0} \overline{C_{1}} + \frac{1}{4} \, A_{0} \overline{E_{1}} \overline{\alpha_{1}} + \frac{1}{4} \, {\left(A_{1} \overline{B_{1}} + A_{0} \overline{B_{2}} + A_{2} \overline{C_{1}}\right)} \alpha_{1} + \frac{1}{4} \, {\left(A_{0} \overline{B_{1}} + A_{1} \overline{C_{1}}\right)} \alpha_{3} + \frac{1}{4} \, A_{3} \overline{B_{1}} + \frac{1}{4} \, A_{2} \overline{B_{2}}\\
	    			& + \frac{1}{4} \, A_{1} \overline{B_{4}} + \frac{1}{4} \, A_{0} \overline{B_{6}} + \frac{1}{4} \, A_{4} \overline{C_{1}} + \frac{1}{4} \, {\left(A_{1} C_{1} + A_{0} C_{2}\right)} \overline{\alpha_{2}}\\
	    			\mu_{10}&=\frac{1}{3} \, A_{1} {\left| A_{1} \right|}^{2} \overline{B_{1}} + \frac{1}{3} \, A_{0} {\left| A_{1} \right|}^{2} \overline{B_{2}} + \frac{1}{3} \, A_{2} {\left| A_{1} \right|}^{2} \overline{C_{1}} + \frac{1}{6} \, A_{0} C_{1} \alpha_{6} + \frac{1}{6} \, A_{0} \alpha_{5} \overline{C_{1}} + \frac{1}{6} \, A_{0} \alpha_{3} \overline{C_{2}} + \frac{1}{6} \, A_{0} B_{1} \overline{\alpha_{2}}\\
	    			& + \frac{1}{6} \, {\left(A_{0} \overline{B_{3}} + A_{1} \overline{C_{2}}\right)} \alpha_{1} + \frac{1}{6} \, {\left(A_{0} \overline{B_{1}} + A_{1} \overline{C_{1}}\right)} \alpha_{4} + \frac{1}{6} \, A_{2} \overline{B_{3}} + \frac{1}{6} \, A_{1} \overline{B_{5}} + \frac{1}{6} \, A_{0} \overline{B_{7}} + \frac{1}{6} \, A_{3} \overline{C_{2}}
	    		\end{align*}
	    		
	    		Then we compute the metric
	    		\footnotesize
	    		\begin{align*}
	    			&e^{2\lambda}=\\
	    			&\begin{dmatrix}
	    			 2 \, A_{0} \overline{A_{0}} & 3 & 3 & 0 \\
	    			2 \, A_{0} \overline{A_{1}} & 3 & 4 & 0 \\
	    			2 \, A_{0} \overline{A_{2}} + \frac{1}{4} \, \overline{A_{0}} \overline{C_{1}} & 3 & 5 & 0 \\
	    			2 \, A_{0} \overline{A_{3}} + \frac{1}{4} \, \overline{A_{1}} \overline{C_{1}} + \frac{1}{6} \, \overline{A_{0}} \overline{C_{2}} & 3 & 6 & 0 \\
	    			\frac{1}{8} \, C_{1} \alpha_{1} \overline{A_{0}} + \frac{1}{8} \, \overline{A_{0}} \overline{C_{1}} \overline{\alpha_{1}} + \frac{1}{64} \, C_{1}^{2} + \frac{1}{4} \, C_{3} \overline{A_{0}} - \frac{1}{32} \, \gamma_{1} \overline{A_{0}} + 2 \, A_{0} \overline{A_{4}} + \frac{1}{4} \, \overline{A_{2}} \overline{C_{1}} + \frac{1}{6} \, \overline{A_{1}} \overline{C_{2}} & 3 & 7 & 0 \\
	    			\nu_{1} & 3 & 8 & 0 \\
	    			\nu_{2} & 3 & 9 & 0 \\
	    			\frac{1}{8} \, A_{0} \overline{C_{1}} & 7 & 1 & 0 \\
	    			\frac{1}{8} \, A_{0} \overline{C_{2}} & 7 & 2 & 0 \\
	    			\frac{1}{8} \, A_{0} C_{1} \alpha_{1} + \frac{1}{8} \, A_{0} \overline{C_{1}} \overline{\alpha_{1}} + \frac{1}{4} \, A_{2} C_{1} + \frac{1}{6} \, A_{1} C_{2} + \frac{1}{4} \, A_{0} C_{3} - \frac{1}{32} \, A_{0} \gamma_{1} + 2 \, A_{4} \overline{A_{0}} + \frac{1}{64} \, \overline{C_{1}}^{2} & 7 & 3 & 0 \\
	    			\nu_3 & 7 & 4 & 0 \\
	    			\nu_4 & 7 & 5 & 0 \\
	    			\nu_5 & 8 & 4 & 0 \\
	    			\frac{1}{10} \, A_{0} \overline{B_{1}} + \frac{1}{8} \, A_{1} \overline{C_{1}} & 8 & 1 & 0 \\
	    			\frac{1}{12} \, A_{0} \alpha_{1} \overline{C_{1}} + \frac{1}{10} \, A_{1} \overline{B_{1}} + \frac{1}{12} \, A_{0} \overline{B_{2}} + \frac{1}{8} \, A_{2} \overline{C_{1}} & 9 & 1 & 0 \\
	    			\frac{1}{14} \, A_{0} \alpha_{3} \overline{C_{1}} + \frac{1}{14} \, A_{0} C_{1} \overline{\alpha_{2}} + \frac{1}{84} \, {\left(6 \, A_{0} \overline{B_{1}} + 7 \, A_{1} \overline{C_{1}}\right)} \alpha_{1} + \frac{1}{10} \, A_{2} \overline{B_{1}} + \frac{1}{12} \, A_{1} \overline{B_{2}} + \frac{1}{14} \, A_{0} \overline{B_{4}} + \frac{1}{8} \, A_{3} \overline{C_{1}} & 10 & 1 & 0 \\
	    			\nu_6 & 11 & 1 & 0 \\
	    			\frac{1}{5} \, A_{0} {\left| A_{1} \right|}^{2} \overline{C_{1}} + \frac{1}{10} \, A_{0} \overline{B_{3}} + \frac{1}{8} \, A_{1} \overline{C_{2}} & 8 & 2 & 0 \\
	    			\frac{1}{6} \, A_{0} {\left| A_{1} \right|}^{2} \overline{B_{1}} + \frac{1}{5} \, A_{1} {\left| A_{1} \right|}^{2} \overline{C_{1}} + \frac{1}{12} \, A_{0} \alpha_{4} \overline{C_{1}} + \frac{1}{12} \, A_{0} \alpha_{1} \overline{C_{2}} + \frac{1}{10} \, A_{1} \overline{B_{3}} + \frac{1}{12} \, A_{0} \overline{B_{5}} + \frac{1}{8} \, A_{2} \overline{C_{2}} & 9 & 2 & 0 \\
	    			\nu_7 & 10 & 2 & 0 \\
	    			\frac{1}{16} \, A_{0} \overline{E_{1}} & 11 & -1 & 0 \\
	    			\frac{1}{18} \, A_{0} \overline{C_{1}} \overline{\alpha_{2}} + \frac{1}{16} \, A_{1} \overline{E_{1}} + \frac{1}{18} \, A_{0} \overline{E_{2}} & 12 & -1 & 0 \\
	    			\frac{1}{20} \, A_{0} \alpha_{7} \overline{C_{1}} + \frac{1}{20} \, A_{0} \alpha_{1} \overline{E_{1}} + \frac{1}{16} \, A_{2} \overline{E_{1}} + \frac{1}{18} \, A_{1} \overline{E_{2}} + \frac{1}{20} \, A_{0} \overline{E_{4}} + \frac{1}{180} \, {\left(9 \, A_{0} \overline{B_{1}} + 10 \, A_{1} \overline{C_{1}}\right)} \overline{\alpha_{2}} & 13 & -1 & 0 \\
	    			\frac{1}{16} \, A_{0} \overline{E_{3}} & 11 & 0 & 0 \\
	    			\frac{1}{9} \, A_{0} {\left| A_{1} \right|}^{2} \overline{E_{1}} + \frac{1}{18} \, A_{0} \alpha_{6} \overline{C_{1}} + \frac{1}{18} \, A_{0} \overline{C_{2}} \overline{\alpha_{2}} + \frac{1}{16} \, A_{1} \overline{E_{3}} + \frac{1}{18} \, A_{0} \overline{E_{5}} & 12 & 0 & 0 \\
	    			\frac{1}{4} \, A_{0} \gamma_{1} + \overline{A_{0}} \overline{E_{1}} & 7 & 3 & 1 \\
	    			\overline{A_{1}} \overline{E_{1}} + \frac{1}{8} \, A_{0} \overline{\gamma_{2}} & 7 & 4 & 1 \\
	    			\frac{1}{4} \, \zeta_{0} \overline{A_{0}} \overline{C_{1}} + \frac{1}{8} \, A_{0} \overline{C_{1}} \overline{\zeta_{0}} + \frac{1}{64} \, {\left(16 \, A_{0} \overline{\alpha_{1}} + 3 \, \overline{C_{1}}\right)} \gamma_{1} + \overline{A_{2}} \overline{E_{1}} + \frac{1}{8} \, A_{0} \overline{\gamma_{3}} & 7 & 5 & 1 
	    			\end{dmatrix}
	    			\begin{dmatrix}
	    			\textbf{(1)}\\
	    			\textbf{(2)}\\
	    			\textbf{(3)}\\
	    			\textbf{(4)}\\
	    			\textbf{(5)}\\
	    			\textbf{(6)}\\
	    			\textbf{(7)}\\
	    			\textbf{(8)}\\
	    			\textbf{(9}\\
	    			\textbf{(10)}\\
	    			\textbf{(11)}\\
	    			\textbf{(12)}\\
	    			\textbf{(13)}\\
	    			\textbf{(14)}\\
	    			\textbf{(15)}\\
	    			\textbf{(16)}\\
	    			\textbf{(17)}\\
	    			\textbf{(18)}\\
	    			\textbf{(19)}\\
	    			\textbf{(20)}\\
	    			\textbf{(21)}\\
	    			\textbf{(22)}\\
	    			\textbf{(23)}\\
	    			\textbf{(24)}\\
	    			\textbf{(25)}\\
	    			\textbf{(26)}\\
	    			\textbf{(27)}\\
	    			\textbf{(28)}
	    			\end{dmatrix}
	    			\end{align*}
	    			\begin{align*}
	    			\begin{dmatrix}
	    			\frac{2}{5} \, A_{0} \gamma_{1} {\left| A_{1} \right|}^{2} + \overline{A_{1}} \overline{C_{1}} \overline{\alpha_{2}} + \frac{1}{5} \, A_{0} \gamma_{4} + \overline{A_{1}} \overline{E_{2}} + \frac{1}{8} \, A_{1} \overline{\gamma_{2}} & 8 & 4 & 1 \\
	    			\frac{1}{4} \, A_{0} C_{1} + 2 \, A_{2} \overline{A_{0}} & 5 & 3 & 0 \\
	    			\frac{1}{4} \, A_{1} C_{1} + \frac{1}{6} \, A_{0} C_{2} + 2 \, A_{3} \overline{A_{0}} & 6 & 3 & 0 \\
	    			\nu_8 & 8 & 3 & 0 \\
	    			\nu_9 & 9 & 3 & 0 \\
	    			\frac{1}{4} \, A_{0} B_{1} + 2 \, A_{2} \overline{A_{1}} & 5 & 4 & 0 \\
	    			\frac{1}{4} \, \alpha_{1} \overline{A_{0}} \overline{C_{1}} + \frac{1}{4} \, A_{0} C_{1} \overline{\alpha_{1}} + \frac{1}{4} \, A_{0} B_{2} + 2 \, A_{2} \overline{A_{2}} + \frac{1}{4} \, \overline{A_{0}} \overline{B_{2}} + \frac{5}{128} \, C_{1} \overline{C_{1}} & 5 & 5 & 0 \\
	    			\nu_{10} & 5 & 6 & 0 \\
	    			\nu_{11} & 5 & 7 & 0 \\
	    			\frac{1}{3} \, A_{0} C_{1} {\left| A_{1} \right|}^{2} + \frac{1}{4} \, A_{1} B_{1} + \frac{1}{6} \, A_{0} B_{3} + 2 \, A_{3} \overline{A_{1}} & 6 & 4 & 0 \\
	    			\nu_{12} & 6 & 5 & 0 \\
	    			\nu_{13} & 6 & 6 & 0 \\
	    			A_{0} E_{1} + \frac{1}{4} \, \gamma_{1} \overline{A_{0}} & 3 & 7 & 1 \\
	    			A_{0} C_{1} \alpha_{2} + A_{0} E_{2} + \frac{1}{4} \, \gamma_{1} \overline{A_{1}} + \frac{1}{10} \, \overline{A_{0}} \overline{\gamma_{2}} & 3 & 8 & 1 \\
	    			A_{0} B_{1} \alpha_{2} + A_{0} E_{1} \overline{\alpha_{1}} + A_{0} C_{1} \overline{\alpha_{7}} + \frac{1}{12} \, \overline{A_{0}} \overline{C_{1}} \overline{\zeta_{0}} + A_{0} E_{4} + \frac{1}{12} \, {\left(2 \, \overline{A_{0}} \overline{\alpha_{1}} + 3 \, \overline{A_{2}}\right)} \gamma_{1} + \frac{1}{8} \, E_{1} \overline{C_{1}} + \frac{1}{10} \, \overline{A_{1}} \overline{\gamma_{2}} + \frac{1}{12} \, \overline{A_{0}} \overline{\gamma_{3}} & 3 & 9 & 1 \\
	    			\nu_{14} & 4 & 7 & 0 \\
	    		\nu_{15} & 4 & 8 & 0 \\
	    			\overline{A_{0}} \overline{C_{1}} \overline{\alpha_{2}} + \frac{1}{4} \, A_{1} \gamma_{1} + \frac{1}{10} \, A_{0} \gamma_{2} + \overline{A_{0}} \overline{E_{2}} & 8 & 3 & 1 \\
	    			\frac{1}{12} \, A_{0} C_{1} \zeta_{0} + \alpha_{7} \overline{A_{0}} \overline{C_{1}} + \alpha_{1} \overline{A_{0}} \overline{E_{1}} + \overline{A_{0}} \overline{B_{1}} \overline{\alpha_{2}} + \frac{1}{12} \, {\left(2 \, A_{0} \alpha_{1} + 3 \, A_{2}\right)} \gamma_{1} + \frac{1}{10} \, A_{1} \gamma_{2} + \frac{1}{12} \, A_{0} \gamma_{3} + \frac{1}{8} \, C_{1} \overline{E_{1}} + \overline{A_{0}} \overline{E_{4}} & 9 & 3 & 1 \\
	    			\frac{1}{8} \, C_{1} \zeta_{0} \overline{A_{0}} + \frac{1}{4} \, A_{0} C_{1} \overline{\zeta_{0}} + A_{2} E_{1} + \frac{1}{64} \, {\left(16 \, \alpha_{1} \overline{A_{0}} + 3 \, C_{1}\right)} \gamma_{1} + \frac{1}{8} \, \gamma_{3} \overline{A_{0}} & 5 & 7 & 1 \\
	    			\frac{1}{16} \, A_{0} \zeta_{0} \overline{C_{1}} + \frac{1}{16} \, \overline{C_{1}} \overline{E_{1}} & 11 & 1 & 1 \\
	    			2 \, A_{1} \overline{A_{0}} & 4 & 3 & 0 \\
	    			2 \, A_{1} \overline{A_{1}} & 4 & 4 & 0 \\
	    			2 \, A_{1} \overline{A_{2}} + \frac{1}{4} \, \overline{A_{0}} \overline{B_{1}} & 4 & 5 & 0 \\
	    			\frac{1}{3} \, {\left| A_{1} \right|}^{2} \overline{A_{0}} \overline{C_{1}} + 2 \, A_{1} \overline{A_{3}} + \frac{1}{4} \, \overline{A_{1}} \overline{B_{1}} + \frac{1}{6} \, \overline{A_{0}} \overline{B_{3}} & 4 & 6 & 0 \\
	    			A_{1} E_{1} + \frac{1}{8} \, \gamma_{2} \overline{A_{0}} & 4 & 7 & 1 \\
	    			\frac{2}{5} \, \gamma_{1} {\left| A_{1} \right|}^{2} \overline{A_{0}} + A_{1} C_{1} \alpha_{2} + A_{1} E_{2} + \frac{1}{5} \, \gamma_{4} \overline{A_{0}} + \frac{1}{8} \, \gamma_{2} \overline{A_{1}} & 4 & 8 & 1 \\
	    			\frac{1}{8} \, C_{1} \overline{A_{0}} & 1 & 7 & 0 \\
	    			\frac{1}{10} \, B_{1} \overline{A_{0}} + \frac{1}{8} \, C_{1} \overline{A_{1}} & 1 & 8 & 0 \\
	    			\frac{1}{12} \, C_{1} \overline{A_{0}} \overline{\alpha_{1}} + \frac{1}{12} \, B_{2} \overline{A_{0}} + \frac{1}{10} \, B_{1} \overline{A_{1}} + \frac{1}{8} \, C_{1} \overline{A_{2}} & 1 & 9 & 0 
	    			\end{dmatrix}
	    				    			\begin{dmatrix}
	    			\textbf{(29}\\
	    			\textbf{(30)}\\
	    			\textbf{(31)}\\
	    			\textbf{(32)}\\
	    			\textbf{(33)}\\
	    			\textbf{(34)}\\
	    			\textbf{(35)}\\
	    			\textbf{(36)}\\
	    			\textbf{(37)}\\
	    			\textbf{(38)}\\
	    			\textbf{(39)}\\
	    			\textbf{(40)}\\
	    			\textbf{(41)}\\
	    			\textbf{(42)}\\
	    			\textbf{(43)}\\
	    			\textbf{(44)}\\
	    			\textbf{(45)}\\
	    			\textbf{(46)}\\
	    			\textbf{(47)}\\
	    			\textbf{(48)}\\
	    			\textbf{(49)}\\
	    			\textbf{(50)}\\
	    			\textbf{(51)}\\
	    			\textbf{(52)}\\
	    			\textbf{(53)}\\
	    			\textbf{(54)}\\
	    			\textbf{(55)}\\
	    			\textbf{(56)}\\
	    			\textbf{(57)}\\
	    			\textbf{(58)}
	    			\end{dmatrix}
	    			\end{align*}
	    			\begin{align*}
	    			\begin{dmatrix}
	    			\frac{1}{14} \, \alpha_{2} \overline{A_{0}} \overline{C_{1}} + \frac{1}{14} \, C_{1} \overline{A_{0}} \overline{\alpha_{3}} + \frac{1}{14} \, B_{4} \overline{A_{0}} + \frac{1}{12} \, B_{2} \overline{A_{1}} + \frac{1}{10} \, B_{1} \overline{A_{2}} + \frac{1}{8} \, C_{1} \overline{A_{3}} + \frac{1}{84} \, {\left(6 \, B_{1} \overline{A_{0}} + 7 \, C_{1} \overline{A_{1}}\right)} \overline{\alpha_{1}} & 1 & 10 & 0 \\
	    			\nu_{16} & 1 & 11 & 0 \\
	    			\frac{1}{16} \, C_{1} \overline{A_{0}} \overline{\zeta_{0}} + \frac{1}{16} \, C_{1} E_{1} & 1 & 11 & 1 \\
	    			\frac{1}{8} \, C_{2} \overline{A_{0}} & 2 & 7 & 0 \\
	    			\frac{1}{5} \, C_{1} {\left| A_{1} \right|}^{2} \overline{A_{0}} + \frac{1}{10} \, B_{3} \overline{A_{0}} + \frac{1}{8} \, C_{2} \overline{A_{1}} & 2 & 8 & 0 \\
	    			\frac{1}{6} \, B_{1} {\left| A_{1} \right|}^{2} \overline{A_{0}} + \frac{1}{5} \, C_{1} {\left| A_{1} \right|}^{2} \overline{A_{1}} + \frac{1}{12} \, C_{2} \overline{A_{0}} \overline{\alpha_{1}} + \frac{1}{12} \, C_{1} \overline{A_{0}} \overline{\alpha_{4}} + \frac{1}{12} \, B_{5} \overline{A_{0}} + \frac{1}{10} \, B_{3} \overline{A_{1}} + \frac{1}{8} \, C_{2} \overline{A_{2}} & 2 & 9 & 0 \\
	    			\nu_{17} & 2 & 10 & 0 \\
	    			\frac{1}{16} \, E_{1} \overline{A_{0}} & -1 & 11 & 0 \\
	    			\frac{1}{18} \, C_{1} \alpha_{2} \overline{A_{0}} + \frac{1}{18} \, E_{2} \overline{A_{0}} + \frac{1}{16} \, E_{1} \overline{A_{1}} & -1 & 12 & 0 \\
	    			\frac{1}{20} \, E_{1} \overline{A_{0}} \overline{\alpha_{1}} + \frac{1}{20} \, C_{1} \overline{A_{0}} \overline{\alpha_{7}} + \frac{1}{180} \, {\left(9 \, B_{1} \overline{A_{0}} + 10 \, C_{1} \overline{A_{1}}\right)} \alpha_{2} + \frac{1}{20} \, E_{4} \overline{A_{0}} + \frac{1}{18} \, E_{2} \overline{A_{1}} + \frac{1}{16} \, E_{1} \overline{A_{2}} & -1 & 13 & 0 \\
	    			\frac{1}{16} \, E_{3} \overline{A_{0}} & 0 & 11 & 0 \\
	    			\frac{1}{9} \, E_{1} {\left| A_{1} \right|}^{2} \overline{A_{0}} + \frac{1}{18} \, C_{2} \alpha_{2} \overline{A_{0}} + \frac{1}{18} \, C_{1} \overline{A_{0}} \overline{\alpha_{6}} + \frac{1}{18} \, E_{5} \overline{A_{0}} + \frac{1}{16} \, E_{3} \overline{A_{1}} & 0 & 12 & 0
	    			\end{dmatrix}
	    			\begin{dmatrix}
	    			\textbf{(59)}\\
	    			\textbf{(60)}\\
	    			\textbf{(61)}\\
	    			\textbf{(62)}\\
	    			\textbf{(63)}\\
	    			\textbf{(64)}\\
	    			\textbf{(65)}\\
	    			\textbf{(66)}\\
	    			\textbf{(67)}\\
	    			\textbf{(68)}\\
	    			\textbf{(69)}\\
	    			\textbf{(70)}
	    			\end{dmatrix}
	    		\end{align*}
	    		\small
	    		where
	    		\begin{align*}
	    			\nu_{1}&=\frac{1}{5} \, C_{2} {\left| A_{1} \right|}^{2} \overline{A_{0}} + \frac{1}{10} \, C_{1} \alpha_{4} \overline{A_{0}} + \frac{1}{10} \, \overline{A_{0}} \overline{C_{1}} \overline{\alpha_{3}} + \frac{9}{320} \, B_{1} C_{1} + \frac{1}{40} \, {\left(4 \, B_{1} \overline{A_{0}} + 5 \, C_{1} \overline{A_{1}}\right)} \alpha_{1} + \frac{1}{4} \, C_{3} \overline{A_{1}} - \frac{1}{32} \, \gamma_{1} \overline{A_{1}} + 2 \, A_{0} \overline{A_{5}}\\
	    			& + \frac{1}{4} \, \overline{A_{3}} \overline{C_{1}} + \frac{1}{6} \, \overline{A_{2}} \overline{C_{2}} + \frac{1}{10} \, \overline{A_{0}} \overline{C_{4}} + \frac{1}{40} \, {\left(5 \, \overline{A_{1}} \overline{C_{1}} + 4 \, \overline{A_{0}} \overline{C_{2}}\right)} \overline{\alpha_{1}} - \frac{1}{100} \, \overline{A_{0}} \overline{\gamma_{2}}\\
	    			\nu_{2}&=\frac{1}{6} \, B_{3} {\left| A_{1} \right|}^{2} \overline{A_{0}} + \frac{1}{5} \, C_{2} {\left| A_{1} \right|}^{2} \overline{A_{1}} + \frac{1}{12} \, C_{1} \beta \overline{A_{0}} + \frac{1}{12} \, C_{2} \overline{A_{0}} \overline{\alpha_{4}} - \frac{1}{144} \, \overline{A_{0}} \overline{C_{1}} \overline{\zeta_{0}} + \frac{1}{80} \, B_{1}^{2} + \frac{5}{192} \, B_{2} C_{1}\\
	    			& + \frac{1}{120} \, {\left(10 \, B_{2} \overline{A_{0}} + 12 \, B_{1} \overline{A_{1}} + 15 \, C_{1} \overline{A_{2}}\right)} \alpha_{1} + \frac{1}{60} \, {\left(5 \, B_{1} \overline{A_{0}} + 6 \, C_{1} \overline{A_{1}}\right)} \alpha_{4} - \frac{1}{288} \, {\left(4 \, \overline{A_{0}} \overline{\alpha_{1}} + 9 \, \overline{A_{2}}\right)} \gamma_{1} + \frac{1}{4} \, C_{3} \overline{A_{2}}\\
	    			& + 2 \, A_{0} \overline{A_{6}} + \frac{1}{256} \, E_{1} \overline{C_{1}} + \frac{1}{4} \, \overline{A_{4}} \overline{C_{1}} + \frac{1}{6} \, \overline{A_{3}} \overline{C_{2}} + \frac{1}{10} \, \overline{A_{1}} \overline{C_{4}} + \frac{1}{12} \, \overline{A_{0}} \overline{C_{5}} + \frac{1}{960} \, {\left(25 \, C_{1}^{2} + 160 \, C_{3} \overline{A_{0}} + 120 \, \overline{A_{2}} \overline{C_{1}} + 96 \, \overline{A_{1}} \overline{C_{2}}\right)} \overline{\alpha_{1}}\\
	    			& + \frac{1}{60} \, {\left(6 \, \overline{A_{1}} \overline{C_{1}} + 5 \, \overline{A_{0}} \overline{C_{2}}\right)} \overline{\alpha_{3}} - \frac{1}{100} \, \overline{A_{1}} \overline{\gamma_{2}} - \frac{1}{144} \, \overline{A_{0}} \overline{\gamma_{3}}\\
	    			\nu_{3}&=\frac{1}{3} \, A_{1} C_{1} {\left| A_{1} \right|}^{2} + \frac{1}{4} \, A_{0} C_{2} {\left| A_{1} \right|}^{2} + \frac{1}{8} \, A_{0} B_{1} \alpha_{1} + \frac{1}{8} \, A_{0} C_{1} \alpha_{4} + \frac{1}{8} \, A_{0} \overline{C_{2}} \overline{\alpha_{1}} + \frac{1}{8} \, A_{0} \overline{C_{1}} \overline{\alpha_{3}} + \frac{1}{4} \, A_{2} B_{1} + \frac{1}{6} \, A_{1} B_{3} + 2 \, A_{4} \overline{A_{1}}\\
	    			& + \frac{5}{192} \, \overline{C_{1}} \overline{C_{2}} + \frac{1}{8} \, A_{0} \overline{C_{4}} + \frac{1}{2} \, \overline{A_{0}} \overline{E_{3}} - \frac{1}{64} \, A_{0} \overline{\gamma_{2}}\\
	    			\nu_{4}&=\frac{1}{3} \, A_{1} B_{1} {\left| A_{1} \right|}^{2} + \frac{1}{4} \, A_{0} B_{3} {\left| A_{1} \right|}^{2} + \frac{1}{8} \, A_{0} B_{1} \alpha_{4} + \frac{1}{8} \, A_{0} C_{1} \beta + \frac{1}{4} \, C_{1} \alpha_{7} \overline{A_{0}} + \frac{1}{4} \, \alpha_{3} \overline{A_{0}} \overline{B_{1}} - \frac{1}{16} \, \zeta_{0} \overline{A_{0}} \overline{C_{1}} + \frac{1}{4} \, C_{2} \overline{A_{0}} \overline{\alpha_{2}}\\
	    			& + \frac{1}{8} \, A_{0} \overline{C_{2}} \overline{\alpha_{3}} - \frac{1}{64} \, A_{0} \overline{C_{1}} \overline{\zeta_{0}} + \frac{1}{4} \, A_{2} B_{2} + \frac{1}{6} \, A_{1} B_{5} + \frac{1}{384} \, {\left(48 \, A_{0} B_{2} + 96 \, \overline{A_{0}} \overline{B_{2}} + 23 \, C_{1} \overline{C_{1}}\right)} \alpha_{1} - \frac{1}{512} \, {\left(16 \, A_{0} \overline{\alpha_{1}} + 3 \, \overline{C_{1}}\right)} \gamma_{1}\\
	    			& + 2 \, A_{4} \overline{A_{2}} + \frac{13}{480} \, C_{2} \overline{B_{1}} + \frac{7}{192} \, C_{1} \overline{B_{2}} + \frac{1}{4} \, \overline{A_{0}} \overline{B_{6}} + \frac{3}{64} \, C_{3} \overline{C_{1}} + \frac{1}{96} \, \overline{C_{2}}^{2} + \frac{1}{8} \, A_{0} \overline{C_{5}} + \frac{1}{2} \, \overline{A_{1}} \overline{E_{3}}\\
	    			& + \frac{1}{384} \, {\left(96 \, A_{2} C_{1} + 64 \, A_{1} C_{2} + 96 \, A_{0} C_{3} + 9 \, \overline{C_{1}}^{2} + 96 \, \overline{A_{0}} \overline{E_{1}}\right)} \overline{\alpha_{1}} + \frac{1}{24} \, {\left(4 \, A_{1} C_{1} + 3 \, A_{0} C_{2}\right)} \overline{\alpha_{4}} - \frac{1}{64} \, A_{0} \overline{\gamma_{3}}\\
	    			\nu_{5}&=\frac{1}{3} \, A_{2} C_{1} {\left| A_{1} \right|}^{2} + \frac{1}{4} \, A_{1} C_{2} {\left| A_{1} \right|}^{2} + \frac{2}{5} \, A_{0} C_{3} {\left| A_{1} \right|}^{2} - \frac{1}{25} \, A_{0} \gamma_{1} {\left| A_{1} \right|}^{2} + \frac{11}{240} \, {\left| A_{1} \right|}^{2} \overline{C_{1}}^{2} + {\left| A_{1} \right|}^{2} \overline{A_{0}} \overline{E_{1}} + \frac{1}{10} \, A_{0} B_{1} \alpha_{3}\\
	    			& + \frac{1}{10} \, A_{0} C_{1} \alpha_{5} + \frac{1}{2} \, \alpha_{6} \overline{A_{0}} \overline{C_{1}} + \frac{1}{2} \, \overline{A_{0}} \overline{C_{2}} \overline{\alpha_{2}} + \frac{1}{10} \, A_{0} \overline{C_{2}} \overline{\alpha_{4}} + \frac{1}{10} \, A_{0} \overline{C_{1}} \overline{\alpha_{5}} + \frac{1}{4} \, A_{3} B_{1} + \frac{1}{6} \, A_{2} B_{3} + \frac{1}{5} \, A_{0} C_{6}\\
	    			& + \frac{1}{40} \, {\left(5 \, A_{1} B_{1} + 4 \, A_{0} B_{3}\right)} \alpha_{1} + \frac{1}{40} \, {\left(5 \, A_{1} C_{1} + 4 \, A_{0} C_{2}\right)} \alpha_{4} - \frac{1}{50} \, A_{0} \gamma_{4} + 2 \, A_{5} \overline{A_{1}} + \frac{11}{480} \, \overline{B_{3}} \overline{C_{1}} + \frac{23}{960} \, \overline{B_{1}} \overline{C_{2}} + \frac{1}{8} \, A_{1} \overline{C_{4}}\\
	    			& + \frac{1}{2} \, \overline{A_{0}} \overline{E_{5}} + \frac{1}{40} \, {\left(4 \, A_{0} \overline{B_{3}} + 5 \, A_{1} \overline{C_{2}}\right)} \overline{\alpha_{1}} + \frac{1}{40} \, {\left(4 \, A_{0} \overline{B_{1}} + 5 \, A_{1} \overline{C_{1}}\right)} \overline{\alpha_{3}} - \frac{1}{64} \, A_{1} \overline{\gamma_{2}}\\
	    			\nu_{6}&=\frac{1}{16} \, A_{0} C_{1} \alpha_{7} - \frac{1}{256} \, A_{0} \zeta_{0} \overline{C_{1}} + \frac{1}{16} \, A_{0} \overline{E_{1}} \overline{\alpha_{1}} + \frac{1}{336} \, {\left(24 \, A_{1} \overline{B_{1}} + 21 \, A_{0} \overline{B_{2}} + 28 \, A_{2} \overline{C_{1}}\right)} \alpha_{1} + \frac{1}{112} \, {\left(7 \, A_{0} \overline{B_{1}} + 8 \, A_{1} \overline{C_{1}}\right)} \alpha_{3}\\
	    			& + \frac{1}{10} \, A_{3} \overline{B_{1}} + \frac{1}{12} \, A_{2} \overline{B_{2}} + \frac{1}{14} \, A_{1} \overline{B_{4}} + \frac{1}{16} \, A_{0} \overline{B_{6}} + \frac{1}{8} \, A_{4} \overline{C_{1}} + \frac{1}{128} \, \overline{C_{1}} \overline{E_{1}} + \frac{1}{112} \, {\left(8 \, A_{1} C_{1} + 7 \, A_{0} C_{2}\right)} \overline{\alpha_{2}}\\
	    			\nu_{7}&=\frac{1}{6} \, A_{1} {\left| A_{1} \right|}^{2} \overline{B_{1}} + \frac{1}{7} \, A_{0} {\left| A_{1} \right|}^{2} \overline{B_{2}} + \frac{1}{5} \, A_{2} {\left| A_{1} \right|}^{2} \overline{C_{1}} + \frac{1}{14} \, A_{0} C_{1} \alpha_{6} + \frac{1}{14} \, A_{0} \alpha_{5} \overline{C_{1}} + \frac{1}{14} \, A_{0} \alpha_{3} \overline{C_{2}} + \frac{1}{14} \, A_{0} B_{1} \overline{\alpha_{2}}\\
	    			& + \frac{1}{84} \, {\left(6 \, A_{0} \overline{B_{3}} + 7 \, A_{1} \overline{C_{2}}\right)} \alpha_{1} + \frac{1}{84} \, {\left(6 \, A_{0} \overline{B_{1}} + 7 \, A_{1} \overline{C_{1}}\right)} \alpha_{4} + \frac{1}{10} \, A_{2} \overline{B_{3}} + \frac{1}{12} \, A_{1} \overline{B_{5}} + \frac{1}{14} \, A_{0} \overline{B_{7}} + \frac{1}{8} \, A_{3} \overline{C_{2}}\\
	    			\nu_{8}&=\frac{1}{5} \, A_{0} {\left| A_{1} \right|}^{2} \overline{C_{2}} + \frac{1}{10} \, A_{0} C_{1} \alpha_{3} + \frac{1}{10} \, A_{0} \overline{C_{1}} \overline{\alpha_{4}} + \frac{1}{4} \, A_{3} C_{1} + \frac{1}{6} \, A_{2} C_{2} + \frac{1}{4} \, A_{1} C_{3} + \frac{1}{10} \, A_{0} C_{4} + \frac{1}{40} \, {\left(5 \, A_{1} C_{1} + 4 \, A_{0} C_{2}\right)} \alpha_{1}\\
	    			& - \frac{1}{32} \, A_{1} \gamma_{1} - \frac{1}{100} \, A_{0} \gamma_{2} + 2 \, A_{5} \overline{A_{0}} + \frac{9}{320} \, \overline{B_{1}} \overline{C_{1}} + \frac{1}{40} \, {\left(4 \, A_{0} \overline{B_{1}} + 5 \, A_{1} \overline{C_{1}}\right)} \overline{\alpha_{1}}\\
	    			\nu_{9}&=\frac{1}{6} \, A_{0} {\left| A_{1} \right|}^{2} \overline{B_{3}} + \frac{1}{5} \, A_{1} {\left| A_{1} \right|}^{2} \overline{C_{2}} - \frac{1}{144} \, A_{0} C_{1} \zeta_{0} + \frac{1}{12} \, A_{0} \beta \overline{C_{1}} + \frac{1}{12} \, A_{0} \alpha_{4} \overline{C_{2}} + \frac{1}{4} \, A_{4} C_{1} + \frac{1}{6} \, A_{3} C_{2} + \frac{1}{4} \, A_{2} C_{3} + \frac{1}{10} \, A_{1} C_{4}\\
	    			& + \frac{1}{12} \, A_{0} C_{5} + \frac{1}{960} \, {\left(120 \, A_{2} C_{1} + 96 \, A_{1} C_{2} + 160 \, A_{0} C_{3} + 25 \, \overline{C_{1}}^{2}\right)} \alpha_{1} + \frac{1}{60} \, {\left(6 \, A_{1} C_{1} + 5 \, A_{0} C_{2}\right)} \alpha_{3} - \frac{1}{288} \, {\left(4 \, A_{0} \alpha_{1} + 9 \, A_{2}\right)} \gamma_{1}\\
	    			& - \frac{1}{100} \, A_{1} \gamma_{2} - \frac{1}{144} \, A_{0} \gamma_{3} + 2 \, A_{6} \overline{A_{0}} + \frac{1}{80} \, \overline{B_{1}}^{2} + \frac{5}{192} \, \overline{B_{2}} \overline{C_{1}} + \frac{1}{256} \, C_{1} \overline{E_{1}} + \frac{1}{120} \, {\left(12 \, A_{1} \overline{B_{1}} + 10 \, A_{0} \overline{B_{2}} + 15 \, A_{2} \overline{C_{1}}\right)} \overline{\alpha_{1}}\\
	    			& + \frac{1}{60} \, {\left(5 \, A_{0} \overline{B_{1}} + 6 \, A_{1} \overline{C_{1}}\right)} \overline{\alpha_{4}}\\
	    			\nu_{10}&=\frac{1}{3} \, {\left| A_{1} \right|}^{2} \overline{A_{0}} \overline{B_{1}} + \frac{1}{4} \, A_{0} \alpha_{2} \overline{C_{1}} + \frac{1}{6} \, \alpha_{4} \overline{A_{0}} \overline{C_{1}} + \frac{1}{4} \, A_{0} B_{1} \overline{\alpha_{1}} + \frac{1}{4} \, A_{0} C_{1} \overline{\alpha_{3}} + \frac{1}{4} \, A_{0} B_{4} + \frac{1}{12} \, {\left(3 \, \overline{A_{1}} \overline{C_{1}} + 2 \, \overline{A_{0}} \overline{C_{2}}\right)} \alpha_{1} + 2 \, A_{2} \overline{A_{3}}\\
	    			& + \frac{1}{4} \, \overline{A_{1}} \overline{B_{2}} + \frac{1}{6} \, \overline{A_{0}} \overline{B_{5}} + \frac{3}{80} \, B_{1} \overline{C_{1}} + \frac{11}{384} \, C_{1} \overline{C_{2}}\\
	    			\nu_{11}&=\frac{1}{3} \, {\left| A_{1} \right|}^{2} \overline{A_{1}} \overline{B_{1}} + \frac{1}{4} \, {\left| A_{1} \right|}^{2} \overline{A_{0}} \overline{B_{3}} + \frac{1}{8} \, C_{2} \alpha_{3} \overline{A_{0}} - \frac{1}{64} \, C_{1} \zeta_{0} \overline{A_{0}} + \frac{1}{8} \, \beta \overline{A_{0}} \overline{C_{1}} + \frac{1}{4} \, A_{0} \alpha_{2} \overline{C_{2}} + \frac{1}{4} \, A_{0} B_{1} \overline{\alpha_{3}} + \frac{1}{8} \, \overline{A_{0}} \overline{B_{1}} \overline{\alpha_{4}}\\
	    			& + \frac{1}{4} \, A_{0} \overline{C_{1}} \overline{\alpha_{7}} - \frac{1}{16} \, A_{0} C_{1} \overline{\zeta_{0}} + \frac{1}{4} \, A_{0} B_{6} + \frac{1}{96} \, C_{2}^{2} + \frac{3}{64} \, C_{1} C_{3} + \frac{1}{2} \, A_{1} E_{3}\\
	    			& + \frac{1}{384} \, {\left(9 \, C_{1}^{2} + 96 \, A_{0} E_{1} + 96 \, C_{3} \overline{A_{0}} + 96 \, \overline{A_{2}} \overline{C_{1}} + 64 \, \overline{A_{1}} \overline{C_{2}}\right)} \alpha_{1} + \frac{1}{24} \, {\left(4 \, \overline{A_{1}} \overline{C_{1}} + 3 \, \overline{A_{0}} \overline{C_{2}}\right)} \alpha_{4} - \frac{1}{512} \, {\left(16 \, \alpha_{1} \overline{A_{0}} + 3 \, C_{1}\right)} \gamma_{1}\\
	    			& + \frac{1}{8} \, C_{5} \overline{A_{0}} - \frac{1}{64} \, \gamma_{3} \overline{A_{0}} + 2 \, A_{2} \overline{A_{4}} + \frac{1}{4} \, \overline{A_{2}} \overline{B_{2}} + \frac{1}{6} \, \overline{A_{1}} \overline{B_{5}} + \frac{7}{192} \, B_{2} \overline{C_{1}} + \frac{13}{480} \, B_{1} \overline{C_{2}} + \frac{1}{384} \, {\left(96 \, A_{0} B_{2} + 48 \, \overline{A_{0}} \overline{B_{2}} + 23 \, C_{1} \overline{C_{1}}\right)} \overline{\alpha_{1}}\\
	    			\nu_{12}&=\frac{1}{3} \, A_{0} B_{1} {\left| A_{1} \right|}^{2} + \frac{1}{4} \, \alpha_{1} \overline{A_{0}} \overline{B_{1}} + \frac{1}{4} \, \alpha_{3} \overline{A_{0}} \overline{C_{1}} + \frac{1}{4} \, C_{1} \overline{A_{0}} \overline{\alpha_{2}} + \frac{1}{6} \, A_{0} C_{1} \overline{\alpha_{4}} + \frac{1}{4} \, A_{1} B_{2} + \frac{1}{6} \, A_{0} B_{5} + 2 \, A_{3} \overline{A_{2}} + \frac{3}{80} \, C_{1} \overline{B_{1}}\\
	    			& + \frac{1}{4} \, \overline{A_{0}} \overline{B_{4}} + \frac{11}{384} \, C_{2} \overline{C_{1}} + \frac{1}{12} \, {\left(3 \, A_{1} C_{1} + 2 \, A_{0} C_{2}\right)} \overline{\alpha_{1}}\\
	    			\nu_{13}&=\frac{1}{3} \, A_{0} B_{2} {\left| A_{1} \right|}^{2} + \frac{1}{3} \, {\left| A_{1} \right|}^{2} \overline{A_{0}} \overline{B_{2}} + \frac{1}{6} \, C_{1} \alpha_{6} \overline{A_{0}} + \frac{1}{6} \, \alpha_{4} \overline{A_{0}} \overline{B_{1}} + \frac{1}{6} \, \alpha_{5} \overline{A_{0}} \overline{C_{1}} + \frac{1}{6} \, A_{0} B_{1} \overline{\alpha_{4}} + \frac{1}{6} \, A_{0} C_{1} \overline{\alpha_{5}} + \frac{1}{6} \, A_{0} \overline{C_{1}} \overline{\alpha_{6}}\\
	    			& + \frac{1}{4} \, A_{1} B_{4} + \frac{1}{6} \, A_{0} B_{7} + \frac{1}{12} \, {\left(3 \, \overline{A_{1}} \overline{B_{1}} + 2 \, \overline{A_{0}} \overline{B_{3}}\right)} \alpha_{1} + \frac{1}{12} \, {\left(2 \, A_{0} \overline{B_{1}} + 3 \, A_{1} \overline{C_{1}}\right)} \alpha_{2} + \frac{1}{12} \, {\left(3 \, \overline{A_{1}} \overline{C_{1}} + 2 \, \overline{A_{0}} \overline{C_{2}}\right)} \alpha_{3} + 2 \, A_{3} \overline{A_{3}}\\
	    			& + \frac{29}{800} \, B_{1} \overline{B_{1}} + \frac{13}{480} \, C_{1} \overline{B_{3}} + \frac{1}{4} \, \overline{A_{1}} \overline{B_{4}} + \frac{1}{6} \, \overline{A_{0}} \overline{B_{7}} + \frac{13}{480} \, {\left(4 \, C_{1} {\left| A_{1} \right|}^{2} + B_{3}\right)} \overline{C_{1}} + \frac{25}{1152} \, C_{2} \overline{C_{2}} + \frac{1}{12} \, {\left(3 \, A_{1} B_{1} + 2 \, A_{0} B_{3}\right)} \overline{\alpha_{1}}\\
	    			& + \frac{1}{12} \, {\left(2 \, B_{1} \overline{A_{0}} + 3 \, C_{1} \overline{A_{1}}\right)} \overline{\alpha_{2}} + \frac{1}{12} \, {\left(3 \, A_{1} C_{1} + 2 \, A_{0} C_{2}\right)} \overline{\alpha_{3}}\\
	    			\nu_{14}&=\frac{1}{3} \, {\left| A_{1} \right|}^{2} \overline{A_{1}} \overline{C_{1}} + \frac{1}{4} \, {\left| A_{1} \right|}^{2} \overline{A_{0}} \overline{C_{2}} + \frac{1}{8} \, C_{2} \alpha_{1} \overline{A_{0}} + \frac{1}{8} \, C_{1} \alpha_{3} \overline{A_{0}} + \frac{1}{8} \, \overline{A_{0}} \overline{B_{1}} \overline{\alpha_{1}} + \frac{1}{8} \, \overline{A_{0}} \overline{C_{1}} \overline{\alpha_{4}} + \frac{5}{192} \, C_{1} C_{2} + \frac{1}{2} \, A_{0} E_{3} + \frac{1}{8} \, C_{4} \overline{A_{0}}\\
	    			& - \frac{1}{64} \, \gamma_{2} \overline{A_{0}} + 2 \, A_{1} \overline{A_{4}} + \frac{1}{4} \, \overline{A_{2}} \overline{B_{1}} + \frac{1}{6} \, \overline{A_{1}} \overline{B_{3}}\\
	    			\nu_{15}&=	\frac{11}{240} \, C_{1}^{2} {\left| A_{1} \right|}^{2} + A_{0} E_{1} {\left| A_{1} \right|}^{2} + \frac{2}{5} \, C_{3} {\left| A_{1} \right|}^{2} \overline{A_{0}} - \frac{1}{25} \, \gamma_{1} {\left| A_{1} \right|}^{2} \overline{A_{0}} + \frac{1}{3} \, {\left| A_{1} \right|}^{2} \overline{A_{2}} \overline{C_{1}} + \frac{1}{4} \, {\left| A_{1} \right|}^{2} \overline{A_{1}} \overline{C_{2}} + \frac{1}{2} \, A_{0} C_{2} \alpha_{2} + \frac{1}{10} \, C_{2} \alpha_{4} \overline{A_{0}}\\
	    			& + \frac{1}{10} \, C_{1} \alpha_{5} \overline{A_{0}} + \frac{1}{10} \, \overline{A_{0}} \overline{B_{1}} \overline{\alpha_{3}} + \frac{1}{10} \, \overline{A_{0}} \overline{C_{1}} \overline{\alpha_{5}} + \frac{1}{2} \, A_{0} C_{1} \overline{\alpha_{6}} + \frac{11}{480} \, B_{3} C_{1} + \frac{23}{960} \, B_{1} C_{2} + \frac{1}{2} \, A_{0} E_{5} + \frac{1}{40} \, {\left(4 \, B_{3} \overline{A_{0}} + 5 \, C_{2} \overline{A_{1}}\right)} \alpha_{1}\\
	    			& + \frac{1}{40} \, {\left(4 \, B_{1} \overline{A_{0}} + 5 \, C_{1} \overline{A_{1}}\right)} \alpha_{3} + \frac{1}{5} \, C_{6} \overline{A_{0}} - \frac{1}{50} \, \gamma_{4} \overline{A_{0}} + \frac{1}{8} \, C_{4} \overline{A_{1}} - \frac{1}{64} \, \gamma_{2} \overline{A_{1}} + 2 \, A_{1} \overline{A_{5}} + \frac{1}{4} \, \overline{A_{3}} \overline{B_{1}} + \frac{1}{6} \, \overline{A_{2}} \overline{B_{3}}\\
	    			& + \frac{1}{40} \, {\left(5 \, \overline{A_{1}} \overline{B_{1}} + 4 \, \overline{A_{0}} \overline{B_{3}}\right)} \overline{\alpha_{1}} + \frac{1}{40} \, {\left(5 \, \overline{A_{1}} \overline{C_{1}} + 4 \, \overline{A_{0}} \overline{C_{2}}\right)} \overline{\alpha_{4}}\\
	    			\nu_{16}&=\frac{1}{16} \, E_{1} \alpha_{1} \overline{A_{0}} + \frac{1}{16} \, \overline{A_{0}} \overline{C_{1}} \overline{\alpha_{7}} - \frac{1}{256} \, C_{1} \overline{A_{0}} \overline{\zeta_{0}} + \frac{1}{128} \, C_{1} E_{1} + \frac{1}{112} \, {\left(8 \, \overline{A_{1}} \overline{C_{1}} + 7 \, \overline{A_{0}} \overline{C_{2}}\right)} \alpha_{2} + \frac{1}{16} \, B_{6} \overline{A_{0}} + \frac{1}{14} \, B_{4} \overline{A_{1}}\\
	    			& + \frac{1}{12} \, B_{2} \overline{A_{2}} + \frac{1}{10} \, B_{1} \overline{A_{3}} + \frac{1}{8} \, C_{1} \overline{A_{4}} + \frac{1}{336} \, {\left(21 \, B_{2} \overline{A_{0}} + 24 \, B_{1} \overline{A_{1}} + 28 \, C_{1} \overline{A_{2}}\right)} \overline{\alpha_{1}} + \frac{1}{112} \, {\left(7 \, B_{1} \overline{A_{0}} + 8 \, C_{1} \overline{A_{1}}\right)} \overline{\alpha_{3}}\\
	    			\nu_{17}&=\frac{1}{7} \, B_{2} {\left| A_{1} \right|}^{2} \overline{A_{0}} + \frac{1}{6} \, B_{1} {\left| A_{1} \right|}^{2} \overline{A_{1}} + \frac{1}{5} \, C_{1} {\left| A_{1} \right|}^{2} \overline{A_{2}} + \frac{1}{14} \, \alpha_{2} \overline{A_{0}} \overline{B_{1}} + \frac{1}{14} \, C_{2} \overline{A_{0}} \overline{\alpha_{3}} + \frac{1}{14} \, C_{1} \overline{A_{0}} \overline{\alpha_{5}} + \frac{1}{14} \, \overline{A_{0}} \overline{C_{1}} \overline{\alpha_{6}} + \frac{1}{14} \, B_{7} \overline{A_{0}}\\
	    			& + \frac{1}{12} \, B_{5} \overline{A_{1}} + \frac{1}{10} \, B_{3} \overline{A_{2}} + \frac{1}{8} \, C_{2} \overline{A_{3}} + \frac{1}{84} \, {\left(6 \, B_{3} \overline{A_{0}} + 7 \, C_{2} \overline{A_{1}}\right)} \overline{\alpha_{1}} + \frac{1}{84} \, {\left(6 \, B_{1} \overline{A_{0}} + 7 \, C_{1} \overline{A_{1}}\right)} \overline{\alpha_{4}}
	    		\end{align*}
	    		\normalsize
	    		There are some non-trivial cancellations which we used to rule out the dual $\textbf{(22)}-\textbf{(67)}$, which comes from the fact that thanks of \eqref{4h4}, we have
    		    \begin{align*}
    		    	\frac{1}{2}\vec{E}_2=\begin{dmatrix}
    		    	-\frac{9}{160} \, C_{1}^{2} & \overline{A_{1}} & -4 & 5 & 0 &\textbf{(3)}\\
    		    	\frac{1}{2} \ccancel{\, A_{0} C_{1} \alpha_{2}} - \frac{17}{160} \ccancel{\, B_{1} C_{1}} & \overline{A_{0}} & -4 & 5 & 0 &\textbf{(4)}\\
    		    	-\frac{1}{80} \, C_{1} \overline{A_{1}} & C_{1} & -4 & 5 & 0 &\textbf{(5)}\\
    		    	\ccancel{C_{1} \alpha_{2} \overline{A_{0}}} - \frac{1}{5} \ccancel{\, E_{1} \overline{A_{1}}} & A_{0} & -4 & 5 & 0 &\textbf{(6)}
    		    	\end{dmatrix}=-\frac{9}{160}\s{\vec{C}_1}{\vec{C}_1}\bar{\vec{A}_1}-\frac{1}{80}\s{\bar{\vec{A}_1}}{\vec{C}_1}\vec{C}_1
    		    \end{align*}
    		    so that
    		    \begin{align*}
    		    	\s{\vec{A}_0}{\vec{E}_2}=\s{\bar{\vec{A}_0}}{\vec{E}_2}=0.
    		    \end{align*}
    		    Furthermore, as by \eqref{4h4}, we have
    		    \begin{align*}
    		    	\frac{1}{2}\vec{E}_3=\begin{dmatrix}
    		    	-\frac{1}{8} \, C_{1} C_{2} + \frac{1}{3} \ccancel{\, A_{1} E_{1}} & \overline{A_{0}} & -3 & 4 & 0 &\textbf{(13)}
    		    	\end{dmatrix}=-\frac{1}{8}\s{\vec{C}_1}{\vec{C}_2}\bar{\vec{A}_0}
    		    \end{align*}
    		    so that (see $\textbf{(24)}-\textbf{(69)}$)
    		    \begin{align*}
    		    	\s{\bar{\vec{A}_0}}{\vec{E}_3}=\s{\vec{A}_1}{\vec{E}_3}=\s{\bar{\vec{A}_1}}{\vec{E}_3}=0
    		    \end{align*}
    		    The relation for $\textbf{(27)}-\textbf{(54)}$ we used is
    		    \begin{align*}
    		    	\s{\vec{A}_0}{\bar{\vec{\gamma}_2}}=0
    		    \end{align*}
    		     and comes from the line $\textbf{(46)}$ of \eqref{4conf}.
	    		The new powers of order $11$ are
	    		\small
	    		\begin{align*}
	    		\begin{dmatrix}
	    		10 & 1 & 0\\
	    		9 & 2 & 0\\
	    		8 & 3 & 0
	    		\end{dmatrix}
	    		\begin{dmatrix}
	    		7 & 4 & 0\\
	    		6 & 5 & 0\\
	    		8 & 3 & 1
	    		\end{dmatrix}
	    		\end{align*}
	    		\normalsize
	    		while the ones of order $12$ are
	    		\small
	    		\begin{align*}
	    		\begin{dmatrix}
	    		13 & -1 & 0\\
	    		12 & 0 & 0\\
	    		11 & 1 & 0\\
	    		10 & 2 & 0\\
	    		9 & 3 & 0\\
	    		8 & 4 & 0
	    		\end{dmatrix}
	    		\begin{dmatrix}
	    		7 & 5 & 0\\
	    		6 & 6 & 0\\
	    		9 & 3 & 1\\
	    	    8 & 4 & 1\\
	    		7 & 5 & 1
	    		\end{dmatrix}
	    		\end{align*}
	    		\normalsize
	    		Remembering that
	    		\begin{align}
	    		e^{2\lambda}=\begin{dmatrix}
	    		1 & 3 & 3 & 0\\
	    		2|{A}_1|^2 & 4 & 4 & 0 \\
	    		\beta & 5 & 5 & 0\\
	    		\alpha_1 & 5 & 3 & 0\\
	    		\alpha_2 & 1 & 8 & 0\\
	    		\alpha_3 & 6 & 3 & 0
	    		\end{dmatrix}
	    		\begin{dmatrix}
	    		\alpha_4 & 5 & 4 & 0\\
	    		\alpha_5 & 6 & 4 & 0\\
	    		\alpha_6 & 8 & 2 & 0\\
	    		\alpha_7 & 9 & 1 & 0\\
	    		\zeta_0 & 7 & 3 & 1\\
	    		\bar{\alpha_1} & 3 & 5 & 0
	    			    		\end{dmatrix}
	    		\begin{dmatrix}
	    		\bar{\alpha_2} & 8 & 1 & 0\\
	    		\bar{\alpha_3} & 3 & 6 & 0\\
	    		\bar{\alpha_4} & 4 & 5 & 0\\
	    		\bar{\alpha_5} & 4 & 6 & 0\\
	    		\bar{\alpha_6} & 2 & 8 & 0\\
	    		\bar{\alpha_7} & 1 & 9 & 0\\
	    		\bar{\zeta_0} & 3 & 7 & 1\\
	    		\end{dmatrix}
	    		\end{align}
	    		we see that there exists
	    		\begin{align*}
	    		\alpha_j\in\mathbb{C},\;\, 8\leq j\leq 19,\quad \zeta_j\in\mathbb{C},\;\, 1\leq j\leq 4,\quad \delta\in \R,
	    		\end{align*}
	    		such that
	    		\begin{align}\label{4endmetric}
	    		&e^{2\lambda}=\begin{dmatrix}
	    		1 & 3 & 3 & 0\\
	    		2|{A}_1|^2 & 4 & 4 & 0 \\
	    		\beta & 5 & 5 & 0\\
	    		\delta & 6 & 6 & 0
	    		\end{dmatrix}
	    		+2\,\Re\left\{\begin{dmatrix}
	    		\alpha_1 & 5 & 3 & 0\\
	    		\alpha_2 & 1 & 8 & 0\\
	    		\alpha_3 & 6 & 3 & 0\\
	    		\alpha_4 & 5 & 4 & 0\\
	    		\alpha_5 & 6 & 4 & 0\\
	    		\alpha_6 & 8 & 2 & 0\\
	    		\alpha_7 & 9 & 1 & 0\\
	    		\alpha_8 & 10 & 1 & 0\\
	    		\alpha_9 & 9 & 2 & 0\\
	    		\alpha_{10} & 8 & 3 & 0\\
	    		\alpha_{11} & 7 & 4 & 0\\
	    		\alpha_{12} & 6 & 5 & 0
	    			    		\end{dmatrix}
	    		\begin{dmatrix}
	    		\alpha_{13} & 13 & -1 & 0\\
	    		\alpha_{14} & 12 & 0 & 0\\
	    		\alpha_{15} & 11 & 1 & 0\\
	    		\alpha_{16} & 10 & 2 & 0\\
	    		\alpha_{17} & 9 & 3 & 0\\
	    		\alpha_{18} & 8 & 4 & 0\\
	    		\alpha_{19} & 7 & 5 & 0\\
	    		\zeta_0 & 7 & 3 & 1\\
	    		\zeta_1 & 8 & 3 & 1\\
	    		\zeta_2 & 9 & 3 & 1\\
	    		\zeta_3 & 8 & 4 & 1\\
	    		\zeta_4 & 7 & 5 & 1
	    		\end{dmatrix}\right\}
	    		\end{align}
    		Comparison \TeX\: (left, we do not write the complex conjugate part of the real part arising in \eqref{4endmetric}) and Sage (right)
    		\small
	    		\begin{align*}
	    			\begin{dmatrix}
	    			1 & 3 & 3 & 0\\
	    			2|{A}_1|^2 & 4 & 4 & 0 \\
	    			\beta & 5 & 5 & 0\\
	    			\delta & 6 & 6 & 0\\
	    			\alpha_1 & 5 & 3 & 0\\
	    			\alpha_2 & 1 & 8 & 0\\
	    			\alpha_3 & 6 & 3 & 0\\
	    			\alpha_4 & 5 & 4 & 0\\
	    			\alpha_5 & 6 & 4 & 0\\
	    			\alpha_6 & 8 & 2 & 0\\
	    			\alpha_7 & 9 & 1 & 0\\
	    			\alpha_8 & 10 & 1 & 0\\
	    			\alpha_9 & 9 & 2 & 0\\
	    			\alpha_{10} & 8 & 3 & 0\\
	    			\alpha_{11} & 7 & 4 & 0\\
	    			\alpha_{12} & 6 & 5 & 0\\
	    			\alpha_{13} & 13 & -1 & 0\\
	    			\alpha_{14} & 12 & 0 & 0\\
	    			\alpha_{15} & 11 & 1 & 0\\
	    			\alpha_{16} & 10 & 2 & 0\\
	    			\alpha_{17} & 9 & 3 & 0\\
	    			\alpha_{18} & 8 & 4 & 0\\
	    			\alpha_{19} & 7 & 5 & 0\\
	    			\zeta_0 & 7 & 3 & 1\\
	    			\zeta_1 & 8 & 3 & 1\\
	    			\zeta_2 & 9 & 3 & 1\\
	    			\zeta_3 & 8 & 4 & 1\\
	    			\zeta_4 & 7 & 5 & 1
	    			\end{dmatrix}
	    			\begin{dmatrix}
	    			1 & 3 & 3 & 0 \\
	    			2 \, {\left| A_{1} \right|}^{2} & 4 & 4 & 0 \\
	    			\beta & 5 & 5 & 0 \\
	    			\delta & 6 & 6 & 0 \\
	    			\alpha_{1} & 5 & 3 & 0 \\
	    			\alpha_{2} & 1 & 8 & 0 \\
	    			\alpha_{3} & 6 & 3 & 0 \\
	    			\alpha_{4} & 5 & 4 & 0 \\
	    			\alpha_{5} & 6 & 4 & 0 \\
	    			\alpha_{6} & 8 & 2 & 0 \\
	    			\alpha_{7} & 9 & 1 & 0 \\
	    			\alpha_{8} & 10 & 1 & 0 \\
	    			\alpha_{9} & 9 & 2 & 0 \\
	    			\alpha_{10} & 8 & 3 & 0 \\
	    			\alpha_{11} & 7 & 4 & 0 \\
	    			\alpha_{12} & 6 & 5 & 0 \\
	    			\alpha_{13} & 13 & -1 & 0 \\
	    			\alpha_{14} & 12 & 0 & 0 \\
	    			\alpha_{15} & 11 & 1 & 0 \\
	    			\alpha_{16} & 10 & 2 & 0 \\
	    			\alpha_{17} & 9 & 3 & 0 \\
	    			\alpha_{18} & 8 & 4 & 0 \\
	    			\alpha_{19} & 7 & 5 & 0 \\
	    			\zeta_{0} & 7 & 3 & 1 \\
	    			\zeta_{1} & 8 & 3 & 1 \\
	    			\zeta_{2} & 9 & 3 & 1 \\
	    			\zeta_{3} & 8 & 4 & 1 \\
	    			\zeta_{4} & 7 & 5 & 1 
	    			\end{dmatrix}
	    			\begin{dmatrix}
	    			\overline{\alpha_{1}} & 3 & 5 & 0 \\
	    			\overline{\alpha_{2}} & 8 & 1 & 0 \\
	    			\overline{\alpha_{3}} & 3 & 6 & 0 \\
	    			\overline{\alpha_{4}} & 4 & 5 & 0 \\
	    			\overline{\alpha_{5}} & 4 & 6 & 0 \\
	    			\overline{\alpha_{6}} & 2 & 8 & 0 \\
	    			\overline{\alpha_{7}} & 1 & 9 & 0 \\
	    			\overline{\alpha_{8}} & 1 & 10 & 0 \\
	    			\overline{\alpha_{9}} & 2 & 9 & 0 \\
	    			\overline{\alpha_{10}} & 3 & 8 & 0 \\
	    			\overline{\alpha_{11}} & 4 & 7 & 0 \\
	    			\overline{\alpha_{12}} & 5 & 6 & 0 \\
	    			\overline{\alpha_{13}} & -1 & 13 & 0 \\
	    			\overline{\alpha_{14}} & 0 & 12 & 0 \\
	    			\overline{\alpha_{15}} & 1 & 11 & 0 \\
	    			\overline{\alpha_{16}} & 2 & 10 & 0 \\
	    			\overline{\alpha_{17}} & 3 & 9 & 0 \\
	    			\overline{\alpha_{18}} & 4 & 8 & 0 \\
	    			\overline{\alpha_{19}} & 5 & 7 & 0 \\
	    			\overline{\zeta_{0}} & 3 & 7 & 1 \\
	    			\overline{\zeta_{1}} & 3 & 8 & 1 \\
	    			\overline{\zeta_{2}} & 3 & 9 & 1 \\
	    			\overline{\zeta_{3}} & 4 & 8 & 1 \\
	    			\overline{\zeta_{4}} & 5 & 7 & 1
	    			\end{dmatrix}
	    		\end{align*}
	    		\normalsize
    		Then we have
    		\small
    		\begin{align*}
    			\h_0=\begin{dmatrix}
    			-\frac{1}{4} & E_{1} & -2 & 8 & 0 \\
    			-\frac{2}{9} & E_{2} & -2 & 9 & 0 \\
    			\frac{1}{36} \, \alpha_{2} & C_{1} & -2 & 9 & 0 \\
    			-\frac{1}{5} & E_{4} & -2 & 10 & 0 \\
    			-\frac{1}{5} \, \overline{\alpha_{1}} & E_{1} & -2 & 10 & 0 \\
    			\frac{1}{20} \, \overline{\alpha_{7}} & C_{1} & -2 & 10 & 0 \\
    			-4 \, \alpha_{2}^{2} + 8 \, \overline{\alpha_{13}} & A_{0} & -2 & 10 & 0 \\
    			-\frac{3}{16} & E_{3} & -1 & 8 & 0 \\
    			-\frac{1}{6} & E_{5} & -1 & 9 & 0 \\
    			-\frac{11}{24} \, {\left| A_{1} \right|}^{2} & E_{1} & -1 & 9 & 0 \\
    			\frac{1}{12} \, \alpha_{2} & C_{2} & -1 & 9 & 0 \\
    			-\frac{1}{24} \, \overline{\alpha_{6}} & C_{1} & -1 & 9 & 0 \\
    			6 \, \overline{\alpha_{14}} & A_{0} & -1 & 9 & 0 \\
    			-\frac{1}{4} & C_{1} & 0 & 4 & 0 \\
    			-\frac{1}{5} & B_{1} & 0 & 5 & 0 \\
    			4 \, \alpha_{2} & A_{0} & 0 & 5 & 0 \\
    			-\frac{1}{6} & B_{2} & 0 & 6 & 0 \\
    			-\frac{1}{6} \, \overline{\alpha_{1}} & C_{1} & 0 & 6 & 0 \\
    			4 \, \overline{\alpha_{7}} & A_{0} & 0 & 6 & 0 \\
    			-\frac{1}{7} & B_{4} & 0 & 7 & 0 \\
    			\frac{5}{14} \, \alpha_{2} & \overline{C_{1}} & 0 & 7 & 0 \\
    			-\frac{1}{7} \, \overline{\alpha_{1}} & B_{1} & 0 & 7 & 0 \\
    			-\frac{1}{7} \, \overline{\alpha_{3}} & C_{1} & 0 & 7 & 0 \\
    			-4 \, \alpha_{2} \overline{\alpha_{1}} + 4 \, \overline{\alpha_{8}} & A_{0} & 0 & 7 & 0 \\
    			-\frac{1}{8} & B_{6} & 0 & 8 & 0 \\
    			-\frac{1}{4} \, \alpha_{1} & E_{1} & 0 & 8 & 0 \\
    			\frac{5}{24} \, \alpha_{2} & \overline{C_{2}} & 0 & 8 & 0 \\
    			-\frac{1}{8} \, \overline{\alpha_{1}} & B_{2} & 0 & 8 & 0 \\
    			-\frac{1}{8} \, \overline{\alpha_{3}} & B_{1} & 0 & 8 & 0 \\
    			\frac{3}{8} \, \overline{\alpha_{7}} & \overline{C_{1}} & 0 & 8 & 0 
    			\end{dmatrix}
    			\begin{dmatrix}
    			-\frac{3}{128} \, \overline{\zeta_{0}} & C_{1} & 0 & 8 & 0 \\
    			-4 \, \alpha_{2} \overline{\alpha_{3}} - 4 \, \overline{\alpha_{1}} \overline{\alpha_{7}} + 4 \, \overline{\alpha_{15}} & A_{0} & 0 & 8 & 0 \\
    			-\frac{1}{8} \, \overline{\zeta_{0}} & C_{1} & 0 & 8 & 1 \\
    			-\frac{1}{8} & C_{2} & 1 & 4 & 0 \\
    			-\frac{1}{10} & B_{3} & 1 & 5 & 0 \\
    			-\frac{9}{20} \, {\left| A_{1} \right|}^{2} & C_{1} & 1 & 5 & 0 \\
    			4 \, \alpha_{2} & A_{1} & 1 & 5 & 0 \\
    			2 \, \overline{\alpha_{6}} & A_{0} & 1 & 5 & 0 \\
    			-\frac{1}{12} & B_{5} & 1 & 6 & 0 \\
    			-\frac{11}{30} \, {\left| A_{1} \right|}^{2} & B_{1} & 1 & 6 & 0 \\
    			-\frac{1}{12} \, \overline{\alpha_{1}} & C_{2} & 1 & 6 & 0 \\
    			-\frac{5}{24} \, \overline{\alpha_{4}} & C_{1} & 1 & 6 & 0 \\
    			4 \, \overline{\alpha_{7}} & A_{1} & 1 & 6 & 0 \\
    			-4 \, \alpha_{2} {\left| A_{1} \right|}^{2} + 2 \, \overline{\alpha_{9}} & A_{0} & 1 & 6 & 0 \\
    			-\frac{1}{14} & B_{7} & 1 & 7 & 0 \\
    			-\frac{13}{42} \, {\left| A_{1} \right|}^{2} & B_{2} & 1 & 7 & 0 \\
    			\frac{3}{7} \, \alpha_{2} & \overline{B_{1}} & 1 & 7 & 0 \\
    			-\frac{1}{14} \, \overline{\alpha_{1}} & B_{3} & 1 & 7 & 0 \\
    			-\frac{1}{14} \, \overline{\alpha_{3}} & C_{2} & 1 & 7 & 0 \\
    			-\frac{6}{35} \, \overline{\alpha_{4}} & B_{1} & 1 & 7 & 0 \\
    			\frac{1}{12} \, {\left| A_{1} \right|}^{2} \overline{\alpha_{1}} - \frac{11}{56} \, \overline{\alpha_{5}} & C_{1} & 1 & 7 & 0 \\
    			\frac{5}{28} \, \overline{\alpha_{6}} & \overline{C_{1}} & 1 & 7 & 0 \\
    			-4 \, \alpha_{2} \overline{\alpha_{1}} + 4 \, \overline{\alpha_{8}} & A_{1} & 1 & 7 & 0 \\
    			-4 \, {\left| A_{1} \right|}^{2} \overline{\alpha_{7}} - 2 \, \alpha_{2} \overline{\alpha_{4}} - 2 \, \overline{\alpha_{1}} \overline{\alpha_{6}} + 2 \, \overline{\alpha_{16}} & A_{0} & 1 & 7 & 0 \\
    			\frac{1}{8} & \gamma_{1} & 2 & 4 & 0 \\
    			-\frac{1}{4} \, \alpha_{1} & C_{1} & 2 & 4 & 0 \\
    			-\overline{\zeta_{0}} & A_{0} & 2 & 4 & 0 \\
    			-\frac{1}{4} \, {\left| A_{1} \right|}^{2} & C_{2} & 2 & 5 & 0 \\
    			\frac{1}{20} & \overline{\gamma_{2}} & 2 & 5 & 0 \\
    			-\frac{1}{5} \, \alpha_{1} & B_{1} & 2 & 5 & 0 
    			    			\end{dmatrix}
    			\end{align*}
    			\begin{align*}
    			\begin{dmatrix}
    			-\frac{1}{4} \, \alpha_{4} & C_{1} & 2 & 5 & 0 \\
    			4 \, \alpha_{2} & A_{2} & 2 & 5 & 0 \\
    			2 \, \overline{\alpha_{6}} & A_{1} & 2 & 5 & 0 \\
    			-\overline{\zeta_{1}} & A_{0} & 2 & 5 & 0 \\
    			-\frac{1}{5} \, {\left| A_{1} \right|}^{2} & B_{3} & 2 & 6 & 0 \\
    			\frac{1}{24} & \overline{\gamma_{3}} & 2 & 6 & 0 \\
    			\frac{1}{10} \, {\left| A_{1} \right|}^{4} + \frac{1}{12} \, \alpha_{1} \overline{\alpha_{1}} - \frac{1}{4} \, \beta & C_{1} & 2 & 6 & 0 \\
    			-\frac{1}{6} \, \alpha_{1} & B_{2} & 2 & 6 & 0 \\
    			-\frac{1}{5} \, \alpha_{4} & B_{1} & 2 & 6 & 0 \\
    			\frac{1}{12} \, \overline{\alpha_{1}} & \gamma_{1} & 2 & 6 & 0 \\
    			-\frac{1}{8} \, \overline{\alpha_{4}} & C_{2} & 2 & 6 & 0 \\
    			-\frac{1}{12} \, \overline{\zeta_{0}} & \overline{C_{1}} & 2 & 6 & 0 \\
    			4 \, \overline{\alpha_{7}} & A_{2} & 2 & 6 & 0 \\
    			-4 \, \alpha_{2} {\left| A_{1} \right|}^{2} + 2 \, \overline{\alpha_{9}} & A_{1} & 2 & 6 & 0 \\
    			-12 \, \overline{\alpha_{1}}^{3} + \overline{\alpha_{1}} \overline{\zeta_{0}} - \overline{\zeta_{2}} & A_{0} & 2 & 6 & 0 \\
    			2 & A_{1} & 3 & 0 & 0 \\
    			-4 \, {\left| A_{1} \right|}^{2} & A_{0} & 3 & 1 & 0 \\
    			\frac{1}{4} & \overline{B_{1}} & 3 & 2 & 0 \\
    			-2 \, \overline{\alpha_{4}} & A_{0} & 3 & 2 & 0 \\
    			-\frac{1}{6} \, {\left| A_{1} \right|}^{2} & \overline{C_{1}} & 3 & 3 & 0 \\
    			\frac{1}{6} & \overline{B_{3}} & 3 & 3 & 0 \\
    			4 \, {\left| A_{1} \right|}^{2} \overline{\alpha_{1}} - 2 \, \overline{\alpha_{5}} & A_{0} & 3 & 3 & 0 \\
    			-\frac{1}{12} \, {\left| A_{1} \right|}^{2} & \overline{C_{2}} & 3 & 4 & 0 \\
    			\frac{3}{64} & \gamma_{2} & 3 & 4 & 0 \\
    			\frac{1}{8} & C_{4} & 3 & 4 & 0 \\
    			-\frac{1}{8} \, \alpha_{1} & C_{2} & 3 & 4 & 0 \\
    			-\frac{1}{4} \, \alpha_{3} & C_{1} & 3 & 4 & 0 \\
    			\frac{1}{8} \, \overline{\alpha_{1}} & \overline{B_{1}} & 3 & 4 & 0 \\
    			-\frac{1}{8} \, \overline{\alpha_{4}} & \overline{C_{1}} & 3 & 4 & 0 \\
    			4 \, {\left| A_{1} \right|}^{2} \overline{\alpha_{3}} + 2 \, \overline{\alpha_{1}} \overline{\alpha_{4}} - 2 \, \overline{\alpha_{11}} & A_{0} & 3 & 4 & 0 \\
    			-\overline{\zeta_{0}} & A_{1} & 3 & 4 & 0 
    			    			\end{dmatrix}
    			\end{align*}
    			\begin{align*}
    			\begin{dmatrix}
    			\frac{1}{8} & \gamma_{2} & 3 & 4 & 1 \\
    			-\frac{1}{10} \, {\left| A_{1} \right|}^{2} & C_{3} & 3 & 5 & 0 \\
    			\frac{2}{25} & \gamma_{4} & 3 & 5 & 0 \\
    			\frac{1}{5} & C_{6} & 3 & 5 & 0 \\
    			\frac{89}{400} \, {\left| A_{1} \right|}^{2} & \gamma_{1} & 3 & 5 & 0 \\
    			-\frac{1}{10} \, \alpha_{1} & B_{3} & 3 & 5 & 0 \\
    			-\frac{1}{5} \, \alpha_{3} & B_{1} & 3 & 5 & 0 \\
    			-\frac{3}{20} \, \alpha_{4} & C_{2} & 3 & 5 & 0 \\
    			\frac{1}{10} \, \alpha_{1} {\left| A_{1} \right|}^{2} - \frac{11}{40} \, \alpha_{5} & C_{1} & 3 & 5 & 0 \\
    			\frac{1}{10} \, \overline{\alpha_{1}} & \overline{B_{3}} & 3 & 5 & 0 \\
    			\frac{1}{10} \, \overline{\alpha_{3}} & \overline{B_{1}} & 3 & 5 & 0 \\
    			-\frac{1}{15} \, \overline{\alpha_{4}} & \overline{C_{2}} & 3 & 5 & 0 \\
    			\frac{1}{4} \, {\left| A_{1} \right|}^{2} \overline{\alpha_{1}} - \frac{3}{20} \, \overline{\alpha_{5}} & \overline{C_{1}} & 3 & 5 & 0 \\
    			4 \, \alpha_{2} & A_{3} & 3 & 5 & 0 \\
    			2 \, \overline{\alpha_{6}} & A_{2} & 3 & 5 & 0 \\
    			-\overline{\zeta_{1}} & A_{1} & 3 & 5 & 0 \\
    			-76 \, {\left| A_{1} \right|}^{2} \overline{\alpha_{1}}^{2} + 2 \, {\left| A_{1} \right|}^{2} \overline{\zeta_{0}} + 2 \, \alpha_{2} \alpha_{3} + 2 \, \overline{\alpha_{3}} \overline{\alpha_{4}} + 2 \, \overline{\alpha_{1}} \overline{\alpha_{5}} + 2 \, \alpha_{1} \overline{\alpha_{6}} - 2 \, \overline{\alpha_{18}} - \overline{\zeta_{3}} & A_{0} & 3 & 5 & 0 \\
    			-\frac{1}{10} \, {\left| A_{1} \right|}^{2} & \gamma_{1} & 3 & 5 & 1 \\
    			\frac{1}{5} & \gamma_{4} & 3 & 5 & 1 \\
    			4 \, {\left| A_{1} \right|}^{2} \overline{\zeta_{0}} - 2 \, \overline{\zeta_{3}} & A_{0} & 3 & 5 & 1 \\
    			4 & A_{2} & 4 & 0 & 0 \\
    			-4 \, \alpha_{1} & A_{0} & 4 & 0 & 0 \\
    			-4 \, {\left| A_{1} \right|}^{2} & A_{1} & 4 & 1 & 0 \\
    			-4 \, \alpha_{4} & A_{0} & 4 & 1 & 0 \\
    			\frac{1}{2} & \overline{B_{2}} & 4 & 2 & 0 \\
    			8 \, {\left| A_{1} \right|}^{4} + 4 \, \alpha_{1} \overline{\alpha_{1}} - 4 \, \beta & A_{0} & 4 & 2 & 0 \\
    			-2 \, \overline{\alpha_{4}} & A_{1} & 4 & 2 & 0 \\
    			\frac{1}{3} & \overline{B_{5}} & 4 & 3 & 0 \\
    			\frac{1}{6} \, {\left| A_{1} \right|}^{2} & \overline{B_{1}} & 4 & 3 & 0 \\
    			-\frac{1}{6} \, \alpha_{4} & \overline{C_{1}} & 4 & 3 & 0 \\
    			4 \, {\left| A_{1} \right|}^{2} \overline{\alpha_{1}} - 2 \, \overline{\alpha_{5}} & A_{1} & 4 & 3 & 0 
    			    			\end{dmatrix}
    			\end{align*}
    			\begin{align*}
    			\begin{dmatrix}
    			8 \, {\left| A_{1} \right|}^{2} \overline{\alpha_{4}} + 4 \, \alpha_{4} \overline{\alpha_{1}} + 4 \, \alpha_{1} \overline{\alpha_{3}} - 4 \, \overline{\alpha_{12}} & A_{0} & 4 & 3 & 0 \\
    			\frac{1}{32} & \gamma_{3} & 4 & 4 & 0 \\
    			\frac{1}{4} & C_{5} & 4 & 4 & 0 \\
    			\frac{1}{6} \, {\left| A_{1} \right|}^{2} & \overline{B_{3}} & 4 & 4 & 0 \\
    			\frac{1}{3} \, {\left| A_{1} \right|}^{4} + \frac{1}{4} \, \alpha_{1} \overline{\alpha_{1}} - \frac{1}{4} \, \beta & \overline{C_{1}} & 4 & 4 & 0 \\
    			\frac{1}{8} \, \alpha_{1} & \gamma_{1} & 4 & 4 & 0 \\
    			-\frac{1}{8} \, \alpha_{3} & C_{2} & 4 & 4 & 0 \\
    			-\frac{1}{12} \, \alpha_{4} & \overline{C_{2}} & 4 & 4 & 0 \\
    			-\frac{1}{32} \, \zeta_{0} & C_{1} & 4 & 4 & 0 \\
    			\frac{1}{4} \, \overline{\alpha_{1}} & \overline{B_{2}} & 4 & 4 & 0 \\
    			4 \, {\left| A_{1} \right|}^{2} \overline{\alpha_{3}} + 2 \, \overline{\alpha_{1}} \overline{\alpha_{4}} - 2 \, \overline{\alpha_{11}} & A_{1} & 4 & 4 & 0 \\
    			-160 \, {\left| A_{1} \right|}^{4} \overline{\alpha_{1}} - 40 \, \alpha_{1} \overline{\alpha_{1}}^{2} + 8 \, {\left| A_{1} \right|}^{2} \overline{\alpha_{5}} + 4 \, \beta \overline{\alpha_{1}} + 4 \, \alpha_{4} \overline{\alpha_{3}} + 2 \, \overline{\alpha_{4}}^{2} + \alpha_{1} \overline{\zeta_{0}} - 4 \, \overline{\alpha_{19}} - \overline{\zeta_{4}} & A_{0} & 4 & 4 & 0 \\
    			-\overline{\zeta_{0}} & A_{2} & 4 & 4 & 0 \\
    			\frac{1}{4} & \gamma_{3} & 4 & 4 & 1 \\
    			-\frac{1}{4} \, \zeta_{0} & C_{1} & 4 & 4 & 1 \\
    			4 \, \alpha_{1} \overline{\zeta_{0}} - 4 \, \overline{\zeta_{4}} & A_{0} & 4 & 4 & 1 \\
    			6 & A_{3} & 5 & 0 & 0 \\
    			-4 \, \alpha_{1} & A_{1} & 5 & 0 & 0 \\
    			-6 \, \alpha_{3} & A_{0} & 5 & 0 & 0 \\
    			-4 \, {\left| A_{1} \right|}^{2} & A_{2} & 5 & 1 & 0 \\
    			-4 \, \alpha_{4} & A_{1} & 5 & 1 & 0 \\
    			12 \, \alpha_{1} {\left| A_{1} \right|}^{2} - 6 \, \alpha_{5} & A_{0} & 5 & 1 & 0 \\
    			\frac{3}{4} & \overline{B_{4}} & 5 & 2 & 0 \\
    			\frac{1}{4} \, \alpha_{1} & \overline{B_{1}} & 5 & 2 & 0 \\
    			\frac{1}{8} \, \overline{\alpha_{2}} & C_{1} & 5 & 2 & 0 \\
    			8 \, {\left| A_{1} \right|}^{4} + 4 \, \alpha_{1} \overline{\alpha_{1}} - 4 \, \beta & A_{1} & 5 & 2 & 0 \\
    			12 \, \alpha_{4} {\left| A_{1} \right|}^{2} + 6 \, \alpha_{3} \overline{\alpha_{1}} + 6 \, \alpha_{1} \overline{\alpha_{4}} - 6 \, \alpha_{12} & A_{0} & 5 & 2 & 0 \\
    			-2 \, \overline{\alpha_{4}} & A_{2} & 5 & 2 & 0 \\
    			\frac{1}{2} & \overline{B_{7}} & 5 & 3 & 0 \\
    			\frac{1}{2} \, {\left| A_{1} \right|}^{2} & \overline{B_{2}} & 5 & 3 & 0 \\
    			\frac{1}{6} \, \alpha_{1} & \overline{B_{3}} & 5 & 3 & 0 
    			    			\end{dmatrix}
    			\end{align*}
    			\begin{align*}
    			\begin{dmatrix}
    			\frac{1}{3} \, \alpha_{1} {\left| A_{1} \right|}^{2} - \frac{1}{4} \, \alpha_{5} & \overline{C_{1}} & 5 & 3 & 0 \\
    			-\frac{1}{8} \, \alpha_{6} & C_{1} & 5 & 3 & 0 \\
    			-112 \, {\left| A_{1} \right|}^{6} - 168 \, \alpha_{1} {\left| A_{1} \right|}^{2} \overline{\alpha_{1}} + 12 \, \beta {\left| A_{1} \right|}^{2} + 6 \, \alpha_{5} \overline{\alpha_{1}} + 6 \, \alpha_{2} \overline{\alpha_{2}} + 6 \, \alpha_{3} \overline{\alpha_{3}} + 6 \, \alpha_{4} \overline{\alpha_{4}} + 6 \, \alpha_{1} \overline{\alpha_{5}} - 6 \, \delta & A_{0} & 5 & 3 & 0 \\
    			4 \, {\left| A_{1} \right|}^{2} \overline{\alpha_{1}} - 2 \, \overline{\alpha_{5}} & A_{2} & 5 & 3 & 0 \\
    			8 \, {\left| A_{1} \right|}^{2} \overline{\alpha_{4}} + 4 \, \alpha_{4} \overline{\alpha_{1}} + 4 \, \alpha_{1} \overline{\alpha_{3}} - 4 \, \overline{\alpha_{12}} & A_{1} & 5 & 3 & 0 \\
    			\frac{1}{2} & \overline{E_{1}} & 6 & 0 & 0 \\
    			8 & A_{4} & 6 & 0 & 0 \\
    			-4 \, \alpha_{1} & A_{2} & 6 & 0 & 0 \\
    			-6 \, \alpha_{3} & A_{1} & 6 & 0 & 0 \\
    			4 \, \alpha_{1}^{2} - \zeta_{0} & A_{0} & 6 & 0 & 0 \\
    			4 & \overline{E_{1}} & 6 & 0 & 1 \\
    			-8 \, \zeta_{0} & A_{0} & 6 & 0 & 1 \\
    			-4 \, {\left| A_{1} \right|}^{2} & A_{3} & 6 & 1 & 0 \\
    			2 & \overline{E_{3}} & 6 & 1 & 0 \\
    			-4 \, \alpha_{4} & A_{2} & 6 & 1 & 0 \\
    			12 \, \alpha_{1} {\left| A_{1} \right|}^{2} - 6 \, \alpha_{5} & A_{1} & 6 & 1 & 0 \\
    			16 \, \alpha_{3} {\left| A_{1} \right|}^{2} + 8 \, \alpha_{1} \alpha_{4} - 8 \, \alpha_{11} & A_{0} & 6 & 1 & 0 \\
    			1 & \overline{B_{6}} & 6 & 2 & 0 \\
    			\frac{1}{2} \, \alpha_{1} & \overline{B_{2}} & 6 & 2 & 0 \\
    			\frac{1}{4} \, \alpha_{3} & \overline{B_{1}} & 6 & 2 & 0 \\
    			\frac{1}{4} \, \alpha_{7} & C_{1} & 6 & 2 & 0 \\
    			-\frac{1}{4} \, \zeta_{0} & \overline{C_{1}} & 6 & 2 & 0 \\
    			\overline{\alpha_{1}} & \overline{E_{1}} & 6 & 2 & 0 \\
    			\frac{3}{8} \, \overline{\alpha_{2}} & C_{2} & 6 & 2 & 0 \\
    			8 \, {\left| A_{1} \right|}^{4} + 4 \, \alpha_{1} \overline{\alpha_{1}} - 4 \, \beta & A_{2} & 6 & 2 & 0 \\
    			12 \, \alpha_{4} {\left| A_{1} \right|}^{2} + 6 \, \alpha_{3} \overline{\alpha_{1}} + 6 \, \alpha_{1} \overline{\alpha_{4}} - 6 \, \alpha_{12} & A_{1} & 6 & 2 & 0 \\
    			-176 \, \alpha_{1} {\left| A_{1} \right|}^{4} + 16 \, \alpha_{5} {\left| A_{1} \right|}^{2} - 44 \, \alpha_{1}^{2} \overline{\alpha_{1}} + 4 \, \alpha_{4}^{2} + 8 \, \alpha_{1} \beta + \zeta_{0} \overline{\alpha_{1}} + 8 \, \alpha_{3} \overline{\alpha_{4}} - 8 \, \alpha_{19} - \zeta_{4} & A_{0} & 6 & 2 & 0 \\
    			-2 \, \overline{\alpha_{4}} & A_{3} & 6 & 2 & 0 \\
    			8 \, \zeta_{0} \overline{\alpha_{1}} - 8 \, \zeta_{4} & A_{0} & 6 & 2 & 1 \\
    			-10 \, \overline{\alpha_{2}} & A_{0} & 7 & -2 & 0 \\
    			-10 \, \alpha_{6} & A_{0} & 7 & -1 & 0 
    			    			\end{dmatrix}
    			\end{align*}
    			\begin{align*}
    			\begin{dmatrix}
    			\frac{1}{2} & \overline{E_{2}} & 7 & 0 & 0 \\
    			10 & A_{5} & 7 & 0 & 0 \\
    			-\frac{3}{4} \, \overline{\alpha_{2}} & \overline{C_{1}} & 7 & 0 & 0 \\
    			-4 \, \alpha_{1} & A_{3} & 7 & 0 & 0 \\
    			-6 \, \alpha_{3} & A_{2} & 7 & 0 & 0 \\
    			10 \, \alpha_{1} \alpha_{3} + 10 \, \overline{\alpha_{1}} \overline{\alpha_{2}} - 10 \, \alpha_{10} - \zeta_{1} & A_{0} & 7 & 0 & 0 \\
    			4 \, \alpha_{1}^{2} - \zeta_{0} & A_{1} & 7 & 0 & 0 \\
    			5 & \overline{E_{2}} & 7 & 0 & 1 \\
    			5 \, \overline{\alpha_{2}} & \overline{C_{1}} & 7 & 0 & 1 \\
    			-8 \, \zeta_{0} & A_{1} & 7 & 0 & 1 \\
    			-10 \, \zeta_{1} & A_{0} & 7 & 0 & 1 \\
    			\frac{5}{2} & \overline{E_{5}} & 7 & 1 & 0 \\
    			5 \, {\left| A_{1} \right|}^{2} & \overline{E_{1}} & 7 & 1 & 0 \\
    			\frac{5}{4} \, \alpha_{6} & \overline{C_{1}} & 7 & 1 & 0 \\
    			\frac{5}{3} \, \overline{\alpha_{2}} & \overline{C_{2}} & 7 & 1 & 0 \\
    			-4 \, {\left| A_{1} \right|}^{2} & A_{4} & 7 & 1 & 0 \\
    			-4 \, \alpha_{4} & A_{3} & 7 & 1 & 0 \\
    			12 \, \alpha_{1} {\left| A_{1} \right|}^{2} - 6 \, \alpha_{5} & A_{2} & 7 & 1 & 0 \\
    			16 \, \alpha_{3} {\left| A_{1} \right|}^{2} + 8 \, \alpha_{1} \alpha_{4} - 8 \, \alpha_{11} & A_{1} & 7 & 1 & 0 \\
    			-92 \, \alpha_{1}^{2} {\left| A_{1} \right|}^{2} + 2 \, \zeta_{0} {\left| A_{1} \right|}^{2} + 10 \, \alpha_{3} \alpha_{4} + 10 \, \alpha_{1} \alpha_{5} + 10 \, \alpha_{6} \overline{\alpha_{1}} + 10 \, \overline{\alpha_{2}} \overline{\alpha_{3}} - 10 \, \alpha_{18} - \zeta_{3} & A_{0} & 7 & 1 & 0 \\
    			-2 \, {\left| A_{1} \right|}^{2} & \overline{E_{1}} & 7 & 1 & 1 \\
    			20 \, \zeta_{0} {\left| A_{1} \right|}^{2} - 10 \, \zeta_{3} & A_{0} & 7 & 1 & 1 \\
    			-12 \, \alpha_{7} & A_{0} & 8 & -2 & 0 \\
    			-10 \, \overline{\alpha_{2}} & A_{1} & 8 & -2 & 0 \\
    			-10 \, \alpha_{6} & A_{1} & 8 & -1 & 0 \\
    			24 \, {\left| A_{1} \right|}^{2} \overline{\alpha_{2}} - 12 \, \alpha_{9} & A_{0} & 8 & -1 & 0 \\
    			\frac{1}{2} & \overline{E_{4}} & 8 & 0 & 0 \\
    			12 & A_{6} & 8 & 0 & 0 \\
    			\frac{1}{2} \, \alpha_{1} & \overline{E_{1}} & 8 & 0 & 0 \\
    			-\alpha_{7} & \overline{C_{1}} & 8 & 0 & 0 \\
    			-\frac{3}{4} \, \overline{\alpha_{2}} & \overline{B_{1}} & 8 & 0 & 0 
    			    			\end{dmatrix}
    			\end{align*}
    			\begin{align*}
    			\begin{dmatrix}
    			-4 \, \alpha_{1} & A_{4} & 8 & 0 & 0 \\
    			-6 \, \alpha_{3} & A_{3} & 8 & 0 & 0 \\
    			10 \, \alpha_{1} \alpha_{3} + 10 \, \overline{\alpha_{1}} \overline{\alpha_{2}} - 10 \, \alpha_{10} - \zeta_{1} & A_{1} & 8 & 0 & 0 \\
    			-16 \, \alpha_{1}^{3} + 24 \, \alpha_{6} {\left| A_{1} \right|}^{2} + 6 \, \alpha_{3}^{2} + \alpha_{1} \zeta_{0} + 12 \, \alpha_{7} \overline{\alpha_{1}} + 12 \, \overline{\alpha_{2}} \overline{\alpha_{4}} - 12 \, \alpha_{17} - \zeta_{2} & A_{0} & 8 & 0 & 0 \\
    			4 \, \alpha_{1}^{2} - \zeta_{0} & A_{2} & 8 & 0 & 0 \\
    			6 & \overline{E_{4}} & 8 & 0 & 1 \\
    			4 \, \alpha_{1} & \overline{E_{1}} & 8 & 0 & 1 \\
    			6 \, \alpha_{7} & \overline{C_{1}} & 8 & 0 & 1 \\
    			6 \, \overline{\alpha_{2}} & \overline{B_{1}} & 8 & 0 & 1 \\
    			-8 \, \zeta_{0} & A_{2} & 8 & 0 & 1 \\
    			-10 \, \zeta_{1} & A_{1} & 8 & 0 & 1 \\
    			12 \, \alpha_{1} \zeta_{0} - 12 \, \zeta_{2} & A_{0} & 8 & 0 & 1 \\
    			-12 \, \alpha_{7} & A_{1} & 9 & -2 & 0 \\
    			14 \, \alpha_{1} \overline{\alpha_{2}} - 14 \, \alpha_{8} & A_{0} & 9 & -2 & 0 \\
    			-10 \, \overline{\alpha_{2}} & A_{2} & 9 & -2 & 0 \\
    			-10 \, \alpha_{6} & A_{2} & 9 & -1 & 0 \\
    			24 \, {\left| A_{1} \right|}^{2} \overline{\alpha_{2}} - 12 \, \alpha_{9} & A_{1} & 9 & -1 & 0 \\
    			28 \, \alpha_{7} {\left| A_{1} \right|}^{2} + 14 \, \alpha_{1} \alpha_{6} + 14 \, \alpha_{4} \overline{\alpha_{2}} - 14 \, \alpha_{16} & A_{0} & 9 & -1 & 0 \\
    			-12 \, \alpha_{7} & A_{2} & 10 & -2 & 0 \\
    			14 \, \alpha_{1} \overline{\alpha_{2}} - 14 \, \alpha_{8} & A_{1} & 10 & -2 & 0 \\
    			16 \, \alpha_{1} \alpha_{7} + 16 \, \alpha_{3} \overline{\alpha_{2}} - 16 \, \alpha_{15} & A_{0} & 10 & -2 & 0 \\
    			-10 \, \overline{\alpha_{2}} & A_{3} & 10 & -2 & 0 \\
    			-18 \, \alpha_{14} & A_{0} & 11 & -3 & 0 \\
    			10 \, \overline{\alpha_{2}}^{2} - 20 \, \alpha_{13} & A_{0} & 12 & -4 & 0
    			\end{dmatrix}
    		\end{align*}
    		
    		Now, we have
    		\footnotesize
    		\begin{align*}
    			&g^{-1}\otimes Q(\h_0)=\\
    			&\begin{dmatrix}
    			-16 \, A_{0} C_{1} \alpha_{1} + 16 \, A_{2} C_{1} + 2 \, A_{1} C_{2} & 0 & 0 & 0 \\
    			-30 \, A_{0} C_{1} \alpha_{3} + 8 \, A_{0} A_{1} \overline{\zeta_{0}} + 30 \, A_{3} C_{1} + 6 \, A_{2} C_{2} - 6 \, {\left(4 \, A_{1} C_{1} + A_{0} C_{2}\right)} \alpha_{1} - A_{1} \gamma_{1} & 1 & 0 & 0 \\
    			\omega_1 & 2 & 0 & 0 \\
    			\omega_2 & 3 & 0 & 0 \\
    			20 \, A_{1} E_{1} & -3 & 4 & 0 \\
    			-48 \, A_{0} E_{1} \alpha_{1} + 48 \, A_{2} E_{1} + 12 \, A_{1} E_{3} & -2 & 4 & 0 \\
    			\frac{9}{80} \, C_{1}^{2} {\left| A_{1} \right|}^{2} - 35 \, A_{0} E_{1} {\left| A_{1} \right|}^{2} - \frac{1}{2} \, A_{0} C_{1} \overline{\alpha_{6}} + \frac{1}{40} \, B_{3} C_{1} - \frac{1}{40} \, B_{1} C_{2} + 20 \, A_{1} E_{2} - \frac{1}{2} \, {\left(19 \, A_{1} C_{1} - A_{0} C_{2}\right)} \alpha_{2} & -3 & 5 & 0 \\
    		    \omega_3 & -3 & 6 & 0 \\
    			\omega_4 & -2 & 5 & 0 \\
    			-\frac{3}{8} \, C_{2} E_{1} - \frac{3}{16} \, C_{1} E_{3} & -5 & 8 & 0 \\
    			\omega_5 & -1 & 4 & 0 \\
    			\omega_6 & 0 & 2 & 0 \\
    			\omega_7 & 1 & 2 & 0 \\
    			\omega_8 & 0 & 3 & 0 \\
    		    \omega_9 & 3 & 0 & 1 \\
    		    \omega_{10} & 0 & 1 & 0 \\
    			\omega_{11} & 1 & 1 & 0 \\
    			\omega_{12} & 2 & 1 & 0 \\
    			\omega_{13} & 4 & -1 & 0 \\
    			\omega_{14} & 5 & -2 & 0 \\
    			\omega_{15} & 4 & -2 & 0 \\
    			-\frac{1}{2} \, C_{1} E_{1} & -6 & 8 & 0 \\
    			-\frac{3}{10} \, B_{1} E_{1} - \frac{5}{9} \, C_{1} E_{2} + \frac{1}{72} \, {\left(5 \, C_{1}^{2} + 432 \, A_{0} E_{1}\right)} \alpha_{2} & -6 & 9 & 0 \\
    			-48 \, A_{0} C_{1} \zeta_{0} + 24 \, C_{1} \overline{E_{1}} & 2 & 0 & 1 \\
    			-480 \, A_{0} \alpha_{2} \overline{E_{1}} + 48 \, {\left(20 \, A_{0}^{2} \alpha_{2} - A_{0} B_{1}\right)} \zeta_{0} + 24 \, B_{1} \overline{E_{1}} & 2 & 1 & 1 \\
    			\omega_{16} & 3 & -2 & 0 \\
    			192 \, A_{0}^{2} \alpha_{2} {\left| A_{1} \right|}^{2} - \frac{48}{5} \, A_{0} B_{1} {\left| A_{1} \right|}^{2} - 3 \, A_{0} C_{1} \overline{\alpha_{4}} - 144 \, A_{0} A_{1} \overline{\alpha_{7}} + 6 \, A_{1} B_{2} + \frac{3}{8} \, C_{1} \overline{B_{1}} & -1 & 2 & 0 \\
    			\omega_{17} & -1 & 3 & 0 \\
    			-96 \, A_{0} A_{1} \zeta_{4} + 48 \, A_{0} \overline{E_{1}} \overline{\alpha_{4}} + 12 \, {\left(8 \, A_{0} A_{1} \overline{\alpha_{1}} - 8 \, A_{0}^{2} \overline{\alpha_{4}} + A_{0} \overline{B_{1}}\right)} \zeta_{0} - 6 \, \overline{B_{1}} \overline{E_{1}} & 5 & -2 & 1 \\
    		    \omega_{18} & 3 & -1 & 0 \\
    			\omega_{19} & 6 & -3 & 0 \\
    			\omega_{20} & 7 & -4 & 0 
    			    			\end{dmatrix}
    			\end{align*}
    			\begin{align*}
    			\begin{dmatrix}
    			\omega_{21} & 6 & -4 & 0 \\
    			-24 \, A_{0} A_{1} \zeta_{0} + 12 \, A_{1} \overline{E_{1}} & 5 & -4 & 0 \\
    			\omega_{22} & 5 & -3 & 0 \\
    			\omega_{23} & 9 & -6 & 0 \\
    			-16 \, A_{0}^{2} \alpha_{1} {\left| A_{1} \right|}^{2} - 8 \, A_{1}^{2} {\left| A_{1} \right|}^{2} + 16 \, A_{0} A_{2} {\left| A_{1} \right|}^{2} - 8 \, A_{0} A_{1} \alpha_{4} & 3 & -3 & 0 \\
    			16 \, A_{0} A_{1} \alpha_{1} {\left| A_{1} \right|}^{2} - 48 \, A_{0}^{2} \alpha_{3} {\left| A_{1} \right|}^{2} - 16 \, A_{1} A_{2} {\left| A_{1} \right|}^{2} + 48 \, A_{0} A_{3} {\left| A_{1} \right|}^{2} - 16 \, A_{1}^{2} \alpha_{4} - 24 \, A_{0} A_{1} \alpha_{5} & 4 & -3 & 0 \\
    			-9 \, A_{0} C_{1} {\left| A_{1} \right|}^{2} - 120 \, A_{0} A_{1} \alpha_{2} + 6 \, A_{1} B_{1} & -1 & 1 & 0 \\
    			-96 \, A_{0}^{2} \zeta_{0} {\left| A_{1} \right|}^{2} + 48 \, A_{0} {\left| A_{1} \right|}^{2} \overline{E_{1}} & 5 & -3 & 1 \\
    			\omega_{24} & 6 & -3 & 1 \\
    			-720 \, A_{0}^{2} \alpha_{7} {\left| A_{1} \right|}^{2} - 1520 \, A_{0} A_{1} {\left| A_{1} \right|}^{2} \overline{\alpha_{2}} - 360 \, A_{0}^{2} \alpha_{4} \overline{\alpha_{2}} + 120 \, A_{0} A_{1} \alpha_{9} - 20 \, {\left(6 \, A_{0}^{2} \alpha_{1} - 5 \, A_{1}^{2} - 6 \, A_{0} A_{2}\right)} \alpha_{6} & 7 & -5 & 0 \\
    			\omega_{25} & 8 & -5 & 0 \\
    			-480 \, A_{0}^{2} {\left| A_{1} \right|}^{2} \overline{\alpha_{2}} + 80 \, A_{0} A_{1} \alpha_{6} & 6 & -5 & 0 \\
    			864 \, A_{0} A_{1} \alpha_{14} & 10 & -7 & 0 \\
    			240 \, A_{0} A_{1} \alpha_{7} - 40 \, {\left(6 \, A_{0}^{2} \alpha_{1} - 5 \, A_{1}^{2} - 6 \, A_{0} A_{2}\right)} \overline{\alpha_{2}} & 7 & -6 & 0 \\
    			-240 \, A_{0}^{2} \alpha_{3} \overline{\alpha_{2}} + 336 \, A_{0} A_{1} \alpha_{8} - 96 \, {\left(4 \, A_{0}^{2} \alpha_{1} - 3 \, A_{1}^{2} - 4 \, A_{0} A_{2}\right)} \alpha_{7} - 16 \, {\left(61 \, A_{0} A_{1} \alpha_{1} - 35 \, A_{1} A_{2} - 15 \, A_{0} A_{3}\right)} \overline{\alpha_{2}} & 8 & -6 & 0 \\
    			-880 \, A_{0} A_{1} \overline{\alpha_{2}}^{2} + 1440 \, A_{0} A_{1} \alpha_{13} & 11 & -8 & 0 \\
    			160 \, A_{0} A_{1} \overline{\alpha_{2}} & 6 & -6 & 0 \\
    			-160 \, A_{0}^{2} \zeta_{0} \overline{\alpha_{2}} + 80 \, A_{0} \overline{E_{1}} \overline{\alpha_{2}} & 9 & -6 & 1 \\
    			6 \, A_{1} C_{1} & -1 & 0 & 0 \\
    			-6 \, A_{1} C_{1} \overline{\zeta_{0}} & -1 & 4 & 1
    			\end{dmatrix}
    		\end{align*}
    		    		\small
    		where

    		Then we have
    		\footnotesize
    		\begin{align*}
    			&-K_g\,\h_0\totimes\h_0=\\
    			&\begin{dmatrix}
    			16 \, A_{1}^{2} {\left| A_{1} \right|}^{2} & 3 & -3 & 0 \\
    			-64 \, A_{0} A_{1} \alpha_{1} {\left| A_{1} \right|}^{2} + 64 \, A_{1} A_{2} {\left| A_{1} \right|}^{2} + 16 \, A_{1}^{2} \alpha_{4} & 4 & -3 & 0 \\
    			\kappa_{1} & 5 & -3 & 0 \\
    			\kappa_{2} & 6 & -3 & 0 \\
    			-4 \, A_{1} C_{1} {\left| A_{1} \right|}^{2} - 80 \, A_{1}^{2} \alpha_{2} & 0 & 1 & 0 \\
    			8 \, A_{0} C_{1} \alpha_{1} {\left| A_{1} \right|}^{2} + 320 \, A_{0} A_{1} \alpha_{1} \alpha_{2} - 8 \, A_{2} C_{1} {\left| A_{1} \right|}^{2} - 2 \, A_{1} C_{2} {\left| A_{1} \right|}^{2} - 320 \, A_{1} A_{2} \alpha_{2} - 4 \, A_{1} C_{1} \alpha_{4} - 40 \, A_{1}^{2} \overline{\alpha_{6}} & 1 & 1 & 0 \\
    			8 \, A_{0} C_{1} {\left| A_{1} \right|}^{4} + 384 \, A_{0} A_{1} \alpha_{2} {\left| A_{1} \right|}^{2} - \frac{16}{5} \, A_{1} B_{1} {\left| A_{1} \right|}^{2} - 4 \, A_{1} C_{1} \overline{\alpha_{4}} - 96 \, A_{1}^{2} \overline{\alpha_{7}} & 0 & 2 & 0 \\
    			\kappa_{3} & 0 & 3 & 0 \\
    			\kappa_{4} & 1 & 2 & 0 \\
    			-40 \, A_{0} C_{1} \alpha_{1} \alpha_{2} + \frac{1}{4} \, C_{1} C_{2} {\left| A_{1} \right|}^{2} - 4 \, A_{1} E_{1} {\left| A_{1} \right|}^{2} + 40 \, A_{2} C_{1} \alpha_{2} + 10 \, A_{1} C_{2} \alpha_{2} + \frac{1}{4} \, C_{1}^{2} \alpha_{4} + 10 \, A_{1} C_{1} \overline{\alpha_{6}} & -2 & 5 & 0 \\
    			\kappa_{5} & 2 & 1 & 0 \\
    			64 \, A_{0}^{2} {\left| A_{1} \right|}^{6} - 16 \, A_{1}^{2} {\left| A_{1} \right|}^{2} \overline{\alpha_{1}} - 96 \, A_{0} A_{1} {\left| A_{1} \right|}^{2} \overline{\alpha_{4}} + 4 \, A_{1} {\left| A_{1} \right|}^{2} \overline{B_{1}} - 24 \, {\left(2 \, {\left| A_{1} \right|}^{2} \overline{\alpha_{1}} - \overline{\alpha_{5}}\right)} A_{1}^{2} & 3 & -1 & 0 \\
    		    \kappa_{6} & 4 & -1 & 0 \\
    			\kappa_{7} & 3 & 0 & 0 \\
    			-128 \, A_{0} A_{1} \zeta_{0} {\left| A_{1} \right|}^{2} + 64 \, A_{1} {\left| A_{1} \right|}^{2} \overline{E_{1}} & 6 & -3 & 1 \\
    			-64 \, A_{0} A_{1} {\left| A_{1} \right|}^{4} + 16 \, A_{1}^{2} \overline{\alpha_{4}} & 3 & -2 & 0 \\
    			\kappa_{8} & 4 & -2 & 0 \\
    			\kappa_{9} & 5 & -2 & 0 \\
    			-64 \, A_{0} A_{1} \alpha_{1} \zeta_{0} + 80 \, A_{1}^{2} \overline{\alpha_{1}} \overline{\alpha_{2}} + 64 \, A_{1} A_{2} \zeta_{0} + 20 \, A_{1}^{2} \zeta_{1} - 20 \, A_{1} \overline{B_{1}} \overline{\alpha_{2}} & 7 & -4 & 0 \\
    			\kappa_{10} & 8 & -5 & 0 \\
    			160 \, A_{0} A_{1} {\left| A_{1} \right|}^{2} \overline{\alpha_{2}} - 40 \, A_{1}^{2} \alpha_{6} & 7 & -5 & 0 \\
    			\frac{1}{4} \, C_{1}^{2} {\left| A_{1} \right|}^{2} + 20 \, A_{1} C_{1} \alpha_{2} & -3 & 5 & 0 \\
    			-48 \, A_{0} C_{1} \alpha_{2} {\left| A_{1} \right|}^{2} - 320 \, A_{0} A_{1} \alpha_{2}^{2} + \frac{2}{5} \, B_{1} C_{1} {\left| A_{1} \right|}^{2} + 16 \, A_{1} B_{1} \alpha_{2} + \frac{1}{4} \, C_{1}^{2} \overline{\alpha_{4}} + 24 \, A_{1} C_{1} \overline{\alpha_{7}} & -3 & 6 & 0 \\
    			-\frac{5}{4} \, C_{1}^{2} \alpha_{2} & -6 & 9 & 0 \\
    			320 \, A_{0} A_{1} \alpha_{1} \overline{\alpha_{2}} - 96 \, A_{1}^{2} \alpha_{7} - 320 \, A_{1} A_{2} \overline{\alpha_{2}} & 8 & -6 & 0 \\
    			\kappa_{11} & 9 & -6 & 0 \\
    			16 \, A_{1}^{2} \zeta_{0} & 6 & -4 & 0 \\
    			-80 \, A_{1}^{2} \overline{\alpha_{2}} & 7 & -6 & 0 \\
    			800 \, A_{0} A_{1} \overline{\alpha_{2}}^{2} & 11 & -8 & 0
    			\end{dmatrix}
    		\end{align*}
    		\small
    		where
    		\begin{align*}
    			\kappa_{1}&=64 \, A_{0}^{2} \alpha_{1}^{2} {\left| A_{1} \right|}^{2} - 80 \, A_{1}^{2} \alpha_{1} {\left| A_{1} \right|}^{2} - 128 \, A_{0} A_{2} \alpha_{1} {\left| A_{1} \right|}^{2} - 96 \, A_{0} A_{1} \alpha_{3} {\left| A_{1} \right|}^{2}\\
    			& - 64 \, A_{0} A_{1} \alpha_{1} \alpha_{4} + 64 \, A_{2}^{2} {\left| A_{1} \right|}^{2} + 96 \, A_{1} A_{3} {\left| A_{1} \right|}^{2} - 24 \, {\left(2 \, \alpha_{1} {\left| A_{1} \right|}^{2} - \alpha_{5}\right)} A_{1}^{2} + 64 \, A_{1} A_{2} \alpha_{4}\\
    			\kappa_{2}&=192 \, A_{0} A_{1} \alpha_{1}^{2} {\left| A_{1} \right|}^{2} + 192 \, A_{0}^{2} \alpha_{1} \alpha_{3} {\left| A_{1} \right|}^{2} + 64 \, A_{0}^{2} \alpha_{1}^{2} \alpha_{4} + 16 \, {\left(4 \, \alpha_{1}^{2} - \zeta_{0}\right)} A_{0} A_{1} {\left| A_{1} \right|}^{2} - 256 \, A_{1} A_{2} \alpha_{1} {\left| A_{1} \right|}^{2}\\
    			& - 192 \, A_{0} A_{3} \alpha_{1} {\left| A_{1} \right|}^{2} - 112 \, A_{1}^{2} \alpha_{3} {\left| A_{1} \right|}^{2} - 192 \, A_{0} A_{2} \alpha_{3} {\left| A_{1} \right|}^{2} - 64 \, A_{0} A_{1} \zeta_{0} {\left| A_{1} \right|}^{2} + 96 \, {\left(2 \, \alpha_{1} {\left| A_{1} \right|}^{2} - \alpha_{5}\right)} A_{0} A_{1} \alpha_{1}\\
    			& - 80 \, A_{1}^{2} \alpha_{1} \alpha_{4} - 128 \, A_{0} A_{2} \alpha_{1} \alpha_{4} - 96 \, A_{0} A_{1} \alpha_{3} \alpha_{4} + 192 \, A_{2} A_{3} {\left| A_{1} \right|}^{2} + 128 \, A_{1} A_{4} {\left| A_{1} \right|}^{2} + 8 \, A_{1} {\left| A_{1} \right|}^{2} \overline{E_{1}}\\
    			& - 32 \, {\left(2 \, \alpha_{3} {\left| A_{1} \right|}^{2} + \alpha_{1} \alpha_{4} - \alpha_{11}\right)} A_{1}^{2} - 96 \, {\left(2 \, \alpha_{1} {\left| A_{1} \right|}^{2} - \alpha_{5}\right)} A_{1} A_{2} + 64 \, A_{2}^{2} \alpha_{4} + 96 \, A_{1} A_{3} \alpha_{4}\\
    			\kappa_{3}&=-448 \, A_{0}^{2} \alpha_{2} {\left| A_{1} \right|}^{4} + \frac{32}{5} \, A_{0} B_{1} {\left| A_{1} \right|}^{4} + \frac{4}{3} \, A_{1} C_{1} {\left| A_{1} \right|}^{2} \overline{\alpha_{1}} + 12 \, A_{0} C_{1} {\left| A_{1} \right|}^{2} \overline{\alpha_{4}} + 448 \, A_{0} A_{1} {\left| A_{1} \right|}^{2} \overline{\alpha_{7}} - \frac{8}{3} \, A_{1} B_{2} {\left| A_{1} \right|}^{2}\\
    			& - \frac{1}{2} \, C_{1} {\left| A_{1} \right|}^{2} \overline{B_{1}} + 80 \, A_{1}^{2} \alpha_{2} \overline{\alpha_{1}} + 224 \, A_{0} A_{1} \alpha_{2} \overline{\alpha_{4}} + 112 \, {\left(\alpha_{2} \overline{\alpha_{1}} - \overline{\alpha_{8}}\right)} A_{1}^{2} + 6 \, {\left(2 \, {\left| A_{1} \right|}^{2} \overline{\alpha_{1}} - \overline{\alpha_{5}}\right)} A_{1} C_{1} - 20 \, A_{1} \alpha_{2} \overline{B_{1}}\\
    			& - \frac{16}{5} \, A_{1} B_{1} \overline{\alpha_{4}}\\
    			\kappa_{4}&=-768 \, A_{0}^{2} \alpha_{1} \alpha_{2} {\left| A_{1} \right|}^{2} + \frac{44}{5} \, A_{1} C_{1} {\left| A_{1} \right|}^{4} + 4 \, A_{0} C_{2} {\left| A_{1} \right|}^{4} + \frac{32}{5} \, A_{0} B_{1} \alpha_{1} {\left| A_{1} \right|}^{2} + 528 \, A_{1}^{2} \alpha_{2} {\left| A_{1} \right|}^{2} + 768 \, A_{0} A_{2} \alpha_{2} {\left| A_{1} \right|}^{2}\\
    			& + 16 \, A_{0} C_{1} \alpha_{4} {\left| A_{1} \right|}^{2} + 192 \, A_{0} A_{1} {\left| A_{1} \right|}^{2} \overline{\alpha_{6}} + 384 \, A_{0} A_{1} \alpha_{2} \alpha_{4} - \frac{32}{5} \, A_{2} B_{1} {\left| A_{1} \right|}^{2} - \frac{8}{5} \, A_{1} B_{3} {\left| A_{1} \right|}^{2} + 8 \, A_{0} C_{1} \alpha_{1} \overline{\alpha_{4}}\\
    			& + 384 \, A_{0} A_{1} \alpha_{1} \overline{\alpha_{7}} + 48 \, {\left(2 \, \alpha_{2} {\left| A_{1} \right|}^{2} - \overline{\alpha_{9}}\right)} A_{1}^{2} + 8 \, {\left(2 \, {\left| A_{1} \right|}^{4} + \alpha_{1} \overline{\alpha_{1}} - \beta\right)} A_{1} C_{1} - \frac{16}{5} \, A_{1} B_{1} \alpha_{4} - \frac{5}{4} \, C_{1}^{2} \overline{\alpha_{2}} - 8 \, A_{2} C_{1} \overline{\alpha_{4}}\\
    			& - 2 \, A_{1} C_{2} \overline{\alpha_{4}} - 384 \, A_{1} A_{2} \overline{\alpha_{7}}\\
    			\kappa_{5}&=-320 \, A_{0}^{2} \alpha_{1}^{2} \alpha_{2} + 8 \, A_{1} C_{1} \alpha_{1} {\left| A_{1} \right|}^{2} + 4 \, A_{0} C_{2} \alpha_{1} {\left| A_{1} \right|}^{2} + 12 \, A_{0} C_{1} \alpha_{3} {\left| A_{1} \right|}^{2} - 16 \, A_{0} A_{1} {\left| A_{1} \right|}^{2} \overline{\zeta_{0}} + 400 \, A_{1}^{2} \alpha_{1} \alpha_{2}\\
    			& + 640 \, A_{0} A_{2} \alpha_{1} \alpha_{2} + 480 \, A_{0} A_{1} \alpha_{2} \alpha_{3} + 8 \, A_{0} C_{1} \alpha_{1} \alpha_{4} - 12 \, A_{3} C_{1} {\left| A_{1} \right|}^{2} - 4 \, A_{2} C_{2} {\left| A_{1} \right|}^{2} + 2 \, A_{1} \gamma_{1} {\left| A_{1} \right|}^{2} + 160 \, A_{0} A_{1} \alpha_{1} \overline{\alpha_{6}}\\
    			& + 6 \, {\left(2 \, \alpha_{1} {\left| A_{1} \right|}^{2} - \alpha_{5}\right)} A_{1} C_{1} - 320 \, A_{2}^{2} \alpha_{2} - 480 \, A_{1} A_{3} \alpha_{2} - 8 \, A_{2} C_{1} \alpha_{4} - 2 \, A_{1} C_{2} \alpha_{4} - 160 \, A_{1} A_{2} \overline{\alpha_{6}}\\
    			\kappa_{6}&=256 \, A_{0} A_{1} {\left| A_{1} \right|}^{6} + 192 \, A_{0}^{2} \alpha_{4} {\left| A_{1} \right|}^{4} + 64 \, A_{0} A_{1} \alpha_{1} {\left| A_{1} \right|}^{2} \overline{\alpha_{1}} + 192 \, A_{0}^{2} \alpha_{1} {\left| A_{1} \right|}^{2} \overline{\alpha_{4}} + 192 \, {\left(2 \, {\left| A_{1} \right|}^{4} + \alpha_{1} \overline{\alpha_{1}} - \beta\right)} A_{0} A_{1} {\left| A_{1} \right|}^{2}\\
    			& - 8 \, A_{0} \alpha_{1} {\left| A_{1} \right|}^{2} \overline{B_{1}} - 64 \, A_{1} A_{2} {\left| A_{1} \right|}^{2} \overline{\alpha_{1}} - 20 \, A_{0} C_{1} {\left| A_{1} \right|}^{2} \overline{\alpha_{2}} - 144 \, A_{1}^{2} {\left| A_{1} \right|}^{2} \overline{\alpha_{4}} - 192 \, A_{0} A_{2} {\left| A_{1} \right|}^{2} \overline{\alpha_{4}}\\
    			& + 96 \, {\left(2 \, {\left| A_{1} \right|}^{2} \overline{\alpha_{1}} - \overline{\alpha_{5}}\right)} A_{0} A_{1} \alpha_{1} + 8 \, A_{2} {\left| A_{1} \right|}^{2} \overline{B_{1}} + 8 \, A_{1} {\left| A_{1} \right|}^{2} \overline{B_{2}} - 16 \, A_{1}^{2} \alpha_{4} \overline{\alpha_{1}} + 480 \, A_{0} A_{1} \alpha_{2} \overline{\alpha_{2}} - 96 \, A_{0} A_{1} \alpha_{4} \overline{\alpha_{4}}\\
    			& - 48 \, {\left(2 \, {\left| A_{1} \right|}^{2} \overline{\alpha_{4}} + \alpha_{4} \overline{\alpha_{1}} + \alpha_{1} \overline{\alpha_{3}} - \overline{\alpha_{12}}\right)} A_{1}^{2} - 96 \, {\left(2 \, {\left| A_{1} \right|}^{2} \overline{\alpha_{1}} - \overline{\alpha_{5}}\right)} A_{1} A_{2} + 10 \, A_{1} C_{1} \alpha_{6} + 4 \, A_{1} \alpha_{4} \overline{B_{1}} + 16 \, A_{1} B_{1} \overline{\alpha_{2}}\\
    			\kappa_{7}&=64 \, A_{0} A_{1} {\left| A_{1} \right|}^{4} \overline{\alpha_{1}} + 128 \, A_{0}^{2} {\left| A_{1} \right|}^{4} \overline{\alpha_{4}} - 8 \, A_{0} {\left| A_{1} \right|}^{4} \overline{B_{1}} - \frac{8}{3} \, A_{1} {\left| A_{1} \right|}^{4} \overline{C_{1}} + 128 \, {\left(2 \, {\left| A_{1} \right|}^{2} \overline{\alpha_{1}} - \overline{\alpha_{5}}\right)} A_{0} A_{1} {\left| A_{1} \right|}^{2} - 16 \, A_{1}^{2} {\left| A_{1} \right|}^{2} \overline{\alpha_{3}}\\
    			& + \frac{8}{3} \, A_{1} {\left| A_{1} \right|}^{2} \overline{B_{3}} - 16 \, A_{1}^{2} \overline{\alpha_{1}} \overline{\alpha_{4}} - 32 \, A_{0} A_{1} \overline{\alpha_{4}}^{2} - 32 \, {\left(2 \, {\left| A_{1} \right|}^{2} \overline{\alpha_{3}} + \overline{\alpha_{1}} \overline{\alpha_{4}} - \overline{\alpha_{11}}\right)} A_{1}^{2} - 4 \, A_{1} C_{1} \zeta_{0} + 4 \, A_{1} \overline{B_{1}} \overline{\alpha_{4}}\\
    			\kappa_{8}&=128 \, A_{0}^{2} \alpha_{1} {\left| A_{1} \right|}^{4} - 96 \, A_{1}^{2} {\left| A_{1} \right|}^{4} - 128 \, A_{0} A_{2} {\left| A_{1} \right|}^{4} - 128 \, A_{0} A_{1} \alpha_{4} {\left| A_{1} \right|}^{2} - 64 \, A_{0} A_{1} \alpha_{1} \overline{\alpha_{4}} - 32 \, {\left(2 \, {\left| A_{1} \right|}^{4} + \alpha_{1} \overline{\alpha_{1}} - \beta\right)} A_{1}^{2}\\
    			& + 20 \, A_{1} C_{1} \overline{\alpha_{2}} + 64 \, A_{1} A_{2} \overline{\alpha_{4}}\\
    			\kappa_{9}&=448 \, A_{0} A_{1} \alpha_{1} {\left| A_{1} \right|}^{4} + 192 \, A_{0}^{2} \alpha_{3} {\left| A_{1} \right|}^{4} + 256 \, A_{0}^{2} \alpha_{1} \alpha_{4} {\left| A_{1} \right|}^{2} - 320 \, A_{1} A_{2} {\left| A_{1} \right|}^{4} - 192 \, A_{0} A_{3} {\left| A_{1} \right|}^{4}\\
    			& + 192 \, {\left(2 \, \alpha_{1} {\left| A_{1} \right|}^{2} - \alpha_{5}\right)} A_{0} A_{1} {\left| A_{1} \right|}^{2} - 176 \, A_{1}^{2} \alpha_{4} {\left| A_{1} \right|}^{2} - 256 \, A_{0} A_{2} \alpha_{4} {\left| A_{1} \right|}^{2} + 64 \, A_{0}^{2} \alpha_{1}^{2} \overline{\alpha_{4}} + 128 \, {\left(2 \, {\left| A_{1} \right|}^{4} + \alpha_{1} \overline{\alpha_{1}} - \beta\right)} A_{0} A_{1} \alpha_{1}\\
    			& - 64 \, A_{0} A_{1} \alpha_{4}^{2} - 40 \, A_{0} C_{1} \alpha_{1} \overline{\alpha_{2}} - 80 \, A_{1}^{2} \alpha_{1} \overline{\alpha_{4}} - 128 \, A_{0} A_{2} \alpha_{1} \overline{\alpha_{4}} - 96 \, A_{0} A_{1} \alpha_{3} \overline{\alpha_{4}} - 48 \, {\left(2 \, \alpha_{4} {\left| A_{1} \right|}^{2} + \alpha_{3} \overline{\alpha_{1}} + \alpha_{1} \overline{\alpha_{4}} - \alpha_{12}\right)} A_{1}^{2} \\
    			&- 128 \, {\left(2 \, {\left| A_{1} \right|}^{4} + \alpha_{1} \overline{\alpha_{1}} - \beta\right)} A_{1} A_{2} + 24 \, A_{1} C_{1} \alpha_{7} + 40 \, A_{2} C_{1} \overline{\alpha_{2}} + 10 \, A_{1} C_{2} \overline{\alpha_{2}} + 64 \, A_{2}^{2} \overline{\alpha_{4}} + 96 \, A_{1} A_{3} \overline{\alpha_{4}}\\
    			\kappa_{10}&=-320 \, A_{0}^{2} \alpha_{1} {\left| A_{1} \right|}^{2} \overline{\alpha_{2}} + 192 \, A_{0} A_{1} \alpha_{7} {\left| A_{1} \right|}^{2} + 304 \, A_{1}^{2} {\left| A_{1} \right|}^{2} \overline{\alpha_{2}} + 320 \, A_{0} A_{2} {\left| A_{1} \right|}^{2} \overline{\alpha_{2}} + 160 \, A_{0} A_{1} \alpha_{1} \alpha_{6} + 160 \, A_{0} A_{1} \alpha_{4} \overline{\alpha_{2}}\\
    			& + 48 \, {\left(2 \, {\left| A_{1} \right|}^{2} \overline{\alpha_{2}} - \alpha_{9}\right)} A_{1}^{2} - 160 \, A_{1} A_{2} \alpha_{6}\\
    			\kappa_{11}&=-320 \, A_{0}^{2} \alpha_{1}^{2} \overline{\alpha_{2}} + 384 \, A_{0} A_{1} \alpha_{1} \alpha_{7} + 400 \, A_{1}^{2} \alpha_{1} \overline{\alpha_{2}} + 640 \, A_{0} A_{2} \alpha_{1} \overline{\alpha_{2}} + 480 \, A_{0} A_{1} \alpha_{3} \overline{\alpha_{2}} + 112 \, {\left(\alpha_{1} \overline{\alpha_{2}} - \alpha_{8}\right)} A_{1}^{2} - 384 \, A_{1} A_{2} \alpha_{7}\\
    			& - 320 \, A_{2}^{2} \overline{\alpha_{2}} - 480 \, A_{1} A_{3} \overline{\alpha_{2}}
    		\end{align*}
    		Gathering these two terms, we obtain
    		\footnotesize
    		\begin{align*}
    			&\mathscr{Q}_{\phi}-\left(\frac{5}{4}|\H|^2\h_0\totimes\h_0+\s{\H}{\h_0}^2\right)=g^{-1}\otimes Q(\h_0)-K_g\,\h_0\totimes\h_0=\\
    			&\begin{dmatrix}
    			-16 \, A_{0} C_{1} \alpha_{1} + 16 \, A_{2} C_{1} + 2 \, A_{1} C_{2} & 0 & 0 & 0 \\
    			-30 \, A_{0} C_{1} \alpha_{3} + 8 \, A_{0} A_{1} \overline{\zeta_{0}} + 30 \, A_{3} C_{1} + 6 \, A_{2} C_{2} - 6 \, {\left(4 \, A_{1} C_{1} + A_{0} C_{2}\right)} \alpha_{1} - A_{1} \gamma_{1} & 1 & 0 & 0 \\
    			\pi_1 & 2 & 0 & 0 \\
    			\pi_2 & 3 & 0 & 0 \\
    			20 \, A_{1} E_{1} & -3 & 4 & 0 \\
    			-48 \, A_{0} E_{1} \alpha_{1} + 48 \, A_{2} E_{1} + 12 \, A_{1} E_{3} & -2 & 4 & 0 \\
    			\frac{29}{80} \, C_{1}^{2} {\left| A_{1} \right|}^{2} - 35 \, A_{0} E_{1} {\left| A_{1} \right|}^{2} - \frac{1}{2} \, A_{0} C_{1} \overline{\alpha_{6}} + \frac{1}{40} \, B_{3} C_{1} - \frac{1}{40} \, B_{1} C_{2} + 20 \, A_{1} E_{2} + \frac{1}{2} \, {\left(21 \, A_{1} C_{1} + A_{0} C_{2}\right)} \alpha_{2} & -3 & 5 & 0 \\
    			\pi_3 & -3 & 6 & 0 \\
    			\pi_4 & -2 & 5 & 0 \\
    			-\frac{3}{8} \, C_{2} E_{1} - \frac{3}{16} \, C_{1} E_{3} & -5 & 8 & 0 \\
    			\pi_5 & -1 & 4 & 0 \\
    			\pi_6 & 0 & 2 & 0 \\
    			\pi_7 & 1 & 2 & 0 \\
    			\pi_8 & 0 & 3 & 0 \\
    			\pi_9 & 3 & 0 & 1 \\
    			\pi_{10} & 0 & 1 & 0 \\
    			\pi_{11} & 1 & 1 & 0 \\
    			\pi_{12} & 2 & 1 & 0 \\
    			\pi_{13} & 4 & -1 & 0 \\
    			\pi_{14} & 5 & -2 & 0 \\
    			\pi_{15} & 4 & -2 & 0 \\
    			-\frac{1}{2} \, C_{1} E_{1} & -6 & 8 & 0 \\
    			-\frac{3}{10} \, B_{1} E_{1} - \frac{5}{9} \, C_{1} E_{2} - \frac{1}{72} \, {\left(85 \, C_{1}^{2} - 432 \, A_{0} E_{1}\right)} \alpha_{2} & -6 & 9 & 0 \\
    			-48 \, A_{0} C_{1} \zeta_{0} + 24 \, C_{1} \overline{E_{1}} & 2 & 0 & 1 \\
    			-480 \, A_{0} \alpha_{2} \overline{E_{1}} + 48 \, {\left(20 \, A_{0}^{2} \alpha_{2} - A_{0} B_{1}\right)} \zeta_{0} + 24 \, B_{1} \overline{E_{1}} & 2 & 1 & 1 \\
    			\pi_{16} & 3 & -2 & 0 \\
    			192 \, A_{0}^{2} \alpha_{2} {\left| A_{1} \right|}^{2} - \frac{48}{5} \, A_{0} B_{1} {\left| A_{1} \right|}^{2} - 3 \, A_{0} C_{1} \overline{\alpha_{4}} - 144 \, A_{0} A_{1} \overline{\alpha_{7}} + 6 \, A_{1} B_{2} + \frac{3}{8} \, C_{1} \overline{B_{1}} & -1 & 2 & 0 \\
    			\pi_{17} & -1 & 3 & 0 \\
    			-96 \, A_{0} A_{1} \zeta_{4} + 48 \, A_{0} \overline{E_{1}} \overline{\alpha_{4}} + 12 \, {\left(8 \, A_{0} A_{1} \overline{\alpha_{1}} - 8 \, A_{0}^{2} \overline{\alpha_{4}} + A_{0} \overline{B_{1}}\right)} \zeta_{0} - 6 \, \overline{B_{1}} \overline{E_{1}} & 5 & -2 & 1 \\
    			\pi_{18} & 3 & -1 & 0 \\
    			\pi_{19} & 6 & -3 & 0 \\
    			\pi_{20} & 7 & -4 & 0 
    			    			\end{dmatrix}
    			\end{align*}
    			\begin{align*}
    			\begin{dmatrix}
    			\pi_{21} & 6 & -4 & 0 \\
    			-24 \, A_{0} A_{1} \zeta_{0} + 12 \, A_{1} \overline{E_{1}} & 5 & -4 & 0 \\
    			\pi_{22} & 5 & -3 & 0 \\
    			\pi_{23} & 9 & -6 & 0 \\
    			-16 \, A_{0}^{2} \alpha_{1} {\left| A_{1} \right|}^{2} + 8 \, A_{1}^{2} {\left| A_{1} \right|}^{2} + 16 \, A_{0} A_{2} {\left| A_{1} \right|}^{2} - 8 \, A_{0} A_{1} \alpha_{4} & 3 & -3 & 0 \\
    			-48 \, A_{0} A_{1} \alpha_{1} {\left| A_{1} \right|}^{2} - 48 \, A_{0}^{2} \alpha_{3} {\left| A_{1} \right|}^{2} + 48 \, A_{1} A_{2} {\left| A_{1} \right|}^{2} + 48 \, A_{0} A_{3} {\left| A_{1} \right|}^{2} - 24 \, A_{0} A_{1} \alpha_{5} & 4 & -3 & 0 \\
    			-9 \, A_{0} C_{1} {\left| A_{1} \right|}^{2} - 120 \, A_{0} A_{1} \alpha_{2} + 6 \, A_{1} B_{1} & -1 & 1 & 0 \\
    			-96 \, A_{0}^{2} \zeta_{0} {\left| A_{1} \right|}^{2} + 48 \, A_{0} {\left| A_{1} \right|}^{2} \overline{E_{1}} & 5 & -3 & 1 \\
    			\pi_{24} & 6 & -3 & 1 \\
    			-720 \, A_{0}^{2} \alpha_{7} {\left| A_{1} \right|}^{2} - 1360 \, A_{0} A_{1} {\left| A_{1} \right|}^{2} \overline{\alpha_{2}} - 360 \, A_{0}^{2} \alpha_{4} \overline{\alpha_{2}} + 120 \, A_{0} A_{1} \alpha_{9} - 60 \, {\left(2 \, A_{0}^{2} \alpha_{1} - A_{1}^{2} - 2 \, A_{0} A_{2}\right)} \alpha_{6} & 7 & -5 & 0 \\
    			\pi_{25} & 8 & -5 & 0 \\
    			-480 \, A_{0}^{2} {\left| A_{1} \right|}^{2} \overline{\alpha_{2}} + 80 \, A_{0} A_{1} \alpha_{6} & 6 & -5 & 0 \\
    			864 \, A_{0} A_{1} \alpha_{14} & 10 & -7 & 0 \\
    			240 \, A_{0} A_{1} \alpha_{7} - 120 \, {\left(2 \, A_{0}^{2} \alpha_{1} - A_{1}^{2} - 2 \, A_{0} A_{2}\right)} \overline{\alpha_{2}} & 7 & -6 & 0 \\
    			-240 \, A_{0}^{2} \alpha_{3} \overline{\alpha_{2}} + 336 \, A_{0} A_{1} \alpha_{8} - 192 \, {\left(2 \, A_{0}^{2} \alpha_{1} - A_{1}^{2} - 2 \, A_{0} A_{2}\right)} \alpha_{7} - 16 \, {\left(41 \, A_{0} A_{1} \alpha_{1} - 15 \, A_{1} A_{2} - 15 \, A_{0} A_{3}\right)} \overline{\alpha_{2}} & 8 & -6 & 0 \\
    			-80 \, A_{0} A_{1} \overline{\alpha_{2}}^{2} + 1440 \, A_{0} A_{1} \alpha_{13} & 11 & -8 & 0 \\
    			160 \, A_{0} A_{1} \overline{\alpha_{2}} & 6 & -6 & 0 \\
    			-160 \, A_{0}^{2} \zeta_{0} \overline{\alpha_{2}} + 80 \, A_{0} \overline{E_{1}} \overline{\alpha_{2}} & 9 & -6 & 1 \\
    			6 \, A_{1} C_{1} & -1 & 0 & 0 \\
    			-6 \, A_{1} C_{1} \overline{\zeta_{0}} & -1 & 4 & 1
    			\end{dmatrix}
    		\end{align*}
    		\small
    		where

    		Finally, the extra terms are (we have two expressions in order to determine in which order we take scalar products)
    		\begin{align}\label{4end0}
    			&\frac{5}{4}|\H|^2\h_0\totimes\h_0=\nonumber\\
    			&\begin{dmatrix}
    			 \frac{5}{4} \, \left(4 \, A_{1}^{2}\right) \left(\frac{1}{4} \, C_{1}^{2}\right) & 2 & 0 & 0 \\
    			\frac{5}{2} \, \left(-8 \, A_{0} A_{1} \alpha_{1}\right) \left(\frac{1}{4} \, C_{1}^{2}\right) + \frac{5}{2} \, \left(8 \, A_{1} A_{2}\right) \left(\frac{1}{4} \, C_{1}^{2}\right) + \frac{5}{2} \, \left(4 \, A_{1}^{2}\right) \left(\frac{1}{4} \, C_{1} C_{2}\right) + \frac{5}{2} \, \left(-\frac{1}{2} \, A_{1} C_{1}\right) \left(\frac{1}{4} \, \overline{C_{1}}^{2}\right) & 3 & 0 & 0 \\
    			\frac{5}{2} \, \left(-\frac{1}{2} \, A_{1} C_{1}\right) \left(\frac{1}{4} \, C_{1}^{2}\right) & -1 & 4 & 0 \\
    			\frac{5}{2} \, \left(4 \, A_{1}^{2}\right) \left(\frac{1}{4} \, B_{1} C_{1}\right) + \frac{5}{2} \, \left(-8 \, A_{0} A_{1} {\left| A_{1} \right|}^{2}\right) \left(\frac{1}{4} \, C_{1}^{2}\right) & 2 & 1 & 0 \\
    			\frac{5}{2} \, \left(4 \, A_{1}^{2}\right) \left(\frac{1}{4} \, C_{1} \overline{C_{1}}\right) & 4 & -2 & 0 \\
    			\frac{5}{2} \, \left(4 \, A_{1}^{2}\right) \left(\frac{1}{4} \, C_{1} \overline{B_{1}}\right) + 5 \, \left(-8 \, A_{0} A_{1} \alpha_{1}\right) \left(\frac{1}{4} \, C_{1} \overline{C_{1}}\right) + 5 \, \left(8 \, A_{1} A_{2}\right) \left(\frac{1}{4} \, C_{1} \overline{C_{1}}\right) + \frac{5}{2} \, \left(4 \, A_{1}^{2}\right) \left(\frac{1}{4} \, C_{2} \overline{C_{1}}\right) & 5 & -2 & 0 \\
    			5 \, \left(-\frac{1}{2} \, A_{1} C_{1}\right) \left(\frac{1}{4} \, C_{1} \overline{C_{1}}\right) & 1 & 2 & 0 \\
    			\frac{5}{2} \, \left(4 \, A_{1}^{2}\right) \left(\frac{1}{4} \, B_{1} \overline{C_{1}}\right) + 5 \, \left(-8 \, A_{0} A_{1} {\left| A_{1} \right|}^{2}\right) \left(\frac{1}{4} \, C_{1} \overline{C_{1}}\right) + \frac{5}{2} \, \left(4 \, A_{1}^{2}\right) \left(\frac{1}{4} \, C_{1} \overline{C_{2}}\right) & 4 & -1 & 0 \\
    			\frac{5}{4} \, \left(4 \, A_{1}^{2}\right) \left(\frac{1}{4} \, \overline{C_{1}}^{2}\right) & 6 & -4 & 0 \\
    			\frac{5}{2} \, \left(4 \, A_{1}^{2}\right) \left(\frac{1}{4} \, \overline{B_{1}} \overline{C_{1}}\right) + \frac{5}{2} \, \left(-8 \, A_{0} A_{1} \alpha_{1}\right) \left(\frac{1}{4} \, \overline{C_{1}}^{2}\right) + \frac{5}{2} \, \left(8 \, A_{1} A_{2}\right) \left(\frac{1}{4} \, \overline{C_{1}}^{2}\right) & 7 & -4 & 0 \\
    			\frac{5}{2} \, \left(-8 \, A_{0} A_{1} {\left| A_{1} \right|}^{2}\right) \left(\frac{1}{4} \, \overline{C_{1}}^{2}\right) + \frac{5}{2} \, \left(4 \, A_{1}^{2}\right) \left(\frac{1}{4} \, \overline{C_{1}} \overline{C_{2}}\right) & 6 & -3 & 0
    			\end{dmatrix}\nonumber\\
    			&=
    			\begin{dmatrix}
    				\frac{5}{4} \, A_{1}^{2} C_{1}^{2} & 2 & 0 & 0 \\
    			-5 \, A_{0} A_{1} C_{1}^{2} \alpha_{1} + 5 \, A_{1} A_{2} C_{1}^{2} + \frac{5}{2} \, A_{1}^{2} C_{1} C_{2} - \frac{5}{16} \, A_{1} C_{1} \overline{C_{1}}^{2} & 3 & 0 & 0 \\
    			-\frac{5}{16} \, A_{1} C_{1}^{3} & -1 & 4 & 0 \\
    			-5 \, A_{0} A_{1} C_{1}^{2} {\left| A_{1} \right|}^{2} + \frac{5}{2} \, A_{1}^{2} B_{1} C_{1} & 2 & 1 & 0 \\
    			\frac{5}{2} \, A_{1}^{2} C_{1} \overline{C_{1}} & 4 & -2 & 0 \\
    			-10 \, A_{0} A_{1} C_{1} \alpha_{1} \overline{C_{1}} + \frac{5}{2} \, A_{1}^{2} C_{1} \overline{B_{1}} + 10 \, A_{1} A_{2} C_{1} \overline{C_{1}} + \frac{5}{2} \, A_{1}^{2} C_{2} \overline{C_{1}} & 5 & -2 & 0 \\
    			-\frac{5}{8} \, A_{1} C_{1}^{2} \overline{C_{1}} & 1 & 2 & 0 \\
    			-10 \, A_{0} A_{1} C_{1} {\left| A_{1} \right|}^{2} \overline{C_{1}} + \frac{5}{2} \, A_{1}^{2} B_{1} \overline{C_{1}} + \frac{5}{2} \, A_{1}^{2} C_{1} \overline{C_{2}} & 4 & -1 & 0 \\
    			\frac{5}{4} \, A_{1}^{2} \overline{C_{1}}^{2} & 6 & -4 & 0 \\
    			-5 \, A_{0} A_{1} \alpha_{1} \overline{C_{1}}^{2} + \frac{5}{2} \, A_{1}^{2} \overline{B_{1}} \overline{C_{1}} + 5 \, A_{1} A_{2} \overline{C_{1}}^{2} & 7 & -4 & 0 \\
    			-5 \, A_{0} A_{1} {\left| A_{1} \right|}^{2} \overline{C_{1}}^{2} + \frac{5}{2} \, A_{1}^{2} \overline{C_{1}} \overline{C_{2}} & 6 & -3 & 0
    			\end{dmatrix}
    			\end{align}
    			and
    			\begin{align}\label{4end1}
    			&\s{\H}{\h_0}^2=\nonumber\\
    			&\begin{dmatrix}
    			\left(A_{1} C_{1}\right) \left(A_{1} C_{1}\right) & 2 & 0 & 0 \\
    			\chi_1 & 3 & 0 & 0 \\
    			\left(A_{1} C_{1}\right) \left(-\frac{1}{8} \, C_{1}^{2}\right) + \left(A_{1} C_{1}\right) \left(-\frac{1}{8} \, C_{1}^{2}\right) & -1 & 4 & 0 \\
    			\left(-2 \, A_{0} C_{1} {\left| A_{1} \right|}^{2}\right) \left(A_{1} C_{1}\right) + \left(-2 \, A_{0} C_{1} {\left| A_{1} \right|}^{2}\right) \left(A_{1} C_{1}\right) + \left(A_{1} B_{1}\right) \left(A_{1} C_{1}\right) + \left(A_{1} B_{1}\right) \left(A_{1} C_{1}\right) & 2 & 1 & 0 \\
    			\left(A_{1} C_{1}\right) \left(A_{1} \overline{C_{1}}\right) + \left(A_{1} C_{1}\right) \left(A_{1} \overline{C_{1}}\right) & 4 & -2 & 0 \\
    			\chi_2 & 5 & -2 & 0 \\
    			\left(-\frac{1}{8} \, C_{1}^{2}\right) \left(A_{1} \overline{C_{1}}\right) + \left(-\frac{1}{8} \, C_{1}^{2}\right) \left(A_{1} \overline{C_{1}}\right) + \left(A_{1} C_{1}\right) \left(-\frac{1}{8} \, C_{1} \overline{C_{1}}\right) + \left(A_{1} C_{1}\right) \left(-\frac{1}{8} \, C_{1} \overline{C_{1}}\right) & 1 & 2 & 0 \\
    		    \chi_3 & 4 & -1 & 0 \\
    			\left(A_{1} \overline{C_{1}}\right) \left(A_{1} \overline{C_{1}}\right) & 6 & -4 & 0 \\
    			\chi_4 & 7 & -4 & 0 \\
    			\left(-2 \, A_{0} {\left| A_{1} \right|}^{2} \overline{C_{1}}\right) \left(A_{1} \overline{C_{1}}\right) + \left(-2 \, A_{0} {\left| A_{1} \right|}^{2} \overline{C_{1}}\right) \left(A_{1} \overline{C_{1}}\right) + \left(A_{1} \overline{C_{1}}\right) \left(A_{1} \overline{C_{2}}\right) + \left(A_{1} \overline{C_{1}}\right) \left(A_{1} \overline{C_{2}}\right) & 6 & -3 & 0
    			\end{dmatrix}\nonumber\\
    			&=\begin{dmatrix}
    			A_{1}^{2} C_{1}^{2} & 2 & 0 & 0 \\
    			-4 \, A_{0} A_{1} C_{1}^{2} \alpha_{1} + 4 \, A_{1} A_{2} C_{1}^{2} + 2 \, A_{1}^{2} C_{1} C_{2} - \frac{1}{4} \, A_{1} C_{1} \overline{C_{1}}^{2} & 3 & 0 & 0 \\
    			-\frac{1}{4} \, A_{1} C_{1}^{3} & -1 & 4 & 0 \\
    			-4 \, A_{0} A_{1} C_{1}^{2} {\left| A_{1} \right|}^{2} + 2 \, A_{1}^{2} B_{1} C_{1} & 2 & 1 & 0 \\
    			2 \, A_{1}^{2} C_{1} \overline{C_{1}} & 4 & -2 & 0 \\
    			-8 \, A_{0} A_{1} C_{1} \alpha_{1} \overline{C_{1}} + 2 \, A_{1}^{2} C_{1} \overline{B_{1}} + 8 \, A_{1} A_{2} C_{1} \overline{C_{1}} + 2 \, A_{1}^{2} C_{2} \overline{C_{1}} & 5 & -2 & 0 \\
    			-\frac{1}{2} \, A_{1} C_{1}^{2} \overline{C_{1}} & 1 & 2 & 0 \\
    			-8 \, A_{0} A_{1} C_{1} {\left| A_{1} \right|}^{2} \overline{C_{1}} + 2 \, A_{1}^{2} B_{1} \overline{C_{1}} + 2 \, A_{1}^{2} C_{1} \overline{C_{2}} & 4 & -1 & 0 \\
    			A_{1}^{2} \overline{C_{1}}^{2} & 6 & -4 & 0 \\
    			-4 \, A_{0} A_{1} \alpha_{1} \overline{C_{1}}^{2} + 2 \, A_{1}^{2} \overline{B_{1}} \overline{C_{1}} + 4 \, A_{1} A_{2} \overline{C_{1}}^{2} & 7 & -4 & 0 \\
    			-4 \, A_{0} A_{1} {\left| A_{1} \right|}^{2} \overline{C_{1}}^{2} + 2 \, A_{1}^{2} \overline{C_{1}} \overline{C_{2}} & 6 & -3 & 0
    			\end{dmatrix}
    		\end{align}
    		\normalsize
    		where
    		\begin{align*}
    			\chi_1&=\left(-2 \, A_{0} C_{1} \alpha_{1}\right) \left(A_{1} C_{1}\right) + \left(-2 \, A_{0} C_{1} \alpha_{1}\right) \left(A_{1} C_{1}\right) + \left(A_{1} C_{1}\right) \left(2 \, A_{2} C_{1}\right) + \left(A_{1} C_{1}\right) \left(2 \, A_{2} C_{1}\right) + \left(A_{1} C_{1}\right) \left(A_{1} C_{2}\right)\\
    			& + \left(A_{1} C_{1}\right) \left(A_{1} C_{2}\right)
    			+ \left(A_{1} \overline{C_{1}}\right) \left(-\frac{1}{8} \, C_{1} \overline{C_{1}}\right) + \left(A_{1} \overline{C_{1}}\right) \left(-\frac{1}{8} \, C_{1} \overline{C_{1}}\right)\\
    			\chi_2&=\left(-2 \, A_{0} \alpha_{1} \overline{C_{1}}\right) \left(A_{1} C_{1}\right) + \left(-2 \, A_{0} \alpha_{1} \overline{C_{1}}\right) \left(A_{1} C_{1}\right) + \left(A_{1} C_{1}\right) \left(A_{1} \overline{B_{1}}\right) + \left(A_{1} C_{1}\right) \left(A_{1} \overline{B_{1}}\right) + \left(-2 \, A_{0} C_{1} \alpha_{1}\right) \left(A_{1} \overline{C_{1}}\right)\\
    			& + \left(-2 \, A_{0} C_{1} \alpha_{1}\right) \left(A_{1} \overline{C_{1}}\right) + \left(2 \, A_{2} C_{1}\right) \left(A_{1} \overline{C_{1}}\right) + \left(2 \, A_{2} C_{1}\right) \left(A_{1} \overline{C_{1}}\right) + \left(A_{1} C_{2}\right) \left(A_{1} \overline{C_{1}}\right) + \left(A_{1} C_{2}\right) \left(A_{1} \overline{C_{1}}\right)\\
    			& + \left(A_{1} C_{1}\right) \left(2 \, A_{2} \overline{C_{1}}\right)
    			+ \left(A_{1} C_{1}\right) \left(2 \, A_{2} \overline{C_{1}}\right)\\
    			\chi_3&=\left(-2 \, A_{0} {\left| A_{1} \right|}^{2} \overline{C_{1}}\right) \left(A_{1} C_{1}\right) + \left(-2 \, A_{0} {\left| A_{1} \right|}^{2} \overline{C_{1}}\right) \left(A_{1} C_{1}\right) + \left(-2 \, A_{0} C_{1} {\left| A_{1} \right|}^{2}\right) \left(A_{1} \overline{C_{1}}\right)\\
    			& + \left(-2 \, A_{0} C_{1} {\left| A_{1} \right|}^{2}\right) \left(A_{1} \overline{C_{1}}\right) + \left(A_{1} B_{1}\right) \left(A_{1} \overline{C_{1}}\right) + \left(A_{1} B_{1}\right) \left(A_{1} \overline{C_{1}}\right) + \left(A_{1} C_{1}\right) \left(A_{1} \overline{C_{2}}\right) + \left(A_{1} C_{1}\right) \left(A_{1} \overline{C_{2}}\right)\\
    			\chi_4&=\left(-2 \, A_{0} \alpha_{1} \overline{C_{1}}\right) \left(A_{1} \overline{C_{1}}\right) + \left(-2 \, A_{0} \alpha_{1} \overline{C_{1}}\right) \left(A_{1} \overline{C_{1}}\right) + \left(A_{1} \overline{B_{1}}\right) \left(A_{1} \overline{C_{1}}\right) + \left(A_{1} \overline{B_{1}}\right) \left(A_{1} \overline{C_{1}}\right)\\
    			& + \left(A_{1} \overline{C_{1}}\right) \left(2 \, A_{2} \overline{C_{1}}\right) + \left(A_{1} \overline{C_{1}}\right) \left(2 \, A_{2} \overline{C_{1}}\right)
    		\end{align*}
    		We need only $\pi_{15}$ for the coefficient in 
    		\begin{align*}
    			z^{\theta_0}\z^{2-\theta_0}\,dz^4=z^{4}\z^{-2}\,dz^4
    		\end{align*}
    		which is
    		\begin{align*}
    	        &\pi_{15}=-\colorcancel{80 \, A_{1}^{2} {\left| A_{1} \right|}^{4}}{blue} - \colorcancel{160 \, A_{0} A_{2} {\left| A_{1} \right|}^{4}}{blue} - \ccancel{16 \, A_{0} A_{1} \alpha_{4} {\left| A_{1} \right|}^{2}} - \ccancel{48 \, A_{0} A_{1} \alpha_{12}} - \ccancel{144 \, A_{0} C_{1} \alpha_{7}} + 2 \, {\left(\ccancel{80 \, A_{0}^{2} {\left| A_{1} \right|}^{4}} + \ccancel{3 \, A_{1} \overline{B_{1}}}\right)} \alpha_{1}\\
    	        & + 6 \, {\left(\ccancel{8 \, A_{0} A_{1} \overline{\alpha_{1}}} + A_{0} \overline{B_{1}}\right)} \alpha_{3} - 6 \, A_{3} \overline{B_{1}} + 6 \, A_{1} \overline{B_{4}} - 15 \, {\left(7 \, A_{1} C_{1} + 3 \, A_{0} C_{2}\right)} \overline{\alpha_{2}} - 48 \, {\left(\ccancel{A_{0} A_{1} \alpha_{1}} + \ccancel{A_{0}^{2} \alpha_{3}} - \colorcancel{A_{1} A_{2}}{blue} - \colorcancel{A_{0} A_{3}}{blue}\right)} \overline{\alpha_{4}}\\
    	        &=6\s{\alpha_3\vec{A}_0-\vec{A}_3}{\bar{\vec{B}_1}}+6\s{\vec{A}_1}{\bar{\vec{B}_4}}-60\bar{\alpha_2}\s{\vec{A}_1}{\vec{C}_1}
    		\end{align*}
    		Now, recall that
    		\begin{align*}
    			\alpha_3=\frac{1}{12}\s{\vec{A}_1}{\vec{C}_1}+2\s{\bar{\vec{A}_0}}{\vec{A}_3}
    		\end{align*}
    		and
    		\begin{align*}
    		    &|\vec{A}_0|^2=\dfrac{1}{2},\qquad
    			\vec{B}_1=-2\s{\bar{\vec{A}_1}}{\vec{C}_1}\vec{A}_0
    		\end{align*}
    		so
    		\begin{align*}
    			&\alpha_2=\frac{1}{10} \, B_{1} \overline{A_{0}} + \frac{1}{8} \, C_{1} \overline{A_{1}}=-\frac{1}{10}\s{\bar{\vec{A}_1}}{\vec{C}_1}+\frac{1}{8}\s{\bar{\vec{A}_1}}{\vec{C}_1}=\frac{1}{40}\s{\bar{\vec{A}_1}}{\vec{C}_1}\\
    			&\s{\alpha_3\vec{A}_0-\vec{A}_3}{\bar{\vec{A}_0}}=\frac{\alpha_3}{2}-\s{\bar{\vec{A}_0}}{\vec{A}_3}=\frac{1}{24}\s{\vec{A}_1}{\vec{C}_1}\\
    			&\s{\alpha_3\vec{A}_0-\vec{A}_3}{\bar{\vec{B}_1}}=-\frac{1}{12}\s{\vec{A}_1}{\vec{C}_1}\s{{\vec{A}_1}}{\bar{\vec{C}_1}}.
    		\end{align*}
    		Now, we have
    		\begin{align*}
    			&\frac{1}{2}\vec{B}_4=\begin{dmatrix}
    			-\frac{7}{48} \, C_{1} \overline{C_{1}} & \overline{A_{1}} & -2 & 3 & 0 &\textbf{(23)}\\
    			A_{0} \alpha_{2} \overline{C_{1}} - \frac{2}{15} \, B_{1} \overline{C_{1}} - \frac{7}{48} \, C_{1} \overline{C_{2}} & \overline{A_{0}} & -2 & 3 & 0 &\textbf{(24)}\\
    			-\frac{1}{48} \, \overline{A_{1}} \overline{C_{1}} & C_{1} & -2 & 3 & 0 &\textbf{(25)}\\
    			\frac{5}{3} \, \alpha_{2} \overline{A_{0}} \overline{C_{1}} + C_{1} \overline{A_{0}} \overline{\alpha_{3}} - \frac{1}{3} \, B_{2} \overline{A_{1}} - \frac{2}{3} \, B_{1} \overline{A_{2}} - C_{1} \overline{A_{3}} + \frac{1}{3} \, {\left(2 \, B_{1} \overline{A_{0}} + 3 \, C_{1} \overline{A_{1}}\right)} \overline{\alpha_{1}} & A_{0} & -2 & 3 & 0 &\textbf{(26)}\\
    			-\frac{1}{24} \, C_{1} \overline{A_{1}} & \overline{C_{1}} & -2 & 3 & 0 &\textbf{(27)}.
    			\end{dmatrix}
    		\end{align*}
    		so there exists $\lambda_1,\lambda_2\in\mathbb{C}$ such that
    		\begin{align*}
    			\vec{B}_4=-\frac{7}{24}|\vec{C}_1|^2\bar{\vec{A}_1}-\frac{1}{24}\bar{\s{\vec{A}_1}{\vec{C}_1}}\vec{C}_1-\frac{1}{12}\s{\bar{\vec{A}_1}}{\vec{C}_1}\bar{\vec{C}_1}+\lambda_1\vec{A}_0+\lambda_2\bar{\vec{A}_0}.
    		\end{align*}
    		and we have
    		\begin{align*}
    			\s{\bar{\vec{A}_1}}{\vec{B}_4}&=-\frac{7}{24}|\vec{C}_1|^2\bar{\s{\vec{A}_1}{\vec{A}_1}}-\frac{1}{24}\bar{\s{\vec{A}_1}{\vec{C}_1}}\s{\bar{\vec{A}_1}}{\vec{C}_1}-\frac{1}{12}\s{\bar{\vec{A}_1}}{\vec{C}_1}\bar{\s{\vec{A}_1}{\vec{C}_1}}\\
    			&=-\frac{7}{24}|\vec{C}_1|^2\bar{\s{\vec{A}_1}{\vec{A}_1}}-\frac{1}{8}\s{\bar{\vec{A}_1}}{\vec{C}_1}\bar{\s{\vec{A}_1}{\vec{C}_1}}.
    		\end{align*}
    		Therefore, we deduce that
    		\begin{align}\label{4end2}
    			\pi_{15}&=6\s{\alpha_3\vec{A}_0-\vec{A}_3}{\bar{\vec{B}_1}}+6\s{\vec{A}_1}{\bar{\vec{B}_4}}-60\bar{\alpha_2}\s{\vec{A}_1}{\vec{C}_1}\nonumber\\
    			&=-\frac{1}{2}\s{\vec{A}_1}{\vec{C}_1}\s{{\vec{A}_1}}{\bar{\vec{C}_1}}-\frac{7}{4}|\vec{C}_1|^2\s{\vec{A}_1}{\vec{A}_1}-\frac{3}{4}\s{\vec{A}_1}{\vec{C}_1}\s{\vec{A}_1}{\bar{\vec{C}_1}}-\frac{3}{2}\s{\vec{A}_1}{\vec{C}_1}\s{\vec{A}_1}{\bar{\vec{C}_1}}\nonumber\\
    			&=-\frac{7}{4}|\vec{C}_1|^2\s{\vec{A}_1}{\vec{A}_1}-\frac{11}{4}\s{\vec{A}_1}{\vec{C}_1}\s{\vec{A}_1}{\bar{\vec{C}_1}}.
    		\end{align}
    		Finally, thanks of \eqref{4end0} and \eqref{4end1}, the coefficient in $z^4\z^{-2}dz^4$ is
    		\begin{align*}
    			\frac{5}{2}|\vec{C}_1|^2\s{\vec{A}_1}{\vec{A}_1}+2\s{\vec{A}_1}{\vec{C}_1}\s{\vec{A}_1}{\bar{\vec{C}_1}}
    		\end{align*}
    		and finally the coefficient in $z^4\z^{-2}\,dz^4$ in $\mathscr{Q}_{\phi}$ is
    		\begin{align}
    			\Omega_0=\frac{3}{4}\left(|\vec{C}_1|^2\s{\vec{A}_1}{\vec{A}_1}-\s{\vec{A}_1}{\vec{C}_1}\s{\vec{A}_1}{\bar{\vec{C}_1}}\right)=0.
    		\end{align}
    		As the coefficient in $\dfrac{\z^4}{z}\log|z|$ in $\mathscr{Q}_{\phi}$ is
    		\begin{align*}
    			\Omega_1=-6 \, A_{1} C_{1} \overline{\zeta_{0}}
    		\end{align*}
    		As $\s{\vec{A}_0}{\vec{\gamma}_1}=0$, and 
    		\begin{align*}
    			\vec{E}_1=-\frac{1}{8}\s{\vec{C}_1}{\vec{C}_1}\bar{\vec{A}_0}
    		\end{align*}
    		so
    		\begin{align*}
    			\zeta_0=\frac{1}{4} \, A_{0} \gamma_{1} + \overline{A_{0}} \overline{E_{1}}=-\frac{1}{16}\s{\vec{C}_1}{\vec{C}_1}.
    		\end{align*}
    		so we finally obtain
    		\begin{align}
    			\Omega_1=\frac{3}{8}\s{\vec{A}_1}{\vec{C}_1}\bar{\s{\vec{C}_1}{\vec{C}_1}}=0
    		\end{align}
    		Thanks of the argument given at the beginning at the section, we are done.

	    		\section{The case where $\theta_0=3$}

	    		Recall the following development valid for all $\theta_0\geq 3$
	    		\begin{align}\label{31dev}
	    		\left\{\begin{alignedat}{1}
	    		\p{z}\phi&=\vec{A}_0z^{\theta_0-1}+\vec{A}_1z^{\theta_0}+\vec{A}_2z^{\theta_0+1}+\frac{1}{4\theta_0}\vec{C}_1z\,\z^{\theta_0}+\frac{1}{8}\bar{\vec{C}_1}z^{\theta_0-1}\z^2+O(|z|^{\theta_0+2-\epsilon})\\
	    		e^{2\lambda}&=e^{2\lambda}=|z|^{2\theta_0-2}\left(1+2|\vec{A}_1|^2|z|^2+2\,\Re\left(\alpha_0z+\alpha_1z^2\right)+O(|z|^{3-\epsilon})\right)\\
	    		\h_0&=2\left(\vec{A}_1-\alpha_0\vec{A}_0-\Big(2|\vec{A}_1|^2-|\alpha_0|^2\Big)\vec{A}_0\z\right)z^{\theta_0-1}+2\left(2\vec{A}_2-\alpha_0\vec{A}_1-(2\alpha_1-\alpha_0^2)\vec{A}_0\right)z^{\theta_0}\\
	    		&-\frac{(\theta_0-2)}{2\theta_0}\vec{C}_1\z^{\theta_0}+O(|z|^{\theta_0+1-\epsilon})\\
	    		\H&=\Re\left(\frac{\vec{C}_1}{z^{\theta_0-2}}\right)+O(|z|^{3-\theta_0-\epsilon})\\
	    		\mathscr{Q}_{\phi}&=(\theta_0-1)(\theta_0-2)\s{\vec{A}_1}{\vec{C}_1}\frac{dz^4}{z}+O(|z|^{-\epsilon}).
	    		\end{alignedat}\right.
	    		\end{align}
	    		Also, recall that by conformality of $\phi$, we have
	    		\begin{align*}
	    		0=\s{\p{z}\phi}{\p{z}\phi}&=\s{\vec{A}_0}{\vec{A}_0}z^{2\theta_0-2}+2\s{\vec{A}_0}{\vec{A}_1}z^{2\theta_0-1}+\left(\s{\vec{A}_1}{\vec{A}_1}+2\s{\vec{A}_0}{\vec{A}_2}\right)z^{2\theta_0}+\frac{1}{2\theta_0}\s{\vec{A}_0}{\vec{C}_1}|z|^{2\theta_0}\\
	    		&+\frac{1}{4}\s{\vec{A}_0}{\bar{\vec{C}_1}}z^{2\theta_0-2}\z^2+O(|z|^{2\theta_0+1-\epsilon})
	    		\end{align*}
	    		so that
	    		\begin{align}\label{3cancel1}
	    		\s{\vec{A}_0}{\vec{A}_0}=\s{\vec{A}_0}{\vec{A}_1}=\s{\vec{A}_0}{\vec{C}_1}=\s{\vec{A}_0}{\bar{\vec{C}_1}}=0.
	    		\end{align}
	    		Now, we take $\theta_0=3$ in \eqref{31dev}, we obtain
	    		\begin{align*}
	    		\left\{\begin{alignedat}{1}
	    		\p{z}\phi&=\vec{A}_0z^{2}+\vec{A}_1z^{3}+\vec{A}_2z^{4}+\frac{1}{12}\vec{C}_1z\z^3+\frac{1}{8}\bar{\vec{C}_1}|z|^4+O(|z|^{5-\epsilon})\\
	    		e^{2\lambda}&=|z|^{4}+2|\vec{A}_1|^2|z|^6+2\,\Re(\alpha_0z^3\z^2+\alpha_1z^4\z^2)+O(|z|^{7-\epsilon})\\
	    		\h_0&=2\left(\vec{A}_1-\alpha_0\vec{A}_0-\Big(2|\vec{A}_1|^2-|\alpha_0|^2\Big)\vec{A}_0\z\right)z^{2}+2\left(2\vec{A}_2-\alpha_0\vec{A}_1-(2\alpha_1-\alpha_0^2)\vec{A}_0\right)z^{3}-\frac{1}{6}\vec{C}_1\z^{3}+O(|z|^{4-\epsilon})\\
	    		\H&=\Re\left(\frac{\vec{C}_1}{z}\right)+O(|z|^{-\epsilon})\\
	    		\mathscr{Q}_{\phi}&=2\s{\vec{A}_1}{\vec{C}_1}\frac{dz^4}{z}+O(|z|^{-\epsilon}).
	    		\end{alignedat}\right.
	    		\end{align*}
	    		And as logarithm will appear in the forthcoming computations, the last column will indicate the power of the logarithm. Furthermore, as this is a new source of possible mistakes (the algorithms will change to include logarithms), we will only start from
	    		\begin{align}\label{30}
	    		\left\{\begin{alignedat}{1}
	    		\p{z}\phi&=\vec{A}_0z^{2}+\vec{A}_1z^{3}+\vec{A}_2z^{4}+\frac{1}{12}\vec{C}_1z\z^3+\frac{1}{8}\bar{\vec{C}_1}|z|^4+O(|z|^{5-\epsilon})\\
	    		e^{2\lambda}&=|z|^{4}+2|\vec{A}_1|^2|z|^6+2\,\Re(\alpha_0z^3\z^2+\alpha_1z^4\z^2)+O(|z|^{7-\epsilon})\\
	    		\H&=\Re\left(\frac{\vec{C}_1}{z}\right)+O(|z|^{-\epsilon})
	    		\end{alignedat}\right.
	    		\end{align}
	    		and recompute $\h_0$.
	    		First, we check that there is no discrepancy in our transcription.
	    		\begin{align}\label{3sage1}
	    		\p{z}\phi=\begin{dmatrix}
	    		1 & A_{0} & 2 & 0 & 0 \\
	    		1 & A_{1} & 3 & 0 & 0 \\
	    		1 & A_{2} & 4 & 0 & 0 \\
	    		\frac{1}{12} & C_{1} & 1 & 3 & 0 \\
	    		\frac{1}{8} & \overline{C_{1}} & 2 & 2 & 0
	    		\end{dmatrix},\quad 
	    		e^{2\lambda}=\begin{dmatrix}
	    		1 & 2 & 2 & 0 \\
	    		2 \, {\left| A_{1} \right|}^{2} & 3 & 3 & 0 \\
	    		\alpha_{0} & 3 & 2 & 0 \\
	    		\overline{\alpha_{0}} & 2 & 3 & 0 \\
	    		\alpha_{1} & 4 & 2 & 0 \\
	    		\overline{\alpha_{1}} & 2 & 4 & 0
	    		\end{dmatrix},\quad
	    		\H=\begin{dmatrix}
	    		\frac{1}{2} & C_{1} & -1 & 0 & 0 \\
	    		\frac{1}{2} & \overline{C_{1}} & 0 & -1 & 0
	    		\end{dmatrix}
	    		\end{align}
	    		which checks with \eqref{30}. Now, we compute the new expression of $\h_0$ as
	    		\begin{align}\label{b1}
	    		\h_0&=\begin{dmatrix}
	    		2 & A_{1} & 2 & 0 & 0 \\
	    		4 & A_{2} & 3 & 0 & 0 \\
	    		-\frac{1}{6} & C_{1} & 0 & 3 & 0 \\
	    		-4 \, {\left| A_{1} \right|}^{2} + 2 \, \alpha_{0} \overline{\alpha_{0}} & A_{0} & 2 & 1 & 0 \\
	    		-2 \, \alpha_{0} & A_{0} & 2 & 0 & 0 \\
	    		-2 \, \alpha_{0} & A_{1} & 3 & 0 & 0 \\
	    		2 \, \alpha_{0}^{2} - 4 \, \alpha_{1} & A_{0} & 3 & 0 & 0
	    		\end{dmatrix}
	    		\end{align}
	    		to be compared with
	    		\begin{align}\label{b2}
	    		\h_0&=2\left(\vec{A}_1-\alpha_0\vec{A}_0-\Big(2|\vec{A}_1|^2-|\alpha_0|^2\Big)\vec{A}_0\z\right)z^{2}+2\left(2\vec{A}_2-\alpha_0\vec{A}_1-(2\alpha_1-\alpha_0^2)\vec{A}_0\right)z^{3}-\frac{1}{6}\vec{C}_1\z^{3}+O(|z|^{4-\epsilon}).
	    		\end{align}
	    		First, there are $7$ different terms (of the form ${\lambda}\vec{\Lambda}z^{\alpha}\z^{\beta}$ for some $\lambda\in\mathbb{C},\vec{\Lambda}\in\C^n,\alpha\in\Z,\beta\in\Z$) in the two expressions, and we easily check that each coefficient coincide between \eqref{b1} and \eqref{b2}.
	    		Now, recall the fundamental equation
	    		\begin{align}\label{3boot1}
	    		\partial\left(\H-2i\vec{L}+\vec{\gamma}_0\log|z|\right)=-|\H|^2\partial\phi-2\,g^{-1}\otimes\s{\H}{\h_0}\otimes\bar{\partial}\phi
	    		\end{align}
	    		If $\vec{Q}\in C^{\infty}(D^2\setminus\ens{0},\C^n)$ is the unique anti-holomorphic free of the equation
	    		\begin{align*}
	    		\partial\vec{Q}=-|\H|^2\partial\phi-2\,g^{-1}\otimes\s{\H}{\h_0}\otimes\bar{\partial}\phi
	    		\end{align*}
	    		we deduce that
	    		\begin{align*}
	    		\partial\left(\vec{H}-2i\vec{L}-\vec{Q}\right)=0,
	    		\end{align*}
	    		so there exists some $\bar{\vec{D}_2}\in \mathbb{C}^n$ such that
	    		\begin{align*}
	    		\H-2i\vec{L}=\frac{\bar{\vec{C}_1}}{z}+\bar{\vec{D}_2}+\vec{Q}+O(|z|^{1-\epsilon}).
	    		\end{align*}
	    		so that as $\vec{L}$ is \emph{real}, one obtains
	    		\begin{align}\label{3Htemp1}
	    		\H=\Re\left(\frac{\vec{C}_1}{z}\right)+\Re\left(\vec{D}_2\right)+\Re\left(\vec{Q}\right)+O(|z|^{1-\epsilon})
	    		\end{align}
	    		and we have
	    		\begin{align}\label{3Htemp2}
	    		\Re(\vec{Q})&=\begin{dmatrix}
	    		2 \, {\left(\ccancel{A_{0} \alpha_{0}} - A_{1}\right)} C_{1} & \overline{A_{0}} & 0 & 0 & 1 \\
	    		2 \, {\left(\ccancel{\overline{A_{0}} \overline{\alpha_{0}}} - \overline{A_{1}}\right)} \overline{C_{1}} & A_{0} & 0 & 0 & 1 \\
	    		{\left(\ccancel{A_{0} \alpha_{0}} - A_{1}\right)} \overline{C_{1}} & \overline{A_{0}} & 1 & -1 & 0 \\
	    		{\left(\ccancel{\overline{A_{0}} \overline{\alpha_{0}}} - \overline{A_{1}}\right)} C_{1} & A_{0} & -1 & 1 & 0
	    		\end{dmatrix}+O(|z|^{1-\epsilon})
	    		=\begin{dmatrix}
	    		-2\, A_{1} C_{1} & \overline{A_{0}} & 0 & 0 & 1 \\
	    		-2 \, \overline{A_{1}} \overline{C_{1}} & A_{0} & 0 & 0 & 1 \\
	    		-A_{1} \overline{C_{1}} & \overline{A_{0}} & 1 & -1 & 0 \\
	    		- \overline{A_{1}} C_{1} & A_{0} & -1 & 1 & 0
	    		\end{dmatrix}+O(|z|^{1-\epsilon})\nonumber\\
	    		&=-2\,\Re\left(\s{\bar{\vec{A}_1}}{\vec{C}_1}\vec{A}_0\frac{\z}{z}\right)-4\,\Re\left(\s{\vec{A}_1}{\vec{C}_1}\bar{\vec{A}_0}\right)\log|z|+O(|z|^{1-\epsilon})
	    		\end{align}
	    		as $\s{\vec{A}_0}{\vec{C}_1}=\s{\bar{\vec{A}_0}}{\vec{C}_1}=0$ by \eqref{3cancel1}. Therefore, if we define
	    		\begin{align}\label{3defconts1}
	    		\left\{
	    		\begin{alignedat}{1}
	    		\vec{C}_2&=\Re\left(\vec{D}_2\right)\in\R^n\\
	    		\vec{B}_1&=-2\s{\bar{\vec{A}_1}}{\vec{C}_1}\vec{A}_0\in\mathbb{C}^n\\
	    		\vec{\gamma}_1&=-\vec{\gamma}_0-4\,\Re\left(\s{\vec{A}_1}{\vec{C}_1}\bar{\vec{A}_0}\right)\in\mathbb{R}^n
	    		\end{alignedat}\right.
	    		\end{align}
	    		we obtain by \eqref{3Htemp1} and \eqref{3Htemp2}
	    		\begin{align}\label{3h2}
	    		\H=\Re\left(\frac{\vec{C}_1}{z}+\vec{B}_1\frac{\z}{z}\right)+\vec{C}_2+\vec{\gamma}_1\log|z|+O(|z|^{1-\epsilon}).
	    		\end{align}
	    		This is also something that we can easily obtain by hand, as we only need the first order development of the tensors. Let us check this.
	    		
	    		First, as for all $\theta_0\geq 3$, we have $\H=O(|z|^{2-\theta_0})$, and $\p{z}\phi=O(|z|^{\theta_0-1})$, we have
	    		\begin{align}\label{osc}
	    		|\H|^2\p{z}\phi=O(|z|^{3-\theta_0}).
	    		\end{align}
	    		Furthermore, we have
	    		\begin{align*}
	    		\h_0&=2\left(\vec{A}_1-\alpha_0\vec{A}_0\right)z^{\theta_0-1}dz^2+O(|z|^{\theta_0}),\qquad
	    		\H=\frac{1}{2}\frac{\vec{C}_1}{z^{\theta_0-2}}+\frac{1}{2}\frac{\bar{\vec{C}_1}}{\z^{\theta_0-2}}+O(|z|^{3-\theta_0-\epsilon})
	    		\end{align*}
	    		so (as $\s{\vec{A}_0}{\vec{C}_1}=\s{\bar{\vec{A}_0}}{\vec{C}_1}=0$)
	    		\begin{align}\label{genpart1}
	    		\s{\H}{\h_0}=\s{\vec{A}_1}{\vec{C}_1}z\,dz+\s{\vec{A}_1}{\bar{\vec{C}_1}}z^{\theta_0-1}\z^{2-\theta_0}\, dz+O(|z|^{2})
	    		\end{align}
	    		Now, as
	    		\begin{align*}
	    		\p{z}\phi=\vec{A}_0z^{\theta_0-1}+O(|z|^{\theta_0}),\quad e^{2\lambda}=|z|^{2\theta_0-2}+O(|z|^{2\theta_0-1})
	    		\end{align*}
	    		we trivially have
	    		\begin{align}\label{genpart2}
	    		e^{-2\lambda}\p{\z}\phi=\bar{\vec{A}_0}z^{1-\theta_0}+O(|z|^{2-\theta_0-\epsilon}).
	    		\end{align}
	    		Finally, by \eqref{genpart1} and \eqref{genpart2}, we have
	    		\begin{align*}
	    		g^{-1}\otimes\s{\H}{\h_0}\otimes\bar{\partial}\phi=\s{\vec{A}_1}{\vec{C}_1}\bar{\vec{A}_0}z^{2-\theta_0}\,dz+\s{\vec{A}_1}{\bar{\vec{C}_1}}\bar{\vec{A}_0}\z^{2-\theta_0}\,dz+O(|z|^{3-\theta_0})
	    		\end{align*}
	    		so we obtain by \eqref{osc} the equation
	    		\begin{align}\label{3boot2}
	    		\partial\left(\vec{H}-2i\vec{L}+\vec{\gamma}_0\log|z|\right)&=-|\H|^2\partial\phi-2\,g^{-1}\otimes\s{\H}{\h_0}\otimes\bar{\partial}\phi\\
	    		&=-2\s{\vec{A}_1}{\vec{C}_1}z^{2-\theta_0}-2\s{\vec{A}_1}{\bar{\vec{C}_1}}\bar{\vec{A}_0}\z^{2-\theta_0}+O(|z|^{3-\theta_0-\epsilon}).
	    		\end{align}
	    		Taking $\theta_0=3$ in \eqref{3boot2} yields
	    		\begin{align*}
	    		\partial\left(\vec{H}-2i\vec{L}+\vec{\gamma}_0\log|z|\right)=-2\s{\vec{A}_1}{\vec{C}_1}\bar{\vec{A}_0}\frac{dz}{z}-2\s{\vec{A}_1}{\bar{\vec{C}_1}}\frac{dz}{\z}+O(|z|^{-\epsilon}).
	    		\end{align*}
	    		so for some $\bar{\vec{D}_2}\in\mathbb{C}^n$, we have
	    		\begin{align*}
	    		\H-2i\vec{L}+\vec{\gamma}_0\log|z|=-4\s{\vec{A}_1}{\vec{C}_1}\bar{\vec{A}_0}\log|z|-2\s{\vec{A}_1}{\bar{\vec{C}_1}}\bar{\vec{A}_0}\frac{z}{\z}+O(|z|^{1-\epsilon})
	    		\end{align*}
	    		which immediately gives by taking the real part the development in \eqref{3h2} thanks of \eqref{3defconts1}.
	    		
	    		Now we check that the transcription was correct :
	    		\begin{align*}
	    		\H=\begin{dmatrix}
	    		\frac{1}{2} & C_{1} & -1 & 0 & 0 \\
	    		\frac{1}{2} & B_{1} & -1 & 1 & 0 \\
	    		1 & C_{2} & 0 & 0 & 0 \\
	    		\end{dmatrix}
	    		\begin{dmatrix}
	    		1 & \gamma_{1} & 0 & 0 & 1 \\
	    		\frac{1}{2} & \overline{C_{1}} & 0 & -1 & 0 \\
	    		\frac{1}{2} & \overline{B_{1}} & 1 & -1 & 0
	    		\end{dmatrix}
	    		\end{align*}
	    		We know from the previous step that we do not have to develop to the next order $\H$ to obtain $\vec{A}_1=0$, and by integrating the equation
	    		\begin{align*}
	    		\p{\z}\left(\p{z}\phi\right)=\frac{e^{2\lambda}}{2}\H
	    		\end{align*}
	    		we obtain for some $\vec{A}_3\in\mathbb{C}^n$
	    		\begin{align*}
	    		\p{z}\phi=\begin{dmatrix}
	    		1 & A_{0} & 2 & 0 & 0 \\
	    		1 & A_{1} & 3 & 0 & 0 \\
	    		1 & A_{2} & 4 & 0 & 0 \\
	    		1 & A_{3} & 5 & 0 & 0 \\
	    		\frac{1}{12} & C_{1} & 1 & 3 & 0 \\
	    		\frac{1}{16} & B_{1} & 1 & 4 & 0 \\
	    		\frac{1}{16} \, \overline{\alpha_{0}} & C_{1} & 1 & 4 & 0 \\
	    		\frac{1}{8} & \overline{C_{1}} & 2 & 2 & 0 
	    		\end{dmatrix}
	    		\begin{dmatrix}
	    		-\frac{1}{36} & \gamma_{1} & 2 & 3 & 0 \\
	    		\frac{1}{6} & C_{2} & 2 & 3 & 0 \\
	    		\frac{1}{12} \, \overline{\alpha_{0}} & \overline{C_{1}} & 2 & 3 & 0 \\
	    		\frac{1}{12} \, \alpha_{0} & C_{1} & 2 & 3 & 0 \\
	    		\frac{1}{6} & \gamma_{1} & 2 & 3 & 1 \\
	    		\frac{1}{8} & \overline{B_{1}} & 3 & 2 & 0 \\
	    		\frac{1}{8} \, \alpha_{0} & \overline{C_{1}} & 3 & 2 & 0.
	    		\end{dmatrix}+O(|z|^{6-\epsilon})
	    		\end{align*}
	    		and by cutting this development to one order less we recover
	    		\begin{align*}
	    		\p{z}\phi=\begin{dmatrix}
	    		1 & A_{0} & 2 & 0 & 0 \\
	    		1 & A_{1} & 3 & 0 & 0 \\
	    		1 & A_{2} & 4 & 0 & 0 \\
	    		\frac{1}{12} & C_{1} & 1 & 3 & 0 \\
	    		\frac{1}{8} & \overline{C_{1}} & 2 & 2 & 0
	    		\end{dmatrix}
	    		\end{align*}
	    		as expected (see \eqref{3sage1}). As $\phi$ is real, we also have by direct integration
	    		\begin{align}
	    		\phi(z)=\begin{dmatrix}
	    		\frac{1}{3} & \overline{A_{0}} & 0 & 3 & 0 \\
	    		\frac{1}{4} & \overline{A_{1}} & 0 & 4 & 0 \\
	    		\frac{1}{5} & \overline{A_{2}} & 0 & 5 & 0 \\
	    		\frac{1}{6} & \overline{A_{3}} & 0 & 6 & 0 \\
	    		\frac{1}{3} & A_{0} & 3 & 0 & 0 \\
	    		\frac{1}{4} & A_{1} & 4 & 0 & 0 \\
	    		\frac{1}{5} & A_{2} & 5 & 0 & 0 \\
	    		\frac{1}{6} & A_{3} & 6 & 0 & 0 \\
	    		\frac{1}{24} & C_{1} & 2 & 3 & 0 \\
	    		\frac{1}{32} & B_{1} & 2 & 4 & 0 
	    		\end{dmatrix}
	    		\begin{dmatrix}
	    		\frac{1}{32} \, \overline{\alpha_{0}} & C_{1} & 2 & 4 & 0 \\
	    		\frac{1}{24} & \overline{C_{1}} & 3 & 2 & 0 \\
	    		-\frac{1}{54} & \gamma_{1} & 3 & 3 & 0 \\
	    		\frac{1}{18} & C_{2} & 3 & 3 & 0 \\
	    		\frac{1}{36} \, \overline{\alpha_{0}} & \overline{C_{1}} & 3 & 3 & 0 \\
	    		\frac{1}{36} \, \alpha_{0} & C_{1} & 3 & 3 & 0 \\
	    		\frac{1}{18} & \gamma_{1} & 3 & 3 & 1 \\
	    		\frac{1}{32} & \overline{B_{1}} & 4 & 2 & 0 \\
	    		\frac{1}{32} \, \alpha_{0} & \overline{C_{1}} & 4 & 2 & 0
	    		\end{dmatrix}
	    		\end{align}
	    		
	    		By conformality, we  have
	    		\begin{align*}
	    		0=\s{\p{z}\phi}{\p{z}\phi}=\begin{dmatrix}
	    		A_{0}^{2} & 4 & 0 & 0 \\
	    		2 \, A_{0} A_{1} & 5 & 0 & 0 \\
	    		A_{1}^{2} + 2 \, A_{0} A_{2} & 6 & 0 & 0 \\
	    		2 \, A_{1} A_{2} + 2 \, A_{0} A_{3} & 7 & 0 & 0 \\
	    		\frac{1}{6} \, A_{0} C_{1} & 3 & 3 & 0 \\
	    		\frac{1}{8} \, A_{0} C_{1} \overline{\alpha_{0}} + \frac{1}{8} \, A_{0} B_{1} & 3 & 4 & 0 \\
	    		\frac{1}{4} \, A_{0} \overline{C_{1}} & 4 & 2 & 0 \\
	    		\frac{1}{6} \, A_{0} \overline{C_{1}} \overline{\alpha_{0}} + \frac{1}{6} \, {\left(A_{0} \alpha_{0} + A_{1}\right)} C_{1} + \frac{1}{3} \, A_{0} C_{2} - \frac{1}{18} \, A_{0} \gamma_{1} & 4 & 3 & 0 \\
	    		\frac{1}{3} \, A_{0} \gamma_{1} & 4 & 3 & 1 \\
	    		\frac{1}{4} \, A_{0} \overline{B_{1}} + \frac{1}{4} \, {\left(A_{0} \alpha_{0} + A_{1}\right)} \overline{C_{1}} & 5 & 2 & 0
	    		\end{dmatrix}
	    		\end{align*}
	    		In particular, we deduce that
	    		\begin{align}\label{3conf1}
	    		\left\{
	    		\begin{alignedat}{1}
	    		&\s{\vec{A}_0}{\vec{A}_0}=\s{\vec{A}_0}{\vec{A}_1}=\s{\vec{A}_0}{\vec{C}_1}=\s{\vec{A}_0}{\bar{\vec{C}_1}}=\s{\vec{A}_0}{\vec{\gamma}_1}=0\\
	    		&\s{{\vec{A}_1}}{\vec{A}_1}+2\s{\vec{A}_0}{\vec{A}_2}=0,\quad \s{\vec{A}_1}{\vec{C}_1}+2\s{\vec{A}_0}{\vec{C}_2}=0.
	    		\end{alignedat}\right.
	    		\end{align}
	    		Remark as $\s{\vec{A}_0}{\vec{A}_0}=0$ that
	    		\begin{align}\label{3true}
	    		0=\s{\vec{A}_0}{\vec{\gamma_1}}=-\bs{\vec{A}_0}{\vec{\gamma}_0+4\,\Re\left(\s{\vec{A}_1}{\vec{C}_1}\bar{\vec{A}_0}\right)}=-\left(\s{\vec{A}_0}{\vec{\gamma}_0}+\s{\vec{A}_1}{\vec{C}_1}\right)
	    		\end{align}
	    		so for a \emph{true} Willmore sphere, we have $\vec{\gamma}_0=0$, and
	    		\begin{align*}
	    		\s{\vec{A}_1}{\vec{C}_1}=0,
	    		\end{align*}
	    		which this proves the holomorphy of the quartic form. Now, let us come back to the general case where no hypothesis is made on $\vec{\gamma}_0$.
	    		
	    		Now, the expansion of the metric is by \eqref{3cancel1} (notice that $\vec{\gamma}_1$ and $\vec{C}_2$ are \emph{real} and $\vec{B}_1\in\mathrm{Span}(\vec{A}_0)$)
	    		\begin{align*}
	    		e^{2\lambda}&=\begin{dmatrix}
	    		2 \, A_{0} \overline{A_{0}} & 2 & 2 & 0 \\
	    		2 \, A_{0} \overline{A_{1}} & 2 & 3 & 0 \\
	    		2 \, A_{0} \overline{A_{2}} + \frac{1}{4} \, \ccancel{\overline{A_{0}} \overline{C_{1}}} & 2 & 4 & 0 \\
	    		\frac{1}{6}\ccancel{ \, C_{1} \alpha_{0} \overline{A_{0}}} + \frac{1}{3} \, C_{2} \overline{A_{0}} - \frac{1}{18} \, \ccancel{\gamma_{1} \overline{A_{0}}} + 2 \, A_{0} \overline{A_{3}} + \frac{1}{12} \, {\left(\ccancel{2 \, \overline{A_{0}} \overline{\alpha_{0}}} + 3 \, \overline{A_{1}}\right)} \overline{C_{1}} & 2 & 5 & 0 \\
	    		\frac{1}{6} \, \ccancel{A_{0} \overline{C_{1}}} & 5 & 1 & 0 \\
	    		\frac{1}{8} \, A_{0} \overline{B_{1}} + \frac{1}{24} \, {\left(\ccancel{3 \, A_{0} \alpha_{0}} + 4 \, A_{1}\right)} \overline{C_{1}} & 6 & 1 & 0 \\
	    		\frac{1}{4} \, \ccancel{A_{0} C_{1}} + 2 \, A_{2} \overline{A_{0}} & 4 & 2 & 0 \\
	    		\frac{1}{6} \, \ccancel{A_{0} \overline{C_{1}}} \overline{\alpha_{0}} + \frac{1}{12} \, {\left(\ccancel{2 \, A_{0} \alpha_{0}} + 3 \, A_{1}\right)} C_{1} + \frac{1}{3} \, A_{0} C_{2} - \frac{1}{18} \, \ccancel{A_{0} \gamma_{1}} + 2 \, A_{3} \overline{A_{0}} & 5 & 2 & 0 \\
	    		\frac{1}{3} \, \ccancel{A_{0} \gamma_{1}} & 5 & 2 & 1 \\
	    		\frac{1}{4} \,\ccancel{ A_{0} C_{1}} \overline{\alpha_{0}} + \frac{1}{4} \,\ccancel{ A_{0} B_{1}} + 2 \, A_{2} \overline{A_{1}} & 4 & 3 & 0 \\
	    		2 \, A_{1} \overline{A_{0}} & 3 & 2 & 0 \\
	    		2 \, A_{1} \overline{A_{1}} & 3 & 3 & 0 \\
	    		\frac{1}{4} \, \alpha_{0} \ccancel{\overline{A_{0}} \overline{C_{1}}} + 2 \, A_{1} \overline{A_{2}} + \frac{1}{4} \, \ccancel{\overline{A_{0}} \overline{B_{1}}} & 3 & 4 & 0 \\
	    		\frac{1}{6} \, \ccancel{C_{1} \overline{A_{0}}} & 1 & 5 & 0 \\
	    		\frac{1}{24} \, {\left(\ccancel{3 \, \overline{A_{0}} \overline{\alpha_{0}}} + 4 \, \overline{A_{1}}\right)} C_{1} + \frac{1}{8} \, B_{1} \overline{A_{0}} & 1 & 6 & 0 \\
	    		\frac{1}{3} \, \ccancel{\gamma_{1} \overline{A_{0}}} & 2 & 5 & 1
	    		\end{dmatrix}\\
	    		&=
	    		\begin{dmatrix}
	    		2 \, A_{0} \overline{A_{0}} & 2 & 2 & 0 \\
	    		2 \, A_{0} \overline{A_{1}} & 2 & 3 & 0 \\
	    		2 \, A_{0} \overline{A_{2}} & 2 & 4 & 0 \\
	    		\frac{1}{3} \, C_{2} \overline{A_{0}} + 2 \, A_{0} \overline{A_{3}} + \frac{1}{12} \, {\left( 3 \, \overline{A_{1}}\right)} \overline{C_{1}} & 2 & 5 & 0 \\
	    		\frac{1}{8} \, A_{0} \overline{B_{1}} + \frac{1}{24} \, {\left(4 \, A_{1}\right)} \overline{C_{1}} & 6 & 1 & 0 \\2 \, A_{2} \overline{A_{0}} & 4 & 2 & 0 \\ 
	    		\frac{1}{12} \, {\left( 3 \, A_{1}\right)} C_{1} + \frac{1}{3} \, A_{0} C_{2} +
	    		2 \, A_{3} \overline{A_{0}} & 5 & 2 & 0 \\
	    		2 \, A_{2} \overline{A_{1}} & 4 & 3 & 0 \\
	    		2 \, A_{1} \overline{A_{0}} & 3 & 2 & 0 \\
	    		2 \, A_{1} \overline{A_{1}} & 3 & 3 & 0 \\
	    		2 \, A_{1} \overline{A_{2}} & 3 & 4 & 0 \\
	    		\frac{1}{24} \, {\left(4 \, \overline{A_{1}}\right)} C_{1} + \frac{1}{8} \, B_{1} \overline{A_{0}} & 1 & 6 & 0 \\
	    		\end{dmatrix}
	    		=	\begin{dmatrix}
	    		2 \, A_{0} \overline{A_{0}} & 2 & 2 & 0 \\
	    		2 \, A_{0} \overline{A_{1}} & 2 & 3 & 0 \\
	    		2 \, A_{0} \overline{A_{2}} & 2 & 4 & 0 \\
	    		\frac{1}{12}\bar{A_1 C_1} + 2 \, A_{0} \overline{A_{3}}  & 2 & 5 & 0 \\
	    		\frac{1}{24} \, \, A_{1} \overline{C_{1}} & 6 & 1 & 0 \\2 \, A_{2} \overline{A_{0}} & 4 & 2 & 0 \\ 
	    		\frac{1}{12} \,  \, A_{1} C_{1}+2 \, A_{3} \overline{A_{0}} & 5 & 2 & 0 \\
	    		2 \, A_{2} \overline{A_{1}} & 4 & 3 & 0 \\
	    		2 \, A_{1} \overline{A_{0}} & 3 & 2 & 0 \\
	    		2 \, A_{1} \overline{A_{1}} & 3 & 3 & 0 \\
	    		2 \, A_{1} \overline{A_{2}} & 3 & 4 & 0 \\
	    		\frac{1}{24} \, \overline{A_{1}} C_{1} & 1 & 6 & 0 \\
	    		\end{dmatrix}.
	    		\end{align*}
	    		as by  \eqref{3cancel1}
	    		\begin{align*}
	    		\left\{
	    		\begin{alignedat}{1}
	    		&\frac{1}{12} \, {\left( 3 \, A_{1}\right)} C_{1} + \frac{1}{3} \, A_{0} C_{2}\\ &\frac{1}{4}\s{\vec{A}_1}{\vec{C}_1}-\frac{1}{6}\s{\vec{A}_1}{\vec{C}_1}=\frac{1}{12}\s{\vec{A}_1}{\vec{C}_1}\\
	    		&\s{\bar{\vec{A}_0}}{\vec{B}_1}=\s{\bar{\vec{A}_0}}{-2\s{\bar{\vec{A}_1}}{\vec{C}_1}}=-\s{\bar{\vec{A}_1}}{\vec{C}_1}\\
	    		&\frac{1}{24} \, {\left(4 \, \overline{A_{1}}\right)} C_{1} + \frac{1}{8} \, B_{1} \overline{A_{0}}=\frac{1}{6}\s{\bar{\vec{A}_1}}{\vec{C}_1}-\frac{1}{8}\s{\bar{\vec{A}_1}}{\vec{C}_1}.
	    		\end{alignedat}
	    		\right.
	    		\end{align*}
	    		Therefore, if $\alpha_0,\alpha_1,\alpha_2,\alpha_3,\alpha_4\in \mathbb{C}$ are defined by
	    		\begin{align}\label{3deflambda1}
	    		\left\{
	    		\begin{alignedat}{1}
	    		\alpha_0&=2\s{\bar{\vec{A}_0}}{\vec{A}_1}\\
	    		\alpha_1&=2\s{\bar{\vec{A}_0}}{\vec{A}_2}\\
	    		\alpha_2&=\frac{1}{24}\s{\bar{\vec{A}_1}}{\vec{C}_1}\\
	    		\alpha_3&=\frac{1}{12}\s{
	    			\vec{A}_1}{\vec{C}_1}+2\s{\bar{\vec{A}_0}}{\vec{A}_3}\\
	    		\alpha_4&=2\s{\bar{\vec{A}_1}}{\vec{A}_2},
	    		\end{alignedat}\right.
	    		\end{align}
	    		we obtain
	    		\begin{align*}
	    		e^{2\lambda}&=|z|^4+2|\vec{A}_1|^2|z|^6+2\,\Re\bigg(\alpha_0z^3\z^2+\alpha_1z^4\z^2+\alpha_2z\z^6+\alpha_3z^5\z^2+\alpha_4z^4\z^3\bigg)+O(|z|^{8-\epsilon})\\
	    		&=|z|^4\left(1+2|\vec{A}_1|^2|z|^2+2\,\Re\bigg(\alpha_0z+\alpha_1z^2+\alpha_2z^{-1}\z^3+\alpha_3z^3+\alpha_4z^2\z\bigg)+O(|z|^{4-\epsilon})\right)\\
	    		&=|z|^{2\theta_0-2}\left(1+2|\vec{A}_1|^2|z|^2+2\,\Re\bigg(\alpha_0z+\alpha_1z^2+\alpha_2z^{2-\theta_0}\z^{\theta_0}+\alpha_3z^3+\alpha_4z^2\z\bigg)+O(|z|^{4-\epsilon})\right)\\
	    		\end{align*}
	    		which is translated as
	    		\begin{align*}
	    		e^{2\lambda}=\begin{dmatrix}
	    		1 & 2 & 2 & 0 \\
	    		2 \, {\left| A_{1} \right|}^{2} & 3 & 3 & 0 \\
	    		\alpha_{0} & 3 & 2 & 0 \\
	    		\alpha_{1} & 4 & 2 & 0 \\
	    		\alpha_{2} & 1 & 6 & 0 \\
	    		\alpha_{3} & 5 & 2 & 0 \\
	    		\alpha_{4} & 4 & 3 & 0 
	    		\end{dmatrix}
	    		\begin{dmatrix}
	    		\overline{\alpha_{0}} & 2 & 3 & 0 \\
	    		\overline{\alpha_{1}} & 2 & 4 & 0 \\
	    		\overline{\alpha_{2}} & 6 & 1 & 0 \\
	    		\overline{\alpha_{3}} & 2 & 5 & 0 \\
	    		\overline{\alpha_{4}} & 3 & 4 & 0.
	    		\end{dmatrix}
	    		\end{align*}
	    		and we obtain
	    		\begin{align*}
	    		\h_0=\begin{dmatrix}
	    		2 & A_{1} & 2 & 0 & 0 \\
	    		4 & A_{2} & 3 & 0 & 0 \\
	    		6 & A_{3} & 4 & 0 & 0 \\
	    		-\frac{1}{6} & C_{1} & 0 & 3 & 0 \\
	    		-\frac{1}{8} & B_{1} & 0 & 4 & 0 \\
	    		-\frac{1}{8} \, \overline{\alpha_{0}} & C_{1} & 0 & 4 & 0 \\
	    		\frac{1}{6} & \gamma_{1} & 1 & 3 & 0 \\
	    		-\frac{1}{6} \, \alpha_{0} & C_{1} & 1 & 3 & 0 \\
	    		\frac{1}{4} & \overline{B_{1}} & 2 & 2 & 0 \\
	    		-4 \, {\left| A_{1} \right|}^{2} + 2 \, \alpha_{0} \overline{\alpha_{0}} & A_{0} & 2 & 1 & 0 \\
	    		-4 \, {\left| A_{1} \right|}^{2} + 2 \, \alpha_{0} \overline{\alpha_{0}} & A_{1} & 3 & 1 & 0 \\
	    		-2 \, \alpha_{0} & A_{0} & 2 & 0 & 0 \\
	    		-2 \, \alpha_{0} & A_{1} & 3 & 0 & 0 \\
	    		-2 \, \alpha_{0} & A_{2} & 4 & 0 & 0 \\
	    		2 \, \alpha_{0}^{2} - 4 \, \alpha_{1} & A_{0} & 3 & 0 & 0 \\
	    		2 \, \alpha_{0}^{2} - 4 \, \alpha_{1} & A_{1} & 4 & 0 & 0 \\
	    		2 \, \alpha_{2} & A_{0} & 0 & 4 & 0 \\
	    		-10 \, \alpha_{0}^{3} + 6 \, \alpha_{0} \alpha_{1} - 6 \, \alpha_{3} & A_{0} & 4 & 0 & 0 \\
	    		8 \, \alpha_{0} {\left| A_{1} \right|}^{2} - 28 \, \alpha_{0}^{2} \overline{\alpha_{0}} + 4 \, \alpha_{1} \overline{\alpha_{0}} - 4 \, \alpha_{4} & A_{0} & 3 & 1 & 0 \\
	    		-8 \, \overline{\alpha_{2}} & A_{0} & 5 & -1 & 0 \\
	    		-8 \, \overline{\alpha_{0}}^{3} & A_{0} & 1 & 3 & 0 \\
	    		4 \, {\left| A_{1} \right|}^{2} \overline{\alpha_{0}} - 26 \, \alpha_{0} \overline{\alpha_{0}}^{2} + 2 \, \alpha_{0} \overline{\alpha_{1}} - 2 \, \overline{\alpha_{4}} & A_{0} & 2 & 2 & 0
	    		\end{dmatrix}
	    		\end{align*}
	    		Notice that there is no logarithm in this expansion, so we can replace in the case $\theta_0=3$
	    		\begin{align*}
	    		\mathscr{Q}_{\phi}=(\theta_0-1)(\theta_0-2)\s{\vec{A}_1}{\vec{C}_1}\frac{dz^4}{z}+O(|z|^{-\epsilon})=2\s{\vec{A}_1}{\vec{C}_1}\frac{dz^4}{z}+O(|z|^{-\epsilon})
	    		\end{align*}
	    		by
	    		\begin{align}
	    		\mathscr{Q}_{\phi}=2\s{\vec{A}_1}{\vec{C}_1}\frac{dz^4}{z}+O(1)
	    		\end{align}
	    		Now, the conservation law associated to the invariance by inversions yields if 
	    		\begin{align*}
	    		&\vec{\alpha}=\partial\H+|\H|^2\partial\phi+2\,g^{-1}\otimes\s{\H}{\h_0}\otimes\bar{\partial}\phi\\
	    		&\vec{\beta}=|\phi|^2\vec{\alpha}-2\s{\phi}{\vec{\alpha}}\phi-g^{-1}\otimes\left(\bar{\partial}|\phi|^2\otimes\h_0-2\s{\phi}{\h_0}\otimes\bar{\partial}\phi\right)
	    		\end{align*}
	    		the identity
	    		\begin{align*}
	    		d\,\Im(\vec{\beta})=d\,\Im\left(|\phi|^2\vec{\alpha}-2\s{\phi}{\vec{\alpha}}\phi-g^{-1}\otimes\left(\bar{\partial}|\phi|^2\otimes\h_0-2\s{\phi}{\h_0}\otimes\bar{\partial}\phi\right)\right)=0
	    		\end{align*}
	    		which is equivalent if $\vec{F}\in C^{\infty}(D^2\setminus\ens{0},\mathbb{C}^n)$ is such that
	    		\begin{align*}
	    		|\phi|^2\vec{\alpha}-2\s{\phi}{\vec{\alpha}}\phi-g^{-1}\otimes\left(\bar{\partial}|\phi|^2\otimes\h_0-2\s{\phi}{\h_0}\otimes\bar{\partial}\phi\right)=\vec{F}(z)dz
	    		\end{align*}
	    		to
	    		\begin{align*}
	    		d\,\Im(\vec{\beta})=\Im\left(\bar{\partial}\vec{\beta}\right)=\Im\left(\p{\z}\vec{F}(z)d\z\wedge dz\right)=-2\,\Re\left(\p{\z}\vec{F}(z)\right)dx_1\wedge dx_2.
	    		\end{align*}
	    		
	    		Then, we compute
	    		\small
	    		\begin{align*}
	    		&\Re(\p{\z}\vec{F}(z))=\\
	    		&\begin{dmatrix}
	    		\frac{1}{6} \, \overline{A_{0}}^{2} & C_{1} & -2 & 5 & 0 \\
	    		\frac{5}{72} \, \overline{A_{0}}^{2} & B_{1} & -2 & 6 & 0 \\
	    		-\frac{1}{24} \, \overline{A_{0}}^{2} \overline{\alpha_{0}} + \frac{1}{3} \, \overline{A_{0}} \overline{A_{1}} & C_{1} & -2 & 6 & 0 \\
	    		\frac{1}{9} \, C_{1} \overline{A_{0}}^{3} \overline{\alpha_{0}} - \frac{1}{9} \, C_{1} \overline{A_{0}}^{2} \overline{A_{1}} - 2 \, \alpha_{2} \overline{A_{0}}^{2} & A_{0} & -2 & 6 & 0 \\
	    		-\frac{1}{24} \, C_{1} \overline{A_{0}} & \overline{A_{1}} & -2 & 6 & 0 \\
	    		\frac{2}{9} \, A_{0} C_{1} \overline{A_{0}} \overline{A_{1}} + 2 \, A_{0} \alpha_{2} \overline{A_{0}} - \frac{1}{72} \, {\left(16 \, A_{0} \overline{A_{0}}^{2} - 3 \, \overline{A_{0}}\right)} C_{1} \overline{\alpha_{0}} + \frac{11}{72} \, B_{1} \overline{A_{0}} & \overline{A_{0}} & -2 & 6 & 0 \\
	    		-\frac{1}{6} \, \overline{A_{0}}^{2} & \gamma_{1} & -1 & 5 & 0 \\
	    		\frac{2}{3} \, A_{0} C_{1} \alpha_{0} \overline{A_{0}}^{2} - 16 \, A_{0} \overline{A_{0}} \overline{\alpha_{0}}^{3} - \frac{2}{3} \, A_{1} C_{1} \overline{A_{0}}^{2} & \overline{A_{0}} & -1 & 5 & 0 \\
	    		16 \, \overline{A_{0}}^{2} \overline{\alpha_{0}}^{3} & A_{0} & -1 & 5 & 0 \\
	    		-2 \, A_{0} \alpha_{0} \overline{A_{0}} + 2 \, A_{1} \overline{A_{0}} & \overline{A_{0}} & 0 & 2 & 0 \\
	    		-2 \, \overline{A_{0}}^{2} & A_{1} & 0 & 2 & 0 \\
	    		2 \, \alpha_{0} \overline{A_{0}}^{2} & A_{0} & 0 & 2 & 0 \\
	    		-2 \, A_{0} \alpha_{0} \overline{A_{0}} + 2 \, A_{1} \overline{A_{0}} & \overline{A_{1}} & 0 & 3 & 0 \\
	    		-4 \, A_{0} {\left| A_{1} \right|}^{2} \overline{A_{0}} - 2 \, A_{0} \alpha_{0} \overline{A_{1}} + 2 \, A_{1} \overline{A_{1}} + 2 \, {\left(2 \, A_{0} \alpha_{0} \overline{A_{0}} - A_{1} \overline{A_{0}}\right)} \overline{\alpha_{0}} & \overline{A_{0}} & 0 & 3 & 0 \\
	    		2 \, \overline{A_{0}}^{2} \overline{\alpha_{0}} - 4 \, \overline{A_{0}} \overline{A_{1}} & A_{1} & 0 & 3 & 0 \\
	    		4 \, {\left| A_{1} \right|}^{2} \overline{A_{0}}^{2} - 4 \, \alpha_{0} \overline{A_{0}}^{2} \overline{\alpha_{0}} + 4 \, \alpha_{0} \overline{A_{0}} \overline{A_{1}} & A_{0} & 0 & 3 & 0 \\
	    		-\frac{11}{36} \, \overline{A_{0}}^{2} & \overline{B_{1}} & 0 & 4 & 0 \\
	    		\lambda_1 & \overline{A_{0}} & 0 & 4 & 0 \\
	    		-2 \, A_{0} \alpha_{0} \overline{A_{0}} + 2 \, A_{1} \overline{A_{0}} & \overline{A_{2}} & 0 & 4 & 0 \\
	    		-4 \, A_{0} {\left| A_{1} \right|}^{2} \overline{A_{0}} - 2 \, A_{0} \alpha_{0} \overline{A_{1}} + 2 \, A_{1} \overline{A_{1}} + 2 \, {\left(2 \, A_{0} \alpha_{0} \overline{A_{0}} - A_{1} \overline{A_{0}}\right)} \overline{\alpha_{0}} & \overline{A_{1}} & 0 & 4 & 0 \\
	    		\lambda_2 & A_{0} & 0 & 4 & 0 \\
	    		-2 \, \overline{A_{0}}^{2} \overline{\alpha_{0}}^{2} + 4 \, \overline{A_{0}} \overline{A_{1}} \overline{\alpha_{0}} + 2 \, \overline{A_{0}}^{2} \overline{\alpha_{1}} - 2 \, \overline{A_{1}}^{2} - 4 \, \overline{A_{0}} \overline{A_{2}} & A_{1} & 0 & 4 & 0 \\
	    		-4 \, \alpha_{0}^{2} \overline{A_{0}}^{2} + 4 \, \alpha_{1} \overline{A_{0}}^{2} + \frac{1}{6} \, C_{1} \overline{A_{0}} & A_{0} & 1 & 2 & 0 \\
	    		4 \, A_{0} \alpha_{0}^{2} \overline{A_{0}} - 4 \, A_{1} \alpha_{0} \overline{A_{0}} - 4 \, A_{0} \alpha_{1} \overline{A_{0}} + 4 \, A_{2} \overline{A_{0}} & \overline{A_{0}} & 1 & 2 & 0 \\
	    		-4 \, \overline{A_{0}}^{2} & A_{2} & 1 & 2 & 0 \\
	    		4 \, \alpha_{0} \overline{A_{0}}^{2} & A_{1} & 1 & 2 & 0 
	    		\end{dmatrix}
	    		\end{align*}
	    		\begin{align*}
	    		\begin{dmatrix}
	    		-\frac{1}{9} \, A_{0} \overline{A_{0}} & B_{1} & 1 & 3 & 0 \\
	    		\lambda_3 & A_{0} & 1 & 3 & 0 \\
	    		-\frac{1}{6} \, A_{0} \overline{A_{0}} \overline{\alpha_{0}} + \frac{1}{6} \, A_{0} \overline{A_{1}} & C_{1} & 1 & 3 & 0 \\
	    		\lambda_4 & \overline{A_{0}} & 1 & 3 & 0 \\
	    		4 \, A_{0} \alpha_{0}^{2} \overline{A_{0}} - 4 \, A_{1} \alpha_{0} \overline{A_{0}} - 4 \, A_{0} \alpha_{1} \overline{A_{0}} - \frac{1}{3} \, A_{0} C_{1} + 4 \, A_{2} \overline{A_{0}} & \overline{A_{1}} & 1 & 3 & 0 \\
	    		4 \, \overline{A_{0}}^{2} \overline{\alpha_{0}} - 8 \, \overline{A_{0}} \overline{A_{1}} & A_{2} & 1 & 3 & 0 \\
	    		8 \, {\left| A_{1} \right|}^{2} \overline{A_{0}}^{2} - 8 \, \alpha_{0} \overline{A_{0}}^{2} \overline{\alpha_{0}} + 8 \, \alpha_{0} \overline{A_{0}} \overline{A_{1}} & A_{1} & 1 & 3 & 0 \\
	    		-2 \, A_{0} \overline{A_{0}} \overline{\alpha_{0}} + 2 \, A_{0} \overline{A_{1}} & A_{0} & 2 & 0 & 0 \\
	    		-2 \, A_{0}^{2} & \overline{A_{1}} & 2 & 0 & 0 \\
	    		2 \, A_{0}^{2} \overline{\alpha_{0}} & \overline{A_{0}} & 2 & 0 & 0 \\
	    		-4 \, A_{0}^{2} \overline{\alpha_{0}}^{2} + 4 \, A_{0}^{2} \overline{\alpha_{1}} + \frac{1}{6} \, A_{0} \overline{C_{1}} & \overline{A_{0}} & 2 & 1 & 0 \\
	    		4 \, A_{0} \overline{A_{0}} \overline{\alpha_{0}}^{2} - 4 \, A_{0} \overline{A_{1}} \overline{\alpha_{0}} - 4 \, A_{0} \overline{A_{0}} \overline{\alpha_{1}} + 4 \, A_{0} \overline{A_{2}} & A_{0} & 2 & 1 & 0 \\
	    		-4 \, A_{0}^{2} & \overline{A_{2}} & 2 & 1 & 0 \\
	    		4 \, A_{0}^{2} \overline{\alpha_{0}} & \overline{A_{1}} & 2 & 1 & 0 \\
	    		\lambda_5 & \overline{A_{0}} & 2 & 2 & 0 \\
	    		\lambda_6 & A_{0} & 2 & 2 & 0 \\
	    		-\frac{1}{4} \, A_{0} \alpha_{0} \overline{A_{0}} + \frac{1}{4} \, A_{1} \overline{A_{0}} & C_{1} & 2 & 2 & 0 \\
	    		-\frac{1}{4} \, A_{0} \overline{A_{0}} \overline{\alpha_{0}} + \frac{1}{4} \, A_{0} \overline{A_{1}} & \overline{C_{1}} & 2 & 2 & 0 \\
	    		-6 \, \alpha_{0}^{2} \overline{A_{0}}^{2} + 6 \, \alpha_{1} \overline{A_{0}}^{2} - \frac{3}{8} \, C_{1} \overline{A_{0}} & A_{1} & 2 & 2 & 0 \\
	    		-6 \, A_{0}^{2} \overline{\alpha_{0}}^{2} + 6 \, A_{0}^{2} \overline{\alpha_{1}} - \frac{3}{8} \, A_{0} \overline{C_{1}} & \overline{A_{1}} & 2 & 2 & 0 \\
	    		-6 \, \overline{A_{0}}^{2} & A_{3} & 2 & 2 & 0 \\
	    		-6 \, A_{0}^{2} & \overline{A_{3}} & 2 & 2 & 0 \\
	    		6 \, \alpha_{0} \overline{A_{0}}^{2} & A_{2} & 2 & 2 & 0 \\
	    		6 \, A_{0}^{2} \overline{\alpha_{0}} & \overline{A_{2}} & 2 & 2 & 0 \\
	    		-2 \, A_{0} \overline{A_{0}} \overline{\alpha_{0}} + 2 \, A_{0} \overline{A_{1}} & A_{1} & 3 & 0 & 0 \\
	    		-4 \, A_{0} {\left| A_{1} \right|}^{2} \overline{A_{0}} - 2 \, A_{0} \alpha_{0} \overline{A_{1}} + 2 \, A_{1} \overline{A_{1}} + 2 \, {\left(2 \, A_{0} \alpha_{0} \overline{A_{0}} - A_{1} \overline{A_{0}}\right)} \overline{\alpha_{0}} & A_{0} & 3 & 0 & 0 \\
	    		2 \, A_{0}^{2} \alpha_{0} - 4 \, A_{0} A_{1} & \overline{A_{1}} & 3 & 0 & 0 \\
	    		4 \, A_{0}^{2} {\left| A_{1} \right|}^{2} - 4 \, {\left(A_{0}^{2} \alpha_{0} - A_{0} A_{1}\right)} \overline{\alpha_{0}} & \overline{A_{0}} & 3 & 0 & 0 
	    		\end{dmatrix}
	    		\end{align*}
	    		\begin{align*}
	    		\begin{dmatrix}
	    		-\frac{1}{9} \, A_{0} \overline{A_{0}} & \overline{B_{1}} & 3 & 1 & 0 \\
	    		\lambda_7 & \overline{A_{0}} & 3 & 1 & 0 \\
	    		-\frac{1}{6} \, A_{0} \alpha_{0} \overline{A_{0}} + \frac{1}{6} \, A_{1} \overline{A_{0}} & \overline{C_{1}} & 3 & 1 & 0 \\
	    		\lambda_8 & A_{0} & 3 & 1 & 0 \\
	    		4 \, A_{0} \overline{A_{0}} \overline{\alpha_{0}}^{2} - 4 \, A_{0} \overline{A_{1}} \overline{\alpha_{0}} - 4 \, A_{0} \overline{A_{0}} \overline{\alpha_{1}} + 4 \, A_{0} \overline{A_{2}} - \frac{1}{3} \, \overline{A_{0}} \overline{C_{1}} & A_{1} & 3 & 1 & 0 \\
	    		4 \, A_{0}^{2} \alpha_{0} - 8 \, A_{0} A_{1} & \overline{A_{2}} & 3 & 1 & 0 \\
	    		8 \, A_{0}^{2} {\left| A_{1} \right|}^{2} - 8 \, {\left(A_{0}^{2} \alpha_{0} - A_{0} A_{1}\right)} \overline{\alpha_{0}} & \overline{A_{1}} & 3 & 1 & 0 \\
	    		-\frac{11}{36} \, A_{0}^{2} & B_{1} & 4 & 0 & 0 \\
	    		\lambda_9 & A_{0} & 4 & 0 & 0 \\
	    		-2 \, A_{0} \overline{A_{0}} \overline{\alpha_{0}} + 2 \, A_{0} \overline{A_{1}} & A_{2} & 4 & 0 & 0 \\
	    		-4 \, A_{0} {\left| A_{1} \right|}^{2} \overline{A_{0}} - 2 \, A_{0} \alpha_{0} \overline{A_{1}} + 2 \, A_{1} \overline{A_{1}} + 2 \, {\left(2 \, A_{0} \alpha_{0} \overline{A_{0}} - A_{1} \overline{A_{0}}\right)} \overline{\alpha_{0}} & A_{1} & 4 & 0 & 0 \\
	    		-8 \, A_{0}^{2} \alpha_{0} {\left| A_{1} \right|}^{2} + 8 \, A_{0} A_{1} {\left| A_{1} \right|}^{2} - 4 \, A_{0}^{2} \alpha_{1} \overline{\alpha_{0}} + 2 \, A_{0}^{2} \alpha_{4} + 2 \, {\left(19 \, A_{0}^{2} \alpha_{0}^{2} - 4 \, A_{0} A_{1} \alpha_{0} + A_{1}^{2} + 2 \, A_{0} A_{2}\right)} \overline{\alpha_{0}} & \overline{A_{0}} & 4 & 0 & 0 \\
	    		-2 \, A_{0}^{2} \alpha_{0}^{2} + 4 \, A_{0} A_{1} \alpha_{0} + 2 \, A_{0}^{2} \alpha_{1} - 2 \, A_{1}^{2} - 4 \, A_{0} A_{2} & \overline{A_{1}} & 4 & 0 & 0 \\
	    		\frac{1}{6} \, A_{0}^{2} & \overline{C_{1}} & 5 & -2 & 0 \\
	    		-\frac{1}{6} \, A_{0}^{2} & \gamma_{1} & 5 & -1 & 0 \\
	    		-16 \, A_{0} \alpha_{0}^{3} \overline{A_{0}} + \frac{2}{3} \, A_{0}^{2} \overline{A_{0}} \overline{C_{1}} \overline{\alpha_{0}} - \frac{2}{3} \, A_{0}^{2} \overline{A_{1}} \overline{C_{1}} & A_{0} & 5 & -1 & 0 \\
	    		16 \, A_{0}^{2} \alpha_{0}^{3} & \overline{A_{0}} & 5 & -1 & 0 \\
	    		\frac{5}{72} \, A_{0}^{2} & \overline{B_{1}} & 6 & -2 & 0 \\
	    		-\frac{1}{24} \, A_{0}^{2} \alpha_{0} + \frac{1}{3} \, A_{0} A_{1} & \overline{C_{1}} & 6 & -2 & 0 \\
	    		\frac{1}{9} \, A_{0}^{3} \alpha_{0} \overline{C_{1}} - \frac{1}{9} \, A_{0}^{2} A_{1} \overline{C_{1}} - 2 \, A_{0}^{2} \overline{\alpha_{2}} & \overline{A_{0}} & 6 & -2 & 0 \\
	    		-\frac{1}{24} \, A_{0} \overline{C_{1}} & A_{1} & 6 & -2 & 0 \\
	    		\frac{2}{9} \, A_{0} A_{1} \overline{A_{0}} \overline{C_{1}} - \frac{1}{72} \, {\left(16 \, A_{0}^{2} \overline{A_{0}} - 3 \, A_{0}\right)} \alpha_{0} \overline{C_{1}} + 2 \, A_{0} \overline{A_{0}} \overline{\alpha_{2}} + \frac{11}{72} \, A_{0} \overline{B_{1}} & A_{0} & 6 & -2 & 0
	    		\end{dmatrix}
	    		\end{align*}
	    		\normalsize
	    		As $\s{\vec{A}_0}{\vec{A}_0}=\s{\vec{A}_0}{\vec{C}_1}=\s{\bar{\vec{A}_0}}{\vec{C}_1}=0$, we have
	    		\begin{align*}
	    		&A_0\bar{A_0}^2C_1\in\ens{\s{\vec{A}_0}{\bar{\vec{A}_0}}\s{\bar{\vec{A}_0}}{\vec{C}_1},\s{\bar{\vec{A}_0}}{\bar{\vec{A}_0}}\s{\vec{A}_0}{\vec{C}_1}}=\ens{0}\\
	    		&A_1\bar{A_0}^2C_1\in \ens{\s{\vec{A}_1}{\bar{\vec{A}_0}}\s{\bar{\vec{A}_0}}{\vec{C}_1},\s{\vec{A}_1}{\vec{C}_1}\bar{\s{\vec{A}_0}{\vec{A}_0}}}=\ens{0},
	    		\end{align*}
	    		and we obtain
	    		\begin{align*}
	    		{\left(A_{0} \alpha_{0} \overline{A_{0}}^{2} - A_{1} \overline{A_{0}}^{2}\right)} C_{1} =0,
	    		\end{align*}
	    		so the coefficient in $\dfrac{\z^{\theta_0+2}}{z}=\dfrac{\z^5}{z}$ in the Taylor expansion is
	    		\begin{align}\label{motif}
	    		\begin{dmatrix}
	    		-\frac{1}{6} \, \ccancel{\overline{A_{0}}^{2}} & \gamma_{1} & -1 & 5 & 0 \\
	    		\frac{2}{3} \ccancel{\, A_{0} C_{1} \alpha_{0} \overline{A_{0}}^{2}} - 16 \, A_{0} \overline{A_{0}} \overline{\alpha_{0}}^{3} - \frac{2}{3} \ccancel{\, A_{1} C_{1} \overline{A_{0}}^{2}} & \overline{A_{0}} & -1 & 5 & 0 \\
	    		\ccancel{16 \, \overline{A_{0}}^{2} \overline{\alpha_{0}}^{3}} & A_{0} & -1 & 5 & 0 
	    		\end{dmatrix}
	    		=-16|\vec{A}_0|^2\bar{\alpha_0}^3\bar{\vec{A}_0}=-8\bar{\alpha_0}^3\bar{\vec{A}_0}=0
	    		\end{align}
	    		as $|\vec{A}_0|^2=\dfrac{1}{2}$. Therefore, we obtain by \eqref{motif}
	    		\begin{align}\label{3yatta1}
	    		\alpha_0=2\s{\bar{\vec{A}_0}}{\vec{A}_1}=0\qquad 
	    		\end{align}

	    		At this point, we will rerun our previous computations by plugging $\alpha_0=0$ from the stage when we integrate 
	    		\begin{align*}
	    		\p{\z}\left(\p{z}\phi\right)=\frac{e^{2\lambda}}{2}\H.
	    		\end{align*}
	    		We obtain
	    		\begin{align}\label{3dev2}
	    		\p{z}\phi&=\begin{dmatrix}
	    		1 & A_{0} & 2 & 0 & 0 \\
	    		1 & A_{1} & 3 & 0 & 0 \\
	    		1 & A_{2} & 4 & 0 & 0 \\
	    		1 & A_{3} & 5 & 0 & 0 \\
	    		\frac{1}{12} & C_{1} & 1 & 3 & 0 \\
	    		\frac{1}{16} & B_{1} & 1 & 4 & 0 \\
	    		\frac{1}{8} & \overline{C_{1}} & 2 & 2 & 0 \\
	    		-\frac{1}{36} & \gamma_{1} & 2 & 3 & 0 \\
	    		\frac{1}{6} & C_{2} & 2 & 3 & 0 \\
	    		\frac{1}{6} & \gamma_{1} & 2 & 3 & 1 \\
	    		\frac{1}{8} & \overline{B_{1}} & 3 & 2 & 0
	    		\end{dmatrix}\;\,
	    		\phi(z)=\begin{dmatrix}
	    		\frac{1}{3} & \overline{A_{0}} & 0 & 3 & 0 \\
	    		\frac{1}{4} & \overline{A_{1}} & 0 & 4 & 0 \\
	    		\frac{1}{5} & \overline{A_{2}} & 0 & 5 & 0 \\
	    		\frac{1}{6} & \overline{A_{3}} & 0 & 6 & 0 \\
	    		\frac{1}{3} & A_{0} & 3 & 0 & 0 \\
	    		\frac{1}{4} & A_{1} & 4 & 0 & 0 \\
	    		\frac{1}{5} & A_{2} & 5 & 0 & 0 \\
	    		\frac{1}{6} & A_{3} & 6 & 0 & 0 \\
	    		\frac{1}{24} & C_{1} & 2 & 3 & 0 \\
	    		\frac{1}{32} & B_{1} & 2 & 4 & 0 \\
	    		\frac{1}{24} & \overline{C_{1}} & 3 & 2 & 0 \\
	    		-\frac{1}{54} & \gamma_{1} & 3 & 3 & 0 \\
	    		\frac{1}{18} & C_{2} & 3 & 3 & 0 \\
	    		\frac{1}{18} & \gamma_{1} & 3 & 3 & 1 \\
	    		\frac{1}{32} & \overline{B_{1}} & 4 & 2 & 0
	    		\end{dmatrix},\;\,
	    		e^{2\lambda}=\begin{dmatrix}
	    		1 & 2 & 2 & 0 \\
	    		2 \, {\left| A_{1} \right|}^{2} & 3 & 3 & 0 \\
	    		\alpha_{1} & 4 & 2 & 0 \\
	    		\alpha_{2} & 1 & 6 & 0 \\
	    		\alpha_{3} & 5 & 2 & 0 \\
	    		\alpha_{4} & 4 & 3 & 0 \\
	    		\overline{\alpha_{1}} & 2 & 4 & 0 \\
	    		\overline{\alpha_{2}} & 6 & 1 & 0 \\
	    		\overline{\alpha_{3}} & 2 & 5 & 0 \\
	    		\overline{\alpha_{4}} & 3 & 4 & 0
	    		\end{dmatrix}\\
	    		\h_0&=\begin{dmatrix}
	    		2 & A_{1} & 2 & 0 & 0 \\
	    		4 & A_{2} & 3 & 0 & 0 \\
	    		6 & A_{3} & 4 & 0 & 0 \\
	    		-\frac{1}{6} & C_{1} & 0 & 3 & 0 \\
	    		-\frac{1}{8} & B_{1} & 0 & 4 & 0 \\
	    		\frac{1}{6} & \gamma_{1} & 1 & 3 & 0 \\
	    		\frac{1}{4} & \overline{B_{1}} & 2 & 2 & 0 \\
	    		-4 \, {\left| A_{1} \right|}^{2} & A_{0} & 2 & 1 & 0 
	    		\end{dmatrix}
	    		\begin{dmatrix}
	    		-4 \, {\left| A_{1} \right|}^{2} & A_{1} & 3 & 1 & 0 \\
	    		-4 \, \alpha_{1} & A_{0} & 3 & 0 & 0 \\
	    		-4 \, \alpha_{1} & A_{1} & 4 & 0 & 0 \\
	    		2 \, \alpha_{2} & A_{0} & 0 & 4 & 0 \\
	    		-6 \, \alpha_{3} & A_{0} & 4 & 0 & 0 \\
	    		-4 \, \alpha_{4} & A_{0} & 3 & 1 & 0 \\
	    		-8 \, \overline{\alpha_{2}} & A_{0} & 5 & -1 & 0 \\
	    		-2 \, \overline{\alpha_{4}} & A_{0} & 2 & 2 & 0
	    		\end{dmatrix}\\
	    		\vec{\alpha}&=\begin{dmatrix}
	    		-\frac{1}{6} \, C_{1}^{2} & \overline{A_{0}} & -3 & 3 & 0 \\
	    		-\frac{1}{2} & C_{1} & -2 & 0 & 0 \\
	    		-\frac{1}{2} & B_{1} & -2 & 1 & 0 \\
	    		-\frac{1}{6} \, C_{1} \overline{C_{1}} & \overline{A_{0}} & -2 & 2 & 0 \\
	    		\frac{1}{2} & \gamma_{1} & -1 & 0 & 0 \\
	    		2 \, A_{1} C_{1} & \overline{A_{0}} & -1 & 0 & 0 \\
	    		2 \, A_{1} C_{1} & \overline{A_{1}} & -1 & 1 & 0 \\
	    		-4 \, A_{0} C_{1} {\left| A_{1} \right|}^{2} + 2 \, A_{1} B_{1} & \overline{A_{0}} & -1 & 1 & 0 \\
	    		\frac{1}{2} & \overline{B_{1}} & 0 & -1 & 0 \\
	    		2 \, A_{1} \overline{C_{1}} & \overline{A_{0}} & 0 & -1 & 0 \\
	    		2 \, A_{1} \overline{C_{1}} & \overline{A_{1}} & 0 & 0 & 0 \\
	    		\frac{1}{4} \, C_{1}^{2} & A_{0} & 0 & 0 & 0 \\
	    		-4 \, A_{0} {\left| A_{1} \right|}^{2} \overline{C_{1}} - 4 \, A_{0} C_{1} \alpha_{1} + 4 \, A_{2} C_{1} + 4 \, A_{1} C_{2} & \overline{A_{0}} & 0 & 0 & 0 \\
	    		4 \, A_{1} \gamma_{1} & \overline{A_{0}} & 0 & 0 & 1 \\
	    		\frac{1}{2} \, C_{1} \overline{C_{1}} & A_{0} & 1 & -1 & 0 \\
	    		-4 \, A_{0} \alpha_{1} \overline{C_{1}} + 2 \, A_{1} \overline{B_{1}} + 4 \, A_{2} \overline{C_{1}} & \overline{A_{0}} & 1 & -1 & 0 \\
	    		\frac{1}{4} \, \overline{C_{1}}^{2} & A_{0} & 2 & -2 & 0
	    		\end{dmatrix}
	    		\end{align}
	    		Then, we have for some $\vec{D}_3\in \mathbb{C}^n$
	    		\begin{align*}
	    			&\H+\vec{\gamma}_0\log|z|=\vec{C}_2+\Re\left(\vec{D}_3z\right)
	    			\\
	    			&+\begin{dmatrix}
	    			-\frac{1}{12} \, C_{1}^{2} & \overline{A_{0}} & -2 & 3 & 0 &\textbf{(1)}\\
	    			-C_{1} \overline{A_{1}} & A_{0} & -1 & 1 & 0 &\textbf{(2)}\\
	    			-\frac{5}{24} \, C_{1} \overline{C_{1}} & \overline{A_{0}} & -1 & 2 & 0 &\textbf{(3)}\\
	    			\ccancel{C_{1} \overline{A_{0}} \overline{\alpha_{1}}} - \frac{1}{2} \ccancel{\, B_{1} \overline{A_{1}}} - C_{1} \overline{A_{2}} & A_{0} & -1 & 2 & 0 &\textbf{(4)}\\
	    			-2 \, A_{1} C_{1} & \overline{A_{0}} & 0 & 0 & 1 &\textbf{(5)}\\
	    			-2 \, \overline{A_{1}} \overline{C_{1}} & A_{0} & 0 & 0 & 1 &\textbf{(6)}\\
	    			-\frac{1}{8} \, \overline{C_{1}}^{2} & \overline{A_{0}} & 0 & 1 & 0 &\textbf{(7)}\\
	    			\ccancel{2 \, C_{1} {\left| A_{1} \right|}^{2} \overline{A_{0}}} + \ccancel{2 \, \overline{A_{0}} \overline{C_{1}} \overline{\alpha_{1}}} - 2 \, C_{2} \overline{A_{1}} + \gamma_{1} \overline{A_{1}} - 2 \, \overline{A_{2}} \overline{C_{1}} & A_{0} & 0 & 1 & 0 &\textbf{(8)}\\
	    			-C_{1} \overline{A_{1}} & A_{1} & 0 & 1 & 0 &\textbf{(9)}\\
	    			-2 \, A_{1} C_{1} & \overline{A_{1}} & 0 & 1 & 1 &\textbf{(10)}\\
	    			\ccancel{4 \, A_{0} C_{1} {\left| A_{1} \right|}^{2}} - \ccancel{2 \, A_{1} B_{1}} & \overline{A_{0}} & 0 & 1 & 1 &\textbf{(11)}\\
	    			-2 \, \gamma_{1} \overline{A_{1}} & A_{0} & 0 & 1 & 1 &\textbf{(12)}\\
	    			-A_{1} \overline{C_{1}} & \overline{A_{0}} & 1 & -1 & 0 &\textbf{(13)}\\
	    			-\frac{1}{8} \, C_{1}^{2} & A_{0} & 1 & 0 & 0 &\textbf{(14)}\\
	    			\ccancel{2 \, A_{0} {\left| A_{1} \right|}^{2} \overline{C_{1}}} + \ccancel{2 \, A_{0} C_{1} \alpha_{1}} - 2 \, A_{2} C_{1} - 2 \, A_{1} C_{2} + A_{1} \gamma_{1} & \overline{A_{0}} & 1 & 0 & 0 &\textbf{(15)}\\
	    			-A_{1} \overline{C_{1}} & \overline{A_{1}} & 1 & 0 & 0 &\textbf{(16)}\\
	    			-2 \, A_{1} \gamma_{1} & \overline{A_{0}} & 1 & 0 & 1 &\textbf{(17)}\\
	    			-2 \, \overline{A_{1}} \overline{C_{1}} & A_{1} & 1 & 0 & 1 &\textbf{(18)}\\
	    			\ccancel{4 \, {\left| A_{1} \right|}^{2} \overline{A_{0}} \overline{C_{1}}} - \ccancel{2 \, \overline{A_{1}} \overline{B_{1}}} & A_{0} & 1 & 0 & 1 &\textbf{(19)}\\
	    			-\frac{5}{24} \, C_{1} \overline{C_{1}} & A_{0} & 2 & -1 & 0 &\textbf{(20)}\\
	    			\ccancel{A_{0} \alpha_{1} \overline{C_{1}}} - \frac{1}{2} \ccancel{\, A_{1} \overline{B_{1}}} - A_{2} \overline{C_{1}} & \overline{A_{0}} & 2 & -1 & 0 &\textbf{(21)}\\
	    			-\frac{1}{12} \, \overline{C_{1}}^{2} & A_{0} & 3 & -2 & 0&\textbf{(22)}
	    			\end{dmatrix}
	    		\end{align*}
	    		Now, recall that
	    		\begin{align*}
	    			\H=\begin{dmatrix}
	    			\frac{1}{2} & C_{1} & -1 & 0 & 0 \\
	    			\frac{1}{2} & B_{1} & -1 & 1 & 0 \\
	    			1 & C_{2} & 0 & 0 & 0 \\
	    			\end{dmatrix}
	    			\begin{dmatrix}
	    			1 & \gamma_{1} & 0 & 0 & 1 \\
	    			\frac{1}{2} & \overline{C_{1}} & 0 & -1 & 0 \\
	    			\frac{1}{2} & \overline{B_{1}} & 1 & -1 & 0
	    			\end{dmatrix}
	    		\end{align*}
	    		and we see that there exists $\vec{C}_3,\vec{B}_2,\vec{E}_1,\vec{\gamma}_2\in \mathbb{C}^n$ such that
	    		\begin{align*}
	    			\H=\Re\left(\frac{\vec{C}_1}{z}+\vec{C}_3z+\vec{B}_1\frac{\z}{z}+\vec{B}_2\frac{\z^2}{z}+\vec{E}_1\frac{\z^3}{z^2}\right)+\vec{C}_2+\vec{\gamma}_1\log|z|+\Re(\vec{\gamma}_2z)\log|z|+O(|z|^{2-\epsilon}).
	    		\end{align*}
	    		Then, we have
	    		\begin{align*}
	    		    \H=\begin{dmatrix}
	    			\frac{1}{2} & C_{1} & -1 & 0 & 0 \\
	    			1 & C_{2} & 0 & 0 & 0 \\
	    			\frac{1}{2} & C_{3} & 1 & 0 & 0 \\
	    			\frac{1}{2} & B_{1} & -1 & 1 & 0 \\
	    			\frac{1}{2} & B_{2} & -1 & 2 & 0 \\
	    			\frac{1}{2} & E_{1} & -2 & 3 & 0 \\
	    			1 & \gamma_{1} & 0 & 0 & 1 
	    			\end{dmatrix}
	    			\begin{dmatrix}
	    			\frac{1}{2} & \gamma_{2} & 1 & 0 & 1 \\
	    			\frac{1}{2} & \overline{C_{1}} & 0 & -1 & 0 \\
	    			\frac{1}{2} & \overline{C_{3}} & 0 & 1 & 0 \\
	    			\frac{1}{2} & \overline{B_{1}} & 1 & -1 & 0 \\
	    			\frac{1}{2} & \overline{B_{2}} & 2 & -1 & 0 \\
	    			\frac{1}{2} & \overline{E_{1}} & 3 & -2 & 0 \\
	    			\frac{1}{2} & \overline{\gamma_{2}} & 0 & 1 & 1
	    			\end{dmatrix}
	    		\end{align*}
	    		and we obtain immediately
	    		\begin{align*}
	    		\frac{1}{2}\vec{B}_2=\begin{dmatrix}
	    		-\frac{5}{24} \, C_{1} \overline{C_{1}} & \overline{A_{0}} & -1 & 2 & 0 &\textbf{(3)}\\
	    		\ccancel{C_{1} \overline{A_{0}} \overline{\alpha_{1}}} - \frac{1}{2} \ccancel{\, B_{1} \overline{A_{1}}} - C_{1} \overline{A_{2}} & A_{0} & -1 & 2 & 0 &\textbf{(4)}
	    		\end{dmatrix}
	    		=-\frac{5}{24}|\vec{C}_1|^2\bar{\vec{A}_0}-\s{\bar{\vec{A}_2}}{\vec{C}_1}\vec{A}_0
	    		\end{align*}
	    		so
	    		\begin{align}\label{3b2}
	    			\vec{B}_2=-\frac{5}{12}|\vec{C}_1|^2\bar{\vec{A}_0}-2\s{\bar{\vec{A}_2}}{\vec{C}_1}\vec{A}_0.
	    		\end{align}
	    		Then, we have
	    		\begin{align}
	    			\frac{1}{2}\vec{E}_1=\begin{dmatrix}
	    			-\frac{1}{12} \, C_{1}^{2} & \overline{A_{0}} & -2 & 3 & 0 &\textbf{(1)}
	    			\end{dmatrix}
	    		\end{align}
	    		so
	    		\begin{align}\label{3e1}
	    			\vec{E}_1=-\frac{1}{6}\s{\vec{C}_1}{\vec{C}_1}\bar{\vec{A}_0}.
	    		\end{align}
	    		Finally, we have
	    		\begin{align*}
	    			\frac{1}{2}\vec{\gamma}_2&=\begin{dmatrix}
	    			-2 \, A_{1} \gamma_{1} & \overline{A_{0}} & 1 & 0 & 1 &\textbf{(17)}\\
	    			-2 \, \overline{A_{1}} \overline{C_{1}} & A_{1} & 1 & 0 & 1 &\textbf{(18)}\\
	    			\ccancel{4 \, {\left| A_{1} \right|}^{2} \overline{A_{0}} \overline{C_{1}}} - \ccancel{2 \, \overline{A_{1}} \overline{B_{1}}} & A_{0} & 1 & 0 & 1 &\textbf{(19)}
	    			\end{dmatrix}\\
	    			&=-2\s{\vec{A}_1}{\vec{\gamma}_1}\bar{\vec{A}_0}-2\bar{\s{\vec{A}_1}{\vec{C}_1}}\vec{A}_1
	    		\end{align*}
	    		so
	    		\begin{align}\label{3gamma2}
	    			\vec{\gamma}_2=-4\s{\vec{A}_1}{\vec{\gamma}_1}\bar{\vec{A}_0}-4\bar{\s{\vec{A}_1}{\vec{C}_1}}\vec{A}_1.
	    		\end{align}
	    		Finally, we obtain by \eqref{3b2}, \eqref{3e1} and \eqref{3gamma2}
	    		\begin{align}\label{3consts2}
	    			\left\{
	    			\begin{alignedat}{1}
	    			\vec{C}_2&=\Re\left(\vec{D}_2\right)\in\R^n\\
	    			\vec{C}_3&=\vec{D}_3-\frac{1}{4}\s{\vec{C}_1}{\vec{C}_1}\vec{A}_0+\left(-4(\s{\vec{A}_2}{\vec{C}_1}+\s{\vec{A}_1}{\vec{C}_2})+2\s{\vec{A}_1}{\vec{\gamma}_1}\right)\bar{\vec{A}_0}-2\s{\vec{A}_1}{\bar{\vec{C}_1}}\bar{\vec{A}_1}\\
	    			\vec{B}_1&=-2\s{\bar{\vec{A}_1}}{\vec{C}_1}\vec{A}_0\\
	    			\vec{B}_2&=-\frac{5}{12}|\vec{C}_1|^2\bar{\vec{A}_0}-2\s{\bar{\vec{A}_2}}{\vec{C}_1}\vec{A}_0.\\
	    			\vec{E}_1&=-\frac{1}{6}\s{\vec{C}_1}{\vec{C}_1}\bar{\vec{A}_0}\\
	    			\vec{\gamma}_1&=-\vec{\gamma}_0-4\,\Re\left(\s{\vec{A}_1}{\vec{C}_1}\bar{\vec{A}_0}\right)\in\mathbb{R}^n\\
	    			\vec{\gamma}_2&=-4\s{\vec{A}_1}{\vec{\gamma}_1}\bar{\vec{A}_0}-4\bar{\s{\vec{A}_1}{\vec{C}_1}}\vec{A}_1.
	    			\end{alignedat}\right.
	    		\end{align}
	    		and
	    		\begin{align*}
	    			\p{z}\phi=\begin{dmatrix}
	    			1 & A_{0} & 2 & 0 & 0 \\
	    			1 & A_{1} & 3 & 0 & 0 \\
	    			1 & A_{2} & 4 & 0 & 0 \\
	    			1 & A_{3} & 5 & 0 & 0 \\
	    			1 & A_{4} & 6 & 0 & 0 \\
	    			\frac{1}{24} & E_{1} & 0 & 6 & 0 \\
	    			\frac{1}{12} & C_{1} & 1 & 3 & 0 \\
	    			\frac{1}{16} & B_{1} & 1 & 4 & 0 \\
	    			\frac{1}{20} & B_{2} & 1 & 5 & 0 \\
	    			\frac{1}{20} \, \overline{\alpha_{1}} & C_{1} & 1 & 5 & 0 \\
	    			\frac{1}{8} & \overline{C_{1}} & 2 & 2 & 0 \\
	    			-\frac{1}{36} & \gamma_{1} & 2 & 3 & 0 \\
	    			\frac{1}{6} & C_{2} & 2 & 3 & 0 \\
	    			\frac{1}{6} & \gamma_{1} & 2 & 3 & 1 \\
	    				\end{dmatrix}
	    			\begin{dmatrix}
	    			-\frac{1}{128} & \overline{\gamma_{2}} & 2 & 4 & 0 \\
	    			\frac{1}{16} & \overline{C_{3}} & 2 & 4 & 0 \\
	    			\frac{1}{8} \, {\left| A_{1} \right|}^{2} & C_{1} & 2 & 4 & 0 \\
	    			\frac{1}{16} \, \overline{\alpha_{1}} & \overline{C_{1}} & 2 & 4 & 0 \\
	    			\frac{1}{16} & \overline{\gamma_{2}} & 2 & 4 & 1 \\
	    			\frac{1}{8} & \overline{B_{1}} & 3 & 2 & 0 \\
	    			-\frac{1}{72} & \gamma_{2} & 3 & 3 & 0 \\
	    			\frac{1}{12} & C_{3} & 3 & 3 & 0 \\
	    			\frac{1}{6} \, {\left| A_{1} \right|}^{2} & \overline{C_{1}} & 3 & 3 & 0 \\
	    			\frac{1}{12} \, \alpha_{1} & C_{1} & 3 & 3 & 0 \\
	    			\frac{1}{12} & \gamma_{2} & 3 & 3 & 1 \\
	    			\frac{1}{8} & \overline{B_{2}} & 4 & 2 & 0 \\
	    			\frac{1}{8} \, \alpha_{1} & \overline{C_{1}} & 4 & 2 & 0 \\
	    			\frac{1}{4} & \overline{E_{1}} & 5 & 1 & 0
	    			\end{dmatrix}
	    		\end{align*}
	    		By conformality, we obtain
	    		\small
	    		\begin{align*}
	    			&0=\s{\p{z}\phi}{\p{z}\phi}=\\
	    			&\begin{dmatrix}
	    			A_{0}^{2} & 4 & 0 & 0 \\
	    			2 \, A_{0} A_{1} & 5 & 0 & 0 \\
	    			A_{1}^{2} + 2 \, A_{0} A_{2} & 6 & 0 & 0 \\
	    			2 \, A_{1} A_{2} + 2 \, A_{0} A_{3} & 7 & 0 & 0 \\
	    			A_{2}^{2} + 2 \, A_{1} A_{3} + 2 \, A_{0} A_{4} & 8 & 0 & 0 \\
	    			\frac{1}{144} \, C_{1}^{2} + \frac{1}{12} \, A_{0} E_{1} & 2 & 6 & 0 \\
	    			\frac{1}{6} \, A_{0} C_{1} & 3 & 3 & 0 \\
	    			\frac{1}{8} \, A_{0} B_{1} & 3 & 4 & 0 \\
	    			\frac{1}{10} \, A_{0} C_{1} \overline{\alpha_{1}} + \frac{1}{10} \, A_{0} B_{2} + \frac{1}{48} \, C_{1} \overline{C_{1}} & 3 & 5 & 0 \\
	    			\frac{1}{4} \, A_{0} \overline{C_{1}} & 4 & 2 & 0 \\
	    			\frac{1}{6} \, A_{1} C_{1} + \frac{1}{3} \, A_{0} C_{2} - \frac{1}{18} \, A_{0} \gamma_{1} & 4 & 3 & 0 \\
	    			\frac{1}{3} \, A_{0} \gamma_{1} & 4 & 3 & 1 \\
	    			\frac{1}{4} \, A_{0} C_{1} {\left| A_{1} \right|}^{2} + \frac{1}{8} \, A_{0} \overline{C_{1}} \overline{\alpha_{1}} + \frac{1}{8} \, A_{1} B_{1} + \frac{1}{64} \, \overline{C_{1}}^{2} + \frac{1}{8} \, A_{0} \overline{C_{3}} - \frac{1}{64} \, A_{0} \overline{\gamma_{2}} & 4 & 4 & 0 \\
	    			\frac{1}{8} \, A_{0} \overline{\gamma_{2}} & 4 & 4 & 1 \\
	    			\frac{1}{4} \, A_{0} \overline{B_{1}} + \frac{1}{4} \, A_{1} \overline{C_{1}} & 5 & 2 & 0 \\
	    			\frac{1}{3} \, A_{0} {\left| A_{1} \right|}^{2} \overline{C_{1}} + \frac{1}{6} \, A_{0} C_{1} \alpha_{1} + \frac{1}{6} \, A_{2} C_{1} + \frac{1}{3} \, A_{1} C_{2} + \frac{1}{6} \, A_{0} C_{3} - \frac{1}{18} \, A_{1} \gamma_{1} - \frac{1}{36} \, A_{0} \gamma_{2} & 5 & 3 & 0 \\
	    			\frac{1}{3} \, A_{1} \gamma_{1} + \frac{1}{6} \, A_{0} \gamma_{2} & 5 & 3 & 1 \\
	    			\frac{1}{4} \, A_{0} \alpha_{1} \overline{C_{1}} + \frac{1}{4} \, A_{1} \overline{B_{1}} + \frac{1}{4} \, A_{0} \overline{B_{2}} + \frac{1}{4} \, A_{2} \overline{C_{1}} & 6 & 2 & 0 \\
	    			\frac{1}{2} \, A_{0} \overline{E_{1}} & 7 & 1 & 0
	    			\end{dmatrix}
	    		\end{align*}
	    		\normalsize
	    		Now, we obtain
	    		\small
	    		\begin{align*}
	    			&e^{2\lambda}=\begin{dmatrix}
	    			2 \, A_{0} \overline{A_{0}} & 2 & 2 & 0 \\
	    			2 \, A_{0} \overline{A_{1}} & 2 & 3 & 0 \\
	    			2 \, A_{0} \overline{A_{2}} + \frac{1}{4} \, \overline{A_{0}} \overline{C_{1}} & 2 & 4 & 0 \\
	    			\frac{1}{3} \, C_{2} \overline{A_{0}} - \frac{1}{18} \, \gamma_{1} \overline{A_{0}} + 2 \, A_{0} \overline{A_{3}} + \frac{1}{4} \, \overline{A_{1}} \overline{C_{1}} & 2 & 5 & 0 \\
	    			\frac{1}{4} \, C_{1} {\left| A_{1} \right|}^{2} \overline{A_{0}} + \frac{1}{8} \, \overline{A_{0}} \overline{C_{1}} \overline{\alpha_{1}} + \frac{1}{3} \, C_{2} \overline{A_{1}} - \frac{1}{18} \, \gamma_{1} \overline{A_{1}} + 2 \, A_{0} \overline{A_{4}} + \frac{1}{4} \, \overline{A_{2}} \overline{C_{1}} + \frac{1}{8} \, \overline{A_{0}} \overline{C_{3}} - \frac{1}{64} \, \overline{A_{0}} \overline{\gamma_{2}} & 2 & 6 & 0 \\
	    			\frac{1}{12} \, A_{0} \overline{E_{1}} & 8 & 0 & 0 \\
	    			\frac{1}{6} \, A_{0} \overline{C_{1}} & 5 & 1 & 0 \\
	    			\frac{1}{8} \, A_{0} \overline{B_{1}} + \frac{1}{6} \, A_{1} \overline{C_{1}} & 6 & 1 & 0 \\
	    			\frac{1}{10} \, A_{0} \alpha_{1} \overline{C_{1}} + \frac{1}{8} \, A_{1} \overline{B_{1}} + \frac{1}{10} \, A_{0} \overline{B_{2}} + \frac{1}{6} \, A_{2} \overline{C_{1}} & 7 & 1 & 0 \\
	    			\frac{1}{4} \, A_{0} C_{1} + 2 \, A_{2} \overline{A_{0}} & 4 & 2 & 0 \\
	    			\frac{1}{4} \, A_{1} C_{1} + \frac{1}{3} \, A_{0} C_{2} - \frac{1}{18} \, A_{0} \gamma_{1} + 2 \, A_{3} \overline{A_{0}} & 5 & 2 & 0 \\
	    			\frac{1}{3} \, A_{0} \gamma_{1} & 5 & 2 & 1 \\
	    			\frac{1}{4} \, A_{0} {\left| A_{1} \right|}^{2} \overline{C_{1}} + \frac{1}{8} \, A_{0} C_{1} \alpha_{1} + \frac{1}{4} \, A_{2} C_{1} + \frac{1}{3} \, A_{1} C_{2} + \frac{1}{8} \, A_{0} C_{3} - \frac{1}{18} \, A_{1} \gamma_{1} - \frac{1}{64} \, A_{0} \gamma_{2} + 2 \, A_{4} \overline{A_{0}} & 6 & 2 & 0 \\
	    			\frac{1}{3} \, A_{1} \gamma_{1} + \frac{1}{8} \, A_{0} \gamma_{2} & 6 & 2 & 1 \\
	    			\frac{1}{4} \, A_{0} B_{1} + 2 \, A_{2} \overline{A_{1}} & 4 & 3 & 0 \\
	    			\frac{1}{3} \, A_{0} C_{1} {\left| A_{1} \right|}^{2} + \frac{1}{6} \, A_{0} \overline{C_{1}} \overline{\alpha_{1}} + \frac{1}{4} \, A_{1} B_{1} + 2 \, A_{3} \overline{A_{1}} + \frac{1}{48} \, \overline{C_{1}}^{2} + \frac{1}{6} \, A_{0} \overline{C_{3}} + \frac{1}{2} \, \overline{A_{0}} \overline{E_{1}} - \frac{1}{36} \, A_{0} \overline{\gamma_{2}} & 5 & 3 & 0 \\
	    			\frac{1}{6} \, A_{0} \overline{\gamma_{2}} & 5 & 3 & 1 \\
	    			\frac{1}{4} \, \alpha_{1} \overline{A_{0}} \overline{C_{1}} + \frac{1}{4} \, A_{0} C_{1} \overline{\alpha_{1}} + \frac{1}{4} \, A_{0} B_{2} + 2 \, A_{2} \overline{A_{2}} + \frac{1}{4} \, \overline{A_{0}} \overline{B_{2}} + \frac{13}{288} \, C_{1} \overline{C_{1}} & 4 & 4 & 0 \\
	    			\frac{1}{3} \, {\left| A_{1} \right|}^{2} \overline{A_{0}} \overline{C_{1}} + \frac{1}{6} \, C_{1} \alpha_{1} \overline{A_{0}} + \frac{1}{48} \, C_{1}^{2} + \frac{1}{2} \, A_{0} E_{1} + \frac{1}{6} \, C_{3} \overline{A_{0}} - \frac{1}{36} \, \gamma_{2} \overline{A_{0}} + 2 \, A_{1} \overline{A_{3}} + \frac{1}{4} \, \overline{A_{1}} \overline{B_{1}} & 3 & 5 & 0 \\
	    			2 \, A_{1} \overline{A_{0}} & 3 & 2 & 0 \\
	    			2 \, A_{1} \overline{A_{1}} & 3 & 3 & 0 \\
	    			2 \, A_{1} \overline{A_{2}} + \frac{1}{4} \, \overline{A_{0}} \overline{B_{1}} & 3 & 4 & 0 \\
	    			\frac{1}{12} \, E_{1} \overline{A_{0}} & 0 & 8 & 0 \\
	    			\frac{1}{6} \, C_{1} \overline{A_{0}} & 1 & 5 & 0 \\
	    			\frac{1}{8} \, B_{1} \overline{A_{0}} + \frac{1}{6} \, C_{1} \overline{A_{1}} & 1 & 6 & 0 \\
	    			\frac{1}{10} \, C_{1} \overline{A_{0}} \overline{\alpha_{1}} + \frac{1}{10} \, B_{2} \overline{A_{0}} + \frac{1}{8} \, B_{1} \overline{A_{1}} + \frac{1}{6} \, C_{1} \overline{A_{2}} & 1 & 7 & 0 \\
	    			\frac{1}{3} \, \gamma_{1} \overline{A_{0}} & 2 & 5 & 1 \\
	    			\frac{1}{3} \, \gamma_{1} \overline{A_{1}} + \frac{1}{8} \, \overline{A_{0}} \overline{\gamma_{2}} & 2 & 6 & 1 \\
	    			\frac{1}{6} \, \gamma_{2} \overline{A_{0}} & 3 & 5 & 1
	    			\end{dmatrix}
	    		\end{align*}
	    		\begin{align*}
	    			&e^{2\lambda}=\\
	    			&\begin{dmatrix}
	    			2 \, A_{0} \overline{A_{0}} & 2 & 2 & 0 &\textbf{(1)}\\
	    			2 \, A_{0} \overline{A_{1}} & 2 & 3 & 0 &\textbf{(2)}\\
	    			2 \, A_{0} \overline{A_{2}} + \frac{1}{4} \, \overline{A_{0}} \overline{C_{1}} & 2 & 4 & 0 &\textbf{(3)}\\
	    			\frac{1}{3} \, C_{2} \overline{A_{0}} - \frac{1}{18} \, \gamma_{1} \overline{A_{0}} + 2 \, A_{0} \overline{A_{3}} + \frac{1}{4} \, \overline{A_{1}} \overline{C_{1}} & 2 & 5 & 0 &\textbf{(4)}\\
	    			\frac{1}{4} \ccancel{\, C_{1} {\left| A_{1} \right|}^{2} \overline{A_{0}}} + \frac{1}{8} \ccancel{\, \overline{A_{0}} \overline{C_{1}} \overline{\alpha_{1}}} + \frac{1}{3} \, C_{2} \overline{A_{1}} - \frac{1}{18} \, \gamma_{1} \overline{A_{1}} + 2 \, A_{0} \overline{A_{4}} + \frac{1}{4} \, \overline{A_{2}} \overline{C_{1}} + \frac{1}{8} \, \overline{A_{0}} \overline{C_{3}} - \frac{1}{64} \, \overline{A_{0}} \overline{\gamma_{2}} & 2 & 6 & 0 &\textbf{(5)}\\
	    			\frac{1}{12} \ccancel{\, A_{0} \overline{E_{1}}} & 8 & 0 & 0 &\ccancel{\textbf{(6)}}\\
	    			\frac{1}{6} \, A_{0} \overline{C_{1}} & 5 & 1 & 0 &\textbf{(7)}\\
	    			\frac{1}{8} \, A_{0} \overline{B_{1}} + \frac{1}{6} \, A_{1} \overline{C_{1}} & 6 & 1 & 0 &\textbf{(8)}\\
	    			\frac{1}{10} \ccancel{\, A_{0} \alpha_{1} \overline{C_{1}}} + \frac{1}{8} \ccancel{\, A_{1} \overline{B_{1}}} + \frac{1}{10} \, A_{0} \overline{B_{2}} + \frac{1}{6} \, A_{2} \overline{C_{1}} & 7 & 1 & 0 &\textbf{(9)}\\
	    			\frac{1}{4} \, A_{0} C_{1} + 2 \, A_{2} \overline{A_{0}} & 4 & 2 & 0 &\textbf{(10)}\\
	    			\frac{1}{4} \, A_{1} C_{1} + \frac{1}{3} \, A_{0} C_{2} - \frac{1}{18} \, A_{0} \gamma_{1} + 2 \, A_{3} \overline{A_{0}} & 5 & 2 & 0 &\textbf{(11)}\\
	    			\frac{1}{3} \, A_{0} \gamma_{1} & 5 & 2 & 1 &\textbf{(12)}\\
	    			\frac{1}{4} \ccancel{\, A_{0} {\left| A_{1} \right|}^{2} \overline{C_{1}}} + \frac{1}{8} \ccancel{\, A_{0} C_{1} \alpha_{1}} + \frac{1}{4} \, A_{2} C_{1} + \frac{1}{3} \, A_{1} C_{2} + \frac{1}{8} \, A_{0} C_{3} - \frac{1}{18} \, A_{1} \gamma_{1} - \frac{1}{64} \, A_{0} \gamma_{2} + 2 \, A_{4} \overline{A_{0}} & 6 & 2 & 0 &\textbf{(13)}\\
	    			\frac{1}{3} \, A_{1} \gamma_{1} + \frac{1}{8} \, A_{0} \gamma_{2} & 6 & 2 & 1 &\textbf{(14)}\\
	    			\frac{1}{4} \, A_{0} B_{1} + 2 \, A_{2} \overline{A_{1}} & 4 & 3 & 0 &\textbf{(15)}\\
	    			\frac{1}{3} \ccancel{\, A_{0} C_{1} {\left| A_{1} \right|}^{2}} + \frac{1}{6} \ccancel{\, A_{0} \overline{C_{1}} \overline{\alpha_{1}}} + \frac{1}{4} \ccancel{\, A_{1} B_{1}} + 2 \, A_{3} \overline{A_{1}} + \frac{1}{48} \, \overline{C_{1}}^{2} + \frac{1}{6} \, A_{0} \overline{C_{3}} + \frac{1}{2} \, \overline{A_{0}} \overline{E_{1}} - \frac{1}{36} \, \ccancel{A_{0} \overline{\gamma_{2}}} & 5 & 3 & 0 &\textbf{(16)}\\
	    			\frac{1}{6} \ccancel{\, A_{0} \overline{\gamma_{2}}} & 5 & 3 & 1 &\ccancel{\textbf{(17)}}\\
	    			\frac{1}{4} \ccancel{\, \alpha_{1} \overline{A_{0}} \overline{C_{1}}} + \frac{1}{4} \ccancel{\, A_{0} C_{1} \overline{\alpha_{1}}} + \frac{1}{4} \, A_{0} B_{2} + 2 \, A_{2} \overline{A_{2}} + \frac{1}{4} \, \overline{A_{0}} \overline{B_{2}} + \frac{13}{288} \, C_{1} \overline{C_{1}} & 4 & 4 & 0 &\textbf{(18)}\\
	    			\frac{1}{3} \, \ccancel{{\left| A_{1} \right|}^{2} \overline{A_{0}} \overline{C_{1}}} + \frac{1}{6} \ccancel{\, C_{1} \alpha_{1} \overline{A_{0}}} + \frac{1}{48} \, C_{1}^{2} + \frac{1}{2} \, A_{0} E_{1} + \frac{1}{6} \, C_{3} \overline{A_{0}} - \frac{1}{36} \ccancel{\, \gamma_{2} \overline{A_{0}}} + 2 \, A_{1} \overline{A_{3}} + \frac{1}{4} \ccancel{\, \overline{A_{1}} \overline{B_{1}}} & 3 & 5 & 0 &\textbf{(19)}\\
	    			2 \, A_{1} \overline{A_{0}} & 3 & 2 & 0 &\textbf{(20)}\\
	    			2 \, A_{1} \overline{A_{1}} & 3 & 3 & 0 &\textbf{(21)}\\
	    			2 \, A_{1} \overline{A_{2}} + \frac{1}{4} \, \overline{A_{0}} \overline{B_{1}} & 3 & 4 & 0 &\textbf{(22)}\\
	    			\frac{1}{12} \ccancel{\, E_{1} \overline{A_{0}}} & 0 & 8 & 0 &\ccancel{\textbf{(23)}}\\
	    			\frac{1}{6} \, C_{1} \overline{A_{0}} & 1 & 5 & 0 &\textbf{(24)}\\
	    			\frac{1}{8} \, B_{1} \overline{A_{0}} + \frac{1}{6} \, C_{1} \overline{A_{1}} & 1 & 6 & 0 &\textbf{(25)}\\
	    			\frac{1}{10} \ccancel{\, C_{1} \overline{A_{0}} \overline{\alpha_{1}}} + \frac{1}{10} \, B_{2} \overline{A_{0}} + \frac{1}{8} \ccancel{\, B_{1} \overline{A_{1}}} + \frac{1}{6} \, C_{1} \overline{A_{2}} & 1 & 7 & 0 &\textbf{(26)}\\
	    			\frac{1}{3} \, \gamma_{1} \overline{A_{0}} & 2 & 5 & 1 &\textbf{(27)}\\
	    			\frac{1}{3} \, \gamma_{1} \overline{A_{1}} + \frac{1}{8} \, \overline{A_{0}} \overline{\gamma_{2}} & 2 & 6 & 1 &\textbf{(28)}\\
	    			\frac{1}{6} \, \ccancel{\gamma_{2} \overline{A_{0}}} & 3 & 5 & 1&\ccancel{\textbf{(29)}}
	    			\end{dmatrix}\\
	    			&=\begin{dmatrix}
	    			2 \, A_{0} \overline{A_{0}} & 2 & 2 & 0 &\textbf{(1)}\\
	    			2 \, A_{0} \overline{A_{1}} & 2 & 3 & 0 &\textbf{(2)}\\
	    			2 \, A_{0} \overline{A_{2}} + \frac{1}{4} \, \overline{A_{0}} \overline{C_{1}} & 2 & 4 & 0 &\textbf{(3)}\\
	    			\frac{1}{3} \, C_{2} \overline{A_{0}} - \frac{1}{18} \, \gamma_{1} \overline{A_{0}} + 2 \, A_{0} \overline{A_{3}} + \frac{1}{4} \, \overline{A_{1}} \overline{C_{1}} & 2 & 5 & 0 &\textbf{(4)}\\
	    			\frac{1}{3} \, C_{2} \overline{A_{1}} - \frac{1}{18} \, \gamma_{1} \overline{A_{1}} + 2 \, A_{0} \overline{A_{4}} + \frac{1}{4} \, \overline{A_{2}} \overline{C_{1}} +\frac{1}{8} \, \overline{A_{0}} \overline{C_{3}}- \frac{1}{64} \, \overline{A_{0}} \overline{\gamma_{2}} & 2 & 6 & 0 &\textbf{(5)}\\
	    			\frac{1}{6} \, A_{0} \overline{C_{1}} & 5 & 1 & 0 &\textbf{(7)}\\
	    			\frac{1}{8} \, A_{0} \overline{B_{1}} + \frac{1}{6} \, A_{1} \overline{C_{1}} & 6 & 1 & 0 &\textbf{(8)}\\
	    			\frac{1}{10} \, A_{0} \overline{B_{2}} + \frac{1}{6} \, A_{2} \overline{C_{1}} & 7 & 1 & 0 &\textbf{(9)}\\
	    			\frac{1}{4} \, A_{0} C_{1} + 2 \, A_{2} \overline{A_{0}} & 4 & 2 & 0 &\textbf{(10)}\\
	    			\frac{1}{4} \, A_{1} C_{1} + \frac{1}{3} \, A_{0} C_{2} - \frac{1}{18} \, A_{0} \gamma_{1} + 2 \, A_{3} \overline{A_{0}} & 5 & 2 & 0 &\textbf{(11)}\\
	    			\frac{1}{3} \, A_{0} \gamma_{1} & 5 & 2 & 1 &\textbf{(12)}\\
	    			 \frac{1}{4} \, A_{2} C_{1} + \frac{1}{3} \, A_{1} C_{2} + \frac{1}{8} \, A_{0} C_{3} - \frac{1}{18} \, A_{1} \gamma_{1} - \frac{1}{64} \, A_{0} \gamma_{2} + 2 \, A_{4} \overline{A_{0}} & 6 & 2 & 0 &\textbf{(13)}\\
	    			\frac{1}{3} \, A_{1} \gamma_{1} + \frac{1}{8} \, A_{0} \gamma_{2} & 6 & 2 & 1 &\textbf{(14)}\\
	    			\frac{1}{4} \, A_{0} B_{1} + 2 \, A_{2} \overline{A_{1}} & 4 & 3 & 0 &\textbf{(15)}\\
	    			2 \, A_{3} \overline{A_{1}} + \frac{1}{48} \, \overline{C_{1}}^{2} + + \frac{1}{6} \, A_{0} \overline{C_{3}} + \frac{1}{2} \, \overline{A_{0}} \overline{E_{1}} & 5 & 3 & 0 &\textbf{(16)}\\
	    			\frac{1}{4} \, A_{0} B_{2} + 2 \, A_{2} \overline{A_{2}} + \frac{1}{4} \, \overline{A_{0}} \overline{B_{2}} + \frac{13}{288} \, C_{1} \overline{C_{1}} & 4 & 4 & 0 &\textbf{(18)}\\
	    			\frac{1}{48} \, C_{1}^{2} + \frac{1}{2} \, A_{0} E_{1} + \frac{1}{6} \, C_{3} \overline{A_{0}} + 2 \, A_{1} \overline{A_{3}} & 3 & 5 & 0 &\textbf{(19)}\\
	    			2 \, A_{1} \overline{A_{0}} & 3 & 2 & 0 &\textbf{(20)}\\
	    			2 \, A_{1} \overline{A_{1}} & 3 & 3 & 0 &\textbf{(21)}\\
	    			2 \, A_{1} \overline{A_{2}} + \frac{1}{4} \, \overline{A_{0}} \overline{B_{1}} & 3 & 4 & 0 &\textbf{(22)}\\
	    			\frac{1}{6} \, C_{1} \overline{A_{0}} & 1 & 5 & 0 &\textbf{(24)}\\
	    			\frac{1}{8} \, B_{1} \overline{A_{0}} + \frac{1}{6} \, C_{1} \overline{A_{1}} & 1 & 6 & 0 &\textbf{(25)}\\
	    			\frac{1}{10} \, B_{2} \overline{A_{0}} +\frac{1}{6} \, C_{1} \overline{A_{2}} & 1 & 7 & 0 &\textbf{(26)}\\
	    			\frac{1}{3} \, \gamma_{1} \overline{A_{0}} & 2 & 5 & 1 &\textbf{(27)}\\
	    			\frac{1}{3} \, \gamma_{1} \overline{A_{1}} + \frac{1}{8} \, \overline{A_{0}} \overline{\gamma_{2}} & 2 & 6 & 1 &\textbf{(28)}\\
	    			\end{dmatrix}
	    		\end{align*}
	    		Then, recall that
	    		\begin{align*}
	    			e^{2\lambda}=\begin{dmatrix}
	    			1 & 2 & 2 & 0 \\
	    			2 \, {\left| A_{1} \right|}^{2} & 3 & 3 & 0 \\
	    			\alpha_{1} & 4 & 2 & 0 \\
	    			\alpha_{2} & 1 & 6 & 0 \\
	    			\alpha_{3} & 5 & 2 & 0 \\
	    			\alpha_{4} & 4 & 3 & 0 
	    			\end{dmatrix}
	    			\begin{dmatrix}
	    			\overline{\alpha_{1}} & 2 & 4 & 0 \\
	    			\overline{\alpha_{2}} & 6 & 1 & 0 \\
	    			\overline{\alpha_{3}} & 2 & 5 & 0 \\
	    			\overline{\alpha_{4}} & 3 & 4 & 0
	    			\end{dmatrix}
	    		\end{align*}
	    		so we need only look at powers of order $2\theta_0+2=8$, which are
	    		\begin{align*}
	    		    (1\;\, 7),\; (7\;\, 1),\; 
	    			(2\;\, 6),\; (6\;\, 2),\; (2\;\, 6\;\, 1),\; (6\;\, 2\;\, 1),\; (3\;\, 5),\; (5\;\, 3),\; (4\;\, 4).	    	
	    		\end{align*}
	    		Therefore, there exists $\alpha_5,\alpha_6,\alpha_7,\zeta_0\in \mathbb{C}$ and $\beta\in \R$ such that (left \TeX\; and right Sage)
	    		\begin{align*}
	    			e^{2\lambda}=\begin{dmatrix}
	    			1 & 2 & 2 & 0 \\
	    			2 \, {\left| A_{1} \right|}^{2} & 3 & 3 & 0 \\
	    			\beta & 4 & 4 & 0 \\
	    			\alpha_{1} & 4 & 2 & 0 \\
	    			\alpha_{2} & 1 & 6 & 0 \\
	    			\alpha_{3} & 5 & 2 & 0 \\
	    			\alpha_{4} & 4 & 3 & 0 \\
	    			\alpha_{5} & 7 & 1 & 0\\
	    			\alpha_{6} & 6 & 2 & 0\\
	    			\alpha_{7} & 3 & 5 & 0\\\
	    			\zeta_{0} & 6 & 2 & 1 \\
	    			\end{dmatrix}
	    			\begin{dmatrix}
	    			\overline{\alpha_{1}} & 2 & 4 & 0 \\
	    			\overline{\alpha_{2}} & 6 & 1 & 0 \\
	    			\overline{\alpha_{3}} & 2 & 5 & 0 \\
	    			\overline{\alpha_{4}} & 3 & 4 & 0 \\
	    			\bar{\alpha_{5}} & 1 & 7 & 0\\
	    			\bar{\alpha_{6}} & 2 & 6 & 0\\
	    			\bar{\alpha_{7}} & 5 & 3 & 0\\
	    			\bar{\zeta_{0}} & 2 & 6 & 1 \\
	    			\end{dmatrix}
	    			\begin{dmatrix}
	    			1 & 2 & 2 & 0 \\
	    			2 \, {\left| A_{1} \right|}^{2} & 3 & 3 & 0 \\
	    			\beta & 4 & 4 & 0 \\
	    			\alpha_{1} & 4 & 2 & 0 \\
	    			\alpha_{2} & 1 & 6 & 0 \\
	    			\alpha_{3} & 5 & 2 & 0 \\
	    			\alpha_{4} & 4 & 3 & 0 \\
	    			\alpha_{5} & 7 & 1 & 0 \\
	    			\alpha_{6} & 6 & 2 & 0 \\
	    			\alpha_{7} & 3 & 5 & 0 \\
	    			\zeta_{0} & 6 & 2 & 1 
	    			\end{dmatrix}
	    			\begin{dmatrix}
	    			\overline{\alpha_{1}} & 2 & 4 & 0 \\
	    			\overline{\alpha_{2}} & 6 & 1 & 0 \\
	    			\overline{\alpha_{3}} & 2 & 5 & 0 \\
	    			\overline{\alpha_{4}} & 3 & 4 & 0 \\
	    			\overline{\alpha_{5}} & 1 & 7 & 0 \\
	    			\overline{\alpha_{6}} & 2 & 6 & 0 \\
	    			\overline{\alpha_{7}} & 5 & 3 & 0 \\
	    			\overline{\zeta_{0}} & 2 & 6 & 1
	    			\end{dmatrix}
	    		\end{align*}
	    		\normalsize
    		Then, by developing the quartic form up to order $3$, the only term involving logarithm terms in the Taylor expansion are
    		\begin{align*}
    			-2\,\zeta_0\s{\vec{A}_1}{\vec{C}_1}z^3\log|z|dz^4 -2\,\bar{\zeta_0}\s{\vec{A}_1}{\vec{C}_1}\z^3\log|z|dz^4
    		\end{align*}
    		so we obtain
    		\begin{align}\label{3extracancel}
    			\zeta_0\s{\vec{A}_1}{\vec{C}_1}=0
    		\end{align}
    		and as
    		\begin{align*}
    			\s{\vec{A}_1}{\vec{\gamma}_1}+2\s{\vec{A}_0}{\vec{\gamma}_2}=0
    		\end{align*}
    		Then, we have
    		\begin{align*}
    			\zeta_0=\frac{1}{3} \, A_{1} \gamma_{1} + \frac{1}{8} \, A_{0} \gamma_{2}=\frac{1}{12}\s{\vec{A}_1}{\vec{\gamma}_1}=-\frac{1}{12}\s{\vec{A}_1}{\vec{\gamma}_0}
    		\end{align*}
    		so we obtain
    		\begin{align*}
    			\s{\vec{A}_1}{\vec{\gamma}_0}\s{\vec{A}_1}{\vec{C}_1}=0
    		\end{align*}
    		so $\s{\vec{A}_1}{\vec{C}_1}=0$ or $\s{\vec{A}_1}{\vec{\gamma}_0}=0$, but the last relation does not imply anything useful, as we have
    		\begin{align*}
    			\s{\vec{A}_0}{\vec{\gamma}_0}+\s{\vec{A}_1}{\vec{C}_1}=0.
    		\end{align*}
    		Then, we have
    		\small
    		\begin{align*}
    		\h_0=
    		\begin{dmatrix}
    		-\frac{1}{6} & E_{1} & -1 & 6 & 0 \\
    		-\frac{1}{6} & C_{1} & 0 & 3 & 0 \\
    		-\frac{1}{8} & B_{1} & 0 & 4 & 0 \\
    		2 \, \alpha_{2} & A_{0} & 0 & 4 & 0 \\
    		-\frac{1}{10} & B_{2} & 0 & 5 & 0 \\
    		-\frac{1}{10} \, \overline{\alpha_{1}} & C_{1} & 0 & 5 & 0 \\
    		2 \, \overline{\alpha_{5}} & A_{0} & 0 & 5 & 0 \\
    		\frac{1}{6} & \gamma_{1} & 1 & 3 & 0 \\
    		-\frac{1}{3} \, {\left| A_{1} \right|}^{2} & C_{1} & 1 & 4 & 0 \\
    		\frac{1}{16} & \overline{\gamma_{2}} & 1 & 4 & 0 \\
    		2 \, \alpha_{2} & A_{1} & 1 & 4 & 0 \\
    		-\overline{\zeta_{0}} & A_{0} & 1 & 4 & 0 \\
    		2 & A_{1} & 2 & 0 & 0 \\
    		-4 \, {\left| A_{1} \right|}^{2} & A_{0} & 2 & 1 & 0 \\
    		\frac{1}{4} & \overline{B_{1}} & 2 & 2 & 0 \\
    		-2 \, \overline{\alpha_{4}} & A_{0} & 2 & 2 & 0 \\
    		-\frac{1}{6} \, {\left| A_{1} \right|}^{2} & \overline{C_{1}} & 2 & 3 & 0 \\
    		\frac{1}{18} & \gamma_{2} & 2 & 3 & 0 \\
    		\frac{1}{6} & C_{3} & 2 & 3 & 0 \\
    		-\frac{1}{6} \, \alpha_{1} & C_{1} & 2 & 3 & 0 \\
    		4 \, {\left| A_{1} \right|}^{2} \overline{\alpha_{1}} - 2 \, \alpha_{7} & A_{0} & 2 & 3 & 0 \\
    		\frac{1}{6} & \gamma_{2} & 2 & 3 & 1 
    		\end{dmatrix}
    		\begin{dmatrix}
    		4 & A_{2} & 3 & 0 & 0 \\
    		-4 \, \alpha_{1} & A_{0} & 3 & 0 & 0 \\
    		-4 \, {\left| A_{1} \right|}^{2} & A_{1} & 3 & 1 & 0 \\
    		-4 \, \alpha_{4} & A_{0} & 3 & 1 & 0 \\
    		\frac{1}{2} & \overline{B_{2}} & 3 & 2 & 0 \\
    		8 \, {\left| A_{1} \right|}^{4} + 4 \, \alpha_{1} \overline{\alpha_{1}} - 4 \, \beta & A_{0} & 3 & 2 & 0 \\
    		-2 \, \overline{\alpha_{4}} & A_{1} & 3 & 2 & 0 \\
    		6 & A_{3} & 4 & 0 & 0 \\
    		-4 \, \alpha_{1} & A_{1} & 4 & 0 & 0 \\
    		-6 \, \alpha_{3} & A_{0} & 4 & 0 & 0 \\
    		-4 \, {\left| A_{1} \right|}^{2} & A_{2} & 4 & 1 & 0 \\
    		\frac{3}{2} & \overline{E_{1}} & 4 & 1 & 0 \\
    		-4 \, \alpha_{4} & A_{1} & 4 & 1 & 0 \\
    		12 \, \alpha_{1} {\left| A_{1} \right|}^{2} - 6 \, \overline{\alpha_{7}} & A_{0} & 4 & 1 & 0 \\
    		-8 \, \overline{\alpha_{2}} & A_{0} & 5 & -1 & 0 \\
    		8 & A_{4} & 5 & 0 & 0 \\
    		-4 \, \alpha_{1} & A_{2} & 5 & 0 & 0 \\
    		-6 \, \alpha_{3} & A_{1} & 5 & 0 & 0 \\
    		4 \, \alpha_{1}^{2} - 8 \, \alpha_{6} - \zeta_{0} & A_{0} & 5 & 0 & 0 \\
    		-8 \, \zeta_{0} & A_{0} & 5 & 0 & 1 \\
    		-10 \, \alpha_{5} & A_{0} & 6 & -1 & 0 \\
    		-8 \, \overline{\alpha_{2}} & A_{1} & 6 & -1 & 0
    		\end{dmatrix}
    		\end{align*}
    		\normalsize
    		
    		Now, we need the next order development of the tensors, and we first have for some $\vec{D}_4\in \mathbb{C}^n$
    		\small
    		\begin{align*}
    			&\H+\vec{\gamma}_0\log|z|=\vec{C}_2+\Re\left(\vec{D}_3z+\vec{D}_4z^2\right)+\\
    			&\begin{dmatrix}
    			-\frac{1}{12} \, C_{1}^{2} & \overline{A_{0}} & -2 & 3 & 0& \textbf{(1)}\\
    			-\frac{7}{96} \, C_{1}^{2} & \overline{A_{1}} & -2 & 4 & 0& \textbf{(2)}\\
    			\frac{1}{2} \, A_{0} C_{1} \alpha_{2} - \frac{13}{96} \, B_{1} C_{1} & \overline{A_{0}} & -2 & 4 & 0 & \textbf{(3)}\\
    			-\frac{1}{48} \, C_{1} \overline{A_{1}} & C_{1} & -2 & 4 & 0 & \textbf{(4)}\\
    			C_{1} \alpha_{2} \overline{A_{0}} - \frac{1}{4} \, E_{1} \overline{A_{1}} & A_{0} & -2 & 4 & 0 & \textbf{(5)}\\
    			-C_{1} \overline{A_{1}} & A_{0} & -1 & 1 & 0 & \textbf{(6)}\\
    			-\frac{5}{24} \, C_{1} \overline{C_{1}} & \overline{A_{0}} & -1 & 2 & 0 & \textbf{(7)}\\
    			C_{1} \overline{A_{0}} \overline{\alpha_{1}} - \frac{1}{2} \, B_{1} \overline{A_{1}} - C_{1} \overline{A_{2}} & A_{0} & -1 & 2 & 0 & \textbf{(8)}\\
    			A_{0} \alpha_{2} \overline{C_{1}} - \frac{1}{3} \, C_{1} C_{2} + A_{1} E_{1} + \frac{1}{36} \, C_{1} \gamma_{1} - \frac{7}{48} \, B_{1} \overline{C_{1}} & \overline{A_{0}} & -1 & 3 & 0 & \textbf{(9)}\\
    			-\frac{1}{6} \, C_{1} \overline{C_{1}} & \overline{A_{1}} & -1 & 3 & 0 & \textbf{(10)}\\
    			-\frac{1}{36} \, \overline{A_{1}} \overline{C_{1}} & C_{1} & -1 & 3 & 0 & \textbf{(11)}\\
    			\frac{4}{3} \, \alpha_{2} \overline{A_{0}} \overline{C_{1}} + C_{1} \overline{A_{0}} \overline{\alpha_{3}} - \frac{1}{3} \, B_{2} \overline{A_{1}} - \frac{2}{3} \, B_{1} \overline{A_{2}} - C_{1} \overline{A_{3}} + \frac{1}{3} \, {\left(2 \, B_{1} \overline{A_{0}} + 3 \, C_{1} \overline{A_{1}}\right)} \overline{\alpha_{1}} & A_{0} & -1 & 3 & 0 & \textbf{(12)}\\
    			-\frac{1}{24} \, C_{1} \overline{A_{1}} & \overline{C_{1}} & -1 & 3 & 0 & \textbf{(13)}\\
    			-\frac{1}{3} \, C_{1} \gamma_{1} & \overline{A_{0}} & -1 & 3 & 1 & \textbf{(14)}\\
    			-2 \, A_{1} C_{1} & \overline{A_{0}} & 0 & 0 & 1 & \textbf{(15)}\\
    			-2 \, \overline{A_{1}} \overline{C_{1}} & A_{0} & 0 & 0 & 1 & \textbf{(16)}\\
    			-\frac{1}{8} \, \overline{C_{1}}^{2} & \overline{A_{0}} & 0 & 1 & 0 & \textbf{(17)}\\
    			2 \, C_{1} {\left| A_{1} \right|}^{2} \overline{A_{0}} + 2 \, \overline{A_{0}} \overline{C_{1}} \overline{\alpha_{1}} - 2 \, C_{2} \overline{A_{1}} + \gamma_{1} \overline{A_{1}} - 2 \, \overline{A_{2}} \overline{C_{1}} & A_{0} & 0 & 1 & 0 & \textbf{(18)}\\
    			-C_{1} \overline{A_{1}} & A_{1} & 0 & 1 & 0 & \textbf{(19)}\\
    			-2 \, A_{1} C_{1} & \overline{A_{1}} & 0 & 1 & 1 & \textbf{(20)}\\
    			4 \, A_{0} C_{1} {\left| A_{1} \right|}^{2} - 2 \, A_{1} B_{1} & \overline{A_{0}} & 0 & 1 & 1 & \textbf{(21)}\\
    			-2 \, \gamma_{1} \overline{A_{1}} & A_{0} & 0 & 1 & 1 & \textbf{(22)}\\
    			-\frac{1}{16} \, \overline{C_{1}}^{2} & \overline{A_{1}} & 0 & 2 & 0 & \textbf{(23)}\\
    			-\frac{1}{8} \, C_{1} \overline{B_{1}} - \frac{1}{4} \, C_{2} \overline{C_{1}} + \frac{1}{16} \, \gamma_{1} \overline{C_{1}} & \overline{A_{0}} & 0 & 2 & 0 & \textbf{(24)}\\
    			-\frac{1}{16} \, \overline{A_{1}} \overline{C_{1}} & \overline{C_{1}} & 0 & 2 & 0 & \textbf{(25)}\\
    			\lambda_{1} & A_{0} & 0 & 2 & 0 & \textbf{(26)}\\
    			C_{1} \overline{A_{0}} \overline{\alpha_{1}} - \frac{1}{2} \, B_{1} \overline{A_{1}} - C_{1} \overline{A_{2}} & A_{1} & 0 & 2 & 0 & \textbf{(27)}
    			\end{dmatrix}
    			\end{align*}
    			\begin{align*}
    			\begin{dmatrix}
    			-2 \, A_{1} C_{1} & \overline{A_{2}} & 0 & 2 & 1 & \textbf{(28)}\\
    			4 \, A_{0} B_{1} {\left| A_{1} \right|}^{2} + 2 \, A_{1} C_{1} \overline{\alpha_{1}} + 2 \, A_{0} C_{1} \overline{\alpha_{4}} - 2 \, A_{1} B_{2} - \frac{1}{12} \, C_{1} \overline{B_{1}} - \frac{5}{12} \, \gamma_{1} \overline{C_{1}} & \overline{A_{0}} & 0 & 2 & 1 & \textbf{(29)}\\
    			4 \, A_{0} C_{1} {\left| A_{1} \right|}^{2} - 2 \, A_{1} B_{1} & \overline{A_{1}} & 0 & 2 & 1 & \textbf{(30)}\\
    			2 \, {\left(\overline{A_{0}} \overline{\alpha_{1}} - \overline{A_{2}}\right)} \gamma_{1} - \frac{1}{2} \, \overline{A_{1}} \overline{\gamma_{2}} & A_{0} & 0 & 2 & 1 & \textbf{(31)}\\
    			-A_{1} \overline{C_{1}} & \overline{A_{0}} & 1 & -1 & 0 & \textbf{(32)}\\
    			-\frac{1}{8} \, C_{1}^{2} & A_{0} & 1 & 0 & 0 & \textbf{(33)}\\
    			2 \, A_{0} {\left| A_{1} \right|}^{2} \overline{C_{1}} + 2 \, A_{0} C_{1} \alpha_{1} - 2 \, A_{2} C_{1} - 2 \, A_{1} C_{2} + A_{1} \gamma_{1} & \overline{A_{0}} & 1 & 0 & 0 & \textbf{(34)}\\
    			-A_{1} \overline{C_{1}} & \overline{A_{1}} & 1 & 0 & 0 & \textbf{(35)}\\
    			-2 \, A_{1} \gamma_{1} & \overline{A_{0}} & 1 & 0 & 1 & \textbf{(36)}\\
    			-2 \, \overline{A_{1}} \overline{C_{1}} & A_{1} & 1 & 0 & 1 & \textbf{(37)}\\
    			4 \, {\left| A_{1} \right|}^{2} \overline{A_{0}} \overline{C_{1}} - 2 \, \overline{A_{1}} \overline{B_{1}} & A_{0} & 1 & 0 & 1 & \textbf{(38)}\\
    			\lambda_2 & A_{0} & 1 & 1 & 0 & \textbf{(39)}\\
    			\lambda_3 & \overline{A_{0}} & 1 & 1 & 0 & \textbf{(40)}\\
    			2 \, A_{0} {\left| A_{1} \right|}^{2} \overline{C_{1}} + 2 \, A_{0} C_{1} \alpha_{1} - 2 \, A_{2} C_{1} - 2 \, A_{1} C_{2} + A_{1} \gamma_{1} & \overline{A_{1}} & 1 & 1 & 0 & \textbf{(41)}\\
    			-A_{1} \overline{C_{1}} & \overline{A_{2}} & 1 & 1 & 0 & \textbf{(42)}\\
    			2 \, C_{1} {\left| A_{1} \right|}^{2} \overline{A_{0}} + 2 \, \overline{A_{0}} \overline{C_{1}} \overline{\alpha_{1}} - 2 \, C_{2} \overline{A_{1}} + \gamma_{1} \overline{A_{1}} - 2 \, \overline{A_{2}} \overline{C_{1}} & A_{1} & 1 & 1 & 0 & \textbf{(43)}\\
    			-C_{1} \overline{A_{1}} & A_{2} & 1 & 1 & 0 & \textbf{(44)}\\
    			-2 \, A_{1} \gamma_{1} & \overline{A_{1}} & 1 & 1 & 1 & \textbf{(45)}\\
    			4 \, A_{0} \gamma_{1} {\left| A_{1} \right|}^{2} - A_{1} \overline{\gamma_{2}} & \overline{A_{0}} & 1 & 1 & 1 & \textbf{(46)}\\
    			-2 \, \gamma_{1} \overline{A_{1}} & A_{1} & 1 & 1 & 1 & \textbf{(47)}\\
    			4 \, \gamma_{1} {\left| A_{1} \right|}^{2} \overline{A_{0}} - \gamma_{2} \overline{A_{1}} & A_{0} & 1 & 1 & 1 & \textbf{(48)}\\
    			-\frac{5}{24} \, C_{1} \overline{C_{1}} & A_{0} & 2 & -1 & 0 & \textbf{(49)}\\
    			A_{0} \alpha_{1} \overline{C_{1}} - \frac{1}{2} \, A_{1} \overline{B_{1}} - A_{2} \overline{C_{1}} & \overline{A_{0}} & 2 & -1 & 0 & \textbf{(50)}\\
    			-\frac{1}{16} \, C_{1}^{2} & A_{1} & 2 & 0 & 0 & \textbf{(51)}\\
    			-\frac{1}{4} \, C_{1} C_{2} + \frac{1}{16} \, C_{1} \gamma_{1} - \frac{1}{8} \, B_{1} \overline{C_{1}} & A_{0} & 2 & 0 & 0 & \textbf{(52)}\\
    			-\frac{1}{16} \, A_{1} C_{1} & C_{1} & 2 & 0 & 0 & \textbf{(53)}\\
    			\lambda_4 & \overline{A_{0}} & 2 & 0 & 0 & \textbf{(54)}\\
    			A_{0} \alpha_{1} \overline{C_{1}} - \frac{1}{2} \, A_{1} \overline{B_{1}} - A_{2} \overline{C_{1}} & \overline{A_{1}} & 2 & 0 & 0 & \textbf{(55)}
    			\end{dmatrix}
    			\end{align*}
    			\begin{align}\label{3hend}
    			\begin{dmatrix}
    			4 \, {\left| A_{1} \right|}^{2} \overline{A_{0}} \overline{B_{1}} + 2 \, \alpha_{4} \overline{A_{0}} \overline{C_{1}} + 2 \, \alpha_{1} \overline{A_{1}} \overline{C_{1}} - \frac{5}{12} \, C_{1} \gamma_{1} - 2 \, \overline{A_{1}} \overline{B_{2}} - \frac{1}{12} \, B_{1} \overline{C_{1}} & A_{0} & 2 & 0 & 1 & \textbf{(56)}\\
    			2 \, {\left(A_{0} \alpha_{1} - A_{2}\right)} \gamma_{1} - \frac{1}{2} \, A_{1} \gamma_{2} & \overline{A_{0}} & 2 & 0 & 1 & \textbf{(57)}\\
    			-2 \, \overline{A_{1}} \overline{C_{1}} & A_{2} & 2 & 0 & 1 & \textbf{(58)}\\
    			4 \, {\left| A_{1} \right|}^{2} \overline{A_{0}} \overline{C_{1}} - 2 \, \overline{A_{1}} \overline{B_{1}} & A_{1} & 2 & 0 & 1 & \textbf{(59)}\\
    			-\frac{1}{12} \, \overline{C_{1}}^{2} & A_{0} & 3 & -2 & 0 & \textbf{(60)}\\
    			-\frac{1}{6} \, C_{1} \overline{C_{1}} & A_{1} & 3 & -1 & 0 & \textbf{(61)}\\
    			C_{1} \overline{A_{0}} \overline{\alpha_{2}} - \frac{7}{48} \, C_{1} \overline{B_{1}} - \frac{1}{3} \, C_{2} \overline{C_{1}} + \frac{1}{36} \, \gamma_{1} \overline{C_{1}} + \overline{A_{1}} \overline{E_{1}} & A_{0} & 3 & -1 & 0 & \textbf{(62)}\\
    			-\frac{1}{36} \, A_{1} C_{1} & \overline{C_{1}} & 3 & -1 & 0 & \textbf{(63)}\\
    			A_{0} \alpha_{3} \overline{C_{1}} + \frac{4}{3} \, A_{0} C_{1} \overline{\alpha_{2}} + \frac{1}{3} \, {\left(2 \, A_{0} \overline{B_{1}} + 3 \, A_{1} \overline{C_{1}}\right)} \alpha_{1} - \frac{2}{3} \, A_{2} \overline{B_{1}} - \frac{1}{3} \, A_{1} \overline{B_{2}} - A_{3} \overline{C_{1}} & \overline{A_{0}} & 3 & -1 & 0 & \textbf{(64)}\\
    			-\frac{1}{24} \, A_{1} \overline{C_{1}} & C_{1} & 3 & -1 & 0 & \textbf{(65)}\\
    			-\frac{1}{3} \, \gamma_{1} \overline{C_{1}} & A_{0} & 3 & -1 & 1 & \textbf{(66)}\\
    			-\frac{7}{96} \, \overline{C_{1}}^{2} & A_{1} & 4 & -2 & 0 & \textbf{(67)}\\
    			\frac{1}{2} \, \overline{A_{0}} \overline{C_{1}} \overline{\alpha_{2}} - \frac{13}{96} \, \overline{B_{1}} \overline{C_{1}} & A_{0} & 4 & -2 & 0 & \textbf{(68)}\\
    			-\frac{1}{48} \, A_{1} \overline{C_{1}} & \overline{C_{1}} & 4 & -2 & 0 & \textbf{(69)}\\
    			A_{0} \overline{C_{1}} \overline{\alpha_{2}} - \frac{1}{4} \, A_{1} \overline{E_{1}} & \overline{A_{0}} & 4 & -2 & 0& \textbf{(70)}
    			\end{dmatrix}
    			\end{align}
    		\normalsize
    		where
    		\small
    		\begin{align*}
    			\lambda_{1}&=B_{1} {\left| A_{1} \right|}^{2} \overline{A_{0}} + 2 \, C_{1} {\left| A_{1} \right|}^{2} \overline{A_{1}} + \frac{3}{2} \, \overline{A_{0}} \overline{C_{1}} \overline{\alpha_{3}} + C_{1} \overline{A_{0}} \overline{\alpha_{4}} - \frac{1}{2} \, {\left(\overline{A_{0}} \overline{\alpha_{1}} - \overline{A_{2}}\right)} \gamma_{1} - 2 \, C_{2} \overline{A_{2}} - \frac{3}{2} \, \overline{A_{3}} \overline{C_{1}} - \frac{1}{2} \, \overline{A_{1}} \overline{C_{3}}\\
    			& + \frac{1}{2} \, {\left(4 \, C_{2} \overline{A_{0}} + 3 \, \overline{A_{1}} \overline{C_{1}}\right)} \overline{\alpha_{1}} + \frac{1}{8} \, \overline{A_{1}} \overline{\gamma_{2}}\\
    			\lambda_{2}&=4 \, C_{2} {\left| A_{1} \right|}^{2} \overline{A_{0}} - 2 \, \gamma_{1} {\left| A_{1} \right|}^{2} \overline{A_{0}} + 4 \, {\left| A_{1} \right|}^{2} \overline{A_{1}} \overline{C_{1}} + C_{1} \alpha_{4} \overline{A_{0}} + C_{1} \alpha_{1} \overline{A_{1}} + 2 \, \overline{A_{0}} \overline{B_{1}} \overline{\alpha_{1}} + 2 \, \overline{A_{0}} \overline{C_{1}} \overline{\alpha_{4}} - \frac{3}{8} \, B_{1} C_{1}\\
    			& - C_{3} \overline{A_{1}} + \frac{1}{2} \, \gamma_{2} \overline{A_{1}} - 2 \, \overline{A_{2}} \overline{B_{1}}\\
    			\lambda_{3}&=4 \, A_{1} C_{1} {\left| A_{1} \right|}^{2} + 4 \, A_{0} C_{2} {\left| A_{1} \right|}^{2} - 2 \, A_{0} \gamma_{1} {\left| A_{1} \right|}^{2} + 2 \, A_{0} B_{1} \alpha_{1} + 2 \, A_{0} C_{1} \alpha_{4} + A_{1} \overline{C_{1}} \overline{\alpha_{1}} + A_{0} \overline{C_{1}} \overline{\alpha_{4}} - 2 \, A_{2} B_{1} - \frac{3}{8} \, \overline{B_{1}} \overline{C_{1}}\\
    			& - A_{1} \overline{C_{3}} + \frac{1}{2} \, A_{1} \overline{\gamma_{2}}\\
    			\lambda_{4}&=A_{0} {\left| A_{1} \right|}^{2} \overline{B_{1}} + 2 \, A_{1} {\left| A_{1} \right|}^{2} \overline{C_{1}} + \frac{3}{2} \, A_{0} C_{1} \alpha_{3} + A_{0} \alpha_{4} \overline{C_{1}} - \frac{3}{2} \, A_{3} C_{1} - 2 \, A_{2} C_{2} - \frac{1}{2} \, A_{1} C_{3} + \frac{1}{2} \, {\left(3 \, A_{1} C_{1} + 4 \, A_{0} C_{2}\right)} \alpha_{1}\\
    			& - \frac{1}{2} \, {\left(A_{0} \alpha_{1} - A_{2}\right)} \gamma_{1} + \frac{1}{8} \, A_{1} \gamma_{2}
    		\end{align*}
    		\normalsize
    		We see that the new powers are
    		\begin{align*}
    			\Re(z^2),\;\, \quad \Re(z^2)\log|z|,\;\, \Re\left( \frac{\z^3}{z}\right),\;\,\Re\left(\frac{\z^3}{z}\right) \log|z|,\;\,\Re\left(\frac{\z^4}{z^2}\right),\;\, |z|^2,\;\, |z|^2\log|z|
    		\end{align*}
    		so there exists $\vec{C}_4,\vec{B}_3,\vec{E}_2,\vec{\gamma}_3,\vec{\gamma}_4\in \mathbb{C}^n$ and $\vec{E}_4,\vec{\gamma}_5\in \R^n$ such that
    		\begin{align*}
    			\H&=\Re\left(\frac{\vec{C}_1}{z}+\vec{C}_3z+\vec{C}_4+\vec{B}_1\frac{\z}{z}+\vec{B}_2\frac{\z^2}{z}+\vec{B}_3\frac{\z^3}{z}+\vec{E}_1\frac{\z^3}{z^2}+\vec{E}_2\frac{\z^4}{z^2}\right)+\vec{C}_2+\vec{B}_4|z|^2\\
    			&+\vec{\gamma}_1\log|z|+\Re\left(\vec{\gamma}_2z+\vec{\gamma}_3z^2+\frac{\z^3}{z}\vec{\gamma}_4\right)\log|z|+\vec{\gamma}_5|z|^2\log|z|+O(|z|^{2-\epsilon}).
    		\end{align*}
    		which gives, once translated in Sage
    		\begin{align*}
    		\begin{dmatrix}
    			\frac{1}{2} & C_{1} & -1 & 0 & 0 \\
    			1 & C_{2} & 0 & 0 & 0 \\
    			\frac{1}{2} & C_{3} & 1 & 0 & 0 \\
    			\frac{1}{2} & C_{4} & 2 & 0 & 0 \\
    			\frac{1}{2} & B_{1} & -1 & 1 & 0 \\
    			\frac{1}{2} & B_{2} & -1 & 2 & 0 \\
    			\frac{1}{2} & B_{3} & -1 & 3 & 0 \\
    			1 & B_{4} & 1 & 1 & 0 \\
    			\frac{1}{2} & E_{1} & -2 & 3 & 0 \\
    			\frac{1}{2} & E_{2} & -2 & 4 & 0 \\
    			1 & \gamma_{1} & 0 & 0 & 1 \\
    			\frac{1}{2} & \gamma_{2} & 1 & 0 & 1 \\
    			\frac{1}{2} & \gamma_{3} & 2 & 0 & 1 
    			\end{dmatrix}
    			\begin{dmatrix}
    			\frac{1}{2} & \gamma_{4} & -1 & 3 & 0 \\
    			1 & \gamma_{5} & 1 & 1 & 1 \\
    			\frac{1}{2} & \overline{C_{1}} & 0 & -1 & 0 \\
    			\frac{1}{2} & \overline{C_{3}} & 0 & 1 & 0 \\
    			\frac{1}{2} & \overline{C_{4}} & 0 & 2 & 0 \\
    			\frac{1}{2} & \overline{B_{1}} & 1 & -1 & 0 \\
    			\frac{1}{2} & \overline{B_{2}} & 2 & -1 & 0 \\
    			\frac{1}{2} & \overline{B_{3}} & 3 & -1 & 0 \\
    			\frac{1}{2} & \overline{E_{1}} & 3 & -2 & 0 \\
    			\frac{1}{2} & \overline{E_{2}} & 4 & -2 & 0 \\
    			\frac{1}{2} & \overline{\gamma_{2}} & 0 & 1 & 1 \\
    			\frac{1}{2} & \overline{\gamma_{3}} & 0 & 2 & 1 \\
    			\frac{1}{2} & \overline{\gamma_{4}} & 3 & -1 & 0
    			\end{dmatrix}
    		\end{align*}
    		Now, by integrating 
    		\begin{align*}
    			\p{\z}\left(\p{z}\phi\right)=\frac{e^{2\lambda}}{2}\H,
    		\end{align*}
    		we obtain for some $\vec{A}_5\in \mathbb{C}^n$
    		\begin{align*}
    			\p{z}\phi=\begin{dmatrix}
    			1 & A_{0} & 2 & 0 & 0 \\
    			1 & A_{1} & 3 & 0 & 0 \\
    			1 & A_{2} & 4 & 0 & 0 \\
    			1 & A_{3} & 5 & 0 & 0 \\
    			1 & A_{4} & 6 & 0 & 0 \\
    			1 & A_{5} & 7 & 0 & 0 \\
    			\frac{1}{24} & E_{1} & 0 & 6 & 0 \\
    			\frac{1}{28} & E_{2} & 0 & 7 & 0 \\
    			\frac{1}{28} \, \alpha_{2} & C_{1} & 0 & 7 & 0 \\
    			\frac{1}{12} & C_{1} & 1 & 3 & 0 \\
    			\frac{1}{16} & B_{1} & 1 & 4 & 0 \\
    			\frac{1}{20} & B_{2} & 1 & 5 & 0 \\
    			\frac{1}{20} \, \overline{\alpha_{1}} & C_{1} & 1 & 5 & 0 \\
    			\frac{1}{24} & B_{3} & 1 & 6 & 0 \\
    			\frac{1}{24} & \gamma_{4} & 1 & 6 & 0 \\
    			\frac{1}{24} \, \alpha_{2} & \overline{C_{1}} & 1 & 6 & 0 \\
    			\frac{1}{24} \, \overline{\alpha_{1}} & B_{1} & 1 & 6 & 0 \\
    			\frac{1}{24} \, \overline{\alpha_{3}} & C_{1} & 1 & 6 & 0 \\
    			\frac{1}{8} & \overline{C_{1}} & 2 & 2 & 0 \\
    			-\frac{1}{36} & \gamma_{1} & 2 & 3 & 0 \\
    			\frac{1}{6} & C_{2} & 2 & 3 & 0 \\
    			\frac{1}{6} & \gamma_{1} & 2 & 3 & 1 \\
    			-\frac{1}{128} & \overline{\gamma_{2}} & 2 & 4 & 0 \\
    			\frac{1}{16} & \overline{C_{3}} & 2 & 4 & 0 
    			\end{dmatrix}
    			\begin{dmatrix}
    			\frac{1}{8} \, {\left| A_{1} \right|}^{2} & C_{1} & 2 & 4 & 0 \\
    			\frac{1}{16} \, \overline{\alpha_{1}} & \overline{C_{1}} & 2 & 4 & 0 \\
    			\frac{1}{16} & \overline{\gamma_{2}} & 2 & 4 & 1 \\
    			-\frac{1}{200} & \overline{\gamma_{3}} & 2 & 5 & 0 \\
    			\frac{1}{20} & \overline{C_{4}} & 2 & 5 & 0 \\
    			\frac{1}{10} \, {\left| A_{1} \right|}^{2} & B_{1} & 2 & 5 & 0 \\
    			\frac{1}{10} \, \overline{\alpha_{1}} & C_{2} & 2 & 5 & 0 \\
    			-\frac{1}{100} \, \overline{\alpha_{1}} & \gamma_{1} & 2 & 5 & 0 \\
    			\frac{1}{20} \, \overline{\alpha_{3}} & \overline{C_{1}} & 2 & 5 & 0 \\
    			\frac{1}{20} \, \overline{\alpha_{4}} & C_{1} & 2 & 5 & 0 \\
    			\frac{1}{20} & \overline{\gamma_{3}} & 2 & 5 & 1 \\
    			\frac{1}{10} \, \overline{\alpha_{1}} & \gamma_{1} & 2 & 5 & 1 \\
    			\frac{1}{8} & \overline{B_{1}} & 3 & 2 & 0 \\
    			-\frac{1}{72} & \gamma_{2} & 3 & 3 & 0 \\
    			\frac{1}{12} & C_{3} & 3 & 3 & 0 \\
    			\frac{1}{6} \, {\left| A_{1} \right|}^{2} & \overline{C_{1}} & 3 & 3 & 0 \\
    			\frac{1}{12} \, \alpha_{1} & C_{1} & 3 & 3 & 0 \\
    			\frac{1}{12} & \gamma_{2} & 3 & 3 & 1 \\
    			-\frac{1}{64} & \gamma_{5} & 3 & 4 & 0 \\
    			\frac{1}{8} & B_{4} & 3 & 4 & 0 \\
    			\frac{1}{4} \, {\left| A_{1} \right|}^{2} & C_{2} & 3 & 4 & 0 \\
    			-\frac{1}{32} \, {\left| A_{1} \right|}^{2} & \gamma_{1} & 3 & 4 & 0 \\
    			\frac{1}{16} \, \alpha_{1} & B_{1} & 3 & 4 & 0 \\
    			\frac{1}{16} \, \alpha_{4} & C_{1} & 3 & 4 & 0 
    			\end{dmatrix}
    			\begin{dmatrix}
    			\frac{1}{16} \, \overline{\alpha_{1}} & \overline{B_{1}} & 3 & 4 & 0 \\
    			\frac{1}{16} \, \overline{\alpha_{4}} & \overline{C_{1}} & 3 & 4 & 0 \\
    			\frac{1}{8} & \gamma_{5} & 3 & 4 & 1 \\
    			\frac{1}{4} \, {\left| A_{1} \right|}^{2} & \gamma_{1} & 3 & 4 & 1 \\
    			\frac{1}{8} & \overline{B_{2}} & 4 & 2 & 0 \\
    			\frac{1}{8} \, \alpha_{1} & \overline{C_{1}} & 4 & 2 & 0 \\
    			-\frac{1}{72} & \gamma_{3} & 4 & 3 & 0 \\
    			\frac{1}{12} & C_{4} & 4 & 3 & 0 \\
    			\frac{1}{6} \, {\left| A_{1} \right|}^{2} & \overline{B_{1}} & 4 & 3 & 0 \\
    			\frac{1}{6} \, \alpha_{1} & C_{2} & 4 & 3 & 0 \\
    			-\frac{1}{36} \, \alpha_{1} & \gamma_{1} & 4 & 3 & 0 \\
    			\frac{1}{12} \, \alpha_{3} & C_{1} & 4 & 3 & 0 \\
    			\frac{1}{12} \, \alpha_{4} & \overline{C_{1}} & 4 & 3 & 0 \\
    			\frac{1}{12} & \gamma_{3} & 4 & 3 & 1 \\
    			\frac{1}{6} \, \alpha_{1} & \gamma_{1} & 4 & 3 & 1 \\
    			\frac{1}{4} & \overline{E_{1}} & 5 & 1 & 0 \\
    			\frac{1}{8} & \overline{B_{3}} & 5 & 2 & 0 \\
    			\frac{1}{8} & \overline{\gamma_{4}} & 5 & 2 & 0 \\
    			\frac{1}{8} \, \alpha_{1} & \overline{B_{1}} & 5 & 2 & 0 \\
    			\frac{1}{8} \, \alpha_{3} & \overline{C_{1}} & 5 & 2 & 0 \\
    			\frac{1}{8} \, \overline{\alpha_{2}} & C_{1} & 5 & 2 & 0 \\
    			\frac{1}{4} & \overline{E_{2}} & 6 & 1 & 0 \\
    			\frac{1}{4} \, \overline{\alpha_{2}} & \overline{C_{1}} & 6 & 1 & 0
    			\end{dmatrix}
    		\end{align*}
    		Then, we obtain by conformality of the immersion $\phi$
    		\footnotesize
    		\begin{align}\label{3conf}
    			&0=\s{\p{z}\phi}{\p{z}\phi}\\
    			&=\begin{dmatrix}
    				A_{0}^{2} & 4 & 0 & 0 \\
    				2 \, A_{0} A_{1} & 5 & 0 & 0 \\
    				A_{1}^{2} + 2 \, A_{0} A_{2} & 6 & 0 & 0 \\
    				2 \, A_{1} A_{2} + 2 \, A_{0} A_{3} & 7 & 0 & 0 \\
    				A_{2}^{2} + 2 \, A_{1} A_{3} + 2 \, A_{0} A_{4} & 8 & 0 & 0 \\
    				2 \, A_{2} A_{3} + 2 \, A_{1} A_{4} + 2 \, A_{0} A_{5} & 9 & 0 & 0 \\
    				\frac{1}{2} \, A_{0} \overline{E_{1}} & 7 & 1 & 0 \\
    				\frac{1}{2} \, A_{0} \overline{C_{1}} \overline{\alpha_{2}} + \frac{1}{2} \, A_{1} \overline{E_{1}} + \frac{1}{2} \, A_{0} \overline{E_{2}} & 8 & 1 & 0 \\
    				\frac{1}{4} \, A_{0} \overline{C_{1}} & 4 & 2 & 0 \\
    				\frac{1}{4} \, A_{0} \overline{B_{1}} + \frac{1}{4} \, A_{1} \overline{C_{1}} & 5 & 2 & 0 \\
    				\frac{1}{4} \, A_{0} \alpha_{1} \overline{C_{1}} + \frac{1}{4} \, A_{1} \overline{B_{1}} + \frac{1}{4} \, A_{0} \overline{B_{2}} + \frac{1}{4} \, A_{2} \overline{C_{1}} & 6 & 2 & 0 \\
    				\frac{1}{4} \, A_{0} \alpha_{3} \overline{C_{1}} + \frac{1}{4} \, A_{0} C_{1} \overline{\alpha_{2}} + \frac{1}{4} \, {\left(A_{0} \overline{B_{1}} + A_{1} \overline{C_{1}}\right)} \alpha_{1} + \frac{1}{4} \, A_{2} \overline{B_{1}} + \frac{1}{4} \, A_{1} \overline{B_{2}} + \frac{1}{4} \, A_{0} \overline{B_{3}} + \frac{1}{4} \, A_{3} \overline{C_{1}} + \frac{1}{4} \, A_{0} \overline{\gamma_{4}} & 7 & 2 & 0 \\
    				\frac{1}{6} \, A_{0} C_{1} & 3 & 3 & 0 \\
    				\frac{1}{6} \, A_{1} C_{1} + \frac{1}{3} \, A_{0} C_{2} - \frac{1}{18} \, A_{0} \gamma_{1} & 4 & 3 & 0 \\
    				\frac{1}{3} \, A_{0} {\left| A_{1} \right|}^{2} \overline{C_{1}} + \frac{1}{6} \, A_{0} C_{1} \alpha_{1} + \frac{1}{6} \, A_{2} C_{1} + \frac{1}{3} \, A_{1} C_{2} + \frac{1}{6} \, A_{0} C_{3} - \frac{1}{18} \, A_{1} \gamma_{1} - \frac{1}{36} \, A_{0} \gamma_{2} & 5 & 3 & 0 \\
    				\lambda_1 & 6 & 3 & 0 \\
    				\frac{1}{3} \, A_{0} \gamma_{1} & 4 & 3 & 1 \\
    				\frac{1}{3} \, A_{1} \gamma_{1} + \frac{1}{6} \, A_{0} \gamma_{2} & 5 & 3 & 1 \\
    				\frac{1}{3} \, {\left(A_{0} \alpha_{1} + A_{2}\right)} \gamma_{1} + \frac{1}{6} \, A_{1} \gamma_{2} + \frac{1}{6} \, A_{0} \gamma_{3} & 6 & 3 & 1 \\
    				\frac{1}{8} \, A_{0} B_{1} & 3 & 4 & 0 \\
    				\frac{1}{4} \, A_{0} C_{1} {\left| A_{1} \right|}^{2} + \frac{1}{8} \, A_{0} \overline{C_{1}} \overline{\alpha_{1}} + \frac{1}{8} \, A_{1} B_{1} + \frac{1}{64} \, \overline{C_{1}}^{2} + \frac{1}{8} \, A_{0} \overline{C_{3}} - \frac{1}{64} \, A_{0} \overline{\gamma_{2}} & 4 & 4 & 0 \\
    				\lambda_2 & 5 & 4 & 0 \\
    				\frac{1}{8} \, A_{0} \overline{\gamma_{2}} & 4 & 4 & 1 \\
    				\frac{1}{2} \, A_{0} \gamma_{1} {\left| A_{1} \right|}^{2} + \frac{1}{4} \, A_{0} \gamma_{5} + \frac{1}{8} \, A_{1} \overline{\gamma_{2}} & 5 & 4 & 1 \\
    				\frac{1}{10} \, A_{0} C_{1} \overline{\alpha_{1}} + \frac{1}{10} \, A_{0} B_{2} + \frac{1}{48} \, C_{1} \overline{C_{1}} & 3 & 5 & 0 \\
    				\lambda_3 & 4 & 5 & 0 \\
    				\frac{1}{120} \, {\left(24 \, A_{0} \overline{\alpha_{1}} + 5 \, \overline{C_{1}}\right)} \gamma_{1} + \frac{1}{10} \, A_{0} \overline{\gamma_{3}} & 4 & 5 & 1 \\
    				\frac{1}{144} \, C_{1}^{2} + \frac{1}{12} \, A_{0} E_{1} & 2 & 6 & 0 \\
    				\frac{1}{12} \, A_{0} \alpha_{2} \overline{C_{1}} + \frac{1}{12} \, A_{0} B_{1} \overline{\alpha_{1}} + \frac{1}{12} \, A_{0} C_{1} \overline{\alpha_{3}} + \frac{1}{12} \, A_{0} B_{3} + \frac{1}{36} \, C_{1} C_{2} + \frac{1}{12} \, A_{1} E_{1} - \frac{1}{216} \, C_{1} \gamma_{1} + \frac{1}{12} \, A_{0} \gamma_{4} + \frac{1}{64} \, B_{1} \overline{C_{1}} & 3 & 6 & 0 \\
    				\frac{1}{36} \, C_{1} \gamma_{1} & 3 & 6 & 1 \\
    				\frac{1}{14} \, A_{0} C_{1} \alpha_{2} + \frac{1}{96} \, B_{1} C_{1} + \frac{1}{14} \, A_{0} E_{2} & 2 & 7 & 0
    			\end{dmatrix}
    		\end{align}  
    		\normalsize
    		where
    		\small
    		\begin{align*}
    			\lambda_{1}&=\frac{1}{3} \, A_{0} {\left| A_{1} \right|}^{2} \overline{B_{1}} + \frac{1}{3} \, A_{1} {\left| A_{1} \right|}^{2} \overline{C_{1}} + \frac{1}{6} \, A_{0} C_{1} \alpha_{3} + \frac{1}{6} \, A_{0} \alpha_{4} \overline{C_{1}} + \frac{1}{6} \, A_{3} C_{1} + \frac{1}{3} \, A_{2} C_{2} + \frac{1}{6} \, A_{1} C_{3} + \frac{1}{6} \, A_{0} C_{4}\\
    			& + \frac{1}{6} \, {\left(A_{1} C_{1} + 2 \, A_{0} C_{2}\right)} \alpha_{1} - \frac{1}{18} \, {\left(A_{0} \alpha_{1} + A_{2}\right)} \gamma_{1} - \frac{1}{36} \, A_{1} \gamma_{2} - \frac{1}{36} \, A_{0} \gamma_{3}\\
    			\lambda_{2}&=\frac{1}{4} \, A_{1} C_{1} {\left| A_{1} \right|}^{2} + \frac{1}{2} \, A_{0} C_{2} {\left| A_{1} \right|}^{2} - \frac{1}{16} \, A_{0} \gamma_{1} {\left| A_{1} \right|}^{2} + \frac{1}{8} \, A_{0} B_{1} \alpha_{1} + \frac{1}{8} \, A_{0} C_{1} \alpha_{4} + \frac{1}{8} \, A_{0} \overline{C_{1}} \overline{\alpha_{4}} + \frac{1}{8} \, A_{2} B_{1} + \frac{1}{4} \, A_{0} B_{4}\\
    			& - \frac{1}{32} \, A_{0} \gamma_{5} + \frac{1}{32} \, \overline{B_{1}} \overline{C_{1}} + \frac{1}{8} \, A_{1} \overline{C_{3}} + \frac{1}{8} \, {\left(A_{0} \overline{B_{1}} + A_{1} \overline{C_{1}}\right)} \overline{\alpha_{1}} - \frac{1}{64} \, A_{1} \overline{\gamma_{2}}\\
    			\lambda_{3}&=\frac{1}{5} \, A_{0} B_{1} {\left| A_{1} \right|}^{2} + \frac{1}{10} \, A_{0} \overline{C_{1}} \overline{\alpha_{3}} + \frac{1}{10} \, A_{0} C_{1} \overline{\alpha_{4}} + \frac{1}{10} \, A_{1} B_{2} - \frac{1}{3600} \, {\left(72 \, A_{0} \overline{\alpha_{1}} + 25 \, \overline{C_{1}}\right)} \gamma_{1} + \frac{1}{48} \, C_{1} \overline{B_{1}} + \frac{1}{24} \, C_{2} \overline{C_{1}} + \frac{1}{10} \, A_{0} \overline{C_{4}}\\
    			& + \frac{1}{10} \, {\left(A_{1} C_{1} + 2 \, A_{0} C_{2}\right)} \overline{\alpha_{1}} - \frac{1}{100} \, A_{0} \overline{\gamma_{3}}
    		\end{align*}
    		\normalsize
    	while  		
    	\footnotesize
    	\begin{align*}
    		e^{2\lambda}=\begin{dmatrix}
    		\frac{1}{12} \, A_{0} \overline{E_{1}} & 8 & 0 & 0 \\
    		\frac{1}{14} \ccancel{\, A_{0} \overline{C_{1}} \overline{\alpha_{2}}} + \frac{1}{12} \ccancel{\, A_{1} \overline{E_{1}}} + \frac{1}{14} \, \ccancel{A_{0} \overline{E_{2}}} & 9 & 0 & 0 \\
    		\frac{1}{6} \, A_{0} \overline{C_{1}} & 5 & 1 & 0 \\
    		\frac{1}{8} \, A_{0} \overline{B_{1}} + \frac{1}{6} \, A_{1} \overline{C_{1}} & 6 & 1 & 0 \\
    		\frac{1}{10} \, A_{0} \alpha_{1} \overline{C_{1}} + \frac{1}{8} \, A_{1} \overline{B_{1}} + \frac{1}{10} \, A_{0} \overline{B_{2}} + \frac{1}{6} \, A_{2} \overline{C_{1}} & 7 & 1 & 0 \\
    		\frac{1}{12} \, A_{0} \alpha_{3} \overline{C_{1}} + \frac{1}{12} \, A_{0} C_{1} \overline{\alpha_{2}} + \frac{1}{60} \, {\left(5 \, A_{0} \overline{B_{1}} + 6 \, A_{1} \overline{C_{1}}\right)} \alpha_{1} + \frac{1}{8} \, A_{2} \overline{B_{1}} + \frac{1}{10} \, A_{1} \overline{B_{2}} + \frac{1}{12} \, A_{0} \overline{B_{3}} + \frac{1}{6} \, A_{3} \overline{C_{1}} + \frac{1}{12} \, A_{0} \overline{\gamma_{4}} & 8 & 1 & 0 \\
    		2 \, A_{0} \overline{A_{0}} & 2 & 2 & 0 \\
    		2 \, A_{1} \overline{A_{0}} & 3 & 2 & 0 \\
    		\frac{1}{4} \, A_{0} C_{1} + 2 \, A_{2} \overline{A_{0}} & 4 & 2 & 0 \\
    		\frac{1}{4} \, A_{1} C_{1} + \frac{1}{3} \, A_{0} C_{2} - \frac{1}{18} \, A_{0} \gamma_{1} + 2 \, A_{3} \overline{A_{0}} & 5 & 2 & 0 \\
    		\frac{1}{4} \, A_{0} {\left| A_{1} \right|}^{2} \overline{C_{1}} + \frac{1}{8} \, A_{0} C_{1} \alpha_{1} + \frac{1}{4} \, A_{2} C_{1} + \frac{1}{3} \, A_{1} C_{2} + \frac{1}{8} \, A_{0} C_{3} - \frac{1}{18} \, A_{1} \gamma_{1} - \frac{1}{64} \, A_{0} \gamma_{2} + 2 \, A_{4} \overline{A_{0}} & 6 & 2 & 0 \\
    		\mu_1 & 7 & 2 & 0 \\
    		\frac{1}{3} \, A_{0} \gamma_{1} & 5 & 2 & 1 \\
    		\frac{1}{3} \, A_{1} \gamma_{1} + \frac{1}{8} \, A_{0} \gamma_{2} & 6 & 2 & 1 \\
    		\frac{1}{15} \, {\left(3 \, A_{0} \alpha_{1} + 5 \, A_{2}\right)} \gamma_{1} + \frac{1}{8} \, A_{1} \gamma_{2} + \frac{1}{10} \, A_{0} \gamma_{3} & 7 & 2 & 1 \\
    		2 \, A_{0} \overline{A_{1}} & 2 & 3 & 0 \\
    		2 \, A_{1} \overline{A_{1}} & 3 & 3 & 0 \\
    		\frac{1}{4} \, A_{0} B_{1} + 2 \, A_{2} \overline{A_{1}} & 4 & 3 & 0 \\
    		\frac{1}{3} \, A_{0} C_{1} {\left| A_{1} \right|}^{2} + \frac{1}{6} \, A_{0} \overline{C_{1}} \overline{\alpha_{1}} + \frac{1}{4} \, A_{1} B_{1} + 2 \, A_{3} \overline{A_{1}} + \frac{1}{48} \, \overline{C_{1}}^{2} + \frac{1}{6} \, A_{0} \overline{C_{3}} + \frac{1}{2} \, \overline{A_{0}} \overline{E_{1}} - \frac{1}{36} \, A_{0} \overline{\gamma_{2}} & 5 & 3 & 0 \\
    		\mu_2 & 6 & 3 & 0 \\
    		\frac{1}{6} \, A_{0} \overline{\gamma_{2}} & 5 & 3 & 1 \\
    		\frac{1}{2} \, A_{0} \gamma_{1} {\left| A_{1} \right|}^{2} + \frac{1}{4} \, A_{0} \gamma_{5} + \frac{1}{6} \, A_{1} \overline{\gamma_{2}} & 6 & 3 & 1 \\
    		2 \, A_{0} \overline{A_{2}} + \frac{1}{4} \, \overline{A_{0}} \overline{C_{1}} & 2 & 4 & 0 \\
    		2 \, A_{1} \overline{A_{2}} + \frac{1}{4} \, \overline{A_{0}} \overline{B_{1}} & 3 & 4 & 0 \\
    		\frac{1}{4} \, \alpha_{1} \overline{A_{0}} \overline{C_{1}} + \frac{1}{4} \, A_{0} C_{1} \overline{\alpha_{1}} + \frac{1}{4} \, A_{0} B_{2} + 2 \, A_{2} \overline{A_{2}} + \frac{1}{4} \, \overline{A_{0}} \overline{B_{2}} + \frac{13}{288} \, C_{1} \overline{C_{1}} & 4 & 4 & 0 \\
    		\mu_3 & 5 & 4 & 0 \\
    		\frac{1}{72} \, {\left(24 \, A_{0} \overline{\alpha_{1}} + 5 \, \overline{C_{1}}\right)} \gamma_{1} + \frac{1}{6} \, A_{0} \overline{\gamma_{3}} & 5 & 4 & 1 \\
    		\frac{1}{6} \, C_{1} \overline{A_{0}} & 1 & 5 & 0 \\
    		\frac{1}{3} \, C_{2} \overline{A_{0}} - \frac{1}{18} \, \gamma_{1} \overline{A_{0}} + 2 \, A_{0} \overline{A_{3}} + \frac{1}{4} \, \overline{A_{1}} \overline{C_{1}} & 2 & 5 & 0 \\
    		\frac{1}{3} \, {\left| A_{1} \right|}^{2} \overline{A_{0}} \overline{C_{1}} + \frac{1}{6} \, C_{1} \alpha_{1} \overline{A_{0}} + \frac{1}{48} \, C_{1}^{2} + \frac{1}{2} \, A_{0} E_{1} + \frac{1}{6} \, C_{3} \overline{A_{0}} - \frac{1}{36} \, \gamma_{2} \overline{A_{0}} + 2 \, A_{1} \overline{A_{3}} + \frac{1}{4} \, \overline{A_{1}} \overline{B_{1}} & 3 & 5 & 0 \\
    		\mu_4 & 4 & 5 & 0 \\
    		\frac{1}{3} \, \gamma_{1} \overline{A_{0}} & 2 & 5 & 1 \\
    		\frac{1}{6} \, \gamma_{2} \overline{A_{0}} & 3 & 5 & 1 \\
    		\frac{1}{72} \, {\left(24 \, \alpha_{1} \overline{A_{0}} + 5 \, C_{1}\right)} \gamma_{1} + \frac{1}{6} \, \gamma_{3} \overline{A_{0}} & 4 & 5 & 1 
    		\end{dmatrix}
    		\end{align*}
    		\begin{align*}
    		\begin{dmatrix}
    		\frac{1}{8} \, B_{1} \overline{A_{0}} + \frac{1}{6} \, C_{1} \overline{A_{1}} & 1 & 6 & 0 \\
    		\frac{1}{4} \, C_{1} {\left| A_{1} \right|}^{2} \overline{A_{0}} + \frac{1}{8} \, \overline{A_{0}} \overline{C_{1}} \overline{\alpha_{1}} + \frac{1}{3} \, C_{2} \overline{A_{1}} - \frac{1}{18} \, \gamma_{1} \overline{A_{1}} + 2 \, A_{0} \overline{A_{4}} + \frac{1}{4} \, \overline{A_{2}} \overline{C_{1}} + \frac{1}{8} \, \overline{A_{0}} \overline{C_{3}} - \frac{1}{64} \, \overline{A_{0}} \overline{\gamma_{2}} & 2 & 6 & 0 \\
    	    \mu_5 & 3 & 6 & 0 \\
    		\frac{1}{3} \, \gamma_{1} \overline{A_{1}} + \frac{1}{8} \, \overline{A_{0}} \overline{\gamma_{2}} & 2 & 6 & 1 \\
    		\frac{1}{2} \, \gamma_{1} {\left| A_{1} \right|}^{2} \overline{A_{0}} + \frac{1}{4} \, \gamma_{5} \overline{A_{0}} + \frac{1}{6} \, \gamma_{2} \overline{A_{1}} & 3 & 6 & 1 \\
    		\frac{1}{10} \, C_{1} \overline{A_{0}} \overline{\alpha_{1}} + \frac{1}{10} \, B_{2} \overline{A_{0}} + \frac{1}{8} \, B_{1} \overline{A_{1}} + \frac{1}{6} \, C_{1} \overline{A_{2}} & 1 & 7 & 0 \\
    	    \mu_6 & 2 & 7 & 0 \\
    		\frac{1}{15} \, {\left(3 \, \overline{A_{0}} \overline{\alpha_{1}} + 5 \, \overline{A_{2}}\right)} \gamma_{1} + \frac{1}{8} \, \overline{A_{1}} \overline{\gamma_{2}} + \frac{1}{10} \, \overline{A_{0}} \overline{\gamma_{3}} & 2 & 7 & 1 \\
    		\frac{1}{12} \, E_{1} \overline{A_{0}} & 0 & 8 & 0 \\
    		\frac{1}{12} \, \alpha_{2} \overline{A_{0}} \overline{C_{1}} + \frac{1}{12} \, C_{1} \overline{A_{0}} \overline{\alpha_{3}} + \frac{1}{12} \, B_{3} \overline{A_{0}} + \frac{1}{12} \, \gamma_{4} \overline{A_{0}} + \frac{1}{10} \, B_{2} \overline{A_{1}} + \frac{1}{8} \, B_{1} \overline{A_{2}} + \frac{1}{6} \, C_{1} \overline{A_{3}} + \frac{1}{60} \, {\left(5 \, B_{1} \overline{A_{0}} + 6 \, C_{1} \overline{A_{1}}\right)} \overline{\alpha_{1}} & 1 & 8 & 0 \\
    		\frac{1}{14} \ccancel{\, C_{1} \alpha_{2} \overline{A_{0}}} + \frac{1}{14} \ccancel{\, E_{2} \overline{A_{0}}} + \frac{1}{12} \ccancel{\, E_{1} \overline{A_{1}}} & 0 & 9 & 0
    		\end{dmatrix}
    	\end{align*}
    \normalsize
    where
    \small
    \begin{align*}
    	\mu_{1}&=\frac{1}{5} \, A_{0} {\left| A_{1} \right|}^{2} \overline{B_{1}} + \frac{1}{4} \, A_{1} {\left| A_{1} \right|}^{2} \overline{C_{1}} + \frac{1}{10} \, A_{0} C_{1} \alpha_{3} + \frac{1}{10} \, A_{0} \alpha_{4} \overline{C_{1}} + \frac{1}{4} \, A_{3} C_{1} + \frac{1}{3} \, A_{2} C_{2} + \frac{1}{8} \, A_{1} C_{3} + \frac{1}{10} \, A_{0} C_{4}\\
    	& + \frac{1}{40} \, {\left(5 \, A_{1} C_{1} + 8 \, A_{0} C_{2}\right)} \alpha_{1} - \frac{1}{450} \, {\left(9 \, A_{0} \alpha_{1} + 25 \, A_{2}\right)} \gamma_{1} - \frac{1}{64} \, A_{1} \gamma_{2} - \frac{1}{100} \, A_{0} \gamma_{3} + 2 \, A_{5} \overline{A_{0}}\\
    	\mu_{2}&=\frac{1}{3} \, A_{1} C_{1} {\left| A_{1} \right|}^{2} + \frac{1}{2} \, A_{0} C_{2} {\left| A_{1} \right|}^{2} - \frac{1}{16} \, A_{0} \gamma_{1} {\left| A_{1} \right|}^{2} + \frac{1}{8} \, A_{0} B_{1} \alpha_{1} + \frac{1}{8} \, A_{0} C_{1} \alpha_{4} + \frac{1}{2} \, \overline{A_{0}} \overline{C_{1}} \overline{\alpha_{2}} + \frac{1}{8} \, A_{0} \overline{C_{1}} \overline{\alpha_{4}} + \frac{1}{4} \, A_{2} B_{1}\\
    	& + \frac{1}{4} \, A_{0} B_{4} - \frac{1}{32} \, A_{0} \gamma_{5} + 2 \, A_{4} \overline{A_{1}} + \frac{7}{192} \, \overline{B_{1}} \overline{C_{1}} + \frac{1}{6} \, A_{1} \overline{C_{3}} + \frac{1}{2} \, \overline{A_{0}} \overline{E_{2}} + \frac{1}{24} \, {\left(3 \, A_{0} \overline{B_{1}} + 4 \, A_{1} \overline{C_{1}}\right)} \overline{\alpha_{1}} - \frac{1}{36} \, A_{1} \overline{\gamma_{2}}\\
    	\mu_{3}&=\frac{1}{3} \, A_{0} B_{1} {\left| A_{1} \right|}^{2} + \frac{1}{4} \, \alpha_{1} \overline{A_{0}} \overline{B_{1}} + \frac{1}{4} \, \alpha_{3} \overline{A_{0}} \overline{C_{1}} + \frac{1}{4} \, C_{1} \overline{A_{0}} \overline{\alpha_{2}} + \frac{1}{6} \, A_{0} \overline{C_{1}} \overline{\alpha_{3}} + \frac{1}{6} \, A_{0} C_{1} \overline{\alpha_{4}} + \frac{1}{4} \, A_{1} B_{2} - \frac{1}{432} \, {\left(24 \, A_{0} \overline{\alpha_{1}} + 5 \, \overline{C_{1}}\right)} \gamma_{1}\\
    	& + 2 \, A_{3} \overline{A_{2}} + \frac{1}{24} \, C_{1} \overline{B_{1}} + \frac{1}{4} \, \overline{A_{0}} \overline{B_{3}} + \frac{5}{72} \, C_{2} \overline{C_{1}} + \frac{1}{6} \, A_{0} \overline{C_{4}} + \frac{1}{2} \, \overline{A_{1}} \overline{E_{1}} + \frac{1}{12} \, {\left(3 \, A_{1} C_{1} + 4 \, A_{0} C_{2}\right)} \overline{\alpha_{1}} - \frac{1}{36} \, A_{0} \overline{\gamma_{3}} + \frac{1}{4} \, \overline{A_{0}} \overline{\gamma_{4}}\\
    	\mu_{4}&=\frac{1}{3} \, {\left| A_{1} \right|}^{2} \overline{A_{0}} \overline{B_{1}} + \frac{1}{6} \, C_{1} \alpha_{3} \overline{A_{0}} + \frac{1}{4} \, A_{0} \alpha_{2} \overline{C_{1}} + \frac{1}{6} \, \alpha_{4} \overline{A_{0}} \overline{C_{1}} + \frac{1}{4} \, A_{0} B_{1} \overline{\alpha_{1}} + \frac{1}{4} \, A_{0} C_{1} \overline{\alpha_{3}} + \frac{1}{4} \, A_{0} B_{3} + \frac{5}{72} \, C_{1} C_{2} + \frac{1}{2} \, A_{1} E_{1}\\
    	& + \frac{1}{12} \, {\left(4 \, C_{2} \overline{A_{0}} + 3 \, \overline{A_{1}} \overline{C_{1}}\right)} \alpha_{1} - \frac{1}{432} \, {\left(24 \, \alpha_{1} \overline{A_{0}} + 5 \, C_{1}\right)} \gamma_{1} + \frac{1}{4} \, A_{0} \gamma_{4} + \frac{1}{6} \, C_{4} \overline{A_{0}} - \frac{1}{36} \, \gamma_{3} \overline{A_{0}} + 2 \, A_{2} \overline{A_{3}} + \frac{1}{4} \, \overline{A_{1}} \overline{B_{2}} + \frac{1}{24} \, B_{1} \overline{C_{1}}\\
    	\mu_{5}&=\frac{1}{2} \, C_{2} {\left| A_{1} \right|}^{2} \overline{A_{0}} - \frac{1}{16} \, \gamma_{1} {\left| A_{1} \right|}^{2} \overline{A_{0}} + \frac{1}{3} \, {\left| A_{1} \right|}^{2} \overline{A_{1}} \overline{C_{1}} + \frac{1}{2} \, A_{0} C_{1} \alpha_{2} + \frac{1}{8} \, C_{1} \alpha_{4} \overline{A_{0}} + \frac{1}{8} \, \overline{A_{0}} \overline{B_{1}} \overline{\alpha_{1}} + \frac{1}{8} \, \overline{A_{0}} \overline{C_{1}} \overline{\alpha_{4}} + \frac{7}{192} \, B_{1} C_{1}\\
    	& + \frac{1}{2} \, A_{0} E_{2} + \frac{1}{24} \, {\left(3 \, B_{1} \overline{A_{0}} + 4 \, C_{1} \overline{A_{1}}\right)} \alpha_{1} + \frac{1}{4} \, B_{4} \overline{A_{0}} - \frac{1}{32} \, \gamma_{5} \overline{A_{0}} + \frac{1}{6} \, C_{3} \overline{A_{1}} - \frac{1}{36} \, \gamma_{2} \overline{A_{1}} + 2 \, A_{1} \overline{A_{4}} + \frac{1}{4} \, \overline{A_{2}} \overline{B_{1}}\\
    	\mu_{6}&=\frac{1}{5} \, B_{1} {\left| A_{1} \right|}^{2} \overline{A_{0}} + \frac{1}{4} \, C_{1} {\left| A_{1} \right|}^{2} \overline{A_{1}} + \frac{1}{10} \, \overline{A_{0}} \overline{C_{1}} \overline{\alpha_{3}} + \frac{1}{10} \, C_{1} \overline{A_{0}} \overline{\alpha_{4}} - \frac{1}{450} \, {\left(9 \, \overline{A_{0}} \overline{\alpha_{1}} + 25 \, \overline{A_{2}}\right)} \gamma_{1} + \frac{1}{3} \, C_{2} \overline{A_{2}} + 2 \, A_{0} \overline{A_{5}}\\
    	& + \frac{1}{4} \, \overline{A_{3}} \overline{C_{1}} + \frac{1}{8} \, \overline{A_{1}} \overline{C_{3}} + \frac{1}{10} \, \overline{A_{0}} \overline{C_{4}} + \frac{1}{40} \, {\left(8 \, C_{2} \overline{A_{0}} + 5 \, \overline{A_{1}} \overline{C_{1}}\right)} \overline{\alpha_{1}} - \frac{1}{64} \, \overline{A_{1}} \overline{\gamma_{2}} - \frac{1}{100} \, \overline{A_{0}} \overline{\gamma_{3}}
    \end{align*}
    \normalsize
    as
    \begin{align*}
    	\frac{1}{2}\vec{E}_2=\begin{dmatrix}
    	-\frac{7}{96} \, C_{1}^{2} & \overline{A_{1}} & -2 & 4 & 0& \textbf{(2)}\\
    	\frac{1}{2} \ccancel{\, A_{0} C_{1} \alpha_{2}} - \frac{13}{96} \ccancel{\, B_{1} C_{1}} & \overline{A_{0}} & -2 & 4 & 0 & \textbf{(3)}\\
    	-\frac{1}{48} \, C_{1} \overline{A_{1}} & C_{1} & -2 & 4 & 0 & \textbf{(4)}\\
    	\ccancel{C_{1} \alpha_{2} \overline{A_{0}}} - \frac{1}{4} \ccancel{\, E_{1} \overline{A_{1}}} & A_{0} & -2 & 4 & 0 & \textbf{(5)}\\
    	\end{dmatrix}=-\frac{7}{96}\s{\vec{C}_1}{\vec{C}_1}\bar{\vec{A}_1}-\frac{1}{48}\s{\bar{\vec{A}_1}}{\vec{C}_1}\vec{C}_1
    \end{align*}
    so
    \begin{align}\label{3defe2}
    	\vec{E}_2=-\frac{7}{48}\s{\vec{C}_1}{\vec{C}_1}\bar{\vec{A}_1}-\frac{1}{24}\s{\bar{\vec{A}_1}}{\vec{C}_1}\vec{C}_1
    \end{align}
    and
    \begin{align*}
    	\s{\vec{A}_0}{\vec{E}_2}=\s{\bar{\vec{A}_0}}{\vec{E}_2}=0.
    \end{align*}
    Therefore, recalling that
    \begin{align*}
    	e^{2\lambda}=\begin{dmatrix}
    	1 & 2 & 2 & 0 \\
    	2 \, {\left| A_{1} \right|}^{2} & 3 & 3 & 0 \\
    	\beta & 4 & 4 & 0 \\
    	\alpha_{1} & 4 & 2 & 0 \\
    	\alpha_{2} & 1 & 6 & 0 \\
    	\alpha_{3} & 5 & 2 & 0 \\
    	\alpha_{4} & 4 & 3 & 0 \\
    	\alpha_{5} & 7 & 1 & 0\\
    	\alpha_{6} & 6 & 2 & 0\\
    	\alpha_{7} & 3 & 5 & 0\\\
    	\zeta_{0} & 6 & 2 & 1 \\
    	\end{dmatrix}
    	\begin{dmatrix}
    	\overline{\alpha_{1}} & 2 & 4 & 0 \\
    	\overline{\alpha_{2}} & 6 & 1 & 0 \\
    	\overline{\alpha_{3}} & 2 & 5 & 0 \\
    	\overline{\alpha_{4}} & 3 & 4 & 0 \\
    	\bar{\alpha_{5}} & 1 & 7 & 0\\
    	\bar{\alpha_{6}} & 2 & 6 & 0\\
    	\bar{\alpha_{7}} & 5 & 3 & 0\\
    	\bar{\zeta_{0}} & 2 & 6 & 1 \\
    	\end{dmatrix}
    \end{align*}
    we see that the only new powers of degree $9$ are
    \begin{align*}
    	\Re(z^8\z),\;\, \Re(z^7\z^2),\;\, \Re(z^7\z^2)\log|z|,\;\, \Re(z^6\z^3),\;\, \Re(z^6\z^3)\log|z|,\;\, \Re(z^5\z^4),\;\, \Re(z^5\z^4)\log|z|,\;\, 
    \end{align*}
    so there exists $\alpha_8,\alpha_9,\alpha_{10},\alpha_{11},\zeta_1,\zeta_2,\zeta_3\in \mathbb{C}$ such that
    \begin{align*}
    	e^{2\lambda}=\begin{dmatrix}
    	1 & 2 & 2 & 0 \\
    	2 \, {\left| A_{1} \right|}^{2} & 3 & 3 & 0 \\
    	\beta & 4 & 4 & 0 \\
    	\alpha_{1} & 4 & 2 & 0 \\
    	\alpha_{2} & 1 & 6 & 0 \\
    	\alpha_{3} & 5 & 2 & 0 \\
    	\alpha_{4} & 4 & 3 & 0 \\
    	\alpha_{5} & 7 & 1 & 0\\
    	\alpha_{6} & 6 & 2 & 0\\
    	\alpha_{7} & 3 & 5 & 0\\
    	\end{dmatrix}
    	\begin{dmatrix}
    	\alpha_{8} & 8 & 1 & 0 \\
    	\alpha_{9} & 7 & 2 & 0\\
    	\alpha_{10} & 6 & 3 & 0\\
    	\alpha_{11} & 5 & 4 & 0\\
    	\zeta_{0} & 6 & 2 & 1 \\
    	\zeta_{1} & 7 & 2 & 1\\
    	\zeta_{2} & 6 & 3 & 1\\
    	\zeta_{3} & 5 & 4 & 1
    	\end{dmatrix}
    	\begin{dmatrix}
    	\overline{\alpha_{1}} & 2 & 4 & 0 \\
    	\overline{\alpha_{2}} & 6 & 1 & 0 \\
    	\overline{\alpha_{3}} & 2 & 5 & 0 \\
    	\overline{\alpha_{4}} & 3 & 4 & 0 \\
    	\bar{\alpha_{5}} & 1 & 7 & 0\\
    	\bar{\alpha_{6}} & 2 & 6 & 0\\
    	\bar{\alpha_{7}} & 5 & 3 & 0
    	\end{dmatrix}
    	\begin{dmatrix}
    	\bar{\alpha_{8}} & 1 & 8 & 0 \\
    	\bar{\alpha_{9}} & 2 & 7 & 0\\
    	\bar{\alpha_{10}} & 3 & 6 & 0\\
    	\bar{\alpha_{11}} & 4 & 5 & 0\\
    	\bar{\zeta_{0}} & 2 & 6 & 1 \\
    	\bar{\zeta_{1}} & 2 & 7 & 1\\
    	\bar{\zeta_{2}} & 3 & 6 & 1\\
    	\bar{\zeta_{3}} & 4 & 5 & 1
    	\end{dmatrix}.
    \end{align*}
    Let us compare this to the \TeX\: version (on the left, while Sage is on the right)
    \begin{align*}
    	e^{2\lambda}=\begin{dmatrix}
    	1 & 2 & 2 & 0 \\
    	2 \, {\left| A_{1} \right|}^{2} & 3 & 3 & 0 \\
    	\beta & 4 & 4 & 0 \\
    	\alpha_{1} & 4 & 2 & 0 \\
    	\alpha_{2} & 1 & 6 & 0 \\
    	\alpha_{3} & 5 & 2 & 0 \\
    	\alpha_{4} & 4 & 3 & 0 \\
    	\alpha_{5} & 7 & 1 & 0\\
    	\alpha_{6} & 6 & 2 & 0\\
    	\alpha_{7} & 3 & 5 & 0\\
    	\alpha_{8} & 8 & 1 & 0 \\
    	\alpha_{9} & 7 & 2 & 0\\
    	\alpha_{10} & 6 & 3 & 0\\
    	\alpha_{11} & 5 & 4 & 0\\
    	\zeta_{0} & 6 & 2 & 1 \\
    	\zeta_{1} & 7 & 2 & 1\\
    	\zeta_{2} & 6 & 3 & 1\\
    	\zeta_{3} & 5 & 4 & 1
    	\end{dmatrix}
    	\begin{dmatrix}
    	\overline{\alpha_{1}} & 2 & 4 & 0 \\
    	\overline{\alpha_{2}} & 6 & 1 & 0 \\
    	\overline{\alpha_{3}} & 2 & 5 & 0 \\
    	\overline{\alpha_{4}} & 3 & 4 & 0 \\
    	\bar{\alpha_{5}} & 1 & 7 & 0\\
    	\bar{\alpha_{6}} & 2 & 6 & 0\\
    	\bar{\alpha_{7}} & 5 & 3 & 0 \\
    	\bar{\alpha_{8}} & 1 & 8 & 0 \\
    	\bar{\alpha_{9}} & 2 & 7 & 0\\
    	\bar{\alpha_{10}} & 3 & 6 & 0\\
    	\bar{\alpha_{11}} & 4 & 5 & 0\\
    	\bar{\zeta_{0}} & 2 & 6 & 1 \\
    	\bar{\zeta_{1}} & 2 & 7 & 1\\
    	\bar{\zeta_{2}} & 3 & 6 & 1\\
    	\bar{\zeta_{3}} & 4 & 5 & 1
    	\end{dmatrix} =\begin{dmatrix}
    	1 & 2 & 2 & 0 \\
    	2 \, {\left| A_{1} \right|}^{2} & 3 & 3 & 0 \\
    	\beta & 4 & 4 & 0 \\
    	\alpha_{1} & 4 & 2 & 0 \\
    	\alpha_{2} & 1 & 6 & 0 \\
    	\alpha_{3} & 5 & 2 & 0 \\
    	\alpha_{4} & 4 & 3 & 0 \\
    	\alpha_{5} & 7 & 1 & 0 \\
    	\alpha_{6} & 6 & 2 & 0 \\
    	\alpha_{7} & 3 & 5 & 0 \\
    	\alpha_{8} & 8 & 1 & 0 \\
    	\alpha_{9} & 7 & 2 & 0 \\
    	\alpha_{10} & 6 & 3 & 0 \\
    	\alpha_{11} & 5 & 4 & 0 \\
    	\zeta_{0} & 6 & 2 & 1 \\
    	\zeta_{1} & 7 & 2 & 1 \\
    	\zeta_{2} & 6 & 3 & 1 \\
    	\zeta_{3} & 5 & 4 & 1 
    	\end{dmatrix}
    	\begin{dmatrix}
    	\overline{\alpha_{1}} & 2 & 4 & 0 \\
    	\overline{\alpha_{2}} & 6 & 1 & 0 \\
    	\overline{\alpha_{3}} & 2 & 5 & 0 \\
    	\overline{\alpha_{4}} & 3 & 4 & 0 \\
    	\overline{\alpha_{5}} & 1 & 7 & 0 \\
    	\overline{\alpha_{6}} & 2 & 6 & 0 \\
    	\overline{\alpha_{7}} & 5 & 3 & 0 \\
    	\overline{\alpha_{8}} & 1 & 8 & 0 \\
    	\overline{\alpha_{9}} & 2 & 7 & 0 \\
    	\overline{\alpha_{10}} & 3 & 6 & 0 \\
    	\overline{\alpha_{11}} & 4 & 5 & 0 \\
    	\overline{\zeta_{0}} & 2 & 6 & 1 \\
    	\overline{\zeta_{1}} & 2 & 7 & 1 \\
    	\overline{\zeta_{2}} & 3 & 6 & 1 \\
    	\overline{\zeta_{3}} & 4 & 5 & 1
    	\end{dmatrix}
    \end{align*}
    Then, we obtain
    \small
    \begin{align*}
    	\h_0=\begin{dmatrix}
    	2 & A_{1} & 2 & 0 & 0 \\
    	4 & A_{2} & 3 & 0 & 0 \\
    	6 & A_{3} & 4 & 0 & 0 \\
    	8 & A_{4} & 5 & 0 & 0 \\
    	10 & A_{5} & 6 & 0 & 0 \\
    	-\frac{1}{6} & C_{1} & 0 & 3 & 0 \\
    	-\frac{1}{8} & B_{1} & 0 & 4 & 0 \\
    	-\frac{1}{10} & B_{2} & 0 & 5 & 0 \\
    	-\frac{1}{10} \, \overline{\alpha_{1}} & C_{1} & 0 & 5 & 0 \\
    	-\frac{1}{12} & B_{3} & 0 & 6 & 0 \\
    	-\frac{1}{12} & \gamma_{4} & 0 & 6 & 0 \\
    	\frac{1}{6} \, \alpha_{2} & \overline{C_{1}} & 0 & 6 & 0 \\
    	-\frac{1}{12} \, \overline{\alpha_{1}} & B_{1} & 0 & 6 & 0 \\
    	-\frac{1}{12} \, \overline{\alpha_{3}} & C_{1} & 0 & 6 & 0 \\
    	\frac{1}{6} & \gamma_{1} & 1 & 3 & 0 \\
    	\frac{1}{16} & \overline{\gamma_{2}} & 1 & 4 & 0 \\
    	-\frac{1}{3} \, {\left| A_{1} \right|}^{2} & C_{1} & 1 & 4 & 0 \\
    	\frac{1}{20} & \overline{\gamma_{3}} & 1 & 5 & 0 \\
    	-\frac{1}{4} \, {\left| A_{1} \right|}^{2} & B_{1} & 1 & 5 & 0 \\
    	\frac{1}{10} \, \overline{\alpha_{1}} & \gamma_{1} & 1 & 5 & 0 \\
    	-\frac{1}{6} \, \overline{\alpha_{4}} & C_{1} & 1 & 5 & 0 \\
    	\frac{1}{4} & \overline{B_{1}} & 2 & 2 & 0 \\
    	\frac{1}{18} & \gamma_{2} & 2 & 3 & 0 \\
    	\frac{1}{6} & C_{3} & 2 & 3 & 0 \\
    	-\frac{1}{6} \, {\left| A_{1} \right|}^{2} & \overline{C_{1}} & 2 & 3 & 0 \\
    	-\frac{1}{6} \, \alpha_{1} & C_{1} & 2 & 3 & 0 \\
    	\frac{1}{6} & \gamma_{2} & 2 & 3 & 1 \\
    	\frac{3}{32} & \gamma_{5} & 2 & 4 & 0 \\
    	\frac{1}{4} & B_{4} & 2 & 4 & 0 \\
    	-\frac{1}{6} \, {\left| A_{1} \right|}^{2} & C_{2} & 2 & 4 & 0 
    	\end{dmatrix}
    	\begin{dmatrix}
    	\frac{43}{144} \, {\left| A_{1} \right|}^{2} & \gamma_{1} & 2 & 4 & 0 \\
    	-\frac{1}{8} \, \alpha_{1} & B_{1} & 2 & 4 & 0 \\
    	-\frac{5}{24} \, \alpha_{4} & C_{1} & 2 & 4 & 0 \\
    	\frac{1}{8} \, \overline{\alpha_{1}} & \overline{B_{1}} & 2 & 4 & 0 \\
    	-\frac{1}{8} \, \overline{\alpha_{4}} & \overline{C_{1}} & 2 & 4 & 0 \\
    	\frac{1}{4} & \gamma_{5} & 2 & 4 & 1 \\
    	-\frac{1}{6} \, {\left| A_{1} \right|}^{2} & \gamma_{1} & 2 & 4 & 1 \\
    	\frac{1}{2} & \overline{B_{2}} & 3 & 2 & 0 \\
    	\frac{1}{36} & \gamma_{3} & 3 & 3 & 0 \\
    	\frac{1}{3} & C_{4} & 3 & 3 & 0 \\
    	\frac{1}{6} \, {\left| A_{1} \right|}^{2} & \overline{B_{1}} & 3 & 3 & 0 \\
    	\frac{1}{6} \, \alpha_{1} & \gamma_{1} & 3 & 3 & 0 \\
    	-\frac{1}{6} \, \alpha_{3} & C_{1} & 3 & 3 & 0 \\
    	-\frac{1}{6} \, \alpha_{4} & \overline{C_{1}} & 3 & 3 & 0 \\
    	\frac{1}{3} & \gamma_{3} & 3 & 3 & 1 \\
    	\frac{3}{2} & \overline{E_{1}} & 4 & 1 & 0 \\
    	\frac{3}{4} & \overline{B_{3}} & 4 & 2 & 0 \\
    	\frac{3}{4} & \overline{\gamma_{4}} & 4 & 2 & 0 \\
    	\frac{1}{4} \, \alpha_{1} & \overline{B_{1}} & 4 & 2 & 0 \\
    	\frac{1}{12} \, \overline{\alpha_{2}} & C_{1} & 4 & 2 & 0 \\
    	2 & \overline{E_{2}} & 5 & 1 & 0 \\
    	\overline{\alpha_{2}} & \overline{C_{1}} & 5 & 1 & 0 \\
    	-\frac{1}{6} & E_{1} & -1 & 6 & 0 \\
    	-\frac{1}{7} & E_{2} & -1 & 7 & 0 \\
    	\frac{1}{42} \, \alpha_{2} & C_{1} & -1 & 7 & 0 \\
    	-4 \, {\left| A_{1} \right|}^{2} & A_{0} & 2 & 1 & 0 \\
    	-4 \, {\left| A_{1} \right|}^{2} & A_{1} & 3 & 1 & 0 \\
    	-4 \, {\left| A_{1} \right|}^{2} & A_{2} & 4 & 1 & 0 \\
    	-4 \, {\left| A_{1} \right|}^{2} & A_{3} & 5 & 1 & 0 
    	\end{dmatrix}
    	\end{align*}
    	\begin{align*}
    	\begin{dmatrix}
    	8 \, {\left| A_{1} \right|}^{4} + 4 \, \alpha_{1} \overline{\alpha_{1}} - 4 \, \beta & A_{0} & 3 & 2 & 0 \\
    	8 \, {\left| A_{1} \right|}^{4} + 4 \, \alpha_{1} \overline{\alpha_{1}} - 4 \, \beta & A_{1} & 4 & 2 & 0 \\
    	-4 \, \alpha_{1} & A_{0} & 3 & 0 & 0 \\
    	-4 \, \alpha_{1} & A_{1} & 4 & 0 & 0 \\
    	-4 \, \alpha_{1} & A_{2} & 5 & 0 & 0 \\
    	-4 \, \alpha_{1} & A_{3} & 6 & 0 & 0 \\
    	2 \, \alpha_{2} & A_{0} & 0 & 4 & 0 \\
    	2 \, \alpha_{2} & A_{1} & 1 & 4 & 0 \\
    	2 \, \alpha_{2} & A_{2} & 2 & 4 & 0 \\
    	-6 \, \alpha_{3} & A_{0} & 4 & 0 & 0 \\
    	-6 \, \alpha_{3} & A_{1} & 5 & 0 & 0 \\
    	-6 \, \alpha_{3} & A_{2} & 6 & 0 & 0 \\
    	-4 \, \alpha_{4} & A_{0} & 3 & 1 & 0 \\
    	-4 \, \alpha_{4} & A_{1} & 4 & 1 & 0 \\
    	-4 \, \alpha_{4} & A_{2} & 5 & 1 & 0 \\
    	-10 \, \alpha_{5} & A_{0} & 6 & -1 & 0 \\
    	-10 \, \alpha_{5} & A_{1} & 7 & -1 & 0 \\
    	4 \, \alpha_{1}^{2} - 8 \, \alpha_{6} - \zeta_{0} & A_{0} & 5 & 0 & 0 \\
    	4 \, \alpha_{1}^{2} - 8 \, \alpha_{6} - \zeta_{0} & A_{1} & 6 & 0 & 0 \\
    	4 \, {\left| A_{1} \right|}^{2} \overline{\alpha_{1}} - 2 \, \alpha_{7} & A_{0} & 2 & 3 & 0 \\
    	4 \, {\left| A_{1} \right|}^{2} \overline{\alpha_{1}} - 2 \, \alpha_{7} & A_{1} & 3 & 3 & 0 \\
    	12 \, \alpha_{1} \overline{\alpha_{2}} - 16 \, \alpha_{8} & A_{0} & 7 & -1 & 0 \\
    	20 \, {\left| A_{1} \right|}^{2} \overline{\alpha_{2}} + 10 \, \alpha_{1} \alpha_{3} - 14 \, \alpha_{9} - \zeta_{1} & A_{0} & 6 & 0 & 0 \\
    	16 \, \alpha_{3} {\left| A_{1} \right|}^{2} + 8 \, \alpha_{1} \alpha_{4} + 8 \, \overline{\alpha_{1}} \overline{\alpha_{2}} - 12 \, \alpha_{10} - \zeta_{2} & A_{0} & 5 & 1 & 0 \\
    	12 \, \alpha_{4} {\left| A_{1} \right|}^{2} + 6 \, \alpha_{3} \overline{\alpha_{1}} + 6 \, \alpha_{1} \overline{\alpha_{4}} - 10 \, \alpha_{11} - \zeta_{3} & A_{0} & 4 & 2 & 0 \\
    	-8 \, \zeta_{0} & A_{0} & 5 & 0 & 1 \\
    	-8 \, \zeta_{0} & A_{1} & 6 & 0 & 1 \\
    	-14 \, \zeta_{1} & A_{0} & 6 & 0 & 1 \\
    	-12 \, \zeta_{2} & A_{0} & 5 & 1 & 1 \\
    	-10 \, \zeta_{3} & A_{0} & 4 & 2 & 1 
    	\end{dmatrix}
    	\end{align*}
    	\begin{align*}
    	\begin{dmatrix}
    	-8 \, \overline{\alpha_{2}} & A_{0} & 5 & -1 & 0 \\
    	-8 \, \overline{\alpha_{2}} & A_{1} & 6 & -1 & 0 \\
    	-8 \, \overline{\alpha_{2}} & A_{2} & 7 & -1 & 0 \\
    	-2 \, \overline{\alpha_{4}} & A_{0} & 2 & 2 & 0 \\
    	-2 \, \overline{\alpha_{4}} & A_{1} & 3 & 2 & 0 \\
    	-2 \, \overline{\alpha_{4}} & A_{2} & 4 & 2 & 0 \\
    	2 \, \overline{\alpha_{5}} & A_{0} & 0 & 5 & 0 \\
    	2 \, \overline{\alpha_{5}} & A_{1} & 1 & 5 & 0 \\
    	-\overline{\zeta_{0}} & A_{0} & 1 & 4 & 0 \\
    	-\overline{\zeta_{0}} & A_{1} & 2 & 4 & 0 \\
    	12 \, \alpha_{1} {\left| A_{1} \right|}^{2} - 6 \, \overline{\alpha_{7}} & A_{0} & 4 & 1 & 0 \\
    	12 \, \alpha_{1} {\left| A_{1} \right|}^{2} - 6 \, \overline{\alpha_{7}} & A_{1} & 5 & 1 & 0 \\
    	-2 \, \alpha_{2} \overline{\alpha_{1}} - 2 \, \overline{\alpha_{8}} & A_{0} & 0 & 6 & 0 \\
    	-4 \, \overline{\alpha_{9}} - \overline{\zeta_{1}} & A_{0} & 1 & 5 & 0 \\
    	4 \, {\left| A_{1} \right|}^{2} \overline{\alpha_{3}} + 2 \, \alpha_{1} \alpha_{2} + 2 \, \overline{\alpha_{1}} \overline{\alpha_{4}} - 6 \, \overline{\alpha_{10}} - \overline{\zeta_{2}} & A_{0} & 2 & 4 & 0 \\
    	8 \, {\left| A_{1} \right|}^{2} \overline{\alpha_{4}} + 4 \, \alpha_{4} \overline{\alpha_{1}} + 4 \, \alpha_{1} \overline{\alpha_{3}} - 8 \, \overline{\alpha_{11}} - \overline{\zeta_{3}} & A_{0} & 3 & 3 & 0 \\
    	-4 \, \overline{\zeta_{1}} & A_{0} & 1 & 5 & 1 \\
    	-6 \, \overline{\zeta_{2}} & A_{0} & 2 & 4 & 1 \\
    	-8 \, \overline{\zeta_{3}} & A_{0} & 3 & 3 & 1
    	\end{dmatrix}
    \end{align*}
    
    We now come to the last expansion of the quartic form.    
    We first have    
    \footnotesize
    \begin{align*}
    	&g^{-1}\otimes Q(\h_0)=\\
    	&\begin{dmatrix}
    	-6 \, A_{0} C_{1} \alpha_{1} + 6 \, A_{2} C_{1} - A_{1} \gamma_{1} & 0 & 0 & 0 \\
    	-10 \, A_{1} C_{1} \alpha_{1} - 12 \, A_{0} C_{1} \alpha_{3} + 4 \, A_{0} \alpha_{1} \gamma_{1} + 12 \, A_{3} C_{1} - 4 \, A_{2} \gamma_{1} + \frac{1}{2} \, A_{1} \gamma_{2} & 1 & 0 & 0 \\
    	\omega_1 & 2 & 0 & 0 \\
    	6 \, A_{1} E_{1} & -2 & 3 & 0 \\
    	\frac{1}{18} \, C_{1}^{2} {\left| A_{1} \right|}^{2} - 10 \, A_{0} E_{1} {\left| A_{1} \right|}^{2} - \frac{10}{3} \, A_{1} C_{1} \alpha_{2} - \frac{1}{3} \, A_{0} \alpha_{2} \gamma_{1} + \frac{1}{6} \, A_{0} C_{1} \overline{\zeta_{0}} + 6 \, A_{1} E_{2} + \frac{1}{48} \, B_{1} \gamma_{1} - \frac{1}{96} \, C_{1} \overline{\gamma_{2}} & -2 & 4 & 0 \\
    	\omega_2 & -1 & 3 & 0 \\
    	\omega_3 & 0 & 2 & 0 \\
    	\omega_4 & 0 & 1 & 0 \\
    	\omega_5 & 1 & 1 & 0 \\
    	\omega_6 & 3 & -1 & 0 \\
    	-16 \, A_{0}^{2} \alpha_{1} \overline{\alpha_{4}} + 16 \, {\left(2 \, {\left| A_{1} \right|}^{4} + \alpha_{1} \overline{\alpha_{1}} - \beta\right)} A_{0} A_{1} + 2 \, A_{0} \alpha_{1} \overline{B_{1}} - \frac{80}{3} \, A_{0} C_{1} \overline{\alpha_{2}} - 8 \, A_{1}^{2} \overline{\alpha_{4}} + 16 \, A_{0} A_{2} \overline{\alpha_{4}} - 2 \, A_{2} \overline{B_{1}} + 2 \, A_{1} \overline{B_{2}} & 2 & -1 & 0 \\
    	48 \, A_{0}^{2} \alpha_{2} {\left| A_{1} \right|}^{2} - 3 \, A_{0} B_{1} {\left| A_{1} \right|}^{2} - \frac{2}{3} \, A_{0} C_{1} \overline{\alpha_{4}} - 40 \, A_{0} A_{1} \overline{\alpha_{5}} + 2 \, A_{1} B_{2} + \frac{1}{12} \, C_{1} \overline{B_{1}} & -1 & 2 & 0 \\
        \omega_7 & 4 & -2 & 0 \\
    	2 \, A_{0} \alpha_{1} \gamma_{2} - 20 \, A_{0} C_{1} \zeta_{0} - 48 \, A_{0} A_{1} \overline{\zeta_{3}} - 2 \, A_{2} \gamma_{2} + 2 \, A_{1} \gamma_{3} & 2 & 0 & 1 \\
    	-32 \, A_{0} A_{1} \alpha_{1} {\left| A_{1} \right|}^{2} - 48 \, A_{0}^{2} \alpha_{3} {\left| A_{1} \right|}^{2} - 16 \, A_{1} A_{2} {\left| A_{1} \right|}^{2} + 48 \, A_{0} A_{3} {\left| A_{1} \right|}^{2} + 24 \, {\left(2 \, \alpha_{1} {\left| A_{1} \right|}^{2} - \overline{\alpha_{7}}\right)} A_{0} A_{1} - 16 \, A_{1}^{2} \alpha_{4} + 6 \, A_{1} \overline{E_{1}} & 3 & -2 & 0 \\
    	-\frac{1}{12} \, C_{1} E_{1} & -4 & 6 & 0 \\
    	-16 \, A_{0}^{2} \alpha_{1} {\left| A_{1} \right|}^{2} - 8 \, A_{1}^{2} {\left| A_{1} \right|}^{2} + 16 \, A_{0} A_{2} {\left| A_{1} \right|}^{2} - 8 \, A_{0} A_{1} \alpha_{4} & 2 & -2 & 0 \\
    	-\frac{8}{3} \, A_{0} C_{1} {\left| A_{1} \right|}^{2} - 32 \, A_{0} A_{1} \alpha_{2} + 2 \, A_{1} B_{1} & -1 & 1 & 0 \\
    	-320 \, A_{0}^{2} \alpha_{5} {\left| A_{1} \right|}^{2} - 480 \, A_{0} A_{1} {\left| A_{1} \right|}^{2} \overline{\alpha_{2}} + 32 \, A_{0}^{2} \alpha_{1} \zeta_{0} - 128 \, A_{0}^{2} \alpha_{4} \overline{\alpha_{2}} - 32 \, A_{1}^{2} \zeta_{0} - 32 \, A_{0} A_{2} \zeta_{0} - 56 \, A_{0} A_{1} \zeta_{1} & 5 & -3 & 0 \\
    	-96 \, A_{0}^{2} \zeta_{0} {\left| A_{1} \right|}^{2} - 72 \, A_{0} A_{1} \zeta_{2} & 4 & -2 & 1 \\
    	-192 \, A_{0}^{2} {\left| A_{1} \right|}^{2} \overline{\alpha_{2}} - 24 \, A_{0} A_{1} \zeta_{0} & 4 & -3 & 0 \\
    	-64 \, A_{0}^{2} \alpha_{1} \overline{\alpha_{2}} + 80 \, A_{0} A_{1} \alpha_{5} + 64 \, A_{1}^{2} \overline{\alpha_{2}} + 64 \, A_{0} A_{2} \overline{\alpha_{2}} & 5 & -4 & 0 \\
    	\omega_8 & 6 & -4 & 0 \\
    	-80 \, A_{0} A_{1} \zeta_{3} & 3 & -1 & 1 \\
    	48 \, A_{0} A_{1} \overline{\alpha_{2}} & 4 & -4 & 0 \\
    	2 \, A_{1} C_{1} & -1 & 0 & 0 \\
    	80 \, A_{0} A_{1} \overline{\zeta_{1}} & 0 & 2 & 1
    	\end{dmatrix}
    \end{align*}
    where
    \begin{align*}
    	\omega_{1}&=-16 \, {\left(2 \, {\left| A_{1} \right|}^{4} + \alpha_{1} \overline{\alpha_{1}} - \beta\right)} A_{0}^{2} {\left| A_{1} \right|}^{2} - 2 \, A_{0} \alpha_{1} {\left| A_{1} \right|}^{2} \overline{C_{1}} + 24 \, {\left(2 \, {\left| A_{1} \right|}^{2} \overline{\alpha_{1}} - \alpha_{7}\right)} A_{0}^{2} \alpha_{1} - 2 \, A_{0} C_{1} \alpha_{1}^{2} + 2 \, A_{1} {\left| A_{1} \right|}^{2} \overline{B_{1}}\\
    	& - 2 \, A_{0} {\left| A_{1} \right|}^{2} \overline{B_{2}} + 2 \, A_{2} {\left| A_{1} \right|}^{2} \overline{C_{1}} + 400 \, A_{0}^{2} \alpha_{2} \overline{\alpha_{2}} - 8 \, A_{0}^{2} \alpha_{4} \overline{\alpha_{4}} + 6 \, {\left(8 \, {\left| A_{1} \right|}^{2} \overline{\alpha_{4}} + 4 \, \alpha_{4} \overline{\alpha_{1}} + 4 \, \alpha_{1} \overline{\alpha_{3}} - 8 \, \overline{\alpha_{11}} - \overline{\zeta_{3}}\right)} A_{0} A_{1}\\
    	& + 12 \, {\left(2 \, {\left| A_{1} \right|}^{2} \overline{\alpha_{1}} - \alpha_{7}\right)} A_{1}^{2} - 24 \, {\left(2 \, {\left| A_{1} \right|}^{2} \overline{\alpha_{1}} - \alpha_{7}\right)} A_{0} A_{2} + \frac{5}{2} \, {\left(4 \, \alpha_{1}^{2} - 8 \, \alpha_{6} - \zeta_{0}\right)} A_{0} C_{1} - 8 \, A_{2} C_{1} \alpha_{1} + 2 \, A_{0} C_{3} \alpha_{1} - 18 \, A_{1} C_{1} \alpha_{3}\\
    	& + 7 \, A_{1} \alpha_{1} \gamma_{1} + 9 \, A_{0} \alpha_{3} \gamma_{1} + \frac{4}{3} \, A_{0} C_{1} \zeta_{0} + A_{0} \alpha_{4} \overline{B_{1}} - A_{1} \alpha_{4} \overline{C_{1}} - 25 \, A_{0} B_{1} \overline{\alpha_{2}} - 32 \, A_{0} A_{1} \overline{\zeta_{3}} + 20 \, A_{4} C_{1} - 2 \, A_{2} C_{3} + 2 \, A_{1} C_{4}\\
    	& + {\left(6 \, A_{0} C_{1} \alpha_{1} - 6 \, A_{2} C_{1} + A_{1} \gamma_{1}\right)} \alpha_{1} - 9 \, A_{3} \gamma_{1} + \frac{3}{2} \, A_{1} \gamma_{3} + 8 \, {\left(2 \, A_{0}^{2} \alpha_{1} {\left| A_{1} \right|}^{2} + A_{1}^{2} {\left| A_{1} \right|}^{2} - 2 \, A_{0} A_{2} {\left| A_{1} \right|}^{2} + A_{0} A_{1} \alpha_{4}\right)} \overline{\alpha_{1}}\\
    	\omega_{2}&=-\frac{16}{5} \, A_{0} C_{1} {\left| A_{1} \right|}^{2} \overline{\alpha_{1}} + 64 \, A_{0}^{2} {\left| A_{1} \right|}^{2} \overline{\alpha_{5}} - \frac{16}{5} \, A_{0} B_{2} {\left| A_{1} \right|}^{2} + 16 \, A_{0}^{2} \alpha_{2} \overline{\alpha_{4}} + 48 \, {\left(\alpha_{2} \overline{\alpha_{1}} + \overline{\alpha_{8}}\right)} A_{0} A_{1} - 16 \, A_{0} E_{1} \alpha_{1} - 2 \, A_{0} \alpha_{2} \overline{B_{1}}\\
    	& - 4 \, A_{1} \alpha_{2} \overline{C_{1}} + 2 \, A_{1} B_{1} \overline{\alpha_{1}} - A_{0} B_{1} \overline{\alpha_{4}} + 2 \, A_{1} B_{3} + 16 \, A_{2} E_{1} - \frac{1}{36} \, C_{1} \gamma_{2} + 2 \, A_{1} \gamma_{4} + \frac{1}{8} \, B_{1} \overline{B_{1}}\\
    	& + \frac{2}{3} \, {\left(4 \, A_{0} C_{1} {\left| A_{1} \right|}^{2} + 48 \, A_{0} A_{1} \alpha_{2} - 3 \, A_{1} B_{1}\right)} \overline{\alpha_{1}}\\
    	\omega_{3}&=-4 \, A_{0} C_{1} {\left| A_{1} \right|}^{4} + 96 \, A_{0} A_{1} \alpha_{2} {\left| A_{1} \right|}^{2} - 12 \, A_{0}^{2} {\left| A_{1} \right|}^{2} \overline{\zeta_{0}} + 72 \, A_{0}^{2} \alpha_{2} \alpha_{4} - 2 \, A_{1} B_{1} {\left| A_{1} \right|}^{2} - 6 \, A_{0} C_{1} \alpha_{1} \overline{\alpha_{1}} + 120 \, A_{0}^{2} \alpha_{1} \overline{\alpha_{5}}\\
    	& + \frac{3}{4} \, A_{0} {\left| A_{1} \right|}^{2} \overline{\gamma_{2}} + 2 \, {\left(2 \, {\left| A_{1} \right|}^{4} + \alpha_{1} \overline{\alpha_{1}} - \beta\right)} A_{0} C_{1} - 6 \, A_{0} B_{2} \alpha_{1} - \frac{9}{2} \, A_{0} B_{1} \alpha_{4} + 10 \, A_{0} A_{1} {\left(4 \, \overline{\alpha_{9}} + \overline{\zeta_{1}}\right)}\\
    	& + \frac{4}{3} \, {\left(4 \, A_{0} C_{1} {\left| A_{1} \right|}^{2} + 48 \, A_{0} A_{1} \alpha_{2} - 3 \, A_{1} B_{1}\right)} {\left| A_{1} \right|}^{2} + 6 \, A_{2} C_{1} \overline{\alpha_{1}} - A_{1} \gamma_{1} \overline{\alpha_{1}} - \frac{4}{3} \, A_{1} C_{1} \overline{\alpha_{4}} + \frac{1}{3} \, A_{0} \gamma_{1} \overline{\alpha_{4}} - 20 \, A_{1}^{2} \overline{\alpha_{5}} - 120 \, A_{0} A_{2} \overline{\alpha_{5}}\\
    	& - 16 \, A_{0} A_{1} \overline{\zeta_{1}} + 6 \, A_{2} B_{2} - \frac{1}{24} \, \gamma_{1} \overline{B_{1}} + \frac{1}{4} \, C_{1} \overline{B_{2}} + {\left(6 \, A_{0} C_{1} \alpha_{1} - 6 \, A_{2} C_{1} + A_{1} \gamma_{1}\right)} \overline{\alpha_{1}} - \frac{1}{2} \, A_{1} \overline{\gamma_{3}}\\
    	\omega_{4}&=96 \, A_{0}^{2} \alpha_{1} \alpha_{2} - \frac{16}{3} \, A_{1} C_{1} {\left| A_{1} \right|}^{2} + \frac{4}{3} \, A_{0} \gamma_{1} {\left| A_{1} \right|}^{2} - 6 \, A_{0} B_{1} \alpha_{1} - 16 \, A_{1}^{2} \alpha_{2} - 96 \, A_{0} A_{2} \alpha_{2} - 4 \, A_{0} C_{1} \alpha_{4} + 8 \, A_{0} A_{1} \overline{\zeta_{0}}\\
    	& + 6 \, A_{2} B_{1} - \frac{1}{2} \, A_{1} \overline{\gamma_{2}}\\
    	\omega_{5}&=-\frac{32}{3} \, A_{0} C_{1} \alpha_{1} {\left| A_{1} \right|}^{2} + 192 \, A_{0} A_{1} \alpha_{1} \alpha_{2} + 192 \, A_{0}^{2} \alpha_{2} \alpha_{3} + \frac{16}{3} \, A_{2} C_{1} {\left| A_{1} \right|}^{2} + 2 \, A_{1} \gamma_{1} {\left| A_{1} \right|}^{2} - \frac{2}{3} \, A_{0} \gamma_{2} {\left| A_{1} \right|}^{2} - 32 \, A_{0}^{2} \alpha_{1} \overline{\zeta_{0}}\\
    	& + 8 \, {\left(2 \, \alpha_{1} {\left| A_{1} \right|}^{2} - \overline{\alpha_{7}}\right)} A_{0} C_{1} - 8 \, A_{1} B_{1} \alpha_{1} - 64 \, A_{1} A_{2} \alpha_{2} - 192 \, A_{0} A_{3} \alpha_{2} - 12 \, A_{0} B_{1} \alpha_{3} - \frac{22}{3} \, A_{1} C_{1} \alpha_{4} + \frac{8}{3} \, A_{0} \alpha_{4} \gamma_{1}\\
    	& + 2 \, {\left(6 \, A_{0} C_{1} \alpha_{1} - 6 \, A_{2} C_{1} + A_{1} \gamma_{1}\right)} {\left| A_{1} \right|}^{2} + 2 \, A_{0} \alpha_{1} \overline{\gamma_{2}} + 32 \, A_{0} A_{2} \overline{\zeta_{0}} - 24 \, A_{0} A_{1} \overline{\zeta_{2}} + 12 \, A_{3} B_{1}\\
    	& + \frac{2}{3} \, {\left(4 \, A_{0} C_{1} {\left| A_{1} \right|}^{2} + 48 \, A_{0} A_{1} \alpha_{2} - 3 \, A_{1} B_{1}\right)} \alpha_{1} + A_{1} \gamma_{5} + 2 \, C_{1} \overline{E_{1}} - 2 \, A_{2} \overline{\gamma_{2}}\\
    	\omega_{6}&=-32 \, A_{0} A_{1} \alpha_{1} \overline{\alpha_{4}} - 48 \, A_{0}^{2} \alpha_{3} \overline{\alpha_{4}} + 8 \, {\left(12 \, \alpha_{4} {\left| A_{1} \right|}^{2} + 6 \, \alpha_{3} \overline{\alpha_{1}} + 6 \, \alpha_{1} \overline{\alpha_{4}} - 10 \, \alpha_{11} - \zeta_{3}\right)} A_{0} A_{1} + 32 \, {\left(2 \, {\left| A_{1} \right|}^{4} + \alpha_{1} \overline{\alpha_{1}} - \beta\right)} A_{1}^{2}\\
    	& - 40 \, A_{0} C_{1} \alpha_{5} - 40 \, A_{0} A_{1} \zeta_{3} + 16 \, {\left(2 \, A_{0}^{2} \alpha_{1} {\left| A_{1} \right|}^{2} + A_{1}^{2} {\left| A_{1} \right|}^{2} - 2 \, A_{0} A_{2} {\left| A_{1} \right|}^{2} + A_{0} A_{1} \alpha_{4}\right)} {\left| A_{1} \right|}^{2} + 6 \, A_{1} \alpha_{1} \overline{B_{1}} + 6 \, A_{0} \alpha_{3} \overline{B_{1}}\\
    	& - \frac{100}{3} \, A_{1} C_{1} \overline{\alpha_{2}} + \frac{64}{3} \, A_{0} \gamma_{1} \overline{\alpha_{2}} - 16 \, A_{1} A_{2} \overline{\alpha_{4}} + 48 \, A_{0} A_{3} \overline{\alpha_{4}} - 6 \, A_{3} \overline{B_{1}} + 6 \, A_{1} \overline{B_{3}} + 6 \, A_{1} \overline{\gamma_{4}}\\
    	\omega_{7}&=12 \, {\left(4 \, \alpha_{1}^{2} - 8 \, \alpha_{6} - \zeta_{0}\right)} A_{0}^{2} {\left| A_{1} \right|}^{2} - 16 \, A_{1}^{2} \alpha_{1} {\left| A_{1} \right|}^{2} - 32 \, A_{0} A_{2} \alpha_{1} {\left| A_{1} \right|}^{2} - 96 \, A_{0} A_{1} \alpha_{3} {\left| A_{1} \right|}^{2} + 32 \, A_{0}^{2} \zeta_{0} {\left| A_{1} \right|}^{2}\\
    	& - 24 \, {\left(2 \, \alpha_{1} {\left| A_{1} \right|}^{2} - \overline{\alpha_{7}}\right)} A_{0}^{2} \alpha_{1} - 24 \, A_{0}^{2} \alpha_{3} \alpha_{4} - 16 \, A_{2}^{2} {\left| A_{1} \right|}^{2} + 96 \, A_{0} A_{4} {\left| A_{1} \right|}^{2} - 48 \, A_{0} A_{1} \overline{\alpha_{1}} \overline{\alpha_{2}} - 144 \, A_{0}^{2} \overline{\alpha_{2}} \overline{\alpha_{4}}\\
    	& + 6 \, {\left(16 \, \alpha_{3} {\left| A_{1} \right|}^{2} + 8 \, \alpha_{1} \alpha_{4} + 8 \, \overline{\alpha_{1}} \overline{\alpha_{2}} - 12 \, \alpha_{10} - \zeta_{2}\right)} A_{0} A_{1} + 36 \, {\left(2 \, \alpha_{1} {\left| A_{1} \right|}^{2} - \overline{\alpha_{7}}\right)} A_{1}^{2} + 24 \, {\left(2 \, \alpha_{1} {\left| A_{1} \right|}^{2} - \overline{\alpha_{7}}\right)} A_{0} A_{2}\\
    	& - 40 \, A_{1} A_{2} \alpha_{4} + 24 \, A_{0} A_{3} \alpha_{4} - 48 \, A_{0} A_{1} \zeta_{2} - 6 \, A_{0} \alpha_{1} \overline{E_{1}} + 18 \, A_{0} \overline{B_{1}} \overline{\alpha_{2}} + 6 \, A_{1} \overline{C_{1}} \overline{\alpha_{2}}\\
    	& + 8 \, {\left(2 \, A_{0}^{2} \alpha_{1} {\left| A_{1} \right|}^{2} + A_{1}^{2} {\left| A_{1} \right|}^{2} - 2 \, A_{0} A_{2} {\left| A_{1} \right|}^{2} + A_{0} A_{1} \alpha_{4}\right)} \alpha_{1} + 6 \, A_{2} \overline{E_{1}} + 12 \, A_{1} \overline{E_{2}}\\
    	\omega_{8}&=-120 \, A_{0}^{2} \alpha_{1} \alpha_{5} - 176 \, A_{0} A_{1} \alpha_{1} \overline{\alpha_{2}} - 48 \, A_{0}^{2} \alpha_{3} \overline{\alpha_{2}} - 40 \, {\left(3 \, \alpha_{1} \overline{\alpha_{2}} - 4 \, \alpha_{8}\right)} A_{0} A_{1} + 100 \, A_{1}^{2} \alpha_{5} + 120 \, A_{0} A_{2} \alpha_{5} + 176 \, A_{1} A_{2} \overline{\alpha_{2}}
    	\\
    	& + 48 \, A_{0} A_{3} \overline{\alpha_{2}}
    \end{align*}
    \normalsize
    Now, we remark as
    \begin{align*}
    	|\H|^2\h_0\totimes\h_0=O(|z|^2),\quad \s{\H}{\h_0}^2=O(|z|^2)
    \end{align*}
    that the development of these two tensors $|\H|^2\h_0\totimes\h_0$ and $\s{\H}{\h_0}^2$ up to order $O(|z|^{3-\epsilon})$ only depends on their first order terms in their Taylor expansion, and as the following expansion is valid at a branch point of multiplicity $\theta_0\geq 3$
    \begin{align*}
    	\H=\Re\left(\frac{\vec{C}_1}{z^{\theta_0-2}}\right)+O(|z|^{3-\epsilon}),\quad \h_0=2\vec{A}_1z^{\theta_0-1}dz^2+O(|z|^{\theta_0})
    \end{align*}
    we can use \eqref{hh2} to obtain
    	\begin{align}\label{hh2}
    \frac{5}{4}|\H|^2\h_0\totimes\h_0&=\begin{dmatrix}
    \left(4 \, A_{1}^{2}\right) \left(\frac{5}{16} \, C_{1}^{2}\right) & 2 & 0 \\
    \left(4 \, A_{1}^{2}\right) \left(\frac{5}{16} \, C_{1} \overline{C_{1}}\right) & \theta_{0} & -\theta_{0} + 2 \\
    \left(4 \, A_{1}^{2}\right) \left(\frac{5}{16} \, C_{1} \overline{C_{1}}\right) & \theta_{0} & -\theta_{0} + 2 \\
    \left(4 \, A_{1}^{2}\right) \left(\frac{5}{16} \, \overline{C_{1}}^{2}\right) & 2 \, \theta_{0} - 2 & -2 \, \theta_{0} + 4
    \end{dmatrix}\\
    \s{\H}{\h_0}^2&=\begin{dmatrix}
    \left(A_{1} C_{1}\right) \left(A_{1} C_{1}\right) & 2 & 0 \\
    \left(A_{1} C_{1}\right) \left(A_{1} \overline{C_{1}}\right) & \theta_{0} & -\theta_{0} + 2 \\
    \left(A_{1} C_{1}\right) \left(A_{1} \overline{C_{1}}\right) & \theta_{0} & -\theta_{0} + 2 \\
    \left(A_{1} \overline{C_{1}}\right) \left(A_{1} \overline{C_{1}}\right) & 2 \, \theta_{0} - 2 & -2 \, \theta_{0} + 4
    \end{dmatrix}
    \end{align}
    which yields for $\theta_0=3$
    \begin{align}\label{3devextra1}
    	&\frac{5}{4}|\H|^2\h_0\totimes\h_0=\frac{5}{4}\s{\vec{A}_1}{\vec{A}_1}\s{\vec{C}_1}{\vec{C}_1}z^2dz^4+\frac{5}{2}|\vec{C}_1|^2\s{\vec{A}_1}{\vec{A}_1}z^{3}\z^{-1}dz^4+\frac{5}{4}\s{\vec{A}_1}{\vec{A}_1}\bar{\s{\vec{C}_1}{\vec{C}_ 1}}z^{4}\z^{-2}dz^4\nonumber\\
    	&\s{\H}{\h_0}^2=\s{\vec{A}_1}{\vec{C}_1}^2z^2dz^4+2\s{\vec{A}_1}{\vec{C}_1}\s{\vec{A}_1}{\bar{\vec{C}_1}}z^{3}\z^{-1}dz^4+\s{\vec{A}_1}{\bar{\vec{C}_1}}^2z^{4}\z^{-2}dz^4
    \end{align}
    Therefore, we need only develop the Gauss curvature, which is
    \footnotesize
    \begin{align*}
    	&-K_g\h_0\totimes\h_0=\\
    	&\begin{dmatrix}
    	16 \, A_{1}^{2} {\left| A_{1} \right|}^{2} & 2 & -2 & 0 \\
    	-64 \, A_{0} A_{1} \alpha_{1} {\left| A_{1} \right|}^{2} + 64 \, A_{1} A_{2} {\left| A_{1} \right|}^{2} + 16 \, A_{1}^{2} \alpha_{4} & 3 & -2 & 0 \\
    	\kappa_1 & 4 & -2 & 0 \\
    	-\frac{8}{3} \, A_{1} C_{1} {\left| A_{1} \right|}^{2} - 32 \, A_{1}^{2} \alpha_{2} & 0 & 1 & 0 \\
    	\frac{16}{3} \, A_{0} C_{1} {\left| A_{1} \right|}^{4} + 160 \, A_{0} A_{1} \alpha_{2} {\left| A_{1} \right|}^{2} - 2 \, A_{1} B_{1} {\left| A_{1} \right|}^{2} - \frac{8}{3} \, A_{1} C_{1} \overline{\alpha_{4}} - 40 \, A_{1}^{2} \overline{\alpha_{5}} & 0 & 2 & 0 \\
    	\frac{16}{3} \, A_{0} C_{1} \alpha_{1} {\left| A_{1} \right|}^{2} + 128 \, A_{0} A_{1} \alpha_{1} \alpha_{2} - \frac{16}{3} \, A_{2} C_{1} {\left| A_{1} \right|}^{2} + \frac{8}{3} \, A_{1} \gamma_{1} {\left| A_{1} \right|}^{2} - 128 \, A_{1} A_{2} \alpha_{2} - \frac{8}{3} \, A_{1} C_{1} \alpha_{4} & 1 & 1 & 0 \\
    	64 \, A_{0}^{2} {\left| A_{1} \right|}^{6} - 16 \, A_{1}^{2} {\left| A_{1} \right|}^{2} \overline{\alpha_{1}} - 96 \, A_{0} A_{1} {\left| A_{1} \right|}^{2} \overline{\alpha_{4}} + 4 \, A_{1} {\left| A_{1} \right|}^{2} \overline{B_{1}} - 24 \, {\left(2 \, {\left| A_{1} \right|}^{2} \overline{\alpha_{1}} - \alpha_{7}\right)} A_{1}^{2} & 2 & 0 & 0 \\
    	-64 \, A_{0} A_{1} {\left| A_{1} \right|}^{4} + 16 \, A_{1}^{2} \overline{\alpha_{4}} & 2 & -1 & 0 \\
    	\kappa_2 & 3 & -1 & 0 \\
    	16 \, A_{1}^{2} \zeta_{0} & 5 & -3 & 0 \\
    	\frac{1}{9} \, C_{1}^{2} {\left| A_{1} \right|}^{2} + \frac{16}{3} \, A_{1} C_{1} \alpha_{2} & -2 & 4 & 0 \\
    	128 \, A_{0} A_{1} \alpha_{1} \overline{\alpha_{2}} - 40 \, A_{1}^{2} \alpha_{5} - 128 \, A_{1} A_{2} \overline{\alpha_{2}} & 6 & -4 & 0 \\
    	-32 \, A_{1}^{2} \overline{\alpha_{2}} & 5 & -4 & 0
    	\end{dmatrix}
    \end{align*}
    \normalsize
    where
    \small
    \begin{align*}
    	\kappa_1&=64 \, A_{0}^{2} \alpha_{1}^{2} {\left| A_{1} \right|}^{2} - 80 \, A_{1}^{2} \alpha_{1} {\left| A_{1} \right|}^{2} - 128 \, A_{0} A_{2} \alpha_{1} {\left| A_{1} \right|}^{2} - 96 \, A_{0} A_{1} \alpha_{3} {\left| A_{1} \right|}^{2} - 64 \, A_{0} A_{1} \alpha_{1} \alpha_{4} + 64 \, A_{2}^{2} {\left| A_{1} \right|}^{2} + 96 \, A_{1} A_{3} {\left| A_{1} \right|}^{2}\\
    	& - 24 \, {\left(2 \, \alpha_{1} {\left| A_{1} \right|}^{2} - \overline{\alpha_{7}}\right)} A_{1}^{2} + 64 \, A_{1} A_{2} \alpha_{4}\\
    	\kappa_2&=128 \, A_{0}^{2} \alpha_{1} {\left| A_{1} \right|}^{4} - 96 \, A_{1}^{2} {\left| A_{1} \right|}^{4} - 128 \, A_{0} A_{2} {\left| A_{1} \right|}^{4} - 128 \, A_{0} A_{1} \alpha_{4} {\left| A_{1} \right|}^{2} - 64 \, A_{0} A_{1} \alpha_{1} \overline{\alpha_{4}} - 32 \, {\left(2 \, {\left| A_{1} \right|}^{4} + \alpha_{1} \overline{\alpha_{1}} - \beta\right)} A_{1}^{2}\\
    	& + \frac{16}{3} \, A_{1} C_{1} \overline{\alpha_{2}} + 64 \, A_{1} A_{2} \overline{\alpha_{4}}
    \end{align*}
    \normalsize
    Finally, we have
    \footnotesize
    \begin{align}\label{part}
    	&g^{-1}\otimes Q(\h_0)-K_g\h_0\totimes\h_0=\\
    	&\begin{dmatrix}
    	-6 \, A_{0} C_{1} \alpha_{1} + 6 \, A_{2} C_{1} - A_{1} \gamma_{1} & 0 & 0 & 0 \\
    	-10 \, A_{1} C_{1} \alpha_{1} - 12 \, A_{0} C_{1} \alpha_{3} + 12 \, A_{3} C_{1} + 4 \, {\left(A_{0} \alpha_{1} - A_{2}\right)} \gamma_{1} + \frac{1}{2} \, A_{1} \gamma_{2} & 1 & 0 & 0 \\
    	\pi_1 & 2 & 0 & 0 \\
    	6 \, A_{1} E_{1} & -2 & 3 & 0 \\
    	\frac{1}{6} \, C_{1}^{2} {\left| A_{1} \right|}^{2} - 10 \, A_{0} E_{1} {\left| A_{1} \right|}^{2} + 2 \, A_{1} C_{1} \alpha_{2} + \frac{1}{6} \, A_{0} C_{1} \overline{\zeta_{0}} + 6 \, A_{1} E_{2} - \frac{1}{48} \, {\left(16 \, A_{0} \alpha_{2} - B_{1}\right)} \gamma_{1} - \frac{1}{96} \, C_{1} \overline{\gamma_{2}} & -2 & 4 & 0 \\
    	\pi_2 & -1 & 3 & 0 \\
    	\pi_3 & 0 & 2 & 0 \\
    	\pi_4 & 0 & 1 & 0 \\
    	\pi_5 & 1 & 1 & 0 \\
        \pi_6 & 3 & -1 & 0 \\
    	\pi_7 & 2 & -1 & 0 \\
    	48 \, A_{0}^{2} \alpha_{2} {\left| A_{1} \right|}^{2} - 3 \, A_{0} B_{1} {\left| A_{1} \right|}^{2} - \frac{2}{3} \, A_{0} C_{1} \overline{\alpha_{4}} - 40 \, A_{0} A_{1} \overline{\alpha_{5}} + 2 \, A_{1} B_{2} + \frac{1}{12} \, C_{1} \overline{B_{1}} & -1 & 2 & 0 \\
    	\pi_8 & 4 & -2 & 0 \\
    	-20 \, A_{0} C_{1} \zeta_{0} - 48 \, A_{0} A_{1} \overline{\zeta_{3}} + 2 \, {\left(A_{0} \alpha_{1} - A_{2}\right)} \gamma_{2} + 2 \, A_{1} \gamma_{3} & 2 & 0 & 1 \\
    	-48 \, A_{0} A_{1} \alpha_{1} {\left| A_{1} \right|}^{2} - 48 \, A_{0}^{2} \alpha_{3} {\left| A_{1} \right|}^{2} + 48 \, A_{1} A_{2} {\left| A_{1} \right|}^{2} + 48 \, A_{0} A_{3} {\left| A_{1} \right|}^{2} - 24 \, A_{0} A_{1} \overline{\alpha_{7}} + 6 \, A_{1} \overline{E_{1}} & 3 & -2 & 0 \\
    	-\frac{1}{12} \, C_{1} E_{1} & -4 & 6 & 0 \\
    	-16 \, A_{0}^{2} \alpha_{1} {\left| A_{1} \right|}^{2} + 8 \, A_{1}^{2} {\left| A_{1} \right|}^{2} + 16 \, A_{0} A_{2} {\left| A_{1} \right|}^{2} - 8 \, A_{0} A_{1} \alpha_{4} & 2 & -2 & 0 \\
    	-\frac{8}{3} \, A_{0} C_{1} {\left| A_{1} \right|}^{2} - 32 \, A_{0} A_{1} \alpha_{2} + 2 \, A_{1} B_{1} & -1 & 1 & 0 \\
    	-320 \, A_{0}^{2} \alpha_{5} {\left| A_{1} \right|}^{2} - 480 \, A_{0} A_{1} {\left| A_{1} \right|}^{2} \overline{\alpha_{2}} - 128 \, A_{0}^{2} \alpha_{4} \overline{\alpha_{2}} - 56 \, A_{0} A_{1} \zeta_{1} + 16 \, {\left(2 \, A_{0}^{2} \alpha_{1} - A_{1}^{2} - 2 \, A_{0} A_{2}\right)} \zeta_{0} & 5 & -3 & 0 \\
    	-96 \, A_{0}^{2} \zeta_{0} {\left| A_{1} \right|}^{2} - 72 \, A_{0} A_{1} \zeta_{2} & 4 & -2 & 1 \\
    	-192 \, A_{0}^{2} {\left| A_{1} \right|}^{2} \overline{\alpha_{2}} - 24 \, A_{0} A_{1} \zeta_{0} & 4 & -3 & 0 \\
    	80 \, A_{0} A_{1} \alpha_{5} - 32 \, {\left(2 \, A_{0}^{2} \alpha_{1} - A_{1}^{2} - 2 \, A_{0} A_{2}\right)} \overline{\alpha_{2}} & 5 & -4 & 0 \\
    	-48 \, A_{0}^{2} \alpha_{3} \overline{\alpha_{2}} + 160 \, A_{0} A_{1} \alpha_{8} - 60 \, {\left(2 \, A_{0}^{2} \alpha_{1} - A_{1}^{2} - 2 \, A_{0} A_{2}\right)} \alpha_{5} - 24 \, {\left(7 \, A_{0} A_{1} \alpha_{1} - 2 \, A_{1} A_{2} - 2 \, A_{0} A_{3}\right)} \overline{\alpha_{2}} & 6 & -4 & 0 \\
    	-80 \, A_{0} A_{1} \zeta_{3} & 3 & -1 & 1 \\
    	48 \, A_{0} A_{1} \overline{\alpha_{2}} & 4 & -4 & 0 \\
    	2 \, A_{1} C_{1} & -1 & 0 & 0 \\
    	80 \, A_{0} A_{1} \overline{\zeta_{1}} & 0 & 2 & 1
    	\end{dmatrix}
    \end{align}
    \normalsize
    where
    \small
    \begin{align*}
    	\pi_1&=32 \, A_{0}^{2} {\left| A_{1} \right|}^{6} + 16 \, A_{0}^{2} \beta {\left| A_{1} \right|}^{2} + 14 \, A_{0} C_{1} \alpha_{1}^{2} + 6 \, A_{1} {\left| A_{1} \right|}^{2} \overline{B_{1}} - 2 \, A_{0} {\left| A_{1} \right|}^{2} \overline{B_{2}} + 2 \, A_{2} {\left| A_{1} \right|}^{2} \overline{C_{1}} + 24 \, A_{0} A_{1} \alpha_{1} \overline{\alpha_{3}}\\
    	& - 18 \, A_{1} C_{1} \alpha_{3} - 20 \, A_{0} C_{1} \alpha_{6} - \frac{7}{6} \, A_{0} C_{1} \zeta_{0} - 48 \, A_{0} A_{1} \overline{\alpha_{11}} - 38 \, A_{0} A_{1} \overline{\zeta_{3}} + 20 \, A_{4} C_{1} - 2 \, A_{2} C_{3} + 2 \, A_{1} C_{4}\\
    	& - 2 \, {\left(A_{0} {\left| A_{1} \right|}^{2} \overline{C_{1}} + 7 \, A_{2} C_{1} - A_{0} C_{3}\right)} \alpha_{1} + {\left(32 \, A_{0} A_{1} \overline{\alpha_{1}} + A_{0} \overline{B_{1}} - A_{1} \overline{C_{1}}\right)} \alpha_{4} - 12 \, {\left(2 \, A_{0}^{2} \alpha_{1} - A_{1}^{2} - 2 \, A_{0} A_{2}\right)} \alpha_{7}\\
    	& + {\left(8 \, A_{1} \alpha_{1} + 9 \, A_{0} \alpha_{3} - 9 \, A_{3}\right)} \gamma_{1} + \frac{3}{2} \, A_{1} \gamma_{3} + 16 \, {\left(3 \, A_{0}^{2} \alpha_{1} {\left| A_{1} \right|}^{2} - 2 \, A_{1}^{2} {\left| A_{1} \right|}^{2} - 4 \, A_{0} A_{2} {\left| A_{1} \right|}^{2}\right)} \overline{\alpha_{1}}\\
    	& + 25 \, {\left(16 \, A_{0}^{2} \alpha_{2} - A_{0} B_{1}\right)} \overline{\alpha_{2}} - 8 \, {\left(6 \, A_{0} A_{1} {\left| A_{1} \right|}^{2} + A_{0}^{2} \alpha_{4}\right)} \overline{\alpha_{4}}\\
    	\pi_{2}&=-\frac{8}{15} \, A_{0} C_{1} {\left| A_{1} \right|}^{2} \overline{\alpha_{1}} + 64 \, A_{0}^{2} {\left| A_{1} \right|}^{2} \overline{\alpha_{5}} - \frac{16}{5} \, A_{0} B_{2} {\left| A_{1} \right|}^{2} - 16 \, A_{0} E_{1} \alpha_{1} + 48 \, A_{0} A_{1} \overline{\alpha_{8}} + 2 \, A_{1} B_{3} + 16 \, A_{2} E_{1}\\
    	& + 2 \, {\left(40 \, A_{0} A_{1} \overline{\alpha_{1}} - A_{0} \overline{B_{1}} - 2 \, A_{1} \overline{C_{1}}\right)} \alpha_{2} - \frac{1}{36} \, C_{1} \gamma_{2} + 2 \, A_{1} \gamma_{4} + \frac{1}{8} \, B_{1} \overline{B_{1}} + {\left(16 \, A_{0}^{2} \alpha_{2} - A_{0} B_{1}\right)} \overline{\alpha_{4}}\\
    	\pi_{3}&=320 \, A_{0} A_{1} \alpha_{2} {\left| A_{1} \right|}^{2} - 12 \, A_{0}^{2} {\left| A_{1} \right|}^{2} \overline{\zeta_{0}} - 8 \, A_{1} B_{1} {\left| A_{1} \right|}^{2} + 2 \, A_{0} C_{1} \alpha_{1} \overline{\alpha_{1}} + \frac{3}{4} \, A_{0} {\left| A_{1} \right|}^{2} \overline{\gamma_{2}} - 6 \, A_{0} B_{2} \alpha_{1} - 2 \, A_{0} C_{1} \beta\\
    	& - 4 \, A_{1} C_{1} \overline{\alpha_{4}} + 40 \, A_{0} A_{1} \overline{\alpha_{9}} - 6 \, A_{0} A_{1} \overline{\zeta_{1}} + 6 \, A_{2} B_{2} + \frac{1}{12} \, {\left(128 \, A_{0} {\left| A_{1} \right|}^{4} + 3 \, \overline{B_{2}}\right)} C_{1} + \frac{9}{2} \, {\left(16 \, A_{0}^{2} \alpha_{2} - A_{0} B_{1}\right)} \alpha_{4}\\
    	& + \frac{1}{24} \, {\left(8 \, A_{0} \overline{\alpha_{4}} - \overline{B_{1}}\right)} \gamma_{1} + 60 \, {\left(2 \, A_{0}^{2} \alpha_{1} - A_{1}^{2} - 2 \, A_{0} A_{2}\right)} \overline{\alpha_{5}} - \frac{1}{2} \, A_{1} \overline{\gamma_{3}}\\
    	\pi_{4}&=-8 \, A_{1} C_{1} {\left| A_{1} \right|}^{2} + \frac{4}{3} \, A_{0} \gamma_{1} {\left| A_{1} \right|}^{2} - 6 \, A_{0} B_{1} \alpha_{1} - 4 \, A_{0} C_{1} \alpha_{4} + 8 \, A_{0} A_{1} \overline{\zeta_{0}} + 6 \, A_{2} B_{1} + 48 \, {\left(2 \, A_{0}^{2} \alpha_{1} - A_{1}^{2} - 2 \, A_{0} A_{2}\right)} \alpha_{2}\\
    	& - \frac{1}{2} \, A_{1} \overline{\gamma_{2}}\\
    	\pi_{5}&=-12 \, A_{2} C_{1} {\left| A_{1} \right|}^{2} - \frac{2}{3} \, A_{0} \gamma_{2} {\left| A_{1} \right|}^{2} - 10 \, A_{1} C_{1} \alpha_{4} - 8 \, A_{0} C_{1} \overline{\alpha_{7}} - 24 \, A_{0} A_{1} \overline{\zeta_{2}} + 12 \, A_{3} B_{1} + \frac{2}{3} \, {\left(38 \, A_{0} C_{1} {\left| A_{1} \right|}^{2} - 15 \, A_{1} B_{1}\right)} \alpha_{1}\\
    	& + 32 \, {\left(11 \, A_{0} A_{1} \alpha_{1} - 6 \, A_{1} A_{2} - 6 \, A_{0} A_{3}\right)} \alpha_{2} + 12 \, {\left(16 \, A_{0}^{2} \alpha_{2} - A_{0} B_{1}\right)} \alpha_{3} + \frac{4}{3} \, {\left(5 \, A_{1} {\left| A_{1} \right|}^{2} + 2 \, A_{0} \alpha_{4}\right)} \gamma_{1} + A_{1} \gamma_{5} + 2 \, C_{1} \overline{E_{1}}\\
    	& + 2 \, {\left(A_{0} \alpha_{1} - A_{2}\right)} \overline{\gamma_{2}} - 32 \, {\left(A_{0}^{2} \alpha_{1} - A_{0} A_{2}\right)} \overline{\zeta_{0}}\\
    	\pi_{6}&=-80 \, A_{1}^{2} {\left| A_{1} \right|}^{4} - 160 \, A_{0} A_{2} {\left| A_{1} \right|}^{4} - 16 \, A_{0} A_{1} \alpha_{4} {\left| A_{1} \right|}^{2} - 80 \, A_{0} A_{1} \alpha_{11} - 40 \, A_{0} C_{1} \alpha_{5} - 48 \, A_{0} A_{1} \zeta_{3} - 28 \, A_{1} C_{1} \overline{\alpha_{2}}\\
    	& + \frac{64}{3} \, A_{0} \gamma_{1} \overline{\alpha_{2}} + 2 \, {\left(80 \, A_{0}^{2} {\left| A_{1} \right|}^{4} + 3 \, A_{1} \overline{B_{1}}\right)} \alpha_{1} + 6 \, {\left(8 \, A_{0} A_{1} \overline{\alpha_{1}} + A_{0} \overline{B_{1}}\right)} \alpha_{3} - 6 \, A_{3} \overline{B_{1}} + 6 \, A_{1} \overline{B_{3}}\\
    	& - 48 \, {\left(A_{0} A_{1} \alpha_{1} + A_{0}^{2} \alpha_{3} - A_{1} A_{2} - A_{0} A_{3}\right)} \overline{\alpha_{4}} + 6 \, A_{1} \overline{\gamma_{4}} \\
    	\pi_{7}&=-32 \, A_{0} A_{1} {\left| A_{1} \right|}^{4} + 16 \, A_{0} A_{1} \alpha_{1} \overline{\alpha_{1}} - 16 \, A_{0} A_{1} \beta + 2 \, A_{0} \alpha_{1} \overline{B_{1}} - \frac{80}{3} \, A_{0} C_{1} \overline{\alpha_{2}} - 2 \, A_{2} \overline{B_{1}} + 2 \, A_{1} \overline{B_{2}}\\
    	& - 8 \, {\left(2 \, A_{0}^{2} \alpha_{1} - A_{1}^{2} - 2 \, A_{0} A_{2}\right)} \overline{\alpha_{4}}\\
    	\pi_{8}&=80 \, A_{0}^{2} \alpha_{1}^{2} {\left| A_{1} \right|}^{2} - 96 \, A_{0} A_{1} \alpha_{3} {\left| A_{1} \right|}^{2} - 96 \, A_{0}^{2} \alpha_{6} {\left| A_{1} \right|}^{2} + 20 \, A_{0}^{2} \zeta_{0} {\left| A_{1} \right|}^{2} + 48 \, A_{2}^{2} {\left| A_{1} \right|}^{2} + 96 \, A_{1} A_{3} {\left| A_{1} \right|}^{2} + 96 \, A_{0} A_{4} {\left| A_{1} \right|}^{2}\\
    	& - 144 \, A_{0}^{2} \overline{\alpha_{2}} \overline{\alpha_{4}} - 72 \, A_{0} A_{1} \alpha_{10} - 54 \, A_{0} A_{1} \zeta_{2} - 2 \, {\left(32 \, A_{1}^{2} {\left| A_{1} \right|}^{2} + 64 \, A_{0} A_{2} {\left| A_{1} \right|}^{2} + 3 \, A_{0} \overline{E_{1}}\right)} \alpha_{1}\\
    	& - 8 \, {\left(A_{0} A_{1} \alpha_{1} + 3 \, A_{0}^{2} \alpha_{3} - 3 \, A_{1} A_{2} - 3 \, A_{0} A_{3}\right)} \alpha_{4} + 6 \, A_{2} \overline{E_{1}} + 12 \, A_{1} \overline{E_{2}} + 6 \, {\left(3 \, A_{0} \overline{B_{1}} + A_{1} \overline{C_{1}}\right)} \overline{\alpha_{2}}\\
    	& + 12 \, {\left(2 \, A_{0}^{2} \alpha_{1} - A_{1}^{2} - 2 \, A_{0} A_{2}\right)} \overline{\alpha_{7}}
    \end{align*}
    \normalsize
    As we will be only interested in $\pi_4$ and $\pi_6$, which correspond respectively to the coefficients 
    \begin{align*}
    	\z\,dz^4,\quad\text{and}\;\, z^{\theta_0}\z^{2-\theta_0}=z^{3}\z^{-1}\,dz^4,
    \end{align*}
    in \eqref{part}, we deduce that we only need to compute $\vec{B}_3$ and $\vec{\gamma_4}$. First, recall that by \eqref{3consts2}
\begin{align}\label{3relend1}
    \left\{
    \begin{alignedat}{1}
    	\vec{C}_2&=\Re\left(\vec{D}_2\right)\in\R^n\\
    	\vec{C}_3&=\vec{D}_3-\frac{1}{4}\s{\vec{C}_1}{\vec{C}_1}\vec{A}_0+\left(-4(\s{\vec{A}_2}{\vec{C}_1}+\s{\vec{A}_1}{\vec{C}_2})+2\s{\vec{A}_1}{\vec{\gamma}_1}\right)\bar{\vec{A}_0}-2\s{\vec{A}_1}{\bar{\vec{C}_1}}\bar{\vec{A}_1}\\
    	\vec{B}_1&=-2\s{\bar{\vec{A}_1}}{\vec{C}_1}\vec{A}_0\\
    	\vec{B}_2&=-\frac{5}{12}|\vec{C}_1|^2\bar{\vec{A}_0}-2\s{\bar{\vec{A}_2}}{\vec{C}_1}\vec{A}_0.\\
    	\vec{E}_1&=-\frac{1}{6}\s{\vec{C}_1}{\vec{C}_1}\bar{\vec{A}_0}\\
    	\vec{\gamma}_1&=-\vec{\gamma}_0-4\,\Re\left(\s{\vec{A}_1}{\vec{C}_1}\bar{\vec{A}_0}\right)\in\mathbb{R}^n\\
    	\vec{\gamma}_2&=-4\s{\vec{A}_1}{\vec{\gamma}_1}\bar{\vec{A}_0}-4\bar{\s{\vec{A}_1}{\vec{C}_1}}\vec{A}_1.
    \end{alignedat}\right.
\end{align}
Therefore, as $\vec{B}_1\in \mathrm{Span}(\vec{A}_0)$, and as by \eqref{3conf}, we have
\begin{align*}
	\s{\vec{A}_0}{\vec{A}_0}=\s{\vec{A}_0}{\vec{A}_1}=\s{\vec{A}_0}{\vec{C}_1}=0,\quad \s{\vec{A}_1}{\vec{A}_1}+2\s{\vec{A}_0}{\vec{A}_2}=0,
\end{align*}
we deduce that
\begin{align}\label{3relend2}
\left\{
\begin{alignedat}{1}
	&\s{\vec{A}_2}{\vec{B}_1}=-2\s{\bar{\vec{A}_1}}{\vec{C}_1}\s{\vec{A}_0}{\vec{A}_2}=\s{\bar{\vec{A}_1}}{\vec{C}_1}\s{\vec{A}_1}{\vec{A}_1}\\
	&\s{\vec{A}_1}{\bar{\vec{\gamma}_2}}=-4|\vec{A}_1|^2\s{\vec{A}_1}{\vec{C}_1}
	\end{alignedat}\right.
\end{align}
so by \eqref{3relend1} and \eqref{3relend2}, we finally obtain
\begin{align*}
	\pi_{4}&=-8 \, A_{1} C_{1} {\left| A_{1} \right|}^{2} + \frac{4}{3} \ccancel{\, A_{0} \gamma_{1} {\left| A_{1} \right|}^{2}} - \ccancel{6 \, A_{0} B_{1} \alpha_{1}} - \ccancel{4 \, A_{0} C_{1} \alpha_{4}} + \ccancel{8 \, A_{0} A_{1} \overline{\zeta_{0}}} + 6 \, A_{2} B_{1} + 48 \, {\left(\ccancel{2 \, A_{0}^{2} \alpha_{1}} - \colorcancel{A_{1}^{2}}{blue} - \colorcancel{2 \, A_{0} A_{2}}{blue}\right)} \alpha_{2}\\
	& - \frac{1}{2} \, A_{1} \overline{\gamma_{2}}\\
	&=-8|\vec{A}_1|^2\s{\vec{A}_1}{\vec{C}_1}+6\s{\bar{\vec{A}_1}}{\vec{C}_1}\s{\vec{A}_1}{\vec{A}_1}-\frac{1}{2}\left(-4|\vec{A}_1|^2\s{\vec{A}_1}{\vec{C}_1}\right)\\
	&=-6\left(|\vec{A}_1|^2\s{\vec{A}_1}{\vec{C}_1}-\s{\bar{\vec{A}_1}}{\vec{C}_1}\s{\vec{A}_1}{\vec{A}_1}\right)\\
	&=0
\end{align*}
    and we obtain the first line of the system in $\s{\vec{A}_1}{\vec{C}_1},\s{\vec{A}_1}{\vec{A}_1}$
    \begin{align}\label{3line1}
    	|\vec{A}_1|^2\s{\vec{A}_1}{\vec{C}_1}=\s{\bar{\vec{A}_1}}{\vec{C}_1}\s{\vec{A}_1}{\vec{A}_1}.
    \end{align}
    Now, notice that thanks of \eqref{3conf}, we have
    \begin{align}\label{3relend3}
    	\s{\vec{A}_0}{\vec{A}_3}+\s{\vec{A}_1}{\vec{A}_2}=0,
    \end{align}
    Then, we have thanks of \eqref{3hend}
    \begin{align*}
    	\frac{1}{2}\vec{B}_3=\begin{dmatrix}
    		A_{0} \alpha_{2} \overline{C_{1}} - \frac{1}{3} \, C_{1} C_{2} + A_{1} E_{1} + \frac{1}{36} \, C_{1} \gamma_{1} - \frac{7}{48} \, B_{1} \overline{C_{1}} & \overline{A_{0}} & -1 & 3 & 0 & \textbf{(9)}\\
    	-\frac{1}{6} \, C_{1} \overline{C_{1}} & \overline{A_{1}} & -1 & 3 & 0 & \textbf{(10)}\\
    	-\frac{1}{36} \, \overline{A_{1}} \overline{C_{1}} & C_{1} & -1 & 3 & 0 & \textbf{(11)}\\
    	\frac{4}{3} \, \alpha_{2} \overline{A_{0}} \overline{C_{1}} + C_{1} \overline{A_{0}} \overline{\alpha_{3}} - \frac{1}{3} \, B_{2} \overline{A_{1}} - \frac{2}{3} \, B_{1} \overline{A_{2}} - C_{1} \overline{A_{3}} + \frac{1}{3} \, {\left(2 \, B_{1} \overline{A_{0}} + 3 \, C_{1} \overline{A_{1}}\right)} \overline{\alpha_{1}} & A_{0} & -1 & 3 & 0 & \textbf{(12)}\\
    	-\frac{1}{24} \, C_{1} \overline{A_{1}} & \overline{C_{1}} & -1 & 3 & 0 & \textbf{(13)}
    	\end{dmatrix}
    \end{align*}
    Therefore, we have
    \begin{align*}
    	\frac{1}{2}\s{\bar{\vec{A}_1}}{\vec{B}_3}&=-\frac{1}{6}|\vec{C}_1|^2\s{\vec{A}_1}{\vec{A}_1}-\frac{1}{36}\s{\vec{A}_1}{\vec{C}_1}\s{\vec{A}_1}{\bar{\vec{C}_1}}-\frac{1}{24}\s{\bar{\vec{A}_1}}{\vec{C}_1}\s{\vec{A}_1}{\vec{C}_1}\\
    	&=-\frac{1}{6}|\vec{C}_1|^2\s{\vec{A}_1}{\vec{A}_1}-\frac{5}{72}\s{\vec{A}_1}{\vec{C}_1}\s{\vec{A}_1}{\bar{\vec{C}_1}}
    \end{align*}
    so
    \begin{align}\label{rel}
    	&\s{\bar{\vec{A}_1}}{\vec{B}_3}=-\frac{1}{3}|\vec{C}_1|^2\s{\vec{A}_1}{\vec{A}_1}-\frac{5}{36}\s{\vec{A}_1}{\vec{C}_1}\s{\vec{A}_1}{\bar{\vec{C}_1}}\nonumber\\
    	&6\s{\bar{\vec{A}_1}}{\vec{B}_3}=-2|\vec{C}_1|^2\s{\vec{A}_1}{\vec{A}_1}-\frac{5}{6}\s{\vec{A}_1}{\vec{C}_1}\s{\vec{A}_1}{\bar{\vec{C}_1}}
    \end{align}
    We also have
    \begin{align*}
    	\frac{1}{2}\vec{\gamma_4}=\begin{dmatrix}
    	-\frac{1}{3} \, C_{1} \gamma_{1} & \overline{A_{0}} & -1 & 3 & 1 & \textbf{(14)}
    	\end{dmatrix}
    \end{align*}
    so
    \begin{align}\label{3defgamma4}
    	\vec{\gamma}_4=-\frac{2}{3}\s{\vec{C}_1}{\vec{\gamma}_1}\bar{\vec{A}_0}=\frac{2}{3}\s{\vec{\gamma}_0}{\vec{C}_1}\bar{\vec{A}_0}
    \end{align}
    and
    \begin{align}
    	\s{\vec{A}_1}{\vec{\gamma}_4}=\s{\bar{\vec{A}_1}}{\vec{\gamma}_4}=0.
    \end{align}
    Also, recall that by \eqref{3deflambda1}, we have
    \begin{align}\label{3deflambda2}
    \left\{
    \begin{alignedat}{1}
    \alpha_0&=2\s{\bar{\vec{A}_0}}{\vec{A}_1}\\
    \alpha_1&=2\s{\bar{\vec{A}_0}}{\vec{A}_2}\\
    \alpha_2&=\frac{1}{24}\s{\bar{\vec{A}_1}}{\vec{C}_1}\\
    \alpha_3&=\frac{1}{12}\s{
    	\vec{A}_1}{\vec{C}_1}+2\s{\bar{\vec{A}_0}}{\vec{A}_3}\\
    \alpha_4&=2\s{\bar{\vec{A}_1}}{\vec{A}_2},
    \end{alignedat}\right.
    \end{align}
    In particular, we deduce as $|\vec{A}_0|^2=\dfrac{1}{2}$ that
    \begin{align}\label{rel2}
         &\s{\alpha_3\vec{A}_0-\vec{A}_3}{\bar{\vec{A}_0}}=\frac{\alpha_3}{2}-\s{\bar{\vec{A}_0}}{\vec{A}_3}=\frac{1}{24}\s{\vec{A}_1}{\vec{C}_1}\nonumber\\
         &\s{\alpha_3\vec{A}_0-\vec{A}_3}{\bar{\vec{B}_1}}=-2\s{\vec{A}_1}{\bar{\vec{C}_1}}\s{\alpha_3\vec{A}_0-\vec{A}_3}{\bar{\vec{A}_0}}=-\frac{1}{12}\s{\vec{A}_1}{\vec{C}_1}\s{\vec{A}_1}{\bar{\vec{C}_1}}\nonumber\\
         &6\s{\alpha_3\vec{A}_0-\vec{A}_3}{\bar{\vec{B}_1}}=-\frac{1}{2}\s{\vec{A}_1}{\vec{C}_1}\s{\vec{A}_1}{\bar{\vec{C}_1}}
    \end{align}
    Now, we have by \eqref{3relend1}, \eqref{3relend2}, \eqref{3relend3}, \eqref{rel} and \eqref{rel2}
    \small
    \begin{align}\label{pi6}
    	\pi_6&=-\colorcancel{80 \, A_{1}^{2} {\left| A_{1} \right|}^{4}}{blue} - \colorcancel{160 \, A_{0} A_{2} {\left| A_{1} \right|}^{4}}{blue} - \ccancel{16 \, A_{0} A_{1} \alpha_{4} {\left| A_{1} \right|}^{2}} - \ccancel{80 \, A_{0} A_{1} \alpha_{11}} - \ccancel{40 \, A_{0} C_{1} \alpha_{5}} - \ccancel{48 \, A_{0} A_{1}} \zeta_{3} - 28 \, A_{1} C_{1} \overline{\alpha_{2}}\nonumber\\
    	& + \frac{64}{3} \ccancel{\, A_{0} \gamma_{1} \overline{\alpha_{2}}} + 2 \, {\left(\ccancel{80 \, A_{0}^{2} {\left| A_{1} \right|}^{4}} + \ccancel{3 \, A_{1} \overline{B_{1}}}\right)} \alpha_{1} + 6 \, {\left(\ccancel{8 \, A_{0} A_{1} \overline{\alpha_{1}}} + A_{0} \overline{B_{1}}\right)} \alpha_{3} - 6 \, A_{3} \overline{B_{1}} + 6 \, A_{1} \overline{B_{3}}\nonumber\\
    	& - 48 \, {\left(\ccancel{A_{0} A_{1} \alpha_{1}} + \ccancel{A_{0}^{2} \alpha_{3}} - \colorcancel{A_{1} A_{2}}{blue} - \colorcancel{A_{0} A_{3}}{blue}\right)} \overline{\alpha_{4}} + \ccancel{6 \, A_{1} \overline{\gamma_{4}}}\nonumber\\
    	&=-28\s{\vec{A}_1}{\vec{C}_1}\bar{\alpha_2}-\frac{1}{2}\s{\vec{A}_1}{\vec{C}_1}\s{\vec{A}_1}{\bar{\vec{C}_1}}-2|\vec{C}_1|^2\s{\vec{A}_1}{\vec{A}_1}-\frac{5}{6}\s{\vec{A}_1}{\vec{C}_1}\s{\vec{A}_1}{\bar{\vec{C}_1}}\nonumber\\
    	&=-\frac{28}{24}\s{\vec{A}_1}{\vec{C}_1}\s{\vec{A}_1}{\bar{\vec{C}_1}}-\frac{1}{2}\s{\vec{A}_1}{\vec{C}_1}\s{\vec{A}_1}{\bar{\vec{C}_1}}-2|\vec{C}_1|^2\s{\vec{A}_1}{\vec{A}_1}-\frac{5}{6}\s{\vec{A}_1}{\vec{C}_1}\s{\vec{A}_1}{\bar{\vec{C}_1}}\\
    	&=-\frac{5}{2}\s{\vec{A}_1}{\vec{C}_1}\s{\bar{\vec{A}_1}}{\vec{C}_1}-2|\vec{C}_1|^2\s{\vec{A}_1}{\vec{A}_1}
    \end{align}
    \normalsize
    Now, thanks of \eqref{3devextra1}, the coefficient in $z^3\z^{-1}dz^4$ in the Taylor expansion of 
    \begin{align}
    	\frac{5}{4}|\H|^2\h_0\totimes\h_0+\s{\H}{\h_0}^2
    \end{align}
    is
    \begin{align}
    	2\s{\vec{A}_1}{\vec{C}_1}\s{\vec{A}_1}{\bar{\vec{C}_1}}+\frac{5}{2}|\vec{C}_1|^2\s{\vec{A}_1}{\vec{A}_1},
    \end{align}
    so the coefficient in $z^{3}\z^{-1}dz^4$ in the Taylor expansion of the \emph{meromorphic} quartic form is
    \begin{align*}
    	\Omega=\pi_6+2\s{\vec{A}_1}{\vec{C}_1}\s{\vec{A}_1}{\bar{\vec{C}_1}}+\frac{5}{2}|\vec{C}_1|^2\s{\vec{A}_1}{\vec{A}_1}=\frac{1}{2}\left(|\vec{C}_1|^2\s{\vec{A}_1}{\vec{A}_1}-\s{\vec{A}_1}{\vec{C}_1}\s{\vec{A}_1}{\bar{\vec{C}_1}}\right)=0,
    \end{align*}
    and we finally recover the system
    \begin{align}\label{3system}
    \left\{\begin{alignedat}{1}
    |\vec{A}_1|^2\s{\vec{A}_1}{\vec{C}_1}&=\s{\bar{\vec{A}_1}}{\vec{C}_1}\s{\vec{A}_1}{\vec{A}_1}\\
    |\vec{C}_1|^2\s{\vec{A}_1}{\vec{A}_1}&=\s{\vec{A}_1}{\bar{\vec{C}_1}}\s{\vec{A}_1}{\vec{C}_1},
\end{alignedat}\right.
\end{align}
Now, as by \eqref{3true}, we have
\begin{align}
	\s{\vec{A}_0}{\vec{\gamma}_0}+\s{\vec{A}_1}{\vec{C}_1}=0,
\end{align}
and for a \emph{true} Willmore disk, we have $\vec{\gamma}_0=0$, so
\begin{align*}
	\s{\vec{A}_1}{\vec{C}_1}=0,
\end{align*}
the meromorphic quartic differential
\begin{align*}
	\mathscr{Q}_{\phi}=2\s{\vec{A}_1}{\vec{C}_1}\frac{dz^4}{z}+O(1)
\end{align*}
is holomorphic, and thanks of \eqref{3system}, we also deduce that the octic form $\mathscr{O}_{\phi}$ is holomorphic thanks of the analysis of chapter \ref{chapteroctic}. This concludes the proof of the case $\theta_0=3$.

    		\normalsize
    		\section{The case where $\theta_0=2$}
    		
    		In this case, as the proof is very short (only two pages), we do not transcript our computer computations here. Notice that we cannot show the holomorphy in general if we assume non-zero first residue $\vec{\gamma}_0$ as such branch points. As for true branch points of multiplicity $\theta_0=2$, the Willmore immersion is smooth, the form $\mathscr{Q}_{\phi}$ and $\mathscr{O}_{\phi}$ are trivially holomorphic.

	    		\normalsize
	    		
	    		\normalsize
	    		\nocite{}
	    		\bibliographystyle{plain}
	    		\bibliography{sagebiblio}
	    		
    		\end{document}